\documentclass[11pt]{book}
\usepackage{hyperref}
\usepackage{amsfonts,amssymb,amsmath,amsthm,cite}
\usepackage{graphicx}
\usepackage[toc,page]{appendix}
\usepackage{nicefrac}
\usepackage[applemac]{inputenc}
\usepackage{amssymb, euscript}
\usepackage[matrix,arrow,curve]{xy}
\usepackage{graphicx}
\usepackage{tabularx}
\usepackage{float}
\usepackage{tikz}
\usepackage{slashed}
\usepackage{mathrsfs}
\usepackage{multirow}
\usepackage{rotating}

\usetikzlibrary{matrix}
\usetikzlibrary{cd}

\usepackage{siunitx}

\usepackage{lmodern}
\usepackage[T1]{fontenc}
\usepackage[babel=true]{microtype}

\usepackage{amsfonts,cite}
\usepackage{graphicx}

\usepackage[applemac]{inputenc}

\usepackage[sc]{mathpazo}
\usepackage{environ}

\DeclareFontFamily{T1}{pzc}{}
\DeclareFontShape{T1}{pzc}{m}{it}{1.8 <-> pzcmi8t}{}
\DeclareMathAlphabet{\mathpzc}{T1}{pzc}{m}{it}

\theoremstyle{plain}
\newtheorem{prop}{Proposition}[section]
\newtheorem{prdf}[prop]{Proposition and Definition}
\newtheorem{lem}[prop]{Lemma}
\newtheorem{cor}[prop]{Corollary}
\newtheorem{thm}[prop]{Theorem}
\newtheorem{theorem}[prop]{Theorem}
\newtheorem{lemma}[prop]{Lemma}
\newtheorem{proposition}[prop]{Proposition}
\newtheorem{corollary}[prop]{Corollary}

\theoremstyle{definition}
\newtheorem{defn}[prop]{Definition}
\newtheorem{empt}[prop]{}
\newtheorem{exm}[prop]{Example}
\newtheorem{rem}[prop]{Remark}

\newtheorem{notn}{Notation}        

\theoremstyle{definition}
\newtheorem{notation}[prop]{Notation}
\newtheorem{definition}[prop]{Definition}
\newtheorem{example}[prop]{Example}
\newtheorem{exercise}[prop]{Exercise}

\newtheorem{axiom}[prop]{Axiom}
\newtheorem{problem}[prop]{Problem}
\newtheorem{remark}[prop]{Remark}

\numberwithin{equation}{section}

\DeclareMathOperator{\Dom}{Dom}              
\newcommand{\Dslash}{{D\mkern-11.5mu/\,}}    

\newcommand{\vertiii}[1]{{\left\vert\kern-0.25ex\left\vert\kern-0.25ex\left\vert #1
		\right\vert\kern-0.25ex\right\vert\kern-0.25ex\right\vert}}
\newcommand{\Ga}{\Gamma}  
\newcommand{\coker}{\mathrm{coker}}                   
\newcommand{\Coo}{C^\infty}                  

\usepackage[sc]{mathpazo}
\newbox\ncintdbox \newbox\ncinttbox 
\setbox0=\hbox{$-$} \setbox2=\hbox{$\displaystyle\int$}
\setbox\ncintdbox=\hbox{\rlap{\hbox
		to \wd2{\hskip-.125em \box2\relax\hfil}}\box0\kern.1em}
\setbox0=\hbox{$\vcenter{\hrule width 4pt}$}
\setbox2=\hbox{$\textstyle\int$} \setbox\ncinttbox=\hbox{\rlap{\hbox
		to \wd2{\hskip-.175em \box2\relax\hfil}}\box0\kern.1em}

\newcommand{\oxyox}{\otimes\cdots\otimes}  
\newcommand{\rep}{\mathfrak{rep}}
\newcommand{\lift}{\mathfrak{lift}}
\newcommand{\desc}{\mathfrak{desc}}
\newcommand{\id}{\mathrm{id}}                
\newcommand{\Id}{\mathrm{Id}}                
\newcommand{\codim}{\mathrm{codim}}          
\renewcommand{\d}{\mathrm{d}}       
\newcommand*{\braket}[2]{\langle#1 {,~} #2\rangle}

\newcommand{\A}{\mathcal{A}}                 

\renewcommand{\a}{\alpha}                    
\newcommand{\B}{\mathcal{B}}                 
\newcommand{\E}{\mathcal{E}}                 
\renewcommand{\b}{\beta}                     

\newcommand{\C}{\mathbb{C}}                  
\newcommand{\cc}{\mathbf{c}}                 
\DeclareMathOperator{\Cl}{C\ell}             
\newcommand{\F}{\mathcal{F}}                 
\newcommand{\G}{\mathcal{G}}                 
\newcommand{\D}{\mathcal{D}}                 
\renewcommand{\H}{\mathcal{H}}               
\newcommand{\half}{\tfrac{1}{2}}             
\newcommand{\hookto}{\hookrightarrow}        
\newcommand{\K}{\mathcal{K}}                 
\renewcommand{\L}{\mathcal{L}}               
\newcommand{\La}{\Lambda}                    
\newcommand{\la}{\lambda}                    
\newcommand{\M}{\mathcal{M}}                 

	\newcommand{\N}{\mathbb{N}}                  
	\newcommand{\nb}{\nabla}                     
	\newcommand{\om}{\omega}                     
	\newcommand{\opp}{{\mathrm{op}}}             
	\newcommand{\ox}{\otimes}                    
	\newcommand{\eps}{\varepsilon}                    
	\newcommand{\Q}{\mathbb{Q}}                  
	\newcommand{\R}{\mathbb{R}}                  
	\newcommand{\sepword}[1]{\quad\mbox{#1}\quad} 
	\newcommand{\set}[1]{\{\,#1\,\}}             

	\renewcommand{\SS}{\mathcal{S}}              
	\DeclareMathOperator{\supp}{\mathfrak{supp}}            
	\newcommand{\T}{\mathbb{T}}                  
	\renewcommand{\th}{\theta}                   
	\DeclareMathOperator{\tr}{tr}                
	\newcommand{\del}{\partial}                  
	\DeclareMathOperator{\tsum}{{\textstyle\sum}} 
	\newcommand{\Z}{\mathbb{Z}}                  
	\newcommand{\7}{\dagger}                     
	\renewcommand{\.}{\cdot}                     
	\renewcommand{\:}{\colon}                    

	\newcommand{\sA}{\mathcal{A}} 
	\newcommand{\sF}{\mathcal{F}}       
	\newcommand{\sG}{\mathcal{G}}       
	\newcommand{\sH}{\mathcal{H}}       
	\newcommand{\sS}{\mathcal{S}}       
	\newcommand{\sU}{\mathcal{U}}       
	\newcommand{\sV}{\mathcal{V}}       
	\newcommand{\sX}{\mathcal{X}}       
	\newcommand{\sY}{\mathcal{Y}}       
	\newcommand{\sZ}{\mathcal{Z}}       
	
	\newcommand{\Om}{\Omega}       

	\newcommand{\bN}{\mathbb{N}}       
	\newcommand{\bZ}{\mathbb{Z}}       
	
	\newcommand{\bydef}{\stackrel{\mathrm{def}}{=}}          

	\newcommand{\al}{\alpha}          
	\newcommand{\bt}{\beta}           
	\newcommand{\dl}{\delta}          
	\newcommand{\ga}{\gamma}          
	\newcommand{\ka}{\kappa}          
	\newcommand{\Th}{\Theta}          
	\renewcommand{\th}{\theta}        
	
	\DeclareMathOperator{\Ind}{Ind}
	\newcommand{\blank}{-}

	\newcommand{\isom}{\cong}          
	\renewcommand{\:}{\colon}           

	\newcommand{\pairing}[2]{(#1\stroke #2)} 
	\def\<#1|#2>{\langle#1\stroke#2\rangle} 
	\def\?#1|#2?{\{#1\stroke#2\}}        
	
	\newcommand{\word}[1]{\quad\text{#1}\enspace} 

	\def\<#1,#2>{\langle#1,#2\rangle}            
	\def\ee_#1{e_{{\scriptscriptstyle#1}}}       
	\def\wick:#1:{\mathopen:#1\mathclose:}       

	\newbox\ncintdbox \newbox\ncinttbox 
	\setbox0=\hbox{$-$}
	\setbox2=\hbox{$\displaystyle\int$}
	\setbox\ncintdbox=\hbox{\rlap{\hbox
			to \wd2{\box2\relax\hfil}}\box0\kern.1em}
	\setbox0=\hbox{$\vcenter{\hrule width 4pt}$}
	\setbox2=\hbox{$\textstyle\int$}
	\setbox\ncinttbox=\hbox{\rlap{\hbox
			to \wd2{\hskip-.05em\box2\relax\hfil}}\box0\kern.1em}
	
	\newcommand{\stroke}{\mathbin|}   
	
	\newcommand{\End}{\mathrm{End}}       
	\newcommand{\Ext}{\mathrm{Ext}}       
	\newcommand{\Hom}{\mathrm{Hom}}       
	\newcommand{\Spin}{\mathrm{Spin}}       
	\newcommand{\Aut}{\mathrm{Aut}}       
	\newcommand{\Homeo}{\mathrm{Homeo}}       
	\newcommand{\Diff}{\mathrm{Diff}}       
	\newcommand{\im}{\mathrm{im}}       

	\newcommand{\SU}{SU}       
	
	\newcommand{\bs}{\begin{split}}
		\newcommand{\es}{\end{split}}
	\newcommand{\be}{\begin{equation}}
		\renewcommand{\ee}{\end{equation}}
	\newcommand{\bea}{\begin{eqnarray}}
		\newcommand{\eea}{\end{eqnarray}}
	\newcommand{\bean}{\begin{eqnarray*}}
		\newcommand{\eean}{\end{eqnarray*}}
	\newcommand{\brray}{\begin{array}}
		\newcommand{\erray}{\end{array}}
	
	\newcommand{\Hsquare}{%
		\text{\fboxsep=-.2pt\fbox{\rule{0pt}{1ex}\rule{1ex}{0pt}}}%
	}

	\title{Noncommutative Geometry of Quantized Coverings}
	
	\author
	{\textbf{Petr R. Ivankov*}\\
		e-mail: * monster.ivankov@gmail.com \\
	}
	
\begin{document}

\maketitle  
\pagestyle{plain}




\chapter*{Introduction}
\paragraph*{}

There are theories of coverings of $C^*$-algebras which can be included into a following list:
\begin{itemize}
	\item coverings of commutative $C^*$-algebras  \cite{clarisson:phd,pavlov_troisky:cov},
	\item coverings of $C^*$-algebras of groupoids and foliations  \cite{ouchi:cov_fol,xiaolu:foli_cov,zhi:cov_group}),
	\item  coverings of noncommutative tori \cite{clarisson:phd,schwieger:nt_cov},
	\item {the double covering of the quantum group $SO_q(3)$} \cite{dijkhuizen:so_doublecov,podles:so_su}. 
\end{itemize}

This work is devoted to a single general theory which includes all theories of this list, i.e. we develop a  system of axioms which can be applied for every item of the list.

A commutative Gelfand-Na\u{i}mark theorem \cite{arveson:c_alg_invt} states a duality correspondence between  locally compact Hausdorff topological spaces and commutative $C^*$-algebras. So  a noncommutative $C^*$-algebra can be regarded as a noncommutative generalization of a topological space. Further development of the noncommutative geometry gives generalizations of following classical geometric and topological notions.
\break
\begin {table}[H]
\caption {Mapping between classical and noncommutative geometry} \label{main_mapping_table} 
\begin{center}
\begin{tabular}{|c|c|}
	\hline
Classical notion & Noncommutative generalization\\
	\hline
	&\\
	Topological space & $C^*$-algebra\\
	&\\
	Measure space & von Neumann algebra\\
	&\\	Riemannian manifold  & Spectral triple\\
	&\\	Topological $K$-theory & $K$-theory of $C^*$-algebras \\
	&\\	Homology and cohomology & Noncommutative homology and cohomology\\
	&\\
	\hline
\end{tabular}
\end{center}
\end {table}

In this book we continue development of the noncommutative geometry, this book contains generalizations of following notions.
\\
\\
\begin {table}[H]
\caption {Mapping between geometry of topological coverings and noncommutative ones} \label{add_mapping_table} 
\begin{center}
\begin{tabular}{|c|c|}
	\hline
	Classical notion & Noncommutative generalization\\
	\hline
	&\\
	Covering & Noncommutative covering\\
	&\\
	Theorem about covering & 	Theorem about covering \\
 of Riemannian manifold & of spectral triple\\
	&\\
	Fundamental group of a space $\pi_1\left(\mathcal X \right)$  & Fundamental group of a $C^*$-algebra $\pi_1\left(A \right)$  \\
	& \\
	Hurewicz homomorphism &  Noncommutative	Hurewicz 
	\\ $\pi_1\left(\sX, x_0 \right)\to K_1\left( \sX\right)  $ & homomorphism $\pi_1\left(A \right)\to K_1\left(A\right)  $\\
	&\\
	
Flat connections given by & Noncommutative flat connections\\
 coverings & given by noncommutative coverings\\
	&\\
Unoriented  Spin$^c$-manifolds  & Unoriented  spectral triples\\
&\\
	\hline
\end{tabular}
\end{center}
\end {table}

\paragraph*{}
There is a set of theories of noncommutative coverings (e.g. \cite{clarisson:phd,schwieger:nt_cov}). In contrary of them the presented here theory gives results which are (almost) equivalent to the classical topological theory. In particular covering spaces of commutative spaces are also commutative. This fact yields pure algebraic definition of the fundamental group   (cf. Theorem \ref{comm_uni_lim_thm}).
\paragraph*{}
The Part \ref{math_part} is devoted to mathematical aspects of the theory.
\paragraph*{}
The Chapter \ref{prel_chap} contains preliminary results. The material of this chapter can be read as needed.
\paragraph*{}
The Chapter \ref{cov_fin_chap} contains the construction of noncommutative    finite-fold coverings of $C^*$-algebras and operator spaces. Sections \ref{cov_fin_bas_sec} - \ref{induced_repr_fin_sec} are basic and needed for the further reading of this book. Other sections are written for  those who are interested in following applications of this theory:
\begin{itemize}
	\item Coverings and strong Morita equivalence.
	\item  Noncommutative path lifting.
	\item Coverings of spectral triples.
	\item Finite noncommutative coverings and flat connections.
	\item Unoriented spectral triples.
\end{itemize}
\paragraph*{}
The Chapter \ref{infinite_covering_chap} is devoted to  noncommutative   infinite coverings of $C^*$-algebras and operator spaces. The Sections \ref{infinite_ca_sec} - \ref{inf_mor_sec} are  basic. The Section \ref{str_cov_glo_sec} is interesting for  those who are interested in coverings of spectral triples.
\paragraph*{}
The Chapter \ref{top_chap} contains applications of described in Chapters  \ref{cov_fin_chap} and \ref{infinite_covering_chap} theory to commutative $C^*$-algebras. It is constructed a given by the Table \ref{add_mapping_table}   natural one to one correspondence between geometry of topological coverings and "noncommutative" ones.
\paragraph*{}
The Chapter \ref{blowing_chap} is devoted  to the construction of  Hausdorff blowing-up. This construction can be applied for obtaining of noncommutative coverings of  $C^*$-algebras with Hausdorff spectrum, and non-Hausdorff one, e.g. $C^*$-algebras of groupoids and foliations (cf. \cite{candel:foliII}).

\paragraph*{}
In the Chapter \ref{stab_chap} it is proven that some properties of noncommutative coverings are stably invariant. 

\paragraph*{} The Chapter \ref{h_chap} is devoted to the generalization of the Hurewicz homomorphism from the fundamental group to $K$-homology.

\paragraph*{} The Chapter \ref{ctr_chap} is devoted to noncommutative coverings of $C^*$-algebras with Hausdorff spectrum. It is proven that the theory of noncommutative coverings of $C^*$-algebras with Hausdorff spectrum contains all ingredients of right row of the Table \ref{add_mapping_table}. 


\paragraph*{} In the Chapter \ref{foliations_chap} we consider noncommutative finite-fold and infinite coverings of $ C^*$-algebras of groupoids and foliations. 

\paragraph*{} The Chapter \ref{nt_chap} is devoted to noncommutative coverings of noncommutative tori.

\paragraph*{} The Chapter \ref{isospectral_chap} is devoted to noncommutative coverings of isospectral deformations. We consider "noncommutative finite-fold coverings" only. The presented in Sections \ref{triple_fin_cov}, \ref{flat_sec} and \ref{unoti_defn_sec} theory of coverings of spectral triples is applied to isospectral deformations.
\paragraph*{} 
The Chapter \ref{qdr_chap} is devoted to  noncommutative coverings of quantum groups.

\paragraph*{} In the Part \ref{phys_part} we consider related to the theoretical physics math models.

\paragraph*{}
In the Chapter \ref{haus_non_haus_dul_sec} we show Hausdorff - non Hausdorff duality, i.e. a $C^*$-algebra with non Hausdorff spectrum can have a covering  $C^*$-algebra with  Hausdorff spectrum. In principle having physical models with 
Hausdorff spectrum one can  construct physical models with 
non Hausdorff spectrum.
\paragraph*{}
The Chapter \ref{wl_chap} is devoted to Wilson lines over noncommutative spaces.

\tableofcontents

\part{Mathematical Theory}\label{math_part}
\chapter{Preliminaries}\label{prel_chap}

\section{Basic facts. Notations}
\paragraph*{}

This research is based on the noncommutative generalizations of both topological spaces and coverings. Following two theorems describe these generalizations.
\begin{thm}\label{gelfand-naimark_thm}\cite{arveson:c_alg_invt} (Commutative Gelfand-Na\u{\i}mark theorem). 
	Let $A$ be a commutative $C^*$-algebra and let $\mathcal{X}$ be the spectrum of A. There is the natural $*$-isomorphism $\gamma:A \xrightarrow{\cong} C_0(\mathcal{X})$.
\end{thm}

\begin{theorem}\label{pavlov_troisky_thm}\cite{pavlov_troisky:cov}
	Suppose $\mathcal X$ and $\mathcal Y$ are compact Hausdorff connected spaces and $p :\mathcal  Y \to \mathcal X$
	is a continuous surjection. If $C(\mathcal Y )$ is a projective finitely generated Hilbert module over
	$C(\mathcal X)$ with respect to the action
	\begin{equation*}
	(f\xi)(y) = f(y)\xi(p(y)), ~ f \in  C(\mathcal Y ), ~ \xi \in  C(\mathcal X),
	\end{equation*}
	then $p$ is a finite-fold  covering.
\end{theorem}

Following table contains  special symbols.
\\
\begin{tabular}{|c|c|}
	\hline
	Symbol & Meaning\\
	\hline
	&\\
	$\hat{A}$ or $A\hat~$ & Spectrum of a  $C^*$- algebra $A$  with the hull-kernel topology \\
	& (or Jacobson topology)\\
	$A_+$  & Cone of positive elements of $C^*$- algebra, i.e. $A_+ \bydef \left\{a\in A \ | \ a \ge 0\right\}$\\
	$A^G$  & Algebra of $G$ - invariants, i.e. $A^G \bydef \left\{a\in A \ | \ ga=a, \forall g\in G\right\}$\\
	$\mathrm{Aut}(A)$ & Group of * - automorphisms of $C^*$- algebra $A$\\
	$A''$  & Enveloping von Neumann algebra  of $A$\\
	
	$B(\H)$ & Algebra of bounded operators on a Hilbert space $\H$\\
	$\mathbb{C}$ (resp. $\mathbb{R}$)  & Field of complex (resp. real) numbers \\
	$C(\mathcal{X})$ & $C^*$- algebra of continuous complex valued \\
	& functions on a compact  space $\mathcal{X}$\\
	$C_0(\mathcal{X})$ & $C^*$- algebra of continuous complex valued functions on a locally \\
	&   compact  topological space $\mathcal{X}$ equal to $0$ at infinity\\
	$C_c(\mathcal{X})$ & Algebra of continuous complex valued functions on a \\
	&  topological  space $\mathcal{X}$ with compact support\\
	$C_b(\mathcal{X})$ & $C^*$- algebra of bounded  continuous complex valued \\
	& functions on a locally compact topological space $\mathcal{X}$ \\
	$G\left( \widetilde{\mathcal{X}}~ |~ \mathcal{X}\right) $ & Group of covering transformations of covering  $\widetilde{\mathcal{X}} \to \mathcal{X}$ \cite{spanier:at}  \\
	$\H$ & Hilbert space \\
	$\mathcal{K}= \mathcal{K}\left(\H \right) $ & $C^*$- algebra of compact operators on the separable Hilbert space $\H$  \\
	$K(A)$ & Pedersen ideal of $C^*$-algebra $A$\\
	$C^*\text{-}\varinjlim$ & $C^*$-inductive limit \\
	$\varprojlim$ & Inverse limits of groups and topological spaces \\
	$M(A)$  & A multiplier algebra of $C^*$-algebra $A$\\
	$\mathbb{M}_n(A)$  & The $n \times n$ matrix algebra over $C^*$-algebra $A$\\
	$\mathbb{N}$  & A set of positive integer numbers\\
	$\mathbb{N}^0$  & A set of nonnegative integer numbers\\
	$\supp ~a$ & Support of $a \in C_b\left( \mathcal X\right)$ \\
	$\rep_x$ or $\rep^A_x$ & An irreducible representation $A \to B\left(\H \right)$  which corresponds\\ &  to a point $x \in \hat A$ of spectrum of $A$ (cf. \eqref{rep_x_eqn}).\\ 
	

	
	$\mathbb{Z}$ & Ring of integers \\
	
	$\mathbb{Z}_n$ & Ring of integers modulo $n$ \\
	$\overline{k} \in \mathbb{Z}_n$ & An element in $\mathbb{Z}_n$ represented by $k \in \mathbb{Z}$  \\
	$X \setminus A$ & Difference of sets  $X \setminus A= \{x \in X \ | \ x\notin A\}$\\
	$|X|$ & Cardinal number of a finite set $X$\\ 
	$\left[x\right]$ & The range projection of element $x$ of a von Neumann algebra.\\ 
	$f|_{A'}$& Restriction of a map $f: A\to B$ to $A'\subset A$, i.e. $f|_{A'}: A' \to B$\\ 
	\hline
\end{tabular}

\section{$C^*$-inductive limits of nonunital algebras}
  \begin{definition}\label{connected_c_a_defn}
	We say that a $C^*$-algebra $A$ is \textit{connected} if it cannot be represented as a direct sum  $A \cong A' \oplus A''$ of nontrivial $C^*$-algebras $A'$ and $A''$.
	
\end{definition}
\begin{definition}\label{connected_comp_defn}
	A connected closed two-sided ideal $A$ of  $C^*$-algebra $B$ is said to be a \textit{connected component of}  $B$ is there is a direct sum $B = A \oplus A'$ of $C^*$-algebras.
\end{definition}
\begin{exercise}\label{conn_comp_exer}
Using the Theorems  \ref{ideal_spectrum_thm} and \ref{dauns_hofmann_thm} prove that there is a natural one-to-one correspondence between connected components of $C^*$-algebra and quasi-components of its spectrum (cf. Definitions \ref{top_quasi_component_defn}, \ref{spectrum_prime_primtive_defn}).
\end{exercise}
\begin{empt}\label{unital_notation_empt}
	Let $A$ be a $*$-algebra. Denote by $A^{\sim}$ the unital $*$-algebra given by
	\be\label{unital_notation_eqn}
	A^\sim\bydef\begin{cases}
		A & \text{if }A \text{ is unital}\\
		A^+ & \text{if }A \text{ is not unital}\\
	\end{cases}
	\ee
	where $A^+$ is given by \eqref{a_plus_eqn}.
	If $A$ is nonunital $C^*$-algebra then 	$A^{\sim}$
	is the minimal unitization of $A$ (cf. Definition \ref{multiplier_min_defn}).
\end{empt}
\begin{definition}\label{principal_non_defn}
	An injective $*$-homomorphism $\phi: A \hookto B$ of $C^*$-algebras is said to be \textit{unital} if it is unital in the sense of the Definition \ref{principal_defn} or can be uniquely extended up to the unital (in sense of the Definition \ref{principal_defn}) $*$-homomorphism $\phi^{\sim}: A^{\sim} \hookto B^{\sim}$ where $A^{\sim}\bydef A^+$ and $B^{\sim} \bydef B^+$ are minimal unitizations of $A$ and $B$ (cf. equation \eqref{unital_notation_eqn} and Definition \ref{multiplier_min_defn}).
\end{definition}
\begin{remark}
	The Definition \ref{principal_non_defn}  is a generalization of the Definition \ref{principal_defn}.
\end{remark}

\begin{empt}\label{inductive_empt}
	Let $\left\{A_\la\right\}_{\la \in \La}$ be a family of $C^*$-algebras where $\La$  denotes an  directed set (cf. Definition \ref{directed_set_defn}). Suppose that for every $\mu, \nu$ with $\mu \le \nu$, there exists the unique unital (in sense of the Definition \ref{principal_non_defn}) injective $*$-homomorphism  $f_{\mu\nu}: A_\mu \hookto A_\nu$ satisfying
	$
	f_{\mu\nu} = f_{\mu\la}\circ f_{\la\nu}$ {where} $\quad \mu < \la < \nu.
	$
	If $C^*$-algebras $A_\la$ are  not unital 	then there are natural unital  (in sense of the Definition \ref{principal_defn}) unique injective $*$-homomorphisms  $f^{\sim}_{\mu\nu}: A^{\sim}_\mu \to A^{\sim}_\nu$ of minimal unitizations (cf. equation \eqref{unital_notation_eqn}). From the Theorem \ref{inductive_lim_thm} it turns out that $C^*$-{inductive limit}  $C^*\text{-}\varinjlim_\La A^{\sim}_{\la}$ of $\left\{A^{\sim}_\la \right\}_{\la\in\La}$. 
\end{empt}
\begin{definition}\label{inductive_lim_non_defn}
	In the situation of \ref{inductive_empt} consider injective $*$-homomorphisms $A_\la \hookto C^*\text{-}\varinjlim_\La A^{\sim}_{\la}$ which are compositions $A_\la \hookto  A_\la^\sim\hookto  C^*\text{-}\varinjlim_\La A^{\sim}_{\la}$.	The $C^*$-norm completion of the union  $\cup_{\la \in \La}A_\la\subset C^*\text{-}\varinjlim_\La A^{\sim}_{\la}$  is said to be the $C^*$-\textit{inductive limit} of $\left\{A_\la \right\}$. It is denoted by $C^*\text{-}\varinjlim_{\la\in\La} A_\la$ or  $C^*\text{-}\varinjlim A_\la$.
\end{definition}

\begin{remark}
	The Definition \ref{inductive_lim_non_defn} is a generalization of the Definition \ref{inductive_lim_defn}. There are  described below generalizations of the Theorem \ref{inductive_lim_thm} and the Proposition  \ref{inductive_lim_prop}.
\end{remark}
\begin{lemma}\label{state_un_lem}
	If $A^+$ is the minimal unitization of a nonunital $C^*$-algebra $A$ and $\Om_A$, $\Om_{A^+}$ are state spaces of $A$ and  $A^+$ respectively then one has
	\be\label{state_un_eqn}
	\Om_A = \left\{\tau^+ \in \Om_{A^+}~|~ \tau^+\left(0\oplus 1 \right)= 0;\quad \mathrm{where} \quad 0\oplus 1 \in A \oplus \C \cong A^+  \right\}.
	\ee
\end{lemma}
\begin{proof}
	Any state $\tau: A\to \C$ induces the  state $\tau^+: A^+\to \C$ given by $\tau^+\left(a\oplus \la  \right)\bydef \tau\left(a \right)+\la$ for any $a \oplus \la \in  A \oplus \C \cong A^+$. Clearly  $\tau^+$ satisfies to  \eqref{state_un_eqn}. 
\end{proof}

\begin{corollary}\label{inductive_lim_state_nor_cor}
	If a nonunital $C^*$-algebra $\widehat A$ is a $C^*$-inductive limit (in sense of the Definition \ref{inductive_lim_non_defn}) of $\left\lbrace A_\la\right\rbrace_{\la \in \La}$, then the
	state space $\Om$ of $\widehat A$ is homeomorphic to the projective limit of the state spaces $\Om_\la$ of $A_\la$. 
\end{corollary}
\begin{proof}
	Let us consider unital  injective $*$-homomorphisms  $f^+_{\mu\nu}: A^+_\mu \to A^+_\nu$ (in the sense of the Definition \ref{principal_defn})  of the minimal unitizations and let $\widehat{A}^+$ be a $C^*$-inductive limit (in sense of the Definition \ref{inductive_lim_non_defn}) of $\left\{A^+_\la\right\}$. From the Theorem \ref{inductive_lim_state_thm} it follows that the state space $\Om_{\widehat{A}^+}$ is the projective limit $\Om^+_\la$ of $A^+_\la$. From the Lemma \ref{state_un_lem} it turns that $\Om_{\widehat{A}} \subset \Om_{\widehat{A}^+}$ and  $\Om_\la \subset \Om^+_\la$ and for any $\la\in \La$. Moreover every state $\widehat{\tau} \in \Om_{\widehat{A}}$ is mapped onto $\tau_\la \in \Om_\la$. It follows that $\Om$ is homeomorphic to the projective limit of the state spaces $\Om_\la$.
\end{proof}
\begin{exercise}\label{inductive_lim_state_exer}
Prove a related to nonunital commutative $C^*$-algebras generalization of the Corollary \ref{inductive_lim_state_cor}.
\end{exercise}

\section{Ext-functor and characters}\label{ext_char_sec}
\paragraph*{}
Let $A \bydef \Z^n /\Gamma$ be a finite Abelian group where $\Gamma\subset \Z^n \subset \Q^n$ is a subgroup. The inclusion $\Gamma\subset \Z^n$ provides the natural homomorphism
$$
 \Hom \left(\Z^n, \Z \right)\to \Hom \left(\Ga, \Z\right)
$$
and 
$$
\Ext\left( A, \Z\right) =  \Hom \left(\Ga, \Z\right)/\Hom \left(\Z^n, \Z \right).
$$
Any $\phi \in  \Hom \left(\Ga, \Z\right)$ can be uniquely extended up to $\psi \in \Hom \left(\Q^n, \Q\right)$. One has a character 
\bean
\varphi : \Q^n \to e^{\frac{i\Q}{2\pi}},\\
x \mapsto e^{\frac{i\psi\left(x \right) }{2\pi}}.
\eean
 The restriction $\varphi|_{\Z^n}: \Z^n \to e^{\frac{i\Q}{2\pi}}$ is a character on $\Z^n$. From $\phi\left(\Ga \right) \subset  \Z$ it follows that $\varphi|_{\Z^n}\left( \Ga\right) = \{1\}$, i.e. $\varphi|_{\Z^n}$ is trivial on $\Ga$. So $\varphi|_{\Z^n}$ yields a homomorphism $\chi_\phi: \Z^n / \Ga \to e^{\frac{i\Q}{2\pi}}$, i.e. $\chi_\phi$ is a character.  If $\phi \in \im ~ \Hom \left(\Z^n, \Z \right)$ then $\phi\left( \Z^n \right) \subset \Z$ and $\varphi|_{\Z^n}$ yields  a trivial character. Otherwise if the character $\varphi|_{\Z^n}$ is trivial then $\phi \in \im ~ \Hom \left(\Z^n, \Z \right)$. So any element of   $\Hom \left(\Ga, \Z\right)/\Hom \left(\Z^n, \Z \right)$ yield a character $\Z^n / \Ga\to e^{\frac{i\Q}{2\pi}}$.
 
Conversely a character $\chi : \Z^n/\Ga \to e^{\frac{i\Q}{2\pi}}$ naturally yields a character $\th : \Z^n \to e^{\frac{i\Q}{2\pi}}$. Let $\left(e_1, ..., e_n\right)\subset \Z^n$ be generators of $\Z^n$ and for any $j = 1, ..., n$ we select $\xi_j \in \Q$ such that $\chi\left(e_j\right)=e^{\frac{i\xi_j}{2\pi}}$. Define 
\bean 
\varphi: \Z^n \to \Q;\\
\forall\left(a_1, ..., a_n\right)\in \Z^n\quad \varphi \left( a_1e_1 + ... + a_1e_n\right)  = a_1\xi_1 + ... + a_1\xi_n
\eean 
From $\chi\left(\Ga\right)= \{1\}$ it follows that $\varphi\left( \Ga\right) = \Z$ i.e. $\varphi|_\Ga\in \Hom\left(\Ga, \Z \right)$.
The reader can prove that the class  $\varphi|_\Ga + \Hom \left(\Ga, \Z \right)\in \Hom \left(\Z^n, \Z \right)/\Hom \left(\Ga, \Z\right)$ does not depend on  a choice of $\xi_1, ..., \xi_n$. Thus  one concludes that $\Ext\left(A, \Z \right)$ is naturally isomorphic to a group of characters of $A= \Z^n/\Ga$.

\section{Representations of hereditary subalgebras}

\begin{empt}\label{hered_repr_p_empt}
	If $\rho: A\hookto B\left(\H \right)$ be a faithful nondegenerate representation then there is a $C^*$-algebra $A''$ such that
	\begin{itemize}
		\item $A''$ is isomorphic to the bicommutant $\rho\left( A\right)''$ in $B\left(\H\right)$ (cf. Definition \ref{commutant_defn}) of $A$,
		\item there is the natural inclusion $A \subset A''$,
		\item one has a natural extension 
		\be\label{rho_ext_eqn}
		\rho'': A'' \hookto B\left(\H\right)
		\ee
		of the representation $\rho$.
	\end{itemize}
	Let $B$ be a hereditary subalgebra $A$, and let  $\left\{u_\la \right\}_{\la \in \La} \subset  B\left(\H \right)_+$ be an increasing  net of positive elements, such that
	\begin{itemize}
		\item $\left\{u_\la \right\}\subset M\left(A \right) \cap M\left( B\right)$.
		\item 
		\be\label{hered_r_u_eqn}
		\forall\la\in\La\quad	\rho\left(u_\la A u_\la \right)\subset \rho\left(B  \right).  
		\ee
		\item There is a limit $\bt$-$\lim_\la u_\la = 1_{ M\left( B\right) }$ with respect to the strict topology of $M\left( B\right)$ (cf. Definition \ref{strict_topology_defn}). 
	\end{itemize}
	The existence ot a net $\left\{u_\la \right\}$ follows from the Lemma \ref{hered_ideal_lem} and the Theorem \ref{left_ideal_thm}. 
	One has
	\be\label{hered_uau_eqn}
	B = \left\{ a \in A \left|~\lim_{\la\in\La} \left\| a - u_\la a u_\la\right\|= 0 \right.\right\}
	\ee
	From the Lemma \ref{increasing_convergent_w_lem} the net $\left\{u_\la \right\}$  is convergent  with respect to the strong topology of $B\left(\H \right)$ (cf. Definition \ref{strong_topology_defn}). If $p \bydef s$-$\lim\rho\left(  u_\la \right) \in B\left( \H\right)$ is a strong limit then $p$ lies in the strong closure of $\rho(A)$ (cf. Definition \ref{strong_topology_defn}). Form the Theorem \ref{von_Neumann_thm} it follows that $p\in \rho\left(A \right)''= \rho''\left(A'' \right)$, i.e.
	\be\label{hered_repr_abc_eqn}
	\exists  p'' \in A''\quad p = \rho''\left(p'' \right).
	\ee
	From  $1_{ M\left( B\right) }= 1^*_{ M\left( B\right) }$ and $1^2_{ M\left( B\right) }=1_{ M\left( B\right) }$ it follows that $p$ is a projector, i.e. $p^*=p$ and $p^2 = p$. From \eqref{hered_uau_eqn} it turns out that 
	\be\label{hered_pap_eqn}
	\forall a \in B \quad  p\rho\left( a\right)p= \rho \left( a\right), \quad p'' a p'' = a,
	\ee
	and  there is a natural representation
	\be\label{hered_rep_eqn}
	B \hookto B\left(p \H \right) 
	\ee
\end{empt}
\begin{remark}\label{hered_dense_rem}
	If $\left\{u_\la \right\}_{\la \in \La} \subset \rho\left(B \right)$ then from $p \bydef s$-$\lim_{\la\in\La} \rho\left( u_\la\right)$ it follows that $p\H$ is a norm completion of $\rho\left(B\right)\H$.
\end{remark}
\begin{lemma}\label{hered_full_lem} 
Under the hypotheses \ref{hered_repr_p_empt} the representation \eqref{hered_rep_eqn} is  faithful.
\end{lemma}
\begin{proof}
	For any $a \in B$ there is $\xi \in \H$ such that $\rho\left(a\right)\xi \neq 0$. From  \eqref{hered_rep_eqn} it follows that $\rho\left( a\right) \xi = p \rho\left( a\right) p\xi \neq 0$. On the other hand $p\xi \in p\H$.
\end{proof}

\begin{lemma}\label{hered_nondegenerate_lem} 
Under the hypotheses \ref{hered_repr_p_empt} the representation \eqref{hered_rep_eqn} is  nondegenerate.
\end{lemma}
\begin{proof}
	For any $\xi \in p\H\setminus\{0\}$ there is $a \in A$ such that $a \xi \neq 0$. From the equation \eqref{four_decompositon_eqn} we can suppose that $a$ is positive. One has $\lim_{\la \in \La}u_\la \xi = \xi$, it follows that  $\lim_{\la \in \La}a u_\la \xi = a\xi$. There is $\la_0\in \La$ such that $a u_{\la_0} \xi \neq 0$, so one has
	\bean
	\left(a u_{\la_0} \xi, a u_{\la_0} \xi \right) = \left(\xi , \left( u_{\la_0} a^*a u_{\la_0}\right) \xi  \right)\neq 0\quad \Rightarrow \quad \left(  u_{\la_0} a^*a u_{\la_0}\right) \xi	\neq 0.
	\eean
	However $ u_{\la_0} a^*a u_{\la_0}\in B$.
\end{proof}
\begin{corollary}\label{hered_representation_cor}
	Let $A$ be a $C^*$-algebra and let $\rho: A\hookto B\left( \H\right)$ be a faithful, nondegenerate representation. If $\left\{u_\la \right\}_{\la\in \La}\subset B$ is an approximate unit of $B$ (cf. Definition \ref{approximate_unit_defn}) and $p \bydef s$-$\lim_{\la\in\La} \rho\left( u_\la\right)$ is a strong limit in $B\left(\H \right)$  (cf. Definition \ref{strong_topology_defn})
	then  there is a faithful, nondegenerate representation
	$$
	B \hookto B\left(p \H \right). 
	$$
\end{corollary}
\begin{lemma}\label{hered_repr_p_lem}
Under the hypotheses \ref{hered_repr_p_empt} if $\pi: A\to B\left(\H \right)$ is an atomic  representation (cf. Definition \ref{atomic_repr_defn}) then the given by \eqref{hered_rep_eqn} representation $ B \hookto B\left(p \H \right)$ is atomic.
\end{lemma}
\begin{proof}
	One can prove this lemma by usage of the proof of the Lemma 
	\ref{hered_repr_lem}.
\end{proof}
\begin{lem}\label{top_ul_lem}
	If $\sX$ is a locally compact Hausdorff space and $\left\{\sU_\la\right\}_{\la \in \La}$ is a family of subsets of $\sX$ such that the union $\cup_{\la \in \La} C_0\left( \sU_\la\right)$ is dense in $C_0\left( \sX\right)$ then $\sX = \cup_{\la \in \La} \sU_\la$.
\end{lem}
\begin{proof}
	If $x \in \sX \setminus \cup_{\la \in \La} \sU_\la$ and $f \in C_0\left( \sX\right)$ is such that $f\left( x\right)= 1$ then
	$$
	\inf_{g \in \cup_{\la \in \La} C_0\left( \sU_\a\right)} \left\|f - g \right\| \ge 1,
	$$
	i.e. $\cup_{\la \in \La} C_0\left( \sU_\a\right)$ is not dense in $C_0\left( \sX\right)$.
	
\end{proof}
\begin{empt}\label{top_hered_empt}
	Let $B\subset C_0\left(\sX \right)$ be a hereditary $C^*$-subalgebra (cf, Definition \ref{hered_defn}),  and let $\left\{u_\la\right\}_{\la \in \La}$ be an approximate unit of $B$ (cf. Definition \ref{approximate_unit_defn}). A set
	\be\label{top_hered_eqn}
	\sU \bydef \bigcup_{\la\in\La}\left\{x \in \sX | u_\la\left(x \right) > 0\right\}
	\ee
	is an open subset of $\sX$ because it is a union of open sets.
\end{empt}

\begin{lemma}\label{top_hered_lem}
	Under the hypotheses \ref{top_hered_empt} one has $B = C_0\left(\sU \right)$.
\end{lemma}
\begin{proof}
	For any $\la\in \La$ denote by $\sU_\la \bydef \left\{x \in \sX | u_\la\left(x \right) > 0\right\}$. The union $\cup_{\la \in \La} C_0\left( \sU_\la\right)$ is dense in $C_0\left( \sU_\la\right)$ so this lemma follows from the \ref{top_ul_lem} one.
	
\end{proof}

\begin{exercise}\label{top_u_net_exer}
	Let $\sX$ be a locally compact, Hausdorff space.\\ Denote by $\Xi\left( \sX\right)\bydef \left\{u_\la\right\}_{\la \in \Xi\left( \sX\right)}\subset C_c\left( \sX\right)_+$ a net of all positive continuous maps such that 
	\bean
	\forall	\la \in \Xi\left( \sX\right)\quad u_\la\left(\sX\right)\subset \left[0,1\right],\\
	\forall	\mu, \nu \in \Xi\left( \sX\right)\quad	\mu \le \nu \quad \Leftrightarrow\quad u_\mu \le u_\nu.
	\eean
	Prove that the net $\Xi\left( \sX\right)$ is an approximate unit of $C_c\left( \sX\right)$. 
\end{exercise}

\begin{definition}\label{top_char_f_defn}
	If $\sX$ is a set and $\sU \subset \sX$  then a map 
	\be\label{top_char_f_eqn}
	\begin{split}
		\chi_\sU : \sX \to \R,\\
		x \mapsto	\begin{cases}
			1 & x \in \sU\\
			0 & x \in \sX\setminus \sU
		\end{cases}
	\end{split}
	\ee
	is said to be a \textit{characteristic function} of $\sU$.
\end{definition}

\begin{empt}\label{top_u_net_rep_empt}
	Let $\sX$ be a locally compact, Hausdorff space, and let $\rho: C_0\left(\sX\right)\to B\left(\H \right)$ be a faithful representation, and let $C_0\left(\sX\right)''\subset  B\left(\H \right)$ be a bicommutant (cf. There is an inclusion $C_0\left(\sX\right)\subset C_0\left(\sX\right)''$ such that the representation $\rho$ can be extended up to $\rho'': C_0\left(\sX\right)''\hookto B\left(\sH\right)$. If $\chi_\sU$ is a characteristic function of a Borel set $\sU$ (cf. Definition \ref{top_char_f_defn}) then one can suppose that $\chi_\sU\in C_0\left(\sX\right)''$.
	Let   $\left\{u_\la\right\}_{\la \in \Xi\left( \sX\right)}\subset C_0\left( \sX\right)_+$ be an explained in the Exercise \ref{top_u_net_exer} net. From the Lemma  \ref{increasing_convergent_w_lem} it turns out that the family $\left\{\rho\left( u_\la\right) \right\}_{\la \in \Xi\left( \sX\right)}$ is convergent with respect to a strong topology of $B\left( \sH\right)$ (cf. Definition \ref{strong_topology_defn}). Moreover the strong limit $s$-$\lim\rho\left(  u_\la \right) \in B\left( \H\right)$ equals to $\rho''\left(\chi_\sU \right)$. 
\end{empt}
\section{Strict and strong limits}

\begin{lemma}\label{lift_mult_lem}
	Let $\pi: A \hookto M\left( \widetilde A\right)$ be an injective *-homomorphism of $C^*$-algebras, and let $\left\{u_\a\right\}_{\a \in \mathscr A}\subset A$ be an approximate unit for $A$ (cf. Definition \ref{approximate_unit_defn}). If
	\be\label{lift_mult_norm_eqn}
	\forall \widetilde a \in \widetilde A \quad \lim_{\a \in \mathscr A}\left\| \widetilde a - \pi\left( u_\a\right)  \widetilde a \right\| = \lim_{\a \in \mathscr A}\left\| \widetilde a -  \widetilde a \pi\left( u_\a\right)\right\| = 0
	\ee
	i.e. $\bt\text{-}\lim_{\a \in \mathscr A} \pi\left( u_\a\right)= 1_{M\left( \widetilde A\right)}$,
	then $\pi$ can be naturally  extended up to an injective *-homomorphism 
	\be\label{lift_mult_eqn}
	\begin{split}
		M\left(\pi \right) : M\left( A\right)  \hookto M\left( \widetilde A\right),\\
		a \mapsto  \bt\text{-}\lim_{\a \in \mathscr A} \pi\left( u_\a  a \right)= \bt\text{-}\lim_{\a \in \mathscr A} \pi\left( a u_\a \right)= \bt\text{-}\lim_{\a \in \mathscr A} \pi\left( u_\a a u_\a \right)
	\end{split}
	\ee
	where $\bt\text{-}\lim_{\a \in \mathscr A}$ means the limit with respect to the strict topology  of $M\left( \widetilde A\right)$ (cf. Definition \ref{strict_topology_defn}).
\end{lemma}
\begin{proof}
	For any $a \in M\left( A\right)$ we define both  maps:
	\bean 
	L_a : \widetilde A \to \widetilde A,\\
	\widetilde a \mapsto  \lim_{\a \in \mathscr A }\widetilde a~\pi \left( u_\a a \right) = \lim_{\a \in \mathscr A }\widetilde a~\pi \left( u_\a a u_\a\right);\\
	R_a : \widetilde A \to \widetilde A,\\
	\widetilde a \mapsto  \lim_{\a \in \mathscr A }\pi \left( a u_\a\right)  \widetilde a = \lim_{\a \in \mathscr A }\pi \left(u_\a a u_\a\right)  \widetilde a 
	\eean
	where the convergence with respect to $C^*$-norm topology is is implied.
	From 
	$$
	\forall \widetilde{a}, \widetilde{b} \in \widetilde{A}\quad \left( \widetilde a \pi \left( u_\a a \right)\right)  \widetilde b=\widetilde a \left( \pi \left( u_\a a \right) \widetilde b\right) 
	$$
	one can deduce that a pair $\left( L_a, R_a\right)$ satisfies to  the equation \eqref{double_centralizer_eqn}, i.e.  $\left( L_a, R_a\right)$ is  a {double centralizer} (cf. Definition \ref{double_centralizer_defn}). From the Remark \ref{double_centralizer_rem} it follows that $\left( L_a, R_a\right)$ yields an element of $M\left( \widetilde{A}\right)$, so one has a natural map $	M\left(\widetilde \pi \right) :	M\left(A \right)  \to M\left( \widetilde A\right)$. 
	We leave to the reader the proof of that $M\left(\pi \right)$ is an injective *-homomorphisms  and $M\left(\pi \right)$ is an extension of  $\pi$.
	
\end{proof}

\begin{lemma}\label{strict_strong_lem}
	Let $\pi: A \to B\left(\H \right)$ be a faithful nondegenerate representation (cf. Definitions \ref{faithful_representation_defn} and \ref{nondegenerate_repr_defn}), and let  $\left\{a_\la\right\}_{\la \in \La} \subset M\left(A\right)$ be a bounded net i.e. 
	$$
	\exists C \in \R_+ \quad \forall \la\in \La \quad \left\| a_\la  \right\|< C. 
	$$
	If the net $\left\{a_\la\right\}_{\la \in \La} \subset M\left(A\right)$ is convergent with respect to the strict topology of $M\left( A\right)$ (cf. Definition \ref{strict_topology_defn})   and a net $\left\{\pi\left( a_\la\right) \right\}_{\la \in \La}\in  B\left(\H \right)$ is convergent with respect to the strong topology of $B\left( \H\right)$ (cf. Definition \ref{strict_topology_defn}) then one has
	\be\label{srict_strong_eqn}
	\pi\left(\bt\text{-}\lim_{\la\in \La} a_\la \right)= s\text{-}\lim_{\la\in \La} \pi\left( a_\la\right)  
	\ee 
	where both  $\bt\text{-}\lim_{\la\in \La}$ and $s\text{-}\lim_{\la\in \La}$ are strict and strong limits respectively.
\end{lemma}
\begin{proof}
	Denote by $a^\bt \bydef  \bt\text{-}\lim_{\la\in \La} a_\la$ and $a^s  \bydef s\text{-}\lim_{\la\in \La} \pi\left( a_\la\right)\in B\left( \H\right)$.
	If the equation \eqref{srict_strong_eqn} is not true then there is $\xi \in \H$ such that
	$$
	\zeta \bydef\pi\left(a_\bt \right)\xi -  a^s\xi\neq 0.
	$$ 
	There is $\la_s \in \La$ such that 
	$$
	\forall \la \in \La \quad \la \ge \la_s \quad \left\| a^s \xi - \pi\left( a_\la\right) \xi  \right\|< \frac{\left\| \zeta  \right\|}{4}.
	$$
	On the other hand from the Definition \ref{strict_topology_eqn} and the Lemma \ref{nondegenerate_repr_lem} it turns out that there are $b \in A$, $\la_\bt \in \La$ and $\eta \in \H$ such that 
	\bean
	\left\| b \eta  - \xi  \right\| <  \frac{\left\| \zeta  \right\|}{4\max\left(  C, \left\| a^\bt  \right\|\right)  },\\
	\forall \la \in \La \quad \la \ge \la_\bt\quad  \left\| a^\bt b - a_\la b   \right\| <  \frac{\left\| \zeta  \right\|}{4 \left\| \eta\right\|}
	\eean
	If $\la_0\in \La$ is such that $\la_0 \ge \la_\bt$ and $\la\ge \la_s$ from the triangle identity it follows that
	\bean
	\left\| \pi\left(a_{\la_0}\right)\xi -  \pi\left(a^\bt\right)\xi  \right\|\le\\\le \left\| \pi\left(a_{\la_0}\right)\xi -  \pi\left(a_{\la_0}b\right)\eta\right\|+ \left\|  \pi\left(a_{\la_0}b\right)\eta- \pi\left(a^\bt b\right)\eta\right\|+\\ 
	+ \left\|   \pi\left(a^\bt b\right)\eta- \pi\left(a^\bt \right)\xi\right\|	< 
	\frac{\left\| \zeta  \right\|}{4}+ \frac{\left\| \zeta  \right\|}{4}+\frac{\left\| \zeta  \right\|}{4}=\frac{3\left\| \zeta  \right\|}{4};\\
	\left\| \zeta  \right\| =  \left\| \pi\left(a^\bt \right)\xi -  a^s\xi  \right\|\le\\\le 
	\left\| \pi\left(a_{\la_0}\right)\xi -  a^s\xi  \right\|+ \left\| \pi\left(a_{\la_0}\right)\xi -  \pi\left(a^\bt\right)\xi  \right\|< \frac{\left\| \zeta  \right\|}{4}+ \frac{3\left\| \zeta  \right\|}{4} = \left\| \zeta  \right\|,
	\eean  
	i.e. there is a contradiction $\left\| \zeta  \right\|< \left\| \zeta  \right\|$. From this contradiction it turns out that $\zeta = 0$ and the equality \eqref{srict_strong_eqn} is true.
\end{proof}


\section{Transitive coverings}\label{top_transitive_section}
\paragraph*{}
The notion of regular covering (cf. Definition \ref{top_regular_defn}) has no any good noncommutative generalization. So  the "regular covering" term  will be replaced by "transitive covering".

\begin{definition}\label{top_transitive_defn}
	A covering $p: \widetilde{\mathcal X}\to {\mathcal X}$ with connected $\widetilde\sX$  is said to be \textit{transitive}  
	if there is a properly discontinuous  group $G$ of homeomorphisms of a topological space $\widetilde\sX$  (cf. Definition \ref{top_properly_disc_defn}) such that $\sX \cong \widetilde\sX/G$ and $p: \widetilde{\mathcal X}\to {\mathcal X}$ is naturally equivalent to the given by the Theorem \ref{top_group_of_covering_transformations_thm} covering $\widetilde{\mathcal X} \to \widetilde{\mathcal X}/G$.
\end{definition}

\begin{remark}
	From the Theorem \ref{top_cov_fact_thm} and the Lemma \ref{top_cov_from_pi1_cor} it turns out that any regular covering is transitive. 
\end{remark}

\begin{lemma}\label{top_trans_tp_cov_lem}
	Let $\sX$ be a  connected, locally connected, locally compact, Hausdorff space, and let $p:	\widetilde{   \mathcal X } \to \mathcal X$ be a  covering of connected spaces. If $g \in \Homeo\left( \widetilde{   \mathcal X }\right)$ is such that $p \circ g = p$ then one has
	$$
	\exists  \widetilde x \in \widetilde\sX \quad  g\widetilde x = \widetilde x \quad  \Leftrightarrow\quad \forall  \widetilde x \in \widetilde\sX \quad  g\widetilde x = \widetilde x.
	$$
\end{lemma}
\begin{proof}
	The implication $\Leftarrow$ is evident. 
	If  $g\widetilde x= \widetilde x$ then there is a connected open neighborhood $\widetilde \sU$ of $\widetilde x$ which mapped homeomorphically onto $\sU \bydef p\left( \widetilde \sU\right)$. The restriction $\left.p\right|_{\widetilde\sU}$ is a homeomorphism from $\sU$. Otherwise there is an open  neighborhood $\widetilde \sV$ of $\widetilde x$ such that $g\widetilde  \sV \subset \widetilde \sU$.
	Taking into account $p \circ g = p$ one has $$\forall \widetilde x \in \widetilde \sU \cap \widetilde  \sV\quad g\widetilde x = \left(\left.p\right|_{\widetilde\sU}\right)^{-1} \circ p \circ g \left( \widetilde x \right) = \left(\left.p\right|_{\widetilde\sU}\right)^{-1} \circ p\left( \widetilde x \right)=  \widetilde x,$$ so a set
	$$
\widetilde{ \mathcal W}	=	\left\{\left.\widetilde x \in \widetilde \sX~\right| g\widetilde x= \widetilde x\right\}
$$
is open. If $\widetilde{ \mathcal X}\setminus \widetilde{ \mathcal W}\neq \emptyset$ and  for any $\widetilde x \in \widetilde{ \mathcal X}\setminus \widetilde{ \mathcal W}$ there is and open neighborhood $\widetilde \sV$ such that $\widetilde{ \mathcal V}\subset\widetilde{ \mathcal X}\setminus \widetilde{ \mathcal W}$ then the set $\widetilde{ \mathcal X}\setminus \widetilde{ \mathcal W}$ is open and $\widetilde\sX$ is a disjoint union of nontrivial open sets. It is impossible since the space $\widetilde{ \mathcal X}$ is connected. So there is $\widetilde x\in \widetilde{ \mathcal X}\setminus \widetilde{ \mathcal W}$ such then for any open neighborhood $\widetilde \sV$ of $\widetilde x$ one has $\widetilde \sV\cap \widetilde{ \mathcal W}\neq \emptyset$. From our previous proof it follows that
\bean
\widetilde{ \mathcal W}\subsetneqq \widetilde{ \mathcal W}\cup \widetilde \sV,\\
\widetilde{ \mathcal W}\cup \widetilde \sV\subset \left\{\left.\widetilde x \in \widetilde \sX~\right| g\widetilde x= \widetilde x\right\},
\eean
so $\widetilde{ \mathcal W}\subsetneqq \widetilde{ \mathcal W}$. It is impossible so $\widetilde{ \mathcal X}\setminus \widetilde{ \mathcal W}=\emptyset$ and $\widetilde{ \mathcal X}=\widetilde{ \mathcal W}$.
\end{proof}
\begin{corollary}\label{top_trans_tp_cov_cor}
	Let $\sX$ be a  connected,  locally compact, Hausdorff space, and let $p:	\widetilde{   \mathcal X } \to \mathcal X$ be a   covering of connected spaces. If there is  $\widetilde x_0\in \widetilde\sX$ such that the map
	\bean 
	G\left(\widetilde{   \mathcal X }~|~{   \mathcal X } \right)\xrightarrow{\approx}p^{-1}p\left(\widetilde x_0 \right);\\
	g \mapsto g\widetilde x_0. 
	\eean
	is bijective	then $p$ is a transitive covering.
\end{corollary}

\begin{corollary}\label{top_trans_cov_factor_cor}
	If $p: \widetilde{\mathcal X}\to {\mathcal X}$  is a {transitive} covering then there is the natural homeomorphism $\sX  \cong \widetilde{\mathcal X}/G\left(\left.\widetilde{\sX}~\right|\sX\right)$.
\end{corollary}
\begin{remark}
	The Corollary \ref{top_trans_cov_factor_cor} can be regarded as a generalization of the Theorem \ref{top_cov_fact_thm}.
\end{remark} 

\begin{theorem}\label{top_covp_cat_thm} 
	Consider a commutative triangle of connected, Hausdorff topological spaces and  continuous maps
	\newline
	\begin{tikzpicture}
		\matrix (m) [matrix of math nodes,row sep=3em,column sep=4em,minimum width=2em]
		{
			\widetilde{\mathcal X}_1 & &\widetilde{\mathcal X}_2\\ 
			& {\mathcal X}\\};
		\path[-stealth]
		(m-1-1) edge node [above] {$p$} (m-1-3)
		(m-1-1) edge node [left]  {$p_1~~$} (m-2-2)
		(m-1-3) edge node [right] {$~~p_2$} (m-2-2);
	\end{tikzpicture}
	\\
where the space $\sX$ is locally connected.	If   both  $p_1$ and $p_2$ are transitive coverings (cf. Definition \ref{top_transitive_defn}) and the map $p$ is surjective then the map $p$ is a transitive covering and there is the natural exact sequence  
	\be\label{top_covp_cat_eqn}
	\{e\}\to G\left( \left.{\widetilde{\sX}}_1~\right|\widetilde{\sX}_2\right)	\to 	G\left( \left.{\widetilde{\sX}}_1~\right|{\sX}\right)\to G\left( \left.{\widetilde{\sX}}_2~\right|{\sX}\right)\to\{e\}
	\ee
	of groups and homomorphisms.	
\end{theorem}

\begin{proof}
	From the Corollary \ref{top_cov_cat_cor} it follows that the map $p$ is a transitive covering. Let both $\widetilde{x}_1\in \widetilde{\sX}_1$ and $g \in G\left( \left.{\widetilde{\sX}}_1~\right|{\sX}\right)$ are such that
	$
p_2\left(g \widetilde{x}_1\right) \neq p_2\left( \widetilde{x}_1\right)	
	$. 
There are connected open neighborhoods $\widetilde\sU'_2$ and $\widetilde\sU''_2$ of $p\left( \widetilde x_1\right) $ and $p\left( \widetilde gx_1\right)$ such that $\widetilde\sU'_2\cap \widetilde\sU''_2 = \emptyset$. Otherwise there is open connected neighborhood $\widetilde\sV_1$ of $\widetilde x_1$ such that $p\left( g\widetilde\sV_1\right)  \subset \widetilde\sU''_2$. It follows that $p\left(  g\widetilde y_1\right) \neq p\left( \widetilde y_1\right)$ for all $\widetilde y \in \widetilde \sU'_1\cap \widetilde\sV_1$. Hence the set
$$
\widetilde{ \mathcal W}_1\bydef \left\{\left.\widetilde y_1 \in \widetilde \sX_1~\right| p\left(  g\widetilde y_1\right) \neq p\left( \widetilde y_1\right)\right\}
$$
is open. If $\widetilde{ \mathcal X}_1\setminus \widetilde{ \mathcal W}_1\neq \emptyset$ and  for any $\widetilde x_1 \in \widetilde{ \mathcal X}_1\setminus \widetilde{ \mathcal W}_1$ there is and open neighborhood $\widetilde \sV_1$ such that $\widetilde{ \mathcal V}_1\subset\widetilde{ \mathcal X}_1\setminus \widetilde{ \mathcal W}_1$ then the set $\widetilde{ \mathcal X}_1\setminus \widetilde{ \mathcal W}_1$ is open and $\widetilde\sX_1$ is a disjoint union of nontrivial open sets. It is impossible since the space $\widetilde{ \mathcal X}_1$ is connected. So there is $\widetilde x'_1\in \widetilde{ \mathcal X}_1\setminus \widetilde{ \mathcal W}_1$ such then for any open neighborhood $\widetilde \sV_1$ of  $\widetilde x'_1$ one has $\widetilde \sV_1\cap \widetilde{ \mathcal W}_1\neq \emptyset$. From our previous proof it follows that
\bean
\widetilde{ \mathcal W}\subsetneqq \widetilde{ \mathcal W}\cup \widetilde \sV,\\
\widetilde{ \mathcal W}\cup \widetilde \sV\subset \left\{\left.\widetilde x \in \widetilde \sX~\right| g\widetilde x= \widetilde x\right\},
\eean
so $\widetilde{ \mathcal W}_1\subsetneqq \widetilde{ \mathcal W}_1$. It is impossible so $\widetilde{ \mathcal X}_1\setminus \widetilde{ \mathcal W}_1=\emptyset$ and $\widetilde{ \mathcal X}_1=\widetilde{ \mathcal W}_1$. From our construction it follows that for any $g \in G\left( \left.{\widetilde{\sX}}_1~\right|{\sX}\right)$ one has
\be\label{top_free_act}
\begin{split}
	\exists  \widetilde x_1 \in \widetilde\sX_1 \quad p \circ g\left( \widetilde x_1\right) \neq p\left( \widetilde x_1\right)  \quad  \Leftrightarrow\quad \forall   \widetilde x_1 \in \widetilde\sX_1 \quad p \circ g\left( \widetilde x_1\right)\neq p\left( \widetilde x_1\right),\\
\exists  \widetilde x_1 \in \widetilde\sX_1 \quad p \circ g\left( \widetilde x_1\right)  = p\left( \widetilde x_1\right)  \quad  \Leftrightarrow\quad \forall   \widetilde x_1 \in \widetilde\sX_1 \quad p \circ g\left( \widetilde x_1\right)= p\left( \widetilde x_1\right) . 
\end{split}
\ee 
	If   $\varphi : \widetilde{\sX}_1\xrightarrow{\approx }  \widetilde{\sX}_1$ is   a homeomorphism such that
	$$
	\forall \widetilde{x}_1\in \widetilde{\sX}_1\quad p\left( \widetilde{x}_1\right) = p\left( \varphi\left(\widetilde{x}_1\right)\right).
	$$
	then
	$$
	\forall \widetilde{x}_1\in \widetilde{\sX}_1\quad p_1\left( \widetilde{x}_1\right) = p_1\left( \varphi\left(\widetilde{x}_1\right)\right).
	$$
	so one has $\varphi \in  G\left( \left.{\widetilde{\sX}}_1~\right|{\sX}\right)$. In result one has
	
	\bean 
	G\left( \left.{\widetilde{\sX}}_1~\right|\widetilde{\sX}_2\right) \bydef \left\{\left.g \in  G\left( \left.{\widetilde{\sX}}_1~\right|{\sX}\right) ~\right|\forall \widetilde x_1\in \widetilde\sX_1\quad p\left(  \widetilde x_1\right) = p\left(  g\widetilde x_1\right) \right\}
	\eean
	and there is an injective homomorphism of groups
	\bean
	G\left( \left.{\widetilde{\sX}}_1~\right|\widetilde{\sX}_2\right)	\hookto 	G\left( \left.{\widetilde{\sX}}_1~\right|{\sX}\right).
	\eean
	For any   $\widetilde x_1\in \widetilde{\sX}_1$
	\be\label{top_guuu_eqn}
	\forall g_2 \in G\left(  \left.{\widetilde{\sX}}_2~\right|{\sX}\right)\quad \exists g_1 \in  G\left( \left.{\widetilde{\sX}}_1~\right|{\sX}\right)\quad  g_2 p\left(\widetilde x_1  \right)= p \left(g_1\widetilde x_1  \right)
	\ee
	and
	\bean
	\forall g', g'' \in G\left( \left.{\widetilde{\sX}}_1~\right|\widetilde{\sX}_2\right)\quad g_2 p\left(\widetilde x_1  \right)= p \left(g'g_1\widetilde x_1  \right)= p \left(g_1 g''\widetilde x_1  \right),\\
	g_1G\left( \left.{\widetilde{\sX}}_1~\right|\widetilde{\sX}_2\right) = G\left( \left.{\widetilde{\sX}}_1~\right|\widetilde{\sX}_2\right)g_1,
	\eean
	or equivalently $G\left( \left.{\widetilde{\sX}}_1~\right|\widetilde{\sX}_2\right)$ is a normal subgroup of $G\left( \left.{\widetilde{\sX}}_1~\right|{\sX}\right)$. From \eqref{top_guuu_eqn} it follows that there ia a surjective homomorphism $G\left( \left.{\widetilde{\sX}}_1~\right|{\sX}\right)\to G\left( \left.{\widetilde{\sX}}_2~\right|{\sX}\right)$. From \eqref{top_free_act} it follows that  for any $\widetilde x_2\in \widetilde\sX_2$ one has $p^{-1}\left( \widetilde x_2\right)=G \left( \left.{\widetilde{\sX}}_1~\right|\widetilde{\sX}_2\right)\widetilde x_1$   where $\widetilde x_1 \in p^{-1}\left( \widetilde x_2\right)$.
So one has $\widetilde{\sX}_2= \widetilde{\sX}_1/G\left( \left.{\widetilde{\sX}}_1~\right|\widetilde{\sX}_2\right)$, it follows that the covering $p$ is transitive.
\end{proof}
\begin{definition}\label{top_fin_cov_defn}
	Let $\sX$ be a connected, locally connected, locally compact, Hausdorff space. Let us consider a   category (cf. Definition \ref{category_defn}) $\mathfrak{FinCov}$-$\sX$ such that
	\begin{itemize}
		\item $\mathfrak{FinCov}$-$\sX$-objects of   are transitive finite-fold coverings $\widetilde\sX \to \sX$ of $\sX$ where the space $\widetilde\sX$ is connected,
		\item a $\mathfrak{FinCov}$-$\sX$-morphism form $p_1:\widetilde\sX_1 \to \sX$ to $p_2:\widetilde\sX_2 \to \sX$ is a surjective continuous map $p_{12}: \widetilde\sX_1 \to \widetilde\sX_2$ such that $p_{12} \circ p_2 = p_1$.
	\end{itemize}
	We say that $\mathfrak{FinCov}$-$\sX$ is the \textit{category of finite coverings} of $\sX$. 	From the Theorem \ref{top_covp_cat_thm} it follows that $p_{12}$ is a transitive covering. Sometimes we write $\widetilde \sX$ instead of a covering $\widetilde \sX\to\sX$ to designate an object of $\mathfrak{FinCov}$-$\sX$.
\end{definition}

\section{Miscellany}

\paragraph{}
Let $X$ be a $C^*$-Hilbert $A$-module. We also write $X_A$ instead $X$ and
   $\left\langle{x},{y}\right\rangle_{X_A}$ instead $\left\langle{x},{y}\right\rangle_{A}$, the meaning depends on context. 
If $\ell^2\left( A\right)$ is the {standard Hilbert $A$-module (cf. Definition \ref{standard_h_m_defn})}
There is the natural $*$-isomorphism
\be\label{triv_eqn}
\K\left( \ell^2\left(  A\right)\right) \cong A \otimes \K.
\ee
Any element $a \in \K\left( \ell^2\left(  A\right) \right)$ corresponds to an infinite matrix
\be\label{inf_m_eqn}
\begin{pmatrix}
	a_{1,1}& \ldots & a_{j,1} &  \dots \\
	\vdots& \ddots & \vdots& \ldots\\
	a_{j,1}& \ldots &a_{j,j} &  \ldots \\
	\vdots& \vdots & \vdots & \ddots \\
\end{pmatrix}\in \K\left( \ell^2\left( A\right) \right).
\ee
The free finitely generated $A$-module $A^n$ is also $C^*$-Hilbert $A$ module with the product 
\be\label{unst_hilb_eqn}
\left\langle\left\{a_j\right\}_{j=1,...,n}, \left\{b_j\right\}_{j=1,...,n}\right\rangle_{A^n}=\sum_{j =1}^n a^*_jb_j.
\ee
The algebra of compact operators is given by
\be\label{comp_fin}
\K\left(A^n \right) = \mathbb{M}_n\left(A \right) 
\ee

\begin{remark}
	Similarly to \eqref{inf_m_eqn} and  \eqref{comp_fin} any compact operator on the separable Hilbert space can be represented by an infinite or a finite matrix
	\be\label{comp_matr_eqn}
	\begin{split}
		\begin{pmatrix}
			x_{1,1}& \ldots & x_{j,1} &  \dots \\
			\vdots& \ddots & \vdots& \ldots\\
			x_{j,1}& \ldots &x_{j,j} &  \ldots \\
			\vdots& \vdots & \vdots & \ddots \\
		\end{pmatrix}\in \K\left( L^2\left( \N \right)\right),\\
		\text{or}\\
		\begin{pmatrix}
			x_{1,1}& \ldots & x_{n,1}  \\
			\vdots& \ddots & \vdots\\
			x_{n,1}& \ldots &x_{n,n} \\
		\end{pmatrix}\in \K\left( L^2\left( \{1,...,n\} \right)\right) = \mathbb{M}_n\left(\C \right) 
	\end{split}
	\ee
\end{remark}

\begin{definition}\label{elementary_defn}
	Let us consider an algebra of finite or infinite matrices given by  \eqref{comp_matr_eqn}. For every $j, k \in \N$ a matrix given by \eqref{comp_matr_eqn} is said to be $j,k$-\textit{elementary} if 
	$$
	x_{pr}= \delta_{pj}\delta_{rk}.
	$$
	Denote by $\mathfrak{e}_{jk}$ the $j,k$-{elementary} matrix. We say that $x= \left( x_1,...,x_n\right) \in \C^n$ or  $x=  \left( x_1,x_2,...\right)\in \ell^2\left(\N \right)$ is $j$-\textit{elementary} if  $x_k = \delta_{jk}$. 	Denote by $\mathfrak{e}_{j}$ the $j$-{elementary} vector.
\end{definition}

\begin{definition}\label{def:rep_Hilbert_module_uni}
Let $A$ be a *-algebra and let $D$ be a $C^*$-algebra. Let $\mathcal{E}$ be a Hilbert $D$-module, and let $A^\sim$ is given by \eqref{unital_notation_eqn}, A \emph{representation} of $A $ on $\mathcal{E} $ is a {representation} of $A^\sim$ on $\mathcal{E}$ in the sense of the Definition \ref{def:rep_Hilbert_module}.
\end{definition}

\begin{empt}
If $A$ is a $C^*$-algebra then one has
	\be\label{four_decompositon_eqn}
	\forall a \in A\quad \exists a_1, a_2, a_3, a_4 \in A_+\quad a=a_1 - a_2 + ia_3 - ia_4.
\ee
\end{empt}

\begin{lemma}\label{pedersen_eps_lem}
Let $\eps > 0$, and let 	 $f_\eps: \R \to \R$ be a continuous function given by 
	\begin{equation}\label{f_eps_eqn}
		f_\eps\left( x\right)  =\left\{
		\begin{array}{c l}
			0 &x \le \eps \\
			x - \eps & x > \eps
		\end{array}\right.
	\end{equation}
If $A$ is a $C^*$-algebra 
then one has
\bea\label{k0_ped_e_eqn}
	K\left( A \right)_0 = \left\{f_\eps \left(a\right) \left|~a \in A_+, \quad \eps > 0 \right.\right\},~~~\\
	\label{kp_ped_e_eqn}
		K\left( A \right)_+ = \left\{a \in A_+ \left|a \le \sum_{j = 1}^n f_{\eps_j}\left(  a_j\right) \quad a_j \in  	K\left( A \right)_0\quad \eps_j > 0\quad j=1,...,n\right.\right\}~~
\eea
where both $	K\left( A \right)_0$ and 	$K\left( A \right)_+$ are given by equations \eqref{pedersen_k0_eqn} and \eqref{pedersen_k_plus_eqn} respectively
\end{lemma}
\begin{proof}
	If  $a\in A_+$ and $\eps > 0$ then $f_\eps \left(a\right)= \phi_\eps\left(a \right)$ where  $\phi_\eps \in K(]0, \infty [$ is given by 
		\bean
		\phi_\eps\left( x\right)  =\left\{
		\begin{array}{c l}
			0 &x \le \eps \\
			x - \eps &  \eps \le x \le \left\| a\right\|\\
		2	\left\| a\right\| - \eps - x & \left\| a\right\| \le x \le 2\left\| a\right\|-\eps\\
	0 & x \ge 2 \left\| a\right\| -\eps\\
		\end{array}\right.
	\eean
It follows that $f_\eps\left(a\right)	\in K\left( A \right)_0$. Conversely  if $a \in	K\left( A \right)_0$ then from \eqref{pedersen_k0_eqn} it turns out that there is $b' \in A_+$ and $\varphi \in   K(]0, \infty [$ such that $a = \varphi\left(b'\right)$. If  $\supp \varphi \subset \left[\eps, c\right]$ and  $\psi \in C_c\left( \R\right)_+$ is given by
 		\bean
 \psi\left( x\right)  =\left\{
 \begin{array}{c l}
 	0 &x \le 0 \\
 x &0 \le x \le \eps \\
 \varphi\left(x\right) + \eps &\eps \le x \le c \\
 \eps + c - x&  c \le x \le c+ \eps\\
 0 & x \ge c + \eps
 \end{array}\right.
 \eean
 then $ \varphi = f_\eps  \circ\psi$. It follows that $a = f_\eps\left(b \right)$ where $b \bydef \psi\left(b'\right)$. So the equation \eqref{k0_ped_e_eqn} is proven. The equation \eqref{kp_ped_e_eqn} is a direct consequence of \eqref{k0_ped_e_eqn} and \eqref{pedersen_k_plus_eqn} ones.
 
\end{proof}

\begin{definition}\label{natural_defn}
	If $ \A$ is  category then a related to $ \A$ construction is said to be \textit{natural} or \textit{functorial} if it is invariant with respect to isomorphisms of $ \A$.
\end{definition}

\chapter{Noncommutative finite-fold coverings}\label{cov_fin_chap}
   \section{Basic definitions}\label{cov_fin_bas_sec}
  \subsection{Coverings of $C^*$-algebras}
   \paragraph{}

 Here the noncommutative generalization of the Theorem \ref{pavlov_troisky_thm} is being discussed.  
   \begin{definition}\label{fin_quasi_defn}
	Let both  $A$ and  $\widetilde{A}$ be  connected $C^*$-algebras (cf. Definition \ref{connected_c_a_defn}), and let $\pi: A \hookto \widetilde{A}$ be an injective $*$-homomorphism of 
	$C^*$-algebras. Let $G$ be a finite  group of *-automorphisms of $\widetilde{A}$ such that 	$\pi\left(A\right) = \widetilde{A}^G\stackrel{\text{def}}{=}\left\{
	\left.a\in \widetilde{A}~\right|~ a = g a;~ \forall g \in G\right\}$.	We say that the triple $\left(A, \widetilde{A}, G \right)$ and/or the quadruple $\left(A, \widetilde{A}, G, \pi \right)$ and/or $*$-homomorphism $\pi: A \hookto \widetilde{A}$   is a \textit{noncommutative finite-fold  quasi-covering}. We write
	\be\label{fin_cov_gr_eqn}
	G\left(\left.\widetilde{A}~\right| {A} \right) \stackrel{\text{def}}{=}  	G.
	\ee
\end{definition}

\begin{lemma}\label{fin_pedersen_lem}
If a quadruple $\left(A, \widetilde{A}, G, \pi \right)$   is a {noncommutative finite-fold  quasi-covering} (cf. Definition \ref{fin_quasi_defn})
then one has
\be\label{fin_pedersen_eqn}
\forall \widetilde a \in K\left(\widetilde{A}\right) \quad \sum_{	g \in G}g \widetilde a \in \pi\left( K\left( A\right)\right)  
\ee
where both  $K\left(\widetilde{A}\right)$ and  $K\left( A\right)$ are Pedersen's ideals of $\widetilde{{A}}$ and $A$ respectively (cf. Definition \ref{pedersen_ideal_defn}). Moreover the map
\be\label{fin_cov_gr_y_eqn}
\begin{split}
K\left(\widetilde{A}\right)\to K\left({A}\right),\\
\widetilde a \mapsto \quad \pi^{-1}\left( \sum_{	g \in G}g \widetilde a \right) 
\end{split}
\ee
is a surjective homomorphism of $A$-$A$ bimodules.
\end{lemma}
\begin{proof}
	From  $a' \bydef \sum_{	g \in G}g \widetilde a'\in \widetilde A^G$ it follows that $a' \in \pi\left(A\right)$.
If $\widetilde a \in K\left(\widetilde{A}\right)_0$ (cf. equation \eqref{pedersen_k0_eqn})  then from the Lemma \ref{pedersen_eps_lem} it follows that there is a positive $\widetilde b \in \widetilde A_+$ and $\eps > 0$  such that $\widetilde a = f_\eps\left(\widetilde b\right)$ where $\eps > 0$ and $f_\eps$ is given by \eqref{f_eps_eqn}. If $b \bydef \pi^{-1}\left( \sum_{	g \in G}g \widetilde b\right) $ then one has: 
\begin{itemize}
	\item $f_\eps \left(b \right) \in  \pi\left( K\left( A\right)\right)_0$,
	\item $ f_\eps  \left(b \right)  \ge a \bydef  \pi^{-1}\left(\sum_{	g \in G}g f_\eps\left(  \widetilde a\right)  \right) $.
\end{itemize}
Taking into account that 
$$
\forall a' \in A_+ \quad a'' \in K\left(A\right)_+ \quad a' \le a'' \quad \Rightarrow \quad a' \in K\left(A\right)_+
$$
one has $a \in K\left( A\right)_+$. Using the equation \eqref{four_decompositon_eqn}  one can prove \eqref{fin_pedersen_eqn} for any $ \widetilde a \in K\left(\widetilde{A}\right)$.

If $a \in K\left( A\right)_0$ then there is a positive $ b \in  A_+$ and $\eps > 0$  such that $ a = f_\eps\left( b\right)$ where $\eps > 0$ and $f_\eps$ is given by \eqref{f_eps_eqn}.  So $\pi\left(a \right)= f_\eps\left(\pi\left(b \right)  \right) \in K\left(\widetilde A \right)$. From $\pi\left(b \right) =\sum_{g\in G}g \pi\left( \frac{1}{\left|G \right|}b\right)$ and $\pi\left(a \right) =\sum_{g\in G} g\pi\left( \frac{1}{\left|G \right|}a\right)$ it follows that the map \eqref{fin_cov_gr_y_eqn} is surjective.
\end{proof}
\begin{lemma}\label{quasi_hareditary_lem}
	
	If $\left(A, \widetilde{A}, G, \pi \right)$ is a noncommutative finite-fold  quasi-covering then $\widetilde{A}$ is a generated by $\pi\left( A\right)$ hereditary subalgebra of $ \widetilde{A}$ (cf. Definition \ref{hered_gen_defn}.)
\end{lemma}
\begin{proof}
	If $\widetilde a \in \widetilde A_+$ is a positive element then from
	$$ 
	\widetilde a < \sum_{	g \in G} g\widetilde a \in \pi\left(A \right) 
	$$
	and the Definition \ref{hered_defn} it turns out that $\widetilde a$ lies in a generated by $\pi\left( A\right)$ hereditary subalgebra of $ \widetilde{A}$.
\end{proof}


\begin{lemma}\label{quasi_approximate_unit_lem}
	If $\left(A, \widetilde{A}, G, \pi \right)$ is a noncommutative finite-fold  quasi-covering then there is an  approximate unit (cf. Definition \ref{approximate_unit_defn} ) of  $\widetilde{A}$ which is contained in $\pi\left( A\right) $.
\end{lemma}
\begin{proof}
	From the Theorem \ref{approximate_unit_thm}  one has an approximate unit  $\left\{\widetilde{u}_\la \right\}_{\la \in \La}$ of  $\widetilde{A}$. For any $g \in G$ and $\widetilde{x}\in  \widetilde{A}$ one has 
\bean
	\lim_\la \left\|\widetilde{x}- \widetilde{x}\left( g\widetilde{u}_\la\right)  \right\| = 	\lim_\la \left\|g\left(g^{-1} \widetilde{x}- \left(g^{-1} \widetilde{x} \right) \widetilde{u}_\la\right)   \right\|=\\=\left\|g\lim_\la\left(g^{-1} \widetilde{x}- \left(g^{-1} \widetilde{x} \right) \widetilde{u}_\la\right)   \right\| =0,
\eean
	i.e.  $\left\{g\widetilde{u}_\la \right\}_{\la \in \La}$ is an approximate unit of  $\widetilde{A}$. 
	If ${v}_\la \bydef \frac{1}{\left|G\right|}\sum_{g\in G} g\widetilde{u}_\la \in \widetilde A^G= \pi\left( A\right) $, then 
	\bean
	\lim_\la \left\| \widetilde{x}- \widetilde{x}v_\la   \right\| = 	\lim_\la \left\|\frac{1}{\left|G\right|}\sum_{g\in G} \left( \widetilde{x}- \widetilde{x}   \left( g\widetilde{u}_\la\right)\right)    \right\|\le \lim_\la \frac{1}{\left|G\right|}\sum_{g\in G}\left\| \widetilde{x}- \widetilde{x}  \left( g\widetilde{u}_\la\right)   \right\| = 0.
	\eean
	Similarly one can prove that $\lim_\la \left\| \widetilde{x}- v_\la \widetilde{x}  \right\|=0$, i.e.   $\left\{{v}_\la \right\}_{\la \in \La}$ is an approximate unit of $\widetilde A$. 
	Otherwise $\left\{{v}_\la \right\}_{\la \in \La}\subset\pi\left(A\right)$.
  \end{proof}
 
\begin{remark}\label{implicit_inclusion_rem}
Since a $*$-homomorphism $\pi: A\hookto \widetilde{A}$ it can be implicitly replaced by an inclusion
\be\label{implicit_inclusion_eqn}
A\subset \widetilde{A}
\ee 
\end{remark}
\begin{lemma}\label{fin_hilbert_mod_lem}
If $\left(A, \widetilde{A}, G, \pi \right)$  is a {noncommutative finite-fold  quasi-covering} then $\widetilde{A}$ is a $C^*$-Hilbert $A$-module with respect to the following $A$-valued product
	\begin{equation}\label{finite_hilb_mod_prod_eqn}
	\left\langle \widetilde a, \widetilde b \right\rangle_{{A}} =\pi^{-1}\left( \sum_{g \in G} g\left( \widetilde a^* \widetilde b\right)\right).  
	\end{equation}
	\end{lemma}
\begin{proof}
	Direct calculation shows that the product \eqref{finite_hilb_mod_prod_eqn} satisfies to conditions (a)-(d) of the Definition \ref{hilbert_module_defn}. Let us prove that $\widetilde{A}$ is closed with respect to the norm \eqref{hilbert_module_norm_eqn}.
	The norm $\left\| \cdot\right\|_H$  of $A$-Hilbert pre-module $\widetilde{A}_A$  is given by
\be\label{c_hilb_norm_eqn}
	\left\| \widetilde{a}\right\|_H\bydef  \sqrt{	\left\langle \widetilde a,\widetilde a \right\rangle_{{A}}} = \sqrt{\left\| \sum_{g \in G} g\left( \widetilde{a}^*\widetilde{a}\right) \right\|_C }
\ee
	where $\left\| \cdot\right\|_C$ is the $C^*$-norm. From the above equation and taking into account $\left\| \widetilde{a}\right\|_C = \sqrt{\left\| \widetilde{a}^*\widetilde{a}\right\|_C}$ it turns out
	\be\label{hilb_c_norm_eqn}
	\begin{split}
	\left\| \widetilde{a}\right\|_C \le \left\| \widetilde{a}\right\|_H,\\
\left\| \widetilde{a}\right\|_H \le \sqrt{\left|G\right|} 	\left\| \widetilde{a}\right\|_C,
	\end{split}
	\ee
	i.e. the norms	$\left\| \cdot\right\|_C$ and 	$\left\| \cdot\right\|_H$ yield equal topologies.
	The algebra $\widetilde{A}$ is closed with respect to $C^*$-norm, hence a pre-$C^*$-Hilbert $A$-module $\widetilde{A}_A$ is closed with respect to the pre-$C^*$-Hilbert $A$-module norm $\left\| \cdot\right\|_H$.
\end{proof}

\begin{definition}\label{hilbert_product_defn}
	If $\left(A, \widetilde{A}, G, \pi \right)$ is a noncommutative finite-fold  quasi-covering then the   given by \eqref{finite_hilb_mod_prod_eqn} product is said to be the \textit{Hilbert product associated} with the noncommutative finite-fold  quasi-covering $\left(A, \widetilde{A}, G, \pi \right)$. The given by the Lemma \ref{fin_hilbert_mod_lem} $C^*$-Hilbert $A$-module $\mathscr L^2\left( \widetilde{A}\right)_A$  is said to be the \textit{$C^*$-Hilbert   associated with the noncommutative finite-fold  quasi-covering} $\left(A, \widetilde{A}, G, \pi \right)$. The $C^*$-Hilbert module  $\mathscr L^2\left( \widetilde{A}\right)$  is $\widetilde A$-$A$-bimodule, so it will be denoted by $_{\widetilde{A}}\mathscr L^2\left( \widetilde{A}\right)_A$ and there is the natural isomorphism
		$$
	_{\widetilde{A}}\mathscr L^2\left( \widetilde{A}\right)_A \cong 	\widetilde A
		$$
		of $\widetilde A$-$A$-bimodules.
\end{definition}
   \begin{definition}\label{fin_pre_defn}
   		Let both  $A$ and  $\widetilde{A}$ are  connected $C^*$-algebras (cf. Definition \ref{connected_c_a_defn}), and let $\pi: A \hookto \widetilde{A}$ be an injective  $*$-homomorphism of 
   	$C^*$-algebras such that following conditions hold:
	\begin{enumerate}
		\item[(a)] If $\Aut\left(\widetilde{A} \right)$ is a group of *-automorphisms of $\widetilde{A}$ then the group  
				\be\nonumber
		G \bydef \left\{ \left.g \in \Aut\left(\widetilde{A} \right)~\right|\forall a \in \pi \left( A\right) \quad ga = a\right\}
		\ee
		is finite.
		\item[(b)] 	$\pi\left( A \right) = \widetilde{A}^G\stackrel{\text{def}}{=}\left\{\left.a\in \widetilde{A}~~\right|\forall g \in G\quad  a = g a\right\}$.
	\end{enumerate}
	We say that the triple $\left(A, \widetilde{A}, G \right)$ and/or the quadruple $\left(A, \widetilde{A}, G, \pi \right)$ and/or $*$-homomorphism $\pi: A \hookto \widetilde{A}$   is a \textit{noncommutative finite-fold  pre-covering}. We write $G\left(\left.\widetilde A~\right|A \right)\bydef G$.
\end{definition}
\begin{remark}
Any noncommutative finite-fold  pre-covering is a noncommutative finite-fold quasi-covering.
\end{remark}

 \begin{remark}\label{unital_rem}
	From the Lemma \ref{quasi_approximate_unit_lem} it follows that the homomorphism $\pi$ in Definitions \ref{fin_quasi_defn} and \ref{fin_pre_defn} is unital in the sense of the Definition \ref{principal_non_defn}.
\end{remark}
\begin{remark}
Indeed conditions (a), (b) of the Definition \ref{fin_pre_defn} state an isomorphism between ordered sets  of groups  and $C^*$-algebras  and it is an analog the fundamental
theorem of Galois theory. (cf. \cite{kurosh:lga} Chapter Six, \S~6 for details).
\end{remark}
\begin{definition}\label{proper_subgroup_fin_defn}
Let $\left(A, \widetilde{A}, G\left(\left.\widetilde A~\right|A \right), \pi \right)$ be  a {noncommutative finite-fold  pre-covering} (cf. Definition \ref{fin_pre_defn}). A normal subgroup $H\subset G\left(\left.\widetilde A~\right|A \right)$ is said to be $\left(A, \widetilde{A}, G\left(\left.\widetilde A~\right|A \right), \pi \right)$-\textit{proper} if a quadruple
\be\label{proper_subgroup_fin_eqn}
\left(A, \widetilde{A}^H, G\left(\left.\widetilde A~\right|A \right)/H, \pi^H \right)
\ee
is a  {noncommutative finite-fold  pre-covering}. The $*$-homomorphism $\pi^H:A\hookto \widetilde{A}^H $ comes from $\pi: A \hookto \widetilde{A}$ and an inclusion  $\pi\left( A\right) \subset \widetilde{A}^H$.
\end{definition}
\begin{lemma}\label{proper_subgroup_fin_lem}
If $\left(A, \widetilde{A}, G\left(\left.\widetilde A~\right|A \right), \pi \right)$ is  a {noncommutative finite-fold  pre-covering} and  $H\subset G\left(\left.\widetilde A~\right|A \right)$ is a subgroup  then a quadruple 
	\be\label{proper_subgroup_finh_eqn}
\left(\widetilde{A}^H, \widetilde{A}, H, \left.\Id_{\widetilde{A}}\right|_{\widetilde{A}^H} \right).
\ee
a {noncommutative finite-fold  pre-covering} (cf. Definition \ref{fin_pre_defn}).
\end{lemma}
\begin{proof}	
	From the Definition \ref{fin_pre_defn} it follows that
	$$
 G\left(\left.\widetilde A~\right|A \right) = 	 \left\{ \left.g \in \Aut\left(\widetilde{A} \right)~\right|~ ga = a;~~\forall a \in \pi \left( A\right) \right\}.
	$$
	So one has
\bean
\forall g \in \Aut\left(\widetilde{A}\right) \quad \forall a^H \in \widetilde A^H\quad   g\widetilde a^H = \widetilde a^H \quad \Rightarrow \quad \forall a \in \pi\left(A \right) \quad g a = a,\\
 \left\{ \left.g \in \Aut\left(\widetilde{A} \right)~\right|~ g\widetilde a^H = \widetilde a^H;~~\forall \widetilde a^H \in \widetilde{A}^H\right\}\subset  G\left(\left.\widetilde A~\right|A \right).	
\eean
	It turns out that
	\bean
	 \left\{ \left.g \in \Aut\left(\widetilde{A} \right)~\right|~ g\widetilde a^H = \widetilde a^H;~~\forall \widetilde a^H \in \widetilde{A}^H\right\}= \\= \left\{ \left.g \in G~\right|~ g\widetilde a^H = \widetilde a^H;~~\forall \widetilde a^H \in \widetilde{A}^H\right\}=H.
	\eean
\end{proof}

\begin{empt}
Let $A \subset \widetilde A$ be an inclusion of $C^*$-algebras such that $\widetilde A$ is a finitely generated left  $A$-module, i.e.
$$
\exists \widetilde a_1, ... \widetilde a_n \quad \widetilde A = A\widetilde a_1 + ...  + A \widetilde a_n.
$$
From $\widetilde A = \widetilde A^*$ it turns out that
$$
\widetilde A = \widetilde a^*_1 A + ...  + \widetilde a^*_n A,
$$
hence one has
\be\label{left_right_eqv_equ}
\begin{split}
\widetilde A \text{ is a finitely generated left  }A \text{-module} \quad \Leftrightarrow \\
\Leftrightarrow\quad \widetilde A \text{ is a finitely generated right }A\text{-module}.
\end{split}
\ee
\end{empt}

   \begin{definition}\label{fin_unital_defn}
   	Let $\left(A, \widetilde{A}, G, \pi \right)$ be a  noncommutative finite-fold  pre-covering (cf. Definition \ref{fin_pre_defn}) (resp. quasi-covering  (cf. Definition \ref{fin_quasi_defn})). Suppose .that both $A$ and  $\widetilde{A}$ are unital $C^*$-algebras. We say that $\left(A, \widetilde{A}, G, \pi \right)$ is an \textit{unital noncommutative finite-fold  covering} (resp. \textit{unital noncommutative finite-fold  quasi-covering}) if $\widetilde{A}$ is a finitely generated $C^*$-Hilbert left  and/or right $A$-module with respect to product given by \eqref{finite_hilb_mod_prod_eqn}.
   \end{definition}
   \begin{remark}
   	Above definition is motivated by the Theorem \ref{pavlov_troisky_thm}.
   \end{remark}
\begin{remark}
	From the equation \eqref{left_right_eqv_equ} it follows that the usage of left  $A$-modules in the Definition \ref{fin_unital_defn} is equivalent to the usage of right ones.
\end{remark}
\begin{remark}\label{fin_cov_kasp_rem}
	From the Lemma \ref{fin_hilbert_mod_lem} and the Corollary \ref{fin_hpro_cor} it turns out that $\widetilde{A}$ is a projective $A$-module.
\end{remark}
   \begin{lem}\label{fin_composition_lem} 
	Let $\left(A, \widetilde{A}, G, \pi \right)$ be an unital noncommutative finite-fold  covering (cf. Definition \ref{fin_unital_defn}). If  $H \subset G$ is a subgroup then 
	\be
		\label{fin_composition2_eqn}
	\left(\widetilde{A}^H, \widetilde{A}, H, \left.\Id_{\widetilde{A}}\right|_{\widetilde{A}^H} \right).
	\ee
is an unital finite-fold noncommutative covering.	Moreover if $H$ a $\left(A, \widetilde{A}, G, \pi \right)$-proper group (cf. Definition \ref{proper_subgroup_fin_defn})  then there is a following unital finite-fold noncommutative covering
	\be\label{fin_composition_eqn}
	\left(A, \widetilde{A}^H, G/H, \rho \right).
	\ee
\end{lem}
\begin{proof}
From  the Lemma \ref{proper_subgroup_fin_lem} it follows that the  quadruple \eqref{fin_composition2_eqn} is a noncommutative finite-fold pre-covering (cf. Definition \ref{fin_pre_defn}). If $\widetilde{A}$ is an $A$-module generated by  $\left\{\widetilde a_1 ,..., \widetilde a_n\right\}\subset \widetilde{A}$ then
$\widetilde{A} = \widetilde a_1 A + ... + \widetilde a_n A$. From $A\subset \widetilde{A}^H$ it  turns out that $\widetilde{A} = \widetilde a_1 \widetilde{A}^H + ... + \widetilde a_n \widetilde{A}^H$, i.e. $\widetilde{A}$ is a finitely generated $\widetilde{A}^H$-module.

If $H$ is a $\left(A, \widetilde{A}, G, \pi \right)$-proper group then from the Definition \ref{proper_subgroup_fin_defn} it follows that the quadruple \ref{fin_composition_eqn}  is a noncommutative finite-fold pre-covering (cf. Definition \ref{fin_pre_defn}).
If $\widetilde a^H \in \widetilde{A}^H$ then $\widetilde a^H = \frac{1}{\left|H \right| }\sum_{g \in H} \widetilde a^H$. Otherwise 
$$
\widetilde a^H = \widetilde a_1 a_n + ...+ \widetilde a_n a_n =  \frac{1}{\left|H \right| }\sum_{g \in H}\left( g\widetilde a_1 a_n + ...+ g\widetilde a_n a_n\right)= \widetilde a^H_1 a_1 + ...+ \widetilde a^H_n a_n
$$
where $\widetilde a^H_j \bydef \frac{1}{\left|H \right| }\sum_{g \in H}g \widetilde a_j\in \widetilde{A}^H$ for any $j = 1,..., n$. So $\widetilde{A}^H$ is an $A$-module  generated by a finite set $\left\{\widetilde a^H_1 ,..., \widetilde a^H_n\right\}\subset \widetilde{A}^H$. 
\end{proof}

\begin{definition}\label{fin_unitization_defn}
   	Let $\left(A, \widetilde{A}, G, \pi \right)$ be a noncommutative finite-fold  pre-covering (resp. quasi-covering) of $C^*$-algebras $A$ and $\widetilde{A}$ such  that following conditions hold:
   	\begin{enumerate}
   		\item[(a)] 
  		there are unitizations $A \hookto B$  and $\widetilde{A} \hookto \widetilde{B}$ (cf. Definition \ref{unitization_defn});
   		\item[(b)] there is a 
   	  unital  noncommutative finite-fold quasi-covering	$\left(B ,\widetilde{B}, G, \widetilde{\pi} \right)$ (cf. Definition \ref{fin_unital_defn}) such that
   	  $\pi = \widetilde{\pi}|_A$ (or, equivalently $\pi\left( A \right)= \widetilde{A}\cap \widetilde \pi\left( B\right)$) and
   	   the action $G \times\widetilde{A} \to \widetilde{A}$ comes from  the $G \times\widetilde{B} \to \widetilde{B}$ one.
   	\end{enumerate}
We say that the triple $\left(A, \widetilde{A}, G \right)$ and/or the quadruple $\left(A, \widetilde{A}, G, \pi \right)$ and/or $*$-homomorphism $\pi: A \hookto \widetilde{A}$ is a
   		 \textit{noncommutative finite-fold covering with unitization}, (resp. \textit{noncommutative finite-fold quasi-covering with unitization}). 
   \end{definition}
   \begin{remark}\label{unital_unitization_covering_rem}
   	If  $\left(A, \widetilde{A}, G, \pi \right)$ is an unital  noncommutative finite-fold covering (cf. Definition \ref{fin_unital_defn}) then both $A$ and $\widetilde{A}$ are unital, so both identical inclusions   $A \hookto A$  and $\widetilde{A} \hookto \widetilde{A}$ are unitizations (cf. Definition \ref{unitization_defn}). Using this fact one can prove that $\left(A, \widetilde{A}, G, \pi \right)$ is a  noncommutative finite-fold covering with unitization. 
   	So any unital noncommutative finite-fold covering is a finite-fold covering with unitization with unitization (cf. Definition \ref{fin_unitization_defn}) 
   	   \end{remark}
   
\begin{lemma}\label{fin_composition_unitization_lem}
		Let $\left(A, \widetilde{A}, G, \pi \right)$ be a noncommutative finite-fold  covering with unitization (cf. Definition \ref{fin_unitization_defn}). If  $H \subset G$ is a subgroup then 
	\be
	\label{fin_composition_uni_2_eqn}
	\left(\widetilde{A}^H, \widetilde{A}, H, \left.\Id_{\widetilde{A}}\right|_{\widetilde{A}^H} \right).
	\ee
	is a finite-fold noncommutative covering with unitization (cf. Definition \ref{fin_unitization_defn}).	Moreover if $H$ a $\left(A, \widetilde{A}, G, \pi \right)$-proper group (cf. Definition \ref{proper_subgroup_fin_defn})  then there is a following  finite-fold noncommutative covering with unitization (cf. Definition \ref{fin_unitization_defn}).
	\be\label{fin_composition_uni_1_eqn}
	\left(A, \widetilde{A}^H, G/H, \rho \right).
	\ee

\end{lemma}
\begin{proof}
	 From the Lemma \ref{proper_subgroup_fin_lem} it follows that the quadruple  \eqref{fin_composition_uni_2_eqn} is a noncommutative  finite-fold pre-covering.
Let $\left(B, \widetilde{B}, G, \widetilde\pi \right)$ be a required by the Definition \ref{fin_unitization_defn} unital noncommutative finite-fold covering, i.e. one has
   	\begin{enumerate}
	\item[(a)] 
	There are unitizations $A \hookto B$  and $\widetilde{A} \hookto \widetilde{B}$ (cf. Definition \ref{unitization_defn}).
	\item[(b)] There is a 
	unital  noncommutative finite-fold quasi-covering	$\left(B ,\widetilde{B}, G, \widetilde{\pi} \right)$ (cf. Definition \ref{fin_unital_defn}) such that
	$\pi = \widetilde{\pi}|_A$ (or, equivalently $\pi\left( A \right)= \widetilde{A}\cap \widetilde \pi\left( B\right)$) and
	the action $G \times\widetilde{A} \to \widetilde{A}$ comes from  the $G \times\widetilde{B} \to \widetilde{B}$ one.
\end{enumerate}
There is a natural inclusion $ \widetilde{A}^H \subset  \widetilde{B}^H$, let us prove that $ \widetilde{A}^H$ is an essential ideal of $\widetilde{B}^H$ (cf. Definition \ref{essential_defn}). Let  $\widetilde b \in \widetilde B^H\setminus \{0\}$ be such that $\widetilde b\widetilde A^H= \{0\}$. Since $\widetilde b\neq 0$ from the Lemma \ref{essential_lem} there is $\widetilde a \in \widetilde A$ such that $\widetilde b \widetilde a \neq 0$. From \eqref{four_decompositon_eqn} it follows that one can assume that $\widetilde a$ is positive. From the equation $\widetilde b \widetilde a = \widetilde b \sqrt{\widetilde a} \sqrt{\widetilde a} \neq 0$ it follows that $\widetilde b \sqrt{\widetilde a}\neq 0$ and $\widetilde b \sqrt{\widetilde a}\left( \widetilde b \sqrt{\widetilde a}\right)^* =  \widetilde b \sqrt{\widetilde a} \sqrt{\widetilde a}\widetilde b^*= \widetilde b\widetilde a\widetilde b^*$ is a nonzero positive element. The following sum
\bean
\sum_{	g \in H} g\left( \widetilde b\widetilde a\widetilde b^*\right) = \widetilde b \left( \sum_{	g \in H} \widetilde a\right)\widetilde  b^*
\eean
of positive nonzero elements is not equals to zero, so one has
$$
 \widetilde b \left( \sum_{	g \in H} \widetilde a\right) \widetilde b^*\neq 0 \quad\Rightarrow\quad\widetilde b \left( \sum_{	g \in H} \widetilde a\right)\neq 0.
$$
On the other hand $\sum_{	g \in H} \widetilde a\in \widetilde A^H$, so $\widetilde b\widetilde A^H\neq \{0\}$. One has a contradiction it follows that
$$
\forall \widetilde b\in \widetilde B^H\quad \widetilde b\widetilde A^H\neq \{0\} \quad \Rightarrow\quad \widetilde b \neq 0.
$$
From the above equation and the Lemma \ref{essential_lem} it follows that $\widetilde A^H$ is an essential ideal of $\widetilde B^H$.
Since $\widetilde B^H$ is unital, an inclusion $\widetilde  A^H\hookto\widetilde B^H$ is an unitization. From the Lemma \ref{fin_composition_lem}  it follows that the quadruple
	\bean
\left(\widetilde{B}^H, \widetilde{B}, H, \left.\Id_{\widetilde{B}}\right|_{\widetilde{B}^H} \right).
\eean
is an unital noncommutative finite-fold covering. From the Definition \ref{fin_unitization_defn} it follows that the quadruple \eqref{fin_composition_uni_2_eqn} is a noncommutative  finite-fold coverings with unitization. 

If $H$ is a $\left(A, \widetilde{A}, G, \pi \right)$-proper group  then from the Definition \ref{proper_subgroup_fin_defn} it follows that the quadruple \eqref{fin_composition_uni_1_eqn} is a noncommutative finite-fold pre-covering (cf. Definition \ref{fin_pre_defn}). From the Lemma \ref{fin_composition_lem}  it follows that the quadruple
\bean
\left(B, \widetilde{B}^H, G/H,\widetilde \rho \right)
\eean
is an unital noncommutative finite-fold covering. From the Definition \ref{fin_unitization_defn} it follows that the quadruple \eqref{fin_composition_uni_1_eqn} is a noncommutative  finite-fold coverings with unitization. 
\end{proof}
\begin{remark}\label{fin_comp_rem}
 If we consider the situation of the Definition \ref{fin_unitization_defn} then since both  inclusions $A \hookto B$, $\widetilde{A} \hookto \widetilde{B}$ correspond to essential ideals then one has
 $$
 \left\{ \left.g \in \Aut\left(\widetilde{A} \right)~\right|~ ga = a;~~\forall a \in A\right\}\cong \left\{ \left.g \in \Aut\left(\widetilde{B} \right)~\right|~ gb = b;~~\forall b \in B\right\},
 $$
 it turns out that
 \be\label{fin_comp_pre_eqn}
 \left(A ,\widetilde{A}, G, {\pi} \right)\text{ is a pre-covering } \quad \Leftrightarrow\quad \left(B ,\widetilde{B}, G, \widetilde{\pi} \right) \text{ is a pre-covering }.
 \ee
 \end{remark}
\begin{remark}\label{fin_ded_rem}
	The equation   $\pi = \widetilde{\pi}|_A$ (or, equivalently $\pi\left( A \right)= \widetilde{A}\cap \widetilde \pi\left( B\right)$) follows from the actions $G \times\widetilde{A} \to \widetilde{A}$ and $G \times\widetilde{B} \to \widetilde{B}$. Really 
	if $ \widetilde a\in \widetilde{A}\cap \pi\left( B\right)$, or   equivalently  $ \widetilde a \in \widetilde \pi\left( A\right)$ then from $\widetilde a\in \widetilde \pi\left( B\right)$ it follows that	$	\widetilde a = \frac{1}{\left|G \right| }\sum_{	g \in G}\widetilde a$, so one has $\widetilde a \in \widetilde A^G = \pi\left(A\right)$. It turns out that $\pi = \widetilde{\pi}|_A$ (or, equivalently $A = \widetilde{A}\cap\widetilde{\pi} \left( B\right)$).
\end{remark}

\begin{empt}\label{act_mult_empt}
If $\left(A, \widetilde{A}, G, \pi \right)$ be {noncommutative finite-fold  quasi-covering} (cf. Definition \ref{fin_quasi_defn}) then there is the action $G \times M\left(\widetilde A \right) \to M\left(\widetilde A\right)$ such that for any $\widetilde a\in M\left( \widetilde A\right)$  (cf. Definition \ref{strict_topology_defn}) one has
\be\label{gma_act_eqn}
g \widetilde a \bydef \bt\text{-}\lim v_\la g \widetilde a
\ee
where $\{v_\la\}\subset A$ is given by  the Lemma  \ref{quasi_approximate_unit_lem} approximate unit of $\widetilde A$ and the strict topology of $M\left( \widetilde A\right)$ is implied. From \eqref{gma_act_eqn} It follows that
\bean
\forall \widetilde a\in M\left(\widetilde A\right)\quad \forall\widetilde b\in A, \quad \forall g\in G \quad g\left(\widetilde a \widetilde b \right)= \bt\text{-}\lim \pi\left(  v_\la\right)  g \left( \widetilde a\widetilde b \right)= \\= \left( \bt\text{-}\lim \pi\left(  v_\la\right)\left(  g\widetilde a\right) \right)\left(g \widetilde b\right)= \left(g\widetilde a\right)\left(g\widetilde b\right),\\
g\left(\widetilde b \widetilde a \right)=\lim\pi\left(  v_\la\right) g\left(  \widetilde b\widetilde a \right)= \left(g\widetilde b\right)\left(\lim\pi\left(  v_\la\right) g\widetilde a\right)= \left(g\widetilde b\right)\left(g\widetilde a\right).
\eean 

\end{empt}
 \begin{lemma}\label{fin_multilier_lem}
	If $\left(A, \widetilde{A}, G, \pi \right)$ is a noncommutative finite-fold  quasi-covering then one has:
	\begin{enumerate}
		\item [(i)] there is the natural inclusion $ M\left(\pi\right) :M\left({A} \right) \hookto M\left(\widetilde{A} \right)$ of $C^*$-algebras;
		\item[(ii)] if an	action of $G$ on $M\left(\widetilde{A} \right)$ is induced by the action  of $G$ on $\widetilde{A}$ (i.e. it is given by \eqref{gma_act_eqn}) then there is the natural $*$-isomorphism
		\begin{equation}\label{mag_ma_acc_eqn}
		\psi: M\left(\widetilde{A} \right)^G \cong  M\left(A \right);
		\end{equation}
	\item[(iii)] if  $\left(A, \widetilde{A}, G, \pi \right)$ is a pre-covering (cf. Definition \ref{fin_pre_defn}) then $\left( M\left({A} \right),  M\left(\widetilde{A}\right)  G, M\left(\pi \right)  \right)$ is a pre-covering. 
	\end{enumerate}
\end{lemma}
\begin{proof} 
	(i)
	From the Lemma \ref{quasi_approximate_unit_lem}  and the Proposition \ref{mult_str_pos_prop} it turns out that $\pi$ extends to the inclusion
	\be\label{mult_map_eqn}
M\left(\pi\right):M\left({A} \right) \hookto M\left(\widetilde{A} \right).
	\ee\\
	(ii) If  $a \in M\left(A \right)$ then there is a net $\left\{a_\la\right\}\bydef\left\{v_\la a\right\}\subset A$ such that $a =\bt$-$\lim a_\la$ with respect to the strict topology of $M\left({A} \right)$ (cf Definition \ref{strict_topology_defn}).	From \eqref{gma_act_eqn} it follows that
	$$
\forall g \in G\left(\left.\widetilde{A}~\right|A  \right)\quad 	g M\left(\pi \right)\left(a \right)= \bt\text{-}\lim g a_\la = \bt\text{-}\lim  a_\la =\varphi\left(a \right) ,
	$$ 
	hence $M\left(\pi \right)\left(M\left(A \right)  \right) \subset  M\left( \widetilde{A}\right)^G$. So there is a natural injective $*$-homomorphism $\phi : M\left(A \right)\hookto  M\left(\widetilde{A} \right)^G$.	On the other hand 	if $a \in M\left(\widetilde{A} \right)^G$ and $b \in \widetilde{A}^G$ then $ab \in \widetilde{A}$ is such that $g\left( ab\right) =\left( ga\right)\left(gb \right) = ab$ for any $g\in G$, i.e. $ab  \in \widetilde{A}^G$. It follows that $M\left( \widetilde{A}\right)^G \widetilde A^G = \widetilde A^G$ and taking into account $\widetilde A^G= A$ one concludes that $M\left( \widetilde{A}\right)^G A = A$. Similarly we can prove that $AM\left( \widetilde{A}\right)^G  = A$. From both  $M\left( \widetilde{A}\right)^GA  = A$ and  $AM\left( \widetilde{A}\right)^G  = A$ it turns out that any $a \in M\left( \widetilde{A}\right)^G$ yields an element of $M\left(A \right)$ (cf. the Definition \ref{multiplier_el_defn} and/or the Remark \ref{double_centralizer_rem}), i.e. there is  a $*$-homomorphism.
	\be\label{mag_mai_eqn}
	\psi:M\left(\widetilde{A} \right)^G \hookto   M\left(A \right).
	\ee
	Clearly $\psi \circ \phi = \Id_{M\left(A \right) }$ and since $\phi$ is injective we conclude that both $\phi$ and $\psi$ are $*$-isomorphisms. \\
	(iii) Both $A$ and $\widetilde A$ are dense in both $M\left(A\right)$ and $M\left(\widetilde A\right)$ with respect to strict topology, so one has
\bean
	G  = 	 \left\{ \left.g \in \Aut\left(\widetilde{A} \right)~\right|~ g\pi\left( a\right)  = \pi\left( a\right);~~\forall a \in A\right\}=\\=  \left\{ \left.g \in \Aut\left(M\left( \widetilde{A}\right)  \right)~\right|~ gM\left(\pi \right)\left( a\right)  =M\left(\pi \right)\left( a\right);~~\forall a \in M\left( A\right) \right\}.
\eean
	 
\end{proof}
\begin{corollary}\label{ind_mult_inv_cor}
	If $\left(A, \widetilde{A}, G, \pi \right)$ is a noncommutative finite-fold  quasi-covering then there is a natural  noncommutative finite-fold  quasi-covering $\left(M\left( A\right) , M\left( \widetilde{A}\right) , G, M\left(\pi \right) \right)$  such that $\pi = \left.\widetilde \pi\right|_A$ and an action $G \times \widetilde A \to \widetilde{A}$ is a restriction of an action  $G \times M\left( \widetilde{A}\right) \to M\left( \widetilde{A}\right)$.
\end{corollary}

\begin{lemma}\label{ind_gmult_act_lem}
	If $\left(A, \widetilde{A}, G, \pi \right)$ is a noncommutative finite-fold  quasi-covering then one has an unique natural isomorphism
\bean
\left\{\left.g \in \Aut\left( \widetilde A\right)\right|\forall a \in A \quad  g\pi\left( a\right)  = \pi\left(a \right) \right\}\cong\\ \cong \left\{\left.g \in \Aut\left(M\left(  \widetilde A\right) \right)\right|\forall a \in M\left( A\right)  \quad  gM\left(\pi \right)\left( a\right) = M\left(\pi \right)\left( a\right)\right\}.
\eean
\end{lemma}
\begin{proof}
	If $\widetilde a \in M\left(  \widetilde A\right)$ then there is a net $\left\{\widetilde a_\la \right\}\subset A$ such that $\widetilde a= \bt$-$\lim_\la \widetilde a_\la$ where $\bt$-$\lim$ means the limit with respect to the strict topology (cf. Definition \ref{strict_topology_defn}) of $M\left( \widetilde{   A} \right)$. Let $g \in \left\{\left.g \in \Aut\left( \widetilde A\right)\right|\forall a \in A \quad  ga = a\right\}$, and let   $\left\{g\widetilde a_\la \right\}\subset \widetilde A$. If  $\widetilde b \in \widetilde A$ and a seminorm $\vertiii{\cdot }_{\widetilde b}$ is given by \eqref{strict_topology_norm_eqn} then one has
	$$
	\vertiii{g\widetilde a_\la}_{\widetilde b} = 	\vertiii{\widetilde a_\la}_{g^{-1}\widetilde b} 
	$$
	From the above equation it follows that the net $\left\{g\widetilde a_\la \right\}\subset \widetilde A$  is convergent with respect to the strict topology (cf. Definition \ref{strict_topology_defn}) of $M\left( \widetilde{   A} \right)$. We define 
	$$
	g \bt\text{-}\lim_\la \widetilde a_\la\bydef \bt\text{-}\lim_\la g\widetilde a_\la 
	$$
	Other details of the proof are left  to the reader.
\end{proof}
\begin{definition}\label{hereditary_lift_defn}
	Let $\left(A, \widetilde{A}, G, \pi \right)$ be a noncommutative finite-fold  quasi-covering (cf. Definition \ref{fin_quasi_defn}). Let $B$ be a hereditary $C^*$-subalgebra of $A$ (cf. Definition \ref{hered_defn}). A generated by $\pi\left( B\right)$ hereditary $C^*$-subalgebra  $\widetilde B$ of  $\widetilde A$ (cf. Definition \ref{hered_gen_defn}) is said to be a \textit{hereditary} $\left(A, \widetilde{A}, G, \pi \right)$-\textit{lift} of $B$. 
\end{definition}
\begin{remark}\label{hered_lift_rem}
If a net $\left\{u_\la\right\}_{\la \in \La} \subset M\left( B\right)_+$ is such that $\bt$-$\lim u_\la = 1_{M\left( B\right)}$ where the strict limit is implied (cf. Definition \ref{strict_topology_defn}) then similarly to \eqref{hered_uau_eqn} one has
\bea
\label{hered_lift_eqn}
\widetilde B = \left\{\left.\widetilde a \in \widetilde A\right|\lim_{\la \in \La}\left\| \widetilde a- \pi\left( u_\la\right)\widetilde a\pi\left( u_\la\right)\right\| =0 \right\}.
\eea 
\end{remark}
\begin{remark}\label{hered_lift_pap_rem}
	If $\widetilde \rho: \widetilde A\hookto B\left(\widetilde \H\right)$ is a faithful nondegenerate representation then similarly to the equation \ref{hered_r_u_eqn} one has
\be	\label{hered_r_u_lift_eqn}
	\forall\la\in\La\quad	\widetilde\rho\left(\pi\left( u_\la\right)\widetilde A\pi\left( u_\la\right) \right)\subset \widetilde\rho\left(\widetilde B  \right). 
\ee	
From the Lemma \ref{increasing_convergent_w_lem} it follows that there is a limit $\widetilde p = s$-$\lim_{\la\in \La} \widetilde\rho\left(\widetilde\pi\left( u_\la\right)  \right)$ with respect to the strong topology of $B\left(\widetilde \H \right)$ (cf. Definition \ref{strong_topology_defn}). Similarly to  the Corollary \ref{hered_representation_cor}  one has  a faithful, nondegenerate representation
\be\label{hered_lift_repr_eqn}
\widetilde B \hookto B\left(\widetilde p \widetilde\H \right). 
\ee

\end{remark}
\begin{remark}\label{hered_lift_au_rem}
	If $\left\{u_\la\right\}_{\la \in \La} \subset B$ be an approximate unit of $B$ then $\left\{\pi\left( u_\la\right) \right\}_{\la \in \La} \subset \widetilde B$ is an approximate unit of $\widetilde B$. 
\end{remark}

\begin{definition}\label{strictly_proper_defn}
	Let $\left(A, \widetilde{A}, G, \pi \right)$ be a noncommutative finite-fold  pre-covering (cf. Definition \ref{fin_pre_defn}) such that $A$ is a connected $C^*$-algebra. 	Let $B\subset A$ 	
	be a hereditary connected $C^*$-subalgebra (cf.  Definitions \ref{hered_defn}, \ref{connected_c_a_defn}), and let $\widetilde B$ be a hereditary $\left(A, \widetilde{A}, G, \pi \right)$-{lift} of $B$ (cf. Definition \ref{hereditary_lift_defn}). We say that $B$ is $\left(A, \widetilde{A}, G, \pi \right)$-\textit{proper} if one has a natural isomorphism 
	\be\label{hered_proper_lift_eqn}	
	\begin{split}
	\left\{ \left.g \in \Aut\left(\widetilde{A} \right)~\right|~ ga = a;~~\forall a \in A\right\}=\\
	= \left\{ \left.g \in \Aut\left(\widetilde{B} \right)~\right|~ gb = b;~~\forall b \in B\right\}.
	\end{split}
	\ee
	Moreover we say that $B$ is $\left(A, \widetilde{A}, G, \pi \right)$-\textit{strictly proper} if the quadruple $\left(B, \widetilde{B}, G, \pi|_B \right)$ is   a covering with unitization (cf. Definition \ref{fin_unitization_defn}).  (The action 	$G \times \widetilde{B}\to \widetilde{B}$, is the restriction on $\widetilde{B}\subset \widetilde{   A}$ of the action $G\times  \widetilde{A}\to \widetilde{A}$).
	
\end{definition}
\begin{lemma}\label{fin_composition_unihereditary_lem}
	Let $\left(A, \widetilde{A}, G, \pi \right)$ be a noncommutative finite-fold pre-covering   (cf. Definition \ref{fin_pre_defn}).  Let $H \subset G$ is  a $\left(A, \widetilde{A}, G, \pi \right)$-proper normal subgroup (cf. Definition \ref{proper_subgroup_fin_defn}). Let $B\subset A$ be a  $\left(A, \widetilde{A}, G, \pi \right)$-hereditary subalgebra. Let both \\$\left(A, \widetilde{A}^H, G\left(\left.\widetilde A~\right|A \right)/H, \pi^H \right)$ and 
	$\left(\widetilde{A}^H, \widetilde{A}, H, \left.\Id_{\widetilde{A}}\right|_{\widetilde{A}^H} \right)$ are given by the Definition \ref{proper_subgroup_fin_defn} and the Lemma \ref{proper_subgroup_fin_lem}. 
	 {noncommutative finite-fold  pre-coverings}. If a hereditary $C^*$-subalgebra $B\subset A$ is $\left(A, \widetilde{A}, G, \pi \right)$-{strictly proper} (cf. Definition \ref{strictly_proper_defn}) and  $\widetilde B$ is a hereditary $\left(A, \widetilde{A}, G, \pi \right)$-{lift} of $B$ (cf. Definition \ref{hereditary_lift_defn})	 
	 then following conditions hold:
\begin{itemize}
	\item[(i)] a $C^*$-algebra $\widetilde B^H$ is a $\left(A, \widetilde{A}^H, G\left(\left.\widetilde A~\right|A \right)/H, \pi^H \right)$-lift of $B$ (cf. Definition \ref{strictly_proper_defn});
	\item[(ii)] a $C^*$-algebra $\widetilde B$ is a $\left(\widetilde{A}^H, \widetilde{A}, H, \left.\Id_{\widetilde{A}}\right|_{\widetilde{A}^H} \right)$-lift of $\widetilde B^H$;
	\item[(iii)] a hereditary $C^*$-subalgebra $B\subset A$  is $\left(A, \widetilde{A}^H, G\left(\left.\widetilde A~\right|A \right)/H, \pi^H \right)$-{strictly proper}, and a hereditary subalgebra $\widetilde B^H\subset \widetilde A^H$ is $\left(\widetilde{A}^H, \widetilde{A}, H, \left.\Id_{\widetilde{A}}\right|_{\widetilde{A}^H} \right)$-strictly proper (cf. Definition \ref{strictly_proper_defn}).
\end{itemize}	 
\end{lemma}
\begin{proof}
	Let $\left\{u_\la\right\}_{\la \in \La} \subset B$ be an approximate unit of $B$.\\
(i) From the Remark \ref{hered_lift_rem} it follows that
\bean
\widetilde B = \left\{\left.\widetilde b \in \widetilde A\right|\widetilde b=\lim_{\la \in \La} \pi\left( u_\la\right)\widetilde b\pi\left( u_\la\right) \right\}.
\eean
On the other hand from
\bean
\widetilde A^H = \left\{\left.\sum_{	g \in H}g \widetilde a~ \right| \widetilde a \in \widetilde A\right\},\\
\widetilde B^H = \left\{\left.\sum_{	g \in H}g \widetilde b~ \right| \widetilde b \in \widetilde B\right\}= \left\{\left.\widetilde b \in \widetilde A^H\right|\widetilde b=\lim_{\la \in \La} \pi\left( u_\la\right)\widetilde b\pi\left( u_\la\right) \right\}
\eean
it follows that
\bean
\widetilde B^H = \left\{\left.\widetilde b \in \widetilde A^H\right|\widetilde b=\lim_{\la \in \La} \pi\left( u_\la\right)\widetilde b\pi\left( u_\la\right) \right\}.
\eean
So from the Remark \ref{hered_lift_rem} it follows that $\widetilde B^H$ is a $\left(A, \widetilde{A}^H, G\left(\left.\widetilde A~\right|A \right)/H, \pi^H \right)$-lift of $B$.\\
(ii) From the Remark \ref{hered_lift_au_rem} it follows that  $\left\{\pi^H\left( u_\la\right) \right\}_{\la \in \La} \subset \widetilde B$ is an approximate unit of $\widetilde B$. If $\widetilde B'$ is a $\left(\widetilde{A}^H, \widetilde{A}, H, \left.\Id_{\widetilde{A}}\right|_{\widetilde{A}^H} \right)$-lift of $\widetilde B^H$ then one has
\bean
\widetilde B'= \left\{\left.\widetilde b \in \widetilde B~\right|\widetilde b= \lim_{\la \in \La} \left.\Id_{\widetilde{A}}\right|_{\widetilde{A}^H}\circ\pi^H\left(u_\la \right)\widetilde b \left.\Id_{\widetilde{A}}\right|_{\widetilde{A}^H}\circ\pi^H\left(u_\la \right)\right\}=\\=\left\{\left.\widetilde b \in \widetilde A\right|\widetilde b=\lim_{\la \in \La} \pi\left( u_\la\right)\widetilde b\pi\left( u_\la\right) \right\}
\eean
and from the Remark \ref{hered_lift_au_rem} it follows that $\widetilde{B}' = \widetilde B$.\\
(iii) Follows from the Lemma \ref{fin_composition_unitization_lem}.
\end{proof}

\begin{definition}\label{fin_defn}
Let  $\left(A, \widetilde{A}, G, \pi \right)$ be a noncommutative finite-fold  pre-covering   (cf. Definition \ref{fin_pre_defn}). The triple $\left(A, \widetilde{A}, G, \pi \right)$ and/or the injective $*$-homomorphism $\pi: A\hookto \widetilde{A}$ is  a \textit{noncommutative finite-fold covering} if there is an indexed by a directed set $\La$ family  $\left\{A_\la \subset A\right\}_{\la\in\La}$ of connected $\left(A, \widetilde{A}, G, \pi \right)$-strictly proper  (cf. the Definition \ref{strictly_proper_defn}) hereditary $C^*$-subalgebras such that following conditions hold:
\begin{enumerate}
	\item[(a)] 
$$
\forall \mu, \nu \in \La \quad \mu \le \nu \quad\Rightarrow\quad A_\mu \subset A_\nu;
$$	
\item[(b)]
a union 
$
\bigcup_{\la\in\La} A_\la
$
is dense in $A$.
\end{enumerate}
  The group $G$ is said to be the \textit{finite covering transformation group} (of $\left(A, \widetilde{A},G, \pi\right)$ ) and we use the following notation
\begin{equation}\label{group_cov_eqn}
	G\left(\left.\widetilde{A}~~\right|~A \right) \stackrel{\mathrm{def}}{=} G.
\end{equation}

\end{definition}

\begin{remark}
	The Definition \ref{fin_defn} is motivated by the Theorem \ref{top_finite_covering_thm}.
\end{remark}
\begin{remark}
Roughly speaking the Definition \ref{fin_defn} is an approximation of any covering by coverings with compact spaces.
\end{remark}
\begin{remark}\label{unitization_covering_rem}
	If  $\left(A, \widetilde{A}, G, \pi \right)$ is a noncommutative finite-fold covering with unitization (cf. Definition \ref{fin_unitization_defn}) then a singleton  family $\left\{A\right\}$ satisfies to hypotheses of the Definition \ref{fin_defn}. 
So any noncommutative finite-fold covering with unitization (cf. Definition \ref{fin_unitization_defn}) is a noncommutative finite-fold covering.
\end{remark}

\begin{lemma}\label{fin_def_lem}
	Under the hypotheses of the Definition \ref{fin_defn}  if for all $\la\in\La$ a $C^*$-subalgebra $\widetilde{A}_\la\subset \widetilde A$ is   a {hereditary} $\left(A, \widetilde{A}, G, \pi \right)$-{lift} of $A_\la$ (cf. Definition \ref{hereditary_lift_defn}) then a union $\cup_{\la \in \La} \widetilde{A}_\la$ is dense in $\widetilde{A}$.
\end{lemma}
\begin{proof}
From the Lemma \ref{quasi_approximate_unit_lem} it follows that for any $\widetilde a \in \widetilde A$ and $\eps > 0$ there is a positive $u \in A_+$ such that
$$
\left\|\widetilde a- \pi\left(u \right)\widetilde a \pi\left(u \right)  \right\| < \frac{\eps}{2}.
$$
On the other hand since the union $\left\{A_\la\right\}_{\la \in \La}$ is dense in $A$ for any $\dl> 0$ there is $\la_\dl \in \La$ and positive $a_{\la_\dl}\in A_{\la_\dl}$ such that
$$
\left\|a_{\la_\dl} - u \right\| <\dl.
$$
So if $x \bydef a_{\la_\dl} - u$ then $\left\|x \right\| <\dl$. It follows that
\bean
\left\|\pi\left(a_{\la_\dl} \right)\widetilde a \pi\left(a_{\la_\dl} \right)- \pi\left(u \right)\widetilde a \pi\left(u \right)  \right\|= \left\|\pi\left(u + x\right)\widetilde a \pi\left(u+x \right) - \pi\left(u \right)\widetilde a \pi\left(u \right)  \right\|= \\
=\left\| \pi\left(x \right) \widetilde a \pi\left(x \right) + \pi\left(x \right) \widetilde a \pi\left(u \right)+ \pi\left(u \right) \widetilde a \pi\left(x \right)\right\|\le \left\|\widetilde a \right\|\left\|x \right\|^2 + 2\left\|\widetilde u \right\|\left\|\widetilde a \right\|\left\|x \right\|.
\eean 
If $\dl$ is such that $\left\|\widetilde a \right\|\dl^2 + 2\left\|\widetilde u \right\|\left\|\widetilde a \right\|\dl < \eps/2$ then 
\bean
\left\|\pi\left(a_{\la_\dl} \right)\widetilde a \pi\left(a_{\la_\dl} \right)- \pi\left(u \right)\widetilde a \pi\left(u \right)  \right\| <  \frac{\eps}{2},
\eean
\be\label{al_dense_eqn} 
\left\|\pi\left(a_{\la_\dl} \right)\widetilde a \pi\left(a_{\la_\dl} \right)- \widetilde a  \right\|< \eps.
\ee 
On the other hand from the Definition \ref{hereditary_lift_defn} and the Lemma \ref{hered_bab_lem} it follows that $\pi\left(a_{\la_\dl} \right)\widetilde a \pi\left(a_{\la_\dl} \right)\in \widetilde A_{\la_\dl}$, so from \eqref{al_dense_eqn} it turns out that the union $\cup_{\la \in \La} \widetilde{A}_\la$ is dense in $\widetilde{A}$.
\end{proof}
\begin{lemma}\label{fin_composition_covering_lem}
	Let $\left(A, \widetilde{A}, G, \pi \right)$ be a noncommutative finite-fold covering (cf. Definition \ref{fin_defn}).  If $H \subset G$ is  a $\left(A, \widetilde{A}, G, \pi \right)$-proper group (cf. Definition \ref{proper_subgroup_fin_defn}) then  both noncommutative finite-fold pre-coverings $\left(A, \widetilde{A}^H, G\left(\left.\widetilde A~\right|A \right)/H, \pi^H \right)$ and \\
	$\left(\widetilde{A}^H, \widetilde{A}, H, \left.\Id_{\widetilde{A}}\right|_{\widetilde{A}^H} \right)$  given by the Definition \ref{proper_subgroup_fin_defn} and the Lemma \ref{proper_subgroup_fin_lem} are noncommutative finite-fold coverings. 
\end{lemma}
\begin{proof}
Let $\left\{A_\la \subset A\right\}_{\la\in\La}$  be the required by the Definition \ref{fin_defn} indexed by a directed set $\La$ family of $\left(A, \widetilde{A}, G, \pi \right)$-{strictly proper} hereditary $C^*$-subalgebras (cf. the Definition \ref{strictly_proper_defn}). From  (iii) of the Lemma \ref{fin_composition_unihereditary_lem} it follows that $A_\la$ is a $\left(A, \widetilde{A}^H, G\left(\left.\widetilde A~\right|A \right)/H, \pi^H \right)$-{strictly proper} hereditary $C^*$-subalgebra for all  $\la\in \La$. Taking into account that a union $\cup_{\la\in \La} A_\la$ is dense in $A$ one concludes that the quadruple $\left(A, \widetilde{A}^H, G\left(\left.\widetilde A~\right|A \right)/H, \pi^H \right)$ a noncommutative finite-fold covering.

If $\widetilde A^H_\la$ is a hereditary $\left(A, \widetilde{A}^H, G\left(\left.\widetilde A~\right|A \right)/H, \pi^H \right)$-lift of $A_\la$ (cf. Definition \ref{hereditary_lift_defn}) then from  the Lemma \ref{fin_def_lem} it turns out that the union $\cup_{\la\in \La} \widetilde A^H_\la$ is dense in $\widetilde A^H$. On the other hand from (iii) 
of the Lemma \ref{fin_composition_unihereditary_lem} it follows that $\widetilde{A}^H_\la$ is a $\left(A, \widetilde{A}^H, G\left(\left.\widetilde A~\right|A \right)/H, \pi^H \right)$-strictly proper hereditary $C^*$-subalgebra. From the Definition \ref{fin_defn} it turns out that the quadruple 
	 $\widetilde{A}^H_\la$ is a $\left(A, \widetilde{A}^H, G\left(\left.\widetilde A~\right|A \right)/H, \pi^H \right)$ is a noncommutative finite-fold covering.
\end{proof}
	\begin{definition}\label{fin_composition_covering_defn}
		Under the hypotheses of  the Lemma \ref{fin_composition_covering_lem} we say that a noncommutative finite-fold covering $\left(A, \widetilde{A}^H, G\left(\left.\widetilde A~\right|A \right)/H, \pi^H \right)$ (resp.
		$\left(\widetilde{A}^H, \widetilde{A}, H, \left.\Id_{\widetilde{A}}\right|_{\widetilde{A}^H} \right)$) is $H$-\textit{regular} (resp. $H$-\textit{singular}).
\end{definition}

\begin{remark}
	If we put $\left\{A_\la\right\}_{\la\in\La} = \left\{A \right\}$  then from the Definition \ref{fin_defn} it follows that   any  noncommutative finite-fold covering with unitization is a  noncommutative finite-fold covering.
\end{remark}

\begin{remark}
	The Definition \ref{fin_defn}  is motivated by the Lemmas  \ref{top_fin_necassary_lem} and \ref{top_fin_sufficient_lem} .
\end{remark}

\begin{definition}\label{composable_coverings_defn}
A pair  $\left( \left(\widetilde{A}', \widetilde{A}'', G'', \pi'' \right), \left(A, \widetilde{A}', G', \pi' \right)\right)$ of noncommutative finite-fold coverings is said to be \textit{composable} if there is a noncommutative finite-fold covering $\left(A, \widetilde{A}'', G, \pi \right)$ satisfying following conditions:
\begin{itemize}
	\item $\pi = \pi' \circ \pi''$,
	\item there is  $H \subset G$ is  a $\left(A, \widetilde{A}, G, \pi \right)$-proper group (cf. Definition \ref{proper_subgroup_fin_defn}), such that both $\left(A, \widetilde{A}', G', \pi' \right)$ and $\left(\widetilde{A}', \widetilde{A}'', G'', \pi'' \right)$ are equivalent to $H$-{regular} and $H$-{singular} coverings (cf. Definition \ref{fin_composition_covering_defn}) respectively.
\end{itemize}
We say that $\left(A, \widetilde{A}'', G, \pi \right)$ is a \textit{composition} of $\left(\widetilde{A}', \widetilde{A}'', G'', \pi'' \right)$ and $\left(A, \widetilde{A}', G', \pi' \right)$.
\end{definition}

\begin{definition}\label{simply_connected_defn}
	A connected  $C^*$-algebra $A$ is said to be \textit{simply connected} if it   has no nontrivial   noncommutative finite-fold coverings.
\end{definition}
\begin{remark}
	The Definition \ref{simply_connected_defn} can be regarded as a generalization of \ref{top_weakly_simply_connected_defn} one (cf. Theorem  \ref{top_finite_covering_thm}  below).
\end{remark}

\begin{definition}\label{fin_weak_defn} 
	A noncommutative finite-fold  pre-covering  $\left(A, \widetilde{A}, G, \pi \right)$  	
	is said to be   a \textit{noncommutative weak finite-fold covering} if there is  a { noncommutative  finite-fold covering} $\left(B, \widetilde{B}, G, \rho \right)$  (cf. Definition \ref{fin_defn} such that
	\begin{enumerate}
		\item [(a)] Both $B$ and $\widetilde{B}$ are essential ideals of $A$ and  $\widetilde{A}$ respectively.
		\item[(b)] $G\widetilde{B}= \widetilde{B}$, and an  action $G\times \widetilde{B}\to \widetilde{B}$ comes from $G\times \widetilde{A}\to \widetilde{A}$ one.
		\item[(c)] $\rho = \pi|_B$.
	\end{enumerate}
	If 	the $C^*$-algebra $\widetilde A$ is  (cf. Definition \ref{connected_c_a_defn}) then $\left(A, \widetilde{A}, G, \pi \right)$  is said to be   a \textit{ noncommutative finite-fold weak covering} or simply a \textit{noncommutative finite-fold weak covering}.
\end{definition}

\begin{definition}\label{fin_category_defn} 
	For any connected $C^*$-algebra $A$ consider a category $\mathfrak{FinCov}$-$A$ (cf. Definition \ref{category_defn}) such that following conditions hold:
	\begin{itemize}
		\item [(i)]  $\mathfrak{FinCov}$-$A$-objects  are  finite-fold noncommutative coverings $\pi: A \hookto \widetilde A$ (or $\left(A, \widetilde{A}, G, \pi \right)$) (cf. Definition \ref{fin_defn}).
		\item[(ii)] A $\mathfrak{FinCov}$-$A$-morphism from $\left(A, \widetilde{A}', G', \pi' \right)$ to $\left(A, \widetilde{A}'', G'', \pi'' \right)$ is an injective *-homomorphism $\pi: \widetilde A''\hookto \widetilde A'$ such that the following diagram 
				\newline
		\begin{tikzcd}
\widetilde A'	& &\arrow[ll,  "\pi" {yshift=10pt}]\widetilde A'' \\
& A \arrow[lu, "\pi'"] \arrow[ru, "\pi''" {xshift=12pt,yshift=-12pt}]&		
\end{tikzcd}
		\\ 	
		
		is commutative.
		\item[(iii)] Compositions of morphisms come from compositions of *-homomorphisms.
	\end{itemize}
	We say that $\mathfrak{FinCov}$-$A$ is the \textit{category of finite-fold coverings} of $A$. 
	Sometimes we write $\widetilde A$ instead of $\pi: A \hookto \widetilde A$ to designate an object of $\mathfrak{FinCov}$-$A$. 
\end{definition}

\begin{definition}\label{fin_composable_cat_defn}
	A full subcategory $\mathscr C\subseteqq\mathfrak{FinCov}$-$A$ (cf. Definition \ref{subcategory_defn}) of $\mathfrak{FinCov}$-$A$ is said to be \textit{composable} then for any $\mathscr C$-morphism $\pi: \widetilde A''\hookto \widetilde A'$ from $\left(A, \widetilde{A}', G', \pi' \right)$ to $\left(A, \widetilde{A}'', G'', \pi'' \right)$  there is a $\left(A, \widetilde{A}', G', \pi' \right)$-proper normal subgroup $H\subset G''$ (cf. Definition \ref{proper_subgroup_fin_defn}) such that there is a *-isomorphism $\widetilde{A}'' \cong \widetilde{A}'^H$ and $\pi$ is equivalent to a natural inclusion $ \widetilde{A}'^H\subset \widetilde{A}'$, i.e. there is a natural noncommutative finite-fold covering 
	$$
	\left(\widetilde A'', \widetilde A', G\left(\left.  \widetilde A'\right| \widetilde A''\right), \pi\right)\cong \left(\widetilde A'^H, \widetilde A', H, \Id_{\widetilde A'}|_{\widetilde A''}\right)
	$$
	(cf. Lemma \ref{fin_composition_covering_lem}).
\end{definition}

\begin{definition}\label{proper_composition_defn}
	The  $C^*$-algebra $A$ is \textit{proper with respect to covering compositions} if for any noncommutative finite-fold coverings  $\left(A, \widetilde{A}', G\left(\left.\widetilde{A}'~\right|{A}  \right), \pi \right)$ and  
	$\left(\widetilde A', \widetilde{A}, G\left(\left.\widetilde{A}~\right|\widetilde {A}'  \right), \widetilde \pi \right)$ and there are a  noncommutative finite-fold covering \\$\left(A, \widetilde{A}, G\left(\left.\widetilde{A}~\right|{A}  \right),\widetilde \pi\circ \pi \right)$ and a $\left(A, \widetilde{A}, G\left(\left.\widetilde{A}~\right|{A}  \right),\widetilde \pi\circ \pi \right)$-proper normal subgroup $H\subset G\left(\left.\widetilde{A}~\right|\widetilde {A}'  \right)$ (cf. Definition \ref{proper_subgroup_fin_defn}) such that there are following equivalences:
	\bean
	\left(A, \widetilde{A}', G\left(\left.\widetilde{A}'~\right|{A}  \right), \pi \right)\cong 	\left(A, \widetilde{A}^H, G\left(\left.\widetilde{A}~\right|{A}  \right)/H, \pi \right),\\
	\left(\widetilde A', \widetilde{A}, G\left(\left.\widetilde{A}~\right|\widetilde {A}'  \right), \widetilde \pi \right)\cong 	\left(\widetilde A^H, \widetilde{A},H,\left. \Id_{\widetilde A}\right|_{\widetilde A^H} \right)
	\eean
	of noncommutative finite-fold coverings.
\end{definition}

\subsection{Coverings of *-algebras} 
\begin{definition}\label{fin_quasi*_defn}
	Let both  $A$ and $\widetilde{A}$ be  connected $*$-algebras (cf. Definition \ref{connected_c_a_defn}), and let $\pi: A \hookto \widetilde{A}$ be an injective, unital (in the sense of the Definition \ref{principal_non_defn}) $*$-homomorphism of 
	$C^*$-algebras. Let $G$ be a finite  group of *-automorphisms of $\widetilde{A}$ such that 	$\pi\left(A\right) = \widetilde{A}^G\stackrel{\text{def}}{=}\left\{
	\left.a\in \widetilde{A}~\right|~ a = g a;~ \forall g \in G\right\}$.	We say that the triple $\left(A, \widetilde{A}, G \right)$ and/or the quadruple $\left(A, \widetilde{A}, G, \pi \right)$ and/or $*$-homomorphism $\pi: A \hookto \widetilde{A}$   is a \textit{noncommutative finite-fold  quasi-covering}. We write
	\be\label{fin_cov_gr*_eqn}
	G\left(\left.\widetilde{A}~\right| {A} \right) \stackrel{\text{def}}{=}  	G.
	\ee
\end{definition}

\begin{definition}\label{fin_pre*_defn}  
	Let $\pi: A \hookto \widetilde{A}$ be an injective  $*$-homomorphism of  connected
	$*$-algebras such that following conditions hold:
	\begin{enumerate}
		\item[(a)] If $\Aut\left(\widetilde{A} \right)$ is a group of *-automorphisms of $\widetilde{A}$ then the group  
		\be\nonumber
		G = \left\{ \left.g \in \Aut\left(\widetilde{A} \right)~\right|~ ga = a;\quad\forall a \in A\right\}
		\ee
		is finite.
		\item[(b)] 	$A \cong \widetilde{A}^G\stackrel{\text{def}}{=}\left\{\left.a\in \widetilde{A}~~\right|~ a = g a;~ \forall g \in G\right\}$.
	\end{enumerate}
	We say that the triple $\left(A, \widetilde{A}, G \right)$ and/or the quadruple $\left(A, \widetilde{A}, G, \pi \right)$ and/or $*$-homomorphism $\pi: A \hookto \widetilde{A}$   is a \textit{noncommutative finite-fold  pre-covering of *-algebras}. If both $A$ and $\widetilde{A}$ are topological $*$-algebras then we assume that the homomorphism $\pi$ is continuous.
\end{definition}
\subsubsection{Coverings of pro-$C^*$-algebras}
\begin{definition}\label{pro_fin_defn} 
	Let  $\left(A, \widetilde{A}, G, \pi \right)$ be a pre-covering (cf. Definition \ref{fin_pre*_defn}) of pro-$C^*$-algebras (cf. Definition \ref{pro_c_defn}) where $\pi$ is regarded as an inclusion. 
		We say that the triple $\left(A, \widetilde{A}, G \right)$ and/or the quadruple $\left(A, \widetilde{A}, G, \pi \right)$ and/or $*$-homomorphism $\pi: A \hookto \widetilde{A}$   is a \textit{noncommutative finite-fold  covering  of pro-$C^*$-algebras} 
		if \\$\left(b\left(A\right) ,  b\left(\widetilde{A}\right) , G, \pi|_{b\left(A\right)} \right)$ is a   noncommutative  finite-fold weak covering (cf. Definition \ref{fin_weak_defn}), where both $b\left({A}\right)$ and $b\left(\widetilde{A}\right)$ are  $C^*$-algebras of bounded elements (cf. the Definition \ref{pro_bound_defn} and the Proposition \ref{pro_bound_prop}).
\end{definition}

\subsubsection{Coverings of bounded operator $*$-algebras}
\begin{empt}\label{fin_oa_empt}
	Let $\left(A, \widetilde{A}', G, \pi \right)$ be  a noncommutative finite-fold  pre-covering of *-algebras (cf. Definition \ref{fin_pre*_defn}). Suppose that both $A$ and $\widetilde{A}'$ are dense *-subalgebras of $C^*$-algebras. Let both  $B$ and  $\widetilde{B}$ be  $C^*$-algebras which are $C^*$-norm completions of $A$ and $\widetilde{A}'$ respectively. Assume that both a $*$-homomorphism $\pi_B: B \hookto\widetilde B$ and an action $G \times \widetilde{B}\to  \widetilde{B}$ come from $\pi: A \hookto\widetilde A'$ and an action $G \times \widetilde{A}'\to  \widetilde{A}'$. If  $\left(B, \widetilde{B}, G, \pi_B \right)$ is a  noncommutative finite-fold covering (cf. Definition \ref{fin_defn}) then we say that *-algebra  $\widetilde{A}'$ is \textit{admissible}.

\end{empt}
	\begin{definition}\label{fin_oa_defn}
	In the above situation if $\widetilde{A}\subset \widetilde{B}$ is a maximal admissible *-algebra 
then we say that  $\left(A, \widetilde{A}, G, \pi \right)$, or $\left(A, \widetilde{A}, G  \right)$, or $\pi$ is a \textit{noncommutative finite-fold covering of bounded operator *-algebras}.
\end{definition}

\begin{remark}
	Any union of simply ordered by inclusion  set of admissible algebras is an admissible algebra.	 From the Zorn's lemma (cf. \ref{zorn_thm}) it follows that  there is one are more maximal admissible algebras. 
\end{remark}

\subsubsection{Coverings of $O^*$-algebras}

\begin{empt}\label{o*fin_empt}
	Let $\left(A, \widetilde{A}, G, \pi \right)$ be a finite-fold pre-covering of *-algebras (cf. Definition \ref{fin_pre*_defn}) such that $\widetilde{A}$ is an $O^*$-algebra (cf. Definition \ref{o*alg_defn}); that is $\widetilde{A}\subset \L^\dagger\left(\D\right)$ (cf. \eqref{l_dag_eqn}) where $\D$ is a dense subspace of a Hilbert space $\H$. Assume that there is an equivariant action $G\times \D\to \D$, i.e.
	\be\label{equ_d_eqn}
\forall \widetilde  a \in \widetilde A\quad	\forall \xi\in \D \quad \forall g\in G \quad g\left( \widetilde  a\xi\right)=   g\left( \widetilde  a\right) g\left( \xi\right).
	\ee 
	 Let $\widetilde{A}_b \bydef \widetilde{A}\cap \L^\dagger\left(\D\right)_b$ (cf. \eqref{o*b_eqn}).
	From the Remark \ref{o*b_rem} it follows that there are natural $C^*$-norms completions $\widetilde B$ and ${B}$  of $\widetilde{A}_b$ and $A_b\bydef \widetilde{A}_b\cap A$ respectively. Suppose  that there is  a  noncommutative weak finite-fold covering
$\left(B, \widetilde{B}, G, \pi_B \right)$
	(cf. Definition  \ref{fin_weak_defn})
	where the action $G\times \widetilde{B} \to \widetilde{B}$ comes from   $G\times \widetilde{A} \to \widetilde{A}$ and ${\pi}_B$ comes from $\pi$.
\end{empt}
\begin{definition}\label{fino*_defn}
	In the  situation \ref{o*fin_empt} we say that $\left(A, \widetilde{A}, G, \pi \right)$, or $\left(A, \widetilde{A}, G  \right)$, or $\pi: A \hookto \widetilde A$ is a \textit{noncommutative finite-fold covering of $O^*$-algebras}.
\end{definition}
\subsubsection{Coverings and unbounded operators on Hilbert modules}\label{fin_chull_sec}
\paragraph*{}
The Definition \ref{fino*_defn} assumes  that both *-algebras $A$ and $\widetilde{A}$ contain  quite enough bounded elements. Here we consider an opposite situation.
\begin{empt}\label{fin_chull_empt}
Let $\left(A, \widetilde{A}, G, \pi_A \right)$ be a {noncommutative finite-fold  pre-covering of *-algebras} (cf. Definition \ref{fin_pre*_defn}), and let $\left(B, \widetilde{B}, G, \pi_B \right)$ be a  {noncommutative finite-fold  covering} (cf. Definition \ref{fin_defn}). 
Suppose that following conditions hold:
\begin{enumerate}
\item[(a)] There are a dense right ideal $\mathfrak{B}$ of $B$ and  a representation  $\left( \mathfrak{B},\mu\right)$ of $A$ on $B$ (cf. Definition \ref{def:rep_Hilbert_module_uni}), where $B$ is regarded as a Hilbert module over $B$, i.e. one has an inclusion  $\mu : A \hookto \End^*_{ B}\left( {\mathfrak B}\right)$  of *-algebras.
\item[(b)] There is  a dense right ideal $\widetilde{\mathfrak{B}}\subset \widetilde B$ , such that $G\widetilde{\mathfrak{B}} = \widetilde{\mathfrak{B}}$ and
\be\label{fin_chull_e_eqn}
{\mathfrak{B}}= \widetilde{\mathfrak{B}}^G \bydef \left\{\left. \widetilde{b}\in \widetilde{\mathfrak{B}}\right| \forall g \in G\quad g \widetilde b = \widetilde b\right\}.
\ee
\item[(c)] There is   a representation  $\left( \widetilde{\mathfrak{B}},\widetilde\mu\right)$ of $\widetilde A$ on $\widetilde B$ (cf. Definition \ref{def:rep_Hilbert_module_uni}),  where $\widetilde B$ is regarded as a Hilbert module over $\widetilde B$, i.e. one has an inclusion  $\widetilde\mu : \widetilde A \hookto \End^*_{ \widetilde B}\left( \widetilde{\mathfrak B}\right)$  of *-algebras.
\end{enumerate}
\end{empt}
\begin{definition}\label{fin_chull_defn}
Consider the situation  \ref{fin_chull_empt}. If following conditions hold:
\begin{enumerate}
	\item [(a)]  $\pi_A: A \hookto \widetilde{A}$ agrees with the natural inclusion $\End^*_{ B}\left( {\mathfrak B}\right)\subset \End^*_{\widetilde B}\left( \widetilde{\mathfrak B}\right)$.
	\item[(b)] The action $G\times  \widetilde A \to \widetilde A$ comes from the natural action $G \times \End^*_{\widetilde B}\left( \widetilde{\mathfrak B}\right) \to\End^*_{\widetilde B}\left( \widetilde{\mathfrak B}\right) $.
\end{enumerate}
	then  we say that $\left(A, \widetilde{A}, G, \pi_A \right)$ is an \textit{associated with $\left(B, \widetilde{B}, G, \pi_B \right)$ noncommutative finite-fold covering of *-algebras}.
\end{definition}

\subsection{Coverings of quasi *-algebras}

\begin{definition}\label{fin_preq_defn}
	Suppose that $\pi: \left(\mathfrak{A}, \mathfrak{A}_0\right) \hookto \left(\widetilde{\mathfrak{A}}, \widetilde{\mathfrak{A}}_0\right)$ is an injective unital $*$-homomorphism of   quasi $*$-algebras (cf. Definitions \ref{quasi_defn} and \ref{quasi_hom_defn})  such that following conditions hold:
	\begin{enumerate}
		\item[(a)] If $\Aut\left(\widetilde{\mathfrak{A}}, \widetilde{\mathfrak{A}}_0\right)$ is a group of *-automorphisms of $\widetilde{\mathfrak{A}}, \widetilde{\mathfrak{A}}_0$ then the group  
		\be\nonumber
		G = \left\{ \left.g \in \Aut\left(\widetilde{\mathfrak{A}}, \widetilde{\mathfrak{A}}_0\right)~\right|~ ga = a;~~\forall a \in \mathfrak A\right\}
		\ee
		is finite.
		\item[(b)] 	$\widetilde{\mathfrak{A}} \cong \widetilde{\mathfrak{A}}^G\stackrel{\text{def}}{=}\left\{\left.a\in \widetilde{\mathfrak{A}}~~\right|~ a = g a;~ \forall g \in G\right\}$.
	\end{enumerate}
	We say that the triple $\left(\left( {\mathfrak{A}}, {\mathfrak{A}}_0\right) , \left( \widetilde{\mathfrak{A}}, \widetilde{\mathfrak{A}}_0\right) , G \right)$ and/or the quadruple\\ $\left(\left( {\mathfrak{A}}, {\mathfrak{A}}_0\right) , \left( \widetilde{\mathfrak{A}}, \widetilde{\mathfrak{A}}_0\right) , G , \pi\right)$ and/or $*$-homomorphism $\pi: \left(\mathfrak{A}, \mathfrak{A}_0\right) \hookto \left(\widetilde{\mathfrak{A}}, \widetilde{\mathfrak{A}}_0\right)$   is a \textit{noncommutative finite-fold pre-covering of quasi $*$-algebras}. 
\end{definition}

\begin{definition}\label{oq*fin_defn}
	Let   $\left(\left( {\mathfrak{A}}, {\mathfrak{A}}_0\right) , \left( \widetilde{\mathfrak{A}}, \widetilde{\mathfrak{A}}_0\right) , G , \pi\right)$ be a {noncommutative finite-fold pre-covering of quasi-$*$-algebras}. If the quadruple  $\left(\mathfrak{A}_0 ,  \widetilde{\mathfrak{A}}_0 , G , \pi|_{\mathfrak{A}_0}\right)$ is a {noncommutative finite-fold covering of $O^*$-algebras} (cf. Definition \ref{fino*_defn}) then we say that  the quadruple $\left(\left( {\mathfrak{A}}, {\mathfrak{A}}_0\right) , \left( \widetilde{\mathfrak{A}}, \widetilde{\mathfrak{A}}_0\right) , G , \pi\right)$ and/or $*$-homomorphism $\pi:  \left(\mathfrak{A}, \mathfrak{A}_0\right) \hookto \left(\widetilde{\mathfrak{A}}, \widetilde{\mathfrak{A}}_0\right)$   is a \textit{noncommutative finite-fold covering of quasi-$*$-algebras}.
\end{definition}

\section{Basic constructions}\label{cov_fin_basc_sec}


 \begin{lemma}\label{frame_lem}
 	If $A$ is an unital $C^*$-algebra, and $P$ is a finitely generated  Hilbert $A$-module then 
 	there is a finite subset $\left\{ \widetilde{a}_1, ..., \widetilde{a}_n, \widetilde{b}_1, ..., \widetilde{b}_n\right\} \subset P$ such that
 	\be\label{frame_eqn}
 	\sum_{j=1}^{n}\widetilde{b}_j\left\rangle \right\langle\widetilde{a}_j =	1_{\End\left( P\right)_A}= 1_{\K\left( P\right)_A},
 	\ee
 where the  Notation \ref{rank_one_notation_eqn} is used.
 \end{lemma}
 
 \begin{proof}
 	If $P$ is generated by $\widetilde{a}_j, ..., \widetilde{a}_n\in P$ as a right $A$-module and 	 $S \in \End^*\left( P\right)_A$ is given by  
 	\be\label{s_matr_eqn}
 	S = \sum_{j=1}^{n}\widetilde{a}_j\left\rangle \right\langle\widetilde{a}_j
 	\ee
 	then $S$ is self-adjoint. From the Corollary 1.1.25  of \cite{jensen_thomsen:kk} it turns out that $S$ is strictly positive. Otherwise $P$ is a finitely generated right $A$-module, so from the Exercise  15.O of \cite{wegge_olsen} it follows that $S$ is invertible. The spectrum of $S$ is contained in 
 	$\R_\eps\bydef \left\{x \in \R| x > \eps\right\}$ where $\eps > 0$ depends on $S$. So if $\phi : \R \to \R$ is given by $z \mapsto \frac{1}{z}$ then   $T\bydef\phi\left(S \right)\in  \End^*\left( P\right)_A $. Moreover one has 
 	$$
 	TS = 1_{\End^*\left( P\right)_A},
 	$$ 
 	it turns out that if $\widetilde{b}_j \bydef T\widetilde{a}_j$
 	then
 	\bean
 	\sum_{j=1}^{n}T\widetilde{a}_j\left\rangle \right\langle\widetilde{a}_j = \sum_{j=1}^{n}\widetilde{b}_j\left\rangle \right\langle\widetilde{a}_j=	1_{\End^*\left( P\right)_A}= 1_{\K\left( P\right)_A}.
 	\eean
 	
 \end{proof}



 \begin{remark}\label{frame_rem}
 	From the Lemma \ref{frame_lem} and the equation \eqref{finite_hilb_mod_prod_eqn}  it follows that if $ \left(A, \widetilde{A}, G, \pi \right)$ is   an {unital noncommutative finite-fold  quasi-covering} (cf. Definition \ref{fin_unital_defn}) then 
 	there is a finite subset $$\left\{ \widetilde{a}_1, ..., \widetilde{a}_n, \widetilde{b}_1, ..., \widetilde{b}_n\right\} \subset \widetilde{A}$$ such that
 	\be\label{frame_rem_eqn}
\forall \widetilde{a} \in \widetilde{A}\quad \widetilde{a}=	\sum_{j=1}^{n}\widetilde{b}_j\left(  \sum_{g \in G} g\left(\widetilde{a}^*_j\widetilde{a} \right) \right).
 	\ee
 \end{remark}
 \begin{empt}\label{frame_l_empt}
 	Let $\left(A, \widetilde{A}, G, \pi \right)$ be a {noncommutative finite-fold covering with unitization}. If $\left(B, \widetilde{B}, G, \phi \right)$ is an unital noncommutative finite-fold covering which satisfies to conditions (a), (b) of the Definition \ref{fin_unitization_defn} then 	there is a finite subset $$\left\{ \widetilde{a}_1, ..., \widetilde{a}_n, \widetilde{b}_1, ..., \widetilde{b}_n\right\} \subset \widetilde{B}$$ such that
 	\be\nonumber
\forall \widetilde{b} \in \widetilde{B}\quad \widetilde{b}=	\sum_{j=1}^{n}\widetilde{b}_j\left(  \sum_{g \in G} g\left(\widetilde{a}^*_j\widetilde{b} \right) \right)
	\ee
 	(cf.\eqref{frame_rem_eqn}). In particular if $\widetilde{a} \in \widetilde{A}$ then
 \be\label{frame_l_eqn}
\forall \widetilde{a} \in \widetilde A	\quad	\widetilde{a}=	\sum_{j=1}^{n}\widetilde{b}_ja_j\quad \text{ where }\quad a_j\bydef  \sum_{g \in G} g\left(\widetilde{a}^*_j\widetilde{a} \right)\in \pi\left( A\right) 
\ee
 	and taking into account 
 	\bean
 	\widetilde{a}^*_j\widetilde{a} \in \widetilde{A},\\
 	\sum_{g \in G} g\left(\widetilde{a}^*_j\widetilde{a} \right)  \in \widetilde{A}^G =A 
 	\eean
 	one has
 	\be\label{fin_comp_eqn}
 	\widetilde{A} = \widetilde{b}_1A+...+\widetilde{b}_n A
 	\ee
 From \eqref{fin_comp_eqn} it follows that
  	\be\label{fin_comp*_eqn}
 \widetilde{A} = A\widetilde{b}_1^*+...+A\widetilde{b}_n^*
 \ee
 Moreover from \ref{four_decompositon_eqn} one can deduce that
   	\be\label{fin_compp_eqn}
 \exists \left\{\widetilde{b}_1, ..., \widetilde{b}_n\right\}\subset  \widetilde{B}\subset M\left( \widetilde A\right)  \quad \widetilde{A} = A\widetilde{b}_1+...+A\widetilde{b}_n =  \widetilde{b}^*_1A+...+\widetilde{b}^*_n A.
 \ee
\end{empt}
\begin{lemma}\label{cov_sep_fin_lem}
	Let $\left(A, \widetilde{A}, G, \pi \right)$ be a {noncommutative finite-fold quasi-covering with unitization}. If $A$ is a separable $C^*$-algebra then $\widetilde{A}$ is a separable $C^*$-algebra.
\end{lemma}
\begin{proof}
	Follows from the equation \eqref{fin_comp_eqn}.
\end{proof}

\begin{lemma}\label{cov_mult_fin_lem}
	Let $\left(A, \widetilde{A}, G, \pi \right)$ be a {noncommutative finite-fold quasi-covering with unitization}. 	Let $\left(B ,\widetilde{B}, G, \widetilde{\pi} \right)$ be an 
	unital  noncommutative finite-fold quasi-covering which satisfies to conditions  (a), (b) of the Definition \ref{fin_unitization_defn}. Then there is $\left\{  \widetilde{b}_1, ..., \widetilde{b}_n\right\} \subset \widetilde{B}$ such that 
	\be\label{ma_fin_eqn}
	M\left(\widetilde{A} \right) = \widetilde{b}_1 M\left(A \right) + ... + \widetilde{b}_n M\left(A \right).
	\ee
	
\end{lemma}
\begin{proof} 
	It is known that for any $\widetilde{a} \in M\left(\widetilde{A} \right)$ there is a net  $\left\{\widetilde{a}_\la\right\} \subset \widetilde{A}$ (cf. Definition \ref{top_net_defn}) such that there is the limit $\widetilde{a}=\bt$-$\lim \widetilde{a}_\la$ with respect to the strict topology (cf. Definition \ref{strict_topology_defn}). From the Remark \ref{frame_rem} turns out that there is finite subset $$\left\{ \widetilde{a}_1, ..., \widetilde{a}_n, \widetilde{b}_1, ..., \widetilde{b}_n\right\} \subset \widetilde{B}$$ such that
	$$
	\widetilde{b}=	\sum_{j=1}^{n}\widetilde{b}_j\left(  \sum_{g \in G} g\left(\widetilde{a}^*_j\widetilde{b} \right) \right)\quad \forall \widetilde{b} \in \widetilde{B} \quad \text{(cf. }\eqref{frame_rem_eqn}\text{)}.
	$$
If 	$\widetilde a \in M\left(\widetilde{A} \right)$ and $\widetilde a= \bt\text{-}\lim \widetilde{a}_\la$ (cf. Definition \eqref{strict_topology_defn}) where $\left\{\widetilde{a}_\la\right\}\subset \widetilde{A}$ then one has 
	$$
	\widetilde{a}_\la=	\sum_{j=1}^{n}\widetilde{b}_j\left(  \sum_{g \in G} g\left(\widetilde{a}^*_j\widetilde{a}_\la \right) \right).
	$$	
 From the above equation it turns out
	\be\label{projective_decomposition_eqn}
	\begin{split}
		\widetilde{a}= \bt\text{-}\lim \sum_{j=1}^{n}\widetilde{b}_j\left(  \sum_{g \in G} g\left(\widetilde{a}^*_j\widetilde{a}_\la \right) \right)= \sum_{j=1}^{n}\widetilde{b}_j\left(  \sum_{g \in G} g\left(\widetilde{a}^*_j\bt\text{-}\lim\widetilde{a}_\la \right) \right)= \\
		= \sum_{j=1}^{n}\widetilde{b}_j\left( \sum_{g \in G} g\left(\widetilde{a}^*_j\widetilde{a} \right) \right)= \sum_{j=1}^{n}\widetilde{b}_ja_j; \text{ where } a_j = \sum_{g \in G} g\left(\widetilde{a}^*_j\widetilde{a} \right)\in M\left(A \right) .
	\end{split}
	\ee 
	The element $\widetilde{a} \in M\left(\widetilde{A} \right)$  is arbitrary, so one has
	\be\nonumber
	M\left(\widetilde{A} \right) = \widetilde{b}_1 M\left(A \right) + ... + \widetilde{b}_n M\left(A \right).
	\ee
	
\end{proof}
\begin{remark}\label{proj_decomp_rem}
	If we consider the equation \eqref{projective_decomposition_eqn} then one has
	\be\label{projective_decompositioneq_eqn}
	\widetilde{a}=0 \quad \Leftrightarrow \quad \forall j =1,...,n \quad a_j = 0.
	\ee
	There are homomorphisms of $B$-modules
	\be\label{tb_bn_eqn}
	\begin{split}
		\iota_B : \widetilde{B}\to B^n,\\
		\widetilde b \mapsto \frac{1}{\left| G\right| } \left( \sum_{g \in G} g\left(\widetilde{a}^*_1\widetilde{b} \right),...,  \sum_{g \in G} g\left(\widetilde{a}^*_n\widetilde{b} \right) \right),\\
		p_{B^n} : B^n \to \widetilde B,\\
		\left(b_1, ..., b_n\right)\mapsto \widetilde{b}_1b_1 + ... + \widetilde{b}_1b_1.
	\end{split}
	\ee\label{tb_bn_p_eqn}
	such that
	\be
	p_{B^n}\circ \iota_{\widetilde B} = \Id_{\widetilde B}.
	\ee
	
\end{remark}
\begin{corollary}\label{cov_mult_fin_cor}
If $\left(A, \widetilde{A}, G, \pi \right)$ is a {noncommutative finite-fold quasi-covering with unitization} (cf. Definition \ref{fin_unitization_defn}) then the quadruple   $\left(M\left(A \right)  , M\left( \widetilde{A}\right) , G, M\left( \pi\right)  \right)$
unital  noncommutative finite-fold quasi-covering (cf. Definition \ref{fin_unital_defn}).
\end{corollary}

\begin{corollary}\label{cov_mult_fin_c_cor}
	If $\left(A, \widetilde{A}, G, \pi \right)$ is a {noncommutative finite-fold covering with unitization} (cf. Definition \ref{fin_unitization_defn}) then there is  an 
	unital  noncommutative finite-fold covering  $\left(M\left(A \right)  ,M\left( \widetilde{A}\right) , G, M\left(\pi \right) \right)$ (cf. Definition \ref{fin_unital_defn}).
\end{corollary}
\begin{proof}
Since both $A \subset M\left(A \right)$ and $\widetilde{A}\subset M\left( \widetilde{A}\right)$ are essential ideals one has
\bean
\left\{\left.g \in \Aut\left( M\left( \widetilde{A}\right)\right)\right| gM\left(\pi \right)\left(a \right)    = M\left(\pi \right)\left(a \right) \left( a\right) \quad \forall a \in M\left(A\right) \right\}\cong \\
\cong \left\{\left.g \in \Aut\left( \widetilde{A}\right)\right| g \pi\left(a \right)    =\pi\left(a \right) \left( a\right) \quad \forall a \in A \right\}.
\eean
It follows that $\left(M\left(A \right)  ,M\left( \widetilde{A}\right) , G, M\left( \pi\right) \right)$ is a noncommutative finite-fold pre-covering (cf. Definition \ref{fin_pre_defn}). Now this corollary follows from \ref{cov_mult_fin_cor} one.
\end{proof}
\begin{remark}
From the Corollary \ref{cov_mult_fin_c_cor} and the equation \eqref{frame_rem_eqn} it follows that if  $\left(A, \widetilde{A}, G, \pi \right)$ is a {noncommutative finite-fold covering with unitization} (cf. Definition \ref{fin_unitization_defn}) then one has

\be\label{cov_mult_fin_c_eqn}
\begin{split}
\exists \left\{ \widetilde{a}_1, ..., \widetilde{a}_n, \widetilde{b}_1, ..., \widetilde{b}_n\right\} \subset   M\left(  \widetilde{A}\right) \quad \forall \widetilde{a} \in M\left( \widetilde{A}\right) \quad \widetilde{a}=	\sum_{j=1}^{n}\widetilde{b}_j a_j \\
\text{where} \quad \forall j\in \left\{1,...,n\right\}\quad a_j\bydef  \sum_{g \in G} g\left(\widetilde{a}^*_j\widetilde{a} \right)\in M\left(\pi \right) \left(  M\left( A\right) \right) .
\end{split}
\ee
\end{remark}
\section{Induced representations}\label{induced_repr_fin_sec}

   \begin{defn}\label{induced_repr_fin_defn}
   	Let $\left(A, \widetilde{A}, G, \pi\right)$ be a finite-fold noncommutative quasi-covering (cf. Definition \ref{fin_quasi_defn}) such that $\widetilde{A}_A\bydef \widetilde A$ is of the $C^*$-Hilbert $A$-module with the given by the equation \eqref{finite_hilb_mod_prod_eqn} product. If $\rho: A \to B\left( \H\right)$ is a representation and  $\widetilde{\rho}=\widetilde{A}-\Ind^A_{\widetilde{A}}\rho: \widetilde{A}\to B\left(\widetilde{\H} \right) $ is given by \eqref{induced_representation_eqn}, i.e. $\widetilde{\rho}$ is the induced representation (cf. Definition \ref{induced_representation_defn}) then we say that $\widetilde{\rho}$ \textit{is induced by the pair} $\left(\rho,\left(A, \widetilde{A}, G, \pi\right)  \right)$.  
   \end{defn}
  \begin{empt}\label{induced_repr_constr_empt}
	If $X=\widetilde{A}\otimes_A \H$ is the algebraic tensor product then according to \eqref{hilb_prod_eqn} there is a sesquilinear $\C$-valued product $\left(\cdot, \cdot\right)_{X}$ on $X$  given by
	\begin{equation}\label{induced_prod_equ}
	\left(\widetilde{a} \otimes \xi, \widetilde{b} \otimes \eta \right)_{X}= 
	\left(\xi, \left\langle \widetilde{a}, \widetilde{b} \right\rangle_{{A}}	 \eta\right)_{\H}
	\end{equation}
	where $ \left(\cdot, \cdot\right)_{\H}$ means the Hilbert space product on $\H$, and $\left\langle \cdot, \cdot \right\rangle_{\widetilde{A}}$ is given by the Definition \ref{hilbert_product_defn}. So $X$ is a pre-Hilbert space. There is a natural map $\widetilde{A} \times \left( \widetilde{A}\otimes_A \H \right)\to \widetilde{A}\otimes_A \H$ given by
	\be\label{ind_act_form_eqn}
	\begin{split}
		\widetilde{A} \times \left( \widetilde{A}\otimes_A \H \right)\to \widetilde{A}\otimes_A \H,\\
		\left( \widetilde{a}, \widetilde{b} \otimes \xi\right)  \mapsto \widetilde{a}\widetilde{b}\otimes \xi.
	\end{split}
	\ee
	The space $\widetilde{\H}$ of the representation  $\widetilde{\rho}=\widetilde{A}-\Ind^A_{\widetilde{A}}\rho: \widetilde{A}\to B\left(\widetilde{\H} \right)$ (cf. \eqref{induced_representation_eqn}) is the Hilbert norm completion of the pre-Hilbert space $X=\widetilde{A}\otimes_A \H$ and the action of $\widetilde{A}$ on $\widetilde{   \H}$ is the completion of the action on $\widetilde{A}\otimes_A \H$ given by \eqref{ind_act_form_eqn}.
\end{empt}

   \begin{lem}\label{induced_faithful_lem}
   	If $A \hookto B\left(\H \right) $ is faithful (cf. Definition \ref{faithful_representation_defn}) then $\widetilde{\rho}: \widetilde{A} \to B\left( \widetilde{\H} \right)$ is faithful. 
   \end{lem}
 \begin{proof}
	If $\widetilde{a} \in \widetilde{A}\setminus \{0\}$ is a nonzero element then
	$$
	a =\left\langle \widetilde{a}~\widetilde{a}^*, \widetilde{a}~\widetilde{a}^*\right\rangle_{{A}} = \sum_{g \in G}g\left(\widetilde{a}^*\widetilde{a}~\widetilde{a}~\widetilde{a}^* \right) \in A_+\setminus\{0\} 
	$$
	is a nonzero positive element. There is $\xi \in \H$ such that $\left( a\xi, \xi\right)_{\H} > 0$.
	However
	\bean
	\left(a \xi, \xi\right)_{\H} = \left(\xi, \left\langle \widetilde{a}~\widetilde{a}^*, \widetilde{a}~\widetilde{a}^*\right\rangle_{{A}} \xi \right)_{\H} = \left( \widetilde{a}\widetilde{a}^*\otimes {\xi}, \widetilde{a}\widetilde{a}^*\otimes {\xi}\right)_{\widetilde{\H}}=\left( \widetilde{a}\widetilde{\xi}, \widetilde{a}\widetilde{\xi}\right)_{\widetilde{\H}}> 0.
	\eean
	where $\widetilde{\xi}\bydef \widetilde{a}^*\otimes \xi \in \widetilde{A}\otimes_A \H \subset \widetilde{\H}$. Hence $\widetilde{a}\widetilde{\xi} \neq 0$.
\end{proof}
\begin{lemma}\label{induced_nondegenerate_lem}	
	If a representation  $\rho: A \hookto B\left(\H \right) $ is nondegenerate then a representation  $\widetilde{\rho}: \widetilde{A} \to B\left( \widetilde{\H} \right)$ is also nondegenerate (cf. Definition \ref{nondegenerate_repr_defn}). 
\end{lemma}
\begin{proof}
	For any $\widetilde \xi \in \widetilde \H$ and any $\eps > 0$ there are $\left\{\widetilde a_1,..., \widetilde a_n\right\}\in \widetilde A$ and $\left\{\xi_1, ..., \xi_n\right\}$ such that
\be\label{nondeg1_eqn}
\left\|\widetilde a_1\otimes \xi_1 + ... + \widetilde a_n\otimes \xi_n - \widetilde \xi\right\| < \frac{\left\| \widetilde \xi\right\|}{3},
\ee	
because $\widetilde A \otimes_A \H$ is dense in $\widetilde \H$. From \eqref{nondeg1_eqn} and the triangle identity it follows that
\be\label{nondeg2_eqn}
\left\|\widetilde a_1\otimes \xi_1 + ... + \widetilde a_n\otimes \xi_n\right\| >\frac{2\left\| \widetilde \xi\right\|}{3},
\ee	
If $\left\{\widetilde u_\la\right\}_{\la \in \La}\subset \widetilde A$ is an approximate unit for $\widetilde A$ (cf. Definition \ref{approximate_unit_defn}) then one has
\be\label{nondeg3_eqn}
\forall j = 1,..., n\quad \lim_{\la\in\La} \widetilde u_\la \widetilde a_j = \widetilde a_j,
\ee	
From \eqref{nondeg3_eqn} it turns out that there is $\widetilde u \in \left\{\widetilde u_\la\right\}_{\la \in \La}$ with
for $\widetilde A$ (cf. Definition \ref{approximate_unit_defn}) then one has
\be\label{nondeg4_eqn}
\left\|\widetilde a_1\otimes \xi_1 + ... + \widetilde a_n\otimes \xi_n - u\left( \widetilde a_1\otimes \xi_1 + ... + \widetilde a_n\otimes \xi_n\right) \right\| < \frac{\left\| \widetilde \xi\right\|}{3},
\ee
Applying the triangle identity to both \eqref{nondeg1_eqn} and \eqref{nondeg4_eqn} one has
\be\label{nondeg5_eqn}
\left\|\widetilde u\left( \widetilde a_1\otimes \xi_1 + ... + \widetilde a_n\otimes \xi_n\right) \right\| > \frac{\left\| \widetilde \xi\right\|}{3},
\ee
From \eqref{nondeg1_eqn} and taking into account $\left\|u \right\|\le 1$ we conclude that 
\be\label{nondeg6_eqn}
\left\|\widetilde u\left( \widetilde a_1\otimes \xi_1 + ... + \widetilde a_n\otimes \xi_n\right)- \widetilde u \widetilde \xi  \right\| < \frac{\left\| \widetilde \xi\right\|}{3},
\ee
The equations \eqref{nondeg5_eqn} and \eqref{nondeg6_eqn} are compatible if and only if $\widetilde u \widetilde \xi\neq 0$, i.e.   the representation $\widetilde{\rho}: \widetilde{A} \to B\left( \widetilde{\H} \right)$ is nondegenerate (cf. Definition \ref{nondegenerate_repr_defn}).
\end{proof}
\begin{empt}\label{g_act_induced_empt}
   	Let $\left(A, \widetilde{A}, G, \pi\right)$ be a  noncommutative finite-fold covering, and let $\rho: A \to B\left(\H \right)$ be a faithful non-degenerate representation, and let  $\widetilde{\rho}: \widetilde{A} \to B\left( \widetilde{\H} \right)$ be induced by the pair $\left(\rho,\left(A, \widetilde{A}, G, \pi\right)  \right)$ (cf. Definition \ref{induced_repr_fin_defn}). There is the natural action of $G\times\widetilde{\H}\to \widetilde{\H}$ which comes from a map
  \be\label{induced_g_act_eqn}
  \begin{split}
  	G\times \left( \widetilde A\otimes_A\H\right)\to\widetilde A\otimes_A\H ,\\ 
\left( 	g,  \left( \widetilde{a} \otimes \xi\right)\right)  \mapsto \left( g\widetilde{a} \right) \otimes \xi; \quad \widetilde{a} \in \widetilde{A}, \quad g \in G, \quad \xi \in \H. 
  \end{split}
      \ee
   	There is the natural orthogonal inclusion 
   	\be\label{hilb_fin_inc_eqn}
   	\H \hookto \widetilde{\H}
   	\ee
   	 induced by inclusions
   	$$
   	A \subset\widetilde{A}, ~~ A \otimes_A \H \subset\widetilde{A} \otimes_A \H.
   	$$
   	The action $G \times \widetilde \H\to \widetilde \H$ is equivariant, i.e.
   	\be\label{induced_equiv_eqn}
   	g \left(\widetilde \rho\left(\widetilde a \right) \widetilde \xi \right) =  \widetilde \rho\left(g\widetilde a \right) \left( g\widetilde \xi \right).
   	\ee
   	From $A = \widetilde A^G$ it follows that
   	\be\label{ind_h_g_eqn}
   	\H = \widetilde \H^G.
   	\ee
   	If $\widetilde{A}$ is an unital $C^*$-algebra then the inclusion \eqref{hilb_fin_inc_eqn} is given by
   	\be\label{hilb_fin_inc_map_eqn}
   	\begin{split}
   	\varphi:		\H \hookto \widetilde{\H},\\
   		\xi \mapsto 1_{\widetilde{A}} \otimes \xi
   	\end{split}
   	\ee
   	where $1_{\widetilde{A}} \otimes \xi\in \widetilde{A} \otimes_A \H$ is regarded as element of $\widetilde{\H}$. 
   The inclusion \eqref{hilb_fin_inc_map_eqn} is not isometric. From 
   $$
   \left\langle 1_{\widetilde{A}}, 1_{\widetilde{A}} \right\rangle = \sum_{g \in G\left(\left.\widetilde{A}~\right|A \right)} g1^2_{\widetilde{A}} = \left|G\left(\left.\widetilde{A}~\right|A \right) \right| 1_{A}
   $$
   it turns out
   \bean
   \left(\xi, \eta \right)_{\H}= \frac{1}{\left|G\left(\left.\widetilde{A}~\right|A \right) \right|}\left( 	\varphi\left(\xi \right), \varphi\left(\eta \right)\right)_{\widetilde{\H}}; ~~ \forall \xi, \eta \in \H 
   \eean
   	Action of $g\in G\left({\widetilde{A}}~|~A \right)$ on $\widetilde{A}$ can be defined by representation as $ \widetilde{\rho}\left( g \widetilde{a}\right)  = g \widetilde{\rho}\left( \widetilde{a}\right) g^{-1}$, i.e.
 $$
   	(g\widetilde{a}) \xi = g\left(\widetilde{a} \left( g^{-1}\xi \right)  \right);~ \forall \xi \in \widetilde{\H}.
   	$$
 \end{empt}
   
   	

\begin{lemma}\label{comact_repr_lem}
   	Let $\left(A, \widetilde{A}, G, \pi\right)$ be a  noncommutative finite-fold covering, let $\rho: A \to B\left(\H \right)$ be a faithful non-degenerate representation, and let  $\widetilde{\rho}: \widetilde{A} \to B\left( \widetilde{\H} \right)$ be induced by the pair $\left(\rho,\left(A, \widetilde{A}, G, \pi\right)  \right)$ (cf. Definition \ref{induced_repr_fin_defn}). If $\widetilde{A}_A$ is regarded as a $C^*$-Hilbert $A$-module (cf. Lemma \ref{fin_hilbert_mod_lem}) then there is a natural representation
   	$$
 \phi_\K  :	\K\left(\widetilde{A}_A \right) \to  B\left( \widetilde{\H} \right).
   	$$
\end{lemma}
\begin{proof}
If $\widetilde a, \widetilde b \in \widetilde A$ then form \eqref{finite_hilb_mod_prod_eqn} it follows that
$$
	\left\langle \widetilde a, \widetilde b \right\rangle_{{A}} =\sum_{g \in G} g\left(\widetilde a^* \widetilde b\right). 
$$
So for any $\widetilde c \in \widetilde A_A$ one has
$\widetilde b\left\rangle \right\langle\widetilde    a ~\left(\widetilde c \right) = \widetilde b\sum_{g \in G} g\left(\widetilde a^* \widetilde c\right)$
   	We define
$$
\forall \widetilde \xi \in \widetilde \H \quad  \phi_\K\left( \widetilde b\left\rangle \right\langle \widetilde   a\right)\left( \widetilde\xi \right)\bydef   \widetilde b\sum_{g \in G} g\left(\widetilde a^* \widetilde \xi\right) \in \widetilde\H. 	
$$
\end{proof}
\begin{definition}\label{state_cov_11_defn}
	Let	$\left(A, \widetilde{A}, G, \pi \right)$ be  a noncommutative finite-fold  quasi-covering. A state $\widetilde{\tau}$ on $\widetilde{A}$ is said to be \textit {invariant}  if $\widetilde{\tau}\left( g\widetilde{a}\right)=\widetilde{\tau}\left( \widetilde{a}\right)$ for all $g \in G$ and $\widetilde{a}\in \widetilde{A}$.
\end{definition}

\begin{lemma}\label{state_cov_11_lem}
	If	$\left(A, \widetilde{A}, G, \pi \right)$ is  a noncommutative finite-fold  quasi-covering then there is a one-to-one correspondence between positive functionals on $A$ and invariant positive functionals on $\widetilde{A}$.
\end{lemma}
\begin{proof}
	If $\tau: A \to \C$ is a positive functional on $A$ then there is the invariant positive functional
	\be\label{state_down_eqn}
	\begin{split}
		\widetilde{\tau}:  \widetilde{A}\to \C,\\
		\widetilde{a}\mapsto \frac{1}{\left|G\right|}\sum_{g \in G} \tau\left(g\widetilde{a} \right),
	\end{split}
	\ee
	i.e. there is a given by $\tau \to \widetilde\tau$ map $\phi$ from the set of positive functionals on $A$ to the set of invariant positive functionals on  $\widetilde{A}$. 
	If $\widetilde{\tau}: \widetilde{A}\to \C$ is an invariant positive functional  then there is positive  $\widetilde{a}\in \widetilde{A}_+$ such that $\widetilde{\tau}\left(\widetilde{a} \right)\neq 0$. So one has
	$$
	\frac{1}{\left|G\right|}\sum_{g \in G} \widetilde{\tau}\left(g\widetilde{a} \right)\neq 0.
	$$
	It follows that the restriction $\tau \bydef \widetilde{\tau}|_A: A\to \C$ is not trivial.  So there is the given by $\widetilde{\tau}\to \tau$ map $\varphi$ from the set of invariant states on $\widetilde{A}$ to states on $A$. The maps $\phi$ and $\varphi$ are mutually inverse.
\end{proof}

\begin{proposition} \label{spectrum_covering_finite_prop}
	If $\left(A, \widetilde{A}, G, \pi \right)$ is  a noncommutative finite-fold  quasi-covering (cf. Definition \ref{fin_quasi_defn})  then following conditions hold:
	\begin{itemize}
		\item [(i)] 
		If $\widetilde{\mathcal X}\bydef {\widetilde{A}}\hat~$ is the spectrum of $\widetilde{A}$ then there is the natural action $G \times \widetilde{\mathcal X} \to \widetilde{\mathcal X}$.

		\item[(ii)]  If $\mathcal X = \hat A$ is the spectrum of $A$ then there is the natural surjective continuous map 
		\be\label{spectrum_cov_eqn}
\begin{split}
		p:\widetilde{\mathcal X} \to {\mathcal X};\\
		\rep^{\widetilde{A}}_{ \widetilde{x}} \mapsto \left.\rep^{\widetilde{A}}_{ \widetilde{x}}\right|_A \cong \rep^A_{p\left(\widetilde{x} \right) }
	\end{split}
\ee
		(cf. Proposition \ref{sur_prop}, and \eqref{rep_x_eqn}), which is $G$-invariant, i.e.
		\be\label{cov_spectrum_inv_eqn}
		p \circ g = p;~~ \forall g \in G,
		\ee
		\item[(iii)] The map $p$ naturally induces the homeomorphism
	\be\label{spectrum_homeo_eqn}
	p_G:\widetilde{\mathcal X}/G \xrightarrow{\approx} {\mathcal X}.
	\ee	
		
	\end{itemize}
\end{proposition}
\begin{proof}
	(i)
	If $\widetilde{\rho}: \widetilde{A} \to B\left(\widetilde{\H} \right)$ is an irreducible representation then  for any $g \in G$ there is an irreducible representation
	\bean
	\widetilde{\rho}_g: \widetilde{A} \to   B\left(\widetilde{\H} \right),\\
	\widetilde{\rho}_g\left(\widetilde{a} \right) = \widetilde{\rho}\left(g\widetilde{a} \right); ~~ \forall  \widetilde{a} \in \widetilde{A}.
	\eean 
	\\
	The map $\widetilde{\rho} \mapsto \widetilde{\rho}_g$ gives a homeomorphism of $\widetilde{\mathcal X}$. Such homeomorphisms give an action $G \times \widetilde{\mathcal X} \to \widetilde{\mathcal X}$, such that $g\widetilde{\rho}= \widetilde{\rho}_g$ for any $g \in G$.\\
	(ii) Let $\widetilde\tau : \widetilde A \to \C$ be a pure state on $\widetilde A$, and let $\tau \bydef \left.\widetilde\tau\right|_{A}: A \to\C$ be its restriction. 
		 If $\widetilde{a} \in \widetilde{A}_+$ is a positive element such that $\widetilde{\tau}\left(\widetilde{a} \right) \neq 0$ then $a = \sum_{g \in G} g \widetilde{a} \in A$ is a positive element such that $\widetilde{\tau}\left(a \right) \neq 0$. It turns out that the restriction $\tau\bydef \left. \widetilde{\tau}\right|_{A}: A \to \C$ is not trivial. Since an inclusion $A \hookto \widetilde A$ is unital in the sense of the Definition \ref{principal_non_defn}  one has $\tau\left(1_{A^\sim } \right) = 1$, i.e. $\tau$ is a state. Let  $q:A \to \C$ be a positive functional such that $q \le \tau$. From the Lemma \ref{state_cov_11_lem} $q$ corresponds to a $G$-invariant positive functional 
		 \bean
		 \widetilde q: \widetilde A \to \C,\\
		 \widetilde a \mapsto \frac{1}{\left|G \right| }q\left( \sum_{	g \in G}g \widetilde a\right). 
		 \eean 
		 Similarly $\tau$ corresponds to $G$-invariant positive state
	 \bean
\widetilde \tau^\oplus: \widetilde A \to \C,\\
\widetilde a \mapsto \frac{1}{\left|G \right| } \tau\left( \sum_{	g \in G}g \widetilde a\right)= \frac{1}{\left|G \right| }\widetilde \sum_{	g \in G}\widetilde \tau\left( g \widetilde a\right). 
\eean 
If 
$
G'' \bydef \left\{g \in G| g \widetilde \tau = \widetilde \tau \right\}, 
$
and $G' = G/G''$ then there is a free transitive  action
$$
G' \times (G\widetilde \tau ) \to (G\widetilde \tau ) 
$$
From our construction it turns out that
\be\label{tau_oplus_eqn}
\begin{split}
\forall \widetilde a \in \widetilde A \quad \widetilde \tau^\oplus \left(\widetilde a\right) =\frac{1}{\left|G' \right| } \sum_{	g \in G'}\widetilde \tau\left( g \widetilde a\right),\\
\widetilde \tau^\oplus = \frac{1}{\left|G' \right| } \sum_{	g \in G'}g \widetilde \tau.
\end{split}
\ee
If $\pi_{\widetilde \tau^\oplus} : \widetilde A \to B\left( \H_{\widetilde \tau^\oplus}\right)$ is an associated with  $\widetilde \tau^\oplus$ representation  (cf. Definition \ref{gns_defn}) then from the Lemma \ref{disj_repr_lem} it follows that $$\H_{\widetilde \tau^\oplus} = \bigoplus_{g \in G'} \H_{g\widetilde \tau}$$ where
$\H_{g\widetilde \tau}$ is a space of the associated with $g\widetilde \tau$ representation $\pi_{g\widetilde \tau} :  \widetilde A \to B\left(\H_{g\widetilde \tau}\right)$. From $\widetilde q \le  \widetilde \tau^\oplus$ and the Lemma \ref{domi_mult_repr_lem} it follows that 
if $\pi_{\widetilde q} : A \to B\left( \H_{\widetilde q}\right)$  is an associated with  ${\widetilde q}$ representation  (cf. Definition \ref{gns_defn})   then  $\pi_{\widetilde q} : A \to B\left( \H_{\widetilde q}\right)$  is spatially 	equivalent to a subrepresentation of $\pi_{\tau^\oplus} : A \to B\left( \H_{\tau^\oplus}\right)$ (cf. Definition \ref{subrepr_defn}, and the Corollary \ref{domin_rep_cor}). From the Definition \ref{subrepr_defn} it follows that $\H_{\widetilde q}$ is a closed $\widetilde A$-invariant subspace of $\bigoplus_{g \in G'} \H_{g\widetilde \tau}$. Taking into account that $\pi_{g\widetilde \tau}$ is an irreducible for any $g \in G'$ we conclude that there is a subset $H \subset G'$ such that $\H_{\widetilde q} = \bigoplus_{g \in H} \H_{g\widetilde \tau}$. From the Theorem \ref{state_repr_thm} and the Lemma \ref{domi_mult_repr_lem} one can deduce that 
$\widetilde q = \sum_{g \in H} b_g \left(g \widetilde{\tau} \right)$ where  $b_g$ is a positive real number for any $g \in H$. So one has $\widetilde q = \sum_{g \in G'} t_g \left(g \widetilde{\tau} \right)$ where  $t_g$ is a non negative real number for any $g \in G'$. On the other hand the state $\widetilde q$ is $G$-invariant, it is possible if and only if there is $t \in \R$ such that $t_g = t$ for each $g\in G'$. It follows that $\widetilde q =t \widetilde \tau^\oplus$  and $q = t\tau$. So the state $\tau$ is pure (cf. Definition \ref{ps_defn}) and  $\widetilde\tau \mapsto \tau$ is 	the  given by \eqref{spectrum_cov_eqn} map $p$.  Let us prove that the map $p$ is continuous.
 From the Theorem \ref{jtop_thm} and the Definition \ref{jtop_defn}	any closed set $\mathcal U \subset  \mathcal X$ corresponds to a closed two-sided ideal $I  \subset A$ such that
	\bean
	\mathcal U  = \left\{ \left.x \in \mathcal X~\right|~  \rep_x\left( I\right)= \{0\}  \right\}.
	\eean 
	The set $p^{-1}\left(\mathcal U \right)$ corresponds to the closed two-sided ideal
	\be\label{ideal_preimage_eqn}
	\widetilde{I}_I= \bigcap_{\substack{\widetilde{x} \in\widetilde{\mathcal X}\\\rep_{\widetilde{x}}\left( I\right) =0}}\ker \rep_{\widetilde{x}},
	\ee 
	so the preimage of   $\mathcal U$, is closed, it turns out that $p$ is continuous. From the Proposition \ref{state_prop} it turns out that $p$ is a surjective map. From the condition (b) of the Definition \ref{fin_pre_defn} it turns out that $ga=a$ for any $a \in A$ and $g \in G$, so one has $\left(g\widetilde{\rho} \right)\left( a\right)= \widetilde{\rho}\left(a \right)$ for any $\widetilde{\rho} \in \widetilde{\mathcal X}$. It follows that $p \circ g = p;~~ \forall g \in G$.
	From the Lemma \ref{state_cov_11_lem} it turns out that $p$ is a bijective map.  \\
	(iii) Any pure state $\phi : A \to \C$ uniquely defines the given by \eqref{state_down_eqn} invariant state $\widetilde{\phi}$. Otherwise $\widetilde{\phi}$  uniquely depends on the $G$-orbit of any pure state $\widetilde{\psi}:\widetilde A \to \C$ such that $\left.\widetilde{\psi}\right|_A = \phi$. So the map $p_G:\widetilde{\mathcal X}/G \xrightarrow{\approx} {\mathcal X}$ is bijective. Open subsets of $\widetilde{\mathcal X}/G$ correspond to open subsets $\widetilde\sU\subset \widetilde\sY$ such that $G\widetilde\sU=\widetilde\sU$. Otherwise such open subset corresponds to the closed two-sided  $G$-invariant  ideal $\widetilde I$, i.e. $G\widetilde I=\widetilde I$. If $I\bydef A\cap \pi^{-1}\left( \widetilde I\right) $ then $\widetilde I$ is the generated by $\pi\left( I\right)$ ideal. There are mutually inverse maps between closed two sided ideals of $A$  an closed two-sided  $G$-invariant  ideals $\widetilde A$ given by
	\bean
	I \mapsto \text{ the generated by }\pi\left( I\right) \text{ ideal of } \widetilde A,\\
	\widetilde I \mapsto  A\cap \pi^{-1}\left( \widetilde I\right)\subset A. 
	\eean
	These maps establishes the one to one correspondence between open subsets of both $\widetilde{\mathcal X}/G $ and $\sX$, i.e. $p_G:\widetilde{\mathcal X}/G \xrightarrow{\approx} {\mathcal X}$ is a homeomorphism.
\end{proof}

\begin{remark}
	From the equation \eqref{spectrum_cov_eqn} it turns out that there is the natural $*$-isomorphism
	\be\label{rep_x_mor_eqn}
	\begin{split}
	p^{\widetilde{A}}_{\widetilde{x}} :\rep^{\widetilde{A}}_{ \widetilde{x}}\left(A \right) \xrightarrow{\cong} \rep^{A}_{p\left( \widetilde{x}\right) }\left(A \right),\\
	 \rep^{\widetilde{A}}_{ \widetilde{x}}\left(a \right) \mapsto \rep^{A}_{p\left( \widetilde{x}\right)}\left(a \right); \quad \forall a \in A.
	\end{split}
	\ee
Taking into account that the atomic representation (cf. Definition \ref{atomic_repr_defn}) is faithful and nondegenerate (cf. Definitions \ref{faithful_representation_defn}, \ref{nondegenerate_repr_defn}), from the Definition \ref{multiplier_el_defn} and the Lemma \ref{fin_multilier_lem} it follows that thee is a natural injective  $*$-homomorphism
	\be\label{rep_x_mor_inj_eqn}
\begin{split}
		\rep^{M\left(A\right)}_{p\left( \widetilde{x}\right) }\left(M\left(A\right) \right)\hookto \rep^{M\left( \widetilde A\right) }_{ \widetilde{x} }\left(M\left( \widetilde A \right) \right).
\end{split}
\ee
	
\end{remark}
\begin{definition}\label{fin_free_defn}
We say that a finite-fold quasi-covering $\left(A, \widetilde{A}, G, \pi \right)$ is \textit{free} if the given by the Proposition \ref{spectrum_covering_finite_prop}   action of $G\times \widetilde\sX \to \widetilde\sX$ is free.
\end{definition}   
\begin{lemma}\label{sub_pure_lem}
Let $\left(A, \widetilde{A}, G, \pi \right)$ be a noncommutative finite-fold  quasi-covering. Let $\widetilde \tau: \widetilde{A}\to \C$ and $\tau \bydef \left. \widetilde{\tau}\right|_A: A \to \C$ be states. If $\rho: A \to \C$ is a pure state such that $\rho\preceq\tau$ (cf. Definition \ref{orth_repr_defn}), then there is a pure state $\widetilde \rho : \widetilde{A}\to \C$ such that 
\bean
\widetilde\rho\preceq\widetilde\tau,\\
\rho = \left. \widetilde{\rho}\right|_A.
\eean 
\end{lemma}
\begin{proof}
From the Proposition \ref{state_prop} it follows that there is a pure state $\widetilde \rho' : \widetilde{A}\to \C$ such that $\rho = \left. \widetilde{\rho}'\right|_A$. On the other hand both $\tau$ and $\rho$ correspond to  given by
\bean
\widetilde\tau^\oplus \bydef \frac{1}{\left| G\right| }\sum_{	g \in G} g\widetilde{\tau},\\
\widetilde\rho^\oplus \bydef \frac{1}{\left| G\right| }\sum_{	g \in G} g\widetilde{\rho}'
\eean 
$G$-invariant states respectively (cf. Definition \ref{state_cov_11_defn} and Lemma \ref{state_cov_11_lem}). From $\rho\preceq\tau$ it follows that $\widetilde\rho^\oplus\preceq\widetilde\tau^\oplus$. Otherwise from the Corollary \ref{domin_rep_cor} it follows that $\widetilde\rho'\preceq\widetilde\rho^\oplus$, so one has $\widetilde\rho'\preceq\widetilde\tau^\oplus$. Since $g\widetilde\rho'$ for all $g \in G$ is an irreducible  then there are two alternatives
\bean
g\widetilde\rho'\perp \widetilde \tau,\\
g\widetilde\rho'\preceq \widetilde \tau.
\eean
If $g\widetilde\rho'\perp \widetilde \tau$  for each $g \in G$ then $\widetilde\rho^\oplus\perp \widetilde\tau^\oplus$. But it is impossible since $\widetilde\rho'\preceq\widetilde\rho^\oplus$. So there is $g \in G$ such that $g\widetilde\rho'\perp \widetilde \tau$.
 If $\widetilde \rho \bydef g^{-1}\widetilde \rho'$ then from $\widetilde \rho' \preceq g \widetilde \tau$ one has
$$
\widetilde{\rho}\preceq\widetilde\tau
$$
 From $\widetilde{\rho}= g^{-1}\widetilde{\rho}'$ it follows that
$$
\left. \widetilde{\rho}\right|_A= \left. \widetilde{\rho}'\right|_A= \rho.
$$
\end{proof}

\begin{empt}
Suppose that $\left(A, \widetilde{A}, G, \pi \right)$ is  a noncommutative finite-fold  quasi-covering.
There is  an isomorphism  $\widetilde{A} \approx A \bigoplus P$ of $A$-modules, so any $\widetilde{a} \in \End\left( \widetilde{A}\right)_A$ can be represented by  four $A$-module homomorphisms $\al : A \to A$, $\bt : P \to A$, $\ga : A \to P$, $\delta : P \to P$ such that $\widetilde{a}$ corresponds to a matrix
	\be\label{matr_unis_rep_eqn}
	M\left(\widetilde{a} \right) = 	\begin{pmatrix}
		\al\left(\widetilde{a} \right) & 	\bt\left(\widetilde{a} \right)\\
		\ga\left(\widetilde{a} \right) & 	\delta\left(\widetilde{a} \right).
	\end{pmatrix}
	\ee
Moreover one has
	\bean
	M\left(\widetilde{a}'\widetilde{a}'' \right)= M\left(\widetilde{a}' \right)M\left(\widetilde{a}'' \right)
	.
	\eean
\end{empt}
\begin{empt}\label{irred_decomp_empt}
	If $\widetilde \rho:\widetilde{A} \to B\left(\widetilde{ \H}\right)$ is an irreducible representation then 
	\be\label{th_eqn}
	\widetilde{ \H} = \rho\left( \widetilde{A} \right) \xi.
	\ee
	If
	\bea\label{hurewicz_eqn}
	\H \bydef \widetilde\rho\left(\al\left(A\right)\right) \xi.
	\eea
	then from the Proposition \ref{spectrum_covering_finite_prop} one has the irreducible representation 
	\bea\label{ah_eqn}
\rho: A \to B\left(\H \right); \quad \rho\left(a \right)\left(   \widetilde \rho\left( a'\right)\xi\right)  \bydef \widetilde \rho\left(aa'\right) \xi\quad \forall a, a' \in A.
\eea
	
	If $p_\H: \widetilde{ \H} \to \widetilde{ \H}$ is the projector onto $\H$ then
	\be\label{al_eqn}
	\begin{split}
\forall \widetilde a\in \widetilde A \quad	\widetilde 	\rho\left(\al\left( \widetilde{a}\right)\right) = p_\H \widetilde \rho\left( \widetilde{a}\right)p_\H,\\ \rho\left(\bt\left( \widetilde{a}\right)\right)= p_\H \widetilde \rho\left( \widetilde{a}\right)\left( 1-p_\H\right),\\
	\widetilde 	\rho\left(\ga\left( \widetilde{a}\right)\right) =\left( 1-p_\H\right) \widetilde \rho\left( \widetilde{a}\right)p_\H,\\ \rho\left(\delta\left( \widetilde{a}\right)\right)= \left( 1-p_\H\right)\widetilde  \rho\left( \widetilde{a}\right)\left( 1-p_\H\right),\\
	\widetilde 	\rho\left(\widetilde{a}\right)= 	
	\begin{pmatrix}
	\widetilde  \rho\left( 	\al\left(\widetilde{a} \right)\right)  & 	\widetilde  \rho\left( \bt\left(\widetilde{a} \right)\right) \\
	\widetilde  \rho\left( 	\ga\left(\widetilde{a} \right)\right)  & 	\widetilde  \rho\left( \delta\left(\widetilde{a} \right)\right)
	\end{pmatrix}.
	\end{split}
	\ee	
	
\end{empt}
 
\begin{definition}\label{reduced_alg_defn}
	Let $\left(A, \widetilde{A}, G, \pi \right)$ be  a noncommutative finite-fold  quasi-covering. A $C^*$-subalgebra
	\be\label{red_incl_eqn}
	\widetilde{A}_{\text{red}}\bydef	\left\{\left. \widetilde a \in \widetilde A\right| 	\rep_{ \widetilde{x}}\left(\widetilde{a} \right) \in  \rep_{ \widetilde{x}}\left(\pi\left( A\right) \right)\quad \forall \widetilde x \in \widetilde{\sX}\right\} \subset \widetilde A
	\ee
	is said to be a \textit{reduced algebra of} $\left(A, \widetilde{A}, G, \pi \right)$.
\end{definition}
\begin{rem}\label{red_alg_inc_rem}
	Clearly one has an inclusion
	\be\label{red_alg_inc_eqn}
	A\subset \widetilde{A}_{\text{red}}.
	\ee
\end{rem}
\begin{rem}\label{red_g_act_rem}
	From $G \widetilde{A}_{\text{red}}= \widetilde{A}_{\text{red}}$ it follows that  there is the natural action
	\be\label{red_g_act_eqn}
	G\times \widetilde{A}_{\text{red}}\to  \widetilde{A}_{\text{red}}.
	\ee
	Moreover from $\widetilde{A}^G = A$ and $A \subset \widetilde{A}_{\text{red}}$ it follows that
	\be\label{red_inv_eqn}
	\widetilde{A}_{\text{red}}^G = A.
	\ee
\end{rem}

\begin{definition}\label{fin_red_defn}
	A noncommutative finite-fold  quasi-covering  $\left(A, \widetilde{A}, G, \pi \right)$ (cf. Definition \ref{fin_quasi_defn}) is said to be \textit{reduced} if $\widetilde{A}=\widetilde{A}_{\text{red}}$ or equivalently if for any irreducible representation $\rho: \widetilde{A}\to B\left(\H \right)$ one has
	\be\label{fin_red_eqn}
	\rho\left( \widetilde A\right) = \rho\left(\pi\left( A\right) \right).
	\ee 
\end{definition}

\begin{remark}
	I do not know whether exists a noncommutative finite-fold  covering which is not reduced.
\end{remark}

\begin{remark}
	I do not know whether exists a noncommutative finite-fold  covering which is not reduced.
\end{remark}

\begin{lemma}\label{covering_atomic_lem} 
	If $\left(A, \widetilde{A}, G, \pi \right)$ is  a  reduced 
	noncommutative finite-fold  quasi-covering (cf. Definitions \ref{fin_quasi_defn} and \ref{fin_red_defn}), $\rho: A \to B\left(\H\right)$ is an atomic representation (cf. Definition \ref{atomic_repr_defn}) and $ \widetilde \rho :\widetilde A \to B\left(\widetilde \H\right)$ is induced by the pair $\left(\rho, \left(A, \widetilde{A}, G, \pi \right)\right)$ (cf. Definition \ref{induced_repr_fin_defn}), then the representation $\widetilde \rho$ is atomic.
\end{lemma}
\begin{proof}
	Let both $\sX$ and $\widetilde \sX$ be spectra of both $A$ and $\widetilde A$. If $x \in \sX$ and $\H_x$ is a space of an irreducible  representation $\rep_x: A \to B\left(\H_x\right)$ then $\H= \oplus_{x\in \sX} \H_x$ where $\oplus$ is a norm completion of an algebraic direct sum (cf. Definition \ref{atomic_repr_defn}). It turns out that $\widetilde \H =\oplus_{x\in \sX} \widetilde \H_x$ where $\widetilde \H_x \bydef \widetilde A\otimes \H_x$. From the equation \eqref{tau_oplus_eqn} it follows that
	$
	\widetilde \H_x = \bigoplus_{\widetilde x \in p^{-1}\left( x\right)} \widetilde \H_{\widetilde x}
	$
	where $p: \widetilde \sX \to \sX$ is given by the Proposition \ref{spectrum_covering_finite_prop} and $\widetilde\H_{\widetilde x}$ a space of an irreducible  representation $\rep_{\widetilde x}: \widetilde A \to B\left(\widetilde \H_{\widetilde x}\right)$. 
	Using the above construction one has
	$$
	\widetilde \H = \bigoplus_{x\in \sX} \widetilde A\otimes \H_x = \bigoplus_{x\in \sX}\left(\bigoplus_{\widetilde x \in p^{-1}\left( x\right)} \widetilde \H_{\widetilde x} \right) =  \bigoplus_{\widetilde x \in \widetilde \sX} \widetilde \H_{\widetilde x}
	$$
where $ \bigoplus_{\widetilde x \in \widetilde \sX}$ means the norm completion of the algebraic direct sum.	Above equation means that $\widetilde \H$ is a Hilbert  space of the atomic representation.
\end{proof}

\begin{empt}\label{hered_lift_repr_empt}
	Let  $\left(A, \widetilde{A}, G, \pi \right)$ be a {noncommutative finite-fold  quasi-covering} (cf. Definition \ref{fin_quasi_defn}), and let $\rho : A \hookto B\left(\H \right)$ be a faithful nondegenerate  representation (cf. Definitions \ref{faithful_representation_defn} and \ref{nondegenerate_repr_defn}). If $\widetilde \rho  : \widetilde{A}\to B\left( \widetilde{\H}\right)$ is induced by a pair $\left(\rho, \left(A, \widetilde{A}, G, \pi \right)\right)$ (cf. Definition \ref{induced_repr_fin_defn}) then $\widetilde \rho$ is faithful nondegenerate representation (cf. Lemmas \ref{induced_faithful_lem} and \ref{induced_nondegenerate_lem}). If $\left\{u_\la \right\}_{\la\in\La} \subset  B\left(\H \right)_+$ is an increasing  net of positive elements explained in \ref{hered_repr_p_empt} then one can define a net $\left\{\rho\left(u_\la \right) \right\}_{\la\in\La}\subset B\left(  \H\right) $ such that then there is a projector $p \in B\left( \H\right)$ such that  $p$ is a strong limit of $\left\{\rho\left(u_\la \right) \right\}_{\la\in\La}$ in $B\left(\H\right)$ and	\bean
\forall a \in 	B  p\quad \rho\left( a\right) p = \rho\left( a\right).
	\eean
	Moreover  a natural representation
	\bean
	B \hookto B\left(p \H \right) 
	\eean
	is faithful and nondegenerate (cf. Corollary \ref{hered_representation_cor}). 
	From the Lemma \ref{fin_multilier_lem} it follows that there is a natural inclusion $M\left(\pi \right) : M\left(A \right) \hookto M\left(\widetilde A \right)$. 
	If  $\widetilde B\subset \widetilde A$ is a hereditary lift  $\left(A, \widetilde{A}, G, \pi \right)$-lift of $B$
	(cf. Definition \ref{hereditary_lift_defn}) then there is a net  $\left\{\widetilde\rho \left( M\left(\pi \right)\left(u_\la \right) \right)  \right\}_{\la\in \La}\subset B\left(\widetilde \H\right)$ 	such that $\left\{M\left(\pi \right)\left(u_\la \right)  \right\}_{\la\in\La}\subset M\left(\widetilde A \right) \cap M\left( \widetilde B\right)$. Similarly to \ref{hered_repr_p_empt} 
	there is  a strong limit $\widetilde p \bydef s-\lim\widetilde \rho\left( M\left(  u_\la\right)  \right) \in B\left( \widetilde \H\right)$ such that   $\widetilde p$ lies in the strong closure  
	 $\widetilde\rho\left(\widetilde A \right)''$ of $\widetilde\rho\left(\widetilde A\right)$ in $B\left(\widetilde \H\right)$.	
Moreover  a natural representation
\be\label{hered_rept_eqn}
\widetilde B \hookto B\left(\widetilde p \widetilde \H \right) 
\ee
is faithful and nondegenerate (cf. Corollary \ref{hered_representation_cor}).
\end{empt}
\begin{empt}\label{spectrum_ff_p_empt}
	Let $\left(A, \widetilde{A}, G, \pi \right)$ be a {noncommutative finite-fold quasi-covering with unitization} (cf. Definition \ref{fin_unitization_defn}). If both $C$ and $\widetilde C$ are centers of $M\left(A \right)$ and $M\left( \widetilde{A}\right)$ respectively then from the Remark \ref{unital_rem} one can deduce that $M\left(\pi \right)\left( C \right)\subset \widetilde C$, so there is an injective inclusion $C\subset \widetilde C$  of commutative $C^*$-algebras. If for both $A$ and $\widetilde A$ a prime spectrum coincides with primitive one (cf. Definition \ref{primitive_prime_spectrum_defn}) then from the Theorem \ref{dauns_hofmann_thm} it follows that  both $C$ and $\widetilde C$  are classes of bounded continuous functions on primitive spectra. From the Proposition \ref{spectrum_covering_finite_prop} it follows that there is an action $G \times \widetilde C \to \widetilde C$ such that $\widetilde C^G = C$. Similarly to the Lemma \ref{fin_hilbert_mod_lem} $\widetilde C$ is a $C^*$-Hilbert $C$-module.
\end{empt}
\begin{lemma}\label{spectrum_ff_p_lem}
	Let  $\left(A, \widetilde{A}, G, \pi \right)$ be a {noncommutative finite-fold quasi-covering with unitization} (cf. Definition \ref{fin_unitization_defn}). If   both $C$ and $\widetilde C$ are centers of $M\left(A \right)$ and $M\left( \widetilde{A}\right)$ then  $\widetilde C$ is a finitely generated $C$-module.
\end{lemma}
\begin{proof}
	If $\widetilde C$ is not a finitely generated $C$-module  then from the Theorem \ref{fin_gen_thm} it follows that  there is an infinite algebraic sum
	$$
	\widetilde C = \sum_{\la\in \La} M_\la,
	$$
	of right $C$-modules, such that for any finite subset $\La_0 \subset \La$ one has
	$$
	\widetilde C \subsetneqq \sum_{\la\in \La_0} M_\la,
	$$
	One has an algebraic sum 
	$$
	M\left(\widetilde A \right)  = \sum_{\la\in \La} M\left( A\right)  M_\la
	$$
	of left  $M\left(A \right)$-nodules such that for any finite subset $\La_0 \subset \La$ one has
	$$
	M\left( \widetilde{A}\right) \neq \sum_{\la\in \La_0} M\left( A\right)  M_\la.
	$$
	From the Theorem \ref{fin_gen_thm} it follows that $M\left(\widetilde A \right)$ is not finitely generated $M\left( A\right)$-module. This fact contradicts with hypotheses of this Lemma.
\end{proof}
\begin{corollary}\label{spectrum_ff_p_cor}
	Let  $\left(A, \widetilde{A}, G, \pi \right)$ be a {noncommutative finite-fold quasi-covering with unitization} (cf. Definition \ref{fin_unitization_defn}). If  the spectrum $\sX$ of $A$ is Hausdorff then one has:
	\begin{enumerate}
		\item[(i)] the spectrum $\widetilde \sX$ of $\widetilde A$ is also Hausdorff,
		\item[(ii)] the given by the Proposition \ref{spectrum_covering_finite_prop} map $p: \widetilde \sX \to \sX$ is a transitive finite-fold covering,
		\item [(iii)] for any $\widetilde x \in \widetilde \sX$ there is an an equality of subgroups
		$$
		\left\{g \in G\left| g\widetilde x = \widetilde x\right.\right\}= \ker\left(  G \to \left\{\left.g \in \mathrm{Homeo}\left(\widetilde\sX \right)\right| \forall \widetilde x \in \widetilde\sX \quad p\left( \widetilde x\right)= p\left( g\widetilde x\right) \right\}\right)
		$$	
		of $G$.
	\end{enumerate}
\end{corollary}

\begin{proof}
	If $C$ is the center of $M\left( \widetilde A\right)$ then $C\cong C_b\left( \sX\right) = C\left(\bt \sX \right)$ where $\bt  \sX$ is the  {Stone-\v{C}ech compactification} of $\sX$ (cf. Definition  \ref{sc_comp_defn}). If $\widetilde C$ is the center of  $M\left( \widetilde{A}\right)$ then there is a compact set $\widetilde \sY$ such that $\widetilde C\cong C\left(\widetilde\sY \right)$. From the Lemma \ref{spectrum_ff_p_lem} it follows that $C\left(\widetilde\sY \right)$ is finitely generated $C\left(\bt \sX \right)$-module, and taking into account the Lemma \ref{pavlov_troisky_thm} the inclusion $ C\left(\bt \sX \right)\subset C\left(\widetilde\sY \right)$ corresponds to a finite-fold covering $M\left( p\right) : \widetilde\sY\to \bt \sX$. \\
	(i) If the spectrum  $\widetilde \sX$ is not Hausdorff then there are $\widetilde x', \widetilde x''\in \widetilde\sX$ such that  $\widetilde x'\neq \widetilde x''$ and
	$$
	\forall \widetilde f \in \widetilde C \quad \widetilde f\left(\widetilde x' \right) = \widetilde f\left(\widetilde x'' \right).
	$$
	If $\widetilde I =\left\{ \left.\widetilde f \in \widetilde C\right| \widetilde f\left(\widetilde x' \right) = \widetilde f\left(\widetilde x' \right) = 0  \right\}$ then there is $\widetilde y_0 \in \widetilde\sY$ such that
	$$
	\widetilde I= \widetilde I_{\widetilde y_0}\bydef \left\{ \left.\widetilde f \in \widetilde C\right| f\left(\widetilde y_0 \right) = 0  \right\}.
	$$
	From the Proposition \ref{spectrum_covering_finite_prop} it follows that $p\left( \widetilde x'\right) = p\left( \widetilde x''\right) =x_0 = M\left( p\right)\left(\widetilde y_0 \right)$.
	If $g \in G$ is such that $g \widetilde y_0 \neq \widetilde y_0$ then $\widetilde x'' \neq g\widetilde x'$, i.e. $G\widetilde x'\neq G\widetilde x''$, so from the Proposition \ref{spectrum_covering_finite_prop} it follows that $p\left( \widetilde x'\right) \neq p\left( \widetilde x''\right)$, i.e. there is a contradiction. So the spectrum $\widetilde \sX$ is Hausdorff.\\
	(ii) There is an inclusion $\widetilde \sX \to \widetilde \sY$
	$$
	\widetilde x \mapsto  \widetilde y \quad \Leftrightarrow \quad \left( \forall  \widetilde f \in \widetilde \sY \quad  \widetilde f\left(\widetilde x \right) = 0  \quad \Rrightarrow \quad \widetilde f\left(\widetilde y \right) = 0 \right).
	$$
	such that $p = M(p)|_{\widetilde x}$ so $p$ is a finite-fold covering. From the homeomorphism \eqref{spectrum_homeo_eqn}  it follows that the covering is transitive.\\
	(iii) Follows from Lemma \ref{top_trans_tp_cov_lem}.
\end{proof}

 \section{Strong Morita equivalence of finite-fold coverings}
 
 \paragraph{}
 The following result about algebraic Morita equivalence of Galois extensions is in fact rephrasing of described in \cite{auslander:galois,mont:hopf-morita} constructions. 
 Let $\left(A, \widetilde{A}, G\right)$ be an unital finite-fold noncommutative covering. Denote by $\widetilde{A} \rtimes G$ a \textit{crossed product}, i.e. $\widetilde{A} \rtimes G$ is a $C^*$-algebra which coincides with a set of maps from $G$ to $\widetilde{A}$ as a set, and operations on  $\widetilde{A}$ are given by
 \begin{equation}\label{crossed_prod_eqn}
 \begin{split}
 \left(a + b \right) \left( g \right)   = a \left( g \right)+ b \left( g \right),\\
 \left(a \cdot b \right) \left( g \right)   = \sum_{g'\in G} a\left( g' \right)\left( g'\left(  b\left(g'^{-1}g \right)\right)\right)  \quad \forall a, b \in \widetilde{A} \rtimes G,~~ \forall g \in G,\\
  a^*\left(g \right) =\left(a\left( g^{-1}\right)  \right)^*\quad \forall a \in \widetilde{A} \rtimes G,~~ \forall g \in G.
 \end{split}
 \end{equation}
 (cf. Equations \eqref{discr_cr_prod_op_eqn}). 
 Let us construct a Morita context $$\left( \widetilde{A} \rtimes G, A,~ _{\widetilde{A} \rtimes G}\widetilde{A}_A, ~ _A\widetilde{A}_{\widetilde{A} \rtimes G}, \widetilde{\tau}, \psi\right) $$ where both $_A\widetilde{A}_{\widetilde{A} \rtimes G}$ and $ _{\widetilde{A} \rtimes G}\widetilde{A}_A$ coincide with $\widetilde{A}$ as $\C$-spaces. The left  and right action of $\widetilde{A} \rtimes G$ on $\widetilde{A}$ is given by
 
 \begin{equation}\label{morita_ta_act_eqn}
 \begin{split}
  \forall \widetilde{a}\in \widetilde{A}\quad \forall a \in \widetilde{A} \rtimes G\quad a \widetilde{a} = \sum_{g \in G} a\left( g\right) \left( g\widetilde{a}\right)   , \quad
 \widetilde{a} a =\sum_{g \in G} g^{-1}\left(\widetilde{a} a\left( g\right)\right). \end{split}
 \end{equation}
The left  (resp. right) action of $A$ on $\widetilde{A}$ we define as left  (resp. right) multiplication by $A$.
 Denote by $\varphi: \widetilde{A} \otimes_{A}  \widetilde{A}  \to \widetilde{A} \rtimes G$ and $\psi: \widetilde{A} \otimes_{\widetilde{A} \rtimes G}  \widetilde{A} \to A$ bilinear maps such that
 \begin{equation}\label{mor_pp_fin_eqn}
 \begin{split}
  \forall ~\widetilde{a} ,\widetilde{b}\in \widetilde{A}, ~~ g \in G \quad \varphi \left(\widetilde{a} \otimes \widetilde{b} \right)\left( g\right) \bydef    \widetilde{a} \left(g\widetilde{b} \right) ,\\
 \psi \left(\widetilde{a} \otimes \widetilde{b} \right) \bydef \sum_{g \in G} g\left(\widetilde{a} \widetilde{b} \right).
 \end{split}
 \end{equation}
 From above equations it follows that 
 \begin{equation}\label{morita_g_eqn}
 \begin{split}
 \forall \widetilde{a},\widetilde{b},\widetilde{c}\in \widetilde{A}\quad\varphi \left(\widetilde{a} \otimes \widetilde{b} \right)\widetilde{c} = \sum_{g \in G} \left( \varphi\left(\widetilde{a} \otimes \widetilde{b} \right)\left( g\right)\right)  g \widetilde{c} =\\
 = \sum_{g \in \G} \widetilde{a} \left( g \widetilde{b} \right) \left( g \widetilde{c}\right) =  \widetilde{a}\sum_{g \in G}  g\left(  \widetilde{b}  \widetilde{c}\right) = \widetilde{a} \psi\left(\widetilde{b}\otimes \widetilde{c} \right),\\
 \widetilde{a}\varphi \left(\widetilde{b} \otimes \widetilde{c} \right) =
 \sum_{g \in G} g^{-1}\left( \widetilde{a} \widetilde{b} g \widetilde{c}\right) =
 \\
 = \left( \sum_{g \in G} g^{-1}\left( \widetilde{a} \widetilde{b} \right) \right)\widetilde{c}= \left( \sum_{g \in G} g\left( \widetilde{a} \widetilde{b} \right) \right)\widetilde{c}=\psi\left(\widetilde{a} \otimes \widetilde{b} \right)\widetilde{c},
 \end{split}
 \end{equation}
 i.e. $\varphi, \psi$ satisfy conditions \eqref{morita_ctx_eqn}, so $\left( \widetilde{A} \rtimes G, A,~ _{\widetilde{A} \rtimes G}\widetilde{A}_A, ~ _A\widetilde{A}_{\widetilde{A} \rtimes G}, \varphi, \psi\right) $ is a Morita context.  Taking into account that the $A$-module $\widetilde{A}_A$ is a finitely generated projective generator and Remark \ref{morita_rem}  one has a following lemma.
 \begin{lem}\label{morita_galois_lem}
 	If $\left(A, \widetilde{A}, G \right)$ is an unital noncommutative finite-fold  covering then  $$\left( \widetilde{A} \rtimes G, A,~ _{\widetilde{A} \rtimes G}\widetilde{A}_A, ~ _A\widetilde{A}_{\widetilde{A} \rtimes G}, \varphi, \psi\right)$$ is an algebraic Morita equivalence. 
 \end{lem}

 \begin{cor}\label{unital_cov_cor}
 	Let $\left(A, \widetilde{A}, G, \pi \right)$ be an unital noncommutative finite-fold  covering. Let us define a structure of Hilbert $\widetilde{A} \rtimes G-A$ bimodule on  $_{\widetilde{A} \rtimes G}\widetilde{A}_A$ given by following products
 	\begin{equation}\label{hilb_gal_eqn}
 	\begin{split}
 	\left\langle\widetilde a, \widetilde b \right\rangle_{\widetilde{A} \rtimes G} = \varphi\left(\widetilde a \otimes  \widetilde b^* \right),\\ 
 	\left\langle \widetilde a,\widetilde b \right\rangle_A = \psi\left(\widetilde a^* \otimes \widetilde b\right).\\
 	\end{split}
 	\end{equation}
 	Following conditions hold:
 	\begin{enumerate}
 		\item [(i)]
 		A bimodule  $_{\widetilde{A} \rtimes G}\widetilde{A}_A$ satisfies the associativity condition (a) of the Definition \ref{strong_morita_defn},
 		
 		\item[(ii)] 
 		\begin{equation*}
 		\begin{split}
 		\left\langle \widetilde{A} , \widetilde{A} \right\rangle_{\widetilde{A} \rtimes G} = \widetilde{A} \rtimes G,\\ 
 		\left\langle \widetilde{A} , \widetilde{A}  \right\rangle_A = A.\\
 		\end{split}
 		\end{equation*}
 	\end{enumerate}
 	It follows that $_{\widetilde{A} \rtimes G}\widetilde{A}_A$  is a ${\widetilde{A} \rtimes G}-A$  equivalence bimodule (cf. Definition \ref{strong_morita_defn}).
 \end{cor}
 \begin{proof}
 	(i) From \eqref{morita_g_eqn} it follows that products $ \left\langle -, - \right\rangle_{\widetilde{A} \rtimes G}$, $~~\left\langle -, - \right\rangle_A$ satisfy condition (a) of the Definition \ref{strong_morita_defn}.\\
 	(ii) From the Lemma  \ref{morita_galois_lem} and the Definition \ref{morita_ctx_defn} of algebraic Morita equivalence  it turns out
 	\begin{equation*}
 	\begin{split}
 	\varphi\left(  \widetilde{A} \otimes_{A}  \widetilde{A}\right)  =\left\langle \widetilde{A} , \widetilde{A} \right\rangle_{\widetilde{A} \rtimes G} = \widetilde{A} \rtimes G,\\
 	\psi\left( \widetilde{A} \otimes_{\widetilde{A} \rtimes G}  \widetilde{A}\right)=\left\langle \widetilde{A} , \widetilde{A} \right\rangle_{A}  = A.
 	\end{split}
 	\end{equation*}
 \end{proof}
 Let us consider the situation of the Corollary \ref{unital_cov_cor}. Denote by $e \in G$ the neutral element. The unity $1_{\widetilde{A} \rtimes G}$ of $\widetilde{A} \rtimes G$ is given by
 \begin{equation}\nonumber
 \begin{split}
 1_{\widetilde{A} \rtimes G}\left(g \right) = \left\{\begin{array}{c l}
 1_{\widetilde{A}}  & g = e \\
 0 & g \neq e
 \end{array}\right..
 \end{split}
 \end{equation}
 From the Lemma \ref{morita_galois_lem} it follows that there are $\widetilde{a}_1,...,\widetilde{a}_n, \widetilde{b}_1,..., \widetilde{b}_n \in \widetilde{A}$ such that
 \begin{equation}\label{1_g_eqn}
 \begin{split}
 1_{\widetilde{A} \rtimes G} = \varphi\left( \sum_{j=1}^n \widetilde{a}_j \otimes\widetilde{ b}_j^*\right) = \sum_{j=1}^n \left\langle \widetilde{a}_j , \widetilde{b}_j \right\rangle_{\widetilde{A} \rtimes G}.
 \end{split}
 \end{equation}
 From the above equation it turns out that for any $g \in G$
 \begin{equation}\label{ga_eqn}
 \varphi\left( \sum_{j=1}^n g\widetilde{a}_j \otimes \widetilde{b}_j^*\right)\left(g' \right) = \left(\sum_{j=1}^n \left\langle g\widetilde{a}_j , \widetilde{b}_j \right\rangle_{\widetilde{A} \rtimes G} \right) \left(g' \right) \left\{\begin{array}{c l}
 1_{\widetilde{A}}  & g' = g \\
 0 & g' \neq g
 \end{array}\right..
 \end{equation}

\begin{corollary}\label{tens_inv_cor}
$\left(A, \widetilde{A}, G, \pi \right)$ be an unital noncommutative finite-fold  covering. The given by  \ref{herm_functor_empt} functor $_A\widetilde{A}_{\widetilde{A}\rtimes G} \otimes_{\widetilde{A}\rtimes G}\left( -\right): \mathbf{Herm}_{\widetilde{A}\rtimes G} \to \mathbf{Herm}_A$ is equivalent to the functor of invariant submodule 
\bean
\left( -\right)^G : \mathbf{Herm}_{\widetilde{A}\rtimes G} \to \mathbf{Herm}_A;\\
\widetilde{   \H} \mapsto \widetilde{   \H}^G=\left\{\widetilde{\xi} \in \widetilde{   \H}~|~ g\widetilde{\xi} = \widetilde{\xi},~ \forall g \in G \right\}.
\eean  
\end{corollary}
\begin{proof} 
If $\xi = \widetilde a \otimes \widetilde\eta \in _A\widetilde{A}_{\widetilde A \rtimes G}\otimes_{A \rtimes G}\widetilde{   \H}$ then $\xi = 1_{\widetilde{A}}\otimes a' \widetilde \eta$ where $a'\in \widetilde A \rtimes G$ is such that 
	\begin{equation}\nonumber
		\begin{split}
			a'\left(g \right) = \left\{\begin{array}{c l}
				\widetilde{a}  & g = e \text{ unity of } G \\
				0 & g \neq e
			\end{array}\right..
		\end{split}
	\end{equation}
It turns out that  ${\xi}\in _A\widetilde{A}_{\widetilde A \rtimes G}\otimes_{A \rtimes G}\widetilde{   \H}$ equals to $1_{\widetilde{A}}\otimes \widetilde{\xi}'$. From $g1_{\widetilde{A}}\otimes \widetilde{\xi}'= 1_{\widetilde{A}}\otimes g\widetilde{\xi}'$ and $g1_{\widetilde{A}}=1_{\widetilde{A}}$ for any $g \in G$ it turns out that $1_{\widetilde{A}}\otimes g\widetilde{\xi}'= 1_{\widetilde{A}} \otimes\widetilde{\xi}'$ and $g\widetilde{\xi}'= \widetilde{\xi}'$ i.e. $\widetilde{\xi}'$ is $G$-invariant.
\end{proof}

 \begin{lem}\label{comp_cov_lem}
 	Let $\left(A, \widetilde{A}, G \right)$ be a noncommutative finite-fold  quasi-covering with unitization (cf. Definition \ref{fin_unitization_defn}). Let us define a structure of Hilbert $\widetilde{A} \rtimes G-A$ bimodule on  $_{\widetilde{A} \rtimes G}\widetilde{A}_A$ given by  products \eqref{hilb_gal_eqn}.
 	Following conditions hold:
 	\begin{enumerate}
 		\item [(i)]
 		A bimodule  $_{\widetilde{A} \rtimes G}\widetilde{A}_A$ satisfies associativity conditions (a) of the Definition \ref{strong_morita_defn}
 		\item[(ii)] 
 		\begin{equation*}
 		\begin{split}
 		\left\langle \widetilde{A} , \widetilde{A}  \right\rangle_A = A,\\
 		\left\langle \widetilde{A} , \widetilde{A} \right\rangle_{\widetilde{A} \rtimes G} = \widetilde{A} \rtimes G.\\ 
 		\end{split}
 		\end{equation*}
 	\end{enumerate}
 	It follows that $_{\widetilde{A} \rtimes G}\widetilde{A}_A$  is a ${\widetilde{A} \rtimes G}-A$  equivalence bimodule.
 \end{lem}
 \begin{proof}(i)
 	From the Definition \ref{fin_unitization_defn} it follows that there are  unital $C^*$-algebras $B$, $\widetilde{B}$  and inclusions 
 	$A \subset B$,  $\widetilde{A}\subset \widetilde{B}$ such that $A$ (resp. $B$) is an essential ideal of $\widetilde{A}$ (resp. $\widetilde{B}$). Moreover there is an unital  noncommutative finite-fold quasi-covering $\left(B ,\widetilde{B}, G \right)$. From the Corollary \ref{unital_cov_cor} it turns out that a bimodule  $_{\widetilde{B} \rtimes G}\widetilde{B}_B$ is a $\widetilde{B} \rtimes G-B$ equivalence bimodule. Both scalar products $ \left\langle -, - \right\rangle_{\widetilde{A} \rtimes G}$, $~~\left\langle -, - \right\rangle_A$ are restrictions of products $ \left\langle -, - \right\rangle_{\widetilde{B} \rtimes G}$, $~~\left\langle -, - \right\rangle_B$, so products $ \left\langle -, - \right\rangle_{\widetilde{A} \rtimes G}$, $~~\left\langle -, - \right\rangle_A$ satisfy to condition (a) of the Definition \ref{strong_morita_defn}.\\
 	(ii) From \eqref{1_g_eqn} it turns out that there are $\widetilde{a}_1,...,\widetilde{a}_n, \widetilde{b}_1,..., \widetilde{b}_n \in \widetilde{B}$ such that
 	\begin{equation*}
 	\begin{split}
 	1_{\widetilde{B} \rtimes G} = \varphi\left( \sum_{j=1}^n \widetilde{a}_j \otimes\widetilde{ b}_j^*\right) = \sum_{j=1}^n \left\langle \widetilde{a}_j , \widetilde{b}_j \right\rangle_{\widetilde{B} \rtimes G}.
 	\end{split}
 	\end{equation*}
 	If $a \in A_+$ is a positive element then there is $x \in A \subset \widetilde{A}$ such that $a = x^*x$ it follows that
 	$$
 	a = \frac{1}{\left|G\right|}\left\langle x, x \right\rangle_A.
 	$$
 	Otherwise $A$ is the $\C$-linear span of positive elements it turns out
 	$$
 	\left\langle \widetilde{A} , \widetilde{A}  \right\rangle_A = A.
 	$$
 	For any positive  $\widetilde{a} \in \widetilde{A}_+$ and any $g \in G$ denote by
 	$y^{\widetilde{a}}_g \in \widetilde{A} \rtimes G $ given by
 	\begin{equation*}
 	y^{\widetilde{a}}_g\left(g' \right) \left\{\begin{array}{c l}
 	\widetilde{a}  & g' = g \\
 	0 & g' \neq g
 	\end{array}\right..
 	\end{equation*}
 	There is $\widetilde{x} \in \widetilde{A}$ such that $\widetilde{x}\widetilde{x}^*= \widetilde{a}$. Clearly $\widetilde{x}\left( g\widetilde{a}_j\right) , \widetilde{x}\widetilde{b}_j \in \widetilde{A}$ and from \eqref{ga_eqn} it follows that
 	$$
 	y^{\widetilde{a}}_g = \sum_{j=1}^n \left\langle \widetilde x \left( g\widetilde{a}_j\right)  , \widetilde x\widetilde{b}_j \right\rangle_{\widetilde{A} \rtimes G}
 	$$
 	The algebra $\widetilde{A} \rtimes G$ is the $\C$-linear span of elements $y^{\widetilde{a}}_g$, so one has $\left\langle \widetilde{A} , \widetilde{A} \right\rangle_{\widetilde{A} \rtimes G} = \widetilde{A} \rtimes G$. 
 	
 \end{proof}
 
 \begin{thm}\label{finite_morita_main_theorem}
 	If $\left(A, \widetilde{A}, G \right)$ is a noncommutative finite-fold covering  then a Hilbert $ \widetilde{A} \rtimes G$-$A$  bimodule  $_{\widetilde{A} \rtimes G}\widetilde{A}_A$ is a $\widetilde{A} \rtimes G-A$ equivalence bimodule.
 \end{thm}
 \begin{proof}
 	From the Definition \ref{fin_defn} there is a family $\left\{\widetilde{A}_\la \subset \widetilde{A} \right\}_{\la \in \La}$ of  hereditary subalgebras of $\widetilde{A}$ such that a union $\bigcup_{\la \in \La} \widetilde{A}_\la$ is a dense subspace of $\widetilde A$ and for any $\la \in \La$ there is a natural  noncommutative finite-fold covering with unitization $\left(\widetilde{A}_\la \bigcap A,   \widetilde{A}_\la, G \right)$. From the Lemma \ref{comp_cov_lem} it turns out that   $\widetilde{A}_\la$ is a $\widetilde{A}_\la \rtimes G-\widetilde{A}_\la\bigcap A$ equivalence bimodule and
 	\begin{equation}\label{ia_eqn}
 	\begin{split}
 	\left\langle  \widetilde{A}_\la ,  \widetilde{A}_\la \right\rangle_{ \widetilde{A}_\la \rtimes G} =  \widetilde{A}_\la \rtimes G,\\ 
 	\left\langle  \widetilde{A}_\la ,  \widetilde{A}_\la\right\rangle_{ \widetilde{A}_\la \bigcap A} =  \widetilde{A}_\la
 	\bigcap A.\\
 	\end{split}
 	\end{equation}
 	for any $\la\in \La$.
 	The union $\bigcup_{\la \in \La} \widetilde{A}_\la$ is a dense subset of $\widetilde{A}$ and $\bigcup_{\la \in \La} \widetilde{A}_\la \bigcap A$ is a dense subset of $A$. So domains of products  $\left\langle -, - \right\rangle_{\widetilde{A}_\la \rtimes G}~$, $~~\left\langle -, - \right\rangle_{A_\la}$ can be extended up to $\widetilde{A} \times \widetilde{A}$ and resulting products satisfy to  (a) of the Definition \ref{strong_morita_defn}. From \eqref{ia_eqn} it turns out that $\left\langle \widetilde{A}  ,  \widetilde{A}  \right\rangle_{ \widetilde{A}  \rtimes G}$ (resp. $\left\langle  \widetilde{A}  ,  \widetilde{A}  \right\rangle_{ A}$) is a dense subset of  $\widetilde{A} \rtimes G$ (resp. $A$), i.e. the condition (b) of the Definition \ref{strong_morita_defn} holds.
 \end{proof}

   \section{Noncommutative unique path lifting}\label{noncomutative_path_liftting}
 
   \paragraph*{}
   This section is rarely used and may be skipped on first reading.
    Here we generalize  the path-lifting property given by the Definition \ref{top_path_lifting_defn}. Because noncommutative geometry has no points it has no a direct generalization of paths, but there is an implicit generalization. Let $\widetilde{\mathcal{X}} \to \mathcal{X}$ be a covering. The following diagram reflects the path lifting problem.
   \newline
   \hspace*{\fill}
   \begin{tikzpicture}
   \matrix (m) [matrix of math nodes,row sep=3em,column sep=4em,minimum width=2em] {
   	& \widetilde{\mathcal{X}} \\
   	I = [0,1]  & \mathcal{X} \\};
   \path[-stealth]
   (m-2-1.east|-m-2-2) edge node [below] {$f$}(m-2-2)
   (m-1-2) edge node [right] {$p$} (m-2-2)
   (m-2-1)      edge [dashed]  node[above] {$f'$} (m-1-2);
   \end{tikzpicture}
   \hspace{\fill}
   \newline
   However above diagram can be replaced with a having noncommutative  generalization diagram.
   \newline
   \hspace*{\fill}
   \begin{tikzpicture}
   \matrix (m) [matrix of math nodes,row sep=3em,column sep=4em,minimum width=2em] {
   	& \mathrm{Homeo}(\widetilde{\mathcal{X}}) \\
   	I = [0,1]  & \mathrm{Homeo}(\mathcal{X}) \\};
   \path[-stealth]
   (m-2-1.east|-m-2-2) edge node [below] {$f$}(m-2-2)
   (m-1-2) edge node [right] {$\al \mapsto \left.\a\right|_{\mathcal X}$} (m-2-2)
   (m-2-1)      edge [dashed]  node[above] {$f'$} (m-1-2);
   \end{tikzpicture}
   \hspace{\fill}
   \newline
   where $\mathrm{Homeo}$ means the group of homeomorphisms with {compact-open} topology (cf \ref{top_comp_open_empt}). Above diagram has the following noncommutative generalization   
  \newline
  \hspace*{\fill}
  \begin{tikzpicture}
  \matrix (m) [matrix of math nodes,row sep=3em,column sep=4em,minimum width=2em] {
  	& \Aut(\widetilde{A}) \\
  	I = [0,1]  & \Aut(A) \\};
  \path[-stealth]
  (m-2-1.east|-m-2-2) edge node [below] {$f$}(m-2-2)
  (m-1-2) edge node [right] {$\al \mapsto \left.\a\right|_{A}$} (m-2-2)
  (m-2-1)      edge [dashed]  node[above] {$f'$} (m-1-2);
  \end{tikzpicture}
  \hspace{\fill}
  \newline
   In the above diagram we require that $\al|_A \in \mathrm{Aut}\left(A\right)$ for any $\al \in \mathrm{Aut}\left(\widetilde{A}\right)$. The diagram means that $f'(t)|_{A} = f(t)$ for any $t \in [0,1]$. Noncommutative generalization of a locally compact space is a $C^*$-algebra, so  the generalization of $\mathrm{Homeo}(\mathcal{X})$ is the group  $\mathrm{Aut}(A)$ of *-automorphisms carries (at least) two different topologies making it into a topological group \cite{thomsem:ho_type_uhf}. The most important is {\it the topology of pointwise norm-convergence} based on the open sets
   \begin{equation*}
   	\left\{\left.\alpha \in \mathrm{Aut}(A) \ \right| \ \|\alpha(a)-a\| < 1 \right\}, \quad a \in A.
   \end{equation*}
   The other topology is the {\it uniform norm-topology} based on the open sets
   \begin{equation}\label{aut_norm_eqn}
   	\left\{\alpha \in \mathrm{Aut}(A) \ \left| \ \sup_{a \neq 0}\ \|a\|^{-1} \|\alpha(a)-a\| < \varepsilon \right. \right\}, \quad \varepsilon > 0
   \end{equation}
   which corresponds to following "norm"
   \begin{equation}\label{uniform_norm_topology_formula_eqn}
   	\|\alpha\|_{\text{Aut}} = \sup_{a \neq 0}\ \|a\|^{-1} \|\alpha(a)-a\| = \sup_{\|a\| =1}\  \|\alpha(a)-a\|.
   \end{equation}
   Above formula does not really means a norm because $\mathrm{Aut}\left(A\right)$ is not a vector space. Henceforth the uniform norm-topology will be considered only.

  \begin{definition}\label{upl_f_defn}
	Let $A \hookto \widetilde{A}$  be an inclusion of $C^*$-algebras. Let $f: [0,1] \to  \mathrm{Aut}\left({A}\right)$ be a continuous function such that $f\left( 0\right)= \Id_A$. If  there is a continuous map $\widetilde{f}: [0,1] \to  \mathrm{Aut}\left(\widetilde{A}\right)$ such that $\left.\widetilde{f}(t)\right|_A = f(t)$ for any $t \in \left[0,1\right]$ and $\widetilde{f}\left( 0\right)  = \Id_{\widetilde{A}}$ then we say that $\widetilde{f}$ is a $\pi$-\textit{lift} of $f$. If a lift $\widetilde{f}$ of $f$ is unique then  a map $\widetilde{f}: [0,1] \to  \mathrm{Aut}\left(\widetilde{A}\right)$ is said to be the \textit{unique $\pi$-lift} of $f$. If $\widetilde{f}$ is the unique $\pi$-lift of $f$ then we denote by 
	\be\label{upl_f_eqn}
	\lift_f \stackrel{\text{def}}{=}\widetilde{f}\left(1 \right).
	\ee
\end{definition}
\begin{empt}\label{lift_group_empt}
	Let $\pi: A \hookto \widetilde{A}$  be an inclusion of $C^*$-algebras. Let both $f_1: [0,1] \to  \mathrm{Aut}\left({A}\right)$, $f_2: [0,1] \to  \mathrm{Aut}\left({A}\right)$ have unique $\pi$-lifts $\widetilde{f}_1$ and  $\widetilde{f}_2$ respectively. If $f_1* f_2: [0,1] \to  \mathrm{Aut}\left({A}\right)$, $\widetilde{f}_1* \widetilde{f}_2: [0,1] \to  \mathrm{Aut}\left(\widetilde{{A}}\right)$ are given by
\bean
\left( f_1* f_2\right) \left(t \right) =\left\{
\begin{array}{c l}
f_2\left( 2t\right)  &0 \le t \le \frac{1}{2}  \\
\\
f_1\left(2\left(t - \frac{1}{2} \right)  \right) f_2\left(1\right) & \frac{1}{2} < t \le 1
\end{array}\right.,\\
\left( \widetilde{f}_1* \widetilde{f}_2\right) \left(t \right) =\left\{
\begin{array}{c l}
	\widetilde{f}_2\left( 2t\right)  &0 \le t \le \frac{1}{2}  \\
	\\ \\
	\widetilde{f}_1\left(2\left(t - \frac{1}{2} \right)  \right) \widetilde{f}_2\left(1\right) & \frac{1}{2} < t \le 1
\end{array}\right.
\eean
then $\widetilde{f}_1* \widetilde{f}_2$ is the unique $\pi$-lift of $f_1*f_2$.
Hence one has $\lift_{f_1* f_2}= \lift_{f_1} \lift_{f_2}$.
It turns out that the set of elements $\lift_f$ is a subgroup of $\Aut\left(\widetilde{A}\right)$.
\end{empt}

\begin{defn}\label{lift_group_defn}
		Let $\pi: A \hookto \widetilde{A}$  be an inclusion of $C^*$-algebras then  the group of elements $\lift_f$ (cf. \ref{lift_group_empt}) is said to be the $\pi$-\textit{lift group}.
\end{defn}
\begin{lemma}\label{lift_commutes_lem}
		Let $A \hookto \widetilde{A}$  be an inclusion of $C^*$-algebras. Let $f: [0,1] \to  \mathrm{Aut}\left({A}\right)$ be a continuous function such that $f\left( 0\right)= \Id_A$ and $\widetilde{f}: [0,1] \to  \mathrm{Aut}\left(\widetilde{A}\right)$ is  the \textit{unique $\pi$-lift} of $f$. If $g \in G=\left\{\left. g \in \Aut\left(\widetilde{A} \right)~\right|~ ga = a;~~\forall a \in A\right\}$ then one has $g\widetilde{f}\left( 1\right) = \widetilde{f}\left( 1\right)g$, i.e. $G$ commutes with the $\pi$-lift group.
\end{lemma}
\begin{proof}
If $\widetilde{f}': [0,1] \to  \mathrm{Aut}\left(\widetilde{A}\right)$ is given by $\widetilde{f}'(t)= g\widetilde{f}(t)g^{-1}$ then $\widetilde{f}'(t)|_A = \widetilde{f}(t)|_A$. From the Definition \ref{upl_f_defn} it turns out $\widetilde{f}'= \widetilde{f}$. It turns out $g\widetilde{f}(1)g^{-1}=\widetilde{f}'(1)=\widetilde{f}(1)$, i.e. one has $\widetilde{f}(1) = g\widetilde{f}(1)g^{-1}$ and, consequently  $\widetilde{f}(1)g = g\widetilde{f}(1)$.
\end{proof}
\begin{lemma}
		Let $\pi: A \hookto \widetilde{A}$  be an inclusion of $C^*$-algebras.
$f: [0,1] \to  \mathrm{Aut}\left({A}\right)$ be a continuous function such that $f\left( 0\right)= \Id_A$. If  there are two different $\pi$-lifts of $f$ then there is a nontrivial continuous  map: $\widetilde{f}:[0,1] \to  \mathrm{Aut}\left(\widetilde{A}\right)$ such that $\widetilde{f}|_A\left( t\right)  = \id_A$ for every $t \in \left[0,1\right]$.
\end{lemma}
\begin{proof}
	If $\widetilde{f}'$ and $\widetilde{f}''$ are different $\pi$-lifts of $f$ then the map $\widetilde{f}:[0,1] \to  \mathrm{Aut}\left(\widetilde{A}\right)$ given by $t \mapsto \widetilde{f}'\left( t\right)\left(  \widetilde{f}''\left( t\right) \right)^{-1}$ satisfies to conditions of this lemma.
\end{proof}
\begin{cor}\label{lift_unique_cor}
	Let $\pi: A \hookto \widetilde{A}$  be an inclusion of $C^*$-algebras, such that the group 
				\be\label{lift_group_eqn}
	G = \left\{ \left.g \in \Aut\left(\widetilde{A} \right)~\right|~ ga = a;~~\forall a \in A\right\}
	\ee
	is finite,  and let   $f: [0,1] \to  \mathrm{Aut}\left({A}\right)$ be such that  $f\left( 0\right)= \Id_A$. If $\widetilde{f}$ is a $\pi$-lift of $f$ then   $\widetilde{f}$ is the unique $\pi$-lift of $f$.
\end{cor}
\begin{proof}
	If there are two different lifts $\widetilde f_1, \widetilde f_2: [0,1] \to  \mathrm{Aut}\left(\widetilde{A}\right)$ of  $f$  then there is a nontrivial continuous  map: $\widetilde{f}\bydef \widetilde{f}_1 \widetilde{f}^{-1}_2:[0,1] \to  G$ such that $\widetilde{f}|_A\left( t\right)  = \id_A$. But it is impossible because the group is $G$ discrete.
\end{proof}
\begin{empt}\label{lift_unique_empt}
Consider the given by the Corollary \ref{lift_unique_cor} situation. If $f: [0,1] \to  \mathrm{Aut}\left({A}\right)$ is such that  $f\left( 0\right)= f\left( 1\right)=\Id_A$, and   $\widetilde{f}$ is the unique $\pi$-lift of $f$, such that $\widetilde{f}\left(0 \right)= \Id_{\widetilde{{A}}}$ then $\widetilde{f}\left(1 \right)\in G$, where $G$ is given by the equation \eqref{lift_group_eqn} group.
\end{empt}


   \section{Coverings of (local) operator spaces}
   \paragraph*{}
   Here we generalize the Section \ref{cov_fin_bas_sec}.
   \begin{definition}\label{operator_space_subunital_defn}
   	Let $Y$ be an unital operator space, and let $C^*_e\left( Y\right)$ be the $C^*$-envelope (cf. Definition \ref{c_env_sp_defn}) of $Y$. 
   	A \textit{sub}-\textit{unital operator space} is a pair $\left(X, Y \right)$ where  $X$ is a subspace of $Y$ such that $X = Y \oplus \C \cdot 1_{C^*_e\left( Y\right)}$ or $X = Y$.
   \end{definition}

   \begin{example}\label{sub_alg_exm} 
   	If $A$ is a $C^*$-algebra and $A^\sim$ is given by \eqref{unital_notation_eqn} then since $A$ is $C^*$-norm closed one has  {sub-unital} operator space $\left(A, A^\sim\right)$.
   \end{example}
   \begin{definition}\label{op_sum_space_defn}
   	If both  $\left(X, Y \right)$ and $\left(\widetilde X, \widetilde Y \right)$ are {sub-unital} operator spaces then a \textit{complete isometry} from $\left(X, Y \right)$ to $\left(\widetilde X, \widetilde Y \right)$  is a complete unital isometry $\pi_Y: Y \hookto \widetilde Y$ such that $\pi_Y\left(X \right) \subset \widetilde X$. We write
   	\be\label{op_sum_space_eqn}
   	\begin{split}
   	\left( \pi_X: X  \hookto \widetilde X\right) \stackrel{\text{def}}{=} \left.\pi_Y\right|_X;\\
   	\left(\pi_X, \pi_Y \right):  \left(X, Y \right) \hookto \left(\widetilde X, \widetilde Y \right).
   	\end{split}
   	\ee	
   \end{definition}

\begin{definition}\label{operator_space_envelope_defn}
   	If a pair   $\left(X, Y \right)$ is a {sub-unital} operator space then the $C^*$-\textit{envelope} of $\left(X, Y \right)$ is the $C^*$-algebra given by
   \be\label{op_su_env_eqn}
   \begin{split}
   	C^*_e\left( X, Y\right) = \\ = \cap \left\{ \left.A \subset 	C^*_e\left( Y\right)~\right| ~ A  \text{ is a } C^*\text{-subalgebra of } C^*_e\left( Y\right) \text{ AND } X \subset A \right\},
   \end{split}
   \ee
   i.e. $C^*_e\left( X, Y\right)$ is the $C^*$-algebra, generated by $X$.
   \end{definition}
   \begin{remark}\label{op_su_env_rem}
   If  a sub-unital  operator space $\left(A,A^\sim \right)$  is given by the  Example \ref{sub_alg_exm}
then the $C^*$-{envelope}  of $\left(A,A^\sim \right)$ is $A$, i.e. $C^*_e\left(A,A^\sim \right)\cong A$. It turns out from the following circumstances:
\begin{itemize}
	\item   $C_e\left(A^\sim \right) = A^\sim$, \item $A$ is $C^*$-norm closed.
\end{itemize}
   \end{remark}
   \begin{definition}\label{fin_op_defn}
   Let both $\left(X, Y \right)$ and $\left(\widetilde X, \widetilde Y \right)$ be {sub-unital} operator spaces and let $\left(\pi_X, \pi_Y \right):  \left(X, Y \right) \hookto \left(\widetilde X, \widetilde Y \right)$ be a {complete isometry} from  $\left(X, Y \right)$ to $\left(\widetilde X, \widetilde Y \right)$.
  Suppose that  following conditions hold:
   	\begin{enumerate}
 	\item[(a)] There is a  finite-fold  noncommutative covering $\left( C^*_e\left( X,  Y \right), C^*_e\left(\widetilde X, \widetilde Y \right), G, \rho \right)$ (cf. Definition \ref{fin_defn}) such that $\pi_Y = \left.\rho\right|_Y$ and $\pi_X = \left.\rho\right|_X$. 
 	\item[(b)] $\widetilde X\subset C^*_e\left(\widetilde X, \widetilde Y \right)$ is the maximal among  $\C$-linear subspaces $\widetilde X' \subset C^*_e\left(\widetilde X, \widetilde Y \right)$  such that $X = \widetilde X' \cap  C^*_e\left( X,  Y \right)$ and $G\widetilde X'= \widetilde X'$.
 	\end{enumerate}
We  say that  $\left(\left(X, Y \right),\left(\widetilde X, \widetilde Y \right), G, \left(\pi_X, \pi_Y \right) \right)$   is a  \textit{noncommutative finite-fold covering} of the {sub-unital} operator space $\left(X, Y\right)$.
\end{definition}
\begin{remark}\label{fin_op_g_rem}
The action $G \times C^*_e\left(\widetilde X, \widetilde Y \right) \to C^*_e\left(\widetilde X, \widetilde Y \right)$ is trivial on $C^*_e\left( X,  Y \right)$, so from (b) of the Definition \ref{fin_op_defn} it follows that there is the natural action $G\times \widetilde X \to \widetilde X$ which is trivial on $X$.
\end{remark}
\begin{definition}\label{op_su_r_space_defn}
	Let $Y$ be a  real operator space (cf. Definitions \ref{real_os_defn}), and let $X\subset Y$ is its subspace. Let both $\C X$ and $\C Y$ are complexifications  (cf. \ref{complexification_empt}) of $X$ and $Y$ respectively.
We say that a pair 	$\left(X, Y \right)$ is a \textit{sub-unital real operator space} if the  pair $\left(\C X, \C Y \right)$  is a sub-unital operator space (cf. Definition \ref{operator_space_subunital_defn}).
\end{definition}
   \begin{definition}\label{op_sum_r_space_defn}
	If both  $\left(X, Y \right)$ and $\left(\widetilde X, \widetilde Y \right)$ are {sub-unital} real operator spaces then \textit{complete isometry} from $\left(X, Y \right)$ to $\left(\widetilde X, \widetilde Y \right)$  is a complete unital isometry $\pi_Y: Y \hookto \widetilde Y$ such that $\pi_Y\left(X \right) \subset \widetilde X$. We write
	\be\label{op_sum_r_space_eqn}
	\begin{split}
		\left( \pi_X: X  \hookto \widetilde X\right) \stackrel{\text{def}}{=} \left.\pi_Y\right|_X;\\
		\left(\pi_X, \pi_Y \right):  \left(X, Y \right) \hookto \left(\widetilde X, \widetilde Y \right).
	\end{split}
	\ee	
\end{definition}
\begin{remark}\label{ros_iso_rem}
Let both  $\left(X, Y \right)$ and $\left(\widetilde X, \widetilde Y \right)$ are {sub-unital} real operator spaces, and let 	$\left(\pi_X, \pi_Y \right):  \left(X, Y \right) \hookto \left(\widetilde X, \widetilde Y \right)$ be a complete isometry.  If $\C X$, $\C Y$,  $\C \widetilde X$, $ \C \widetilde Y$ are complexifications of $X$, $Y$,  $\widetilde X$, $ \widetilde Y$ then from the Theorem \ref{ros_contr_iso_thm} it turns out that the map  $\left(\C\pi_X, \C\pi_Y \right):\left( \C X,\C Y\right)  \hookto  \left( \C \widetilde X, \C \widetilde Y\right)$ is a complete isometry.
\end{remark}
\begin{definition}\label{r_comp_defn}
	Under the hypotheses of the Remark \ref{ros_iso_rem} we say that $\left( \C X,\C Y\right)  \hookto  \left( \C \widetilde X, \C \widetilde Y\right)$ is the complexification of  $\left(\pi_X, \pi_Y \right):  \left(X, Y \right) \hookto \left(\widetilde X, \widetilde Y \right)$.
\end{definition}

\begin{definition}\label{fin_rop_defn}
	Let both $\left(X, Y \right)$ and $\left(\widetilde X, \widetilde Y \right)$ be {sub-unital} real operator spaces and let $\left(\pi_X, \pi_Y \right):  \left(X, Y \right) \hookto \left(\widetilde X, \widetilde Y \right)$ be a {complete isometry} from  $\left(X, Y \right)$ to $\left(\widetilde X, \widetilde Y \right)$. Let  $\left(\C\pi_X, \C\pi_Y \right):\left( \C X,\C Y\right)  \hookto  \left( \C \widetilde X, \C \widetilde Y\right)$ is complexification of  $\left(\pi_X, \pi_Y \right):  \left(X, Y \right) \hookto \left(\widetilde X, \widetilde Y \right)$. We say that   $\left(\left(X, Y \right),\left(\widetilde X, \widetilde Y \right), G, \left(\pi_X, \pi_Y \right) \right)$   is a  \textit{noncommutative finite-fold covering of the sub-unital real operator space} $\left(X, Y\right)$ if $$\left(\left(\C X, \C Y \right),\left(\C \widetilde X, \C \widetilde Y \right), G, \left(\C\pi_X, \C\pi_Y \right) \right)$$   is a  {noncommutative finite-fold covering} of the {sub-unital} operator space $\left(\C X, \C Y\right)$ (cf. Definition \ref{fin_op_defn}).
\end{definition}


\begin{empt}\label{loc_op_env_empt} 
Let $X$ be a local operator space (cf. Definition \ref{loc_op_sp_defn}) and let $\varphi : X \hookto B$ be a complete isometry into pro-$C^*$-algebra (cf. \ref{complete_loc_maps_empt} and the Remark \ref{loc_op_abs_conc_rem}). Let $b\left( X\right) \bydef b\left(B \right) \cap X$. If $C_e\left( b\left( X\right)\right)$  {envelope } $C^*$-{algebra} of $b\left( X\right)$ then there is the surjective *- homomorphism $ b\left(B \right)\to C_e\left( b\left( X\right)\right)$. Let $b\left( I\right) \bydef \ker\left(  b\left(B \right)\to C_e\left( b\left( X\right)\right)\right)$. If
$$
I \bydef \bigcap_{\substack{ J \text{ is an ideal of } B\\ J\cap b\left(B \right)\supset b\left( I\right)}}J
$$
then $A \bydef B/I$ is a pro-$C^*$-algebra and there is the complete isometry $X \hookto A$.
\end{empt}
\begin{definition}\label{loc_op_env_defn}
In the described in \ref{loc_op_env_empt} we say that $A$ is the \textit{envelope pro}-$C^*$-\textit{algebra} of $X$. We denote it by $C^*_e\left( X\right)$. 
\end{definition}
\begin{remark} Under the hypotheses \ref{loc_op_env_empt},  \ref{loc_op_env_defn}  there is the natural surjective $*$-homomorphism $B\to C^*_e\left( X\right)$. So $C^*_e\left( X\right)$ is unique up to $*$-isomorphism.
\end{remark}
\begin{definition}\label{bounded_loc_op_el_defn}
An element $a$ of local operator space $X$ is said to be \textit{bounded} if $a \in X \cap b\left(C^*_e\left(X\right) \right)$.  Denote by $b\left( X\right)$ the $\C$-space of bounded elements.  
\end{definition}

\begin{rem}\label{bounded_loc_op_el_rem}
For any local operator space $X$ the space $b\left( X\right)$ is an operator space and $ b\left(C^*_e\left(X\right)\right) $ is the $C^*$-envelope of $b\left(X\right)$. 
\end{rem}

\begin{remark}
Any complete isomerty $\phi: X \to Y$ of local operator spaces (cf. \ref{complete_loc_maps_empt}) induces an injective $*$-homomorphism $C^*_e\left(Y\right)\hookto C^*_e\left(X\right)$  of {envelope pro}-$C^*$-{algebras} (cf. Definition \ref{loc_op_env_defn}).
\end{remark}

  \begin{definition}\label{op_loc_su_space_defn}
	Let $Y$ be an unital local operator space, and let $C^*_e\left( Y\right)$ be the envelope pro-$C^*$-algebra (cf. Definition \ref{loc_op_env_defn}) of $Y$. 
	A \textit{sub-unital local operator space} is a pair $\left(X, Y \right)$ where  $X$ is a  subspace of $Y$ such that $Y = X \oplus \C \cdot 1_{C^*_e\left( Y\right)}$ or $X = Y$.
\end{definition}
 \begin{definition}\label{op_loc_sum_space_defn}
	If both  $\left(X, Y \right)$ and $\left(\widetilde X, \widetilde Y \right)$ are {sub-unital} local operator spaces then \textit{complete isometry} from $\left(X, Y \right)$ to $\left(\widetilde X, \widetilde Y \right)$  is a complete unital isometry $\pi_Y: Y \hookto \widetilde Y$ (cf. \ref{complete_loc_maps_empt}) such that $\pi_Y\left(X \right) \subset \widetilde X$. We write
	\be\label{op_loc_sum_space_eqn}
	\begin{split}
		\left( \pi_X: X  \hookto \widetilde X\right) \stackrel{\text{def}}{=} \left.\pi_Y\right|_X;\\
		\left(\pi_X, \pi_Y \right):  \left(X, Y \right) \hookto \left(\widetilde X, \widetilde Y \right).
	\end{split}
	\ee	
\end{definition}
\begin{definition}\label{op_loc_su_env_defn}
	If a pair   $\left(X, Y \right)$ is {sub-unital} local operator space then the pro-$C^*$-\textit{envelope} of $\left(X, Y \right)$ is the pro-$C^*$-algebra given by
	\be\label{op_loc_su_env_eqn}
	\begin{split}
		C^*_e\left( X, Y\right) = \\ = \cap \left\{ \left.A \subset 	C^*_e\left( Y\right)~\right| ~ A  \text{ is a } \text{pro-}C^*\text{-subalgebra of } C^*_e\left( Y\right) \text{ AND } X \subset A \right\}.
	\end{split}
	\ee
\end{definition}

   \begin{definition}\label{fin_loc_op_defn}
	Let both $\left(X, Y \right)$ and $\left(\widetilde X, \widetilde Y \right)$ be {sub-unital} local operator spaces and let $\left(\pi_X, \pi_Y \right):  \left(X, Y \right) \hookto \left(\widetilde X, \widetilde Y \right)$ be a {complete isometry} from  $\left(X, Y \right)$ to $\left(\widetilde X, \widetilde Y \right)$.
	If  the following conditions hold:
	\begin{enumerate}
		\item[(a)] There is a finite-fold  noncommutative covering $\left( C^*_e\left( X,  Y \right), C^*_e\left(\widetilde X, \widetilde Y \right), G, \rho \right)$ of pro-$C^*$-algebras (cf. Definition \ref{pro_fin_defn}) such that $\pi_Y = \left.\rho\right|_Y$ and $\pi_X = \left.\rho\right|_X$. 
		\item[(b)] If $\widetilde X' \subset C^*_e\left(\widetilde X, \widetilde Y \right)$ is a $\C$-linear space  such that $X = \widetilde X' \cap  C^*_e\left( X,  Y \right)$ and $G\widetilde X'= \widetilde X'$  then $\widetilde X' \subseteq \widetilde X$.
	\end{enumerate}
	Then we  say that  $\left(\left(X, Y \right),\left(\widetilde X, \widetilde Y \right), G, \left(\pi_X, \pi_Y \right) \right)$   is a  \textit{noncommutative finite-fold covering} of the {sub-unital} operator space $\left(X, Y\right)$.
\end{definition}

\begin{definition}\label{loc_op_re_sp_defn}
The \textit{real local operator space}  is an complete  $\R$-space with seminorms which satisfy to conditions M1, M2 of the Definition \ref{loc_op_sp_defn}.
\end{definition}
\begin{remark}\label{loc_op_re_sp_rem}
Similarly to the Remark \ref{loc_op_sp_rem} any abstract real operator local space $V$ is an inverse limit of real operator spaces.
\end{remark}
\begin{remark}
 Similarly to definition \ref{fin_loc_op_defn} one can define a finite-fold noncommutative covering for local real operator spaces.
\end{remark}

\section{Coverings of spectral triples}\label{triple_fin_cov}

\paragraph*{} 
 If $\left( \A, \H, D, J\right)$ is a spectral triple then from the Axioms \ref{regularity_axiom}, \ref{finiteness_axiom} it follows that:
 \begin{itemize}
 	\item the space  
 	$
 \H^\infty\bydef	\bigcap_{k\in\bN} \Dom D^k
 	$
 	is dense in $\H$;
 	\item one has
 	\bean
 \A\subset \L^\dagger \left( \H^\infty \right), \\
 D \in 	\L^\dagger \left( \H^\infty \right)
 	\eean
 	(cf. equation \eqref{l_dag_eqn}).
 \end{itemize}

\begin{empt}\label{triple_conn_lift_empt} 
Let $\left(\A, \widetilde \A, G\left( \left.\widetilde \A~\right| \A\right), \pi  \right)$ be a {noncommutative finite-fold covering of  bounded operator *-algebras} (cf. Definition \ref{fin_oa_defn}) such that  $\widetilde \A$ is a finitely generated projective $\A$-module. 
 If both $A$ and $\widetilde A$ are $C^*$-norm completions of $\A$ and  $\widetilde \A$ then from the Definition  \ref{fin_oa_defn} it follows that there is  a  unital noncommutative finite-fold covering $\left(A, \widetilde A, G\left( \left.\widetilde A\right| A\right), \widetilde \pi  \right)$ (cf. Definition \ref{fin_unital_defn}) such that $\pi= \widetilde\pi|_{\A}$. Similarly to the equation \eqref{finite_hilb_mod_prod_eqn} one has a sesquilinear product
 $\left\langle \cdot , \cdot  \right\rangle_{\widetilde A} : \widetilde A\times \widetilde A \to \A$ such that
\bean
 \forall \widetilde a, \widetilde b \in \widetilde A\quad \left\langle\widetilde a, \widetilde b \right\rangle_{\widetilde A}\bydef\frac{1}{\left| G\right| } \sum_{	g \in G}g\left(\widetilde a^* \widetilde b \right),\\
 \left\langle\widetilde \A, \widetilde \A \right\rangle_{\widetilde A}= \A.
\eean
   Denote by 
\bean
\rho : A \hookto B\left(\H \right),\\
\eean
a natural representations of $A$. If the natural representation $\widetilde \rho: \widetilde A \to B\left(\widetilde \H \right)$ is induced by a pair $\left(\rho, \left(A, \widetilde A, G\left( \left.\widetilde A\right| A\right), \widetilde \pi  \right) \right)$ (cf. Definition \ref{induced_repr_fin_defn}) then we have
\be\label{st_h_defn_eqn}
\widetilde \H \bydef \widetilde A\ox_A \H
\ee
because $ \widetilde A$ is a finitely generated $A$-module.  Moreover the  Hilbert space $\widetilde A \ox_A \H$ has  the following scalar product
\bean
\left(\widetilde a \ox \xi, \widetilde b \ox \eta \right)_{\widetilde{\H}}\bydef \left(\xi,  \left\langle a, b \right\rangle_{\widetilde A} \eta \right)_\H,
\eean
If operators $\widetilde J, \widetilde J^\7 : \widetilde{\H}\to \widetilde{\H}$ are given by
\bea\label{cc_eqn}
\widetilde J\left( \widetilde a \otimes \xi \right) \bydef \widetilde a^* \otimes J\xi,\\
\label{ccc_eqn}
\widetilde J^\7\left( \widetilde a \otimes \xi \right) \bydef \widetilde a^* \otimes J^\7\xi
\eea
then on has 
\bean
\widetilde J^\dagger\widetilde J= \widetilde J\widetilde J^\dagger = 1_{B\left( \widetilde \H\right) },\\
\left(\widetilde J\left( \widetilde a \otimes \xi \right), \widetilde J\left( \widetilde b \otimes \eta \right)\right)_{\widetilde\H}= \left( \widetilde a^*\ox J\xi ,  \widetilde b^*\ox J \eta \right)_{\widetilde\H}= \left(J\xi, \left\langle \widetilde a^*, \widetilde b^* \right\rangle J\eta \right)_{\H}  =
\\
= \left(\widetilde J \left( 1_{\widetilde A}\ox \xi \right), \widetilde J\left(  \left\langle \widetilde b, \widetilde a \right\rangle \otimes \eta\right)    \right)_{\widetilde\H} =  \left(\widetilde J^\dagger\widetilde J \left( 1_{\widetilde A}\ox \xi \right),  \left\langle \widetilde a, \widetilde b \right\rangle \otimes \eta   \right)_{\widetilde\H}=\\=\left(\xi,  \left\langle a, b \right\rangle_{\A} \eta \right)_\H= \left(\widetilde a \ox \xi, \widetilde b \ox \eta \right)_{\widetilde{\H}}.
\eean 
From the  equation\eqref{cc_eqn} and from
$$
\forall z \in \C \quad \widetilde J z \left(\widetilde a \ox \xi \right) = \widetilde J  \left(\widetilde a \ox z \xi \right)= \widetilde a^* \ox Jz\xi= \widetilde a^* \ox z^* J\xi= z^* \widetilde a^* \ox J\xi= z^* \widetilde J  \left(\widetilde a \ox \xi \right) 
$$
it follows that $\widetilde J$ is an antiunitary operator. Similarly one can prove that the operator $\widetilde J^\dagger$ is also  antiunitary.
If $\xi \in \H$ and $\widetilde a, \widetilde b, \widetilde c\in  \widetilde A$ then one has
\bean
\forall \widetilde a, \widetilde b, \widetilde c \quad \widetilde a \left(\widetilde J \widetilde b^* J^\7 \right) \left(\widetilde c \ox \xi \right)=\widetilde a \widetilde J \left( \widetilde b^*\widetilde c^*\ox \widetilde J^\7 \xi\right) = \widetilde a \widetilde c \widetilde b \otimes JJ^\7\xi= \widetilde a \widetilde c \widetilde b \otimes \xi,
\eean
\be\label{jbj_eqn}
\begin{split}
\left(\widetilde J \widetilde b^* J^\7 \right)\widetilde a\left(\widetilde c \ox \xi \right)=\left(\widetilde J \widetilde b^* J^\7 \right)\left(\widetilde a\widetilde c \ox \xi \right)= \widetilde J \widetilde b^*\left(  \widetilde c^* \widetilde a^*\otimes J^\7\xi \right) = \\= \widetilde J \left(   \widetilde b^*\widetilde c^* \widetilde a^*\otimes J^\7\xi \right)
=  \widetilde a \widetilde c \widetilde b \otimes JJ^\7\xi= \widetilde a \widetilde c \widetilde b \otimes\xi.
\end{split}
\ee
So a following condition holds
$$
\widetilde a \left(\widetilde J \widetilde b^* J^\7 \right)=\left(\widetilde J \widetilde b^* J^\7 \right)\widetilde a,
$$
or equivalently
\be\label{st_comn_wt_eqn}
\forall \widetilde a,\widetilde b \in \widetilde \A\quad \left[ \widetilde a,\widetilde J\widetilde b^*\widetilde J^\dagger\right] = 0,
\ee
i.e. one has a situation of the Definition \ref{df:spt-real_defn}. 
 Let $\Om^1_D$ be given by \eqref{om_d_eqn}, and let  
\be\label{nabla_cov_herm_eqn}
\nabla : \widetilde \A\to \widetilde \A\otimes_\A \Om^1_D.
\ee
be  a connection (cf. Definition \ref{connection_defn}).
 We suppose that $\nabla$ is $G$-equivariant, i.e.
\be\label{conn_equ_conn}
\forall g \in G \quad \nabla \widetilde a = \sum_{j=1}^n \widetilde a_j \otimes \om_j \quad \Rightarrow \quad \nabla \left( g\widetilde a\right)  = \sum_{j=1}^n g\widetilde a_j \otimes \om_j. 
\ee
Let us prove that there is an action
$$
\left( \widetilde \A\otimes_\A \Om^1_D\right)  \times \left(\widetilde \A\otimes_\A \H^\infty\right)\to \left(\widetilde \A\otimes_\A \H^\infty\right)
$$
such that for all $\widetilde a \otimes b\left[D, a\right]\in \widetilde \A\ox_\A \Om^1_D$ and $\widetilde c\otimes\xi \in\widetilde \A\ox_\A  \H^\infty$ one has
\be\label{om_act_eqn}
\begin{split}
 \left(\left(\widetilde a \otimes b\left[D, a\right]\right) , \left(\widetilde c\otimes\xi\right)\right)\mapsto \widetilde a \widetilde c \otimes b\left[D, a\right]\xi.
\end{split}
\ee
The  action \eqref{om_act_eqn} is correct if for all $c \in\A$ one has
\bean
\left( \widetilde a  c \otimes b\left[D, a\right]\right)\left(  \widetilde c\otimes\xi\right)   =\left( \widetilde a  \otimes cb\left[D, a\right]\right) \left(  \widetilde c\otimes\xi\right),
\eean
or, equivalently
\be\label{corr_conn_act_eqn}
\widetilde a c\widetilde c \otimes b\left[D, a\right]\xi = \widetilde a \widetilde c \otimes cb\left[D, a\right]\xi.
\ee
From \eqref{jbj_eqn} it follows that
\bean
\widetilde a c\widetilde c \otimes b\left[D, a\right]\xi= \left(\widetilde J \widetilde c^*\widetilde J^\7 \right) \left(\widetilde a c \otimes b\left[D, a\right]\xi \right)= \\ =\left(\widetilde J \widetilde c^*\widetilde J^\7 \right)\left(\widetilde a \otimes cb\left[D, a\right]\xi \right) 
=\widetilde J \widetilde c^* \widetilde a^* \otimes  J^\dagger cb\left[D, a\right]\xi =\\= \widetilde a \widetilde c \otimes J J^\dagger cb\left[D, a\right]\xi=\widetilde a \widetilde c \otimes cb\left[D, a\right]\xi.
\eean
So the condition \eqref{corr_conn_act_eqn} holds and the given by \eqref{om_act_eqn} action is correct. 
We define an operator $\widetilde D\in \L^\dagger \left( \widetilde \H^\infty\right) $  such that 
\be\label{wtd_eqn}
\begin{split}
\widetilde D \left(\widetilde a \otimes \xi  \right)\bydef \left(\nabla \widetilde a \right)\xi +  \widetilde a\otimes D\xi = \sum_{j=1}^n \widetilde a_j \otimes \om_j \xi+  \widetilde a\otimes D\xi \\ \text{where}\quad \nabla \widetilde a = \sum_{j=1}^n\widetilde a_j \otimes \om_j.
\end{split}
\ee
We also require that the operator $\widetilde D$ is self-adjoint, i.e.
$$
\forall \widetilde a \otimes \xi, \widetilde b \otimes \eta\in  \widetilde \A\otimes_\A \H^\infty\quad \left( \widetilde D\left( \widetilde a \otimes \xi \right), \widetilde b \otimes \eta \right)_{\widetilde \H}   = \left(  \widetilde a \otimes \xi ,\widetilde D\left( \widetilde b \otimes \eta \right)\right)_{\widetilde \H}. 
$$
If $\nabla \widetilde a = \sum_{j=1}^n\widetilde a_j \otimes \om'_j$ and $\nabla \widetilde b = \sum_{j=1}^m\widetilde b_j \otimes \om''_j$ then one has
\bean
\left( \widetilde D\left( \widetilde a \otimes \xi \right), \widetilde b \otimes \eta \right)_{\widetilde \H}   = \left( \sum_{j=1}^n \widetilde a_j \otimes \om'_j \xi + \widetilde a \otimes D \xi , \widetilde b \otimes \eta \right)_{\widetilde \H}=\\= \sum_{j=1}^n \left(\om'_j \xi,  \left\langle \widetilde a_j, \widetilde b \right\rangle_{\widetilde A} \eta \right)_\H  + \left(D\xi, \left\langle \widetilde a, \widetilde b \right\rangle_{\widetilde A} \eta \right)_\H, \\
\left(  \widetilde a \otimes \xi ,\widetilde D\left( \widetilde b \otimes \eta \right)\right)_{\widetilde \H}= \sum_{j=1}^m\left(   \xi ,\left\langle \widetilde a, \widetilde b_j \right\rangle_{\widetilde A} \om''_j \eta \right)_{\H}+ \left(\xi, \left\langle \widetilde a, \widetilde b \right\rangle_{\widetilde A} D \eta \right)_\H.
\eean
It follows that if 
\be\label{herma_conn_eqn}
\begin{split}
 \sum_{j=1}^n \left(\om'_j \xi,  \left\langle \widetilde a_j, \widetilde b \right\rangle_{\widetilde A} \eta \right)_\H  + \left(D\xi, \left\langle \widetilde a, \widetilde b \right\rangle_{\widetilde A} \eta \right)_\H =\\=\sum_{j=1}^m\left(   \xi ,\left\langle \widetilde a, \widetilde b_j \right\rangle_{\widetilde A} \om''_j \eta \right)_{\H}+ \left(\xi, \left\langle \widetilde a, \widetilde b \right\rangle_{\widetilde A} D \eta \right)_\H 
\end{split}
\ee
then the given by \eqref{wtd_eqn} operator $\widetilde{D}$ is self-adjoint.
\end{empt}
\begin{definition}\label{herma_conn_defn}
If a connection $\nabla$ satisfies to the equation \eqref{herma_conn_eqn} then we say that $\nabla$ is  \textit{Hermitian} (cf. Definition \ref{herm_conn_defn}).
\end{definition}
Henceforth we require that the connection  $\nabla$  is Hermitian. 
From the equations \eqref{nabla_cov_herm_eqn} and \eqref{wtd_eqn} it follows that
\be\label{dta_eqn}
\forall \widetilde{a}\in \widetilde{\A}\quad \exists \widetilde{a_1},..., \widetilde{a_n}\in \widetilde{\A}\quad \exists a_1, ..., a_n \in \A\quad \left[ \widetilde D, \widetilde a\right]= \sum_{j=1}^n \widetilde{a_j} \otimes \left[ D, a_j\right]
\ee
If  $\widetilde{b}, \widetilde{c}\in \widetilde{\A}$ and $\xi \in \H^\infty$ then for any $j = 1,..., n$  one has
\bean
\widetilde J\widetilde b^*\widetilde J^\dagger \left( \widetilde a_j \otimes \left[  D,  a_j\right]\right) \left(\widetilde c \ox \xi \right) = \widetilde J\widetilde b^*\widetilde J^\dagger\left(\widetilde a_j \widetilde c \ox \left[  D,  a_j\right]\xi \right) = \\=\widetilde J \widetilde b^*\widetilde c^*\widetilde a^*_j \ox J^\dagger \left[  D,  a_j\right]\xi
= \widetilde a_j\widetilde c \widetilde b \otimes JJ^\dagger \left[  D,  a_j\right]\xi= \widetilde a_j\widetilde c \widetilde b \otimes \left[  D,  a_j\right]\xi, \\
 \left( \widetilde a_j \otimes \left[  D,  a_j\right]\right)\widetilde J\widetilde b^*\widetilde J^\dagger \left(\widetilde c \ox \xi \right)= \left( \widetilde a_j \otimes \left[  D,  a_j\right]\right)\widetilde J\widetilde b^*  \left(\widetilde c^* \ox J^\dagger \xi \right)=\\ \left( \widetilde a_j \otimes \left[  D,  a_j\right]\right]\left( \widetilde c \widetilde b\ox JJ^\dagger \xi\right) =\widetilde a_j\widetilde c \widetilde b \otimes \left[  D,  a_j\right]=
 \widetilde J\widetilde b^*\widetilde J^\dagger \left( \widetilde a_j \otimes \left[  D,  a_j\right]\xi\right) \left(\widetilde c \ox \xi \right)
\eean
and taking into account \eqref{dta_eqn} we obtain the following equation
\bea
\label{first_order_wt_eqn}
\forall\widetilde a, \widetilde b \in \widetilde\A\quad \left[\left[ \widetilde D, \widetilde a\right] , \widetilde J\widetilde b^*\widetilde J^\dagger\right]=  0.
\eea
(cf. equation \eqref{fist_order_eqn})
\begin{exercise}\label{st_ax_exer}
Prove that  for all $\widetilde a, \widetilde b \in \widetilde{\A}$ a following condition hold.
\bea
 \label{first_order_wtb_eqn}
  \left[\widetilde a, \left[ \widetilde D, J\widetilde b^*\widetilde J^\dagger \right]\right]  = 0.
\eea
 (cf. equation \eqref{first_order_dual_eqn}).
\end{exercise}
From our construction it follows that
\be\label{dirac_om_fin_lift_eqn}
\nabla \widetilde a = \sum_{j = 1}^n \widetilde a_j \otimes \om_j \quad\Rightarrow \quad\left[\widetilde D, \widetilde a\right]=  \sum_{j = 1}^n \widetilde a_j \otimes \om_j\in \widetilde\A \otimes_\A \Om^1_D.
\ee
\begin{defn}\label{spectral_triple_fin_lift_defn} 	In the above situation the spectral  triple $\left( \widetilde{\A}, \widetilde{\H}, \widetilde{D}, \widetilde J\right)$  is said to be the  $\left(A, \widetilde{A}, G, \pi \right)$-\textit{lift} of $\left( {\A}, {\H}, {D}, J\right)$.  
\end{defn} 
\begin{definition}\label{ext_st_alg_defn}
	If $\A_D$ is a *-subalgebra of $\L^\dagger \left(\D \right)$ generated by a union $\A\cup \left\{D\right\}$ then we say that $\A_D$ is the \textit{extended algebra of the spectral triple} $\left( \A, \H, D\right)$.
\end{definition}

\begin{exercise}\label{str_ad_exer}
Let both $\A_D$ and  $\widetilde\A_{\widetilde D}$ are {extended algebras of the spectral triples} $\left( \A, \H, D\right)$ and $\left(\widetilde \A, \widetilde\H, \widetilde D\right)$. Prove that there is a {noncommutative finite-fold covering of $O^*$-algebras} 
$\left(\A_{D}, \widetilde \A_{\widetilde D}, G, \psi  \right)$ (cf. Definition \ref{fino*_defn}), such that $\pi = \psi|_\A$ and $\widetilde D = \psi\left(D \right)$.
\end{exercise}

\section{Finite noncommutative coverings and flat connections}\label{flat_sec}

 \subsection{Basic construction}\label{nc_flat_sec}
 \paragraph*{}
 Let $\left( \A, \H, D\right)$ be a spectral triple, let  $\left( \widetilde{\A}, \widetilde{\H}, \widetilde{D}\right)$ be the $\left(A, \widetilde{A}, G \right)$-lift of $\left( \A, \H, D \right)$ (cf. Definition \ref{spectral_triple_fin_lift_defn}). Let $V = \C^n$ and with left  action of $G$, i.e. there is a linear representation $\rho: G \to GL\left(\C,n\right)$. Let $\widetilde{\mathcal E} \bydef  \widetilde\A \otimes \C^{n} \approx \widetilde{\A}^n$ be a free module over $\widetilde{\A}$, so $\widetilde{\mathcal E}$ is a projective finitely generated $\A$-module (because $\widetilde{ \A}$ is a finitely generated projective $\A$-module). Let $\widetilde{\nabla} : \widetilde{\mathcal E} \to \widetilde{\mathcal E} \otimes_{\widetilde{\A}} \Om^1_{\widetilde{D}}$ be the trivial flat connection (cf. Definition \ref{triv_conn_defn} and Remark \ref{triv_conn_rem}). 
 From \eqref{dirac_om_fin_lift_eqn} it turns out that $\Om^1_{\widetilde{D}} = \widetilde{\A}\otimes_{\A}\Om^1_{D}$. therefore  the connection $\widetilde{\nabla} : \widetilde{\mathcal E} \to \widetilde{\mathcal E} \otimes_{\widetilde{\A}} \Om^1_{\widetilde{D}}$ can be regarded as a map $\nabla':\widetilde{\mathcal E} \to \widetilde{\mathcal E} \otimes_{\widetilde{\A}} \widetilde{\A} \otimes_{\A} \Om^1_{\widetilde{D}}= \widetilde{\mathcal E} \otimes_{\A}\Om^1_{{D}}$, i.e. one has a connection
 \be\nonumber
 \nabla':\widetilde{\mathcal E} \to  \widetilde{\mathcal E} \otimes_{\A}\Om^1_{{D}}.
 \ee
 From $\widetilde{\nabla} \circ \widetilde{\nabla} |_{\widetilde{\mathcal E}}=0$ it turns out that  $\nabla' \circ \nabla' |_{\widetilde{\mathcal E}}=0$, i.e. $\nabla'$ is flat. There is the action of $G$ on $\widetilde{\mathcal E}= \widetilde{\A} \otimes \C^n$ given by 
 \be
 g\left( \widetilde{a}\otimes  x\right)   = g \widetilde{a} \otimes g x; ~~ \forall g \in G,~ \widetilde{a} \in \widetilde{\A}, ~ x \in \C^n.
 \ee
 Denote by
 \be
 \mathcal E = \widetilde{\mathcal E}^G = \left\{\widetilde{\xi} \in  \widetilde{\mathcal E}~|~ G\widetilde{\xi} = \widetilde{\xi}\right\}
 \ee
 Clearly $\mathcal E$ is an $\A$-$\A$-bimodule. For any $\widetilde{\xi} \in \widetilde{\mathcal E}$ there is the unique decomposition
 \be
 \begin{split}
 	\widetilde{\xi} = \xi + \xi_\perp, \\
 	\xi = \frac{1}{\left|G\right|}\sum_{g \in G} g \widetilde{\xi},\\
 	\xi_\perp = \widetilde{\xi} - \xi. 
 \end{split}
 \ee
 From the above decomposition it turns out the direct sum $\widetilde{\mathcal E} = \widetilde{\mathcal E}^G \bigoplus {\mathcal E}_\perp$ of $\A$-modules. So  $\mathcal E = \widetilde{\mathcal E}^G$ is a projective finitely generated $\A$-module, it follows that there is an idempotent $e \in \End_{\A}{\widetilde{\mathcal E}}$ such that $\mathcal E = e \widetilde{\mathcal E}$. The Proposition \ref{conn_prop} gives the canonical connection 
 \be\label{nc_flat_conn}
 \nabla : \mathcal E \to \mathcal E \otimes_{\A} \Om^1_D
 \ee
 which is defined by the connection $\nabla':\widetilde{\mathcal E} \to  \widetilde{\mathcal E} \otimes_{\A}\Om^1_{{D}}$ and the idempotent $e$.
 
 	\begin{lemma}\label{flat_lem}
 	If $\nabla:  \mathcal E \to \mathcal E \otimes_{\A} \Om^1$ is the trivial connection and $e \in \End_{\A}\left( \E\right)$ is an idempotent then the given by \eqref{idem_conn} connection  
 	\bean
 	\nabla_e: e\mathcal E \to e\mathcal E \otimes \Om^1;\\
 	\xi \mapsto \left(e \otimes 1 \right) \nabla \xi	
 \eean
 is flat.
 \end{lemma} 
 \begin{proof}
 	
 	If we consider the given by the equation  \eqref{eop_into_eqn} map
 	 		\bean
 	\nabla:  \mathcal E \otimes_{\A} \Om \to \mathcal E \otimes_{\A} \Om, \\	\nabla\left(\xi \otimes \om \right) \bydef \nabla\left(\xi \right) \om + \xi \otimes d\om; \quad\forall \xi \in \mathcal E, ~ \om \in \Om. 
 	\eean
 	or equivalently 
 	\bean
 	\nabla \bydef \Id_{\mathcal E} \otimes d
 	\eean
 	then one has
 	$$
 	\nabla_e\bydef	\left(e \otimes 1 \right)\left(\Id_{\mathcal E} \otimes d \right).
 	$$
 	From
 	$$
 	\left(e \otimes 1 \right)\left(\Id_{\mathcal E} \otimes d \right) \circ \left(e \otimes 1 \right)\left(\Id_{\mathcal E} \otimes d \right) = e \otimes d^2 =e \otimes 0 = 0	
 	$$
 	it turns out that $\nabla_e \circ \nabla_e = 0$, i.e. $\nabla_e$ is flat.
 \end{proof}
 \begin{remark}
 	The notion of the trivial connection is an algebraic version of geometrical canonical connection explained in the Appendix \ref{geom_flat_subsec}.
 \end{remark}
   From the Lemma \ref{flat_lem} it turns out that the given by \eqref{nc_flat_conn} connection $\nabla$ is flat.
 \begin{definition}
 	We say that  $\nabla$ is a \textit{flat connection induced by} noncommutative covering $\left(A, \widetilde{A}, G\right)$ and the	linear representation $\rho: G \to GL\left(\C,n\right)$, or we say the $\nabla$ \textit{comes from the representation} $\rho: G \to GL\left(\C,n\right)$.
 \end{definition}
 \break

 \subsection{Mapping between geometric and algebraic constructions}
 
 \paragraph{}
 The geometric (resp. algebraic) construction of flat connection is explained in the Appendix \ref{geom_flat_subsec} (resp. Section \ref{nc_flat_sec}). The following table gives a mapping between these constructions.
 \newline
 	\begin {table}[H]
 \caption {The mapping between geometric and algebraic notions} \label{geo_alg_table}
 \begin{tabular}{|c|c|c|}
 	\hline
 	&Geometry & Algebra\\
 	\hline
 	&	&\\
 	1&Riemannian manifold $M$.  & Spectral triple $\left(\Coo\left(M \right), L^2\left(M, \sS \right), \slashed D   \right)$. \\ & & \\
 	2&Topological covering $\widetilde{M} \to M$. & Noncommutative covering, \\ & & $\left(C\left(M \right), C\left(\widetilde{M} \right), G\left( \left.\widetilde{M}~\right|M\right)   \right),$  \\
 	& & given by the Theorem \ref{pavlov_troisky_thm}. \\ & & \\
 	3&Natural structure of Riemannian  & Triple $\left(\Coo\left(\widetilde{M} \right), L^2\left(\widetilde{M}, \widetilde{\sS} \right), \widetilde{\slashed D}   \right)$ is the\\
 	& manifold on the covering space $\widetilde{M}$.&  $\left(C\left(M \right), C\left(\widetilde{M} \right), G\left( \left.\widetilde{M}~\right|M\right)   \right)$ -lift
 	\\
 	&  & of $\left(\Coo\left(M \right), L^2\left(M, \sS \right), \slashed D   \right)$.\\ & & \\
 	4&Group homomorphism   & Action $ G\left( \left.\widetilde{M}~\right|M\right) \times \C^n \to \C^n$\\
 	&$ G\left( \left.\widetilde{M}~\right|M\right) \to GL\left(n, C \right)$ & \\ & & \\
 	5&Trivial bundle $\widetilde{M}\times \C^n$. & Free module $\Coo\left(\widetilde{M} \right) \otimes \C^n$. \\ & & \\
 	6&Canonical flat connection on $\widetilde{M}\times \C^n$ & Trivial flat connection on $\Coo\left(\widetilde{M} \right) \otimes \C^n$\\ & & \\
 	7&Action of $G\left( \left.\widetilde{M}~\right|M\right)$ on $\widetilde{M}\times \C^n$  & Action of $G\left( \left.\widetilde{M}~\right|M\right)$ on  $\Coo\left(\widetilde{M} \right) \otimes \C^n$\\ & & \\
 	8&Quotient space  & Invariant module    \\
 	& $P = \left( \widetilde{M}\times \C^n\right)/G\left( \left.\widetilde{M}~\right|M\right).$& $\mathcal E =  \left( \Coo\left( \widetilde{M}\right) \otimes \C^n\right)^{G\left( \left.\widetilde{M}~\right|M\right)}$ \\ & & \\
 	9&Geometric flat connection on $P$ & Algebraic flat connection on $\mathcal E$.\\ & & \\
 	\hline
 \end{tabular}
\end{table}
\section{Unoriented spectral triples}\label{unoti_defn_sec}
 \paragraph*{}
 Let $M$ be an unoriented  Riemannian manifold, and let $\widetilde{M} \to M$ be a two-fold covering by an oriented  Riemannian manifold $\widetilde{M}$ which admits a Spin$^c$ structure (cf. Definition \ref{spin_str_defn}). There is an action of $\Z_2 \times \widetilde{M} \to \widetilde{M}$ such that $M \cong \widetilde{M}/ \Z_2$. These 
 considerations inspire the following construction.
 
 \begin{empt}\label{unoriented_empt}
  Let $\left(\A, \widetilde \A, \Z_2, \pi  \right)$ be a {noncommutative finite-fold covering of  bounded operator *-algebras} (cf. Definition \ref{fin_oa_defn}), and let $\left(\widetilde \A, \widetilde \H, \widetilde D, \widetilde J  \right)$ be a spectral triple. Let $\rho: A \to B\left( \H\right)$ be a faithful representation.  Suppose that following conditions hold:
 \begin{enumerate}
 	\item[(a)] If both $A$ and $\widetilde A$ are $C^*$-norm completions of $\A$ and $\widetilde \A$ respectively, $\left(A, \widetilde A, \Z_2, \widetilde\pi  \right)$ is an unital noncommutative finite-fold covering (cf. Definition \ref{fin_unital_defn}) which corresponds to  $\left(\A, \widetilde \A, \Z_2, \pi  \right)$ then  the natural representation  $\widetilde\rho:\widetilde A \hookto B\left(\widetilde\H\right)$ is induced by the pair
 	$$
 	\left( \rho, \left(A, \widetilde A, \Z_2, \widetilde\pi  \right)\right) 
 	$$
 	(cf. Definition \ref{induced_repr_fin_defn}).
 	\item[(b)] If $\Z_2\times\widetilde \H\to   \widetilde \H$ is an explained in \ref{g_act_induced_empt} action then
 	$$
 \Z_2\times\Dom \widetilde D = \Dom \widetilde D.
 	$$
 \item[(c)]	
 The operator $\widetilde D$ is $\Z_2$-invariant i.e.
 $$
 \forall \widetilde \xi \in \Dom \widetilde D \quad \forall g \in \Z_2 \quad \widetilde D\left(  g \widetilde \xi \right) = g \left( \widetilde D\widetilde \xi\right).
 $$
 \item[(d)] The  The operator $\widetilde J$ is $\Z_2$-invariant i.e.
 $$
 \forall \widetilde \xi \in \sH \quad \forall g \in \Z_2 \quad \widetilde J\left(  g \widetilde \xi \right) = g \left( \widetilde J\widetilde \xi\right).
 $$
 \item[(e)] If
 	$$
 \cc = \tsum_j (a_j^0 \ox b_j^0) \ox a_j^1 \oxyox a_j^n 
 \in Z_n(\sA, \sA \ox \sA^\opp),
 $$
 is a required by an Axiom \ref{orientation_st_ax} Hochschild $n$-cycle and an element $g \in \Z_2$ is not trivial then one has
 \bean
 \tsum_j  ga_j^0 (J gb_j^{0*} J^{-1}) \,[D,ga_j^1] \dots [D,ga_j^n]= - \tsum_j a_j^0 (J b_j^{0*} J^{-1}) \,[D,a_j^1] \dots [D,a_j^n]
 \eean
 (cf. Equation \ref{eq:vol-cond}).
  \end{enumerate}

 \end{empt}
\begin{exercise}\label{unoriented_exer}
	Under the hypotheses \ref{unoriented_empt} there is a natural inclusion $\H \subset\widetilde \H$ (cf. equation \eqref{hilb_fin_inc_eqn}).
Prove  that $\widetilde D \left( \H \cap \Dom \widetilde D \right)\subset \H$, i.e. the restriction $D \bydef \left. \widetilde D\right|_{\H}$ can be regarded as an unbounded operator on $\H$. 
\end{exercise}
 \begin{definition}\label{unoriented_defn} 
 	Under the hypotheses \ref{unoriented_empt} we say that
 $\left( \A,  \H,  D, J \bydef\left.\widetilde J\right|_{\H} \right)$ is an \textit{unoriented  spectral triple}.
\end{definition}

 \chapter{Noncommutative infinite coverings}\label{infinite_covering_chap}

 \section{Coverings of $C^*$-algebras}\label{infinite_ca_sec}
 \paragraph*{}
 
 This section contains  noncommutative generalizations of infinite coverings.
 
 \subsection{Basic constructions}
 \begin{definition}\label{g_category_defn}
 	Let $ G$ be a residually finite group (cf. Definition \ref{residually_finite_defn}). 
 	Let $\left\{G_\la\right\}_{\la\in \La}$ be the indexed by a set $\La$ family of all finite factor groups of $ G$. If we define an order on $\La$ by the following way
 	\be\label{top_group_order_eqn}
 	\mu \ge \nu \quad \Leftrightarrow \quad \text{there is the natural homomorphism} \quad G_\mu \to G_\nu
 	\ee
 	with a natural commutative diagram
 	\newline
 	\begin{tikzcd}
 		G \arrow[rr]\arrow[rd] & & G_\nu\\
 		& G_\mu\arrow[ru] &
 	\end{tikzcd}
 	\\
 	then $\La$ becomes an directed set (cf. Definition \ref{directed_set_defn}). We say that the set $\La$ is $  G$-\textit{set} and the given by $\La$ pre-ordering category (cf. Definition \ref{preordercat_defn}) is $  G$-\textit{category}  $  G$-{category} denoted by $\La_{ G }$.
 	There is the natural functor from $\La_{ G }$ to the category of finite groups and surjective homomorphisms. The element $\la_{\min }\in \La$ which corresponds to the trivial factor-group $G_{\la_{\min}}= \{e\}$ is said to be \textit{minimal}.
 \end{definition}
 \begin{remark}\label{g_category_rem}
 If $ G$ is a residually finite group (cf. Definition \ref{residually_finite_defn}) and $\widehat{G}\bydef \varprojlim_{\la\in \La}G_\la$ is its profinite completion (cf. Example \ref{profinite_exm}) then ${G}$-category $\La_{ G }$ is naturally equivalent to $\widehat G$-category $\La_{ \widehat G }$.
 \end{remark}
 
 \begin{empt}\label{infinite_quasicovering_empt}
 	If  $\widehat{G}$ is a profinite group (cf. Example \ref{profinite_exm})  then 
 	  $\widehat{G} \bydef \varprojlim_{\la \in \La}  G_\la$ is an inverse limit of finite groups (cf. Definition \ref{group_inv_lim_defn}). According to the Definition \ref{inverse_limit_defn} the set $\La$ is directed. Indeed $\La$ is the $\widehat{G}$-set (cf. Definition \ref{g_category_defn}). Let $\overline A$ be a $C^*$-algebra with an action $\widehat{G}\times \overline A\to \overline A$ such that any $ g \in \widehat G$ yields an $*$-automorphism of $\overline A$.	  
 	  Suppose that for any element $\overline a \in K\left(\overline A \right)$ of the Pedersen's ideal of $\overline A$ (cf. Definition \ref{pedersen_ideal_defn}) a series 
 	\be\label{infinite_covering_basic_eqn}
 	\sum_{	g \in \widehat{G}}g \overline a
 	\ee
 	is convergent with respect to the strict topology of $M\left(\overline A\right)$ (cf. Definition \ref{strict_topology_defn}). For any subset $G \subset \widehat{G}$ a series is $	\sum_{	g \in {G}}g \overline a$ is also convergent with respect to the strict topology of $M\left(\overline A\right)$. For any $\la \in \La$ denote by $A_\la$ a generated by elements
 	\be\label{basic_cov_cl_eqn}
 	a_\la =\bt\text{-} \sum_{	g \in \ker\left( \widehat{G}\to G_\la\right) }g \overline a
 	\ee
 	$C^*$-subalgebra of $M\left(\overline A\right)$, where  $\bt\text{-} \sum$ means a convergence  with respect to the strict topology of $M\left(\overline A\right)$. 
 \end{empt}
 \begin{definition}\label{infinite_quasicovering_defn} 	Under the hypotheses  \ref{infinite_quasicovering_empt}  $\la_{\mathrm{min}}\in \La$ is the minimal element and $A\bydef A_{\la_{\min}}$ (cf. Definition \ref{g_category_defn} then we say that the triple $\left( A, \overline A, \widehat G\right)$ is an  \textit{infinite quasi-covering}. We say that $A_\la$ is the $\la$-\textit{descent of} $\overline A$. The natural injective $*$-homomorphism $\lift_\la: A_\la \hookto M\left(\overline A \right)$ is the $\la$-\textit{lift}.
 \end{definition}
 
 \begin{lemma}\label{infinite_quasi_covering_lem}
 	Under the hypotheses  \ref{infinite_quasicovering_empt} all $\mu, \nu \in \La$ such that $\nu\ge\mu$ there is a natural noncommutative finite-fold quasi-covering $\left(A_\mu, A_\la, G_\nu/ G_\mu, \pi^\mu_\nu\right)$ (cf. Definition \ref{fin_quasi_defn}).
 \end{lemma}
 \begin{proof}
 	If $\left\{g_1, ..., g_n\right\}\in G_\nu$ is a set of representatives of   $G_\nu/ G_\mu$ then for any $j = 1,...,n$ there is $\widehat g_j \in \widehat G$ such that $g_j$ is the image of $\widehat g_j$. If   	$$
 	a_\mu =\bt\text{-} \sum_{	g \in \ker\left( \widehat{G}\to G_\mu\right) }g \overline a\in  A_\mu
 	$$
 	and $b \bydef \sum_{j=1}^n \widehat g_j \overline a$ then
 	$$
 	b_\nu \bydef \sum_{	g \in \ker\left( \widehat{G}\to G_\nu\right) }g \overline b\in  A_\nu
 	$$
 	Clearly $b_\nu = a_\mu$ where the natural inclusions $A_\mu \subset M\left(\overline A\right)$ and $A_\mu \subset M\left(\overline A\right)$ are implied. So there is the natural inclusion $A_\mu \subset A_\nu$ or equivalently one has an injective $*$-homomorphism $\pi^\mu_\nu : A_\mu \hookto A_\nu$. The proof of 
 	$$
 	\pi\left(A_\mu\right) = A_\nu^{ G_\nu/ G_\mu}\stackrel{\text{def}}{=}\left\{
 	\left.a\in A_\nu~\right|~ a = g a;\quad \forall g \in  G_\nu/ G_\mu\right\}
 	$$
 	is left  to the reader. 
 \end{proof}
 \begin{empt}\label{inf_cat_empt}
 	Under the hypotheses  of the Lemma \ref{infinite_quasi_covering_lem} there is a category $\mathfrak{S}$ such that
 	\begin{itemize}
 		\item $\mathfrak{S}$-objects are $C^*$-algebras $A_\la$ where $\la\in \La$,
 		\item $\mathfrak{S}$-morphisms are natural injective $*$-homomorphisms $\pi^\mu_\nu : A_\mu \hookto A_\nu$
 		such that the diagram
 		\\
 		\begin{tikzcd}
 		A_\mu \arrow[rr, "\pi_{\mu\nu}"] \arrow[rd, "\desc_{\mu}"]&& A_\nu\arrow[ld, "\desc_\nu"]\\
 		& M\left( \overline A\right) &
 		\end{tikzcd}
 		\\
 		is commutative.
 	\end{itemize}
 	
 \end{empt}
 \begin{notation}\label{inf_cat_not}
 	Under the hypotheses   \ref{inf_cat_empt}    we use the following notation
 	\be\label{inf_cat_not_eqn}
 	A \bydef A_{\la_{\text{min}}}
 	\ee
 	where $\la_{\text{min}}\in \La$ is the minimal element,
 	\be\label{inf_cat_not_pi_eqn}
 	\pi_\la \bydef \pi^{\la_{\text{min}}}_\la : A \hookto A_\la.
 	\ee
 	
 \end{notation}
 
   \begin{lemma}\label{infinite_pre_covering_lem}
 	Under the hypotheses \ref{infinite_quasicovering_empt} one has
 	\be\label{a_la_eqn}
 	\overline a \in K\left(\overline A \right)\quad  	a_\la \bydef 	\bt~\text{-}\sum_{g \in \ker\left(\widehat{G}\to  G_\la \right)}g \overline a \in  K\left( A_\la\right).	
 	\ee
 \end{lemma}
 \begin{proof}
 	If $\overline a \in K\left( \overline A\right)_+$ then $\overline a \le \sum_{j = 1}^n\overline a_j$ where  $\overline a_j \in  K\left( \overline A\right)_0$ for any $j = 1, ..., n$. From the Lemma \ref{pedersen_eps_lem} it follows that there is a positive element   $\overline b_j \in \overline A_+$ and $\eps_j > 0$ such that  $\overline a = f_{\eps_j}\left( \overline b\right)$ where $f_\eps$ is given by \eqref{f_eps_eqn}. If $\overline c_j \bydef f_{\eps_j/2}\left( \overline b_j\right)$ then  $\overline a_j = f_{\eps_j/2}\left( \overline c_j\right)$. Let denote
 	\bean
 	a^j_\la \bydef \bt~\text{-}\sum_{g \in \ker\left(\widehat{G}\to   G_\la \right)}g \overline a_j, \\
 	c^j_\la \bydef \bt~\text{-}\sum_{g \in \ker\left(\widehat{G}\to   G_\la \right)}g \overline c_j 
 	\eean
 	and let us prove that $a^j_\la \le f_{\eps_j/2}\left( c^j_\la\right)$. If we consider a faithful representation then $\overline A \to B\left( \overline \H\right)$ then for any $\overline \xi \in \overline \H$ such that $\left\| \overline \xi\right\|= 1$ one has
 	\bea\label{subordinated_1_eqn}
 	\left(\overline \xi, f_{\eps_j/2}\left( \overline c_j\right)\overline \xi \right) = \max\left(0,  \left(\overline \xi, \overline c_j\overline \xi \right)- \frac{\eps_j}{2} \right),\\\label{subordinated_2_eqn}
 	\left(\overline \xi, c^j_\la \overline \xi \right)= \sum_{ g \in \ker\left(\widehat{G}\to   G_\la\right)} \left(\overline \xi,\left(  g \overline c_j\right) \overline \xi \right),\\
 	\label{subordinated_3_eqn} \left(\overline \xi, a^j_\la \overline \xi \right)= \sum_{ g \in \ker\left(\widehat{G}\to   G_\la\right)} \max\left(0,  \left(\overline \xi,\left(g  \overline c_j\right) \overline \xi \right)- \frac{\eps_j}{2} \right) 
 	\eea
 	From 
 	\bean
 	\max\left(0, \sum_{ g \in \ker\left(\widehat{G}\to   G_\la\right)} \left(\overline \xi,\left(g  \overline c_j\right) \overline \xi \right)- \frac{\eps_j}{2} \right)\ge \\\ge\sum_{ g \in \ker\left(\widehat{G}\to   G_\la\right)} \max\left(0,  \left(\overline \xi,\left(g  \overline c_j\right) \overline \xi \right)- \frac{\eps_j}{2} \right)
 	\eean
 	and taking into account \eqref{subordinated_1_eqn}-\eqref{subordinated_3_eqn} one has
 	$$
 	\left(\overline \xi, a^j_\la \overline \xi \right)\le \left(\overline \xi, f_{\eps_j/2}\left(  c^j_\la\right) \overline \xi \right)
 	$$
 	So $a^j_\la \le f_{\eps_j/2}\left( c^j_\la\right)$ and taking into account $f_{\eps_j/2}\left(c^j_\la\right)\in K\left(A_\la \right)_0$ we conclude that $a^j_\la\in K\left(A_\la \right)_+$. Using the Definition \ref{pedersen_ideal_defn} one can prove this Lemma for all $\overline a \in K\left( \overline A\right)$.
 \end{proof}

 \begin{definition}\label{infinite_desc_defn}
 	Under the hypotheses of the Lemma \ref{infinite_pre_covering_lem}  for all $\la\in\La$ there is a  natural   homomorphism of $A_\la$-$A_\la$-bimodules given by
 	\be\label{inf_desc_eqn}
 	\begin{split}
 		\desc_{\la} : K\left(\overline A \right) \to K\left(A_\la \right),\\
 		\overline a \mapsto\bt\text{-} \sum_{	g \in \ker\left( \widehat{G}\to G_\la\right) }g \overline a
 	\end{split}
 	\ee
 	where  $\bt\text{-} \sum$ means the convergence with respect to the strict topology of $M\left(\overline A\right)$ (cf. Definition \ref{strict_topology_defn}).
 	The given by the equation  \eqref{inf_desc_eqn}
 	homomorphism   $\desc_{\la}$ of $A_\la$-$A_\la$-bimodules is the $\la$-\textit{descent}. If $\la_{\min}\in \La$ is the minimal element (cf. Definition \ref{g_category_defn}) then a given by	
 	$$
 	\desc \bydef \desc_{\la_{\min}}: K\left(\overline A \right) \to K\left(A \right)
 	$$  homomorphism   of $A$-$A$-bimodules is the \textit{minimal descent}.
 \end{definition}

 \begin{definition}\label{algebraical_finite_covering_category_defn} 
 	The  given by \ref{inf_cat_empt} category is said to be \textit{algebraical finite covering category} if one has:
 	\begin{enumerate}
 		\item [(a)] 
 		any $\mathfrak{S}$-morphism $\pi^\mu_\nu : A_\mu \hookto A_\nu$ is a noncommutative finite-fold  covering (cf. Definition \ref{fin_defn});
 		\item[(b)] for all $\la \in \La$ is the $\la$-{descent} $\desc_{\la} : K\left(\overline A \right) \to K\left(A_\la \right)$  (cf. Definition \ref{infinite_desc_defn}) is surjective, i.e. $\desc_{\la} \left(  K\left(\overline A \right)\right) = K\left(A_\la \right)$, moreover we require that $\desc_{\la} \left(  K\left(\overline A \right)_+\right) = K\left(A_\la \right)_+$. 
 	\end{enumerate}
 		
 	We write
 	\be\label{algebraical_finite_covering_category_eqn}
 	\mathfrak{S}\bydef \left\{\left\{A_\la\right\}_{\la\in \La}, \left\{\pi^\mu_\nu : A_\mu \hookto A_\nu\right\}_{\substack{\mu, \nu \in \La\\\mu \le \nu}}\right\}
 	\ee
 Moreover the given by the Definition \ref{infinite_quasicovering_defn} infinite quasi-covering $\left( A, \overline A, \widehat G\right)$ is said to be a \textit{pre}-\textit{covering of the algebraical finite covering category}  $\mathfrak{S}$.
 \end{definition}
  
\begin{definition}\label{algebraical_reduced_finite_covering_category_defn}
 	A given by the equation \eqref{algebraical_finite_covering_category_eqn}  algebraical finite covering category 	(cf. Definition \ref{algebraical_finite_covering_category_defn}  is \textit{reduced} if for  any $\la\in \La$ the given by \eqref{inf_cat_not_pi_eqn} injective $*$-homomorphism $
 	\pi_\la : A \hookto A_\la$ is a reduced  noncommutative finite-fold covering (cf. Definitions  \ref{fin_defn} and \ref{fin_red_defn}).
 \end{definition}
 \begin{empt}\label{infinite_covering_empt} Let  $\mathfrak{S}\bydef\left\{\left\{A_\la\right\}_{\la\in \La}, \left\{\pi^\mu_\nu : A_\mu \hookto A_\nu\right\}_{\substack{\mu, \nu \in \La\\\mu \le \nu}}\right\}$ be an algebraical finite covering category (cf. Definition \ref{algebraical_finite_covering_category_defn}). Consider a category $\mathfrak{Cov}\left(\mathfrak{S} \right)$ such that
 	\begin{itemize}
 		\item $\mathfrak{Cov}\left(\mathfrak{S} \right)$-objects are pre-coverings of $\mathfrak{S}$ (cf. Definition  \ref{algebraical_finite_covering_category_eqn}),
 		\item A $\mathfrak{Cov}\left(\mathfrak{S} \right)$-morphism from $\left(A, \overline A', \widehat{G}\right)$ to  $\left( A, \overline A'', \widehat{G}\right)$ is an injective \\$*$-homomorphism $\phi: \overline A' \hookto  \overline A''$  such that
 		\bea\label{infinite_covering_g_eqn}
 		\forall	\overline a' \in \overline A'\quad \forall g \in \widehat{G}\quad g\phi\left( \overline a'\right) = \phi\left(g \overline a'\right),~~~\\
 		\label{infinite_covering_eqn}\forall	\overline a' \in \overline K\left(\overline A'\right) \quad \forall \la\in \La \quad \bt\text{-} \sum_{	g \in \ker\left( \widehat{G}\to G_\mu\right) }g \overline a'= \bt\text{-} \sum_{	g \in \ker\left( \widehat{G}\to G_\mu\right) }g \phi\left( \overline a'\right).~~~ 
 		\eea
 	\end{itemize} 
 \end{empt}
 \begin{definition}\label{disconnected_infinite_noncommutative_covering_defn}
 	The terminal object $\left(A, \overline A, \widehat{G}\right)$ (cf. Definition \ref{terminal_object_defn}) of the category $\mathfrak{Cov}\left(\mathfrak{S} \right)$ (cf. \ref{infinite_covering_empt}) is the  \textit{disconnected infinite noncommutative covering} of $\mathfrak{S}=\left\{\left\{A_\la\right\}_{\la\in \La}, \left\{\pi^\mu_\nu : A_\mu \hookto A_\nu\right\}_{\substack{\mu, \nu \in \La\\\mu \le \nu}}\right\}$. 
 	We also  denote it by the triple
 	\be\label{disconnected_infinite_noncommutative_covering_eqn}
 	\left(A, \overline{A},  \widehat{G}\right).
 	\ee
 \end{definition}
  \begin{lemma}\label{uni_dicsonnected_lem}
 	If a disconnected infinite noncommutative covering (cf. Definition \ref{disconnected_infinite_noncommutative_covering_defn}) exists then it is unique up to $*$-isomorphism. 
 \end{lemma}
 \begin{proof}
 	If both $\overline{A}'$ and $\overline{A}''$ are disconnected infinite noncommutative coverings then there are injective $*$-homomorphism $\phi' : \overline{A}'\hookto\overline{A}''$ and $\phi'': \overline{A}''\hookto \overline{A}'$. From the Lemma \ref{bt_lim_lem} it follows that any $\overline a' \in \overline A'$ is given by
 	$$
 	\overline a' =\bt\text{-}\lim_{\la \in \La} a'_\la
 	$$
 	where 
 	$$
 	a'_\la \bydef 	\bt~\text{-}\sum_{g \in \ker\left(\widehat{G}\to G_\la\right)}g \overline a' \in \phi\circ \varphi_\la\left( K\left( A_\la\right) \right) 
 	$$
 	where $\bt\text{-}\lim$ is the limit with respect to the strict topology of $M\left( \overline{A}\right)$ (cf. \eqref{bt_lim_eqn}). On the other hand from \eqref{bt_lim_eqn} it follows that $$
 	\phi'\left(\overline a' \right)= \bt\text{-}\lim_{\la \in \La} a'_\la.
 	$$ 
 	Similarly one can prove that
 	$$
 	\phi'' \circ \phi'\left( \overline a'\right) = \bt\text{-}\lim_{\la \in \La} a'_\la =  \overline a',
 	$$
 	i.e. $\phi'' \circ \phi' = \Id_{\overline A'}$. It turns out that both $\phi'$ and $\phi''$ are $*$-isomorphisms.
 \end{proof}

 \begin{lemma}\label{bt_lim_lem}
 	Under the hypotheses of the Definition  \ref{infinite_quasicovering_defn}  one has
 	\be\label{bt_lim_eqn}
 	\forall \overline a \in K\left( \overline A\right)  \quad \overline a =\bt\text{-}\lim_{\la \in \La} a_\la= \bt\text{-}\lim_{\la \in \La} \desc_\la\left(\overline a \right) 
 	\ee
 	(cf. \eqref{inf_desc_eqn})
 	where $a_\la$ is given by \eqref{basic_cov_cl_eqn} $\bt\text{-}\lim$ is the limit with respect to the strict topology of $M\left( \overline{A}\right)$ (cf. Definition \ref{strict_topology_defn}).
 \end{lemma}
 \begin{proof}
 	If $\la = \la_{\min}$ then from \eqref{basic_cov_cl_eqn} it follows that there is $a \in A = A_{\la_{\min}}$ such that
 	$$
 	a  = \bt\text{-}\sum_{g \in \widehat{G}}g \overline a
 	$$
 	If $\overline b \in \overline A$ and $e \in \widehat{G}$ is a unity of  $\widehat{G}$  then for any $\eps > 0$ there is a finite subset $G_0\subset  \widehat{G}$ with $e \in G_0$ such that for any finite set $G \subset \widehat G$ with $G_0\subset G$ following condition holds
 	$$
 	\left\| \sum_{g \in G}\left( g \overline a \right) \overline b - a \overline b \right\| <\eps.
 	$$
 	It turns out that 
 	for any finite set $G'\subset  \widehat{G} \setminus G_0$ one has
 	$$
 	\left\| \sum_{g \in G'}\left( g \overline a\right)  \overline b \right\| <\eps.
 	$$
 	Since the set $G'$ is finite there is $\la_\eps\in \La$ such that 
 	$$
 	\ker\left(\widehat{G}\to G_{\la_e}\right)\cap G' = \{e\}
 	$$
 	It follows that
 	$$
 	\forall \la \ge \la_\eps \quad \left\| a_\la \overline b - \overline a \overline b\right\| < \eps.
 	$$
 	The above equation means the following
 	$$
 	\lim_{\la\in \La} a_\la \overline b - \overline a \overline b
 	$$
 	where the $C^*$-norm limit is implied. Similarly we can prove that $\lim_{\la\in \La} \overline b a_\la - \overline b \overline a$, so one has \eqref{bt_lim_eqn}.
 	
 \end{proof}
  \begin{corollary}\label{infinite_multiplier_cor}
 	Under the hypotheses of the Definition  \ref{infinite_quasicovering_defn} for all $\la\in \La$  the $\la$-lift $\lift_\la: A_\la \hookto M\left(\overline A \right)$  (cf. Definition \ref{infinite_quasicovering_defn})    can be uniquely extended up to the natural injective $*$-homomorphism  
 	\be\label{inf_mult_mult_hom_eqn}
 	M\left( \lift_\la \right) :	M\left(A_\la \right)  \hookto M\left( \overline A\right).
 	\ee
 \end{corollary}
 \begin{proof}
 	If $\overline{a}	\in K\left(  \overline A\right)_+ $ is a positive element of Pedersen's ideal (cf. Definition \ref{pedersen_ideal_defn}) then $\overline{b}\bydef \sqrt{\overline{a}}	\in K\left(  \overline A\right)_+$. From the Lemma \ref{bt_lim_lem} it follows that is
 	$$
 	b_\la \bydef  \desc_{\la}\left(\overline b \right) \in K\left(A_\la \right)_+ 
 	$$
 (cf. Definition \ref{infinite_desc_defn})	then there is the $C^*$-norm limit.
 	$$
 	\overline{a} = \lim_{\la\in \La}b_\la \overline{b}.
 	$$
 	where $\desc_{\la}$ is the $\la$-descent (cf. Definition \ref{infinite_desc_defn}). 
 	One has
 	\be \label{inf_la_conv_eqn}
 	\begin{split}
 		\forall \la', \la'' \quad  \la' \le \la''\quad 	\left(  b_{\la'} - \overline b\right) \overline{b}^2\left(  b_{\la'} - \overline b\right)\ge \left(  b_{\la''} -  \overline b\right) \overline{b}^2\left(  b_{\la''} - \overline b\right)\quad \Rightarrow \\ \Rightarrow\left\|   \left(  b_{\la'} - \overline b\right)\overline{b} \right\|\ge \left\|   \left(  b_{\la''} - \overline b\right)\overline{b} \right\|.
 	\end{split}
 	\ee
 	If $\eps > 0$ then there is $\la_\eps \in \La$ such that $\left\|  b_{\la_\eps} \overline b   -  \overline b\overline b \right\|\le {\eps}/{3}$ end taking into account \eqref{inf_la_conv_eqn} we conclude that
 	\be\label{inf_la_bb_eqn}
 	\forall\la \in \La \quad  \la \ge \la_\eps \quad  \left\|  b_\la \overline b   -  \overline b\overline b \right\|\le \frac{\eps}{3}.
 	\ee
 	If  $\left\{u_\a\right\}_{\a \in \mathscr A}$ be an approximate unit of $A$ (cf. Definition \ref{approximate_unit_defn}), then from 
 	\bean
 	\forall \la', \la'' \quad  \la' \le \la''\quad b_{\la'} - \overline b\ge \overline b _{\la''} - \overline b \\
 	\forall \a', \a''\in  \mathscr A  \quad \a' \le \a'' \quad \left(  1_{M\left( \overline A\right) }-u_{\a'}\right) \le \left(  1_{M\left( \overline A\right) }-u_{\a''}\right)
 	\eean
 	it follows that
 	\bean
 	\forall \la', \la'' \quad  \la' \le \la''\quad\Rightarrow\quad  \forall \a', \a''\in  \mathscr A  \quad \a' \le \a''\quad \Rightarrow\\
 	\Rightarrow\quad  \left(  1_{M\left( \overline A\right) }-u_{\a'}\right)\left(  b_{\la'} - \overline b\right) \overline{b}^2\left(  b_{\la'} - \overline b\right)\left(  1_{M\left( \overline A\right) }-u_{\a'}\right)\ge\\
 	\ge \left(  1_{M\left( \overline A\right) }-u_{\a''}\right)\left(  b_{\la''} - \overline b\right) \overline{b}^2\left(  b_{\la''} - \overline b\right)\left(  1_{M\left( \overline A\right) }-u_{\a''}\right).
 	\eean
 	So one has 
 	\be \label{inf_in_eqn}
 	\begin{split}
 		\forall \la', \la'' \in \La\quad  \la' \le \la''\quad\Rightarrow\quad  \forall \a', \a''\in  \mathscr A  \quad \a' \le \a''\quad \Rightarrow\\
 		\Rightarrow\quad \left\| \left(  1_{M\left( \overline A\right) }-u_{\a''}\right)\left(  b_{\la''} - \overline b\right) \overline{b}\right\|\le \left\| \left(  1_{M\left( \overline A\right) }-u_{\a'}\right)\left(  b_{\la'} - \overline b\right) \overline{b}\right\|\quad \Leftrightarrow\\
 		\Leftrightarrow \left\| u_{\a''}\left(\overline b^2-b_{\la''}\overline b \right)-\left(\overline b^2-b_{\la''}\overline b \right)   \right\|\le  \left\| u_{\a'}\left(\overline b^2-b_{\la'}\overline b \right)-\left(\overline b^2-b_{\la'}\overline b \right)   \right\|
 	\end{split}
 	\ee
 	Similarly one can proof 
 	\be \label{inf_inn_eqn}
 	\begin{split}
 		\forall \la\in\La \quad   \forall \a', \a''\in  \mathscr A  \quad \a' \le \a''\quad\Rightarrow\quad \left\| u_{\a''}\left(b_{\la}-\overline b \right)  \right\|\le \left\| u_{\a'}\left(\overline b-b_{\la} \right)  \right\|
 	\end{split}
 	\ee
 	There is $\a_\eps \in \mathscr A$ such that 
 	$$
 	\left\| u_{\a_\eps}  b_{\la_\eps } \overline b   -   b_{\la_\eps } \overline b  \right\| \le \frac{\eps}{3\left\| \overline b \right\|}\quad\Rightarrow \quad \left\| u_{\a_\eps}  b_{\la_\eps } \overline b   -   b_{\la_\eps } \overline b  \right\| \le \frac{\eps}{3}
 	$$
 	
 	From the equations \eqref{inf_la_conv_eqn} and \eqref{inf_inn_eqn} it follows that 
 	\be\label{inf_bb_eqn}
 	\forall\la \in \La \quad  \la \ge \la_\eps \quad\Rightarrow\quad  \forall \a\in  \mathscr A  \quad \a\ge  \a_{\eps} \quad\Rightarrow\quad\left\| u_{\a}  b_{\la } \overline b   -   b_{\la } \overline b  \right\|\le \frac{\eps}{3}
 	\ee
 	
 	From both \eqref{inf_la_bb_eqn} and \eqref{inf_in_eqn} it follows that 
 	
 	\be\label{inf_mn_eqn}
 	\begin{split}
 		\forall \la\in \La \quad  \la_\eps \le \la\quad\Rightarrow\quad  \forall \a\in  \mathscr A  \quad \a_\eps  \le \a\quad \Rightarrow\\
 		\Rightarrow \left\| u_{\a}\left(\overline b^2-b_{\la}\overline b \right)-\left(\overline b^2-b_{\la}\overline b \right)   \right\|\le  \left\| u_{\a_\eps}\left(\overline b^2-b_{\la_\eps}\overline b \right)-\left(\overline b^2-b_{\la_\eps}\overline b \right)   \right\|\Rightarrow\\
 		\Rightarrow \left\| u_{\a}\left(\overline b^2-b_{\la}\overline b \right)-\left(\overline b^2-b_{\la}\overline b \right)   \right\|\le\\ \le \left\| u_{\a_\eps}\right\| \left\| \overline b^2-b_{\la_\eps}\overline b   \right\|+\left\| \overline b^2-b_{\la_\eps}\overline b   \right\|\le \frac{2}{3}.
 	\end{split}
 	\ee
 	From the triangle inequality it follows that
 	\bean
 	\forall\la \in \La \quad    \forall \a\in  \mathscr A\quad  \left\|  u_{\a}\overline{b}^2- \overline{b}^2 \right\|\le \\
 	\left\|\left(  u_{\a}\overline{b}^2- \overline{b}^2\right)+ \left( b_{\la}  \overline{b} - u_{\a}b_{\la} \overline{b}\right)  \right\|+ \left\| b_{\la} \overline{b} - u_{\a}b_{\la} \overline{b}  \right\|,
 	\eean
 	and taking into account the equations \eqref{inf_bb_eqn} \eqref{inf_mn_eqn} one has
 	$$
 	\forall \a\in  \mathscr A  \quad \a_\eps  \le \a\quad \Rightarrow \left\|  u_{\a}\overline{b}^2- \overline{b}^2 \right\| < \eps.
 	$$
 	Taking into account that $\overline{a}=\overline{b}^2$ we conclude that 
 	$$
 	\lim_{\a \in \mathscr A}\left\|  u_{\a}\overline{a}- \overline{a} \right\|= 0.
 	$$ 
 	Since Pedersen's ideal $K\left(\overline A\right)$ is dense in $\overline A$ one has
 	$$
 	\forall \overline a \in \overline A_+\quad \lim_{\a \in \mathscr A}\left\|  u_{\a}\overline{a}- \overline{a} \right\|= 0,
 	$$
 	and taking into account \eqref{four_decompositon_eqn} we conclude
 	$$
 	\forall \overline a \in \overline A\quad \lim_{\a \in \mathscr A}\left\|  u_{\a}\overline{a}- \overline{a} \right\|= 0.
 	$$
 	Similarly one can prove that
 	$$
 	\forall \overline a \in \overline A\quad\forall\la \in \La\quad  \lim_{\a \in \mathscr A}\left\|  \overline{a}u_{\a}- \overline{a} \right\|= 0.
 	$$
 	Now this corollary follows from the Lemma \ref{lift_mult_lem}.
 \end{proof}

\begin{remark}\label{ininite_covering_inductive_rem}
 	Since the set $\La$ is directed there is a $C^*$-inductive limit 
 	\bean
 	\widehat A \bydef C^*\text{-}\varinjlim_{\la \in \La} A_\la	
 	\eean
 	in the sense of the Definition \ref{inductive_lim_non_defn}.  Indeed $\widehat A$ is the $C^*$-norm completion of the union
 	$$
 	\bigcup_{\la \in \La} A_\la	\subset M\left(\overline A \right) 	
 	$$
 	It turns out that there is an inclusion
 	\be\label{ininite_covering_inductive_eqn}
 	\widehat A \bydef C^*\text{-}\varinjlim_{\la \in \La} A_\la\subset  M\left(\overline A \right).
 	\ee
 \end{remark}
 \begin{corollary}\label{bt_lim_cor}
 	Under the above hypotheses  both if $\widehat A$ is given by \eqref{ininite_covering_inductive_eqn} then both sets
 	\bean
 	\widehat A  \overline A\bydef \left\{\widehat a  \overline a\left| \widehat a\in \widehat A \quad \overline{a}\in  \overline A\right.\right\},\\
 	\overline A  \widehat A\bydef \left\{ \overline a  \widehat a\left| \widehat a\in \widehat A \quad \overline{a}\in  \overline A\right.\right\}
 	\eean
 	are dense in $\overline A$.
 \end{corollary}
 \begin{proof}
 	Let $\overline a \in  K\left( \overline{A}\right)_+$ be a positive element and $\overline b\bydef \sqrt{\overline a}$. If $b_\la  A_\la$ is such that
 	$$
 	b_\la =	\bt-\sum_{g \in \ker\left(\widehat{G}\to G_\la\right)}g \overline a 
 	$$
 	(cf. equation \eqref{basic_cov_cl_eqn}. It follows that
 	for all $\eps > 0$ there is $\la_\eps$ such that 
 	$$
 	\forall \la \in \La \quad \la \ge \la_\eps \quad \Rightarrow \quad \left\| b_\la \overline b - \overline b~\overline b\right\| = \left\| b_\la \overline b - \overline a\right\| < \eps.
 	$$
 	Taking into account that $K\left( \overline{A}\right)$ is dense in $\overline{A}$ and $b_\la\in \widehat A$ and  we conclude that $\widehat A  \overline A$ is dense in $\overline A$. Similarly one can prove that  $\overline A\widehat A$ is dense in $\overline A$.
 \end{proof}
 \begin{corollary}\label{bt_lim_a_cor}
 	Under the above hypotheses if $A$ is given by \eqref{inf_cat_not} then both sets
 	\bean
 	A  \overline A\bydef \left\{ a  \overline a\left| \widehat a\in  A \quad \overline{a}\in  \overline A\right.\right\},\\
 	\overline A   A\bydef \left\{ \overline a   a\left|  a\in  A \quad \overline{a}\in  \overline A\right.\right\}
 	\eean	
 	are dense in $\overline A$.
 \end{corollary}
 \begin{proof}
 	From the Lemma \ref{quasi_approximate_unit_lem} it follows that $\pi_\la\left( A\right)A_\la$ (cf. equation \eqref{inf_cat_not_pi_eqn}) is dense in $A_\la$  for any $\la \in \La$. If follows that a set
 	$$
 	A\bigcup_{\la\in \La} A_\la\bydef \left\{\left.a a^\cup\right| a \in A\quad a^\cup \in  \bigcup_{\la\in \La} A_\la\right\}
 	$$
 	is dense in $\bigcup_{\la\in \La} A_\la$ where the natural inclusion $A \subset A_\la$ is implied. Since $\bigcup_{\la\in \La} A_\la$ is dense in $\widehat{A}$ the set $\varphi\left(  A \right)  \widehat{A}$ is dense in $ \widehat{A}$. From the Corollary \ref{bt_lim_cor} it follows that  both $ A \overline A$ is dense in $\overline A$. Similarly one can prove that   $\overline A A $ is dense in $\overline A$ 
 \end{proof}

\begin{empt}
 	Let	 $\left(A, \overline{A},\widehat{G} \right)$ be a  disconnected infinite noncommutative covering of $\mathfrak{S}=\left\{\left\{A_\la\right\}_{\la\in \La}, \left\{\pi^\mu_\nu : A_\mu \hookto A_\nu\right\}_{\substack{\mu, \nu \in \La\\\mu \le \nu}}\right\}$ (cf. Definition \ref{disconnected_infinite_noncommutative_covering_defn}). If $\widetilde A$ is a connected component (cf. Definition \ref{connected_comp_defn}) of $\overline{A}$, i.e. $\overline{A} = \widetilde A \oplus \widetilde A^\perp$, and
 	\be\label{infinite_covering_transformation_group_eqn}
 	G\left(\left.\widetilde{A}~\right| A\right)\bydef 
 	\left\{\left. g \in \widehat{G}\right| \forall \widetilde a^\perp \in \widetilde A^\perp \quad g \widetilde a^\perp= \widetilde a^\perp\right\}
 	\ee
 	then there is a natural action
 	\be\label{gta_act_eqn} 
 	G\left(\left.\widetilde{A}~\right| A\right)\times \widetilde{A} \to \widetilde{A}.
 	\ee
 	The $\la_{\min}$-lift  $\la_{\min}: A \hookto M\left(\overline A \right)$ (cf. Definition \ref{infinite_quasicovering_defn}) induces a $*$-homomorphism
 	\be\label{inj_to_m_eqn}
 	\widetilde{\lift} : A \to M\left(\widetilde A \right)
 	\ee
 	Similarly for all $\la\in \La$ one can obtain a natural $*$-homomorphism
 	\be\label{inj_to_ml_eqn}
 		\widetilde{\lift}_\la: A_\la \to M\left(\widetilde A \right)
 	\ee

 \end{empt}
 \begin{definition}\label{infinite_lift_defn}
 The given by the equations \eqref{inj_to_m_eqn} and \eqref{inj_to_ml_eqn} injective *-homomorphisms are  \textit{lift} and $\la$-\textit{lift} respectively.
 \end{definition}
 \begin{remark}\label{dense_inf_rem} 
 Form the Corollary \ref{bt_lim_a_cor} it follows that the subsets $\widetilde {\lift}\left( A\right) \widetilde A$, $~\widetilde A\widetilde {\lift}\left( A\right)$, $~\widetilde {\lift}_\la\left( A_\la\right) \widetilde A$ and $\widetilde A\widetilde {\lift}_\la\left( A_\la\right)$ are dense in $\widetilde A$.
 \end{remark}

 \begin{remark}
The given by \eqref{algebraical_finite_covering_category_eqn} {algebraical finite covering category} (cf. Definition \ref{algebraical_finite_covering_category_defn}) is a subcategory of the  {category $\mathfrak{FinCov}$-$A$  of finite-fold coverings} of $A$ (cf. Definition \ref{subcategory_defn}) of (cf. Definition \ref{fin_category_defn}).
 \end{remark}
 
\begin{definition}\label{good_defn}
A  disconnected infinite noncommutative covering 	$\left(A, \overline{A},\widehat{G}\right)$ be of $\mathfrak{S}$ (cf. Definition \ref{disconnected_infinite_noncommutative_covering_defn}) is \textit{good} if  following conditions hold:
 	\begin{enumerate}
 		\item[(a)] if both $\widetilde{A}'$ and $\widetilde{A}''$ are  {connected components} of $\overline A$ then there is  $g \in \widehat{G}$ such that $g \widetilde{A}'= \widetilde{A}''$,
 		\item [(b)] if $\widetilde A$ is a   connected component of $\overline{A}$ if  $G\left(\left.\widetilde{A}~\right| A\right)$ is given by \eqref{infinite_covering_transformation_group_eqn} then for any $\la \in \La$ the restriction $h_\la|_{\widetilde A}$ is an epimorphism, i. e. $h_\la\left(G\left(\left.\widetilde{A}~\right| A\right) \right) = G\left(\left. A_\la~\right|~A \right)$.
 	\end{enumerate}
 \end{definition}
 \begin{remark}\label{cond_b_ref_rem}
 	If $\overline  \sX$ is a spectrum of $\overline{A}$ then the condition (b) of the Definition \ref{good_defn}  can be reformulated by a following way. If both $\widetilde{\sX}'$ and $\widetilde{\sX}''$ are  {quasi-components} of $\overline \sX$  (cf. Definition \ref{top_quasi_component_defn}) then there is  $g \in \widehat{G}$ such that $g \widetilde{\sX}'= \widetilde{\sX}''$ (cf. Exercise \ref{conn_comp_exer}).
 \end{remark}
 \begin{lemma}
 	 If $\left(A, \overline{A},\widehat{G}\right)$ is a good  disconnected infinite noncommutative covering of $\mathfrak{S}$ (cf. Definitions \ref{disconnected_infinite_noncommutative_covering_defn}  and \ref{good_defn}) then both  given by  	 \ref{infinite_lift_defn}  lift and $\la$-lift are injective.
\end{lemma}
\begin{proof}
Let	$\widetilde {\lift} : A \to M\left(\widetilde A \right)$ be the given by the Definition  \ref{infinite_lift_defn} lift. If $\widetilde {\lift}$ is not injective then there is $a \in A\setminus \{0\}$ such that $\widetilde {\lift}\left( a\right)  \widetilde A = \{0\}$. For any $g \in \widehat G$ from $g\left(\widetilde {\lift}\left( a\right)\widetilde A \right)= \widetilde {\lift}\left( a\right)g \widetilde A $  it follows that $\widetilde {\lift}\left( a\right) \left( g  \widetilde A \right) = \{0\}$ . Now from the condition (a) of the Definition \ref{good_defn} the operator $\widetilde {\lift}\left( a\right)$ acts trivially on every connected component of $\overline A$. It turns out that $a  \overline A = \{0\}$, however it contradicts with the inclusion $a \subset M\left(\overline A \right)$. So $\widetilde {\lift}$ is injective. Similarly one can prove that $\la$-lift (cf. Definition \ref{infinite_lift_defn})  is also injective.
\end{proof}
 
 \begin{definition}\label{infinite_noncommutative_covering_defn}
 	If $\left(A, \overline{A},\widehat{G}\right)$ is a good  disconnected infinite noncommutative covering of $\mathfrak{S}=\left\{\left\{A_\la\right\}_{\la\in \La}, \left\{\pi^\mu_\nu : A_\mu \hookto A_\nu\right\}_{\substack{\mu, \nu \in \La\\\mu \le \nu}}\right\}$ (cf. Definitions \ref{disconnected_infinite_noncommutative_covering_defn}  and \ref{good_defn}) then a connected component $\widetilde{A} \subset \overline{A}$ (cf. Definition \ref{connected_comp_defn}) is said to be the \textit{inverse noncommutative limit of $\mathfrak{S}=\left\{\left\{A_\la\right\}_{\la\in \La}, \left\{\pi^\mu_\nu : A_\mu \hookto A_\nu\right\}_{\substack{\mu, \nu \in \La\\\mu \le \nu}}\right\}$}. The given by the equation \eqref{infinite_covering_transformation_group_eqn} group $G\left(\left.\widetilde{A}~\right| A\right)$  is said to be the \textit{covering transformation group}.  The triple
 	\be\label{infinite_noncommutative_covering_eqn}
 	\left(A, \widetilde{A}, G\left(\left.\widetilde{A}~\right| A\right)\right)
 	\ee
 	 is said to be the  \textit{infinite noncommutative covering} or the  \textit{covering inverse  limit} of  $$\mathfrak{S}=\left\{\left\{A_\la\right\}_{\la\in \La}, \left\{\pi^\mu_\nu : A_\mu \hookto A_\nu\right\}_{\substack{\mu, \nu \in \La\\\mu \le \nu}}\right\}.$$ 
 	We will use the following notation 
 	\begin{equation}\label{inf_lim_not_eqn}
 		\begin{split}
 			\varprojlim \mathfrak{S}\stackrel{\mathrm{def}}{=}\widetilde{A}.
 		\end{split}
 	\end{equation}	
 \end{definition}
\begin{remark}
The motivation of the Definition \ref{infinite_noncommutative_covering_defn} is the Theorem \ref{top_main_thm}.
\end{remark}
\subsection{Equivariant representations}
\subsubsection{Faithful representations}
\begin{definition}\label{equivariant_representation_defn}
If $\widehat{A}\bydef C^*\text{-}\varinjlim_{\la \in \La} A_\la$ (cf. Remark \ref{ininite_covering_inductive_rem}) then a faithful (cf. Definition \ref{faithful_representation_defn}) and nondegenerate (cf. Definition \ref{nondegenerate_repr_defn}) representation ${\pi}: \widehat{A} \hookto B\left( \H\right)$ (cf. Definition \ref{faithful_representation_defn}) is said to be \textit{equivariant} if there is the action  $\widehat{G}\times {\H} \to {\H}$ such that
	\be\label{equivt_eqn}
	g \left(\widehat{a}{\xi} \right) = \left(g\widehat{a} \right) \left(g{\xi} \right); \quad \forall g \in \widehat{G},~ \widehat{a} \in \widehat{A},~ {\xi}\in {\H}.
	\ee
\end{definition}
\begin{remark}\label{equiv_rem}
	If 	a faithful representation ${\pi}: \widehat{A} \hookto B\left( \H\right)$ is {equivariant} then one has:
	\be\label{equivb_eqn}
	\left( ga\right)\xi \bydef  g\left( a \left( g^{-1}\xi\right) \right);\quad \forall g \in \widehat{G},~ a \in \widehat A,~ {\xi}\in {\H}.
	\ee
Moreover the equation \eqref{equivb_eqn} uniquely defines an action $\widehat G \times \widehat A \to \widehat A$.
\end{remark}
\begin{lemma}\label{equiv_lem}
	Following conditions hold:
	\begin{enumerate}
		\item [(i)] the universal  representation $\widehat{\pi}: \widehat{A} \hookto B\left( \widehat\H\right)$ (cf. Definition \ref{universal_rep_defn}) is equivariant,
		\item [(ii)] the atomic  representation ${\pi}_a: \widehat{A} \hookto B\left(\widehat \H_a\right)$ (cf. Definition \ref{atomic_repr_defn}) is equivariant.
	\end{enumerate}
\end{lemma}
\begin{proof}
	(i) If $S$ is the state space of $\widehat{A}$ then there is the natural action $\widehat{G}\times S \to S$ given by
	$$
	\left(gs \right) \left(\widehat a\right)\stackrel{\text{def}}{=} s\left(ga \right)
	$$
	where $s \in S$, $g \in \widehat{G}$ and $\widehat{a}\in \widehat{A}$
	Let $s \in S$, and let $L^2\left(\widehat A, s\right)$ is the Hilbert space of the representation $\pi_s:\widehat{A} \to B\left(L^2\left(\widehat A, s\right)\right)$ which corresponds to $s$ (cf. \ref{gns_constr_sec}). If $f_s: \widehat A\to L^2\left(\widehat A, s\right)$ is the given by \eqref{from_a_to_l2_eqn} natural $\widehat{A}$-linear map then the $\widehat{A}$-module $f_s\left( \widehat A\right) $ is dense in $L^2\left(\widehat A, s\right)$. Since $$
	\widehat \H\bydef \bigoplus_{s \in S} L^2\left(\widehat A, s\right)
	$$
	the $\C$-linear span of given by
	$$
	\xi^s_{\widehat{a}}=\left(0, ..., \underbrace{f_s\left( \widehat{a}\right) }_{s^{\text{th}}-\text{place}}, ..., 0\right)\in \bigoplus_{s \in S} L^2\left(\widehat A, s\right)=		\widehat \H
	$$
	elements	is dense in  $\widehat \H$. Hence if for $g \in  \widehat G$ we define
	$$
	g \xi^s_{\widehat{a}}\stackrel{\text{def}}{=}	\left(0, ..., \underbrace{f_{gs}\left( g\widehat{a}\right) }_{gs^{\text{th}}-\text{place}}, ..., 0\right)\in \widehat \H
	$$
	then   $g \xi$ can be uniquely defined for all $\xi \in \widehat \H$ i.e. there is the  natural  action $\widehat{G}\times \widehat{\H} \to \widehat{\H}$. Using the equality
	$$
	\widehat a_1 f_s\left(\widehat a_2\right)= f_s\left(\widehat a_1\widehat a_2\right)
	$$
	one can prove that the action $\widehat{G}\times \widehat{\H} \to \widehat{\H}$ satisfies to \eqref{equivt_eqn}.\\
	(ii) The replacement of the "state" word in (i) with "pure state" one gives (ii).
\end{proof}

\begin{proposition}\label{infinite_spectrum_limit_prop} 
	If $\left(A, \overline A, \widehat{G}\right)$ is a {pre}-{covering of the algebraical finite covering category}  
		\bean
	\mathfrak{S}\bydef \left\{\left\{A_\la\right\}_{\la\in \La}, \left\{\pi^\mu_\nu : A_\mu \hookto A_\nu\right\}_{\substack{\mu, \nu \in \La\\\mu \le \nu}}\right\}
	\eean
(cf. Definition \ref{algebraical_finite_covering_category_defn}) then one has:
\begin{enumerate}
	\item[(i)] the  spectrum of $\widehat A\bydef \widehat C^*$-$\varinjlim_{\la \in \La}A_\la$ (cf. Definition \ref{spectrum_prime_primtive_defn}) is the inverse limit  (cf. Definition \ref{top_inverse_limit_defn}) of spectra of $A_\la$.
	\item[(ii)] there is the natural one-to-one map between spectra of $\overline A$ and $\widehat A$.
\end{enumerate}	
  \end{proposition}
\begin{proof}
	(i)
 	Denote by $\sX_\la$ and $\widehat \sX$ the spectra of $A_\la$ and $\widehat A$ respectively. From the Proposition  \ref{spectrum_covering_finite_prop}   it turns out that for any $\mu, \nu \in \La$ such that $\mu > \nu$ the inclusion  $\pi^\nu_{ \mu}: A_\nu \to A_\mu$ induces the surjective continuous map $p^\nu_\mu :  \sX_\mu \to \sX_\nu$. Otherwise if $\sY_\la$ and $\widehat \sY$ are the state spaces of $A_\la$ and $\widehat A$ then one has  $\sX_\la\subset \sY_\la$ and $\widehat \sX\subset \widehat  \sY$. Moreover for any $\mu, \nu \in \La$ such that $\mu > \nu$ there is the natural  continuous surjective map $t^\nu_\mu :  \sY_\mu \to \sY_\nu$ such that
 	$$
 	p^\nu_\mu = \left.t^\nu_\mu\right|_{\sX_\mu}
 	$$
 	From the Corollary \ref{inductive_lim_state_nor_cor} and the Remark \ref{unital_rem} it follows that the state space $\widehat \sY = \varprojlim \sY_\la$ of $\widehat A$ is the inverse limit  (cf. Definition \ref{top_inverse_limit_defn}) of state spaces of  $A_\la$, i.e. for any $
 	\la \in \La$ there is the natural map $\widehat t_\la: \widehat \sY \to \sY_\la$ such that the following diagram  
 	\newline
 	\begin{tikzpicture}
 		\matrix (m) [matrix of math nodes,row sep=3em,column sep=4em,minimum width=2em]
 		{
 			&  \widehat \sY &  \\
 			\sY_\mu &  & \sY_\nu \\};
 		\path[-stealth]
 		(m-1-2) edge node [left] {$\widehat t_\mu~~$} (m-2-1)
 		(m-1-2) edge node [right] {$~~~\widehat t_\nu$} (m-2-3)
 		(m-2-1) edge node [above] {$t^\nu_\mu$}  (m-2-3);
 		
 	\end{tikzpicture}
 	\newline
 	is commutative.	Similarly there is the inverse limit $\widehat \sX = \varprojlim \sX_\la$  (cf. Definition \ref{top_inverse_limit_defn}) and the following  diagram 
 	\newline
 	\begin{tikzpicture}
 		\matrix (m) [matrix of math nodes,row sep=3em,column sep=4em,minimum width=2em]
 		{
 			&  \widehat \sX &  \\
 			\sX_\mu &  & \sX_\nu \\};
 		\path[-stealth]
 		(m-1-2) edge node [left] {$\widehat p_\mu~~$} (m-2-1)
 		(m-1-2) edge node [right] {$~~~\widehat p_\nu$} (m-2-3)
 		(m-2-1) edge node [above] {$p^\nu_\mu$}  (m-2-3);
 	\end{tikzpicture}
 	\newline
 	is commutative.

 	Suppose $\left(x_\la\right) \in \prod \sX_\la$  corresponds to an element of the inverse limit $\widehat x \in\varprojlim \sX_\la$, and a state $\widehat\tau: \widehat A \to \C$. For any $\la\in \La$ let $\tau_\la: A_\la \to \C$ be a corresponding $x_\la$ pure state.  Let $\widehat p: \widehat A \to \C$ be a positive functional such that $\widehat p < \widehat\tau$. If $\widehat a \in \widehat A$ then there is $\left(a_\la\right)_{\la\in \La} \in \prod A_\la$ such that $\widehat a = \lim_{\la \in \La}a_\la$. It follows that $\widehat \tau\left(\widehat a \right) = \lim \tau_\la\left( a_\la\right)$ and $\widehat p\left( \widehat a\right) = \lim \widehat p\left( a_\la\right)$. Since the state $\tau_\la$ is pure there is $t \in \left(0, 1\right]$ such that $\left.\widehat p\right|_{A_\la}= t \tau_{\la}$ for all $\la\in \La$. It follows that 
 	\bean
 	\widehat p\left( \widehat a\right) = \lim \widehat p\left( a_\la\right) = \lim t \tau_{\la} \left(a_\la \right)= t \lim  \tau_{\la} \left(a_\la \right)= t \widehat \tau \left(\widehat a\right),
 	\eean 
 	i.e. the state $\widehat\tau$ is pure.

 	Conversely let $\widehat\tau : \widehat A \to \C$ be a pure state and suppose that it correspond to an element  $\left\{y_\la\right\} \in \prod \sY_\la$  which represents  to an element of the inverse limit $\varprojlim \sY_\la$. Assume that the state $\tau_\la: A_\la\to \C$ corresponds to $y_\la$ for all $\la\in\La$. Suppose that $\mu\in \La$ is such that $\tau_\mu$  is not pure. Let $\left\{\la_j\right\}_{j \in \N}\subset \La$ be an increasing cofinal sequence (cf. Definition \ref{cofinal_defn}) such that $\la_1 = \mu$. Since a state $\tau_{\la_1}=\tau_\mu$ is not pure there are pure states $\rho'_{\la_1}, \rho''_{\la_1}: :A_{\la_1}\to\C$ such that
 	\bean
 	\rho'_{\la_1}\perp \rho''_{\la_1},\\
 	\rho'_{\la_1}\preceq \tau_{\la_1},\\
 	\rho''_{\la_1}\preceq \tau_{\la_1},\\
 	\eean
 	Now using induction by $j \in \N$  we find pure states $	\rho'_{\la_j},	\rho'_{\la_j}:A_{\la_j}\to\C$ for all $j\in \N$. If $	\rho'_{\la_{j-1}},	\rho'_{\la_{j-1}}:A_{\la_{j}-1}\to\C$ are known then using the Lemma \ref{sub_pure_lem} one can define  $\rho'_{\la_j},	\rho'_{\la_j}$ such that
 	\bean
 	\rho'_{\la_j}\perp \rho''_{\la_j},\\
 	\rho'_{\la_j}\preceq \tau_{\la_j},\\
 	\rho''_{\la_j}\preceq \tau_{\la_j},\\
 	\left.\rho'_{\la_j}\right|_{A_{\la_{j-1}}}=\rho'_{\la_{j-1}},\\
 	\left.\rho''_{\la_j}\right|_{A_{\la_{j-1}}}=\rho''_{\la_{j-1}}.
 	\eean
 	There are two sequences $\left\{x'_{\la_j}\in \sX_{\la_j}\right\}$ and $\left\{x''_{\la_j}\in \sX_{\la_j}\right\}$. The set $\left\{\la_j\right\}_{j \in \N}$ 
 	is cofinal, so using maps $p^\nu_\mu$ one can construct sets $\left\{x'_{\la}\in \sX_{\la}\right\}_{\la\in\La}$,
 	$\left\{x''_{\la}\in \sX_{\la}\right\}_{\la\in\La}$ such that
 	\bean
 	p^\nu_\mu\left(x'_\nu \right) = x'_\mu,\\
 	p^\nu_\mu\left(x''_\nu \right) = x''_\mu.
 	\eean
 	It follows that one has two elements $\widehat x', \widehat x''\in \widehat\sX$ such that
 	\bean
 	\widehat p_\la\left( \widehat x'\right) = x'_\la,\\
 	\widehat p_\la\left( \widehat x''\right) = x''_\la.
 	\eean
 	From the first part of the proof of this theorem it follows that both $\widehat x'$ and $\widehat x''$ correspond to pure states $\widehat \tau', \widehat \tau'': \widehat A\to \C$. Otherwise from our construction one has
 	\bean
 	\widehat \tau'\perp \widehat \tau'',\\
 	\widehat \tau' \preceq \widehat \tau,\\
 	\widehat \tau'' \preceq \widehat \tau.
 	\eean
 	From the above equation it turns out that the state $\widehat \tau$ is not pure however it contradicts with our assumption. So $y_\la$ corresponds to a pure state for all $\la\in\La$.\\
 	  	(ii) If $\overline \rho: \overline A \to B\left(\H \right)$ is irreducible representation  then from the inclusion $A_\la \subset M\left(\overline A \right)$ and the Definition \ref{multiplier_el_defn} one has a natural representation $\rho_\la : A_\la \to B\left(\H \right)$. From the condition (c)  of the Definition \ref{algebraical_finite_covering_category_defn} it follows that $\rho_\la$ is not trivial. If $\rho_{\la_0} : A_\la \to B\left(\H_\la \right)$ is not irreducible there are two irreducible representations $\rho'_{\la_0} : A_\la \to B\left(\H'_{\la_0} \right)$   and $\rho''_{\la_0} : A_\la \to B\left(\H''_{\la_0} \right)$ such that $\rho'_{\la_0} \perp \rho''_{\la_0}$ and $\rho'_{\la_0}, \rho''_{\la_0}\preceq \rho_{\la_0}$. If $\la \ge \la_0$ and $p: \sX_\la \to \sX_{\la_0}$ is a given by the Proposition \ref{spectrum_covering_finite_prop} and then there are are two irreducible representations $\rho'_{\la} : A_\la \to B\left(\H'_\la \right)$   and $\rho''_{\la} : A_\la \to B\left(\H''_\la \right)$ such that $\rho'_{\la} \perp \rho''_{\la}$ and $\rho'_{\la}, \rho''_{\la}\preceq \rho_{\la}$ and \bean 
 	  	p\left(\rho'_{\la} \right)= \rho'_{\la_0},\\
 	  	 p\left(\rho''_{\la} \right)= \rho''_{\la_0}
 	  	\eean
 	If for  all $\la \ge \la_0$ the representations $\rho'_{\la}$ and $\rho''_{\la}$ correspond to states $\tau'_{\la}$ and $\tau''_{\la}$ then there are two states
 	\bean
 	\overline \tau', \overline \tau '' : \overline A \to \C,\\
\forall \overline a \in K\left(\overline  A\right)\quad  	\overline \tau' \left(\overline a \right) \bydef \lim_{\substack{\la \in \La\\ \la \ge \la_0}}\tau'_{\la}\left(a_\la  \right),\quad 
 	\overline \tau'' \left(\overline a \right) \bydef \lim_{\substack{\la \in \La\\ \la \ge \la_0}}\tau''_{\la}\left(a_\la \right) 
 	\eean 
 	where $a_\la$ is given by the equation \eqref{basic_cov_cl_eqn}.
 	If both $\overline \rho'$ and  $\overline \rho  ''$ a corresponding to  $\overline \tau', \overline \tau ''$ representations then one has
 	\bean
 	\overline \rho'\perp \overline \rho  '',\\
 	\overline \rho', \overline \rho  ''\preceq \overline \rho
 	\eean 
 	i.e. the representation $\overline \rho$ is not irreducible. It is a contradiction, so for all $\la \in \La$ there is a unique representation $\rho_\la: A_\la \to B\left( \H_\la\right)$ which corresponds to $\overline \rho$. If $\rho_\la$ comes from a state $\tau_{\la}: A_\la \to \C$, such that$$
 	 \forall \overline a \in K\left(\overline  A\right)\quad  	\overline \tau \left(\overline a \right) \bydef \lim_{{\la \in \La}}\tau_{\la}\left(a_\la  \right).
 	 $$
 	According to the Definition \ref{inverse_limit_defn} a net $\left\{\rho_\la\right\}$ yields the unique element of the set theoretical limit $\varprojlim_{\la\in \La} \sX_\la$.
 \end{proof}

 \begin{lemma}\label{inverse_lim_h_lem}
 	Let 	 $\left(A, \overline A, \widehat{G}\right)$ be a reduced  {pre}-{covering of the algebraical finite covering category}  
 	\bean
 	\mathfrak{S}\bydef \left\{\left\{A_\la\right\}_{\la\in \La}, \left\{\pi^\mu_\nu : A_\mu \hookto A_\nu\right\}_{\substack{\mu, \nu \in \La\\\mu \le \nu}}\right\}
 	\eean
 	(cf. Definitions \ref{algebraical_finite_covering_category_defn}  and \ref{algebraical_reduced_finite_covering_category_defn}).
Let $\widehat p :  \widehat\sX \to\sX= \sX_{\la_{\mathrm{min}}}$ be the natural continuous map from the spectrum of the inverse limit to the spectrum of $A$ (cf. Definition \ref{top_inverse_limit_defn}). Suppose that for all $\la\in \La$ a noncommutative finite-fold covering $A \hookto A_\la$ (cf. Definition \ref{fin_defn}) is reduced (cf. Definition \ref{fin_red_defn}). 
 	If $\widehat{x }\in \widehat{\sX }$ is a point and both $\rep_{  \widehat{x }}: \widehat A \to B\left( \H_{\widehat{x}}\right)$ and  $\rep_{ \widehat{p} \left( \widehat{x }\right) }:  A \to B\left( \H_{\widehat{p} \left( \widehat{x }\right)}\right)$ are corresponding to $\widehat{x }$  and $\widehat{p} \left( \widehat{x }\right)$ irreducible representations then there are natural isomorphisms  
 	\be\label{inverse_lim_hs_eqn}
 	\begin{split}
 		\rep_{  \widehat{x }}\left(\widehat A \right)= \rep_{\widehat{p} \left( \widehat{x}\right)  } \left( A\right),\\
 		\H_{\widehat{x}}\cong 	\H_{\widehat{p} \left( \widehat{x }\right)}
 	\end{split}
 	\ee
 	of $C^*$-algebras and Hilbert spaces respectively.
 \end{lemma}
 \begin{proof}
 	For any $\widehat a\in \widehat A$ there is a net $\left\{a_\la \in A_\la\right\}_{\la \in \La}$ such that there is a $C^*$-norm limit $\widehat a = \lim_{\la \in \La} a_\la$. It follows that 	$\rep_{  \widehat{x }}\left(\widehat a \right)= \lim_{\la \in \La}\rep_{  \widehat{x }}\left(a_\la \right)$. On the other hand since a covering $A\hookto A_\la$ is reducible one has $\rep_{  \widehat{x }}\left(a_\la \right)\in \rep_{  \widehat{x }}\left(A \right)$, so $\rep_{  \widehat{x }}\left(\widehat a \right)\in \rep_{  \widehat{x }}\left(A \right)$ for all $\widehat x \in \widehat \sX$.
 	If $\widehat x$ corresponds to a state $\tau_{\widehat x}$ then $\H_{\widehat x}$ is a completion of $\widehat A/	\mathcal{I}_{\tau_{\widehat x}}$ where 
 	$$
 	\mathcal{I}_{\tau_{\widehat x}}\bydef \left\{\left.\widehat a \in \widehat A \ \right| \ \tau_{\widehat x}(\widehat a^*\widehat a)=0\right\}
 	$$
 	(cf. equation \ref{tau_ideal_eqn}). If $f_{\tau_{\widehat x}}: \widehat A \to \H_{\widehat x}$ is given by the equation \eqref{from_a_to_l2_eqn} then $f_{\tau_{\widehat x}}\left(  \widehat A\right)$ is a dense subset of $\H_{\widehat x}$. On the other hand since $\bigcup A_\la$ is dense in $ \widehat A $ the space $f_{\tau_{\widehat x}}\left(  \bigcup A_\la\right)$  a dense subset of $\H_{\widehat x}$. However from the equation \eqref{red_incl_eqn} it turns out that $f_{\tau_{\widehat x}}\left(  A_\la\right)= f_{\tau_{\widehat x}}\left(   A\right)$ for every $\la\in \La$ it follows that $f_{\tau_{\widehat x}}\left(  A\right)$ is dense in $\H_{\widehat x}$. We leave to the reader the proof of the following isomorphism
 	$$
 	f_{\tau_{\widehat x}}\left(   A\right)\cong f_{\tau_{\widehat p\left( \widehat x\right) }}\left(  A\right).
 	$$
 	Left and right parts of the above bijection are dense subsets of both $\H_{\widehat{x}}$ and	$\H_{\widehat{p} \left( \widehat{x }\right)}$, so one can obtain the natural isomorphism \eqref{inverse_lim_hs_eqn}.
 \end{proof}
 
 \begin{corollary}\label{inverse_lim_h_cor}
 	If  	 $\left(A, \overline A, \widehat{G}\right)$ is  a reduced  {pre}-{covering of the algebraical finite covering category}  
 \bean
 \mathfrak{S}\bydef \left\{\left\{A_\la\right\}_{\la\in \La}, \left\{\pi^\mu_\nu : A_\mu \hookto A_\nu\right\}_{\substack{\mu, \nu \in \La\\\mu \le \nu}}\right\}
 \eean
 (cf. Definitions \ref{algebraical_finite_covering_category_defn} and \ref{algebraical_reduced_finite_covering_category_defn}) and the Hilbert space of the atomic representation (cf. Definition \ref{atomic_repr_defn})  of both $\widehat A$ and $\overline A$ is the Hilbert norm completion of the given by
 	\be\label{inverse_lim_h_eqn}
 	\bigoplus_{ \widehat x \in \widehat \sX} 	\H_{\widehat{x }}\cong 	\bigoplus_{ \widehat x \in \widehat \sX} 	\H_{\widehat{p} \left( \widehat{x }\right)}
 	\ee 
 	algebraic direct sum.
 \end{corollary}
  \begin{lemma}\label{inverse_faithful_nondegenetate_lem}
  	If 		 $\left(A, \overline A, \widehat{G}\right)$ is  a  {pre}-{covering of the algebraical finite covering category}  
  	\bean
  	\mathfrak{S}\bydef \left\{\left\{A_\la\right\}_{\la\in \La}, \left\{\pi^\mu_\nu : A_\mu \hookto A_\nu\right\}_{\substack{\mu, \nu \in \La\\\mu \le \nu}}\right\}
  	\eean
  	(cf. Definitions \ref{algebraical_finite_covering_category_defn})
 and  $\pi:\overline A \hookto B\left(\overline \H\right)$ is faithful nondegenerate  representation (cf. Definitions \ref{faithful_representation_defn}, \ref{nondegenerate_repr_defn}).  and $\overline a \in K\left(\overline A \right)_+$ is given by the equation \eqref{bt_lim_eqn}, i.e.
 	  	\bean
 	\overline a =\bt\text{-}\lim_{\la \in \La} a_\la\quad a_\la \in K\left(A_\la \right)_+ 
 	 \eean 
 	 (cf. equation \eqref{bt_lim_eqn})
 	 then one has
 	 \be\label{inverse_atomic_lim_h_eqn}
 	 \pi \left(\overline a  \right) = s\text{-}\lim_{{\la \in \La}}\pi\left(a_\la \right) 
 	 \ee
 	 where $s\text{-}\lim_{{\la \in \La}}$ means the convergence with respect to the strong topology of $B\left(\overline\H \right)$ (cf. Definition \ref{strong_topology_defn}). 
 \end{lemma}
 \begin{proof} 
 	From the Lemma  \ref{increasing_convergent_w_lem} it follows that the decreasing net $\left\{\pi\left(a_\la \right) \right\}_{\la \in \La}$ is convergent with respect to the strong topology of $B\left(\overline\H \right)$. Now this lemma follows from the \ref{strict_strong_lem} one.
 \end{proof}
 \begin{empt}\label{infinite_atomic_empt}
	If 	 $\left(A, \overline A, \widehat{G}\right)$ is a   {pre}-{covering of the algebraical finite covering category}  
\bean
\mathfrak{S}\bydef \left\{\left\{A_\la\right\}_{\la\in \La}, \left\{\pi^\mu_\nu : A_\mu \hookto A_\nu\right\}_{\substack{\mu, \nu \in \La\\\mu \le \nu}}\right\}
\eean
(cf. Definitions \ref{algebraical_finite_covering_category_defn} and \ref{algebraical_reduced_finite_covering_category_defn}) and $\pi_a : \overline A\hookto B\left(\overline \H_a \right)$ is the atomic representation (cf. Definition \ref{atomic_repr_defn}) then from the Corollary \ref{inverse_lim_h_cor} and the Lemma \ref{inverse_faithful_nondegenetate_lem} it follows that any $\mathfrak{Cov}\left(\mathfrak{S} \right)$-morphism $\phi: \overline A' \hookto  \overline A''$ (cf. \ref{infinite_covering_empt}) can be included into the following commutative diagram
\newline
\begin{tikzcd}
\overline A'\arrow[rr, "\phi"] \arrow[rd, "\pi_{\overline A'}"]&& \overline A''\arrow[ld, "\pi_{\overline A''}"]\\
& B\left(\overline \H_a \right) &
\end{tikzcd}
\\
where both $\pi_{\overline A'}$ and $\pi_{\overline A'}$ are atomic representations of $\overline A'$ and $\overline A''$ respectively.
 \end{empt}

 \begin{lemma}\label{sum_hilb_rep_lem} 
 	Let $\left\{\H_\la\right\}_{\la\in\La}$ be a family of Hilbert spaces, and let $\H$ be the norm completion of the algebraic direct sum $\oplus_{\la\in\La}\H_\la$. The bounded net $\left\{a_\a\right\}_{\a\in \mathscr A} \in B\left( \H\right)$ is strongly convergent (cf. Definition \ref{strong_topology_defn}) if and only if for all $\la\in\La$ and  $\zeta \in \H$ given by
 	$$
 	\zeta =\left(0,...~, \underbrace{\zeta_\la}_ {{\la}^{\text{th}}-\mathrm{place}} ,~...,0\right)\in \bigoplus_{\la\in\La}\H_\la
 	$$
 	the net $\left\{a_\a\zeta\right\} \in \H$ is norm convergent. 
 \end{lemma}
 \begin{proof}
 	If 	$\left\{a_\a\right\}_{\a\in \mathscr A} \in B\left( \H\right)$ is strongly convergent then  the net $\left\{a_\a\zeta\right\} \in \H$ is norm convergent for all $\zeta\in \H$.
 	
 	Conversely suppose that $\left\| a_\a\right\| < C$ for all $\a\in \mathscr A$.
 	If	$\xi= \left(\xi_\la\right)_{\la\in \La} \in \H$ 
 	and $\eps > 0$. There is a finite subset   $\La_0\subset \La$ and  $\eta = \left\lbrace \eta_{\la}\in \H_{\la}\right\rbrace_{\la\in \La}\in \oplus_{\la\in\La}\H_\la$ such that
 	\bean
 	\begin{split}
 		\eta_{\la}= \begin{cases}
 			\xi_{\la}& \la\in \La_0\\
 			0& \la\notin \La_0\\
 		\end{cases},\\ \left\|\xi- \eta\right\|< \frac{\eps} {2C}~.	\end{split}
 	\eean
 	For any $\la \in \La_0$ there is $\a_\la \in \mathscr A$ such that
 	\be\label{cunv_2_un_eqn}
 	\bt , \ga \ge \a_\la \quad \Rightarrow \quad \left\|\left(a_\bt- a_{\ga}\right)\xi_{\la}\right\|< \frac{\eps} {2\left|\La_0\right|C}~.
 	\ee
 	If $\a_{\max} \ge a_\la$ for every $\la\in\La_0$ then  one has 
 	\bean
 	\bt , \ga \ge \a_{\max} \quad \Rightarrow \quad \left\|\left(a_\bt- a_{\ga}\right)\xi\right\|< \eps.
 	\eean
 \end{proof}
 
 \begin{lemma}\label{direct_sum_strong_lim_lem}
 	Let $\left\{\H_\la\right\}_{\la\in\La}$ be a family of Hilbert spaces, and let $\H$ be the norm completion of the algebraic direct sum $\oplus_{\la\in\La}\H_\la$. Let $\left\{a_\a\right\}_{\a\in\mathscr A}\in B\left(\H\right)$ be net such that:
 	\begin{itemize}
 		\item there is $C > 0$ such that $\left\|a_\a \right\|<C$ for all $\a\in\mathscr A$,
 		\item $a_\a \H_\la \subset \H_\la$ where the natural inclusion $\H_\la\subset \H$ is implied.
 		\item	for all $\la \in \La$ there is a limit $\lim_{\a \in \mathscr A}a_\a|_{\H_\la} = a_\la$ with respect to the strong topology of $B\left( \H_\a\right)$ (cf. Definition \ref{strong_topology_defn}). 		
 	\end{itemize} 
 	Then the net $\left\{a_\a\right\}$ is convergent  with respect to the strong topology of $B\left( \H\right)$ (cf. Definition \ref{strong_topology_defn}). If
 	$\xi \in \H$ corresponds to a family
 	\be\nonumber
 	\left\{\xi_{\la}\in \H_{\la}\right\}_{\widehat \la\in \La}.
 	\ee
 	then $\left(\lim_{\a \in \mathscr A } a_\a\right) \xi$ is represented by a family
 	$$
 	\left\{ \left( \lim_{\a \in \mathscr A } a_\a|_{\H_\la}\right) \xi_{\la}\in \H_{\la}\right\}_{\la\in \La}
 	$$
 	where the strong limit  $\lim_{\a \in \mathscr A } a_\a$ is implied (cf. Definition \ref{strong_topology_defn}).
 \end{lemma}\label{sum_act_repr_lem}
 \begin{proof}
 	From the Lemma \ref{sum_hilb_rep_lem}  it follows that $\left\{a_\a\right\}$ is convergent  with respect to the strong topology of $B\left( \H\right)$. If $\xi= \left(\xi_\la\right)_{\la\in \La} \in \H$ 
 	and $\eps > 0$ then there is a finite subset   $\La_0\subset \La$ and  $\eta \in \oplus_{\la\in\La}\H_\la$ such that
 	\be\label{cunv_unp_eqn}
 	\eta_{\la}= \begin{cases}
 		\xi_{\la}& \la\in \La_0\\
 		0& \la\notin \La_0\\
 	\end{cases}\text{ AND }	\eta =\left(\eta_{\la}\in \H_{\la}\right)_{\la\in \La}	 \quad \Rightarrow \quad \left\|\xi- \eta\right\|< \frac{\eps} {2C}.
 	\ee
 	For any $\la \in \La_0$ there is $\a_\la \in \mathscr A$ such that
 	\be\label{cunv_2p_un_eqn}
 	\bt , \ga \ge \a_\la \quad \Rightarrow \quad \left\|\left(a_\bt- a_{\ga}\right)\xi_{\la}\right\|< \frac{\eps} {2\left|\La_0\right|C}
 	\ee
 	If $\a_{\max} \ge a_\la$ for every $\la\in\La_0$ then  one has 
 	\bean
 	\bt , \ga \ge \a_{\max} \quad \Rightarrow \quad \left\|\left(a_\bt- a_{\ga}\right)\xi\right\|< \eps.
 	\eean
  	If  $\eta =\left\lbrace\eta_\la \right\rbrace  = \left(\lim_{\a \in \mathscr A } a_\a\right) \xi$ then from  $a_\a \H_\la \subset \H_\la$ it follows that $a_\a \xi_\la = a_\a|_{\H_\la} \xi_{\la}\in \H_{\la}$ and $\eta_\la =  \left( \lim_{\a \in \mathscr A } a_\a|_{\H_\la}\right) \xi_{\la}$.
 \end{proof}
 \begin{empt}\label{atomic_faithful_empt}
 	Under the situation of the Definition  \ref{algebraical_finite_covering_category_eqn}. Let $\widehat \pi : \widehat A \to B\left(\widehat{\H}\right)$ be a representation. From the Lemma \ref{bt_lim_eqn} it turns out that any $\overline a \in K\left( \overline A\right) $ is given by
 	$$
 	\overline a =\bt\text{-}\lim_{\la \in \La} a_\la
 	$$
 	where $\bt\text{-}\lim$ is the limit with respect to the strict topology of $M\left( \overline{A}\right)$ (cf. Definition \ref{strict_topology_defn}).	If 
 	$a_\la = a^1_\la - a^2_\la + ia^3_\la -i a^4_\la$ where $a^1_\la, a^2_\la, a^3_\la, a^4_\la \in A_\la$ are positive elements (cf. \eqref{four_decompositon_eqn}, then for $j = 1,..., 4$ one has a decreasing net $\left\{\widehat\pi\left( a^j_\la\right)\right\} $. From the Proposition \ref{increasing_convergent_w_lem} all these  nets are convergent with respect to the strong topology of $B\left(\widehat{\H}\right)$ (cf. Definition \ref{strong_topology_defn}). It follows the net $\{\widehat\pi\left( a_\la\right) \}$ is convergent with respect to the strong topology of $B\left(\widehat{\H}\right)$.
 	We define a map
 	\be\label{ext_repr_eqn}
 	\begin{split}
 		ext(\widehat \pi) : K\left( \overline{A}\right) \to B \left(\widehat \H\right),\\
 		\overline a \mapsto s\text{-}\lim_{\la\in\La}\widehat\pi\left( a_\la\right)
 	\end{split}
 	\ee
 	where $ s\text{-}\lim$ means the limit with respect to the strong topology of $B\left(\widehat{\H}\right)$.
 	From the equations 
 	\bean
 	\bt\textit{-}\lim_{\la\in\La}\left(  a'_\la + a''_\la\right)  =\left(\bt\textit{-}\lim_{\la\in\La} a'_\la \right)+\left(\bt\textit{-}\lim_{\la\in\La} a''_\la \right),\\
 	s\text{-}\lim_{\la\in\La}\left( \widehat\pi\left( a'_\la\right)+\widehat\pi\left( a''_\la\right)\right) =\left(  s\text{-}\lim_{\la\in\La}\widehat\pi\left( a'_\la\right)\right)+ \left(  s\text{-}\lim_{\la\in\La}\widehat\pi\left( a''_\la\right)\right) ,\\
 	\bt\textit{-}\lim_{\la\in\La} a'_\la a''_\la =\left(\bt\textit{-}\lim_{\la\in\La} a'_\la \right)\left(\bt\textit{-}\lim_{\la\in\La} a''_\la \right),\\
 	s\text{-}\lim_{\la\in\La}\widehat\pi\left( a'_\la\right)\widehat\pi\left( a''_\la\right)=\left(  s\text{-}\lim_{\la\in\La}\widehat\pi\left( a'_\la\right)\right) \left(  s\text{-}\lim_{\la\in\La}\widehat\pi\left( a''_\la\right)\right) ,\\
 	\left( \bt\textit{-}\lim_{\la\in\La} a_\la\right)^* = \bt\textit{-}\lim_{\la\in\La} a^*_\la ,\\
 	\left(s\text{-}\lim_{\la\in\La}\widehat\pi\left( a_\la\right)\right)^* = s\text{-}\lim_{\la\in\La}\widehat\pi\left( a_\la\right)^*
 	\eean
 	It follows that the map \eqref{ext_repr_eqn} is a $*$-homomorphism. Since $ K\left( \overline A\right)$ is dense in $\overline A$ one has a representation
 	\be\label{ext_repro_eqn}
 	ext(\widehat \pi) : \overline{A}\to B \left(\widehat \H\right).
 	\ee
 \end{empt}
 
\begin{lemma}\label{infinite_faithful_nondegererate_lem}
Under the situation of the Definition \ref{algebraical_finite_covering_category_defn}  if an equivariant representation $\widehat \pi : \widehat A  \bydef C^*$-$\varinjlim_{\la \in \La} A_\la \hookto B\left(\widehat{\H}\right)$ (cf. Definition \ref{equivariant_representation_defn}) is faithful and non-degenerate (cf. Definitions \ref{faithful_representation_defn} and \ref{nondegenerate_repr_defn}) then the given by \eqref{ext_repro_eqn} representation is faithful and non-degenerate.
\end{lemma}
\begin{proof}
The notation of the Definition  \ref{algebraical_finite_covering_category_eqn} is used. If $\overline a \in  K\left( \overline A\right)$ then there and 	$a \bydef 	\bt\text{-}\sum_{g\in \widehat{G}}g \overline a \in  K\left( A\right)$ then there is $b \in A$ such that $ab\neq 0$. From 	$a \bydef 	\bt\text{-}\sum_{g\in \widehat{G}}g \overline a$ it follows that there is a finite subset $G_0\subset  \widehat{G}$ such that $\left( \sum_{g\in  {G}_0}g \overline a\right) b \neq 0$.  Moreover there is $g \in G_0$ such that $ (g\overline a)b \neq 0$. The representation $\widehat \pi$ is faithful, so there is $\xi \in \widehat{\H}$ such that $(g\overline a)b \xi \neq 0$, so one has $g^{-1}\left((g\overline a)b \xi  \right) = \overline a b g^{-1}\xi \neq 0$ (cf. equations \eqref{equivt_eqn} and \eqref{equivt_eqn}). It turns out that $\overline a \widehat{\H} \neq \{0\}$. But $K\left( \overline A\right) $ is dense in $\overline A$, so  $\overline a \widehat{\H} \neq \{0\}$ for all  $\overline a \in \overline A$, i.e. the given by \eqref{ext_repro_eqn} representation is faithful. Since $A_\la$ is generated by given by \eqref{basic_cov_cl_eqn} elements
  	\bean
a_\la =\bt\text{-} \sum_{	g \in \ker\left( \widehat{G}\to G_\la\right) }g \overline a
\eean there is  $\overline a \in K\left( \overline A\right) $ such that $\left( \bt\text{-}\sum_{g \in \ker\left(\widehat{G}\to G_\la\right)}g \overline a\right) \xi \neq 0$. It turns out that there is a finite subset $G_0\subset  \ker\left(\widehat{G}\to G_\la\right)$ such that $\left( \sum_{g\in  {G}_0}g \overline a\right) \xi \neq 0$, and one has $g \in G_0$ satisfying $g \overline a \xi \neq 0$, i.e. $K\left( \overline A\right)\xi \neq 0$. So the given by \eqref{ext_repro_eqn} representation is nondegenerate.
\end{proof}

\begin{lemma}\label{strict_strong_1_lem}
Let	$\left( A, \overline A, \widehat{G}\right)$ be a {pre}-{covering} of the category $\mathfrak{S}$ (cf. Definition  \ref{algebraical_finite_covering_category_eqn}). If  $\overline\pi: \overline{A}: \to B\left(\overline \H \right)$ is a faithful (cf. Definition \ref{faithful_representation_defn}), nondegenerate (cf. Definition \ref{nondegenerate_repr_defn}) representation. If  $\overline{a}\in K\left( \overline{A}\right) $,  $\la \in \La$ and  $a_\la\in K\left(A_\la \right) $ is the sum of convergent with respect  to the strict topology of $M\left( \overline A\right)$ (cf. Definition \ref{strict_topology_defn}) series 
	$$
\bt\text{-}\sum_{g \in \ker\left(\widehat{G}\to G_\la\right)}g \overline a  	
$$
 then 
$\overline\pi\left(a_\la\right)$ is the sum of the convergent  with respect to the strong topology of 	 $B\left(\overline \H \right)$ (cf. Definition \ref{strong_topology_defn}) series
		$$
s\text{-}	\sum_{g \in \ker\left(\widehat{G}\to G_\la\right)}\overline\pi\left( g \overline a  \right)	
	$$
\end{lemma}
\begin{proof}
	If $\overline a \in K\left( \overline A\right)_+$ is positive then this Lemma follows from the \ref{strict_strong_lem} one. General case is a consequence of the equation \eqref{four_decompositon_eqn}.

\end{proof}
\begin{lemma}\label{strict_strong_2_lem}
Let	$\left( A,\overline A, \widehat{G}\right)$ be a {pre}-{covering} of the category $\mathfrak{S}$ (cf. Definition  \ref{algebraical_finite_covering_category_eqn}). If $\overline\pi: \overline{A}: \to B\left(\overline \H \right)$ be a faithful (cf. Definition \ref{faithful_representation_defn}), nondegenerate (cf. Definition \ref{nondegenerate_repr_defn}) representation.  Under the hypotheses of the Lemma \ref{bt_lim_lem} if one has a given by \eqref{bt_lim_eqn} strict limit
	$$
\forall \overline a \in K\left( \overline A\right)  \quad \overline a =\bt\text{-}\lim_{\la \in \La}a_\la$$
	 (cf. Definition \ref{strict_topology_defn}) then there is a limit
	 $$
\overline \pi \left( \overline a \right) = s\text{-}\lim_{\la\in \La }	 \pi \left(  a_\la \right)
	 $$
with respect to the strong topology of 	 $B\left(\overline \H \right)$ (cf. Definition \ref{strong_topology_defn}).
\end{lemma}
\begin{proof}
The proof of this lemma is similar the \ref{strict_strong_1_lem} one.
\end{proof}
\subsubsection{Induced representations}
\paragraph{}
	If  $\left( A, \overline A, \widehat{G}\right)$ is a {quasi}-{covering of algebraical finite covering category}  $$\mathfrak{S}\bydef \left\{\left\{A_\la\right\}_{\la\in \La}, \left\{\pi^\mu_\nu : A_\mu \hookto A_\nu\right\}_{\substack{\mu, \nu \in \La\\\mu \le \nu}}\right\}$$ (cf. Definition   \ref{infinite_quasicovering_defn}) then from \eqref{basic_cov_cl_eqn} it follows that there is a sesquilinear $A$-valued product
	\be\label{inf_hilb_product_eqn}
	\begin{split}
		\left\langle \cdot, \cdot  \right\rangle_A : K\left( \overline{A}\right)\times K\left( \overline{A}\right)  \to K\left( A\right) ,\\
		\left(\overline{a}, \overline{b}\right) \mapsto \bt\text{-}\sum_{g \in \widehat{G}}g \left( \overline a^* \overline{b}\right)= \desc\left(\overline a^* \overline{b} \right) 
	\end{split}
	\ee
	where the above series is convergent with respect to the strict topology (cf. Definition \ref{strict_topology_defn}) of $M\left( \overline{A}\right)$ and $\desc$ is the minimal descent (cf. Definition \ref{infinite_desc_defn}). The  Pedersen's ideal  $K\left( \overline{A}\right)$ possesses a structure of $\overline{A}$-$A$-bimodule which comes from the structure of $\overline{A}$-$A$-bimodule $\overline{A}$. A completion  of  $K\left( \overline{A}\right) $ with respect to  a norm
	\begin{equation}\label{inf_hilb_norm_eqn}
		\begin{split}
			\left\| \overline{a}\right\| = \sqrt{\left\| \left\langle \overline{a}, \overline{a}  \right\rangle_A\right\|}
		\end{split}
	\end{equation}
	will be  denoted by $\mathscr L^2\left( \overline{A}\right)$.

	The natural action $M\left( \overline{A}\right) \times K\left( \overline{A}\right) \to  K\left( \overline{A}\right)$ can be extended up to the action $M\left( \overline{A}\right)\times \mathscr L^2\left( \overline{A}\right) \to \mathscr L^2\left( \overline{A}\right)$ since $ K\left( \overline{A}\right)$ is a dense $\overline{A}$-module of $ \mathscr L^2\left( \overline{A}\right)$, so there is a homomorphism 
	\be\label{mult_adj_eqn}
\varphi : M\left( \overline{A}\right) \hookto \End^*_A \left(\mathscr L^2\left( \overline{A}\right)\right) 
\ee
 where  $\End^*_A \left(\mathscr L^2\left( \overline{A}\right)\right)$ is $C^*$-algebra of adjointable operators (cf. Definition \ref{adjointable_operator_defn}).

\begin{empt}\label{infinite_hilb_mod_empt}
	If $\left(A, \overline{A},\widehat{G}\right)$ is good (cf. Definition \ref{good_defn}), and a triple $\left(A, \widetilde{A}, G\left(\left.\widetilde{A}~\right| A\right)\right)$  the  {infinite noncommutative covering} of $\mathfrak{S}$ (cf. Definition \ref{infinite_noncommutative_covering_defn}) then  one can define a product
	\be\label{inf_hilb_prod_eqn}
	\begin{split}
		\left\langle \cdot, \cdot  \right\rangle_A : K\left( \widetilde{A}\right) \times K\left( \widetilde{A}\right) \to K\left(A \right) 
	\end{split}
	\ee
	which is a restriction of the \eqref{inf_hilb_product_eqn} one. We suppose that  $K\left( \widetilde{A}\right)$ possesses a structure of $\widetilde{A}$-$A$-bimodule which comes from the structure of $\widetilde{A}$-$A$-bimodule $\widetilde{A}$. Then a completion $ \mathscr L^2\left( \widetilde{A}\right)$ of  $ K\left( \widetilde{A}\right)$ with respect to  a norm
	\begin{equation}\label{inf_hilb_tnorm_eqn}
		\begin{split}
			\left\| \widetilde{a}\right\| = \sqrt{\left\| \left\langle \widetilde{a}, \widetilde{a}  \right\rangle_A\right\|}
		\end{split}
	\end{equation}
\end{empt}

\begin{definition}\label{infinite_hilb_mod_defn} 
	Under the hypotheses \ref{infinite_hilb_mod_empt} we say that the $C^*$-Hilbert $A$-module $\mathscr L^2\left( \widetilde{A}\right) $ is the \textit{$C^*$-Hilbert module   associated with the  infinite noncommutative covering} $\left(A, \widetilde{A}, G\left(\left.\widetilde{A}~\right| A\right)\right)$ of $\mathfrak{S}$ (cf. Definition \ref{infinite_noncommutative_covering_defn}). The $C^*$-Hilbert module  $\mathscr L^2\left( \widetilde{A}\right)$  is $\widetilde A$-$A$-bimodule, so  we use following notations for it
	\bean
	\mathscr L^2\left( \widetilde{A}\right)_A\bydef \mathscr L^2\left( \widetilde{A}\right), \quad 
	_{\widetilde{A}}\mathscr L^2\left( \widetilde{A}\right)_A\bydef \mathscr L^2\left( \widetilde{A}\right).
	\eean
\end{definition}

\begin{remark}
Similarly to \eqref{mult_adj_eqn} one has an injective $*$-homomorphism
\be\label{mult_adjt_eqn}
\varphi : M\left( \widetilde{A}\right) \hookto \End^*_A \left(\mathscr L^2\left( \widetilde{A}\right)\right) 
\ee
\end{remark}

\begin{remark}\label{infinite_hilb_inv_rem}
	The involution on $\widetilde{A}$ induces an involution on $\mathscr L^2\left( \widetilde{A}\right)$, we denote it by
	\be\label{infinite_hilb_inv_eqn}
	x \mapsto x^*
	\ee
\end{remark}

\begin{defn}\label{induced_repr_inf_defn}
	Let 
	$
	\left(A, \widetilde{A}, G\left(\left.\widetilde{A}~\right| A\right)\right)
	$
	be the  {infinite noncommutative covering}  of  $$\mathfrak{S}=\left\{\left\{A_\la\right\}_{\la\in \La}, \left\{\pi^\mu_\nu : A_\mu \hookto A_\nu\right\}_{\substack{\mu, \nu \in \La\\\mu \le \nu}}\right\}.$$ 
	(cf. Definition \ref{infinite_noncommutative_covering_defn}). 	 Suppose that 	  $_{\widetilde{A}}\mathscr L^2\left( \widetilde{A}\right)_A$ the {$C^*$-Hilbert   associated the  infinite noncommutative covering} $\left(A, \widetilde{A}, G\left(\left.\widetilde{A}~\right| A\right)\right)$.  
	If $\rho: A \to B\left( \H\right)$ is a representation and   $\widetilde{\rho}=_{\widetilde{A}}\mathscr L^2\left( \widetilde{A}\right)_A-\Ind^A_{\widetilde{A}}\rho: \widetilde{A}\to B\left(\widetilde{\H} \right) $ is given by \eqref{induced_representation_eqn}, i.e. $\widetilde{\rho}$ is the induced representation (cf. Definition \ref{induced_representation_defn}). 
	The representation  $\widetilde{\rho}:\widetilde{A} \to  B\left( \widetilde{\H}\right)$ is said to be \textit{induced} by  $\rho$. We also say that  $\widetilde{\rho}$ is  \textit{induced by a pair} $\left( \rho, \left( A, \widetilde{A}, G\left(\left.\widetilde{A}~\right|A \right) \right)\right) $. 
\end{defn}

\begin{remark}

The given by  \ref{inf_hilb_prod_eqn}  pairing can be uniquely extended up to the following product
\be\label{inf_hilb_f_prod_eqn}
\begin{split}
	\left\langle \cdot, \cdot  \right\rangle_A : \mathscr L^2\left( \widetilde{A}\right)_A ~\otimes_A \mathscr L^2\left( \widetilde{A}\right)_A \to A,\\
	\left\langle \widetilde{a}, \widetilde{b}  \right\rangle_A  \bydef\bt\text{-} \sum_{g \in  G\left(\left.\widetilde{A}_\pi~\right|A \right)} g \left(\widetilde{a}^* \widetilde{b}  \right) \in A 
\end{split}
\ee
where the above series is convergent with respect to the strict topology (cf. Definition \ref{strict_topology_defn}) of $M\left( \widetilde{A}\right)$.
\end{remark}
\begin{rem}\label{ap_act_hilb_rem} 
	If  $\widetilde{\rho}:\widetilde{A} \to  B\left( \widetilde{\H}\right)$ is   {induced} by $\left( \rho, \left( A, \widetilde{A}, G\left(\left.\widetilde{A}~\right|A \right) \right)\right) $ then $\widetilde{\H}$ is the completion of the pre-Hilbert space $ \mathscr L^2\left( \widetilde{A}\right) \otimes_A \H$ with given by
	\be\label{comp_hilb_eqn}
	\left(\widetilde{a} \otimes \xi, \widetilde{b} \otimes \eta \right)_{\widetilde{\H}} = \left( \xi, \left\langle \widetilde{a} , \widetilde{b} \right\rangle_A \eta\right)_{\H} \quad \forall \widetilde{a} \otimes \xi,~ \widetilde{b} \otimes \eta \in \mathscr L^2\left( \widetilde{A}\right) \otimes_A \H
	\ee
	scalar product (cf. equation \eqref{hilb_prod_eqn}). The action $\widetilde{A}\otimes \widetilde \H\to \widetilde \H$ corresponds to the completion of the following action
	\be\label{comp_hilb_actf_eqn}
	\begin{split}
		\widetilde{A} \times	\left( \mathscr L^2\left( \widetilde{A}\right) \otimes_A \H \right) \to  \mathscr L^2\left( \widetilde{A}\right) \otimes_A \H  ;\\
		\left(\widetilde a, \widetilde b\otimes \xi \right) \mapsto \widetilde a\widetilde b\otimes \xi.
	\end{split}
	\ee
	If $\widetilde \H$ is the Hilbert norm completion of $K\left(\widetilde{A}\right) \otimes \H$ then the action \eqref{comp_hilb_actf_eqn} uniquely defines an action
	\be\label{comp_hilbh_actf_eqn}
	\begin{split}
		\widetilde{A} \times	\widetilde \H \to \widetilde \H,
	\end{split}
	\ee 
	so there is the inclusion
	\be\label{comp_hilbr_actf_eqn}
	\begin{split}
		\widetilde{A} \subset B\left( \widetilde \H\right).
	\end{split}
	\ee
	
\end{rem}

\begin{lemma}\label{faithful_inf_lem}
	Under the hypotheses of the Definition \ref{induced_repr_inf_defn} if $\rho$ is faithful, then  $\widetilde{\rho}$ is faithful.
\end{lemma}
\begin{proof}
	If $\widetilde\rho$ is not faithful then there is a positive $\widetilde a \in \widetilde{A}\setminus\{0\}$ such that
	$$
	\forall \widetilde \xi \in \widetilde{\H} \quad \left(\widetilde a \widetilde \xi, \widetilde a\widetilde \xi \right)_{\widetilde\H} = 0.
	$$
	If $f_\eps$ is given by \eqref{f_eps_eqn} then there is $0 < \delta < 1$ such that $f_\delta\left( \widetilde a\right)\neq 0$.  Moreover from $f_\delta\left( \widetilde a\right) <  \widetilde a$ it turns out that
	$$
	\forall \widetilde \xi \in \widetilde{\H} \quad \left(f_\delta\left( \widetilde a\right) \widetilde \xi, f_\delta\left( \widetilde a\right) \widetilde \xi \right)_{\widetilde\H}\le \left(\widetilde a \widetilde \xi, \widetilde a\widetilde \xi \right)_{\widetilde\H} = 0.
	$$
	From the Lemma \ref{pedersen_eps_lem} one has $f_\delta\left( \widetilde a\right)\in  \mathscr L^2\left( \widetilde{A}\right)$ and  if $\widetilde \xi \bydef  f_\delta\left( \widetilde a\right)\otimes  \xi \in  \mathscr L^2\left( \widetilde{A}\right)\otimes_A \H \subset \widetilde \H$ then from \eqref{hilb_prod_eqn} it follows that
	\be\label{inf_ind_c_eqn}
	\begin{split}
		\forall \xi \in \H\quad 	\left(f_\delta\left( \widetilde a\right) \left( f_\delta\left( \widetilde a\right)\otimes \xi \right) , f_\delta\left( \widetilde a\right) \left( f_\delta\left( \widetilde a\right)\otimes \xi \right) \right)_{\widetilde\H}=\\= \left( \xi, \left\langle f^2_\delta\left( \widetilde a\right), f^2_\delta\left( \widetilde a\right)\right\rangle_A \xi\right)_\H=0. 
	\end{split}
	\ee
	On the other hand 
	$$
	a_\delta\bydef	\left\langle f^2_\delta\left( \widetilde a\right), f^2_\delta\left( \widetilde a\right)\right\rangle_A = \sum_{g \in G\left(\left.\widetilde{A}~\right|A \right)} f^4_\delta\left( \widetilde a\right)\in A_+\setminus \{0\}.
	$$
	Since $\rho$ if faithful there is $\xi \in \H$ such that $a_\delta \xi \neq 0$. So one has
	$$
	\left( \xi, \left\langle f^2_\delta\left( \widetilde a\right), f^2_\delta\left( \widetilde a\right)\right\rangle_A \xi\right)_\H= \left( \xi, a_\delta \xi\right)_\H > 0.
	$$
	Above equation contradicts with \eqref{inf_ind_c_eqn} so for any positive $\widetilde a \in \widetilde{A}\setminus\{0\}$ there is $\widetilde\xi \in \widetilde{\H}$ such that $\left(\widetilde a \widetilde \xi, \widetilde a\widetilde \xi \right)_{\widetilde\H} > 0$, i.e. $\widetilde{\rho}$ is faithful.
\end{proof}
\begin{rem}\label{a_act_hilb_rem} 
	There is the action of $G\left(\left.\widetilde{A}~\right|A \right)$ on $\widetilde{\H}$ which comes from the natural action of $G\left(\left.\widetilde{A}~\right|A \right)$ on the $\widetilde{A}$-bimodule $K\left( \widetilde{A}\right) $. If the representation $\widetilde A \to 	B\left( \widetilde{\H} \right)$ is faithful, so an action of  $G\left(\left.\widetilde{A}~\right|A \right)$ on $\widetilde{A}$ is given by
	\be\label{ind_h_ind_a_eqn}
	\left( g  \widetilde a\right) \xi =   g \left(  \widetilde a  \left(  g^{-1}\widetilde\xi\right) \right); ~ \forall  g  \in {G}, ~ \forall\widetilde a  \in \widetilde{A}, ~\forall\widetilde \xi \in \widetilde{\H}.
	\ee
\end{rem}

\begin{lemma}\label{induced_nondegenerate_inf_lem}
	Under the hypotheses of the Definition \ref{induced_repr_inf_defn}	if  the representation  $\rho: A \hookto B\left(\H \right) $ is nondegenerate (cf. Definition \ref{nondegenerate_repr_defn}) then the induced representation  $\widetilde{\rho}: \widetilde{A} \to B\left( \widetilde{\H} \right)$ (cf. Definition \ref{induced_repr_inf_defn}) is also nondegenerate. 
\end{lemma}
\begin{proof}
	This lemma can be proven as well as \ref{induced_nondegenerate_lem} one.
\end{proof}

\begin{empt}\label{hurewicz_n_to_h_constr_empt}
	Let us fix $\la\in \La$, and 
	let $\H_\la$ be a Hilbert completion of $A_{\la} \otimes_A \H$ which is constructed in the section \ref{induced_repr_fin_sec}. From  the Remark \ref{dense_inf_rem} it follows that both  $ \widetilde{\lift}_\la\left(A_\la \right)\widetilde{A}$ and $\widetilde{A}~ \widetilde{\lift}_\la\left(A_\la \right)$ are dense subsets of $	\widetilde{A}$. Taking into account $	K\left(\widetilde{A}\right)\widetilde{A} =	K\left(\widetilde{A}\right)$ one has
	a natural inclusion of pre-Hilbert spaces
	\begin{equation}\label{tensor_n_eqn}
		\begin{split}
			K\left(\widetilde{A}\right) \otimes_A \H = \left( K\left(\widetilde{A}\right) \otimes_{A_{\la}}   A_{\la}\right)  \otimes_A \H=\\
			=K\left(\widetilde{A}\right) \otimes_{A_{\la}} \left( A_{\la} \otimes_A \H\right)\hookto
			K\left(\widetilde{A}\right)\otimes_{A_{\la}} \H_{\la}.
		\end{split}
	\end{equation}
	such that $K\left(\widetilde{A}\right) \otimes_A \H$ is dense in $K\left(\widetilde{A}\right)\otimes_{A_{\la}} \H_{\la}$ with respect to pre-Hilbert norm, i.e. the inclusion \ref{tensor_n_eqn} induces the isomorphism of Hilbert completions. So for all $\la\in\La$ there is a dense inclusion
	\be\label{kla_dens_eqn}
	K\left(\widetilde{A}\right)\otimes_{A_{\la}} \H_{\la}\hookto \widetilde \H. 
	\ee
\end{empt}

\section{Coverings of *-algebras}

\subsection{Coverings of pro-$C^*$-algebras}
\paragraph*{} Here we discuss infinite coverings of operator *-algebras which can contain unbounded operators. If $A$ is pro-$C^*$-algebra then denote by $b\left(A \right)\subset A$ the $C^*$-algebra  of bounded operators (cf. the Definition \ref{pro_bound_defn} and the Proposition \ref{pro_bound_prop}).
\begin{empt}\label{comp_pro_empt}
	Let $\La$ be a  directed set (cf. Definition \ref{directed_set_defn}) such that there is the unique minimal element $\la_{\min} \in \La$.
	Let $A$ be a pro-$C^*$-algebra.  Let us consider a set noncommutative finite-fold coverings  pro-$C^*$-algebras $\mathfrak{S}=\left\{ \pi_{\la}:A \hookto A_{\la}\right\}_{\la \in \La}$ (cf. Definition \ref{pro_fin_defn}) indexed by  $\La$ such that $A_{\la_{\mathrm{min}}}= A$, and $\pi_{\la_{\mathrm{min}}}= \Id_A$. Suppose that there is a  category 	
	\be\label{comp_pro_eqn}
	\begin{split}
		\mathfrak{S} \bydef \left(\left\{\pi_\la: A \hookto A_\la \right\}_{\la \in \La}, \left\{\pi^\mu_\nu: A_\mu \hookto A_\nu\right\}_{\substack{\mu, \nu \in \La\\ \nu > \mu}}\right),\\
		\text{or simply} \quad 	\mathfrak{S} \bydef \left(\left\{\pi_\la: A \hookto A_\la \right\}, \left\{\pi^\mu_\nu\right\}\right).
	\end{split}
	\ee
	of noncommutative finite-fold coverings  pro-$C^*$-algebras (cf. Definition \ref{pro_fin_defn}). 
	Suppose that  $\left(B, \widetilde{B}, G\left(\left.\widetilde{B}~\right| B\right)\right)$ is 
	{infinite noncommutative covering} of 
\be\label{pro_cs_cat_eqn}	
	\mathfrak{S}_B=\left\{\left\{B_\la\right\}_{\la\in \La}, \left\{\widetilde \pi^\mu_\nu : B_\mu \hookto B_\nu\right\}_{\substack{\mu, \nu \in \La\\\mu \le \nu}}\right\}
\ee	
	 (cf. Definition \ref{infinite_noncommutative_covering_defn}).
Assume that for each $\la\in \La$ a $C^*$-algebra  $B_\la \subset b\left(A_\la \right)$ is an essential ideal of the $C^*$-algebra $b\left(A_\la \right) \subset A_\la$ of bounded elements of $A$, such that 
	\bean
\widetilde \pi^\mu_\nu =  \left.\pi^\mu_\nu\right|_{B_\mu}.
	\eean
\end{empt}

\begin{definition}\label{comp_pro_defn}
	The the given by  \eqref{comp_pro_eqn} set $\mathfrak{S}$
	is said to be an \textit{algebraical  finite covering category of pro-$C^*$-algebras}. 
	We write $\mathfrak{S} \in$ pro-$\mathfrak{FinAlg}$.
\end{definition}
\begin{empt}\label{inv_pro_lim_empt}
Let $\mathfrak{S}$ algebraical  finite covering category of pro-$C^*$-algebras, and let $\mathfrak{S}_B$ be a given by \eqref{comp_pro_eqn}  { algebraical  finite covering category}. 	
	For any $\la \in \La$ there is a given by the Definition  \ref{infinite_lift_defn} $\la$-lift
		$\widetilde{\lift}_\la :B_\la \to M\left( \widetilde{B} \right) $. For all $\la\in \La$  any seminorm $p_\a: M\left( \widetilde{B} \right)\to\R$ induces a seminorm
		\be\label{pro_ind_semi_eqn}
		 p^\la_\a\bydef p_\a \circ \widetilde{\lift}_\la: B_\la \to\R
		 \ee
		 A $\widetilde G$-invariant  family $\left\{p_\a\right\}_{\a \in \mathscr A}$  of $C^*$-seminorms on $ M\left( \widetilde{B} \right)$ is said to be \textit{admissible} if for all $\la\in \La$ a pro-$C^*$ algebra $A_\la$  is a completion of $B_\la$  with respect to the family of given by \eqref{pro_ind_semi_eqn} $C^*$-seminorms. For example the equation  \eqref{seminorms_inv_lim} yields an admissible family of $C^*$-seminorms. A completion of $\widetilde B$ with respect to an admissible family of seminorms is said to be an \textit{admissible} pro-$C^*$-algebra.
\end{empt}

\begin{definition}\label{inv_pro_lim_defn}
In the described in \ref{inv_pro_lim_empt} situation we say that the triple
$\left(A, \widetilde{A}, G\right)$ is an 
 \textit{inverse noncommutative limit of the category $\mathfrak{S}$ of pro-$C^*$-algebras} if $\widetilde{A}$ is a maximal among admissible pro-$C^*$-algebras.
\end{definition}

\subsection{Coverings of bounded operator $*$-algebras}
\begin{empt}\label{inf_oa_empt}
	Let $\La$ be a directed set (cf. Definition \ref{directed_set_defn}) such that there is a unique minimal element $\la_{\min} \in \La$.
	Let $A$ be a dense *-subalgebra of $C^*$-algebra $B$.  Let us consider a family of   noncommutative finite-fold coverings of bounded operator *-algebras $\mathfrak{S}=\left\{ \pi_{\la}:A \hookto A_{\la}\right\}_{\la \in \La}$ (cf. Definition \ref{fin_oa_defn}) indexed by  $\La$ such that $A_{\la_{\mathrm{min}}}= A$, and $\pi_{\la_{\mathrm{min}}}= \Id_A$. Suppose that there is a  category 	
	\be\label{comp_oa_eqn}
	\begin{split}
		\mathfrak{S} \bydef \left(\left\{\pi_\la: A \hookto A_\la \right\}_{\la \in \La}, \left\{\pi^\mu_\nu: A_\mu \hookto A_\nu\right\}_{\substack{\mu, \nu \in \La\\ \nu > \mu}}\right),\\
		\text{or simply} \quad \mathfrak{S} \bydef \left(\left\{\pi_\la: A \hookto A_\la \right\}, \left\{\pi^\mu_\nu\right\}\right).
	\end{split}
	\ee
	of noncommutative finite-fold covering of bounded operator *-algebras (cf. Definition \ref{fin_oa_defn}). 	Suppose that  $\left(B, \widetilde{B}, G\left(\left.\widetilde{B}~\right| B\right)\right)$ is 
	{infinite noncommutative covering} of 
	\be\label{bound_cs_cat_eqn}	
	\mathfrak{S}_B=\left\{\left\{B_\la\right\}_{\la\in \La}, \left\{\widetilde \pi^\mu_\nu : B_\mu \hookto B_\nu\right\}_{\substack{\mu, \nu \in \La\\\mu \le \nu}}\right\}
	\ee	
	(cf. Definition \ref{infinite_noncommutative_covering_defn}).
	For each $\la\in \La$  denote by  $B_\la$ the $C^*$-norm completion of $A_\la$ and $\pi^\mu_\nu = \left.\widetilde\pi^\mu_\nu\right|_{A_\mu}$ 
\end{empt}
\begin{definition}\label{comp_oa_defn}
	The the given by  \eqref{comp_pro_eqn} set $\mathfrak{S}$
	is said to be an \textit{algebraical  finite covering category of  bounded operator *-algebras}. 
\end{definition}

\begin{definition}\label{inf_cov_oa_sub_defn}
	Under the hypotheses \ref{inf_oa_empt} we say that  $G$-invariant *-subalgebra $\widetilde A'$ of $\widetilde B$ is \textit{subordinated  to}  $\mathfrak{S}$ if for all $
	\la\in\La$ one has
	\bea\label{inf_cov_oa_sub1_eqn}
\widetilde{\lift}_\la\left( 	A_\la\right)  \widetilde A'\subset \widetilde A',\\
	\label{inf_cov_oa_sub2_eqn}		\widetilde A'\widetilde{\lift}_\la\left( 	A_\la\right)\subset \widetilde A',\\
	\label{inf_cov_oa_sub3_eqn}		\forall \widetilde b \in K\left(\widetilde B \right)
	\quad \forall \widetilde a \in \widetilde A'\quad \exists  \widetilde a', \widetilde a'' \in \bigcup_{\la\in\La} \widetilde{\lift}_\la\left( 	A_\la\right)\quad \widetilde b \widetilde a = \widetilde b \widetilde a', \quad \widetilde a \widetilde b = \widetilde a'' \widetilde b.
	\eea
	where $K\left(\widetilde B \right)$ is the Pedersen's ideal of $\widetilde B$ (cf. Definition \ref{pedersen_ideal_defn}) and $\widetilde{\lift}_\la$ is the $\la$-lift (cf. Definition \ref{infinite_lift_defn}).
\end{definition}
\begin{remark}\label{inf_cov_oa_sub_rem}
The equation \eqref{inf_cov_oa_sub1_eqn} implies following inclusions:
\bean
\widetilde A'\subset \widetilde B\subset M\left(\widetilde B \right) ,\\
\forall \la \in \La \quad A_\la \subset B_\la \subset M\left(\widetilde B \right) \quad \text{(cf. (a) of the Definition \ref{good_defn})}.
\eean
\end{remark}

\begin{remark}
	Any union of simply ordered by inclusion  set of { subordinated   to}  $\mathfrak{S}$ *-subalgebras is a { subordinated   to}  $\mathfrak{S}$ *-subalgebra.	 From the Zorn's lemma (cf. \ref{zorn_thm}) it follows that  there are  \ref{spec_loc_lim_defn} maximal  {subordinated   to}  $\mathfrak{S}$ *-subalgebras. 
\end{remark}

\begin{definition}\label{inf_cov_oa_defn}
	Consider the situation of the Definition \ref{inf_cov_oa_sub_defn}.
	If $\widetilde A\subset \widetilde B$  is a maximal {subordinated   to}  $\mathfrak{S}$ *-subalgebra then
	we say that the triple
	$
	\left(A,  \widetilde A, G \right)$ is the  \textit{inverse noncommutative limit} or \textit{infinite noncommutative covering}  of $\mathfrak{S}$.
\end{definition}
\begin{remark}\label{inf_cov_oa_rem}
	The set of maximal  {subordinated   to}  $\mathfrak{S}$ algebras should not be a singleton, thus an operator algebra  $\widetilde A$ in the Definition \ref{inf_cov_oa_defn}  should not be unique.
\end{remark}

\subsection{Coverings of $O^*$-algebras}

\begin{definition}\label{comp_o*_defn}
	Let $\La$ be a  directed set (cf. Definition \ref{directed_set_defn}) such that there is the unique minimal element $\la_{\min} \in \La$.
	Let $A$ be an $O^*$-algebra (cf. Definition \ref{o*alg_defn}).  Let us consider a set noncommutative finite-fold coverings of  $O^*$-algebras $\mathfrak{S}=\left\{ \pi_{\la}:A \hookto A_{\la}\right\}_{\la \in \La}$ (cf. Definition \ref{fino*_defn}) indexed by  $\La$ such that $A_{\la_{\mathrm{min}}}= A$, and $\pi_{\la_{\mathrm{min}}}= \Id_A$.  Suppose that $\mathfrak{S}$ is a subcategory of $\mathfrak S$, such that:
	\begin{enumerate}
		\item [(a)] any object of $\mathfrak{S}$ corresponds to a class of isomorphic objects of $\mathfrak{S}$,
		\item[(b)]  for any $\nu, \mu\in \La$ such that $\nu >\mu$ there is the unique morphism from $A \hookto A_{\nu}$ to $A \hookto A_{\mu}$ which is a noncommutative finite-fold covering $\pi^\mu_\nu:  A_{\mu}\hookto A_\nu$ of  $O^*$-algebras  such that $\pi_\nu =\pi^\mu_\nu \circ \pi_\mu$,
		\item[(c)] all morphisms of   $\mathfrak{S}$ are described in (b).
	\end{enumerate}
We say that $\mathfrak{S}$ is an \textit{algebraic finite covering category of $O^*$-algebras}.   We write $\mathfrak{S} \in  O^*$-$ \mathfrak{FinAlg}$, and use the following notation
\be\label{comp_pt_o*_eqn}
\begin{split}
	\mathfrak{S} = \left(\left\{ A_\la \right\}_{\la \in \La}, \left\{\pi^\mu_\nu: A_\mu \hookto A_\nu\right\}_{\substack{\mu, \nu \in \La\\ \nu > \mu}}\right),\\
	\text{or simply} \quad \mathfrak{S} = \left(\left\{\pi_\la: A \hookto A_\la \right\}, \left\{\pi^\mu_\nu\right\}\right).
\end{split}
\ee
\end{definition}
	
\begin{remark}
	The category \ref{comp_pt_o*_eqn}
 is equivalent to the pre-order category given by $\La$ (cf. the Definition \ref{preordercat_defn} and the Remark \ref{pre_order_rem}).
\end{remark}

\begin{empt}\label{comp_pt_o*_empt}
Under the hypotheses of the Definition \ref{comp_o*_defn} morphisms of $\mathfrak{S}$ can be regarded as inclusions of *-algebras.
 Let  $\left(B, \overline{B}, \widehat G\right)$  be a  disconnected infinite noncommutative covering of 
 \be\label{bl_cat_eqn}
 \mathfrak{S}_B\bydef \left\{\left\{B_\la\right\}_{\la\in \La}, \left\{\widetilde \pi^\mu_\nu : B_\mu \hookto B_\nu\right\}_{\substack{\mu, \nu \in \La\\\mu \le \nu}}\right\}
 \ee 
 (cf.  Definition \ref{disconnected_infinite_noncommutative_covering_defn}) where  $\widehat G \bydef \varprojlim G\left(\left.A_\la\right|A\right)$ is an inverse limit of groups (cf. Definition \ref{group_inv_lim_defn}).	
  	The union $\widehat A \bydef\cup_{\la\in \La} A_\la$ is a *-algebra.  There is the natural action $\widehat G\times \widehat A\to \widehat A$.
  	Suppose that there is an equivariant  faithful nondegenerate representation (cf. Definitions \ref{equivariant_representation_defn}, \ref{faithful_representation_defn} and \ref{nondegenerate_repr_defn}) $\widehat \pi: \widehat B \bydef C^*$-$\lim_{\la\in\La} B_\la\hookto B \left(\widehat\H \right) $ (cf. Definition \ref{equivalent_representation_defn}) such that the given by the equation \eqref{ext_repro_eqn} representation
  	$$
  ext\left(\widehat \pi \right) : \overline B \to 	 B \left(\widehat\H \right)
  	$$
  	is faithful and nondegenerate (cf. Lemma \ref{infinite_faithful_nondegererate_lem}).  Suppose that there is a $\widehat G$-invariant dense subspace $\widehat \D\subset \widehat\H$ such that there is an  injective $*$-homomorphism   
	\be\label{pi_o*_eqn}
	\pi: \widehat A\hookto \L^\dagger\left(\widehat\D \right)\quad \text{where } \L^\dagger\left(\widehat\D \right)  \text{ is given by \eqref{l_dag_eqn}}
	\ee
such that
	\be\label{o*_dinv_eqn}
\begin{split}
	\widehat G\widehat \D= \widehat \D,\\
\forall \widehat \xi\in\widehat \D\quad \forall \widehat a\in\widehat A\quad	 \forall g\in\widehat G\quad g\left(\widehat a\widehat\xi \right)= \left(g\widehat a\right) \left(g\widehat\xi\right),\\
\forall \widehat \xi\in\widehat \D\quad \forall \widehat b\in\widehat B\cap \L^\dagger\left(\widehat\D \right)\quad \widehat{\pi}\left( \widehat b\right) \widehat \xi = {\pi}\left( \widehat b\right) \widehat \xi
\end{split}
\ee
For any $\la\in \La$ we denote by $A^b_\la$ the $C^*$-norm completion of $A_\la \cap \L^\dagger\left(\widehat\D \right)_b$ (cf. equation \eqref{o*b_eqn}), and assume that $B_\la$ is an essential ideal of $A^b_\la$. Suppose then both $*$-homomorphisms $\pi^\mu_\nu: A_\mu \hookto A_\nu$ ans $ \widetilde\pi^\mu_\nu: B_\mu \hookto B_\nu$  coincide on the intersection $A_\la \cap \L^\dagger\left(\widehat\D \right)_b$,
	Assume that $	\mathfrak{S}_B$ is good (cf. Definition \ref{good_defn}) and  $\left(B,  \widetilde{B}, G\right)$  is an  {infinite noncommutative covering} of  $\mathfrak{S}$. (cf. Definition \ref{infinite_noncommutative_covering_defn}). If $\widetilde{\H}$ is the Hilbert norm completion of $\widetilde{B}\widehat{\H}$  and $\widehat\H= \widetilde\H \oplus \H^\perp$  is the orthogonal sum then  $\widetilde{B}\H^\perp= \{0\}$.  
If 
$\widetilde \D \bydef \widehat{\D}\cap \widetilde\H$ then 
there is the natural inclusion $\L^\dagger\left( \widetilde \D \right)\subset \L^\dagger\left( \widehat \D \right)$. Also there are actions $G\times \widetilde\H\to \widetilde\H$ and $G\times \widetilde\D\to \widetilde\D$
which come from the action $\widehat G\times \widehat \H\to \widehat \H$.
A $G$-invariant  *-subalgebra $\widetilde A'\subset \L^\dagger\left( \widetilde \D \right)$ is said to be \textit{admissible} if it satisfies to the following conditions: 
	\bea\label{inf_cov_o*_sub1_eqn}
A_\la \widetilde A'\subset \widetilde A',\\
\label{inf_cov_o*_sub2_eqn}		\widetilde A'A_\la\subset \widetilde A',\\
\label{inf_cov_o*_sub3_eqn}		\forall \widetilde b \in K\left(\widetilde B\right) \cap \L^\dagger\left( \widetilde \D \right)_b 
~~ \forall \widetilde a \in \widetilde A'~~ \exists  \widetilde a', \widetilde a'' \in \bigcup_{\la\in\La} A_\la~~ \widetilde b \widetilde a = \widetilde b \widetilde a', ~ \widetilde a \widetilde b = \widetilde a'' \widetilde b.
\eea
where $K\left(\widetilde B \right)$ is the Pedersen's ideal of $\widetilde B$ (cf. Definition \ref{pedersen_ideal_defn}) and $\L^\dagger\left( \widetilde \D \right)_b$ is given by the equation  \eqref{o*b_eqn}.

\end{empt}

\begin{definition}\label{inv_o*_lim_defn}    
Consider the described in \ref{comp_pt_o*_empt} situation and suppose that there exists admissible *-algebra $\widetilde A\subset \L^\dagger\left( \widetilde \D \right)$ such that any -admissible *-algebra is a subalgebra of  $\widetilde A$. 
We say that the triple
	$\left(A, \widetilde{A}, G\right)$ is the 
\textit{ inverse noncommutative limit of $\mathfrak{S}$}. We also say that 	$\left(A, \widetilde{A}, G\right)$ is the 
	\textit{ inverse noncommutative limit of $O^*$-algebras}.
\end{definition}
\begin{remark}
	From the equations \eqref{inf_cov_o*_sub1_eqn}-\eqref{inf_cov_o*_sub3_eqn}  it turns out that the  union  is $\bigcup_{\la\in \La} A_\la$ is an admissible so there is the natural inclusion
	\be\label{adm_inc_eqn}
	\bigcup_{\la\in \La} A_\la\subset \widetilde A.
	\ee 
\end{remark} 
\subsection{Coverings and unbounded operators on Hilbert modules}\label{unb_hilb_sec}

\paragraph*{}
This section is an infinite counterpart of the Section \ref{fin_chull_sec}.
\begin{empt}\label{comp_chull_empt}
	Let $A$ be a *-algebra and let $B$ be a $C^*$-algebra which is a Hilbert $B$-module. Let $\mathfrak B\subset B$ be a dense right ideal, and let $A \hookto \End^*_{ { B}}\left( {\mathfrak B}\right)$ be a {representation} of $A $ on $B$ (cf. Definition \ref{def:rep_Hilbert_module_uni}). Suppose that  $\left(B, \widetilde{B}, G\left(\left.\widetilde{B}~\right| B\right)\right)$ is 
	{infinite noncommutative covering} of 
	\be\label{hm_cs_cat_eqn}	
	\mathfrak{S}_B=\left\{\left\{B_\la\right\}_{\la\in \La}, \left\{\widetilde \pi^\mu_\nu : B_\mu \hookto B_\nu\right\}_{\substack{\mu, \nu \in \La\\\mu \le \nu}}\right\}
	\ee	
	(cf. Definition \ref{infinite_noncommutative_covering_defn}).
	Let 	
	\be\label{comp_pt_rep_a_eqn}
	\begin{split}
		\mathfrak{S}_A \bydef \left(\left\{\pi_\la: A \hookto A_\la \right\}_{\la \in \La}, \left\{\pi^\mu_\nu: A_\mu \hookto A_\nu\right\}_{\substack{\mu, \nu \in \La\\ \nu > \mu}}\right).
	\end{split}
	\ee
be a category of $*$-algebras such that for any $\la, \mu, \nu\in \La$ the triples $$\left(A, A_\la, G\left(\left. A_\la\right| A\right)\bydef G\left(\left. B_\la\right| B\right), \pi_\la \right)$$ and $\left(A_\mu, A_\nu, G\left(\left. A_\mu\right| A_\nu\right)\bydef G\left(\left. B_\mu\right| B_\nu\right), \pi^\mu_\nu \right)$ are associated with $\left(B, B_\la, G\left(\left. B_\la\right| B\right), \rho_\la \right)$ and $\left(B_\mu, B_\nu,G\left(\left. B_\mu\right| B_\nu\right), \rho^\mu_\nu \right)$ noncommutative finite-fold covering of *-algebras (cf. Definition \ref{fin_chull_defn}.  Suppose that there is a dense $G\left(\left. \widetilde B\right| B\right)$-equivariant right  ideal $\widetilde{\mathfrak B}'\subset \widetilde B$. If $\la\in\La$ and $\mathfrak B_\la\in B_\la$ is a dense right ideal which corresponds to a covering $\left(A, A_\la, G\left(\left. A_\la\right| A\right), \pi_\la \right)$ (cf. Definition \ref{fin_chull_defn}) then we assume that there is an action   $\widetilde{\mathfrak B}'\times \mathfrak B_\la \to \widetilde{\mathfrak B}'$ which comes from the product $\widetilde{ B}\times B_\la\to\widetilde{ B}$. Suppose that for all $\la\in \La$ there is a  right action $ \widetilde{\mathfrak B}'\times A_\la \to \widetilde{\mathfrak B}'$ such that
	\be\label{hilb_act_eqn}
	\forall a \in A_\la \quad \forall b \in \mathfrak B_\la \quad \forall \widetilde b \in \widetilde {\mathfrak B}' \quad  \widetilde b\left(ba \right) =\left(  \widetilde bb\right) a. 
	\ee
	A $G\left(\left.\widetilde B\right| B\right)$-equivariant *-subalgebra $\widetilde A' \subset \End^*_{\widetilde { B}}\left(\widetilde {\mathfrak B}'\right)$ is said to be \textit{admissible} if for all $\widetilde a'\in\widetilde A'$ following conditions hold:
	\bea		\label{hilb_adm_eqn}
	\forall \widetilde a'\in\widetilde A'\quad  \forall \widetilde b \in K\left(\widetilde B \right)\cap \widetilde{\mathfrak B}'\quad  \exists a^\cup \in \bigcup_{\la\in\La} A_\la \quad \widetilde b \widetilde a' =   \widetilde ba^\cup.
	\eea
\end{empt}

\begin{definition}\label{main_ch_defn}
	Consider the described in \ref{comp_chull_empt} situation,
	suppose that there exists an admissible *-algebra $\widetilde A\subset  \End^*_{\widetilde { B}}\left(\widetilde {\mathfrak B}'\right)$ such that any admissible *-algebra is a subalgebra of  $\widetilde A$. We  say that the triple $\left(A, \widetilde{A}, G\left(\left.\widetilde A \right| A\right)\right)$	is an \textit{infinite noncommutative covering of 	$\mathfrak{S}_A$ associated with $\mathfrak{S}_B$}. 
\end{definition}

\begin{definition}\label{b_core_defn}
	The completion  $\widetilde {\mathfrak B}$ of  $\widetilde {\mathfrak B}'$ with respect to graph topology (cf. Definition \ref{def:rep_Hilbert_module}) given by the family of seminorms
	\be
	\forall \widetilde a \in \widetilde A\quad \widetilde \xi\in \widetilde {\mathfrak B}'\quad \left\| \widetilde \xi \right\|_{\widetilde a} \bydef\left\langle {\xi, {\pi(1_{\widetilde A^\sim}+\widetilde a^*\widetilde a)\xi}}\right\rangle ^{\nicefrac12}_{{\widetilde B}} 
	\ee
	is said to be the \textit{core of} $\left(A, \widetilde{A}, G\left(\left.\widetilde A \right| A\right)\right) $.
\end{definition}

\section{Coverings of quasi *-algebras}

\begin{definition}\label{comp_qo*_defn}
	Let $\La$ be a  directed set (cf. Definition \ref{directed_set_defn}) such that there is the unique minimal element $\la_{\min} \in \La$.
	Let $\left( \mathfrak A, \mathfrak A_0 \right) $ be a quasi-$O^*$-algebra  (cf. Definitions \ref{quasi_defn}, \ref{o*alg_defn}).  Let us consider a family noncommutative finite-fold coverings of  quasi-$O^*$-algebras $\mathfrak{S}=\left\{ \pi_{\la}:\left( \mathfrak A, \mathfrak A_0 \right) \hookto \left( \mathfrak A^\la, \mathfrak A_0^\la \right)\right\}_{\la \in \La}$ (cf. Definition \ref{oq*fin_defn}) indexed by  $\La$ such that $\left( \mathfrak A^{{\la_{\mathrm{min}}}}, \mathfrak A_0^{{\la_{\mathrm{min}}}} \right)= \left( \mathfrak A, \mathfrak A_0 \right)$, and $\pi_{\la_{\mathrm{min}}}= \Id_{\left( \mathfrak A, \mathfrak A_0 \right)}$.  Suppose that there is a subcategory $\mathfrak{S}$ of the category $\mathfrak{S}$ which is equivalent to the pre-ordering category $\La$ (cf. Definition \ref{preordercat_defn}).
	We say that $\mathfrak{S}$ is an \textit{algebraic finite covering category of quasi-$*$-algebras}.   We write $\mathfrak{S} \in  QO^*$-$ \mathfrak{FinAlg}$, and use the following notation
	\be\label{comp_pt_qo*_eqn}
	\begin{split}
		\mathfrak{S} =\\= \left(\left\{\pi_\la: \left( \mathfrak A, \mathfrak A_0 \right) \hookto \left( \mathfrak A^\la, \mathfrak A_0^\la \right) \right\}_{\la \in \La}, \left\{\pi^\mu_\nu: \left( \mathfrak A^\nu, \mathfrak A^\nu_0 \right) \hookto \left( \mathfrak A^\mu, \mathfrak A_0^\mu \right)\right\}_{\substack{\mu, \nu \in \La\\ \nu > \mu}}\right),\\
		\text{or simply} \quad \mathfrak{S} =\left(\left\{\pi_\la: \left( \mathfrak A, \mathfrak A_0 \right) \hookto \left( \mathfrak A^\la, \mathfrak A_0^\la \right) \right\}_{\la \in \La}, \left\{\pi^\mu_\nu: \right\}_{\substack{\mu, \nu \in \La\\ \nu > \mu}}\right).
	\end{split}
	\ee
\end{definition}
\begin{empt}\label{comp_qo*_empt}
Consider a given by \eqref{comp_pt_qo*_eqn} algebraic finite covering category $\mathfrak{S}$ of quasi-$*$-algebras  (cf. Definition \ref{comp_qo*_defn}). If invective $*$-homomorphisms are regarded as inclusions then there is a union $\widehat{\mathfrak A}_0\bydef\bigcup_{\la \in \La} \mathfrak A_0^\la$ which is a *-algebra. Suppose that
the category
\be\label{comp_pt_qo*0_eqn}
\begin{split}
	\mathfrak{S}_0\bydef \left(\left\{\left.\pi_\la\right|_{\mathfrak A_0}: \mathfrak A_0 \hookto \mathfrak A_0^\la \right\}_{\la \in \La}, \left\{\left.\pi^\mu_\nu\right|_{\mathfrak A^\mu_0}: \mathfrak A_0^\mu \hookto \mathfrak A_0^\nu\right\}_{\substack{\mu, \nu \in \La\\ \nu > \mu}}\right)
\end{split}
\ee
satisfies to the conditions  described in \ref{comp_pt_o*_empt}. If 	$\left(\mathfrak A_0, \widetilde{\mathfrak A}_0, G\right)$ is a good 
{ inverse noncommutative limit of $\mathfrak{S}_0$} then from the Definition \ref{inv_o*_lim_defn} it follows that following conditions hold.
\begin{enumerate}
	\item[(i)] There is a good {algebraical  finite covering category} 
	\be\label{blu_cat_eqn}
	\mathfrak{S}_B=\left\{ \left( B ,  B_{\la} , G\left(\left.\mathfrak A_0^\la~\right|\mathfrak A_0 \right) , \overline \pi_{\la}|_{\mathfrak A_0}\right) \right\}_{\la \in \La}
	\ee
	(cf. Definitions \ref{algebraical_finite_covering_category_defn}  and \ref{good_defn}), and there is an  {infinite noncommutative covering} of  $\mathfrak{S}_B$. (cf. Definition \ref{infinite_noncommutative_covering_defn}). \be\label{be_eqn} 
	\left(B, \widetilde{B}, G_{\widehat\pi}\right)
	\ee 
\item[(ii)] There is an equivariant faithful representation $\widehat \pi: \widehat B \bydef C^*$-$\lim_{\la\in\La} B_\la\hookto B \left(\widehat\H \right) $ (cf. Definition \ref{equivalent_representation_defn}) such that the given by the equation \eqref{ext_repro_eqn} representation
$$
ext\left(\widehat \pi \right) : \overline B \to 	 B \left(\widehat\H \right)
$$
is faithful.  
\item[(iii)]There is a $\widehat G$-invariant dense subspace $\widehat \D\subset \widehat\H$ such that there is an  injective $*$-homomorphism   
	\be\label{pi_o*f_eqn}
\pi: \widehat{\mathfrak A}_0\hookto \L^\dagger\left(\widehat\D \right) \text{where } \L^\dagger\left(\widehat\D \right)  \text{ is given by \eqref{l_dag_eqn}}.
\ee
\end{enumerate}
\end{empt}
	
\begin{defn}\label{comp_pt_qo*_defn}
Consider the situation described in \ref{comp_qo*_empt}.  The union $\bigcup_{\la\in \La} \left( \mathfrak A^\la, \mathfrak A_0^\la \right)$ is a quasi $*$-algebra. 
If $\left( \widetilde{\mathfrak A}', \widetilde{\mathfrak A}_0\right)$ is a quasi *-algebra  such that
\begin{enumerate}
	\item[(a)]  There is an  action $G\times \left( \widetilde{\mathfrak A}', \widetilde{\mathfrak A}_0\right)\to \left( \widetilde{\mathfrak A}', \widetilde{\mathfrak A}_0\right)$ which complies with the natural action $ G\times  \widetilde{\mathfrak A}_0\to  \widetilde{\mathfrak A}_0$.
	\item[(b)] There is an injective $*$-homomorphism   $\pi:\bigcup_{\la\in \La} \left( \mathfrak A^\la, \mathfrak A^\la_0\right)\hookto \left( \widetilde{\mathfrak A}', \widetilde{\mathfrak A}_0\right)$ such that the restriction $
	\pi|_{\bigcup_{\la\in \La}  \mathfrak A^\la_0}$ is the natural inclusion $ \bigcup_{\la\in \La}  \mathfrak A^\la_0\hookto \widetilde{\mathfrak A}_0$.
	
	\item[(c)] Following conditions hold:
\bea
\label{inv_qo*_rzero_eqn}
\forall\widetilde a\in \widetilde{\mathfrak A}'\quad \left( K\left(\widetilde{B} \right)\cap  \widetilde{\mathfrak A}_0 \right) \widetilde{a}= \{0\}\quad \Rightarrow\quad \widetilde{a}= 0, \\
\label{inv_qo*_lzero_eqn}
\forall\widetilde a\in \widetilde{\mathfrak A}'\quad \widetilde{a}\left( K\left(\widetilde{B} \right)\cap  \widetilde{\mathfrak A}_0 \right) = \{0\}\quad \Rightarrow\quad \widetilde{a}= 0, \\
\label{inv_qo*_rlim_eqn}
\forall  \widetilde b \in \widetilde{\mathfrak A}_0\cap K\left(\widetilde{B} \right)~\forall \widetilde a \in \widetilde{\mathfrak A}' ~\exists a',a'' \in    \bigcup_{\la\in \La} \mathfrak A_\la ~~ \widetilde b\widetilde a=\widetilde b a', ~~  \widetilde a\widetilde b=a''\widetilde b.
\eea
where $K\left(\widetilde B \right)$ is the Pedersen's ideal of $\widetilde B$ (cf. Definition \ref{pedersen_ideal_defn}) .
	
\end{enumerate}
then $\left( \widetilde{\mathfrak A}', \widetilde{\mathfrak A}_0\right)$ is said to be \textit{admissible}.
\end{defn}

\begin{definition}\label{inv_qo*_lim_defn}
Consider the situation \ref{comp_qo*_empt}.	Suppose that there exists an admissible quasi *-algebra $\left( \widetilde{\mathfrak A}, \widetilde{\mathfrak A}_0\right)$ such that any admissible  quasi *-algebra  $\left( \widetilde{\mathfrak A}', \widetilde{\mathfrak A}_0\right)$ is a subalgebra of  $\left( \widetilde{\mathfrak A}, \widetilde{\mathfrak A}_0\right)$. We say that the triple
	$\left(\left( {\mathfrak A}, {\mathfrak A}_0 \right), \left( \widetilde{\mathfrak A}, \widetilde{\mathfrak A}_0 \right), G_{\widehat{\pi}}\right)$ is the 
	${\pi}$-\textit{ inverse noncommutative limit of $\mathfrak{S}$}.
\end{definition}

\section{Coverings of (local) operator spaces}

\paragraph*{}
Here we consider a generalization of the discussed in the Section \ref{infinite_ca_sec} construction.

\begin{empt}\label{comp_op_pt_empt}
	Consider an algebraical  finite covering category (cf. Definition \ref{algebraical_finite_covering_category_defn}) given by
	$$
	\mathfrak{S} = \left(\left\{\pi_\la: A \hookto A_\la \right\}_{\la \in \La}, \left\{\pi^\mu_\nu: A_\mu \hookto A_\nu\right\}_{\substack{\mu, \nu \in \La\\ \nu > \mu}}\right),
	$$
	and suppose that $\mathfrak{S}$ is good (cf. Definition \ref{good_defn}). Assume that for every $\la \in \La$  there is a sub-unital operator space $\left(X_\la, Y_\la \right)$ (cf. Definition \ref{operator_space_subunital_defn}) such that $A_\la \cong C^*_e\left(X_\la, Y_\la \right)$ (cf. equation \ref{op_su_env_eqn}) suppose that any $\pi_\la$ and $\pi^\mu_\nu$ corresponds to a noncommutative finite-fold coverings of operator spaces (cf. Definition \ref{fin_op_defn})
\be\label{comp_op_ob_pt_eqn}
\left(\left(X, Y \right) , \left(\left(X_\la, Y_\la \right) \right), G\left(\left.A_\la\right|A \right), \left(\pi^\la_{X_\la}, \pi^\la_{Y_\la}\right)\right)
\ee and   
	\be\label{comp_op_mo_pt_eqn}
	\left(\left(X_\mu, Y_\mu \right) , \left(\left(X_\nu, Y_\nu \right) \right), G\left(\left.A_\mu\right|A_\nu \right), \left(\pi^{\mu\nu}_{X_\mu}, \pi^{\mu\nu}_{Y_\mu}\right)\right)
	\ee
	respectively.
\end{empt}
\begin{definition}\label{comp_op_pt_defn}
	Let 	$\mathfrak{S}_{\text{op}}$ be a category such that its objects and morphisms are given by \eqref{comp_op_ob_pt_eqn} and \eqref{comp_op_mo_pt_eqn} respectively. We say that $\mathfrak{S}_{\text{op}}$ 
	is an \textit{algebraical  finite covering category of operator spaces}. We write
	\be\label{comp_op_pt_eqn}
	\mathfrak{S}_{\text{op}}\in 	\mathfrak{OSp}_{\text{op}}.
	\ee
\end{definition}
\begin{definition}\label{spec_op_defn}
	Consider the situation of \ref{comp_op_pt_empt} and \ref{comp_op_pt_defn}. 	Let $\left(A, \widetilde{A}, G\left(\left.  \widetilde{A}~\right| A\right)\right)$ be a an infinite noncommutative covering of $\mathfrak{S}$ (cf. Definition \ref{infinite_noncommutative_covering_defn}).	  An element $\widetilde x  \in  \widetilde{A}$ is said to be \textit{subordinated} to $\mathfrak{S}_{\text{op}}$ if there is a net $\left\{x_\la \in X_\la\right\}_{\la \in \La}\subset \widehat{A}$ such that
	\be\label{spec_op_eqn}
	\widetilde x = \bt \text{-}\lim_{\la \in \La } x_\la
	\ee
	where the $\bt \text{-}\lim$ implies the  limit   with respect to the strict topology of $M\left(\widetilde A \right)$  (cf. Definition \ref{strong_topology_defn}).
\end{definition}
\begin{definition}\label{spec_lim_defn}
	Let $\Xi_{\text{op}} \subset \widetilde{A}$ is the space of {subordinated} to $\mathfrak{S}_{\text{op}}$ elements. The $C^*$-norm completion of the $\C$-linear space $\Xi_{\text{op}}$ is said to be the \textit{inverse noncommutative limit} of 	$\mathfrak{S}_{\text{op}}$.
\end{definition}
\begin{remark}
	Since the {inverse noncommutative limit} is a $C^*$-norm closed subspace of $C^*$-algebra $\widetilde A$ it has the natural structure of the operator space.
\end{remark}
\begin{definition}\label{comp_rop_pt_defn}
Similarly to \ref{comp_op_pt_empt} consider a category of noncommutative finite-fold coverings of real sub-unital operator spaces (cf. Definition \ref{fin_rop_defn})
\bean
\left(\left(X, Y \right) , \left(\left(X_\la, Y_\la \right) \right), G_\la, \left(\pi^\la_{X_\la}, \pi^\la_{Y_\la}\right)\right), \\
\left(\left(X_\mu, Y_\mu \right) , \left(\left(X_\nu, Y_\nu \right) \right), G^\mu_\nu, \left(\pi^{\mu\nu}_{X_\mu}, \pi^{\mu\nu}_{Y_\mu}\right)\right)
\eean
real operator spaces. Denote it by $	\mathfrak{S}_{\text{op}}$.  We say that 
is an \textit{algebraical  finite covering category of real operator spaces} if the  given by 
\bean
\left(\left(\C X, \C Y \right) , \left(\left(\C X_\la, \C Y_\la \right) \right), G_\la, \left(\C \pi^\la_{X_\la}, \C \pi^\la_{Y_\la}\right)\right),\\ 
\left(\left(\C X_\mu, \C Y_\mu \right) , \left(\left(\C X_\nu, \C Y_\nu \right) \right), G^\mu_\nu, \left(\C \pi^{\mu\nu}_{X_\mu}, \C \pi^{\mu\nu}_{Y_\mu}\right)\right)
\eean
complexification $\C	\mathfrak{S}_{\text{op}}$ of $	\mathfrak{S}_{\text{op}}$ is algebraical  finite covering category of  real operator spaces (cf Definition \ref{comp_op_pt_defn}).  Let $\widetilde{A}$ be the inverse noncommutative limit of 	$\mathfrak{S}$ if there is a net $\left\{x_\la \in X_\la\right\}_{\la \in \La}\subset \widehat{A}$ such that
\be\label{spec_rop_eqn}
\widetilde x =\bt\text{-} \lim_{\la \in \La } x_\la
\ee
	where the $\bt \text{-}\lim$ implies the  limit   with respect to the strict topology of $M\left(\widetilde A \right)$  (cf. Definition \ref{strong_topology_defn}).
\end{definition}
\begin{definition}\label{spec_r_lim_defn}
Let $	\mathfrak{S}_{\text{op}}$.  be a {algebraical  finite covering category of real operator spaces}. Similarly to the Definition an element $\widetilde x  \in  \widetilde{A}$ is said to be \textit{subordinated} to $\mathfrak{S}_{\text{op}}$. The $C^*$-norm completion of the $\C$-linear space $\Xi_{\text{op}}$ is said to be the \textit{inverse noncommutative limit} of 	$\mathfrak{S}_{\text{op}}$.
\end{definition}

\begin{lemma}\label{spec_r_lim_lem}
Let $	\mathfrak{S}_{\mathrm{op}}$.  be an {algebraical  finite covering category of real operator spaces} and let $	\C\mathfrak{S}_{\mathrm{op}}$ be its complexification. If both $\widetilde X$ and $\widetilde X'$ are {inverse noncommutative limits} of both	$\mathfrak{S}_{\mathrm{op}}$ and $\C\mathfrak{S}_{\mathrm{op}}$ then $\widetilde X'$ there is the complex conjugation $\widetilde X'$ such that $\widetilde X'$ is the complexification of $\widetilde X$, i.e. $\widetilde X'\cong\C\widetilde X$.
\end{lemma}
\begin{proof}
Form \ref{complexification_empt} it follows that for any $\la\in \La$ there is the conjugation $x \mapsto \overline x$ of $X_\la$. The equation \eqref{spec_op_eqn} yields the natural conjugation $\widetilde x \mapsto \overline{\widetilde x}$ on $\widetilde X'$. Direct check shows that
\be\label{spec_r_lim_eqn}
\widetilde X= \left\{\left. \widetilde x\in\widetilde X'\right| \widetilde x = \overline{\widetilde x}\right\},
\ee
it turns out that  $\widetilde X'\cong\C\widetilde X$.
\end{proof}

\begin{empt}\label{comp_op_loc_pt_empt}
	Consider the situation \ref{inv_pro_lim_empt}. Assume that for every $\la \in \La$  there is a sub-unital local operator space $\left(X_\la, Y_\la \right)$ (cf. Definition \ref{op_loc_su_space_defn}) such that $A_\la \cong C^*_e\left(X_\la, Y_\la \right)$ suppose that any $\pi_\la$ and $\pi^\mu_\nu$ corresponds to a noncommutative finite-fold coverings of local operator spaces \be\label{comp_lop_ob_pt_eqn}
	\left(\left(X, Y \right) , \left(\left(X_\la, Y_\la \right) \right), G\left(\left.A_\la\right|A \right), \left(\pi^\la_{X_\la}, \pi^\la_{Y_\la}\right)\right)
	\ee and   
	\be\label{comp_lop_mo_pt_eqn}
	\left(\left(X_\mu, Y_\mu \right) , \left(\left(X_\nu, Y_\nu \right) \right), G\left(\left.A_\mu\right|A_\nu \right), \left(\pi^{\mu\nu}_{X_\mu}, \pi^{\mu\nu}_{Y_\mu}\right)\right)
	\ee
	
		respectively.
\end{empt}
\begin{definition}\label{comp_loc_op_pt_defn}
		Let 	$\mathfrak{S}_{\text{lop}}$ be a category such that its objects and morphisms are given by \eqref{comp_lop_ob_pt_eqn} and \eqref{comp_lop_mo_pt_eqn} respectively.
	Consider the situation \ref{comp_op_loc_pt_empt}. We say that 
	is an \textit{algebraical  finite covering category of local operator spaces}. We write
	\be\label{comp_loc_op_pt_eqn}
	\mathfrak{S}_{\text{lop}}\in 	\mathfrak{LOSp}_{\text{lop}}.
	\ee
\end{definition}
\begin{definition}\label{spec_loc_op_defn}
	Consider the situation of the Definition \ref{comp_loc_op_pt_defn}. 	Let $\widetilde{A}$ be the inverse noncommutative limit of 	$\mathfrak{S}$. Let $\widehat{A}= C^*$-$\varinjlim_{\la \in \La} A_\la$ and $\pi_a:b\left( \widehat{A}\right)  \to B\left(\H_a \right)$ be the atomic representation. An element $\widetilde x  \in  b\left( \widetilde{A}\right) $ is said to be \textit{subordinated} to $\mathfrak{S}_{\text{op}}$ if there is a net $\left\{x_\la \in b\left( X_\la\right) \right\}_{\la \in \La}\subset \widehat{A}$ such that
	\be\label{spec_loc_op_eqn}
	\widetilde x = \lim_{\la \in \La } x_\la
	\ee
	where the $\bt \text{-}\lim$ implies the  limit   with respect to the strict topology of $M\left(\widetilde A \right)$  (cf. Definition \ref{strong_topology_defn}).
\end{definition}
\begin{definition}\label{spec_loc_lim_defn}
	Consider the situation of the Definition \ref{spec_loc_op_defn}.
	Let $\Xi_{\text{op}} \subset \widetilde{A}$ is the space of {subordinated} to $\mathfrak{S}_{\text{op}}$ elements. Let $\widetilde X_{b}$ ba $C^*$-norm completion of the $\C$-linear space $\Xi_{\text{op}}$.	The maximal local operator subspace $\widetilde{X}\subset \widetilde{A}$ such that $b\left(\widetilde X \right) = \widetilde X_{b}$  said to be the \textit{inverse noncommutative limit} of 	$\mathfrak{S}_{\text{op}}$.
\end{definition}
\begin{remark} 
	Similarly to the Definition \ref{spec_loc_lim_defn}  one can define inverse noncommutative limits for local real operator spaces.
\end{remark}

\section{Universal  coverings and  fundamental groups}\label{inf_fund_g_sec}

\paragraph{}
Very natural choice of fundamental group is proposed in \cite{milne:etale}, where the fundamental group of algebraic manifold is an inverse limit of finite covering groups  (cf. Definition \ref{group_inv_lim_defn}). However this theory does not yield the fundamental group, it provides  the profinite completion of it only. Another  construction of the noncommutative fundamental group is an application of the described in this section theory. A discussed here  construction yields a noncommutative  fundamental group which  in case of commutative $C^*$-algebras is finer that the profinite completion (cf. Theorem \ref{comm_uni_lim_thm}).

\subsection{Universal coverings of $C^*$-algebras}\label{uni_sec}
\paragraph{} Here a noncommutative generalization of the universal covering space (cf. Definition \ref{top_universal_covering_defn}) is discussed.
\begin{definition}\label{fundamental_group_nc_defn} 
	Let $A$ be a connected $C^*$-algebra, and let 
	$
	\left(A, \widetilde{A}, G\left(\left.\widetilde{A}~\right| A\right)\right)
	$
	be the  {infinite noncommutative covering} (cf. Definition \ref{infinite_noncommutative_covering_defn}) of $$\mathfrak{S}=\left\{\left\{A_\la\right\}_{\la\in \La}, \left\{\pi^\mu_\nu : A_\mu \hookto A_\nu\right\}_{\substack{\mu, \nu \in \La\\\mu \le \nu}}\right\}$$ such that $A = A_{\la_{\mathrm{min}}}$.
	Suppose that $\mathfrak{S}$ contains \textit{all} classes of isomorphisms of  noncommutative finite-fold coverings of $A$ (cf. Definition \ref{fin_defn}). Then the triple $\left(A, \widetilde{A}, G\left(\left.\widetilde{A}~\right|{A}  \right) \right)$ of $\mathfrak{S}$ is said to be the \textit{universal  covering} 
	of $A$.  The group $G\left(\left.\widetilde{A}~\right|{A}  \right)$  is said to be the \textit{fundamental group} of $A$. We use the following notation
	\be\label{fg_bas_eqn}
	\pi_1\left(A \right) \stackrel{\mathrm{def}}{=} G\left(\left.\widetilde{A}~\right|{A} \right).
	\ee
\end{definition}

\begin{definition}\label{fundamental_group_nc_p_defn} 
	Let $P$ be a property of noncommutative finite-fold coverings.
	Let $A$ be a $C^*$-algebra, and let 
	$
	\left(A, \widetilde{A}, G\left(\left.\widetilde{A}~\right| A\right)\right)
	$
	be the  {infinite noncommutative covering} (cf. Definition \ref{infinite_noncommutative_covering_defn}) of $$\mathfrak{S}=\left\{\left\{A_\la\right\}_{\la\in \La}, \left\{\pi^\mu_\nu : A_\mu \hookto A_\nu\right\}_{\substack{\mu, \nu \in \La\\\mu \le \nu}}\right\}$$ such that $A = A_{\la_{\mathrm{min}}}$.
	Suppose that $\mathfrak{S}$ contains \textit{all} classes of isomorphisms of noncommutative finite-fold coverings of $A$ (cf. Definition \ref{fin_defn}) which possess the property $P$. Assume that for all $\mu, \nu \in \La$ such that $\mu \le \nu$ the finite-fold noncommutative cornering  $\pi^\mu_\nu : A_\mu \hookto A_\nu$ possesses the property $P$. 
	Then the triple $\left(A, \widetilde{A}, G\left(\left.\widetilde{A}~\right|{A}  \right) \right)$ of $\mathfrak{S}$ is said to be the $P$-\textit{universal  covering} 
	of $A$.  The group $G\left(\left.\widetilde{A}~\right|{A}  \right)$  is said to be the $P$-\textit{fundamental group} of $A$. We use the following notation
	\be\label{fg_bas_p_eqn}
	\pi^P_1\left(A \right) \stackrel{\mathrm{def}}{=} G\left(\left.\widetilde{A}~\right|{A} \right).
	\ee
\end{definition}
\begin{example}
There is a property $P_{\mathrm{untz}}$ of f noncommutative finite-fold coverings such that
\be\label{unitization_p_eqn}
\begin{split}
\left(A, \widetilde A, G, \pi \right) \in P_{\mathrm{untz}}\quad \Leftrightarrow \\ \Leftrightarrow \left(A, \widetilde A, G, \pi \right) \text{ is a covering with unitization (cf. Definition \ref{fin_unitization_defn})}.
\end{split}
\ee
\end{example}

\begin{exercise}\label{fundamental_group_mor_exer}
	Consider the situation of the Definition \ref{fundamental_group_nc_defn}. 
	Suppose that the $C^*$-algebra $A$ is {proper with respect to covering compositions} (cf. Definition \ref{proper_composition_defn}). Prove following statements.
	\begin{enumerate}
		\item There is a natural functor (cf. Definition \ref{functor_defn}) $\pi_1$ from $\mathfrak{S}$ to the category of groups.
		\item  For any morphism $\pi: \widetilde{A}' \hookto \widetilde{A}''$ of the category $\mathfrak{S}$ there is the natural isomorphism
		$$
		G\left( \left.\widetilde{A}'' \right|\widetilde{A}'\right) \cong \pi_1\left(\widetilde{A}' \right) / \pi_1\left(\pi \right)\left(\widetilde{A}'' \right).  
		$$
	\end{enumerate}
	
\end{exercise}
\begin{remark}
	The Exercise \ref{fundamental_group_mor_exer} can be regarded as a generalization of the Theorem \ref{top_fundamental_group_mor_thm}. (cf. also the Exercise \ref{top_fundamental_group_mor_exer}).
\end{remark}

\begin{definition}\label{fundamental_group_nc_stable_defn} Similarly to \eqref{fundamental_group_nc_defn}  consider an  {algebraical  finite covering}  category
\be\label{stable_fg_eqn}
\mathfrak{S}_\K = \left( \left\{ \rho_\la:A\otimes \K \hookto A^\K_\la\right\}_{\la\in\La},\left\{\rho^\nu_\mu\right\}_{\substack{ \mu, \nu \in \La \\ \mu \ge \nu}} \right).
\ee
If 	$\left(A\otimes\K, \widetilde{A}_\K, G\left(\left.\widetilde{A}_\K~\right|{A}\otimes\K  \right) \right)$ of $\mathfrak{S}$ is an {universal  covering}
of $A\otimes\K$.
	The group $G\left(\left.\widetilde{A}_\K~\right|{A}\otimes \K \right)$  is said to be the \textit{stable fundamental group}.
\end{definition}

\subsection{Universal coverings of operator spaces}
\paragraph{} A generalization of the Section \ref{uni_sec} is being discussed below. 
\begin{definition}\label{fg_bas_os_defn} 
	Let $\left(X, Y\right)$ be a {sub-unital} operator space (cf. Definition \ref{operator_space_subunital_defn}). Let us consider the family $\left(\left(X, Y \right) , \left(\left(X_\la, Y_\la \right) \right), G\left(\left.A_\la\right|A \right)\right)$ of \textit{all} noncommutative finite-fold coverings of $\left(X, Y\right)$.
	Denote by $\mathfrak{S} = \left\{ \pi_\la:A \hookto A_\la\right\}_{\la \in \La}$ and suppose that there is   an  {algebraical  finite covering category of operator spaces}  category $\mathfrak{S}$  (cf. Definition \ref{comp_op_pt_defn}) given by
		\bean
	\begin{split}
		\mathfrak{S}_{\text{op}}
		= \left(\left\{\left( \pi_{X_\la}, \pi_{Y\la}\right) : \left(X, Y \right) \to \left(X_\la, Y_\la \right)^{~} \right\}_{\substack{\la\in \La\\ \substack{~\\~}}},\right. \\
		\left.\left\{\left( \pi^{\nu}_{X_\mu}, \pi^{\nu}_{Y_\mu}\right) : \left(X_\mu, Y_\mu \right) \to \left(X_\nu, Y_\nu \right) \right\}_{\substack{\mu, \nu \in \La\\ \nu > \mu}}\right).
	\end{split}
	\eean
	The  {inverse noncommutative limit} of $\mathfrak{S}_{\text{op}}$ \ref{spec_lim_defn} is said to be  the \textit{universal  covering} of $\left(X, Y \right)$. If 	$C^*_e\left( X, Y\right)$ is the $C^*$-\textit{envelope} of $\left(X, Y \right)$ then the fundamental group  $\pi_1\left( C^*_e\left( X, Y\right), \left\{\pi^\nu_\mu\right\} \right) $ is said to be the \textit{fundamental group} of $\left(X, Y \right)$. We write
		\be\label{fg_bas_os_eqn}
	\pi_1\left(\left(X, Y \right) , \left\{\pi^\nu_\mu\right\} \right)\bydef	\pi_1\left(C^*_e\left(X, Y \right) , \left\{\pi^\nu_\mu\right\} \right)
	\ee
	
\end{definition}

\begin{remark}\label{fg_bas_os_rem}
If $\widetilde X$ is the universal covering of $\left(X, Y\right)$ and $\widetilde{C^*_e\left(X, Y \right)}$ is the universal covering of $C^*_e\left(X, Y \right)$ then there is the natural inclusion $\widetilde X \subset \widetilde{C^*_e\left(X, Y \right)}$. Otherwise if $\pi_1\left(C^*_e\left(X, Y \right) , \left\{\pi^\nu_\mu\right\} \right) \times C^*_e\left(X, Y \right)\to C^*_e\left(X, Y \right)$ is the natural action then $\pi_1\left(C^*_e\left(X, Y \right) , \left\{\pi^\nu_\mu\right\} \right) \widetilde X = \widetilde X$, hence from  \eqref{fg_bas_os_eqn} it follows the existence of the natural action
 \be\label{fg_bas_act_os_eqn}
 \pi_1\left(\left(X, Y \right) , \left\{\pi^\nu_\mu\right\} \right)\times 	\widetilde X \to \widetilde X.
 \ee
\end{remark}
\begin{remark}
	The both  notions of the universal covering and the fundamental group of a $C^*$-algebra are specializations of the universal covering and the fundamental group of an operator space.
\end{remark}
\begin{remark}\label{real_op_rem}
Similarly to the Definition \ref{fg_bas_os_defn} one can define the universal covering and the fundamental group of real operator space.
\end{remark}

\subsection{Fundamental group and noncommutative closed paths}
\paragraph*{}
Topological fundamental group are given by classes of closed paths (cf. Remark \ref{top_homotopy_group_rem}). Sometimes elements of noncommutative fundamental group correspond to described it the Section \ref{noncomutative_path_liftting} noncommutative closed paths.
\begin{exercise}
Suppose that there is   an  {algebraical  finite covering} category $\mathfrak{S} = \left( \left\{ \pi_\la:A \hookto A_\la\right\}_{\la\in\La},\left\{\pi^\nu_\mu\right\}_{\substack{ \mu, \nu \in \La \\ \mu \ge \nu}} \right)$ which is good. Let  $f: [0,1] \to  \mathrm{Aut}\left({A}\right)$ be such that  $f\left( 0\right)= f\left( 1\right)=\Id_A$. From \ref{lift_unique_empt} it turns out that for all $\la\in\La$ there is $g_\la \in G\left(\left. A_\la\right| A\right)$ which corresponds to $f$.
\begin{enumerate}
\item Prove that 
$
\forall \mu, \nu\in\La \quad \mu > \nu \quad\Rightarrow\quad  g_\nu = h^\mu_\nu\left(g_\mu \right) 
$ where the homomorphism  $ h^\mu_\nu: G\left(\left. A_\mu\right| A\right)\to G\left(\left. A_\nu\right| A\right)$.
\item Prove that $f$ corresponds to the unique $x_f \in 	\pi_1\left(A, \left\{\pi^\nu_\mu\right\}_{\substack{ \mu, \nu \in \La \\ \mu \ge \nu}} \right)$ (cf. equation \ref{fg_bas_eqn}).
\end{enumerate}
\end{exercise}
\section{Strong Morita equivalence of infinite coverings}\label{inf_mor_sec}

\paragraph{}  Let  $ \mathscr L^2\left( \widetilde{A}\right) = \mathscr L^2\left( \widetilde{A}\right)_A=_{\widetilde{A}}\mathscr L^2\left( \widetilde{A}\right)_A$ {$C^*$-Hilbert module   associated with the  infinite noncommutative covering}  $\left(A, \widetilde{A}, G\left(\left.\widetilde{A}~\right| A\right)\right)$ (cf. Definition \ref{infinite_hilb_mod_defn})
Any pair $\left(\widetilde{a}, g\right) \in \widetilde{A}\times G$ yields a $*$-endomorphism 
\be\label{hm_inc_eqn}
\begin{split}
	h_{\widetilde{a}, g}\in \End^*_A\left( \mathscr L^2\left( \widetilde{A}\right) \right) ,\\
	\xi \mapsto\widetilde{a} \left(g\xi  \right)\quad \xi\in  \mathscr L^2\left( \widetilde{A}\right). 
\end{split}
\ee
If $C_c\left(G, \widetilde{A}\right)$ is an explained in the Appendix \ref{discr_cr_prod_sec} *-algebra then the equation \eqref{hm_inc_eqn} yields an inclusion  $C_c\left(G, \widetilde{A}\right)\subset \End^*_A\left( \mathscr L^2\left( \widetilde{A}\right) \right)$. Henceforth we suppose that:
\begin{enumerate}
	\item [(a)] $C_c\left(G, \widetilde{A}\right) \subset \K\left( \mathscr L^2\left( \widetilde{A}\right)_A \right)$,
	\item[(b)] $C_c\left(G, \widetilde{A}\right)$ is dense *-subalgebra of $\K\left( \mathscr L^2\left( \widetilde{A}\right)_A \right)$.
\end{enumerate}
In this case  $\K\left( \mathscr L^2\left( \widetilde{A}\right)_A \right)$ is a completion of $C_c\left(G, \widetilde{A}\right) $ with respect to the  $C^*$-norm given by the inclusion (a). Similarly to the Appendix \ref{discr_cr_prod_sec} we use the following notation $\widetilde A\rtimes_cG\bydef \K\left( \mathscr L^2\left( \widetilde{A}\right)_A \right)$. Indeed  $\widetilde A\rtimes_cG$ is a completion of  $C_c\left(G, \widetilde{A}\right)$ with respect to a following $C^*$-norm
$$
\forall a \in C_c\left(G, \widetilde{A}\right)\quad \left\| a\right\|\bydef   \left\| \ka\left(a \right)\right\|
$$
where $\ka: C_c\left(G, \widetilde{A}\right)\hookto \K\left( \mathscr L^2\left( \widetilde{A}\right)_A \right)$ is the natural inclusion.
There is the natural left  action $\left( \widetilde A\rtimes_cG\right) \times  \mathscr L^2\left( \widetilde{A}\right)_A\to  \mathscr L^2\left( \widetilde{A}\right)_A$. 
\begin{definition}\label{cross_cov_defn}
	In the above situation we say that $\widetilde A\rtimes_cG$ is an \textit{crossed product associated with the covering} $\left(A, \widetilde{A}, G\right)$.
\end{definition}
In the above situation one can define analogs of equations \eqref{mor_pp_fin_eqn}.
\be\label{mor_pp_inf_eqn}
	\begin{split}
\varphi:\mathscr L^2\left( \widetilde{A}\right)\otimes_{A}  \mathscr L^2\left( \widetilde{A}\right)   \to \widetilde{A} \rtimes_c G,\quad		\varphi \left(\widetilde{a} \otimes \widetilde{b} \right)\left( g\right) =    \widetilde{a} \left(g\widetilde{b} \right); \\
\psi: \mathscr L^2\left( \widetilde{A}\right)  \otimes_{\widetilde{A} \rtimes_c G}  \mathscr L^2\left( \widetilde{A}\right)  \to A, \quad		\psi \left(\widetilde{a} \otimes \widetilde{b} \right) = \sum_{g \in G} g\left(\widetilde{a} \widetilde{b} \right), 
\end{split}
\ee
Unlike \eqref{mor_pp_fin_eqn} the group $G$ is not finite. However from the equation \eqref{inf_hilb_f_prod_eqn} and the dense inclusion  $C_c\left(G, \widetilde{A}\right) \subset \K\left(X_A \right)$ it turns out that both  $\varphi$ and $\psi$ are well defined.
In result one has the following pairings:
	\begin{equation}\label{hilb_inf_eqn}
	\begin{split}
\left\langle \cdot, \cdot \right\rangle_{\widetilde{A} \rtimes_c G}: \mathscr L^2\left( \widetilde{A}\right) \otimes_{A}  \mathscr L^2\left( \widetilde{A}\right)\to \widetilde{A} \rtimes_c G,\quad \left\langle \widetilde a, \widetilde b \right\rangle_{\widetilde{A} \rtimes_c G} \bydef \varphi\left(\widetilde a \otimes \widetilde b^* \right);\\ 
\left\langle \cdot, \cdot \right\rangle_{\widetilde{A}}: \mathscr L^2\left( \widetilde{A}\right) \otimes_{\widetilde{A} \rtimes_c G} \mathscr L^2\left( \widetilde{A}\right)\to A, \quad \left\langle\widetilde a,\widetilde b \right\rangle_A \bydef \psi\left( \widetilde a^* \otimes \widetilde b\right)
\end{split}
\end{equation}
Using isomorphism $\widetilde{A} \rtimes_c G \cong \K \left(\mathscr L^2\left( \widetilde{A}\right)_A \right)$ the pairing $\left\langle \cdot, \cdot \right\rangle_{\widetilde{A} \rtimes_c G}$ can be defined by the following  way 
 	\begin{equation}\label{hilb_infz_eqn}
 	\begin{split}
 		\left\langle \xi, \eta \right\rangle_{\widetilde{A} \rtimes_c G}\bydef \xi\left\rangle \right\langle \eta \bydef \zeta \mapsto \eta	\left\langle \xi, \zeta \right\rangle_{A}\in \K \left(\mathscr L^2\left( \widetilde{A}\right)_A \right)\cong  \widetilde{A} \rtimes_c G\quad \forall\zeta \in \mathscr L^2\left( \widetilde{A}\right)_A.
 	\end{split}
 \end{equation}
 
\begin{definition}\label{allows_morita_defn}
In the  situation of the Definition \ref{cross_cov_defn}
	we say that an infinite noncommutative covering $\left(A, \widetilde{A}, G\right)$ \textit{allows strong Morita equivalence} if the linear span of $\left\langle \mathscr L^2\left( \widetilde{A}\right), \mathscr L^2\left( \widetilde{A}\right) \right\rangle_{\widetilde{A} \rtimes_c G}$  is dense in $A$.
\end{definition}
\begin{exercise}\label{allows_morita_exer}
	Let  $\left(A, \widetilde{A}, G\right)$ be an infinite noncommutative covering which {allows strong Morita equivalence}.  Prove that $\mathscr L^2\left( \widetilde{A}\right)$ is an  {$A$-$ \widetilde A\rtimes_cG$-equivalence
		bimodule} (cf. Definition \ref{strong_morita_defn}). 
\end{exercise} 

 \section{Coverings of spectral triples}\label{str_cov_glo_sec}

\begin{defn}\label{spectral_triple_weakly_coh_defn}
	Let  $\left(\sA, \sH, D, J\right)$ be a real  spectral triple (cf. Section \ref{df:spt-real_defn}), and let $A$ be the $C^*$-norm completion of $\A$ with the natural representation $A \to B\left( \H\right)$. Let
	\begin{equation}\label{cov_sec_triple_eqn}
\mathfrak{S}\bydef \left\{\left\{A_\la\right\}_{\la\in \La}, \left\{\pi^\mu_\nu : A_\mu \hookto A_\nu\right\}_{\substack{\mu, \nu \in \La\\\mu \le \nu}}\right\}
	\end{equation}
	be a good  algebraic  finite covering category (cf.  Definitions \ref{algebraical_finite_covering_category_defn} and \ref{good_defn}). 
	Suppose that for any $\la \in \La$ there is a spectral triple $\left(\sA_\la, \sH_\la, D_\la, J_\la\right)$, such that:
	\begin{enumerate}
		\item[(a)] $\left(\sA_\la, \sH_\la, D_\la, J_\la\right)$ is the  $\left(A, A_\la, G\left(\left.A_\la~\right|A \right), \pi_\la  \right)$-lift of $\left(\sA, \sH, D, J\right)$ (cf. Definition \ref{spectral_triple_fin_lift_defn}),
		\item[(b)] $A_\la$ is the $C^*$-norm completion of $\A_\la$,
		\item[(c)] there is an infinite noncommutative covering $\left(A, \widetilde{A}, G\left(\left.\widetilde{A}~\right| A\right)\right)$ of 	$\mathfrak{S}$ (cf. Definition \ref{infinite_noncommutative_covering_defn}),
		\item[(d)] for any $\mu > \nu$  the spectral triple $\left(\sA_\mu, \sH_\mu, D_\mu, J_\mu\right)$ is a $\left(A_\nu, A_\mu, G\left(\left.A_\mu~\right|A_\nu \right)  \right)$-lift of $\left(\sA_\nu, \sH_\nu, D_\nu, J_\mu\right)$ (cf. Definition \ref{spectral_triple_fin_lift_defn}).
	\end{enumerate}
	We say that
	\begin{equation}\label{spectral_triple_sec_eqn}
		\begin{split}
			\mathfrak{S}_{\left(\sA, \sH, D, J\right)}= \left\{\left(\sA_\la, \sH_\la, D_\la, J_\la\right)\right\}_{\la \in \La}
		\end{split}
	\end{equation}
	is a \textit{coherent set of spectral triples}. We write $\mathfrak{S}_{\left(\sA, \sH, D, J\right)} \in \mathfrak{CohTriple}$.
\end{defn}

Let $\mathscr L^2\left( \widetilde{A}\right) $ be the {$C^*$-Hilbert   associated the  infinite noncommutative covering} $\left(A, \widetilde{A}, G\left(\left.\widetilde{A}~\right| A\right)\right)$ of $\mathfrak{S}$ (cf. Definition \ref{infinite_hilb_mod_defn}). Let $\rho: A \hookto B\left(\H\right)$ be the natural representation, and let $\widetilde\rho: \widetilde A \hookto B\left(\widetilde \H\right)$  be a representation {induced} by $\left( \rho, \left( A, \widetilde{A}, G\left(\left.\widetilde{A}~\right|A \right) \right)\right)$ (cf. Definition \ref{induced_repr_inf_defn}), where $\widetilde \H$ is a norm completion of a pre-Hilbert space $ \mathscr L^2\left( \widetilde{A}\right)\otimes_{A_{}} \H$. Using the given by \eqref{infinite_hilb_inv_eqn} involution, similarly to the equations \eqref{cc_eqn}, \eqref{ccc_eqn} one can define maps 
\be\label{cc_inf_defn}
\begin{split}
	\widetilde	J : \mathscr L^2\left( \widetilde{A}\right)\otimes_{A_{}} \H\to \mathscr L^2\left( \widetilde{A}\right)\otimes_{A_{}} \H,\\
	\widetilde J\left( \widetilde a \otimes \xi \right) \bydef \widetilde a^* \otimes J\xi;\\
	\widetilde	J^\dagger : \mathscr L^2\left( \widetilde{A}\right)\otimes_{A_{}} \H\to \mathscr L^2\left( \widetilde{A}\right)\otimes_{A_{}} \H,\\
	\widetilde J^\dagger \left( \widetilde a \otimes \xi \right) \bydef \widetilde a^* \otimes J^\dagger\xi,
\end{split}
\ee
which can be extended up to anti-isometries $\widetilde J: \widetilde \H\xrightarrow{\approx }\widetilde \H$ and $\widetilde J^\dagger: \widetilde \H\xrightarrow{\approx }\widetilde \H$ such that $\widetilde J^\dagger\widetilde J= \widetilde J\widetilde J^\dagger= \Id_{\widetilde \H}$. 
From \eqref{kla_dens_eqn} it turns out that for all $\la\in \La$ there is a dense inclusion 	
\bean
\mathscr L^2\left( \widetilde{A}\right)\otimes_{A_{\la}} \H_{\la}\hookto \widetilde \H. 
\eean
From the equations  \eqref{first_order_wt_eqn}, \eqref{first_order_wtb_eqn} it turns out that
\be\label{first_order_la_eqn}
\begin{split}
	\forall a_\la, b_\la \in \A_\la\quad  \left[b_\la, \left[ D_\la, J_\la a^*_\la J^\dagger_\la\right] \right]= 0 \quad \Rightarrow \\
	\Rightarrow b_\la \left[ D_\la, J_\la a^*_\la J^\dagger_\la\right]=  \left[ D_\la, J_\la a^*_\la J^\dagger_\la\right]b_\la.
\end{split}
\ee
If  $a_\la \in \A_\la$ then one can define right action of an operator $\left[D_\la, Ja^*_\la J^\dagger\right]$  on the space 	$\mathscr L^2\left( \widetilde{A}\right)\otimes_{A_{\la}} \H_{\la}$ such that
\be\label{fist_order_la_eqn}
\left(  \widetilde  a\ox\xi\right)\left[D_\la, J_\la a^*_\la J^\dagger_\la\right]= \widetilde  a\ox \left[D_\la, J_\la a^*_\la J^\dagger_\la\right] \xi
\ee
For any $b_\la \in \A_\la$ one has $\widetilde  ab_\la\ox\xi= \widetilde  a\ox b_\la\xi$, so one needs proof that 
$$
\left( \widetilde  ab_\la\ox\xi\right) \left[D_\la, J_\la a^*_\la J_\la^\dagger\right]= \left( \widetilde  a\ox b_\la\xi\right) \left[D_\la, J_\la a^*_\la J_\la^\dagger\right].
$$
However from the equation \eqref{first_order_la_eqn} it follows that
\bean
\left(  \widetilde  ab_\la \ox\xi\right)\left[D_\la, J_\la a^*_\la J^\dagger_\la\right]= \widetilde  ab_\la \ox \left[D_\la, J_\la a^*_\la J^\dagger_\la\right] \xi= \widetilde  a\ox b_\la  \left[D_\la, J_\la a^*_\la J^\dagger_\la\right]\xi =\\
=  \widetilde  a\ox   \left[D_\la, J_\la a^*_\la J^\dagger_\la\right]b_\la\xi= \left( \widetilde  a\ox b_\la\xi\right) \left[D_\la, J_\la a^*_\la J_\la^\dagger\right],
\eean 
i.e. an action is correctly defined. 
From $J^\dagger_\la\left[D_\la, J_\la a^*_\la J^\dagger_\la\right]^*J_\la= \left[D_\la, a_\la \right]$ and 
the  action \eqref{fist_order_la_eqn} one can construct an inclusion
\be\label{da_la_eqn}
\left[D_\la, \A_\la\right]\subset B\left( \widetilde\H\right).
\ee
If $\pi^s_\la: \A_\la \hookto B\left( \H^{2^s}_\la \right)$ is given by  \eqref{s_diff1_repr_equ} and \eqref{s_diff_repr_equ} then from \eqref{da_la_eqn} it turns out that there is the inclusion
\be\label{ps_inc_eqn}
\pi^s_\la\left( \A_\la\right) \hookto B\left(\widetilde\H^{2^s} \right). 
\ee

\begin{defn}\label{smooth_el_defn}
	Let us consider the  situation of the Definition \ref{spectral_triple_weakly_coh_defn}.
	Let $\Om^1_D$ be the {module of differential forms associated} with the spectral triple  $\left( \A, \H, D, J\right)$ (cf. Definition \ref{ass_cycle_defn}).
	An element $\widetilde{a}  \in \widetilde{A}$ is said to be $\mathfrak{S}_{\left(\sA, \sH, D, J\right)}$- \textit{smooth} if the following conditions hold:
	\begin{enumerate}
		\item[(a)] The element $\widetilde{a}$ lies in $\mathscr L^2\left( \widetilde{A}\right)$.
		\item[(b)] For any $\la \in \La$ the series 
		$$
		a_\la=	\sum_{g \in \ker\left(  G\left(\left.\widetilde{A}\right|A \right)\to G\left(\left.A_\la\right|A \right)\right) } ~ g \widetilde{a} 
		$$
		is convergent with respect to strict topology of $M\left(\widetilde{A}\right)$ and
		$a_\la \in \A_\la$ where the inclusion $\A_\la\subset A_\la\subset M\left(\widetilde{A}\right)$ is implied.
		\item[(c)]  If for any $\la \in \La$ and $s \in \N$ the representation $\pi^s_\la: \A_\la \hookto B\left( \H^{2^s}_\la \right)$ is given by  \eqref{s_diff1_repr_equ}, \eqref{s_diff_repr_equ} 
		and  
		inclusion $\pi^s_\la\left( \A_\la\right) \subset B\left(\widetilde\H^{2^s} \right)$ (cf.  \eqref{ps_inc_eqn})  is implied then  the net $\left\{\pi^s_\la\left(a_\la \right)\right\}_{\la \in\La}$  is convergent with respect to the strong topology of  $B\left(\widetilde{\H}^{2^s} \right)$ (cf. Definition \ref{strong_topology_defn}).	
		\item[(d)] If for any $\la\in \La$ the given by \eqref{da_la_eqn}  inclusion $
		\left[D_\la, \A_\la\right]\subset B\left( \widetilde\H\right)
		$ is implied then the net $\left\{\left[D_\la, a_\la\right]\right\}_{\la \in\La}$ is convergent with respect to the strong topology of  $B\left(\widetilde{\H}^{2^s} \right)$   (cf. Definition \ref{strong_topology_defn}).	.
	\end{enumerate}
	Denote by 
	\be\label{a_s_eqn}	
	\widetilde{a}^s\bydef \lim_{\la\in\La}  \pi^s_\la\left(a_\la \right) \in B\left(\widetilde\H^{2^s} \right)
	\ee	
	in sense the strong convergence of $B\left(\widetilde\H^{2^s} \right)$, and denote by $\widetilde{W}^\infty$ the space of smooth elements. 
\end{defn}
\begin{empt}\label{smooth_alg_empt}
	There is a subalgebra $\widetilde{A}_{\text{smooth}} \subset \widetilde{A}$ generated by smooth elements. For any $s > 0$ there is a seminorm $\left\| \cdot \right\|_s$  on $\widetilde{A}_{\text{smooth}}$ given by
	\begin{equation}\label{smooth_seminorms_eqn}
		\left\| \widetilde{a}\right\|_s \bydef \left\|\widetilde{a}^s \right\|=\lim_{\la\in\La}\left\| \pi^s_\la\left(a_\la \right) \right\|.
	\end{equation}
\end{empt}
\begin{defn}\label{smooth_alg_defn}
	The completion of $\widetilde{A}_{\text{smooth}}$ in the topology induced by the seminorms $\left\| \cdot \right\|_s$ is said to be a \textit{smooth algebra} of the coherent set \eqref{spectral_triple_sec_eqn} of spectral triples. This algebra is denoted by $\widetilde{\A}$. 
	We say that the set of spectral triples is \textit{weakly good} if $\widetilde{\A}$ is dense in $\widetilde{A}$. 
\end{defn}
\begin{empt}\label{dirac_inf_constr_empt}
	For any $\widetilde{a}  \in \widetilde{W}^\infty$ and $\xi \in \H^\infty$ we define
	\begin{equation}\label{inf_lift_D_eqn}
		\widetilde D\left(  \widetilde{a}  \ox \xi \right) \bydef \lim_{\la \in \La}\left[D_\la, a_\la \right]\left(  \widetilde{a}  \ox \xi \right)+  \widetilde{a}  \ox D\xi 
	\end{equation}

\end{empt}
\begin{defn}\label{reg_triple_defn}
	Let $\mathfrak{S} \in \mathfrak{FinAlg}$ is given by \eqref{cov_sec_triple_eqn}.    Let $\left(A, \widetilde{A}, G\left(\left.\widetilde{A}~\right|A \right)\right) $  be an $\widehat \pi$-infinite noncommutative covering   of $\mathfrak{S}$.
	Let \eqref{spectral_triple_sec_eqn} be a good coherent set of spectral triples. 
	Let $\widetilde{ D}$ be given by \eqref{inf_lift_D_eqn}.
	We say that $\left( \widetilde{\A}, \widetilde{\H}, \widetilde{D}, \widetilde J\right)$ is a $\left(A, \widetilde{A}, G\left(\left.\widetilde{A}~\right|A \right)\right)$-\textit{lift} of $\left(\A, \H, D, J\right)$.  Also we say that  $\left( \widetilde{\A}, \widetilde{\H}, \widetilde{D}, \widetilde J\right)$ is the \textit{limit} of the coherent set	$\mathfrak{S}_{\left(\sA, \sH, D, J\right)}= \left\{\left(\sA_\la, \sH_\la, D_\la, J_\la \right)\right\}_{\la \in \La}$.
\end{defn}

\chapter{Coverings of commutative $C^*$-algebras}\label{top_chap}

\section{Partitions of unity and compact subsets}

\begin{lemma}\label{top_a_u_lem}
	If $\mathcal X$  is a locally compact, 
	Hausdorff space then for any $x_0 \in  \mathcal X$ and any open neighborhood $\mathcal U\subset\mathcal X$ of $x_0$  there is a continuous function $a: \mathcal X \to \left[0,1 \right]$  such that following conditions hold:
	\begin{itemize}
		\item $a\left(\sX\right) = \left[0,1\right]$,
		\item 	$a\left( x_0\right) = 1$,
		\item $\supp a \subset \mathcal U$,
		\item	there is an open neighborhood $\mathcal V\subset\mathcal U$ of $x_0$ which satisfies to the following condition
		\be\label{com_a_u_eqn}
		a\left(\mathcal V \right)= \{1\}.
		\ee 
	\end{itemize}
	If $\sX$ is locally connected (cf. Definition \ref{top_locally_connected_defn}) then one $a$ can be selected such that the support of $a$ (cf. Definition \ref{top_support_defn}) is connected.
\end{lemma}
\begin{proof}
	From the Exercise \ref{top_completely_regular_exer} it turns out that that $\sX$ is completely regular (cf. Definition \ref{top_completely_regular_defn}), i.e. there is a  continuous function $b: \mathcal X \to \left[0,1 \right]$  such that $b\left( x_0\right) = 1$ and $b\left(\mathcal X \setminus \mathcal U \right)= \{0\}$. A set 
	$\sU' \bydef \left\{x \in \sX | a'\left(x \right) \neq 0\right\}$ is open. 
	The set $\mathcal V = \left\{ x \in \mathcal X~|~c\left(x \right)> \frac{2}{3}  \right\}$ is open.
	If $f: \R \to \R$ is a continuous function given by
	\be\nonumber
f\left(t\right)=	\left\{\begin{array}{c l}
		0 & t \le \frac{1}{3}\\
		3t - 1 & \frac{1}{3} < t \le \frac{2}{3}\\
		1 & t > \frac{2}{3}
	\end{array}
	\right.
	\ee
	then $a \bydef f\left(c \right): \mathcal X \to \left[0,1 \right]$ satisfies to all conditions of this lemma. If $\sX$ is locally connected and $\sU''$ is a component of $x_0$ is $\sU'$ (cf. Definition \ref{top_connected_component_defn}) then   from the Exercise \ref{top_loc_conn_exer} it follows that $\sU''$ is a quasi-component of $\sU'$ (cf. Definition \ref{top_quasi_component_defn}). If we define $a' \in C_c\left(\sX \right)$ by following way
	$$
	a'\left(x \right) \bydef \begin{cases}
	a\left( x\right)& x \in \sU''\\
	0& x \notin \sU''	\end{cases}
	$$ 
\end{proof}
then the support $\supp a'$ is the closure of $\sU''$. From the Theorem \ref{top_connected_closure_thm} it follows that the set $\supp a'$ is connected.
\begin{definition}\label{top_stump_defn}
Let $\mathcal X$ be a locally compact,  
Hausdorff space, $x_0\in \sX$ and $f$ is given by the Lemma \ref{top_a_u_lem} then we denote by
\be\label{top_fx_eqn}
{f}_{{x}_0} \stackrel{\mathrm{def}}{=} 	{f}.
\ee
 We say that ${f}_{{x}_0}$ is an $x_0$-\textit{stump}. Denote by $\mathfrak{Stumps}_{x_0}$ the set of $x_0$-stumps.
\end{definition}
\begin{corollary}\label{com_a_u_cor}
	Let $\mathcal X$ be a locally compact,   
	Hausdorff space. For any $x_0 \in  \mathcal X$, and any open neighborhood $\mathcal U$ of $x_0$ there is an open neighborhood  $\mathcal V$ of $x_0$ such that the closure of $\mathcal V$ is a subset of $\mathcal U$.
\end{corollary}
\begin{proof}
	If $a$ satisfies to the Lemma \ref{top_a_u_lem} then set $\mathcal V = \left\{x \in \mathcal X~|~ a\left(x \right)>0 \right\}$ is open and the closure of $\mathcal V$ is a subset of $\mathcal U$.
\end{proof}

\begin{definition}\label{top_stump_cov_defn}
	Let $p: \widetilde{\sX} \to {\sX}$ be a covering. Let $\widetilde{x}_0\in   \widetilde{\sX}$, and let $\widetilde{\sU}$ be an open neighborhood of $\widetilde{x}_0$ such that the restriction $\left.p\right|_{\widetilde{\sU}}:\widetilde{\sU}\xrightarrow{\approx}\sU \bydef p\left(\widetilde{\sU} \right)$ is a homeomorphism. Since $\widetilde{\sX}$ is locally compact and Hausdorff there is $	\widetilde{f} \in C_c\left(\widetilde{\sX} \right)$ and open subset  $\widetilde{\sV}$ such that $\widetilde{x}_0\in  \widetilde{\sV} \subset {\widetilde{\sU}}$, $~\widetilde{f}\left(\widetilde{\sV} \right)= 1$,  $~\widetilde{f}\left(\widetilde{\sX} \right)= \left[0,1\right]$ and $\supp 	\widetilde{f} \subset \widetilde{\sU}$. We write 
\be\label{top_tfx_eqn}
\widetilde{f}_{\widetilde{x}_0} \stackrel{\mathrm{def}}{=} 	\widetilde{f},
\ee
and we say that $\widetilde{f}_{\widetilde{x}_0}$ is a $p$-$\widetilde x_0$-\textit{stump}.
\end{definition}
\begin{empt}\label{top_fx_empt} 	
Let $\sX$ be a locally compact, Hausdorff space, and let $\sY \subset\sX$ be a compact subset. For any $x \in \sY$ we select an $x$-stump $f_x\in C_c\left(\sX\right)$ (cf. Definition \ref{top_stump_defn}) such that
\begin{itemize}
	\item a set $\supp f_x$ is compact for all $x \in \sY$,
	\item there is an open neighborhood $\sV_x$ of $x$ such that $f_x\left(\sV_x\right)=\{1\}$.
\end{itemize}
A family 
\be\label{top_vx_eqn}
\left\{\sV_x\right\}_{x\in\sX}.
\ee
is such that $\sY \subset \cup_{x \in \sY} \sV_x$.
The set $
\sU_x \bydef \left\{\left. x\in\sX\right| f_x\left( x\right)> 0 \right\}
$ is dense in $\supp f_x$ so $\supp f_x$ is connected (cf. Theorem \ref{top_connected_closure_thm}).  
Moreover $\sV_x \subset \sU_x$ and there is a finite set $\left\{x_1,..., x_n \right\}\subset \sX$ such that $\sY \subset \bigcup_{j=1}^n \sV_{x_j}$, because the set $\sY$ is compact.
If $f \bydef f_{x_1} + ... +  f_{x_n}$ then $f\left(x \right) \ge 1$ for all $x \in \sY$. If 
$$
f_j \bydef \frac{f_{x_j}}{\max\left(1, f\right)}\in C_c\left(\sX \right)_+ 
$$
then
\be\label{top_cfs_eqn}
\sum_{j = 1}^n f_j\left( x\right) = 1\quad\forall x \in \sY.
\ee
\end{empt}
\begin{defn}\label{top_covering_sum_defn}
	Let $\sY\subset\sX$ is a compact subset of a locally compact,  Hausdorff space $\sX$. The  $f \bydef \sum_{j = 1}^n f_j$.
	is said to be a \textit{covering sum for} $\sY$. We also say that the covering sum \eqref{top_cfs_eqn} is \textit{dominated} by the family 
	$\left\{\sV_x\right\}_{x\in\sX}$ (cf. equation \ref{top_vx_eqn}). 
	\end{defn}
\begin{lemma}\label{top_compact_nigthbohood_lem}
	If $\sY\subset\sX$ is a compact subset of a locally compact,  Hausdorff space then there is an open subset $\sU \subset\sX$ with compact closure, such that $\sY \subset \sU$.
\end{lemma}
\begin{proof}
If  $f \bydef \sum_{j = 1}^n f_j\left( x\right)$ is a covering sum of $\sY$  then one has
\begin{itemize}
	\item a set $\sU \bydef \left\{x \in \sX | f\left( x\right)> 0 \right\}$ is open,
	\item the closure of $\sU$ is compact because it equals to  a finite union $\bigcup_{j=1}^n \supp f_j$ of compact sets.
\end{itemize}
\end{proof}

\begin{remark}\label{top_covering_sum_compact_rem}
The support of the covering sum is compact since $\supp f_x$ is compact. 	
\end{remark}
\begin{remark}
Sometimes instead of \eqref{top_cfs_eqn} we use the following notation for a covering  sum
\be\label{top_cfsa_eqn}
\sum_{\a\in \mathscr A} a_\a
\ee
where $\mathscr A$ is a finite set.
\end{remark}
\begin{remark}\label{top_smooth_part_unity_rem} 
	Similarly to the Proposition \ref{top_smooth_part_unity_prop}  if $\sX$ is a paracompact smooth manifold then for every compact set $\sY$ there is covering sum 
	$
	\sum_{j = 1}^n f_j\left( x\right) = 1\quad\forall x \in \sY
	$
	such that $f_j$ is smooth for all $j = 1, ..., n$.
\end{remark}

\begin{definition}\label{top_covering_sum_subordinated_defn}
	Let $p:\widetilde\sX \to \sX$ be a covering. Let both $\sY\subset\sX$ and $\widetilde\sY\subset\widetilde\sX$ are compact subsets. A covering sum $\sum_{j = 1}^n f_j$ for $\sY$ (cf. Definition \ref{top_covering_sum_defn}) is said to be \textit{subordinated  to} $p$ if for all $j=1,..., n$ the set $\left\{\left.x\in\sX\right|f_j\left( x\right)>  0 \right\}$ is evenly covered by $p$ (cf. Definition \ref{top_covering_defn}). A covering sum $\sum_{j = 1}^n \widetilde f_j$ for $\widetilde \sY$ is said to be \textit{subordinated to} $p$ if for all $j=1,..., n$ the set $\widetilde\sU_j\bydef \left\{\left.\widetilde x\in\widetilde\sX\right|\widetilde f_j\left(\widetilde x\right)>  0 \right\}$ is homeomorphically mapped onto $p\left(\widetilde\sU_j\right)$. 
\end{definition}
\begin{definition}\label{top_coveing_triple_defn}
	Let $p:\widetilde\sX \to \sX$ be a covering. Suppose that there are connected open subsets $\widetilde \sU, \widetilde \sV \subset \widetilde \sX$, and a continuous map $\widetilde{s}: \widetilde \sX \to \left[0,1\right]$ such that
	\begin{enumerate}
		\item[(a)] the closure of $\widetilde \sV$ is compact,
		\item[(b)] $\widetilde{   \mathcal V}$ is mapped homeomorphically onto $\sV \bydef p\left(\widetilde{   \mathcal V} \right)$,
		\item[(c)] $\supp \widetilde s \subset \widetilde{   \mathcal V}$,
		\item[(d)] $\widetilde s\left(\widetilde \sU \right) = \{1\}$.
	\end{enumerate}
	We say that $\left(\widetilde \sU, \widetilde \sV, \widetilde s\right)$ is a \textit{covering triple for} $p$.
\end{definition}
\begin{remark}\label{top_covering_sum_subordinated_rem}
From the Lemma \ref{top_a_u_lem} it follows that if the space $\sX$ is Hausdorff, locally compact, and locally connected then for any  $\widetilde x \in \widetilde\sX$ there is a {covering triple $\left(\widetilde \sU, \widetilde \sV, \widetilde s\right)$  for} $p$ such that $\widetilde x \in \widetilde \sX$.
\end{remark}

\begin{empt}\label{ctr_unity_empt}
	If $\sX$ is a compact Hausdorff space and $p: \widetilde\sX \to \sX$ is a transitive finite-fold covering then $\widetilde\sX$ is compact. Let
	\be\nonumber
	1_{C\left(\widetilde{\mathcal X} \right) }= \sum_{\widetilde\a \in \widetilde{\mathscr A} }\widetilde{a}_{\widetilde{\a}}\quad \widetilde{a}_{\widetilde{\a}}\in C\left(\widetilde \sX\right)_+
	\ee
	be a subordinated  to $p$ covering sum for $\widetilde{\mathcal X}$ (cf. Definitions \ref{top_covering_sum_defn}, \ref{top_covering_sum_subordinated_defn} and the equation\eqref{top_cfsa_eqn}). If $\widetilde{e}_{\widetilde{\a}}\bydef \sqrt{\widetilde{a}_{\widetilde{\a}}}\in C\left(\widetilde \sX\right)_+$ then one has
	\be\label{top_ea_eqn}
	\begin{split}
		1_{C\left(\widetilde{\mathcal X} \right) }= \sum_{\widetilde\a \in \widetilde{\mathscr A} }\widetilde{e}^2_{\widetilde{\a}},\\
		\widetilde{e}_{\widetilde{\a}}	\left(g \widetilde{e}_{\widetilde{\a}} \right) = 0 \quad \text{for any nontrivial}\quad g\in G\left(\left.\widetilde \sX \right|\sX \right).
	\end{split}
	\ee
	We also write
	\be\label{top_finite_covering_basis_eqn}
	\begin{split}
		\left\{\widetilde{e}_{\widetilde{\a}}\right\}_{\a  \in \widetilde{\mathscr A}}= \left\{\widetilde{e}_1, ..., \widetilde{e}_n\right\}\subset C\left(\widetilde \sX \right)_+ ,\\
		1_{C\left(\widetilde{\mathcal X} \right) }= \sum_{j=1 }^n\widetilde{e}^2_{\widetilde{j}},\\
		\forall j =1,...,n \quad \widetilde{e}_{j}	\left(g \widetilde{e}_{j} \right) = 0 \quad \text{for any nontrivial}\quad g\in G\left(\left.\widetilde \sX \right|\sX \right).
	\end{split}
	\ee	
\end{empt}

\begin{lemma}\label{top_tietze_ext_lem}
If $\sX$ is a Hausdorff,  locally compact space, and $\sY \subset \sX$ is a compact subspace then the restriction map
\bean
C_c\left(\sX \right) \to C\left(\sY \right) ,\\
f \mapsto f|_{\widetilde\sY}
\eean
is surjective.
\end{lemma}
\begin{proof}
	Let $f \in C\left(\sY \right)$.
If $f' \in C_c\left( \sX\right)$ is a covering sum of $\sY$ (cf. Definition \ref{top_covering_sum_defn}) then from the Theorems \ref{comp_normal_thm} and \ref{tietze_ext_thm} it follows that there if $f'' \in C\left(\supp f'\right)$ such that $ff''|_\sY = f$. However $ff''f' \in C_c\left(\sX \right)$ and  $ff''f'|_\sY = f$.
\end{proof}
\begin{lemma}\label{top_compact_preimage_lem}
	If $p: \widetilde \sX\to \sX$ is a finite-fold covering, 
	and a set $\sY \subset \sX$ is compact then the preimage $p^{-1}\left(\sY \right)$ is also compact.
\end{lemma}
\begin{proof}
	For any $x \in \sY$ we select an open neighborhood $\sU_x$ with compact closure which is evenly covered by $p$ (cf. Definition \ref{top_covering_defn}). From the Corollary \ref{com_a_u_cor}	it follows that there are $\mathcal W_x \subset \sV_x \subset \sU_x$ such that $\mathcal W_x$ is an open neighborhood of $x$ and the set $\subset \sV_x$ is closed. The set $\sV_x$ is compact because it is a subset of the compact closure of $\sU_x$. There is $\left\{x_1, ..., x_n\right\}$ such that $\sY \subset \mathcal W_{x_1}\cup ...\cup \mathcal W_{x_n}\subset \mathcal V_{x_1}\cup ...\cup \mathcal V_{x_n}$. For any $j = 1,..., n$ the set $p^{-1}\left(\sV_j \right)$ is a finite disjoint union of homeomorphic  $\sV_{x_j}$ sets. So $p^{-1}\left(\sV_j \right)$ is compact since $\sV_{x_j}$ is compact. The finite union $\bigcup_{j=1}^n p^{-1}\left(\sV_{x_j} \right)$ of compact sets is compact. From $p^{-1}\left( \sY \right) \subset \bigcup_{j=1}^n p^{-1}\left(\sV_{x_j} \right)$ and taking into account that $p^{-1}\left( \sY \right)$ is closed we conclude that  $p^{-1}\left( \sY \right)$ is compact.
	
\end{proof}
\begin{corollary}\label{top_compact_cc_c0_lem}
	If $p: \widetilde \sX\to \sX$ is a finite-fold covering then there are natural injective $*$-homomorphisms:
	\bea\label{top_compact_cc_eqn}
	C_c\left( p\right) : C_c\left(  \sX\right) \hookto C_c\left( \widetilde \sX\right),\\
	\label{top_compact_cc_c0_eqn}
	C_0\left( p\right) : C_0\left(  \sX\right) \hookto C_0\left( \widetilde \sX\right).
	\eea
\end{corollary}
\begin{proof}
	If $f \in C_c\left(\sX \right)$ then we define  $C_c\left( p\right)\left(f\right)\left(\widetilde x \right) \bydef f\circ p\left(\widetilde{x}\right)$. Since the map $p$ is continuous the map $C_c\left( p\right)\left(f\right)$ is a continuous map from  $\widetilde \sX$ to $\C$. From then Lemma \ref{top_compact_preimage_lem} it follows that $C_c\left( p\right)\left(f\right)\in C_c\left( \widetilde \sX\right)$. It is clear that $C_c\left( p\right)$ is injective $*$-homomorphism. The equation \eqref{top_compact_cc_c0_eqn} follows from \eqref{top_compact_cc_eqn} and the Definition \ref{c_c_closure_defn}, i.e. $C_0\left( p\right)$ can be obtained from $C_c\left( p\right)$ by topological completion.
\end{proof}

\begin{lemma}\label{top_cov_comp_lem}
	Consider a commutative diagram 
	\newline
	\begin{tikzpicture}
		\matrix (m) [matrix of math nodes,row sep=3em,column sep=4em,minimum width=2em]
		{
			\widetilde{\mathcal X}_1 & &\widetilde{\mathcal X}_2\\ 
			& {\mathcal X}\\};
		\path[-stealth]
		(m-1-1) edge node [above] {$p$} (m-1-3)
		(m-1-1) edge node [left]  {$p_1~~$} (m-2-2)
		(m-1-3) edge node [right] {$~~p_2$} (m-2-2);
	\end{tikzpicture}
	\\
	with connected Hausdorff spaces and continuous maps. If both $p_1$ and $p_2$ are finite-fold coverings and the map $p$ is surjective then $p$ is also finite-fold covering.
\end{lemma}
\begin{proof}
	Let $\widetilde x_2 \in \widetilde \sX_2$ and $x = p_2\left(\widetilde x_2 \right)$. Select an open neighborhood $\sU'$ evenly covered by $p_1$ and $p_2$ (cf. Definition \ref{top_covering_defn}). Let $\left\{\widetilde x^1_1, ...,\widetilde x^n_1\right\}\bydef p_2^{-1}\left(\widetilde x_2 \right)$, and let us select for any $j=1,...,n$ an open neighborhood $\widetilde{\sV}^j_1$ of $\widetilde x^j_1$ such that
	$$
	j \neq k \quad \Rightarrow \quad \widetilde{\sV}^j_1\cap \widetilde{\sV}^k_1=\emptyset,
	$$
	so one has
	$$
	\bigcup_{j=1}^n  \widetilde{\sV}^j_1 = \bigsqcup_{j=1}^n  \widetilde{\sV}^j_1
	$$
	A finite intersection 
	$$
	\sU \bydef \sU' \bigcap \left(\bigcap_{j=1}^n p\left(\widetilde{\sV}^j_1 \right)  \right) 
	$$
	is open and  evenly covered by $p_1$. Now one has
	$$
	p^{-1}_1\left( \sU\right) \bigcap \left(\bigsqcup_{j=1}^n  \widetilde{\sV}^j_1 \right) = \bigsqcup_{j=1}^n\left(p^{-1}_1\left( \sU\right)\bigcap   \widetilde{\sV}^j_1 \right) 
	$$
	and for all $j = 1,..., n$ a set $p^{-1}_1\left( \sU\right)\bigcap   \widetilde{\sV}^j_1$ is mapped homeomorphically onto $\sU$. If $ \widetilde\sU_2 \bydef p\left(\bigsqcup_{j=1}^n\left(p^{-1}_1\left( \sU\right)\bigcap   \widetilde{\sV}^j_1 \right) \right)$ then 
	\begin{itemize}
		\item $\widetilde\sU_2$ is an open neighborhood of $\widetilde x_2$,
		\item $p^{-1}\left( \widetilde\sU_2\right) = \bigsqcup_{j=1}^n\left(p^{-1}_1\left( \sU\right)\bigcap   \widetilde{\sV}^j_1 \right)$,
		\item for all $j = 1,..., n$ a set $p^{-1}_1\left( \sU\right)\bigcap   \widetilde{\sV}^j_1$ is mapped homeomorphically onto $\widetilde\sU_2$.
	\end{itemize} 
\end{proof}

\section{Weak fundamental groups}
\paragraph{}
The notion of the fundamental group (cf. Remark \ref{top_homotopy_group_rem}) is relevant to path connected spaces (cf. Definition \ref{top_path_connected_defn}). Here we consider a generalized notion which can be applied for more general connected spaces.

\begin{definition}\label{top_weak_path_defn}
	Let $\sX$ be a topological space.	A \textit{weak path on} $\sX$ is a finite sequence $\left(\sU_1, ..., \sU_n\right)$ of open  connected subsets of $\sX$ such that $\sU_j \cap \sU_{j+1} \neq \emptyset$ for any $j = 1, ..., n - 1$.
\end{definition}
\begin{definition}\label{top_u_path_defn}
	Let $\sX$ be a connected,  locally compact, Hausdorff space.
	Let $\mathfrak{A}\bydef\left\{\sU_\a\right\}_{\a \in \mathscr A}$ be a family of connected open subsets of $\sX$, such that $\sX = \cup~ \sU_\a$ and  a closure of $\sU_\a$ is compact for all $\a\in \mathscr A$. A weak path $\left(\sU_{\a_1},...,\sU_{\a_n}\right)$ (cf. Definition \ref{top_weak_path_defn}) on $\sX$  is said to be an  $\mathfrak{A}$-\textit{path} if  $\sU_j \in \mathfrak{A}$ for all $j = 1,...,n-1$.
\end{definition}

\begin{empt}\label{top_path_inv_empt}
	There is an involution * on the set of paths such that
	\be\label{top_path_inv_eqn}
	\left(\sU_{1},...,\sU_{n}\right)^*\bydef \left(\sU_{n},...,\sU_{1}\right).
	\ee
\end{empt}
\begin{definition}\label{top_path_comp_defn}
	If both $\mathfrak p'=\left(\sU'_{1},...,\sU'_{k}\right)$ and $\mathfrak p''=\left(\sU''_{1},...,\sU_{l}''\right)$ are weak paths such that $\sU'_{k} \cap \sU''_{1}\neq\emptyset$ then we say that a pair $\left( \mathfrak p', \mathfrak p''\right)$ is \textit{composable}, and the path $\mathfrak p=\left(\sU'_{1},...,\sU'_{k},\sU''_{_1},...,\sU_{l}'' \right)$  is said to be the \textit{composition} of $\mathfrak p'$ and $\mathfrak p''$. We write
	\be\label{top_path_comp_eqn}
	\mathfrak p'\circ \mathfrak p''\bydef\mathfrak p.
	\ee
\end{definition}
\begin{remark}\label{top_path_ass_rem}
	The given by the Definition \ref{top_path_comp_defn} composition is associative, i.e.  if both pairs $\left( \mathfrak p', \mathfrak p''\right)$ and $\left( \mathfrak p'', \mathfrak p'''\right)$ are composable then one has
	\be\label{top_path_ass_eqn}
	\left( \mathfrak p'\circ \mathfrak p''\right)  \circ \mathfrak p''' = \mathfrak p'\circ \left( \mathfrak p'' \circ \mathfrak p'''\right) .
	\ee 
\end{remark}

If $p: \widetilde{\sX} \to \sX$ is a covering and $\widetilde{\mathfrak p}\bydef \left(\widetilde\sU_1, ..., \widetilde\sU_n\right)$ is weak path on $\widetilde\sX$ such that $\widetilde\sU_j$ is homeomorphically mapped onto $p\left( \widetilde\sU_j\right)$ for all $j \in\left\{1,...,n\right\}$
then $\mathfrak p\bydef \left(p\left( \widetilde\sU_1\right) , ..., p\left( \widetilde\sU_n\right) \right)$ is a weak path on $\sX$.
\begin{defn}\label{top_path_desc_defn}
	In the above situation we say that $\mathfrak p$ is a $p$-\textit{descent} of $\widetilde{\mathfrak p}$. We write
	\be\label{top_path_desc_eqn}
	\mathfrak p\bydef \desc_p\left( \widetilde{\mathfrak p}\right). 
	\ee
\end{defn}
\begin{definition}\label{top_gen_path_defn}
	Let $\sX$ be a connected, locally connected, locally compact, Hausdorff space.
	Let $\mathfrak{A}\bydef\left\{\sU_\a\right\}_{\a \in \mathscr A}$ be a family of connected open subsets of $\sX$, such that $\sX = \cup~ \sU_\a$ and  a closure of $\sU_\a$ is compact for all $\a\in \mathscr A$. An ordered finite subset $\left(\sU_{\a_1},...,\sU_{\a_n}\right)$ of the family $\left\{\sU_\a\right\}$ is said to be an  $\left\{\sU_\a\right\}$-\textit{path} or simply \textit{path} if  $\sU_{\a_j}\cap \sU_{\a_{j+1}}\neq \emptyset$ for all $j = 1,...,n-1$. We say that a path $\left(\sU_{\a_1},...,\sU_{\a_n}\right)$ is \textit{closed} if $\sU_{\a_1}=\sU_{\a_n}$.
\end{definition}

\begin{lemma}\label{top_gen_path_lem}
Under hypotheses  of the Definition \ref{top_u_path_defn}, for any $\sU_{\a'}, \sU_{\a''} \in \mathfrak{A}$ there is a $\left\{\sU_\a\right\}$-{path} $\left(\sU_{\a_1},...,\sU_{\a_n}\right)$ such that $\sU_{\a_1}=\sU_{\a'}$ and $\sU_{\a_n}=\sU_{\a''}$ .
\end{lemma}

\begin{proof}
	Denote by
	\bean
	\mathfrak{A}_{\a'}\bydef \left\{\sU_{\a''}\in \mathfrak A\left|\exists \text{ path }\left(\sU_{\a_1},...,\sU_{\a_n}\right)\quad\a_1 = \a',\quad \a_n = \a'' \right.\right\},\\
	\mathscr A'\bydef \left\{\left.\a \in \mathscr A\right| \sU_{\a} \in\mathfrak{A}_{\a'} \right\}.
	\eean
	If $\mathscr A'=\mathscr A$ then this lemma is proven. If not   and 
	$$
	\left( \bigcup_{\bt \in\mathscr A'}\sU_{\bt}\right) \bigcap \left( \bigcup_{\a \in\mathscr A\setminus\mathscr A'}\sU_{\a}\right)\neq \emptyset,
	$$
	then there is  $\bt \in\mathscr A'$ and $\a \in\mathscr A\setminus\mathscr A'$ such that $\sU_{\bt}\cap \sU_{\a}\neq \emptyset$. According to the definition of $\mathfrak{A}_{\a'}$ there is a path $\left(\sU_{\a_1},...,\sU_{\a_n}\right)$ such that $\a_1 = \a',\quad \a_n = \bt$. So one has a path $\left(\sU_{\a_1},...,\sU_{\a_n}, \sU_{\a} \right)$, it follows that $\a \in \mathscr A'$ and one has a contradiction. It turns out that 
	\bean
	\left( \bigcup_{\bt \in\mathscr A'}\sU_{\bt}\right) \bigcap \left( \bigcup_{\a \in\mathscr A\setminus\mathscr A'}\sU_{\a}\right)= \emptyset.
	\eean
	Both unions $\bigcup_{\bt \in\mathscr A'}\sU_{\bt}$ and  $\bigcup_{\a \in\mathscr A\setminus\mathscr A'}\sU_{\a}$ of open sets are open, so one has a disjoint union
	$$
	\sX = \left( \bigcup_{\bt \in\mathscr A'}\sU_{\bt}\right) \bigsqcup \left( \bigcup_{\a \in\mathscr A\setminus\mathscr A'}\sU_{\a}\right),
	$$
	i.e.  $\sX$ is not connected. It contradicts with our assumption.	From this contradiction it turns out that $\mathscr A' = \mathscr A$.  
\end{proof}

\begin{corollary}\label{top_connected_union_cor}
	Under the hypotheses of the Definition \ref{top_u_path_defn}  then there is a  directed set $\La$ and  family $\left\{\sU_\la\right\}_{\la\in\La}$ of connected open subsets of $\sX$ such that
	\begin{itemize}
		\item 
		$$
		\sX = \bigcup_{\la\in \La}\sU_\la,
		$$
		\item a closure of $\sU_\la$ is compact for each $\la\in \La$,
		\item
		\be\label{top_lmit_eqn}
		\forall \mu, \nu \in \La  \quad \sU_\mu \cap \sU_\nu \neq \emptyset,
		\ee
		\item
		\be\label{top_lmord_eqn}
		\forall \mu, \nu \in \La \quad \mu \le \nu\quad\Leftrightarrow \quad \sU_\mu \subset \sU_\nu.
		\ee
	\end{itemize}
\end{corollary}
\begin{proof}
	Let $\mathfrak{A}\bydef\left\{\sU_\a\right\}_{\a \in \mathscr A}$ be a given by the Definition \ref{top_u_path_defn} family.
	Select $\a_0 \in \mathscr A$. For any $\a\in \mathscr A$ we define  a path $\left\{\sU_\a\right\}$-{path} $\left(\sU_{\a_1},...,\sU_{\a_n}\right)$ such that $\sU_{\a_1}=\sU_{\a_0}$ and $\sU_{\a_n}=\sU_{\a}$ (cf Lemma \ref{top_gen_path_lem}). For all $\a \in \mathscr A$ a closure of a finite union $\sV_{\a}\bydef \bigcup \sU_{\a_j}$ is compact because a closure of $\sU_{\a_j}$ is compact for every $j = 1,..., n$. From $\sU_{\a_j}\cap \sU_{\a_{j+1}}\neq \emptyset$ for all $j=1,...,n-1$ it turns out that $\sV_{\a}$ is connected, and taking into account $\bigcup_{\a\in \mathscr A}\sU_\a = \sX$ one has 
	\bean
	\forall \a\in \mathscr A \quad  \sU_\a\subset \sV_\a\quad \Rightarrow\quad \bigcup_{\a\in \mathscr A}\sU_\a\subset \bigcup_{\a\in \mathscr A}\sV_\a \quad \Rightarrow\quad\bigcup_{\a\in \mathscr A}\sV_\a= \sX.
	\eean 
	If $\left\{\sU_\la\right\}_{\la \in \La}\subset \sU_\la$ is a family of all finite unions  $$\sU_\la\bydef \bigcup_{j=1}^n \sV_{\a_j}$$ then from $\sU_{\a_0} \subset \sV_{\a}$ it follows that $\sU_\la$ is connected for all $\la \in \La$ and the condition \ref{top_lmit_eqn} holds.
	The condition \eqref{top_lmord_eqn} follows from $\sU_\a \subset \sU_\la$ for all $\la\in \La$.
	A union of two finite unions is also finite union, so one has
	$$
	\forall \mu, \nu \in \La\quad \exists \la \in \La\quad \sU_\mu \cup \sU_\nu = \sU_\la.
	$$
	So if we consider the given by \eqref{top_lmord_eqn} pre-order then $\La$ is a directed set.	
\end{proof}
\begin{corollary}\label{top_compact_union_cor}
	If $\sX$ is a locally compact, locally connected, Hausdorff space then any compact subset of $\sX$ is a subset of connected compact subset of $\sX$.
\end{corollary}
\begin{proof}
	If $\sV$ is a compact set and $\left\{\sU_\la\right\}_{\la\in\La}$ is given by the Corollary \ref{top_connected_union_cor} then from $\sV \subset\sX= \bigcup_{\la\in \La}\sU_\la$ it follows that there is a finite subset $\left\{\la_1, ..., \la_n\right\}\subset \La$ such that $\sV \subset\sX= \bigcup^n_{j =1}\sU_{\la_j}$. The set $\La$ is directed so there is $\la\in\La$ such that $\la \ge \la_j$ for all $j = 1,...,n$. If follows that $\sV \subset \sU_\la$. According to the Corollary the closure $\mathcal W$ if $\sU_\la$ is compact. From the Theorem \ref{top_connected_closure_thm} it follows that the set $\mathcal W$ is connected.
\end{proof}
\
\begin{definition}\label{top_weakly_simply_connected_defn} 
A connected space $\sX$ is said to be \textit{weakly simply connected} if for any covering $p: \widetilde \sX\to \sX$ the space $\sX$ is evenly covered by $p$ (cf. Definition \ref{top_covering_defn})
\end{definition}
\begin{remark}
Any connected, path connected, simply connected  space is weakly simply connected.
\end{remark}

	\begin{defn}\label{top_weakly_semi1_defn}
		An open connected subset $\sU\subset \sX$ of a topological space is said to be \textit{semilocally proper} if it is evenly covered by any transitive covering $\widetilde{\sX}\to\sX$ (cf. Definition \ref{top_covering_defn}).
	A space $\sX$ is said to be \textit{weakly semilocally 1-connected} if for  every point $x\in \sX$ there is a semilocally proper neighborhood. A weak path $\left(\sU_1, ..., \sU_n\right)$ (cf. Definition \ref{top_weak_path_defn}) is \textit{weakly semilocally 1-connected} if $\sU_j$ is semilocally proper for every $j = 1,..., n$.
\end{defn} 
\begin{exercise}
	Prove that any connected, path connected, locally path connected, semilocally 1-connected (cf. Definition \ref{top_semi1_defn}),  space is weakly semilocally 1-connected.
\end{exercise}

\begin{defn}\label{top_x0_path_defn}
Let $\left(\sX, x_0\right)$ be a pointed  space (cf. Definition \ref{top_pointed_defn}) such that $\sX$ is locally connected (cf. Definition \ref{top_locally_connected_defn}), weakly semilocally 1-connected space. A semilocally 1-connected path $\left(\sU_1, ..., \sU_n\right)$ on $\sX$ is said to be $x_0$-\textit{path} if $x_0\in\sU_1$.   Denote by $\mathfrak{Paths} \left( x_0, x\right)$ a set of all $x_0$-paths such that 
	\be\label{top_closed_path_eqn}
	\left( \sU_1, ..., \sU_n\right)\in \mathfrak{Paths} \left( x_0, x\right) \quad \Leftrightarrow \quad x_0 \in \sU_1 \quad x \in \sU_n.
\ee
An $x_0$-{path} $\left(\sU_1, ..., \sU_n\right)$ is \textit{closed} if $x_0\in \sU_n$, i.e.
	\be\nonumber
\left( \sU_1, ..., \sU_n\right)\in \mathfrak{Paths} \left( x_0, x_0\right).
\ee

\end{defn}
\begin{empt}\label{top_lconn_part}
	Let $\sX$ be a connected, locally connected, locally compact, Hausdorff space. Let $p: \widetilde{   \mathcal X }\to \sX$ is a transitive covering. There is a  set $\left\{\sU_\a\right\}_{\a \in \mathscr A}$  set of open connected subsets of $\sX$ such that following conditions hold:
	\begin{itemize}
		\item $\sX = \cup_{\a \in \mathscr A} \sU_\a$.
		\item The closure of $\sU_{\a}$ is a compact set for all $\a \in \mathscr A$.
		\item $\sU_\a$ is evenly covered by $p$ (cf. Definition \ref{top_covering_defn}).
	\end{itemize}
	For any $\a \in \mathscr A$ we select an open connected set $\widetilde \sU_\a \in \widetilde\sX$. One has
	$$
	\widetilde \sX=\bigcup_{\left(g, \a\right)\in  G\left(\left.\widetilde{\sX}~\right|\sX\right)\times \mathscr A} 	g \widetilde \sU_\a
	$$
\end{empt}

\begin{definition}\label{top_lconn_defn}
	
	Under hypotheses \ref{top_lconn_part} we say that the family $\left\{\sU_\a\right\}_{\a \in \mathscr A}$ is \textit{compliant with} $p$, and the family
	$\left\{g \widetilde \sU_\a\right\}_{\left(g, \a\right)\in  G\left(\left.\widetilde{\sX}~\right|\sX\right)\times \mathscr A}$ is the $p$-\textit{lift} of $\left\{\sU_\a\right\}_{\a \in \mathscr A}$. We say that  family $\left\{\sU_\a\right\}_{\a \in \mathscr A}$ is \textit{compliant with coverings} if it is complaint with any finite-fold covering of $\sX$.
\end{definition}
\begin{remark}\label{top_lconn_rem}
	Sometimes we denote by $\widetilde{ \mathscr A}\bydef G\left(\left.\widetilde{\sX}~\right|\sX\right)\times \mathscr A$ and
	\be\label{top_lconn_eqn}	
	\left\{\widetilde \sU_{\widetilde \a }\right\}_{\widetilde{ \a}\in \widetilde{ \mathscr A}}\bydef \left\{g \widetilde \sU_\a\right\}_{\left(g, \a\right)\in  G\left(\left.\widetilde{\sX}~\right|\sX\right)\times \mathscr A}.
	\ee
\end{remark}

\begin{empt}\label{top_path_comp_empt}
	Let $p: \left(\widetilde{\sX}, \widetilde{x_0}\right)\to \left(\sX, x_0\right)$ be a pointed  covering such that both $\widetilde{\sX}$ and $\sX$ are connected, locally connected, locally compact, Hausdorff spaces. Let $\mathfrak A\bydef \left\{\sU_\a\right\}_{\a\in \mathscr A}$ be a family of all connected, evenly covered by $p$ open sets with compact closure, and let $\widetilde{\mathfrak A}$ be a $p$-lift of $\mathfrak A$ (cf. Definition \ref{top_lconn_defn}). Select $\widetilde\sU_0 \in \widetilde{\mathfrak A}$ such that $\widetilde x_0\in \widetilde\sU_0$. Let $\sU_0 = p\left(\widetilde\sU_0 \right)$. From the Lemma \ref{top_gen_path_lem} it follows that for all $g \in G\left(\left.\widetilde{   \mathcal X }~\right|{   \mathcal X } \right)$ there is a $\widetilde{\mathfrak A}$-path $\widetilde{\mathfrak p}_g=\left(\widetilde\sU_{\widetilde\a_1},...,\widetilde\sU_{\widetilde\a_n}\right)$ such that
	\be\label{top_path_compl_eqn}
	\widetilde\sU_{\widetilde\a_1}= \widetilde\sU_{0}; \quad g \widetilde x_0 \in  \widetilde\sU_{\widetilde\a_n}.
	\ee
	A path which satisfies to  \eqref{top_path_compl_eqn} is not unique. One has
	\be\label{top_path_comp_lem_eqn}
	\begin{split}
		\widetilde{\mathfrak p}_g=\left(\widetilde\sU_{\widetilde\a_1},...,\widetilde\sU_{\widetilde\a_n}\right)\quad \Rightarrow p\left( \widetilde x_0\right) \in  p\left( \widetilde\sU_{\widetilde\a_1}\right)\cap p\left( \widetilde\sU_{\widetilde\a_n}\right)
	\end{split}
	\ee
	If $g', g'' \in G\left(\left.\widetilde{   \mathcal X }~\right|{   \mathcal X } \right)$ then from \eqref{top_path_comp_lem_eqn} it follows that there is the composition $\desc_p \widetilde{\mathfrak p}_{g'}\circ \desc_p \widetilde{\mathfrak p}_{g''}$ (cf. Definition \ref{top_path_comp_defn}).
\end{empt}
\begin{definition}\label{top_closed_path_corr_defn}
	Under the hypotheses  \ref{top_path_comp_empt} we say that the closed path $\widetilde{\mathfrak p}_g$ \textit{corresponds to} $g \in G\left(\left.\widetilde{   \mathcal X }~\right|{   \mathcal X } \right)$. A closed path is said to be $p$-\textit{contractible} if it corresponds to the trivial element of $G\left(\left.\widetilde{   \mathcal X }~\right|{   \mathcal X } \right)$.
\end{definition}
\begin{rem}\label{top_inv_path_rem}
	If a weak path $\widetilde{\mathfrak p}$ {corresponds to} $g \in G\left(\left.\widetilde{   \mathcal X }~\right|{   \mathcal X } \right)$ then its involution $\widetilde{\mathfrak p}^*$ {corresponds to} $g^{-1} \in G\left(\left.\widetilde{   \mathcal X }~\right|{   \mathcal X } \right)$.
\end{rem}

\begin{lemma}\label{top_lift_comp_lem}
Under hypotheses \ref{top_path_comp_empt} and suppose that $\widetilde{\mathfrak p}_{g'} = \left( \widetilde\sU_{\widetilde\a'_1},..., \widetilde\sU_{\widetilde\a'_k}\right)$, $\widetilde{\mathfrak p}_{g''} = \left( \widetilde\sU_{\widetilde\a''_1},..., \widetilde\sU_{\widetilde\a''_l}\right)$. If
	\bean
	\widetilde{\mathfrak p} = \lift^{\widetilde\sU_0}_p\left( \desc_p \widetilde{\mathfrak p}_{g'}\circ \desc_p \widetilde{\mathfrak p}_{g''}\right),\\ 
	\widetilde{\mathfrak p}= \left( \widetilde\sU_{\widetilde\a_1},..., \widetilde\sU_{\widetilde\a_{k + l}}\right)
	\eean
	then $ \widetilde\sU_{\widetilde\a_{k + l}}= g'g''\widetilde \sU_0$.
\end{lemma}
\begin{proof}
	From $\widetilde\sU_{\widetilde\a'_k}= g' \widetilde\sU_0$ it follows that
	$$
	\widetilde{\mathfrak p} = \left( \widetilde\sU_{\widetilde\a'_1},..., \widetilde\sU_{\widetilde\a'_k}, g'\widetilde\sU_{\widetilde\a''_1},..., g'\widetilde\sU_{\widetilde\a''_l}\right)
	$$
	and taking into account $g'\widetilde\sU_{\widetilde\a''_1}= \widetilde\sU_{\widetilde\a'_k}= g' \widetilde\sU_0$, $\quad\widetilde\sU_{\widetilde\a''_l}=g''\widetilde\sU_0$ one has $ \widetilde\sU_{\widetilde\a_{k + l }}= g'g''\widetilde \sU_0$.
\end{proof}

\begin{empt}\label{top_g_homo_empt}
	In the situation  \ref{top_path_comp_empt}	consider a commutative triangle of connected topological spaces and  transitive coverings.
	\newline
	\begin{tikzpicture}
		\matrix (m) [matrix of math nodes,row sep=3em,column sep=4em,minimum width=2em]
		{
			\widetilde{\mathcal X}' & &\widetilde{\mathcal X}''\\ 
			& {\mathcal X}\\};
		\path[-stealth]
		(m-1-1) edge node [above] {$p$} (m-1-3)
		(m-1-1) edge node [left]  {$p'~~$} (m-2-2)
		(m-1-3) edge node [right] {$~~p''$} (m-2-2);
	\end{tikzpicture}
	\\
	Suppose that $\sX$ is a  locally connected, locally compact, Hausdorff space. Consider a described in \ref{top_lconn_part} family $\left\{\sU_\a\right\}_{\a \in \mathscr A}$ of open subsets of $\sX$. Let both  $\left\{\widetilde \sU_{\widetilde \a' }\right\}_{\widetilde{ \a'}\in \widetilde{\mathscr A}'}~$, $\left\{\widetilde \sU_{\widetilde \a'' }\right\}_{\widetilde{ \a''}\in \widetilde{\mathscr A}''}~$ be the  $p'$-lift and $p''$-lift of $\left\{\sU_\a\right\}_{\a \in \mathscr A}$ respectively (cf. Definition \ref{top_lconn_defn} and equation \eqref{top_lconn_eqn}). Let $\widetilde \sU'_0 \in \left\{\widetilde \sU_{\widetilde \a' }\right\}$, and $\widetilde x_0 \in \widetilde \sU'_0$. Denote by $ \widetilde \sU''_0\bydef p\left(\widetilde \sU'_0\right)\subset \widetilde \sX''$, $  \sU_0\bydef p'\left(\widetilde \sU'_0\right)\subset \sX$, $\widetilde x''_0\bydef p\left(\widetilde x'_0\right)\subset \widetilde \sU''_0$, $  x_0\bydef p'\left(\widetilde x'_0\right)\subset \sU_0$. Both $p'$ and $p''$ are transitive coverings, hence from the Corollary \ref{top_trans_tp_cov_cor} it follows that there are bijecive maps
	\bean 
	G\left(\left.\widetilde{   \mathcal X }'~\right|~{   \mathcal X } \right)\xrightarrow{\approx}p'^{-1}p'\left(\widetilde x'_0 \right),\\
	g \mapsto g\widetilde x'_0;\\
	G\left(\left.\widetilde{   \mathcal X }''~\right|~{   \mathcal X } \right)\xrightarrow{\approx}p''^{-1}p''\left(\widetilde x''_0 \right),\\
	g \mapsto g\widetilde x''_0. 
	\eean
	On the other hand the map $p$ induces the surjective map $p'^{-1}p'\left(\widetilde x'_0 \right)\to p''^{-1}p''\left(\widetilde x''_0 \right)$ so one has the surjective map
	\be\label{top_g_homo_eqn}
	h :  G\left(\widetilde{   \mathcal X }'~|~{   \mathcal X } \right)\to G\left(\widetilde{   \mathcal X }''~|~{   \mathcal X } \right)
	\ee
\end{empt}
\begin{definition}\label{top_path_corr_defn}
	Under hypotheses  \ref{top_path_comp_empt} we say that the closed path $\widetilde{\mathfrak p}_g$ \textit{corresponds to} $g \in G\left(\left.\widetilde{   \mathcal X }~\right|{   \mathcal X } \right)$. 
\end{definition}

\begin{lemma}\label{top_homo_lem}
Under hypotheses \ref{top_g_homo_empt} following condition holds:
	\begin{itemize}
		\item [(i)] 	The given by \eqref{top_g_homo_eqn} surjective map $ h :  G\left(\widetilde{   \mathcal X }'~|~{   \mathcal X } \right)\to G\left(\widetilde{   \mathcal X }''~|~{   \mathcal X } \right)$ is a homomorphism of groups.
		\item[(ii)] If we consider the natural inclusion $ G\left(\widetilde{   \mathcal X }'~|~{   \mathcal X }'' \right)\subset  G\left(\widetilde{   \mathcal X }'~|~{   \mathcal X } \right)$ then $ G\left(\widetilde{   \mathcal X }'~|~{   \mathcal X }'' \right)$ is a normal subgroup of $ G\left(\widetilde{   \mathcal X }'~|~{   \mathcal X } \right)$ (cf. Definition \ref{normal_subgroup_defn}).
	\end{itemize}
\end{lemma}
\begin{proof} (i)
	Let $g_1, g_2 \in G\left(\widetilde{   \mathcal X }'~|~{   \mathcal X } \right)$ and let both $\widetilde{\mathfrak p}'_{g_1}= \left( \widetilde\sU'_{\widetilde\a'_{1,1}},..., \widetilde\sU'_{\widetilde\a'_{1, k}}\right)$, $\widetilde{\mathfrak p}'_{g_2}= \left( \widetilde\sU'_{\widetilde\a'_{2,1}},..., \widetilde\sU'_{\widetilde\a'_{2, l}}\right)$ are $\left\{\widetilde \sU'_{\widetilde \a' }\right\}$-paths which correspond to  $g_1$ and $g_2$ respectively (cf. Definition \ref{top_path_corr_defn}). Clearly both $\desc_p\widetilde{\mathfrak p}'_{g'_1}$ and $\desc_p \widetilde{\mathfrak p}'_{g'_2}$ are $\left\{\widetilde \sU_{\widetilde \a'' }\right\}$-paths which correspond to  $h\left( g_1\right) $ and $h\left( g_2\right)$ respectively. From the Lemma  \ref{top_lift_comp_lem} it follows that both 
	\bean
	\widetilde{\mathfrak p}' \bydef \lift^{\widetilde\sU'_0}_{p'}\left( \desc_{p'} \widetilde{\mathfrak p}_{g_1}\circ \desc_{p'} \widetilde{\mathfrak p}_{g_2}\right),\\
	\widetilde{\mathfrak p}''\bydef \lift^{\widetilde\sU''_0}_{p''}\left( \desc_{p'} \widetilde{\mathfrak p}_{g_1}\circ \desc_{p'} \widetilde{\mathfrak p}_{g_2}\right)
	\eean
	are $\left\{\widetilde \sU'_{\widetilde \a' }\right\}$ and $\left\{\widetilde \sU''_{\widetilde \a'' }\right\}$-paths which correspond to $g_1 g_2$ and $h\left( g_1\right)h\left( g_2\right)$ respectively. It means that that   $\widetilde{\mathfrak p}'= \left( \widetilde\sU'_{\widetilde\a'_{1}},..., \widetilde\sU'_{\widetilde\a'_{k + l }}\right)$ then $\widetilde\sU'_{\widetilde\a'_{k + l }}= g_1g_2\widetilde \sU'_0$. On the other hand  $\widetilde{\mathfrak p}''= \desc_p \widetilde{\mathfrak p}'= \left( p\left( \widetilde\sU'_{\widetilde\a'_{1}}\right) ,..., p\left( \widetilde\sU'_{\widetilde\a'_{k + l }}\right) \right)$, so $p\left( \widetilde\sU'_{\widetilde\a'_{k + l}}\right)= h\left(g_1g_2 \right) \widetilde\sU''_0$, i.e. $\widetilde{\mathfrak p}''$ corresponds to $h\left(g_1g_2 \right)$. It follows that $h\left(g_1g_2 \right)=h\left(g_1 \right)h\left(g_2 \right)$, i.e. $h$ is a homomorphism.\\
	(ii) Follows from $G\left(\widetilde{   \mathcal X }''~|~{   \mathcal X }' \right)= \ker h$.
\end{proof}
\begin{rem}
The proof of the Lemma \ref{top_homo_lem} is an alternative version of the proof of the Theorem \ref{top_covp_cat_thm}.
\end{rem}

\begin{definition}\label{top_contactible_path_defn}
Let 
$\left(\sX, x_0\right)$ be  a pointed  space.  A closed  semilocally 1-connected path $\mathfrak p \bydef \left(\sU_1, ..., \sU_n\right)\in \mathfrak{Paths} \left( x_0, x_0\right) $ (cf. Definitions \ref{top_weak_path_defn} and \ref{top_x0_path_defn}) is said to be \textit{contractible} if  for any  pointed transitive covering  $p: \left(\widetilde\sX, \widetilde x_0\right)\to\left(\sX, x_0\right)$ it is $p$-contractible (cf. Definition \ref{top_closed_path_corr_defn}).
\end{definition}
\begin{rem}\label{top_path_inv_rem}
If $\mathfrak p$ is a closed semilocally 1-connected $x_0$-path then both compositions $\mathfrak p^* \circ \mathfrak p$ and $\mathfrak p \circ \mathfrak p^*$ are contractible (cf. the Remark \ref{top_inv_path_rem} and the Lemma \ref{top_homo_lem}).
\end{rem}
\begin{definition}\label{top_path_equiv_defn}
We say that both  weakly semilocally 1-connected  $x_0$-paths $\mathfrak p_1,\mathfrak p_2\in \mathfrak{Paths} \left( x_0, x\right)$ are \textit{homotopically equivalent} if  a composition $\mathfrak p_1 \circ \mathfrak p^*_2\in \mathfrak{Paths} \left( x_0, x_0\right)$ is contractible. We write  $\mathfrak p_1 \sim_{\left( x_0, x\right)} \mathfrak p_2$.
\end{definition}
\begin{lemma}\label{top_path_lift_lem}
	If ${\mathfrak p}=\left(\sU_{\a_1},...,\sU_{\a_n}\right)$ is a  $\left\{\sU_\a\right\}$-{path}  then for any $\widetilde \sU_{\widetilde \a' }\in \left\{\widetilde \sU_{\widetilde \a }\right\}$ such that $p\left(\widetilde \sU_{\widetilde \a' } \right) = \sU_{\a_1}$ there is the unique  $\left\{\widetilde \sU_{\widetilde \a }\right\}$-path $\widetilde{\mathfrak p}$ such that ${\mathfrak p}=\desc_p~\widetilde{\mathfrak p}$ and the first element of $\widetilde{\mathfrak p}$ is $\widetilde \sU_{\widetilde \a' }$.
\end{lemma}
\begin{proof}
	From the condition of this lemma the first element of $\widetilde{\mathfrak p}$ is $\widetilde \sU_{\widetilde \a' }$, other elements will be constructed by the induction.
	Suppose that we already have $j$ elements $\widetilde\sU_{\widetilde\a_1},...,\widetilde\sU_{\widetilde\a_j}$ of the path $\widetilde{\mathfrak p}$. If $\widetilde\sU' \in \left\{\widetilde \sU_{\widetilde \a }\right\}$ is such that $p\left( \widetilde\sU'\right) = \sU_{\a_{j+1}}$ then from $\sU_{\a_j}\cap \sU_{\a_{j+1}}\neq \emptyset$ it turns out that there is the unique $g \in G\left(\left.\widetilde{   \mathcal X }~\right|{   \mathcal X } \right)$ such that $g\widetilde\sU'\cap \widetilde\sU_{\widetilde\a_j}\neq\emptyset$. We set $g\widetilde\sU'$ as $j+1^{\text{th}}$ element of $\widetilde{\mathfrak p}$, i.e. $\sU_{\widetilde\a_{j+1}}\bydef g\widetilde\sU'$.  In result the path
	$$
	\widetilde{\mathfrak p} \bydef \left(\widetilde\sU_{\widetilde\a_1},...,\widetilde\sU_{\widetilde\a_n}\right)
	$$
	satisfies to the requirements of this lemma. The above  construction is unique so $\widetilde{\mathfrak p}$ is uniquely defined by the conditions of this lemma.
\end{proof}

\begin{definition}\label{top_path_lift_defn}
	The given by the Lemma \ref{top_path_lift_lem} path $\widetilde{\mathfrak p}$ is said to be the $p$-$\sU_{\widetilde \a' }$-\textit{lift} of $\mathfrak p$. We write
	\be\label{top_path_lift_eqn}
	\lift_p^{\sU_{\widetilde \a' }}~\mathfrak p\bydef\widetilde{\mathfrak p}.
	\ee 
\end{definition}

\begin{empt}\label{top_pp_lift_empt}
	For any  $\mathfrak p = \left( \sU_1, ..., \sU_n\right)\in \mathfrak{Paths} \left( x_0, x\right)$ (cf. the equation \eqref{top_closed_path_eqn}) and any  pointed  transitive covering  $p: \left( \widetilde{\sX}, \widetilde x_0\right) \to \left( \sX, x_0\right) $ there is  a unique $p$-$\widetilde x_0$-lift $$\lift_p^{\widetilde x_0}= \left( \widetilde\sU_1, ...,\widetilde \sU_n\right)$$ (cf. the Definition \eqref{top_path_lift_defn}). There is a unique $\widetilde x \in \widetilde \sU_n$ such that $p\left( \widetilde x\right) = x$.
\end{empt}
\begin{definition}\label{top_pp_lift_defn}
Under the hypotheses \ref{top_pp_lift_empt} we say that $\widetilde x$ is the $p$-$\left(\widetilde x_0, x \right)$-\textit{lift} of $\mathfrak p$. We write
\be\label{top_pp_lift_eqn}
\lift_p^{\left( \widetilde x_0, x\right) }\left( \mathfrak p \right) \bydef\widetilde x.
\ee
\end{definition}
\begin{exercise}\label{top_path_equiv_exer}
Prove following statements:
\begin{enumerate}
	\item Both  weakly semilocally 1-connected  $x_0$-paths $\mathfrak p_1,\mathfrak p_2\in \mathfrak{Paths} \left( x_0, x\right)$ are {homotopically equivalent} $\mathfrak p_1 \sim_{\left( x_0, x\right)} \mathfrak p_2$ (cf. Definition \ref{top_path_equiv_defn}) if and only if for any pointed  transitive covering $p: \left( \widetilde{\sX}, \widetilde x_0\right) \to \left( \sX, x_0\right)$ one has
\be\label{top_pe_lift_eqn}
\lift_p^{\left( \widetilde x_0, x\right) }\left( \mathfrak p_1 \right) =\lift_p^{\left( \widetilde x_0, x\right) }\left( \mathfrak p_2 \right).
\ee
(cf. the equation \eqref{top_pp_lift_eqn}).
\item The  given by the Definition \ref{top_path_equiv_defn} relation $\sim_{\left( x_0, x\right)}$ is an equivalence relation (cf.  Definition \ref{equivalence_relation_defn}).
\end{enumerate}
\end{exercise}
\begin{remark}
	Any pair of closed $x_0$-paths is composable (cf. Definitions \ref{top_path_comp_defn} and \ref{top_closed_path_corr_defn}). Moreover if $\mathfrak p \in  \mathfrak{Paths} \left( x_0, x_0\right)$ then the involution $\mathfrak p^*$ of $\mathfrak p$ also lies in $\mathfrak{Paths} \left( x_0, x_0\right)$.
\end{remark}

\begin{notation}\label{top_path_equiv_not}
From the Exercise \ref{top_path_equiv_exer} it follows that there is a quotient set $\mathfrak{Paths} \left( x_0, x\right)/\sim_{\left( x_0, x\right)}$ (cf. Definition \ref{top_path_equiv_exer}). We use the following notation $$\left[\mathfrak{Paths} \left( x_0, x\right)\right]\bydef \mathfrak{Paths} \left( x_0, x\right)/\sim_{\left( x_0, x\right)}.$$ For any $\mathfrak p\in  \mathfrak{Paths} \left( x_0, x\right)$ we denote by $\left[\mathfrak p\right]\in \left[\mathfrak{Paths} \left( x_0, x\right)\right]$ a represented by $\mathfrak p$ element.
\end{notation}

\begin{empt}
	If both  $\mathfrak p'$ and  $\mathfrak p''$ are contractible then from the Exercise \ref{top_path_equiv_exer} it follows that they are homotopically equivalent, i.e.  $\mathfrak p'\sim_{\left( x_0, x\right)}\mathfrak p''$, or $\left[\mathfrak p'\right]= \left[\mathfrak p''\right]$. Denote by $e \in \left[\mathfrak{Paths} \left( x_0, x_0\right)\right]$ a represented by contractible path element. For any pair $\left[\mathfrak p'\right], \left[\mathfrak p''\right]\in  \left[\mathfrak{Paths} \left( x_0, x_0\right)\right]$ we define the composition $\left[\mathfrak p'\right]\circ\left[\mathfrak p''\right]\bydef\left[\mathfrak p'\circ\mathfrak p''\right]$ where the composition  $\mathfrak p'\circ\mathfrak p''$ is given by the Definition \eqref{top_path_comp_defn}. From the Exercise \ref{top_path_equiv_exer} it turns out that the product $\left[\mathfrak p'\right]\circ\left[\mathfrak p''\right]$ does not depend on representatives of both $\left[\mathfrak p'\right]$ and $\left[\mathfrak p''\right]$. Using \eqref{top_path_ass_eqn}
one can prove that the $\circ$-operation on $\left[\mathfrak{Paths} \left( x_0, x_0\right)\right]$ is associative.	
	From the Definition \ref{top_path_equiv_defn} it follows that for all $\left[\mathfrak p\right]\in \left[\mathfrak{Paths} \left( x_0, x_0\right)\right]$ one has
	\bean
	\left[\mathfrak p\right]\circ  	\left[\mathfrak p^*\right]= 	\left[\mathfrak p^*\right]\circ  	\left[\mathfrak p\right]= e,\\
	\left[\mathfrak p\right]\circ e = e \circ \left[\mathfrak p\right] = \left[\mathfrak p\right].
	\eean
	In result the set $\left[\mathfrak{Paths} \left( x_0, x_0\right)\right]$ becomes a group with respect to the following operations of a multiplication 
	$$
\left( \left[\mathfrak p'\right],  	\left[\mathfrak p''\right]\right) \mapsto 	\left[\mathfrak p'\circ\mathfrak p''\right]
	$$
	and the inversion
$$
\left[\mathfrak p\right]\mapsto \left[\mathfrak p^*\right].
$$	
\end{empt}

\begin{definition}\label{top_weak_fundamental_group_defn}
If
$\left(\sX, x_0\right)$ is  a pointed  space, such that $\sX$ is a connected, locally connected (cf. Definition \ref{top_locally_connected_defn}) and weakly semilocally 1-connected space
then the described above group is said to be a \textit{weak fundamental group} of $\sX$. We denote it by $\pi_1^{\mathrm{w}}\left(\sX, x_0 \right) $.
\end{definition}
\begin{exercise}\label{top_fundamental_group_mor_exer}
	Let  $p: \left( \widetilde{\sX}, \widetilde{x}_0\right) \to \left(\sX, x_0\right)$ be a pointed  covering where $\sX$ is a connected, locally connected  $q_1\left(\widetilde g'_1 \widetilde g''_1 \right)= \widetilde g'_1 \widetilde g''_1\widetilde x_1$and weakly semilocally 1-connected space. Prove that the $p$-descent (cf. Definition  \ref{top_path_desc_defn}) yields an injective homomorphism $\pi_1^{\mathrm{w}}\left(p \right): \pi_1^{\mathrm{w}}\left(\widetilde\sX, \widetilde x_0 \right)\hookto \pi_1^{\mathrm{w}}\left(\sX, x_0 \right)$. Prove a generalization of the Theorem \ref{top_fundamental_group_mor_thm}.
\end{exercise}
\begin{exercise}
	Let $\sX$ be a connected, locally path-connected  and  semilocally 1-connected (cf. Definition \ref{top_semi1_defn}) Hausdorff space. Prove that there is a natural isomorphism
	\be
	\pi_1^{\mathrm{w}}\left(\sX, x_0 \right)\cong \pi_1\left(\sX, x_0 \right)
	\ee
	between the weak fundamental group, and the explained in the Remark \ref{top_homotopy_group_rem} one.
\end{exercise}
Let $\left(\sX, x_0\right)$ be  a pointed space, such that $\sX$ is a connected, locally connected (cf. Definition \ref{top_locally_connected_defn}) and weakly semilocally 1-connected space. 
Let $\widetilde \sX$ be a set of pairs $\left(x, \left[\mathfrak p\right]_{\left( x_0, x\right)}\right)$ where $$\left[\mathfrak p\right]_{\left( x_0, x\right)}\in  \mathfrak{Paths} \left( x_0, x\right)/\sim_{\left( x_0, x\right)}.$$ There is a natural map 
\be
\begin{split}
\widetilde p: \widetilde \sX\to \sX,\\
\left(x, \left[\mathfrak p\right]_{\left( x_0, x\right)}\right)\mapsto x.
\end{split}
\ee
For any $\mathfrak p = \left( \sU_1, ..., \sU_n\right)\in \mathfrak{Paths} \left( x_0, x\right)$ we define a set 
\be\label{top_u_path_eqn}
\widetilde \sU_{\mathfrak p} \bydef \left\{\left.  \left( x',  \left[\mathfrak p\right]_{\left( x_0, x'\right)}\right)\in \widetilde{\sX} \right|x' \in \sU_n\right\}.
\ee
\begin{lemma}\label{top_uni_top_lem}
If $\sX$ is a connected,  weakly semilocally 1-connected (cf. Definition \ref{top_weakly_semi1_defn}) space then a family of sets given by the equation \eqref{top_u_path_eqn} a {basis} for a topology on $\widetilde\sX$. 
\end{lemma}
\begin{proof}
One needs check conditions (a) and (b) of the Definition \ref{top_base_defn}. The condition (a) is evident, let us check the condition (b). Let both ${\mathfrak p'}\bydef \left( \sU'_1, ..., \sU'_{n'}\right)$ and ${\mathfrak p''}\bydef \left( \sU''_1, ..., \sU''_{n''}\right)$ be paths  such that there is $\widetilde x \in \widetilde \sU_{\mathfrak p'}\cap \widetilde \sU_{\mathfrak p''}$. It follows   that $\widetilde p\left(\widetilde x \right) \in \sU'_{n'}\cap \sU''_{n''}$. From 
\eqref{top_u_path_eqn} it turns out that
\be\label{top_up_path_eqn}
\left[\mathfrak p'\right]_{\left( x_0, \widetilde p\left(\widetilde x \right)\right)}= \left[\mathfrak p''\right]_{\left( x_0, \widetilde p\left(\widetilde x \right)\right)}
\ee
From the Definition \ref{top_weakly_semi1_defn} it turns out that there is a weakly 1-connected neighborhood $\sU$ of $p\left(\widetilde x \right)$ which is a subset of $\sU'_{n'}\cap \sU''_{n''}$. If ${\mathfrak p}\bydef \left( \sU'_1, ..., \sU'_{n'}, \sU\right)$ then from \eqref{top_up_path_eqn} one has
\be\label{top_upp_path_eqn}
\left[\mathfrak p\right]_{\left( x_0, \widetilde p\left(\widetilde x \right)\right)}=\left[\mathfrak p'\right]_{\left( x_0, \widetilde p\left(\widetilde x \right)\right)}= \left[\mathfrak p''\right]_{\left( x_0, \widetilde p\left(\widetilde x \right)\right)}.
\ee
From \eqref{top_upp_path_eqn} it follows that $\widetilde x \in \widetilde \sU_{\mathfrak p}$ and  $\widetilde\sU_{\mathfrak p}\subset \widetilde\sU_{\mathfrak p'}\cap \widetilde\sU_{\mathfrak p''}$.
\end{proof}
\begin{lemma}
If we consider the given by the Lemma \ref{top_uni_top_lem} topology on $\widetilde \sX$ then the map $p: \widetilde \sX \to \sX$ is a covering.
\end{lemma}
\begin{proof}
For any  $x \in \sX$  there is an open weakly 1-connected neighborhood $\sU$.  If $\widetilde x \in p^{-1}\left(x\right)$ then $\widetilde{x}= \left( x,  \left[\mathfrak p\right]_{\left( x_0, x\right)}\right)$. If $\mathfrak p= \left( \sU_1, ..., \sU_{n}\right)$ and $\mathfrak p'= \left( \sU_1, ..., \sU_{n}, \sU\right)$ and $\widetilde \sU_{\mathfrak p'}$ is given by  \eqref{top_u_path_eqn} then $p\left( \widetilde \sU_{\mathfrak p'}\right)= \sU$ and $\widetilde \sU_{\mathfrak p'}$ is mapped homeomorphically onto $\sU$. So for all $\widetilde x \in p^{-1}\left(x\right)$ one can find an open neighborhood $\widetilde \sU_{\widetilde x}$ which is is mapped homeomorphically onto $\sU$, i.e. $\sU$ is evenly covered by $p$.
\end{proof}
Let $p': \widetilde \sX'\to \sX$ is a covering then we define a map $\widetilde p': \widetilde \sX \to \widetilde \sX'$ by a following way. For any $\widetilde{x}= \left( x,  \left[\mathfrak p\right]_{\left( x_0, x\right)}\right)$ we define a unique $p'$-lift  $\left( \widetilde\sU'_1, ..., \widetilde\sU'_{n}\right)$ of $\mathfrak p = \left( \sU_1, ..., \sU_{n}\right)$. There is a unique $\widetilde x' \in \widetilde \sX'$ such that $\left\{\widetilde x'\right\}= p^{-1}\left( x\right) \cap \widetilde\sU'_{n}$. So one can define a map 
\be\label{top_from_uni_defn}
\begin{split}
\widetilde p': \widetilde \sX \to \widetilde \sX',\\
\widetilde x \mapsto  \widetilde x'.
\end{split}
\ee
\begin{exercise}\label{top_uni_cov_exer}
Prove that the given by \eqref{top_from_uni_defn} map is a covering.
\end{exercise}
From the Exercise \ref{top_uni_cov_exer} it follows that $\widetilde \sX$ is a universal covering space of $\sX$ (cf. Definition \ref{top_universal_covering_defn}).
So one has the following Theorem.
\begin{theorem}\label{top_simply_con_cov_thm}If $\sX$ is a connected, locally connected (cf. Definition \ref{top_locally_connected_defn}), weakly semilocally 1-connected space then there is an universal covering space for $\sX$.
\end{theorem}

\begin{remark}
The Theorem \ref{top_simply_con_cov_thm} is a generalization of the Lemma \ref{top_simply_con_cov_lem}.
\end{remark}
\begin{exercise}\label{top_weak_covering_iso_exer}
Let $\left(\sX, x_0\right)$ be  a  pointed space, such that $\sX$ is a connected, locally connected (cf. Definition \ref{top_locally_connected_defn}) and weakly semilocally 1-connected space, and let $p: \widetilde\sX\to\sX$ be an universal covering. Prove that there is a following natural isomorphism
$$
\pi^{\text{w}}_1\left(\sX, x_0\right)\cong G\left(\left. \widetilde \sX\right|\sX \right)  
$$
where $\pi^{\text{w}}_1\left(\sX, x_0\right)$ is a weak fundamental group (cf. Definition \ref{top_weak_fundamental_group_defn}) and $G\left(\left. \widetilde \sX\right|\sX \right)$ is a group of covering transformations (cf. Definition \ref{top_group_of_covering_transformations_defn}). 

\end{exercise}
\begin{lemma}\label{top_lex_square_lem}
The ordered square $I^2_o$ (cf. \ref{top_lex_square_empt}) is weakly simply connected (cf. Definition \ref{top_weakly_simply_connected_defn}).
\end{lemma}

\begin{proof}
It is known that the space $I^2_o$ is locally connected (cf. Definition \ref{top_locally_connected_defn}). 
Let $p: \widetilde \sX \to I^2_o$ be a covering. 
For any $x \in I^2_o$ we select an open connected neighborhood $\sU_x$ evenly covered by $p$ (cf. Definition \ref{top_covering_defn}). 
From the Lemma \ref{top_gen_path_lem} it follows that there is a $\left\{\sU_x\right\}_{x \in \sX}$-path $\left( \sU_1, ..., \sU_{n}\right)$ 
(cf. Definition \ref{top_u_path_defn})  from $\left(0,0\right)$ to $\left(1,1\right)$. 
If there is $x \in I^2_o\setminus\cup_{j = 1}^n\sU_j$ then both sets $I^2_o\setminus \{x\}$ and $\cup_{j = 1}^n\sU_j$ are not connected (cf. \cite{counter_topology}). 
However the set $\cup_{j = 1}^n\sU_j$ is connected, so one has $I^2_o=\cup_{j = 1}^n\sU_j$. Let $\widetilde x_0\in \widetilde \sX$, and let us construct a continuous map $\phi_{\widetilde x_0} : I^2_o \to \widetilde \sX$ such that $\phi_{\widetilde x_0} \circ p\left(\widetilde x_0 \right)= \widetilde x_0$ and $p \circ \phi_{\widetilde x_0} = \Id_{I^2_o }$. If $p\left(\widetilde x_0 \right)\in \sU_j$ and $\widetilde\sU_j$ is a connected component of $p^{-1}\left(\sU_j \right)$ such that  $\widetilde x_0\in \widetilde \sU_j$ then there is a natural homeomorphism $\phi_j:\sU_j\to \widetilde\sU_j$.  There is $\widetilde x_{j+1}\in \widetilde \sU_j$ such $p\left( \widetilde x_{j+1}\right)  \in \sU_{j+1}$. If $\widetilde \sU_{j+1}$ is a connected component of $p^{-1}\circ \left( \sU_{j+1} \right)$ such that $\widetilde x_{j+1} \in \widetilde \sU_{j+1}$ so there is a homeomorphism $\phi_{j+1}:\sU_{j+1}\to \widetilde\sU_{j+1}$. Going on one can construct homeomorphism  $\phi_{m}:\sU_{m}\to \widetilde\sU_{m}$ for any $m = j+2,..., n$. Similarly there is homeomorphism  $\phi_{m}:\sU_{m}\to \widetilde\sU_{m}$ for any $m = j-2,..., 1$. One can proof that
$$
\forall l, k \in 1,..., n \quad \forall x \in \sU_{j}\cap \sU_{j}\quad \phi_{l}\left(x \right). =\phi_{k}\left(x \right) 
$$
From the above equation it follows that there is a continuous map $\phi_{\widetilde x_0}: I^2_o \to \widetilde \sX$ such that
$$
\forall x \in \sU_j \quad \phi_{\widetilde x_0}\left(x \right) = \phi_j\left( x\right).
$$
If $\widetilde{  \mathcal X}_{\widetilde{x}_0} \bydef  \phi_{\widetilde x_0}\left( \sX\right)$ then this set is homeomorphic to $\sX$. From  $\widetilde{  \mathcal X}_{\widetilde{x}_0}= \cup_{j=1}^n\widetilde\sU_{j}$ it follows that the set  $\widetilde{  \mathcal X}_{\widetilde{x}_0}$ is open since it is a union of open sets. If $\widetilde x'_0 \in \widetilde{  \mathcal X}\setminus \widetilde{  \mathcal X}_{\widetilde{x}_0}$ then there are
\begin{itemize}
	\item $j \in 1,..., n $ such that $p\left( \widetilde x\right) \in \sU_j$,
	\item a  connected  component $\widetilde \sU'_j$ of $p^{-1}\left(\widetilde \sU_j \right)$ such that $\widetilde x'_0 \subset \widetilde \sU'_j$.
\end{itemize}
The set  $\widetilde \sU'_j$ is a neighborhood of  $\widetilde x'_0$ such that $\widetilde \sU'_j\cap \widetilde{  \mathcal X}_{\widetilde{x}_0} = \emptyset$, or equivalently $\widetilde \sU'_j \subset \widetilde{  \mathcal X}\setminus \widetilde{  \mathcal X}_{\widetilde{x}_0}$. It turns out that the set $\widetilde{  \mathcal X}\setminus \widetilde{  \mathcal X}_{\widetilde{x}_0}$ is open and $\widetilde{  \mathcal X}_{\widetilde{x}_0}$ is closed. Thus $\widetilde{  \mathcal X}_{\widetilde{x}_0}$ is a quasi-component of $\widetilde{  \mathcal X}$ which is mapped homeomorphically onto $\sX$. This lemma follows from the arbitrary choice of the point $\widetilde x_0$.
\end{proof}

\begin{exercise}\label{top_lex_square_exer}
Let  $I^2_o$ be ordered square (cf. \ref{top_lex_square_empt}). Denote by $x_0 \bydef (0,0), x_1 \bydef (1,1)\in I^2_o$. Let $\sX \bydef I^2_o/\sim_\sX$ where  $\sim_\sX$ is a minimal equivalence relation such that $x_0 \sim_\sX x_1$.
\begin{enumerate}
	\item Prove that $\sX$ is a connected, locally connected (cf. Definition \ref{top_locally_connected_defn}), compact, Hausdorff space, but $\sX$ is not path connected (cf. \cite{counter_topology}).
	\item Similarly to the proof of the Lemma \ref{top_lex_square_lem} prove that the space $\sX$ is weakly semilocally 1-connected (cf. Definition \ref{top_weakly_semi1_defn}).
\end{enumerate}
\end{exercise}
\begin{exercise}\label{top_lex_square1_exer}
Let us use the notation of the Exercise \ref{top_lex_square_exer}. Let $\sX_j$ be a homeomorphic to $I^2_o$ space for any $j\in \Z$. Let $x^j_0 , x^j_1 \in \sX_j$ be points which correspond to $x_0 \bydef (0,0), x_1 \bydef (1,1)\in I^2_o$. Let $\sim_{\widetilde \sX}$ be a minimal equivalence relation on $\bigsqcup_{j\in \Z} \sX_j$ such that
$$
\forall j \in \Z \quad x^j_1 \sim_{\widetilde \sX}  x^{j+1}_0.
$$
Let $\widetilde{\sX}\bydef \bigsqcup_{j\in \Z} \sX_j/ \sim_{\widetilde \sX}$.
	\begin{enumerate}
		\item Prove that $\widetilde\sX$ is a connected, locally connected (cf. Definition \ref{top_locally_connected_defn}), locally compact, Hausdorff space, but $\widetilde\sX$ is not path connected (cf. \cite{counter_topology}).
	\item Similarly to the proof of the Lemma \ref{top_lex_square_lem} prove that the space $\widetilde\sX$ is weakly simply connected (cf. Definition \ref{top_weakly_simply_connected_defn}).
	\end{enumerate}
\end{exercise}

\begin{exercise}\label{top_lex_square2_exer}
	Let us use the notation of the Exercises \ref{top_lex_square_exer}, \ref{top_lex_square1_exer}.   for any $n \in \Z$ define a homeomorphism  $\varphi_n^\sqcup$ of $\sqcup\sX_j$ which homeomorphically  maps $\sX_j$ onto $\sX_{j + n}$. Prove following statements:
	\begin{enumerate}
		\item for any $n\in \Z$ the homeomorphism $\varphi_n^\sqcup$ induces a homeomorphism $\varphi_n$ of $\widetilde \sX$,
		\item a family of $\left\{\varphi_n\right\}$ induces an action $\Z\times \widetilde \sX \to \widetilde{\sX}$ such that $\Z$ is a properly discontinuous group of homeomorphisms (cf. Definition \ref{top_properly_disc_defn}),
		\item there is a natural isomorphism $\sX = \widetilde\sX /\Z$.
	\end{enumerate}
\end{exercise}
\begin{example}
	Let us use the notation of the Exercises \ref{top_lex_square_exer}, \ref{top_lex_square1_exer}. From the Theorem \ref{top_cov_fact_thm} it follows that there is a covering $p: \widetilde\sX\to \sX$. Since $\widetilde\sX$ is simply 1-connected the coveting is universal. If  $n \in \Z$ and $\widetilde x_n \bydef x^n_0 / \sim_{\widetilde \sX}\in \widetilde \sX$ then $p\left(\widetilde x_n \right) = x_0 \bydef p\left(\widetilde x_0 \right)$ for all $n \in \Z$. One can regard $p$ as a pointed  covering $\left(\widetilde\sX, \widetilde x_0 \right)\to  \left(\sX, x_0 \right)$. If $\mathfrak p = \left( \sU_1, ..., \sU_k\right)$ is a semilocally 1-connected closed $x_0$-path on $\sX$ and  $\widetilde{\mathfrak p} = \left( \widetilde\sU_1, ..., \widetilde\sU_k\right)$ is its $p$-lift, then there is a unique $n \in \Z$ such that $\widetilde x_n \in \widetilde\sU_k$. One can proof that  the map
	$$
\mathfrak	p\mapsto n
	$$
	yields a homomorphism $h: \pi^{\text{w}}_1\left(\sX \right)  \to \Z$. Since the covering $p$ is universal a path $\mathfrak p$ is contractible if and only  $\lift_p \mathfrak p$ is contractible, i.e. $h\left(e \right) = 0$ where $e$ is a unity of $\pi^{\text{w}}_1\left(\sX \right)$. So the map $h$ is injective. 	From the Lemma \ref{top_gen_path_lem} it follows that for all $k \in \Z$ there is a path
	$\widetilde {\mathfrak p} = \left(\widetilde \sU_1, ..., \widetilde\sU_k\right)$ such that
	\begin{itemize}
		\item $\widetilde x_0\in \widetilde \sU_1, \quad \widetilde x_n \subset \widetilde \sU_k$, 
		\item the set $\widetilde \sU_j$ is weakly 1-connected.
	\end{itemize}
	If $\mathfrak p = \desc_p \widetilde {\mathfrak p}$ is a $p$ descent  of $\widetilde {\mathfrak p}$ (cf. Definition \ref{top_path_desc_defn}) then  the path $\mathfrak p$ is closed and
	$h\left(\left[\mathfrak p\right]\right)= n$. So the homomorphism $h$ is surjective. In result one has an isomorphism $\pi^{\text{w}}_1\left(\sX \right)\cong\Z$. However the explained in the Remark \ref{top_homotopy_group_rem} fundamental group $\pi_1\left(\sX \right)$ is trivial since $\sX$ is not path connected.
\end{example}

\section{Vector bundles and spaces of sections}
\paragraph*{} This chapter is not exclusively devoted to commutative $C^*$-algebras. It contains a more general theory, where $C^*$-algebras are represented by locally trivial  bundles. Homogeneous $C^*$-algebras and $C^*$-algebras with continuous trace   can be represented as locally trivial bundles. Any commutative $C^*$-algebra is also a vector bundle with one-dimensional fibers. So this theory is used in this chapter and the chapters \ref{ctr_chap},  \ref{foliations_chap}. 
\begin{remark}\label{com_inf_rem}
	The section \ref{top_vb_sub_sub} describes the theory of vector bundles with  finite-dimensional fibers. However many results of the theory can be applied to  vector bundles where any fiber is a more general Banach space. Below the word \textit{vector bundle} means that a fiber may be arbitrary Banach space. 
\end{remark} 
\begin{empt}\label{top_bundle_sec_empt}
	Let $\sX$ be a locally compact Hausdorff space and   $\xi= \left( E, \pi, \mathcal X\right)$ is a vector bundle (cf. Definition \ref{top_vb_defn} and Remark \ref{com_inf_rem}) such that any fiber is a Banach $\C$-space. The space $\Ga\left( E, \pi, \mathcal X\right)$ of sections of $\left( E, \pi, \mathcal X\right)$ is a $\C$-space. For any $a \in \Ga\left( E, \pi, \mathcal X\right)$ the given by
	\be\label{top_normb_eqn}
	\begin{split}
		\text{norm}_a : \sX \to \R,\\
		x \mapsto \left\| a_x\right\| 
	\end{split}
	\ee
	map is continuous. Denote by
	\be\label{top_ub_eqn}
	\begin{split}
		\Ga_c\left( E, \pi, \mathcal X\right)\bydef \left\{\left.a \in \Ga\left( E, \pi, \mathcal X\right)\right| \text{norm}_a \in C_c\left(\sX\right) \right\},\\
		\Ga_0\left( E, \pi, \mathcal X\right)\bydef \left\{\left.a \in \Ga\left( E, \pi, \mathcal X\right)\right| \text{norm}_a \in C_0\left(\sX\right) \right\},\\
		\Ga_b\left( E, \pi, \mathcal X\right)\bydef\left\{\left.a \in \Ga\left( E, \pi, \mathcal X\right)\right| \text{norm}_a \in C_b\left(\sX\right) \right\}.\\
	\end{split}
	\ee
	There is the norm  $\left\|a\right\|\bydef\left\| \text{norm}_a\right\|$ on the above spaces such that both $\Ga_0\left( E, \pi, \mathcal X\right)$ and $\Ga_b\left( E, \pi, \mathcal X\right)$ are Banach spaces and $\Ga_c\left( E, \pi, \mathcal X\right)$ is dense in $\Ga_0\left( E, \pi, \mathcal X\right)$. If $f: \sY \to \sX$ is a continuous map and $\left(E\times_{\sX}\sY, \rho, \mathcal Y\right)$ is the {inverse image} of $\left( E, \pi, \mathcal X\right)$ by $f$ (cf. Definition \ref{vb_inv_img_funct_defn}) then the natural projection $E\times_{\sX}\sY\to E$ yields the natural $\C$-linear maps
	\be\label{top_uba_eqn}
	\begin{split}
		\Ga\left( E, \pi, \mathcal X\right)\hookto \Ga\left(E\times_{\sX}\sY, \rho, \mathcal Y\right),\\
		\Ga_b\left( E, \pi, \mathcal X\right)\hookto \Ga_b\left(E\times_{\sX}\sY, \rho, \mathcal Y\right)\\
	\end{split}
	\ee
	such that the second map is an isometry. Consider a category $\mathfrak{Top}$-$\sX$ such that
	\begin{itemize}
		\item Objects of are continuous maps $\sY \to \sX$.
		\item Morphism  $\sY' \to \sX$ to $\sY'' \to \sX$ is a continuous map $f: \sY' \to \sY''$ such that the following diagram
		\newline
		\begin{tikzpicture}
		\matrix (m) [matrix of math nodes,row sep=3em,column sep=4em,minimum width=2em]
		{
			{\mathcal Y'} & &{\mathcal Y''}\\ 
			& {\mathcal X}\\};
		\path[-stealth]
		(m-1-1) edge node [above] {$f$} (m-1-3)
		(m-1-1) edge node [left]  {} (m-2-2)
		(m-1-3) edge node [right] {} (m-2-2);
		\end{tikzpicture}
		\\
		is commutative.
	\end{itemize}
	For any vector bundle $\left( E, \pi, \mathcal X\right)$ equations \eqref{top_ub_eqn} and \eqref{top_uba_eqn} yield two contravariant functors $\Ga$ and $\Ga_b$ (cf. Definitions \ref{functor_defn} and \ref{functor_contravariant_defn}) form  $\mathfrak{Top}$-$\sX$ to the category of vector spaces given such that 
	\be\label{top_f_eqn}
	\begin{split}
		\Ga\left(\sY \right) \bydef  \Ga\left(E\times_{\sX}\sY, \rho, \mathcal Y\right),\\
		\Ga_b\left(\sY \right) \bydef  \Ga_b\left(E\times_{\sX}\sY, \rho, \mathcal Y\right),\\
		\Ga\left(f \right)\bydef \left(\Ga\left(E\times_{\sX}\sY'', \rho, \mathcal Y''\right)\to \Ga\left(E\times_{\sX}\sY', \rho, \mathcal Y'\right) \right),
		\\
		\Ga_b\left(f \right)\bydef \left(\Ga_b\left(E\times_{\sX}\sY'', \rho, \mathcal Y''\right)\to \Ga_b\left(E\times_{\sX}\sY', \rho, \mathcal Y'\right) \right)  
	\end{split}
	\ee
	and $\Ga_b$ is a functor to the category of Banach spaces and isometries.
\end{empt}

\section{Continuity structures lift and descent}\label{top_lift_desc_sec}
\paragraph*{}
Let $\sX$ be a locally compact, Hausdorff space; and  for each $x$ in $\sX$, let $A_x$ be a (complex) Banach space. Let us consider a continuity structure $\sF$ for $\sX$ {and the} $\left\{A_x\right\}_{x \in \sX}$ (cf. Definition \ref{operator_fields_continuity_defn}). Consider a $\C$-space given by
\be\label{top_c_sec_eqn}
C\left(\sX, \left\{A_x\right\}, \sF \right) \bydef \left\{\left.\left\{a_x\right\}_{x \in \sX}\in \prod_{x \in \sX} A_x\right|\left\{a_x\right\} \text{ is a continuous section} \right\}
\ee 
(cf. Definition \ref{op_cont_fields_defn}).
\begin{remark}\label{top_c_mod_rem}
	From the Lemma \ref{op_cont_module_lem} it follows that $C\left(\sX, \left\{A_x\right\}, \sF \right)$ is $C_0\left(\sX\right)$-module.
\end{remark}

\begin{lemma}\label{top_c_cont_str_lem}
	\begin{enumerate}
		\item[(i)] If $\sF$ is a continuity structure  for $\sX$ {and the} $\left\{A_x\right\}$ (cf. Definition \ref{operator_fields_continuity_defn}), then $C\left(\sX, \left\{A_x\right\}, \sF \right)$ is a continuity structure  for $\sX$ {and the} $\left\{A_x\right\}$.
		\item[(ii)]
		\be
		C\left(\sX, \left\{A_x\right\}, \sF \right)= C\left(\sX, \left\{A_x\right\}, C\left(\sX, \left\{A_x\right\}, \sF \right) \right).
		\ee
	\end{enumerate}
\end{lemma} 
\begin{proof}
	(i) One needs check conditions (a)-(c) of the Definition \ref{operator_fields_continuity_defn}.
	\begin{enumerate}
		\item [(a)] Follows from the Lemma \ref{op_cont_con_lem}.
		\item[(b)] Follows from the inclusion $\sF \subset C\left(\sX, \left\{A_x\right\}, \sF \right)$.
		\item[(c)] If $A_x$ is a $C^*$-algebra for all $x \in \sX$ then from (c) of the Definition \ref{operator_fields_continuity_defn} it turns out that $\sF$ is closed under multiplication and involution. Select $x_0 \in \sX$ and $\eps > 0$. For all $a \in C\left(\sX, \left\{A_x\right\}, \sF \right)$ there is an open neighborhood $\sU$ of $x_0$ and a family $\left\{a'_x\right\}\in \sF$ such that $\left\| a_x- a'_x\right\| < \eps$ for all $x \in \sU$. It follows that $\left\| a^*_x- a'^*_x\right\| < \eps$ for all $x \in \sU$, and taking into account $ a'^*_x \in \sF$ we conclude that $a^* \in C\left(\sX, \left\{A_x\right\}, \sF \right)$. If    $a, b \in C\left(\sX, \left\{A_x\right\}, \sF \right)$. There is $\dl > 0$ such that $\delta^2+\delta\left( \left\| a\right\|+ \left\| b\right\|\right) < \eps$  there are open neighborhoods $\sU', \sU''$ of $x_0$ and families $ \left\{a'_x\right\}\in \sF$  such that $\left\| a_x- a'_x\right\| < \delta$ for all $x \in \sU'$ and $\left\| b_x- b'_x\right\| < \delta$ for all $x \in \sU''$. It turns out that $\left\|a_x b_x- a'_xb'_x\right\| < \eps$ for all $x \in \sU' \cap \sU''$. It follows that the family $\left\{a_xb_x\right\}$ lies in 	$C\left(\sX, \left\{A_x\right\}, \sF \right)$.
	\end{enumerate}
	(ii) From $\sF \subset C\left(\sX, \left\{A_x\right\}, \sF \right)$ it follows that $$C\left(\sX, \left\{A_x\right\}, \sF \right)\subset C\left(\sX, \left\{A_x\right\}, C\left(\sX, \left\{A_x\right\}, \sF \right)\right).$$
	Let $a \in C\left(\sX, \left\{A_x\right\}, C\left(\sX, \left\{A_x\right\}, \sF \right)\right) $, let $x_0 \in \sX$ and $\eps > 0$. There is an open neighborhood  $\sU'$ of $x_0$ and $a' \in C\left(\sX, \left\{A_x\right\}, \sF \right)$ such that $\left\|a_x  - a'_x\right\| < \frac{\eps}{2}$ for all $x \in \sU'$. Otherwise there is an open neighborhood  $\sU''$ of $x_0$ and $a'' \in \sF$ such that $\left\|a'_x  - a''_x\right\| < \frac{\eps}{2}$ for every $x \in \sU''$. It follows that $\left\|a_x  - a''_x\right\| < \eps$ for every $x \in \sU'\cap\sU''$, hence one has $a \in C\left(\sX, \left\{A_x\right\}, \sF \right)$. In result we conclude that
	$$ C\left(\sX, \left\{A_x\right\}, C\left(\sX, \left\{A_x\right\}, \sF \right)\right)\subset C\left(\sX, \left\{A_x\right\}, \sF \right).$$
\end{proof} 
\begin{remark}
	The Lemma \ref{top_c_cont_str_lem} is a generalization of the Lemma \ref{top_full_oaf_lem}.
\end{remark}
Any vector bundle $\left( E, \pi, \mathcal X\right)$ (cf. Definition \ref{top_vb_defn} and Remark \ref{com_inf_rem}) yields a family of fibers $\left\{E_x\right\}_{x \in \sX}$.

For any $a \in C\left(\sX, \left\{A_x\right\}, \sF \right)$ denote by 
\bea
\label{top_norm_a_eqn}
\text{norm}_a:  \sX \to \R, \quad x \mapsto \left\|a_x \right\|.
\eea

Let us define $\left\|\cdot\right\|: C\left(\sX, \left\{A_x\right\},\sF\right)\to \left[0, \infty\right)$ given by
\be\label{top_norm_eqn}
\left\|a\right\|\bydef \left\|\text{norm}_a\right\|.
\ee

\begin{definition}\label{top_s_bounded_defn}
	The space $C_b\left(\sX, \left\{A_x\right\},\sF\right)$ given by
	\be\label{top_s_bounded_eqn}
	C_b\left(\sX, \left\{A_x\right\},\sF\right)\bydef \left\{\left.a \in C\left(\sX, \left\{A_x\right\},\sF\right)\right|\left\|a\right\| < \infty\right\}.
	\ee
	is said to be the \textit{space of bounded continuous sections}.
\end{definition}

\begin{lemma}\label{top_bounded_norm_lem}
	One has
	\be\nonumber
	a \in C_b\left(\sX, \left\{A_x\right\},\sF\right) \Leftrightarrow  \mathrm{norm}_a \in C_b\left(\sX \right). 
	\ee
\end{lemma}
\begin{proof}
	From the Lemma 	\ref{op_cont_con_lem} it follows that the map $x \mapsto \left\|a_x \right\|$ is continuous, from the equation \eqref{top_s_bounded_eqn} it follows that $x \mapsto \left\|a_x \right\|$ is bounded.
\end{proof}

\begin{definition}\label{top_cc_supp_defn}
	If $a\in  C\left(\sX, \left\{A_x\right\}, \sF \right)$ is presented by the family
	$\left\{a_x \in A_x\right\}$ then the closure of $\left\{\left. x \in \sX~\right|\left\| a_x\right\| > 0\right\}$ is said to be the \textit{support} of $a$. The support of $a$ coincides with the support of the $\mathrm{norm}_a$,  hence we write
	\be\label{top_support_eqn}
	\supp a  \stackrel{\text{def}}{=} \supp \mathrm{norm}_a.
	\ee
\end{definition}
\begin{definition}\label{top_cc_c0_defn}
	The $C_b\left({\sX}\right) $-module 
	\be
	C_c\left(\sX, \left\{A_x\right\}, \sF \right) = \left\{\left.a \in C_b\left(\sX, \left\{A_x\right\}, \sF \right)~\right| \supp a \text{ is compact } \right\}
	\ee
	is said to be the \textit{compactly supported submodule}. The norm completion of $C_c\left(\sX, \left\{A_x\right\}, \sF \right)$ is said to be \textit{converging to zero submodule} and we denote it by $C_0\left(\sX, \left\{A_x\right\}, \sF \right)$. We also use the following notation
	\bea
	\label{top_cc_eqn}
	C_c\left(A \right) \bydef C_c\left(\sX, \left\{A_x\right\}, A \right),\\
	\label{top_c0_eqn}
	C_0\left(A \right) \bydef C_0\left(\sX, \left\{A_x\right\}, A \right),\\
	\label{top_cb_eqn}
	C_b\left(A \right) \bydef C_b\left(\sX, \left\{A_x\right\}, A \right),\\
	\label{top_c_eqn}
	C\left(A \right) \bydef C\left(\sX, \left\{A_x\right\}, A \right),\\
	\label{top_res_eqn}
	\left.C_0\left(\sX, \left\{A_x\right\}, \sF  \right) \right|_{\sU}\stackrel{\text{def}}{=}\left\{\left.a \in C_0\left(\sX, \left\{A_x\right\}, \sF  \right)~\right|~ \supp a \subset \sU \right\}.
	\eea
	
\end{definition}

\begin{remark}\label{top_cc0_mod_rem}
	Both 	$C_c\left(\sX, \left\{A_x\right\}, A \right)$ and $ C\left(\sX, \left\{A_x\right\}, A \right)$ are $C\left(\sX\right)$-modules (cf. Remark \ref{top_c_mod_rem}). All spaces  $C_c\left(\sX, \left\{A_x\right\}, A \right)$, $C_0\left(\sX, \left\{A_x\right\}, A \right)$, $C_b\left(\sX, \left\{A_x\right\}, A \right)$ and $C\left(\sX, \left\{A_x\right\}, A \right)$ are $C_b\left(\sX\right)$-modules
\end{remark}
\begin{remark}
	The equation \eqref{top_c0_eqn} complies with the explained in \ref{top_cs_not_empt} notation.
		\end{remark}

\begin{remark}
	If $C_0\left(\sX, \left\{A_x\right\}, \sF  \right)$ is a $C^*$-algebra with Hausdorff spectrum then the notation \eqref{top_res_eqn} complies  with \eqref{open_ideal_eqn} one.
\end{remark}

\begin{lemma}\label{top_norm_c0cc_lem}
	Following conditions hold
	\bea
	\label{top_ccn_eqn}
C_c\left(\sX, \left\{A_x\right\}, A \right)= \left\{\left.a \in 	C_b\left(A \right)~\right| ~\mathrm{norm}_a \in 	C_c\left(\sX\right)\right\},\\
	\label{top_cc0_eqn}
C_0\left(\sX, \left\{A_x\right\}, A \right) = \left\{\left.a \in 	C_b\left(A \right)~\right| ~\mathrm{norm}_a \in 	C_0\left(\sX\right)\right\}
	\eea
	(cf. equations \eqref{top_norm_a_eqn}, \eqref{top_c0_eqn},  \eqref{top_cb_eqn}).
\end{lemma}
\begin{proof}
	The equation \eqref{top_ccn_eqn} is evident. If $a \in C_0\left( A\right)$ then there is a $C^*$-norm convergent net $\left\{a_\a \in 	C_c\left(A \right)\right\}_{\a \in \mathscr A}$ such that $a = \lim a_\a$. For any $\a, \bt \in \mathscr A$ one has
	\be\label{top_ab_norm_eqn}
	\left\|\text{norm}_{a_\a} - \text{norm}_{a_\bt} \right\|\le \left\|{a_\a} - {a_\bt} \right\|.
	\ee
	From \eqref{top_ab_norm_eqn} it turns out that the net $\left\{\text{norm}_{a_\a}\right\}$ is uniformly convergent. From $\text{norm}_{a_\a} \in C_c\left( \sX\right)$ and the Definition \ref{top_cc_c0_defn} it turns out $\lim_\a \text{norm}_{a_\a} \in C_0\left(\sX \right)$. Conversely let $a \in C_0\left( A\right)_+$ be a positive element such that $\mathrm{norm}_a \in 	C_0\left(\sX\right)$. If $f_\eps$ is given by \eqref{f_eps_eqn} then one has
	\bean
\left\|f_{\eps/2}\left(a \right)  - a\right\| < \eps,\\
\mathrm{norm}_{f_{\eps/2}\left(a \right)}= f_{\eps/2}\left(\mathrm{norm}_a \right),	
	\eean 
	so $a = \lim_{\eps\to 0}  f_{\eps/2}\left(a \right)$. On the other hand
from  the Lemma \ref{pedersen_eps_lem} it follows that $f_{\eps/2}\left(\mathrm{norm}_a \right)\in K\left(C_0\left(\sX \right)  \right)= C_c\left(\sX \right)$.	
\end{proof}
\begin{corollary}\label{top_c0_cont_str_cor}
	One has
	\be\label{top_c0_cont_strc_eqn}
	\sX \mathrm{~is~compact~} \Rightarrow	C_c\left(A \right)= 	C_0\left(A \right)= 	C_b\left(A \right)= 	C\left(A \right).
	\ee	
\end{corollary}

\begin{lemma}\label{top_bundle_cs_ex_lem}
	If $\left( E, \pi, \mathcal X\right)$ is a vector bundle (cf. Definition \ref{top_vb_defn}) then following conditions hold.
	\begin{enumerate}
		\item [(i)] The space 	$\Ga\left( E, \pi, \mathcal X\right)$ of continuous sections (cf. \ref{top_bundle_sec_empt}) 
		is a continuity structure  for $\sX$ {and the} $\left\{E_x\right\}_{x \in \sX}$ (cf. Definition \ref{operator_fields_continuity_defn}).
		\item[(ii)] There are isomorphisms
		\bea\label{top_cg_eqn}
		C\left(\sX , \left\{E_x\right\}, \Ga\left( E, \pi, \mathcal X\right) \right) \cong \Ga\left( E, \pi, \mathcal X\right),\\
		\label{top_cgc_eqn}
		C_c\left(\sX , \left\{E_x\right\}, \Ga\left( E, \pi, \mathcal X\right) \right) \cong \Ga_c\left( E, \pi, \mathcal X\right),\\
		\label{top_cg0_eqn}
		C_0\left(\sX , \left\{E_x\right\}, \Ga\left( E, \pi, \mathcal X\right) \right) \cong \Ga_0\left( E, \pi, \mathcal X\right),\\
		\label{top_cgb_eqn}
		C_b\left(\sX , \left\{E_x\right\}, \Ga\left( E, \pi, \mathcal X\right) \right) \cong \Ga_b\left( E, \pi, \mathcal X\right)
		\eea
		of are isomorphisms of $C_0\left(\sX \right)$-modules (cf. equation \eqref{top_f_eqn}).
	\end{enumerate}
	
\end{lemma}
\begin{proof} (i)
	One needs check conditions (a)-(c) of the Definition \ref{operator_fields_continuity_defn}.\\
	(a)	
	Let $x_0 \in \sX$ be any point, and select an open neighborhood $\sU$ of $x_0$ such that $\left.\left( E, \pi, \mathcal X\right)\right|_{\sU}$ is trivial (cf.  \ref{trivial_vb_empt} and Definition \ref{top_vb_defn}), i.e. $\left.E\right|_{\sU}\cong \sU \times E_{x_0}$. Any section $s\in\Ga\left( E, \pi, \mathcal X\right)$ corresponds to a continuous  map $s: \sX \to E$ such that $\pi \circ s = \Id_\sX$, so the map
	\bean
	s|_\sU : \sU \to \sU\times E_{x_0},\\
	x \mapsto\left(x, u\left(x \right) \right)
	\eean 
	is continuous at $x_0$. It follows that both maps
	\bean
	x \mapsto u\left(x \right),\\
	x \mapsto \left\| u\left(x \right)\right\| 
	\eean 
	are continuous. Since a point $x_0$ is arbitrary the map $ x \mapsto \left\| s\left(x \right)\right\|$ is continuous on $\sX$. \\
	(b) If $x_0 \in \sX$ and $e_{x_0} \in  E_{x_0}$ then from the Definition \ref{top_vb_defn} it follows that $\left.E\right|_{\sU}\cong \sU \times E_{x_0}$. It turns out that there is a section
	\bean
	s|_\sU : \sU \to \sU\times E_{x_0},\\
	x \mapsto\left(x, e_{x_0} \right),
	\eean 
	i.e. there is a section $s: \sX\to E$ such that $s\left( x_0\right) = e_{x_0}$.\\
	(c) A proof of this property we leave to the reader.\\
	(ii)  A proof  we leave to the reader.
\end{proof}

\begin{lemma}\label{top_c0_cont_str_lem}

	\begin{enumerate}
		\item[(i)] If $\sF$ is a continuity structure  for $\sX$ {and the} $\left\{A_x\right\}$ (cf. Definition \ref{operator_fields_continuity_defn}), then   $C_c\left(\sX, \left\{A_x\right\}, \sF \right)$, $~~C_0\left(\sX, \left\{A_x\right\}, \sF \right)$ and $C_b\left(\sX, \left\{A_x\right\}, \sF \right)$ are  continuity structures  for $\sX$ {and the} $\left\{A_x\right\}$.
		\item[(ii)]
		\be\label{top_c0_cont_str_eqn}
		\begin{split}
			C\left(\sX, \left\{A_x\right\}, \sF \right)= C\left(\sX, \left\{A_x\right\}, C_c\left(\sX, \left\{A_x\right\}, \sF \right) \right)=\\
			C\left(\sX, \left\{A_x\right\}, C_0\left(\sX, \left\{A_x\right\}, \sF \right) \right)= C\left(\sX, \left\{A_x\right\}, C_b\left(\sX, \left\{A_x\right\}, \sF \right) \right).
		\end{split}
		\ee
	\end{enumerate}
\end{lemma}

\begin{proof}
	(i) One needs check conditions (a)-(c) of the Definition \ref{operator_fields_continuity_defn}.
	\begin{enumerate}
		\item [(a)] Follows from the Lemma \ref{op_cont_con_lem}.
		\item[(b)] 	Let $x\in \sX$, and let $f_{x}\in C_c\left( \sX\right)$ be a $x$-stump (cf. Definition \ref{top_stump_defn}).
		For any $a \in \sF$ one has $f_{x} a \in  C_c\left(\sX, \left\{A_x\right\}, \sF \right)$, and from $f_{x}\left( x\right) =1$ it turns out
		\bean
		A\bydef	\left\{\left. a_x \in A_x\right|\exists a' \in \sF \quad a_x = a'_x  \right\}\subset \\ \subset A_c\bydef	\left\{\left. a_x \in A_x\right|\exists a' \in C_c\left(\sX, \left\{A_x\right\}, \sF \right) \quad a_x = a'_x  \right\}.
		\eean
	 If
		\bean
		A_0\bydef	\left\{\left. a_x \in A_x\right|\exists a' \in C_0\left(\sX, \left\{A_x\right\}, \sF \right) \quad a_x = a'_x   \right\},\\
		A_b\bydef	\left\{\left. a_x \in A_x\right|\exists a' \in C_b\left(\sX, \left\{A_x\right\}, \sF \right) \quad a_x = a'_x   \right\}
		\eean 
		then from
		$$
		A \subset A_c \subset A_0 \subset A_b \subset A_x
		$$
		and taking into account that $A$ is dense in $A_x$,
		it turns out that $A_c$, $A_0$ and $A_b$ are dense in $A_x$.
		\item[(c)] The proof is similar to the proof of (i) (c)  of the Lemma \ref{top_c_cont_str_lem}.
	\end{enumerate}
	(ii) From the inclusions
	$$
	C_c\left(\sX, \left\{A_x\right\}, \sF \right) \subset C_0\left(\sX, \left\{A_x\right\}, \sF \right) \subset C_b\left(\sX, \left\{A_x\right\}, \sF \right) \subset C\left(\sX, \left\{A_x\right\}, \sF \right).
	$$
	it follows that 
	\be\label{top_c0_cont_str_j_lem}
	\begin{split}
		C\left(\sX, \left\{A_x\right\}, C_c\left(\sX, \left\{A_x\right\}, \sF \right) \right)\subset 	C\left(\sX, \left\{A_x\right\}, C_0\left(\sX, \left\{A_x\right\}, \sF \right)\right)\subset \\ \subset  
		C\left(\sX, \left\{A_x\right\}, C_b\left(\sX, \left\{A_x\right\}, \sF \right) \right)\subset 	C\left(\sX, \left\{A_x\right\}, \sF \right).
	\end{split}
	\ee 
	Let $a \in C\left(\sX, \left\{A_x\right\}, C\left(\sX, \left\{A_x\right\}, \sF \right)\right) $, let $x_0 \in \sX$ and $\eps > 0$. There is an open neighborhood  $\sU'$ of $x_0$ and $a' \in C\left(\sX, \left\{A_x\right\}, \sF \right)$ such that $\left\|a_x  - a'_x\right\| < \eps$ for all $x \in \sU'$. Otherwise if $f_{x_0}\in C_c\left(\sX\right)$ is an $x_0$-stump  (cf. Definition \ref{top_stump_defn}) then there is an open neighborhood $\sU''$ of $x_0$ such that $f_{x_0} \left( \sU''\right)= \{1\}$. If  and $a''\bydef f_{x_0}a' \in  C_c\left(\sX, \left\{A_x\right\}, \sF \right)$ then  $\left\|a_x - a''_x\right\| < \eps$ for all $x \in \sU' \cap\sU''$. It follows that $a \in 	C\left(\sX, \left\{A_x\right\}, C_c\left(\sX, \left\{A_x\right\}, \sF \right) \right)$, so one has
	$$
	C\left(\sX, \left\{A_x\right\}, \sF \right) \subset C\left(\sX, \left\{A_x\right\}, C_c\left(\sX, \left\{A_x\right\}, \sF \right) \right),
	$$
	and taking into account \eqref{top_c0_cont_str_j_lem} one  obtains \eqref{top_c0_cont_str_eqn}.
\end{proof}

\begin{empt}
	If	$\sF$ is  a continuity structure  for $\sX$ {and the} $\left\{A_x\right\}$ then from the Lemma \ref{op_cont_module_lem} it follows that $C_0\left(\sX \right)\sF \subset C_0\left(\sX, \left\{A_x\right\}, \sF \right)$.
\end{empt}

\begin{lemma}\label{top_cs_nc_lem}
	If $\sF$ is a continuity structure for $\sX$ {and the} $\left\{A_x\right\}$, such that $\sF \subset C_0\left( \sX, \left\{A_x\right\}, \sF\right)$
then a space $C_0\left(\sX \right)\sF$ is dense in 	$C_0\left(\sX, \left\{A_x\right\}, \sF \right)$ with respect to the norm topology.
\end{lemma}
\begin{proof}
	Let $a \in C_0\left(\sX, \left\{A_x\right\}, \sF \right)$ be represented by the family $\left\{a_x \in A_x\right\}$, and let $\eps > 0$. For any $x' \in \sX$ we select an open neighborhood $\sU_{x'} \subset \sX$ and $a^{x'} \in \sF$  represented by the family $\left\{a'_x \in A_x\right\}$ such that $\left\| a_x  - a^{x'}_x\right\|< \eps$ for all $x \in \sU_{x'}$.  The  set $\sU \bydef \left\{\left. x \in \sX~\right| \left\|  a_x   \right\| \le \eps\right\}$ is compact so there is a  dominated by the family $\left\{\sU_{x'}\right\}_{x'\in \sX}$ finite covering sum 
	$
	\sum_{j=1}^n f_{x_j}
	$
	for $\sU$ (cf. Definition \ref{top_covering_sum_defn}), such that  $\sum_{j=1}^n f_{x_j}\left(x \right) = 1$ for each $x \in \sU$. From  $\left\| a_x  - a^{x'}_x\right\|< \eps$ for all $x \in \sU_{x'}$ it turns out  that
	$$
	\left\| \sum_{j=1}^n f_{x_j}\left( x\right)  a^{x_j}_x - a_x\right\|<  \eps \quad \forall x \in \sU.
	$$
	and taking into account 
	\bean
	x \in \sX \setminus \sU\quad  \Rightarrow \quad 0\le \sum_{j=1}^n f_{x_j}\left( x\right)\le 1 \quad \text{ AND } \\\
	\text{ AND }\quad  \left\| a_x  \right\|< \eps \quad \text{ AND } \quad  \left\|a_x - a^{x'}_x \right\|< \eps  \quad  \Rightarrow\\
	\Rightarrow		\left\| \sum_{j=1}^n f_{x_j}\left( x\right)  a^{x_j}_x - a_x\right\|< \eps
	\eean
	one has
	$$
	\left\| \sum_{j=1}^n f_{x_j}\left( x\right)  a^{x_j}_x - a_x\right\|< \eps \quad \forall x \in \sX \quad \Leftrightarrow \quad \left\|  \sum_{j=1}^n f_{x_j}  a^{x_j} - a\right\|< \eps.
	$$	
	From $\sum_{j=1}^n f_{x_j}  a^{x_j}\in C_0\left(\sX \right)\sF$ it turns out that the space $C_0\left(\sX \right)\sF$ is dense in 	$C_0\left(\sX, \left\{A_x\right\}, \sF \right).$	
\end{proof}
\begin{remark}
The Lemma \ref{top_cs_nc_lem} is a generalization of \ref{top_full_oaf_lem} one.
\end{remark}
\begin{corollary}\label{top_cl_cc_iso_cor}
	Let $\sF$ be a continuity structure for $\sX$ {and the} $\left\{A_x\right\}$, such that $\sF \subset C_0\left( \sX, \left\{A_x\right\}, \sF\right)$.
	If  $C_0\left(\sX \right)\sF\subset \sF$  and $\sF$ is norm closed then the the natural inclusion $\sF \hookto C_0\left(\sX, \left\{A_x\right\}, \sF \right)$ is the $C_0\left(\sX \right)$-isomorphism
	$
	\sF \cong C_0\left(\sX, \left\{A_x\right\}, \sF \right)
	$. 
\end{corollary}
\begin{remark}\label{top_cl_cc_iso_rem}
	From the Lemma \ref{top_cs_nc_lem} it follows that $C_0\left(\sX, \left\{A_x\right\}, \sF \right)$ is a minimal norm closed subspace of $C_b\left(\sX, \left\{A_x\right\}, \sF \right)$ which contains all products $fa$ where $a \in \sF$ and $f \in C_0\left(\sX\right)$.
\end{remark}

 
\begin{definition}\label{top_total_defn}
	For any family of Banach spaces $\left\{A_x\right\}_{x \in \sX}$ a union $\widehat{\left\{A_x\right\}}\bydef \cup_{x \in \sX} A_x$ is said to be a \textit{total set} of $\left\{A_x\right\}_{x \in \sX}$. The natural map 
	\bean
	\widehat{\left\{A_x\right\}} \to\sX,\\
	a_x \mapsto x
	\eean
	is said to be the \textit{total projection}.
\end{definition}
\begin{empt} Let us consider a continuity structure $\sF$ for $\sX$ {and the} $\left\{A_x\right\}_{x \in \sX}$ (cf. Definition \ref{operator_fields_continuity_defn}).
	Let $\sU \subset \sX$ be an open subset and let $s \in C\left(\sX, \left\{A_x\right\}, \sF \right)$. For any $\eps > 0$ denote by
	\be\label{top_cs_neigh_eqn}
	\mathfrak{O}\left(\sU, s, \eps \right) \bydef \left\{\left.a_x \in \widehat{\left\{A_x\right\}}\right| x\in \sU  \quad\left\|a_x - s_x \right\|< \eps \right\}
	\ee
\end{empt}
\begin{lemma}
	A collection of given by the equation \eqref{top_cs_neigh_eqn} subsets satisfies to (a) and (b) of the Definition \ref{top_base_defn}.
\end{lemma}
\begin{proof}
	Condition (a) is evident, let us proof (b). If $a^0_{x_0} \in \widehat{\left\{A_x\right\}}$, and  
	$$\mathfrak{O}\left(\sU', s', \eps' \right), \mathfrak{O}\left(\sU'', s'', \eps'' \right)\subset \widehat{\left\{A_x\right\}}$$ 
	are given by \eqref{top_cs_neigh_eqn} sets such that $a^0_{x_0}\in \mathfrak{O}\left(\sU', s', \eps' \right)\cap\mathfrak{O}\left(\sU'', s'', \eps'' \right)$ then one has
	\bean
		\delta ' \bydef\eps'-\left\|s'_{x_0}-a^0_{x_0} \right\|>0,\\
		\delta '' \bydef\eps''-\left\|s''_{x_0}-a^0_{x_0} \right\|> 0.
	\eean 
If $\eps< \min\left(\delta',\delta''\right)$ then there is $s \in C\left(\sX, \left\{A_x\right\}, \sF \right)$ such that  $\left\|a^0_{x_0}-s_{x_0} \right\|<\eps/3$, and one has
	\bean
	\left\|s'_{x_0}-s_{x_0} \right\|<\eps'- \frac{\eps}{3},\\
	\left\|s''_{x_0}-s_{x_0} \right\|<\eps''-\frac{\eps}{3}.
	\eean 
From the above equations it follows that
\bean
\forall b_{x_0}\in A_{x_0}\quad \left\|b_{x_0}-s_{x_0} \right\|< \frac{\eps}{3} \quad\Rightarrow\quad \left\|b_{x_0}-s'_{x_0} \right\|<\eps'-\frac{\eps}{3},\quad \left\|b_{x_0}-s''_{x_0} \right\|<\eps''-\frac{\eps}{3},
\eean 
it turns out that there are open neighborhoods 	$\sV', \quad \sV''$ of $x_0$ which satisfy to following conditions
\be\label{top_o_eqn}
\begin{split}
x \in \sV' \quad\Rightarrow\quad\left(\left\|b_{x}-s_{x} \right\|< \frac{\eps}{3}\quad\Rightarrow\quad \left\|b_{x}-s'_{x} \right\|<\eps'\right),\\
x \in \sV'' \quad\Rightarrow\quad\left(\left\|b_{x}-s_{x} \right\|< \frac{\eps}{3}\quad\Rightarrow\quad \left\|b_{x}-s''_{x} \right\|<\eps''\right). 
\end{split}
\ee
From the equation \eqref{top_o_eqn} it turns out that
$$
\mathfrak{O}\left(\sU\bydef\sV'\cap\sV''\cap\sU'\cap\sU'', s, \eps/3 \right)\subset \mathfrak{O}\left(\sU', s', \eps' \right)\cap \mathfrak{O}\left(\sU'', s'', \eps'' \right).
$$
\end{proof}
\begin{definition}\label{top_total_sp_defn}
	A topological $\mathscr T\left(\sX, \left\{A_x\right\}, \sF \right)$ space such that:
	\begin{itemize}
		\item the space  $\mathscr T\left(\sX, \left\{A_x\right\}, \sF \right)$ which coincides with the total set $ \widehat{\left\{A_x\right\}}\bydef \cup_{x \in \sX} A_x$ as a set;
		\item the topology of $\mathscr T\left(\sX, \left\{A_x\right\}, \sF \right)$ is generated by given by \eqref{top_cs_neigh_eqn} sets (cf. Definition \ref{top_base_defn}).
	\end{itemize}
is said to be the \textit{total space for} $\left(\sX, \left\{A_x\right\}, \sF \right)$.
\end{definition}
\begin{remark}
	The total projection $\mathscr T\left(\sX, \left\{A_x\right\}, \sF \right) \to\sX$ (cf. Definition \ref{top_total_defn}) is a continuous map.
\end{remark}
\begin{lemma}\label{top_js_lem} Following conditions hold.
	\begin{enumerate}
		\item [(i)] If $p:\mathscr T\left(\sX, \left\{A_x\right\}, \sF \right) \to\sX$  is the	total projection  (cf. Definition \ref{top_total_defn}) and $j : \sX\to \mathscr T\left(\sX, \left\{A_x\right\}, \sF \right)$ is  a continuous map  such that $p\circ j = \Id_\sX$ then the family $\left\{j\left(x\right)\in A_x\right\}$ is  continuous with respect to $\sF$ (cf. Definition \ref{op_cont_fields_defn}).
		\item[(ii)] 	 Conversely any continuous with respect to $\sF$ family $\left\{a_x\right\}$  corresponds to a continuous map 
		\be\label{top_js_eqn}
		\begin{split}
			j : \sX\to \mathscr T\left(\sX, \left\{A_x\right\}, \sF \right),\\
			x \mapsto a_x.
		\end{split}
		\ee
	\end{enumerate}
	
\end{lemma}

\begin{proof}
	(i) If $\eps > 0$ and $x_0\in\sX$ then from (b) of the Definition \ref{operator_fields_continuity_defn} it follows that there is $a \in \sF$ such that $\left\| a_{x_0}- j\left( x_0\right) \right\| < \eps/2$. If $\sU'$ is an open neighborhood of $x_0$ then 
	from $j\left( x_0\right) \in \mathfrak{O}\left(\sU', a, \eps \right)$, and since $j$ is continuous it follows that $\sU\bydef j^{-1} \left(  \mathfrak{O}\left(\sU', a, \eps \right)\right)$ is open. So for all $x\in \sU$ one has $\left\| a_{x}- j\left( x\right) \right\|< \eps$, i.e. the family $\left\{j\left(x\right)\right\}$ is continuous (cf. Definition \ref{op_cont_fields_defn}).\\
	(ii)  Suppose $\eps > 0$, $x_0\in\sX$ and $j$ is given by \eqref{top_js_eqn}. From the condition (b) of the Definition \ref{operator_fields_continuity_defn} it follows that there is a $b \in \sF$   such that $\left\| a_{x_0}- b_{x_0} \right\| < \eps$, i.e. $j\left( x_0\right) \in \mathfrak{O}\left(\sU', b, \eps \right)$. From the Lemma \ref{op_cont_con_lem} it follows that there is an open neighborhood $\sU$ of $x_0$ such that $\left\| a_{x}- b_{x} \right\|<\eps$ for all $x \in \sU$. So one has $\sU\cap \sU'\subset j^{-1}\left( \mathfrak{O}\left(\sU', b, \eps \right)\right)$, i.e. $j$ is continuous at $x_0$. Similarly one can prove that $j$ is continuous at any point of $\sX$.
\end{proof}
\begin{definition}\label{ctr_restr_defn}
	For any subset $\sV\subset \sX$ 
	denote by  $C\left(\sV, \left\{A_x\right\}, \sF \right)$ the $\C$-space of continuous maps $j:\sV\to \mathscr T\left(\sX, \left\{A_x\right\}, \sF \right)$ such that if $p$ is the total projection then $p\circ j= \Id_\sV$. The space  $C\left(\sV, \left\{A_x\right\}, \sF \right)$ is said to be the $\sV$-\textit{restriction} of $C\left(\sX, \left\{A_x\right\}, \sF \right)$. The natural $\C$-linear map $C\left(\sX, \left\{A_x\right\}, \sF \right)\to C\left(\sV, \left\{A_x\right\}, \sF \right)$ is said to be the \textit{restriction map}.
\end{definition}
\begin{empt}
	Similarly to the equation \eqref{top_norm_a_eqn} any $a\in C\left(\sV, \left\{A_x\right\}, \sF \right)$ yields a continuous map 
	\be\label{top_rn_eqn}
	\text{norm}_a:  \sV \to \R, \quad x \mapsto \left\|a_x \right\|.
	\ee
	Similarly to \eqref{top_norm_eqn} define $\left\|\cdot\right\|: C\left(\sV, \left\{A_x\right\},\sF\right)\to \left[0, \infty\right)$ given by
	\be\label{top_r_n_eqn}
	\left\|a\right\|\bydef \left\|\text{norm}_a\right\|.
	\ee
	Denote by
	\be\label{top_r_b0_eqn}
	\begin{split}
		C_b\left(\sV, \left\{A_x\right\},\sF\right) \bydef \left\{\left.a \in C\left(\sV, \left\{A_x\right\},\sF\right)\right|\text{norm}_a \in C_b\left( \sV\right) \right\},\\
		C_0\left(\sV, \left\{A_x\right\},\sF\right) \bydef \left\{\left.a \in C\left(\sV, \left\{A_x\right\},\sF\right)\right|\text{norm}_a \in C_0\left( \sV\right) \right\}.
	\end{split}
	\ee
\end{empt}
\begin{remark}\label{top_cv_nc_rem}
	Both given by the equations \eqref{top_r_b0_eqn} spaces are  closed with respect to the given by \eqref{top_r_n_eqn} norm.
\end{remark}
\begin{remark}
	If $\sV$ is compact then one has $	C\left(\sV, \left\{A_x\right\},\sF\right)=	C_b\left(\sV, \left\{A_x\right\},\sF\right)=	C_0\left(\sV, \left\{A_x\right\},\sF\right)$.
\end{remark}

\begin{lemma}\label{top_ax_k_lem}
	Let $\sX$ be a locally compact,  Hausdorff space; and for each $x$ in $\sX$, let $A_x$ be a (complex) Banach space. Let us consider a continuity structure $\sF$ for $\sX$ {and the} $\left\{A_x\right\}_{x \in \sX}$ (cf. Definition \ref{operator_fields_continuity_defn}). If $\sY \subset \sX$ is compact then the restriction map $C_0\left(\sX, \left\{A_x\right\}, \sF \right) \to C\left(\sY, \left\{A_x\right\}, \sF \right)$ (cf. Definition \ref{ctr_restr_defn}) is surjective.
\end{lemma}
\begin{proof}
	One needs prove that for any $a\in C\left(\sY, \left\{A_x\right\}, \sF \right)$ there is $b\in C_0\left(\sX, \left\{A_x\right\}, \sF \right)$ such that
\be\label{top_ak_eqn}
a = \left. b\right|_\sY.
\ee
Our proof has two parts:
\begin{enumerate}
	\item [(i)] For all $a\in C\left(\sY, \left\{A_x\right\}, \sF \right)$ and $\eps > 0$ there is $b\in C_0\left(\sX, \left\{A_x\right\}, \sF \right)$ such that
	\be\label{top_abk_eqn}
	\begin{split}
		\left\|a - \left.b\right|_\sY \right\| < \eps,\\
		\left\| b\right\| <  \left\|a \right\|+  \eps.
	\end{split}
	\ee
	\item[(ii)] Looking for $b$ which satisfies to \eqref{top_ak_eqn}.
\end{enumerate}
	(i) For any $x \in \sY$ we select $s^x \in C\left(\sX, \left\{A_x\right\}, \sF \right)$ such that $\left\|s^x_x - a_x\right\|< {\eps}/{2}$. There is an open neighborhood $\sU'_x\subset\sX$ of $x$ with compact closure such that \be\label{ctr_up_eqn}
\left\|s^x_{x'} \right\| < \eps+ \left\|a \right\|\quad \forall x' \in \sU'_x.
\ee
On the other hand there is an open  neighborhood $\sU''_x\subset\sX$ of $x$ such that
\be\label{ctr_b_eqn}
\left\|s^x_y- a_y \right\| < \eps\quad \forall y \in \sY\cap\sU''_x.
\ee
If $\sU_x \bydef \sU'_x \cap\sU''_x$ then $\sY \subset \bigcup_{x \in \sY} \sU_x$.
Let $\sum_{j = 1}^n f_j$ be a covering sum for $\sY$ dominated by the family 
$\left\{\sU_x\right\}_{x\in\sX}$ (cf. Definition \ref{top_covering_sum_defn}), such that $\supp f_j \subset \sU_{x_j}$ and $\sU_{x_j}\in \left\{\sU_x\right\}_{x\in\sX}$. If $b \bydef \sum_{j=1}^n f_js^{x_j}\in C_0\left(\sX, \left\{A_x\right\}, \sF \right)$ then from the equations \eqref{ctr_up_eqn}, \eqref{ctr_b_eqn} and the Definition \ref{top_covering_sum_defn} it turns out that $b$ satisfies to conditions \eqref{top_abk_eqn}.\\
(ii)
Let $\eps > 0$. Using induction we will construct sequences $\left\{a_n\right\}_{n\in \N} \subset  C\left(\sY, \left\{A_x\right\}, \sF \right)$ and $\left\{b_n\right\}_{n\in \N} \subset  C_0\left(\sX, \left\{A_x\right\}, \sF \right)$ such that 
\be\label{top_abn_eqn}
\begin{split}
	\left\|a_n\right\|< \frac{\eps}{2^{n-1}} ,\\ 
	n > 1 \quad\Rightarrow \quad	\left\|b_n\right\|< \frac{3\eps}{2^{n-1}},\\
\left\|a - \left.\left( \sum_{n = 1}^m b_n\right) \right|_\sV \right\| < \frac{\eps}{2^{m-1}}.
\end{split}
\ee
From \eqref{top_abk_eqn}  it follows that there is  $b_1\bydef b\in C\left(\sX, \left\{A_x\right\}, \sF \right)$ such that
\bean
\begin{split}
	\left\|a - \left.b_1\right|_\sY \right\| < \eps,\\
	\left\| b_1\right\| <  \left\|a \right\|+  \eps.
\end{split}
\eean
	If  $a_1 \bydef a - \left.b_1\right|_\sY\in   C\left(\sY, \left\{A_x\right\}, \sF \right)$  then one has $\left\|a_1 \right\|< \eps$. Applying \eqref{top_abk_eqn} once again (replacing $\eps$ with $\eps/2$) one can obtain  $b_2\in C_0\left(\sX, \left\{A_x\right\}, \sF \right)$ such that
\bean
\begin{split}
	\left\|a_1 - \left.b_2\right|_\sV \right\| < \frac{\eps}{2},\\
	\left\| b_2\right\| <  \left\|a_1 \right\|+ \frac{\eps}{2} < \frac{3\eps}{2},\\
\left\| a -	(b_1 + b_2)|_\sV\right\|  = \left\| a_1 - b_2|_\sV\right\|< \frac{\eps}{2}. 
\end{split}
\eean
	If  $a_2\bydef a_1 - \left.b_2\right|_\sV$ then  $a_1, a_2, b_1, b_2$ satisfy to the equations \eqref{top_abn_eqn}. If $a_n$ and $b_n$ are known then using \eqref{top_abk_eqn} one can find $b_{n+1}$ such that
\be\label{top_abnp_eqn}
\begin{split}
	\left\|a_n - \left.b_{n+1}\right|_\sY \right\| < \frac{\eps}{2^n},\\
	\left\| b_{n+1}\right\| <  \left\|a_n \right\|+ \frac{\eps}{2^{n}} <   \frac{3\eps}{2^{n}},\\
\left\|a - \left.\left( \sum_{n = 1}^m b_n\right) \right|_\sV \right\| < \frac{\eps}{2^{m-1}}.
\end{split}
\ee
From $	\left\| b_{n+1}\right\|<    \frac{3\eps}{2^{n}}$ it follows that the series $b\bydef \sum_{n = 1}^\infty b_n$ is norm convergent. Otherwise from \eqref{top_abn_eqn}, \eqref{top_abnp_eqn} it follows that
$$
\left\|a - \left.\left( \sum_{n = 1}^m b_n\right) \right|_\sV \right\| < \frac{\eps}{2^{m-1}},
$$
hence one has $a = \left.b\right|_\sY$.

\end{proof}

\begin{remark}
The Lemma \ref{top_ax_k_lem} is an analog of the Tietze extension theorem \ref{tietze_ext_thm}. The method of proof of the Lemma \ref{top_ax_k_lem} is similar to the explained in \cite{munkres:topology} one.
\end{remark}
Below in this section we assume that all topological spaces are  connected and locally connected (cf. Definition \ref{top_locally_connected_defn}).

\begin{empt}\label{top_lift_empt}
	If $p: \widetilde{\sX} \to \sX$ is a covering then the family $\left\{A_x\right\}$ naturally induces the family $\left\{\widetilde{A}_{\widetilde{x}}\right\}_{\widetilde{x}\in \widetilde{\sX}}$ such that for any ${\widetilde{x}}\in {\widetilde{\sX}}$ there is the isomorphism 
	\be\label{top_ct_iso_eqn}
	c_{\widetilde{x}}:\widetilde{A}_{\widetilde{x}} \cong {A}_{\left( p\left( \widetilde{x}\right) \right) }.
	\ee
\end{empt}
\begin{definition}\label{top_lift_defn}
	Let $\sF$  be a	continuity structure for $\sX$ {and the} $\left\{A_x\right\}_{x \in \sX}$, let $p: \widetilde{\sX} \to \sX$ be a covering. For any $a \in C\left(\sX, \left\{A_x\right\}, \sF \right)$ we define 
	the family 
	\be\label{comm_lift_eqn}
	\left\{\widetilde{a}_{\widetilde{x}}=c^{-1}_{\widetilde{x}}\left(a_{p\left( \widetilde{x}\right) } \right)\in\widetilde{A}_{\widetilde{x}} \right\}_{\widetilde{x}\in\widetilde{\sX}}
	\ee
	is said to be the  $p$-\textit{lift} of $a$ and we write 
	\be\label{comm_lift_d_eqn}
	\left\{\widetilde{a}_{\widetilde{x}}\right\}\bydef \lift_p\left( a\right) .
	\ee
\end{definition}
\begin{lemma}\label{top_lift_lem}
	Let $\sF$ be a	continuity structure  for $\sX$ {and the} $\left\{A_x\right\}_{x \in \sX}$, let $p: \widetilde{\sX}\to \sX$ be a covering. The space
	\be\label{top_lift_eqn}
	\widetilde{   \sF} = \left\{\left.\lift_p\left( a\right)~\right|~a \in C\left(\sX, \left\{A_x\right\}_{x \in \sX}, \sF  \right)\right\} 
	\ee
	is continuity structure  for $\widetilde{   \sX}$ {and the} $\left\{\widetilde{   A}_{\widetilde{x}}\right\}_{\widetilde{x} \in \widetilde{   \sX}}$.	
\end{lemma}
\begin{proof}
	
	Check (a) - (c) of the Definition \ref{operator_fields_continuity_defn}. \\
	(a) For every  $\widetilde x_0 \in \widetilde \sX$ there is an open neighborhood $\widetilde \sU$ such that the restriction $\left.p\right|_{\widetilde\sU}$ is injective.	 For any $a \in C\left(\sX, \left\{A_x\right\}_{x \in \sX}, \sF  \right)$ the function $\text{norm}_a$ is continuous at $x_0= p\left(\widetilde x_0\right)$. Otherwise for all $\widetilde x \in \widetilde\sU$ one has
	$
	\text{norm}_{\lift_p\left(a \right) }\left( \widetilde x\right)= \text{norm}_a\left(p\left( \widetilde x\right)  \right)$, hence  $\text{norm}_{\lift_p\left(a \right) }$ is  continuous at $\widetilde{x}_0$.\\
	(b) The subspace 
	$$
	\left\{\left. a_{x_0} \in A_{x_0}~\right|~ a \in C\left(\sX, \left\{A_x\right\}, \sF  \right) \right\}
	$$
	is dense in $A_{x_0}$ so the space 
	$$
	\left\{\left. \lift_p\left( a\right)_{\widetilde{x}_0} \in \widetilde{A}_{\widetilde{x}_0}~\right|~ a \in C\left(\sX, \left\{A_x\right\}, \sF  \right) \right\}
	$$
	is dense in $\widetilde{A}_{\widetilde{x}_0}$.\\
	(c) If $A_x$ is a $C^*$-algebra for all $x \in \sX$ then from the condition (c) of the Definition \ref{operator_fields_continuity_defn} it follows that $\sF$ is  closed under pointwise multiplication and involution, i.e. for any $a, b \in \sF$ represented by families $\left\{a_x\right\}$ and $\left\{b_x\right\}$ respectively the families $\left\{a^*_x\right\}$ and $\left\{a_xb_x\right\}$ lie in $\sF$. If $\left\{\widetilde{a}^*_{\widetilde{x}}\right\}= \lift_p\left( a^*\right)$ and $\left\{\widetilde{a}_{\widetilde{x}}\widetilde{b}_{\widetilde{x}}\right\}= \lift_p\left( ab\right)$ then from \eqref{comm_lift_eqn} it turns out that $\lift_p\left( a^*\right), \lift_p\left(ab\right) \in \widetilde\sF$.
\end{proof}
\begin{definition}\label{top_lift_main_defn}
	Let $\widetilde{   \sF}$ be given by \eqref{top_lift_eqn} continuity structure. The space of continuous sections (with respect to) $\widetilde{\sF}$ vector fields is said to be $p$-\textit{lift} of $C\left(\sX, \left\{A_x\right\}, \sF \right)$. We write
	\bea\label{top_cont_lift_eqn}
	\lift_p\left[ C\left(\sX, \left\{A_x\right\}, \sF \right) \right]  \stackrel{\text{def}}{=}	C\left(\widetilde{\sX}, \left\{\widetilde{A}_{\widetilde{x}}\right\}, \widetilde{\sF} \right),\\
	\label{top_cont_s_lift_eqn}
	\lift_p\left[  \sF\right]  \stackrel{\text{def}}{=}	C\left(\widetilde{\sX}, \left\{\widetilde{A}_{\widetilde{x}}\right\}, \widetilde{\sF} \right)
	\eea
\end{definition}
Following lemma is a direct consequence of the Lemma \ref{top_lift_lem} and the Definition \ref{top_lift_main_defn}.
\begin{lemma}\label{top_lift_composition_lem}.
	If $\sF$ is a continuity structure  for $\sX$ {and the} $\left\{A_x\right\}$ then one has:
	\begin{enumerate}
		\item [(i)] For any covering  $p :\widetilde{\sX} \to \sX$ there is the natural injective $C_b\left(\sX \right)$-linear  map
		\be\label{top_glo_lift_h_eqn}
		\lift_p: C\left(\sX, \left\{A_x\right\}, \sF \right) \hookto \lift_p\left[C\left(\sX, \left\{A_x\right\}, \sF \right) \right]. 
		\ee
		\item[(ii)] 	 If $p' :\sX' \to \sX$ and $p'' :\sX'' \to \sX'$ are coverings then one has
		\bea
		\label{top_lift_composition_eqn}
		\lift_{p' \circ p''}\left[ C\left(\sX, \left\{A_x\right\}, \sF\right)  \right] = 	\lift_{p'}\left[	\lift_{p''}\left[C\left(\sX, \left\{A_x\right\}, \sF \right) \right]  \right],\\
		\label{top_lift_composition_a_eqn} 	\forall a \in  C\left(\sX, \left\{A_x\right\}, \sF \right) \quad \lift_{p' \circ p''}\left(a\right) = 	\lift_{p'}\circ	\lift_{p''}\left(a\right).
		\eea
	\end{enumerate}
	
\end{lemma}
\begin{remark}\label{top_top_lift_rem}
Sometimes we write $\lift_p^{\text{top}}$ instead of $\lift_p$ for distinguishing  this lift with the given by the sheaf theory one (cf. equation  \eqref{comm_sh_lift_desc_l_eqn}).
\end{remark}
The following Lemma is evident.
\begin{lemma}\label{top_lift_bounded_lem}
	For any $a \in C_b\left(\sX, \left\{A_x\right\}, \sF \right)$ one has
	\be\label{top_lift_bounded_eqn}
	\left\|a\right\|= \left\|\lift_p\left( a\right) \right\|,
	\ee
	hence the map \eqref{top_glo_lift_h_eqn} induces the norm preserving map
	\be\nonumber
	\lift_p:  C_b\left(\sX, \left\{A_x\right\}, \sF \right)\hookto C_b\left(\widetilde{\sX}, \left\{\widetilde{A}_{\widetilde{x}}\right\}, \widetilde{\sF} \right),
	\ee 
	i.e. one has the isometry
	\be\label{top_glo_lift_bh_eqn}
	\lift_p:  C_b\left(\sX, \left\{A_x\right\}, \sF \right)\hookto C_b\left(\lift_p\left[C\left(\sX, \left\{A_x\right\}, \sF \right) \right] \right).
	\ee 
\end{lemma}

\begin{lemma}\label{top_compact_c0_lem}
	If $p: \widetilde{\sX} \to \sX$ is a finite-fold covering then one has
	\be\label{top_compact_c0_eqn}
	\begin{split}
		\lift_p\left( C_c\left(\sX, \left\{A_x\right\}, \sF \right)\right)  \subset C_c\left(\widetilde{\sX}, \left\{\widetilde{A}_{\widetilde{x}}\right\}, \widetilde{\sF} \right)=\\= C_c\left(\lift_p\left[ C\left(\sX, \left\{A_x\right\}, \sF \right) \right] \right), \\
	\lift_p\left( 	C_0\left(\sX, \left\{A_x\right\}, \sF \right)\right)  \subset C_0\left(\widetilde{\sX}, \left\{\widetilde{A}_{\widetilde{x}}\right\}, \widetilde{\sF} \right)=\\= C_0\left(\lift_p\left[ C\left(\sX, \left\{A_x\right\}, \sF \right) \right] \right). \\
	\end{split}
	\ee
\end{lemma}
\begin{proof}
	If $a \in C_c\left(\sX, \left\{A_x\right\}, \sF \right)$ then $\supp a$ is compact. Since $p$ is a finite-fold covering $\supp \lift_p\left( a\right)= p^{-1}\left( \supp a\right) \subset \widetilde{\sX}$ is compact (cf. Lemma \ref{top_compact_preimage_lem}), hence $\lift_p\left( a\right)\in C_c\left(\widetilde{\sX}, \left\{\widetilde{A}_{\widetilde{x}}\right\}, \sF \right)$. For any $a \in C_0\left(\sX, \left\{A_x\right\}, \sF \right)$ there is a net $\left\{a_\a\in C_c\left(\sX, \left\{A_x\right\}, \sF \right)\right\}_{\a \in \mathscr A}$ such that there is a norm limit $\lim_{\a}a_\a = a$. Otherwise for any $\al, \bt\in \mathscr A$ one has $\left\|a_\a - a_\bt\right\|= \left\|\lift_p\left( a_\a\right)  - \lift_p\left( a_\bt\right) \right\|$ it follows that the net $\left\{\lift_p\left( a_\a\right)\in C_c\left(\widetilde{\sX}, \left\{\widetilde{A}_{\widetilde{x}}\right\}, \widetilde\sF \right)\right\}$ is convergent and $\lim_{\a} \lift_p\left( a_\a\right) = \lift_p\left( a\right)$.
\end{proof}

\begin{remark}\label{top_compact_c0_rem}
	Similarly to the Remark \ref{top_cl_cc_iso_rem} $C_0\left(\lift_p\left[ C\left(\sX, \left\{A_x\right\}, \sF \right) \right] \right)$
	is the  minimal norm closed subspace of $C_b\left(\lift_p\left[ C\left(\sX, \left\{A_x\right\}, \sF \right) \right] \right)$ which contains all products $\widetilde f \lift_p\left( a\right) $ where $a \in C_0\left(\sX, \left\{A_x\right\}, \sF \right)$ and $\widetilde f \in C_0\left(\widetilde\sX\right)$.
\end{remark}
\begin{exercise}
Let $M$ be a differential manifold (cf. Definition \ref{diff_mani_defn}), let $p: \widetilde M \to M$ be a covering (cf. Definition \ref{top_covering_defn}). Suppose that $\widetilde M$ has a given by the Proposition \ref{top_cov_mani_prop} structure of differential manifold. Prove that
\be\label{top_lift_smooth_eqn}
\lift_p\left( \Coo\left( M\right)\cap C_b\left( M\right) \right) \subset  \Coo\left(\widetilde M\right)\cap C_b\left(\widetilde M\right).
\ee

\end{exercise}
\begin{empt}\label{top_norm_sub_empt} 
	If $B \subset 	C_0\left(A \right) = C_0\left(\sX, \left\{A_x\right\}, A \right)$  is a norm closed $\C$-space then for any $x \in \sX$ there is a normed $\C$-subspace
	\be
	\mathring{B}_x = \left\{\left. b_x \in A_x\right|~ \exists b' \in B \quad b_x = b'_x  \right\}
	\ee
	Denote by $B_x$ the norm closure of $\mathring{B}_x$ and suppose $B_x \neq \{0\}$ for any $x \in \sX$.  
\end{empt}
\begin{lemma}\label{top_norm_sub_lem} 
Under the hypotheses \ref{top_norm_sub_empt} the $\C$-space $B$ satisfies to the conditions  (a),(b) of the Definition \ref{operator_fields_continuity_defn}.
\end{lemma}
\begin{proof}

	\begin{enumerate}
	\item[(a)] If $b \in B$ then $b \in C_0\left(\sX, \left\{A_x\right\}, A \right)$, so from (a) the Definition \ref{operator_fields_continuity_defn} it follows that the map  $x \mapsto \left\| b_x\right\|$  is continuous on $\sX$.
	\item[(b)] $B_x$ the norm closure of $\mathring{B}_x$, so $\mathring{B}_x$ is dense in $B_x$.

	\end{enumerate}
\end{proof}

\begin{definition}\label{top_norm_sub_cc_defn}
	Under the hypotheses \ref{top_norm_sub_empt} the family $\left\{B_x\right\}_{x \in \sX}$ is the $B$-\textit{restriction} of  $\left\{A_x\right\}_{x \in \sX}$.
\end{definition}

\begin{empt}\label{top_cstr_res_empt}
From the Corollary \ref{top_cl_cc_iso_cor} it turns out that there is the natural isomorphism
	\be\label{top_cstr_res_eqn}
	B \cong C_0\left(\sX, \left\{B_x\right\}, B \right).
	\ee
	
\end{empt}

\begin{definition}\label{top_norm_sub_defn}
	Under the hypotheses \ref{top_cstr_res_empt} we define
	\bea
	\label{top_sub_ccn_eqn}
	C_c\left(B \right) = B \cap C_c\left( A\right)= \left\{\left.b \in 	B~\right| ~\mathrm{norm}_b \in 	C_c\left(\sX\right)\right\}\\
	\label{top__subcc0_eqn}
	C_0\left(B \right) = B \cap C_0\left( A\right)= \left\{\left.b \in 	B~\right| ~\mathrm{norm}_b \in 	C_0\left(\sX\right)\right\}
	\eea
	There are following inclusions
	\be
	C_c\left(B \right) \subset C_0\left(B \right) = B
	\ee
	such that $C_c\left(B \right)$ is dense in $C_0\left(B \right)$.
\end{definition}
\begin{lemma}\label{top_sub_eq_lem}
	Let us consider a continuity structure $A$ for $\sX$ {and the} $\left\{A_x\right\}_{x \in \sX}$ such that $A \cong C_0\left(\sX, \left\{A_x\right\}, A \right)$. Let both $B', B'' \subset A$ are norm closed $\C$-subspaces. If both $\left\{B'_x\right\}$ and $\left\{B''_x\right\}$ are  $B'$ and $B''$-restrictions  of $\left\{A_x\right\}$ (cf. Definition \ref{top_norm_sub_cc_defn}) then one has
	\be\nonumber
	\forall x \in \sX \quad B''_x \subset B'_x \quad \Rightarrow \quad B'' \subset B'.
	\ee
\end{lemma}
\begin{proof}
	From the Corollary \ref{top_cl_cc_iso_cor} it follows that
	\bean
	B' = C_0\left(\sX, \left\{B'_x\right\}, B' \right).
	\eean
	Select $b'' \in B''$ and $x_0 \in \sX$. For any $\eps > 0$ there is $b' = \left\{b'_x\right\} \in C_0\left(\sX, \left\{B'_x\right\}, B' \right)$ such that $\left\|b'_{x_0} - b''_{x_0} \right\| < \frac{\eps}{2}$. The map $x \mapsto \left\|b'_{x} - b''_{x} \right\|$ is continuous (cf. Lemma \ref{op_cont_con_lem}), it turns out that there is an open neighborhood $\sU$ of $x_0$ such that $\left\|b'_{x} - b''_{x} \right\| < \eps$ for every $x \in \sU$. If follows that $b''$ is continuous with respect to $B'$ (cf. Definition \ref{op_cont_fields_defn}). Taking into account  that  $\text{norm}_{b''} \in C_0\left(\sX\right)$ it turns out that $b'' \in C_0\left(\sX, \left\{B'_x\right\}, B' \right)= B'$, hence $B'' \subset B'$ (cf. Corollary \ref{top_cl_cc_iso_cor}).
\end{proof}
\begin{corollary}\label{top_sub_eq_cor}
	Let us consider a continuity structure $A$ for $\sX$ {and the} $\left\{A_x\right\}_{x \in \sX}$ such that $A \cong C_0\left(\sX, \left\{A_x\right\}, A \right)$. Let both $B', B'' \subset A$ be  norm closed subspaces. If both $\left\{B'_x\right\}$ and $\left\{B''_x\right\}$ are  $B'$ and $B''$-restrictions  of $\left\{A_x\right\}$ (cf. Definition \ref{top_norm_sub_cc_defn}) then one has
	\be\nonumber
	\forall x \in \sX \quad B''_x = B'_x \quad \Rightarrow \quad B'' = B'.
	\ee
\end{corollary}

\begin{empt}
	Let $\sF$ be a	continuity structure for $\sX$ {and the} $\left\{A_x\right\}_{x \in \sX}$, let $p: \widetilde{\sX}\to \sX$ be a covering. Let $B \subset C_0\left(\sX, \left\{A_x\right\}, \sF \right)$ is a norm closed $C_0\left(\widetilde{\sX} \right)$-module such that there is  the $B$-{restriction} $\left\{B_x\right\}_{x \in \sX}$ of  $\left\{A_x\right\}_{x \in \sX}$ (cf. Definition \ref{top_norm_sub_cc_defn}). If $p: \widetilde{\sX}\to \sX$ is a finite-fold covering then there are  $p$-lifts  $\left\{\widetilde{A}_{\widetilde{x}}\right\}_{\widetilde{x}\in \widetilde{\sX}}$ and  $\left\{\widetilde{B}_{\widetilde{x}}\right\}_{\widetilde{x}\in \widetilde{\sX}}$  of both $\left\{A_x\right\}_{x \in \sX}$ and $\left\{B_x\right\}_{x \in \sX}$ respectively (cf. \ref{top_lift_defn}). For any $x \in \sX$ there is the natural inclusion $B_x \subset A_x$, which yields the inclusion
	\be\label{top_bx_ax_eqn}
	\widetilde{B}_{\widetilde{x}} = c^{-1}_{\widetilde x}\left(B_x \right) \subset \widetilde{A}_{\widetilde{x}} \quad \forall \widetilde{x}\in {\widetilde{\sX}}\quad c_{\widetilde x}\text{ is given by }\eqref{top_ct_iso_eqn}.
	\ee
	From \eqref{top_bx_ax_eqn} and $B \subset  C_0\left(\sX, \left\{A_x\right\}, \sF \right)$ for any $b \in B$ one has
	$$
	\lift_p\left(b \right) \in \lift_p\left[ C_0\left(\sX, \left\{A_x\right\}, \sF \right)\right]. 
	$$
	Taking into account the Definition \ref{top_lift_main_defn} and equations 
	\eqref{top_ccn_eqn}, \eqref{top_cc0_eqn}, \eqref{top_lift_bounded_eqn}
	one has 
	\be\label{top_lift_inc_eqn}
	\begin{split}
		\lift_p\left[C_0\left(\sX, \left\{B_x\right\}, B \right)\right] \subset \lift_p\left[C_0\left(\sX, \left\{A_x\right\}, \sF \right)\right],\\
		C_c\left( \lift_p\left[C_0\left(\sX, \left\{B_x\right\}, B \right)\right]\right)  \subset C_c\left( \lift_p\left[C_0\left(\sX, \left\{A_x\right\}, \sF \right)\right]\right) ,\\
		C_0\left( \lift_p\left[C_0\left(\sX, \left\{B_x\right\}, B \right)\right]\right)  \subset C_0\left( \lift_p\left[C_0\left(\sX, \left\{A_x\right\}, \sF \right)\right]\right) ,\\
		C_b\left( \lift_p\left[C_0\left(\sX, \left\{B_x\right\}, B \right)\right]\right)  \subset C_b\left( \lift_p\left[C_0\left(\sX, \left\{A_x\right\}, \sF \right)\right]\right) .\\
	\end{split}
	\ee	 
	The space $	\lift_p\left(B\right)$ is a {continuity structure for} $\widetilde\sX$ {and the} $\left\{\widetilde B_{\widetilde x}\right\}$ (cf. Definition \ref{operator_fields_continuity_defn}).	
\end{empt}

\begin{definition}\label{top_lift_sub_defn}
	In the above situation the space $\widetilde B$ of {continuous (with respect to $\lift_p\left(B\right)$)} vector fields $\left\{b_{\widetilde x}\in \widetilde B_{\widetilde x}\right\}_{{\widetilde x} \in \widetilde{\sX}}$ (cf. Definition \ref{op_cont_fields_defn}) is said to be the $p$-\textit{lift} of $B$. We write
	\be\label{top_lift_sub_eqn}
	\lift_{p}\left[B\right]\bydef \widetilde B.
	\ee
\end{definition}

\begin{remark}
	From the Definition \ref{top_lift_sub_defn} it turns out that there are injective $C_b\left(\sX \right)$-linear maps
	\be\label{top_lift_sub_inc_eqn}
	\begin{split}
		\lift_p: B \hookto \lift_p\left[B \right], \\
		\lift_p:	C_b\left(B \right)  \hookto C_b\left( \lift_p\left[B \right]\right).
	\end{split}
	\ee 
	Moreover if $p$ is a finite-fold covering there are following inclusions
	\be\label{top_lift_sub_ff_eqn}
	\begin{split}
		\lift_p: C_c\left( B\right) \hookto C_c\left(\lift_p\left[ B\right]  \right);\\
		\lift_p: C_0\left( B\right) \hookto C_0\left(\lift_p\left[ B\right]  \right)
	\end{split}
	\ee
	(cf. Definition \ref{top_norm_sub_defn}).
\end{remark}

\begin{remark}
	Similarly to \eqref{top_lift_composition_eqn} one has
	\bea
	\label{top_lift_composition_sub_eqn}
	\lift_{p' \circ p''}\left[ B\right] = 	\lift_{p'}\left[	\lift_{p''}\left[B \right]\right] .
	\eea
\end{remark}

\begin{definition}\label{top_lift_desc_defn}
	Let $\sF$ 	continuity structure for $\sX$ {and the} $\left\{A_x\right\}_{x \in \sX}$, let $p: \widetilde{\sX}\to\sX$ be a covering. Let $\widetilde{\sU}\subset \widetilde{\sX}$ be an open subset such that the restriction $\left.p\right|_{\widetilde{\sU}}$ is injective. Let $\sU\bydef p\left( \widetilde{\sU}\right)$. If $a \in C\left(\sX, \left\{A_x\right\}, \sF \right)$ is such that $\supp a \subset \sU$ (cf. Definition \ref{top_support_defn}) and $\widetilde{a}'= \lift_p\left( a\right)$ then 
	the family 
	\be\label{comm_lift_desc_eqn}
	\widetilde{a}_{\widetilde{x}}=\left\{\begin{array}{c l}
	\widetilde{a}'_{\widetilde{x}}=	c^{-1}_{\widetilde x}\left( {a}_{p\left( \widetilde{x} \right) }\right)    &{\widetilde{x}}\in \widetilde{\sU} \\
		0 &{\widetilde{x}}\notin \widetilde{\sU}\\
	\end{array}\right. \quad c_{\widetilde x}\text{ is given by }\eqref{top_ct_iso_eqn}
	\ee 
	is said to be the $p$-$\widetilde{\sU}$-\textit{lift} or simply the $\widetilde{\sU}$-\textit{lift} of $a$. Otherwise we say that $a$ is the $p$-\textit{descent} of $ \widetilde{a}$. We write
	\be\label{top_lift_desc_eqn}
	\begin{split}
		\widetilde{a}\stackrel{\text{def}}{=} \lift^p_{\widetilde{\sU}}\left(a \right)  \text{ or simply } \widetilde{a}\stackrel{\text{def}}{=} \lift_{\widetilde{\sU}}\left(a \right), \\
		a \stackrel{\text{def}}{=} \desc_{p} \left(\widetilde{a}  \right)  \text{ or simply } a \stackrel{\text{def}}{=} \desc \left(\widetilde{a}  \right).
	\end{split}
	\ee
	If $\widetilde{\sV}$ is the closure of $\widetilde{\sU}$ then we define  the \textit{closed}-$p$-$\widetilde{\sU}$-\textit{lift}   of $a$ by the following way
	\be\label{top_lift_desc_closed_eqn}
\begin{split}
\lift^p_{\widetilde{\sV}}\left(a \right)  \bydef \lift^p_{\widetilde{\sU}}\left(a \right).
\end{split}
\ee
	
\end{definition}
\begin{remark}
	if $a = \desc_{p} \left(\widetilde{a}  \right)$ then $a$ is represented by the family $\left\{a_x\right\}_{x \in \sX}$ such that
	\be\label{top_desc_eqn}
	a_x = \begin{cases}
		c_{\widetilde x}\left( 	\widetilde{a}_{\widetilde{x}}\right)  & x \in \sU ~\text{AND} ~ \widetilde x \in \widetilde\sU~\text{AND} ~p\left( \widetilde{x}\right) = x \\
		0 &  x \notin \sU
	\end{cases}	\quad 
	\ee 
	where	$c_{\widetilde x}$ \text{ is given by }\eqref{top_ct_iso_eqn}.
\end{remark}
\begin{remark}\label{comm_lift_desc_l_rem}
	From the equation \eqref{comm_lift_eqn} it turns out 
	\begin{equation}\label{comm_lift_desc_l_eqn}
	\begin{split}
	a = \mathfrak{desc}_{p}\left( \mathfrak{lift}^p_{\widetilde{\mathcal U}}\left(a \right)\right) ,\\
	\widetilde{a} = \mathfrak{lift}^p_{\widetilde{\mathcal U}}\left( \mathfrak{desc}_{p}\left(\widetilde{a} \right)\right). 
	\end{split}
	\end{equation}	
\end{remark}
\begin{defn}\label{top_compactly_supported_descent_defn}
	Let $\sF$ 	continuity structure for $\sX$ {and the} $\left\{A_x\right\}_{x \in \sX}$, and let $p: \widetilde \sX \to \sX$ be a transitive covering such that the group  of covering transformations $G\left(\left.\widetilde\sX\right| \sX  \right)$ (cf. Definition \ref{top_group_of_covering_transformations_defn}) is residually finite (cf. Definition \ref{residually_finite_defn}). If $\widetilde{a}\in  C_c\left( \lift_p\left[C_0\left(\sX, \left\{A_x\right\}, \sF \right)\right]\right)$ then $\supp \widetilde{a}$ (cf. Equation \eqref{top_support_eqn}) is compact. From the 
Theorem  \ref{top_compact_thm} 	it follows that there is a finite-fold transitive covering $p': \widetilde \sX' \to \sX$ and  a transitive covering $\widetilde p': \widetilde \sX \to \widetilde \sX'$ such that $p = p' \circ \widetilde p'$ and $\supp \widetilde{a}$ is mapped homeomorphically onto $\widetilde p'\left(\supp \widetilde{a} \right)$. If
\be\label{top_desc_compact_eqn}
\desc^c_p\left(\widetilde a \right) \bydef \sum_{	g \in G\left( \left. \widetilde \sX'\right| \sX\right) } g \desc_{\widetilde p'} \left(\widetilde a \right)\in C_c\left(\sX, \left\{A_x\right\}, \sF \right)
\ee
then we say that $\desc^c_p$ is \textit{compactly supported} $p$-\textit{descent}.
\end{defn}
\begin{remark}\label{top_desc_compact_rem}
	Under the hypotheses of the Definition \ref{top_compactly_supported_descent_defn} the  map  \eqref{top_desc_compact_eqn} induces a homomorphism 
	$$
\desc^c_p : C_c\left( \lift_p\left[C_0\left(\sX, \left\{A_x\right\}, \sF \right)\right]\right)\to  C_0\left(\sX, \left\{A_x\right\}, \sF \right)
	$$
	of $C_0\left(\sX \right)$-modules. 
\end{remark}

\begin{empt}
 Let $\sX$ be a locally compact Hausdorff space and 	there is a positive functional 
\be\label{top_tau_eqn}
\begin{split}
\tau : C_c\left( \sX\right) \to \C,\\
a \mapsto  \int_{\sX} a~ d\mu.
\end{split}
\ee
(cf. Theorem \ref{meafunc_thm}) such that $a > 0  \Rightarrow \tau\left(a \right)>0$. If $p: \widetilde \sX\to \sX$ is a covering then one can define a functional

\be\label{top_tau_l_eqn}
\begin{split}
\lift_p\tau: C_c\left( \widetilde\sX\right) \to \C,\\
\widetilde a \mapsto\tau\left(\desc^c_p\left(\widetilde a \right) \right)
\end{split}
\ee
such that $\widetilde a > 0 \quad \Rightarrow \quad \lift_p\tau\left(\widetilde a \right) > 0$.	
From the Theorem \ref{meafunc_thm} one has a measure $\lift_p\mu$ such that
\be
\lift_p\tau \left( \widetilde a\right)=  \int_{\widetilde\sX} \widetilde a~ d~\lift_p\mu
\ee
\end{empt}
\begin{definition}\label{top_lift_measure_defn}
In the described above situation we say:
\begin{itemize}
	\item the functional $\lift_p\tau$ is the $p$-\textit{lift of the functional}  $\tau$,
		\item the measure $\lift_p\mu$ is the $p$-\textit{lift of the measure}  $\mu$.
	\end{itemize}
\end{definition}

\begin{lemma}\label{top_lift_bundle_lem}
	Let  $\left( E, \pi, \mathcal X\right)$ be a  vector bundle with the family of fibers $\left\{E_x\right\}_{x \in \sX}$ and let $\Ga\left( E, \pi, \mathcal X\right)$ be a space of sections which is a continuous structure for $\sX$ and $\left\{E_x\right\}_{x \in \sX}$ (cf. Lemma \ref{top_bundle_cs_ex_lem}). If $p: \widetilde{\sX} \to \sX$ is a covering and $\left( E\times_\sX\widetilde\sX,\rho, \widetilde\sX\right)$ is the {inverse image} of $\left( E, \pi, \mathcal X\right)$ by $p$ (cf. Definition \ref{vb_inv_img_funct_defn}) then there is the natural $\C$-isomorphism
	$$
	\lift_p\left[C\left(\sX, \left\{E_x\right\}, \Ga\left( E, \pi, \mathcal X\right) \right) \right]\cong \Ga\left( E\times_\sX\widetilde\sX,\rho, \widetilde\sX\right).
	$$
\end{lemma}
\begin{proof}
	From the  equations \eqref{top_uba_eqn} and \eqref{top_cg_eqn} it follows that there is the natural inclusion
	$$
	\lift_p\left(C\left(\sX, \left\{E_x\right\}, \Ga\left( E, \pi, \mathcal X\right) \right) \right)\subset \Ga\left( E\times_\sX\widetilde\sX,\rho, \widetilde\sX\right).
	$$
From the Definition \ref{top_lift_main_defn} it turns out that
	$$
\lift_p\left[C\left(\sX, \left\{E_x\right\}, \Ga\left( E, \pi, \mathcal X\right) \right) \right] \subset C\left(\widetilde\sX, \left\{\left(E\times_\sX\widetilde\sX\right)_{\widetilde x}\right\}_{\widetilde{x}\in \widetilde\sX}, \Ga\left( E\times_\sX\widetilde\sX,\rho, \widetilde\sX\right) \right). 	
	$$
	However from \eqref{top_cg_eqn} it follows that 
	$$
	C\left(\widetilde\sX, \left\{\left(E\times_\sX\widetilde\sX\right)_{\widetilde x}\right\}_{\widetilde{x}\in \widetilde\sX}, \Ga\left( E\times_\sX\widetilde\sX,\rho, \widetilde\sX\right) \right)\cong \Ga\left( E\times_\sX\widetilde\sX,\rho, \widetilde\sX\right),
	$$
	so one has
	$$
	\lift_p\left[C\left(\sX, \left\{E_x\right\}, \Ga\left( E, \pi, \mathcal X\right) \right) \right]\subset \Ga\left( E\times_\sX\widetilde\sX,\rho, \widetilde\sX\right).	
	$$
	Let $\widetilde s \in \Ga\left( E\times_\sX\widetilde\sX,\rho, \widetilde\sX\right)$ and let $\widetilde x_0 \in \widetilde \sX$. If   $\widetilde{f}_{\widetilde{x}_0}$ is a $p$-$\widetilde x_0$-{stump} (cf. Definition \ref{top_stump_cov_defn}) then there is an open subset $\widetilde \sU$ such that $\supp \widetilde{f}_{\widetilde{x}_0}\subset \widetilde\sU$ and $\widetilde\sU$ is mapped homeomorphically onto $\sU \subset p\left(\widetilde\sU\right)$. There is a natural isomorphism
	$$
\phi :	\Ga\left( E|_\sU, \pi|_{E|_\sU}, \sU \right)\xrightarrow{\approx} 	\Ga\left( \left.\left( E\times_\sX\widetilde\sX\right)\right|_{\widetilde\sU}, \rho|_{E|\left( E\times_\sX\widetilde\sX\right) |_{\widetilde\sU}}, \widetilde\sU\right)
	$$
	If $\widetilde s' = \widetilde{f}_{\widetilde{x}_0}\widetilde s$ and $\widetilde s'=\phi\left( s\right)$ then $s$ can be extended to the global section
\bean
s \in \Ga\left( E, \pi, \mathcal X\right),\\
s\left(x\right)= \begin{cases}
s'\left(x\right)& x \in \sU\\
0 & x\notin \sU
\end{cases}.
\eean
If $\widetilde \sV$ is an open neighborhood of $\widetilde{x}_0$ such that $\widetilde{f}_{\widetilde{x}_0}\left(\widetilde \sV\right)= \{1\}$ then
$$
\forall \widetilde x \in \widetilde{\sV}\quad \lift_p\left( s\right) \left(\widetilde x \right) = \widetilde s\left( \widetilde x\right). 
$$
It follows that $\widetilde s$ is {continuous (with respect to $\lift_p\left( C\left(\sX, \left\{E_x\right\}, \Ga\left( E, \pi, \mathcal X\right) \right)\right) $)} at $\widetilde x_0$ (cf. Definition \ref{op_cont_fields_defn}). A point $\widetilde x_0$ is arbitrary so  $\widetilde s$ is {continuous} on $\widetilde \sX$, and taking into account Definition \ref{top_lift_main_defn} one has $\widetilde s \in \lift_p\left[C\left(\sX, \left\{E_x\right\}, \Ga\left( E, \pi, \mathcal X\right) \right) \right]$, it turns out that
$$
\Ga\left( E\times_\sX\widetilde\sX,\rho, \widetilde\sX\right)\subset \lift_p\left[C\left(\sX, \left\{E_x\right\}, \Ga\left( E, \pi, \mathcal X\right) \right) \right].
$$
\end{proof}

 \section{Continuity structures and $C^*$-algebras}
\subsection{Basic definitions}
 \paragraph{}
 From the Lemma \ref{oa_haus_alg_lem} it follows that any $C^*$-algebra $A$ with Hausdorff spectrum $\sX$ is a full algebra of operator fields (cf. Definition \ref{full_algebra_operator_fields_defn}, i.e.
\be\label{top_a_ax_eqn}
 A \cong C_0\left( \sX, \left\{\rep_x  \left(  A\right) \right\}_{x\in \sX}, A \right).
\ee
 where  $C_0\left( \sX, \left\{\rep_x  \left(  A\right) \right\}_{x\in \sX}, A \right)$ is the {converging to zero submodule} (cf. Definition \ref{top_cc_c0_defn})

\begin{lemma} 
	In the above situation   $C_b\left(\sX, \left\{\rep_x  \left(  A\right) \right\}_{x\in \sX}, A \right)$ is a $C^*$-algebra.
\end{lemma}
\begin{proof}
	Follows from the condition (c) of the Definition \ref{operator_fields_continuity_defn} and the Lemma  \ref{op_cont_module_lem}.
\end{proof}

\begin{lemma}\label{top_mult_inc_b_lem}
	In the above there is the natural inclusion
	$$
	C_b\left(\sX, \left\{ A_{x}\right\}, A  \right)  \hookto M\left(C_0\left(\sX, \left\{ A_{x}\right\}, A \right)\right) 
	$$
	of $C^*$-algebras.
\end{lemma}	
\begin{proof}
	If $a \in 	C_0\left(\sX, \left\{ A_{x}\right\}, A  \right)$ and $b \in 	C_b\left(\sX, \left\{ A_{x}\right\}, A  \right)$ then $ab  \in 	C_b\left(\sX, \left\{ A_{x}\right\}, A  \right)$. Form the Lemma
	\ref{top_norm_c0cc_lem} it follows that $\text{norm}_a\in C_0\left(\sX \right)$ and from  the Lemma \ref{top_bounded_norm_lem} it follows that $\text{norm}_a, \text{norm}_{ab}\in C_b\left(\sX \right)$. Taking into account $\text{norm}_{ab} \le \text{norm}_{a}\text{norm}_{b}$ we conclude that  $\text{norm}_{ab}\in C_0\left(\sX \right)$, hence from the Lemma \ref{top_norm_c0cc_lem} one has $ab \in C_0\left(\sX, \left\{ A_{x}\right\}, A  \right)$. Similarly we prove that $ba \in C_0\left(\sX, \left\{ A_{x}\right\}, A  \right)$, so $b \in M\left(C_0\left(\sX, \left\{ A_{x}\right\}, A  \right)\right)$.
\end{proof}	

\begin{empt}\label{top_a_res_rem}
	If $A = C_0\left(\sX, \left\{ A_{x}\right\}, A  \right)$ is a $C^*$-algebra with Hausdorff spectrum and $\sV\subset\sX$ is any subset then denote by $C_0\left(\sV, \left\{ A_{x}\right\}, A  \right)$ the $\sV$-{restriction} of $C_0\left(\sX, \left\{A_x\right\}, A \right)$ (cf. Definition \ref{ctr_restr_defn}). The space $C_0\left(\sV, \left\{ A_{x}\right\}, A  \right)$ has the natural structure  of *-algebra with $C^*$-norm. It is norm closed (cf. Remark \ref{top_cv_nc_rem}), so $C_0\left(\sV, \left\{ A_{x}\right\}, A  \right)$ is a $C^*$-algebra.
\end{empt}
\begin{lemma}\label{top_ax_sp_lem}
Under the hypotheses \ref{top_a_res_rem} if a subset   $\sY \subset \sX$ is compact then the following conditions hold.
\begin{enumerate}
	\item [(i)] If $\left.A\right|^\sY$ is given by the equation \eqref{closed_ideal_eqn} then there is the natural $*$-isomorphism $\left.A\right|^\sY\cong C_0\left(\sY, \left\{ A_{x}\right\}, A  \right)$  where  $C_0\left(\sY, \left\{ A_{x}\right\}, A  \right)$ is 
	an $\sY$-restriction  of  $C_0\left(\sX, \left\{ A_{x}\right\}, A  \right)$   (cf. Definition \ref{ctr_restr_defn}).
	\item[(ii)] The spectrum of $\left.A\right|^\sY$ coincides with $\sY$.
\end{enumerate}
\end{lemma}
\begin{proof}
(i) From the equation \eqref{closed_ideal_eqn} it follows that there is a natural surjective $*$-homomorphism $p':A \to\left.A\right|^\sY$. From the Lemma \ref{top_ax_k_lem} it turns out that is a natural surjective $*$-homomorphism $p':A \to C_0\left(\sY, \left\{ A_{x}\right\}, A  \right)$. From $\ker p'=\ker p'' = \left.A\right|_{\sX\setminus\sY}$ (cf. equation \eqref{open_ideal_eqn}) it follows that $\left.A\right|^\sY\cong C_0\left(\sY, \left\{ A_{x}\right\}, A  \right)$.\\
(ii)  From the Theorem \ref{ctr_as_field_thm} it follows that the spectrum of $C_0\left(\sY, \left\{ A_{x}\right\}, A  \right)$ equals to $\sY$, so from (i) of this lemma and the Lemma  \ref{oa_haus_alg_lem} it follows that the spectrum of $\left.A\right|^\sY$ coincides with $\sY$.
\end{proof}

\begin{lemma}\label{top_lift_c_alg_lem}
	If $\sF$ is a continuity structure for  $\sX$ {and the} $\left\{A_x\right\}_{x \in \sX}$, such that $A_x$ is a $C^*$-algebra for all $x \in \sX$ and $p: \widetilde{\sX}\to\sX$ is a covering then
	\begin{enumerate}
		\item[(i)] the map
		$\lift_p: C\left(\sX, \left\{A_x\right\}, \sF \right) \hookto \lift_p \left[\sX, \left\{A_x\right\}, \sF \right]$ is an injective $*$-homomorphism of *-algebras,
		\item[(ii)] the restriction  $$\left.\left.\lift\right._p\right|_{C_b\left(\sX, \left\{A_x\right\}, \sF \right)}:  C_b\left(\sX, \left\{A_x\right\}, \sF \right) \hookto C_b\left( \lift_p \left[\sX, \left\{A_x\right\}, \sF \right]\right)$$
		is an injective  $*$-homomorphism of $C^*$-algebras,
		\item[(iii)] if $p$ is a finite-fold covering then there is 	 an injective $*$-homomorphism of $C^*$-algebras
	\be\label{top_fin_lift_eqn}
		\left.\left.\lift\right._p\right|_{C_0\left(\sX, \left\{A_x\right\}, \sF \right)}:  C_0\left(\sX, \left\{A_x\right\}, \sF \right) \hookto C_0\left( \lift_p \left[\sX, \left\{A_x\right\}, \sF \right]\right)
		\ee
	\end{enumerate}
\end{lemma}
\begin{proof}(i)
	If 
	\bean
	{a} \in C\left({\sX},  \left\{{A}_{{x}}\right\}_{{x} \in {\sX}}, {\sF}\right),\quad 
	{b} \in C\left({\sX},  \left\{{A}_{{x}}\right\}_{{x} \in {\sX}}, {\sF}\right)
	\eean
	then $a$, $b$ and $ab$ are represented by families
	$$
	\left\{{a}_{{x}}\in {A}_{{x}}\right\}, \quad \left\{{b}_{{x}}\in {A}_{{x}}\right\}, \quad \left\{a_x{b}_{{x}}\in {A}_{{x}}\right\}.
	$$
	Otherwise $\lift_p\left(a\right), \lift_p\left(a\right), \lift_p\left(ab\right)\in \lift_p \left[\sX, \left\{A_x\right\}, \sF \right]$ are represented by families
	\bean
	\left\{\widetilde{a}_{\widetilde{x}}= {a}_{p\left( \widetilde{x}\right)}\in {A}_{p\left( \widetilde{x}\right) }\right\}_{\widetilde{x}\in \widetilde{\sX}}, \quad \left\{\widetilde{b}_{\widetilde{x}}= {b}_{p\left( \widetilde{x}\right)}\in {A}_{p\left( \widetilde{x}\right) }\right\}_{\widetilde{x}\in \widetilde{\sX}},\\
	\left\{\widetilde{a}_{\widetilde{x}}\widetilde{b}_{\widetilde{x}}= {a}_{p\left( \widetilde{x}\right)}{b}_{p\left( \widetilde{x}\right)}= (ab)_{p\left( \widetilde{x}\right)}\in {A}_{p\left( \widetilde{x}\right) }\right\}_{\widetilde{x}\in \widetilde{\sX}},
	\eean
	hence one has
	\be\label{top_mult_lift_eqn}
	\lift_p\left(ab\right)= \lift_p\left(a\right)\lift_p\left(b\right).
	\ee
	Element $a^*$ is represented by the family $\left\{{a}^*_{{x}}\in {A}_{{x}}\right\}$ and  $\lift_p\left(a\right)$ is represented by the family $\left\{\widetilde{a}^*_{\widetilde{x}}= {a}^*_{p\left( \widetilde{x}\right)}\in {A}_{p\left( \widetilde{x}\right) }\right\}_{\widetilde{x}\in \widetilde{\sX}}$ it turns out that
	\be\label{top_star_lift_eqn}
	\lift_p\left(a^*\right)= \lift_p\left(a\right)^*.
	\ee
	\\
	(ii)
	Follows from (i) and \eqref{top_glo_lift_bh_eqn}.\\
	(iii) Follows from (i) and \eqref{top_compact_c0_eqn}.
\end{proof}

\begin{lemma}\label{top_lift_d_lem}
	Let $\sF$ be a	continuity structure  for $\sX$ {and the} $\left\{A_x\right\}_{x \in \sX}$, let $p: \widetilde{\sX}\to \sX$ be a transitive covering.
	\begin{enumerate}
		\item [(i)] There is the natural action $$G\left(\left.\widetilde{\sX}~\right| \sX \right) \times \lift_p\left[C\left(\sX, \left\{A_x\right\}, \sF \right) \right] \to \lift_p\left[C\left(\sX, \left\{A_x\right\}, \sF \right)\right] $$ which yields actions  
		\bean
		G\left(\left.\widetilde{\sX}~\right| \sX \right) \times C_c\left( \lift_p\left[C\left(\sX, \left\{A_x\right\}, \sF \right) \right]\right)\to C_c\left( \lift_p\left[C\left(\sX, \left\{A_x\right\}, \sF \right) \right]\right) , \\
		G\left(\left.\widetilde{\sX}~\right| \sX \right) \times C_0\left( \lift_p\left[C\left(\sX, \left\{A_x\right\}, \sF \right) \right]\right)\to C_0\left( \lift_p\left[C\left(\sX, \left\{A_x\right\}, \sF \right) \right]\right), \\ 
		G\left(\left.\widetilde{\sX}~\right| \sX \right) \times C_b\left( \lift_p\left[C\left(\sX, \left\{A_x\right\}, \sF \right) \right]\right)\to C_b\left( \lift_p\left[C\left(\sX, \left\{A_x\right\}, \sF \right) \right]\right).
		\eean
		\item[(ii)] If $p': \widetilde{\sX}' \to \widetilde{\sX}$ is the transitive covering then
		\be\label{top_g_inv_eqn}
		\begin{split}
			G\left(\left.\widetilde{\sX}'~\right| \sX \right) 
			\lift_{p'}\left( \lift_p\left[C\left(\sX, \left\{A_x\right\}, \sF \right) \right] \right) =
			\\
			= \lift_{p'}\left( \lift_p\left[C\left(\sX, \left\{A_x\right\}, \sF \right) \right] \right),\\
			G\left(\left.\widetilde{\sX}'~\right| \sX \right) 
			\lift_{p'}\left(C_b\left(  \lift_p\left[C\left(\sX, \left\{A_x\right\}, \sF \right) \right] \right) \right) =\\
			=\lift_{p'}\left(C_b\left(  \lift_p\left[C\left(\sX, \left\{A_x\right\}, \sF \right) \right] \right) \right).
		\end{split}
		\ee
		Moreover if $p'$ is a finite-fold covering then one has
		\be\label{top_g_c0_inv_eqn}
		\begin{split}
			G\left(\left.\widetilde{\sX}'~\right| \sX \right) 
			\lift_{p'}\left(C_c\left(  \lift_p\left[C\left(\sX, \left\{A_x\right\}, \sF \right) \right] \right) \right) =\\
			=\lift_{p'}\left(C_c\left(  \lift_p\left[C\left(\sX, \left\{A_x\right\}, \sF \right) \right] \right) \right),\\
			G\left(\left.\widetilde{\sX}'~\right| \sX \right) 
			\lift_{p'}\left(C_0\left(  \lift_p\left[C\left(\sX, \left\{A_x\right\}, \sF \right) \right] \right) \right) =\\
			=\lift_{p'}\left(C_0\left(  \lift_p\left[C\left(\sX, \left\{A_x\right\}, \sF \right) \right] \right) \right).
		\end{split}
		\ee
		
		\item[(iii)] If $A_x$ is a $C^*$-algebra for any $x \in \sX$ then any $g\in 	G\left(\left.\widetilde{\sX}~\right| \sX \right)$ yields  automorphisms of involutive algebras $\lift_p\left[C\left(\sX, \left\{A_x\right\}, \sF \right) \right]$, $C_c\left( \lift_p\left[C\left(\sX, \left\{A_x\right\}, \sF \right) \right]\right)$, $C_0\left( \lift_p\left[C\left(\sX, \left\{A_x\right\}, \sF \right) \right]\right)$ and $C_b\left( \lift_p\left[C\left(\sX, \left\{A_x\right\}, \sF \right) \right]\right)$.
		
	\end{enumerate} 
	
\end{lemma}
\begin{proof} 
	(i)	
	For every $g \in G\left(\left.\widetilde{\sX}~\right| \sX \right)$ from \eqref{top_ct_iso_eqn} it turns out that
	$$
	\widetilde{A}_{\widetilde{x}} \cong {A}_{ p\left( \widetilde{x}\right) } \cong 	\widetilde{A}_{g\widetilde{x}}
	$$
	If $\widetilde{a} \in \lift_p\left[C\left(\sX, \left\{A_x\right\}, \sF \right) \right]$ corresponds to a family $\left\{\widetilde{a}_{\widetilde{x}}\right\}_{\widetilde{x}\in \widetilde{\sX}}$ then we define $g\widetilde{a}$ such that it is given by the family $\left\{\widetilde{a}_{g\widetilde{x}}\right\}_{\widetilde{x}\in \widetilde{\sX}}$. The given by  \eqref{top_lift_eqn} continuous structure $\widetilde{\sF}$ is  $G\left(\left.\widetilde{\sX}~\right| \sX \right)$- invariant, it turns out that  $\left\{\widetilde{a}_{g\widetilde{x}}\right\}_{\widetilde{x}\in \widetilde{\sX}}$ is {continuous (with respect to  $\widetilde{\sF}$)} (cf. Definition \ref{op_cont_fields_defn}). So for any $g \in G\left(\left.\widetilde{\sX}~\right| \sX \right)$ one has an isomorphism
	\bean
	\lift_p\left[C\left(\sX, \left\{A_x\right\}, \sF \right) \right] \xrightarrow{\approx} \lift_p\left[C\left(\sX, \left\{A_x\right\}, \sF \right)\right],\\
	\widetilde{a} \mapsto g\widetilde{a}.
	\eean
	Any $g \in G\left(\left.\widetilde{\sX}~\right| \sX \right)$ is in fact a homeomorphism it follows that $\supp \widetilde{a}$ is homeomorphic to $\supp g\widetilde{a}$.
	In particular if $\supp \widetilde{a}$ is compact then $\supp g\widetilde{a}$ is also compact, so one has
	$$
	G\left(\left.\widetilde{\sX}~\right| \sX \right)C_c\left( \lift_p\left[C\left(\sX, \left\{A_x\right\}, \sF \right) \right]\right) = C_c\left( \lift_p\left[C\left(\sX, \left\{A_x\right\}, \sF \right) \right]\right).
	$$
	Taking into account that $C_0\left( \lift_p\left[C\left(\sX, \left\{A_x\right\}, \sF \right) \right]\right)$ is the norm completion of $C_c\left( \lift_p\left[C\left(\sX, \left\{A_x\right\}, \sF \right) \right]\right)$ we conclude 
	$$
	G\left(\left.\widetilde{\sX}~\right| \sX \right)C_0\left( \lift_p\left[C\left(\sX, \left\{A_x\right\}, \sF \right) \right]\right) = C_0\left( \lift_p\left[C\left(\sX, \left\{A_x\right\}, \sF \right) \right]\right).
	$$
	From $\left\| \widetilde{a}\right\|= \left\| g\widetilde{a}\right\|$ one concludes
	$$
	G\left(\left.\widetilde{\sX}~\right| \sX \right)C_b\left( \lift_p\left[C\left(\sX, \left\{A_x\right\}, \sF \right) \right]\right) = C_b\left( \lift_p\left[C\left(\sX, \left\{A_x\right\}, \sF \right) \right]\right).
	$$\\
	(ii)
	If  $\widetilde{a} \in  \lift_p\left[C\left(\sX, \left\{A_x\right\}, \sF \right)\right]$ then $\widetilde{a}$ corresponds to a family $\left\{\widetilde{a}_{\widetilde{x}}\in \widetilde{A}_{\widetilde{x}}\right\}_{\widetilde{x}\in \widetilde{\sX}}$ where $\widetilde{A}_{\widetilde{x}}\cong A_{p\left( \widetilde{x}\right) }$ for each $\widetilde{x}\in \widetilde{\sX}$. The element $\widetilde{a}'=\lift_p\left(\widetilde{a}\right) \in \lift_{p \circ p'}\left[C\left(\sX, \left\{A_x\right\}, \sF \right) \right]$ corresponds to the family
	$
	\left\{
	\widetilde{a}'_{\widetilde{x}'}= \widetilde{a}_{p'\left( \widetilde{x}\right)} \in \widetilde{A}'_{\widetilde{x}} 
	\right\}_{\widetilde{x}'\in \widetilde{\sX}'}
	$.
	For all $g' \in G\left(\left.\widetilde{\sX}'~\right| \sX\right)$ the element $g'\widetilde{a}'$ corresponds to the family 
	\be\label{top_gpg_eqn}
	\left\{
	\widetilde{a}'_{g\widetilde{x}'}= \widetilde{a}_{p'\left( h\left(g' \right) \widetilde{x}\right)} \in \widetilde{A}'_{\widetilde{x}} 
	\right\}
	\ee
	where $h: 	G\left(\left.\widetilde{\sX}'~\right| \sX \right) \to 	G\left(\left.\widetilde{\sX}~\right| \sX \right)$ is the natural surjective homomorphism of covering groups induced by the transitive covering $p'$.  From \eqref{top_gpg_eqn} it follows that
	$$
	g'\widetilde{a}'=g' \lift_p\left(\widetilde{a}\right)= \lift_p\left(h\left(g' \right) \widetilde{a}\right)\in \lift_p\left[C\left(\sX, \left\{A_x\right\}, \sF \right)\right].
	$$
	and taking into account \eqref{top_lift_bounded_eqn} one obtains \eqref{top_g_inv_eqn}. 
	If $\widetilde{a} \in C_c\left( \lift_p\left[C\left(\sX, \left\{A_x\right\}, \sF \right) \right]\right)$ then $\supp \widetilde{a}\in \widetilde{\sX}$ is compact. Moreover if $p'$ is a finite-fold covering then $\supp \lift_{p'}\left(\widetilde{a} \right)\in  \widetilde{\sX}'$ is compact. Since any $g' \in 	G\left(\left.\widetilde{\sX}'~\right| \sX \right)$ is a homeomorphism $\supp g'\lift_{p'}\left(\widetilde{a} \right)= g'\supp \lift_{p'}\left(\widetilde{a} \right)$ is also compact so one has
	\be\label{top_gpgp_eqn}
	g'\lift_{p'}\left(\widetilde{a} \right) \in C_c\left( \lift_p\left[C\left(\sX, \left\{A_x\right\}, \sF \right) \right]\right)
	\ee
	If $\widetilde{b} \in C_0\left( \lift_p\left[C\left(\sX, \left\{A_x\right\}, \sF \right) \right]\right)$ then there is a net $\left\{\widetilde{b}_\a\in C_c\left( \lift_p\left[C\left(\sX, \left\{A_x\right\}, \sF \right) \right]\right)\right\}$ such that $\widetilde{b} = \lim_{\a}\widetilde{b}_\a$. Taking into account \eqref{top_gpgp_eqn} one has
	\be\label{top_gpgpg_eqn}
	g'\lift_{p'}\left(\widetilde{b} \right) = \lim_{\a}  g'\lift_{p'}\left(\widetilde{b}_\a \right)\in C_0\left( \lift_p\left[C\left(\sX, \left\{A_x\right\}, \sF \right) \right]\right).
	\ee
	The equations \eqref{top_gpgp_eqn} and \eqref{top_gpgpg_eqn} yield \eqref{top_g_c0_inv_eqn}.\\
	(iii) If  $\widetilde{a}, \widetilde{b} \in \lift_p\left[C\left(\sX, \left\{A_x\right\}, \sF \right) \right]$ correspond to  $\left\{\widetilde{a}_{\widetilde{x}} \in \widetilde{A}_{\widetilde{x}} 
	\right\}_{\widetilde{x} \in \widetilde{\sX}}$ and $\left\{\widetilde{b}_{\widetilde{x}} \in \widetilde{A}_{\widetilde{x}} 
	\right\}_{\widetilde{x} \in \widetilde{\sX}}$ respectively, then the product $\widetilde{a} \widetilde{b}$ corresponds to the families  $\left\{\widetilde{a}_{\widetilde{x}}\widetilde{b}_{\widetilde{x}} \right\}$. For any $g \in G\left(\left.\widetilde{\sX}'~\right| \sX\right)$ elements $g\widetilde{a}$, $g\widetilde{b}$, $\left(g \widetilde{a}\right) \left( g \widetilde{b}\right)$ correspond to
	$\left\{\widetilde{a}_{g\widetilde{x}} 	\right\}$,	$\left\{\widetilde{b}_{g\widetilde{x}} 	\right\}$, $\left\{\widetilde{a}_{g\widetilde{x}}\widetilde{b}_{g\widetilde{x}}\right\}$. Otherwise $g\left(\widetilde{a} \widetilde{b} \right)$ corresponds to  $\left\{\widetilde{a}_{g\widetilde{x}}\widetilde{b}_{g\widetilde{x}}\right\}$, so one has $\left(g \widetilde{a}\right) \left( g \widetilde{b}\right)= g\left(\widetilde{a} \widetilde{b} \right)$. The elements $\widetilde{a}^*$, $g\widetilde{a}^*$ correspond to the families $\left\{\widetilde{a}^*_{\widetilde{x}} \in \widetilde{A}_{\widetilde{x}} 
	\right\}$, $\left\{\widetilde{a}^*_{g\widetilde{x}} \in \widetilde{A}_{\widetilde{x}} 
	\right\}$ and taking into account that $g\widetilde{a}$ corresponds to $\left\{\widetilde{a}_{g\widetilde{x}}\right\}$ we conclude that $\left( g \widetilde{a}\right)^*=  g \widetilde{a}^*$. Thus $g$ is an automorphism of the involutive algebra $\lift_p\left[C\left(\sX, \left\{A_x\right\}, \sF \right) \right]$, and taking into account  (i) of this Lemma we conclude that  $g$ is  an automorphism of the involutive algebras $C_c\left( \lift_p\left[C\left(\sX, \left\{A_x\right\}, \sF \right) \right]\right)$, $C_0\left( \lift_p\left[C\left(\sX, \left\{A_x\right\}, \sF \right) \right]\right)$ and $C_b\left( \lift_p\left[C\left(\sX, \left\{A_x\right\}, \sF \right) \right]\right)$.
\end{proof}

\begin{lemma}\label{top_res_iso_lem}
	Let $\sX$ be a connected, locally compact, Hausdorff space, and let $\sF$ be	continuity structure for $\sX$ {and the} $\left\{A_x\right\}_{x \in \sX}$ (cf. Definition \ref{operator_fields_continuity_defn}) where $A_x$ is a $C^*$-algebra for every $x \in \sX$. Let $p: \widetilde{\sX} \to \sX$ be a covering. Let $\widetilde{\sU}$ be an open set such that the restriction $\left.p\right|_{\widetilde{\sU}}$ is injective and $\sU = p\left(\widetilde{\sU}\right)$. If the closure of $\sU$ is compact, both $\left.C_0\left(\sX, \left\{A_x\right\}, \sF \right)\right|_{\sU}$ and $\left.C_0\left( \lift_p\left[ C\left(\sX, \left\{A_x\right\}, \sF \right)\right] \right) \right|_{\widetilde{\sU}}$ are given by  \eqref{top_res_eqn} then there is a natural injective $*$-isomorphism
	$$
	\left.C_0\left(\sX, \left\{\rep_x\left( A\right) \right\}_{x \in \sX},A \right)\right|_{\sU} \cong \left.C_0\left( \lift_p\left[ C\left(\sX, \left\{\rep_x\left( A\right) \right\}_{x \in \sX},A \right)\right] \right) \right|_{\widetilde{\sU}}.
	$$
\end{lemma}
\begin{proof}
	There are two $\C$-linear maps 
	\bean
	\lift^p_{\widetilde{\sU}}: \left.C_0\left(\sX, \left\{\rep_x\left( A\right) \right\}_{x \in \sX},A \right)\right|_{\sU} \to  \left.C_0\left( \lift_p\left[ C\left(\sX, \left\{\rep_x\left( A\right) \right\}_{x \in \sX},A \right)\right] \right) \right|_{\widetilde{\sU}},\\
	\desc_p:   \left.C_0\left( \lift_p\left[ C\left(\sX, \left\{\rep_x\left( A\right) \right\}_{x \in \sX},A \right)\right] \right) \right|_{\widetilde{\sU}}\to \left.C_0\left(\sX, \left\{\rep_x\left( A\right) \right\}_{x \in \sX},A \right)\right|_{\sU}
	\eean
	and from \eqref{top_lift_desc_eqn} it turns out that the maps are mutually inverse. Moreover from \eqref{comm_lift_desc_hom_eqn} it follows that both maps are $*$-homomorphisms.
\end{proof}

\begin{lemma}\label{top_mult_inc_l_lem}
	Let $\sX$ be a connected, locally compact, Hausdorff space, and let $\sF$ be	continuity structure for $\sX$ {and the} $\left\{A_x\right\}_{x \in \sX}$ (cf. Definition \ref{operator_fields_continuity_defn}) where $A_x$ is a $C^*$-algebra for every $x \in \sX$. If $
	A\bydef C_0 \left(\sX, \left\{\rep_x\left( A\right) \right\}_{x \in \sX},A \right)$ then one has:
	\begin{enumerate}
		\item[(i)]	if $p: \widetilde{\sX} \to \sX$ is a transitive covering then there is the natural injective $*$-homomorphism
		\be\label{top_mult_inc_eqn}
		M\left(A\left( p\right)  \right) : M\left( 	A_0\left(\sX \right)\right)  \hookto M\left( A_0\left(\widetilde{\sX}\right)\right),
		\ee
		\item[(ii)] there is the natural action $G\left(\left.\widetilde{\sX}~\right|\sX\right) \times M\left( A_0\left(\widetilde{\sX}\right)\right) \to M\left( A_0\left(\widetilde{\sX}\right)\right)$,
		\item[(iii)] 	if $p: \widetilde{\sX} \to \sX$ is a transitive finite-fold covering then 	$M\left(A\left( p\right)  \right)$ induces the following $*$-isomorphism
		\be\label{top_mult_iso_eqn}
		M\left( 	A_0\left(\sX \right)\right)  \cong M\left( A_0\left(\widetilde{\sX}\right)\right)^{G\left(\left.\widetilde{\sX}~\right|\sX\right)} .
		\ee
		
	\end{enumerate}
	
\end{lemma}

\begin{proof}
	(i)
	From the Theorem \ref{cross_mult_thm} it turns out than any $a \in M\left(A  \right)$ corresponds to the strictly continuous section $\left\{a_x \in M\left( \rep_x\left( A\right)\right)  \right\}_{x \in \sX}$ (cf. Definition \ref{ctr_crooss_alg_defn}). The section defines a section 
	$$\left\{\widetilde a_{\widetilde x} = a_{p\left(\widetilde x\right)}\in  M\left(\rep_{\widetilde x}\left( \widetilde A\right) \right) = M\left(\rep_{p\left(\widetilde x\right) }\left(A \right)\right) \right\}_{\widetilde x \in \widetilde\sX}$$ . Let $\widetilde{x}_0 \in \widetilde \sX$
	be any point consider an open neighborhood $\widetilde\sU$ of $\widetilde{x}_0$ such that the restriction $\left.p\right|_{\widetilde{\sU}}:\widetilde{\sU}\xrightarrow{\approx}\sU = p\left(\widetilde{\sU} \right)$
	is a injective. If $\widetilde{f}_{\widetilde{x}_0}$ is a $p$-$\widetilde x_0$-stump (cf. Definition  \ref{top_stump_defn}) such that  there is an open neighborhood $\widetilde\sV$ of $\widetilde{x}_0$ such that $\widetilde{f}_{\widetilde{x}_0}\left(\widetilde\sV\right)= 1$. If $\widetilde{c}\in A_0\left(\widetilde{\sX}\right) = C_0\left(\lift_p\left[C\left( \sX, \left\{\rep_x\left( A\right) \right\}, A\right) \right]\right)$, then $\supp \widetilde{f}_{\widetilde{x}_0}\widetilde{c} \subset \widetilde{\sU}$. Let $c = \desc_p \left(\widetilde{f}_{\widetilde{x}_0} \widetilde{c} \right)$ be the $p$-descent (cf. Definition \ref{top_lift_desc_defn}). From the Theorem  \ref{cross_mult_thm} it follows that for any $\eps$ there is an open $\sV'$ neighborhood  of $x_0 = p\left(\widetilde{x}_0\right)$ and an element $b \in \sF$ such that $\left\|c_x \left( a_x - b_x\right)   \right\|+\left\|\left( a_x - b_x\right)c_x    \right\|< \eps$ for every $x$ in $\mathcal V'$. One can suppose $\mathcal V' \subset \sV$. If $f_{x_0}= \desc_p\left( \widetilde{f}_{\widetilde{x}_0}\right)$, $\widetilde a = \lift^p_{\widetilde{\sU}}\left(f_{x_0}a\right)$, $\widetilde b = \lift^p_{\widetilde{\sU}}\left(f_{x_0}b\right)$ (cf. equations \eqref{top_lift_desc_eqn}) then
	\bean
	\forall \widetilde x \in p^{-1}\left(  \mathcal V'\right)  \cap \widetilde{\sU}\quad	\left\|\widetilde c_{\widetilde x} \left(\widetilde a_{\widetilde x}- \widetilde b_{\widetilde x}\right)   \right\|+\left\|\left(\widetilde a_{\widetilde x}- \widetilde b_{\widetilde x}\right) \widetilde c_{\widetilde x}    \right\|=\\=\left\| c_{p(\widetilde x)} \left( a_{p(\widetilde x)}-  b_{p(\widetilde x)}\right)   \right\|+\left\|\left( a_{p(\widetilde x)}-  b_{p(\widetilde x)}\right)  c_{p(\widetilde x)}    \right\| < \eps,
	\eean
	i.e. the section $\left\{\widetilde a_{\widetilde x}\right\}_{\widetilde x\in \widetilde \sX}$ is strictly continuous (cf. Definition \ref{ctr_crooss_alg_defn}). From the  Theorem \ref{cross_mult_thm} it turns out that 
	$$
	\widetilde a \bydef \left\{\widetilde a_{\widetilde x}\right\}\in M\left( C_0\left(\lift_p\left[C\left(\sX, \left\{\rep_x\left( A\right) \right\}_{x \in \sX},A \right) \right]\right)\right) \cong  M\left( A_0\left(\widetilde{\sX}\right)\right).
	$$
	The map $a \mapsto 	\widetilde a$ is the required by this lemma injective $*$-homomorphism
	from  $M\left( 	A_0\left(\sX \right)\right)$ to $M\left( A_0\left(\widetilde{\sX}\right)\right)$.\\
	(ii) If $\widetilde a\in  M\left( A_0\left(\widetilde{\sX}\right)\right)$ represented by the section $\left\{\widetilde a_{\widetilde x}\in M\left( A_{p\left(\widetilde x \right) }\right) \right\}_{\widetilde x \in \widetilde \sX }$ and $g \in G\left(\left.\widetilde{\sX}~\right|\sX\right)$ then we define $g\widetilde a\in  M\left( A_0\left(\widetilde{\sX}\right)\right)$ as a  represented by a strictly continuous section $\left\{\widetilde a_{g\widetilde x}\in M\left( A_{p\left(\widetilde x \right) }\right) \right\}_{\widetilde x \in \widetilde \sX }$ element.
	\\
	(iii) Every  $\widetilde a\in  M\left( A_0\left(\widetilde{\sX}\right)\right)^{G\left(\left.\widetilde{\sX}~\right|\sX\right)}$ represented by a section $\left\{\widetilde a_{\widetilde x}\in M\left( A_{p\left(\widetilde x \right) }\right) \right\}_{\widetilde x \in \widetilde \sX }$ such that $\widetilde a_{\widetilde x }= \widetilde a_{g\widetilde x }$ for each $g \in G\left(\left.\widetilde{\sX}~\right|\sX\right)$. If turns out that  $p\left(\widetilde x' \right)= p\left(\widetilde x'' \right)\Rightarrow \widetilde a_{\widetilde x'}= \widetilde a_{\widetilde x''}$,  hence there is the section  $\left\{ a_{ x}\in M\left( A_{p\left( x \right) }\right) \right\}_{ x \in \sX }$ such that $\widetilde a_{\widetilde x} = a_{p\left(\widetilde x\right)}$ for every $\widetilde x \in \sX$. It turns out
	$$
	M\left( 	A_0\left(\sX \right)\right) \left( a\right) = \widetilde a,
	$$
	i.e. one has the surjective $*$-homomorphism  
	$$
	M\left( 	A_0\left(\sX \right)\right) \to M\left( A_0\left(\widetilde{\sX}\right)\right)^{G\left(\left.\widetilde{\sX}~\right|\sX\right)}.
	$$
	However the $*$-homomorphism  $M\left( 	A_0\left(\sX \right)\right) \to M\left( A_0\left(\widetilde{\sX}\right)\right)$ is injective and $$M\left( A_0\left(\widetilde{\sX}\right)\right)^{G\left(\left.\widetilde{\sX}~\right|\sX\right)}\subset M\left( A_0\left(\widetilde{\sX}\right)\right)$$ hence the map  $M\left( 	A_0\left(\sX \right)\right) \to M\left( A_0\left(\widetilde{\sX}\right)\right)^{G\left(\left.\widetilde{\sX}~\right|\sX\right)}$ is injective, so one has the $*$-isomorphism  $M\left( 	A_0\left(\sX \right)\right)  \cong M\left( A_0\left(\widetilde{\sX}\right)\right)^{G\left(\left.\widetilde{\sX}~\right|\sX\right)}$.
\end{proof}

\subsection{Finite covering functor}\label{top_functor_sec}
\begin{empt}\label{top_cs_funct_empt}
If $A$ is a $C^*$-algebra with Hausdorff spectrum $\sX$ then 
	\bean
	A \cong C_0\left( \sX, \left\{\rep_x  \left(  A\right) \right\}_{x\in \sX}, A \right).
	\eean
	(cf. equation \eqref{top_a_ax_eqn})
	Let $\sX$ be a locally compact,  Hausdorff space called the base space; and for each $x$ in $\sX$, let $A_x$ be a $C^*$-algebra. 
	Let us consider the category  $\mathfrak{FinCov}$-$\sX$ given by the Definition \ref{top_fin_cov_defn}. From \ref{ctr_crooss_alg_empt} it turns out that $C_0\left( \sX, \left\{A_x\right\}, \sF \right)$ is a $C^*$-algebra. If $p$ is a finite-fold covering then from the Lemma \ref{top_lift_c_alg_lem} it follows that $\lift_p$ induces the injective $*$-homomorphism 
	$$
	C_0\left(\sX, \left\{A_x\right\}, \sF \right) \hookto C_0\left(\lift_p\left[\sX, \left\{A_x\right\}, \sF \right] \right) 
	$$
of $C^*$-algebras. Let $\mathfrak{FinCov}$-$\sX$ be the {finite covering category} of $\sX$ (cf. Definition \ref{top_fin_cov_defn}). If $\widetilde \sX$ is $\mathfrak{FinCov}$-$\sX$-object (cf. Definition \ref{top_fin_cov_defn}) then we Denote by
\be\label{top_c0_ob_eqn}
A_0\left(\widetilde \sX  \right) \bydef C_0\left(\lift_p\left[\sX, \left\{A_x\right\}, \sF \right] \right)
\ee
 If $p: \widetilde\sX_1 \to \widetilde\sX_2$ is a $\mathfrak{FinCov}$-$\sX$ morphism, then
		\be\label{top_c0p_ob_eqn}
		\begin{split}
		A_0\left(p \right) \stackrel{\mathrm{def}}{=} \left.\lift_{p}\right|_{A_0\left(\widetilde\sX_2 \right)}	: A_0\left(\widetilde\sX_2 \right) \hookto A_0\left(\widetilde\sX_1 \right) 
		\end{split}
				\ee	
			where $\left.\lift_{p}\right|_{A_0\left(\widetilde\sX_2 \right)}$ is given by the equation \ref{top_fin_lift_eqn}
\end{empt}
\begin{definition}\label{top_functor_c_algebra_defn}
	The described in \ref{top_cs_funct_empt} contravariant functor $A_0$   from the category  the category  $\mathfrak{FinCov}$-$\sX$ to the category of $C^*$- algebras and $*$-homomorphisms is said to by the \textit{finite covering functor associated with} $A$.
\end{definition}
\begin{definition}\label{top_lift_a_f_defn}
Under the hypotheses of the Definition \ref{top_functor_c_algebra_defn} we say that both  the $C^*$-algebra  $A_0\left(\widetilde\sX  \right)$ and the *-homomorphism  is $	A_0\left(p \right)  	: A_0\left(\widetilde\sX_2 \right) \hookto A_0\left(\widetilde\sX_1 \right)$ is  the $p$-\textit{lift} of $A$.
\end{definition}
\begin{definition}\label{top_cs_functa_b_defn}
	$A$ be a $C^*$-algebra with connected, locally connected, locally compact, Hausdorff spectrum $\sX$
	If $p: \widetilde{\sX} \to \sX$ is a transitive covering then we use the following notation for $C^*$-algebras and their injective *-homomorphisms:
	\be\label{top_cb_defna_eqn}
	\begin{split}
		A_c\left(\sX \right)\stackrel{\mathrm{def}}{=} C_c\left(\sX, \left\{\rep_x\left( A\right) \right\}_{x \in \sX},A \right),\\
		A_0\left(\sX \right) \bydef A =C_0\left(\sX, \left\{\rep_x\left( A\right) \right\}_{x \in \sX},A \right),\\
		A_b\left(\sX \right)\stackrel{\mathrm{def}}{=} C_b\left(\sX, \left\{\rep_x\left( A\right) \right\}_{x \in \sX},A \right),
		\\ 
		A_c\left(\widetilde{\sX} \right)\stackrel{\mathrm{def}}{=} C_c\left( \lift_p\left[C\left(\sX, \left\{\rep_x\left( A\right) \right\}_{x \in \sX},A \right)\right] \right),
		\\
		A_0\left(\widetilde{\sX} \right)\bydef C_0\left( \lift_p\left[C\left(\sX, \left\{\rep_x\left( A\right) \right\}_{x \in \sX},A \right)\right] \right),
		\\
		A_b\left(\widetilde{\sX} \right)\stackrel{\mathrm{def}}{=} C_b\left( \lift_p\left[C\left(\sX, \left\{\rep_x\left( A\right) \right\}_{x \in \sX},A \right)\right] \right),
		\\
		A_b\left(p \right) \stackrel{\mathrm{def}}{=} \left.\left.\lift\right._p\right|_{A_0\left(\sX \right)} : A_0\left(\sX \right)  \hookto A_b\left(\widetilde{\sX} \right),
		\\
		A_b\left(p \right) \stackrel{\mathrm{def}}{=} \left.\left.\lift\right._p\right|_{A_b\left(\sX \right)} : A_b\left(\sX \right)  \hookto A_b\left(\widetilde{\sX} \right)
	\end{split}
	\ee
	(cf. (ii) of the Lemma \ref{top_lift_c_alg_lem}). 
\end{definition}

\begin{remark}
	If $A_x$ is a $C^*$-algebra for any $x \in \sX$ then for any $a, b \in C_0\left(\sX, \left\{A_x\right\}, \sF \right)$ and $\widetilde{a}, \widetilde{b} \in C_0\left( \lift_p\left[C\left(\sX, \left\{A_x\right\}, \sF \right)\right]\right)$ one has 
	\begin{equation}\label{comm_lift_desc_hom_eqn}
	\begin{split}
	\supp \widetilde{a} \cap \supp \widetilde{b} \subset \widetilde{\sU} \Rightarrow \mathfrak{desc}_{p}\left( \widetilde{a}\widetilde{b} \right) = \mathfrak{desc}_{p}\left( \widetilde{a} \right)\mathfrak{desc}_{p}\left( \widetilde{b} \right),\\
	\supp a \cap \supp {b} \subset {\sU} \Rightarrow \mathfrak{lift}^p_{\widetilde{\mathcal U}}\left(ab\right)=\mathfrak{lift}^p_{\widetilde{\mathcal U}}\left(a\right)\mathfrak{lift}^p_{\widetilde{\mathcal U}}\left(b\right),\\
	\supp \widetilde{a} \subset \widetilde{\sU} \Rightarrow \mathfrak{desc}_{p}\left( \widetilde{a}^* \right)= \mathfrak{desc}_{p}\left( \widetilde{a}\right)^*,\\
	\supp a \subset {\sU} \Rightarrow \mathfrak{lift}^p_{\widetilde{\mathcal U}}\left(a^*\right)
	=\mathfrak{lift}^p_{\widetilde{\mathcal U}}\left(a\right)^*
	\end{split}
	\end{equation}	
where $\widetilde \sU \subset \widetilde \sX$ is an open subset homeomorphically mapped onto $\sU \bydef p\left(\widetilde \sU\right)$.
\end{remark}
\begin{empt} Consider a special case of 	
	the Definition  \ref{top_functor_c_algebra_defn}
		Let $\sX$ be a locally compact Hausdorff space. If we consider a family $\left\{\C_x\right\}_{x \in \sX}$ of Banach spaces each of which is isomorphic to $\C$ then $C_0\left(\sX\right)$ is  a continuity structure  for $\sX$ {and the} $\left\{\C_x\right\}$. Moreover if $p: \widetilde{\sX} \to \sX$ is a finite fold covering then one has 
\be\label{top_cocc_eqn}
\begin{split}
	C_b\left(\widetilde{\sX}\right) \cong C_b\left(\lift_{p}\left[ C_0\left( \sX, \left\{\C_x\right\},C_0\left(\sX\right)\right)\right]\right),\\ C_0\left(\widetilde{\sX}\right) \cong C_0\left(\lift_{p}\left[ C_0\left( \sX, \left\{\C_x\right\},C_0\left(\sX\right)\right)\right]\right).
\end{split}
\ee
 
\end{empt}
\begin{definition}\label{top_c_funct_defn}    
	Consider the finite covering functor associated with $$A = \left( \sX, \left\{\C_x\right\},C_0\left(\sX\right)\right)$$ (cf. Definition  \ref{top_functor_c_algebra_defn}) from the category  $\mathfrak{FinCov}$-$\sX$ to the category of $C^*$- algebras and $*$-homomorphisms is said to be the \textit{finite covering algebraic functor}. We denote this functor by $C_0$. 
 For any finite fold covering $p: \widetilde \sX \to \sX$ the functor $C_0$ yields an injective $*$-homomorphism
	\be\label{top_c0p_eqn}
	C_0\left(p\right): 	C_0\left(\sX\right)\hookto	C_0\left(\widetilde\sX\right).
	\ee
\end{definition}
\begin{exercise}\label{top_act_exer}
If  $p: \widetilde \sX \to \sX$ is a transitive covering of compact spaces. From the Theorem \ref{pavlov_troisky_thm} it follows that there is an action $C\left(\sX \right) \times C\left(\widetilde \sX \right)\to C\left(\widetilde \sX \right)$. Prove that if the covering is transitive, the action is given by
$$
\left(a,b \right)\mapsto C_0\left(p\right)\left(a\right)b.
$$
\end{exercise}
\subsection{Convergent series}

\begin{lemma}\label{comm_lift_desc_sum_lem}
	Let $\sX$ be,  locally compact, Hausdorff space, and let $\sF$ be	continuity structure for $\sX$ {and the} $\left\{A_x\right\}_{x \in \sX}$ (cf. Definition \ref{operator_fields_continuity_defn}) where $A_x$ is a $C^*$-algebra for every $x \in \sX$. Let 	$p: \widetilde{\mathcal X}\to\mathcal{X}$ be a transitive covering. Suppose that  $\mathcal U\subset \mathcal X$  is a connected open subset evenly {covered} by $\widetilde{\mathcal U}\subset \widetilde{\mathcal X}$. Let $a \in C_c\left(\sX, \left\{\rep_x\left( A\right) \right\}_{x \in \sX},A \right)$ be such that $\supp a \subset \sU$. If 	$\widetilde{a} = \mathfrak{lift}_{\widetilde{\mathcal U}}\left(a \right)$ 
	then following conditions hold:
	\begin{enumerate}
		\item [(i)] The series 	
		\be\nonumber
		\sum_{g \in G\left(\left.\widetilde{\sX}~\right|\sX\right) } g\widetilde{a}.
		\ee 
		is convergent in the strict topology of $M\left(C_0\left( \lift_p\left[\left(\sX, \left\{\rep_x\left( A\right) \right\}_{x \in \sX},A \right)\right]\right)  \right)$ (cf. Definition \ref{strict_topology_defn}),
		\item[(ii)] 
		\be\label{comm_desc_sum_eqn}
		\begin{split}
			\bt\text{-}			\sum_{g \in G\left(\left.\widetilde{\sX}~\right|\sX\right) } g\widetilde{a} = a = \lift_p \circ \desc_p\left( \widetilde{a}\right) \in 
			\\
			\in C_b\left( \lift_p\left[C\left(\sX, \left\{\rep_x\left( A\right) \right\}_{x \in \sX},A \right)\right] \right)\subset  \\ \subset  M\left(C_0\left( \lift_p\left[C\left(\sX, \left\{\rep_x\left( A\right) \right\}_{x \in \sX},A \right)\right] \right)\right).
		\end{split}
		\ee
		where both $ \lift_p$ and $\desc_p$ are given by Definitions \ref{top_lift_main_defn} and \ref{top_lift_desc_defn} respectively.
		
	\end{enumerate}
	
\end{lemma}

\begin{proof} (i)
	If $\widetilde{b} \in C_0\left( \lift_p\left[C\left(\sX, \left\{\rep_x\left( A\right) \right\}_{x \in \sX},A \right)\right] \right)$ and $\eps > 0$ then from the Definition \ref{c_c_compact_defn} it turns out there is a compact set $\widetilde{\mathcal V}\subset \widetilde{\mathcal X}$ such that $\left\| \widetilde{b} \left(\widetilde{x}  \right) \right\| < \frac{\eps}{2\left\|\widetilde{a}\right\|}$ for any $\widetilde{x} \in \widetilde{\mathcal X} \setminus \widetilde{\mathcal V}$. Let us prove that the set
	$$
	G_0=\left\{\left. g \in G\left(\left.\widetilde{\sX}~\right|\sX\right)\right| \widetilde \sV \cap g \widetilde\sU \neq \emptyset \right\}
	$$
	is finite.
	If $G_0$ is not finite then there is an infinite closed set $\widetilde\sX' = \left\{\widetilde x_g \right\}_{g \in G_0}$ such that $\widetilde x_g \in \widetilde \sV \cap g\widetilde \sU$. If $g'\neq g''$ then both $g' \widetilde\sU$ and $g'' \widetilde\sU$ are open neighborhoods of both $\widetilde x_{g'}$ and $\widetilde x_{g''}$ such that $g' \widetilde\sU\cap g'' \widetilde\sU = \emptyset$. It means that the set $\widetilde\sX'$ is discrete. From $\widetilde\sX' \subset \widetilde\sV$ it follows that $\widetilde\sX'$ is compact, however any infinite discrete set is not compact. From this contradiction we conclude that the set  $G_0$ is finite. For any $G' \subset  G\left(\left.\widetilde{\sX}~\right|\sX\right)$ such that $G_0\subset G'$ a following condition holds
	$$
	\left\|\widetilde{b}\left(\sum_{g \in G' } g\widetilde{a}-\sum_{g \in G_0 } g\widetilde{a} \right)  \right\| +
	\left\|\left(\sum_{g \in G' } g\widetilde{a}-\sum_{g \in G_0 } g\widetilde{a} \right)\widetilde{b}  \right\| < \eps.
	$$
	The	above equation means that the series is convergent in the strict topology of $M\left(C_0\left( \lift_p\left[C\left(\sX, \left\{\rep_x\left( A\right) \right\}_{x \in \sX},A \right)\right] \right) \right)$.\\
	(ii)  If $\widetilde{b} \in C_0\left( \lift_p\left[C\left(\sX, \left\{\rep_x\left( A\right) \right\}_{x \in \sX},A \right)\right] \right)$ and $\widetilde{x}_0 \in \widetilde{\mathcal X}$ then $\left( \widetilde{a}\widetilde{b}\right)_{\widetilde{x}_0} = {a}_{p\left(  \widetilde{x}_0\right) } \widetilde{b}_{\widetilde{x}_0}$.  If ${a}_{ p\left(  \widetilde{x}_0\right) }  = 0$ then $\left( \widetilde{a}\widetilde{b}\right)_{\widetilde{x}_0} = 0$. Otherwise there is $\widetilde{x}'_0\in \widetilde{\mathcal U}$ such that $p\left(\widetilde{x}'_0 \right)= p\left(\widetilde{x}_0 \right)$. From \eqref{comm_lift_eqn} it turns out $\widetilde{a}_{\widetilde{x}'_0} = a_{ p\left(\widetilde{x}'_0 \right)}$. The covering $p$ is transitive, so there is the unique $g \in  G\left(\left.\widetilde{\sX}~\right|\sX\right)$ such that  $\widetilde{x}_0 = g\widetilde{x}'_0$. So one has  $\left( g \widetilde{a}\right)_{\widetilde{x}_0}= a_{p\left(\widetilde{x}_0 \right)} $. If $g' \neq g$ then $g\widetilde{\mathcal U} \bigcap \widetilde{\mathcal U} = \emptyset$. It turns out that $\left( g' \widetilde{a}\right)_{\widetilde{x}_0}= 0$ for any $g'\neq g$. Hence one has
	$$
	\left(\sum_{g' \in G\left(\left.\widetilde{\sX}~\right|\sX\right) } g'\widetilde{a} \right)_{\widetilde{x}_0} = 	\left(g\widetilde{a} \right)_{\widetilde{x}_0}= c^{-1}_{\widetilde x_0} \left( {a}_{ p\left(  \widetilde{x}_0\right)} \right) \text{ where } c_{\widetilde x_0}\text{ is given by }\eqref{top_ct_iso_eqn},
	$$  
	it follows that
	$$
	\left( \left(\sum_{g' \in G\left(\left.\widetilde{\sX}~\right|\sX\right) } g'\widetilde{a} \right)\widetilde{b}\right)_{\widetilde{x}_0}= 	\left(\widetilde{b} \left(\sum_{g' \in G\left(\left.\widetilde{\sX}~\right|\sX\right) } g'\widetilde{a} \right)\right)_{\widetilde{x}_0} = \widetilde{b}_{\widetilde{x}_0}  {a}_{ p\left(  \widetilde{x}_0\right)}.
	$$
	The above equation means that 
	\bean
	\sum_{g \in G\left(\left.\widetilde{\sX}~\right|\sX\right) } g\widetilde{a} = a =\lift_p \circ \desc_p\left( \widetilde{a}\right)\in\\\in M\left(C_0\left( \lift_p\left[C\left(\sX, \left\{\rep_x\left( A\right) \right\}_{x \in \sX},A \right)\right]\right)  \right)
	.
	\eean
\end{proof}

\begin{corollary}\label{comm_lift_desc_sum_cor}
	Let $\sX$ be a connected, locally compact, Hausdorff space, and let $\sF$ be	continuity structure for $\sX$ {and the} $\left\{A_x\right\}_{x \in \sX}$ (cf. Definition \ref{operator_fields_continuity_defn}) where $A_x$ is a $C^*$-algebra for every $x \in \sX$. Let 	$p: \widetilde{\mathcal X}\to\mathcal{X}$ be a transitive covering such that the group $G\left(\left.\widetilde{\sX}~\right|\sX\right)$ is residually finite (cf. Definition \ref{residually_finite_defn}).  If $\widetilde a \in C_c\left( \lift_p\left[C\left(\sX, \left\{\rep_x\left( A\right) \right\}_{x \in \sX},A \right)\right] \right)$
then following conditions hold:
\begin{enumerate}
	\item [(i)] The series 	
	\be\label{top_comact_series_eqn}
	\sum_{g \in G\left(\left.\widetilde{\sX}~\right|\sX\right) } g\widetilde{a}.
	\ee 
	is convergent in the strict topology of $M\left(C_0\left( \lift_p\left[\left(\sX, \left\{\rep_x\left( A\right) \right\}_{x \in \sX},A \right)\right]\right)  \right)$ (cf. Definition \ref{strict_topology_defn}),
	\item[(ii)] 
	\be\label{comm_desc_c_sum_eqn}
	\begin{split}
		\bt\text{-}			\sum_{g \in G\left(\left.\widetilde{\sX}~\right|\sX\right) } g\widetilde{a} = \lift_p \circ \desc^c_p\left( \widetilde{a}\right) \in 
		\\
		\in C_b\left( \lift_p\left[C\left(\sX, \left\{\rep_x\left( A\right) \right\}_{x \in \sX},A \right)\right] \right)\subset  \\ \subset  M\left(C_0\left( \lift_p\left[C\left(\sX, \left\{\rep_x\left( A\right) \right\}_{x \in \sX},A \right)\right] \right)\right).
	\end{split}
	\ee
	where both $ \lift_p$ and $\desc^c_p$ are given by Definitions \ref{top_lift_main_defn} and \ref{top_compactly_supported_descent_defn} respectively.
\end{enumerate}
\end{corollary}
\begin{proof}
(i) A set $\supp \widetilde{a}$ (cf. Equation \eqref{top_support_eqn}) is compact. From the 
Theorem \ref{top_compact_thm} 	it follows that there is a finite-fold transitive covering $p': \widetilde \sX' \to \sX$ and  a transitive covering $\widetilde p': \widetilde \sX \to \widetilde \sX'$ such that $p = p' \circ \widetilde p'$ and $\supp \widetilde{a}$ is mapped homeomorphically onto $\widetilde p'\left(\supp \widetilde{a} \right)$. 
From the Lemma \ref{comm_lift_desc_sum_lem} it follows that 
\begin{itemize}
		\item The series 	
\be\nonumber
\sum_{g \in G\left(\left.\widetilde{\sX}~\right|\widetilde{\sX}'\right) } g\widetilde{a}.
\ee 
is convergent in the strict topology of $M\left(C_0\left( \lift_p\left[\left(\sX, \left\{\rep_x\left( A\right) \right\}_{x \in \sX},A \right)\right]\right)  \right)$ (cf. Definition \ref{strict_topology_defn}),
\item
\be\label{comm_desc_cc_sum_eqn}
\begin{split}
	\bt\text{-}			\sum_{g \in G\left(\left.\widetilde{\sX}~\right|\sX\right) } g\widetilde{a} = a = \lift_{\widetilde p'} \circ \desc_{\widetilde p'}\left( \widetilde{a}\right) \in 
	\\
	\in C_b\left( \lift_{\widetilde p'}\left[C\left(\sX, \left\{\rep_x\left( A\right) \right\}_{x \in \sX},A \right)\right] \right)\subset  \\ \subset  M\left(C_0\left( \lift_{\widetilde p'}\left[C\left(\sX, \left\{\rep_x\left( A\right) \right\}_{x \in \sX},A \right)\right] \right)\right).
\end{split}
\ee
where both $ \lift_p$ and $\desc_p$ are given by Definitions \ref{top_lift_main_defn} and \ref{top_lift_desc_defn} respectively.
\end{itemize}
Now one can prove that the series \eqref{top_comact_series_eqn} is convergent to the following finite sum
\be\label{top_desc_p_compact_eqn}
\lift_p\left( \sum_{	g \in G\left( \left. \widetilde \sX'\right| \sX\right) } g \desc_{\widetilde p'} \left(\widetilde a \right)\in C_c\left(\sX, \left\{A_x\right\}, \sF \right)\right) 
\ee
where the strict topology of $M\left(C_0\left( \lift_p\left[\left(\sX, \left\{\rep_x\left( A\right) \right\}_{x \in \sX},A \right)\right]\right)  \right)$ is implied.\\
(ii) Follows from the equations  \eqref{top_desc_compact_eqn} and \eqref{top_desc_p_compact_eqn}.
\end{proof}

\begin{empt}\label{comm_lift_desc_sum_strong_empt}
	Suppose that $\K= \K\left(\H \right) $ is a simple  $C^*$-algebra (cf. Definition \ref{simple_ca_defn}) of compact operators where $\H = \C$ or $\H = \C^n$ or $\H =  \ell^2\left(\N \right)$. If $A$ is  $C^*$-algebra such that the spectrum $\sX$ of $A$ is Hausdorff then $A$ is a $CCR$ algebra (cf. Remark \ref{ctr_open_res_rem}). There is a family $\left\{\rep_x\left(A \right) \right\}_{x\in \sX}$ such that $\rep_x\left(A \right)$ is simple $C^*$-algebra for all $x\in\sX$. From the Lemma \ref{hausdorff_spectrum_lem} it follows that $A$ is  a {continuity structure for} $\sX$ {and the} $\left\{A_x\right\}$
	(cf. Definition \ref{operator_fields_continuity_defn}). Otherwise from the Lemma
	 \ref{oa_haus_alg_lem} it turns out that
	\be\label{top_acs_eqn}
	A \cong  C_0\left(\sX, \left\{A_x\bydef \rep_x\left(A \right)\right\}_{x \in \sX}, A \right)
	\ee
	where $A_x \cong \K\left(\H \right)$ for all $x \in \sX$. If $p: \widetilde \sX \to \sX$ is a covering then  
	\be
	A_0\left(\widetilde\sX \right) \stackrel{\text{def}}{=} C_0\left(\lift_p\left[C_0\left(\sX, \left\{A_x\bydef \rep_x\left(A \right)\right\}_{x \in \sX}, A \right)\right]\right)
	\ee
	$A_0\left(\widetilde\sX \right) \stackrel{\text{def}}{=} C_0\left({\widetilde\sX},  \left\{\widetilde{A}_{\widetilde{x}}\right\}, \widetilde{\sF}\right)$ is a  $C^*$-algebra such that $\widetilde{ \sX}$ is the spectrum of   $A_0\left(\widetilde\sX \right)$.
	If $\widetilde a \in A_0\left(\widetilde\sX \right)$) is represented by the family $\left\{\widetilde a_{\widetilde x}\right\}_{\widetilde x\in \widetilde\sX}$  then one has
	\be
	\begin{split}
		\widetilde a_{\widetilde x}= \rep_{\widetilde x}\left( \widetilde a\right) 
	\end{split}
	\ee
	where $\rep_{ \widetilde{x}}: A_0\left( \widetilde\sX\right) \to B\left(\H_{\widetilde{x}}\right)$ is the irreducible representation which corresponds to  $\widetilde{x}$ (cf. equation \ref{rep_x_eqn}).	If $\pi_a: A_0\left( \widetilde\sX\right) \hookto B\left(\H_a \right)$ is the atomic representation (cf. Definition \ref{atomic_repr_defn}) then $\H_a = \bigoplus_{\widetilde x\in\widetilde\sX}\H_{\widetilde{x}}$ and $\pi_a = \bigoplus_{\widetilde x\in\widetilde\sX}\rep_{ \widetilde{x}}$.
	Otherwise the representation $\pi_a$ can be uniquely extended up to $\pi_a: A_b\left( \widetilde\sX\right) \to B\left(\H_{a}\right)$. 
	From \eqref{top_cb_defna_eqn} it turns out there is the injective $*$-homomorphism $\lift_p:A_0\left(\sX \right) \hookto A_b\left(\widetilde\sX \right)$  so  $A_0\left(\sX \right)$ can be regarded as a subalgebra of $A_b\left(\widetilde\sX \right)$.
\end{empt}
\begin{exercise}\label{comm_lift_desc_sum_exer}
	Consider the situation of the Definition \ref{top_lift_desc_defn}. Prove following statements:
	\begin{enumerate}
		\item The series
		$$
		\sum_{	g \in G\left( \left. \widetilde\sX\right| \sX\right) }g\widetilde a  = \lift_p\circ \desc_p\left(\widetilde a \right).
		$$
		is convergent with respect to the strict topology (cf. Definition \ref{strict_topology_defn}) of $M\left( C_0\left(\widetilde\sX\right)\right)$.
		\item 	If $\widetilde f\in C_0\left(\widetilde\sX\right)_+$ a positive element is such that $\supp  \widetilde f\subset \widetilde \sU$ then one has
		\be\label{top_desc_fin_eqn}
		\widetilde f= \sqrt{ \widetilde f}\sqrt{ \desc_p\left( \widetilde f\right) }.
		\ee
		\item 
		If $f \in C_0\left(\sX \right)$ is such that $\supp f \subset \sU$ then for all $\xi \in C_b\left(\sX, \left\{A_x\right\}, \sF \right)$
		one has
		\be\label{top_lift_pa_eqn}
		\lift^p_{\widetilde{\sU}}\left(f \xi \right)=  \lift^p_{\widetilde{\sU}}\left(f \right)\lift_p\left( \xi \right)
		\ee
		
	\end{enumerate}
\end{exercise}
\
\begin{exercise}
		Consider the situation \ref{comm_lift_desc_sum_strong_empt}.  Let $\widetilde a \in A_c\left(\widetilde\sX \right)$ (cf. \eqref{top_cb_defna_eqn}). Prove following statements:
		\begin{enumerate}
		\item [(i)] the series 	
\bean
\sum_{g \in G\left(\left.\widetilde{\sX}~\right|\sX\right) } \pi_a\left( g\widetilde{a}\right),
\eean 
is convergent in the strong  topology of $B\left(\H_a\right)$  (cf. Definition \ref{strong_topology_defn} ),
\item[(ii)] 
\bean
\begin{split}
	\sum_{g \in G\left(\left.\widetilde{\sX}~\right|\sX\right) } \pi_a\left( g\widetilde{a}\right)  = \pi_a\left(  \desc^c_p\left( \widetilde{a}\right)\right)
\end{split}
\eean
where $\desc^c_p$  is a compactly supported $p$-descent (cf. Definition \ref{desc_c_sheaf_defn}). 
		\end{enumerate}
\end{exercise}

\section{Sheaves lift and descent}
\subsection{Inverse images lift and descent}
\paragraph*{}
If we consider the situation of the Definition \ref{sheaf_inv_im_defn} then the {inverse image}  $f^{-1}\mathscr G$ on $\sX$ is the sheaf  associated to the presheaf  $\sU \mapsto \lim_{\sV \supseteq f\left(\sU\right)} \mathscr G\left(\sV\right)$, where $\sU$ is any open set in $\sX$, and the limit is taken over all open sets $\sV$ of $\sX$ containing $f\left(\sU\right)$. If the map $f$ is surjective then one has  $\lim_{\sV \supseteq f\left(\sX\right)} \mathscr G\left(\sV\right)\cong \mathscr  G\left(\sX\right)$ it follows that there is the injective homomorphism 
\be\label{top_cs_sheaf_inc_eqn}
f^{-1}_{\mathscr  G}:\mathscr  G\left(\sX\right)\hookto f^{-1}\mathscr  G\left(\sY\right)
\ee
of Abelian groups.

\begin{exercise}\label{top_sheaf_lift_exer}
	Let $p:\widetilde \sX\to \sX$ be a covering (cf. Definition \ref{top_covering_defn}), and let $\mathscr F$ be a sheaf on $\sX$. Prove following statements:
	\begin{enumerate}
		\item For all $\widetilde x\in\widetilde\sX$ there is a natural isomorphism of stalks (cf. Definition \ref{sheaf_stalk_defn})
		\be\label{top_st_iso_eqn}
		p_{\widetilde x}:	\mathscr  F_{p\left( \widetilde x\right) }\cong p^{-1}\mathscr  F_{\widetilde x}.
		\ee
		\item The given by \eqref{top_cs_sheaf_inc_eqn} homomorphism
		$$
		p^{-1}_{\mathscr  F}:\mathscr  F\left(\sX\right)\hookto p^{-1}\mathscr  F\left(\widetilde \sX\right)
		$$ can be described by the following way. 
		Elements $\xi \in \mathscr  F\left(\sX\right)$ and $\widetilde \xi\in  p^{-1}\mathscr  F\left(\widetilde \sX\right)$ can be represented by global continuous sections $s : \sX \to \mathrm{Sp\acute{e}}\left(\mathscr F \right)$ and $\widetilde s : \widetilde\sX \to \mathrm{Sp\acute{e}}\left(p^{-1}\mathscr F \right)$ of les \'espaces etal\'e (cf. Exercise \ref{sheaf_etale_exer}). We put
		\be\label{top_shaf_cx_eqn}
		\widetilde \xi = p^{-1}_{\mathscr  F}\xi \quad \Leftrightarrow\quad  \forall \widetilde x \in \widetilde \sX \quad p_{\widetilde x}\left( \widetilde s \left(\widetilde x\right)\right) = s\left({p\left( \widetilde x\right) }\right).
		\ee
		\item If both $p': \sX'\to \sX$ and $p'' :\sX''\to \sX'$ are coverings
		then one has
		\bea\label{top_sh1_comp_eqn}
		(p'\circ p'')^{-1}\mathscr  F = p''^{-1}\left( p'^{-1}\mathscr  F\right),\quad\\
		\label{top_sh2_comp_eqn}  
		(p''\circ p')_{\widetilde x}= p''_{\widetilde x}\circ p'_{p''\left( \widetilde x\right) }~~:\quad \mathscr  F_{(p''\circ p')\left( \widetilde x\right) }\cong (p''\circ p')^{-1}\mathscr  F_{\widetilde x}\quad \text{(cf. \eqref{top_st_iso_eqn})}\quad
		\eea
		
	\end{enumerate}
	
\end{exercise}
\begin{definition}\label{top_sheaf_lift_defn}
	Under the hypotheses of the Exercise  \ref{top_sheaf_lift_exer} we say that $\widetilde{\xi}$ is the $p$-\textit{lift} of $\xi$. We also use the following notation
	\be\label{top_sheaf_lift_eqn}
	\begin{split}
		\lift_p:\mathscr  F\left(\sX\right)\hookto p^{-1}\mathscr  F\left(\widetilde \sX\right),\\
		\widetilde{\xi} = \lift_p\left( \xi\right) .
	\end{split}
	\ee
\end{definition}

\begin{remark}
	The Definition \ref{top_sheaf_lift_defn} complies with the equation \eqref{top_glo_lift_h_eqn}.
\end{remark}
\begin{remark}\label{top_sheaf_lift_rem}
	Sometimes we write $\lift_p^{\text{sheaf}}$ for distinguishing  this lift with the given by  \eqref{top_glo_lift_h_eqn} one (cf. the Remark \ref{top_top_lift_rem}) and the Exercise \ref{top_sheaf_corr_exer} below).
\end{remark}
\begin{remark}
	If both $p': \sX'\to \sX$ and $p'' :\sX''\to \sX'$ are coverings then from \eqref{top_sh2_comp_eqn} it follows that
	\be\label{top_sh_lift_comp_eqn}
	\lift_{p''\circ p'} = \lift_{p''} \circ \lift_{p'},
	\ee
	or equivalently
		\be\label{top_sh_inv_comp_eqn}
	\left( p''\circ p'\right)^{-1}  = p''^{-1}\circ p'^{-1}.
	\ee
\end{remark}
Following definition is an analog of the Definition \ref{top_lift_desc_defn}.
\begin{definition}\label{top_lift_sh_desc_defn}
	Let $\mathscr F$ 	be a sheaf on  $\sX$ (cf. Definition \ref{sheaf_defn}), and let $p: \widetilde{\sX}\to\sX$ be a covering. Let $\widetilde{\sU}$ be an open subset such that the restriction $\left.p\right|_{\widetilde{\sU}}$ is injective. Let $\sU\bydef p\left( \widetilde{\sU}\right)$. If $a \in\mathscr F\left(\sX\right)$ is such that $\supp a \subset \sU$ and $\widetilde{a}'= \lift_p\left( a\right)$. Suppose that $\widetilde{a}'$ is represented by a continuous section $\widetilde s' : \widetilde\sX \to \mathrm{Sp\acute{e}}\left(p^{-1}\mathscr F \right)$ (cf. Exercise \ref{sheaf_etale_exer}). There is a continuous section $\widetilde s : \widetilde\sX \to \mathrm{Sp\acute{e}}\left(p^{-1}\mathscr F \right)$ such that
	\be\label{comm_sh_lift_desc_eqn}
	\widetilde{s}\left( \widetilde{x}\right) =\left\{\begin{array}{c l}
		\widetilde{s}'\left( \widetilde{x}\right)    &{\widetilde{x}}\in \widetilde{\sU} \\
		0 &{\widetilde{x}}\notin \widetilde{\sU}\\
	\end{array}\right.
	\ee
	Corresponding element $\widetilde \xi \in  p^{-1}\mathscr F\left(\widetilde\sX\right)$
	is said to be the $p$-$\widetilde{\sU}$-\textit{lift} or simply the $\widetilde{\sU}$-\textit{lift} of $\xi$. Otherwise we say that $\xi$ is the $p$-\textit{descent} of $ \widetilde{\xi}$. We write
	\be\label{top_sh_lift_desc_eqn}
	\begin{split}
		\widetilde{\xi}\stackrel{\text{def}}{=} \lift^p_{\widetilde{\sU}}\left(\xi \right)  \text{ or simply } \widetilde{\xi}\stackrel{\text{def}}{=} \lift_{\widetilde{\sU}}\left(\xi \right), \\
		\xi \stackrel{\text{def}}{=} \desc_{p} \left(\widetilde{\xi}  \right)  \text{ or simply } \xi \stackrel{\text{def}}{=} \desc \left(\widetilde{\xi}  \right).
	\end{split}
	\ee
\end{definition}
\begin{remark}
	if $\xi = \desc_{p} \left(\widetilde{\xi}  \right)$ then $\xi$ is represented by the section  $s: \sX \to \mathrm{Sp\acute{e}}\left(\mathscr F \right)$ such that
	\be\label{top_sh_desc_eqn}
	s_x = \begin{cases}
		p^{-1}_{\widetilde x}\left( 	\widetilde{s}_{\widetilde{x}}\right)  & x \in \sU ~\text{AND} ~ \widetilde x \in \widetilde\sU~\text{AND} ~p\left( \widetilde{x}\right) = x \\
		0 &  x \notin \sU
	\end{cases}	\quad 
	\ee 
	where	$p_{\widetilde x}$ \text{ is given by }\eqref{top_st_iso_eqn} and the continuous section 
	$\widetilde s: \widetilde \sX\to  \mathrm{Sp\acute{e}}\left(p^{-1}\mathscr F\right) \left( \widetilde \sX \right)$ represents $\widetilde \xi$.
\end{remark}
\begin{remark}
	From the equation \eqref{top_shaf_cx_eqn} it turns out that
	\begin{equation}\label{comm_sh_lift_desc_l_eqn}
		\begin{split}
			\xi = \mathfrak{desc}_{p}\left( \mathfrak{lift}^p_{\widetilde{\mathcal U}}\left(\xi \right)\right) ,\\
			\widetilde{\xi} = \mathfrak{lift}^p_{\widetilde{\mathcal U}}\left( \mathfrak{desc}_{p}\left(\widetilde{\xi} \right)\right). 
		\end{split}
	\end{equation}	
	The system of equations \eqref{comm_sh_lift_desc_l_eqn} is an analog of \eqref{comm_lift_desc_l_eqn} one.
\end{remark}

\begin{exercise}
	Consider the situation of the Definition \ref{top_lift_sh_desc_defn}. Suppose that  $ \mathscr End\left(\mathscr F \right) \bydef \mathscr Hom \left(\mathscr F,\mathscr F \right)$ is the {sheaf of local morphisms} (cf. Definition \ref{sheaf_hom_defn}).
	Let $p: \widetilde \sX\to\sX$ be a covering, and let  $\widetilde\sU\subset\widetilde \sX$ be an open subset homeomorphically mapped onto $\sU \bydef p\left(\widetilde\sU\right)$. 
	Prove that
	\be\label{top_ldx_eqn}
	\begin{split}
		\forall D \in \mathscr End\left(\mathscr F \right)\quad \forall a \in \mathscr End\left(\mathscr F \right)\quad \supp a\subset \sU\quad\Rightarrow\quad \lift^p_{\widetilde\sU}\left( \left[ D, a\right]\right) = \\= \left[ \lift_p\left(D \right) ,\lift^p_{\widetilde \sU}(a)\right].
	\end{split}
	\ee
\end{exercise}

\subsection{Modules and sheaves}
\paragraph*{}
	Let $\sX$ be a locally compact,  Hausdorff space, and let $R \subset C_0\left(\sX\right)$ be a dense *-subalgebra. Let $A$ be a left  $R$-module such that
	$$
\forall a \in A \quad a = 0 \quad \Leftrightarrow \quad \left\{r \in R | ra = 0\right\}	= \{0\}.
	$$
	For  every open subset $\sU\subset\sX$ we define  an equivalence relation $\sim_\sU$ (cf. Definition \ref{equivalence_relation_defn}) on $A$ such that
\be\label{top_sheaf_e_eqn}
	a' \sim_\sU a'' \quad \Leftrightarrow\quad \left( \forall r \in R \quad \supp r \subset \sU \quad \Rightarrow\quad ra' = r a''\right). 
\ee
From the above equation it turns out that
$$
\sV \subset \sU \quad \Rightarrow \quad \left( a' \sim_\sU a'' \quad \Rightarrow\quad  a' \sim_\sV a''\right).
$$
It follows that one has a natural projection
\be\label{top_puv_eqn}
\rho_{\sU \sV}:A/\sim_\sU \to A/\sim_\sV.
\ee
One can proof maps $\rho_{\sU \sV}$  satisfy to the Definition \ref{presheaf_defn}, i.e. we have a presheaf  on $\sX$.
\begin{definition}\label{top_x_sheaf_defn}
A sheaf associated with the described above presheaf (cf. Proposition and Definition \ref{sheaf_prdf}) is said to be the $A$-\textit{sheaf}. Henceforth we denote the $A$-{sheaf} by $\mathscr S^A$.
\end{definition}
Any $s \in A$ corresponds to a global section of $\mathscr S^A$, it follows that there is an inclusion of Abelian groups
\be\label{top_x_sheaf_eqn}
A\subset \mathscr S^A\left( \sX\right). 
\ee

\begin{remark}
Since $R$ is an $R$-module, there is an $R$-sheaf $\mathscr S^R$.
One has a ringed space $\left(\sX, \mathscr S^R\right)$ (cf. Definition \ref{sheaf_ringed_space_defn}).
\end{remark}
Let $R\subset C_0\left(\sX\right)$ be a dense *-subalgebra. Similarly to the equation \eqref{unital_notation_eqn} one has an unital *-algebra $R^\sim$ which is an $R$-module. There is a natural inclusion $R^\sim \subset C_b\left(\sX\right)$.
\begin{definition}\label{top_soft_r_defn}
	In the above situation we suppose that the space $\sX$ is paracompact (cf. Definition \ref{top_paracompact_defn}).
If 	$\Phi$   is a  {paracompactifying} family of supports (cf. Definition \ref{phi_supp_defn}) 	
	then we say that a dense *-algebra $R\subset C_0\left(\sX\right)$ is $\Phi$-\textit{soft} if following conditions hold:
	\begin{enumerate}
		\item [(a)] a sheaf $\mathscr S^{R^\sim}$ is $\Phi$-\textit{soft} (cf. Definition \ref{soft_sheaf_defn}),
		\item [(b)] one has $R= C_0\left(\sX\right) \cap \mathscr S^{R^\sim}\left( \sX\right) $.
	\end{enumerate}
\end{definition}
\begin{exercise}\label{top_soft_exer}
Let $R$ be a $c$-soft *-algebra (cf. Definitions \ref{top_soft_r_defn} and \ref{soft_sheaf_defn}). Prove that for any $x\in \sX$ and any  neighborhood $\sU$ of $x$
there is a map $f \in R\cap C_c\left(\sX\right)$ such that
\begin{enumerate}
	\item $\supp f \subset \sU$.
	\item $f\left(\sX\right)= \left[0,1\right]$.
	\item There is an open neighborhood $\sV$ such that $f\left(\sV\right)= \left[1\right]$.
\end{enumerate}

\end{exercise}
\begin{definition}\label{top_stump_soft_defn}
Under the hypotheses of the Exercise \ref{top_soft_exer}  we say that $f$ is an $\left(R, x\right)$-\textit{stump}.
\end{definition}

\begin{exercise}\label{top_stump_soft_y_exer}
Let $R$ be a $c$-soft *-algebra (cf. Definitions \ref{top_soft_r_defn}. Prove that  if $\sY \subset \sX$ is a compact subset then there is $f \in R \cap C_c\left( \sX\right)$ such that $f\left(\sX \right) = \left[0, 1\right]$ and $f\left(\sY \right) = \{1\}$ (cf. Theorem \ref{sheaf_a_soft_thm}).
\end{exercise}
\begin{definition}\label{top_stump_soft_y_defn}
	Under the hypotheses of the Exercise \ref{top_stump_soft_y_exer} we say that $f$ is an $\left(R, \sY\right)$-\textit{stump}.
\end{definition}
\begin{remark}\label{top_stump_soft_y_rem}
Under the hypotheses \ref{top_stump_soft_y_defn} for any open set $\sU\subset\sX$ such that $\sY \subset \sU$ there is an $\left(R, \sY\right)$-{stump} $f$ such that $\supp f \subset\sU$.
\end{remark}
\begin{exercise}\label{top_soft_c_exer}
Let $R \subset C_0\left(\sX\right)$ be a $c$-soft *-subalgebra.
Prove that if $p: \widetilde \sX \to \sX$ is a covering then an intersection $p^{-1}\mathscr S^{R}\left( \widetilde\sX\right) \cap C_0\left(\widetilde \sX\right)$ is a $c$-soft *-subalgebra of $C_0\left(\widetilde \sX\right)$.
\end{exercise}
\begin{exercise}\label{top_soft_p_exer}
	If $R \subset C_0\left(\sX\right)$ is a $c$-soft *-subalgebra and $ p: \widetilde \sX \to \sX$ is a covering then for all $\widetilde x \in \widetilde\sX$ there  is $\widetilde f \in C_c\left( \widetilde \sX\right)\cap p^{-1}\mathscr S^{R}\left( \widetilde\sX\right)$ such that following conditions hold:
\begin{enumerate}
	\item[(a)] there is an open subset $\widetilde \sU\subset\widetilde\sX$ such that $\supp \widetilde f\subset \widetilde\sU$ and $\widetilde\sU$ is mapped homeomorphically onto $p\left(\widetilde\sU\right)$,
	\item[(b)] $\widetilde f\left(\sX \right) = \left[0, 1\right]$,
	\item[(c)] there is an open neighborhood $\widetilde \sV$ of $\widetilde x$ such that $\widetilde f\left(\widetilde \sV \right) = \{1\}$.
\end{enumerate}
\end{exercise}

\begin{definition}\label{top_stump_soft_p_defn}
	Under the hypotheses of the Exercise \ref{top_soft_p_exer} we say that $f$ is an $\left(p, R, \widetilde x\right)$-\textit{stump}.
\end{definition}

\begin{exercise}\label{smooth_soft_exer}
Let $M$ be a smooth paracompact manifold. Using the Proposition \ref{top_smooth_part_unity_prop} prove that $\Coo_0\left(M\right)\bydef \Coo\left(M\right)\cap C_0\left(M\right)\subset C_0\left(M\right)$ is a $c$-soft *-algebra (cf. Definition \ref{top_soft_r_defn}).
\end{exercise}
\subsection{Sheaves and continuity structures}
\paragraph{}
	Let $\sX$ be a locally compact,  Hausdorff space.	
	If $\sF$ is a continuity structure  for  $\sX$ {and the} $\left\{A_x\right\}_{x \in \sX}$ (cf. Definition \ref{operator_fields_continuity_defn}) then from the Lemma \ref{op_cont_module_lem} one has following $C_0\left(\sX \right)$-modules
	\be\label{op_cont_module_eqn}
	\begin{split}
	C\left({\sX}, \left\{{A}_{{x}}\right\}, {\sF}\right),\\
C_b\left({\sX}, \left\{{A}_{{x}}\right\}, {\sF}\right),\\
C_0\left({\sX}, \left\{{A}_{{x}}\right\}, {\sF}\right),\\
C_c\left({\sX}, \left\{{A}_{{x}}\right\}, {\sF}\right).
\end{split}
	\ee
	(cf. equation \eqref{top_c_sec_eqn} and Definition \ref{top_s_bounded_defn}, \ref{top_cc_c0_defn}). If the space $\sX$ is paracompact, then $C_0\left( \sX\right)$ is a $c$-soft *-subalgebra of  $C_0\left( \sX\right)$  (cf. Definitions \ref{top_soft_r_defn}). We leave to the reader the proof of following equations:
		\bea\label{top_cont_sh_eqn}
		\mathscr S^{C\left({\sX}, \left\{{A}_{{x}}\right\}, {\sF}\right)}= 	\mathscr S^{C_b\left({\sX}, \left\{{A}_{{x}}\right\}, {\sF}\right)}=
		\mathscr S^{C_0\left({\sX}, \left\{{A}_{{x}}\right\}, {\sF}\right)}= 	\mathscr S^{C_c\left({\sX}, \left\{{A}_{{x}}\right\}, {\sF}\right)},\\
		\label{top_cont_shc_eqn}
		\mathscr S^{C\left({\sX}, \left\{{A}_{{x}}\right\}, {\sF}\right)}\left(\sX\right)=  	C\left({\sX}, \left\{{A}_{{x}}\right\}, {\sF}\right)		
		\eea 
where the notation of the Definition \ref{top_x_sheaf_defn} is used.

\begin{lemma}\label{top_sheaf_lift_lem}
	
	Let $\sX$ be a locally compact,  Hausdorff space.
	Let  $\sF$ be a continuity structure  for  $\sX$ {and the} $\left\{A_x\right\}_{x \in \sX}$. If $p: \widetilde{\sX} \to \sX$ is  a covering and 
	$$
	C\left(\widetilde{\sX}, \left\{\widetilde{A}_{\widetilde{x}}\right\},~ \widetilde{\sF} \right) = \lift_p\left[ C\left(\sX, \left\{A_x\right\}, ~\sF \right) \right]
	$$
	is the $p$-lift of $C\left(\sX, \left\{A_x\right\}, ~\sF \right)$ (cf. Definition \ref{top_lift_main_defn}) then the ${C\left(\widetilde{\sX}, \left\{\widetilde{A}_{\widetilde{x}}\right\},~ \widetilde{\sF} \right)}$-sheaf $\mathscr S^{C\left(\widetilde{\sX}, \left\{\widetilde{A}_{\widetilde{x}}\right\},~ \widetilde{\sF} \right)}$ (cf. Definition \ref{top_x_sheaf_defn}) is the $p$-inverse image (cf. Definition \ref{sheaf_inv_im_defn}) of the ${C\left({\sX}, \left\{{A}_{{x}}\right\}, ~{\sF} \right)}$-sheaf $\mathscr S^{C\left({\sX}, \left\{{A}_{{x}}\right\}, ~{\sF} \right)}$   i.e.
	$$
	\mathscr S^{C\left(\widetilde{\sX}, \left\{\widetilde{A}_{\widetilde{x}}\right\}, \widetilde{\sF} \right)}= p^{-1}\mathscr S^{C\left({\sX}, \left\{{A}_{{x}}\right\}, {\sF} \right)}.
	$$
\end{lemma}
\begin{proof}
	From the Definition \ref{sheaf_inv_im_defn} it follows that $p^{-1}\left( \mathscr  C\left(\sX, \left\{A_x\right\}, \sF \right)\right)$ is the associated to the presheaf  $\mathscr  G$ given by $\widetilde \sU \mapsto \lim_{\sU\supseteq p\left(\widetilde \sU\right)} \mathscr  S^{\left(\sX, \left\{A_x\right\}, \sF \right)}\left( \sU\right) $, where $\widetilde \sU$ is any open set in $\widetilde \sX$, and the limit is taken over all open sets $\sU$ of $\sX$ containing $p\left(\widetilde \sU\right)$. In particular if $\widetilde \sU$ is homeomorphically mapped onto $p\left(\widetilde \sU\right)$ then $\mathscr  G\left( \widetilde \sU\right)\cong \mathscr  S^{\left(\sX, \left\{A_x\right\}, \sF \right)}\left(p\left(\widetilde \sU\right) \right)$. Let $\widetilde s \in p^{-1}\mathscr  S^{\left(\sX, \left\{A_x\right\}, \sF \right)}\left(\widetilde \sU \right) $.
	From (2) of \ref{sheaf_empt} it follows that   for any $\widetilde x \in \widetilde\sU$ there is an open neighborhood  $\widetilde\sV$ of $\widetilde x$ such that $\widetilde\sV \subset \widetilde\sU$ and the element $\widetilde t \in \mathscr  S^{\left(\sX, \left\{A_x\right\}, \sF \right)}$ such that  for all $\widetilde y \in \widetilde\sV$ the germ $\widetilde t_{\widetilde y}$ of $\widetilde t$ at $\widetilde y$ is equal to $\widetilde s\left(\widetilde y\right)$. For any $\widetilde x$ we select an  neighborhood $\widetilde\sV_{\widetilde x}$ of $\widetilde x$ such that   $\widetilde\sV_{\widetilde x}$ satisfies to (2) of \ref{sheaf_empt} and $\widetilde \sV_{\widetilde x}$ is homeomorphically mapped onto $p\left(\widetilde \sV_{\widetilde x}\right)$. From  $\mathscr  G\left( \widetilde \sV_x\right)\cong \mathscr  S^{\left(\sX, \left\{A_x\right\}, \sF \right)}\left(p\left(\widetilde \sV_x\right) \right)$ it follows that the restriction $\left.\widetilde s\right|_{\widetilde \sV_x}$ corresponds to element $s_{\widetilde x} \in \mathscr  S^{\left(\sX, \left\{A_x\right\}, \sF \right)}\left( p\left(\widetilde \sV_x\right)\right)$. Otherwise $s_{\widetilde x}$ corresponds to a {continuous (with respect to $\sF$)} family $\left\{a_y \in A_y \right\}_{y \in p\left(\widetilde \sV_x\right) }$. So if $c_{\widetilde y}: \widetilde A_{\widetilde y}\to A_{p\left(\widetilde y \right) }$ is given by \eqref{top_ct_iso_eqn} 	for any $\widetilde y \in \widetilde \sV_x$  then there is  {continuous (with respect to $\widetilde \sF$)}
	family $\left\{\widetilde a_{\widetilde y} = c^{-1}\left(a_{p\left( \widetilde y\right) } \right)  \right\}_{\widetilde y\in \widetilde\sV_{\widetilde x} }$. Let $ \widetilde x',  \widetilde x'' \in  \widetilde\sU$. If both  $\left\{\widetilde a'_{\widetilde y'}  \right\}_{\widetilde y'\in \widetilde\sV_{\widetilde x'} }$  and  $\left\{\widetilde a''_{\widetilde y''}  \right\}_{\widetilde y''\in \widetilde\sV_{\widetilde x''} }$ are families corresponding to both $\left.\widetilde s\right|_{\widetilde \sV_{x'}}$ and $\left.\widetilde s\right|_{\widetilde \sV_{x''}}$ 
	then from $\left.\left.\widetilde s\right|_{\widetilde \sV_{x'}}\right|_{\widetilde \sV_{x'}\cap \widetilde \sV_{x''} }=\left.\left.\widetilde s\right|_{\widetilde \sV_{x''}}\right|_{\widetilde \sV_{x'}\cap \widetilde \sV_{x''} }$ it follows that $\widetilde a'_{\widetilde y}= \widetilde a''_{\widetilde y}$ for all $\widetilde y \in  \widetilde \sV_{x'}\cap \widetilde \sV_{x''}$. So there is the single family $\left\{\widetilde a_{\widetilde y}  \right\}_{\widetilde y\in \widetilde\sU}$ which is the combination of  {continuous (with respect to $\widetilde \sF$)} families $\left\{\widetilde a_{\widetilde y}  \right\}_{\widetilde y\in \widetilde\sV_{\widetilde x} }$ 
	The family $\left\{\widetilde a_{\widetilde y}  \right\}_{\widetilde y\in \widetilde\sU}$ is {continuous (with respect to $\widetilde \sF$)} so it corresponds the element of $\mathscr S^{C \left(\widetilde{\sX}, \left\{\widetilde{A}_{\widetilde{x}}\right\}, \widetilde{\sF} \right)}\left(\widetilde\sU \right) $.
\end{proof}

\begin{exercise}\label{top_top_shift_iso_exer}
	Under the hypotheses of the Lemma \ref{top_sheaf_lift_lem} one has
	isomorphisms
	\be\label{top_sh_cong_eqn}
	\begin{split}
			\mathscr S^{C\left({\sX}, \left\{{A}_{{x}}\right\}, {\sF}\right)}\left(\sX\right)\cong  	C\left({\sX}, \left\{{A}_{{x}}\right\}, {\sF}\right),	\\	
\mathscr S^{C\left({\widetilde \sX}, \left\{{\widetilde A}_{{\widetilde x}}\right\}, {\widetilde \sF}\right)}\left(\widetilde \sX\right)\cong  	C\left({\widetilde\sX}, \left\{{A}_{{\widetilde x}}\right\}, {\widetilde\sF}\right),		
	\end{split}
	\ee
	(cf. \eqref{top_cont_shc_eqn}). On the other hand the   equations \eqref{top_glo_lift_h_eqn} and   \eqref{top_cs_sheaf_inc_eqn} yield injective $\C$-linear maps
	\bean
	\lift_p^{\text{sheaf}}: \mathscr S^{C\left({\sX}, \left\{{A}_{{x}}\right\}, {\sF}\right)}\left(\sX\right) \hookto\mathscr S^{C\left({\widetilde \sX}, \left\{{\widetilde A}_{{\widetilde x}}\right\}, {\widetilde \sF}\right)}\left(\widetilde \sX\right),\\
\lift_p^{\text{top}}:	C\left({\sX}, \left\{{A}_{{x}}\right\}, {\sF}\right)\hookto C\left({\widetilde\sX}, \left\{{A}_{{\widetilde x}}\right\}, {\widetilde\sF}\right)
	\eean
	(cf. Remarks \ref{top_top_lift_rem}, \ref{top_sheaf_lift_rem}).
Prove that the following diagram
	\newline
\begin{tikzpicture}
	\matrix (m) [matrix of math nodes,row sep=3em,column sep=4em,minimum width=2em]
	{
	\mathscr S^{C\left({\sX}, \left\{{A}_{{x}}\right\}, {\sF}\right)}\left(\sX\right)  & C\left({\sX}, \left\{{A}_{{x}}\right\}, {\sF}\right) \\		
p^{-1}\mathscr S^{C\left({\sX}, \left\{{A}_{{x}}\right\}, {\sF}\right)}\left(\sX\right)=	\mathscr S^{C\left({\widetilde \sX}, \left\{{\widetilde A}_{{\widetilde x}}\right\}, {\widetilde \sF}\right)}\left(\widetilde \sX\right)	&  C\left({\widetilde\sX}, \left\{{A}_{{\widetilde x}}\right\}, {\widetilde\sF}\right) \\};
	\path[-stealth]
	(m-1-1) edge node [above] {$\cong$} (m-1-2)
	(m-2-1) edge node [above]  {$\cong$} (m-2-2)
	(m-1-1) edge node [right] {$\lift_p^{\text{sheaf}}$} (m-2-1)
	(m-1-2) edge node [right] {$\lift_p^{\text{top}}$} (m-2-2);
\end{tikzpicture}
\\
is commutative.
\end{exercise}

\begin{exercise}\label{top_sheaf_corr_exer}
	Under the hypotheses of the Exercise \ref{top_top_shift_iso_exer} suppose that the space $\sX$ is paracompact (cf. Definition \ref{top_paracompact_defn}). Let  $\Phi$   is a  {paracompactifying} family of supports (cf. Definition \ref{phi_supp_defn}). Let $R \subset C_0\left(\sX\right)$ be a $\Phi$-soft *-subalgebra, and let $A \subset 	\mathscr S^{C\left({\sX}, \left\{{A}_{{x}}\right\}, {\sF}\right)}\left(\sX\right)$ be an $R$-submodule. If $\mathscr S^{A}$ is $A$-sheaf (cf. Definition \ref{top_x_sheaf_defn}), then from \eqref{top_sh_cong_eqn} one has  natural inclusions
	\bean
	\mathscr S^{A}\left(\sX\right) \subset \mathscr S^{C\left({\sX}, \left\{{A}_{{x}}\right\}, {\sF}\right)}\left(\sX\right),\\
	p^{-1}\mathscr S^{A}\left(\widetilde\sX\right) \subset  C\left({\widetilde\sX}, \left\{{A}_{{\widetilde x}}\right\}, {\widetilde\sF}\right).\\
	\eean 
	Using above inclusions prove following statements:
	\begin{enumerate}
		\item The given by the equation \eqref{top_sheaf_lift_eqn} map $
		\lift_p^{\text{sheaf}}:\mathscr S^{A}\left(\sX\right)\to p^{-1}\mathscr S^{A}\left(\widetilde\sX\right)$
		(cf. Remark \ref{top_sheaf_lift_rem}) can be regarded as restriction  of the given by \eqref{top_glo_lift_h_eqn} map.
		$$
		\lift_p^{\text{top}}:	C\left({\sX}, \left\{{A}_{{x}}\right\}, {\sF}\right)\hookto C\left({\widetilde\sX}, \left\{{A}_{{\widetilde x}}\right\}, {\widetilde\sF}\right),
		$$	
		(cf. Remark \ref{top_top_lift_rem}), i.e. $\lift_p^{\text{sheaf}}= \left.\lift_p^{\text{top}}\right|_{\mathscr S^{A}\left(\sX\right)}$.
		\item Similarly the given  by the Definition \ref{top_lift_sh_desc_defn} $ \lift^p_{\widetilde{\sU}}$ and $\desc_{p}$ can be regarded as restrictions of the given by the Definition \ref{top_lift_desc_defn} maps $ \lift^p_{\widetilde{\sU}}$ and $\desc_{p}$.
		\item If $\desc^c_p$ is a compactly supported $p$-descent (cf. Definition \ref{desc_c_sheaf_defn}) then one has
		\be\label{top_cdesc_soft_eqn}
		\begin{split}
		\desc^c_p\left(C_c\left({\widetilde\sX}, \left\{{A}_{{\widetilde x}}\right\}, {\widetilde\sF}\right)\cap 	p^{-1}\mathscr S^{A}\left(\widetilde\sX\right) \right) =\\= C_c\left({\sX}, \left\{{A}_{{ x}}\right\}, {\sF}\right)\cap 	\mathscr S^{A}\left(\sX\right).
		\end{split}
		\ee
		
	\end{enumerate}
\end{exercise}

\begin{remark}
	Henceforth we will use the unique notation $\lift_p$ for both $\lift_p^{\text{top}}$ and $\lift_p^{\text{sheaf}}$. Similarly the unique notations  
	$ \lift^p_{\widetilde{\sU}}$, $\desc_{p}$ and $\desc^c_p$ will be used for objects given by the Definitions \ref{top_lift_sh_desc_defn} and \ref{top_lift_desc_defn}. 
\end{remark}
\begin{exercise}  
	Let $R\subset C_0\left(\sX \right)$ be a $c$-soft algebra, and let $X$ be an $R$ module. Let $\mathscr S^X$ be the $X$-sheaf (cf. Definition \ref{top_x_sheaf_defn}). Denote by
	$$
	\mathscr S^X\left(\sX \right)_c \bydef \left\{\left. a \in \mathscr S^X\left(\sX\right) \right| \supp a \text{ is compact }\right\}
	$$
	where $\supp$ is given by the equation \eqref{sheaf_supp_eqn}. Prove that if $p: \widetilde \sX \to \sX$ is a covering then similarly to 	\eqref{top_compactly_supported_descent_defn} there is an $R$-module homomorphism
	\be\label{desc_c_sheaf_eqn}
	\desc^c_p: p^{-1}\mathscr S^X\left(\widetilde \sX \right)_c\to 	\mathscr S^X\left(\sX \right)_c.
	\ee
\end{exercise}

\begin{definition}\label{desc_c_sheaf_defn}
	The given by \eqref{desc_c_sheaf_eqn} $R$-linear homomorphism $\desc^c_p$ is said to be a \textit{compactly supported} $p$-\textit{descent}
\end{definition}
\begin{remark}
The Definition \ref{desc_c_sheaf_defn} is compliant with the Definition \ref{top_compactly_supported_descent_defn}.
\end{remark}

\begin{exercise}\label{top_pre_sheaf_exer} 
Let $\sX$ be a locally compact, paracompact, Hausdorff space, and let $R \subset C_0\left(\sX \right)$ be a $c$-soft *-subalgebra (cf. Definition  \ref{top_soft_r_defn}). 	
Let $\mathscr F$ be an $\mathscr S^R$ -{module} (cf. Definition \ref{sheaf_of_modules_den})
Let $p: \widetilde\sX \to \sX$ be a transitive covering.
 Let $\widetilde{\mathscr F}\bydef p^{-1}\mathscr F$ be an inverse image (cf. Definition \ref{sheaf_inv_im_defn}), and let $\mathscr F\left( \sX\right) \subset  \widetilde{\mathscr F}\left( \widetilde\sX\right)$ be a natural inclusion of $\mathscr S^R$-modules. Let $\Aut_{ \mathscr S^R\left( \sX\right)}\left( \widetilde{\mathscr F}\left( \widetilde\sX\right)\right)$ be a group of  $ \mathscr S^R\left( \sX\right)$-module automorphisms.
 Prove following statements:
 \begin{enumerate}
 	\item [(a)] There is a natural isomorphism of groups
 	$$
 	\left\{ \left.g \in \Aut_{\mathscr S^R\left( \sX\right)}\left( \widetilde{\mathscr F}\left( \widetilde\sX\right)\right)~\right|~ ga = a;~~\forall a \in \mathscr F\left( \sX\right) \right\}\cong G\left(\left. \widetilde\sX \right| \sX \right).
 	$$
 	\item[(b)] The above isomorphism of groups yields a natural  action $G\left(\left. \widetilde\sX \right| \sX \right) \times \widetilde{\mathscr F}\left( \widetilde\sX\right) \to \widetilde{\mathscr F}\left( \widetilde\sX\right)$ such that
 	$$
\mathscr F\left( \sX\right) 	= \widetilde{\mathscr F}\left( \widetilde\sX\right)^{G\left(\left. \widetilde\sX \right| \sX \right)}\bydef \left\{\left.a\in \widetilde{\mathscr F}\left( \widetilde\sX\right)~~\right|~ a = g a;~ \forall g \in G\left(\left. \widetilde\sX \right| \sX \right)\right\}
 	$$
 	
 	 \item[(c)]If $\mathscr F$ is a sheaf of *-algebras then any $g \in G\left(\left. \widetilde\sX \right| \sX \right)$ corresponds to an *-automorphism.
 \end{enumerate}

\end{exercise}

\begin{exercise}\label{top_hom_sub_exer} Let $\sX$ be a locally compact, paracompact, Hausdorff space, and let $R \subset C_0\left(\sX \right)$ a $c$-soft   *-subalgebra  \ref{top_soft_r_defn}). Let both $F$, $G$ be $R$-modules and let $\mathscr S^F$ and $\mathscr S^G$ be $F$-sheaf and $G$-sheaf respectively (cf. Definition \ref{top_x_sheaf_defn}).
	There is a sheaf  of local morphisms  $\mathscr Hom \left(\mathscr S^F,
	\mathscr S^G \right)$ (cf. Definition \ref{sheaf_hom_defn}).
Prove following statements.
\begin{enumerate}
	\item The sheaf $\mathscr Hom \left(\mathscr S^F,
	\mathscr S^G \right)$ has a natural structure of $\mathscr S^R$-module (cf. Definition \ref{sheaf_of_modules_den}).
	\item There is a natural inclusion
	\be\label{top_hom_eqn}
\begin{split}
	p^{-1}_{\mathscr Hom \left(\mathscr S^F,
		\mathscr S^G \right)} : \mathscr Hom \left(\mathscr S^F,
	\mathscr S^G \right)\left(\sX \right) \hookto\\\hookto \mathscr Hom \left(p^{-1}\mathscr S^F,
	p^{-1} \mathscr S^G \right)\left(\widetilde \sX \right).
\end{split}
\ee
	\item If both $p_1: \widetilde \sX_1 \to \sX$	and  $p_2: \widetilde \sX_2 \to \widetilde \sX_1$	are coverings and $p= p_2 \circ p_1$ then for every $D \in\mathscr Hom \left(\mathscr S,
	\mathscr G \right)$ one has
	\be\label{top_dirac_comp_eqn}
	p^{-1} D = p_2^{-1} \left(  p_1^{-1} D \right).
	\ee
\end{enumerate}
\end{exercise}

\begin{definition}\label{top_smooth_inv_im_defn}
	Let us	consider the situation of the 
	Exercise \ref{top_hom_sub_exer}, and suppose that $D \in \mathscr Hom \left(\mathscr S^F,
	\mathscr S^G \right)\left( M\right)$, $\xi \in \mathscr S^F \left( M\right)$. We say that both  $p^{-1}_{\mathscr Hom \left(\mathscr S^F,
		\mathscr S^G \right)}\left(  D\right) \in \mathscr Hom \left(p^{-1}\mathscr S^F,
	p^{-1} \mathscr S^G \right)\left(\widetilde M\right)$ and $p^{-1}_{\mathscr S^F}\left(  \xi\right) \in \mathscr S^F \left(\widetilde M\right)$ (cf. equation \eqref{top_cs_sheaf_inc_eqn})  are the $p$-\textit{inverse images} of $D$ and $\xi$ respectively. We write 
	\be\label{top_smooth_inv_im_eqn}
	\begin{split}
		p^{-1} D \stackrel{\mathrm{def}}{=} p^{-1}_{\mathscr Hom \left(\mathscr S^F,
			\mathscr S^G \right)}\left(  D \right),\\
		p^{-1} \xi \stackrel{\mathrm{def}}{=} p^{-1}_{\mathscr S^F }\left(  \xi \right).
	\end{split}
	\ee
\end{definition}

\begin{empt}\label{top_smooth_sh_empt}
	Let $M$ be a smooth manifold, $F$ and $G$ as smooth bundles $R = \Coo_0\left( M\right)$ (cf. Definition \ref{top_sm_bundle_defn}) and let both $\Ga^\infty\left( M, F\right)$, $\Ga^\infty\left( M, G\right)$   are spaces of smooth sections (cf. Definition \ref{top_sm_sec_defn}).	Any differential operator  $D: \Ga^\infty\left( M, F\right)\to \Ga^\infty\left( M, G\right)$ (cf. Definition \ref{do_man_defn}) can be regarded as global section of 
	$
	\mathscr Hom \left(\mathscr S^{\Ga^\infty\left( M, F\right)},
	\mathscr S^{\Ga^\infty\left( M, G\right)} \right),
	$
	i.e.
	$$
	D\in \mathscr Hom \left(\mathscr S^{\Ga^\infty\left( M, F\right)},
	\mathscr S^{\Ga^\infty\left( M, G\right)} \right)\left( M\right). 
	$$
where $\mathscr Hom$ means a sheaf of local morphisms (cf. Definition \ref{sheaf_hom_defn}).
 If  	$D\left(M, E_1, E_2 \right)$ is a space  differential operators (cf. Definition \ref{do_man_defn}) then there is a $D\left(M, E_1, E_2 \right)$-sheaf $\mathscr S^{D\left(M, E_1, E_2 \right)}$ (cf. Definition \ref{top_x_sheaf_defn}) which is an $\mathscr S^{\Coo_0\left(M\right)}$-module (cf. Definition \ref{sheaf_of_modules_den}). We denote it by $\mathscr D\left(M, E_1, E_2\right)$.
	
\end{empt}

\begin{exercise}\label{top_smooth_sh_exer} Consider the situation \ref{top_smooth_sh_empt}. Prove following statements:
	\begin{enumerate}
		\item[(i)] There is the natural inclusion $$\mathscr D\left(M, E_1, E_2\right)\subset \mathscr Hom \left(\mathscr S^{\Ga^\infty\left( M, F\right)},
		\mathscr S^{\Ga^\infty\left( M, G\right)} \right).$$
		\item[(ii)]  If $p: \widetilde M \to M$ is a covering (cf. Definition \ref{top_covering_defn}) and the structure of smooth manifold $\widetilde M$ is given by the Proposition then there is a unique  smooth bundle $\widetilde E\to \widetilde M$ such that there is the natural isomorphism
		$$
			p^{-1} \mathscr S^{\Ga^\infty\left( M, E\right)}\cong  \mathscr  S^{\Ga^\infty\left(\widetilde M, \widetilde E\right)}.
		$$
		Moreover if we forgot a smooth structure then $\widetilde E$ is the {inverse image} of $E$ in the sense of the Definition \ref{vb_inv_img_funct_defn}.
		\item[(iii)] If we use the above notation then there is the following isomorphism of sheaves
		$$
			p^{-1} \mathscr D\left(M, E_1, E_2\right)\cong  \mathscr \mathscr \mathscr D\left(M, \widetilde E_1, \widetilde E_2\right).
		$$
	\end{enumerate} 

\end{exercise}

\begin{empt}\label{top_h_sp_empt}
Let $M$ be a smooth manifold and there is a positive functional 
\bean
\tau : C_c\left( M\right) \to \C,\\
a \mapsto \int_{M} a~ d\mu.
\eean
(cf. Theorem \ref{meafunc_thm}) such that $a > 0  \Rightarrow \tau\left(a \right)>0$. Suppose that $E$ is a smooth vector bundle (cf. Definition \ref{top_sm_bundle_defn}) with a sesquilinear form $\varphi: E\times_M  E\to \C$ (cf. Definition \ref{top_herm_bundle_form_defn}). Let $\Ga^\infty_c\left( M, E\right) \bydef  \Ga_c\left( M, E\right)\cap \Ga^\infty\left( M, E\right)$ is a $\Coo\left( M\right)$-module of smooth sections (cf. Definition \ref{top_sm_sec_defn}) with compact support (cf. Definitions \ref{top_cc_supp_defn} and \ref{phi_supp_defn}) then there is a  sesquilinear product
\bean
\left\langle \cdot, \cdot \right\rangle_c: \Ga_c^\infty\left( M, E\right)\times \Ga_c^\infty\left( M, E\right) \to C_c\left( M\right)
\eean
which comes from \eqref{top_ggc_eqn}. On the other hand there is a product
\be\label{top_prod_eqn}
\begin{split}
\left(\cdot, \cdot\right): \Ga_c^\infty\left( M, E\right)\times \Ga_c^\infty\left( M, E\right)\to \C,\\
\left(a, b \right) \mapsto\tau\left( \left\langle a, b \right\rangle_c\right),
\end{split}
\ee
i.e. $\Ga_c^\infty\left( M, E\right)$ is a pre-Hilbert space. If $L^2\left(M, E,  \mu\right)$ is the Hilbert norm completion of $\Ga_c^\infty\left( M, E\right)$ then there is a representation
\be\label{top_h_sp_eqn}
C_0\left( M \right) \hookto B\left( L^2\left(M, E,  \mu\right)\right) 
\ee
which is faithful because $a > 0 \Rightarrow \tau\left(a \right)>0$.
 If $P \in D\left(M,E\right): \Ga^\infty\left( M, E\right)\to \Ga^\infty\left( M, E\right)$ is s differential operator (cf. Definition \ref{do_man_defn} and Remark \ref{do_man_rem}) then from $P\left( \Ga_c^\infty\left( M, E\right)\right) \subset \Ga_c^\infty\left( M, E\right)$ it follows that
$$
P \in 	\L\left(\Ga_c^\infty\left( M, E\right)\right)
$$
(cf. Equation \eqref{l_dag_eqn}). If 
\be\label{top_diff*_eqn}
D^*\left(M,E\right) \bydef \left\{\left.P \in D\left(M,E\right)\cap  \L^\dagger\left( \Ga_c^\infty\left( M, E\right)\right)\right|P^*\in D\left(M,E\right) \right\}
\ee
then $D^*\left(M,E\right)$ is an $O^*$-algebra (cf. Definition \ref{o*alg_defn}).
\end{empt}
\begin{definition}\label{top_diff*_defn}
The given by the equation \eqref{top_diff*_eqn} $O^*$-algebra  $D^*\left(M,E\right)$ is said to be an \textit{algebra of adjointable differential operators}.
\end{definition}
\begin{remark}
The $O^*$-algebra  $D^*\left(M,E\right)$ is a module over a $c$-soft *-algebra  $\Coo_0\left(M \right)$ (cf. Definition \ref{top_soft_r_defn} and Exercise \ref{smooth_soft_exer}).
\end{remark}
\begin{lemma}\label{top_diff*_lem}
If $\mathscr S^{D^*\left(M,E\right)}$ is the $D^*\left(M,E\right)$-sheaf (cf. Definition \ref{top_x_sheaf_defn}) then the natural inclusion
$$
D^*\left(M,E\right)\subset \mathscr S^{D^*\left(M,E\right)}\left( M\right) 
$$
is an isomorphism.
\end{lemma}
\begin{proof}
Let $s \in\mathscr S^{D^*\left(M,E\right)}\left( M\right)$ be a global section. Let us define a corresponding to $s$ operator $A(s)\in \L\left(\Ga^\infty_c\left( M, E\right)  \right)$  where $\L\left(\cdot \right) $ is given by the equation \eqref{l_dag_eqn}.   If $\xi \in \Ga^\infty_c\left( M, E\right)$ then $\sY \bydef \supp \xi$ is compact. There is a family of open subsets $\left\{\sU_\a\subset \sX\right\}_{\a\in\mathscr A}$ such that for all $\a \in\mathscr A$ a restriction $s|_\a$ is represented by an operator $A_\a \in D^*\left(M,E\right)$.  There is a family of open subsets $\left\{\sU_\a\subset \sX\right\}_{\a\in\mathscr A}$ such that for all $\a \in\mathscr A$ a restriction $s|_\a$ is represented by an operator $A_\a \in D^*\left(M,E\right)$. Let $\sum_{\a \in \mathscr A_0}f_\a$ be a {covering sum} for $\sY$ {dominated} by the family 
$\left\{\sU_\a\right\}$ (cf. Definition \ref{top_covering_sum_defn})
We suppose that $f_\a\in \Coo\left( M\right)$ for all $\a\in  \mathscr A_0$ (cf. Remark \ref{top_smooth_part_unity_rem}). We define 
$$
A(s) \xi \bydef \sum_{\a \in \mathscr A_0}f_\a A_\a\xi
$$
We leave to the reader the full proof of that the above equation yields an operator  in $\L\left(\Ga^\infty_c\left( M, E\right)  \right)$ and it is a differential operator (cf. Definition \ref{do_man_defn}). Similarly we define an adjoint
$$
A(s)^* \xi \bydef \sum_{\a \in \mathscr A_0} A^*_\a f_\a\xi.
$$
We leave to the reader proof of following statements:
\begin{itemize}
	\item The operator $A(s)^*$ is adjoint to $A(s)$, i.e. $A(s)\in \L^\dagger\left(\Ga^\infty_c\left( M, E\right)  \right)$  where $\L^\dagger\left(\cdot \right) $ is given by the equation \eqref{l_dag_eqn}.
	\item $A(s)^*$ is a differential operator (cf. Definition \ref{do_man_defn}).
	\item $A(s)$ is the representative of the global section $s \in\mathscr S^{D^*\left(M,E\right)}\left( M\right)$.
\end{itemize}

\end{proof}
\begin{exercise}\label{top_diff*_exer}
Consider the situation of the Lemma \ref{top_diff*_lem}. Let $p: \widetilde M\to M$ be a covering, and let $\widetilde M$ has a smooth structure given by the Proposition \ref{top_cov_mani_prop}. Let $\widetilde E$ be the inverse image of $E$ by $p$ (cf. Definition \ref{vb_inv_img_funct_defn}). We consider Hilbert spaces $L^2\left( M, E, \mu \right)$, $L^2\left( \widetilde M, \widetilde E, \lift_p~\mu \right)$ where $\lift_p~\mu$ is the $p$-lift of $\mu$ (cf. Definition \ref{top_lift_measure_defn}) and their dense subspaces $\Ga_c^\infty\left( M, E\right)\subset L^2\left( M, E, \mu \right)$, $\Ga_c^\infty\left(\widetilde M,\widetilde E\right)\subset L^2\left(\widetilde M,\widetilde E, \lift_p~\mu \right)$. 
Prove following statements:
\begin{enumerate}
	\item $\widetilde E$ has a natural structure of a smooth bundle (cf. Definition \ref{top_sm_bundle_defn}).
	\item If $\mathscr S^{ D^*\left( M,  E\right)}_x$ is a space of stalks at $x$ (cf. Definition \ref{sheaf_stalk_defn}) then there is a natural involution
	$$
	* : \mathscr S^{ D^*\left( M,  E\right)}_x\xrightarrow{\approx}\mathscr S^{ D^*\left( M,  E\right)}_x.
	$$
	\item If both $D^*\left(M,E\right)$ and $D^*\left(\widetilde M, \widetilde E\right)$ are algebras of adjointable differential operators then there is an isomorphism of sheaves
$$
\mathscr S^{ D^*\left(\widetilde M, \widetilde E\right)}\cong p^{-1}\mathscr S^{D^*\left( M,  E\right)}
$$
where both $\mathscr S^{ D^*\left( M,  E\right)}$ and $\mathscr S^{ D^*\left(\widetilde M, \widetilde E\right)}$  are  ${ D^*\left( M,  E\right)}$ and ${ D^*\left(\widetilde M, \widetilde E\right)}$-sheaves (cf. Definition \ref{top_x_sheaf_defn}).
	
\end{enumerate}
\end{exercise}
\begin{empt}
There is an injective $\C$-homomorphism
$$
\lift_p: \mathscr S^{ D^*\left( M,  E\right)} \left( M\right) \hookto \mathscr S^{ D^*\left(\widetilde M, \widetilde E\right)}\left(\widetilde M\right)
$$
(cf. Definition \ref{top_sheaf_lift_defn}). From the Lemma \ref{top_diff*_lem} it turns out that this homomorhism is equivalent to a $\C$-homomorphism
\be\label{top_diff_*alg_lift_eqn}
D^*\left(p, E \right) : D^*\left(M, E\right)\hookto D^*\left(\widetilde M, \widetilde E\right).
\ee 
The proof of that $D^*\left(p, E \right)$ is a $*$-homomorphism is left  to the reader.  
\end{empt}

\begin{example}\label{top_spinor_smooth_exm}
 Suppose that $M$ is a Riemannian manifold, and $S\to M$ is a spinor bundle (cf. Section \ref{spin_mani_sec}). 
 There is a $\Coo_0\left( M\right)$-module $	\Ga^\infty\left(M, S \right)\subset 	\Ga\left(M, S \right)$ of smooth sections (cf. Definition \ref{top_smooth_m_defn}), and there is an $\Ga^\infty\left(M, S \right)$-sheaf $\mathscr S^{\Ga^\infty\left(M, S \right)}$ (cf. Definition \ref{top_x_sheaf_defn}). Similarly to \ref{top_h_sp_empt} we define a Hilbert space $L^2\left(M, S, \mu\right)$ which is a completion of $\Ga^\infty\left(M, S \right)$ and a faithful representation 
 $$
 C\left( M \right) \hookto B\left( L^2\left(M, S,  \mu\right)\right) 
 $$
 (cf. equation \eqref{top_h_sp_eqn}). If	\bean
 \left(\Coo\left( M\right), L^2\left(M, S \right), \Dslash   \right).
 \eean
is a given by the equation \eqref{comm_sp_tr_eqn} spectral triple then $\Dslash$ can be regraded as an operator
$$
\Dslash: \Ga^\infty\left(M, S \right)\to \Ga^\infty\left(M, S \right).
$$
 The above equation can be rewritten by the following way
 \be
 \Dslash \in \mathscr End \left( \mathscr S^{\Ga^\infty\left(M, S \right)} \right) \left(M\right) \bydef \mathscr Hom \left( \mathscr S^{\Ga^\infty\left(M, S \right)},  \mathscr S^{\Ga^\infty\left(M, S \right)} \right) \left( M\right) 
 \ee
 where $\mathscr Hom$ means the {sheaf  of local morphisms} (cf. Definition \ref{sheaf_hom_defn}).
 Suppose that $p: \widetilde M \to M$ is a finite fold covering. If both $\widetilde{\mathscr G}^\infty$ and $\widetilde S$ are $p$-inverse images of ${\mathscr G}^\infty$ and $S$ respectively (cf. Definitions \ref{vb_inv_img_funct_defn} and \ref{sheaf_inv_im_defn}).
 then there is an operator  
  \be\label{top_d_shaef_eqn}
p^{-1}\Dslash \in \mathscr End \left(p^{-1}\mathscr S^{\Ga^\infty\left(M, S \right)} \right) \left(\widetilde M\right)
 \ee
 such that
 \be\label{top_d_shaef_m_eqn}
p^{-1}\Dslash  : \Ga^\infty\left(\widetilde M,\widetilde S \right)\to \Ga^\infty\left(\widetilde M,\widetilde S \right).
 \ee

\end{example}

\begin{exercise}\label{top_distr_exer}
If $M$ be a smooth manifold then   a space of distribution densities  $\Coo_c\left(M \right)'$ (cf. Definition \ref{top_distr_dens_def}) is a $\Coo_0\left( M\right)$-module. Let   $\mathscr D_M\bydef \mathscr S^{\Coo_c\left(M \right)'}$ be the $\Coo_c\left(M \right)'$-sheaf (cf. Definition \ref{top_x_sheaf_defn}).
	Prove following statements:
	\begin{enumerate}
		\item[(i)] If $p: \widetilde M \to M$ is a covering and the given by the Proposition  \ref{top_cov_mani_prop} structure of a smooth manifold $\widetilde M$  is implied then $\mathscr D_{\widetilde M}= p^{-1}\mathscr D_M$ where $p^{-1}\mathscr D_M$  means the inverse image (cf. Definition \ref{sheaf_inv_im_defn}).
		\item[(ii)] Using above statements prove that there is the natural injective homomorphism
		\be\label{top_distr_eqn}
		\Coo_c\left(p \right)':\Coo_c\left(M \right)'\hookto \Coo_c\left(\widetilde M \right)'
		\ee
		of $\Coo\left(M \right)$-modules.
		\item[(iii)] If the covering $p$ is transitive then one has:
		\begin{enumerate}
			\item [(a)] There is a natural isomorphism of groups
			$$
			\left\{ \left.g \in \Aut_{\Coo\left(M\right)}\left( \Coo_c\left(\widetilde M \right)'\right) ~\right|~ ga = a;~~\forall a \in \mathscr F\left( \sX\right) \right\}\cong G\left(\left. \widetilde M \right| M \right).
			$$
			\item[(b)] The above isomorphism of groups yields a natural  action $G\left(\left. \widetilde M \right| M \right) \times \Coo_c\left(\widetilde M \right)' \to \Coo_c\left(\widetilde M \right)'$ such that
			\bean
			\Coo_c\left( M \right)' 	= \left( \Coo_c\left(\widetilde M \right)'\right) ^{G\left(\left. \widetilde M \right|M  \right)}\bydef\\\bydef  \left\{\left.a\in \Coo_c\left(\widetilde M \right)'~~\right|~ a = g a;\quad \forall g \in G\left(\left. \widetilde M \right| M \right)\right\}.
			\eean
			
		\end{enumerate}
		
	\end{enumerate}
\end{exercise}

\section{Lifts of continuous functions and bundles}

\begin{lemma}\label{top_comp_sum_lem}
		Let $\widetilde{   \mathcal X }$ be a  connected, locally compact, Hausdorff space. Let $p: \widetilde{   \mathcal X }\to \sX$ be a transitive covering. Let us use the notation  \eqref{top_cb_defna_eqn}.  Any $\widetilde a \in A_c\left( \widetilde{\sX}\right)$ can be represented by the following way
			$$
	\widetilde a =	\widetilde a_1 + ... + \widetilde a_m 
		$$
		where  $ \widetilde a_j\in A_c\left( \widetilde{   \mathcal X }\right)$ and the set $\widetilde \sU_j = \left\{\left.\widetilde x \in \widetilde{   \mathcal X}\right|\widetilde a_j\left(\widetilde x\right)\neq 0 \right\}$ is homeomorphically mapped onto $p\left( \widetilde \sU_j\right)$   for all $j = 1,..., m$. Moreover if	$\widetilde a$ is positive then one can select positive $\widetilde a_1,..., \widetilde a_m$. 
\end{lemma}
\begin{proof}
If $\sum_{j=1}^m\widetilde f_j$ is  a subordinated to $p$ covering  sum for $\supp \widetilde  a$ (cf. Definitions  \ref{top_covering_sum_defn}  and  \ref{top_covering_sum_subordinated_defn}) then one has
$$
	\widetilde a =  \widetilde f_{1}	\widetilde a 	  + ... +  \widetilde f_{m}	\widetilde a 
$$
or, equivalently
$$
	\widetilde a =	\widetilde a_1 + ... + \widetilde a_m 
	$$
	where $\widetilde a_j \bydef  \widetilde f_{j} \widetilde a 	$ for every $j=1,...,m$. If $\widetilde a$ is positive then $\widetilde a_j$ is positive for each $j=1,...,m$. 
\end{proof}
	\begin{empt}
	Let $\sX$ be a locally compact, Hausdorff space; and for each $x$ in $\sX$, let $A_x$ be a (complex) Banach space. Let us consider a continuity structure $\sF$ for $\sX$ and $\left\{A_x\right\}$ (cf. Definition \ref{operator_fields_continuity_defn}). If $p: \widetilde \sX \to \sX$ is a covering then $C_c\left(\widetilde \sX\right)$ is a  $C_0\left( \sX\right)$-module.  Let $\widetilde f \in C_c\left( \widetilde \sX\right)$ and let 
	$\sum_{l=1}^m \widetilde a_l$ be a subordinated  to  $p$ covering sum for $\supp\widetilde f$ (cf. Definitions \ref{top_covering_sum_defn} and \ref{top_covering_sum_subordinated_defn}). Without loss of generality  one can suppose that for all $l=1,...,m$ there is an open subset $\widetilde \sU_l$ such that $ \supp \widetilde a_l\subset  \widetilde \sU_l$ and $ \widetilde \sU_l$ is homeomorphically mapped onto $p\left( \widetilde \sU_l \right)$. If $\xi \in C_c\left( \lift_p\left[C\left(\sX, 
	\left\{A_x\right\}, \sF \right)\right]\right)$ then we define
\be\label{top_xi_eqn}
\widetilde \xi \bydef \sum_{l=1}^m \widetilde f\sqrt{\widetilde a_l} \lift^p_{\widetilde{\sU}_l}\left(\desc_p\left( \sqrt{\widetilde a_l}\right) \xi \right) 
\ee
where both $\lift^p_{\widetilde{\sU}}$ and $\desc_p$ are given by the Definition \ref{top_lift_desc_defn}. Element $\widetilde \xi$ is represented by a family 
	\be\label{top_tensor_lcfam_eqn}
\left\{ \widetilde f\left(\widetilde x \right)  c^{-1}_{\widetilde x}\left( \xi_{p\left( \widetilde x\right) }\right)  \right\}_{{\widetilde x}\in {\widetilde \sX}}\quad  \text{ where } c_{\widetilde x}\text{ is given by } \eqref{top_ct_iso_eqn}.
\ee
because $\sum_{j=1}^n  \widetilde a_j\left( \widetilde x\right)=1$  for every $\widetilde x\in  \supp \widetilde f$. From the equation \eqref{top_tensor_lcfam_eqn} it follows that $\widetilde\xi$ uniquely depends on the pair $\left(\widetilde f, \xi \right) \in C_c\left(\widetilde \sX\right)\times C\left(\sX, 
\left\{A_x\right\}, \sF \right)$, so one has a bilinear map
\bean
\varphi: C_c\left(\widetilde \sX\right)\times C\left(\sX, 
\left\{A_x\right\}, \sF \right)\to C_c\left( \lift_p\left[C\left(\sX, 
\left\{A_x\right\}, \sF \right)\right]\right) ,\\
\left(\widetilde f, \xi \right)\mapsto \widetilde\xi.
\eean
The reader can easily state that $\varphi\left( \widetilde ff, \xi\right) = \varphi\left( \widetilde f, f\xi\right)$ for all $f\in C_0\left( \sX\right)$ it turns out that $\varphi$ yields  the linear map 
	\be\label{top_tensor_lccompact_iso_eqn}
\begin{split}
	\phi: C_c\left(\widetilde \sX\right)\otimes_{C_0\left( \sX\right)}C\left(\sX, 
	\left\{A_x\right\}, \sF \right) \to C_c\left( \lift_p\left[C\left(\sX, 
	\left\{A_x\right\}, \sF \right)\right]\right). \end{split}
\ee
\end{empt}

\begin{lemma}\label{top_tensor_lcompact_lem}
	The given by \eqref{top_tensor_lccompact_iso_eqn} map $\phi$ is the isomorphism of left  $C_0\left(\widetilde \sX\right)$-modules.
\end{lemma}
\begin{proof}  
If $\left(\widetilde f, \xi \right) \in C_c\left(\widetilde \sX\right)\times C\left(\sX, 
\left\{A_x\right\}, \sF \right)$ then from \eqref{top_xi_eqn} it follows that
$$
\varphi\left(\widetilde f, \xi \right)  = \phi\left(\widetilde f\otimes \xi \right)= \sum_{l=1}^m \widetilde f\sqrt{\widetilde a_l} \lift^p_{\widetilde{\sU}_l}\left(\desc_p\left( \sqrt{\widetilde a_l}\right) \xi \right)  
$$
where $\lift^p_{\widetilde{\sU}_l}$ is $p$-$\widetilde{\sU}_l$-{lift} (cf. Definition \ref{top_lift_desc_defn}). If $\widetilde a \in C_0\left(\widetilde \sX\right)$ then from the above equation it turns out that
\bean
\phi\left(\widetilde a\left( \widetilde f\otimes \xi \right) \right)= \sum_{l=1}^m \widetilde a\widetilde f\sqrt{\widetilde a_l} \lift^p_{\widetilde{\sU}_l}\left(\desc_p\left( \sqrt{\widetilde a_l}\right) \xi \right)=\\
  =\widetilde a \sum_{l=1}^m \widetilde f\sqrt{\widetilde a_l} \lift^p_{\widetilde{\sU}_l}\left(\desc_p\left( \sqrt{\widetilde a_l}\right) \xi \right)= \widetilde a\phi\left( \widetilde f\otimes \xi \right),
 \eean
i.e. $\phi$ is a homomorphism of left  $C_0\left(\widetilde \sX\right)$-modules because \\ $C_c\left(\widetilde \sX\right)\otimes_{C_0\left( \sX\right)}C\left(\sX, 
\left\{A_x\right\}, \sF \right)$ is a generated by elements $\widetilde f\otimes \xi$ Abelian group.\\ 
Let $\sum_{j =1}^n \widetilde f_j \otimes \xi_j$ be any element in  $C_c\left(\widetilde \sX\right)\otimes_{C_0\left( \sX\right)}C\left(\sX, 
\left\{A_x\right\}, \sF \right)$. 
The finite union $\widetilde \sY \bydef \cup_{j=1}^n \supp \widetilde f_j$ of compact sets is compact. Let 
$\sum_{l=1}^m \widetilde a_m$ be a subordinated  o  $p$ covering sum for $\widetilde \sY$ (cf. Definitions \ref{top_covering_sum_defn} and \ref{top_covering_sum_subordinated_defn}). Without loss of generality  one can suppose that for all $l=1,...,m$ there is an open subset $\widetilde \sU_l \subset \widetilde\sX$ such that $\supp \widetilde a_l\subset  \widetilde \sU_l$ and $ \widetilde \sU_l$ is homeomorphically mapped onto $p\left( \widetilde \sU_l \right)$. 
From the  equation \eqref{top_xi_eqn} it turns out that
	\be\label{top_tensore_lcfam_eqn}
	\begin{split}
		\sum_{j=1}^n \widetilde f_j \otimes \xi_j= \sum_{l=1}^m \sum_{j=1}^n \widetilde{e}_{l}\widetilde{e}_{l}\widetilde f_j \otimes \xi_j= \sum_{l=1}^m\sum_{j=1}^n \widetilde{e}_{l} \desc_p\left(\widetilde e_j \right)\widetilde f_j \otimes \xi_j = \\
		\sum_{l=1}^m \widetilde{e}_{l}\otimes \desc_p\left(\widetilde{e}_{l}\sum_{j=1}^n\widetilde f_j   \right)\xi_j 
	\end{split}
	\ee
	where $\widetilde e_l\bydef\widetilde a_l^{~\nicefrac{1}{2}}$. 
	If $\phi\left(\sum_{j=1}^n \widetilde f_j \otimes \xi_j \right) = 0$ then from \eqref{top_tensor_lcfam_eqn}  it follows that\\ $\sum_{j=1}^n \widetilde f_{j}\left(\widetilde x \right)  c^{-1}_{\widetilde x}\left( \xi_{j p\left(\widetilde x \right)}\right)=0$ for all $\widetilde x \in \widetilde\sX$. It turns out that
	\be\label{top_t_eqn}
	\sum_{l=1}^m\sum_{j=1}^n \widetilde{e}_{l}\left(\widetilde x \right)\widetilde f_{j}\left(\widetilde x \right)  c^{-1}_{\widetilde x}\left( \xi_{j p\left(\widetilde x \right) }\right)=0 \quad \widetilde x \in \widetilde\sX.
	\ee  
	Since $\widetilde{e}_{l}\left(\widetilde x \right)\in \R$ and $\widetilde{e}_{l}\left(\widetilde x \right)\ge 0$ from \eqref{top_t_eqn} it follows that $\widetilde{e}_{l}\left(\widetilde x \right)\sum_{j=1}^n\widetilde f_{j}\left(\widetilde x \right)  c^{-1}_{\widetilde x}\left( \xi_{j p\left(\widetilde x \right) }\right)=0$ for all $l = 1, ..., m$ and $\widetilde x \in \widetilde\sX$. 
On the other hand family $\left\{\widetilde{e}_{l}\left(\widetilde x \right)\sum_{j=1}^n\widetilde f_{j}\left(\widetilde x \right)  c^{-1}_{\widetilde x}\left( \xi_{j p\left(\widetilde x \right) }\right)\right\}$ corresponds to
$
\lift^p_{\widetilde \sU_j}\left(  \desc_p\left( \widetilde{e}_{l}\sum_{j=1}^n\widetilde f_{j}\right) \xi_j\right), 
$ so one has $\lift^p_{\widetilde \sU_l}\left(  \desc_p\left( \widetilde{e}_{l}\sum_{j=1}^n\widetilde f_{j}\right) \xi_j\right)=0$. Since $\lift^p_{\widetilde \sU_j}$ is an isomorphism one concludes that $$\desc_p\left( \widetilde{e}_{l}\sum_{j=1}^n\widetilde f_{j}\right)\xi_j=0.
$$
Substitution of the above equation into the right part of \eqref{top_tensore_lcfam_eqn} yields the following 
$
	\sum_{j=1}^n \widetilde f_j \otimes \xi_j=0
$.
It turns out that the map $\phi$ is injective. 
	For all $\widetilde \xi \in C_c\left( \lift_p\left[C\left(\sX, 
	\left\{A_x\right\}, \sF \right)\right]\right)$ the set $\supp \widetilde \xi$ is compact. If  $
	\sum_{l=1 }^k\widetilde{a}'_l
	$
	is a subordinated  o $p$ covering sum of compact set $\supp\widetilde \xi$ (cf. Definitions \ref{top_covering_sum_defn} and  \ref{top_covering_sum_subordinated_defn}) then one has
	$$
	\widetilde \xi= \sum_{l=1 }^k\widetilde{a}'_l \widetilde \xi= \sum_{l=1 }^k\widetilde{e}'_l\widetilde{e}'_l \widetilde \xi=\sum_{l=1 }^k\widetilde{e}'_l\desc_p\left( \widetilde{e}'_l  \widetilde \xi \right).
	$$
		where $\widetilde e'_l\bydef\sqrt{\widetilde a'_l}$. From the above equation it turns out that
	
	\bean
	\widetilde \xi = \phi\left( \sum_{l=1}^k \widetilde{e}'_{l}\otimes \desc_p\left( \widetilde{e}'_{l} \widetilde \xi \right) \right), 
	\eean
	it turns out that the map $\phi$ is surjective.
\end{proof}

\section{Geometrical lifts of spectral triples}\label{top_cov_sp_tr_sec}
\paragraph*{} Let $\left(C^{\infty}\left( M\right) , L^2\left( M, S\right) ,\slashed D, J\right)$ be a real commutative spectral triple (cf. the Definition \ref{df:spt-real_defn} and the equation \eqref{comm_sp_tr_eqn}), and  let
 $p:\widetilde{M} \to M$  be a finite-fold covering. From the Proposition \ref{comm_cov_mani_prop} it follows that $\widetilde{M}$ has natural structure of the Riemannian manifold. If $TM$ is a tangent bundle then from the Lemma \ref{top_bundle_cs_ex_lem} it follows that $\Ga\left( M, TM\right)$ is a continuity structure for $M$ and the family $\left\{T_xM\right\}_{x \in M}$ of tangent spaces (cf. Definition \ref{operator_fields_continuity_defn}) such that
$$
\Ga\left( M, TM\right)=C\left(M, \left\{T_xM\right\}, \Ga\left( M, TM\right) \right). 
$$ 
It is known \cite{kobayashi_nomizu:diff_geom} that the tangent bundle $T\widetilde{M}$ of $\widetilde{M}$ is the inverse image of $TM$ (cf. Definition \ref{vb_inv_img_funct_defn}). From the Theorem \ref{top_lift_bundle_lem} it follows that $\Ga\left(\widetilde{M}, T\widetilde{M} \right)$ is the $p$-lift of  $C\left(M, \left\{T_xM\right\}, \Ga\left( M, TM\right) \right)$ (cf. Definition \ref{top_lift_main_defn}), i.e. 
\be\label{top_tf_lift_eqn}
\Ga\left(\widetilde{M}, T\widetilde{M} \right)\cong \lift_p\left[C\left(M, \left\{T_xM\right\}, \Ga\left( M, TM\right) \right)\right]
\ee
The described in the Section \ref{comm_sp_tr_sec} tensor fields on $\widetilde{M}$ are $p$-lifts of  corresponding tensor fields on $M$.  
In particular the metric tensor $g$ on $M$ corresponds to an element of the module $\Ga\left(M, TM \right)\otimes_{C\left( M\right) }\Ga\left(M, TM \right)$. The metric tensor $\widetilde g$ of $\widetilde M$ is given by $\widetilde g = \lift_p\left( g\right)$ (cf. Definition \ref{top_lift_defn}). Similarly  the given by \eqref{top_vol_eqn} volume elements $v$ and $\widetilde v$ on $M$ and $\widetilde M$ satisfy to the equation 
\be\label{top_vol_lift_eqn}
\widetilde v = \lift_p\left(v \right)\quad \text{(cf. Definitions \ref{top_lift_defn} and \ref{top_lift_measure_defn})}.
\ee
Indeed $L^2\left( M, S\right)= L^2\left( M, S, \mu\right)$ where $\mu$ corresponds to the volume element $v$.  
If both 	$\C\ell(M) \to M$ and 	$\C\ell(\widetilde M) \to \widetilde M$ are described in \ref{spin_mani_sec} bundles of $C^*$- algebras then $\Ga\left(\widetilde M, \C\ell(\widetilde M) \right)\cong \lift_p\left[C\left( M, \left\{C\ell(M)_x\right\}, \Ga\left( M, \C\ell( M) \right)\right) \right]$. If  $\widetilde{S} = p^*S$ the inverse image of the spinor bundle $S$ (cf. Definition \ref{vb_inv_img_funct_defn}) then
$$
\Ga\left(\widetilde M, \widetilde S \right)\cong  \lift_p\left[C\left( M, \left\{S_x\right\}, \Ga\left( M, S \right)\right) \right]
$$
For any $x \in M$ there is the natural isomorphism $\C\ell(M)_x \cong \End_\C\left(S_x \right)$, so taking into account $\C\ell( \widetilde M)_{ \widetilde x}\cong \C\ell( M)_{p\left(  \widetilde x\right) }$ and $\End_\C\left( \widetilde S_{ \widetilde x} \right)\cong \End_\C\left(  S_{ p\left( \widetilde x\right)}  \right)$ for any $ \widetilde x \in \widetilde  M$    there is the natural isomorphism $\C\ell( \widetilde M)_{ \widetilde x} \cong \End_\C\left( \widetilde S_{ \widetilde x} \right)$. It follows that the bundle $\widetilde S$ yields the  Spin$^c$-structure on $T\widetilde M$. (cf. Definition \ref{spin_str_defn}). The given by \eqref{top_vol_lift_eqn} volume form yields the measure $\widetilde \mu$ on $\widetilde M$, in result one has the Hilbert space
$$
\widetilde \H \bydef L^2\left(\widetilde M, \widetilde S, \widetilde\mu\right).
$$
Moreover $\widetilde\mu\bydef \lift_p\mu$ is the $p$-lift of $\mu$ (cf. Definition \ref{top_lift_measure_defn}).  If  $C: \Ga\left(M,S \right)\cong \Ga\left(M,S \right)$ be a given be the Proposition  \ref{pr:charge-conj_prop} isomorphism of Abelian groups  then $J : L^2\left( M,  S, \mu\right)\cong L^2\left( M,  S, \mu\right)$ is an extension of $C$ (cf. Remark \ref{pr:charge-conj_rem}). On the other hand there is $\Ga\left(M,S \right)$ is a $C\left(M\right)$-module so one can define a $\Ga\left(M,S \right)$-sheaf $\mathscr S^{\Ga\left(M,S \right)}$ (cf. Definition \ref{top_x_sheaf_defn}). The operator $C$ can be regarded as a global section of a  sheaf of local morphisms (cf. Definition \ref{sheaf_hom_defn}), i.e.
$$
C \in \mathscr Hom \left(\mathscr S^{\Ga\left(M,S \right)}, \mathscr S^{\Ga\left(M,S \right)}\right)\left(M\right).
$$
The bundle $\widetilde S$ is an inverse image of $S$ (cf. Definition \ref{vb_inv_img_funct_defn}) so a sheaf\\ $\mathscr Hom \left(\mathscr S^{\Ga\left(\widetilde M, \widetilde S \right)}, \mathscr S^{\Ga\left(\widetilde M,\widetilde S \right)}\right)$ is 
an inverse image of $\mathscr Hom \left(\mathscr S^{\Ga\left(M,S \right)}, \mathscr S^{\Ga\left(M,S \right)}\right)$ in the sense of the Definition \ref{sheaf_inv_im_defn}. According to the Definition \ref{top_sheaf_lift_defn} one can define $p$-lift 
\be\label{top_wtc_eqn}
\widetilde C \bydef \lift_p\left(C \right) \in \mathscr Hom \left(\mathscr S^{\Ga\left(\widetilde M, \widetilde S \right)}, \mathscr S^{\Ga\left(\widetilde M,\widetilde S \right)}\right)\left(\widetilde M \right). 
\ee
The operator $\widetilde C$ similarly to $C$ corresponds to an isomorphism $\Ga\left(\widetilde M, \widetilde S \right)\cong \Ga\left(\widetilde M, \widetilde S\right)$ of Abelian groups. 
An extension of $\widetilde C$ yields an antiunitary operator $\widetilde J : L^2\left(\widetilde M, \widetilde S, \widetilde\mu\right)\cong  L^2\left(\widetilde M, \widetilde S, \widetilde\mu\right)$.
Using the Definition \ref{top_smooth_m_defn} one can construct the  $\Coo\left(  M\right)$-module $\Ga^\infty\left(  M,  S\right)$ and  $\Coo\left( \widetilde M\right)$-module $\Ga^\infty\left( \widetilde M, \widetilde S\right)$   of smooth sections. Also there is  the Spin$^c$ connection $\nabla^S$ (cf. \eqref{spin_conn_eqn}). Similarly to the Equation \eqref{comm_dirac_eqn} one can define the Dirac operator
\bean
\widetilde \Dslash:\Ga^\infty\left( \widetilde M, \widetilde S\right) \to \Ga^\infty\left( \widetilde M, \widetilde S\right) 
\eean
which can be regarded us an unbounded operator on $\widetilde \sH$. 
 So one has a quadruple $\left( \Coo\left(\widetilde M \right), L^2\left(\widetilde M, \widetilde S, \widetilde\mu\right), \widetilde \Dslash, \widetilde J \right)$.
 \begin{exercise}
 	Prove that the quadruple $\left( \Coo\left(\widetilde M \right), L^2\left(\widetilde M, \widetilde S, \widetilde\mu\right), \widetilde \Dslash, \widetilde J \right)$ is a commutative real spectral triple.
 \end{exercise}
\begin{definition}\label{top_geom_lift_defn}
	The constructed above spectral triple $\left( \Coo\left(\widetilde M \right), L^2\left(\widetilde M, \widetilde S, \widetilde\mu\right), \widetilde \Dslash, \widetilde J \right)$ is said to be the \textit{geometrical} $p$-\textit{lift} of the spectral triple $\left( \Coo\left( M \right), L^2\left( M,  S, \mu\right),  \Dslash, J \right)$.
\end{definition}
\begin{exercise}\label{top_geom_dirac_lift_exer}
Prove that 
\be\label{top_geom_dirac_lift_eqn}
  \widetilde \Dslash = p^{-1}\Dslash \in \mathscr End \left( {\mathscr S}^{\Ga^\infty\left( \widetilde M, \widetilde S\right) } \right) \left(\widetilde M\right)
\ee
(cf. Example \ref{top_spinor_smooth_exm}).
\end{exercise}

\section{Finite-fold coverings}
\paragraph*{} This section contains a generalization of the Theorem \ref{pavlov_troisky_thm} which can be applied to finite-fold coverings of locally compact spaces.

\subsection{Coverings of $C^*$-algebras}
\begin{lemma}\label{top_comm_lem}
Let $p: \widetilde\sX \to \sX$ be a transitive covering  $\widetilde\sX$ and $\sX$. If the $*$-homomorphism  $C_0\left( p\right) : C_0\left( \sX\right)\hookto C_0\left( \widetilde\sX\right)$ is given by \eqref{top_c0p_eqn} 
 then the quadruple $$\left(C_0\left( \sX\right),C_0\left( \widetilde\sX\right), G\left(\left.\widetilde{\sX}~\right|\sX\right), C_0\left(p \right)   \right)$$   is a  noncommutative finite-fold pre-covering (cf. Definition \ref{fin_pre_defn}).  
\end{lemma}
\begin{proof}
There is a one to one correspondence between homeomorphisms of  $\widetilde\sX$ and $*$-automorphisms of $C_0\left( \widetilde\sX\right)$ it follows that
$$
\left\{ \left.g \in \Aut\left(C_0\left( \widetilde\sX\right) \right)~\right|\forall a \in C_0\left(p \right)  \left( C\left( \sX\right)\right) \quad ga = a\right\}\cong G\left(\left.\widetilde{\sX}~\right|\sX\right)
$$
i.e. the quadruple $\left(C_0\left( \sX\right),C_0\left( \widetilde\sX\right), G\left(\left.\widetilde{\sX}~\right|\sX\right), C_0\left(p \right)   \right)$
is an  noncommutative finite-fold pre-covering (cf. Definition \ref{fin_pre_defn}). \end{proof}
\begin{corollary}\label{top_comm_cor}
If $\left(C\left( {\sX}\right), \widetilde{A}, G, \pi \right)$ is an unital noncommutative  finite-fold covering (cf. Definition \ref{fin_unital_defn}) and 
\be\label{top_comm_eqn}
H  \bydef \ker\left(  G \to \left\{\left.g \in \mathrm{Homeo}\left(\widetilde\sX \right)\right| \forall \widetilde x \in \widetilde\sX \quad p\left( \widetilde x\right)= p\left( g\widetilde x\right) \right\}\right)
\ee
then one has:
\begin{enumerate}
	\item [(i)] the spectrum $\widetilde \sX$ of $\widetilde A$ is Hausdorff,
	\item[(ii)] the given by the Proposition  \ref{spectrum_covering_finite_prop} map $p : \widetilde \sX\to \sX$ is a transitive finite-fold covering with $G\left(\left.\widetilde{\sX}~\right|\sX\right)\cong G/H$.
	\item[(iii)] The normal subgroup $H \subset G$ is $\left(C\left( {\sX}\right), \widetilde{A}, G, \pi \right)$-proper (cf. Definition \ref{proper_subgroup_fin_defn})
\end{enumerate}
\end{corollary}
\begin{proof}
(i) Follows from the Corollary \ref{spectrum_ff_p_cor}.\\
(ii) Follows from the equation \eqref{top_comm_eqn} and the Corollary \ref{spectrum_ff_p_cor}.\\
(iii) From \eqref{top_comm_eqn} it follows that $\widetilde{A}^H\cong \left(\widetilde\sX \right)$, however from the Lemma \ref{top_comm_lem} it follows that
$$
\left(C_0\left( \sX\right),\widetilde{A}^H, G/H, C_0\left(p \right)   \right)=\left(C_0\left( \sX\right),C_0\left( \widetilde\sX\right), G\left(\left.\widetilde{\sX}~\right|\sX\right),\pi^H  \right)
$$
is a {noncommutative finite-fold  pre-covering}.
\end{proof}
\begin{lemma}\label{comm_ccr_lem}
	If $\left(C\left( {\sX}\right), \widetilde{A}, G, \pi \right)$ is an unital noncommutative  finite-fold quasi-covering (cf. Definition \ref{fin_unital_defn}).
	If $\widetilde \sX$ is  the spectrum of $\widetilde{A}$  then the dimension of  any irreducible representation of  $\widetilde{A}$  is finite. Moreover  there is $N \in \N$ such that for any ${x} \in {\mathcal X}$ there is $m \le N$ which satisfies to the following condition  $$\mathfrak{rep}_{{x}}\left(\widetilde{A} \right) \approx \mathbb{M}_m\left(\C \right).$$
	
\end{lemma}

\begin{proof}
	From the Definition \ref{principal_non_defn} and the  Remark \ref{unital_rem}  it follows that $\pi \left(1_{C\left(\sX \right) }\right)= 1_{\widetilde A}$. On the other hand there are $\widetilde{a}_1, ..., \widetilde{a}_n \in \widetilde{A}$ such that
	$$
	\widetilde{A}= C\left(\sX \right) \widetilde{a}_1 +...+  C\left(\sX \right)\widetilde{a}_n
	$$
	(cf. Definition \ref{fin_unital_defn}).
	It follows that
	$$
	\forall \widetilde x \in \widetilde \sX \quad \rep_{\widetilde x}\left(\widetilde A \right) = 	 \rep_{\widetilde x}\left(A \right) \rep_{\widetilde x}\left(\widetilde a_1 \right)+...+  \rep_{\widetilde x}\left(A \right) \rep_{\widetilde x}\left(\widetilde a_n \right)
	$$
	For any $x \in \sX$ and $a \in C\left(\sX\right)$ the $\rep_x\left(a \right)= a\left(x\right)  \rep_x\left(1_{C\left( x\right) } \right)$, so 
	$$
	\forall 	\widetilde x \in \widetilde \sX\quad \rep_{\widetilde x}\left(\pi\left( a \right) \right)= a\left( p\left( \widetilde x\right) \right)  \rep_{\widetilde x}\left( 1_{{\widetilde A}}\right) 
	$$
	where the map $p: \widetilde\sX \to \sX$ is given by the Proposition \ref{spectrum_covering_finite_prop}.	 It turns out that
	$$
	\forall 	\widetilde x \in \widetilde \sX \quad \forall \ka \in  \rep_{\widetilde x}\left(\widetilde A \right)\quad \exists k_1,..., k_n \in \C  \quad	\ka= k_1\rep_{\widetilde x}\left(\widetilde a_1 \right)+ k_n\rep_{\widetilde x}\left(\widetilde a_n \right),
	$$
	i.e. the $C^*$-algebra   $\rep_{\widetilde x}\left(\widetilde A \right)$ is a finite-dimensional $\C$-space. So for all $x \in \sX$ there is $m \in \N$ such that $\rep_x\left(\widetilde A \right)\approx \mathbb{M}_m\left(\C \right)$. A $\C$-dimension of $\mathbb{M}_m\left(\C \right)$ equals to $m^2$, so $m < N$ where $N \ge \sqrt{n}$.
\end{proof}
\begin{corollary}\label{comm_comm_cor}
	Let $\left(C\left( {\sX}\right), \widetilde{A}, G, \pi \right)$ be an unital noncommutative  finite-fold covering (cf. Definition \ref{fin_unital_defn}) and $p :  \widetilde \sX \to \sX$ be a given by the Proposition \ref{spectrum_covering_finite_prop} map from the spectrum $\widetilde \sX$ of $\widetilde A$ to the spectrum of $A$. Then the space $\widetilde \sX$ is Hausdorff, the map $p$ is a finite-fold transitive covering and the quadruple  $\left(C\left( {\sX}\right), \widetilde{A}, G, \pi \right)$ is equivalent to the $\left(C\left( \sX\right),C\left( \widetilde\sX\right), G\left(\left.\widetilde{\sX}~\right|\sX\right), C_0\left(p \right)   \right)$ one  where $C_0\left(p \right) : C\left({\mathcal X}\right) \hookto C\left( \widetilde{\mathcal X}\right)$ 	comes from the  finite covering algebraic functor (cf. Definition \ref{top_c_funct_defn}).
\end{corollary}
\begin{proof}   	
	If $\widetilde \sX$ is the spectrum of $\widetilde{A}$ and then from the Corollary  \ref{top_comm_cor} and the Lemma \ref{fin_composition_lem} it follows that if 
	$$
	H  \bydef \ker\left(  G \to \left\{\left.g \in \mathrm{Homeo}\left(\widetilde\sX \right)\right| \forall \widetilde x \in \widetilde\sX \quad p\left( \widetilde x\right)= p\left( g\widetilde x\right) \right\}\right)	
	$$ then there are $H$-regular and $H$-singular unital noncommutative finite-fold coverings 
	\bean
	\left(C\left( \sX\right),\widetilde{A}^H, G/H, C_0\left(p \right)   \right)=\left(C\left( \sX\right),C\left( \widetilde\sX\right), G\left(\left.\widetilde{\sX}~\right|\sX\right),\pi^H  \right),\\
	\left(C\left( \widetilde\sX\right), \widetilde A, H, \left.\Id_{\widetilde A}\right|_{C\left( \widetilde\sX\right)} \right) 
	\eean   
	From the Lemma \ref{comm_ccr_lem} it follows that there is $n\in \N$ such that the dimension of  $\rep_{\widetilde x}\left(\widetilde A \right)\le n$. 	
	Suppose that  $n> 1$.
	From the Proposition \ref{less_n_pi_prop} it follows that there is an open subset $\widetilde \sU'\subset \widetilde \sX$ such that 
	$\forall \widetilde x \in \widetilde \sU'\quad  \dim \rep_{\widetilde x}\left( \widetilde A\right) = n$.
	So the ideal $\left.\widetilde{A}\right|_{\widetilde{\sU}'}$ (cf. equation  \eqref{open_ideal_eqn}) is a homogeneous $C^*$-algebra (cf. Definition \ref{ctr_homo_defn}). From the Remark \ref{ctr_homo_rem} it turns out that there is an open subset $\widetilde \sU \subset \widetilde \sU'$ such that
	$$
	\left.\widetilde{A}\right|_{\widetilde{\sU}}\cong C_0\left( \widetilde\sU \right)\otimes \mathbb{M}_n\left(\C \right). 	
	$$
	From the Remark \ref{unital_rem} and the  Definition \ref{principal_non_defn}  it follows that $\pi \left(1_{C\left(\sX \right) }\right)= 1_{\widetilde A}$, so 
	\be\label{top_a_rep_eqn}
	\begin{split}
		\forall \widetilde h \in C\left(\widetilde \sX \right) \quad \widetilde x \in \widetilde \sU \quad \rep_{\widetilde x}\left(\widetilde h \right)= \widetilde h\left( \widetilde x  \right) 1_{\rep_{\widetilde x}\left(\widetilde A \right) }=\\=\begin{pmatrix}
			\widetilde h\left(\widetilde x\right)   & 0 &\ldots & 0\\
			0& \widetilde h\left(\widetilde x\right)  &\ldots&0 \\
			\vdots& \vdots &\ddots&\vdots\\
			0& 0 & \ldots &\widetilde h\left(\widetilde x\right)  
		\end{pmatrix}.
	\end{split}
	\ee	
	where the map $p: \widetilde\sX \to \sX$ is given by the Proposition \ref{spectrum_covering_finite_prop}.	Any $\widetilde{a} \in \left.\widetilde{A}\right|_{  \widetilde \sU }$ can be represented by a continuous function $\widetilde \sU \to \mathbb{M}_N\left(\C \right)$. Let  $\psi, \vartheta :\widetilde \sU\to \mathbb{M}_n\left(\C \right)$ be given by
	\bean
	\psi\left( { x}\right) = 	
	\begin{pmatrix}
		f\left( { x}\right) & 0 &\ldots & 0\\
		0& 0 &\ldots&0 \\
		\vdots& \vdots &\ddots&\vdots\\
		0& 0 & \ldots &0
	\end{pmatrix},\\
	\vartheta\left( { x}\right) = 	
	\begin{pmatrix}
		0 & 	f\left( { x}\right) &\ldots & 0\\
		0& 0 &\ldots&0 \\
		\vdots& \vdots &\ddots&\vdots\\
		0& 0 & \ldots &0
	\end{pmatrix}.
	\eean
	Let both $\widetilde{p}^\parallel, \widetilde b \in  \left.\widetilde{A}\right|_{{   \widetilde \sU }}$ be represented by both $\psi, \vartheta$ respectively. From $\left.\widetilde{A}\right|_{{   \widetilde \sU }} \subset \widetilde{A}$ it turns out that  both $\widetilde{p}^\parallel$ and $\widetilde b$ can be regarded as elements of $\widetilde{A}$, i.e. $\widetilde{p}^\parallel, ~\widetilde b \in \widetilde{A}$. There are $\widetilde x_0 \in \widetilde\sU$ and a continuous positive $\varphi \in C_0\left( \widetilde\sU\right)_+$ such that $\varphi \left(\widetilde {x}_0 \right)  = 1$ and $\supp \varphi \subset {   \mathcal U }$. Denote by
	$\widetilde{p}^\perp \stackrel{\text{def}}{=} 1_{A }- \widetilde{p}^\parallel\in  \widetilde{A}$. 
	If $\epsilon \in \R$ and $\widetilde{u}_\epsilon \bydef p^\parallel e^{i\epsilon \varphi} + p^\perp\in   \widetilde{A}$ then for any ${x} \in {\mathcal U}$ following condition holds
	\be\label{top_rep_eqn}  
	\begin{split}
		\forall \widetilde{x}\in \widetilde{\sU}\quad 	\rep_{\widetilde{x}}\left( \widetilde{u}_\epsilon\right) = \left\{
		\begin{array}{c l}
			\begin{pmatrix}
				e^{i\epsilon\varphi\left(\widetilde{x} \right)}  & 0 &\ldots & 0\\
				0& 1 &\ldots&0 \\
				\vdots& \vdots &\ddots&\vdots\\
				0& 0 & \ldots &1
			\end{pmatrix} &\widetilde x \in   \widetilde\sU \\ \\
			\rep_{\widetilde{x}}\left(1_{\widetilde A}\right)	 & \widetilde x \in \widetilde\sX \setminus  \widetilde\sU
		\end{array}\right..\end{split}
	\ee
	Otherwise for any $\widetilde {x} \in \widetilde\sU$ one has $\rep_{\widetilde{x}}\left( \widetilde{u}_\epsilon  \right)\rep_{{x}}\left( \widetilde{u}^*_\epsilon  \right) = \rep_{{x}}\left( \widetilde{u}^*_\epsilon  \right) \rep_{\widetilde {x}}\left( \widetilde{u}_\epsilon  \right) = 1_{B\left(\H_{ \widetilde x} \right) }$ where $\H_{ \widetilde x}$ is a space of an irreducible representation $\rep{\widetilde {x}}:\widetilde A\to B\left(\H_{ \widetilde {x}}\right)$. It follows that $\rep_{\widetilde{x}}\left( \widetilde{u}_\epsilon  \right)\rep_{\widetilde {x}}\left( \widetilde{u}^*_\epsilon  \right) = \rep_{{x}}\left( \widetilde{u}^*_\epsilon  \right) \rep_{\widetilde{x}}\left( \widetilde{u}_\epsilon  \right) = 1_{B\left(\H_{\widetilde {x}} \right) }$  for any $\widetilde {x} \in \widetilde {\mathcal X}$.  Hence one has $\widetilde{u}_\epsilon^*\widetilde{u}_\epsilon= \widetilde{u}_\epsilon \widetilde{u}_\epsilon^*=  1_{  \widetilde{A} }$ i.e. $\widetilde{u}_\epsilon$ is unitary.  There is an internal *-automorphism $\chi_\epsilon \in \Aut\left(\widetilde{A} \right)$ given by $\widetilde{a}  \mapsto \widetilde{u}_\epsilon^*\widetilde{a}\widetilde{u}_\epsilon$.  From
	$$
	\rep_{{x}_0}\left( \widetilde{u}^*_\epsilon\right)  \vartheta\left(x_0 \right)  \rep_{{x}_0}\left( \widetilde{u}_\epsilon\right) = \begin{pmatrix}
		0 & 	e^{-i\epsilon} &\ldots & 0\\
		0& 0 &\ldots&0 \\
		\vdots& \vdots &\ddots&\vdots\\
		0& 0 & \ldots &0
	\end{pmatrix}
	$$
	it turns out that
	\be\label{top_infe_eqn}
	\begin{split}
		\forall\epsilon_1 ,\epsilon_2\in [0, 2\pi)\quad \epsilon_1 \neq \epsilon_2   \Rightarrow \rep_{ {x}_0}\left(\chi_{\epsilon_1}\left( \widetilde{b}\right)  \right) \neq \rep_{ {x}_0}\left(\chi_{\epsilon_2}\left( \widetilde{b}\right)  \right)\Rightarrow \\\Rightarrow \quad \chi_{\epsilon_1} \neq \chi_{\epsilon_2}.
	\end{split}
	\ee
	From the equation \eqref{top_a_rep_eqn} it follows that 
	\bean
	\forall \widetilde h \in C\left(\widetilde \sX \right) \quad \forall \widetilde x \in   \widetilde{\mathcal X} \quad \rep_{\widetilde{x}}\left( \widetilde{u}_\epsilon  \right)\rep_{\widetilde{x}}\left( \widetilde h\right)\rep_{\widetilde{x}}\left( \widetilde{u}^*_\epsilon  \right)= \rep_{\widetilde{x}}\left( a\right).
	\eean 
	It means that $\chi_{\epsilon}\left( a\right) = a$ for all $a \in C\left(\sX\right)$
	hence from \eqref{top_infe_eqn} it turns out that the  group
	\be\nonumber
	H = \left\{ g \in \Aut\left(\left.\widetilde{A} \right)~\right|~ g \left.\Id_{\widetilde A}\right|_{C\left( \widetilde\sX\right)}\left( \widetilde h \right) =  \left.\Id_{\widetilde A}\right|_{C\left( \widetilde\sX\right)}\left( \widetilde h \right)\quad\forall\widetilde h \in C\left(\widetilde \sX\right)\right\}
	\ee
	is not finite. So $G$ is also infinite  i.e. one has the contradiction with (a) of the Definition \ref{fin_pre_defn}. From this contradiction it turns out that $\dim\rep_{ {x}}\left(\widetilde{A}\right)= 1$ for any ${x} \in {\sX}$, it means that $\widetilde{A}$ is commutative. So $\widetilde{A} = C\left( \widetilde  \sX\right)$. 
\end{proof}

The algebraic construction of nonunital finite-fold coverings requires the following definition.
\begin{definition}\label{top_covering_compactification_defn}
	A   covering $p: \widetilde{   \mathcal X } \to \mathcal X$ is said to be a \textit{ covering with compactification} if there are compactifications ${   \mathcal X } \hookto {   \mathcal Y }$ and $\widetilde{   \mathcal X } \hookto \widetilde{   \mathcal Y }$ such that:
	\begin{itemize}
		\item there is a covering $\widetilde{p}:\widetilde{   \mathcal Y }\to {   \mathcal Y }$,
		\item the covering $p$ is the restriction of $\widetilde{p}$, i.e. $p = \widetilde{p}|_{\widetilde{   \mathcal X }}$.
	\end{itemize}
\end{definition}
\begin{example}
	
	Let $g: S^1 \to S^1$ be an $n$-fold covering of a circle. Let $\mathcal X \cong \widetilde{\mathcal X} \cong S^1 \times \left[0,1\right)$. 
	The map
	\begin{equation*}
		\begin{split}
			p: \widetilde{   \mathcal X } \to \mathcal X,\\
			p = g \times \Id_{\left[0,1\right)}
		\end{split}
	\end{equation*}
	is an $n$-fold covering. If $\mathcal Y \cong \widetilde{\mathcal Y} \cong S^1 \times \left[0,1\right]$ then a compactification $\left[0,1\right) \hookto \left[0,1\right]$ induces  compactifications  $\mathcal X \hookto\mathcal Y$, $\widetilde{   \mathcal X } \hookto \widetilde{   \mathcal Y }$. The map
	\begin{equation*}
		\begin{split}
			\widetilde{p}: \widetilde{   \mathcal Y } \to \mathcal Y,\\
			\widetilde{p} = g \times \Id_{\left[0,1\right]}
		\end{split}
	\end{equation*}
	is a covering  such that $\widetilde{p}|_{\widetilde{   \mathcal X }}=p$. So if $n > 1$ then $p$ is  a nontrivial covering with  compactification.\\
\end{example}
\begin{example}
	Let $\mathcal X  = \C \setminus \{0\}$ be a complex plane with punctured 0, which is parameterized by the complex-valued variable $z\in \C$. Let  $\mathcal X \hookto\mathcal Y$ be any compactification. If both $\left\{z'_\a \in \mathcal X\right\}$, 	$\left\{z''_\a \in \mathcal X\right\}$ are  nets such that $\lim_{\a}\left|z'_\a\right|=\lim_\a\left|z''_\a\right| = 0$ then form $\lim_{\a}\left|z'_\a-z''_\a\right|= 0$ it turns out 
	\begin{equation}\label{x_0_eqn}
		x_0 = \lim_{\a} z'_\a = \lim_{\a} z''_\a \in \mathcal Y.
	\end{equation}
	If $\widetilde{   \mathcal X } \cong \mathcal X$ then for any $n \in \N$ such that $n > 1$ there is an $n$-fold covering 
	\begin{equation*}
		\begin{split}
			p: \widetilde{   \mathcal X } \to \mathcal X,\\
			z \mapsto z^n.
		\end{split}
	\end{equation*}
	If both  $\mathcal X \hookto\mathcal Y$, $\widetilde{   \mathcal X } \hookto \widetilde{   \mathcal Y }$ are compactifications, and $\widetilde{p}: \widetilde{   \mathcal Y } \to \mathcal Y$ is a covering such that $\widetilde{p}|_{\widetilde{   \mathcal X }} = p$ then from \eqref{x_0_eqn} it turns out $\widetilde{p}^{-1}\left(x_0 \right)= \left\{\widetilde{x}_0\right\}$ where $\widetilde{x}_0$ is the unique point such that following conditions hold:
	\begin{equation*}
		\begin{split}
			\widetilde{x}_0 = \lim_{\a} \widetilde{z}_\a \in \widetilde{   \mathcal Y },\\
			\lim_{\a}\left|\widetilde{z}_\a\right|= 0.
		\end{split}
	\end{equation*}
	It turns out $\left|\widetilde{p}^{-1}\left(x_0 \right) \right|=1$. However $\widetilde{p}$ is an $n$-fold covering and if $n >1$ then  $\left|\widetilde{p}^{-1}\left(x_0 \right) \right|=n>1$. It contradicts with $\left|\widetilde{p}^{-1}\left(x_0 \right) \right|=1$, and from the contradiction it turns out that for any $n > 1$ the map  $p$ is not a covering with  compactification.
\end{example}

The following lemma supplies the quantization of coverings with compactification.

\begin{lemma}\label{comm_fin_lem}  
	Let $\mathcal X$ be a locally compact,  Hausdorff space and there is an injective $*$-homomorphism $\pi: C_0\left( \mathcal X\right) \hookto \widetilde{A}$ of $C^*$-algebras	then following conditions are equivalent:
	\begin{enumerate}
		\item [(i)] the quadruple $\left(C_0\left(\mathcal  X \right), \widetilde{A}, G  , \pi  \right)$ is a  noncommutative finite-fold covering with unitization,
		\item[(ii)] the $C^*$-algebra $\widetilde{A}$ is commutative, moreover if $\mathcal X$ (resp. $\widetilde{\mathcal X}$) is the spectrum of $C_0\left( \sX\right) $ (resp. $\widetilde{A}$) then a given by the Proposition \ref{spectrum_covering_finite_prop} map  $p:\widetilde{\mathcal X}\to {\mathcal X}$ is a transitive topological finite-fold covering with compactification (cf. Definition \ref{top_covering_compactification_defn}) such that the quadruple $\left(C_0\left(\mathcal  X \right), \widetilde{A}, G, \pi   \right)$ is naturally equivalent to $$\left(C_0\left(\mathcal  X \right), C_0\left(\widetilde{\mathcal X} \right), G\left(\left. \widetilde{\sX}~\right|~ {\sX}\right) , C_0\left( p\right)     \right)$$
		where $*$-homomorphism $C_0\left(p \right): C_0\left(\mathcal  X \right)\to C_0\left(\widetilde{\mathcal X} \right)$ is $C_0$ is a finite covering algebraic functor (cf. Definition \ref{top_c_funct_defn}).  
	\end{enumerate}
\end{lemma}

\begin{proof}
	
	(i)$\Rightarrow$(ii)  	  From the Definition \ref{fin_unitization_defn} it turns out that following conditions hold:
	\begin{enumerate}
		\item[(a)] 
		there are unital $C^*$-algebras $B$, $\widetilde{B}$  and inclusions 
		$C_0\left(\mathcal  X \right) \subset B$,  $\widetilde{A}\subset \widetilde{B}$ such that both  $C_0\left(\mathcal  X \right) $ (resp. $ \widetilde{A}$) are an essential ideals of both $B$ (resp. $\widetilde{B}$),
		\item[(b)] there is an 
		unital  noncommutative finite-fold covering $\left(B ,\widetilde{B}, G, \widetilde{\pi} \right)$ such that $\pi = \widetilde{\pi}|_{C_0\left(\mathcal  X \right)}$ and the action $G \times\widetilde{A} \to \widetilde{A}$ is induced by $G \times\widetilde{B} \to \widetilde{B}$.
	\end{enumerate}
	
	The unital  $C^*$-algebra $B$ is a subalgebra  of a commutative multiplier algebra $M\left(C_0\left(\sX \right) \right)$ (cf. Definition \ref{multiplier_defn}), so the $C^*$-algebra $B$ is also commutative. From the Theorem \ref{gelfand-naimark_thm} it follows that there is a compact Hausdorff space $\sX$ such that $B = C\left(\sX\right)$. From the Corollary \ref{comm_comm_cor} it follows that there is a compact Hausdorff space $\widetilde{\sX}$ and a (topological) transitive finite-fold covering $\widetilde{p}: \widetilde{\sX}\to\sX$ such that the quadruple $\left(B ,\widetilde{B}, G, \widetilde{\pi} \right)$ is equivalent to 
	$$\left(C\left({\mathcal Y} \right) , C\left( \widetilde{\mathcal Y}\right) , G\left(\left.\widetilde{\sY}~\right|\sY\right) , C_0\left(\widetilde p \right)  \right).$$
	Both $C_0\left(\mathcal  X \right)$ and $C_0\left(\widetilde \sX\right)$ are essential ideals, of $C\left(\mathcal  Y \right)$ and $\widetilde C\left(\widetilde \sY\right)$, hence from the Example \ref{comm_ess_exm} it follows that both $\sX$ and $\widetilde \sX$ are open dense subsets of both $\sX$ and $\widetilde \sX$. It means that both  inclusions $\sX\hookto\sX$ and $\widetilde \sX\hookto \widetilde \sX$ are compactifications. 
	From $\pi = \widetilde{\pi}|_{C_0\left(\mathcal  X \right)}$  it follows that $\widetilde{p}^{-1}\left({\sX}\right) = \widetilde{\sX}$ so $p = \left.\widetilde{p}\right|_{\widetilde{\sX}}: \widetilde{\sX}\to \sX$ is a covering. From transitivity of action  $G \times \widetilde{\mathcal Y} \to \widetilde{\mathcal Y}$ and the condition $G \times C_0\left( \widetilde{\mathcal X}\right)=  C_0\left( \widetilde{\mathcal X}\right)$ it turns out that the action  $G \times \widetilde{\mathcal X} \to \widetilde{\mathcal X}$ is transitive. From $p = \left.\widetilde{p}\right|_{\widetilde{\sX}}$ it follows that $p$ is a covering with compactification.\\
	(ii)$\Rightarrow$(i)
	The map  $p:\widetilde{\mathcal X}\to {\mathcal X}$ is a transitive topological finite-fold covering with compactification, hence there are compactifications $\sX \hookto \sX$,  $\widetilde{\mathcal X} \hookto \widetilde{\mathcal Y}$ and the transitive finite-fold covering $\widetilde{p}:\widetilde{\mathcal Y}\to \sX$ such that $p = \left.\widetilde{p}\right|_{\widetilde{\mathcal X}}$.  Firstly we prove that the quadruple  $\left(C_0\left(\mathcal  X \right), C_0\left(  \widetilde{   \mathcal X }\right) ,G\left(\left. \widetilde{   \mathcal X } ~\right| {   \mathcal X }\right),   \pi\right)$ is a finite-fold noncommutative pre-covering, i.e. it satisfies to the conditions  (a)  and (b)  of the Definition \ref{fin_pre_defn}.
	\begin{enumerate}
		\item[(a)]
		From 
		\be\nonumber
		G = \left\{ \left.g \in \Aut\left(C_0\left(\widetilde{   \mathcal X }\right) \right)~\right|~ ga = a;\quad\forall a \in C_0\left({   \mathcal X }\right) \right\}\cong G\left(\left. \widetilde{   \mathcal X } ~\right| {   \mathcal X }\right),
		\ee
		and taking into account that $p:\widetilde{\mathcal X}\to {\mathcal X}$ is a (topological) finite-fold covering one concludes that the group 	$G\cong G\left(\left. \widetilde{   \mathcal X } ~\right| {   \mathcal X }\right)$ is finite.
		\item[(b)] From $\sX \cong \widetilde{   \mathcal X }/G\left(\left. \widetilde{   \mathcal X } ~\right| {   \mathcal X }\right)$ it follows that 	$C_0\left({   \mathcal X }\right) \cong C_0\left(\widetilde{   \mathcal X }\right)^{G\left(\left. \widetilde{   \mathcal X } ~\right| {   \mathcal X }\right)}$.
	\end{enumerate}
	Secondly we prove that 
	$\left(C_0\left(\mathcal  X \right), C_0\left(\widetilde{\mathcal X} \right), G\left( \widetilde{\mathcal X}~|~ {\mathcal X}\right) , C_0\left( p\right)     \right)$
	satisfies to the conditions (a), (b) of the Definition \ref{fin_unitization_defn}.
	\begin{enumerate}
		\item [(a)] Both maps $\sX \hookto \sY$,  $\widetilde{\mathcal X} \hookto \widetilde{\mathcal Y}$ are compactifications, hence both inclusions  $C_0\left(\sX \right)\hookto C\left(\sY \right)$ and   $C_0\left( \widetilde{\mathcal X}\right) \hookto C\left( \widetilde{\mathcal Y}\right) $ correspond to essential ideals (cf. Example \ref{comm_ess_exm}). It means that the inclusions are unitizations.
		\item[(b)] From the Theorem \ref{pavlov_troisky_thm} it follows that $C_0\left( \widetilde{\mathcal Y}\right)$ is a finitely generated projective module, it follows that the quadruple
		$$
		\left(C\left( {\mathcal Y}\right) , C\left( \widetilde{\mathcal Y}\right) , G, C_0\left(  \widetilde{p}\right) \right)
		$$
		is an  unital noncommutative finite-fold     covering (cf. Definition  \ref{fin_unital_defn}). From $p = \left.\widetilde{p}\right|_{\widetilde{\mathcal X}}$ it follows that $\left.\widetilde{\pi}\right|_{C_0\left(\sX \right)} = \pi$. The action  $G \times \widetilde{\mathcal X} \to   \widetilde{\mathcal X}$ is induced by the action $G \times \widetilde{\mathcal Y} \to   \widetilde{\mathcal Y}$, hence  the action $G \times C_0\left( \widetilde{\mathcal X}\right) \to  C_0\left( \widetilde{\mathcal X}\right)$ is induced by the action $G \times C\left( \widetilde{\mathcal Y}\right) \to  C\left( \widetilde{\mathcal Y}\right)$. 
	\end{enumerate}
\end{proof}
\begin{lemma}\label{top_conn_fam_lem}  
	Let $p:\widetilde \sX \to \sX$ be a transitive finite-fold covering where both $\sX$ and $\widetilde \sX$ are connected, locally connected, locally compact, Hausdorff spaces. There is  $\left\{\sU_\la\subset \sX\right\}_{\la\in \La}$  a family of connected open subsets such that:
	\begin{enumerate}
		\item[(a)] for any $\la \in \La$ the closure of $\sU_\la$ is compact,
		\item[(b)] for every $\la \in \La$ the preimage $p^{-1}\left(\sU_\la\right)$ is connected,
		\item[(c)] $\sX = \cup_{\la\in\La} \sU_\la$,
		\item[(d)] for all $\la\in\La$ there are  isomorphisms $ G\left(\left.\widetilde\sX\right| \sX\right)\cong  G\left(\left.\widetilde\sU_\la\right| \sU_\la\right)$.
	\end{enumerate}
\end{lemma}
\begin{proof}
	Since $\widetilde \sX$ is locally compact there is a family $\left\{\widetilde \sV_\a\subset \widetilde\sX\right\}_{\a\in \mathscr A}$ of open sets such that 
	\begin{enumerate}
		\item the closure of $\widetilde \sV_\a$ is compact.
		\item $\widetilde\sX = \cup_{\a\in \mathscr A} \widetilde \sV_\a$,
		\item for each $\a \in \mathscr A$ a set $\widetilde \sV_\a$ is mapped homeomorphically onto $p\left(\widetilde \sV_\a\right)$.
	\end{enumerate}
	Let $\widetilde x_0\in \widetilde \sX$ be any point and $\a_0\in \mathscr A$ be such that $\widetilde x_0\in \widetilde \sV_{ \a_0}$. For any nontrivial $g \in G\left(\left.\widetilde\sX\right| \sX\right)$ we select $\a_g\in \mathscr A$ such that $g\widetilde x_0 \in \widetilde \sV_{ \a_g}$. From the Lemma \ref{top_gen_path_lem} it follows that for each  nontrivial $g \in G\left(\left.\widetilde\sX\right| \sX\right)$ there is $\left\{\widetilde\sV_\a\right\}$-{path} $\left(\widetilde\sV_{\a^g_1},...,\widetilde\sV_{\a^g_{n_g}}\right)$ (cf. Definition \ref{top_gen_path_defn}) such that $\widetilde\sV_{\a^g_1}=\widetilde\sV_{\a_0}$ and $\widetilde\sV_{\a^g_{n_q}}=\widetilde\sV_{\a_g}$. Since the union  $$\widetilde \sU' \bydef \bigcup_{j = 1, ..., n_g} \widetilde \sV_{\a^g_{j_g}}$$ is finite its closure is compact. 
	Similarly the closure of the finite union  $\widetilde\sU\bydef \bigcup_{g \in G\left(\left.\widetilde\sX\right| \sX\right)}g \widetilde \sU'$ is compact. For all $g',g'' \in G\left(\left.\widetilde\sX\right| \sX\right)$ one 
	has $\widetilde x_0 \in g' \widetilde \sU'\cap g'' \widetilde \sU'$ it follows that the set  $\widetilde \sU$ is connected. We put $\widetilde \sU = \widetilde \sU_{\la_{\text{min}}}$ where $\la_{\text{min}}\in \La$ is the minimal element. For any $\la\in \La$ there is a finite subset $\mathscr A_0 \subset \mathscr A$ such that a union    $\widetilde\sU'_\la \bydef \widetilde\sU' \bigcup \cup_{\a \in \mathscr A_0}\widetilde \sV_\a$ is connected. The union is finite and any element of the union has a compact closure, so the closure of $\widetilde\sU'_\la$ is compact. A finite union  $\widetilde\sU_\la\bydef \bigcup_{G\left(\left.\widetilde\sX\right| \sX\right)}\widetilde\sU'_\a$ is  connected and have a compact closures. Let $\sU_\la \bydef p\left( \widetilde \sU'_\la\right)=  p\left( \widetilde \sU_\la\right)$ ad let us prove that the family satisfies to conditions (a)-(d) of this lemma.
	\begin{enumerate}
		\item[(a)] For each $\la\in \La$ the set $\sU_\la$ is a finite union of a sets with compact closures.
		\item[(b)] One has $p^{-1}\left(\sU_\la\right)= \widetilde \sU_\la$ and according to our construction the set $\widetilde \sU_\la$ is connected.
		\item[(c)] If $x \in \sX$ is any point then there is  $\widetilde x \in \widetilde \sX$ such that $x = p\left(\widetilde x\right)$. From the Lemma \ref{top_gen_path_lem} it follows that there is  $\left\{\widetilde\sV_\a\right\}$-{path} $\left(\widetilde\sV_{\a_0},...,\widetilde\sV_{\a_n}\right)$ such that $\widetilde x \in \widetilde\sV_{\a_n}$. On the other hand there is $\la_{\widetilde x} \in \La$ such that $\widetilde\sU'_{\la_{\widetilde x}}= \widetilde \sU' \cup \widetilde \sV_{\a_0}\cup ...\cup  \widetilde \sV_{\a_n}$. It turns out that $\widetilde x \in \widetilde\sU'_{\la_{\widetilde x}}$ and $x \in \sU_{\la_{\widetilde x}}$. One has $\sX = \cup_{\la\in\La} \sU_\la$ because the point $x$ is arbitrary.
		\item[(d)] If $\widetilde x_0 \in \widetilde\sU_\la$ then one has
		$$
		p^{-1}\circ p\left(\widetilde x_0 \right) \cong G\left(\left.\widetilde\sX\right| \sX\right)\cong  G\left(\left.\widetilde\sU_\la\right| \sU_\la\right).
		$$
	\end{enumerate}
\end{proof}

\begin{lemma}\label{top_fin_sufficient_lem}
	Let $\sX$ be  a connected, locally connected (cf. Definition \ref{top_locally_connected_defn}), locally compact, Hausdorff space. If $p: \widetilde \sX \to \sX$  is finite-fold transitive covering and a $*$-homomorphism  $C_0\left(p \right):  C_0\left(  {\sX }\right) \hookto C_0\left(  \widetilde{\sX }\right)$ is given by the equation  \ref{top_c0_eqn} then the quadruple
	$$
	\left(C_0\left(\sX \right), C_0\left(  \widetilde{\sX }\right) ,G\left(\left. \widetilde{   \mathcal X} ~\right| {   \mathcal X }\right),   C_0\left(p \right)\right)
	$$
	is a noncommutative finite-fold covering (cf. Definition \ref{fin_defn}).
\end{lemma}
\begin{proof}
	If $\left\{\sU_\la\right\}_{\la\in\La}$ is a family of given by the Lemma \ref{top_conn_fam_lem} open subsets then one has:
	\begin{itemize}
		\item 
		$$
		\sX = \bigcup_{\la\in \La}\sU_\la,
		$$
		\item a closure of $\sU_\la$ is compact for each $\la\in \La$,
		\item
		\bean
		\forall \mu, \nu \in \La  \quad \sU_\mu \cap \sU_\nu \neq \emptyset,
		\eean
		\item
		\bean
		\forall \mu, \nu \in \La \quad \mu \le \nu\quad\Leftrightarrow \quad \sU_\mu \subset \sU_\nu,
		\eean
		\item For each $\la\in \La$ both sets $\sU_\la$ and if $\widetilde \sU_\la\bydef p^{-1}\left(\sU_\la\right)$ are connected (cf. Definition \ref{top_connected_defn}), so both $C^*$-algebras $C_0\left(\sU_\la\right)$ and $C_0\left(\sU_\la\right)$ are connected $C^*$-algebras (cf. Definition \ref{connected_c_a_defn}).
	\end{itemize}
	If $\sV_\la$ is the closure of $\sU_\la$ then the preimage $\widetilde \sV_\la \bydef p^{-1}\left(\sV_\la  \right)$ is compact  (cf. Lemma \ref{top_compact_preimage_lem}). There is a transitive covering $\widetilde \sV_\la\to  \sV_\la$, so $\widetilde \sU_\la\to \sU_\la$ is a covering with compactification (cf. Definition \ref{top_covering_compactification_defn}). From the Lemma \ref{comm_fin_lem}  it follows that 
	$$
	\left(C_0\left(\mathcal  U_\la \right), C_0\left(  \widetilde{   \mathcal U}_\la \right) ,G\left(\left. \widetilde{   \mathcal X} ~\right| {   \mathcal X }\right),   \left.C_0\left(p \right)\right|_{C_0\left(\mathcal  U_\la \right)}\right) 	
	$$
	is a noncommutative finite-fold covering with unitization (cf. Definition \ref{fin_unitization_defn}), i.e. the hereditary subalgebra $C_0\left(\mathcal  U_\la \right)\subset C_0\left( \sX \right)$ is  $\left(C_0\left(\sX \right), C_0\left(  \widetilde{\sX }\right) ,G\left(\left. \widetilde{   \mathcal X} ~\right| {   \mathcal X }\right),   C_0\left(p \right)\right)$-{strictly proper} (cf. Definition \ref{strictly_proper_defn}). If $f \in C_0\left( \sX\right)$ then from the Definition \ref{c_c_closure_defn} it follows that for any $\eps > 0$ there is $f' \in C_c\left( \sX\right)$ such that $\left\| f - f'\right\| < \eps$. However the support $\supp f'$ of $f'$ (cf. Definition \ref{top_support_defn}) is compact, so from $\bigcup_{\la\in\La}\sU_\la = \sX$ it follows that there is $\la_0 \in \La$ such that $\supp f' \subset \sU_{\la_0}$. If turns out that $f' \in C_0\left(\sU_{\la_0} \right)\subset \bigcup_{\la\in\La} C_0 \left( \sU_\la\right)$.
\end{proof}

\begin{lemma}\label{top_fin_necassary_lem} 
	If   $\sX$ is a connected,  locally connected, locally compact, Hausdorff space, a quadruple $\left(C_0\left(\mathcal  X \right), \widetilde{A}, G,    \pi\right)$ is a noncommutative finite-fold covering (cf. Definition \ref{fin_defn}) and $p: \widetilde{   \mathcal X } \to \sX$ is the given by the Proposition \ref{spectrum_covering_finite_prop} continuous map from the spectrum $\widetilde \sX$ of $\widetilde A$ to the spectrum $\sX$ of $A$ then following conditions hold:
	\begin{enumerate}
		\item [(i)] $\widetilde A$ is a commutative $C^*$-algebra and $\widetilde A \cong C_0\left( \widetilde \sX\right)$, 
		\item[(ii)] the map $p: \widetilde{   \mathcal X } \to \sX$ is a transitive finite-fold covering, 
		\item[(iii)]  the quadruple   $\left(C_0\left(\mathcal  X \right), \widetilde{A}, G,    \pi\right)$ is equivalent to $$\left(C_0\left(\mathcal  X \right), C_0\left(  \widetilde{   \mathcal X }\right) ,G\left(\left. \widetilde{   \mathcal X } ~\right| {   \mathcal X }\right),   C_0\left(p \right)\right)$$ one,  where   $C_0$ is a finite covering algebraic functor (cf. Definition \ref{top_c_funct_defn}).
	\end{enumerate}
\end{lemma}
\begin{proof}
	Let $\left\{A_\la \subset C_0\left\{\sX\right\}\right\}_{\la\in\La}$ be a required by the Definition \ref{fin_defn} family of \\ $\left(C_0\left(\mathcal  X \right), \widetilde{A}, G,    \pi\right)$-{strictly proper} $C^*$-subalgebras (cf. the Definition \ref{strictly_proper_defn}). Let us fix $\la\in\La$. Since $A_\la$ is a hereditary subalgebra of $C_0\left(\mathcal  X \right)$ from the Theorem \ref{gelfand-naimark_thm} and  Proposition \ref{hered_spectrum_prop} it follows that there is an open subset $\sU_\la\subset \sX$ such that $A_\la\cong C_0\left(\sU_\la\right)$. Let $\widetilde A_\la\subset \widetilde A$ be   a {hereditary} $\left(C_0\left(\mathcal  X \right), \widetilde{A}, G, \pi \right)$-{lift} of $A_\la\cong C_0\left(\sU_\la\right)$ (cf. the Definition \ref{hereditary_lift_defn}). From the Definitions (\ref{hereditary_lift_defn} and \ref{strictly_proper_defn}) it follows that $\left(C_0\left(\mathcal  U_\la \right), \widetilde{A}_\la, G,    \pi|_{C_0\left(\mathcal  U_\la \right)}\right)$ is a noncommutative finite-fold covering with unitization. From the Lemma \ref{comm_fin_lem}  it follows that the $C^*$-algebra $\widetilde{A}_\la$ is commutative, i.e. $\widetilde{A}_\la\cong C_0\left(\widetilde\sU_\la\right)$ and a restriction $\left.p\right|_{\widetilde\sU_\la}: \widetilde\sU_\la\to\sU_\la$ is a transitive finite-fold covering.\\
	(i)
	Form  Lemma \ref{fin_def_lem} it follows that a  union  $\cup \widetilde{A}_\la \cong \cup  C_0\left(\widetilde\sU_\la\right)$ is dense in  $\widetilde A$, so $\widetilde A$ is a commutative $C^*$-algebra. Taking into account the Theorem \ref{gelfand-naimark_thm} one has the natural $*$-isomorphism $$\widetilde A \cong C_0\left( \widetilde \sX\right).$$
	(ii) If  $x_0\in   \sX\setminus \bigcup_{\la\in \La}\sU_\la$  is any point then there if $f \in C_0\left( \sX\right)$ such that 
	\bean
	f\left(x_0  \right) = 1,\\
	\forall f'\in \bigcup_{\la\in \La}C_0\left( \sU_\la\right) \quad \left\|f'-f \right\| \ge 1,
	\eean 
	i.e. the union $\bigcup_{\la\in \La}C_0\left( \sU_\la\right)$ is not dense in $C_0\left( \sX\right)$. This fact contradicts with the Definition \ref{fin_defn}, so there is $\la_0\in \La$ such that $x_0\in \sU_{\la_0}$. However the restriction  $\left.p\right|_{\widetilde\sU_{\la_0}}: \widetilde\sU_{\la_0}\to\sU_{\la_0}$ is a transitive finite-fold covering, so there is an open neighborhood $\sV \subset \sU_{\la_0}$ of $x_0$ evenly covered by $\left.p\right|_{\widetilde\sU_{\la_0}}$ (cf. Definition \ref{top_covering_defn}). However  $\sV$ is an open subset of $\sX$ evenly covered by $p$, so $p$ is a covering (cf. Definition \ref{top_covering_defn}). From $\sX = \widetilde \sX/G$ it follows that the covering is transitive.\\
	(iii) We already know that $\widetilde A \cong C_0\left( \widetilde \sX\right)$ and $G \cong G\left(\left. \widetilde{   \mathcal X } ~\right| {   \mathcal X }\right)$. From the Lemma  \ref{comm_fin_lem} it follows that the noncommutative finite-fold covering with unitization
	$\left(C_0\left(\mathcal  U_\la \right), \widetilde{A}_\la, G,    \pi|_{C_0\left(\mathcal  U_\la \right)}\right)$ is equivalent to 
	$$
	\left(C_0\left(\mathcal  U_\la \right), C_0\left( \widetilde{\sU}_\la\right) ,G\left(\left. \widetilde{   \mathcal X } ~\right| {   \mathcal X }\right),    C_0\left(p|_{\mathcal  U_\la} \right)=\left.C_0\left(p \right)\right|_{C_0\left(\mathcal  U_\la \right)} \right)$$
	i.e. $ \pi|_{C_0\left(\mathcal  U_\la \right)}=\left.C_0\left(p \right)\right|_{C_0\left(\mathcal  U_\la \right)}$
 So one has $\pi = C_0\left(p \right)$ since the union $\bigcup_{\la\in \La}C_0\left( \sU_\la\right)$ is not dense in $C_0\left( \sX\right)$. 
\end{proof}

\begin{theorem}\label{top_finite_covering_thm} 
Following conditions hold.
\begin{enumerate}
	\item [(i)] 	If $\sX$ is a connected, locally connected (cf. Definition \ref{top_locally_connected_defn}), locally compact, Hausdorff space,  then    the {category $\mathfrak{FinCov}$-$C_0\left(\sX \right)$ of finite-fold coverings} of $C_0\left(\sX \right)$ (cf. Definition \ref{fin_category_defn}) is composable (cf. Definition \ref{fin_composable_cat_defn})
	\item[(ii)] The given by the Definition \ref{top_c_funct_defn}  finite covering algebraic functor   $C_0:\mathfrak{FinCov}$-$\sX \cong \mathfrak{FinCov}$-$C_0\left(\sX \right)$ is the natural equivalence of the categories (cf. Definition \ref{category_equivalence_definition}).
	
\end{enumerate}
\end{theorem}

\begin{proof}
	 From the Lemmas    \ref{top_fin_sufficient_lem} and \ref{top_fin_necassary_lem} it follows that the functor $C_0$ yields a one-to-one correspondence between objects of $\mathfrak{FinCov}$-$\sX$ and objects of  $\mathfrak{FinCov}$-$C_0\left(\sX \right)$ given by
	$$
	\left(p:\widetilde{\sX} \to \sX \right) \mapsto 	C_0\left( p\right) :   C_0\left( \sX\right)    \hookto C_0\left( \widetilde{\sX}\right)
	$$
	where $C_0\left(p\right)$ is given by the Definition \ref{top_c_funct_defn}.\\
	(i)  Any $\mathfrak{FinCov}$-$C_0\left(\sX \right)$-morphism from $C_0\left( p'\right) :   C_0\left( \sX\right)    \hookto C_0\left( \widetilde{\sX}'\right)$ to $C_0\left( p''\right) :   C_0\left( \sX\right)    \hookto C_0\left( \widetilde{\sX}''\right)$ is an injective $*$-homomorphism $\pi: C_0\left( \widetilde{\sX}''\right) \hookto C_0\left( \widetilde{\sX}'\right)$ 
	such that the following diagram 
	\newline
	\begin{tikzpicture}
		\matrix (m) [matrix of math nodes,row sep=3em,column sep=4em,minimum width=2em]
		{
			C_0\left( \widetilde{\sX}'\right)  & &C_0\left( \widetilde{\sX}''\right) \\ 
			& C_0\left( {\sX}\right)  & \\};
		\path[-stealth]
		(m-1-3) edge node [above] {$\pi$} (m-1-1)
		(m-2-2) edge node [left]  {$C_0\left( p'\right)~~$} (m-1-1)
		(m-2-2) edge node [right] {$~~C_0\left( p''\right)$} (m-1-3);
	\end{tikzpicture}
	\\ 	
	is commutative. Since $\pi$ is injective it corresponds to a continuous surjective map $p: \widetilde{\sX}'\to \widetilde{\sX}''$ such that the following diagram 
		\newline
	\begin{tikzpicture}
		\matrix (m) [matrix of math nodes,row sep=3em,column sep=4em,minimum width=2em]
		{
			\widetilde{\mathcal X}' & &\widetilde{\mathcal X}''\\ 
			& {\mathcal X}\\};
		\path[-stealth]
		(m-1-1) edge node [above] {$p$} (m-1-3)
		(m-1-1) edge node [left]  {$p'~~$} (m-2-2)
		(m-1-3) edge node [right] {$~~p''$} (m-2-2);
	\end{tikzpicture}
	\\
From the Theorem \ref{top_covp_cat_thm}	it follows that the map $p$ is a transitive covering, and taking into account the Lemma \ref{top_fin_sufficient_lem} we conclude that the quadruple 
	$$
\left(C_0\left(\widetilde\sX'' \right), C_0\left(  \widetilde{\sX }'\right) ,G\left(\left. \widetilde\sX' ~\right| \widetilde\sX''\right),   C_0\left(p \right)\right)$$
	is a noncommutative finite-fold covering (cf. Definition \ref{fin_defn}). It turns out that the category $\mathfrak{FinCov}$-$C_0\left(\sX \right)$ satisfies to the Definition \ref{fin_composable_cat_defn} .\\
	(ii) We already know that  the functor $C_0$ yields a one-to-one correspondence between objects of $\mathfrak{FinCov}$-$\sX$ and objects of  $\mathfrak{FinCov}$-$C_0\left(\sX \right)$. From the proof of (i) it follows that any $\mathfrak{FinCov}$-$C_0\left(\sX \right)$ morphism yields $\mathfrak{FinCov}$-$\sX$-morphism. Conversely from the Lemma    \ref{top_fin_sufficient_lem} it turns out that any $\mathfrak{FinCov}$-$\sX$ morphism uniquely defines a $\mathfrak{FinCov}$-$C_0\left(\sX \right)$-morphism. These correspondences between morphisms of categories  $\mathfrak{FinCov}$-$\sX$ and $\mathfrak{FinCov}$-$C_0\left(\sX \right)$ are mutually inverse.

\end{proof}
\begin{remark}\label{comm_comm_rem}  
The Theorem  \ref{top_finite_covering_thm}  yields a pure algebraic definition of the fundamental group (cf. the Theorem \ref{comm_uni_lim_thm} and the Theorem \ref{comm_uni_lim_thm}). Alternative theories (e.g. described in \cite{clarisson:phd}) which allow coverings of commutative $C^*$-algebras by noncommutative ones do not have this property.
\end{remark}

\begin{corollary}\label{top_consents_comp_lem} 
	Let  $\mathcal X$ is a locally compact, connected, locally connected, Hausdorff space. If $\sX$ then the $C^*$-algebra $C_0\left(\sX \right)$ is {proper with respect to covering compositions} (cf. Definition \ref{proper_composition_defn}).
\end{corollary}

\begin{corollary}\label{top_fin_cor}
	Let  $\mathcal X$ be a locally compact, connected, locally connected, Hausdorff space, $\widetilde{   \mathcal X }$ and  $p: \widetilde{   \mathcal X } \to \sX$ is a finite-fold transitive topological covering with connected $\widetilde \sX$. If the map  $\lift_p : C_b\left(\sX\right)\hookto C_b\left( \widetilde \sX\right)$ is given by  the Lemma \ref{top_lift_bounded_lem} 	then the quadruple
	$$
\left(C_b\left(\mathcal  X \right), C_b\left(  \widetilde{   \mathcal X }\right) ,G\left(\left. \widetilde{   \mathcal X } ~\right| {   \mathcal X }\right),   \lift_p \right)
	$$
is a connected noncommutative  finite-fold weak covering (cf. Definition \ref{fin_weak_defn}).
\end{corollary}
\begin{proof}
	It is well known that both $C_0\left(  \sX\right)\subset C_b\left(  \sX\right)$ and $C_0\left( \widetilde \sX\right)\subset C_b\left( \widetilde \sX\right)$ are essential ideals. Proof of conditions (b) and (c) of the Definition \ref{fin_weak_defn} is left  to the reader.
\end{proof}

\subsection{Coverings of *-algebras}

\subsubsection{Coverings of pro-$C^*$-algebras}
\paragraph{}

If $\sX$ is a locally compact, Hausdorff space then a *-algebra $Cont\left(\sX \right)$ of continuous complex-valued functions on $\sX$ is a completion of $C_0\left(\sX\right)$ with respect to a family of $C^*$-seminorms
\bean
p_\sV : \sX \to \R,\\
a \mapsto \max_{x \in \sV}\left|a\left(x\right)\right|
\eean
where $\sV \subset\sX$ is a compact subset. So $Cont\left(\sX \right)$ is a pro-$C^*$-algebra (cf. Definition \ref{pro_c_defn}) such that $C_0\left(\sX\right)$ is an essential ideal (cf. Definition \ref{essential_defn}) of the $C^*$-algebra  $b\left(Cont\left(\sX \right) \right)$ of bounded elements (cf. Definition \ref{pro_bound_defn}).
\begin{exercise}\label{top_pro_fin_exer} Let $\sX$ be a locally compact, Hausdorff space, and let $p: \widetilde{\sX}\to \sX$ be a finite-fold transitive covering.
	Prove following statements:
	\begin{itemize}
		\item Using a notation \eqref{top_c_sec_eqn} and the Definition \ref{operator_fields_continuity_defn} proof that there are natural $*$-isomorphisms
		$Cont\left(\sX\right) \cong  C\left(\sX, \left\{\C_x\right\}, C_0\left(\sX \right)  \right)$ and\\ $Cont\left(\widetilde \sX\right) \cong  C\left(\widetilde\sX, \left\{\C_{\widetilde x}\right\}, C_0\left( \widetilde\sX \right)  \right)$.
		\item Using a notation \eqref{top_cont_lift_eqn} prove that there is a natural $*$-isomorphisms \\$ C\left(\widetilde\sX, \left\{\C_{\widetilde x}\right\}, C_0\left( \widetilde\sX \right)  \right)\cong \lift_p\left[C\left(\sX, \left\{\C_x\right\}, C_0\left(\sX \right)  \right)\right]$.
		\item The given by the Definition \eqref{top_lift_defn} map
		$$
		\lift_p :  C\left(\sX, \left\{\C_x\right\}, C_0\left(\sX \right)  \right)\hookto  C\left(\widetilde\sX, \left\{\C_{\widetilde x}\right\}, C_0\left( \widetilde\sX \right)  \right)
		$$
		yields a natural injective  $*$-homomorphism
		\be\label{top_pro_hom_eqn}
		Cont\left( p\right)\bydef \lift_p : Cont\left(\sX\right)\hookto Cont\left(\widetilde\sX\right)
		\ee
		of pro-$C^*$-algebras.
		\item There is a natural isomorphism of groups
			\bean
		\left\{ \left.g \in \Aut\left(Cont\left(\widetilde\sX\right) \right)~\right|~ ga = a;~~\forall a \in Cont\left(\sX\right)\right\}\cong G\left(\left.\widetilde\sX\right| \sX \right)
		\eean
		such that
		$$
		Cont\left(\widetilde\sX\right)^{G\left(\left.\widetilde\sX\right| \sX \right)}\cong Cont\left(\sX\right).
		$$
			
	\end{itemize}
\end{exercise}
\begin{theorem}\label{top_pro_fin_thm}
	Let $\sX$ be a connected, locally connected, locally compact, Hausdorff space. If $\sX$ consents with noncommutative coverings,  $p:\widetilde\sX \to \sX$ is a finite-fold, transitive covering and   $
	Cont\left(  p\right): Cont\left(\sX\right)\hookto Cont\left(\widetilde\sX\right)$ is given by the equation \eqref{top_pro_hom_eqn}  then the quadruple 
	\be\label{top_pro_fin_eqn}
	\left(Cont\left(\sX\right), Cont\left(\widetilde\sX\right), G\left(\left.\widetilde\sX\right| \sX \right) , Cont\left( p\right) \right)
	\ee
	is a {noncommutative finite-fold  covering  of pro-$C^*$ algebras} (cf. Definition \ref{pro_fin_defn}).
\end{theorem}
\begin{proof}
	From the Exercise \ref{top_pro_fin_exer} it follows that 	the quadruple \eqref{top_pro_fin_eqn} is a pre-covering of $*$-algebras (cf. Definition \ref{fin_pre*_defn}).	
Both $Cont\left(\sX\right)$ and  $Cont\left(\widetilde\sX\right)$ are pro-$C^*$-algebras of complex-valued continuous maps on $\sX$ and $\widetilde\sX$ respectively. So we have $b\left(Cont\left(\sX\right)\right)= C_b\left(\sX\right)$ and $b\left(Cont\left(\widetilde \sX\right)\right)= C_b\left(\widetilde\sX\right)$. From the Corollary \ref{top_fin_cor} it follows that the quadruple
$$
\left(C_b\left(\mathcal  X \right), C_b\left(  \widetilde{   \mathcal X }\right) ,G\left(\left. \widetilde{   \mathcal X } ~\right| {   \mathcal X }\right),    \left.Cont\left( p\right)\right|_{C_b\left(\mathcal  X \right)} \right)
$$
is a connected noncommutative  finite-fold weak covering (cf. Definition \ref{fin_weak_defn}). Thus the quadruple \eqref{top_pro_fin_eqn} matches  to the Definition \ref{pro_fin_defn}.
\end{proof}

\subsubsection{Coverings of bounded operator $*$-algebras}

\begin{theorem}\label{top_oa_cov_thm} 
	Let $\sX$ be a connected, locally connected,  locally compact, paracompact, Hausdorff  space. If $\sX$ consents with noncommutative coverings and $R \subset C_0\left(\sX\right)$ is a $c$-soft  *-subalgebra (cf. Definition \ref{top_soft_r_defn}  then	following conditions hold:
	\begin{enumerate}
		\item [(i)] if 	$\left(R, \widetilde A, G, \pi\right)$ noncommutative finite-fold covering of bounded operator *-algebras  (cf. Definition \ref{fin_oa_defn}) then there is a finite fold transitive topological covering  $p: \widetilde \sX \to \sX$ such that there is a natural dense inclusion  $\widetilde A \subset C_0\left( \widetilde \sX\right)$ (cf. Definition \ref{top_x_sheaf_defn}) and
		$G \cong G\left(\left.\widetilde \sX~ \right| \sX\right)$ and $\pi= \left.C_0\left(p \right)\right|_{R }$ where $C_0\left(p \right)$ is given by the equation \eqref{top_c0p_eqn},
		\item[(ii)] if $\mathscr S^R$ is the $R$-sheaf (cf. Definition \ref{top_x_sheaf_defn}) and $p^{-1}\mathscr S^R$ is its inverse image (cf. Definition \ref{sheaf_inv_im_defn}) then a quadruple	$$
	\left( 	R, C_0\left(\widetilde\sX \right)\cap p^{-1} \mathscr S^R\left(\widetilde\sX \right), G\left(\left.\widetilde \sX~ \right| \sX\right), \left.C_0\left(p \right)\right|_{R}\right)
		$$
		finite-fold covering of bounded operator *-algebras  (cf. Definition \ref{fin_oa_defn}).
	\end{enumerate}
\end{theorem}
\begin{proof}
	(i)
	The $C^*$-norm completion of $R$ is $C_0\left(\sX\right)$. If a  quadruple $\left(C_0\left(\sX \right), \widetilde{B}, G,    \pi\right)$ is a connected noncommutative finite-fold covering (cf. Definition \ref{fin_defn}  then   there is a finite-fold transitive topological covering and $p: \widetilde{    \sX }\to \sX$ such that $\left(C_0\left(\mathcal  X \right), \widetilde{B}, G,    \pi\right)$ is equivalent to $$\left(C_0\left(\sX \right), C_0\left(  \widetilde{\sX}\right) , G\left(\left.\widetilde \sX~ \right| \sX\right),    C_0\left( p\right) \right)$$ (cf. the Theorem \ref{top_fin_necassary_lem}. On the other hand from \ref{fin_oa_empt} it follows that there is a dense inclusion  $\widetilde A\subset C_0\left(  \widetilde{\sX}\right)$. Conditions $G \cong G\left(\left.\widetilde \sX~ \right| \sX\right)$ and $\pi= \left.C_0\left(p \right)\right|_{R }$ follow from the construction \ref{fin_oa_empt}.\\
	(ii)	
	We leave to the reader a proof of that  $ C_0\left(\widetilde\sX \right)\cap p^{-1} \mathscr S^R\left(\widetilde\sX \right)$ is an admissible subalgebra, i.e. it satisfies to conditions of \ref{fin_oa_empt}. If  $ C_0\left(\widetilde\sX \right)\cap p^{-1} \mathscr S^R\left(\widetilde\sX \right)$ is not maximal admissible *-algebra then there is an admissible algebra  $\widetilde A' \subset C_0\left(\widetilde \sX\right)$  such that $ C_0\left(\widetilde\sX \right)\cap p^{-1} \mathscr S^R\left(\widetilde\sX \right)\subsetneqq \widetilde A'$. If $\widetilde a \in \widetilde A' \setminus  C_0\left(\widetilde\sX \right)\cap p^{-1} \mathscr S^R\left(\widetilde\sX \right)$ then there is $\widetilde x_0 \in \widetilde \sX$ such that $\widetilde a_{\widetilde x_0} \notin p^{-1} \mathscr S^R_{\widetilde x_0}$ where $\widetilde a_{\widetilde x_0}$ is the stalk at $\widetilde x_0$ (cf. Definition \ref{sheaf_stalk_defn}). If $\widetilde f\in p^{-1}\mathscr S^{R}\left( \widetilde \sX\right) \cap C_0\left(\widetilde \sX\right)$ be an $\left(p, R, \widetilde x\right)$-{stump} (cf. Definition \ref{top_stump_soft_p_defn}) then $\widetilde a_{\widetilde x_0}=\left(\widetilde f\widetilde a \right)_{\widetilde x_0}\notin p^{-1} \mathscr S^R_{\widetilde x_0}$.
	One has 
\bean
a \bydef \sum_{	g \in G\left(\left.\widetilde{\sX}~\right|\sX\right)} g\left(\widetilde f \widetilde a \right) \in \widetilde A'  \cap C_0\left(\sX\right),\\
a = \desc_p\left( \widetilde f\widetilde a\right)\in \widetilde A'^{G\left(\left.\widetilde{\sX}~\right|\sX\right)}
\eean	
where $\desc_p$ is given by the Definitions \ref{top_lift_desc_defn}, \ref{top_lift_sh_desc_defn}. From the isomorphism \eqref{top_st_iso_eqn} and $\left( \widetilde f\widetilde a\right)_{\widetilde x_0} \notin p^{-1} \mathscr S^R_{\widetilde x_0}$ it turns out that
	 $a_{p\left( \widetilde x_0\right)}  \notin \mathscr S^R_{p\left( \widetilde x_0\right) }$, i.e. $a \notin R$. On the other hand from \ref{fin_oa_empt} it follows that $$\left(R, \widetilde A', G, \pi\right)$$ is a noncommutative finite-fold  pre-covering of *-algebras, and taking into account the condition  (b) of the Definition \ref{fin_pre*_defn} one has $A'^{G\left(\left.\widetilde{\sX}~\right|\sX\right)}=R$ so   $a \in R$. From this contradiction it turns out that $ C_0\left(\widetilde\sX \right)\cap p^{-1} \mathscr S^R\left(\widetilde\sX \right)$ is a maximal admissible algebra.
\end{proof}

\subsubsection{Coverings of $O^*$-algebras}\label{top_pd_fin_sec}

\paragraph*{}
If $M$ be a smooth manifold then  the product $T \bydef M \times \C$ is also a smooth manifold, because $\C$ has the natural smooth structure. So one has a smooth bundle $T \to M$ (cf. Definition \ref{top_sm_bundle_defn}). Every fiber of $T$ is naturally homeomorphic to $\C$, so having a sesquilinear product
\be\label{top_ses_one_eqn}
\C\times\C \to \C,\\
\left(z', z'' \right) \mapsto \overline{z}'z'',
\ee
one can define a sesquilinear map   $T\times_M T \to \C$  (cf. Definition \ref{top_herm_bundle_form_defn}). From   \eqref{top_h_sp_eqn}  it follows that there is a faithful representation 
	\bean
\rho :C_0\left( M \right) \hookto B\left( L^2\left(M, T,  \mu\right)\right).
\eean
If  $p: \widetilde M\to M$ is a topological transitive finite-fold covering with connected  $\widetilde M$  then the space  $\widetilde M$ has a natural structure of smooth manifold given by the Proposition \ref{top_cov_mani_prop}. From the Theorems \ref{top_fin_necassary_lem} and \ref{top_fin_sufficient_lem}  it follows that $$
\left(C_0\left(M \right), C_0\left(  \widetilde{   M }\right) ,G\left(\left. \widetilde{  M } ~\right| {   M }\right),   C_0\left(p \right)  \right)
$$
the connected  noncommutative finite-fold covering (cf. Definition \ref{fin_defn}). If a representation $\widetilde \rho: C_0\left(  \widetilde{   M }\right) \to B\left(\widetilde \H \right)$ is induced be the pair  $$\left(\rho,\left(C_0\left(M \right), C_0\left(  \widetilde{   M }\right) ,G\left(\left. \widetilde{  M } ~\right| {   M }\right),   C_0\left(p \right)  \right)  \right)$$ (cf. Definition \ref{induced_repr_fin_defn}) then it is faithful (cf. Lemma \ref{induced_faithful_lem}) and equivariant (cf. equation \eqref{induced_equiv_eqn}).
We leave to the reader the proof of a natural isomorphism $\widetilde \H \cong L^2\left(\widetilde M,  \widetilde T,   \widetilde\mu\right)$ where $\widetilde T$ is an inverse image of $T$ by $p$ (cf. Definition \ref{vb_inv_img_funct_defn}) and $\widetilde\mu\bydef \lift_p\mu$ is the $p$-lift of $\mu$ (cf. Definition \ref{top_lift_measure_defn}). 
Below we will prove that the given by the equation \eqref{top_diff_*alg_lift_eqn} $*$-homomorphism 
\bean
D^*\left(p,T \right) : D^*\left(M, T\right)\hookto D^*\left(\widetilde M, \widetilde T\right)
\eean
is  a {noncommutative finite-fold covering of $O^*$-algebras}
(cf. Definition \ref{fino*_defn}).  We leave to the reader the prove of that an action $G\left(\left. \widetilde{  M } ~\right| {   M }\right)\times \widetilde{  M } \to \widetilde{  M } $  naturally induces an action $G\left(\left. \widetilde{  M } ~\right| {   M }\right)\times D^*\left(\widetilde M, \widetilde T\right)\to D^*\left(\widetilde M, \widetilde T\right)$.
\begin{empt}
Here we consider a related to the $*$-homomorphism $D^*\left(p,T \right) : D^*\left(M, T\right)\hookto D^*\left(\widetilde M, \widetilde T\right)$ specialization of the construction \ref{o*fin_empt}.
 Denote by
\be\label{top_o*_triple_eqn}
\left(A, \widetilde{A}, G, \pi \right)\bydef \left(  D^*\left(M,T\right), D^*\left(\widetilde M, \widetilde T\right), G\left(\left. \widetilde{  M } ~\right| {   M }\right), D^*\left(p,T \right) \right),
\ee
We leave to the reader the prove of that the quadruple \eqref{top_o*_triple_eqn} noncommutative finite-fold  pre-covering of *-algebras (cf. Definition \ref{fin_pre*_defn}). If
\bean
\H \bydef L^2\left(\widetilde M,  \widetilde T,   \widetilde\mu\right), \\
\D \bydef \Ga^\infty_c\left(\widetilde M,  \widetilde T \right) 
\eean
then from the construction \ref{top_h_sp_empt} it turns out the $\D$ is a dense subspace of $\H$. From 
\eqref{top_diff*_eqn} one has an inclusion $\widetilde{A} \bydef D^*\left(\widetilde M, \widetilde T\right)\subset \L^\dagger\left(\D\bydef \Ga^\infty_c\left(\widetilde M,  \widetilde T \right)\right)$ of $O^*$-algebras.
The action $G\times \widetilde M\to \widetilde M$ induces an action $G\times \D\to \D$ which is equivariant (cf. equation \ref{equ_d_eqn}). If $\mathcal L^\dagger\left(\D\right)_b$ is given by \eqref{o*b_eqn} then the intersection $\widetilde{A}_b \bydef D^*\left(\widetilde M, \widetilde T\right)\cap \L^\dagger\left(\D\right)_b$ contains bounded operators only. However any bounded differential operator has order 0 (cf. Definition \ref{do_man_order_defn}). It follows that 
\be\label{top_ddound_eqn}
\begin{split}
\widetilde{A}_b \bydef \widetilde{A}\cap \L^\dagger\left(\D\right)_b \cong D^*\left(\widetilde M, \widetilde T\right)\cap \L^\dagger\left(\D\right)_b\cong \Coo\left(\widetilde M \right)\cap C_b\left(\widetilde M \right),\\
A_b \bydef  \widetilde{A}_b \cap A \cong \Coo\left( M \right)\cap C_b\left( M \right)
\end{split}
\ee
If both $B$ and $\widetilde B$ are $C^*$-norm completions of both $\widetilde{A}_b \cong \Coo\left(\widetilde M \right)\cap C_b\left(\widetilde M \right)$ and $A_b  \cong \Coo\left( M \right)\cap C_b\left( M \right)$ then one has $B\cong C_b\left(M \right)$ and $\widetilde B\cong C_b\left(\widetilde M \right)$. However from the Corollary \ref{top_fin_cor} it turns out that
$$
\left(B, \widetilde{B}, G, \pi_B \right)\cong \left(C_b\left(M \right), C_b\left(  \widetilde{ M }\right) ,G\left(\left. \widetilde{   M } ~\right| {   M }\right),   \lift_p \right)
$$
is  a connected noncommutative weak finite-fold covering. From this construction one has the following theorem. 
\end{empt}

\begin{theorem}\label{top_d_fin_thm}
In the described above situation the quadruple
$$
\left( D^*\left(M,T\right), D^*\left(\widetilde M, \widetilde T\right) ,G\left(\left. \widetilde{   M } ~\right| {   M }\right),D^*\left(p,T \right)\right)
$$
is a {noncommutative finite-fold covering of $O^*$-algebras}
(cf. Definition \ref{fino*_defn}).
\end{theorem}
Similarly the reader can prove the following lemma.

\begin{lemma}\label{top_s_fin_lem}
The quadruple
$$
\left(\Coo\left(M \right),\Coo\left(  \widetilde{   M }\right)  ,G\left(\left. \widetilde{   M } ~\right| {   M }\right), \Coo\left(p \right)\bydef \lift_p|_{\Coo\left(M \right)}\right)
$$
is a {noncommutative finite-fold covering of $O^*$-algebras}
(cf. Definition \ref{fino*_defn}).
\end{lemma}
\begin{remark}
The Lemma \ref{top_s_fin_lem} is a noncommutative counterpart of the Proposition \ref{top_cov_mani_prop}. 
\end{remark}
\subsubsection{Coverings and unbounded operators on Hilbert modules}\label{top_fin_chull_sec}
\paragraph{Algebras of unbounded functions}

Let both $\sX$ be  and $\widetilde\sX$ be connected, locally connected, locally compact, H Hausdorff spaces, and let $p:\widetilde\sX\to \sX$ be a transitive covering. Suppose that $\sX$ consents with noncommutative coverings.
If both $A\bydef Cont\left( \sX\right) $ and $\widetilde A\bydef Cont\left(\widetilde  \sX\right)$  are *-algebras of (unbounded) $\C$-valued continuous maps then $p$ induces an injective $*$-homomorphism
$$
\pi_A : A\hookto \widetilde A.
$$
If $G \bydef G\left(\left.\widetilde \sX~\right|\sX \right)$ then the natural action   $G \times \widetilde \sX \to \widetilde \sX$ induces an action $G\times \widetilde A \to \widetilde A$. We leave to the reader the proof of that a quadruple $\left(  A ,\widetilde{ A}, G\left(\left.\widetilde \sX~\right|\sX \right), \pi_A\right)$ is a  {noncommutative finite-fold  pre-covering of *-algebras} (cf. Definition \ref{fin_pre*_defn}).
If
\bean
B\bydef  C_0\left( \sX\right),\quad
\widetilde B\bydef  C_0\left( \widetilde\sX\right),\\
\mathfrak B\bydef  C_c\left( \sX\right),\quad
\widetilde{\mathfrak B}\bydef  C_c\left( \widetilde\sX\right).
\eean 
then from the Lemma \ref{top_fin_necassary_lem}
 it turns out that
$$
\left( B, \widetilde B, G, \pi_B  \right) \bydef \left( C_0\left( \sX\right),  C_0\left( \widetilde\sX\right), G\left(\left.\widetilde \sX~\right|\sX \right), C_0\left(p \right)  \right) 
$$
is a connected {noncommutative finite-fold  covering} (cf. Definition \ref{fin_defn}).
\begin{empt}\label{top_fin_hilb_empt}
	Let us prove that above objects satisfy to conditions (a)-(c) of \ref{fin_chull_empt}.\\

	(a) An algebra  $C_c\left(\sX\right)$ is a dense right ideal of $C_0\left(\sX \right)$ so from  $B\bydef  C_0\left( \sX\right)$ and 
	$\mathfrak B\bydef  C_c\left( \sX\right)$ if follows that $\mathfrak B$ is a dense right ideal of $B$. On the other hand from 
	$$
	\forall a \in Cont\left( \sX\right)\quad \forall b \in C_c\left( \sX\right) \quad a b \in C_c\left( \sX\right)
	$$
	it follows that there is an action $Cont\left( \sX\right)\times  C_c\left( \sX\right)\to C_c\left( \sX\right)$ or using the above notation one has $A\times \mathfrak B\to\mathfrak B$. This action is equivalent to an inclusion  $\mu : A \hookto \End^*_{ B}\left( {\mathfrak B}\right)$  of *-algebras.\\
	(b) Similarly to (a) one has a dense right ideal $\widetilde{\mathfrak{B}}\bydef C_c\left( \widetilde \sX\right) \subset \widetilde B\bydef C_0\left( \widetilde \sX\right)$. Any homeomorphism of $\widetilde \sX$ maps compact subsets of $\widetilde \sX$ onto compact ones. It turns out that $G\left(\left.\widetilde \sX~\right|\sX \right)C_c\left( \widetilde \sX\right) = C_c\left( \widetilde \sX\right)$, or using the above notation  one has $G\widetilde{\mathfrak{B}} = \widetilde{\mathfrak{B}}$. Any $G\left(\left.\widetilde \sX~\right|\sX \right)$-invariant element of $C_c\left( \widetilde \sX\right)$ lies in $C_c\left(\sX \right)$ and vice versa, it turns out that
	$$
	C_c\left(\sX \right)= C_c\left( \widetilde{\sX}\right) ^{G\left(\left.\widetilde \sX~\right|\sX \right)} \bydef \left\{\left. \widetilde{b}\in C_c\left( \widetilde{\sX}\right) \right| \forall g \in G\left(\left.\widetilde \sX~\right|\sX \right)\quad g \widetilde b = \widetilde b\right\}.
	$$
	The above equation is a specialization of \eqref{fin_chull_e_eqn}.\\
(c) Can be proven similarly to (a).
\end{empt}

\begin{theorem}\label{top_fin_hilb_thm}
In the above situation  the quadruple $\left(  A ,\widetilde{ A}, G \bydef G\left(\left.\widetilde \sX~\right|\sX \right), \pi_A\right)$ is  
	an associated with $\left(C_0\left( \sX\right),C_0\left( \widetilde\sX\right) ,G\left(\left.\widetilde \sX~\right|\sX \right), C_0\left(p \right) \right) $  noncommutative finite-fold covering of *-algebras in the sense of the  Definition \ref{fin_chull_defn}.
	\end{theorem}
\begin{proof}
The construction \ref{top_fin_hilb_empt} is a specialization of  \ref{fin_chull_empt} one. We leave to the reader prove of conditions (a) and (b)  of the Definition \ref{fin_chull_defn}.
\end{proof}
\begin{remark}
There    are completions of  ${\mathfrak{B}}$ and $\widetilde{\mathfrak{B}}$ with respect to graph topology (cf. Definition \ref{def:rep_Hilbert_module}) given by norms
\bean
\forall a \in A \quad \forall b  \in {\mathfrak{B}} \quad  \left\| b\right\|_a \bydef {\braket{b}{\pi(1+a^*a)b}}^{\nicefrac12},\\
\forall \widetilde a \in\widetilde A \quad \forall\widetilde b  \in  \widetilde{\mathfrak{B}} \quad  \left\| \widetilde b\right\|_{\widetilde a} \bydef {\braket{\widetilde b}{\widetilde\pi(1+\widetilde a^*\widetilde a)\widetilde b}}^{\nicefrac12}.
\eean
In the above proof the ideals ${\mathfrak{B}}'$ and $\widetilde{\mathfrak{B}}'$ can be replaced by their completions.
\end{remark}

\paragraph{Polynomials in one variable} \label {top_fin_hull_sec}
Let $\sX$ be a locally compact,  connected, locally connected, Hausdorff space, and let $R\subset C_0\left(\sX \right)$ be a $c$-soft *-algebra (cf. Definition \ref{top_soft_r_defn}). 
Let  $R\left[t\right]$ be a commutative *-algebra of polynomials of $t$ such that $t= t^*$. Clearly $R\left[t\right]$ is an $R$-module. Let $\mathscr S^{R\left[t\right]}$ be an $R\left[t\right]$-sheaf (cf. Definition \ref{top_x_sheaf_defn}).
 Let 
\be\label{top_poly_a_defn}
A \bydef \mathscr S^{R\left[t\right]}\left(\sX\right)
\ee
 be a *-algebra of global sections.
If $B \bydef C_0\left(\sX \times \R\right)$ 
and $\mathfrak B'\bydef C_c\left(\sX \times \R\right)$
then $\mathfrak B'$ is a dense right ideal of $B$. Moreover if $T: \R \to \R$ is an identical map then there is a representation
\bean
R\left[t\right]\hookto \End_B\left(\mathfrak B' \right),\\
\forall a \in R\quad n \in \N^0 \quad \forall b \in \mathfrak B' \quad a t^n \mapsto \left(b \mapsto  b aT^n\right). 
\eean
which can be extended up to a representation $\left( \mathfrak{B}',\mu'\right)$ of $A$ on $B$ (cf. Definition \ref{def:rep_Hilbert_module_uni}), i.e. one has
$$
\mu': A\hookto \End_B\left(\mathfrak B' \right).
$$
If $\mathfrak B$ is the completion of $\mathfrak B'$ with respect to the graph topology (cf. Definition \ref{def:rep_Hilbert_module} ) then there is the representation  $\mu:  A \hookto \End_B\left(\mathfrak B \right)$ such that $\mathfrak B$ is a core (cf. Definition \ref{def:rep_Hilbert_module}).
 If  $p:\widetilde \sX \to \sX$ is a finite-fold transitive covering with connected $\widetilde \sX$ and $G \bydef G\left(\left.\widetilde \sX \right| \sX  \right)$ then there is the natural finite-fold transitive covering $\widetilde \sX\times\R \to \sX\times\R$ such that $G \cong  G\left(\left.\widetilde \sX\times\R \right| \sX \times\R \right)$. 
 From the 
 Lemma \ref{top_fin_necassary_lem} 
  it follows that one has a natural  {noncommutative finite-fold  covering} $\left(B, \widetilde{B}\bydef C_0\left(\widetilde \sX \times \R\right) , G, \pi_B \right)$).
If $\mathscr S^R$ is an $R$-sheaf (cf. Definition \ref{top_x_sheaf_defn}) and $p^{-1}\mathscr S^{R}$ is its inverse image (cf. Definition \ref{sheaf_inv_im_defn}) then from the Exercise \ref{top_soft_c_exer} it turns out that
an intersection $\widetilde R \bydef p^{-1}\mathscr S^{R}\left( \widetilde\sX\right) \cap C_0\left(\widetilde \sX\right)$ is a $c$-soft *-subalgebra. We leave to the reader the proof of the natural isomorphism
$$
p^{-1}\mathscr S^{R\left[t\right]}\cong \mathscr S^{\widetilde R\left[t\right]}.
$$
Denote by $\widetilde A \bydef \mathscr S^{\widetilde R\left[t\right]}\left(\widetilde \sX\right)$.
\begin{exercise}\label{top_fin_hull_exer}
Prove following statements.
	\begin{enumerate}
				\item There exists a natural pre-covering   covering of *-algebras  $\left(A, \widetilde{A}, G, \pi_A \right)$   (cf. Definition \ref{fin_pre*_defn}).
		\item The quadruple  $\left(A, \widetilde{A}, G, \pi_A \right)$ is an {associated with $\left(B, \widetilde{B}, G, \pi_B \right)$  noncommutative finite-fold  covering of *-algebras} (cf. Definition \ref{fin_chull_defn}).
		\item The quadruple $\left(A, \widetilde{A}, G, \pi_A \right)$ can be regarded as a noncommutative finite-fold covering of $O^*$-algebras (cf. Definition \ref{fino*_defn}).
	\end{enumerate}
\end{exercise}
If both $tA\subset A$ and $t\widetilde A\subset \widetilde A$ are ideals then a quadruple  $\left(tA, t\widetilde{A}, G, \pi_A \right)$ is a noncommutative finite-fold  pre-covering of *-algebras (cf. Definition \ref{fin_pre*_defn}  ). Let both $B^T$ and $\widetilde B^T$ be $C^*$-norm completion of $T {\mathfrak B}$ and $T\widetilde {\mathfrak B}$ respectively.
\begin{exercise}
Prove that the quadruple  $\left(tA, t\widetilde{A}, G, \left.\pi_A\right|_{tA} \right)$ is an {associated with $\left(B^T, \widetilde{B}^T, G, \left.\pi_B\right|_{B^T} \right)$  noncommutative finite-fold  covering of *-algebras} (cf. Definition \ref{fin_chull_defn}).

\end{exercise}
\begin{remark}
Both ideals  $tA$ and $t\widetilde A$ do not contain nontrivial bounded elements, so the quadruple  $\left(tA, t\widetilde{A}, G, \left.\pi_A\right|_{tA} \right)$ is not a noncommutative finite-fold covering of $O^*$-algebras.
\end{remark}

\subsection{Coverings of quasi *-algebras}
\paragraph*{}
	Let $M$ be a smooth manifold, $\Coo_c\left( M\right)\bydef C_c\left(M\right)\cap \Coo\left(M\right)$ and $\Coo_c\left(M \right)'$ is a space of distribution densities (cf. Definition \ref{top_distr_dens_def}). If $p: \widetilde M \to M$ is a finite-fold transitive covering then there is the natural injective homomorphism
	$$
	\Coo_c\left(p \right)':\Coo_c\left(M \right)'\hookto \Coo_c\left(\widetilde M \right)'
$$
	of $\Coo\left(M \right)$-modules (cf. equation \eqref{top_distr_eqn}). 		From the Remark 
	\ref{top_distr_dens_q_rem} it follows that there are following quasi *-algebras
	\bean
	\left(\mathfrak{A}, \mathfrak{A}_0\right)\bydef \left(\Coo_c\left( M \right)' , \Coo\left( M \right)\right),\\
	\left(\widetilde{\mathfrak{A}}, \widetilde{\mathfrak{A}}_0\right)\bydef \left(\Coo_c\left(\widetilde M \right)' , \Coo\left( \widetilde M \right)\right)
	\eean
	(cf. Definition \ref{qousi_star_defn}).
	\begin{exercise}\label{top_quasi_exer}
	Prove following statements:
	\begin{enumerate}
			\item[(i)] Using the equation \eqref{top_distr_eqn} and the Lemma \ref{top_s_fin_lem} prove that there is an injective $*$-homomorphism 
			$$
	\left( \Coo_c\left(p \right)', \Coo\left( p\right)  \right) :		\left(\mathfrak{A}, \mathfrak{A}_0\right)\hookto \left(\widetilde{\mathfrak{A}}, \widetilde{\mathfrak{A}}_0\right)
			$$
			(cf. Definition \ref{quasi_hom_defn}).
			\item[(ii)] Using the Exercise \ref{top_distr_exer} and  the Lemma \ref{top_s_fin_lem} define an action
			$$
	 G\left(\left. \widetilde M \right| M \right)\times \left(\widetilde{\mathfrak{A}}, \widetilde{\mathfrak{A}}_0\right)\to \left(\widetilde{\mathfrak{A}}, \widetilde{\mathfrak{A}}_0\right)	
			$$
		and prove that the quadruple
		$$
		\left(\left( {\mathfrak{A}}, {\mathfrak{A}}_0\right) , \left( \widetilde{\mathfrak{A}}, \widetilde{\mathfrak{A}}_0\right) ,  G\left(\left. \widetilde M \right| M \right) , \left( \Coo_c\left(p \right)', \Coo\left( p\right)  \right)\right)
		$$
	is a {noncommutative finite-fold pre-covering of quasi $*$-algebras}	(cf. Definition \ref{oq*fin_defn}).
	\end{enumerate}
	\end{exercise}

\begin{lemma}\label{top_q_fin_lem}
	The quadruple 
	$$
	\left(\left(\Coo_c\left(M \right)', \Coo\left(M \right)    \right), \left(\Coo_c\left(\widetilde M \right)', \Coo\left( \widetilde M \right)    \right), G\left(\left. \widetilde M\right| M\right),\left( \Coo_c\left(p \right)', \Coo\left( p\right)\right)\right) 
	$$
	is a	{noncommutative finite-fold covering of quasi $*$-algebras} (cf. Remark \ref{top_distr_dens_q_rem} and Definition \ref{oq*fin_defn}).
\end{lemma}
\begin{proof}
	From the Exercise \ref{top_quasi_exer} it follows that the quadruple 
a {noncommutative finite-fold pre-covering of quasi $*$-algebras}. From the Lemma \ref{top_s_fin_lem} it follows that the quadruple
$$
\left(\Coo\left(M \right),\Coo\left(  \widetilde{   M }\right)  ,G\left(\left. \widetilde{   M } ~\right| {   M }\right), \Coo\left(p \right)\right)
$$
is a noncommutative finite-fold covering of $O^*$-algebras.
\end{proof}

\subsection{Coverings of operator spaces}
\begin{empt}\label{top_r_fin_empt}
	Let $\sX$ be a connected, locally connected,  locally compact, Hausdorff space, and let
	\be\label{top_sim_eqn}
	\sX^\sim \bydef \begin{cases}
		\sX & \sX \text{ is compact}\\
		\text{the one point compactification of } \sX & \sX \text{ is not compact}
	\end{cases} 
	\ee (cf. Definition \ref{top_comp_defn}). Suppose  $\sX$ is locally connected.
	From the Example \ref{op_real_exm} it turns out that he space $C\left(\sX^\sim, \R \right)$ of all real-valued continuous functions is a real operator space. If $C_0\left(\sX, \R \right)$ is an algebra of vanishing at infinity real-valued functions then the pair $$\left(C_0\left(\sX, \R \right), C\left(\sX^\sim, \R \right) \right)$$ is a sub-unital real operator space (cf. the Definitions \ref{operator_space_subunital_defn}, \ref{op_su_r_space_defn} and the Example \ref{sub_alg_exm}). If $p: \widetilde{\sX} \to \sX$ is a finite-fold transitive covering then there are the following complexifications
	\be\label{top_complex_eqn}
	\begin{split}
		\C C_0\left(\sX, \R \right)\cong C_0\left(\sX\right),\\
		\C C\left(\sX^\sim, \R \right)\cong C\left(\sX^\sim \right),\\
		\C C_0\left(\widetilde\sX, \R \right)\cong C\left(\widetilde\sX\right),\\
		\C C\left(\widetilde\sX^\sim, \R \right)\cong C\left(\sX^\sim \right).
	\end{split}
	\ee
 The $C^*$-envelopes (cf. Definition \ref{operator_space_envelope_defn}) of   	 $\left(C_0\left(\sX^\sim \right), C\left(\sX^\sim \right) \right)$ and $\left(C_0\left(\widetilde\sX \right), C\left(\widetilde\sX^\sim\right) \right)$ are equal to 
	\bean
	C^*_e\left(C_0\left(\sX \right), C\left(\sX^\sim \right) \right) = C_0\left(\sX \right),\\
	C^*_e\left(C_0\left(\widetilde\sX \right), C\left(\widetilde\sX^\sim \right) \right) = C_0\left(\widetilde\sX \right)
	\eean
	respectively (cf. Remark \ref{op_su_env_rem}).
	If $G \bydef G\left(\left. \widetilde{   \mathcal X } ~\right| {   \mathcal X }\right)$ then 	
	from the Lemma \ref{top_fin_necassary_lem}
	 it follows that
	$$
	\left(\left(C_0\left(\sX \right), C\left(\sX^\sim \right) \right),\left(C_0\left(\sX \right), C\left(\sX^\sim \right)\right) , G, \left(\pi_{C_0\left(\sX \right)}, \pi_{C\left(\sX^\sim \right)} \right) \right)
	$$   
	is a  {noncommutative finite-fold covering} of the {sub-unital} operator space $\left(C_0\left(\sX \right), C\left(\sX^\sim \right) \right)$. Taking into account the equation \eqref{top_complex_eqn} and the Definition \ref{fin_rop_defn} one can prove that there is a  {noncommutative finite-fold covering 
		\be\label{top_real_fin_eqn}
		\left(\left(C_0\left(\sX, \R \right), C\left(\sX^\sim , \R\right) \right),\left(C_0\left(\sX, \R \right), C\left(\sX^\sim \right), \R\right) , G, \left(\pi, \pi^\sim \right)\right) 
		\ee
		of the sub-unital real operator space} and  $\left(C_0\left(\sX, \R \right), C\left(\sX^\sim , \R\right) \right)$ where  both maps $\pi:C_0\left(\sX, \R \right)\hookto C_0\left(\widetilde\sX, \R \right)$ and $\pi^\sim:C\left(\sX^\sim, \R \right)\hookto C\left(\widetilde\sX^\sim, \R \right)$  are induced by the covering $p$.
\end{empt}

\subsection{Induced representations}\label{comm_induced_finite_sec}

\paragraph*{} 
Let  $\sX$ be a connected, locally connected, locally compact, Hausdorff space. If 
 and  $E \to \sX$ is a locally trivial bundle with a {sesquilinear form}  $\varphi: E \times_\sX E \to \C$ (cf. Definitions \ref{top_vb_defn} and \ref{top_herm_bundle_form_defn}) then there are  given by \eqref{top_ggc_eqn} and \eqref{top_prod_eqn} pairings
\bean
	\left\langle \cdot, \cdot \right\rangle_c:  \Ga_c\left(\sX, E\right)\times \Ga_c\left(\sX, E\right) \to C_c\left( \sX\right),\\
 	\left(\cdot, \cdot\right): \Ga_c\left(\sX, E\right)\times \Ga_c\left( \sX, E\right)\to \C,\\
\left(a, b \right) \mapsto\tau\left( \left\langle a, b \right\rangle_c\right)
\eean
where $\tau$ is a positive functional 
\bean
\tau : C_c\left( \sX\right) \to \C,\\
a \mapsto \int_{\sX} a~ d\mu.
\eean
(cf. Theorem \ref{meafunc_thm}) such that $a > 0  \Rightarrow \tau\left(a \right)>0$.
If $L^2\left(\sX, E,  \mu\right)$ is the Hilbert norm completion of $\Ga_c\left( \sX, E\right)$ then similarly  to \eqref{top_h_sp_eqn} one can construct a faithful representation
 \be\label{top_l2_xe_eqn}
\rho: C_0\left(\sX \right) \hookto L^2\left(\sX, E,  \mu\right).
 \ee
 On the other hand from the Lemma \ref{top_bundle_cs_ex_lem}  it follows that $\Ga\left(\sX, E\right)$ is {continuity structure for} $\sX$ {and the} $\left\{E_x\right\}_{x \in \sX}$ (cf. Definition \ref{operator_fields_continuity_defn}), such that
\bean
\Ga\left(\sX, E\right) \cong C\left(\sX, \left\{E_x\right\},  \Ga\left(\sX, E\right)\right),\\
\Ga_c\left(\sX, E\right) \cong C_c\left(\sX, \left\{E_x\right\},  \Ga\left(\sX, E\right)\right). 
\eean
If $p: \widetilde{\sX} \to \sX$ is a   transitive covering then from the Lemma \ref{top_tensor_lcompact_lem} it follows that there is the natural isomorphism
\be\label{top_lift_g_eqn}
\begin{split}
\phi: C_c\left(\widetilde \sX\right)\otimes_{C_0\left(\sX\right)} C_c\left(\sX,\left\{E_x\right\},  \Ga\left(\sX, E\right)\right) \xrightarrow{\approx}\\
\xrightarrow{\approx} C_c\left( \lift_p\left[C\left(\sX, \left\{E_x\right\},  \Ga\left(\sX, E\right)\right)\right]\right) 
\end{split}
\ee
of $C\left(\widetilde{\sX} \right)$-modules. From the Lemma \ref{top_lift_bundle_lem} it follows that \\$\lift_p\left[C\left(\sX, \left\{E_x\right\},  \Ga\left(\sX, E\right)\right)\right]\cong \Ga\left(\widetilde{\sX}, \widetilde{E} \right)$ where $\widetilde{E}$ is the {inverse image} of $E$ by $p$ (cf. Definition \ref{vb_inv_img_funct_defn}). Hence the equation \eqref{top_lift_g_eqn} can be rewritten by the following way
\be\label{top_lift_gs_eqn}
\phi: C_c\left(\widetilde \sX\right)\otimes_{C_0\left(\sX\right)} \Ga_c\left(\sX, E\right) \xrightarrow{\approx} \Ga_c\left(\widetilde{\sX}, \widetilde{E} \right).
\ee
If $p$ is a finite-fold transitive covering then from the Lemma \ref{state_cov_11_lem} it follows that  there is a functional
\be\label{top_wt_eqn}
\begin{split}
	\widetilde	\tau : C_c\left(\widetilde \sX \right)\to \C,\\
	\widetilde	a \mapsto \tau\left( \sum_{g \in G\left(\left.\widetilde \sX~\right|\sX \right) }g\widetilde a\right)=   \int_{\widetilde \sX} \widetilde a~ d\widetilde\mu 
\end{split}
\ee
where a sum $\sum_{g \in G\left(\left.\widetilde \sX~\right|\sX \right) }g\widetilde a$ is regarded as an element of $C\left(  \sX\right)$. The measure $\widetilde\mu\bydef \lift_p\mu$ is the $p$-lift of $\mu$ (cf. Definition \ref{top_lift_measure_defn}). 
Similarly to the above construction we have a Hilbert space $L^2\left(\widetilde{\sX}, \widetilde{E} \right)= L^2\left(\widetilde{\sX}, \widetilde{E}, \widetilde\mu \right)$ and a representation 
\be\label{comm_bundle_repr_eqn}
C_0\left(\widetilde \sX\right) \to B\left( L^2\left( \widetilde \sX, \widetilde E\right)\right).
\ee

\begin{lemma}\label{top_hilb_homo_lem}
	Following conditions hold:
	\begin{itemize}
		\item [(i)] the given by \eqref{top_lift_gs_eqn} map $\phi: C_c\left(\widetilde \sX\right)\otimes_{C_0\left(\sX\right)} \Ga_c\left(\sX, E\right) \xrightarrow{\approx} \Ga_c\left(\widetilde{\sX}, \widetilde{E} \right)$ can be extended up to an injective homomorphism 
		$$
	\varphi:	C_0\left(\widetilde\sX\right)\otimes_{C_0\left(\sX\right)} L^2\left(\sX, E\right) \hookto L^2\left(\widetilde{\sX}, \widetilde{E} \right)
		$$
		of left  $C_0\left(\widetilde\sX\right)$-modules,
		\item[(ii)] the image of $C_0\left(\widetilde \sX\right)\otimes_{C_0\left(\sX\right)} L^2\left(\sX, E\right)$ is dense in $L^2\left(\widetilde{\sX}, \widetilde{E} \right)$.
	\end{itemize}
\end{lemma}

\begin{proof}(i) If $\sum_{j = 1}^n \widetilde a_j \otimes \eta_j\in 	C_0\left(\widetilde \sX\right)\otimes_{C_0\left(\sX\right)} \Ga_c\left(\sX, E\right)$ then since the covering $p$ is finite-fold, the set $p^{-1}\left( \supp \eta_j\right)$ is compact for any $j =1, ..., n$. If $\widetilde b_j\in C_c\left(\widetilde \sX \right) $ is a covering sum for $p^{-1}\left( \supp \eta_j\right)$  (cf. Definition \ref{top_covering_sum_defn}) the  reader can proof that
	$$
\forall j =1,...,n \quad	\widetilde a_j \otimes \eta_j = 	\widetilde a_j 	\widetilde b_j \otimes \eta_j 
	$$
	So the map $\phi$ can be extended up to
	\bean
	\overline \phi:  C_0\left(\widetilde \sX\right)\otimes_{C_0\left(\sX\right)} \Ga_c\left(\sX, E\right) \xrightarrow{\approx} \Ga_c\left(\widetilde{\sX}, \widetilde{E} \right),\\
	\sum_{j = 1}^n \widetilde a_j \otimes \eta_j \mapsto \sum_{j=1}^n \phi \left( \widetilde a_j 	\widetilde b_j \otimes \eta_j \right). 
	\eean 
	Let $\sum_{j = 1}^n \widetilde a_j \otimes \xi_j\in 	C_0\left(\widetilde \sX\right)\otimes_{C_0\left(\sX\right)} L^2\left(\sX, E\right)$ be any element. From the above construction it turns out that $L^2\left( \sX, E\right)$ is the Hilbert norm completion of the pre-Hilbert space $\Ga_c\left(\sX, E\right)$  with the given by $\left(a, b \right) \mapsto\tau\left( \left\langle a, b \right\rangle_c\right)$ scalar product. Hence for any $j =1,..., n$ there is a net $\left\{\xi_{j\a}\right\}\subset \Ga_c\left(\sX, E\right)$ such that $\xi_j = \lim_\a \xi_{j\a}$ where we mean the convergence with respect to the Hilbert norm $\left\| \cdot\right\|_{L^2\left(\sX, E\right)}$. If $C = \max_{j = 1,...,n}\left\|\widetilde a_j\right\|$ then there is $\a_0$ such that
	$$
	\a \ge \a_0 \quad \Rightarrow \left\| \xi_j - \xi_{j\a}\right\|_{L^2\left(\sX, E\right)} < \frac{\eps}{nC}~~,
	$$
	so one has
$$
	\a \ge \a_0 \quad \Rightarrow \left\| \overline \phi\left( \sum_{j = 1}^n \widetilde a_j \otimes \xi_{j\a_0} - \sum_{j = 1}^n \widetilde a_j \otimes \xi_{j\a}\right) \right\|_{L^2\left(\widetilde \sX, \widetilde E\right)} < \eps.
	$$
From the above equation it follows that the net $\left\{\overline \phi\left(\sum_{j = 1}^n \widetilde a_j \otimes \xi_{j\a}\right) \right\}_\a\subset L^2\left(\widetilde \sX, \widetilde E\right)$ satisfies to the Cauchy condition, so it is convergent with respect to the topology of $L^2\left(\widetilde \sX, \widetilde E\right)$  (cf. \ref{top_cauchy_empt} and the Remark \ref{top_cauchy_rem}). We define
	$$
	\varphi\left(\sum_{j = 1}^n \widetilde a_j \otimes \xi_j \right) \bydef \lim_{\a}\overline\phi \left( \sum_{j = 1}^n \widetilde a_j \otimes \xi_{j\a}\right). 
	$$
	\\
	(ii) From the Lemma \ref{top_tensor_lccompact_iso_eqn} it follows that the image of $C_0\left(\widetilde \sX\right)\otimes_{C\left(\sX\right)} L^2\left(\sX, E\right)$ in $L^2\left(\widetilde \sX, \widetilde E\right)$ contains $\Ga_c\left(\widetilde{\sX}, \widetilde{E} \right)$, however $\Ga_c\left(\widetilde{\sX}, \widetilde{E} \right)$ is dense in $L^2\left(\widetilde{\sX}, \widetilde{E} \right)$. It follows that  $\phi\left(  C_0\left(\widetilde \sX\right)\otimes_{C_0\left(\sX\right)} L^2\left(\sX, E\right)\right)$ is dense in $L^2\left(\widetilde{\sX}, \widetilde{E} \right)$. 
\end{proof}
\begin{remark}
	The $C_0\left(\widetilde \sX\right)$-module structure of 
	$C_0\left(\widetilde \sX\right)\otimes_{C_0\left(\sX\right)} L^2\left(\sX, S\right)$ is given by an action
	\be\label{top_tens_a_eqn}
	\begin{split}
		C_0\left(\widetilde \sX\right)\times \left( C_0\left(\widetilde \sX\right)\otimes_{C_0\left(\sX\right)} L^2\left(\sX, E\right)\right)\to \left( C_0\left(\widetilde \sX\right)\otimes_{C_0\left(\sX\right)} L^2\left(\sX, E\right)\right), \\
		\widetilde b \left(\widetilde a \otimes \xi \right)  \mapsto \widetilde b\widetilde a \otimes \xi.
	\end{split}
	\ee
\end{remark}
\begin{remark}
	Indeed if $\widetilde \sX$ is compact then $C_0\left(\widetilde \sX\right)\cong C\left(\widetilde \sX\right)$ and the image of $C\left(\widetilde \sX\right)\otimes_{C_0\left(\sX\right)} L^2\left(\sX, E\right)$ coincides with $L^2\left(\widetilde{\sX}, \widetilde{E} \right)$. However this fact is not used here.
\end{remark}
If $\xi, \eta \in \Ga_c\left(\sX, E \right)$ then since $p$ is a finite-fold covering the set\\ $\widetilde \sV\bydef  p^{-1}\left(\supp \xi \cup \supp\eta  \right)$ is compact.  Let    $\sum_{\a \in \mathscr A} \widetilde a_\a$ is a covering sum for $\widetilde \sV$  {subordinated  to} $p$ (cf. Definition \ref{top_covering_sum_subordinated_defn}), and $\widetilde e_\a \bydef \sqrt {\widetilde a_\a}$. Denote by $\widetilde \sU_\a$ an open set such that $ \supp\widetilde a_\a\subset \widetilde\sU_\a$ and the restriction $p|_{\widetilde \sU_\a}: \widetilde \sU_\a \to p\left(\widetilde \sU_\a\right)$ is a homeomorphism. Let $e_\a \bydef \desc_p\left( \widetilde e_\a\right)$ (cf. Definition \ref{top_lift_desc_defn}). The specialization of the  given by the equation \eqref{induced_prod_equ} scalar product $\left(\cdot, \cdot \right)_{\text{ind}}$ on  $ C_0\left(\widetilde \sX\right)) \otimes L^2\left({\sX},{E}, \mu\right)$ satisfies to the following equation 
\begin{equation}\label{comm_ind_l2l_eqn}
	\begin{split}
		\left( \widetilde{a} \otimes \xi ,    \widetilde{b} \otimes \eta\right)_{\text{ind}}= \left(\xi, \left\langle \widetilde{a}, \widetilde{b}\right\rangle_{C\left(\widetilde{\sX} \right) } \eta \right)_{L^2\left({\sX},{E} \right)}=\\=\sum_{\a\in \widetilde{\mathscr A}}\left(\xi, \left\langle \widetilde{a}_{\a}\widetilde{a}, \widetilde{b}\right\rangle_{C\left(\widetilde{\sX} \right) } \eta \right)_{L^2\left({\sX},{E} \right)}
		= \sum_{\a\in {\mathscr A}}\left(\xi, \mathfrak{desc}_p\left( \widetilde{a}_{\a} \widetilde{a}^* \widetilde{b}\right)  \eta \right)_{L^2\left({\sX},{E} \right)}=
		\\
		=
		\sum_{\a\in {\mathscr A}}~ \int_\sX \left(\xi_x, \mathfrak{desc}_p\left( \widetilde{a}_{\a} \widetilde{a}^* \widetilde{b}\right)  \eta_x \right)_x d\mu=
		\\
		=
		\sum_{\a\in {\mathscr A}}~ \int_M \left(\mathfrak{desc}_p\left( \widetilde{e}_{\a} \widetilde{a}\right)\xi_x, \mathfrak{desc}_p\left( \widetilde{e}_{\a}  \widetilde{b}\right)  \eta_x \right)_x d\mu=
		\\
		=\sum_{\a\in {\mathscr A}}~  \int_{\widetilde{\sX}}\left(\widetilde{a} \lift_{\widetilde{\sU}_{\a}}\left({e}_{\a}\xi \right)_{\widetilde{x}},  \widetilde{b} ~\lift_{\widetilde{\sU}_{\a}}\left({e}_{\a}\eta \right)_{\widetilde{x}}\right)_{\widetilde{x}}d \widetilde{   \mu} =
		\\
		=\sum_{\a\in {\mathscr A}}~  \int_{\widetilde{\sX}}\widetilde{a}_{\a}\left(\widetilde{a} \lift_{\widetilde{\sU}_{\a}}\left(\xi \right)_{\widetilde{x}},  \widetilde{b} ~\lift_{\widetilde{\sU}_{\a}}\left(\eta \right)_{\widetilde{x}}\right)_{\widetilde{x}}d \widetilde{   \mu} =
		\\
		=  \int_{\widetilde{\sX}}\left( \widetilde{a} \lift_p\left(\xi \right)_{\widetilde{x}},  \widetilde{b} ~\lift_p\left(\eta \right)_{\widetilde{x}}\right)_{\widetilde{x}}d \widetilde{   \mu}
		= \left( \widetilde{a}~ \lift_p\left(\xi \right),  \widetilde{b} ~\lift_p\left(\eta \right)\right)_{L^2\left(\widetilde{\sX},\widetilde{E} \right)} =\\
		=  \left(\phi\left(  \widetilde{a}\otimes\xi\right) ,  \phi \left( \widetilde{b}\otimes\eta\right)\right)_{L^2\left(\widetilde{\sX},\widetilde{E} \right)}
	\end{split} 
\end{equation}
where $\phi$ is given by \eqref{top_lift_gs_eqn} and $\widetilde\mu\bydef \lift_p\mu$ is the $p$-lift of $\mu$ (cf. Definition \ref{top_lift_measure_defn}).  The $C_0\left(\sX \right)$-module $\Ga_c\left(\sX, E \right)$ is a dense subspace of  $L^2\left({\sX},{E} \right)$, so from the equation \eqref{comm_ind_l2l_eqn} it follows that
\be\label{comm_ind_l2_eqn}
\begin{split}
\forall \xi, \eta \in L^2\left(\widetilde{\sX},\widetilde{E} \right)\quad \forall \widetilde{a}, \widetilde{b} \in C_0\left( \widetilde{\sX}\right) 	\quad	\left( \widetilde{a} \otimes \xi ,    \widetilde{b} \otimes \eta\right)_{\text{ind}}=\\=  \left(\phi\left(  \widetilde{a}\otimes\xi\right) ,  \phi \left( \widetilde{b}\otimes\eta\right)\right)_{L^2\left(\widetilde{\sX},\widetilde{E} \right)}.
\end{split}
\ee
The equation \eqref{comm_ind_l2_eqn} means that $\left(\cdot, \cdot \right)_{\text{ind}}= \left(\cdot, \cdot \right)_{L^2\left(\widetilde{\sX},\widetilde{E} \right)}$, and taking into account the dense inclusion $C_0\left(\widetilde{\sX} \right) \otimes_{C_0\left(\sX \right) } \Ga\left( \sX, E\right) \subset L^2\left(\widetilde{\sX},\widetilde{E} \right)$ with respect to the Hilbert norm of $L^2\left(\widetilde{\sX},\widetilde{E}\right) $. One has a following lemma.
\begin{lemma}\label{comm_ind_lem}
	If the representation  $\widetilde{\rho}: C_0\left(\widetilde{\sX}\right) \to  B\left( \widetilde{   \H}\right)  $ is induced by the pair $$\left(C_0\left(\sX \right)\to B\left(  L^2\left(\sX,S \right)\right) ,\left(C_0\left(\sX \right) , C_0\left( \widetilde{\sX}\right) , G\left( \left.\widetilde{\sX}~\right|\sX\right) \right)  \right)$$ (cf. Definition \ref{induced_repr_fin_defn}) then following conditions hold:
	\begin{enumerate}
		\item[(i)] there is the homomorphism of Hilbert spaces $\widetilde{   \H}\cong L^2\left(\widetilde{\sX},\widetilde{E} \right)$,
		\item[(ii)] the representation $\widetilde{\rho}$ is given by the natural action of $C_0\left(\widetilde{\sX}\right)$ on $ L^2\left(\widetilde{\sX},\widetilde{E} \right)$ (cf.   \eqref{comm_bundle_repr_eqn}).
	\end{enumerate}
\end{lemma}
\begin{proof}
	(i) Follows from the above construction (cf. equation \eqref{comm_ind_l2_eqn}).\\
	(ii) 	The natural  action $C_0\left(\widetilde\sX\right)\times \left( C_0\left(\widetilde \sX\right)\otimes_{C_0\left(\sX\right)} L^2\left(M, E\right)\right)\to C_0\left(\widetilde \sX\right)\otimes_{C_0\left(\sX\right)} L^2\left(\sX, E\right)$ is given by the equation \eqref{top_tens_a_eqn}. However the formula \eqref{top_tens_a_eqn} is a specialization of the general equation \eqref{ind_act_form_eqn} which yields an induced representations.
	
\end{proof}

\subsection{Coverings of spectral triples}\label{comm_d_sec}
\paragraph*{}
Consider the situation of the Section \ref{comm_induced_finite_sec}. Moreover assume  that $\sX = M$ where $M$ is a compact Riemannian manifold and $E= S$ is a smooth spinor bundle (cf. Section \ref{spin_mani_sec}). Let  $\left(\Coo\left(M \right), L^2\left( M, S\right), \slashed D , J\right)$ be an explained in the Appendix \ref{comm_sp_tr_sec} commutative spectral triple. If $p: \widetilde M \to M$ is a transitive finite-fold covering then $\widetilde M$ has a (unique) structure of Riemannian manifold given by the Proposition \ref{comm_cov_mani_prop}. From the Lemmas \ref{top_fin_necassary_lem}   and \ref{top_fin_sufficient_lem}   it follows that the quadruple
\be\label{top_am_eqn}
\left(A, \widetilde A   , G\left(\left. \widetilde A \right|A \right) , \pi \right)\bydef \left(C\left(M \right), C\left(\widetilde M \right)   , G\left(\left. \widetilde M \right|M \right) , C\left(p \right) \right)  
\ee
is an (unital) noncommutative finite-fold covering. If  $\widetilde S$ is  the {inverse image} of $S$ by $p$  (cf. Definition \ref{vb_inv_img_funct_defn}) then from the Lemma \ref{comm_ind_lem} it follows that the representation $\widetilde{\rho}$ is given by the natural action of $C\left(\widetilde{M}\right)$ on $ L^2\left(\widetilde{M},\widetilde{S} \right)$ is induced by the pair $$\left(C\left(M \right)\to B\left(  L^2\left(M,S \right)\right) ,\left(C\left(M \right) , C\left( \widetilde{M}\right) , G\left( \left.\widetilde{M}~\right|M\right) \right)  \right)$$ (cf. Definition \ref{induced_repr_fin_defn}). Both spaces $\Ga^\infty\left( M, S\right)$ and  $\Ga^\infty\left( \widetilde M, \widetilde S\right)$ of smooth sections (cf.  Definition \ref{top_sm_sec_defn}) are  $\Coo\left(M \right)$-module and $\Coo\left(\widetilde M \right)$-module respectively. We leave to the reader the proof of the following isomorphism
$$
\mathscr S^{\Ga^\infty\left( \widetilde M, \widetilde S\right)}\cong p^{-1}\mathscr S^{\Ga^\infty\left(  M,  S\right)}
$$
where both $\mathscr S^{\Ga^\infty\left( \widetilde M, \widetilde S\right)}$ and
$\mathscr S^{\Ga^\infty\left(  M,  S\right)}$ are  $\Ga^\infty\left( \widetilde M, \widetilde S\right)$-sheaf and $\Ga^\infty\left(  M,  S\right)$-sheaf respectively (cf. Definition \ref{top_x_sheaf_defn}) and $p^{-1}$ means an inverse image of sheaves (cf. Definition \ref{sheaf_inv_im_defn}).  
From \eqref{top_d_shaef_m_eqn} it follows that there is an operator 
$$
\widetilde \Dslash  : \Ga^\infty\left(\widetilde M,\widetilde S \right)\to \Ga^\infty\left(\widetilde M,\widetilde S \right).
$$
such that 
\be\label{top_d_shaef_mt_eqn}
\begin{split}
	\widetilde \Dslash \bydef p^{-1}\Dslash \in \mathscr End \left(\mathscr S^{\Ga^\infty\left( \widetilde M, \widetilde S\right)} \right) \left(\widetilde M\right)\bydef\\\bydef \mathscr Hom \left(\mathscr S^{\Ga^\infty\left( \widetilde M, \widetilde S\right)} , \mathscr S^{\Ga^\infty\left( \widetilde M, \widetilde S\right)}\right) \left(\widetilde M\right).
\end{split}
\ee
If $\H^\infty $ and $\widetilde\H^\infty $ are given by
\bean
\H^\infty \bydef	\bigcap_{k\in\N} \Dom \slashed D^k\subset L^2\left(  M, S\right),\\
\widetilde \H^\infty \bydef	\bigcap_{k\in\N} \H^\infty om \widetilde{\slashed D}^k\subset L^2\left( \widetilde M,\widetilde S\right)
\eean 
then
\bean
\H^\infty = \Ga^\infty\left( M, S\right)=\mathscr S^{\Ga^\infty\left(  M,  S\right)} \left(  M\right),\\
\widetilde \H^\infty = \Ga^\infty\left( \widetilde M, \widetilde S\right)=\mathscr S^{\Ga^\infty\left( \widetilde M, \widetilde S\right)}\left( \widetilde M\right)
\eean  
and one has
\bean 
\Coo\left(M \right)\subset \L^\dagger \left(\H^\infty  \right), \quad
\slashed D \in 	\L^\dagger \left(\H^\infty  \right),\\
\Coo\left(\widetilde M \right)\subset \L^\dagger \left(\widetilde \H^\infty  \right), \quad
\widetilde{\slashed D} \in 	\L^\dagger \left(\widetilde \H^\infty  \right)
\eean
where both $\L^\dagger \left( \H^\infty  \right)$ and $\L^\dagger \left(\widetilde \H^\infty  \right)$ are given by the equation \eqref{l_dag_eqn}.
\begin{exercise}\label{top_cc_exer}
Prove that the specialization of the given by the equation \eqref{cc_eqn} operator $\widetilde J$ is an extension of the given by \eqref{top_wtc_eqn} operator $\widetilde C$.
\end{exercise}
\begin{empt}
	Here we consider a related to the spectral triple 
	$$
	\left(\A, \H ,D, J \right) \bydef	\left(\Coo\left(M \right), L^2\left( M, S\right), \slashed D ,J \right)
	$$
	specialization of the construction \ref{triple_conn_lift_empt}. 
	From the Exercise \ref{smooth_soft_exer} it turns out that $\Coo\left(M \right)\subset C\left( M\right)$ is a $c$-soft *-algebra (cf. Definition \ref{top_soft_r_defn}). From the  Theorem \ref{top_oa_cov_thm} it follows that the quadruple	$$
	\left(\A, \widetilde \A   , G\left(\left. \widetilde \A \right|\A \right) , \pi \right)\bydef  \left(\Coo\left(M\right), \Coo\left(\widetilde M \right), G\left(\left.\widetilde M~ \right| M\right), \left.C_0\left(p \right)\right|_{\Coo\left(M\right)}\right)
	$$
	is a finite-fold covering of bounded operator *-algebras  (cf. Definition \ref{fin_oa_defn}).	 
	We leave to the reader the proof of that   $\widetilde \A$ is a left  and right finitely generated projective $\A$-module.  Both $A$ and $\widetilde A$ are $C^*$-norm completions of $\A$ and  $\widetilde \A$ respectively, so from the Definition  \ref{fin_oa_defn} it follows that there is  an unital connected noncommutative finite-fold covering $\left(A, \widetilde A, G\left( \left.\widetilde A\right| A\right), \widetilde \pi  \right)$ (cf. Definition \ref{fin_unital_defn} and the equation \eqref{top_am_eqn}) such that $\pi= \widetilde\pi|_{\A}$. If 
	$$
	\left(  \rho : A \to B\left(\H \right)\right) \bydef \left(C\left( M\right)\to B\left(L^2\left( M, S\right) \right)  \right) 
	$$
	is a representation  which corresponds to the spectral triple $\left( \A, \H, D\right)$ then from the Lemma \ref{comm_ind_lem} it follows that the induced by the pair $\left(\rho, \left(A, \widetilde A, G\left( \left.\widetilde A\right| A\right), \widetilde \pi  \right)\right)$ representation $\widetilde \rho : \widetilde A \to B\left(\widetilde{   \H} \right)$ (cf. Definition \ref{induced_repr_fin_defn}) is given by
	$$
	\widetilde \rho :C\left(\widetilde M \right)\hookto B\left( L^2\left(\widetilde{M},\widetilde{S} \right)\right).    
	$$
	The space $M$ is compact so there is is a subordinated  to  $p:\widetilde{M}\to{M}$ covering sum for $\widetilde M$ (cf. Definitions \ref{top_covering_sum_defn}, \ref{top_covering_sum_subordinated_defn}).
	\be\label{top_ta_eqn}
	\sum_{\widetilde \a \in \widetilde {\mathscr A}  } \widetilde{a}_{\widetilde\a} = 1_{C\left(\widetilde{M} \right) }.
	\ee
	Similarly to the Proposition \ref{top_smooth_part_unity_prop}   one can assume that $\widetilde{a}_{\widetilde\a}\in \Coo\left( \widetilde{M}\right)$, so one has
	\be\label{top_te_eqn}
	\begin{split}
		\widetilde{e}_{\widetilde\a}\bydef \sqrt{\widetilde{a}_{\widetilde\a}}\in \Coo\left( \widetilde{M}\right),\\
		\sum_{ \widetilde\a \in \widetilde{ \mathscr A} } \widetilde{e}^2_{\widetilde\a} =\sum_{ \widetilde\a \in \widetilde{ \mathscr A} } \widetilde{e}_{\widetilde\a} \desc_p\left(\widetilde{e}_{\widetilde\a} \right)= 1_{C\left(\widetilde{M} \right) }
	\end{split}
	\ee 
	where $\desc_p$ is given by the Definition \ref{top_lift_desc_defn} $p$-descent.	
	For any $\widetilde{a} \in \Coo\left(\widetilde{M} \right)$ a following condition holds
	\be\label{top_ph_ps_eqn}
	\widetilde{a} = \sum_{\widetilde{\a}\in \widetilde{\mathscr A}} \widetilde{\ga}_{\widetilde{\a}}= \sum_{\widetilde{\a}\in \widetilde{\mathscr A}} \widetilde{\phi}_{\widetilde{\a}}\widetilde{\psi}_{\widetilde{\a}}; \text{ where } \widetilde{\phi}_{\widetilde{\a}}\bydef \widetilde{e}_{\widetilde{\a}},~ \widetilde{\psi}_{\widetilde{\a}}\bydef \widetilde{e}_{\widetilde{\a}}\widetilde{a},~ \widetilde{\ga}_{\widetilde{\a}}\bydef\widetilde{a}_{\widetilde{\a}}\widetilde{e}^2_{\widetilde{\a}}.
	\ee
	Let $\Om^1_{\slashed D}$ be the {module of differential forms associated}  with the given by  \eqref{comm_sp_tr_eqn} spectral triple  $\left(\Coo\left({M} \right), L^2\left(M,S \right), \slashed D\right)$ (cf. Definition \ref{ass_cycle_defn}), and denote by
	\be\label{top_descs_eqn}
	{\phi}_{\widetilde{\a}} \bydef \mathfrak{desc} \left(\widetilde{\phi}_{\widetilde{\a}} \right),\quad {\psi}_{\widetilde{\a}} \bydef \mathfrak{desc} \left(\widetilde{\psi}_{\widetilde{\a}} \right), \quad {\ga}_{\widetilde{\a}} \bydef \mathfrak{desc} \left(\widetilde{\ga}_{\widetilde{\a}} \right)\in \Coo\left({M} \right)
	\ee
	Let us define a $\C$-linear map
	\be\label{top_nabla_eqn}
	\begin{split}
		{\slashed{\nabla}}: \Coo\left(\widetilde{M} \right) \to \Coo\left(\widetilde{M} \right) \otimes_{\Coo\left({M} \right)}\Om^1_{\slashed D},\\
		\widetilde{a} \mapsto \sum_{\widetilde{\a}\in \widetilde{\mathscr A}}\left(  \widetilde{\phi}_{\widetilde{\a}} \otimes \left[\slashed D, {\psi}_{\widetilde{\a}}  \right] +\widetilde{\psi}_{\widetilde{\a}} \otimes \left[\slashed D, {\phi}_{\widetilde{\a}}   \right] \right).
	\end{split}
	\ee
	For any $a \in \Coo\left({M} \right)$ from $\left[\slashed D, {\psi}_{\widetilde{\a}} a  \right] = {\psi}_{\widetilde{\a}} \left[\slashed D, a  \right] + \left[\slashed D, {\psi}_{\widetilde{\a}}   \right]a$, $~~\widetilde{\phi}_{\widetilde{\a}} \otimes {\psi}_{\widetilde{\a}}$ = $\widetilde{\phi}_{\widetilde{\a}}  \widetilde{\psi}_{\widetilde{\a}}\otimes 1_{\Coo\left( M\right) }$ and the equation  \eqref{comm_matr_x_eqn} it follows that
	\bean
	{\slashed{\nabla}}\left(\widetilde{a}a \right)=  \sum_{\widetilde{\a}\in \widetilde{\mathscr A}}\left(  \widetilde{\phi}_{\widetilde{\a}} \otimes \left[\slashed D, {\psi}_{\widetilde{\a}} a  \right] +\widetilde{\psi}_{\widetilde{\a}}a \otimes \left[\slashed D, {\phi}_{\widetilde{\a}}   \right] \right) =
	\\
	= \sum_{\widetilde{\a}\in \widetilde{\mathscr A}}\left(  \widetilde{\phi}_{\widetilde{\a}} \otimes \left[\slashed D, {\psi}_{\widetilde{\a}}  \right]a+\widetilde{\phi}_{\widetilde{\a}} \otimes {\psi}_{\widetilde{\a}}\left[\slashed D,  a  \right] +\widetilde{\psi}_{\widetilde{\a}} \otimes a \left[\slashed D, {\phi}_{\widetilde{\a}}   \right]\right) =
	\\
	= \sum_{\widetilde{\a}\in \widetilde{\mathscr A}}\left(  \widetilde{\phi}_{\widetilde{\a}} \otimes \left[\slashed D, {\psi}_{\widetilde{\a}}  \right]a+\widetilde{\phi}_{\widetilde{\a}} \widetilde{\psi}_{\widetilde{\a}}\otimes \left[\slashed D,  a  \right] +\widetilde{\psi}_{\widetilde{\a}} \otimes  \left[\slashed D, {\phi}_{\widetilde{\a}}   \right]a\right) =
	\\
	= \sum_{\widetilde{\a}\in \widetilde{\mathscr A}}\left(  \widetilde{\phi}_{\widetilde{\a}} \otimes \left[\slashed D, {\psi}_{\widetilde{\a}}   \right] +\widetilde{\psi}_{\widetilde{\a}} \otimes \left[\slashed D, {\phi}_{\widetilde{\a}}  \right] \right) a+
	\sum_{\widetilde{\a}\in \widetilde{\mathscr A}}\widetilde{\phi}_{\widetilde{\a}}\widetilde{\psi}_{\widetilde{\a}}\otimes\left[ \slashed D,a\right]=
	\widetilde{\slashed{\nabla}}\left( \widetilde{a}\right) a + \widetilde{a}\otimes  \left[ \slashed D,a\right]. 
	\eean
	The comparison of the above equation and \eqref{conn_triple_eqn} states that $\widetilde{\slashed{\nabla}}$ is a connection. If $g \in  G\left( \left.\widetilde{M}~\right|M\right)$ then from $\mathfrak{desc}\left(g\widetilde{\phi}_{\widetilde{\a}}\right) = {\phi}_{\widetilde{\a}} $ and $\mathfrak{desc}\left(g\widetilde{\psi}_{\widetilde{\a}}\right) = {\psi}_{\widetilde{\a}} $ it follows that
	\bean
	{\slashed{\nabla}}\left(g\widetilde{a} \right)=  {\slashed{\nabla}}\left(\sum_{\widetilde{\a}\in \widetilde{\mathscr A}} \left( g\widetilde{\phi}_{\widetilde{\a}}\right) \left( g\widetilde{\psi}_{\widetilde{\a}}\right)\right)=\\=\sum_{\widetilde{\a}\in \widetilde{\mathscr A}}\left( g \widetilde{\phi}_{\widetilde{\a}} \otimes \left[\slashed D, {\psi}_{\widetilde{\a}}   \right] +g\widetilde{\psi}_{\widetilde{\a}} \otimes \left[\slashed D, {\phi}_{\widetilde{\a}}   \right] \right),  
	\eean
	i.e. ${\slashed{\nabla}}$ is  $G\left( \left.\widetilde{M}~\right|M\right)$-{equivariant} (cf. equation \eqref{conn_equ_conn}).  We leave to the reader the proof of that the connection $\slashed{\nabla}$ is Hermitian  (cf. Definition \ref{herma_conn_defn}).
	If $\widetilde D: \widetilde \A\otimes_\A \H^\infty  \to \widetilde \A\otimes_\A \H^\infty $ is given by \eqref{wtd_eqn} then
	\be\label{top_dirac_om_fin_lift_eqn}
	\slashed \nabla \widetilde a = \sum_{j = 1}^n \widetilde a_j \otimes \om_j \quad\Rightarrow \quad\left[\widetilde D, \widetilde a\right]=  \sum_{j = 1}^n \widetilde a_j \otimes \om_j\in \widetilde\A \otimes_\A \Om^1_D.
	\ee
	
\end{empt}
\begin{exercise}
	Using \eqref{top_lift_gs_eqn} prove that there is a $\Coo\left(\widetilde M \right)$-module isomorphism
	\be\label{top_phinf_eqn}
	\begin{split}
		\phi_\infty :\Coo\left(\widetilde M \right)\otimes \Ga^\infty\left(M, S \right)  \to \Ga^\infty\left(\widetilde M, \widetilde S \right),\\
		\sum_{j=1}^n \widetilde a_j \otimes \xi_j\mapsto  \sum_{j=1}^n\widetilde a_j \lift_p\left(\xi_j\right)
	\end{split}
	\ee
\end{exercise}
A following table reflects the mapping between general theory and its commutative specialization\\ \\
\begin{tabular}{|c|c|c|}
	\hline
	&General theory & Commutative specialization\\ 
	\hline
	&	&\\
	Hilbert spaces & $\H$  and $\widetilde\H$ &  $L^2\left(M, S \right)$ and $L^2\left(\widetilde M, \widetilde S \right)$\\ & & \\
	Pre-$C^*$-algebras	& $\A$ and $\widetilde\A$ & $\Coo\left(M \right)$ and $\Coo\left(\widetilde M \right)$   \\  & & \\
Charge conjugation 	& $J$ and $\widetilde J$ & $C$ and $\widetilde C$ (cf. Proposition \ref{pr:charge-conj_prop})  \\  & & \\
	Dirac operators & $D$ & $\Dslash$  \\
	& $\widetilde{D}$ & ?\\  & & \\
	&$\H^\infty\bydef \bigcap_{n =1}^\infty \Dom D^n\subset \H$ & $\Ga^\infty\left(M, \sS \right)= \bigcap_{n =1}^\infty \Dom \Dslash^n$ \\& & \\
	\hline
\end{tabular}
\\
\\
\\
If $\slashed\nabla \widetilde{a}= \sum_{j = 1}^m\widetilde{a}_j \otimes \om_j\in \Coo\left(\widetilde M \right)\otimes_{\Coo\left( M \right)}\Om^1_{\Dslash}$ then for all $\xi \in \Ga^\infty\left(M, S \right)$ denote by
$$
\widetilde{a}_\xi = \phi_\infty\left(\sum_{j = 1}^m\widetilde{a}_j \otimes \om_j\left( \xi\right)  \right)= \sum_{j=1}^m\widetilde a_j \lift_p\left(\om_j\left( \xi\right)\right)\in \Ga^\infty\left(\widetilde M, \widetilde S \right)
$$
where $\phi_\infty$ is given by the Equation \eqref{top_phinf_eqn}. Taking into account the equation \eqref{top_nabla_eqn} one has
$$
\widetilde{a}_\xi = \sum_{\widetilde{\a}\in \widetilde{\mathscr A}}\left(  \widetilde{\phi}_{\widetilde{\a}} \lift_p\left(  \left[\slashed D, {\psi}_{\widetilde{\a}}  \right]\xi \right)  +\widetilde{\psi}_{\widetilde{\a}} \lift_p\left( \left[\slashed D, {\phi}_{\widetilde{\a}}   \right]\xi\right)  \right).
$$
The partition of unity \eqref{top_ta_eqn} compliant with  the covering $p:\widetilde{M}\to{M}$, so for any $\widetilde{\a}\in \widetilde{\mathscr A}$ there is an open subset $\widetilde{\sU}_{\widetilde \a}\subset \widetilde M$ such that 
\begin{itemize}
	\item $\supp \widetilde{a}_{\widetilde \a} \subset \widetilde{\sU}_{\widetilde \a}$.
	\item The restriction $\left.p\right|_{\widetilde{\sU}_{\widetilde \a}}$ is injective. 
\end{itemize}
From \eqref{top_ph_ps_eqn} it follows that $\supp \widetilde{\phi}_{\widetilde \a},\quad \supp \widetilde{\psi}_{\widetilde \a} ,\quad \supp \widetilde{\ga}_{\widetilde \a} \subset \widetilde{\sU}_{\widetilde \a}$. From the Equations \eqref{comm_lift_desc_l_eqn}  and  \eqref{top_descs_eqn} it follows that
$$
\widetilde{\phi}_{\widetilde \a}= \lift^p_{\widetilde{\mathcal U}_{\widetilde \a}}\left(\phi_{\widetilde \a} \right), \quad  \widetilde{\psi}_{\widetilde \a}= \lift^p_{\widetilde{\mathcal U}_{\widetilde \a}}\left(\psi_{\widetilde \a} \right), \quad  \widetilde{\ga}_{\widetilde \a}= \lift^p_{\widetilde{\mathcal U}_{\widetilde \a}}\left(\ga_{\widetilde \a} \right),
$$ hence taking into account \eqref{top_lift_pa_eqn} one has
\bean
\widetilde{a}_\xi = \sum_{\widetilde{\a}\in \widetilde{\mathscr A}}\left(  \lift^p_{\widetilde{\mathcal U}_{\widetilde \a}}\left(\phi_{\widetilde \a} \right) \lift_p\left(  \left[\slashed D, {\psi}_{\widetilde{\a}}  \right]\xi \right)  +\lift^p_{\widetilde{\mathcal U}_{\widetilde \a}}\left(\psi_{\widetilde \a} \right) \lift_p\left( \left[\slashed D, {\phi}_{\widetilde{\a}}   \right]\xi\right)  \right)=
\\
= \sum_{\widetilde{\a}\in \widetilde{\mathscr A}}\left(  \lift^p_{\widetilde{\mathcal U}_{\widetilde \a}}\left(\phi_{\widetilde \a}  \left[\slashed D, {\psi}_{\widetilde{\a}}  \right]\xi \right)  +\lift^p_{\widetilde{\mathcal U}_{\widetilde \a}}\left(\psi_{\widetilde \a}  \left[\slashed D, {\phi}_{\widetilde{\a}}   \right]\xi\right)  \right).
\eean
Otherwise from \eqref{comm_matr_x_eqn} it turns out that ${\psi}_{\widetilde{\a}}  \left[\slashed D, {\phi}_{\widetilde{\a}}   \right] =   \left[\slashed D, {\phi}_{\widetilde{\a}}   \right]{\psi}_{\widetilde{\a}}$, so it follows that
\bean
\widetilde{a}_\xi = 
\sum_{\widetilde{\a}\in \widetilde{\mathscr A}}\left(  \lift^p_{\widetilde{\mathcal U}_{\widetilde \a}}\left(\phi_{\widetilde \a}  \left[\slashed D, {\psi}_{\widetilde{\a}}  \right]\xi \right)  +\lift^p_{\widetilde{\mathcal U}_{\widetilde \a}}\left(\psi_{\widetilde \a} \left[\slashed D, {\phi}_{\widetilde{\a}}   \right]\xi\right)  \right)=
\\
=
\sum_{\widetilde{\a}\in \widetilde{\mathscr A}}\left(  \lift^p_{\widetilde{\mathcal U}_{\widetilde \a}}\left(\phi_{\widetilde \a}  \left[\slashed D, {\psi}_{\widetilde{\a}}  \right]\xi  +\psi_{\widetilde \a} \left[\slashed D, {\phi}_{\widetilde{\a}}   \right]\xi\right)  \right)=
\\
=\sum_{\widetilde{\a}\in \widetilde{\mathscr A}}\left(  \lift^p_{\widetilde{\mathcal U}_{\widetilde \a}}\left(  \left[\slashed D, \phi_{\widetilde \a} {\psi}_{\widetilde{\a}}  \right]\xi \right)    \right) = \sum_{\widetilde{\a}\in \widetilde{\mathscr A}}\left(  \lift^p_{\widetilde{\mathcal U}_{\widetilde \a}}\left(  \left[\slashed D, \ga_{\widetilde \a}  \right]\xi \right)    \right) .
\eean
If $ p^{-1}\left(\Dslash \right)\in  \mathscr End\left( \mathscr  G^\infty\right)$ is the   $p$-{inverse image} of $\Dslash$ (cf. Definition \ref{top_smooth_inv_im_defn}) then taking into account  \eqref{top_ldx_eqn} one has
\bean
\forall \widetilde{\a}\in \widetilde{\mathscr A}\quad \lift^p_{\widetilde{\mathcal U}_{\widetilde \a}}\left(  \left[\slashed D, \ga_{\widetilde \a}  \right]\xi \right) =\left[p^{-1}\left(\Dslash \right) , \widetilde{\ga}_{\widetilde{\a}}\right]\lift_p\left(\xi \right)
\eean
From the above equation it follows that
\bean
\widetilde{a}_\xi 
=\sum_{\widetilde{\a}\in \widetilde{\mathscr A}}
\left[p^{-1}\left(\Dslash \right) , \widetilde{\ga}_{\widetilde{\a}}\right]\lift_p\left(\xi \right) =\left[p^{-1}\left(\Dslash \right) , \sum_{\widetilde{\a}\in \widetilde{\mathscr A}}\widetilde{\ga}_{\widetilde{\a}}\right]\lift_p\left(\xi \right)=
\\
= \left[p^{-1}\left(\Dslash \right) , \widetilde{a}\right] \lift_p\left(\xi \right). 
\eean

If $\widetilde{D}$ is the specialization of the given by \eqref{wtd_eqn} operator then
\bean
\phi_\infty\left(\widetilde{D}\left(\widetilde{a}\otimes \xi \right)  \right) = \widetilde{a}_\xi + \phi_\infty\left( \widetilde{a} \otimes D\xi\right) = 
\\
= \left[p^{-1}\left(\Dslash \right) , \widetilde{a}\right] \lift_p\left(\xi \right)+ \widetilde{a}\lift_p\left(\Dslash\xi \right)=
\\
= \left[p^{-1}\left(\Dslash \right) , \widetilde{a}\right] \lift_p\left(\xi \right)+ \widetilde{a}p^{-1}\left(\Dslash \right)\lift_p\left(\xi \right)=  p^{-1}\Dslash   \widetilde{a} \lift_p\left(\xi \right)=
\\
= p^{-1}\left(\Dslash \right)\phi_\infty\left( \widetilde{a}\otimes\xi\right),
\eean
or, shortly
\be\label{comm_liftd_comp_eqn} 
\phi_\infty\left(\widetilde{D}\left(\widetilde{a}\otimes \xi \right)  \right) = p^{-1}\Dslash \phi_\infty\left( \widetilde{a}\otimes\xi\right),
\ee
However $\phi_\infty$ is the isomorphism, it can be regarded as the identical map, so the equation \eqref{comm_liftd_comp_eqn} is equivalent to 
\be\label{comm_dxi_eqn}
\widetilde{D}\left(\widetilde{a}\otimes \xi \right)   = p^{-1}\left(\Dslash \right)\left( \widetilde{a}\otimes\xi\right),
\ee in result one has
\be\label{comm_d_eqn}
\widetilde{D} =  p^{-1}{\slashed D}.
\ee
In more detail the above formula has the following meaning
\bean
\begin{split}
	\widetilde{\slashed{\nabla}}\left(\widetilde{a} \right)= \sum_{j=1}^m \widetilde{a}_j\otimes  \left[\Dslash, a_j \right]~ \Rightarrow~ \sum_{j=1}^m \widetilde{a}_j\otimes  \left[\Dslash, a_j \right]\xi + \widetilde a \otimes \Dslash\xi =p^{-1}{\slashed D}\left(\widetilde a \otimes \xi \right) 
\end{split}
\eean
where  $p^{-1}\Dslash$ is the   $p$-{inverse image} of $\Dslash$ (cf. Definition \ref{top_smooth_inv_im_defn}). From the equations \eqref{top_geom_dirac_lift_eqn}, \eqref{comm_d_eqn} and the Exercise \ref{top_cc_exer} one has a following theorem.

\begin{thm}\label{comm_fin_sp_tr_thm} 
	In the above situation the  $\left(C\left(M \right) , C\left( \widetilde{M}\right) , G\left(\left.\widetilde M\right| M\right), C\left(p \right)  \right)$-{lift} of $\left(C^{\infty}\left( M\right) , L^2\left( M, S\right) ,\slashed D, J\right)$ (cf. Definition \ref{spectral_triple_fin_lift_defn}) coincides with the given by the Definition \ref{top_geom_lift_defn} geometrical $p$-lift $\left(\Coo\left( \widetilde M\right) , L^2\left(\widetilde M, \widetilde S\right) ,\widetilde \Dslash \cong p^{-1} \Dslash, \widetilde J \right)$.
\end{thm}

\begin{remark}
	The Theorem \ref{comm_fin_sp_tr_thm} is a noncommutative analog of the Proposition \ref{comm_cov_mani_prop}.
\end{remark}

\subsection{Unoriented spectral triples}\label{comm_sp_tr} 

\paragraph*{}
If $M$ is an unoriented   Riemannian manifold (cf.  Definitions \ref{ori_man_defn} and \ref{riemann_mani_defn}) then there is a transitive covering $p: \widetilde M \to M$ such that:
\begin{itemize}
	\item $G\left( \left.\widetilde{M}\right| M\right) \cong \Z_2$,
	\item $\widetilde{M}$ has a structure of Riemannian manifold which comes from the covering $p$ (cf. Proposition \ref{comm_cov_mani_prop}), 
	\item the manifold $\widetilde M$ is oriented.
\end{itemize}
Let $\C\ell(M) \to M$ be a given by \eqref{st_cliffod_eqn} Clifford  bundle.
\begin{exercise}\label{top_unori_exer}
Suppose that there is a complex vector bundle $S\to M$ with a sesquilinear form (cf. Definition \ref{top_herm_bundle_form_defn}) such that following conditions hold:
\begin{itemize}
	\item for all $x \in M$ there is a $*$-isomorphism $\C\ell_x(M) \simeq \End(S_x)$ (cf. Equation \eqref{st_spinor_eqn}),
	\item there  is an antilinear endomorphism $C : \Ga\left(M, S \right)\to \Ga\left(M, S \right)$ which satisfies to the Proposition   \ref{pr:charge-conj_prop}.
\end{itemize}
Prove that there is a natural unoriented spectral triple
\be\label{top_unori_eqn}
\left( \Coo\left( M\right), L^2\left(M, S \right), \Dslash, C  \right).
\ee
\end{exercise}

\begin{remark}
	If a tangent bundle $TM$ of an unoriented  Riemannian manifold $M$ is trivial then $M$ satisfies to the conditions of the Exercise \ref{top_unori_exer}. In particular one has an  unoriented spectral triple
	$$
	\left( \Coo\left( \R P^3\right), L^2\left(\R P^3, S \right), \Dslash, C  \right)
	$$
	where $\R P^3$ is a real projective space.
\end{remark}
\section{Infinite coverings}

\subsection{Inverse limits of coverings in topology}\label{top_inf_to_sec}
\paragraph{}
Here we discuss a related to commutative $C^*$-algebras specialization of the explained in the Chapter \ref{infinite_covering_chap} theory.
\begin{empt}\label{top_fin_cov_cat_empt}  Let 
	Let $\widetilde \sX$ be a topological space with an action $  G\times \widetilde \sX\to \widetilde\sX$ of residually finite group $  G$ (cf. Definition \ref{residually_finite_defn}) of properly discontinuous  group of homeomorphisms (cf. Definition \ref{top_properly_disc_group_defn}). Let $\sX \bydef \widetilde\sX/   G$ and $p: \widetilde \sX\to\sX$ be a natural covering. 	For any finite factor group $G_\la =  G/ H_\la$ we define a space $\sX_\la \bydef \widetilde \sX/ H_\la$. Then there is a category of topological spaces and finite-fold transitive coverings given by
	\be\label{top_g_x_cat_eqn}
	\mathfrak{S}_p \bydef \left\{\left\{\sX_\la\right\}_{\la \in \La}, \left\{p^\mu_\nu:\sX_\mu\to \sX_\nu\right\}_{\substack{\mu,\nu \in \La\\\mu\ge\nu}}\right\}.
	\ee
	(cf. \eqref{top_category_fin_eqn})
	Usage of the {finite covering algebraic functor} $C_0$ (cf. Definition \ref{top_c_funct_defn}) yields a category  given by
	\be\label{top_g_a_x_cat_eqn}
	\begin{split}
		\mathfrak{S}_{C_0\left(p\right) } \bydef \\
		\left\{ \left\{ C_0\left( p_\la\right)  :C_0\left( \mathcal{X}\right)  \hookto C_0\left( \mathcal{X}_\la\right) \right\}, \left\{ C_0\left( p^\mu_\nu\right)  :C_0\left( \mathcal{X}_\mu\right)  \hookto C_0\left( \mathcal{X}_\nu\right) \right\}  \right\}
	\end{split}
	\ee
	Both categories are equivalent to $  G$-{category} (cf. Definition \ref{g_category_defn})
\end{empt}
\begin{definition}\label{top_fin_cov_cat_defn}
	We say that both categories  \eqref{top_g_x_cat_eqn} and \eqref{top_g_a_x_cat_eqn} are the \textit{topological} $\widetilde \sX$-$  G$-\textit{category} and the \textit{algebraic} $\widetilde \sX$-$  G$-\textit{category} respectively.
\end{definition}

	Let $\widetilde p : \widetilde \sX \to \sX$ be a transitive covering with connected,  locally compact, Hausdorff space  $\widetilde \sX$ such that a group  $G\left(\left. \widetilde\sX~\right|\sX \right)$ is residually finite (cf. Definition \ref{residually_finite_defn}). Let $\La$ be a $G\left(\left. \widetilde\sX~\right|\sX \right)$-set (cf. Definition \ref{g_category_defn}), and let $\left\{ G_\la\right\}_{\la \in \La}$ be a family of all finite factor groups of  $G\left(\left. \widetilde\sX~\right|\sX \right)$. According to the Definition \ref{g_category_defn} the set $\La$ is directed.
If  $\widehat{G} \bydef \varprojlim_{\la \in \La}  G_\la$ is an inverse limit of finite groups (cf. Definition \ref{group_inv_lim_defn}) then the group  $\widehat{G}$ is profinite (cf. Example \ref{profinite_exm}).  If $\sX_\la \bydef \widetilde \sX/ \ker\left( G\left(\left. \widetilde\sX~\right|\sX \right)\to G_\la\right)$ then there is an inverse limit $\widehat \sX = \varprojlim_{\la \in \La} \sX_\la$ (cf. Definition \ref{top_inverse_limit_defn}). There is a natural continuous map $\widetilde{\widehat p}: \widetilde \sX  \to \widehat \sX$. If we consider a  {final} {with respect to the family of maps} $\left\{g \circ \widehat p\right\}_{g\in \widehat{G}}$ topology on $\widehat \sX$ (cf. Definition \ref{top_final_defn}) then we obtain a topological space $\overline \sX$.

\begin{lemma}\label{top_disconnected_lem}
	Under the above hypotheses  the following conditions hold.
\begin{enumerate}
		\item[(i)] If $\left\{g_\iota G\left(\left. \widetilde\sX~\right|\sX \right)\right\}_{\iota \in I}$   is  a set of all left  cosets of $G\left(\left. \widetilde\sX~\right|\sX \right)$ in $\widehat{G}$ (cf. Definition \ref{group_coset_defn}) then there is a natural homeomorphism 
\be\label{top_disconnected_eqn}
	\overline \sX \cong \bigsqcup_{\iota\in I} g_\iota \widetilde\sX.
\ee
	
	\item[(ii)] The natural map $\widetilde{\widehat p}: \widetilde \sX  \to \widehat \sX$ yields a natural inclusion $\widetilde \sX \subset \overline \sX$ such that $\widetilde \sX$ is a quasi-component (cf. Definition \ref{top_quasi_component_defn}) of $\overline \sX$.
	\item[(iii)] For any  a quasi-component $\widetilde \sX' \subset \overline \sX$ there is $g \in \widehat{G}$ such that $\widetilde \sX' = g \widetilde \sX$.
	\item[(iv)]  For any $\la \in \La$ the natural surjective map $\widehat p_\la: \widehat \sX \to \sX_\la$ yields a covering  $\overline p_\la: \overline \sX \to \sX_\la$ such that $\sX_\la\cong \overline\sX/ \ker \left(\widehat{G} \to G_\la \right)$ (cf. Theorem \ref{top_group_of_covering_transformations_thm}). 
	\item[(v)] There is a natural bijective continuous map $\overline{\widehat p}: \overline \sX \to \widehat \sX$.
\end{enumerate}	
\end{lemma}

\begin{proof}(i)
	From the Definition \ref{top_inverse_limit_defn} it follows that there is the natural continuous map  $\widetilde{\widehat p}: \widetilde \sX  \to \widehat \sX$.
	For any $\widehat x \in \widehat \sX$ there is $\widetilde x \in \widetilde{\sX}$ such that $\widetilde p\left(\widetilde x \right) = \widehat p \left( \widehat x\right)$. So there is $g \in  \widehat{G}$ such that $\widehat x = g \widetilde{\widehat p}\left(\widetilde x \right) $. There is a left  coset $g_\iota G\left(\left. \widetilde\sX~\right|\sX \right)$ such that $g \in g_\iota G\left(\left. \widetilde\sX~\right|\sX \right)$. If $\widetilde x' \in \widetilde{\sX}$ such that $\widetilde p\left(\widetilde x' \right) = \widehat p \left( \widehat x\right)=\widetilde p\left(\widetilde x \right)$ then there is $g' \in G\left(\left. \widetilde\sX~\right|\sX \right)$ such that $\widetilde x = g'\widetilde x'$ and $\widehat x = g g' \widetilde{\widehat p}\left(\widetilde x' \right)$ where  $g g' \in g_\iota G\left(\left. \widetilde\sX~\right|\sX \right)$. So any $\widehat x \in \widehat \sX$ uniquely defines a left  coset. From this fact one has a set theoretic disjoint union
	$$
 \overline \sX = \bigsqcup_{\iota \in I}	g_\iota \widetilde{\widehat p}\left(  \widetilde \sX\right).
	$$
	where the topological disjoint union is implied.\\
(ii) Follows from the proof of (i).\\
(iii) Follows from the proof of (i).\\	
(iv) Follows from (i).\\
(v) For any $\le\in \La$ there is a covering $\overline \sX \to \sX_\la$ so from  $\widehat \sX \bydef \varprojlim_{\la \in \La} \sX_\la$ one has the unique continuous map $\overline{\widehat p }: \overline \sX \to \widehat{\sX}$. This map is bijective since  $\overline \sX$ coincides with $\widehat \sX$ as a set.
\end{proof}

\begin{definition}\label{top_disconnected_defn}
	Under the  hypotheses of the Lemma \ref{top_disconnected_lem} we say that the map $\overline p: \overline \sX  \to  \sX$  is the \textit{disconnected covering of} $p : \widetilde \sX \to \sX$. The {topological} $\overline \sX$-$\widehat  G$-{category} $\mathfrak{S}_p$ is the \textit{finite covering category of} $p : \widetilde \sX \to \sX$.
	Write
	\be\label{top_category_fin_eqn}
	\mathfrak{S}_p \bydef \left\{\left\{\sX_\la\right\}_{\la \in \La}, \left\{p^\mu_\nu:\sX_\mu\to \sX_\nu\right\}_{\substack{\mu,\nu \in \La\\\mu\ge\nu}}\right\}.
	\ee
	We say that $p : \widetilde \sX \to \sX$ is the \textit{covering inverse limit of} $\mathfrak{S}_p$
	and we write
	\be\label{top_disconnected_lim_eqn}
\widetilde \sX \bydef \varprojlim \mathfrak{S}_p
\ee
\end{definition}
\begin{remark}
	The category $\mathfrak{S}_p$ is a subcategory of  category $\mathfrak{FinCov}$-$\sX$ of finite coverings of $\sX$ (cf. Definition \ref{top_fin_cov_defn}).
\end{remark}
\begin{remark}
Indeed $\widetilde \sX$ is an  inverse limit if we consider a category having coverings as morphisms.
\end{remark}

\begin{theorem}\label{top_compact_thm} 
Let $\overline \sX$ be a Hausdorff, locally compact topological space with an action $\widehat  G\times \overline \sX\to \overline\sX$ (cf. Definition \ref{residually_finite_defn}) of properly discontinuous  group $\widehat  G$ of homeomorphisms (cf. Definition \ref{top_properly_disc_group_defn}). Assume that the group $\widehat  G$ is  residually finite (cf. Definition \ref{residually_finite_defn}). Suppose that $\sX \cong \overline \sX/\widehat G$ and   $p: \overline \sX \to \sX$ is  a natural transitive covering (cf. Theorem \ref{top_group_of_covering_transformations_thm}).  For any compact subset  $\overline \sY \subset \overline \sX$ there are both a finite-fold transitive covering $\widetilde p: \widetilde \sX\to \sX$ and a transitive covering $\overline p:\overline \sX \to \widetilde \sX$ such that: 
\begin{itemize}
	\item $p = \widetilde p \circ \overline p$,
	\item the set $\overline \sY$ is homeomorphically mapped onto $\widetilde \sY \bydef \overline p \left( \overline \sY\right)$.
\end{itemize}	
\end{theorem}

\begin{proof}
For any $\overline{x}\in \overline{  \sY}$ there is  a connected open neighborhood $\overline{   \mathcal U }_{\overline{x}}$ which is mapped homeomorphically onto ${p}\left(\overline{   \mathcal U }_{\overline{x}} \right)$ i.e. ${p}\left(\overline{   \mathcal U }_{\overline{x}} \right)$ is evenly covered by $p$ (cf. Definition \ref{top_covering_defn}). One has $\overline{   \mathcal Y } \subset \cup_{\overline{x} \in \overline{   \mathcal U }} \overline{   \mathcal U }_{\overline{x}}$. Since $\overline{   \mathcal Y }$ is compact there is a finite set $\left\{ \overline{   \mathcal U }_{\overline{x}_1}, ...,  \overline{   \mathcal U }_{\overline{x}_n}\right\} \subset  \left\{ \overline{   \mathcal U }_{\overline{x}}\right\}_{\overline{x }\in \overline{   \mathcal U }}$ such that $\overline{   \mathcal Y } \subset \cup_{j=1}^n \overline{   \mathcal U }_{\overline{x}_j}$. Write $\left\{ \overline{   \mathcal U }_{1}, ...,  \overline{   \mathcal U }_{n}\right\}\stackrel{\text{def}}{=}\left\{ \overline{   \mathcal U }_{\overline{x}_1}, ...,  \overline{   \mathcal U }_{\overline{x}_n}\right\}$. Let  $J$ be a maximal subset of  $\left\{ 1,..., n\right\}^2$ such that for all $\left(j,k\right) \in J$ one has
\be\label{top_compact_eqn}
\begin{split}
	j > k,\\
	\overline{   \mathcal U }_{j} \cap \overline{   \mathcal U }_{k} = \emptyset,\\
	{p}\left(\overline{   \mathcal U }_{j} \right)\cap {p}\left( \overline{   \mathcal U }_{k}\right) \neq \emptyset.	
\end{split}
\ee	
Then for every $\left(j,k\right) \in J$  there are $\overline{x}_j \in \overline{   \mathcal U }_{j}$ and $\overline{x}_k \in \overline{   \mathcal U }_{k}$ such that ${p}\left(\overline{x}_j \right) = {p}\left(\overline{x}_k \right)$. Since $\overline{x}_j \neq \overline{x}_k$ and the covering $p$ is transitive (cf. Definition \ref{top_transitive_defn}) there is a nontrivial $g_{jk} \in \overline G$ such that $\overline{x}_j = g_{jk} \overline{x}_k$, so   $\overline{   \mathcal U }_j \cap g_{jk} \overline{   \mathcal U }_k\neq \emptyset$. Since the group $\overline G$ is residually finite (cf. Definition \ref{residually_finite_defn}) there is a homomorphism $h_{jk}:\overline G\to G_{jk}$ such that:
\begin{itemize}
	\item $G_{jk}$ is a finite group,
	\item $h_{jk}\left(g_{jk} \right)\in G_{jk}$ is not trivial.
\end{itemize}
If $H_{jk}\bydef \ker h_{jk}$, and $\sX_{jk} \bydef \overline{\sX}/ H_{jk}$ then there are natural coverings $\overline p_{jk}: \overline \sX\to  \sX_{jk}$ and $\widetilde p_{jk}:  \sX_{jk}\to \sX$
such that
\be\label{top_compactp_eqn}
\begin{split}
	G\left(\left. \sX_{jk}\right| \sX  \right)\cong G_{jk},\\
	\widetilde p_{jk}\circ \overline p_{jk} = p,\\
	\overline p_{jk}\left( \overline{\sU}_j\right) \cap  h_{jk}\left(g_{jk} \right)	\overline p_{jk}\left( \overline{\sU}_k\right)= \overline p_{jk}\left(  \overline{\sU}_j \cap g_{jk}  \overline{\sU}_k\right)  \neq \emptyset
\end{split}
\ee
Using  the Lemmas \ref{top_cov_lem} and \ref{top_om_lem} and taking into account the Theorem \ref{top_bicont_thm} one can prove that natural maps 
\bean
\overline{\sU}_j \xrightarrow{\cong} 	\overline p_{jk}\left(\overline{\sU}_j \right),\\ 
\overline p_{jk}\left(\overline{\sU}_j \right)  \xrightarrow{\cong} 	\overline p\left(\overline{\sU}_j \right), \\
\overline{\sU}_k \xrightarrow{\cong} 	\overline p_{jk}\left(\overline{\sU}_k \right),\\ 
\overline p_{jk}\left(\overline{\sU}_k \right)  \xrightarrow{\cong} 	\overline p\left(\overline{\sU}_k \right), 
\eean
are homeomorphisms. Assume that  $\overline p_{jk}\left(\overline{\sU}_k \right)\cap  h_{jk}\left( g_{jk}\right) \overline p_{jk}\left( \overline{\sU}_k\right)\neq \emptyset$ and $\widetilde x \in \overline p_{jk}\left(\overline{\sU}_k \right)\cap  h_{jk}\left( g_{jk}\right) \overline p_{jk}\left( \overline{\sU}_k\right)$. Then one has
\bean
h_{jk}\left( g^{-1}_{jk}\right)\left(\widetilde x \right)\in \overline p_{jk}\left( \overline{\sU}_k\right),\\
h_{jk}\left( g^{-1}_{jk}\right)\widetilde x  \neq \widetilde x,\\	
\widetilde p_{jk}\left(  h_{jk}\left( g^{-1}_{jk}\right)\widetilde x \right) = 	\widetilde p_{jk} \left( \widetilde x \right),
\eean
i.e.  the map $\overline p_{jk}\left(\overline{\sU}_k \right)  \xrightarrow{\cong} 	\overline p\left(\overline{\sU}_k \right)$ is not bijective, so there is a contradiction. From this contradiction one concludes that the assumption $\overline p_{jk}\left(\overline{\sU}_k \right)\cap  h_{jk}\left( g_{jk}\right) \overline p_{jk}\left( \overline{\sU}_k\right)\neq \emptyset$ is not true, i.e. one has
\be\label{top_uiniu_eqn}
\overline p_{jk}\left(\overline{\sU}_k \right)\cap  h_{jk}\left( g_{jk}\right) \overline p_{jk}\left( \overline{\sU}_k\right)=\emptyset.
\ee
If we assume   $\overline p_{jk}\left( \overline{\sU}_j\right) \cap   	\overline p_{jk}\left( \overline{\sU}_k\right)\neq \emptyset$ then there is $\widetilde x \in \overline p_{jk}\left( \overline{\sU}_j\right) \cap   	\overline p_{jk}\left( \overline{\sU}_k\right)$	 and
\bean
\widetilde p_{jk}\left(\widetilde x \right) \in \widetilde p_{jk}\circ \overline p_{jk}\left( \overline{\sU}_j\right) \cap   \widetilde p_{jk}\circ 	\overline p_{jk}\left( \overline{\sU}_k\right)= \widetilde p_{jk}\circ \overline p_{jk}\left( \overline{\sU}_j\right) \cap   \widetilde p_{jk} \left( 	\overline p_{jk}\left( g_{jk} \overline{\sU}_k\right)\right) =\\
=  \widetilde p_{jk}\left( 	\overline p_{jk}\left(\overline{\sU}_j\right) \cap  h_{jk}\left( g_{jk}\right) p_{jk}\left( \overline{\sU}_k\right)\right) .
\eean
From the above equation it follows that there is $\widetilde x' \in \overline p_{jk}\left(\overline{\sU}_j\right) \cap  h_{jk}\left( g_{jk}\right) p_{jk}\left( \overline{\sU}_k\right)$ such that $\widetilde p_{jk}\left( \widetilde x'\right) = \widetilde p_{jk}\left( \widetilde x\right)$. Taking into account \eqref{top_uiniu_eqn} one has $\widetilde x' \neq \widetilde x$, however it is impossible because the map $\overline p_{jk}\left(\overline{\sU}_j \right)  \xrightarrow{\cong} 	\overline p\left(\overline{\sU}_j \right)$ is bijective, i.e. there is a contradiction. So the assumption   $\overline p_{jk}\left( \overline{\sU}_j\right) \cap   	\overline p_{jk}\left( \overline{\sU}_k\right)\neq \emptyset$ is not true and  
\be\label{top_compactppp_eqn}	\overline p_{jk}\left( \overline{\sU}_j\right) \cap   	\overline p_{jk}\left( \overline{\sU}_k\right)= \emptyset.\ee  The group	$H_{jk}$ is a normal subgroup of $\overline G$ for any pair $\left(j, k\right)\in J$,
The index   $\left( \overline G : H_{jk} \right) = \left|G_{jk} \right| $ is finite (cf. Section \ref{grops_sec}). The set $J = \left\{\left(j_1, k_1 \right),...,  \left(j_m, k_m \right) \right\} \subset \left\{1,..., n\right\}^2$ is finite. If
$$
H \bydef \bigcap_{\left(j, k \right)\in J} H_{j k}
$$
then index   $\left( \overline G: H\right)$ is finite (cf Example \ref{profinite_exm}). So if $\widetilde \sX \bydef \overline \sX / H$ then there  are transitive  coverings $\widetilde p: \widetilde \sX\to \sX$ and  $\overline p:\overline \sX \to \widetilde \sX$  such that $p = \widetilde p \circ \overline p$. From $\left( \overline G: H\right)= \left|G\left(\left.\widetilde\sX\right| \sX  \right)\right| $ it follows that $G\left(\left.\widetilde\sX\right| \sX  \right)$ is a finite group and  $\widetilde p: \widetilde \sX\to \sX$ is a finite-fold covering. Using \eqref{top_compactppp_eqn} one has
\be\label{top_compactpp_eqn}
\overline{\sU}_j \cap  \overline{\sU}_k= \emptyset \quad \Leftrightarrow \quad 	\overline p\left( \overline{\sU}_j\right) \cap \overline p\left( \overline{\sU}_k\right)= \emptyset.
\ee
Let $\overline y', \overline y'' \in \overline \sY$ be such that $\overline p\left(\overline y'\right)= \overline p\left(\overline y''\right)= \widetilde y$.
If $\overline y' \in \overline \sU_j$ and $\overline y'' \in \overline \sU_k$  then from \eqref{top_compactpp_eqn} it follows that $\overline y',\overline y''\in \overline \sU_j\cap \overline \sU_k$, so  $\overline y'= \overline y''$, i.e. $\overline p$ yields a bijective map  
$\th: \overline \sY\xrightarrow{\cong}\widetilde \sY \bydef \overline p \left( \overline \sY\right)$. From the Lemmas \ref{top_cov_lem} and \ref{top_om_lem} it follows that the map $\th$ is open (cf. Definition \ref{top_bicont_defn}). Taking into account the Theorem \ref{top_bicont_thm} one concludes that $\th$ is a homeomorphism.
\end{proof}

\begin{lemma}\label{top_compact_la_lem}
Let $\overline p:\overline\sX \to \sX$ be a transitive covering, and let $\left(\overline \sU, \overline \sV, \overline s\right)$ is a {covering triple for} $\overline p$ (cf. Definition \ref{top_coveing_triple_defn})
Under the hypotheses of the Theorem \ref{top_compact_thm} there is $\la_{\left(\widetilde \sU, \widetilde \sV, \widetilde s\right)} \in \La$ such that for any  $\la \ge \la_{\left(\overline \sU, \overline \sV, \overline s\right)} \in \La$ such that  the closure of  $\overline \sV$ is mapped homeomorphically onto  its image in $\sX_\la$.
\end{lemma}
\begin{proof}
The proof of this lemma can be obtained from the Theorem \ref{top_compact_thm}.
\end{proof}

\begin{empt}\label{top_fin_empt}
	Let $G$ be a group, and let $\left\{h: G \to G_\la\right\}_{\la\in \La}$ is a set of all homomorphisms from $G$ to finite groups. If 
	\be\label{top_inf_eqn}
	\mathfrak{ResInf}\left( G\right)\bydef \bigcap_{\la \in \La} \ker h_\la
	\ee
	then $\mathfrak{ResInf}\left( G\right)\subset G$ be a normal subgroup $G$. We define 
	\be\label{top_fin_eqn}
	\mathfrak{ResFin}\left( G\right)\bydef G / 
	\mathfrak{ResInf}\left( G\right).
	\ee
	A given by the equation group is residually finite (cf. Definition \ref{residually_finite_defn}).
\end{empt}
\begin{empt}\label{top_fin_sec_empt}
	Let $ \sX $ be a connected,  locally compact, Hausdorff  pointed space, and let $\widetilde{p}: \widetilde{\sX},  \to \sX$ be a     transitive covering such that $\widetilde{\sX}$ is connected (cf. Definition \ref{top_pointed_defn}). Consider a family of all surjective homomorphisms $\left\{h_\la:G\left( \left.{\widetilde{\sX}}~\right|\sX\right)\to G_\la\right\}$ onto finite groups. If $H_\la \bydef \ker h_\la$ and 
\be\label{top_xh_eqn}
	\sX_\la \bydef \widetilde \sX / H_\la
\ee 
then there is a  finite-fold transitive covering $p_\la :  \sX_\la \to \sX$ such that
	\begin{itemize}
	\item there is a natural sequence of transitive coverings
	\be\label{top_fin_sec_eqn}
	\widetilde{\sX}\xrightarrow{\widetilde p_\la}  {\sX}_\la \xrightarrow{ p_\la} \sX;
	\ee
	\item one has a natural isomorphism of groups $G\left( \left.{{\sX}_\la}~\right|\sX\right)\cong G_\la$,
	\item The set $\La$ is ordered by following way
$\mu \ge \nu$ if there is the following commutative diagram
	\newline
\begin{tikzcd}
	& \widetilde \sX \arrow[ld, "\widetilde p_\mu"] \arrow[rd, "\widetilde p_\nu"] &\\
\sX_\mu \arrow[rr, "p^\mu_\nu"] \arrow[rd, "p_\mu"]	& &\sX_\nu \arrow[ld, "p_\nu"] \\
	& \sX&
\end{tikzcd}
\newline
where all arrows are transitive coverings.
\end{itemize}

Conversely if a finite-fold transitive covering satisfies to the equation \eqref{top_fin_sec_eqn} then there is a natural surjective homomorphism $G\left( \left.{\widetilde{\sX}}~\right|\sX\right)\to G\left( \left.{{\sX}_\la}~\right|\sX\right)$. It follows that there is an 1-1 correspondence between surjective homomorphism and coverings which satisfy to the equation \eqref{top_fin_sec_eqn}. One can proof that
$$
\bigcap_{\la\in\La} H_\la = \mathfrak{ResInf}\left( G\left( \left.{\widetilde{\sX}}~\right|\sX\right)\right)
$$
(cf. \eqref{top_inf_eqn}). If 
\be\label{top_x_fr_eqn}
\widetilde \sX_{\mathrm{res~fin}}\bydef \widetilde \sX / \mathfrak{ResInf}\left( G\left( \left.{\widetilde{\sX}}~\right|\sX\right)\right)
\ee 
then one has:
\begin{itemize}
	\item for any $\la\in \La$ there is a sequence of coverings
	\be\label{top_fin_c_sec_eqn}
 \widetilde{\sX}_{\mathrm{res~fin}}\to {\sX}_\la\to  {\sX};
	\ee
	\item if $\mathfrak{ResFin}$ is given by \eqref{top_fin_eqn} then there is a natural isomorphism.
	\be\label{top_fin_g_sec_eqn}
	G\left( \left.{\widetilde{\sX}_{\mathrm{res~fin}}}~\right|\sX\right)\cong \mathfrak{ResFin}\left( G\left( \left.{\widetilde{\sX}}~\right|\sX\right)\right).
	\ee 	
\end{itemize}

\end{empt}

\begin{lemma}\label{top_associded_lem}
	Under the hypotheses \ref{top_fin_sec_empt} the natural covering $\widetilde{\sX}_{\mathrm{res~fin}}\to   {\sX}$ is the {covering inverse limit of} 
of the category
	\bean
\mathfrak{S}_p \bydef \left\{\left\{\sX_\la\right\}_{\la \in \La}, \left\{p^\mu_\nu:\sX_\mu\to \sX_\nu\right\}_{\substack{\mu,\nu \in \La\\\mu\ge\nu}}\right\}.
\eean
(cf. Definition \ref{top_disconnected_defn}).
\end{lemma}
\begin{proof}
Follows from the construction \ref{top_fin_sec_empt}.
\end{proof}

\subsection{Infinite coverings and algebraic constructions}
\paragraph{}
Here we explain  general algebraic constructions which can be applied to the coverings of commutative $C^*$-algebras and another ones discussed in Chapters \ref{blowing_chap}, \ref{ctr_chap} and \ref{foliations_chap}.
Let  $\sX$ be a connected, locally coconnected,  locally compact, Hausdorff space, and let $p: \widetilde \sX\to \sX$ a  be a covering.  Suppose that  $\widetilde{G}$ is a residually finite group (cf. Definition \ref{residually_finite_defn})  of properly discontinuous homeomorphisms of $\widetilde \sX$  (cf. Definition \ref{top_properly_disc_group_defn}) such that there is the natural homeomorphism $\sX = \widetilde\sX /  \widetilde G$. Suppose   $\widetilde{G}$ is residually finite and suppose that $\left\{G_\la\right\}_{\la \in \La}$ is a indexed by a directed set family of finite factor-groups of  $\widetilde G$ such that for all $\mu, \nu \in \La$ one has $\mu \ge \nu$ if and only if there is a following natural commutative diagram
	\newline
	\begin{tikzcd}
		\widetilde G \arrow[rr]\arrow[rd] & & G_\nu\\
		& G_\mu\arrow[ru] &
	\end{tikzcd}
	\\
	of groups and surjective homomorphisms.	For any $\la\in \La$ denote by $\sX_\la \bydef \widetilde \sX / \ker \left(\widetilde G\to G_\la \right)$. For any $\la\in \La$ there are natural coverings   $p_\la: \sX_\la \to \sX$ and $\widetilde p_\la: \widetilde \sX \to \sX$ such that $p_\la$ is finite-fold. For any $\la\in \La$ there are  the $p$-lift  $C_0\left( p_\la\right): C_0\left( \sX\right) \hookto C_0\left(\sX_\la \right)$ (cf. Definition \ref{top_lift_a_f_defn}) and the given by the Definition \ref{top_cs_functa_b_defn} injective $*$-homomorphism $C_b\left(\widetilde p_\la\right):  C_0\left(\sX_\la \right)\hookto C_b\left(\widetilde \sX \right)$. There is an  {algebraic} $\overline  \sX$-$  \widehat G$-{category} 
		\be\label{top_category_fin_a0_eqn}
	\begin{split}
		\mathfrak{S}_{C_0\left(p\right) } \bydef \\
		\left\{ \left\{ C_0\left( p_\la\right)  :C_0\left( \mathcal{X}\right)  \hookto C_0\left( \mathcal{X}_\la\right) \right\}, \left\{C_0\left( p^\mu_\nu\right)  :C_0\left( \mathcal{X}_\mu\right)  \hookto C_0\left( \mathcal{X}_\nu\right) \right\}  \right\}
	\end{split}
	\ee
	(cf. Definition \ref{top_fin_cov_cat_defn})
	
\begin{definition}\label{top_discr_defn}
	If $\sX$ is a topological space, and  $\sX_{\mathrm{discr}}$ equals to $\sX$ as a set but the topology of $\sX_{\mathrm{discr}}$ is discrete
	then the $C^*$-algebra $C\left(\sX_{\mathrm{discr}} \right)$ is the \textit{discontinuous extension of} $C_0\left(  \sX\right)$.
\end{definition}
	If  $C_0\left(\widetilde\sX \right)$ be given by the Definition \ref{top_cs_functa_b_defn} then from the Theorem \ref{dauns_hofmann_thm} it follows that there is an action
	$$
	C_b\left(\widetilde \sX \right) \times C_0\left(\widetilde\sX \right) \to C_0\left(\widetilde\sX \right).
	$$
	This action can be naturally extended up to the action
	$$
	C_b\left(\widetilde \sX \right) \times C\left(\widetilde\sX_{\mathrm{discr}} \right) \to C\left(\widetilde\sX_{\mathrm{discr}} \right)
	$$
	where $C\left(\widetilde\sX_{\mathrm{discr}}\right) $ is the \textit{discontinuous extension of} $C_0\left(\widetilde  \sX\right)$.	For any $\la\in \La$ there is the natural inclusion $C_b\left( \sX_\la\right) \subset C_b\left(\widetilde \sX \right)$ so there is the natural action
	\be\label{ctr_inf_a_eqn}
	C_b\left( \sX_\la \right) \times C\left(\widetilde\sX_{\mathrm{discr}} \right) \to C\left(\widetilde\sX_{\mathrm{discr}} \right).
	\ee
	For any $\la\in\La$ there is  the given by \eqref{top_cb_defna_eqn} inclusion 	$C_b\left(\widetilde p_\la \right): C_0\left(\sX_\la \right)  \hookto C_b\left(\widetilde{\sX} \right)$ it follows that for any $\widetilde{x}\in \widetilde \sX$ there is a $*$-homomorphism 
	\be\label{ctr_inf_rho_defn}
	\rho^\la_{\widetilde x} : C_0\left(\sX_\la \right) \to \rep_{  \widetilde{x }}\left(  C_0\left(\widetilde\sX \right)\right)\cong \C
	\ee 
	where $\H_{ \widetilde x}\cong \C$ is the Hilbert space of the representation $ \rep_{  \widetilde{x }}$.

\begin{empt}\label{top_long_p_empt}
Suppose that  $\widetilde a\in  K\left( C\left(\widetilde \sX_{\mathrm{discr}} \right) \right) _+$ is  a positive element such that
\be\label{top_atomic_cc_eqn}
a_\la \bydef s\text{-}	\sum_{g \in \ker\left(\widehat{G}\to  G\left( \left.\sX_\la~\right|\sX\right)\right) }g \widetilde a \in  C_c\left(\sX_\la\right)
\ee 
where the point-wise convergence of the series is implied and $C_c\left(\sX_\la\right)$ is regarded as a $C^*$-subalgebra of $C\left(\widetilde\sX_{\mathrm{discr}} \right)$.
Denote by $\widetilde b^1 \bydef \widetilde a$. Let $\widetilde x_0\in \widetilde \sX$ be such that $\frac{2}{3}\left\|\widetilde a\right\|<  \widetilde a \left( \widetilde x_0\right) \le \left\| \widetilde a\right\|$. From  the equation \eqref{top_atomic_cc_eqn} it follows that there is $\la_0\in \La$ such that $\frac{2}{3}\left\| \widetilde b^1\right\|<b^1_{\la_0}\left(\widetilde{p}_{\la_0}\left( \widetilde x_0\right)  \right) < \frac{4}{3}\left\| \widetilde b^1\right\|$ where $b^1_{\la_0}=s\text{-}\sum_{g \in \ker\left(\widehat{G}\to G_{\la_0}\right)}g \widetilde b^1\in   C_c\left(\sX_{\la_0}\right)$ and the point-wise convergence of the series is implied. There is an open connected neighborhood $\sU'_{\la_0}\subset \sX_{\la_0}$ of $\widetilde{p}_{\la_0}\left( \widetilde x_0\right)$ such that $\frac{2}{3}\left\|  \widetilde b^1\right\|<b^1_{\la_0} \left(  x_{\la_0}\right)  < \frac{4}{3}\left\|  \widetilde b^1\right\|$ for all  $ x_{\la_0}\in  \sU'_{\la_0}$. Let $\widetilde{a}^1 \bydef f_\eps \left(\widetilde{b}_1\right)$ where $f_\eps$ is given by the equation \eqref{f_eps_eqn} and $\eps =  \frac{2}{3}\left\|  \widetilde b^1\right\|$. If $a^1_{\la_0}= \sum_{g \in \ker\left(\widehat{G}\to G_{\la_0}\right)}g \widetilde a^1\in K\left(C_0\left( \sX_{\la_0}\right)\right)$ then from $\frac{2}{3}\left\| \widetilde b^1\right\| >\widetilde a^1 \left( \widetilde x_0\right)$ it follows that
$0<  \widetilde a^1 \left( \widetilde x_0\right)< \frac{1}{3}\left\|  \widetilde b^1\right\|$ and 
$0< a^1_{\la_0}\left(\widetilde{p}_{\la_0}\left(\widetilde x_0 \right)  \right)< b^1_{\la_0}\left(\widetilde{p}_{\la_0}\left(\widetilde x_0 \right) \right) - \frac{2}{3}\left\| \widetilde b^1\right\|< \frac{2}{3}\left\| \widetilde b^1\right\|$. So   there is an open  neighborhood $\sU''_{\la_0}\subset \sX_{\la_0}$ of $\widetilde{p}_{\la_0}\left( \widetilde x_0\right)$ such that  $0<a^1_{\la_0} \left(  x_{\la_0}\right)  < \frac{2}{3}\left\| \widetilde b^1\right\|$ 
for all $x_{\la_0}\in \sU''_{\la_0}$. If $\mu \in \La$ is such that $\mu > \la_0$ and  $a^1_{\mu}=\sum_{g \in \ker\left(\widehat{G}\to G_{\mu}\right)}g \widetilde a^1$ then one has
$$
a^1_{\la_0} = \sum_{	g \in G\left( \left. \sX_\mu\right| \sX_{\la_0}\right) } a^1_{\mu}
$$
From the above equation it follows that
\bean
p^\mu_{ \la_0}\left(\left\{x_\mu \in \sX_{\mu} |  a^1_{\mu}\left(x_\mu \right)  > 0\right\} \right) = \left\{x_{ \la_0} \in \sX_{ \la_0} |  a^1_{ \la_0}\left( x_{\la_0}\right)  > 0\right\}, 
\eean
so if $\sU_{\la_0} \bydef \sU'_{\la_0}\cap  \sU''_{\la_0}\cap \left\{x_{ \la_0} \in \sX_{ \la_0} |  a^1_{ \la_0}\left( x_{\la_0}\right) > 0\right\}$ and $$\sU_\mu \bydef \left\{\left.x_\mu \in \sX_{\mu} \right|  a^1_{\mu}\left( x_\mu\right)  > 0\right\}\cap \left( p^\mu_{ \la_0}\right)^{-1}\left(\sU_{\la_0}\right)$$
then $p^\mu_{ \la_0}\left(\sU_\mu  \right)= \sU_{\la_0}$, i.e. one has a surjective continuous map
$$
\th\bydef \left.p^\mu_{ \la_0} \right|_{\sU_\mu}:\sU_\mu\to \sU_{\la_0}
$$
If $\th$ is not a one-to-one map then there are $x^1_\mu, x^2_\mu\in \sU_\mu$ such that $x^1_\mu\neq x^2_\mu$ and $ p^\mu_{ \la_0}\left(x^1_\mu\right)=  p^\mu_{ \la_0}\left(x^2_\mu\right)$. From $a^1_{ \la_0}\left(x^1_\mu \right) , a^1_{ \la_0}\left(x^2_\mu \right) > 0$ it follows that there are $\widetilde{x}_1, \widetilde{x}_2 \in \widetilde{\sX}$ such that $\widetilde a^1 \left(\widetilde{x}_1\right),  \widetilde a^1 \left(\widetilde{x}_2\right) > 0$, so one has $ \widetilde b^1 \left(\widetilde{x}_1\right), \widetilde b^1 \left(\widetilde{x}_2\right) > \frac{2}{3}\left\| \widetilde b^1\right\|$.
It follows that $b^1_\mu\left(x^1_\mu \right), b^1_\mu\left(x^2_\mu \right)  > \frac{2}{3}\left\| b^1\right\|$ so $b^1_{\la_0}\left( p^\mu_{ \la_0}\left(x^1_\mu\right) \right)  > b^1_\mu\left(x^1_\mu \right)+ b^1_\mu\left(x^2_\mu \right) > \frac{4}{3}\left\| \widetilde b^1\right\|$. It contradicts with $\frac{2}{3}\left\| \widetilde b^1\right\|<b^1_{\la_0}\left( p^\mu_{ \la_0}\left(x^1_\mu\right) \right) < \frac{4}{3}\left\| \widetilde b^1\right\|$. From this contradiction turns out that a map $\th$ is one to one. From the Lemmas \ref{top_cov_lem} and \ref{top_om_lem} it follows that the map $\th$ is open (cf. Definition \ref{top_bicont_defn}). Taking into account the Theorem \ref{top_bicont_thm} one concludes that $\th$ is a homeomorphism. From this fact it follows that
\be\label{top_a_desc}
\forall x_\mu \in \sU_\mu\quad a^1_{\mu }\left( x_\mu\right) = a^1_{\la_0 }\left( p^\mu_{ \la_0\left(x_\mu \right) }\right). 
\ee
If $\widetilde{\sU}\subset \widetilde \sX$ is a connected component of $\widetilde{p}^{-1}_{\la_0}\left(\sU_{{\la_0}} \right)$ such that $\widetilde x_0 \in \widetilde\sU$ then from the equation \eqref{top_a_desc} it follows that
\be\label{top_lift_al_eqn}
\begin{split}
	\forall \widetilde x \in \widetilde \sU\quad  \widetilde a^1\left( \widetilde x\right) = a^1_{\la_0 }\left(\widetilde p_{ \la_0\left(\widetilde x\right) }\right),\\
	\left.\widetilde a^1\right|_{\widetilde\sU} = \lift^{ \widetilde{p}_{\la_0}}_{\widetilde{\sU}}\left( a^1_{\la_0 }\right) 
\end{split}
\ee
where $\lift^{ \widetilde{p}_{\la_0}}_{\widetilde{\sU}}$ is given by the Definition \ref{top_lift_desc_defn} $\widetilde{p}_{\la_0}$-${\widetilde{\sU}}$-lift.
So the restriction $\left. \widetilde a^1_{\la_0 }\right|_{\widetilde \sU}$ is continuous, i.e.  $\widetilde a^1$ is continuous at $\widetilde{x}_0$. Since a point $\widetilde x_0$ is arbitrary one has  $\widetilde a^1\in C_b\left( \widetilde{\sX}\right) $. One has $\left\|\widetilde{a}^1 \right\| =   \frac{1}{3}\left\|\widetilde{a}\right\|$.
If $\widetilde{b}^2 \bydef \widetilde{b}^1 - \widetilde{a}^1$ then one can proof following:
\bean
\widetilde{b}^2 \ge 0,\\
\left\|\widetilde{b}^2 \right\| \le   \frac{2}{3}\left\|\widetilde{a}\right\|. 
\eean 
So if $\eps' =\left( \frac{2}{3}\right)^2\left\|\widetilde{a}\right\| $ and $\widetilde{a}^2= f_{\eps'}\left( \widetilde{b}_2\right)$ then similarly to the above constriction one has $\widetilde a^2 \in C_b\left( \widetilde{\sX}\right)$ and
$\left\|\widetilde{a}^2 \right\| =    \frac{2}{3}\frac{1}{3}\left\|\widetilde{a}\right\|$. Going on one has a sequence $\widetilde{a}^1,..., \widetilde{a}^n, ... $ such that $ \widetilde a^n \in C_b\left( \widetilde{\sX}\right)$ and $\left\|\widetilde{a}^n \right\| =    \left( \frac{2}{3}\right)^{n-1} \frac{1}{3}\left\|\widetilde{a}\right\|$.
From  $\left\|\widetilde{a}^n \right\| =    \left( \frac{2}{3}\right)^{n-1} \frac{1}{3}\left\|\widetilde{a}\right\|$ if follows than one has $C^*$-norm convergent series
\be\label{top_ser_eqn}
\widetilde{a}' \bydef	\sum_{n = 1}^\infty \widetilde  a^n.
\ee
the sum of the series lies in $C_b\left( \widetilde{\sX}\right)$
On the other hand from
$$ 
\forall n \in \N \quad \left\|\widetilde{a}- \widetilde{a}^n\right\|\le \frac{1}{3}\left\|\widetilde{a}\right\|  \left( \frac{2}{3}\right)^{n-1}
$$	
it follows that  $\widetilde a= \widetilde a'$ and o $\widetilde a\in C_b\left( \widetilde{\sX}\right)\subset   C\left(\widetilde\sX_{\mathrm{discr}} \right)$.

\end{empt}

\begin{thm}\label{top_discrete_p_thm}
	Under the hypotheses \ref{top_long_p_empt} one has $\widetilde{a}\in C_c\left( \widetilde{\sX}\right)$.
\end{thm}
\begin{proof}
	From the Lemma \ref{pedersen_eps_lem} and the equality \eqref{four_decompositon_eqn} we can suppose that there is $\dl > 0$ such $\widetilde a = f_\dl \left( \widetilde b\right)$ where $\widetilde b \in C_b\left(\widetilde \sX \right)_+$ and $f_\dl$ is given by   the equation \eqref{f_eps_eqn}. If $\widetilde a' \bydef f_{\dl/2}\left( \widetilde b\right)$ then  $\widetilde a = f_{\dl/2}\left( \widetilde a'\right)$. 
	Let
	$$
	a \bydef \sum_{g\in \widehat G} g \widetilde a'\in C_c\left( \sX\right) 
	$$
	Let $\widetilde a^1 \bydef f_\eps\left( \widetilde a\right)$ where $f_\eps$ is given by \eqref{f_eps_eqn} and $\eps \bydef \frac{2}{3}\lVert \widetilde a\rVert$ and denote by $a^1 \bydef \sum_{g\in \widehat G} \widetilde a^1\in C_c\left( \sX\right)$. For any $\widetilde x \in \supp \widetilde a^1$ we select an open neighborhood $\widetilde \sU_{\widetilde x}\subset \widetilde \sX$ such that:
	\begin{itemize}
		\item $\widetilde \sU_{\widetilde x}$ is mapped homeomorphically onto $ \sU_{\widetilde x}\bydef \widetilde p\left(\widetilde \sU_{\widetilde x} \right)$,
		\item there is $\la_{\widetilde x} \in \La$ such that  
		\be\label{top_cond_eqn}
		\widetilde x \in \widetilde \sU_{\widetilde x}\quad \widetilde a^1\left(\widetilde x \right) = a_{\la_{\widetilde x}} \left( \widetilde p_{\la_{\widetilde x}}\left(\widetilde x \right)  \right) 	
		\ee 
		where 
		$$
		a_{\la_{\widetilde x}} \bydef  \sum_{g\in \ker \left( \widehat G \to G\left( \sX_{\la{\widetilde x}}| \sX\right)\right)  } \widetilde a^1.	
		$$
	\end{itemize}
	The existence of $\widetilde \sU_{\widetilde x}$ follows from the equation \eqref{top_lift_al_eqn}. From  $a^1\le a$ it follows that $\supp a^1$ is compact. From this fact and taking into account that $\supp a^1\subset \bigcup_{\widetilde x \in \supp \widetilde a^1}\sU_{\widetilde x}$ we conclude that there is a finite set $\left\{\widetilde x_1, ..., \widetilde x_n \right\}\subset \supp \widetilde a^1$ such that
	$$
	\supp a^1 \subset  \sU_{\widetilde x_1}\cup ...\cup \sU_{\widetilde x_n}.
	$$
	The set $\La$ is directed, so there is $\widetilde \la \in \La$ such that for any $\widetilde x_j \in \left\{\widetilde x_1, ..., \widetilde x_n \right\}$ one has $\widetilde \la\ge \la_{\widetilde x_j}$. If $x_{\widetilde \la} \in \widetilde p_{\widetilde \la}\left(\supp \widetilde a^1 \right)$ then $p_{\widetilde \la  }\left( x_{\widetilde \la} \right)\in \supp a^1$. For any  $x_{\widetilde \la} \in \widetilde p_{\widetilde \la}\left(\supp \widetilde a^1 \right)$ there is $\widetilde x_j \in \left\{\widetilde x_1, ..., \widetilde x_n \right\}$ such that $p_{\widetilde \la  }\left( x_{\widetilde \la} \right)\in \sU_{\widetilde x_j}$. If $x_{\la_{\widetilde x_j}}\bydef p^{\widetilde \la}_{\la_{\widetilde x_j}}\left(  x_{\widetilde \la}\right)$ then from \eqref{top_cond_eqn} it follows that there is the unique $\widetilde x \in \supp \widetilde a^1$ such that $\widetilde p_{\la_{\widetilde x_j}}\left(\widetilde x  \right) = x_{\la_{\widetilde x_j}}$. Thus there is   the unique point  $\widetilde x\in\supp \widetilde a^1$ such that $\widetilde p_{\widetilde \la}\left(\widetilde x  \right) = x_{\widetilde\la}$, so the restriction $\widetilde p_{\widetilde \la}|_{\supp \widetilde a^1}: \supp \widetilde a^1 \to \widetilde p_{\widetilde \la}\left(\supp \widetilde a^1 \right)$ is a bijective map.
	From the Lemmas \ref{top_cov_lem} and \ref{top_om_lem} and  the Theorem \ref{top_bicont_thm} it follows that the restriction is a homeomorphism. On the other hand 
	\bean
	p_{\widetilde \la}\left(\supp \widetilde a^1 \right)= \supp a^1_{\widetilde \la},\\
	\text{where} \quad 	
	a_{\widetilde\la} \bydef  \sum_{g\in \ker \left( \widehat G \to G\left( \sX_{\widetilde \la}| \sX\right)\right)  } \widetilde a^1.	
	\eean 
	From \eqref{a_la_eqn} it follows that $ a^1_{\widetilde \la} \in K\left(C_0\left(\sX_{\widetilde \la} \right)  \right) = C_c\left(\sX_{\widetilde \la} \right)$ 
	so the set $\supp a^1_{\widetilde \la}$ is compact. It turns out that the set $\supp \widetilde a^1$ is compact 
	and $\widetilde{a}^1\in C_c\left(\widetilde{\sX} \right)$. If $\widetilde{a}^2, ..., \widetilde{a}^n, ...$ are explained in \ref{top_long_p_empt} elements then similarly to the above proof one has $\widetilde{a}^n \in  C_c\left(\widetilde{\sX} \right)$. On the other hand one has a $C^*$-norm convergent series
	\bean
	\widetilde a'= \sum_{n = 1}^\infty  \widetilde a^n
	\eean
	(cf. Equation \eqref{top_ser_eqn}) so from the Definition \ref{c_c_closure_defn} it follows that $ \widetilde a \in C_0\left(  \widetilde \sX\right)$. From $\widetilde a= f_{\dl/2}\left( \widetilde a'\right)$ it follows  that $\widetilde a \in C_c\left(  \widetilde \sX\right)$.
\end{proof}

\begin{empt}\label{top_cov_empt}
  Let $p: \widetilde{\sX}\to \sX$ be a transitive covering such that 
   $\sX \cong \widetilde{\sX}/G$ where $G$ is  properly discontinuous group of homeomorphisms (cf. Definition \ref{top_properly_disc_group_defn}). Suppose that the space $\sX$ is locally connected (cf. Definition \ref{top_locally_connected_defn}). Let
	$$
	\left\{\left(\sU_\a, \sV_\a, s_\a\right)\right\}_{\a \in \mathscr A}
	$$
	be a family of triples such that 
	\begin{itemize}
		\item for any $\a \in \mathscr A$ there is  a covering  triple  $\left(\widetilde \sU, \widetilde \sV, \widetilde s\right)$ for $p$ (cf. Definition \ref{top_coveing_triple_defn}) with
		\bean
		\sU \bydef p\left( \widetilde \sU \right),\\
		\sV \bydef p\left( \widetilde \sV \right),\\
		 s_\a = \desc_p\left(\widetilde s \right), 
		\eean 
		where $\desc_p$ is the $p$-descent (cf. Definition \ref{top_lift_desc_defn}),
		\item  $\sX = \bigcup_{\a \in \mathscr A} \sU_{{\a}}$.
	\end{itemize}
	The existence of the family follows from the Remark \ref{top_coveing_triple_defn}.
	 For any $\a\in\mathscr A$ we select  a connected open set $\widetilde \sV_\a\subset \widetilde \sX$ which is mapped homeomorphically onto $\sV_\a$.  Let $\left(\widetilde \sU_\a, \widetilde \sV_\a, \widetilde s_\a\right)$ be a triple  such that  $\widetilde \sU_\a \bydef \widetilde\sV_\a \cap p^{-1}\left(\sU_\a \right)$ and $\widetilde s_\a \bydef \lift^p_{\widetilde{\sV}_\a}\left( s_\a\right)$ where $\lift^p_{\widetilde{\sV}_\a}$ is $p$-$\widetilde{\sV}_\a$-lift (cf. Definition \ref{top_desc_eqn}). If $\widetilde{\mathscr A}\bydef G\times  \mathscr A$ then for any $\widetilde \a = \left(g, \a \right)\in \widetilde{\mathscr A}$. If  we denote by
	 $$
\left(\widetilde \sU_{\widetilde\a}, \widetilde \sV_{\widetilde\a}, \widetilde s_{\widetilde\a}\right)\bydef \left(g\widetilde \sU_{\a},g \widetilde \sV_\a, g\widetilde s_\a\right). 
	$$ 
There are following objects:
	\begin{itemize}
		\item  the natural action $G\times \widetilde{\mathscr A}\to \widetilde{\mathscr A}$,
		\item the natural projection $p_{\widetilde{\mathscr A}}: \widetilde{\mathscr A} \to {\mathscr A}$.
	\end{itemize}
Moreover one has $\widetilde \sX = \bigcup_{\widetilde \a \in \widetilde{\mathscr A}} \widetilde \sU_{\widetilde \a}$.
	
\end{empt}

\begin{definition}\label{top_cov_defn}
	Under the hypotheses \ref{top_cov_empt} we say that a pair $$\left(\	\left\{\left(\sU_\a, \sV_\a, s_\a\right)\right\}_{\a \in \mathscr A}, \left\{\left(\widetilde \sU_{\widetilde\a}, \widetilde \sV_{\widetilde\a}, \widetilde s_{\widetilde\a}\right)\right\}_{\widetilde \a \in \widetilde{\mathscr A}} \right)$$ is a $p$-\textit{covering}. Both the action $G\times \widetilde{\mathscr A}\to \widetilde{\mathscr A}$ and the map  $p_{\widetilde{\mathscr A}}: \widetilde{\mathscr A} \to {\mathscr A}$ are $p$-\textit{action} and $p$-\textit{projection} respectively.
\end{definition}

\subsection{Coverings $C^*$-algebras and operator spaces}

\paragraph*{}
This section supplies a purely algebraic  analog of the topological construction given by the Subsection \ref{top_inf_to_sec}.
\begin{empt}\label{top_inf_sufficienf_empt}
Let $p:\widetilde  \sX \to \sX$ be a transitive covering such that:
\begin{itemize}
	\item $\sX$ is locally connected, locally compact, Hausdorff space,
	\item the space $\widetilde{\sX}$ is connected,
	\item the  covering group $G\left(\left.\widetilde{\sX}~\right|\sX\right)$ (cf. Definition  \ref{top_group_of_covering_transformations_defn}) is residually finite (cf. Definition  \ref{residually_finite_defn}).
\end{itemize}
Let $\left\{G_\la\right\}_{\la\in \La}$ be the indexed by the directed set $\La$  a family  of all finite factor-groups of  $G\left(\left.\widetilde{\sX}~\right|\sX\right)$ (cf. Definition \ref{g_category_defn}). 
 Let $\overline p: \overline \sX  \to  \sX$  be the {disconnected covering of} $p : \widetilde \sX \to \sX$, and let
 	\bean
 \mathfrak{S}_p \bydef \left\{\left\{\sX_\la\right\}_{\la \in \La}, \left\{p^\mu_\nu:\sX_\mu\to \sX_\nu\right\}_{\substack{\mu,\nu \in \La\\\mu\ge\nu}}\right\}.
 \eean 
	{finite covering category of} $p : \widetilde \sX \to \sX$ (cf. Definition \ref{top_disconnected_defn}).

\be\label{top_x_g_eqn}		\begin{split}
		\mathfrak{S}_{C_0\left(p\right) } \bydef \\
		\left\{ \left\{ C_0\left( p_\la\right)  :C_0\left( \mathcal{X}\right)  \hookto C_0\left( \mathcal{X}_\la\right) \right\}, \left\{ C_0\left( p^\mu_\nu\right)  :C_0\left( \mathcal{X}_\mu\right)  \hookto C_0\left( \mathcal{X}_\nu\right) \right\}  \right\}
	\end{split}
	\ee 
be the 	{algebraic} $\widetilde \sX$-$  G$-{category} (cf. Definition \ref{top_fin_cov_cat_defn}).
 If $\widehat G \bydef \varprojlim G_\la$ is the profinite completion (cf. Example \ref{profinite_exm}) then
 $\mathfrak{S}_{C_0\left(p\right)}$ {algebraic} $\overline \sX$-$ \widehat G$-{category}
\end{empt}
\begin{lemma}\label{top_inf_sufficienf_lem}
Under the hypotheses \ref{top_inf_sufficienf_empt} the triple $\left( C_0\left(\sX \right) ,  C_0\left(\overline\sX \right), \widehat G\right)$ is a {pre}-\textit{covering of the algebraical finite covering category}  $\mathfrak{S}_{C_0\left(p\right)}$ (cf. Definition \ref{algebraical_finite_covering_category_defn}).
\end{lemma}
\begin{proof}
Firstly we prove that  $\left( C_0\left(\sX \right) ,  C_0\left(\overline\sX \right), \widehat G\right)$ is  an  {infinite quasi-covering} (cf. Definition \ref{infinite_quasicovering_defn}. From the Corollary \ref{comm_lift_desc_sum_cor} it follow that for any $\la \in \La$ and $\overline{a }\in K\left(C_0\left(\overline\sX \right) \right)= C_c\left(\overline\sX \right)$ a series 
$$
 	a_\la =\bt\text{-} \sum_{	g \in \ker\left( \widehat{G}\to G_\la\right) }g \overline a
$$
is convergent with respect to the strict topology of $M\left(C_0\left(\overline\sX \right)\right) =  C_b\left(\overline\sX \right)$. So if $a_\la$ can be regarded as an element of $C_0\left(\overline\sX \right)$. On the other hand $a_\la$ is invariant with respect to action of $H_\la \bydef \ker\left( \widehat{G}\to G_\la\right)$, so $a_\la$ can be regarded as an element of $C_b\left( \sX_\la\right)$ where $\sX_{\la}\bydef \overline \sX / H_\la$ (cf. the equation \eqref{top_xh_eqn}). Indeed one has
\be\label{top_a_l_eqn}
\desc_{\la}\left(\overline a \right)  = \desc^c_{\overline p_\la}\left(\overline a \right)
\ee
where $\desc_\la$  is the $\la$-{descent} (cf. Definition \ref{infinite_desc_defn}) and $\desc^c_{\overline p_\la}\left(\overline a \right)$ is the  compactly supported $\overline p_\la$ descent. From the Theorem  \ref{top_compact_img_thm} it follows that $\desc^c_{\overline p_\la}\left(\overline a \right)\in C_c\left(\sX_\la \right)$ and one can proof that any element of  $ C_c\left(\sX_\la \right)$ is a finite linear combination of given by the equation \eqref{top_a_l_eqn} elements. The $C^*$-norm completion of  $ C_c\left(\sX_\la \right)$ is $C_0\left(\sX_\la \right)$ so the $C^*$-algebra $C_0\left(\sX_\la \right)$ is the $\la$-descent of $\overline{A}$ (cf. Definition \ref{infinite_quasicovering_defn})
 From the Theorem \ref{top_finite_covering_thm} it follows that any $*$-homomorphism $ C_0\left( p^\mu_\nu\right)  :C_0\left( \mathcal{X}_\mu\right)  \hookto C_0\left( \mathcal{X}_\nu\right)$ is a noncommutative finite-fold covering i.e. under hypotheses of this lemma the condition (a) of the Definition \ref{algebraical_finite_covering_category_defn}   is  satisfied. We leave to the reader an elementary  proof of that the   triple $\left( C_0\left(\sX \right) ,  C_0\left(\overline\sX \right), \widehat G\right)$ satisfies to the condition (b) of the Definition \ref{algebraical_finite_covering_category_defn} 
\end{proof}

\begin{lemma}\label{top_discontinuous_lem}
	Under the hypotheses \ref{top_inf_sufficienf_empt} the triple $\left( C_0\left(\sX \right) ,  C_0\left(\overline\sX \right), \widehat G\right)$ is the {disconnected infinite noncommutative covering} of the given by \eqref{top_x_g_eqn} algebraical finite covering category $\mathfrak{S}_{C_0\left(p\right)}$ (cf. Definitions \ref{algebraical_finite_covering_category_defn} and \ref{disconnected_infinite_noncommutative_covering_defn}).
\end{lemma}
\begin{proof}
If 	$\left( C_0\left(\sX \right) ,  \overline A', \widehat G\right)$ is a  {pre}-{covering of the algebraical finite covering category}  $\mathfrak{S}_{C_0\left(p\right)}$
then from the Corollary \ref{inverse_lim_h_cor} one can deduce that there is the natural inclusion
$$
\phi_{\overline A'}: \overline A' \hookto C\left(\overline \sX_{\mathrm{discr}}\right)	\subset B\left( \overline \H_a\right) 
$$
	where $C \left( \overline \sX_{\mathrm{discr}}\right)$ is the {discontinuous extension of} $C_0\left( \overline \sX\right)$ (cf. Definition \ref{top_discr_defn}, and
where $\overline \H_a$ is a Hilbert space of the atomic representation (cf. Definition \ref{atomic_repr_defn}) $C^*$-$\varinjlim_{\la \in \La} C_0\left(\sX_\la \right) \hookto  B\left( \overline \H_a\right)$ of the inductive limit in the sense of the Definition \ref{inductive_lim_non_defn}. From the Theorem \ref{top_discrete_p_thm} it follows that
$$
\forall \overline a \in \overline A'\quad \phi_{\overline A'}\left( \overline a\right) \in C_c\left(\overline \sX\right)
$$
From this circumstance it follows that $\overline A'\subset C_0\left( \overline \sX\right)$.  Taking into account the Lemma \ref{top_inf_sufficienf_lem} we conclude that $\left( C_0\left(\sX \right) ,  C_0\left(\overline\sX \right), \widehat G\right)$ is the terminal object of the category $\mathfrak{Cov}\left(\mathfrak{S}_{C_0\left(p\right)} \right)$ (cf. Definition \ref{disconnected_infinite_noncommutative_covering_defn}), i.e. $\left( C_0\left(\sX \right) ,  C_0\left(\overline\sX \right), \widehat G\right)$ is the {disconnected infinite noncommutative covering} of $\mathfrak{S}_{C_0\left(p\right)}$.
\end{proof}

 	\begin{empt}\label{comm_transitive_empt} 
 		From the Lemma \ref{top_discontinuous_lem}	 it follows that  the disconnected infinite noncommutative covering of $	\mathfrak{S}_p$  (cf. Definition \ref{disconnected_infinite_noncommutative_covering_defn}) is isomorphic to $C_0\left( \overline{   \mathcal X }\right)$. 
		From (iv)  of the Lemma \eqref{top_disconnected_lem} one has
	\be\nonumber
\widehat G= \varprojlim_{\la \in \La} G\left(\left.\sX_\la \right|\sX\right).
	\ee
The natural action $\widehat G\times C_0\left( \overline{   \mathcal X }\right)\to   C_0\left( \overline{   \mathcal X }\right)$ comes from the action 
	$ \widehat G\times  \overline{   \mathcal X} \to C_0\left( \overline{   \mathcal X }\right).
$
 From the equation \eqref{top_disconnected_eqn} it follows that any quasi-component of $\overline \sX$ is homeomorphic to $\widetilde{\mathcal X}$. 
	If $\widetilde A \subset C_0\left(\overline{\mathcal X}\right)$ is a connected component of $C_0\left(\overline{\mathcal X}\right)$ (cf. Definition \ref{connected_comp_defn}) then from the Exercise \ref{conn_comp_exer} it turns out that
 $$
\widetilde A \cong  C_0\left(\widetilde{\mathcal X}\right)\cong \left\{\left.\widetilde a \in C_0\left(\overline{\mathcal X}\right)\right|\widetilde a\left( \overline \sX \setminus \widetilde\sX' \right)= \{0\}  \right\}.
 $$
Moreover one has $C_0\left(\overline{\mathcal X}\right)= C_0\left(\widetilde {\mathcal X}\right)\oplus C_0\left(\widetilde {\mathcal X}^\perp \right)$ where $\widetilde {\mathcal X}^\perp \bydef \overline{\mathcal X} \setminus \widetilde {\mathcal X}$. Similarly to  equations \eqref{infinite_covering_transformation_group_eqn}  we define
	
\be\label{top_conn_group_eqn}
\begin{split}
G\left(\left. C_0\left(\widetilde{\mathcal X}\right)~\right| C_0\left(\sX\right) \right)\bydef \\\bydef\left\{\left.g \in  \widehat G~\right| g \widetilde a^\perp = \widetilde a^\perp\quad  \forall \widetilde a^\perp \in C_0\left( \widetilde {\mathcal X}^\perp\right)  \right\}\cong \\\cong \left\{\left.g\in  \widehat G\right| g \widetilde x^\perp = \widetilde x^\perp\quad  \forall \widetilde x^\perp \in \widetilde {\mathcal X}^\perp \right\}\cong G\left(\left. C_0\left(\widetilde{\mathcal X}\right)~\right| C_0\left(\sX\right)\right) .
\end{split}
\ee	

and taking into account $\widetilde {\mathcal X} = \widetilde {\mathcal X}\sqcup \bydef \overline{\mathcal X}^\perp$   one has
	\be\label{top_tga_eqn}
G\left(\left. C_0\left(\widetilde{\mathcal X}\right)~\right| C_0\left(\sX\right) \right)	\cong  G\left(\left. \widetilde\sX~\right|\sX \right).
\ee
If both $\widetilde{\mathcal X}', \widetilde{\mathcal X}'' \subset \overline{\mathcal X}$ are quasi-components   (cf. Definition \ref{top_connected_component_defn}) then from the Lemma \ref{top_disconnected_lem} it follows that there is $g \in G\left(\left.\overline{\sX} \right|\sX\right)$ such that $\widetilde{\mathcal X}''= g\widetilde{\mathcal X}'$. So one has $C_0\left(\widetilde{\sX}'' \right) = gC_0\left(\widetilde{\sX}' \right)$. For any $\la\in \La$ the natural homomorphism
$$
 G\left(\left. \widetilde\sX~\right|\sX \right) \to  G\left(\left. \sX_\la~\right|\sX \right)
$$
is surjective, so one has the surjective homomorphism
$$
G\left(\left. C_0\left( \widetilde\sX\right) ~\right|C_0\left( \sX \right) \right) \to  G\left(\left. C_0\left( \sX_\la\right) ~\right|C_0\left( \sX\right)  \right).
$$
From the above construction one can deduce the following theorem.	
	
\end{empt}

\begin{theorem}\label{top_main_thm}
	Under the hypotheses  of the Definition \ref{top_disconnected_defn} which include
	\begin{itemize}
		\item the {disconnected  $\overline p: \overline \sX  \to  \sX$ covering of} $p : \widetilde \sX \to \sX$,
		\item the finite covering category 	$\mathfrak{S}_p \bydef \left\{\left\{\sX_\la\right\}_{\la \in \La}, \left\{p^\mu_\nu:\sX_\mu\to \sX_\nu\right\}\right\}$ of $p : \widetilde \sX \to \sX$,
	\end{itemize}
following conditions hold.
\begin{enumerate}
	\item [(i)]  The given by \eqref{top_cb_defna_eqn} algebraic finite covering category (cf. Definition \ref{algebraical_finite_covering_category_defn}) $\mathfrak{S}_{C_0\left(p \right) }$ is good  (cf. Definition \ref{good_defn}) and the triple $$\left(C_0\left(\mathcal{X}\right), C_0\left(\varprojlim  \mathfrak{S}_p\right),G\left(\left.\varprojlim  \mathfrak{S}_p~\right| \mathcal X\right)\right)$$ is  the  {infinite noncommutative covering} of $\mathfrak{S}_{C_0\left( p\right) }$ (cf. Definition \ref{infinite_noncommutative_covering_defn}), where $\varprojlim  \mathfrak{S}_{p}$ is covering inverse limit of  given by  \eqref{top_category_fin_eqn} category (cf. Definition \ref{top_disconnected_defn}.
	\item[(ii)] There are  isomorphisms:
	\begin{itemize}
		\item $\varprojlim  \mathfrak{S}_{C_0\left(p\right)} \approx C_0\left(\varprojlim  \mathfrak{S}_{p}\right)$.
		\item $G\left(\left.\varprojlim  \mathfrak{S}_{C_0\left(p\right)}~\right| C_0\left(\mathcal X\right)\right) \approx G\left(\left.\varprojlim  \mathfrak{S}_{p}~\right| \mathcal X\right)$ 
\end{itemize}
		where the notation of the Definition \eqref{inf_lim_not_eqn} is used.

\end{enumerate}

\end{theorem}
\begin{empt}\label{top_r_inf_empt}
Under the hypotheses  of the Definition \ref{top_disconnected_defn} one has  noncommutative finite-fold coverings of operator spaces (cf. Definition \ref{fin_op_defn})
\be\label{top_comp_cat_eqn}
\begin{split}
\left(\left(C_0\left(\sX\right), C\left(\sX^\sim\right) \right) , \left(C_0\left(\sX_\la\right), C\left(\sX^\sim_\la\right) \right), G_\la, \left(\pi_\la, \pi^\sim_\la\right) \right),\\\left(\left(C_0\left(\sX_\mu\right), C\left(\sX^\sim_\mu\right) \right) , \left(C_0\left(\sX_\nu\right), C\left(\sX^\sim_\nu\right) \right) G_{\mu\nu}, \left(\pi_{\mu\nu}, \pi^{\sim}_{\mu\nu}\right) \right)
\end{split}
\ee
 where the $ ^\sim$ operation is given by \eqref{top_sim_eqn}, and 
 $G_\la\bydef  G\left(\left.\sX_\la\right|\sX \right)$ and  $G_{\mu\nu}\bydef G\left(\left.\sX_\mu\right|\sX_\nu \right)$.
 On the other hand from  the given by \eqref{top_real_fin_eqn}
 noncommutative finite-fold covering  the sub-unital real operator spaces (cf. Definition \ref{fin_rop_defn})  one has a family  finite-fold coverings  the sub-unital real operator spaces
\be\label{top_real_fam_eqn}
\begin{split}
\left(\left(C_0\left(\sX, \R\right), C\left(\sX^\sim, \R\right) \right) , \left(C_0\left(\sX_\la, \R\right), C\left(\sX^\sim_\la, \R\right) \right), G_\la, \left(\rho_\la, \rho^\sim_\la\right) \right),\\ 
\left(\left(C_0\left(\sX_\mu, \R\right), C\left(\sX^\sim_\mu, \R\right) \right) , \left(C_0\left(\sX_\nu, \R\right), C\left(\sX^\sim_\nu, \R\right) \right) G_{\mu\nu}, \left(\rho_{\mu\nu}, \rho^{\sim}_{\mu\nu}\right) \right)
\end{split}
\ee
 such that
\be\label{top_creal_fam_eqn}
\begin{split}
C_0\left( \sX\right) = \C C_0\left( \sX, \R\right),\quad
C_0\left( \sX^\sim\right) = \C C_0\left( \sX^\sim , \R\right),\\ 
C_0\left( \sX_\la\right) = \C C_0\left( \sX_\la, \R\right),\quad 
C_0\left( \sX^\sim_\la\right) = \C C_0\left( \sX^\sim_\la, \R\right),\\
\pi_\la= \C\rho_\la, \quad  \pi_{\mu\nu}= \C\rho_{\mu\nu},\quad  \pi^\sim_\la= \C\rho^\sim_\la, \quad  \pi^\sim_{\mu\nu}= \C\rho^\sim_{\mu\nu}.
\end{split}
\ee
(cf. \eqref{top_complex_eqn}). From \eqref{top_creal_fam_eqn} it turns out that the given by \eqref{top_comp_cat_eqn} noncommutative finite-fold coverings of operator spaces is the complexification of the family \eqref{top_comp_cat_eqn}. It follows that the given by \eqref{top_creal_fam_eqn}  family  $\mathfrak{S}_{\text{op}}$ is a algebraical  finite covering category of real operator spaces (cf. Definition \ref{comp_rop_pt_defn}).
\end{empt}
\begin{empt}\label{top_mu_empt} 
	Denote by 
		\bean
	\begin{split}
		\mathfrak{S}_{ C_0\left(p \right) }  = \\
		=\left(\left\{ C_0\left( p_\la\right) :  C_0\left(\sX \right)  \hookto  C_0\left(\sX_\la \right)\right\}, \left\{ C_0\left( p_\nu^\mu\right) :  C_0\left(\sX_\nu \right) \hookto  C_0\left(\sX_\mu \right)\right\}\right).
	\end{split}
	\eean
	the given by the Corollary \ref{top_cb_defna_eqn}
algebraical  finite covering category (cf. Definition \ref{algebraical_finite_covering_category_defn}).
	Let $\tau: C_c\left(\sX \right)\to \C$ be a positive functional given by the Theorem \eqref{meafunc_eqn}, i.e. there exists a Borel
	measure $\mu$ on $\sX$  such that
	\be\label{top_meafunc_eqn}
	\tau\left(a \right) = \int_{\sX} a~ d\mu \quad \forall a\in C_c\left(\sX \right).
	\ee
	Suppose that $\mu\left( \sU\right) > 0$ for any  nonempty open set $\sU\subset\sX$, or equivalently 
	$$
	\forall a \in C_c\left(\widetilde \sX\right)\quad a > 0 \quad \Rightarrow \tau\left(a \right) > 0.
	$$

	Let $\overline{\mathcal X}$ be to be the {disconnected inverse limit} of $\mathfrak{S}_{p}$	(cf. Definition \ref{top_disconnected_defn}), and let $\overline p: \overline{\mathcal X}\to\sX$, $\overline p_\la: \overline{\mathcal X}\to\sX_\la$ be natural coverings. From the Lemma \ref{state_cov_11_lem} it follows that for all $\la\in \La$ there is a functional
	\be\label{top_wtl_eqn}
	\begin{split}
		\tau_\la : C_c\left(\sX_\la \right)\to \C,\\
		a_\la \mapsto \tau\left( \sum_{g \in G\left(\left.\sX_\la~\right|\sX \right) }ga\right) 
	\end{split}
	\ee
	where a $\sum_{g \in G\left(\left.\sX_\la~\right|\sX \right)}ga$ is regarded as an element of $C_c\left( \sX\right)$. From the Theorem \ref{meafunc_thm} it turns out that $\tau_\la$ is given by 	
	\bean
	\tau_\la\left(a_\la \right) = \int_{\sX_\la} a_\la~ d{\mu_\la} \quad \forall a_\la \in C_c\left(\sX_\la \right).
	\eean
	where $\mu_\la$ is  a Borel measure $\sX_\la$. Indeed for all $\la\in\La$ the measure $\mu_\la\bydef  \lift_{p_\la}\mu$ is the $p_\la$-lift of $\mu$ (cf. Definition \ref{top_lift_measure_defn}). For all $\la\in \La$ similarly to equations \eqref{tau_prod_eqn} and \eqref{mu_prod_eqn} one can define the $\C$-valued product $\left(\cdot, \cdot \right):C_c\left( \sX_\la\right) \times C_c\left( \sX_\la\right) \to\C$ given by
	\bean
	\left(a_\la, b_\la \right) \bydef \tau_\la\left( a^*_\la, b_\la\right) \cong  \int_{\sX_\la} a^*_\la b_\la~ d{\mu_\la} 
	\eean
	so $C_c\left( \sX_\la\right)$ becomes a pre-Hilbert space.
	Denote by $L^2\left(\sX_\la, \mu_\la \right)$ the Hilbert norm completion of $C_c\left( \sX_\la\right)$.  Let us define a functional
	\bean
	\overline \tau: C_c\left(\overline  \sX\right) \to\C,\\
	\overline a \mapsto \int_{\overline\sX}  \overline a_\la~ d \overline \mu
	\eean
	where $\overline \mu \bydef  \lift_{\overline p}\mu$ is the $\overline p$-lift of $\mu$ (cf. Definition \ref{top_lift_measure_defn}). 
	One has
	\be\label{top_ot_eqn}
	\overline\tau\left(\overline a \right)= \tau_\la\left(\desc_{\overline p_\la}\left( \overline  a\right)\right)\quad\la\ge \la_{\supp \overline  a}
	\ee
	
	Let us define the $\C$-valued product $\left(\cdot, \cdot \right):C_c\left( \overline\sX\right) \times C_c\left(\overline  \sX\right) \to\C$.
	The product is given by
	\bean
	\left(\overline  a, \overline  b\right)\bydef \overline\tau \left( \overline  a^*\overline  b\right).
	\eean
	Thus  $C_c\left(\overline  \sX\right)$ becomes a pre-Hilbert space denote by $L^2\left(\overline  \sX , \overline \mu\right)$ its norm completion.  A natural action $C_b\left(\overline  \sX \right)  \times C_c\left(\overline  \sX \right)\to C_c\left(\overline  \sX \right)$ yields an action $C_b\left(\overline  \sX \right)  \times  L^2\left(\overline  \sX , \overline \mu\right)\to L^2\left(\overline  \sX , \overline \mu\right)$ which gives  a faithful representation	\be\label{top_phi_mu_eqn}
	\overline	\pi_\tau :  C_b\left(\overline  \sX \right)\hookto B\left(L^2\left(\overline  \sX , \overline \mu\right) \right).
	\ee
	On the other hand from if $\widehat{C_0\left(\sX \right) }\bydef C^*\text{-}\lim_{\la\in \La} C_0\left(  \sX_\la \right)$ then from \eqref{ininite_covering_inductive_eqn}  one has an inclusion $\widehat{C_0\left(\sX \right) }\hookto M\left(C_0\left(\overline  \sX \right) \right) = C_b\left(\overline  \sX \right)$ so the representation $\overline	\pi_\tau$ yields a representation 	\be\label{top_pi_mu_eqn}
	\widehat	\pi_\tau\bydef\left.\overline	\pi_\tau \right|_{\widehat{C_0\left(\sX \right) }}: \widehat{C_0\left(\sX \right) }\hookto B\left(L^2\left(\overline  \sX , \overline \mu\right) \right).
	\ee
\end{empt}
\begin{exercise}\label{top_pi_mu_exer}
	Prove following statements.
	\begin{itemize}
		\item The representation \ref{top_phi_mu_eqn} is faithful and nondegenerate (cf. Definitions \ref{faithful_inf_lem} and \ref{nondegenerate_repr_defn}).
		\item The representation \eqref{top_pi_mu_eqn} is equivariant (cf. Definition \ref{equivariant_representation_defn}.
	\end{itemize}
\end{exercise}

 \subsection{Universal  coverings and  fundamental groups}
\paragraph*{}

A following theorem gives universal  coverings of commutative $C^*$-algebras. 

\begin{theorem}\label{comm_uni_lim_thm}
	If $\widetilde p: \widetilde \sX \to \sX$ is a transitive covering such that:
	\begin{enumerate}
		\item[(a)] both $\sX$ and $\widetilde \sX$ are connected, locally compact, Hausdorff spaces,
		\item[(b)] For any transitive finite-fold covering $p': \sX' \to \sX$ there is a transitive covering $\widetilde p': \widetilde \sX \to \sX'$  with $\widetilde p= p'\circ \widetilde p'$.
	\end{enumerate}
then one has:
	\begin{enumerate}
		\item [(i)]  	a triple 
	\be\label{comm_uni_lim_eqn}
\left(C_0\left(\sX\right), C_0\left( \widetilde{\sX}_{\mathrm{res~fin}}\right),G\left( \left.{\widetilde{\sX}}_{\mathrm{res~fin}}~\right|\sX\right)\right)	
	\ee	
 with given by \eqref{top_x_fr_eqn}	$\widetilde{\sX}_{\mathrm{res~fin}}$, is an universal covering of $C_0\left( {\sX}\right)$ (cf. Definition \ref{fundamental_group_nc_defn});
	\item[(ii)]	there is a  natural  
	group isomorphism 
	\be\label{comm_fg_iso_eqn}
	\pi_1\left( C_0\left(\sX \right)\right)\cong G\left( \left.{\widetilde{\sX}}_{\mathrm{res~fin}}~\right|\sX\right) 
	\ee
	where  $\pi_1\left( C_0\left(\sX \right)\right)$ 
	is the {fundamental group}	of   $ C_0\left(\sX \right)$  (cf. Definition \ref{fundamental_group_nc_defn}.
	
	\end{enumerate}

\end{theorem}

\begin{proof}
	(i) Denote by $p:  \widetilde \sX_{\mathrm{res~fin}}\to \sX$ the natural covering.
	From the construction \ref{top_fin_sec_empt} it follows that the group $G\left( \left.{\widetilde{\sX}}_{\mathrm{res~fin}}~\right|\sX\right)$ is residually finite (cf. the equation \eqref{top_x_fr_eqn} and the Definition \ref{residually_finite_defn}).  From the hypothesis (b) of this theorem it turns out that for any transitive finite-fold covering $p': \sX' \to \sX$ there is a transitive covering $\widetilde p'_{\text{fin}}: \widetilde \sX_{\mathrm{res~fin}} \to \sX'$  with $ p= p'\circ \widetilde p'$. From the Lemma \ref{top_associded_lem} it follows that 
	 the natural covering $\widetilde{\sX}_{\mathrm{res~fin}}\to   {\sX};$ is the {covering inverse limit of} 
		of the category
		\bean
		\mathfrak{S}_p \bydef \left\{\left\{\sX_\la\right\}_{\la \in \La}, \left\{p^\mu_\nu:\sX_\mu\to \sX_\nu\right\}_{\substack{\mu,\nu \in \La\\\mu\ge\nu}}\right\}.
		\eean
which includes all finite-fold coverings of $\sX$. Taking into account the Lemmas  \ref{top_fin_necassary_lem}  and \ref{top_fin_sufficient_lem}
 the category $\mathfrak{S}_p$ contains any noncommutative finite-fold covering of $C_0\left(\sX\right)$.
From the   Definition \ref{fundamental_group_nc_defn} it follows that the triple \eqref{comm_uni_lim_eqn} is an universal covering of $C_0\left( {\sX}\right)$.\\
(ii) Follows from (i).

\end{proof}

\begin{corollary}\label{comm_uni_lim_cor}
	Let $\mathcal X$ be a connected, weakly semilocally 1-connected (cf. Definition \ref{top_weakly_semi1_defn}), locally compact, locally connected,  Hausdorff space. If  $\widetilde \sX_{\mathrm{res~fin}}$ is a given by the equation \eqref{top_x_fr_eqn} space, and   $\widetilde{p}:\widetilde{\sX}_{\mathrm{res~fin}}\to \mathcal X$ is the natural covering (cf. \eqref{top_fin_c_sec_eqn}) then one has:
	\begin{enumerate}
		\item [(i)]  	a triple 
		\be\label{comm_uniw_lim_eqn}
		\left(C_0\left(\sX\right), C_0\left( \widetilde{\sX}_{\mathrm{res~fin}}\right),G\left( \left.{\widetilde{\sX}}_{\mathrm{res~fin}}~\right|\sX\right)\right)	
		\ee	
with given by \eqref{top_x_fr_eqn}	$\widetilde{\sX}_{\mathrm{res~fin}}$	is an universal covering of $C_0\left( {\sX}\right)$ (cf. Definition \ref{fundamental_group_nc_defn}),
		\item[(ii)]	there is a  natural  
		group isomorphism 
		\be\label{comm_fg_iso_eqn!}
		\pi_1\left( C_0\left(\sX \right)\right)\cong \mathfrak{ResFin}\left( 	\pi_1^{\mathrm{w}}\left( \sX, x_0\right)\right) 
		\ee
		where  $\pi_1\left( C_0\left(\sX \right)\right)$ 
		is the {fundamental group}	of $ C_0\left(\sX \right)$  (cf. Definition \ref{fundamental_group_nc_defn} and $\pi_1^{\mathrm{w}}\left( \sX, x_0\right)$ is a weak fundamental group (cf. Definition \ref{top_weak_fundamental_group_defn}) and $\mathfrak{ResFin}$ is given by the equation \eqref{top_fin_eqn}.
		
	\end{enumerate}
\end{corollary}
\begin{proof}
	From the Exercise \ref{top_uni_cov_exer} it follows that there is the universal covering of $\sX$, so this corollary is a consequence of the Theorem \ref{comm_uni_lim_thm}.
\end{proof}

\begin{exercise}\label{comm_uni_lim_exe}
	Let $\mathcal X$ be a connected, locally path connected, weakly semilocally 1-connected (cf. Definition \ref{top_weakly_semi1_defn}), locally compact,  Hausdorff space then $C_0\left( \widetilde{\sX}_{\mathrm{res~fin}}, \R\right)$ is the universal covering of the real sub-unital operator space $\left( C_0\left( {\sX}, \R\right),C_0\left( {\sX}^\sim, \R\right)\right)$ (cf. Remark \ref{real_op_rem}). Moreover there is an isomorphism of fundamental groups
	$$
	\pi_1\left(\left( C_0\left( {\sX}, \R\right),C_0\left( {\sX}^\sim, \R\right)\right)\right) \cong \mathfrak{ResFin}\left( \pi_1^{\text{w}}\left( \sX\right)\right) . 
	$$
\end{exercise}

\subsection{Coverings of *-algebras}

\subsubsection{Coverings of pro-$C^*$-algebras}
\paragraph*{}
Let $p : \widetilde \sX \to \sX$ is the {covering inverse limit of} a category 
$$
	\mathfrak{S}_p \bydef \left\{\left\{\sX_\la\right\}_{\la \in \La}, \left\{p^\mu_\nu:\sX_\mu\to \sX_\nu\right\}_{\substack{\mu,\nu \in \La\\\mu\ge\nu}}\right\}.
$$
(cf. Definition \ref{top_disconnected_defn}).
 From the equation   \ref{top_pro_hom_eqn}  it follows that  there are injective homomorphisms of pro-$C^*$-algebras 
\bean
Cont\left( p_\la\right): Cont\left(\sX\right)\hookto Cont\left(\sX_\la\right),\\
Cont\left( p^\mu_\nu\right): Cont\left(\sX_\nu\right)\hookto Cont\left(\sX_\mu\right).
\eean
From the Theorem \ref{top_pro_fin_thm} it follows that these $*$-homomorphisms are {noncommutative finite-fold  coverings  of pro-$C^*$-algebras} (cf. Definition \ref{pro_fin_defn}). So the pair 
\be\label{top_pro_cat_eqn}
\mathfrak{S}_{Cont\left(\sX \right) }\bydef \left\{ \left\{Cont\left( p_\la\right)\right\}_{\la\in \La}, \left\{Cont\left( p^\mu_\nu\right)\right\}_{\mu,\nu\in\La}\right\}
\ee
is a algebraical  finite covering category of pro-$C^*$-algebras (cf. Definition \ref{comp_pro_defn}).
\begin{empt}\label{top_pro_empt}
	Let us obtain the $\mathfrak{S}_{Cont\left(\sX \right) }$-specialization of the given by \ref{inv_pro_lim_empt} construction. If for any $\la\in \La$ the set $b\left(Cont\left(\sX_\la\right) \right)$ is a $C^*$-algebra of bounded elements of $Cont\left(\sX_\la\right)$   (cf. Definition \ref{pro_bound_defn} and Proposition \ref{pro_bound_prop}) then $ b\left(Cont\left(\sX_\la\right) \right)= C_b\left(\sX_\la\right)$. Moreover $B_\la \bydef C_0\left(\sX_\la\right)$ is an essential ideal of $C_b\left(\sX_\la\right)$.  If $B \bydef C_0\left( \sX\right)$ then from the Theorem \ref{top_main_thm}
	one has a good   { algebraical  finite covering category} (cf. Definition \ref{good_defn})
	\be\label{top_sb_pro_eqn}
	\begin{split}
		\mathfrak{S}_{B } \bydef \left(\left\{C_0\left( p_\la\right):B  \hookto B_\la  \right\}_{\la \in \La},\left\{C_0\left( p^\mu_\nu\right): B_\nu \hookto B_\mu \right\}_{\substack{\mu, \nu \in \La~~\nu > \mu}}\right).
	\end{split}
	\ee
	If $\widetilde \sX$ is a topological inverse limit  of $\mathfrak{S}_{\left(\sX, x_0 \right) }$ (cf. Definition \ref{top_disconnected_defn}) then from the Theorem \ref{top_main_thm} it follows that the triple 
	\bean
	\left(B , \widetilde{B}, \widetilde G\right)\bydef \left(C\left(\sX \right), C_0\left( \widetilde\sX\right), G\left(\left.\widetilde\sX\right|\sX\right)  \right) 
	\eean
	is an infinite noncommutative covering  of $\mathfrak{S}_{B}= \mathfrak{S}_{C_0\left(\sX \right) }$ (cf. Definition \ref{infinite_noncommutative_covering_defn}). One has $M\left(\widetilde B \right)\cong C_b\left(\widetilde\sX \right)$, so for any $C^*$-seminorm of $p_\a : C_b\left(\widetilde\sX \right)\to \R$ there is a closed subset $\widetilde{\sV}_\a \subset \widetilde\sX$ such that
	$$
	\forall \widetilde a \in C_b\left(\widetilde\sX \right) \quad p_\a\left(\widetilde a \right)  \bydef \sup_{\widetilde x \in \widetilde \sV}\left|\widetilde a\left(\widetilde x \right) \right|.
	$$
\end{empt}
\begin{exercise}\label{top_pro_exer}
	Prove following statements:
	\begin{enumerate}
		\item A $\widetilde G$-invariant family $\left\lbrace \widetilde \sV_\a\subset \widetilde \sX\right\rbrace_{\a \in \mathscr A}$ corresponds to an admissible family of $C^*$-seminorms  (cf. \ref{inv_pro_lim_empt}) if for all $\la\in \La$ a family   $\left\{\widetilde p_\la\left( \widetilde \sV_\a\right)\right\}_{\a\in \mathscr A}$ contains all  compact sets of $\sX_\la$;
		
		\item Any admissible pro-$C^*$-algebra (cf. \ref{inv_pro_lim_empt}) 
		is a *-subalgebra of pro-$C^*$-algebra  $Cont\left( \widetilde \sX\right)$;
		\item If a set  $\left\{\widetilde p_\a\right\}$ of $\widetilde G$-invariant  $C^*$-seminorms corresponds to a  family $\left\{\widetilde \sV_\a\right\}$ of all compact subsets of $\widetilde \sX$, then $\left\{\widetilde p_\a\right\}$ is admissible and a completion of $C_0\left(\widetilde \sX \right)$ with respect to the family of $C^*$-seminorms $\left\{\widetilde p_\a\right\}$ equals to  $Cont\left( \widetilde \sX\right)$.
	\end{enumerate}
\end{exercise}

\begin{theorem}\label{top_pro_thm}
	If we consider the situation \ref{top_pro_empt} then the triple 
	$$
	\left(Cont\left(\sX \right), Cont\left( \widetilde{\sX}\right) ,G\left(\left.\widetilde \sX\right| \sX \right) \right)  
	$$
	is an inverse noncommutative limit of the given by \eqref{top_pro_cat_eqn} {algebraical  finite covering category $\mathfrak{S}_{Cont\left(\sX \right) }$ of pro-$C^*$-algebras}  (cf. Definition \ref{inv_pro_lim_defn}).
\end{theorem}
\begin{proof}
	Follows from the Exercise \ref{top_pro_exer}.
\end{proof}

\subsubsection{Coverings of bounded operator $*$-algebras}
\begin{empt}
Let $p : \widetilde \sX \to \sX$ is the {covering inverse limit of} a category 
$$
\mathfrak{S}_p \bydef \left\{\left\{\sX_\la\right\}_{\la \in \La}, \left\{p^\mu_\nu:\sX_\mu\to \sX_\nu\right\}_{\substack{\mu,\nu \in \La\\\mu\ge\nu}}\right\}.
$$
(cf. Definition \ref{top_disconnected_defn}).
	Let $R\subset C_0\left(\sX \right)$ be a $c$-soft *-subalgebra (cf. Definition \ref{top_soft_r_defn}), and for any $\la\in \La$ let $R_\la\bydef  C_0\left(\sX_\la \right)\cap p^{-1}_\la \mathscr S^{R}\left( \sX_\la\right)$ and $\widetilde R\bydef  C_0\left(\widetilde \sX \right)\cap \widetilde p^{-1} \mathscr S^{R}\left(\widetilde \sX\right)$ where a notation of the Definition \ref{top_x_sheaf_defn} is used.	
	The Theorem \ref{top_oa_cov_thm} enables us to construct a algebraical  finite covering category of  bounded operator *-algebras (cf. Definition \ref{comp_oa_defn}) given by
	\be\label{top_inf_cov_oae_eqn}
	\begin{split}
		\mathfrak{S}_{R } \bydef
		\left(\left\{ \left.C_0\left( p_\la\right)\right|_R :  R  \hookto  R_\la\right\}, \left\{ \left.C_0\left( p^\mu_{ \nu}\right)\right|_{R_\nu} : R_\nu \hookto R_\mu\right\}\right).
	\end{split}
	\ee
	The formula \eqref{top_inf_cov_oae_eqn} is a specialization of the given by \eqref{comp_oa_eqn} one. Moreover from the  the Theorem \ref{top_main_thm} one has a 	good	{algebraical  finite covering category} (cf. Definitions \ref{algebraical_finite_covering_category_defn}  and \ref{good_defn}) given by
	\be\label{top_inf_cov_oa_eqn}
	\begin{split}
	\mathfrak{S}_{C_0\left(p\right) } \bydef \\
	\left\{ \left\{ C_0\left( p_\la\right)  :C_0\left( \mathcal{X}\right)  \hookto C_0\left( \mathcal{X}_\la\right) \right\}, \left\{ C_0\left( p^\mu_\nu\right)  :C_0\left( \mathcal{X}_\mu\right)  \hookto C_0\left( \mathcal{X}_\nu\right) \right\}  \right\}
\end{split}
	\ee
	
	The formula \eqref{top_inf_cov_oa_eqn} can be regarded as a specialization of \eqref{algebraical_finite_covering_category_eqn} one. From the Theorem \ref{top_main_thm} it follows that there is a transitive covering $~~\widetilde p:\widetilde \sX \to \sX$ such that 
	$
	\left(C_0\left(\sX\right), C_0\left(\widetilde \sX\right)  G\left(\left.\widetilde \sX \right|  \sX\right)\right) 
	$
	is  an infinite noncommutative covering of $\mathfrak{S}_{C_0\left(\sX\right) }$ (cf. Definition \ref{infinite_noncommutative_covering_defn}).
	So we have all described in  \ref{inf_oa_empt} ingredients. 
\end{empt}
\begin{lemma}\label{top_inf_cov_oa_lem}
	Consider the above situation.  If $\widetilde A$ is a *-algebra of continuous maps from $\widetilde \sX$ to $\C$  such that
	\be\label{top_inf_cov_oap_eqn}
	\forall \widetilde b \in C_c\left( \widetilde \sX \right) \quad \forall \widetilde a \in  \widetilde A \quad \exists a^\cup \in \bigcup_{\la\in \La}~ \lift_{ \widetilde{p}_{\la}}\left(  R_\la\right) \quad
	\widetilde b  \widetilde a =  \widetilde b a^\cup
	\ee
	where inclusions $\lift_{ \widetilde{p}_{\la}} : C_0\left(\sX_\la  \right)\hookto C_b\left(\widetilde \sX\right)$ are  given by the Definition \ref{top_lift_defn} injective $*$-homomorphisms then there is an  inclusion
	$$
	\widetilde A \subset \mathscr S^{\widetilde R}\left( \widetilde \sX\right) 
	$$
	where $\mathscr S^{\widetilde R}$ is the $\widetilde R$-sheaf (cf. Definition \ref{top_x_sheaf_defn}).
\end{lemma}

\begin{proof}
	Algebra $\widetilde A$ contains continuous maps only, it turns out that $\mathscr S^{\widetilde A}_{\widetilde x}\subset \mathscr S^{C_0\left( \widetilde \sX\right)}_{\widetilde x}$ for any $\widetilde x \in \widetilde \sX$ where $\mathscr S^{\widetilde A}_{\widetilde x}$ means the space of stalks at $\widetilde x$ (cf. Definitions \ref{sheaf_stalk_defn} and \ref{top_x_sheaf_defn}). On the other hand 
	from $\widetilde R \subset C_0\left( \widetilde \sX\right)$ it follows  that $\mathscr S^{\widetilde R}_{\widetilde x}\subset \mathscr S^{C_0\left( \widetilde \sX\right)}_{\widetilde x}$ for any $\widetilde x \in \widetilde \sX$. We would like prove that 
	\be\label{top_ra_sh_eqn}
	\forall \widetilde x \in \widetilde \sX	\quad \mathscr S^{\widetilde A}_{\widetilde x} \subset \mathscr S^{\widetilde R}_{\widetilde x}
	\ee
	If $\widetilde a_{\widetilde x} \in \mathscr S^{\widetilde A}_{\widetilde x}$ then $\widetilde a_{\widetilde x}$ can be represented by an element $\widetilde a\in\widetilde A$. From the Exercise \ref{top_soft_c_exer} it follows that $\widetilde R$ is a $c$-soft *-algebra (cf. Definition \ref{top_soft_r_defn}). If $\widetilde b \in  \widetilde R$ is a $\left(\widetilde p, R, \widetilde x\right)$-{stump} (cf. Definition \ref{top_stump_soft_p_defn})
	then from  $\widetilde b \in C_c\left(\widetilde \sX \right)$ and \eqref{top_inf_cov_oap_eqn} it follows that there is $a^\cup \in 	\bigcup_{\la\in \La}~ \lift_{ \widetilde{p}_{\la}}\left(  R_\la\right)$ such that
	$$
	\widetilde b \widetilde a= \widetilde b a^\cup
	$$
	(cf. equation \eqref{top_inf_cov_oap_eqn}). On the other hand if $\widetilde \sV$ is such that $\widetilde x \subset \widetilde \sV$ and $\widetilde b\left( \widetilde \sV\right)= \{1\} $
	then one has
	$$
	\left.\left( \widetilde b \widetilde a\right)\right|_{\widetilde \sV } = \left.\widetilde a\right|_{\widetilde \sV }.
	$$
	It follows that $\left( \widetilde b \widetilde a\right)_{\widetilde x}=  \widetilde a_{\widetilde x}\in \mathscr S^{\widetilde A}_{\widetilde x}$. On the other hand there is $\la_0\in \La$  and $a^{\la_0}\in R_{\la_0}$ such that $a^\cup = \lift_{p_{\la_0}}\left(a^{\la_0}\right)$ and
	$$
	\left.\left( \widetilde b \widetilde a\right)\right|_{\widetilde \sV } = \left.\lift_{p_{\la_0}}\left(a^{\la_0}\right)\right|_{\widetilde \sV }.
	$$
	It turns out that $\widetilde a_{\widetilde x}= \left( \widetilde p_{\widetilde x}^\la\right)^{-1} \left(a^\la_{\widetilde p^\la\left( \widetilde x\right) } \right) $ 
	where $\widetilde p^\la : \widetilde \sX \to \sX_\la$ is the natural covering and the isomorphism $\widetilde p_{\widetilde x}^\la$ corresponds to the equation \eqref{top_st_iso_eqn}. 
	On the other hand from  $a^{\la_0}\in R_{\la_0}$ it follows that $a^\la_{\widetilde p^\la\left( \widetilde x\right) }\in p_\la^{-1}\left( \mathscr S^R \right)_{p^\la\left(\widetilde x \right)}$ where $p^\la: \sX_\la \to \sX$ is the natural covering.
	From the equations \eqref{top_sh1_comp_eqn} and 
	\eqref{top_sh2_comp_eqn} it follows that 
	$$
	\widetilde a_{\widetilde x}\in \widetilde p^{-1}\mathscr S^R_{\widetilde x}= \mathscr S^{\widetilde R}_{\widetilde x}.
	$$
	From the above equation it follows that $\widetilde a\in \mathscr S^{\widetilde R}\left( \widetilde \sX\right)$. But $\widetilde a$ is an arbitrary element of $\widetilde A$ so one has
	$$
	\widetilde A \subset\mathscr S^{\widetilde R}\left( \widetilde \sX\right) 
	$$
	
\end{proof}
\begin{theorem}\label{top_inf_cov_oa_thm}
	In the situation of the Lemma \ref{top_inf_cov_oa_lem} the triple
	$$
	\left(R,\widetilde R, G\left(\left.\widetilde \sX \right|  \sX\right)\right) 
	$$ is an {infinite noncommutative covering}  of $\mathfrak{S}_{R}$ in the sense of the  Definition \ref{inf_cov_oa_defn}.
\end{theorem}

\begin{proof}
	Firstly we prove that $\widetilde R$   {is subordinated to}  $	\mathfrak{S}_{R} $ (cf. Definition \ref{inf_cov_oa_sub_defn}), i.e. it satisfies to equations \eqref{inf_cov_oa_sub1_eqn}-\eqref{inf_cov_oa_sub3_eqn}. From
	\bean
	\left( 	p^{-1}_\la \mathscr S^{R}\right)\left(\sX_\la\right) \left( \widetilde p^{-1} \mathscr S^{R}\right)\left( \widetilde \sX\right)  \subset \widetilde p^{-1} \mathscr S^{R}\left( \widetilde \sX\right) ,\\
	\left( \widetilde p^{-1} \mathscr S^{R}\right)\left( \widetilde \sX\right)\left( 	p^{-1}_\la \mathscr S^{R}\right)\left( \sX_\la\right)  \subset \widetilde p^{-1} \mathscr S^{R}\left( \widetilde \sX\right) 
	\eean 
	it follows that 
	\bean
	R_\la \widetilde R \subset  \widetilde R,\\
	\widetilde R \subset R_\la  \widetilde R
	\eean
	it means that 	$\widetilde R$ satisfies to the equations \eqref{inf_cov_oa_sub1_eqn}, \eqref{inf_cov_oa_sub2_eqn}. If  $\widetilde b \in K\left( C_0\left( \widetilde \sX\right)\right)\cong  C_c\left( \widetilde \sX\right)$ then $\supp \widetilde b$ is compact.  Let $\widetilde f \in  \widetilde R$ be a $\left(R, \supp \widetilde b\right)$-{stump} (cf. Definition \ref{top_stump_soft_y_defn}). Using the Lemma  \ref{top_compact_la_lem}  one can assume that there is an open subset $\widetilde \sU \subset \widetilde \sX$ and $\la_{\widetilde\sU}\in\La$ such that the natural covering $\widetilde p_{\la_{\widetilde\sU}} : \widetilde\sX \to \sX_{\la_{\widetilde\sU}}$ homeomorphically maps $\widetilde\sU$ onto $p_{\la_{\widetilde\sU}}\left(\widetilde\sU\right)$. For all $\widetilde a\in \widetilde R$ one has
	\bean
	a^\cup \bydef \desc_{p_{\la_{\widetilde\sU}}}\left( \widetilde f \widetilde a\right) \in R_{\la_{\widetilde\sU}}\subset \bigcup_{\la\in\La} R_\la
	\eean
	One the other hand from our construction one has
	$$
	\widetilde a \widetilde b = \widetilde b \widetilde a^\cup = a\widetilde b = \widetilde ba^\cup.
	$$
	The above equation is a specialization of \eqref{inf_cov_oa_sub3_eqn}, it follows that $\widetilde R$ is {subordinated  to}  $\mathfrak{S}_R$ (cf. Definition \ref{inf_cov_oa_sub_defn}). On the other hand from the Lemma \ref{top_inf_cov_oa_lem}
	it turns out that any {subordinated to}  $\mathfrak{S}_R$ *-algebra is a *-subalgebra of $\widetilde R$, i.e. $\widetilde R$ is the unique maximal {subordinated to}  $\mathfrak{S}_R$ *-algebra.
\end{proof}

\subsubsection{Coverings of $O^*$-algebras}\label{top_inf_o_sec}
\paragraph*{}

Let $M$ be a finite-dimensional smooth manifold, and let $p: \widetilde M \to M$ be a transitive covering such that the covering group $G\left(\left. {\widetilde{M}}~\right|~M\right)$ (cf. Definition  \ref{top_group_of_covering_transformations_defn} is residually finite (cf. Definition \ref{residually_finite_defn})
Let 
\be\label{top_secm_p_eqn}
\begin{split}
	\mathfrak{S}_p \bydef \left\{\left\{M_\la\right\}_{\la \in \La}, \left\{p^\mu_\nu:M_\mu\to M_\nu\right\}_{\substack{\mu,\nu \in \La\\\mu\ge\nu}}\right\}.
\end{split}
\ee
be  a finite covering category of  $p: \widetilde M \to M$ (cf. Definition \ref{top_disconnected_defn}).
The product $T \bydef M \times \C$ is also a smooth manifold, because $\C$ has the natural smooth structure. So one has a smooth bundle $T \to M$ (cf. Definition \ref{top_sm_bundle_defn}) having a given by \eqref{top_ses_one_eqn} sesquilinear product.
From the Theorem \ref{top_d_fin_thm} it turns out that there is  an {algebraic finite covering category of $O^*$-algebras} (cf. Definition \ref{comp_o*_defn}) given by
\be\label{top_comp_pt_o*_eqn}
\begin{split}
	\mathfrak{S}_{D^*\left(M, T \right) } = \left(\left\{\pi_\la: D^*\left( M, T\right)  \hookto D^*\left( M_\la, T_\la \right) \right\}, \left\{\pi^\mu_\nu\right\}\right).
\end{split}
\ee
such that  $D^*\left( M , T\right)$ and  $D^*\left( M_\la, T_\la \right)$ are a *-algebras of differential operators $P:\Coo\left( M, T\right) \to \Coo\left( M, T\right)$ and $P_\la:\Coo\left( M_\la, T_\la\right) \to \Coo\left( M_\la, T_\la\right)$ on $M$ and $M_\la$ respectively (cf. Definition \ref{do_man_defn}) 
which satisfy to the condition \eqref{top_diff*_eqn}  for all $\la\in\La$. Note that for every $\la\in\La$ the smooth structure of $M_\la$ comes from the smooth structure on $M$ and the
covering $M_\la\to M$. Injective $*$-homomorphisms $\pi_\la: D^*\left( M, T \right)\hookto D^*\left( M_\la, T_\la \right)$ and 
$\pi^\nu_\mu: D^*\left( M_\nu, T_\nu \right)\hookto D^*\left( M_\mu, T_\mu \right)$ come from the equation
\eqref{top_diff_*alg_lift_eqn}.  If $\overline M$ is the disconnected covering  of $p : \widetilde M \to M$ (cf. Definition \ref{top_disconnected_defn}) then there are natural coverings $\overline p: \overline M\to M$, $\overline p_\la: \overline M\to M_\la$ which induce the structure of smooth manifold on $\overline M$. So for all $\la \in \La$ there is a natural injective $*$-homomorphism $D^*\left( M_\la, T_\la \right)\hookto D^*\left(\overline M, \overline T \right)$. If these injective $*$-homomorphisms are regarded as inclusions and $\widehat A \bydef \bigcup_{\la \in \La}D^*\left( M_\la, T_\la \right)$ then there is the following inclusion
\be\label{top_dinc_eqn}
\widehat A\hookto D^*\left(\overline M, \overline T \right)
\ee
Let $\tau : C_0\left(M \right) \to \C$ be a faithful state  and let $L^2\left(\overline M, \overline \mu\right)$ be explained in  \eqref{top_mu_empt} Hilbert space. 
Then 
\be\label{top_tdd_eqn}
\widehat \D\bydef  \Coo_c\left(\overline M \right)\cong \Ga^\infty_c\left(\overline M, \overline T \right) 
\ee
is a dense subspace of $L^2\left(\overline M, \overline \mu\right)$. There is the natural inclusion
\be D^*\left( \overline M, \overline T\right) \subset \L^\dagger\left(\widehat\D \right)
\ee 
(cf. equation \eqref{l_dag_eqn}) and taking into account  \eqref{top_comp_pt_o*_eqn}  one has an inclusion 
\be\label{top_pi_eqn}
\pi: \widehat A\hookto \L^\dagger\left(\widehat\D \right).
\ee
(cf. the condition \eqref{top_diff*_eqn}).
If $\widehat G \bydef \varprojlim G\left(\left.M_\la\right|M\right)$ is the inverse limit (cf. Definition \ref{group_inv_lim_defn}) then the natural action $\widehat G\times \overline M\to \overline M$ induces the action $\widehat G\times L^2\left(\overline M, \overline \mu\right)\to L^2\left(\overline M, \overline \mu\right)$ which satisfies to \eqref{o*_dinv_eqn}. Similarly to the equation \ref{top_ddound_eqn} the bounded part of $D^*\left( M_\la, T_\la\right)$ equals to $\Coo\left(M_\la \right) \cap C_b\left(M_\la \right)$ where we natural inclusion $\Coo \left(M_\la \right)\subset  D\left( M_\la\right)$ is implied. The $C^*$-norm completion of  $\Coo\left(M_\la \right) \cap C_b\left(M_\la \right)$ equals to $C_b\left(M_\la \right)$ and $B_\la \bydef C_0\left(M_\la \right)$ is an essential ideal of $C_b\left(M_\la \right)$. Similarly to  \eqref{bl_cat_eqn}
one has  an algebraical  finite covering category
\be\label{top_bl_cat_eqn}
\begin{split}
	\mathfrak{S}_b\bydef\left\{ \left( B ,  B_{\la} , G\left(\left.M_\la~\right|M \right) , \overline\pi_{\la}\right) \right\}_{\la \in \La}= \\ =\left\{ \left(C_0\left(M \right)  ,  C_0\left( M_{\la}\right)  , G\left(\left.M_\la~\right|M \right) , \overline\pi_{\la}\right) \right\}_{\la \in \La}.
\end{split}
\ee
If $\widehat B\bydef C^*$-$\varinjlim B_\la$ then from the Exercise \ref{top_pi_mu_exer} it turns out that there is the natural given by \eqref{top_pi_mu_eqn} faithful, nondegenerate, equivariant representation $\pi_\tau: \widehat{B}\hookto B\left( L^2\left(\overline M, \overline \mu\right)\right)$ (cf. Definitions \ref{faithful_representation_defn}, \ref{nondegenerate_repr_defn} and \ref{equivariant_representation_defn}). The triple $\left(C_0\left(M\right), C_0\left(\widetilde M\right),G\left(\left.\widetilde M~\right| M\right)\right)$ is  the  infinite noncommutative covering of $\mathfrak{S}_b$ (cf. Definition \ref{infinite_noncommutative_covering_defn}. If $\widetilde \H$ is the norm completion of $C_0\left(\widetilde M\right)L^2\left(\overline M, \overline \mu\right)$ then one has
\be
\begin{split}
	\widetilde \H \cong L^2\left( \widetilde M, \widetilde\mu\bydef \overline\mu|_{\widetilde M}\right) ,\\
	\widetilde \D\bydef \widehat\D \cap \widetilde \H \cong \Coo_c\left( \widetilde M\right) 
\end{split}
\ee
Above construction is a version of the described in \ref{comp_pt_o*_empt} general construction.
\begin{theorem}\label{top_inf_o*_thm}
	In the described above situation if $\pi$ is given by \eqref{top_pi_eqn} then the triple
	$\left(D^*\left(M , T\right) , D^*\left( \widetilde{M}, \widetilde T\right) , G\left(\left.\widetilde{M}\right|M\right)\right) $ is the 
	${\pi}$- inverse noncommutative limit (cf. Definition \ref{inv_o*_lim_defn}) of the given by
	\eqref{top_comp_pt_o*_eqn} {algebraic finite covering category of $O^*$-algebras} $\mathfrak{S}_{D^*\left( M, T\right) }$  (cf. Definition \ref{comp_o*_defn}) 
\end{theorem}
\begin{proof}
	Let $\widetilde A' \in \L^\dagger\left(\Coo_c\left(\widetilde M \right)\cong  \Ga^\infty_c\left(\widetilde M, \widetilde T \right) \right)$ be an admissible  algebra (cf. Definition \ref{comp_pt_o*_empt}). Let $\widetilde a'\in \widetilde A'$ be any element. If $\widetilde a'$ is not a differential operator then there is $\widetilde \xi \in \Coo_c\left(\widetilde M \right)$ such that $\supp \widetilde a\widetilde\xi \subsetneqq \supp \widetilde\xi$ (cf. Definition \ref{do_man_defn}). Both  $\supp \widetilde a\widetilde\xi$ and  $\supp \widetilde\xi$ a closures of open sets, so there is an open set $\widetilde \sU\subset \supp \widetilde a\widetilde\xi$ such that $\widetilde\sU \cap \supp \widetilde \xi = \emptyset$. If $\widetilde x_0 \in \widetilde\sU$ and $\widetilde b$ is an $\left(\Coo_0\left(\widetilde M\right) , \widetilde x_0\right)$-{stump} (cf. Definition \ref{top_stump_soft_defn}) such that $\supp \widetilde b\subset \widetilde \sU$ then 
	$$
	\widetilde b \in K\left(C_0\left(\widetilde M \right)  \right) \cap \L^\dagger\left(C_0\left(\widetilde M \right)  \right)_b.
	$$
	(cf. equation \eqref{o*b_eqn}).
	From the condition \ref{inf_cov_o*_sub3_eqn} it turns out that
	\be\label{top_o_1_eqn}
	\exists \la \in \La \quad \exists a_\la \in D^*\left( M_\la, T_\la\right)\quad \widetilde b \widetilde a = \widetilde b \lift_{\widetilde p_\la}  \left( a_\la\right)
	\ee
	where $\lift_{\widetilde p_\la}$ is given by the Definition \ref{top_sheaf_lift_defn}. From the Lemma \ref{top_diff_*alg_lift_eqn} it follows that $\lift_{\widetilde p_\la}$ is a differential operator so $\supp \left( \lift_{\widetilde p_\la}\left( \widetilde a\right) \widetilde\xi\right)  \subset \supp \widetilde \xi$ and
	\be\label{top_o_2_eqn}
	\widetilde b \lift_{\widetilde p_\la}\left(a_\la\right) \widetilde\xi = 0.
	\ee
	On the other hand one has
	\be\label{top_o_3_eqn}
	\widetilde b \widetilde a \widetilde \xi \neq 0.
	\ee
	There is a contradiction between an equation \eqref{top_o_2_eqn} and \eqref{top_o_3_eqn} one. From this contradiction it turns out that any $\widetilde a\in \widetilde A'$ is a differential operator.
	We leave to the reader a proof of following statements:
	\begin{itemize}
		\item If $\widetilde a\in \widetilde A'$ then $\widetilde a\in  D^*\left( \widetilde{M}, \widetilde T\right)$.
		\item *-algebra $D^*\left( \widetilde{M}, \widetilde T\right)$ is admissible, i.e. it satisfies to conditions \eqref{inf_cov_o*_sub1_eqn}- \eqref{inf_cov_o*_sub3_eqn}
	\end{itemize}
\end{proof}

Similarly to \eqref{top_comp_pt_o*_eqn} one has 
{algebraic finite covering category of $O^*$-algebras} (cf. Definition \ref{comp_o*_defn}) given by
\be\label{top_comps_pt_o*_eqn}
\begin{split}
	\mathfrak{S}_{\Coo\left(M \right) } = \left(\left\{\pi_\la: \Coo\left( M\right)  \hookto \Coo\left( M_\la \right) \right\}, \left\{\pi^\mu_\nu\right\}\right).
\end{split}
\ee
Similarly to the Theorem \ref{top_inf_o*_thm}
one can proof the following Theorem
\begin{theorem}\label{top_inf_o*_m_thm}
	In the described above situation the triple\\
	$\left(\Coo\left(M \right) ,\Coo\left( \widetilde{M}\right) , G\left(\left.\widetilde{M}\right|M\right)\right) $ is the 
	${\pi}|_{\bigcup_{\la \in \La}\Coo\left( M_\la \right)}$-\textit{ inverse noncommutative limit of $\mathfrak{S}_{\Coo\left(M \right) }$} an {algebraic finite covering category of $O^*$-algebras} (cf. Definition \ref{comp_o*_defn}) given by
	\eqref{top_comps_pt_o*_eqn}
\end{theorem}
The proof of the Theorem \ref{top_inf_o*_m_thm} is left  to the reader.
\begin{remark}
	The Theorem \ref{top_inf_o*_m_thm} is a noncommutative counterpart of the Proposition \ref{top_cov_mani_prop}. 
\end{remark}
\subsubsection{Infinite coverings and unbounded operators on Hilbert modules}
\paragraph{}
Here we consider two  specializations of described in the Section \ref{unb_hilb_sec} theory.

\paragraph{Algebras of unbounded functions}
Let $p : \widetilde \sX \to \sX$ is the {covering inverse limit of} a category 
$$
\mathfrak{S}_p \bydef \left\{\left\{\sX_\la\right\}_{\la \in \La}, \left\{p^\mu_\nu:\sX_\mu\to \sX_\nu\right\}_{\substack{\mu,\nu \in \La\\\mu\ge\nu}}\right\}.
$$
(cf. Definition \ref{top_disconnected_defn}), and let \be\nonumber
\begin{split}
	\mathfrak{S}_{C_0(p)} \bydef \left(\left\{\rho_\la: B \hookto B_\la \right\}_{\la \in \La}, \left\{\rho^\mu_\nu: B_\mu \hookto B_\nu\right\}_{\substack{\mu, \nu \in \La\\ \nu > \mu}}\right)\bydef \\ \left(\left\{C_0\left( p_\la\right) : C_0\left(\sX \right)  \hookto C_0\left(\sX_\la \right)\right\}, \left\{C_0\left( p_\nu^\mu\right) : C_0\left(\sX_\nu \right) \hookto C_0\left(\sX_\mu \right)\right\}\right).
\end{split}
\ee
be an algebraical  finite covering category (cf. Definition \ref{algebraical_finite_covering_category_defn}), given by \eqref{top_x_g_eqn}
From the Theorem \ref{top_main_thm} it follows that  a triple 
$\left(C_0\left(\sX\right), C_0\left(\widetilde \sX\right),G\left(\left.\widetilde \sX~\right| \mathcal X\right)\right)
$ is  an  {infinite noncommutative covering} of $\mathfrak{S}_{C_0\left(p\right)}$ (cf. Definition \ref{infinite_noncommutative_covering_defn}). This triple we also denote by $\left(B, \widetilde B,G\left(\left.\widetilde B~\right| B\right)\right)
$.  For all $\la\in \La$ denote by $A_\la\bydef Cont\left(\sX_\la \right) $ a *-algebra of continuous (possibly unbounded) maps from $\sX_\la$ to $\C$.  From the Theorem \ref{top_fin_hilb_thm} it turns out that	
\bean
\begin{split}
	\mathfrak{S}_A \bydef \left(\left\{\pi_\la: A \hookto A_\la \right\}_{\la \in \La}, \left\{\pi^\mu_\nu: A_\mu \hookto A_\nu\right\}_{\substack{\mu, \nu \in \La\\ \nu > \mu}}\right)
\end{split}
\eean
is a category of $*$-algebras such that for any $\la, \mu, \nu\in \La$ the triples\\ $\left(A, A_\la, G\left(\left. A_\la\right| A\right)\bydef G\left(\left. B_\la\right| B\right), \pi_\la \right)$ and $\left(A_\mu, A_\nu, G\left(\left. A_\mu\right| A_\nu\right)\bydef G\left(\left. B_\mu\right| B_\nu\right), \pi^\mu_\nu \right)$ are associated with $\left(B, B_\la, G\left(\left. B_\la\right| B\right), \rho_\la \right)$ and $\left(B_\mu, B_\nu,G\left(\left. B_\mu\right| B_\nu\right), \rho^\mu_\nu \right)$ noncommutative finite-fold covering of *-algebras (cf. Definition \ref{fin_chull_defn}).  The  Pedersen's ideal $\widetilde{\mathfrak B}\bydef K\left( \widetilde B\right)= C_c\left(\widetilde \sX \right) $ (cf. Definition \ref{pedersen_ideal_defn} and equation \eqref{peder_c0_eqn}) is dense in ${\widetilde B}$.

\begin{exercise}\label{top_hilb_adm_exer}
	Prove following statements:
	\begin{itemize}
		\item The right  $\widetilde B$-module $ \widetilde{\mathfrak B}'$ satisfies to all conditions described in \ref{comp_chull_empt}.
		\item If $\widetilde A\bydef Cont\left(\widetilde\sX\right)$ is a commutative *-algebra of (possibly unbounded) $\C$-valued continuous maps from  $\widetilde{ {\sX}}$ to $\C$ then $\widetilde A$ is $G\left(\left.\widetilde B\right| B\right)$-equivariant. Moreover there is the natural inclusion $\widetilde A \subset \End^*_{\widetilde { B}}\left(\widetilde {\mathfrak B}'\right)$
		which satisfies to the equations \eqref{hilb_act_eqn} and \eqref{hilb_adm_eqn}, i.e. $\widetilde A$ is admissible.
	\end{itemize}
\end{exercise}
\begin{lemma}
	In the described above situation 
	the triple $\left(A, \widetilde{A}, G\left(\left.\widetilde A \right| A\right)\cong G\left(\left.\widetilde B \right| B\right)\right)$ is a {unique infinite noncommutative covering of $\mathfrak{S}_A$ associated with $\mathfrak{S}_B$} (cf. Definition \ref{main_ch_defn}).
\end{lemma}
\begin{proof}
	If $\widetilde A'$ is admissible (cf. \ref{comp_chull_empt}) and  $\widetilde a'\in \widetilde A'$. then from the equation \eqref{hilb_adm_eqn} it follows that
	\be\label{top_hilb_adm_eqn}
	\forall \widetilde b \in C_c\left(\widetilde \sX \right)\quad \exists\la\in\La\quad \exists a_\la \in C_0\left(\sX_\la \right) \quad \widetilde b \widetilde a' =\widetilde b \lift_{\widetilde p_\la} \left( a_\la \right)  
	\ee
	where 	$\lift_{\widetilde p_\la}$ is the ${\widetilde p_\la}$-lift (cf. Definition \ref{top_lift_defn}). Let both $\mathscr S^{C_0\left(\sX \right) }$ and $\mathscr S^{C_0\left(\widetilde \sX \right) }$ be the ${C_0\left(\sX \right) }$-sheaf and the ${C_0\left(\widetilde \sX \right) }$-sheaf respectively (cf. Definition \ref{top_x_sheaf_defn}). Let $\left\{\widetilde \sU_\a\subset \widetilde\sX\right\}$ be a family of open subsets with compact closure  such that $\widetilde\sX = \cup \widetilde\sU_\a$. For any $\sU_\a$ there is a $\left(C_0\left(\widetilde \sX\right) , \widetilde \sY_\a\right)$-{stump} $\widetilde b_\la$ (cf. Remark \ref{top_stump_soft_y_rem})  where $\widetilde \sY_\a$ is the closure of $\widetilde \sU_\a$. From \eqref{top_hilb_adm_eqn} one has
	$$
	\exists\la\in\La\quad \exists a_\la \in C_0\left(\sX_\la \right) \quad \widetilde b_\a\widetilde a'  =   \widetilde b_\a\lift_{\widetilde p_\la} \left( a_\la \right)
	$$
	and taking into account $\lift_{\widetilde p_\la} \left( a_\la \right)  \widetilde b_\a\in C_0\left(\widetilde\sX\right)$ we conclude that 
	$\widetilde a' \widetilde b_\a\in C_0\left(\widetilde\sX\right)$. 
	The restriction $\left. \widetilde a' \widetilde b_\a\right|_{\widetilde \sU_\a}$ yields  a section 
	$\widetilde s_\a \in \mathscr S^{C_0\left(\widetilde \sX \right) }\left(\widetilde\sU_\a \right)$. From (4) of the Definition \ref{sheaf_defn} it follows that a family $\left\{\widetilde s_\a\right\}$ of local sections  yields a global section $\widetilde s \in \mathscr S^{C_0\left(\widetilde \sX \right) }\left(\widetilde\sX \right)$ which correspond to element of $\widetilde f \in Cont\left( \widetilde \sX\right)$. On the other hand there is a natural inclusion $Cont\left( \widetilde \sX\right)\subset   \End^*_{\widetilde { B}}\left(\widetilde {\mathfrak B}'\right)$. We leave to the user the proof of the equality
	$$
	\forall \widetilde b' \in \widetilde {\mathfrak B}'\quad \widetilde b'\left( \widetilde f - \widetilde a'\right) = 0,
	$$
	or equivalently
	$$
	\widetilde f - \widetilde a'=0\in \End^*_{\widetilde { B}}\left(\widetilde {\mathfrak B}'\right)
	$$
	It follows that $\widetilde A \subset Cont\left(\widetilde\sX\right)$. We leave to the reader the proof of that  $Cont\left(\widetilde\sX\right)$ satisfies to all conditions of the Definition \ref{main_ch_defn}.
\end{proof}

\paragraph{Polynomials in one variable}
Let $p : \widetilde \sX \to \sX$ be the {covering inverse limit of} a category 
$$
\mathfrak{S}_p \bydef \left\{\left\{\sX_\la\right\}_{\la \in \La}, \left\{p^\mu_\nu:\sX_\mu\to \sX_\nu\right\}_{\substack{\mu,\nu \in \La\\\mu\ge\nu}}\right\}.
$$
(cf. Definition \ref{top_disconnected_defn}), and let 
\be\nonumber
\begin{split}
	\mathfrak{S}_{C_0(p)} \bydef \left(\left\{\rho_\la: B \hookto B_\la \right\}_{\la \in \La}, \left\{\rho^\mu_\nu: B_\mu \hookto B_\nu\right\}_{\substack{\mu, \nu \in \La\\ \nu > \mu}}\right)\bydef \\ \left(\left\{C_0\left( p_\la\right) : C_0\left(\sX \right)  \hookto C_0\left(\sX_\la \right)\right\}, \left\{C_0\left( p_\nu^\mu\right) : C_0\left(\sX_\nu \right) \hookto C_0\left(\sX_\mu \right)\right\}\right).
\end{split}
\ee
If $\widetilde{\sY}\bydef \widetilde \sY \times \R$, 
 $~\sY \bydef  \sX \times \R $, $~ \sY_\la \bydef \times \R $ then $\widetilde p : \widetilde \sY \to \sY$ be the {covering inverse limit of} a category 
 $$
 \mathfrak{S}_{\widetilde p} \bydef \left\{\left\{\sY_\la\right\}_{\la \in \La}, \left\{\widetilde p^\mu_\nu:\sX_\mu\to \sX_\nu\right\}_{\substack{\mu,\nu \in \La\\\mu\ge\nu}}\right\}.
 $$
 (cf. Definition \ref{top_disconnected_defn}). Let
 \be\nonumber
 \begin{split}
 	\mathfrak{S}_{C_0\left(\widetilde p \right) } \bydef \left(\left\{\rho_\la: B \hookto B_\la \right\}_{\la \in \La}, \left\{\rho^\mu_\nu: B_\mu \hookto B_\nu\right\}_{\substack{\mu, \nu \in \La\\ \nu > \mu}}\right)\bydef \\ \left(\left\{C_0\left(\widetilde p_\la\right) : C_0\left(\sX \right)  \hookto C_0\left(\sX_\la \right)\right\}, \left\{C_0\left(\widetilde p_\nu^\mu\right) : C_0\left(\sX_\nu \right) \hookto C_0\left(\sX_\mu \right)\right\}\right).
 \end{split}
 \ee
 be a given by \eqref{top_x_g_eqn} algebraical  finite covering category (cf. Definition \ref{algebraical_finite_covering_category_defn})
\be\label{top_sb_eqn}
\begin{split}
	\mathfrak{S}_B \bydef \left(\left\{\rho_\la: B \hookto B_\la \right\}_{\la \in \La}, \left\{\rho^\mu_\nu: B_\mu \hookto B_\nu\right\}_{\substack{\mu, \nu \in \La\\ \nu > \mu}}\right).
\end{split}
\ee
From the Theorem \ref{top_main_thm} it follows that $\mathfrak{S}_{C_0\left(\widetilde p \right) }$ is good  (cf. Definition \ref{good_defn}). Moreover and there is a space $\widetilde \sY$  and a transitive covering $\widetilde p: \widetilde \sY\to  \sY$ such that a triple 
$\left(C_0\left(\sY\right), C_0\left(\widetilde \sY\right),G\left(\left.\widetilde \sY~\right| \mathcal Y\right)\right)
$ is  an  { infinite noncommutative covering} of $\mathfrak{S}_{B}$ (cf. Definition \ref{infinite_noncommutative_covering_defn}). This triple we also denote by $\left(B, \widetilde B,G\left(\left.\widetilde B~\right| B\right)\right)
$.
If  $A$ be a given by \eqref{top_poly_a_defn} *-algebra, and  If $\mathfrak B$ is the completion of $C_c\left(\sY\right)$ with respect to the graph topology (cf. Definition \ref{def:rep_Hilbert_module} ) then there is the representation  $\mu:  A \hookto \End_B\left(\mathfrak B \right)$.

\begin{exercise}\label{top_ch_exer}
	\begin{enumerate}
		\item Using the Theorem \ref{top_fin_hilb_thm}  construct an  {algebraic finite covering category of $*$-algebras}  given by
		\be\label{top_sa_eqn}
		\begin{split}
			\mathfrak{S}_A \bydef \left(\left\{\pi_\la: A \hookto A_\la \right\}_{\la \in \La}, \left\{\pi^\mu_\nu: A_\mu \hookto A_\nu\right\}_{\substack{\mu, \nu \in \La\\ \nu > \mu}}\right)
		\end{split}
		\ee
		such that for any $\la, \mu, \nu\in \La$ the triples $\left(A, A_\la, G\left(\left. A_\la\right| A\right)\bydef G\left(\left. B_\la\right| B\right), \pi_\la \right)$ and $\left(A_\mu, A_\nu, G\left(\left. A_\mu\right| A_\nu\right)\bydef G\left(\left. B_\mu\right| B_\nu\right), \pi^\mu_\nu \right)$ are associated with $\left(B, B_\la, G\left(\left. B_\la\right| B\right), \rho_\la \right)$ and $\left(B_\mu, B_\nu,G\left(\left. B_\mu\right| B_\nu\right), \rho^\mu_\nu \right)$ noncommutative finite-fold covering of *-algebras (cf. Definition \ref{fin_chull_defn}.
		\item Prove that 	in the described above situation 
		the triple $\left(A, \widetilde{A}, G\left(\left.\widetilde A \right| A\right)\cong G\left(\left.\widetilde B \right| B\right)\right)$ is an {infinite noncommutative covering of  $\mathfrak{S}_A$  associated with $\mathfrak{S}_B$} (cf. Definition \ref{main_ch_defn}).
		
	\end{enumerate}
	
\end{exercise}

\begin{empt}\label{top_at_empt}
	Consider the above situation. One can define *-algebra $tA$ and $B^T$ as well as it is explained in the section \ref{top_fin_hull_sec}. Similarly to  \eqref{top_sa_eqn} and \eqref{top_sb_eqn} one can define
	\bean
	\begin{split}
		\mathfrak{S}_{B^T} = \left(\left\{\rho_\la: B^T \hookto B^T_\la \right\}_{\la \in \La}, \left\{\rho^\mu_\nu: B^T_\mu \hookto B^T_\nu\right\}_{\substack{\mu, \nu \in \La\\ \nu > \mu}}\right),\\
		\mathfrak{S}_{tA} \bydef \left(\left\{\pi_\la: tA \hookto tA_\la \right\}_{\la \in \La}, \left\{\pi^\mu_\nu: tA_\mu \hookto tA_\nu\right\}_{\substack{\mu, \nu \in \La\\ \nu > \mu}}\right).
	\end{split}
	\eean
	
\end{empt}
\begin{exercise}\label{top_chull_exer}
	Prove that in the described above situation 
	the triple
	\\$\left(tA, t\widetilde{A}, G\left(\left.t\widetilde A \right| tA\right)\cong G\left(\left.\widetilde B^T \right| B^T\right)\right)$ is an {infinite noncommutative covering of  $\mathfrak{S}_{tA}$  associated with $\mathfrak{S}_{B^T}$} (cf. Definition \ref{main_ch_defn}).
\end{exercise}
\begin{remark}
	The given by the Exercise \ref{top_ch_exer} inverse noncommutative limit can be regarded as an inverse noncommutative limit of $O^*$-algebras (cf. Definition \ref{inv_o*_lim_defn}). 	The given by the Exercise \ref{top_chull_exer} inverse noncommutative limit can not be regarded as an inverse noncommutative limit of $O^*$-algebras because $tA$ does not have nontrivial bounded elements.
\end{remark}

\subsection{Coverings of   quasi $*$-algebras}
	\paragraph{}
	Let $M$ be a finite-dimensional smooth manifold, and let $p: \widetilde M \to M$ be a transitive covering such that the covering group $G\left(\left. {\widetilde{M}}~\right|~M\right)$ (cf. Definition  \ref{top_group_of_covering_transformations_defn} is residually finite (cf. Definition \ref{residually_finite_defn})
	Let 
	\be\label{top_secmm_p_eqn}
	\begin{split}
		\mathfrak{S}_p \bydef \left\{\left\{M_\la\right\}_{\la \in \La}, \left\{p^\mu_\nu:M_\mu\to M_\nu\right\}_{\substack{\mu,\nu \in \La\\\mu\ge\nu}}\right\}.
	\end{split}
	\ee
	be  a finite covering category of  $p: \widetilde M \to M$ (cf. Definition \ref{top_disconnected_defn}).
From the Corollary \ref{top_q_fin_lem} it turns out that there  is a
a family of {noncommutative finite-fold coverings of quasi $*$-algebras} (cf. Remark \ref{top_distr_dens_q_rem} and Definition \ref{oq*fin_defn}) given by
\be\label{top_comp_pt_qo*_eqn}
\begin{split}
 \mathfrak{S}_{\left(\Coo_c\left(M \right)', \Coo\left(M \right)    \right) } =\\=\left(\left\{ \left(\Coo_c\left(M_\la \right)', \Coo\left(M_\la \right)    \right) \right\}_{\la \in \La}, \left\{\pi^\mu_\nu: \right\}_{\substack{\mu, \nu \in \La\\ \nu > \mu}}\right).
\end{split}
\ee

\begin{lemma}
The given by \eqref{top_comp_pt_qo*_eqn} family $ \mathfrak{S}_{\left(\Coo_c\left(M \right)', \Coo\left(M \right)    \right) }$  is a \textit{algebraic finite covering category of quasi-$*$-algebras} (cf. Definition \ref{comp_qo*_defn}).
\end{lemma}
\begin{proof}
	In this lemma the pair $\left( \mathfrak A, \mathfrak A_0 \right)$ (resp. $\left( \mathfrak A^\la, \mathfrak A^\la_0 \right)$) of the Definition \ref{comp_qo*_defn}  corresponds to $\left(\Coo_c\left(M \right)', \Coo\left(M \right)\right) $ (resp. $\left(\Coo_c\left(M_\la \right)', \Coo\left(M_\la \right)\right)$). The given by \eqref{comp_pt_qo*_eqn} {algebraic finite covering category of quasi-$*$-algebras} corresponds to  $\mathfrak{S}_{\left(\Coo_c\left(M \right)', \Coo\left(M \right)    \right) }$.
Now one needs check that the family $ \mathfrak{S}_{\left(\Coo_c\left(M \right)', \Coo\left(M \right)    \right) }$ satisfies to   the Definition \ref{comp_qo*_defn}. However it directly follows from the given by \eqref{top_comp_pt_qo*_eqn} definition of $
	\mathfrak{S}_{\left(\Coo_c\left(M \right)', \Coo\left(M \right)    \right) }$.
\end{proof}	
\begin{empt}\label{top_comp_qo*_empt}
	Here we consider a relevant to $ \mathfrak{S}_{\left(\Coo_c\left(M \right)', \Coo\left(M \right)    \right) }$ specialization of the described in \ref{comp_qo*_empt} construction.		
	
	The following table gives a mapping between notations of the  \ref{comp_qo*_empt} and this specialization.
	\newline
	\begin {table}[H]
	\caption {The mapping between notations of  \ref{comp_qo*_empt} and this specialization} \label{top_m_table}
	\begin{tabular}{|c|c|c|}
		\hline
		& \ref{comp_qo*_empt} & This specialization\\
		\hline
		&	&\\
		1& $\left( \mathfrak A, \mathfrak A_0 \right) $  &  $\left(\Coo_c\left(M \right)', \Coo\left(M \right)\right) $ \\ & & \\
		2& $\left( \mathfrak A^\la, \mathfrak A^\la_0 \right) $  &  $\left(\Coo_c\left(M_\la \right)', \Coo\left(M_\la \right)\right) $ \\ & & \\
				3& $\widehat \H$  &  $L^2\left(\overline M, \overline \mu\right)$ (cf. \ref{top_mu_empt}) \\ & & \\
				4& $\widehat D$  &  $\Coo_c\left(\overline M \right)  $ (cf. \eqref{top_tdd_eqn}) \\ & & \\
				5& $	\pi: \widehat{\mathfrak A}_0\hookto \L^\dagger\left(\widehat\D \right)$ (cf. \eqref{pi_o*f_eqn})  &  $\pi: \bigcup_{\la\in\La}\Coo_c\left( M_\la \right)\hookto \L^\dagger\left(\Coo_c\left(\overline M \right) \right) $  \\ & & \\
				6& $\mathfrak{S}_0$ (cf. \eqref{comp_pt_qo*0_eqn})  &  $\mathfrak{S}_{\Coo\left(M \right) } $ (cf. \eqref{top_comps_pt_o*_eqn}) \\ & & \\
				7& $\left\{ \left( B ,  B_{\la} , G\left(\left.\mathfrak A_0^\la~\right|\mathfrak A_0 \right) , \overline \pi_{\la}|_{\mathfrak A_0}\right) \right\}_{\la \in \La}$   &  $	\mathfrak{S}_b$ (cf. \eqref{top_bl_cat_eqn}) \\ & (cf. \eqref{blu_cat_eqn}) & \\
				 & & \\
							8& $\left(B, \widetilde{B}_{\pi_\tau}, G_{{\pi_\tau}}\right) $ (cf. \eqref{be_eqn})  &$\left(C_0\left(M \right),  C_0\left(\widetilde M \right), G\left(\left. \widetilde M \right| M \right)  \right)$  \\ & & \\
										\hline
	\end{tabular}
\end{table}
Using the Table \ref{top_m_table} one can proof that described here situation is a specialization of \ref{comp_qo*_empt}.
\end{empt}

\begin{theorem}\label{top_inf_q*_thm}
	In situation \ref{top_comp_qo*_empt} the triple
	$$
\left(\left(\Coo_c\left(M \right)', \Coo\left(M \right)    \right), \left(\Coo_c\left(\widetilde M \right)', \Coo\left( \widetilde M \right)    \right), G\left(\left. \widetilde M\right| M\right)\right)
$$
 is the 
	${\pi}$- inverse noncommutative limit  (cf. Definition \ref{inv_qo*_lim_defn}) of  the given by \eqref{top_comp_pt_qo*_eqn} {algebraic finite covering category of quasi $*$-algebras} $\mathfrak{S}_{\left(\Coo_c\left(M \right)', \Coo\left(M \right)    \right)  }$.
\end{theorem}
\begin{proof}

	Firstly note that the specific versions of the equations \eqref{inv_qo*_rzero_eqn}-\eqref{inv_qo*_rlim_eqn} are presented below
	\bea
	\label{top_inv_qo*_rzero_eqn}
	\forall\widetilde a\in \widetilde{\mathfrak A}'\quad \Coo_c\left(  \widetilde{M}\right) \widetilde{a}= \{0\}\quad \Rightarrow\quad \widetilde{a}= 0,~\quad\quad \\
	\label{top_inv_qo*_lzero_eqn}
	\forall\widetilde a\in  \widetilde{\mathfrak A}'\quad \widetilde{a}\Coo_c\left(  \widetilde{M}\right) = \{0\}\quad \Rightarrow\quad \widetilde{a}= 0,~ \quad\quad\\
	\label{top_inv_qo*_rlim_eqn}\forall  \widetilde b \in \Coo_c\left(\widetilde M \right) ~\forall \widetilde a' \in \widetilde{\mathfrak A}' ~\exists a',a'' \in    \bigcup_{\la\in \La} \Coo_c\left( M_\la \right) ~~ \widetilde b\widetilde a=\widetilde b a', ~~  \widetilde a\widetilde b=a''\widetilde b. \quad\quad
	\eea
	Let $\left\{\widetilde f_\a\right\}_{\a\in \mathscr A}\subset C_c\left(\widetilde{M} \right)_+\cap \Coo\left(M \right)$ be a net of all positive, smooth, compactly supported maps such that 
\bean
	\forall \a \in \mathscr A\quad\widetilde f_\a \left(\widetilde{M}\right)\subset [0,1],\\
	\forall \bt, \ga \in \mathscr A\quad \bt\le\ga \quad \Leftrightarrow \quad \widetilde f_\bt \le \widetilde f_\ga.
\eean 
From the Exercise \ref{smooth_soft_exer} it follows that $\Coo\left(\widetilde M \right)\cap C_0\left( \widetilde M\right)$ is a $c$-soft *-algebra. Taking into account the Exercise \ref{top_stump_soft_y_exer} we conclude that for any compact subset $\widetilde \sY$ of $\widetilde M$ there is $\a_{\widetilde\sY} \in  \mathscr A$ such that $f_{\a}\left({\widetilde\sY}\right) = \{1\}$ for all $\a \ge\a_{\widetilde\sY}$.  If $\left( \widetilde{\mathfrak A}', \Coo\left( \widetilde{M}\right)\right) $ is an admissible quasi *-algebra (cf. Definition \ref{comp_pt_qo*_defn}) then from the equation \ref{top_inv_qo*_rlim_eqn}  it follows that
\bean
\forall \widetilde a\in \widetilde{\mathfrak A}'\quad \forall \a\in \mathscr A\quad	\widetilde f_\a\widetilde a \in \Coo_c\left( \widetilde{M}\right)\left( \bigcup_{\la\in\La} \Coo_c\left( {M}_\la\right)'\right).
\eean
	On the other hand from $\bigcup_{\la\in\La} \Coo_c\left( {M}_\la\right)'\subset \Coo_c\left( \widetilde{M}\right)'$ it turns out that $\widetilde f_\a \widetilde a\in \Coo_c\left( \widetilde{M}\right)'$.  For any  $\widetilde b \in \Coo_c\left( \widetilde{M}\right)$ there is $\a_{\supp \widetilde b}$ such that $f_{\a_{\supp \widetilde b}}\left( \supp \widetilde b\right)= \{1\}$. One has
	\be\label{top_sw_lim_eqn}
\a\ge \a_{\supp b}\quad \Rightarrow \quad \widetilde f_{\a} \ge \widetilde f_{\a_{\supp \widetilde b}}  \quad \Rightarrow \quad	\left\langle \widetilde f_\a\widetilde a , \widetilde b \right\rangle = 	\left\langle \widetilde f_{\a_{\supp \widetilde b}} \widetilde a , \widetilde b \right\rangle
\ee
where $\left\langle \cdot , \cdot \right\rangle: \Coo_c\left( \widetilde{M}\right)'\times \Coo_c\left( \widetilde{M}\right)\to \C$ is the natural pairing.
The equation \eqref{top_sw_lim_eqn} means that the net $ \left\{ \left\langle \widetilde f_\a\widetilde a , \widetilde b \right\rangle\right\} \subset\C$ is convergent, so the net  $ \left\{ \widetilde f_\a \widetilde a\right\} \subset\Coo_c\left( \widetilde{M}\right)'$ is convergent with respect to the weak topology of $\Coo_c\left( \widetilde{M}\right)'$ (cf. Definition \ref{w_topology_defn}). Thus there is a map
\bean
\Psi : \widetilde{\mathfrak A}'\to \Coo_c\left( \widetilde{M}\right)',\\
\widetilde a\mapsto\lim_{\a\in \mathscr A} \widetilde f_\a \widetilde a
\eean 
 If $\widetilde a\in \ker \Psi$ then $\Coo_c\left( \widetilde{M}\right)\widetilde a= \{0\}$ so from \eqref{inv_qo*_rzero_eqn} (or equivalently \eqref{top_inv_qo*_rzero_eqn}) it follows that $\widetilde a=0$. So one has $\ker \Psi = \{0\}$, $\Psi$ is an invective map and any admissible *- algebra $\left( \widetilde{\mathfrak A}', \Coo\left( \widetilde{M}\right)\right) $ is a quasi *-subalgebra  of 
$\left( \Coo\left( \widetilde{M}\right)', \Coo\left( \widetilde{M}\right)\right)$.  Thus if we prove that  $\left( \Coo_c\left( \widetilde{M}\right)', \Coo\left( \widetilde{M}\right)\right)$ is admissible then the proof of this lemma would be completed (cf. Definition \ref{inv_qo*_lim_defn}). It means that $\left( \Coo\left( \widetilde{M}\right)', \Coo\left( \widetilde{M}\right)\right)$ satisfies to the conditions (a)-(c) of the Definition \ref{comp_pt_qo*_defn}.
\begin{enumerate}
	\item[(a)] The action $G\left(\left. \widetilde M\right| M\right)\times \widetilde M\to \widetilde M$ naturally induces the action $G\left(\left. \widetilde M\right| M\right)\times \left( \Coo_c\left( \widetilde{M}\right)', \Coo\left( \widetilde{M}\right)\right)\to \left( \Coo_c\left( \widetilde{M}\right)', \Coo\left( \widetilde{M}\right)\right)$.
	\item[(b)] The existence of injective $*$-homomorphism 
	$$
	\varphi: \bigcup_{\La \in \la} \left( \Coo_c\left( {M}_\la\right)', \Coo\left( {M}_\la\right)\right)\hookto \left( \Coo_c\left( \widetilde{M}\right)', \Coo\left( \widetilde{M}\right)\right) 
	$$ 
	follows from the natural inclusions $ \left( \Coo_c\left( {M}_\la\right)', \Coo\left( {M}_\la\right)\right)\hookto \left( \Coo_c\left( \widetilde{M}\right)', \Coo\left( \widetilde{M}\right)\right)$ for all $\la\in\La$.
	\item[(c)]  The specific versions of equations \eqref{inv_qo*_rzero_eqn} -\eqref{inv_qo*_rlim_eqn} are presented below

		\bea
	\label{top_inv_qo*_rzerop_eqn}
	\forall\widetilde a\in \Coo_c\left(  \widetilde{M}\right)'\quad \Coo_c\left(  \widetilde{M}\right) \widetilde{a}= \{0\}\quad \Rightarrow\quad \widetilde{a}= 0,~\quad\quad \\
	\label{top_inv_qo*_lzerop_eqn}
	\forall\widetilde a\in \Coo_c\left(  \widetilde{M}\right)'\quad \widetilde{a}\Coo_c\left(  \widetilde{M}\right) = \{0\}\quad \Rightarrow\quad \widetilde{a}= 0,~ \quad\quad\\
	\label{top_inv_qo*_rlimp_eqn}\forall  \widetilde b \in \Coo_c\left(\widetilde M \right) ~\forall \widetilde a \in \Coo_c\left(  \widetilde{M}\right)' ~\exists a',a'' \in    \bigcup_{\la\in \La} \Coo_c\left( M_\la \right) ~~ \widetilde b\widetilde a=\widetilde b a', ~~  \widetilde a\widetilde b=a''\widetilde b. ~\quad\quad 
	\eea
	so one needs check equations
 \eqref{top_inv_qo*_rzerop_eqn} -\eqref{top_inv_qo*_rlimp_eqn}.
Suppose that $\widetilde a\in \Coo_c\left(  \widetilde{M}\right)'$ is such that $\Coo_c\left( {M}_\la\right) \widetilde{a}= \{0\}$. For all  $\widetilde b\in \Coo_c\left(  \widetilde{M}\right)$, there is $ \widetilde{c} \in \Coo_c\left(  \widetilde{M}\right)_+$ such that $\widetilde{c}\left( \supp \widetilde{b}\right) = \{1\}$. It follows that 
$$
\left\langle \widetilde a, \widetilde b \right\rangle = \left\langle \widetilde a, \widetilde{c}\widetilde b \right\rangle = \left\langle \widetilde a \widetilde{c},\widetilde b \right\rangle,
$$
so if $\widetilde a \widetilde{c}=0$ then $\left\langle \widetilde a, \widetilde b \right\rangle= 0$. Thus the equation \eqref{top_inv_qo*_rzerop_eqn} is proven. The equation \eqref{top_inv_qo*_lzerop_eqn}  is equivalent to \eqref{top_inv_qo*_rzerop_eqn}. 

Consider the equation  \eqref{top_inv_qo*_rlimp_eqn}. If  $\widetilde c\in \Coo_c\left(\widetilde M \right)$ is a $\left(\Coo_0\left( \widetilde M\right), \supp \widetilde b  \right)$-stump (cf. Remark \ref{top_stump_soft_y_rem}) then $\widetilde b = \widetilde b\widetilde c= \widetilde c\widetilde b$. From the Lemma \ref{top_compact_la_lem}  it follows that there is $\la_0\in \La$ and an open set $\widetilde \sU\subset \widetilde \sX$ such that
\begin{itemize}
	\item $\supp \widetilde s \subset \widetilde \sU$,
	\item $\widetilde \sU$ is mapped homeomorphically onto $\sU_{\la_0}\bydef \widetilde p_{\la_0}\left(\widetilde \sU \right)$ where $\widetilde p_{\la_0}: \widetilde M \to M_{\la_0}$ is the natural covering.
\end{itemize}
If $a'\bydef  \desc_{p_{\la_0}}\left(\widetilde c \widetilde a \right)\in \Coo_c\left(M_{\la_0}\right)'$ where $\desc_{\la_0}$ is a given by the Definition  \ref{top_lift_sh_desc_defn} $p_{\la_0}$-descent then one has
$$
\widetilde b\widetilde a= \widetilde a\widetilde b= \widetilde ba'= a'\widetilde b.
$$
The above equation is equivalent to \eqref{top_inv_qo*_rlimp_eqn} one.
\end{enumerate}

\end{proof}

\subsection{Strong Morita equivalence of infinite coverings}

	\paragraph*{} Here we consider a specialization of the explained in the Section \ref{inf_mor_sec} theory.
\begin{lemma}\label{top_all_lem}
If $\mathscr L^2\left(C_0\left( \widetilde \sX\right)  \right) _{C_0\left({\sX}\right)}$ is the
 {$C^*$-Hilbert module  associated the  infinite noncommutative covering} 
 $$
 \left(C_0\left({\sX}\right), C_0\left({\sX}\right), G\left(\left.\widetilde{\sX}~\right| \sX\right)\right)
 $$ of $\mathfrak{S}$ (cf. Definition \ref{infinite_noncommutative_covering_defn}).
given by the Theorem \ref{top_main_thm}
then  $\mathscr L^2\left(C_0\left( \widetilde \sX\right)  \right)_{C_0\left({\sX}\right)}$ satisfies to the 
conditions (a), (b)  of the Section \ref{inf_mor_sec}.
\end{lemma}
\begin{proof}
(a) Let $\left(\widetilde{a}, g\right) \in C_0\left(\mathcal{X}\right)\times G\left(\widetilde{\sX}~|~ \mathcal X\right)$ and let $h_{\widetilde{a}, g}\in \End^*_{C_0\left(\sX \right) }\left(\mathscr L^2\left(C_0\left( \widetilde \sX\right)  \right) \right)$ is given by \eqref{hm_inc_eqn} i.e. 
\be\label{topo_hm_inc_eqn}
\begin{split}
	h_{\widetilde{a}, g}\left(\xi \right) 
	 = g\left(\xi \widetilde{a} \right) \quad \forall \xi\in \mathscr L^2\left(C_0\left( \widetilde \sX\right)  \right)
\end{split}
\ee
For any $\eps > 0$ the set
$$
\widetilde\sY\bydef \left\{ \widetilde x\in \widetilde\sX \left|\left| \widetilde a\left(\widetilde x\right) \right|\ge \frac{\eps}{2} \right.\right\}
$$
is compact.
If  $ \widetilde p: \widetilde\sX \to \sX$ is the natural covering and  $\sum_{j = 1}^n f_j$  is a covering sum for $\widetilde\sY$ {subordinated  to} $\widetilde p$ (cf. Definition \ref{top_covering_sum_subordinated_defn}). then one has
$$
\left\|  \widetilde a -  \widetilde a\sum_{j = 1}^n f_j\right\| < \eps.
$$
From the above equation it turns out that
\be\label{top_hg_e_eqn}
\left\|  	h_{\widetilde{a}, g}-\sum_{j = 1}^n h_{f_j \widetilde a, g}\right\| < \eps.
\ee
On the other hang $h_{f_j \widetilde a, g}=\theta_{\sqrt{f_j}, g\left( \sqrt{f_j}\widetilde{a}\right)  } = \sqrt{f_j}\left\rangle \right\langle g\left( \sqrt{f_j}\widetilde{a}\right)$ is a rank-one operator (cf. equation \eqref{rank_one_eqn}). 
The number $\eps$ is arbitrary small, so from the equation  \eqref{top_hg_e_eqn} and the Definition \ref{compact_a_operator_defn} one concludes that 	$h_{\widetilde{a}, g}$ is compact, i.e. 
$h_{\widetilde{a}, g}\in \K\left(\mathscr L^2\left(C_0\left( \widetilde \sX\right)  \right)_{C_0\left( \sX\right) }  \right)$. Since $C_c\left(G\left(\widetilde{\sX}~|~ \mathcal X\right), C_0\left( \widetilde{\sX}\right) \right)$ is generated  by operators $h_{\widetilde{a}, g}$ one has $C_c\left(G\left(\widetilde{\sX}~|~ \mathcal X\right), C_0\left( \widetilde{\sX}\right) \right) \subset \K\left(\mathscr L^2\left(C_0\left( \widetilde \sX\right)  \right)_{C_0\left( \sX\right) } \right)$.\\
(b)  If $\widetilde a, \widetilde b \in K\left(C_0\left( \widetilde{\sX}\right)  \right) = C_c\left( \widetilde{\sX}\right)$ then both sets $\supp\widetilde a$ and  $\supp\widetilde b$ are compact. Let  $ \widetilde p: \widetilde\sX \to \sX$ be the natural covering and  $\sum_{j = 1}^n f_j$  is a covering sum for $\supp \widetilde a$ {subordinated  to} $\widetilde p$ (cf. Definition \ref{top_covering_sum_subordinated_defn}).
If $j \in \left\{1,..., n\right\}$ and the set
$$
G_j \bydef \left\{\left. g \in G\left(\widetilde{\sX}~|~ \mathcal X\right) \right| g\supp f_j \cap \supp \widetilde b\neq \emptyset \right\}
$$
is not finite then the compact set $\supp \widetilde b$ contains an infinite disjoint union
$$
\bigsqcup_{g \in G_j} g\supp f_j \cap \supp \widetilde b
$$
of compact subsets. It is impossible, so $G_j$ is  finite for all  $j=1,..., n$.
If $\ka_j \bydef f_j \widetilde a \left\rangle \right\langle \widetilde b\in  \K\left(\mathscr L^2\left(C_0\left( \widetilde \sX\right)  \right)_{C_0\left( \sX\right) }  \right)$ is a rank-one operator (cf. equation \eqref{rank_one_eqn}) then $a \left\rangle \right\langle b = \sum_{j=1}^n \ka_j$. For any $\widetilde c \in \K\left(\mathscr L^2\left(C_0\left( \widetilde \sX\right)  \right)_{C_0\left( {\sX}\right)}  \right)$ one has
$$
\ka_j \widetilde c = \sum_{g \in G\left(\widetilde{\sX}~|~ \mathcal X\right)} g\left(  \widetilde c f_j\widetilde a^*\right) \widetilde b=  \sum_{g \in G_j} g\left(  \widetilde cf_j\widetilde a^*\right) \widetilde b
$$
On the other hand from
$$
g_j\left(  \widetilde cf_j\widetilde a^*\right) \widetilde b = h_{f_j\widetilde a^*g_j^{-1}\widetilde b, g_j}\left(\widetilde c \right) 
$$
it follows that
$$
\ka_j = \sum_{g \in G_j}h_{f_j\widetilde a^*g_j^{-1}\widetilde b, g_j}\in C_c\left(G\left(\widetilde{\sX}~|~ \mathcal X\right), C_0\left( \widetilde{\sX}\right) \right)
$$
and taking into account $\widetilde a \left\rangle \right\langle \widetilde b = \sum_{j=1}^n \ka_j$ one has $\widetilde a \left\rangle \right\langle \widetilde b\in  C_c\left(G\left(\widetilde{\sX}~|~ \mathcal X\right), C_0\left( \widetilde{\sX}\right) \right)$. The $\C$-space  $K\left(C_0\left( \widetilde{\sX}\right)  \right)$  is dense in $\mathscr L^2\left(C_0\left( \widetilde \sX\right)  \right)_{C_0\left( {\sX}\right)}$, so for any \\$\xi, \eta \in \mathscr L^2\left(C_0\left( \widetilde \sX\right)  \right)_{C_0\left( {\sX}\right)}$ and $\eps > 0$ there are $\widetilde a, \widetilde b \in K\left(C_0\left( \widetilde{\sX}\right)  \right)$ such that
$$
\left\|  \widetilde a \left\rangle \right\langle \widetilde b-\xi \left\rangle \right\langle \eta\right\|< \eps.
$$
Taking into account $\widetilde a \left\rangle \right\langle \widetilde b\in  C_c\left(G\left(\widetilde{\sX}~|~ \mathcal X\right), C_0\left( \widetilde{\sX}\right) \right)$ one concludes that the rank-one operator $\xi \left\rangle \right\langle \eta$ lies in the $C^*$-norm completion of $C_c\left(G\left(\widetilde{\sX}~|~ \mathcal X\right), C_0\left( \widetilde{\sX}\right) \right)$.
On the other hand from the Definition \ref{compact_a_operator_defn} it follows that $\K\left(\mathscr L^2\left(C_0\left( \widetilde \sX\right)  \right)_{C_0\left( \sX\right) }  \right)$ is the $C^*$-norm completion of the algebraical linear span of rank-one operators, so\\ $C_c\left(G\left(\widetilde{\sX}~|~ \mathcal X\right), C_0\left( \widetilde{\sX}\right) \right)$ is dense in $\K\left(\mathscr L^2\left(C_0\left( \widetilde \sX\right)  \right)_{C_0\left( \sX\right) }  \right)$.
\end{proof}
\begin{exercise}\label{top_all_exer}
Prove that a linear span of \\ $\left\langle \mathscr L^2\left(C_0\left( \widetilde \sX\right)  \right)_{C_0\left( \sX\right)}, \mathscr L^2\left(C_0\left( \widetilde \sX\right)  \right)_{C_0\left( \sX\right) } \right\rangle_{C_0\left( \widetilde\sX\right)\rtimes_c G\left(\widetilde{\sX}~|~ \mathcal X\right)}$  is dense in $C_0\left( \sX\right)$, i.e.
$$\left(C_0\left(\mathcal{X}\right), C_0\left(\widetilde{\sX}\right),G\left(\widetilde{\sX}~|~ \mathcal X\right)\right)$$ allows strong Morita equivalence  (cf. Definition \ref{allows_morita_defn}).  
\end{exercise}
\begin{remark}
The Exercise \ref{top_all_exer} corresponds to the Proposition \ref{top_allows_morita_prop} (cf. Remark \ref{top_allows_morita_rem}).
\end{remark}

\subsection{Induced representations}\label{comm_induced_infinite_sec}

\paragraph*{} 

Similarly to the section \ref{comm_induced_finite_sec} consider  is a connected,  locally compact, Hausdorff space $\sX$  and a locally trivial bundle  $E \to \sX$ with a {sesquilinear form}
and  is a locally trivial bundle with a {sesquilinear form}  $\varphi: E \times_\sX E \to \C$ (cf. Definitions \ref{top_vb_defn} \ref{top_herm_bundle_form_defn}). Suppose that there is a measure $\mu$ and a given by \eqref{top_l2_xe_eqn} representation 
\be\label{comm_m_bundle_inf_repr_eqn}
\rho: C_0\left(\sX \right) \hookto L^2\left(\sX, E,  \mu\right).
\ee

Let $p: \widetilde \sX \to \sX$ be an infinite transitive covering, such that the group $G\left(\left.\widetilde \sX   ~\right| \sX\right)$ is residually finite. From the Lemma \ref{top_associded_lem}  and the  Theorem \ref{top_main_thm} it turns out that the triple  $\left(C_0\left( \sX\right) , C_0\left( \widetilde \sX\right) G\left(\left.\widetilde \sX   ~\right| \sX\right)\right)$ is an infinite noncommutative covering (cf. Definition \ref{infinite_noncommutative_covering_defn}). Let both $\widetilde\tau \bydef \lift_p\left(\tau \right)$ and $\widetilde \mu$ be the $p$-lift of $\tau$ and $\mu$ respectively (cf. Definition \ref{top_lift_measure_defn}) 
Define the $\C$-valued product
\be
\begin{split}
	\left( \cdot, \cdot \right): \Ga_c\left(\widetilde \sX, \widetilde E\right)\times \Ga_c\left( \widetilde \sX, \widetilde E\right) \to \C,\\
	\left(\widetilde\xi,\widetilde \eta \right) \mapsto \widetilde \tau\left( \left\langle\widetilde\xi, \widetilde\eta \right\rangle_c \right) 
\end{split}
\ee
where the pairing $	\left\langle \cdot, \cdot \right\rangle_c:  \Ga_c\left(\widetilde \sX,\widetilde E\right)\times \Ga_c\left(\widetilde \sX, \widetilde E\right)\to C_c\left(\widetilde \sX \right) $ is given by \eqref{top_ggc_eqn}.
This product yields the structure of pre-Hilbert space on $\Ga_c\left(\widetilde \sX,\widetilde E\right)$. Denote by $ L^2\left(\widetilde  \sX,\widetilde E\right)\cong  L^2\left( \widetilde \sX, \widetilde E, \widetilde \mu\right)$ the Hilbert norm completion of $\Ga_c\left(\widetilde \sX,\widetilde E\right)$. The natural action $C_0\left(\widetilde E\right) \times \Ga_c\left(\widetilde \sX,\widetilde E\right)\to \Ga_c\left(\widetilde \sX,  \widetilde E\right)$ induces  the natural representation
\be\label{comm_bundlei_repr_eqn}
C_0\left(\widetilde \sX\right) \to B\left( L^2\left(\widetilde  \sX, \widetilde E, \widetilde \mu\right)\right). 
\ee
Thee is an  induced by $\left(\rho, \left(C_0\left( \sX\right) , C_0\left( \widetilde \sX\right) G\left(\left.\widetilde \sX   ~\right| \sX\right)\right)\right)$ representation $\widetilde \rho: C_0\left( \widetilde \sX\right)\hookto B\left( \widetilde \H\right)$ (cf. Definition \ref{induced_repr_inf_defn}), such that $\widetilde \H$ is the Hilbert norm completion of
$$
K\left( C_0\left( \widetilde \sX\right)\right)\otimes_{C_0\left( {\sX}\right)} L^2\left(  \sX, E\right)
$$
where $K\left( C_0\left( \widetilde M\right)\right)$ is Pedersen's ideal of $C_0\left( \widetilde M\right)$ (cf. Definition \ref{pedersen_ideal_defn} and Remark \ref{ap_act_hilb_rem}). On the other hand 	it is proven in \cite{pedersen:mea_c} that the Pedersen's ideal $K\left(C_0\left(\widetilde{\sX} \right)\right) $ of $C_0\left(\widetilde{\sX} \right)$ coincides with $C_c\left(\widetilde{\sX} \right)$.  If
$$
\Ga_c\left(\widetilde M, \widetilde S \right)\bydef C_c\left( \lift_p\left[C\left(M, \left\{S_x\right\}, \Ga\left(M, S\right)\right)\right]\right)  
$$
then from the Lemma \ref{top_tensor_lcompact_lem} it follows that there is the given by \eqref{top_tensor_lccompact_iso_eqn} isomorphism
$$
C_c\left(\widetilde \sX\right)\otimes_{C\left(\sX\right)}\Ga\left(\sX,E \right) \xrightarrow{\approx}\Ga_c\left(\widetilde \sX, \widetilde E \right)
$$
of left  $C_0\left(\widetilde \sX \right)$-modules. On the other hand  there is the dense with respect to the Hilbert norm inclusion $\Ga_c\left(\widetilde \sX, \widetilde E \right)\subset L^2\left(\widetilde \sX, \widetilde E \right)$, hence there is  the homomorphism
\be\label{top_lift_g_inf_eqn}
\phi: C_c\left(\widetilde \sX\right)\otimes_{C_0\left(\sX\right)}\Ga\left(\sX,E \right) \to L^2\left(\widetilde \sX, \widetilde E\right)
\ee
of left  $C_0\left(\widetilde \sX \right)$-modules.
\begin{lemma}\label{top_hilb_inf_lem}
	Following conditions hold:
	\begin{itemize}
		\item [(i)] The map \eqref{top_lift_g_inf_eqn} can be extended up to the following homomorphism 
		$$
		C_c\left(\widetilde \sX\right)\otimes_{C_0\left(\sX\right)} L^2\left(\sX, E\right) \xrightarrow{\approx} L^2\left(\widetilde{\sX}, \widetilde{E} \right)
		$$
		of left  $C_0\left(\widetilde M\right)$-modules.
		\item[(ii)] The image of $C_c\left(\widetilde \sX\right)\otimes_{C\left(\sX\right)} L^2\left(\sX, E\right)$ is dense in $L^2\left(\widetilde{\sX}, \widetilde{E} \right)$.
	\end{itemize}
\end{lemma}
\begin{proof}(i)
	Let $\sum_{j = 1}^n \widetilde a_j \otimes \xi_j\in 	C_c\left(\widetilde \sX\right)\otimes_{C_0\left(\sX\right)} L^2\left(\sX, E\right)$ be any element. The space $L^2\left( \sX, E, \mu\right)$ is the Hilbert norm completion of $\Ga\left(\sX, E\right)$. Hence for any $j =1,..., n$ there is a net $\left\{\xi_{j\a}\right\}\subset \Ga\left(\sX, E\right)$ such that $\xi_j = \lim_\a \xi_{j\a}$ where we mean the convergence with respect to the Hilbert norm $\left\| \cdot\right\|_{L^2\left(\sX, E, \mu\right)}$. If $C = \max_{j = 1,...,n}\left\|\widetilde a_j\right\|$ then there is $\a_0$ such that
	$$
	\a \ge \a_0 \quad \Rightarrow \left\| \xi_j - \xi_{j\a}\right\|_{L^2\left(M, \sS\right)} < \frac{\eps}{nC},
	$$
	so one has
	$$
	\a \ge \a_0 \quad \Rightarrow \left\| \sum_{j = 1}^n \widetilde a_j \otimes \xi_{j\a_0} - \sum_{j = 1}^n \widetilde a_j \otimes \xi_{j\a}\right\|_{L^2\left(\widetilde M, \widetilde\sS\right)} < \eps.
	$$
	The above equation it follows that the net $\left\{\sum_{j = 1}^n \widetilde a_j \otimes \xi_{j\a}\right\}_\a$ satisfies to the Cauchy condition, so it is convergent with respect to the topology of $L^2\left(\widetilde \sX, \widetilde E, \widetilde \mu\right)$.\\
	(ii) From the Lemma \eqref{top_tensor_lcompact_lem} it follows that the image of $C_c\left(\widetilde M\right)\otimes_{C\left(\sX\right)} L^2\left(\sX, E\right)$ in $L^2\left(\widetilde \sX, \widetilde E\right)$ contains $\Ga_c\left(\widetilde{\sX}, \widetilde{E} \right)$, however $\Ga_c\left(\widetilde{\sX}, \widetilde{E} \right)$ is dense in $L^2\left(\widetilde{\sX}, \widetilde{E} \right)$. It follows that  $\phi\left(  C_c\left(\widetilde \sX\right)\otimes_{C\left(\sX\right)} L^2\left(\sX, E\right)\right)$ is dense in $L^2\left(\widetilde{\sX}, \widetilde{E} , \widetilde \mu\right)$. 
\end{proof}

\begin{empt} 
	If $\widetilde{a} \otimes \xi, \widetilde{b} \otimes \eta \in C_c\left(\widetilde{\sX} \right) \otimes_{C\left(\sX \right) } \Ga\left( \sX, E\right) \subset L^2\left(\widetilde{M},\widetilde{S} \right)$  then $\supp \widetilde{a}$ is compact, so there exists a subordinated  to $\widetilde p$ {covering sum} $	\sum_{\widetilde \a\in \widetilde{\mathscr A}} \widetilde a_{\widetilde\a}$ for $\supp \widetilde{a}$ (cf. Definitions \ref{top_covering_sum_defn} and \ref{top_covering_sum_subordinated_defn}) 
	where a given by \eqref{top_cfsa_eqn} notation  is used. If $\widetilde e_\a \bydef \sqrt{\widetilde a_\a}$ and $e_\a = \desc_{\widetilde p}\left( \widetilde e_\a\right) $ means $\widetilde p$-descent (cf. Definition \ref{top_lift_desc_defn}) then from \eqref{top_desc_fin_eqn} it turns out that
	$$
	\forall\a\in \mathscr A \quad \widetilde a_\a = \widetilde e_\a e_\a.
	$$
	
	The given by the equation \eqref{comp_hilb_eqn} scalar product $\left(\cdot, \cdot \right)_{\text{ind}}$ on  $C_c\left(\widetilde{\sX} \right) \otimes_{C\left(\sX \right) } \Ga\left(\sX, E\right)$ satisfies to the following equation
	\begin{equation}\label{comm_ind_l2_lc_eqn}
		\begin{split}
			\left( \widetilde{a} \otimes \xi ,    \widetilde{b} \otimes \eta\right)_{\text{ind}}= \left(\xi, \left\langle \widetilde{a}, \widetilde{b}\right\rangle_{C\left(\widetilde{M} \right) } \eta \right)_{L^2\left({M},{S}, \mu \right)}=\\=\sum_{\widetilde{\a}\in \widetilde{\mathscr A}_{\supp \widetilde{a}}}\left(\xi, \left\langle \widetilde{a}_{\widetilde{\a}}\widetilde{a}, \widetilde{b}\right\rangle_{C\left(\widetilde{M} \right) } \eta \right)_{L^2\left({M},{S}, \mu \right)}=
			\\
			= \sum_{\widetilde{\a}\in \widetilde{\mathscr A}_{\supp \widetilde{a}}}\left(\xi, \mathfrak{desc}\left( \widetilde{a}_{\widetilde{\a}} \widetilde{a}^* \widetilde{b}\right)  \eta \right)_{L^2\left({M},{S} \right)}=
			\\
			=
			\sum_{\widetilde{\a}\in \widetilde{\mathscr A}_{\supp \widetilde{a}}}~ \int_M \left(\xi_x, \mathfrak{desc}\left( \widetilde{a}_{\widetilde{\a}} \widetilde{a}^* \widetilde{b}\right)  \eta_x \right)_x d\mu=
			\\
			=
			\sum_{\widetilde{\a}\in \widetilde{\mathscr A}_{\supp \widetilde{a}}}~ \int_M \left(\mathfrak{desc}\left( \widetilde{e}_{\widetilde{\a}} \widetilde{a}\right)\xi_x, \mathfrak{desc}\left( \widetilde{e}_{\widetilde{\a}}  \widetilde{b}\right)  \eta_x \right)_x d\mu=
			\\
			=\sum_{\widetilde{\a}\in \widetilde{\mathscr A}_{\supp \widetilde{a}}}~  \int_{\widetilde{M}}\left(\widetilde{a} \lift_{\widetilde{\sU}_{\widetilde{\a}}}\left({e}_{\widetilde{\a}}\xi \right)_{\widetilde{x}},  \widetilde{b} ~\lift_{\widetilde{\sU}_{\widetilde{\a}}}\left({e}_{\widetilde{\a}}\eta \right)_{\widetilde{x}}\right)_{\widetilde{x}}d \widetilde{   \mu} =
			\\
			=\sum_{\widetilde{\a}\in \widetilde{\mathscr A}_{\supp \widetilde{a}}}~  \int_{\widetilde{M}}\widetilde{a}_{\widetilde{\a}}\left(\widetilde{a} \lift_{\widetilde{\sU}_{\widetilde{\a}}}\left(\xi \right)_{\widetilde{x}},  \widetilde{b} ~\lift_{\widetilde{\sU}_{\widetilde{\a}}}\left(\eta \right)_{\widetilde{x}}\right)_{\widetilde{x}}d \widetilde{   \mu} =
			\\
			=\int_{\widetilde{M}}\left( \widetilde{a} \lift_p\left(\xi \right)_{\widetilde{x}},  \widetilde{b} ~\lift_p\left(\eta \right)_{\widetilde{x}}\right)_{\widetilde{x}}d \widetilde{   \mu}
			= \left( \widetilde{a}~ \lift_p\left(\xi \right),  \widetilde{b} ~\lift_p\left(\eta \right)\right)_{L^2\left(\widetilde{M},\widetilde{S}, \widetilde \mu \right)} =\\
			=  \left(\phi\left(  \widetilde{a}\otimes\xi\right) ,  \phi \left( \widetilde{b}\otimes\eta\right)\right)_{L^2\left(\widetilde{M},\widetilde{S}, \widetilde\mu \right)}
		\end{split} 
	\end{equation}
	where $\phi$ is given by \eqref{top_tensor_lccompact_iso_eqn}.
	The equation \eqref{comm_ind_l2_lc_eqn} means that $\left(\cdot, \cdot \right)_{\text{ind}}\cong \left(\cdot, \cdot \right)_{L^2\left(\widetilde{M},\widetilde{S} \right)}$, and taking into account the dense inclusion $C\left(\widetilde{\sX} \right) \otimes_{C\left(\sX \right) } \Ga_c\left(\sX, E\right) \subset L^2\left(\widetilde{\sX},\widetilde{E}, \widetilde\mu \right)$ with respect to the Hilbert norm of $L^2\left(\widetilde{\sX},\widetilde{E}, \widetilde\mu\right) $
	one concludes that the space of induced representation coincides with $L^2\left(\widetilde{\sX},\widetilde{E}, \widetilde\mu\right)$. It means that induced representation $C\left(\widetilde{\sX}\right) \times L^2\left(\widetilde{\sX},\widetilde{E} \right) \to L^2\left(\widetilde{\sX},\widetilde{E} \right)$ is given by \eqref{comm_bundle_repr_eqn}. So one has the following lemma.
\end{empt}
\begin{lemma}\label{comm_ind_lc_lem}
	If the representation  $\widetilde{\rho}: C\left(\widetilde{\sX}\right) \to  B\left( \widetilde{   \H}\right)  $ is induced by the pair $$\left(\rho:C\left(\sX \right)\to B\left(  L^2\left(\sX,E \right)\right) ,\left(C\left(\sX \right) , C_0\left( \widetilde{\sX}\right) , G\left( \left.\widetilde{\sX}~\right|\sX\right) \right)  \right)$$ (cf. Definition \ref{induced_repr_inf_defn}) then following conditions hold:
	\begin{enumerate}
		\item[(a)] There is the homomorphism of Hilbert spaces $\widetilde{   \H}\cong L^2\left(\widetilde{\sX},\widetilde{E}, \widetilde\mu \right)$,
		\item[(b)] The representation $\widetilde{\rho}$ is given by the natural action of $C_0\left(\widetilde{\sX}\right)$ on $ L^2\left(\widetilde{\sX},\widetilde{E} , \widetilde \mu\right)$ (cf. \ref{comm_bundlei_repr_eqn}).
	\end{enumerate}
\end{lemma}
\begin{proof}
	(a) Follows from \eqref{comm_ind_l2_lc_eqn},\\
	(b) From the Lemma \ref{top_hilb_homo_lem}  it follows that the map 		$	C_c\left(\widetilde M\right)\otimes_{C\left(\sX\right)} L^2\left(\sX, E, \mu\right) \xrightarrow{\approx} L^2\left(\widetilde{\sX}, \widetilde{E}, \widetilde \mu \right)$ is the homomorphism of $C_0\left(\widetilde \sX\right)$ modules, so the given by \eqref{comm_bundle_repr_eqn} $C_0\left(\widetilde \sX\right)$-action coincides with the $C_0\left(\widetilde \sX\right)$-action  the given by \eqref{comp_hilb_actf_eqn}.
\end{proof}

\subsection{Coverings of spectral triples}
\paragraph*{} Let $\left(C^{\infty}\left( M\right) , L^2\left( M, S, \mu\right) ,\slashed D, J\right)$ be a commutative spectral triple (cf. equation \eqref{comm_sp_tr_eqn}), and let $\widetilde p:\widetilde{M} \to M$ be an infinite regular covering with a connected $\widetilde M$ such that  the {covering group} $G\left( {\widetilde{M}}~|~M\right)$ (cf. Definition \ref{top_group_of_covering_transformations_defn}) is residually finite (cf. Definition \ref{residually_finite_defn}). There is a  {finite covering category of} $\widetilde p : \widetilde M \to M$.
\be\label{top_m_cat_eqn}
\mathfrak{S}_{\widetilde p} = \left\{p_\la:\left(  M_\la, x^\la_0 \right) \to \left( M,{x_0}\right) \right\}
\ee
 (cf. Definition \ref{top_disconnected_defn}) such that  the {(topological) inverse limit} of $\mathfrak{S}_p$ is naturally homeomorphic  to $\widetilde{M}$.  
 From the Proposition \ref{comm_cov_mani_prop} it follows that $\widetilde{M}$ and $M_\la$ ($\forall \la \in \La$)  have natural structures of the Riemannian manifolds. From the Theorem \ref{top_main_thm} it follows that one has a good algebraical  finite covering category (cf. Definitions \ref{algebraical_finite_covering_category_defn} and \ref{good_defn})
 $$
 \mathfrak{S}_{C\left(M\right)_0\left(\widetilde p \right) } = \left\{C_0\left( p_\la\right) :C\left(M\right)\hookto\left(  M_\la\right)  \right\}
 $$
 such that the triple 
 $$ \left(C_0\left(M\right), C_0\left(\widetilde M\right),G\left(\left.\widetilde M~\right| M\right)\right)
 $$ is  the  { infinite noncommutative covering} of $\mathfrak{S}_{C\left(M\right)_0\left(\widetilde p \right) }$ (cf. Definition \ref{infinite_noncommutative_covering_defn}).  Denote by $\widetilde{S}\stackrel{\mathrm{def}}{=} \widetilde p^*S$ and $S_\la \bydef p^*_\la S $  (for all $\la \in \La$) inverse images of the Spin$^{c}$ bundle (cf. \ref{vb_inv_img_funct_defn}). Similarly to \eqref{top_diff_*alg_lift_eqn} we define the following operators
	\be\label{top_dlift_eqn}
\begin{split}
	 \Dslash_\la \bydef D^*\left(p_\la, S \right)\left(  \Dslash\right)  \in D^*\left( M_\la,  S_\la\right),\\
\widetilde 	 \Dslash \bydef  D^*\left(\widetilde p, S \right)\left(  \Dslash\right) \in D^*\left(\widetilde M, \widetilde S\right).
\end{split}
\ee
(cf. equations  \eqref{top_d_shaef_eqn} and \eqref{top_d_shaef_m_eqn})
Otherwise both $\Dslash_\la$ and $\widetilde \Dslash$ can be regarded as $\C$-linear operators, i.e.
\bean
 \Dslash_\la :\Ga^\infty\left( M_\la,  \sS_\la \right)\to\Ga^\infty\left( M_\la,  \sS_\la \right),\\
\widetilde \Dslash :\Ga^\infty_c\left(\widetilde M,  \widetilde\sS \right)\to\Ga^\infty_c\left(\widetilde M,  \widetilde\sS \right) .
\eean
Both $\Ga^\infty\left( M_\la,  \sS_\la \right)$ and $\Ga^\infty_c\left(  \widetilde M,  \widetilde \sS \right)$ are dense in $L^2\left( M_\la,  \sS_\la, \mu_\la \bydef \lift_{p_\la}\left( \mu\right)  \right)$ and $L^2\left(  \widetilde M,  \widetilde \sS, \widetilde \mu \bydef \lift_{\widetilde p}\left(\mu \right)  \right)$ where $\lift$ means the lift of the measure (cf. Definition \ref{top_lift_measure_defn}). It follows that  both  operators  $\Dslash_\la$ and $\widetilde\Dslash$ are  unbounded operators on $L^2\left( M_\la,  \sS_\la \right)$ and $L^2\left(  \widetilde M,  \widetilde \sS \right)$. Using the equation \eqref{cc_eqn} one can define an anti-isometry operator $J_\la : \H_\la \to \H_\la$ for any $\la\in\La$. There is a family 
\begin{equation}\label{comm_triple_seq_eqn}
	\begin{split}
		\mathfrak{S}_{\left(C^{\infty}\left( M\right) , L^2\left( M, S\right) ,\slashed D,J\right) } = \left\{
		\left(C^{\infty}\left( M_\la\right) , L^2\left( M_\la, S_\la\right) , \slashed D_\la, J_\la\right)
		\right\}_{\la \in \La}
	\end{split}
\end{equation}
of spectral triples. 
\begin{lemma}
The given by the equation \eqref{comm_triple_seq_eqn} family $	\mathfrak{S}_{\left(C^{\infty}\left( M\right) , L^2\left( M, S\right) ,\slashed D, J\right) }$ is the coherent set of spectral triples (cf. Definition \ref{spectral_triple_weakly_coh_defn}).
\end{lemma}

\begin{proof}
One needs check that $	\mathfrak{S}_{\left(C^{\infty}\left( M\right) , L^2\left( M, S\right) ,\slashed D\right) }$ satisfies to conditions (a)-(d) of the Definition \ref{spectral_triple_weakly_coh_defn}.\\
	(a) Follows from the Theorem \ref{comm_fin_sp_tr_thm}.\\
	(b)  $C\left( M_\la\right)$ is the $C^*$-norm completion of $\Coo\left( M_\la\right)$ for all $\la\in\La$.\\
	(c)	 From the Theorem \ref{top_main_thm} it turns out that the if  $\widehat \pi: C^*$-$\varinjlim_{\la \in \La} C\left( M_\la\right)  \hookto B\left(\widehat \H\right) $ is an atomic representation (cf. Definition \ref{atomic_repr_defn}) it turns out  that an {algebraical  finite covering category} $
	\mathfrak{S}_{C\left(M\right)}=
	\left\{C(M)\hookto C(M_\la)\right\}_{\la \in \La}
	$ is $\widehat\pi$-good  (cf. Definition \ref{good_defn}).\\
	(d)  Follows from the Theorem \ref{comm_fin_sp_tr_thm}.
\end{proof}
\begin{empt} Here a specialization of the explained in the Section \ref{str_cov_glo_sec} construction is described. 	
	If $\rho: C\left(M\right) \hookto B\left(L^2\left( M, S, \mu\right)\right)$ a faithful representation which corresponds to a spectral triple $\left(C^{\infty}\left( M\right) , L^2\left( M, S, \mu\right) ,\slashed D, J\right)$ (cf. Definition \ref{sp_tr_defn_sec}) then $\rho$ can be constructed as the completion of the action
	$$
	C\left(M \right) \times \Ga\left(M, S \right) \to \Ga\left(M, S \right)
	$$
	because a $\Ga\left(M, S \right)$ is a dense $	C\left(M \right)$-module of $L^2\left( M, S, \mu\right)$. From the Lemma \ref{comm_ind_lc_lem} it turns out that is induced by the pair $$\left(\rho,C\left(M \right)\to B\left(  L^2\left(M,S \right)\right) ,\left(C\left(M \right) , C_0\left( \widetilde{M}\right) , G\left( \left.\widetilde{M}~\right|M\right) \right)  \right)$$ representation (cf. Definition \ref{induced_repr_inf_defn}) is equivalent to the natural representation $\widetilde \rho: C_0\left(\widetilde{M} \right) \hookto B\left(L^2\left(\widetilde{M},\widetilde{S}, \widetilde\mu \right)\right)$ where $\widetilde\mu \bydef \lift_{\widetilde  p}\left(\mu \right)$ is the $\widetilde p$-lift of $\mu$ (cf. Definition \ref{top_lift_measure_defn}). 	We define an anti-isometry $\widetilde J: (L^2\left(\widetilde M, \widetilde S, \widetilde \mu\right)\to L^2\left(\widetilde M, \widetilde S, \widetilde \mu\right) $ which is a unique extension of the given by \eqref{cc_inf_defn} anti-linear operator. Similarly to the equations \eqref{da_la_eqn} and \eqref{ps_inc_eqn} we define inclusions
	\bea
	\left[\Dslash_\la, \Coo\left( M_\la\right) \right]\subset B\left( L^2\left(\widetilde M, \widetilde S, \widetilde \mu\right)\right),\\
		\pi^s_\la: \Coo\left( M_\la\right) \hookto B\left(L^2\left(\widetilde M, \widetilde S, \widetilde \mu\right)^{2^s} \right).
	\eea
	
	On the other hand the representation  $\widetilde \rho$ can be regarded as a completion of the action
	$$
	C_0\left(\widetilde M \right) \times \Ga_c\left(\widetilde M,\widetilde S \right) \to \Ga_c\left(\widetilde M, \widetilde S \right).
	$$ (cf. Lemma \ref{top_hilb_inf_lem}).
\end{empt}
\begin{empt}\label{top_smooth_a_empt}
	If   $\widetilde a \in K\left( C_0\left(\widetilde M \right) \right) \cong C_c\left(\widetilde M\right)$ then $\widetilde a$ satisfies to the condition (b)    of the Definition \ref{smooth_el_defn} then
from  the Lemma \ref{top_compact_la_lem}  there are an open subset $\widetilde{   \mathcal U }$ of $\widetilde M$ and $\la_{\widetilde \sU} \in \La$
	  such that
	\begin{itemize}
		\item $\supp \widetilde a \subset \widetilde \sU$,
		\item   for  any  $\la \ge \la_{\widetilde \sU}$ the set $\widetilde{   \mathcal U }$  is mapped homeomorphically onto 
		$\widetilde{p}_\la\left(\widetilde{   \mathcal U } \right) \subset M_\la$ where $\widetilde{p}_\la:   \widetilde{   M}   \to  M_\la$ is the natural covering.
	\end{itemize}
\end{empt}
\begin{lemma}\label{top_smooth_a_lem}
	Under the hypotheses \ref{top_smooth_a_empt} the element $\widetilde a$ satisfies to the condition (b) of the Definition \ref{smooth_el_defn} if and only if $\widetilde a\in \Coo_c\left(\widetilde M\right)\bydef C_c\left(\widetilde M\right)\cap \Coo\left(\widetilde M\right)$.
\end{lemma}	
\begin{proof}
If $\widetilde a\in \Coo_c\left(\widetilde M\right)\bydef C_c\left(\widetilde M\right)\cap \Coo\left(\widetilde M\right)$ then from the Lemma \ref{comm_lift_desc_sum_lem} one has
$$
 		a_\la=	\sum_{g \in \ker\left(  G\left(\left.\widetilde{A}_\pi\right|A \right)\to G\left(\left.A_\la\right|A \right)\right) } ~ g \widetilde{a
 		} = \desc^c_{\widetilde p_{\la}}\left(\widetilde a \right) 
$$
where $\desc^c_{\widetilde p_{\la}}$ is a compactly supported $\widetilde p_{\la}$-descent (cf. Definition  \ref{top_compactly_supported_descent_defn}). Taking into account that $\Coo_0\left( \widetilde M\right)$ is a $c$-soft *-algebra (cf. Definition \ref{top_soft_r_defn}) from \eqref{top_cdesc_soft_eqn} it follows that
$$
a_\la = \desc^c_{\widetilde p_{\la}}\left(\widetilde a \right)\in \Coo\left(\widetilde M \right).
$$
There is an open subset $\widetilde \sU \subset \widetilde \sX$ with compact closure such that $\supp \widetilde a \subset \widetilde \sU$. If $\widetilde p_\la: \widetilde \sX \to \sX$ is a natural covering then from the Lemma \ref{top_compact_la_lem}  one has  $\la_{\widetilde \sU} \in \La$ such that for all $\la \ge \la_{\widetilde \sU}$ the set  $\widetilde p_\la\left(\widetilde \sU \right)$ is evenly covered by $\widetilde{p}_\la$. 
If $\la \ge \la_{\widetilde \sU}$ then from the Lemma  \ref{comm_lift_desc_sum_lem} one has
$$
a_\la=	\sum_{g \in \ker\left(  G\left(\left.\widetilde{A}_\pi\right|A \right)\to G\left(\left.A_\la\right|A \right)\right) } ~ g \widetilde{a
} = \desc_{\widetilde p_{\la}}\left(\widetilde a \right) 
$$
where $\desc_{\widetilde p_{\la}}$ is a $\widetilde p_{\la}$-descent (cf. Definition \ref{top_lift_desc_defn}). From the equation \eqref{top_lift_desc_eqn} it follows that $\widetilde{a}$ is a   $\widetilde p_{\la}$-$\widetilde{\sU}$-{lift} of $\widetilde a$, i.e. 
\be\label{top_lift_exer_eqn}
\widetilde{a}= \lift^{\widetilde p_{\la}}_{\widetilde{\sU}}\left(a_\la \right)
\ee
From n$^\text{0}$ 2 of the Exercise \ref{top_sheaf_corr_exer} it follows that $\widetilde a\in \Coo\left( \widetilde M\right)$. 

\end{proof}
 
\begin{lemma}\label{top_smooth_lem}
An element $\widetilde a \in K\left( C_0\left(\widetilde M \right) \right) \cong C_c\left(\widetilde M\right)$  is smooth  (cf.  Definition \ref{smooth_el_defn})  if and only if $\widetilde a\in \Coo_c\left(\widetilde M\right)=C_c\left(\widetilde M\right)\cap \Coo\left(\widetilde M\right)$.
\end{lemma}
\begin{proof}
	If $\widetilde a \in K\left( C_0\left(\widetilde M \right) \right) \cong C_c\left(\widetilde M\right)$ is smooth then  from the Lemma \ref{top_smooth_a_lem} it follows that $\widetilde a \in \Coo_c\left(\widetilde M\right)$.
	
	Conversely suppose that $\widetilde a\in \Coo_c\left(\widetilde M\right)$. One needs check that $\widetilde a$ satisfies to conditions (a)-(d) of the Definition \ref{smooth_el_defn}.
	 	\begin{enumerate}
		\item[(a)] Evident.
		\item[(b)] Follows from the Lemma \ref{top_smooth_a_lem}
		\item[(c)] From $\widetilde a\in C_c\left(\widetilde M\right)$ it follows that $\supp \widetilde a$ is compact. From the Lemma \ref{top_compact_la_lem} it follows that for there is $\la_{\supp \widetilde a}\in \La$ such that $\supp \widetilde a$ is mapped homeomorphically onto $\widetilde p_{\la_{\supp \widetilde a}}\left(\supp \widetilde a\right)$.  If $G_\la \bydef \ker\left(G\left(\left.\widetilde M\right| M\right)\to G\left(\left.M_\la \right| M\right)\right) $ and 
	\bean
\forall \la \in \La \quad \la \ge 	\la_{\supp \widetilde a} \quad \forall g_1, g_2 \in G_\la \quad 	g_1 \neq g_2 \quad \Rightarrow \\
\Rightarrow g_1 \supp \widetilde a\cap g_2 \supp \widetilde a= \emptyset
	\eean
		If $\widetilde\xi\in L^2\left(\widetilde M, \widetilde S, \widetilde \mu\right)^{2^s}$  then 
		$$
	\left\| \widetilde \xi \right\|^2= \sum_{g \in G_{\la{\supp \widetilde a}}}
\left\|\chi_{\supp g \widetilde a} \widetilde \xi \right\|^2$$
where $ \chi_{\supp g\widetilde a}$ is a characteristic function of $\supp g\widetilde a$  (cf. Definition \ref{top_char_f_defn}).  If $e\in G\left(\left.\widetilde M\right| M\right)$  is the  neutral  element then for any $\delta > 0$ there is a finite subset $S'\subset  G\left(\left.\widetilde M\right| M\right)\setminus\{e\}$ such that
 $$
 \sum_{g \in G_{\la_{\supp \widetilde a}}\setminus \left(S' \cup \{e\}\right)}
 \left\| \chi_{\supp g\widetilde a} \widetilde \xi \right\|^2 < \delta.
 $$ 
 So for any $\eps > 0$ one can select a finite subset $S\subset G_{\la_{\supp \widetilde a}}\setminus\{e\}$ such that
 \be\label{top_se_eqn}
 \left\| \sum_{g \in G_{\la_{\supp \widetilde a}}\setminus \left(S \cup \{e\}\right)} \chi_{\supp g\widetilde a} \widetilde \xi \right\|< \frac{\eps }{ \left\|\pi^s\left( \widetilde a\right)  \right\|}.
\ee 
 Since the set $S$ is finite and 
$$
\bigcap_{\la \in \La}G_\la = \{e\}
$$
 there is $\la_S \in \La$ such that $\la_S \ge \la_{\supp \widetilde a}$ and
 $$
\la \ge \la_S \quad \Rightarrow \quad S \cap  G_\la= \emptyset.
 $$
 From our construction it follows that
\bean
\la \ge \la_S~~\Rightarrow~~ \quad\left\|\pi^s\left( \widetilde a\right) \widetilde \xi - \sum_{	g \in G_\la }\pi^s\left( g\widetilde a\right) \widetilde \xi \right\|= \left\|\pi^s\left( \widetilde a\right) \widetilde \xi - \pi^s\left( a_\la\right) \widetilde \xi \right\|< \eps,
\eean
i.e. the net $\left\{\pi^s_\la\left(a_\la \right)\right\}_{\la \in\La}$  is convergent with respect to the strong topology of 
 $B\left(L^2\left(\widetilde M, \widetilde S, \widetilde \mu\right)^{2^s}\right)$.		
 \item[(d)] 	The space $\Ga_c^\infty\left(\widetilde M, \widetilde S \right)$ of compactly supported is dense in $L^2\left(\widetilde M, \widetilde S, \widetilde \mu\right)$ so we should prove that a net 
\bean
\left\lbrace \left[D_\la, \sum_{	g \in G_\la}\widetilde a_\la \right]\widetilde\eta\right\rbrace_{\la\in\La}
\eean  
 for any $\widetilde\eta\in \Ga_c^\infty\left(\widetilde M, \widetilde S \right)$. Similarly to the equation \ref{top_se_eqn} for all $\eps$ one can find a finite subset  $T\subset G_{\la_{\supp \widetilde a}}\setminus\{e\}$
\bean
 \left\| \sum_{g \in G_{\la_{\supp \widetilde a}}\setminus \left(T \cup \{e\}\right)} \chi_{\supp g\widetilde a} \widetilde \eta \right\|< \frac{\eps }{ \left\|\left[\Dslash_{\la_{\supp \widetilde a}}, a_{\la_{\supp \widetilde a}} \right]\right\|}.
 \eean 
There is $\la_T \in \La$ such that $\la_T \ge \la_{\supp \widetilde a}$ and
 $$
 \la \ge \la_T \quad \Rightarrow \quad T \cap  G_\la= \emptyset.
 $$
 From our construction it follows that
 \bean
 \la \ge \la_T~~\Rightarrow~~ \quad\left\|\left[\Dslash_{\la}, a_{\la} \right]\widetilde \eta-\left[\Dslash_{\la_{T}}, a_{\la_{T}} \right]\widetilde \eta \right\|< \eps,
 \eean
 i.e. the net $\left\{\left[\Dslash_{\la}, a_{\la} \right]\right\}_{\la\in\La}$  is convergent with respect to the strong topology of 
 $B\left(L^2\left(\widetilde M, \widetilde S, \widetilde \mu\right)\right)$. If  $\mathscr S^{\Ga_c^\infty\left(\widetilde M, \widetilde S \right)}$ is a $\Ga_c^\infty\left(\widetilde M, \widetilde S \right)$-sheaf (cf. Definition \ref{top_x_sheaf_defn}) then $\widetilde\eta \in \mathscr S^{\Ga_c^\infty\left(\widetilde M, \widetilde S \right)}\left(\widetilde M \right) $. Moreover if   $\mathscr S^{\Ga_c^\infty\left( M_{\la_{\supp \widetilde a}},  S_{\la_{\supp \widetilde a}} \right)}$ is a $\Ga_c^\infty\left( M_{\la_{\supp \widetilde a}},  S_{\la_{\supp \widetilde a}} \right)$-sheaf then on has
 \bean
 \Dslash_{\la_{\supp \widetilde a}}\in \mathscr Hom \left( \mathscr S^{\Ga_c^\infty\left( M_{\la_{\supp \widetilde a}},  S_{\la_{\supp \widetilde a}} \right)},  \mathscr S^{\Ga_c^\infty\left( M_{\la_{\supp \widetilde a}},  S_{\la_{\supp \widetilde a}} \right)}\right)\left(M_{\la_{\supp \widetilde a}} \right)  
 \eean
 where  $\mathscr Hom$ is a sheaf of local morphisms (cf. Definition \ref{sheaf_hom_defn}). If $\widetilde p_{\la_{\supp \widetilde a}}: \widetilde M \to M_{\la_{\supp \widetilde a}}$ is a natural covering then from our construction it follows that
  \be\label{top_llla_eqn}
\lim_{\la \in \La} \left[\Dslash_{\la}, a_{\la} \right]\widetilde \eta = \left[ p^{-1}_{\la_{\supp \widetilde a}}\left(\Dslash_{\la_{\supp \widetilde a}} \right) , \widetilde a \right]\widetilde \eta
 \ee
 where the map
\bean
\widetilde p^{-1}_{\la_{\supp \widetilde a}} : \mathscr Hom \left( \mathscr S^{\Ga_c^\infty\left( M_{\la_{\supp \widetilde a}},  S_{\la_{\supp \widetilde a}} \right)},  \mathscr S^{\Ga_c^\infty\left( M_{\la_{\supp \widetilde a}},  S_{\la_{\supp \widetilde a}} \right)}\right)\left(M_{\la_{\supp \widetilde a}} \right)\hookto\\\hookto \mathscr Hom\left( \Ga_c^\infty\left(\widetilde M, \widetilde S \right), \Ga_c^\infty\left(\widetilde M, \widetilde S \right)\right) 
\eean
is the  $\widetilde p_{\la_{\supp \widetilde a}}$-inverse image (cf. Definition \ref{top_smooth_inv_im_defn}).
	\end{enumerate}

\end{proof}

\begin{remark}
From the Lemma \ref{top_smooth_lem} and the Definition \ref{smooth_alg_defn} it follows that  {smooth algebra} $\widetilde\A$ of the coherent set \eqref{comm_triple_seq_eqn} of spectral triples  is a completion of $\Coo_c\left(\widetilde M\right)$ in the topology induced by the seminorms $\left\| \cdot \right\|_s$ (cf. equation \ref{smooth_seminorms_eqn}). Since $\Coo_c\left(\widetilde M\right)$ is dense in $\Coo_c\left(\widetilde M\right)$ the coherent set \eqref{comm_triple_seq_eqn} of spectral triples is weakly good (cf. Definition \ref{smooth_alg_defn})
\end{remark}
\begin{empt}
	Our construction can be expressed by  the following mapping.
	\\
	\\
	\begin{tabular}{|c|c|c|}
		\hline
		&General theory & The specialization\\ 
		\hline
		&	&\\
		Hilbert spaces & $\H$  and $\widetilde\H$ &  $L^2\left(M, \sS \right)$ and $L^2\left(\widetilde M, \widetilde \sS \right)$\\ & & \\
		Pre-$C^*$-algebra	& $\A$  & $\Coo\left(M \right)$   \\  & & \\
		Pedersen's ideal	& $K\left(\widetilde A\right) $  & $C_c\left(\widetilde M \right)$   \\  & & \\
		The space of  	& $\widetilde{W}^\infty$  & $\Coo_c\left(\widetilde M \right)$   \\ smooth elements & & \\& & \\
Charge conjugation 	& $J$ and $\widetilde J$ & $C$ and $\widetilde C$  \\  & & \\		Dirac operators & $D$ & $ \Dslash$  \\
		& $\widetilde{D}$ & ?\\  & & \\
		&$\H^\infty\bydef \bigcap_{n =1}^\infty \Dom D^n\subset \H$ & $\Ga^\infty\left(M, \sS \right)= \bigcap_{n =1}^\infty \Dom \Dslash^n$ \\& & \\
		\hline
	\end{tabular}
	\\
	\\
To complete it one should find an explicit formula for $\widetilde D$. From \eqref{top_llla_eqn}
	it follows that 
	\bean
	\lim_{\la \in \La} \left[\Dslash_{\la}, a_{\la} \right] = \left[\widetilde p^{-1}_{\la_{\supp \widetilde a}}\left(\Dslash_{\la_{\supp \widetilde a}} \right) , \widetilde a \right].
	\eean
	If $p_{\la_{\supp \widetilde a}}: M_{\la_{\supp \widetilde a}}\to M$ is the natural covering then form the equation \eqref{comm_d_eqn} it follows that $\Dslash_{\la_{\supp \widetilde a}}= p^{-1}_{\la_{\supp \widetilde a}}\left( \Dslash\right)$ and taking into account \eqref{top_sh_inv_comp_eqn} we have 
	$\widetilde p^{-1}_{\la_{\supp \widetilde a}}\left(\Dslash_{\la_{\supp \widetilde a}}\right) = \widetilde p^{-1}\left(\Dslash \right)$ where $\widetilde p : \widetilde M \to M$ is the natural covering. In result one has
\be\label{top_ds_eqn}
\left[\widetilde p^{-1}\left(\Dslash \right) , \widetilde a \right]=\lim_{\la\in\La }\left[\Dslash_{\la}, a_{\la} \right].
\ee
On the other hand if $\widetilde D$ is a specialization of 
the given by \eqref{inf_lift_D_eqn} operator then one has
$$
\left[\widetilde D, \widetilde a \right]= \lim_{\la\in\La }\left[\Dslash_{\la}, a_{\la} \right],
$$
and taking into account \eqref{top_ds_eqn}  we conclude that
 \be\label{top_lift_inf_dirac}
 \widetilde D=\lift_{\widetilde p}\left(  \Dslash\right).
 \ee
  
\end{empt}
A following theorem is a result of the above constructions. 
\begin{theorem}\label{comm_inf_sp_tr_thm}
	The triple  $\left(  \widetilde{\A} , L^2\left(\widetilde M, \widetilde S \right)  , \widetilde{\Dslash}\bydef {\widetilde p^{-1}}\left( \Dslash\right), \widetilde J \right)$ is the\\ $\left(C\left( M\right)  C_0\left( \widetilde{M}\right) , G\left(\left.\widetilde{M}~\right|M \right)\right)$-{lift} of $\left( \Coo\left( {M}\right) , L^2\left( M,  S \right)  , {\Dslash}, J\right)$ (cf. Definition \ref{reg_triple_defn}).
\end{theorem}
\begin{remark}
	The Theorem \ref{comm_inf_sp_tr_thm} is a noncommutative analog of the Proposition \ref{comm_cov_mani_prop}.
\end{remark}

\chapter{Hausdorff blowing-up and coverings}\label{blowing_chap}
\paragraph*{}
If the spectrum of  $C^*$-algebra $A$ is Hausdorff then $A$ possesses  a family of good properties. For example noncommutative coverings of $C^*$-algebras with Hausdorff spectra correspond to coverings of their spectra (cf. Chapters \ref{top_chap} and \ref{ctr_chap}). Sometimes the spectrum of $C^*$-algebra can be replaced by its Hausdorff blowing-up and some results concerning $C^*$-algebras with  Hausdorff spectrum can be applied to $C^*$-algebras with Hausdorff blowing-up. We prove that some noncommutative coverings of $C^*$-algebras correspond to coverings of Hausdorff blowing-up. The results of this chapter are used in both $C^*$-algebras with Hausdorff spectrum and algebras of foliations having  non - Hausdorff spectrum. 
The "blowing-up" word is inspired by following reasons:
\begin{itemize}
	\item sometimes there is  the natural partially defined  surjective  map from  Hausdorff blowing-up to the spectrum, 
	\item  in the algebraic geometry   "blowing-up" means  excluding of singular points  (cf. \cite{hartshorne:ag}).
\end{itemize}

\section{Basic definitions}

\paragraph*{}
 Before the  definition of blowing-up we consider an illustrative example.
For all $\th\in \R_+$ there  is the natural action of group  $G\bydef \Z\th/ \Z$ on the circle
$\sX\bydef S^1\cong \R/\Z$. If $\th$ is rational then the space $\sX/G$ is Hausdorff, otherwise it is not Hausdorff. How can we obtain $\sX/G$ algebraically? The Table \ref{hnh_alg_table} shows different methods of obtaining of the spectra.
\newline
\begin {table}[H]
\caption {Comparison of an invariant algebra and a crossed product} \label{hnh_alg_table}
\begin{tabular}{|c|c|c|}
	\hline
	&$\sX/G$ is Hausdorff & $\sX/G$ is not Hausdorff\\
	\hline
	&&\\
	$C\left( \sX\right)^G$ 	& $C\left( \sX/G\right)$	&$\C$\\
	& & \\
	Spectrum of $C\left( \sX\right)^G$ 	& $ \sX/G$	& A point\\
	& & \\
	\hline
	&&\\
	$C\left( \sX\right)\rtimes G$ 	& $C\left( \sX\right)\rtimes G\neq\C$	&$C\left( \sX\right)\rtimes G\neq\C$\\
	& & \\
	Spectrum of $C\left( \sX\right)\rtimes G$ 	& $ \sX/G$	& $\sX/G$\\
	& & \\
	
	\hline
\end{tabular}
\end{table}
Thus the $C\left( \sX\right)^G$ invariant $C^*$-algebra does not always yield $\sX/G$ as a spectrum but the crossed product $C\left( \sX\right)\rtimes G$ does. 
Moreover since $C\left( \sX/G\right)$ can be trivial it does not know the map $\sX \to \sX/G$  but $C\left( \sX\right)\rtimes G$ knows it, because there  is an inclusion
$C\left(\sX \right) \subset C\left( \sX\right)\rtimes G$. If $\sX/G$ is Hausdorff then sometimes a covering of $\sX/G$ can by find by following way:
\begin{enumerate}
\item find a covering $\widetilde\sX \to\sX$ with an action $G\times \widetilde\sX \to \widetilde\sX$,
\item then obtain  the covering $\widetilde\sX / G \to \sX/G$.
\end{enumerate}
However if $\sX/G$ is not Hausdorff this construction does not work. 
However there is an algebraic version of it containing two steps:
\begin{enumerate}
\item find a covering $\widetilde\sX \to\sX$ with an action $G\times \widetilde\sX \to \widetilde\sX$, and invective $*$-homomorphism $C\left( \sX\right)\hookto  C\left( \widetilde\sX\right)$,
\item then obtain a noncommutative covering  $C\left( \sX\right)\rtimes G\hookto C\left( \widetilde\sX\right)\rtimes G $.
\end{enumerate}
The algebraic construction can be applied if $\sX/G$ is not Hausdorff.  The map $\sX \to \sX/G$ is yields a specialization of discussed below Hausdorff blowing up.


\begin{definition}\label{blowing_defn}
Any  commutative $C^*$-subalgebra $C\cong C_0\left( \sY\right) \subset M\left(A \right)$ is said to be   \textit{Hausdorff blowing-up} of $A$ if  both sets
\be\label{blowing_eqn}
\begin{split}
C_c\left( \sY\right)A \bydef \left\{fa| f \in C_c\left( \sY\right)\quad a \in A \right\},\\
AC_c\left( \sY\right) \bydef \left\{af| f \in C_c\left( \sY\right)\quad a \in A \right\}
\end{split}
\ee
are dense in $A$.
\end{definition}
\begin{remark}\label{blowing_rem}
$C_c\left( \sY\right)A$ is dense in $A$ if and only if $AC_c\left( \sY\right)$ is dense in $A$ (cf. equations \eqref{blowing_eqn}).
\end{remark}
\begin{example}\label{blowing_hausdorff_exm}	If $A$ is a  given by \eqref{cross_mult_eqn} $C^*$-algebra with Hausdorff spectrum $\sX$  then from the Theorems \ref{oa_haus_prim_thm} and  \ref{dauns_hofmann_thm} one has the natural inclusion $C_0\left(\sX \right) \hookto M\left( A\right)$, which is Hausdorff blowing-up.
\end{example}
\begin{example}\label{blowing_almost_hausdorff_exm}
Under the hypotheses \ref{spectrum_ff_p_empt} there are Hausdorff blowing-ups:
\bean
C \hookto M\left( A\right),\\
\widetilde C \hookto M\left(\widetilde A \right) 
\eean
where the notation of \ref{spectrum_ff_p_empt} is used.
\end{example}
\begin{definition}\label{blowing_ideals_au_ua_defn}
	Let  $ C_0\left( \sY\right)\subset  M\left( A\right) $ be  Hausdorff blowing-up of $A$ (cf. Definition \ref{blowing_defn}), and let $\sU \subset \sY$ be an open subset. Both   left and right  closed ideals $A_\sU$  and $_\sU A$ of $A$ (cf. Definition \ref{ideal_left_right_defn}) generated by sets 	$AC_0\left( \sU\right)$ and $C_0\left( \sU\right)A$ are the \textit{left} $\sU$-\textit{ideal} and the \textit{right} $\sU$-\textit{ideal} respectively. A hereditary $C^*$-subalgebra of $A$
\be\label{blowing_hereditary_u_eqn} 
\begin{split}
	_\sU A_\sU \bydef		_\sU A\cap  A_\sU = A^*	_\sU \cap  A_\sU\\ (\text{cf. Definition \ref{hered_defn} and the Lemma \ref{hered_bab_lem}}).
\end{split}
\ee	
is the $\sU$-\textit{subalgebra}.
	
\end{definition}

\begin{lemma}\label{blowing_compact_lem}
	If $C\cong C_0\left( \sY\right) \subset M\left(A \right)$ is    {Hausdorff blowing-up} of $A$  (cf. Definition \ref{blowing_defn}	for any $a \in A$ and $\eps > 0$ following conditions hold:
	\begin{enumerate}
		\item[(i)] there is a positive $f \in C_c\left( \sY\right)_+$ with 
		\be\label{blowing_compact_eqn}
		\begin{split}
			\left\| f  \right\| \le 1,\\
			\left\| a - af  \right\|< \eps,\\
			\left\| a - fa \right\|< \eps,\\
			\left\| a - faf  \right\|< \eps,
		\end{split}
		\ee
		\item[(ii)] there are an open subset $\sU \subset \sY$ with compact closure and $b \in ~_\sU A_\sU$ such that
		\be\label{blowing_compact_b_eqn}
		\begin{split}
			b \in  ~_\sU A_\sU,\\
			\left\| a - b \right\|< \eps.
		\end{split}
		\ee
	\end{enumerate}
\end{lemma}
\begin{proof}(i)
	If $g \in C_c\left( \sY\right)$ then $g = fg = gf$ for any positive $f \in C_c\left( \sY\right)_+$ such that
	\bean
	\left\| f \right\|= 1,\\
	f\left(\supp g \right) = 1.
	\eean 
	If $f' \in C_c\left( \sY\right)$ and $c' \in A$ such that
	$\left\|a-f'c' \right\| < \eps/4$ (cf. equation \eqref{blowing_eqn}) then there for any positive $f_1$ such that $\left\| f_1 \right\|= 1$, $~f_1\left(\supp f' \right) = 1$ one has
	\bean
	\left\|f_1\left( a-f'c \right) \right\|\le \left\|f_1 \right\|\left\|a-f'c \right\|\le \frac{\eps}{4},\\
	\left\|a-f_1a \right\| < \left\|f_1a - f_1f'c\right\|+ \left\|a - f_1f'c\right\|\le \frac{\eps}{2}
	\eean
	Similarly 	If $f'' \in C_c\left( \sY\right)$ and $c''/ \in A$ such that
	$\left\|a-f'c'' \right\| < \eps/4$ then for any positive for any positive $f_2$ such that $\left\| f_2 \right\|= 1$, $~f_1\left(\supp f' \right) = 1$ one has
	\bean
	\left\|a-a f_2 \right\| <  \frac{\eps}{2}
	\eean
	If $f = \max\left(f_1, f_2 \right)$ then $f\left(\supp f' \cup \supp f'' \right)= \{1\}$ 
	and
	\bean
	\begin{split}
		\left\| f  \right\| \le 1,\\
		\left\| a - af  \right\|< \frac{\eps}{2},\\
		\left\| a - fa \right\|< \frac{\eps}{2},\\
	\end{split}
	\eean
	On the other hand
	\bean
	\left\| a - faf \right\|\le \left\| a - fa \right\|+ \left\| fa - faf \right\|< \frac{\eps}{2}+ 	\left\| f\right\|\left\|a - fa\right\|< \eps.
	\eean
	(ii) The set $\sU \bydef \left\{y \in \sY | f\left(y\right)\neq 0\right\}\subset \sY$ is open and closure of $\sU$ is the compact set $\supp f$. Moreover $f a f \in ~_\sU A_\sU$ and  $\left\|a - faf\right\|< \eps$.
	
\end{proof}
\begin{empt}\label{blowing_compact_supp_empt}
If $\left\{\sU_\la\right\}_{\la\in\La}$ is a given by the Corollary	\ref{top_connected_union_cor} family of open subsets of $\sY$ with compact closures then $\sY = \cup_{\la \in \La} \sU_\la$. Since the family is directed the union $A_c \bydef \cup_{\la \in \La} \sU_\la$ is a $*$-subalgebra of $A$. From the Lemma \ref{blowing_compact_lem} it follows that $A_c$ is dense in $A$.
\end{empt}
\begin{definition}\label{blowing_compact_supp_defn}
The given by \ref{blowing_compact_supp_empt} dense $*$-subalgebra $A_c \subset A$ is said to be the \textit{algebra of compactly supported elements}.
\end{definition}

\begin{empt} 

We leave to the reader a proof of following equations:
\be\label{blowing_su_eqn}
\begin{split}
	_\sU A \bydef \left\{ a \in A~ |~ \forall f \in C_0\left( \sY\right)\quad f\left(\sU \right)= \{0\}  \quad \Rightarrow \quad fa = 0\right\} ,\\
	A_\sU \bydef \left\{ a \in A~ |~ \forall f \in C_0\left( \sY\right)\quad  f\left(\sU \right)= \{0\} \quad \Rightarrow \quad af = 0\right\}
\end{split}
\ee

\bea\label{blowing_sue1_eqn}
	\sU'\cap \sU'' =\emptyset\quad \Rightarrow\quad A_{\sU'}~_{\sU''}A= \{0\}.
\eea

From \eqref{blowing_su_eqn} it follows that
\be\label{blowing_su_inc_eqn}
\sU'\subset \sU'' \quad \Rightarrow\quad _{\sU'}A \subset~  _{\sU'}A ~\text{ AND } ~A_{\sU'}\subset A_{\sU''}~\text{ AND }~  _{\sU'}A _{\sU'}\subset _{\sU''}A _{\sU''}.
\ee
\end{empt}

\begin{definition}\label{blowing_support_defn}
 	If $C_0\left( \sY\right) \subset M\left(A \right)$ is    {Hausdorff blowing-up} of $A$  (cf. Definition \ref{blowing_defn}),  $a \in A$ and
	$
	\sU_a \bydef\bigcap 
	\left\{\left.{\sU} \subset \sX\right| a\in~_\sU A_{\sU} \right\}
	$
	then the  closure $\sV_a$  of $\sU_a$ is said to be the \textit{support} of $a$. We write $\supp a \bydef \sV_a$.
\end{definition}
\begin{remark}
	The Definition \ref{blowing_support_defn} complies with the Definition  \ref{top_support_defn} (cf. Example \ref{blowing_hausdorff_exm}).
\end{remark}

\begin{lemma}\label{blowing_pedersen_compact_lem}
If  $C_0\left( \sY\right)\hookto M\left( A\right)$ is Hausdorff blowing-up and $a\in A$ belongs to the Pedersen's ideal $K\left(A \right)$ (cf. Definition \ref{pedersen_ideal_defn}) then the support of $a$ (cf. Definition \ref{blowing_support_defn}) is compact, i.e.
\be\label{blowing_pedersen_compact_eqn}
K\left(A \right)\subset A_c
\ee 
where $A_c \subset A$ is the {algebra of compactly supported elements} (cf. Definition \ref{blowing_compact_supp_defn}).
\end{lemma}
\begin{proof}
If $a \in K\left(A\right)_0$ (cf \eqref{pedersen_k0_eqn}) then from the Lemma  \eqref{pedersen_eps_lem} it follows that there is $\eps > 0$ and $b \in A_+$ such that $a = f_\eps \left( b\right)$ where $f_\eps$ is given by \eqref{f_eps_eqn}. On the other hand  there is a positive element  $c\in  A_+$  such that $\left\|c - b \right\|  < \eps/2$ and $\supp c$ is compact (cf. Definition  \ref{blowing_defn}).
 If $a \le c$ does not hold and $\rho: A \hookto B\left(\H \right)$ is a faithful  nondegenerate representation (cf. Definitions \ref{faithful_representation_defn}, \ref{nondegenerate_repr_defn}) then there is $\xi \in \H$ such that
\bean\label{blowing_ac_eqn}
\forall \xi \in \H \quad \left( \xi, \rho\left( a \right) \xi\right) > \left( \xi, \rho\left( c \right) \xi\right)
\eean
From $\left\|c - b \right\|  < \eps/2$ if follows that
\be\label{blowing_acc_eqn}
\forall \xi \in \H  \quad \left\| \xi \right\| = 1\quad \Rightarrow\quad  \left|\left( \xi, \rho\left( b \right) \xi\right)-\left( \xi, \rho\left( c \right) \xi\right) \right| < \eps/2
\ee
On the other hand from $a =f_\eps \left( b\right)> 0$ it follows that there is $\xi \in \H$ and $\la \in \R_+$ such that
\bean
\left\| \xi \right\| = 1,\\
\rho\left(a \right) \xi = \la \xi,\\
\rho\left(b \right)\xi = \left( \la+\eps \right) \xi,\\
\rho\left(a \right)\xi = \la  \xi,\\
\left( \xi, \rho\left( b \right) \xi\right) - \left( \xi, \rho\left( a \right) \xi\right)= \eps
\eean
and taking into account \eqref{blowing_acc_eqn} one has
$$
\left( \xi, \rho\left( c \right) \xi\right)- \left( \xi, \rho\left( a \right) \xi\right)> \frac{\eps}{2}
$$
Above condition contradicts with \eqref{blowing_ac_eqn} so $a \le c$.
If $\supp a \subsetneqq \supp c$ then there is a nonempty  open set  $\sU \subset \supp a\setminus \supp c$. For any $f \in C_0\left(\sU \right)\setminus \{0\}$ one has
\bean
faf^* > 0,\\
fcf^* = 0.
\eean
However it is impossible since $a \le c$, so $\supp a \subsetneqq \supp c$ is not true and $\supp a \subset\supp c$. Thus the set $\supp a$ is a closed subset of the compact set $\supp c$ therefore $\supp a$ is compact. Using this fact and the Definition \ref{pedersen_ideal_defn} we conclude that $\supp a$ is compact for any $a \in K\left( A\right)$. 
\end{proof}
\begin{remark}
	The Lemma \ref{blowing_pedersen_compact_lem} can be regarded as a generalization of the equation  \eqref{peder_c0_eqn}.
\end{remark}

\begin{empt}\label{blowing_strict_empt} 
	Suppose that $\left\{u_\la \in C_0\left( \sU\right) \right\}_{\la\in\Xi\left( \sU\right)}$ where $\Xi\left( \sU\right)$ is given by the Exercise \ref{top_u_net_exer} net and let $a \in ~_\sU A_\sU$. For all $a  \in~_\sU A_\sU$  there is $a' \in A$ and $f',f'' \in C_0\left( \sU\right) $ such that $\left\|a - f'a'f'' \right\|\le \frac{\eps}{3}$. On the other hand there is $\la_0 \in \La$  such that $\left\|f' - u_\la f' \right\| < \frac{\eps}{3\left\|a'f'' \right\|}$ for all $\la\ge\la_0$. From the triangle inequality it follows that
	\bean
	\la\ge\la_0\Rightarrow\left\|a - u_\la a \right\|\le \\\le \left\|a - f'a'f''\right\|+ \left\|u_\la \left( a - f'a'f''\right)  \right\|+\left\|f'a'f'' - u_\la f'a'f''\right\| < \eps.
	\eean
	Similarly there is $\la'_0\in \La$ such that $\la\ge\la'_0\Rightarrow\left\|a - au_\la  \right\| < \eps$.
	In result one has the following limit
	\be\label{blowing_strict_eqn}
	\bt\text{-}\lim u_\la = 1_{M\left(_\sU A_\sU\right) }
	\ee
	with respect to the strict topology of $M\left( _\sU A_\sU\right)$ (cf. Definition \ref{strict_topology_defn}).
\end{empt}
\begin{lemma}\label{blowing_multiplier_lem}
	Let $A$ be a $C^*$-algebra (cf. Definition \ref{connected_c_a_defn}),  and let $C_0\left(\sY\right)\hookto M\left( A\right)$  be Hausdorff blowing-up (cf. Definition \ref{blowing_defn}). If $\sU \subset \sY$ is an open subset and $_\sU A_\sU\subset A$ is the $\sU$-subalgebra (cf. Definition \ref{blowing_ideals_au_ua_defn}) then there is a natural inclusion
	$$
	C_b\left(\sU\right)\subset M\left(_\sU A_\sU\right).
	$$
\end{lemma}
\begin{proof}
	If  $f\in C_b\left(\sY\right)$ and $\left\{u_\la \in C_0\left( \sU\right) \right\}_{\la\in\Xi\left( \sU\right)}$ where $\Xi\left( \sX\right)$ is given by the Exercise \ref{top_u_net_exer} net, then from $fu_\la \in C_0\left(\sU\right)$ it turns out that $fu_\la a \in _\sU A_\sU$ for all $\la\in\La$. From
	$\left\| fu_\la a\right\| \le \left\| f\right\|\left\|a\right\|$ it follows that a net $\left\{fu_\la a  \right\}_{\la\in\La}\subset ~_\sU A_\sU$ is $C^*$-norm convergent. So one has a pairing
	\bean
	C_b\left(\sU\right)\times~ _\sU A_\sU\to~ _\sU A_\sU,\\
	(f, a)\mapsto fa \bydef \lim_{\la \in \La} fu_\la a.
	\eean
	Similarly to the equation \eqref{double_centralizer_eqn} there is a map
	\bean 
L: ~	_\sU A_\sU \to ~_\sU A_\sU,\\
a \mapsto fa.
	\eean 
	We leave to the reader details of the proof of the existence of a  double centralizer $\left(L, R\right)$ (cf. Definition \ref{double_centralizer_defn}) which corresponds to $f$. The Proposition \ref{dc_prop} yields a map $C_b\left(\sU\right)\to M\left(_\sU A_\sU\right)$. On can prove that this map is an injective $*$-homomorphism.
\end{proof}

\begin{lemma}\label{blowing_uau_conn_lem}
	If $A$ is a connected $C^*$-algebra (cf. Definition \ref{connected_c_a_defn}) and $C_0\left(\sY\right)\hookto M\left( A\right)$ is Hausdorff blowing-up then an open subset $\sU \subset \sY$ is connected if and only if  \ref{blowing_ideals_au_ua_defn} $\sU$-algebra  $C^*$-subalgebra $_\sU A_\sU\subset A$ is  connected $C^*$-algebra (cf. Definition \ref{connected_c_a_defn}. 
\end{lemma}
\begin{proof}
	If $_\sU A_\sU$ is not connected then one has $_\sU A_\sU= A_1 \oplus A_2$. It follows that  $M\left(_\sU A_\sU\right) = M\left( A_1\right)\oplus M\left(A_2\right)$. From the inclusion $C_0\left(\sU\right)\hookto M\left( _\sU A_\sU\right)$ it follows that
	$C_0\left(\sU \right)= M\left( A_1\right)\cap C_0\left(\sU \right)\oplus M\left( A_2\right)\cap C_0\left(\sU \right)$. Since $C_0\left(\sU \right)$ is a connected $C^*$-algebra one has  $M\left( A_1\right)\cap C_0\left(\sU \right)= \{0\}$ or $M\left( A_2\right)\cap C_0\left(\sU \right)=\{0\}$.  Suppose that $M\left( A_1\right)\cap C_0\left(\sU \right)=  C_0\left(\sU \right)$. Any $a \in~ _\sU A_\sU$ can be uniquely represented by $a = a_1 + a_2$ where $a_1 \in A_1$ and $a_2 \in A_2$. 	  	If  $\left\{u_\la \in C_0\left( \sU\right) \right\}_{\la\in\Xi\left( \sY\right)}$  is a given by the Exercise \ref{top_u_net_exer}  net then $u_\la a_2 = 0$ for all $\la\in\Xi\left( \sY\right)$. On the other hand from \eqref{blowing_strict_eqn} it follows that
	$$
a = \lim_{\la\in\Xi\left( \sU\right)}u_\la a = \lim_{\la\in\Xi\left( \sU\right)}u_\la a_1 + 0 = a_1,	
	$$
i.e. $_\sU A_\sU = A_1$. If	$\sU$ is not connected then there is a disjoint union $\sU = \sU_1 \sqcup \sU_2$ of nonempty sets. From the equations  \eqref{blowing_hereditary_u_eqn} and \eqref{blowing_sue1_eqn} one can deduce that there is a direct sum $_\sU A_\sU =~ _{\sU_1}A_{\sU_1}\oplus ~_{\sU_2} A_{\sU_2}$ of nontrivial $C^*$-algebras, i.e. $_\sU A_\sU$ is not connected.
\end{proof}
\begin{corollary}\label{blowing_connected_comp_cor}
	If $A$ is a connected $C^*$-algebra (cf. Definition \ref{connected_c_a_defn}) and $C_0\left(\sY\right)\hookto M\left( A\right)$ is Hausdorff blowing-up   then any connected component of $A$ (cf. Definition \ref{connected_comp_defn}) is given by
	$$
	_\sU A_\sU = ~_\sU A = A_\sU
	$$
	where $\sU \subset \sY$ is a connected component of $\sY$ (cf. Definition \ref{top_connected_component_defn}).
\end{corollary}

\begin{lemma}\label{blowing_ajk_lem}
	Let $A$ be a $C^*$-algebra with Hausdorff blowing-up $C_0\left(\sY\right)\to M\left(A \right)$. If a family $\left\{\sU_\a\right\}_{\a \in \mathscr A}$ of open subsets of $\sY$ is such that $\sY = \bigcup_{\a \in \mathscr A}\sU_\a$ then any element $a \in A_c$ of  {algebra of compactly supported elements} (cf. Definition \ref{blowing_compact_supp_defn}) is given by
\be\label{blowing_ajk_eqn}
a = \sum_{\substack{j=1\\ k = 1}}^n f_ja f_k \quad \forall j = 1,..., n\quad \exists \a \in \mathscr A\quad \supp f_j \subset \sU_\a
\ee
\end{lemma}
\begin{proof}
From the lemma \ref{blowing_pedersen_compact_lem} it follows that the support $\supp a$ (cf. Definition \ref{blowing_support_defn}) of $a$ is compact. If $f \bydef \sum_{j=1}^n f_j$ is a dominated to the family $\left\{\sU_\a\right\}_{\a \in \mathscr A}$ covering sum for $\supp a$ (cf. Definition \ref{top_covering_sum_defn}) then one has
$$
a = faf = \sum_{\substack{j=1\\ k = 1}}^n f_ja f_k.
$$
\end{proof}

\section{Hausdorff blowing-up and noncommutative  coverings}
\subsection{Lift and descent}
\paragraph{}
Here we suppose that we have Hausdorff blowing-us and investigate conditions when they enable us to construct noncommutative coverings.
\begin{empt}\label{blowing_lift_empt}
	Let $A$ be a $C^*$-algebra, and let $C_0\left(\sY \right)\subset  M\left(A \right)$ be  Hausdorff blowing-up (cf. Definition \ref{blowing_defn}). Let $q: \widetilde \sY\to \sY$ be a transitive covering such that $\sY = \widetilde \sY/G$ where $G$ is properly   discontinuous,  residually finite group  of homeomorphisms of $\widetilde \sY$ (cf. Definitions \ref{top_properly_disc_defn} and \ref{residually_finite_defn}). There is $C_0\left(\sY \right)$-valued product
	\bean
	\left\langle\cdot, \cdot  \right\rangle_{C_0\left(\sY \right)}: C_c \left(\widetilde \sY \right) \times C_c \left(\widetilde \sY \right)\to C_0\left(\sY \right),\\
	\left(\widetilde a, \widetilde b \right) \mapsto \desc^c_q\left(\widetilde a^* \widetilde b \right) 
	\eean
	where $\desc^c_q$ is  a {compactly supported} $q$-{descent} \ref{top_compactly_supported_descent_defn}. So $ C_c \left(\widetilde \sY\right)$ is a pre-Hilbert $C_0\left(\sY \right)$-module (cf. Definition \ref{hilbert_module_defn}). If $\mathscr L^2\left(C_0 \left(\widetilde \sY \right) \right)$ is a completion of $C_c \left(\widetilde \sY \right)$ with respect to given by \eqref{hilbert_module_norm_eqn} norm then $\mathscr L^2\left(C_0 \left(\widetilde \sY \right) \right)$ is a $C^*$-Hilbert $C_0\left(\sY \right)$-module. If $\mathscr L^2\left(C_0 \left(\widetilde \sY \right) \right)\otimes_{C_0\left(\sY \right)} A$ is an algebraic tensor product then there is an $A$-valued product.
	\bean
	\left\langle\cdot, \cdot  \right\rangle: \left( \mathscr L^2\left(C_c \left(\widetilde \sY \right) \right)\otimes_{C_0\left(\sY \right)} A\right)\times \left( \mathscr L^2\left(C_c \left(\widetilde \sY \right) \right)\otimes_{C_0\left(\sY \right)} A\right) \to A,\\
	\left(\left(\widetilde f \otimes a \right) , \left(\widetilde g  \otimes b \right)  \right) \mapsto a^*\left\langle\widetilde f, \widetilde g \right\rangle_{C_0\left(\sY \right)} b.
	\eean
	so $\mathscr L^2\left(C_0 \left(\widetilde \sY \right) \right)\otimes_{C_0\left(\sY \right)} A$ becomes a pre-Hilbert $A$-module. We denote by $\mathscr L^2\left(\widetilde \sY \right) _A$ the $C^*$-Hilbert $A$-module which is a completion  of  $\mathscr L^2\left(C_0 \left(\widetilde \sY \right) \right)\otimes_{C_0\left(\sY \right)} A$ with respect to given by \eqref{hilbert_module_norm_eqn} norm. There are both a left  action 	$C_0 \left(\widetilde \sY \right)\times \mathscr L^2\left(\widetilde \sY \right) _A\to \mathscr L^2\left(\widetilde \sY \right) _A$ and the right action $\mathscr L^2\left(\widetilde \sY \right) _A\times A \to \mathscr L^2\left(\widetilde \sY \right) _A$ such that $\mathscr L^2\left(\sY\right)_A$ is $C_0 \left(\widetilde \sY \right)$-$A$-bimodule. These actions induce injective homomorphisms of $C^*$-algebras
	\bea\label{blowing_a_h_eqn}
	\varphi_A : A \hookto \End^*_A\left( \mathscr L^2\left(\widetilde \sY \right)_A\right), \\
	\label{blowing_c_h_eqn}
	\varphi_C : C_0 \left(\widetilde \sY \right) \hookto \End^*_A\left( \mathscr L^2\left(\widetilde \sY \right)_A\right)
	\eea
	where $\End^*_A\left( \mathscr L^2\left(\widetilde \sY \right)_A\right)$ is the $C^*$-algebra of adjointable  endomorphisms of  the $C^*$-Hilbert $A$-module $\mathscr L^2\left(\widetilde \sY \right)_A$ (cf. Definition \ref{adjointable_operator_defn}).
	There is a homomorphism 
	\be\label{blowing_tensor_eqn}
	\begin{split}
		\varphi: C_0\left( \widetilde \sY \right)\otimes_{C_0\left(\sY\right) } K\left( A\right) \to \End^*_A\left( \mathscr L^2\left(\widetilde \sY \right)_A\right),\\
		\sum_{j=1}^n 	\widetilde f_j \otimes a_j \mapsto 	\sum_{j=1}^n 	\varphi_C\left( \widetilde f_j\right) \varphi_A\left(a_j \right)
	\end{split}
	\ee
	of right $C_0\left( \widetilde\sY\right)$-$A$-bimodules, where 	$\varphi_A, ~	\varphi_C$ are given by  \eqref{blowing_a_h_eqn}, \eqref{blowing_c_h_eqn} and the algebraic tensor product is implied. Since both $C_c\left(\widetilde \sY \right)$ and   $K\left(A \right)$ are dense in $C_0\left(\widetilde \sY \right)$ and  $A$ respectively the -homomorphism \eqref{blowing_tensor_eqn} can be uniquely extended up to a  $C_0\left( \widetilde\sX\right)$-$A$-bimodule homomorphism 
	\be\label{blowing_tensoar_eqn}
	\begin{split}
		\varphi': C_c\left( \widetilde \sY \right)\otimes_{C_0\left(\sY\right) } A \to \End^*_A\left( \mathscr L^2\left(\widetilde \sY \right)_A\right)
	\end{split}
	\ee
	such that the $C^*$ closure of $\varphi\left(  C_c\left( \widetilde \sY \right)\otimes_{C_0\left(\sY\right) } K\left( A\right)\right)$ coincides with\\  	$\varphi'\left(  C_c\left( \widetilde \sY \right)\otimes_{C_0\left(\sY\right) } A\right) $ one.
\end{empt}

\begin{definition}\label{blowing_lift_hm_defn}
	Under the hypotheses \ref{blowing_lift_empt}  we say that the $C^*$-Hilbert $A$-module  $\mathscr L^2\left(\widetilde \sY \right)_A$ is the $\widetilde \sY$-$A$-\textit{module}.
\end{definition}
\begin{definition}\label{blowing_a_regular_defn}
	Under the hypotheses  \ref{blowing_lift_empt}
if 	$C_0\left( \widetilde\sY\right)$-$A$-bimodule
	\be\label{blowing_prods_eqn}
	\varphi\left(C_c\left( \widetilde \sY \right)\otimes_{C_0\left(\sY\right) } K\left( A\right) \right)\subset  \End^*_A\left( \mathscr L^2\left(\widetilde \sY \right)_A\right)
	\ee	
	(cf. equation \eqref{blowing_tensor_eqn}) is a $*$-subalgebra of $\End^*_A\left( \mathscr L^2\left(\widetilde \sY \right)_A\right)$
	then  we say that the covering $q: \widetilde{\sY}\to \sY$ is $A$-\textit{regular}.  
\end{definition}

\begin{lemma}\label{blowing_lift_eq_lem}
If $C_0\left(\sY \right)\subset  M\left(A \right)$ is  Hausdorff blowing-up (cf. Definition \ref{blowing_defn}) then a covering $q: \widetilde{\sY}\to \sY$ is $A$-{regular} (cf. Definition \ref{blowing_a_regular_defn}) if and only if for any $\widetilde a \in \varphi\left(C_c\left( \widetilde \sY \right)\otimes_{C_0\left(\sY\right) } K\left( A\right) \right)$ there is $\widetilde f \in C_c\left(\widetilde\sY\right)$ such that
\be\label{blowing_lift_eq_eqn}
\forall \widetilde a \in \varphi\left(C_c\left( \widetilde \sY \right)\otimes_{C_0\left(\sY\right) } K\left( A\right) \right) \quad \exists  \widetilde f \in C_c\left(\widetilde\sY\right) \quad \varphi_A\left(  \widetilde a\right)  =\varphi_A\left(  \widetilde a\right) \varphi_C\left(\widetilde f \right). 
\ee 
\end{lemma}
\begin{proof}
If $\widetilde a \in \varphi\left(C_c\left( \widetilde \sY \right)\otimes_{C_0\left(\sY\right) } K\left( A\right) \right)$ then under the hypotheses of the Definition \ref{blowing_a_regular_defn} one has $\widetilde a^* \in \varphi\left(C_c\left( \widetilde \sY \right)\otimes_{C_0\left(\sY\right) } K\left( A\right) \right)$ and
\bean
\widetilde a^* = \varphi\left(\sum_{j = 1}^n \widetilde f_j \otimes a_j \right),\\
\widetilde a = \sum_{j = 1}^n \varphi_A \left(a^*_j\right)\varphi_C\left(\widetilde f_j^* \right)   
\eean
The finite union $\widetilde \sV \bydef \bigcup_{j = 1}^n \supp \widetilde f_j$ is compact. If $\widetilde f\in C_c\left(\widetilde\sY \right) $ is a covering sum for $\widetilde \sV$ (cf. Definition \ref{top_covering_sum_defn}) then $\varphi_A\left(  \widetilde a\right) = \varphi_A\left(  \widetilde a\right) \varphi_C\left(f \right)$.
\end{proof}
\begin{definition}\label{blowing_lift_defn}
	If 	$C_0\left( \sY\right)\hookto M\left( A\right) $ is  Hausdorff blowing-up (cf. Definition \ref{blowing_defn}), and   $q: \widetilde{\sY}\to \sY$ is an   $A$-{regular} covering (cf. Definition \ref{blowing_a_regular_defn}), then the $C^*$-norm completion $A_0\left(\widetilde\sY\right)$ of $\varphi\left(C_c\left( \widetilde \sY \right)\otimes_{C_0\left(\sY\right) } K\left( A\right) \right)$ (cf. equation \eqref{blowing_prods_eqn}) is the $q$-\textit{lift of} $A$.
\end{definition}
\begin{lemma}\label{blowing_lift_constr_lem}
	Let 	$C_0\left( \sY\right)\hookto M\left( A\right) $ be  Hausdorff blowing-up (cf. Definition \ref{blowing_eqn}), and let   $q: \widetilde{\sY}\to \sY$ be a  transitive, $A$-{regular} covering (cf. Definition \ref{blowing_a_regular_defn}). If $A_0\left(\widetilde\sY\right)$ is the  $q$-\textit{lift of} $A$ (cf. Definition \ref{blowing_lift_defn}) then there is a natural injective $*$-homomorphism
	\be\label{blowing_lift_eqn}
	\begin{split}
		A_b\left(q \right) : A \hookto M\left( A_0\left(\widetilde\sY\right)\right).
	\end{split}
	\ee
\end{lemma}
\begin{proof}
	If  $a \in A$ and $\widetilde{a} \in A_0\left(\widetilde\sY\right)$ then for any $\eps > 0$ there is $\widetilde{a}' \in \varphi\left(C_c\left( \widetilde \sY \right)\otimes_{C_0\left(\sY\right) } K\left( A\right) \right)$ (cf. equation \eqref{blowing_prods_eqn}) such that
\be\label{blowing_le_eqn}
	\left\| \widetilde{a} - \widetilde{a}'\right\| < \frac{\eps}{\left\|a \right\| }.
\ee
	From the Lemma \ref{blowing_lift_eq_lem} it follows that there is a positive function $\widetilde f \in C_c\left(\widetilde{\sY} \right)_+$ such that 
	$$
	\widetilde{a}' = \varphi_C\left(\widetilde f \right)\widetilde{a}'\varphi_C\left(\widetilde f \right)
	$$
	If $ \left\{\widetilde u_\la\right\}_{\la \in \Xi\left( \sY\right)}\subset C_c\left( \sY\right)_+$ is a given by the Exercise \ref{top_u_net_exer} approximate unit of $C_0\left(\widetilde{\sY} \right)$ then there is $\la' \in \La$ such that 
	$$\forall \la \in \La	\quad \la \ge \la' \quad \Rightarrow \quad\widetilde u_\la\left(\supp \widetilde f \right)= \{1\}.$$ 
	It turns out that  
\bean
\forall \la' \in \La	\quad \la \ge \la' \quad \Rightarrow \quad   \varphi_C\left(  \widetilde{u}_\la\right)   \widetilde a' =   \widetilde a' \quad \Rightarrow \\
\Rightarrow\quad \varphi_A\left( a\right)  \varphi_C\left(  \widetilde{u}_\la\right) \widetilde{a}'= \varphi_A\left( a\right)  \varphi_C\left(  \widetilde{u}_{\la'}\right) \widetilde{a}'\in  \varphi\left(C_c\left( \widetilde \sY \right)\otimes_{C_0\left(\sY\right) } K\left( A\right) \right).
\eean
Taking into account the Lemma \ref{blowing_lift_eq_lem}  one has a positive function $\widetilde f'' \in C_c\left(\widetilde{\sY} \right)_+$ such that 
	$$
\forall \la' \in \La	\quad \varphi_A\left( a\right)  \varphi_C\left(  \widetilde{u}_{\la'}\right) \widetilde{a}'= \varphi_C\left(  \widetilde{f}''\right) \varphi_A\left( a\right)  \varphi_C\left(  \widetilde{u}_{\la'}\right) \widetilde{a}'
	$$. 
There is $\la'' \in \La$ such that $$\forall \la \in \La	\quad \la \ge \la' \quad \Rightarrow \quad\widetilde u_\la\left(\supp \widetilde f'' \right)= \{1\}.$$
So if $\la_0\in \La$ is such that $\la_0 \ge \la'$ and $\la_0\ge \la''$ then
\bean
\forall \la \in \La	\quad \la \ge \la_0 \quad \Rightarrow\quad \varphi_A\left( a\right)  \varphi_C\left(  \widetilde{u}_\la\right) \widetilde{a}'= \varphi_C\left(  \widetilde{f}''\right)\varphi_A\left( a\right)  \varphi_C\left(  \widetilde{u}_{\la}\right) \widetilde{a}'=\\=  \varphi_C\left(  \widetilde u_\la\right)\varphi_A\left( a\right)  \varphi_C\left(  \widetilde{u}_{\la}\right) \widetilde{a}'= \varphi_C\left(  \widetilde u_{\la_0}\right)\varphi_A\left( a\right)  \varphi_C\left(  \widetilde{u}_{\la_0}\right) \widetilde{a}'.
\eean
Taking into account $\left\|\widetilde u_{\la_0} \right\|= 1$ and the inequality \eqref{blowing_le_eqn} on has
$$
\forall \la \in \La	\quad \la \ge \la_0 \quad \Rightarrow\quad \left\| \varphi_C\left(  \widetilde u_{\la}\right)\varphi_A\left( a\right)  \varphi_C\left(  \widetilde{u}_{\la}\right) \widetilde{a}- \varphi_C\left(  \widetilde u_{\la_0}\right)\varphi_A\left( a\right)  \varphi_C\left(  \widetilde{u}_{\la_0}\right) \widetilde{a} \right\|< \eps
$$
Since the number $\eps$ can be arbitrary small we conclude that the net
$$\left\{\varphi_C\left(  \widetilde u_{\la}\right)\varphi_A\left( a\right)  \varphi_C\left(  \widetilde{u}_{\la}\right) \widetilde{a}\right\}_{\la\in \La} $$ is $C^*$-norm convergent. We define a map
\bean
L_a :  A_0\left(\widetilde\sY\right)\to A_0\left(\widetilde\sY\right),\\
\widetilde{a}\mapsto \lim_{\la\in \La} \varphi_C\left(  \widetilde u_{\la}\right)\varphi_A\left( a\right)  \varphi_C\left(  \widetilde{u}_{\la}\right) \widetilde{a}
\eean
Similarly there is a map
\bean
R_a :  A_0\left(\widetilde\sY\right)\to A_0\left(\widetilde\sY\right),\\
\widetilde{a}\mapsto \lim_{\la\in \La}\widetilde{a} \varphi_C\left(  \widetilde u_{\la}\right)\varphi_A\left( a\right)  \varphi_C\left(  \widetilde{u}_{\la}\right). 
\eean
From the equality
$$
\forall \la \in \La\quad \forall \widetilde b',  \widetilde b'' \quad \widetilde b'\left( \varphi_C\left(  \widetilde u_{\la}\right) \varphi_A\left( a\right)  \varphi_C\left(  \widetilde{u}_{\la}\right) \widetilde b''\right) = \left( \widetilde b'\varphi_C\left(  \widetilde u_{\la}\right)\varphi_A\left( a\right)  \varphi_C\left(  \widetilde{u}_{\la}\right)\right)  \widetilde b''
$$	
we conclude that a pair  $\left(L_a, R_a\right)$ of maps  is a {double centralizer} (cf. Definition \ref{double_centralizer_defn}).From the Remark \ref{double_centralizer_rem} it follows that there is a map 
\bean
A \to M\left( A_0\left(\widetilde\sY\right)\right),\\
a \mapsto \text{ the corresponding to }\left(L_a, R_a\right)\text{ multiplier}.
\eean 
We leave to the reader an elementary proof of that the map is an injective $*$-homomorphism.	
	\end{proof}

\begin{definition}\label{blowing_lift_hom_defn}  If	$C_0\left( \sY\right)\hookto M\left( A\right) $ is  Hausdorff blowing-up (cf. Definition \ref{blowing_defn}), and    $q: \widetilde{\sY}\to \sY$ is an   $A$-{regular} covering (cf. Definition \ref{blowing_a_regular_defn}) then the given by the Lemma \eqref{blowing_lift_constr_lem} $*$-homomorphism $	A_b\left(q \right) : A \hookto M\left( A_0\left(\widetilde\sY\right)\right)$
	is the  $q$-\textit{lift} of $A$. The difference between both  Definitions  \ref{blowing_lift_defn} and \ref{blowing_lift_hom_defn} of $q$-lift will be explained as we go along.
\end{definition}

\begin{lemma}\label{blowing_mult_iclusion_lem}
	If	$C_0\left( \sY\right)\hookto M\left( A\right) $ is  Hausdorff blowing-up (cf. Definition \ref{blowing_defn}), and    $q: \widetilde{\sY}\to \sY$ is an   $A$-{regular} covering (cf. Definition \ref{blowing_a_regular_defn}) then the  $q$-{lift} $	A_b\left(q \right) : A \hookto M\left( A_0\left(\widetilde\sY\right)\right)$ of $A$ (cf. Definition \ref{blowing_lift_hom_defn}) can be uniquely extended up to an injective $*$-homomorphism
	\be\label{blowing_lift_m_eqn}
	M\left( A_b\left(q \right)\right)  : M\left( A\right)  \hookto M\left( A_0\left(\widetilde\sY\right)\right).
	\ee 
\end{lemma}
\begin{proof}
	Let $\left\{u_\la\right\}_{\la \in \La}$ be an approximate unit of $A$ (cf. Definition \ref{approximate_unit_defn}), and let $\widetilde a \in A_0\left(\widetilde\sY\right)$. Form the Definition \ref{blowing_lift_defn} it follows that for any $\eps > 0$ there is $\widetilde a' \in  \varphi\left(C_c\left( \widetilde \sY \right)\otimes_{C_0\left(\sY\right) } K\left( A\right) \right)$  such that
	$$
	\left\| \widetilde a - \widetilde a'\right\| 	< \frac{\eps}{3}.
	$$
	On the other hand there are two  families  $\left\{\widetilde f_1, ..., \widetilde f_n\right\}\subset C_c\left(\widetilde \sY \right)$ and $\left\{a_1, ..., a_n\right\}$ such that 
	$$
	\widetilde a' = \varphi\left( \sum_{j=1}^n \widetilde f_j \otimes a_j\right) \sum_{j=1}^n \phi_C\left( \widetilde f_j\right) \phi_A\left(a_j \right). 	
	$$
	On the other hand for any $j=1,...,n$ there is $\la_j \in \la$ such that
	$$
	\forall \la\in \La \quad \la \ge \la_j \quad \Rightarrow \quad 	\left\|a_j - a_ju_\la \right\| < \frac{\eps}{3n \left\|\widetilde f_j\right\|}
	$$
	If $\la_0\in \La$ is such that $\la_O \ge \la_j$ for every $j=1,...,n$ then from the triangle identity it follows that
	\bean
	\forall \la\in \La \quad \la \ge \la_0 \quad  \left\|\widetilde a - \widetilde a\varphi_A\left( u_\la\right) \right\|\le \\\le  \left\| \widetilde a - \widetilde a'\right\|+\left\|\left( \widetilde a - \widetilde a'\right) \varphi_A\left( u_\la\right)\right\|+\left\|\widetilde a' - \widetilde a'\varphi_A\left( u_\la\right)\right\|<\\< \frac{2\eps}{3} + \sum_{j = 1}^n \left\|\widetilde f_j\right\|\left\|a_j - a_ju_\la \right\|< \eps,
	\eean 
	so one has
	$
	\lim_{{\la \in \La}}\left\|\widetilde a - \widetilde a\varphi_A\left( u_\la\right) \right\|= 0
	$. Similarly one can prove that \\ $\lim_{{\la \in \La}}\left\|\widetilde a - \varphi_A\left( u_\la\right)\widetilde a \right\|= 0$. This Lemma becomes  a consequence of the \ref{lift_mult_lem} one.
\end{proof}

\begin{lemma}\label{blowing_blowing_lem}
	Let	$C_0\left( \sY\right)\hookto M\left( A\right) $ be  Hausdorff blowing-up (cf. Definition \ref{blowing_defn}), and let  $q: \widetilde{\sY}\to \sY$ be an   $A$-{regular} covering (cf. Definition \ref{blowing_a_regular_defn}).
	If   $A_0\left(\widetilde\sY\right)$  is the $q$-{lift of} $A$ then the natural inclusion $C_0\left( \widetilde\sY\right)\hookto M\left( A_0\left(\widetilde\sY\right)\right)$  is Hausdorff  blowing-up (cf. Definition \ref{blowing_defn}). 
\end{lemma}
\begin{proof}
Form the Definition \ref{blowing_lift_defn} it follows that for any $\eps > 0$ there is\\ $\widetilde a' \in  \varphi\left(C_c\left( \widetilde \sY \right)\otimes_{C_0\left(\sY\right) } K\left( A\right) \right)$  such that
$$
\left\| \widetilde a - \widetilde a'\right\| 	< \eps .
$$
On the other hand there are two  families  $\left\{\widetilde f_1, ..., \widetilde f_n\right\}\subset C_c\left(\widetilde \sY \right)$ and $\left\{a_1, ..., a_n\right\}$ such that 
$$
\widetilde a' = \varphi\left( \sum_{j=1}^n \widetilde f_j \otimes a_j\right) \sum_{j=1}^n \phi_C\left( \widetilde f_j\right) \phi_A\left(a_j \right). 	
$$
A finite union $\widetilde \sV \bigcup_{j=1}^n  \supp \widetilde f_j$ of compact sets is compact. If $\widetilde f \in C_c\left( \widetilde \sY\right)$ is a covering sum for  $\widetilde \sV$ (cf. Definition \ref{top_covering_sum_defn}) then $\widetilde a' = \widetilde f \widetilde a'$ and 
\bean
\left\| \widetilde a -\widetilde f \widetilde a'\right\| 	< \eps, 
\eean
i.e. a set $C_c\left( \widetilde\sY\right)A_0\left(\widetilde\sY\right) \bydef \left\{\widetilde f\widetilde a \left| \widetilde f \in C_c\left( \widetilde\sY\right)\quad\widetilde a \in A_0\left(\widetilde\sY\right)\right. \right\}$ is dense in $A_0\left(\widetilde\sY\right)$. Taking into account the Remark \ref{blowing_rem} we conclude that $C_0\left( \widetilde\sY\right)\hookto M\left( A_0\left(\widetilde\sY\right)\right)$  is Hausdorff  blowing-up.
\end{proof}
Let 	$C_0\left( \sY\right)\hookto M\left( A\right) $ be  Hausdorff blowing-up (cf. Definition \ref{blowing_defn}), and  let $q: \widetilde{\sY}\to \sY$ be  an   $A$-{regular} covering (cf. Definition \ref{blowing_a_regular_defn}). Let  
$$
\left(\	\left\{\left(\sU_\a, \sV_\a, s_\a\right)\right\}_{\a \in \mathscr A}, \left\{\left(\widetilde \sU_{\widetilde\a}, \widetilde \sV_{\widetilde\a}, \widetilde s_{\widetilde\a}\right)\right\}_{\widetilde \a \in \widetilde{\mathscr A}} \right)$$ 
be  a $q$-{covering} (cf. Definition \ref{top_cov_defn}).  If 
$$
\widetilde a \bydef \varphi\left(\widetilde f \otimes a \right) \in \varphi\left(C_c\left( \widetilde \sY \right)\otimes_{C_0\left(\sY\right) } K\left( A\right) \right)
$$ 
then from the Lemma \ref{blowing_lift_eq_lem} it turns out that  there is $\widetilde f'\in  C_c\left( \widetilde \sY \right)$ such that
 $$
 \widetilde a = \widetilde a \varphi_C\left(\widetilde f' \right)= \varphi_C\left(\widetilde f \right)\varphi_A\left( a\right) \varphi_C\left(\widetilde f' \right)
 $$
From the Lemma \ref{blowing_ajk_lem}  it follows that 
 \bean
 a = \sum_{\substack{j=1\\ k = 1}}^n f_ja f_k \quad \forall j = 1,..., n\quad \exists \a \in \mathscr A\quad \supp f_j \subset \sU_\a
 \eean
 Select $j \in \left\{1,..., n\right\}$ and  $\a_j\in \mathscr A$ such that $\supp f_j \subset \sU_{\a_j}$. Under the hypotheses of \ref{top_cov_empt} select $\widetilde \a_j\in\widetilde{\mathscr A}$ such that $p_{\widetilde{\mathscr A}}\left(\widetilde\a_j \right) = \a_j$ where $p_{\widetilde{\mathscr A}}: \widetilde{\mathscr A} \to {\mathscr A}$ is the  $p$-{projection} (cf. Definition \ref{top_cov_defn}).
  The set
 $$
 G_j \bydef \left\{\left.g \in G\right| g \widetilde \sU_{\widetilde a_j} \cap \left( \supp \widetilde f\cup \supp  \widetilde f'   \right)\neq \emptyset  \right\}
. $$
 is finite since the union $\supp \widetilde f\cup \supp  \widetilde f'$ is compact. From our construction it follows that
\be\label{blowing_lift_eqp_eqn}
 \widetilde a = \sum_{\substack{j=1\\ k = 1}}^n\varphi_C\left(\sum_{g' \in G_j} \lift^q_{g' \widetilde \sU_{\widetilde a_j} }\left(f_j \right) \widetilde f \right) \varphi_A\left( a\right)  \varphi_C\left(\sum_{g'' \in G_k} \lift^q_{g'' \widetilde \sU_{\widetilde a_k} }\left(f_k \right)\widetilde f'\right) 
\ee
where $\lift^q_{\widetilde \sU}$ is a $q$-$\widetilde \sU$-lift (cf. Definition \ref{top_lift_desc_defn}). From \eqref{blowing_lift_eqp_eqn} one can deduce that any element  $\widetilde a \in \varphi\left(C_c\left( \widetilde \sY \right)\otimes_{C_0\left(\sY\right) } K\left( A\right) \right)$ can be represented by following finite sum, i. e.
\be\label{blowing_tajk_eqn}
\begin{split}
\widetilde a=  \sum_{\substack{j=1\\ k = 1}}^m\varphi_C\left(\widetilde f'_j \right) \varphi_A\left(a_{jk}\right)\varphi_C\left(\widetilde f''_k \right)\quad a_{jk} \in K\left( A\right), \\ 
 \forall j,k \in \left\{1,...,m\right\}\quad \exists \widetilde \sU'_j , \widetilde \sU''_k \in \left\{\widetilde \sU_{\widetilde \a}\right\}\quad \supp \widetilde f'_j\in \widetilde \sU'_j \quad \supp \widetilde f''_k\in \widetilde \sU''_k.
 \end{split}
\ee

\begin{lemma}\label{blowing_left_ideal_lem} 
	Let 	$C_0\left( \sY\right)\hookto M\left( A\right) $ be  Hausdorff blowing-up (cf. Definition \ref{blowing_defn}), and let  $q: \widetilde{\sY}\to \sY$ be  an   $A$-{regular} covering (cf. Definition \ref{blowing_a_regular_defn}) with $\sY \cong \widetilde{\sY}/G$ where $G$ is a properly disconnected group of homeomorphisms of $\widetilde{\sY}$ (cf. Definition \ref{top_properly_disc_defn}). If $\left(\widetilde \sU, \widetilde \sV, \widetilde s\right)$ is   a {covering triple for} $q$ (cf. Definition \ref{top_coveing_triple_defn}), and  $\sU \bydef q\left(\widetilde\sU \right)$. 
	then there is a natural isomorphism of left $A$-modules
	\be
	\begin{split}
		\pi_{A_0\left(\widetilde\sY \right)_{\widetilde \sU}}:  A_{\sU}\xrightarrow{\cong} A_0\left(\widetilde\sY \right)_{\widetilde \sU}\\
			a \mapsto A_0\left(q \right)\left(   a \right) \widetilde s
		\end{split} 
		\ee
		where	 both $A_{\sU}$ and $A_0\left(\widetilde\sY \right)_{\widetilde \sU}$ are left ${\sU}$- and $\widetilde \sU$-ideals respectively (cf. Definition \ref{blowing_ideals_au_ua_defn}). The isomorphism $\pi_{A_0\left(\widetilde\sY \right)_{\widetilde \sU}}$ preserves $C^*$-norm.
	\end{lemma}
	\begin{proof}
If $a\in A_{\sU}$ and $\widetilde a \bydef \pi_{A_0\left(\widetilde\sY \right)_{\widetilde \sU}}\left(a \right) \in  A_0\left(\widetilde\sY \right)_{\widetilde \sU}$ then $a^*a \in  ~_{\sU}A_{\sU}$ and $\widetilde a^*\widetilde{a} \in ~_{\widetilde \sU}A_0\left(\widetilde\sY \right)_{\widetilde \sU}$ where the notation of the Definition \ref{blowing_ideals_au_ua_defn} is used. If $\widetilde b \in A_0\left(\widetilde\sY \right)$ then from the Lemma \ref{blowing_blowing_lem} and the Definition \ref{blowing_a_h_eqn} it follows that for any $\eps > 0$ there is $\widetilde \varphi_\eps \in C_c\left(\widetilde \sY \right)$ such that
$$
\left\|\widetilde b -  \widetilde \varphi_\eps \widetilde b \right\| < \frac{\eps}{2\left\|a^*a \right\|}.
$$
Since both sets $\supp \widetilde s$ and $\supp \widetilde \varphi_\eps$ are compact the set
$$
G_0\bydef \left\{g \in G \left|\supp \widetilde s\cap \supp \widetilde \varphi_\eps \neq \emptyset \right.\right\}
$$
is finite. If follows that for any finite subset $G'\subset G$ one has
$$
G_0\subset G' \quad \Rightarrow \quad \sum_{g \in G_0} \left( g\widetilde a^* \widetilde a\right) \widetilde \varphi_\eps\widetilde b= \sum_{g \in G'} \left( g \widetilde a^*\widetilde a\right) \widetilde \varphi_\eps\widetilde b,
$$
it turns out that
\bean
G_0\subset G' \quad \Rightarrow \quad \left\|\sum_{g \in G'} \left( g \widetilde a^*\widetilde a\right) \widetilde b - a^*a \widetilde b\right\|\le \left\| a^*a  \widetilde \varphi_\eps \widetilde b  - a^*a \widetilde b\right\|+\\+ \left\|\sum_{g \in G'} \left( g \widetilde a^*\widetilde a\right)\left( \widetilde \varphi_\eps  \widetilde b - \widetilde b\right) \right\|+\left\|\sum_{g \in G'} \left( g \widetilde a^*\widetilde a\right) \widetilde \varphi_\eps \widetilde b - a^*a \widetilde \varphi_\eps \widetilde b\right\|< \\< \frac{\eps}{2\left\|a^*a \right\|}\left\|a^*a \right\|+ \frac{\eps}{2\left\|a^*a \right\|}\left\|a^*a \right\|+ 0 = \eps.
\eean
So the series $\sum_{g \in G} \left( g \widetilde a^*\widetilde a\right) \widetilde b$ is $C^*$-norm convergent. Similarly one can prove that the series  $\sum_{g \in G} \widetilde b\left( g \widetilde a^*\widetilde a\right) $ is also $C^*$-norm convergent, so one has
$$
a^*a = \b\text{-}\sum_{g \in G} \widetilde b\left( g \widetilde a^*\widetilde a\right)
$$
where the sum of the series implies the strict topology of $A_0\left(\widetilde\sY \right)$. So $\left\|a^*a \right\| > 0 \quad \Rightarrow \quad \left\|\widetilde a^*\widetilde a \right\|> 0$ and $a \neq 0\quad  \Rightarrow  \quad \widetilde a = 0$, i.e. the homomorphism $\pi_{A_0\left(\widetilde\sY \right)_{\widetilde \sU}}$ is injective. If $\widetilde \pi: A_0\left(\widetilde\sY \right)\hookto B\left( \widetilde \H\right) $ is a faithful, nondegenerate, $G$-equivariant  representation  (cf. Definitions \ref{faithful_representation_defn}, \ref{nondegenerate_repr_defn} and \ref{equivariant_representation_defn}) and $\chi$ is a characteristic function of $\widetilde \sU$ (cf. Definition \ref{top_char_f_defn}) then $\chi$ naturally yields a projector $p_\chi \in B\left(\widetilde{   \H} \right)$. Moreover 
\bean
\widetilde \pi \left(a^* a\right) \left( \widetilde \H \ominus \left(\bigoplus_{g\in G} g p_\chi \widetilde \H\right) \right) = \{0\},\\
\forall g \in G \quad \widetilde \pi \left( a\right) g p_\chi \widetilde \H\subset g p_\chi \widetilde \H
\eean
where $\bigoplus$ is the norm completion of the algebraic direct sum. For any $\eps > 0$ there is
$$
\xi = \sum_{g \in G} \xi_g \quad \quad \frac{\left\|  \widetilde \pi \left(a^* a\right) \xi\right\|^2}{\left\|   \xi\right\|^2}=\frac{\sum_{g \in G}\left\|  \widetilde \pi \left(a^* a\right) \xi_g\right\|^2  }{\sum_{g \in G}\left\|   \xi_g\right\|^2}> \left\| a^*a\right\|^2 - \eps, \quad \xi_g \in p_\chi \widetilde \H
$$
On the other hand for all $g \in G$ one has
\bean
\widetilde \pi \left(a^* a\right) \xi_g = \widetilde \pi \left(g\left( \widetilde a^* \widetilde a\right)\right)  \xi_g,\\
\frac{\left\| \widetilde \pi \left(g\left( \widetilde a^* \widetilde a\right)\right) \xi_g\right\|^2  }{\left\|   \xi_g\right\|^2}\le \left\| g\left( \widetilde a^* \widetilde a\right) \right\|^2= \left\|  \widetilde a^* \widetilde a \right\|^2.
\eean
So one has
$$
\frac{\left\|  \widetilde \pi \left(a^* a\right) \xi\right\|^2}{\left\|   \xi\right\|^2}\le \frac{\sum_{g \in G}\left\|  \widetilde a^* \widetilde a \right\|^2\left\|   \xi_g\right\|^2  }{\sum_{g \in G}\left\|   \xi_g\right\|^2}\le  \left\|  \widetilde a^* \widetilde a \right\|^2,
$$
i.e. $ \left\| a^*a\right\|^2 - \eps \le \left\|  \widetilde a^* \widetilde a \right\|^2$. On the other hand from $a^*a > \widetilde a^* \widetilde a$ it follows that $\left\| a^*a\right\|^2 \ge  \left\|  \widetilde a^* \widetilde a \right\|^2$ it follows that $\left\| a^*a\right\|=\left\|  \widetilde a^* \widetilde a \right\|$. So one has $\left\| a\right\|= \sqrt{\left\| a^*a\right\|}= \sqrt{\left\| \widetilde a^*\widetilde a\right\|}= \left\| \widetilde a\right\|$, i.e. $\pi_{A_0\left(\widetilde\sY \right)_{\widetilde \sU}}$ preserves the $C^*$-norm. Let us prove that $\pi_{A_0\left(\widetilde\sY \right)_{\widetilde \sU}}$ is surjective. If $\widetilde a \in  A_0\left(\widetilde\sY \right)$ then for any $\eps > 0$ 
there are sets $\left\{\widetilde f'_1,..., \widetilde f'_m\right\}\in C_c\left( \widetilde{\sY}\right)$ and $\left\{a_1, ..., a_m\right\}\in A$ such that
$$
\left\|  \widetilde a-\sum_{j=1}^m \varphi_A\left( a_j\right)\varphi_C\left( \widetilde f'_j\right)  \right\|< \frac{\eps}{2}.
$$
If $\left\{\widetilde u_\a\right\}_{\a \in \mathscr A}$ id sn approximate unit for $C_0\left( \widetilde \sU\right)$  then there is $\a_\eps \in \mathscr A$ such that
$$
\left\|  \widetilde a-\widetilde a \widetilde u_\a   \right\|< \frac{\eps}{2}
$$.
One has 
\bean
\left\|  \widetilde a-\sum_{j=1}^m \varphi_A\left( a_j\right)\varphi_C\left( \widetilde f'_j\right)\widetilde u_\a  \right\|\le \left\| \left(  \widetilde a-\sum_{j=1}^m \varphi_A\left( a_j\right)\varphi_C\left( \widetilde f'_j\right)\right) \widetilde u_\a  \right\|+\\
\left\|  \widetilde a-\widetilde a \widetilde u_\a   \right\| < \eps
\eean 
On the other hand if 
$$
a = \sum_{j=1}^{m} a_j\desc_q\left( \widetilde f'_j \widetilde u_\a\right) 
$$
where $\desc_q$ is $q$-descent (cf. Definition \ref{top_lift_desc_defn}) then
\bean
a \in  A_{\sU},\\
\sum_{j=1}^m \varphi_A\left( a_j\right)\varphi_C\left( \widetilde f'_j\right)\widetilde u_\a= \pi_{A_0\left(\widetilde\sY \right)_{\widetilde \sU}}.
\eean 
So the image $\pi_{A_0\left(\widetilde\sY \right)_{\widetilde \sU}}\left( A_{\sU}\right)$ is dense in $A_0\left(\widetilde\sY \right)_{\widetilde \sU}$, taking into account that  $A_{\sU}$ is $C^*$-norm closed we conclude that the homomorphism $\pi_{A_0\left(\widetilde\sY \right)_{\widetilde \sU}}$ is bijective.

	\end{proof}

\begin{lemma}\label{blowing_descent_lem}

	Let 	$C_0\left( \sY\right)\hookto M\left( A\right) $ be  Hausdorff blowing-up (cf. Definition \ref{blowing_defn}), and let  $q: \widetilde{\sY}\to \sY$ be  an   $A$-{regular} covering (cf. Definition \ref{blowing_a_regular_defn}). If $G$ is a properly discontinuous group of homeomorphisms \ref{top_properly_disc_group_defn} of $\widetilde\sY$ with $\sY \cong \widetilde \sY/G$, and 
		\be\label{blowing_com_eqn}
	A_c\left(\widetilde\sY \right) \bydef \left\{\left.\widetilde a \in A_0\left(\widetilde\sY \right)\right|\supp \widetilde a \mathrm{~is ~compact}  \right\}
	\ee
	is the {algebra of compactly supported elements} (cf. Definition \ref{blowing_compact_supp_defn}) then for any $\widetilde a \in  A_c\left(\widetilde\sY \right)$ the series 
	$$
a  \bydef 	\bt\text{-} \sum_{	g \in {G}}\widetilde a
	$$
	is convergent with respect to the strict topology of $M\left( A_0\left(\widetilde\sY \right)\right)$ (cf. Definition \ref{strict_topology_defn}) and $a \in A_0\left(\sY \right)$ 
\end{lemma}	
\begin{proof}
	If a pair 
	$
	\left(\	\left\{\left(\sU_\a, \sV_\a, s_\a\right)\right\}_{\a \in \mathscr A}, \left\{\left(\widetilde \sU_{\widetilde\a}, \widetilde \sV_{\widetilde\a}, \widetilde s_{\widetilde\a}\right)\right\}_{\widetilde \a \in \widetilde{\mathscr A}} \right)$ 
	is   a $q$-{covering} (cf. Definition \ref{top_cov_defn})
  then there is a covering sum $\widetilde f =\sum_{j=1}^n \widetilde{f}_j$ for  $\supp \widetilde{a}$ dominated by the family $\left\{\widetilde \sU_{\widetilde \a}\right\}$ (cf. Definition \ref{top_covering_sum_defn}). It follows that
   
 \be\label{blowing_convt_eqn}
 \widetilde a =\widetilde a\pi_C\left(  \widetilde f\right)  = \sum_{j=1}^n \widetilde a\pi_C\left(  \widetilde f_j\right) 
 \ee 
 Suppose that $\supp \widetilde f_j \subset \widetilde \sU$ where $ \widetilde \sU \in \left\{\widetilde \sU_{\widetilde\a}\right\}$. Let   $\sU_{\a_j} \bydef q\left(\widetilde \sU \right) \in \left\{ \sU_{\a}\right\}$, and let
 both $A_{\sU_{\a_j}}$ and $A_0\left(\widetilde\sY \right)_{\widetilde \sU}$ are left ${\sU_{\a_j}}$- and $\widetilde \sU$-ideals respectively (cf. Definition \ref{blowing_ideals_au_ua_defn}).
 There a  given by the Lemma \ref{blowing_left_ideal_lem} isomorphism
\bean
\pi_{A_0\left(\widetilde\sY \right)_{\widetilde \sU}}:  A_{\sU_{\a_j}}\to A_0\left(\widetilde\sY \right)_{\widetilde \sU}\\
	a \mapsto A_0\left(q \right)\left(   a \right) \widetilde s
\eean
of left $A$-modules. If $\widetilde a \in A_0\left(\widetilde\sY \right)$ then similarly to \eqref{blowing_convt_eqn} one has
 \be\label{blowing_convtt_eqn}
\widetilde a =\pi_C\left(  \widetilde f\right) \widetilde a\pi_C\left(  \widetilde f\right)  =  \sum_{\substack{j=1\\k=1}}^n \widetilde  a_{kj}\quad \text{where}\quad  \widetilde  a_{kj} \bydef  \pi_C\left(  \widetilde f_k\right)\widetilde  a\pi_C\left(  \widetilde f_j\right) 
\ee 
If $a_{kj}\bydef \pi^{-1}_{A_0\left(\widetilde\sY \right)_{\widetilde \sU}}\left( \widetilde a_{kj}\right) \in A_{\sU_{\a_j}}$ then
$$
 \widetilde a_{kj} = \widetilde \pi_C\left(\widetilde s_{\a_k}\right) A_b\left( q\right) \left(a_{jk} \right) \pi_C\left(\widetilde s_{\a_j}\right)
$$
If $\widetilde b \in A_0\left(\widetilde\sY \right)$ then from the Lemma \ref{blowing_blowing_lem} and the Definition \ref{blowing_a_h_eqn} it follows that for any $\eps > 0$ there is $\widetilde \varphi_\eps \in C_c\left(\widetilde \sY \right)$ such that
$$
\left\|\widetilde b -  \widetilde \varphi_\eps \widetilde b \right\| < \frac{\eps}{2\left\|a_{kj} \right\|}
$$
Since both sets $\supp \widetilde s$ and $\supp \widetilde \varphi_\eps$ are compact the set
$$
G_0\bydef \left\{g \in G \left|\supp \widetilde s\cap \supp \widetilde \varphi_\eps \neq \emptyset \right.\right\}
$$
is finite. If follows that for any finite subset $G'\subset G$ one has
$$
G_0\subset G' \quad \Rightarrow \quad \sum_{g \in G_0} \left( g \widetilde a_{kj}\right) \widetilde \varphi_\eps\widetilde b= \sum_{g \in G'} \left( g \widetilde a_{kj}\right) \widetilde \varphi_\eps\widetilde b,
$$
it turns out that
\bean
G_0\subset G' \quad \Rightarrow \quad \left\|\sum_{g \in G'} \left( g  \widetilde a_{kj}\right) \widetilde b - A_b\left( q\right) \left(a_{jk} \right) \widetilde b\right\|\le \left\|A_b\left( q\right) \left(a_{jk} \right)  \widetilde \varphi_\eps \widetilde b  - A_b\left( q\right) \left(a_{jk} \right) \widetilde b\right\|+\\+ \left\|\sum_{g \in G'} \left( g \widetilde a_{jk}\right)\left( \widetilde \varphi_\eps  \widetilde b - \widetilde b\right) \right\|+\left\|\sum_{g \in G'} \left( g\widetilde a_{jk} \right) \widetilde \varphi_\eps \widetilde b - A_b\left( q\right) \left(a_{jk} \right) \widetilde \varphi_\eps \widetilde b\right\|< \\< \frac{\eps}{2\left\|a_{jk}  \right\|}\left\|a_{jk}\right\|+ \frac{\eps}{2\left\||a_{jk} \right\|}\left\||a_{jk} \right\|+ 0 = \eps.
\eean
So the series $\sum_{g \in G} \left( g A_b\left( q\right) \left(a_{jk} \right)\right) \widetilde b$ is $C^*$-norm convergent. Similarly one can prove that the series  $\sum_{g \in G} \widetilde b\left( g A_b\left( q\right) \left(a_{jk} \right)\right) $ is also $C^*$-norm convergent, so one has
$$
a^*a = \b\text{-}\sum_{g \in G} \widetilde b\left( g \widetilde a^*\widetilde a\right)
$$
where the sum of the series implies the strict topology of $A_0\left(\widetilde\sY \right)$. 
\end{proof}

\begin{definition}\label{blowing_descent_compactly_supported_defn} 
  A given by the Lemma \ref{blowing_descent_lem} $A$-$A$-bimodule homomorphism 
	\be\label{blowing_compactly_supported_descent_eqn}
	\begin{split}
\widetilde a\mapsto \bt\text{-} \sum_{	g \in {G}}\widetilde a
	\end{split} 
	\ee
 is the \textit{compactly supported $q$-descent}.
\end{definition}
\begin{definition}\label{blowing_descent_defn}
	If $\widetilde a \in A\left(\widetilde\sY \right)$ is such that $\supp \widetilde a$ is mapped homeomorphically onto $q\left( \supp \widetilde a\right)$ then we say that $\desc^c_q\left(\widetilde a\right)$ is  $q$-\textit{descent} of $\widetilde a$. We denote it by $\desc_q\left(\widetilde a\right)$.
\end{definition}
\begin{remark}
	Both  Definitions \ref{blowing_descent_compactly_supported_defn} and  \ref{blowing_descent_defn} are generalizations  of \ref{top_compactly_supported_descent_defn} and \ref{top_lift_desc_defn} ones.
\end{remark}
\begin{lemma}\label{blowing_descent_pedersen_lem} 
If	$C_0\left( \sY\right)\hookto M\left( A\right) $ is   Hausdorff blowing-up (cf. Definition \ref{blowing_defn}), and   $q: \widetilde{\sY}\to \sY$ is   an   $A$-{regular} covering (cf. Definition \ref{blowing_a_regular_defn}) then one has
	\be\label{blowing_compactly_supported_descen_oedersen_eqn}
	\begin{split}
	 K\left(A \right) \subset 	\desc^c_q \left(  K\left( A_0\left(\widetilde\sY \right)\right)\right) 
	\end{split} 
	\ee	
\end{lemma}
\begin{proof}
If  $a \in K\left(A \right)$ then  $\supp a $ is compact. 	If a pair \\
$
\left(\	\left\{\left(\sU_\a, \sV_\a, s_\a\right)\right\}_{\a \in \mathscr A}, \left\{\left(\widetilde \sU_{\widetilde\a}, \widetilde \sV_{\widetilde\a}, \widetilde s_{\widetilde\a}\right)\right\}_{\widetilde \a \in \widetilde{\mathscr A}} \right)$ 
is   a $q$-{covering} (cf. Definition \ref{top_cov_defn})
then there is a covering sum $ f =\sum_{j=1}^n {f}_j$ for  $\supp \widetilde{a}$ dominated by the family $\left\{ \sU_{ \a}\right\}$ (cf. Definition \ref{top_covering_sum_defn}). It follows that
\be\label{blowing_convlt_eqn}
a =\widetilde af  = \sum_{j=1}^n a f_j.
\ee 
For any $j = 1,...,n$ one has $b_j \bydef  a f_j\in K\left(A \right)$. Let $\supp f_j \subset \sU_{\a_j}$ where $\a_j\in \mathscr A$, and  let $\widetilde  \a_j \in \widetilde{\mathscr A}$ is such that $q\left( \sU_{\widetilde\a_j}\right) = \sU_{\a_j}$. If $\widetilde b_j =\varphi_A\left(b_j \right)  \varphi_C\left( \widetilde s_{\widetilde\a_j}\right)$ then from the proof of the Lemma \ref{blowing_descent_lem} it follows that $b_j = \desc^c_p\left(\widetilde b_j\right)$. If  $b_j\in K\left( A\right)_+$ then  $b_j \le  \sum_{l= 1} f_{\eps_l}\left(b^l_j \right)$ where $b^l_j \in  A_+$ (cf. Lemma \ref{pedersen_eps_lem}). We can suppose that $b^l_j=b^l_js_{\a_j}$. If $\widetilde b^l_j\bydef \varphi_A\left(b_j \right)  \varphi_C\left( \widetilde s_{\widetilde\a_j}\right)$ then $f_{\eps_j}\left( \widetilde b^l_j\right) = \varphi_A\left(f_{\eps_j}\left( b^l_j \right) \right)  \varphi_C\left( \widetilde s_{\widetilde\a_j}\right)$ and $\widetilde b_j \le\sum_{l= 1}^m f_{\eps_l}\left(\widetilde b^l_j \right)$, i.e. $\widetilde b_j \in  K\left( A_0\left(\widetilde\sY \right)\right)$. The  completion of the proof is left to the reader.
\end{proof}

\begin{exercise}\label{blowing_hausdorff_exer}
	Let $A$ be a $C^*$-algebra with Hausdorff, connected spectrum $\sX$, and let $C_0\left( \sX\right)\subset M\left(A \right)$ be Hausdorff blowing-up (cf. Example \ref{blowing_hausdorff_exm}). Let $p: \widetilde{\sX}\to \sX$ be a transitive covering with residually finite group $G\left(\left.\widetilde \sX\right|\sX \right)$ of covering transformations. Prove following statements.
	\begin{enumerate}
		\item There is $p$-lift $A_b\left(p\right)\bydef \pi : A \hookto	M\left( A_0\left( \widetilde \sX\right)\right)$ (cf. Definition \ref{blowing_lift_defn}  ).
		\item The $C^*$-algebra    $A_0\left( \widetilde \sX\right)$ given by the Definition \ref{top_cs_functa_b_defn} is a specialization of   given by the Definition \ref{blowing_lift_defn}   one.
		\item The Definition \ref{top_compactly_supported_descent_defn} is a specialization of the  \ref{blowing_descent_compactly_supported_defn}.
	\end{enumerate}
\end{exercise}

\begin{empt}
	Both Lemmas \ref{blowing_pedersen_compact_lem} and  \ref{blowing_descent_lem}  \ref{foli_fibration_exm} yield a sesquilinear pairing 
	\be\label{blowing_pa_eqn}
	\begin{split}
		\left\langle\cdot, \cdot \right\rangle : K\left(A_0\left( \widetilde \sY\right) \right)\times K\left(A_0\left( \widetilde \sY\right) \right) \to A,\\
		\left\langle \widetilde a, \widetilde b \right\rangle \bydef \bt\text{-}	\sum_{	g \in G\left( \left. \widetilde \sY\right| \sY\right) }g \left( \widetilde a^* \widetilde b\right) 
	\end{split}
	\ee
	where the convergence with respect to the strict topology of $M\left(A_0\left( \widetilde \sY\right) \right)$ (cf. Definition \ref{strict_topology_defn}) is implied. Thus $K\left(A_0\left( \widetilde \sY\right)\right)$ becomes a pre-Hilbert $A$-module (cf. Definition \ref{hilbert_module_defn}).
\end{empt}
\begin{definition}\label{blowing_hilb_defn}
	The $C^*$-Hilbert $A$-module which is a completion of $K\left(A_0\left( \widetilde \sY\right)\right)$ with respect to norm
	$$
	\left\|\widetilde a\right\|\bydef \sqrt{\left\| \left\langle \widetilde a, \widetilde a \right\rangle \right\|} 
	$$
	(cf. \ref{hilbert_module_norm_eqn})  we denote by $\mathscr{L}^2\left(A_0\left( \widetilde \sY\right) \right)$. 
\end{definition}
\begin{lemma}\label{blowing_hilb_iso_lem}
If 	$C_0\left( \sY\right)\hookto M\left( A\right) $ is   Hausdorff blowing-up (cf. Definition \ref{blowing_defn}), and   $q: \widetilde{\sY}\to \sY$ is  an   $A$-{regular} covering (cf. Definition \ref{blowing_a_regular_defn}) then there is a natural isomorphism 
	\be\label{blowing_hilb_iso_eqn}
	\mathscr{L}^2\left(A_0\left( \widetilde \sY\right) \right)\cong \mathscr L^2\left(\widetilde \sY \right)_A
	\ee
	where $\mathscr L^2\left(\widetilde \sY \right)_A$ is $\widetilde \sY$-$A$-{module} (cf. Definition \ref{blowing_lift_hm_defn})
\end{lemma}
\begin{proof}
	If $\widetilde a \in K\left( A_0\left( \widetilde \sY\right)\right)$ the from \eqref{blowing_pedersen_compact_lem} it follows that 
	$$
	\exists a \in A \quad \exists \widetilde f\in C_c\left(\widetilde\sY \right)\quad  \widetilde a= \varphi_C\left( \widetilde f\right)\varphi_A\left(a \right),
	$$
	so there is an inclusion
	$$
	K\left( A_0\left( \widetilde \sY\right)\right)\subset C_c\left( \widetilde \sY \right)\otimes_{C_0\left(\sY\right) } A
	$$
	and taking into account that $\mathscr L^2\left(\widetilde \sY \right)_A$ is a completion of $K\left( A_0\left( \widetilde \sY\right)\right)\subset C_c\left( \widetilde \sY \right)\otimes_{C_0\left(\sY\right) } A$ one has
	$$
	\mathscr{L}^2\left(A_0\left( \widetilde \sY\right) \right)\subset\mathscr L^2\left(\widetilde \sY \right)_A
	$$
	Let $\eps > 0$. For any  $\widetilde a \in  A_0\left( \widetilde \sY\right)$ and $\dl > 0$ there is $\widetilde a' \in  A_0\left( \widetilde \sY\right)$ such that
	$$
	x \bydef  a -  a'; \quad \left\| x \right\|< \dl. 
	$$
	One has
	$$
	\left\|\widetilde f \otimes a - \widetilde f \otimes a'\right\|_H=  \left\|\widetilde f x\right\|_H \le \left\|x\right\|_C \left\|\desc^c_q\left(\widetilde f^*\widetilde f \right) \right\|_C \le \dl \left\|\desc^c_q\left(\widetilde f^*\widetilde f \right) \right\|_C
	$$
	where $ \left\|\cdot\right\|_H$ is a norm of $C^*$-Hilbert $A$-module $\mathscr L^2\left(\widetilde \sY \right)_A$,  $ ~\left\|\cdot\right\|_H$ is the $C^*$-norm and $\desc^c_q$ is compactly supported $q$-descent (cf. Definition \ref{top_compactly_supported_descent_defn}).  If $\dl \bydef \frac{\eps }{\left\|\desc^c_q\left(\widetilde f^*\widetilde f \right) \right\|_C}$ then $\left\|\widetilde f \otimes a - \widetilde f \otimes a'\right\|_H < \eps$ so $K\left( A_0\left( \widetilde \sY\right)\right)$ is dense in $\mathscr L^2\left(\widetilde \sY \right)_A$ and there is an isomorphism \eqref{blowing_hilb_iso_eqn}.
\end{proof}

\begin{lemma}\label{blowing_fin_lem}
	Let 	$C_0\left( \sY\right)\hookto M\left( A\right) $ be   Hausdorff blowing-up (cf. Definition \ref{blowing_defn}). If  $q: \widetilde{\sY}\to \sY$ is   an   $A$-{regular}  (cf. Definition \ref{blowing_a_regular_defn}) finite-fold  transitive covering and 
	\be\label{blowing_fin_eqn}
	A_b\left(q\right) : A \hookto	M\left( A_0\left( \widetilde \sY\right)\right)  
	\ee
	is a $q$-{lift} of $A$ (cf. Definition \ref{blowing_lift_defn}  ) then one has $A_b\left(q\right) \left(  A\right) \subset  A_0\left( \widetilde \sY\right)$. 
\end{lemma}
\begin{proof}
According to   the Lemma \ref{blowing_lift_constr_lem}   it follows that $A_0\left( \widetilde \sY\right)$ is the $C^*$-norm completion of 
$$
\varphi\left(C_c\left( \widetilde \sY \right)\otimes_{C_0\left(\sY\right) } K\left( A\right) \right)
$$
For any $\widetilde a \in 	\varphi\left(C_c\left( \widetilde \sY \right)\otimes_{C_0\left(\sY\right) } K\left( A\right) \right)$ there is $\widetilde f\in C_c\left( \widetilde \sY\right)_+$ such that 
$$
\varphi_C\left(\widetilde f \right) \widetilde a.
$$
on the other hand for any $a \in  K\left( A\right)$ there is $f \in C_c\left(\sY \right)_+$ such that $f a = a$.  Since a covering $q$ is finite-fold one has 	$A_b\left(q\right)\left( a\right) =  \varphi\left(C_c\left( q\right)\left(f \right) \otimes a \right)\in \varphi\left(C_c\left( \widetilde \sY \right)\otimes_{C_0\left(\sY\right) } K\left( A\right) \right)\subset  A_0\left( \widetilde \sY\right)$ where $C_c\left( q\right)$ is given by \eqref{top_compact_cc_eqn}.
Since $K\left(A \right)$ is dense in $A$  one has  
$$
\forall a \in A \quad A_b\left(q\right) \left( a\right) \subset  A_0\left( \widetilde \sY\right),
$$
or equivalently the equation \eqref{blowing_fin_eqn} holds.
\end{proof}

\begin{definition}\label{blowing_finite_lift_defn}
	The given by the Lemma \ref{blowing_fin_lem} injective $*$-homomorphism
	\be	\label{blowing_lift_fin_eqn}
	A_0\left(q\right) : A \hookto	 A_0\left( \widetilde \sY\right)  
	\ee
	is a \textit{finite}-$q$-\textit{lift} or simply a $q$-\textit{lift} of $A$. 	
\end{definition}
\begin{lemma}\label{blowing_proper_group_lem}
	Let 	$C_0\left( \sY\right)\hookto M\left( A\right) $ be  Hausdorff blowing-up (cf. Definition \ref{blowing_defn}), and let  $q: \widetilde{\sY}\to \sY$ be  an   $A$-{regular} covering (cf. Definition \ref{blowing_a_regular_defn}).
	If the space $\sY$ is locally connected (cf. Definition \ref{top_locally_connected_defn}) and $\widetilde \sY$ is connected then	there is a natural isomorphism of groups
	\be\label{blowing_group_action_eqn}
	G \cong \left\{\left. g \in \Aut\left( M\left( A_0\left( \widetilde \sY\right)\right)  \right) \right| \forall a \in A_b\left(q \right) \left(A \right)  \quad g a = a \right\}.
	\ee
\end{lemma}
\begin{proof}
	Clearly there is the natural inclusion 
	$$
	G \subset \left\{\left. g \in \Aut\left( M\left( A_0\left( \widetilde \sY\right)\right)  \right) \right| \forall a \in A_b\left(q \right) \left(A \right)  \quad g a = a \right\}.
	$$
	If $\widetilde \sU\subset \widetilde \sY$ is connected set  homeomorphically mapped onto $\sU\bydef q\left(\widetilde \sU \right)$ then any hereditary $C^*$-subalgebra of $A_0\left( \widetilde \sY\right)$ which is isomorphic to $_\sU A_\sU$ as 
	$_\sU A_\sU-~_\sU A_\sU$-bimodule equals to $_{g'\widetilde \sU}A_0\left( \widetilde \sY\right)_{g'\widetilde \sU}$  where $g' \in G$. For any 
	$$
	g \in  \left\{\left. g \in \Aut\left( M\left( A_0\left( \widetilde \sY\right)\right)  \right) \right| \forall a \in A_b\left(q \right) \left(A \right)  \quad g a = a \right\}
	$$
	denote by $g_{\widetilde \sU} \in  G$ such that $g\left(  _{\widetilde \sU}A_0\left( \widetilde \sY\right)_{\widetilde \sU}\right) = ~_{g_{\widetilde \sU}\widetilde \sU}A_0\left( \widetilde \sY\right)_{g_{\widetilde \sU}\widetilde \sU}$. If both $\widetilde \sU', \widetilde \sU'' \subset \widetilde \sY$ are  is connected set  homeomorphically mapped onto $ q\left(\widetilde \sU' \right)$ and $q\left(\widetilde \sU'' \right)$ sets  such that $\widetilde \sU'\cap \widetilde \sU'' \neq \emptyset$ then there are $\widetilde a' \in ~_{\widetilde \sU'}A_0\left( \widetilde \sY\right)_{\widetilde \sU'}$ and \\ $\widetilde a'' \in ~_{\widetilde \sU''}A_0\left( \widetilde \sY\right)_{\widetilde \sU''}$ with $\widetilde a'\widetilde a''\neq 0$.
	Since $g$ is an automorphism 
	$$
	g \left(\widetilde a'\widetilde a'' \right)= \left( g_{\widetilde \sU'}\widetilde a'\right) \left( g_{\widetilde \sU''}\widetilde a''\right)\neq 0
	$$
	If $g_{\widetilde \sU'} \neq g_{\widetilde \sU''}$ then $g_{\widetilde \sU'}\widetilde \sU'\cap g_{\widetilde \sU''}\widetilde \sU''= \emptyset$ and $g \left(\widetilde a'\widetilde a'' \right)= \left( g_{\widetilde \sU'}\widetilde a'\right) \left( g_{\widetilde \sU''}\widetilde a''\right)= 0$. From this contradiction we conclude that $g_{\widetilde \sU'} = g_{\widetilde \sU''}$. The space $\widetilde\sY$ is connected and using the Theorem \ref{zorn_thm} we conclude that there is $g' \in G$ such that  $g_{\widetilde \sU}= g'$ for any connected  $\widetilde \sU$ homeomorphically mapped onto $ q\left(\widetilde \sU \right)$. Using this fact one can prove that if $A_0\left( \widetilde \sY\right)_{\widetilde \sU}$ is left $\widetilde \sU$-ideal  (cf. Definition \ref{blowing_ideals_au_ua_defn}) then  $g A_0\left( \widetilde \sY\right)_{\widetilde \sU} = g'A_0\left( \widetilde \sY\right)_{\widetilde \sU}$ is a left $g_\iota\widetilde \sU$-ideal. If $\widetilde a\in K\left( A_0\left( \widetilde \sY\right)\right)$ is an element of the Pedersen's ideal (cf. Definition \ref{pedersen_ideal_defn}) then the support $\supp\widetilde a$ of $\widetilde a$ is compact. There is a finite family $\left\{\widetilde \sU_1, ..., \widetilde \sU_n\right\}$ of open connected subsets of $\widetilde \sY$ such that: 
	\begin{itemize}
		\item $\supp\widetilde a\subset \bigcup_{j=1}^n \widetilde\sU_j$,
		\item for each $j=1,...,n$ the set $\widetilde \sU_j \sY$ is  homeomorphically mapped onto $ q\left(\widetilde \sU_j \right)$.
	\end{itemize}
	If $\widetilde f \bydef \sum_{j = 1}^n \widetilde f_j$ is a {dominated} by the family {covering sum} for $\supp\widetilde a$ (cf Definition \ref{top_covering_sum_defn}) then 
	$$
	\widetilde a = \widetilde a \widetilde f =\sum_{j = 1}^n  \widetilde a \widetilde f_j.
	$$
	From $ \widetilde a  \widetilde f_j\in A_0\left( \widetilde \sY\right)_{\widetilde \sU}$ it follows that $g\left(  \widetilde a'\widetilde f_j\right) = g'\left(  \widetilde a'\widetilde f_j\right)$ and  
	$$
	g\widetilde a= g \left( \widetilde a \widetilde f\right)  = g' \widetilde a
	$$
	Since $K\left( A_0\left( \widetilde \sY\right)\right)$ is dense in $M\left( A_0\left( \widetilde \sY\right)\right)$ with respect to the strict topology of $M\left( A_0\left( \widetilde \sY\right)\right)$ (cf. Definition \ref{strict_topology_defn}) one has
	$$
	\forall \widetilde a \in M\left( A_0\left( \widetilde \sY\right)\right) \quad g\widetilde a= g' \widetilde a
	$$
\end{proof}

\begin{lemma}\label{blowing_lift_composition_lem}
Let 	$C_0\left( \sY\right)\hookto M\left( A\right) $ be  Hausdorff blowing-up (cf. Definition \ref{blowing_defn}), and let  $q: \widetilde{\sY}\to \sY$ be  an   $A$-{regular} covering (cf. Definition \ref{blowing_a_regular_defn}). Suppose that  $\widetilde G$ is a {properly discontinuous} \ref{top_properly_disc_group_defn} group of homeomorphisms of $\widetilde \sY$ such that $\sY \cong \widetilde \sY / \widetilde G$.
	Let 	 	$A_b\left(q\right) : A \hookto	M\left( A_0\left( \widetilde \sY\right)\right)  
	$
	be a $q$-{lift} (cf. Definition \ref{blowing_lift_defn}  ), and let $\phi:\widetilde G\to G$ be a homomorphism onto a finite group $G$. If both  $q': \sY' \to \sY$ and $\widetilde q': \widetilde\sY \to \sY'$ are a natural transitive  coverings with
	$\sY' \bydef \widetilde \sY/ \ker\left( \widetilde G\to G \right)$ then one has:
	\begin{enumerate}
		\item [(i)] the covering $q'$ is  $A$-{regular} covering $q'$ is $A$-regular, so  there is a {finite}-$q'$-{lift} 	$A_0\left(q'\right) : A \hookto	 A_0\left(  \sY'\right)$ (cf. Definition \ref{blowing_finite_lift_defn}),
		\item[(ii)]  the covering $\widetilde q'$ is  $A_0\left(  \sY'\right)$-{regular}, so there is a $\widetilde q'$-{lift} 	$A_b\left(\widetilde q'\right) : A_0\left(  \sY'\right)\hookto	M\left(  A_0\left( \widetilde \sY\right)\right) $ (cf. Definition \ref{blowing_lift_hom_defn}).
	\end{enumerate} 
\end{lemma}

\begin{proof}
	(i)	Let both $\mathscr L^2\left(\widetilde \sY \right)_A$ and $\mathscr L^2\left(\sY' \right)_A$ be explained in \ref{blowing_lift_empt} $C^*$-Hilbert $A$-modules. 	A natural action $C_0\left( \sY'\right) \times C_0\left(\widetilde \sY\right)$ yields an action $C_0\left( \sY'\right) \times \mathscr L^2\left(\widetilde \sY \right)_A\to \mathscr L^2\left(\widetilde \sY \right)_A$. Using these actions one can obtain  inclusions 
	$$
	C_0\left( \sY'\right)\subset	M\left(  A_0\left( \widetilde \sY\right)\right)= 	M\left(  A'_0\left( \widetilde \sY\right)\right) 
	$$
	If $A'$ is an algebraically generated by by elements
	\be\label{blowing_fp_eqn}
	f'a;\quad a \in K\left(A \right) \quad f' \in C_c\left(q' \right) 
	\ee
	right $A$-submodule of $M\left(  A'_0\left( \widetilde \sY\right)\right)$ then $\mathscr L^2\left(\widetilde \sY \right)_A$ is the completion of $A'$ with respect to following norm
	$$
	\left\| f'a \right\|_H =\sqrt{\left\| \sum_{	g \in G}g\left( ({f'a})^*f'a\right)\right\|_C  }  
	$$
	where $\left\| f'a \right\|_C$ is the $C^*$-norm of $	M\left(  A'_0\left( \widetilde \sY\right)\right)$. From 
	$$
	\left\| f'a \right\|_H \le \left\|G \right| \left\| f'a \right\|_C,
	\left\| f'a \right\|_C \le \left\| f'a \right\|_H
	$$
	it follows that $\mathscr L^2\left(\widetilde \sY \right)_A$ coincides with the $C^*$-norm completion of $A'$.
	If $a \in\left( A\right)$ and $h$ is a covering  sum for $\supp a$ (cf. Definition \ref{top_covering_sum_defn}) then
	$$
	f' a = f' (a h) = f' a f''; \quad f'' =  C_c\left(q' \right)\left(h \right)\in C_c\left( \sY'\right)  
	$$
	where $C_c\left(q' \right):  C_c\left( \sY\right)\hookto C_c\left( \sY'\right)$ is given by \eqref{top_compact_cc_eqn}.

	Consider a family $\left\{{\sU}_\a\subset  \sY\right\}_{\a\in \mathscr A}$ of connected open sets with compact closures such that $\sU_\a$ is evenly covered $q$  for all $\a\in \mathscr A$ and $\sY = \bigcup_{\a \in \mathscr A}\sU_\a$. For any $\a\in\mathscr A$ we select  a connected open set $\widetilde \sU_\a\subset \widetilde \sY$  which is mapped homeomorphically onto $\sU_\a$. If $\widetilde{\mathscr A}\bydef G\times  \mathscr A$ then there are the natural action $G\times \widetilde{\mathscr A}\to \widetilde{\mathscr A}$ and the projection $p: \widetilde{\mathscr A} \to {\mathscr A}$.  If 	
	$$
	\sum_{j=1}^n f'_j
	$$
		
	is  a covering sum  for $\supp f' \cup \supp f''$ dominated by the family $\left\{{\sU}'_{\a'}\subset  \sY'\right\}_{\a'\in \mathscr A'}$ 
	(cf. Definition  \ref{top_covering_sum_defn}) then
	\bean
	f' a\in A' \quad \Rightarrow \quad f' a f'' = \sum_{\substack{j=1\\ k=1}}^n f'_j a' f''_k = \sum_{\substack{j=1\\ k=1}}^n h'_j a_{jk} h''_k, \\ \text{where} \quad  a_{jk} \bydef \sum_{	g \in G}g\left(\sqrt{f'_j}f' a f''\sqrt{f'_k}\right),\quad  h'_j = f'\sqrt{f'_j},\quad  h''_k = f''\sqrt{f''_k}   
	\eean
	So any element $a' \in A'$ can represented by the following way
	\be\label{blowing_jkp_eqn}
	a' = \sum_{\substack{j=1\\ k=1}}^n h'_j a_{jk} h''_k
	\ee
	where for all $j, k$ there are $\sU'_1, \sU'_2 \in \left\{{\sU}'_{\a'}\subset  \sY'\right\}_{\a'\in \mathscr A'}$ such that $\supp h'_j \subset \sU'_1$ and  $\supp h''_k \subset \sU'_2$. We leave to the reader proof of that the set of given by \eqref{blowing_jkp_eqn} elements is closed with respect addition, subtraction, multiplication and $*$-operation, i.e. $A'$ is a $*$-algebra.  So the completion $A'_0\left(\sY' \right)$ is a $C^*$-algebra, i.e. a covering $q': \sY' \to \sY$ is $A$-regular (cf. Definition \ref{blowing_a_regular_defn}). From the Lemma    \ref{blowing_fin_lem} it follows that  there is a {finite}-$q'$-{lift} 	$A_0\left(q'\right) : A \hookto	 A_0\left(  \sY'\right)$\\
	(ii) If $X \bydef  C_c \left(\widetilde \sY \right) \otimes_{C_0\left(\sY \right) K\left( A\right)  }$ then  there is a $A$-valued product
	\bean
	\left\langle\cdot , \cdot  \right\rangle_A : X \times X \to A,\\
	\left\langle a' \otimes \widetilde f' ,  a'' \otimes \widetilde f''  \right\rangle_A \bydef a'^* \desc_q (\widetilde f'^* \widetilde f'') a''
	\eean
	From the construction \ref{blowing_lift_empt}  it follows that $\mathscr L^2\left(\widetilde \sY \right)_A$ id a completion of $X$ with respect to the norm
	$$
	\left\| \widetilde{f} \otimes a\right\|_A \bydef  \sqrt{\left\| a^*\left( \bt\text{-}\sum_{g \in G}f^*f\right)  a \right\|}
	$$
	where $\left\|\left\langle\widetilde{a} , \widetilde{a}  \right\rangle_A \right|$ is the $C^*$-norm of $C^*$ algebra $A$.
	There is an  $A_0\left(  \sY'\right)$-valued product
	\bean
	\left\langle\cdot , \cdot  \right\rangle_{A_0\left(  \sY'\right)} : X \times X \to A_0\left(  \sY'\right),\\
	\left\langle a' \otimes \widetilde f' ,  a'' \otimes \widetilde f''  \right\rangle_{A_0\left(  \sY'\right)} \bydef a^*\left(\bt\text{-}\sum_{g \in \ker\left( G\to G\right) }f'^*f'' \right) a
	\eean
	with the following norm
	$$
	\left\| \widetilde{f} \otimes a\right\|_{A_0\left(  \sY'\right)} \bydef  \sqrt{\left\| a^*\left( \bt\text{-}\sum_{g g \in \ker\left( G\to G\right)}f^*f\right)  a \right\|}
	$$
	From the above equations it turns out that
	\bean
	\left\| \widetilde{f} \otimes a\right\|_{A_0\left(  \sY'\right)} \le \left\| \widetilde{f} \otimes a\right\|_{A},\\
	\left\| \widetilde{f} \otimes a\right\|_{A}\le \left|G \right| \left\| \widetilde{f} \otimes a\right\|_{A_0\left(  \sY'\right)},
	\eean 
	so the completion of $\mathscr L^2\left(\widetilde \sY \right)_A$ with respect to the norm $\left\|\cdot\right\|_{A}$ coincides with the completion of $X$ with respect to the norm $\left\|\cdot\right\|_{A_0\left(  \sY'\right)}$. Taking into account the product $\left\langle\cdot , \cdot  \right\rangle_{A_0\left(  \sY'\right)}$ we conclude that $\mathscr L^2\left(\widetilde \sY \right)_A$ is a $C^*$-Hilbert $A_0\left(  \sY'\right)$-module. From $A \subset A_0\left(  \sY'\right)$ it follows that the algebraically generated by elements 
	\be\label{blowing_far_eqn}
	\widetilde f  a, \quad \widetilde f \in C_c \left( \widetilde \sY\right), \quad a \in  A
	\ee
	right $A$-submodule of $M\left(A_0\left(\widetilde \sY \right) \right)$ is a $\C$-subspace of the algebraically generated by elements 
	\be\label{blowing_farf_eqn}
	\widetilde f  a', \quad \widetilde f \in C_c \left( \widetilde \sY\right), \quad a' \in  A_0\left(\widetilde \sY \right).
	\ee
	right $A_0\left(\widetilde \sY \right)$-submodule of $M\left(A_0\left(\widetilde \sY \right) \right)$
	However for any $a' \in A_0\left( \sY' \right)$ and $\eps > 0$  there a sets $\left\{a_1, ..., a_n\right\}\subset A$ and $\left\{f'_1, ..., f'_n\right\}\subset C_c\left( \sY' \right)$ such that
	\be
	\left\|a' - \sum_{j=1}^nf'_j a_j   \right\|< \frac{\eps}{\left\|\widetilde f \right\|}.
	\ee
	So one has
	\be\label{blowing_farr_eqn}
	\left\|\widetilde f a' - \sum_{j=1}^n \widetilde ff'_j a_j   \right\|< \eps.
	\ee
	Since $\sum_{j=1}^n \widetilde ff'_j a_j$ belongs to generated by elements \eqref{blowing_far_eqn} module one has $\sum_{j=1}^n \widetilde ff'_j a_j\in A_0\left(\widetilde \sY \right)$ and taking into account \eqref{blowing_farr_eqn} one has $\widetilde f a'\in A_0\left(\widetilde \sY \right)$. So the closure of generated by \eqref{blowing_farf_eqn} module equals to  $A_0\left(\widetilde \sY \right)$. So the covering $\widetilde q': \widetilde{\sY}\to \sY'$ is $A_0\left( \sY' \right)$-{regular} (cf. Definition \ref{blowing_a_regular_defn}). From  the Lemma \ref{blowing_lift_constr_lem}  it follows that there is a  $\widetilde q'$-lift  $	A_b\left(\widetilde q'\right) :A_0\left( \sY' \right) \hookto	M\left( A_0\left( \widetilde \sY\right)\right)$.
\end{proof}

\subsection{Induced representations}

\paragraph{}

	Under the hypotheses \ref{blowing_lift_empt} let  $\mathscr L^2\left(\widetilde \sY \right)_A$ be the $\widetilde \sY$-$A$-module (cf. Definition \ref{blowing_lift_hm_defn}). If $\rho : A \to B\left(\H \right)$  there is a  given by \eqref{induced_representation_eqn} representation
	\be\label{blowing_induced_eqn}
		\widetilde{\rho}\bydef\mathscr L^2\left(\widetilde \sY \right)_A\text{-}\Ind^A_{{A}_0\left(\widetilde\sY \right) }\rho: A_0\left(\widetilde\sY \right)\to B\left(\widetilde{\H} \right)
		\ee
		is given by \eqref{induced_representation_eqn}, i.e. $\widetilde{\rho}$
		is the induced representation (cf. Definition \ref{induced_representation_defn}). 
\begin{lemma}\label{blowing_induced_representation_lem}
	If  $A_b\left(q \right) :A \hookto	M\left( A_0\left( \widetilde \sY\right)\right)$ is $q$-lift (cf. Definition \ref{blowing_lift_defn}  ) and $\rho : A\to B\left( \H\right)$ is a faithful, nondegenerate  representation (cf. Definitions  \ref{faithful_representation_defn} and \ref{nondegenerate_repr_defn}) then the given by \eqref{blowing_induced_eqn} induced representation $\widetilde{\rho}: A_0\left(\widetilde\sY \right)\to B\left(\widetilde{\H} \right)$ if faithful.
\end{lemma}
\begin{proof}
If $\widetilde{\rho}$ is  faithful if and only if for any positive element $\widetilde a \in K\left( A_0\left(\widetilde\sY \right)\right)_+$ of Pedersen's ideal (cf. Definition \ref{pedersen_ideal_defn}) of  $A_0\left(\widetilde\sY \right)$ one has $\widetilde{\rho}\left( \widetilde{a}\right) \neq 0$. From the Lemma \ref{blowing_pedersen_compact_lem} it follows that the support $\supp \widetilde a$ is compact. If $\sum_{j = 1}^n \widetilde f_j$ is a covering sum for $\supp \widetilde a$   {subordinated to} $q$ (cf. Definition \ref{top_covering_sum_subordinated_defn}) then one has
$$
\widetilde a = \sum_{\substack{j = 1\\k=1}}^n \widetilde f_j \widetilde a\widetilde f_k.
$$
The sum is not zero if and only is there is $j \in \left\{1,..., n\right\}$ such that $\widetilde f_j \widetilde a\widetilde f_j\neq 0$. There is and open subset $\widetilde\sU$ homeomorphically mapped onto $q\left(\widetilde\sU \right)$ such that $\supp  \widetilde f_j \subset \widetilde \sU$. If $\left\{\widetilde u_\a\right\}\subset C_0\left( \widetilde \sU\right)$ is an approximate unit for  $C_0\left( \widetilde \sU\right)$ (cf. Definition \ref{approximate_unit_defn}) then there is $ \widetilde u\in \left\{\widetilde u_\a\right\}$ such that $\widetilde u \widetilde f_j \widetilde a\widetilde f_j \widetilde u\neq 0$. If $a \bydef \desc_q\left( \widetilde f_j\widetilde a{\widetilde f_j}\right) \in A$ where $\desc_q$ means $q$-descent (cf. Definition \ref{blowing_descent_defn}) then $\widetilde u \widetilde f_j \widetilde a\widetilde f_j \widetilde u=\widetilde ua \widetilde u$. Moreover from the Lemma \ref{blowing_descent_lem} it follows that
\bean
\sum_{	g \in G\left(\left.\widetilde \sY\right|\sY \right)}g\left(  \widetilde ua \widetilde u\right) = 	\desc_q\left( \widetilde u\right)  a ~	\desc_q\left(\widetilde u\right)\quad \Rightarrow \\ \Rightarrow \desc_q\left( \widetilde u\right)  a ~	\desc_q\left(\widetilde u\right)>   \widetilde ua \widetilde u\quad \Rightarrow\quad \desc_q\left( \widetilde u\right)  a ~	\desc_q\left(\widetilde u\right)> 0\quad\Rightarrow\\ \Rightarrow \exists \xi \in \H \quad\desc_q\left( \widetilde u\right)  a ~	\desc_q\left(\widetilde u\right)\xi \neq 0.
\eean 
On the other hand there is an inclusion $\phi : C_c\left( \widetilde \sY\right)\otimes_{C_0\left( \sY\right) } \H \hookto   \widetilde \H$. If $\xi \in \H$ is such that  $\desc_q\left( \widetilde u\right)  a ~	\desc_q\left(\widetilde u\right)\xi > 0$ and $\widetilde \xi \bydef \phi \left({\widetilde u} \otimes \xi \right)$ then from  \eqref{hilb_prod_eqn} it follows that
\bean
\left(\left(\widetilde f_j \widetilde a\widetilde f_j \right) \widetilde \xi, \left(\widetilde f_j \widetilde a\widetilde f_j \right) \widetilde \xi \right)_{\widetilde \H}>\left( \desc_q\left( \widetilde u\right)  a ~	\desc_q\left(\widetilde u\right)\xi,\desc_q\left( \widetilde u\right)  a ~	\desc_q\left(\widetilde u\right)\xi\right)>  0,
\eean
From $\widetilde{a} \ge \left(\widetilde f_j \widetilde a\widetilde f_j \right)$ it follows that
$$
\left( \widetilde a \widetilde \xi, \widetilde a \widetilde \xi \right)_{\widetilde \H}\ge \left(\left(\widetilde f_j \widetilde a\widetilde f_j \right) \widetilde \xi, \left(\widetilde f_j \widetilde a\widetilde f_j \right) \widetilde \xi \right)_{\widetilde \H}> 0,
$$
so $ \widetilde a \widetilde \xi \neq 0$ and $\widetilde{\rho}\left( \widetilde{a}\right) \neq 0$.

Let us prove that the representation $\widetilde \rho$ is nondegenerate. If $\widetilde \xi = \sum_{j-1}^n \widetilde{f}_j \otimes \xi_j\in C_c\left( \widetilde \sY\right)\otimes_{C_0\left( \sY\right) } \H$ then 
$$
\widetilde \xi \neq 0 \quad \Leftrightarrow \quad \left(\widetilde \xi , \widetilde \xi\right)_{\widetilde{\H}} =\sum_{\substack{j=1\\k=1}}^{n}\left(\xi_j, \left\langle \widetilde f_j, \widetilde f_k \right\rangle_A \xi_k  \right)_\H > 0.
$$
The union $\widetilde\sZ \bydef \bigcup_{j=1} \supp \widetilde f_j$ is compact and if $\widetilde f$ is a covering sum of $\widetilde\sZ$ (cf. Definition \ref{top_covering_sum_defn}) then $\widetilde f \widetilde f_j = \widetilde f_j$ for any $j=1,..., n$. If $\left\{u_\a\right\}_{\a \in \mathscr A}\subset A$ is an approximate unit of $A$ then for all $\xi \in \H$ one has 
$$
\lim_{\a \in \mathscr A}u_\a \xi = \xi.
$$
So for any $j,k \in \left\{1,..., n\right\}$ there is $\a_{jk} \in \mathscr A$ such that
\be\label{groupoid_inen2_eqn}
\a \ge \a_{jk}\quad \Rightarrow\quad \left\|\left(\xi_j, \left\langle \widetilde f_j, \widetilde f_k \right\rangle_A \xi_k  \right)- \left(u_{\a}\xi_j, \left\langle \widetilde f_j, \widetilde f_k \right\rangle_A u_{\a} \xi_k  \right) \right\| < \frac{\left(\widetilde \xi , \widetilde \xi\right)_{\widetilde{\H}}}{n^2}
\ee
If $\a_0$ is such that $\a_0 \ge \a_{jk}$ for all  $j,k \in \left\{1,..., n\right\}$ and $\widetilde a \bydef \varphi'\left( \widetilde f\otimes u_{\a_0}\right)$ where $\varphi'$ is given by \eqref{blowing_tensoar_eqn}. From \eqref{groupoid_inen2_eqn} and the triangulate identity it follows that
$$
\left| \left(\widetilde \xi , \widetilde \xi\right)_{\widetilde{\H}}- \left(\widetilde a\widetilde \xi , \widetilde a\widetilde \xi\right)_{\widetilde{\H}}\right| < \left(\widetilde \xi , \widetilde \xi\right)_{\widetilde{\H}}
$$
so 
$$
\left(\widetilde a\widetilde \xi , \widetilde a\widetilde \xi\right)_{\widetilde{\H}}> 0  \quad \Rightarrow\quad \widetilde a \widetilde \xi\neq 0.
$$
The subspace of elements $\widetilde \xi = \sum_{j-1}^n \widetilde{f}_j \otimes \xi_j$ is dense in $\widetilde\H$ so for all $\widetilde\xi \in \widetilde\H$ there is $\widetilde a \in A_0\left(\widetilde\sY \right)$ such that $\widetilde a \widetilde\xi \neq 0$.

\end{proof}

\subsection{Finite-fold coverings}
\paragraph{}
Here we establish correspondence between noncommutative finite-fold covering and coverings of Hausdorff blowing-ups.
\begin{theorem}\label{blowing_sufficient_covering_thm}
	Let $C_0\left(\sY \right) \hookto M\left(A \right)$ be a Hausdorff blowing-up (cf. Definition \ref{blowing_defn}) with  locally connected $\sY$, and let 
	 $q: \widetilde\sY \to \sY$  be a finite-fold transitive $A$-regular covering (cf. Definition \ref{blowing_a_regular_defn}) with connected $\widetilde \sY$. If	$A_0\left(q\right) : A \hookto	 A_0\left( \widetilde \sY\right)$ is a {finite}-$q$-{lift} (cf. Definition \ref{blowing_finite_lift_defn})  then the quadruple $$\left(A,	 A_0\left( \widetilde \sY\right), G\left(\left.\widetilde \sY\right|\sY \right), A_0\left(q\right)\right)$$ is a noncommutative finite-fold covering (cf. Definition \ref{fin_defn}).
\end{theorem}

\begin{proof}
	From the Lemma \ref{blowing_proper_group_lem} it follows that $\left(A,	 A_0\left( \widetilde \sY\right), G\left(\left.\widetilde \sY\right|\sY \right), A_0\left(q\right)\right)$ is a noncommutative finite-fold pre-covering (cf. Definition \ref{fin_pre_defn}). From the proof of the Lemma \ref{top_fin_sufficient_lem} it turns out that there is a family $\left\{\sU_\la\subset \sY\right\}_{\la\in\La}$ indexed by directed set $\La$ such that
	\begin{enumerate}
		\item[(i)] for all $\la\in \La$ a quadruple $$\left(C_0\left(\sU_\la \right), C_0\left( \widetilde{\sU}_\la\bydef p^{-1}\left(\sU_\la \right) \right) ,  G\left(\left.\widetilde \sY\right|\sY \right), C_0\left( q\right)|_{C_0\left(\sU_\la \right)}  \right)$$ is  a noncommutative covering with unitization (cf. Definition \ref{fin_unitization_defn}),
		\item[(ii)] ${\sY}= \bigcup_{\la\in\La} \sU_\la$.
	\end{enumerate}
	  On the other hand from the Lemma \ref{cov_mult_fin_lem} it turns out that $$
	  \left(C_b\left(\sU_\la \right), C_b\left( \widetilde{\sU}_\la\right) ,  G\left(\left.\widetilde \sY\right|\sY \right), M\left( C_0\left( q\right)|_{C_0\left(\sU_\la \right)} \right)  \right)$$ 
	  is an unital noncommutative finite-fold covering (cf. Definition \ref{fin_unital_defn}). For all $\la\in \La$ denote by $A_\la \bydef ~_{\sU_\la}A_{\sU_\la}$ and $\widetilde A_\la \bydef ~_{\widetilde\sU_\la} A_0\left( \widetilde \sY\right)_{\widetilde \sU_\la}$. From our construction it follows that for any $\la \in \La$ there are inclusions $C_b\left(\sU_\la \right)\subset M\left( A_\la\right)$ and $C_b\left(\widetilde\sU_\la \right)\subset M\left( \widetilde A_\la\right)$. Let  both $B_\la$ and $\widetilde B_\la$  be $C^*$-subalgebras of  $M\left( A_\la\right)$ and  $M\left(\widetilde A_\la\right)$   generated by unions 	$A_\la \cup C_b\left(\sU_\la \right)$ and $\widetilde A_\la \cup C_b\left(\widetilde \sU_\la \right)$  respectively. Since  $\left(C_b\left(\sU_\la \right), C_b\left( \widetilde{\sU}_\la\right) ,  G\left(\left.\widetilde \sY\right|\sY \right), M\left( C_0\left( q\right)|_{C_0\left(\sU_\la \right)} \right)  \right)$ is an unital noncommutative finite-fold covering  then similarly to \eqref{top_finite_covering_basis_eqn} there is a finite subset $\left\{\widetilde{e}_1, ..., \widetilde{e}_n\right\}\subset C_b\left( \widetilde{\sU}_\la\right)_+$ such that
	\be\label{blowing_finite_covering_basis_eqn}
	\begin{split}
	1_{C_b\left(\widetilde{\mathcal Y} \right) }= \sum_{j=1 }^n\widetilde{e}^2_{\widetilde{j}},\\
	\forall j =1,...,n \quad \widetilde{e}_{j}	\left(g \widetilde{e}_{j} \right) = 0 \quad \text{for any nontrivial}\quad g\in G\left(\left.\widetilde \sX \right|\sX \right).
\end{split}
\ee	
From \eqref{blowing_finite_covering_basis_eqn} it follows that 
$$
\widetilde{b}\in \widetilde B_\la \quad \Rightarrow \quad \widetilde{b} = \sum_{j=1 }^n \widetilde e_j b_j; \quad b_j \bydef \sum_{g \in G\left(\left.\widetilde \sY\right|\sY \right)} g\left(\widetilde e_j  \widetilde{b}\right)\in B, 
$$
so one has
$$
\widetilde{B} = \widetilde{e}_1 B + ... +  \widetilde{e}_n B,
$$
i.e. $\widetilde{B}$ is a finitely generated $B$-module. So if $\pi_{B_\la}:B_\la \hookto\widetilde{B}_\la$ is the natural $*$-homomorphism the quadruple 
$
\left(B_\la, \widetilde{B}_\la,  G\left(\left.\widetilde \sY\right|\sY \right), \pi_{B_\la} \right) 
$
is an unital noncommutative finite-fold covering. Since both inclusions $A_\la \subset B_\la$ and $\widetilde A_\la \subset \widetilde B_\la$ are unitizations the quadruple $
\left(A_\la, \widetilde{A}_\la,  G\left(\left.\widetilde \sY\right|\sY \right), \pi_{B_\la}|_{A_\la} \right) 
$
is  a  noncommutative finite-fold covering with unitization, so $A_\la$ is  $\left(A,	 A_0\left( \widetilde \sY\right), G\left(\left.\widetilde \sY\right|\sY \right), A_0\left(q\right)\right)$-strictly proper  (cf. the Definition \ref{strictly_proper_defn}). 
If $ a \in K\left( A\right) $ is an element of Pedersen's ideal (cf. Definition \ref{pedersen_ideal_defn})  then $\supp  a$ is compact. From ${\sY}= \bigcup_{\la\in\La} \sU_\la$ it follows that there is $\la_a \in \La$ such that $\supp  a\subset \sU_{\la_a}$ and $ a \in A_{\la_a}$. The union $\bigcup_{\la \in \La} A_\la$ is dense in $A$ because $K\left( A\right)$ is dense in $A$.
\end{proof}
\begin{empt}
Under the conditions of the Theorem \ref{blowing_sufficient_covering_thm} 
\end{empt}
\begin{empt}\label{blowing_ext_empt}
	Let $\left( A, \widetilde A, G, \pi\right)$ be a noncommutative finite-fold quasi-covering (cf. Definition \ref{fin_quasi_defn}), and let 
 $M\left(\pi\right) : M\left( A\right)   \hookto M\left( \widetilde A\right) $ be given by the Lemma \ref{fin_multilier_lem}  injective $*$- homomorphism, and let $C_0\left( \sY\right) \subset M\left( A\right)$ be Hausdorff blowing-up (cf. Definition \ref{blowing_defn}). Let $\widetilde \sX$ be the spectrum of   $\widetilde A$.  For any $\widetilde x \in\widetilde \sX$ the $C^*$-algebra 
	$$
	\widetilde{C}_{\widetilde{x}}\bydef \left\{c = \rep_{  \widetilde{x }}\left( a\right)\left|  a \in C_0\left(\sY \right)\right. \right\}\subset \rep_{  \widetilde{x }}\left(M\left( \widetilde A\right) \right)  
	$$
	is commutative. So a $C^*$-algebra
	$$
	\widetilde{C}\bydef \left\{\left.\widetilde c \in M\left( \widetilde A\right)\right| \rep_{  \widetilde{x }}\left( \widetilde c\right) \in \widetilde{C}_{\widetilde{x}} \right\}\subset M\left( \widetilde A\right)
	$$
	is also commutative. There is a locally compact Hausdorff space $\widetilde\sY$ such that
	\be\label{blowing_ext_eqn}
	\widetilde{C}\cong C_0\left(\widetilde \sY\right)
	\ee	
	If for any $\widetilde a\in \widetilde A$ and $\eps > 0$ there is $a \in A$   such that
	\be\label{blowing_dens_cond}
	\left\|\widetilde a - \widetilde a \pi\left(a \right) \right\| < \eps 
	\ee
	then 
	\bean
	\begin{split}
		\widetilde A 	C_c\left(  \sY\right)\bydef \left\{ \widetilde a f \left|  f \in C_c\left( \sY\right)\quad \widetilde a \in \widetilde A \right.\right\},\\
	\end{split}
	\eean 
	are dense in $A$.
	There is an injective *-homomorphism
	\be\label{blowing_cb_eqn}
	C_b\left( \sY\right) \hookto C_b\left( \widetilde \sY\right) 
	\ee
\end{empt}

\begin{definition}\label{blowing_ext_defn}
	Under the hypotheses \ref{blowing_ext_empt} we say that Hausdorff blowing-up is \textit{induced by noncommutative quasi-covering} $\left( A, \widetilde A, G, \pi\right)$.
\end{definition}

\begin{lemma}\label{blowing_unitization_lem} 
	Let  $\left(A, \widetilde{A}, G, \pi \right)$  be noncommutative finite-fold covering with unitization (cf. Definition  \ref{fin_unitization_defn}), and let $C\left(\sY\right)\to M\left( A\right)$ be Hausdorff blowing-up with compact $\sY$. Let $\rho: C\left( \sY\right) \hookto C\left( \widetilde \sY\right)$  be given by the equation \eqref{blowing_cb_eqn}. The action $G \times M\left( \widetilde{A}\right) \to M\left( \widetilde{A}\right)$ induces an action  $G\times C\left(\widetilde \sY \right) \to C\left(\widetilde \sY \right)$. Then $C\left(\widetilde \sY \right)$ is a finitely generated $C\left(\sY \right)$-module. Moreover the natural covering $q: \widetilde \sY \to \sY$ is $A$-regular (cf. Definition \ref{blowing_a_regular_defn}).  
\end{lemma}
\begin{proof}
	From the Corollary \ref{cov_mult_fin_c_cor} it follows that there is  an 
	unital  noncommutative finite-fold covering  $\left(M\left(A \right)  ,M\left( \widetilde{A}\right) , G, \widetilde{\pi} \right)$ (cf. Definition \ref{fin_unital_defn}). If $C\left(\widetilde \sY \right)$ is not a finitely generated $C\left(\sY \right)$ then from the Theorem \ref{fin_gen_thm} it follows that  there is an infinite algebraic sum
	$$
	C\left(\widetilde \sY \right) = \sum_{\la\in \La} M_\la,
	$$
	of right $C\left(\sY \right)$-modules, such that for any finite subset $\La_0 \subset \La$ one has
	$$
	C\left(\widetilde \sY \right) \neq \sum_{\la\in \La_0} M_\la,
	$$
	One has an algebraic sum 
	$$
	M\left(\widetilde A \right)  = \sum_{\la\in \La} M\left( A\right)  M_\la
	$$
	of left  $M\left(A \right)$-nodules such that for any finite subset $\La_0 \subset \La$ one has
	$$
	M\left( \widetilde{A}\right) \neq \sum_{\la\in \La_0} M\left( A\right)  M_\la.
	$$
	From 
	$$
\forall \widetilde a \in \widetilde A \quad \widetilde a = \widetilde a \varphi_C\left( 1_{C\left(\widetilde\sY \right) }\right) 	
	$$
	and the Lemma \ref{blowing_lift_eq_lem} it follows that the natural covering $q: \widetilde \sY \to \sY$ is $A$-regular (cf. Definition \ref{blowing_a_regular_defn}).
\end{proof}
\begin{remark}\label{blowing_unitization_rem}
	From the Theorem \ref{pavlov_troisky_thm} and the Lemma  \ref{blowing_unitization_lem}  it follows that there is a natural transitive finite-fold covering $q: \widetilde \sY\to \sY$.
\end{remark}

\begin{lemma}\label{blowing_compact_unitization_lem}
		Let $C\left(\sY \right) \hookto M\left(A \right)$ be a Hausdorff blowing-up (cf. Definition \ref{blowing_defn}) with compact, locally connected space $\sY$, and let 
	$q: \widetilde\sY \to \sY$  be a finite-fold  $A$-regular covering (cf. Definition \ref{blowing_a_regular_defn})
	Let   both $\sY$ and $\widetilde \sY$ are connected, compact, Hausdorff spaces.
	If $A_0\left(q\right) : A \hookto	 A_0\left( \widetilde \sY\right)$ is a {finite}-$q$-{lift} (cf. Definition \ref{blowing_finite_lift_defn}) then the quadruple
	\be\label{blowing_pre_eqn}
	\left(A,	 A_0\left( \widetilde \sY\right), A_0\left(q\right), G\left(\left.\widetilde \sY\right|\sY \right) \right)
	\ee
	is a noncommutative finite-fold covering with unitization (cf. Definition \ref{fin_unitization_defn})
\end{lemma}

\begin{proof}
	From the Theorem \ref{blowing_sufficient_covering_thm} it follows that the quadruple \eqref{blowing_pre_eqn} a noncommutative finite-fold pre-covering  (cf. Definition \ref{fin_pre_defn}). There is a given by the equation \eqref{mult_map_eqn} injective $*$-homomorphism $M\left(\pi\right): M\left(A \right) \hookto M\left( A_0\left( \widetilde \sY\right)\right)$. 
	If both $B$ and $\widetilde B$ are generated by $A \cup C\left( \sY\right)$ and $\widetilde A \cup C\left( \widetilde \sY\right)$ subalgebras of $M\left(A \right)$ and $M\left(A_0\left( \widetilde \sY\right) \right)$ then the triple $\left(B, \widetilde B, G, M\left(\pi \right)|_B  \right)$  is a noncommutative finite-fold   pre-covering.
	If  $\left\{\widetilde{e}_1, ..., \widetilde{e}_n\right\}\subset C\left(\widetilde \sY \right)_+$ satisfies to \eqref{top_finite_covering_basis_eqn} then one has
	\bean
	\forall \widetilde b \in \widetilde B\quad \widetilde b = \sum_{j=1}^n \widetilde e_j b_j, \quad b_j = \sum_{g \in G} g\left(\widetilde e_j \widetilde b \right)\in B;\\
	\widetilde B = \widetilde e_1 B + ... + \widetilde e_n B,
	\eean
	i.e.  $\widetilde B$ is a finitely generated $B$-module. So $\left(B, \widetilde B, G, M\left(\pi \right)|_B  \right)$ is an unital noncommutative finite-fold  covering (cf. Definition \ref{fin_unital_defn}).
\end{proof}

	Let $C\left(\sY \right) \hookto M\left(A \right)$ be a Hausdorff blowing-up (cf. Definition \ref{blowing_defn}) with compact,  locally connected $\sY$.
If $\left(A, \widetilde A,  G, \pi\right)$  is a noncommutative finite-fold covering with unitization (cf. Definition \ref{fin_unitization_defn}) 
then from the Lemma  \ref{blowing_unitization_lem}  and the Remark \ref{blowing_unitization_rem} it follows that there is a finite-fold covering $q: \widetilde{\sY}\to \sY$ with Hausdorff blowing-up $C\left( \widetilde{\sY} \right) \hookto M\left( \widetilde A\right)$. From the Lemma \ref{blowing_compact_unitization_lem}  it follows that  a quadruple 
\be\label{blowing_c_fin_eqn}
\left(A, A_0\left(\widetilde\sY \right) , G\left(\left.\widetilde\sY \right|\sY \right) , A_0\left(q \right) \right)
\ee
is a noncommutative finite-fold covering with unitization.

\begin{lemma}\label{blowing_ah_lem}
Under the above hypotheses if 
	\be\label{blowing_ah_eqn}
H \bydef \left\{ g \in G \left| \forall \widetilde a \in  A_0\left(\widetilde \sY \right) \quad g\widetilde a=\widetilde a \right.\right\}
\ee
then $\widetilde A^H = A_0\left(\widetilde \sY \right)$. 
\end{lemma}
\begin{proof}
	There is a natural action $G\left(\left.\widetilde\sY \right|\sY \right)  \times \widetilde A^H \to \widetilde A^H$ which comes from the action $G \times \widetilde A\to \widetilde A$. From the Theorem \ref{pavlov_troisky_thm} it follows that there   is a finite subset $\left\{\widetilde f_1, ..., \widetilde f_n\right\} \subset C\left( \widetilde \sY\right)_+$ such that 
	\begin{itemize}
		\item $\widetilde f_1+ ...+ \widetilde f_n = 1_{C\left( \widetilde \sY\right) }$,
		\item for all $j=1,...,n$ a support  $\supp \widetilde f_j$ is homeomorphically mapped onto $q\left( \supp\widetilde f_j\right)$.
	\end{itemize}
	It follows that any $\widetilde a \in \widetilde A^H$ is given by $\widetilde a = \widetilde a\widetilde f_1+...+\widetilde a\widetilde f_n$.	There are $a_j$ and $a^{\sqrt{~}}_j \in A$ such that $\sum_{	g \in  G\left(\left.\widetilde\sY\right|\sY\right)}\widetilde a \widetilde f_j = A_0\left( q\right)\left(  a_j\right)  $ and $\sum_{	g \in  G\left(\left.\widetilde\sY\right|\sY\right)}\widetilde a \sqrt{\widetilde f_j} = A_0\left( q\right)\left(  a^{\sqrt{~}}_j\right)$. If $\widetilde a_j \bydef  A_0\left( q\right)\left(  a^{\sqrt{~}}_j\right)\sqrt{\widetilde f_j}$ then $\sum_{	g \in  G\left(\left.\widetilde\sY\right|\sY\right)}\widetilde a_j = a_j$. Let 
	$$
	\widetilde \sU_j \bydef \left\{\left. \widetilde x \in \widetilde\sY \right| \widetilde f_j \left(\widetilde x_j \right) > 0 \right\},
	$$
	and let be $\widetilde A^H_{\widetilde \sU_j}$ is the $C^*$-norm completion of $\widetilde A^HC_0\left(\widetilde \sU_j \right)$. One has 
	\bean
	\widetilde a_j \in \widetilde A^H_{\widetilde \sU_j},\\
		\widetilde a\widetilde f_j \in \widetilde A^H_{\widetilde \sU_j}, 
	\eean 
	and $\widetilde A^H_{\widetilde \sU_j}\cap \widetilde A^H_{g\widetilde \sU_j}= \{0\}$  for any nontrivial $g \in  G\left(\left.\widetilde\sY\right|\sY\right)$. From this circumstance it follows that
	$$
	\widetilde a_j - \widetilde a\widetilde f_j\neq 0 \quad \Rightarrow \quad \sum_{	g \in  G\left(\left.\widetilde\sY\right|\sY\right)}g\left( \widetilde a_j - \widetilde a\widetilde f_j\right) \neq 0.
	$$
	But it is impossible since
	$$
	\widetilde a_j - \widetilde a\widetilde f_j\neq 0 \quad \Rightarrow \quad \sum_{	g \in  G\left(\left.\widetilde\sY\right|\sY\right)}g\left( \widetilde a_j - \widetilde a\widetilde f_j\right) = A_0\left( q\right) \left( a_j\right) - A_0\left( q\right) \left( a_j\right) = 0.
	$$	
	This contradiction proves that $\widetilde a_j = \widetilde a\widetilde f_j$ and $\widetilde a = \widetilde a_1+...+\widetilde a_n \in A_0\left( \widetilde\sY\right) $.

\end{proof}
\begin{empt}

From the Lemma \ref{blowing_compact_unitization_lem} it follows that 
	the group
	\bean
	H \bydef \left\{\left. g \in G \right| \forall a \in  A_0\left(q \right)\left( A\right) \quad ga=a \right\}
	\eean
	is $\left(A,  \widetilde A, G, \pi\right)$-proper (cf. Definition \ref{proper_subgroup_fin_defn}). From the Lemma \ref{fin_composition_unitization_lem} it follows that there are $H$-{regular}  $H$-{singular} noncommutative finite-fold coverings with unitization (cf. Definition \ref{fin_unitization_defn})
	\bea\label{blowing_reg_sing_m_eqn}
	\left(A, A_0\left(\widetilde\sY \right) , G\left(\left.\widetilde\sY \right|\sY \right) \cong G/H, A_0\left(q \right) \right),\\
	\label{blowing_reg_sing_eqn}
	\left( A_0\left(\widetilde\sY \right)\cong \widetilde A^H , \widetilde A,H ,\left.\Id_{\widetilde A} \right|_{ A_0\left(\widetilde\sY \right)}\right) 	
	\eea
	(cf. Definition \ref{fin_composition_covering_defn}).  From the Corollary \ref{cov_mult_fin_c_cor} it follows that 		
	\be\label{blowing_frame_l_eqn}
	\begin{split}
		\exists \left\{ \widetilde{a}_1, ..., \widetilde{a}_n, \widetilde{b}_1, ..., \widetilde{b}_n\right\} \subset   M\left(  \widetilde{A}\right) \quad \forall \widetilde{a} \in \widetilde{A}\quad \widetilde{a}=	\sum_{j=1}^{n}\widetilde{b}_j a_j \\
		\forall j\in \left\{1,...,n\right\}\quad a_j\bydef  \sum_{g \in H} g\left(\widetilde{a}^*_j\widetilde{a} \right)\in M\left( \left.\Id_{\widetilde A} \right|_{ A_0\left(\widetilde\sY \right)}\right) M\left(  A_0\left(\widetilde\sY \right)\right) .
	\end{split}
	\ee
	(cf. the equation \eqref{cov_mult_fin_c_eqn})	
	
	 There are Hausdorff blowing-up  $C\left( \widetilde\sY\right)  \subset  M\left(A_0\left(\widetilde \sY\right) \right)$ and $C\left( \widetilde\sY\right)  \subset M\left(\widetilde A\right)$  one.

\end{empt}
\begin{lemma}
	Under the above hypotheses if a transitive covering $q': \widetilde\sY' \to \widetilde\sY$ is $A_0\left(\widetilde \sY\right)$-regular then it is  $\widetilde A$-regular.
\end{lemma}
\begin{proof}
	Similarly to \eqref{blowing_tensor_eqn} there are homomorphisms 
	\bean
		\phi': C_0\left( \widetilde \sY' \right)\otimes_{C_0\left(\widetilde \sY\right) } K\left(  A_0\left(\widetilde \sY \right)\right)   \to \End^*_A\left( \mathscr L^2\left(\widetilde \sY' \right)_{\left(  A_0\left(\widetilde \sY \right)\right)_0 \left(  \widetilde \sY'\right)}\right),\\
		\phi'': C_0\left( \widetilde \sY' \right)\otimes_{C_0\left(\widetilde \sY\right) } K\left(\widetilde A\right) \to \End^*_A\left( \mathscr L^2\left(\widetilde \sY' \right)_{\widetilde A}\right)
	\eean 
	of right $C_0\left( \widetilde \sY' \right)$-modules.
	If $$
	\widetilde a \in\phi''\left(  C_0\left( \widetilde \sY' \right)\otimes_{C_0\left(\widetilde \sY\right) } K\left(\widetilde A\right)\right) 
	$$ is a positive then $\widetilde a^\oplus  \bydef \sum_{	g \in {H}}\widetilde a\in K\left( A_0\left(\widetilde \sY\right)\right)$ is also positive. Since $\widetilde a^\oplus$ is $H$-invariant one can suppose that 
	$$
\widetilde a^\oplus\in \phi'\left( C_0\left( \widetilde \sY' \right)\otimes_{C_0\left(\widetilde \sY\right) } K\left(\left(  A_0\left(\widetilde \sY \right)_0\left(  \widetilde \sY'\right)\right)_0\right) \right).
$$	
From the Lemma \ref{blowing_lift_eq_lem} it follows that there is a positive $\widetilde f' \in C_c \left( \widetilde \sY'\right)$  such that
$$
\widetilde a^\oplus = \widetilde f'\widetilde a^\oplus\widetilde f'
$$
and taking into account $\widetilde a < \widetilde a^\oplus$ one has
$$
\widetilde a = \widetilde f'\widetilde a\widetilde f'
$$
Now from the Lemma \ref{blowing_lift_eq_lem} it follows that  $q': \widetilde\sY' \to \widetilde\sY$ is $A_0\left(\widetilde \sY\right)$-regular then it is  $\widetilde A$-regular.
\end{proof}
\begin{lemma}\label{blowing_h_act_lem}
Under the above hypotheses the natural action 
\bean
\begin{split}
	H \times \phi''\left( C_0\left( \widetilde \sY' \right)\otimes_{C_0\left(\widetilde \sY\right) } K\left(\widetilde A\right)\right) \to\phi''\left( C_0\left( \widetilde \sY' \right)\otimes_{C_0\left(\widetilde \sY\right) } K\left(\widetilde A\right)\right),\\
	\left(g, \phi''\left( \widetilde f \otimes a \right) \right) \mapsto \phi''\left( \widetilde f \otimes a \right).
\end{split}
\eean
yields an action 
\be\label{blowing_h_act_eqn}
	H \times \widetilde A_0\left(\widetilde \sY'\right) \to \widetilde A_0\left(\widetilde \sY'\right)
\ee
such that any $g \in H$ is a $*$-automorphism.
\end{lemma}
\begin{proof}
If $\widetilde a \in \phi''\left( C_0\left( \widetilde \sY' \right)\otimes_{C_0\left(\widetilde \sY\right) } K\left(\widetilde A\right)\right)$ then from the equation \eqref{blowing_tensor_eqn} one has
$$
\widetilde a = 	\sum_{j=1}^n  	\varphi_C\left( \widetilde f'_j\right) \varphi_A\left(a_j \right)
$$
and
$$
\forall g \in H \quad g \widetilde a = \sum_{j=1}^n  	\varphi_C\left( \widetilde f'_j\right) \varphi_A\left(g a_j \right)
$$
If $\widetilde b \in \phi''\left( C_0\left( \widetilde \sY' \right)\otimes_{C_0\left(\widetilde \sY\right) } K\left(\widetilde A\right)\right)$ from the Lemma \ref{blowing_lift_eq_eqn} one can deduce that 
$$
\widetilde b = 	\sum_{l=1}^m  	 \varphi_A\left(b_l \right) \varphi_C\left( \widetilde f''_l\right)
$$
and
$$
\forall g \in H \quad g \widetilde b = \sum_{l=1}^m  	 \varphi_A\left(g b_l \right) \varphi_C\left( \widetilde f''_l\right).
$$
From the above equations one has
\bean
 \left( g  \widetilde a\right)  \left( g \widetilde b\right) = \sum_{\substack{j=1\\l=1}}^{\substack{j=n\\l=m}}\varphi_C\left( \widetilde f'_j\right) \varphi_A\left(ga_j \right) \varphi_A\left(g b_l \right) \varphi_C\left( \widetilde f''_l\right)=\\= \sum_{\substack{j=1\\l=1}}^{\substack{j=n\\l=m}}\varphi_C\left( \widetilde f'_j\right) \varphi_A\left(g\left( a_j  b_l\right)  \right) \varphi_C\left( \widetilde f''_l\right)= g\left(  \widetilde a\widetilde b\right). 
\eean 
Similarly one can prove that  for each $\widetilde a,  \widetilde b\in \phi''\left( C_0\left( \widetilde \sY' \right)\otimes_{C_0\left(\widetilde \sY\right) } K\left(\widetilde A\right)\right)$ one has
\bean
g \widetilde a^* = \left(g \widetilde a \right)^*,\\
g \widetilde a + g \widetilde b = g\left( \widetilde a\widetilde b\right).
\eean 
Since $\phi''\left( C_0\left( \widetilde \sY' \right)\otimes_{C_0\left(\widetilde \sY\right) } K\left(\widetilde A\right)\right)$ id dense in $\widetilde A_0\left(\widetilde \sY'\right)$ we obtain the action \eqref{blowing_h_act_eqn}.

\end{proof}
From the Lemma \ref{blowing_h_act_lem} it follows that $A_0\left(\widetilde \sY'\right) = \widetilde A_0\left(\widetilde \sY'\right)^H$, a noncommutative finite-fold quasi-covering
\be\label{blowing_inf_eqn}
\left(A_0\left(\widetilde\sY' \right), \widetilde A_0\left(\widetilde\sY' \right), H, \left.\Id_{\widetilde A_0\left(\widetilde\sY' \right)} \right|_{ A_0\left(\widetilde\sY' \right)} \right) 
\ee
(cf. Definition \ref{fin_quasi_defn})
 which can be naturally extended uo to an injective $*$-homomorphism 
$$
M\left( \left.\Id_{\widetilde A_0\left(\widetilde\sY' \right)} \right|_{ A_0\left(\widetilde\sY' \right)}\right)  : M\left(  A_0\left(\widetilde \sY'\right)\right)\hookto M\left( \widetilde A_0\left(\widetilde \sY'\right)\right).
$$
(cf. Lemma \ref{fin_multilier_lem}). From the Lemma \ref{blowing_mult_iclusion_lem} it follows that there are natural inclusions
\bean
\varphi' : M\left(  A_0\left(\widetilde \sY\right)\right) \hookto M\left(  A_0\left(\widetilde \sY'\right)\right) ,\\
\varphi'': M\left( \widetilde A\right) \hookto M\left( \widetilde A_0\left(\widetilde \sY'\right) \right).
\eean 

(cf. the equations \ref{blowing_reg_sing_eqn})
Otherwise the Lemma \ref{fin_multilier_lem} yields an inclusion 
\bean
M\left(\left.\Id_{\widetilde A} \right|_{ A_0\left(\widetilde\sY \right)} \right)  : M\left(  A_0\left(\widetilde \sY\right)\right) \hookto M\left(  \widetilde A\right) 
\eean
such that a following diagram
\\
\begin{tikzcd}
M\left(  A_0\left(\widetilde \sY\right)\right) \arrow[rrrr, "M\left(\left.\Id_{\widetilde A} \right|_{ A_0\left(\widetilde\sY \right)}\right)"]\arrow[d , "\varphi'"]&&&& M\left(  \widetilde A\right)\arrow[d, "\varphi''"]\\
M\left(  A_0\left(\widetilde \sY'\right)\right)\arrow[rrrr, "M\left( \left.\Id_{\widetilde A_0\left(\widetilde\sY' \right)} \right|_{ A_0\left(\widetilde\sY' \right)}\right)  "]  &&&&M\left( \widetilde A_0\left(\widetilde \sY'\right)\right)  
\end{tikzcd}
\\
is commutative. If $\left\{ \widetilde{a}_1, ..., \widetilde{a}_n, \widetilde{b}_1, ..., \widetilde{b}_n\right\} \subset M\left( \widetilde{A}\right)$ is the  given by the equation \eqref{blowing_frame_l_eqn}  set then
from the equation \eqref{blowing_frame_l_eqn} it follows that
\be\label{blowing_frame_ll_eqn}
\begin{split}
\forall \widetilde{a} \in M\left( \widetilde A_0\left(\widetilde \sY'\right)\right) 	\quad	\widetilde{a}=	\sum_{j=1}^{n}\varphi''\left( \widetilde{b}_j\right) a_j\\ \text{ where }  a_j\bydef  \sum_{g \in G} g\left(\varphi''\left( \widetilde{a}^*_j\right) \widetilde{a} \right)\in M\left( \left.\Id_{\widetilde A_0\left(\widetilde\sY' \right)} \right|_{ A_0\left(\widetilde\sY' \right)}\right)\left( M\left(  A_0\left(\widetilde\sY \right)	\right) \right).
\end{split}
\ee
where the natural inclusion $ A_0\left(\widetilde\sY \right)	\subset  A_0\left(\widetilde\sY' \right)$ is implied 
	
\begin{lemma}\label{blowing_inf_lem}
Under the  above hypotheses the given by \eqref{blowing_inf_eqn} noncommutative finite-fold quasi-covering (cf. Definition \ref{fin_quasi_defn}) 
is a noncommutative finite-fold covering with unitization (cf. Definition \ref{fin_unitization_defn}).
\end{lemma}
\begin{proof}
 From the Corollary \ref{ind_mult_inv_cor} it follows that 
\be\label{fin_quasi_m_defn}
 \left(M\left( A_0\left(\widetilde\sY' \right)\right) ,M\left(  \widetilde A_0\left(\widetilde\sY' \right)\right) , H, M\left( \left.\Id_{\widetilde A_0\left(\widetilde\sY' \right)} \right|_{ A_0\left(\widetilde\sY' \right)}\right)  \right) 
\ee
 is noncommutative finite-fold quasi-covering. If $\left\{\widetilde u_\la\right\}_{\la\in \La}\subset C_0\left( \widetilde\sY\right)$ is an approximate unit for $ C_0\left( \widetilde\sY\right)$ (cf. Definition \ref{approximate_unit_defn}) then
\bean
  \varphi''\left( M\left(  \widetilde A\right)\right)= \left\{\left.a \in  M\left(  \widetilde A_0\left(\widetilde\sY' \right)\right)\right| a = \b\text{-}\lim  a\varphi_C\left( \widetilde u_\la\right) \right\}
\eean
 where $\bt$-$\lim$ implies the convergence with respect to the strict topology of $M\left(  \widetilde A_0\left(\widetilde\sY' \right)\right)$ (cf. Definition \ref{strict_topology_defn}).
 If 
 \bean
H' \bydef\\ \left\{g \in \Aut\left(M\left(  \widetilde A_0\left(\widetilde\sY' \right)\right) \right)  \left| \forall a\in  M\left( \left.\Id_{\widetilde A_0\left(\widetilde\sY' \right)} \right|_{ A_0\left(\widetilde\sY' \right)}\right)\left( M\left(  A_0\left(\widetilde\sY \right)	\right) \right) ~ ga = a\right.\right\} 
\eean
then
\bean
\forall a \in \varphi''\left( M\left(  \widetilde A\right)\right) \quad \forall g \in H' \quad ga =  g\left(  \b\text{-}\lim  a\varphi_C\left( \widetilde u_\la\right)\right)=\\ \left( ga\right) \b\text{-}\lim  \varphi_C\left( \widetilde u_\la\right) \in  \varphi''\left( M\left(  \widetilde A\right)\right).
\eean 
or equivalently
$$
H' \varphi''\left( M\left(  \widetilde A\right)\right) =  \varphi''\left( M\left(  \widetilde A\right)\right).
$$
It follows that there is the natural action $H' \times  M\left(  \widetilde A\right)\to  M\left(  \widetilde A\right)$ such that for all $ga = a$ for any $a \in M\left(A_0\left(\widetilde \sY \right)  \right)$. Since the quadruple  \eqref{blowing_reg_sing_eqn} is a noncommutative finite-fold covering with unitization one has $H' \cong H$, i.e. the quadruple is a noncommutative finite-fold pre-covering (cf. Definition \ref{fin_pre_defn})
From the equation \eqref{blowing_frame_ll_eqn} it follows that $M\left( \widetilde A_0\left(\widetilde \sY'\right)\right)$ is a finitely generated  $M\left(  A_0\left(\widetilde \sY'\right)\right)$-module. So the quadruple \eqref{fin_quasi_m_defn} is an unital noncommutative finite-fold covering (cf, Definition \ref{fin_unital_defn}), and the quadruple \eqref{blowing_inf_eqn} is   a noncommutative finite-fold covering with unitization

\end{proof}
\begin{theorem}\label{blowing_sufficient_fin_thm}
	Let $C\left(\sY \right) \hookto M\left(A \right)$ be a Hausdorff blowing-up (cf. Definition \ref{blowing_defn}) with compact, locally connected space $\sY$, and let 
$q': \widetilde\sY' \to \sY$  be an  $A$-regular covering (cf. Definition \ref{blowing_a_regular_defn}), with connected $\widetilde \sY'$. If the $q'$-lift $A_0\left( \widetilde\sY'\right)$ of $A$ (cf. Definition \ref{blowing_lift_defn}) is simply connected (cf. Definition \ref{simply_connected_defn}) then for any finite-fold covering with unitization 
\be\label{blowing_sufficient_fin_eqn}
\left(A, \widetilde A, G, \pi \right) 
\ee
there is an $A$-regular  finite-fold covering $q : \widetilde \sY \to \sY$ such that that the triple \eqref{blowing_sufficient_fin_eqn} is equivalent to a
$$
\left(A, A_0\left( \widetilde\sY\right) , G\left(\left.\widetilde \sY \right| \sY\right) , A_0\left(q \right) \right)
$$
one.

\end{theorem}
\begin{proof}
From the Lemma \ref{blowing_ah_lem} it follows that there is a given by \eqref{blowing_ah_eqn} normal subgroup $H\subset G$ such that there given by the equations \eqref{blowing_reg_sing_m_eqn}, \eqref{blowing_reg_sing_eqn} regular and singular noncommutative finite-fold coverings
\bean
\left(A, A_0\left(\widetilde\sY \right) , G\left(\left.\widetilde\sY \right|\sY \right) \cong G/H, A_0\left(q \right) \right),\\
\left( A_0\left(\widetilde\sY \right)\cong \widetilde A^H , \widetilde A,H ,\left.\Id_{\widetilde A} \right|_{ A_0\left(\widetilde\sY \right)}\right) 	
\eean
If $H$ is not trivial then the Lemma \ref{blowing_inf_lem} yields a nontrivial noncommutative finite-fold covering 
\bean
\left(A_0\left(\widetilde\sY' \right), \widetilde A_0\left(\widetilde\sY' \right), H, \left.\Id_{\widetilde A_0\left(\widetilde\sY' \right)} \right|_{ A_0\left(\widetilde\sY' \right)} \right) 
\eean
(cf. the equation \eqref{blowing_inf_eqn}). Under the hypothesis of this theorem it is impossible, since $A_0\left(\widetilde\sY' \right)$ is simply connected. So the group $H$ is trivial,  and the quadruple \eqref{blowing_sufficient_fin_eqn} is equivalent to the $\left(A, A_0\left( \widetilde\sY\right) , G\left(\left.\widetilde \sY \right| \sY\right) , A_0\left(q \right) \right)$ one.
\end{proof}

\subsection{Finite covering functor}

\paragraph{}
Here we consider a generalization of the functor $A_0$  described in the Section \ref{top_functor_sec}.
\begin{definition}\label{blowing_reg_cov_cat_defn}
	If $A$ is a connected $C^*$-algebra with  Hausdorff blowing-up and $C_0\left(\sY \right)\hookto M\left(A \right)$ then we consider a full subcategory $\mathfrak{FinCov}\text{-}\sY\text{-}A\subseteqq \mathfrak{FinCov}\text{-}\sY$ (cf. Definition \ref{subcategory_defn}) {finite covering category} of $\sY$ (cf. Definition \ref{top_fin_cov_defn}), such that any $\mathfrak{FinCov}\text{-}\sY\text{-}A$ object is $A$-regular. We say that $\mathfrak{FinCov}\text{-}\sY\text{-}A$ is the $A$-\textit{regular finite covering category} of $\sY$.  Sometimes we write $\widetilde \sY$ instead of a covering $\widetilde \sY\to\sY$ to designate an object of $\mathfrak{FinCov}$-$\sY$-$A$.
	\end{definition}

\begin{definition}\label{blowing_alg_reg_cov_cat_defn}
Let $A$ be a connected $C^*$-algebra with  Hausdorff blowing-up and $C_0\left(\sY \right)\hookto M\left(A \right)$. Consider a category $\mathfrak{FinCov}\text{-}A\text{-}\sY$ such that 
\begin{itemize}
	\item any $\mathfrak{FinCov}\text{-}A\text{-}\sY$-object is
	a {finite}-$q$-{lift} $A_0\left(q\right) : A \hookto	 A_0\left( \widetilde \sY\right)$ of $A$ (cf. Definition \ref{blowing_finite_lift_defn})
	\item  any $\mathfrak{FinCov}\text{-}A\text{-}\sY$-morphism from   $A_0\left(q'\right) : A \hookto	 A_0\left( \widetilde \sY'\right)$ to  $A_0\left(q''\right) : A \hookto	 A_0\left( \widetilde \sY''\right)$ is an injective $*$-homomorphism $\pi: A_0\left( \widetilde \sY''\right)\hookto A_0\left( \widetilde \sY'\right)$ such that the following diagram
				\newline
\begin{tikzcd}
	 A_0\left( \widetilde \sY'\right)	& &\arrow[ll,  "\pi" {yshift=10pt}] A_0\left( \widetilde \sY''\right) \\
	& A \arrow[lu, "A_0\left(q'\right)"] \arrow[ru, "A_0\left(q''\right)" {xshift=30pt,yshift=-12pt}]&		
\end{tikzcd}
\\ 	

is commutative.
		\item Compositions of morphisms come from compositions of *-homomorphisms.

\end{itemize}
	We say that $\mathfrak{FinCov}$-$A$-$\sY$ is the \textit{category of finite-fold} $\sY$-\textit{coverings} of $A$.
	Sometimes we write $A_0\left( \widetilde \sY\right)$ instead of $A_0\left(q\right) : A \hookto	 A_0\left( \widetilde \sY\right)$ to designate an object of $\mathfrak{FinCov}\text{-}A\text{-}\sY$.
\end{definition}
\begin{theorem}\label{blowing_functor_thm}
	If $A$ is a connected $C^*$-algebra with  Hausdorff blowing-up and $C_0\left(\sY \right)\hookto M\left(A \right)$ such that the space $\sY$ is locally connected, then one has:
	\begin{enumerate}
		\item [(i)]  the {category $\mathfrak{FinCov}$-$A$-$\sY$ of finite-fold} $\sY$-{coverings} of $A$ (cf. Definition \ref{blowing_alg_reg_cov_cat_defn}) is a subcategory (cf. Definition \ref{subcategory_defn}) of    the {category $\mathfrak{FinCov}$-$A$ of finite-fold coverings} of $A$ (cf. Definition \ref{fin_category_defn}),
		\item[(ii)]the given by the  \ref{blowing_finite_lift_defn} 	{finite}-$q$-{lifts} (e.g. $A_0\left(q\right) : A \hookto	 A_0\left( \widetilde \sY\right)$) of $A$  yield a full functor (cf. Definition \ref{funct_full_faithfull_defn})
		\be\label{blowing_functor_y_eqn}
				\mathfrak{FinCov}\text{-}\sY\text{-}A\xrightarrow{A_0}\mathfrak{FinCov}\text{-}A\text{-}\sY
		\ee
		$A$-{regular finite covering category} of $\sY$ (cf. Definition \ref{blowing_reg_cov_cat_defn}) to the category $\mathfrak{FinCov}\text{-}A\text{-}\sY$.
	\end{enumerate}
\end{theorem}
\begin{proof}
(i) Follows from the Theorem \ref{blowing_sufficient_covering_thm}.\\
(ii) Both  maps   $\widetilde \sY\mapsto A_0\left( \widetilde \sY\right)$ and $\mapsto A_0\left( \widetilde \sY\right)$defines a functor on objects.  From the Theorem \ref{top_covp_cat_thm} it follows that any $\mathfrak{FinCov}\text{-}\sY\text{-}A$-morphism $q^1_2:\widetilde \sY_1 \to \widetilde \sY_2$ is  a transitive finite fold covering and there is a normal subgroup $H\subset G\left(\left.\widetilde \sY\right|\sY \right)$ such that  $\widetilde \sY_2\cong \widetilde \sY_1/H$. From the  Lemma  \ref{blowing_lift_composition_lem} it turns out  that the covering $q^1_2:  \widetilde \sY\to \widetilde \sY/ H$ is $A_0\left(\widetilde \sY_1/ H = \sY_2\right)$-regular. We define 
\be\label{blowing_functor_eqn}
A_0\left(q^1_2 \right) \bydef A_0\left(\widetilde  \sY_2\right)_0\left( q^1_2\right):  A_0\left(\widetilde \sY_2\right)\hookto A_0\left(\widetilde \sY_1\right)
\ee
where  $A_0\left(\widetilde  \sY_2\right)_0\left( q^1_2\right):A_0\left(\widetilde \sY_2\right)\hookto A_0\left(\widetilde \sY_1\right)$ is the $q^1_2$-lift. The $*$-homomorphism $A_0\left(q^1_2 \right)$ is injective so $A_0$ maps arrows of $\mathfrak{FinCov}\text{-}\sY\text{-}A$ to arrows of $\mathfrak{FinCov}\text{-}A\text{-}\sY$. If both $q':\widetilde \sY_1 \to \widetilde \sY_2$ and $q'':\widetilde \sY_1 \to \widetilde \sY_2$ are two different $\mathfrak{FinCov}\text{-}\sY\text{-}A$-arrows then from the Theorem \ref{top_covp_cat_thm} it follows that there is a nontrivial $g\in G\left( \left.{\widetilde{\sY}}_1~\right|{\widetilde{\sY}}_2\right)$ such that $q'' = g\circ q'$. It turns out that $A_0\left(q'' \right)= gA_0\left(q' \right)\neq A_0\left(q'' \right)$ where the given by the Lemma \ref{blowing_proper_group_lem} nontrivial action $G\left( \left.{\widetilde{\sY}}_1~\right|{\widetilde{\sY}}_2\right)\times  A_0\left(\widetilde \sY_1\right)\to  A_0\left(\widetilde \sY_1\right)$ is implied. So the functor $A_0$ is full.
\end{proof}

\begin{definition}\label{blowing_functor_defn}
	The given by the Theorem \ref{blowing_functor_thm} functor $A_0$ is said to by $\sY$-$A$-\textit{functor}.
\end{definition}
\begin{remark}\label{blowing_functor_rem}
From (i) of the Theorem \ref{blowing_functor_thm} it follows that   $\sY$-$A$-functor can be regarded as a functor
$$
	\mathfrak{FinCov}\text{-}\sY\text{-}A\xrightarrow{A_0}\mathfrak{FinCov}\text{-}A
$$
 from $A$-{regular finite covering category} of $\sY$ to the category of finite-fold coverings of $A$ (cf. Definition \ref{fin_category_defn})
\end{remark}
\begin{empt}\label{blowing_reg_cov_cat_empt}
	If $A$ is a $C^*$-algebra with  Hausdorff blowing-up and $C_0\left(\sY \right)\hookto M\left(A \right)$. If  the transitive finite-fold covering  covering $q:\widetilde \sY\to \sY$ is $A$-regular (cf. Definition \ref{blowing_a_regular_defn}), and $H\subset G\left(\left.\widetilde \sY\right|\sY \right)$ is a normal subgroup  then from the Lemma  \ref{blowing_lift_composition_lem} it follows that the natural transitive covering $\widetilde \sY/ H \to \sY$ is $A$-regular. Moreover from the Lemma  \ref{blowing_lift_composition_lem} it turns out that that the covering $\widetilde \sY\to \widetilde \sY/ H$ is $A_0\left(\widetilde \sY/ H \right)$ regular.
	\end{empt}

\subsection{Infinite coverings}
\begin{theorem}\label{blowing_sufficient_covering_inf_thm} 
	Let $C_0\left( \sY\right) \hookto M\left( A\right) $ be    Hausdorff blowing-up (cf. the Definition \ref{blowing_defn}), and
	let $q:\widetilde \sY \to \sY$ be an $A$-regular (cf. Definition \ref{blowing_a_regular_defn})  transitive covering with connected $\widetilde \sY$ and residually finite group  $G\left(\left.\widetilde \sY\right|\sY \right)$ of covering transformations.  If  $\overline q: \overline \sY  \to  \sY$  is {disconnected covering of} $q : \widetilde \sY \to \sY$ (cf. Definition \ref{top_disconnected_defn}) and if $\widehat{G}$ is the profinite completion of  $G\left(\left.\widetilde \sY\right|\sY \right)$ (cf. Example \ref{profinite_exm}) then 	following conditions hold:
	\begin{enumerate}
		\item[(i)] the covering  $\overline q$ is $A$-regular (cf. Definition \ref{blowing_a_regular_defn}, i.e. there is a $\overline q$-{lift} 	$A_b\left(\overline q\right): A \hookto	M\left( A_0\left( \overline \sY\right)\right)$ of $A$ (cf. Definition \ref{blowing_lift_defn}  ,
		\item[(ii)] if the space if $\sY$ is  locally connected  (cf. Definition \ref{top_locally_connected_defn}) then the  {finite covering category of} $q : \widetilde \sY \to \sY$ (cf. Definition \ref{top_disconnected_defn}) given by
		\be\label{blowing_category_fin_eqn}
	\mathfrak{S}_p \bydef \left\{\left\{\sY_\la\right\}_{\la \in \La}, \left\{q^\mu_\nu:\sY_\mu\to \sY_\nu\right\}_{\substack{\mu,\nu \in \La\\\mu\ge\nu}}\right\}
	\ee
 yields an algebraical finite covering category 
\be\label{blowing_algebraical_finite_covering_category_eqn}
\mathfrak{S}_{A_0\left(q \right) }\bydef \left\{\left\{A_0\left( \sY_\la\right) \right\}_{\la\in \La}, \left\{A_0\left( q^\nu_\mu\right)  : A_0\left( \sY_\mu\right) \hookto A_0\left( \sY_\nu\right)\right\}_{\substack{\mu, \nu \in \La\\\mu \le \nu}}\right\}
\ee
(cf. Definition \ref{algebraical_finite_covering_category_defn}) such that 	a triple $\left(A, A_0\left(\overline\sY \right), \widehat{G} \right)$  is a {pre}-{covering of algebraical finite covering category}  $\mathfrak{S}_{A_0\left(q \right) }$ (cf. Definition \ref{algebraical_finite_covering_category_defn}), 
\item[(iii)] if the triple $\left(A, A_0\left(\overline\sY \right), \widehat{G} \right)$   is the  \textit{disconnected infinite noncommutative covering} of  $\mathfrak{S}_{A_0\left(q \right) }$ (cf. Definition \ref{disconnected_infinite_noncommutative_covering_defn}) then it is good (cf. Definition \ref{good_defn}) and a triple $\left(A, A_0\left(\widetilde\sY \right), G\left(\left.\widetilde \sY\right|\sY \right) \right)$ the  {infinite noncommutative covering} of $\mathfrak{S}_{A_0\left(q \right) }$ (cf. Definition \ref{infinite_noncommutative_covering_defn}).
	\end{enumerate}
\end{theorem}

\begin{proof}
	According to the construction of the Section \ref{top_inf_to_sec} one has $\overline G= \varprojlim_{\la \in \La} G_\la$ where $\left\{G_\la\right\}_{\la \in \La}$ is a set of all finite factor-groups of $G\left(\left.\widetilde \sY\right|\sY \right)$, i.e. the group $G\left(\left.\widetilde \sY\right|\sY \right)$ is profinite (cf. Example \ref{profinite_exm}). From the Lemma \ref{top_disconnected_lem} it follows that if  $\left\{g_\iota G\left(\left. \widetilde\sY~\right|\sY \right)\right\}_{\iota \in I}$   is  a set of all left  cosets of $G\left(\left. \widetilde\sY~\right|\sY \right)$ in $\widehat{G}$ (cf. Definition \ref{group_coset_defn}) then there is a natural homeomorphism 
\be\label{blowing_disconnected_eqn}
\overline \sY \cong \bigsqcup_{\iota\in I} g_\iota \widetilde\sY.
\ee

	(i)
	If $\overline A$ is a $C^*$-norm completion of an algebraic direct sum 
	$$
	\bigoplus _{\iota\in I} A_0\left( g_\iota \widetilde\sY\right) 
	$$
	then there is  natural Hausdorff blowing-up $C_0\left(\overline \sY \right) \hookto \overline A$.
	An action $\overline G\times \overline \sY\to \overline \sY$ yields an action $G\left(\left.\overline \sY\right|\sY \right)\times \overline A\to \overline A$. We leave to the reader that the natural $*$-homomorphism $A \hookto M\left( \overline A\right)$ satisfies to the Definition \ref{blowing_lift_defn}  .\\
(ii)	
Firstly we prove that the triple  $\left( A ,  A_0\left(\overline\sY \right), \widehat G\right)$ is  an  {infinite quasi-covering} (cf. Definition \ref{infinite_quasicovering_defn}). If $\overline a \in K\left(A_0\left( \overline \sY\right)  \right)$ then from the  Lemmas \ref{blowing_blowing_lem} and \ref{blowing_pedersen_compact_lem} it follows that $\overline a \in A_c\left( \overline \sY\right)$ where $A_c\left(\widetilde\sY \right)$  is given by the equation \eqref{blowing_com_eqn}. If $\overline q_\la: \overline \sY \to \sY_\la$ is the natural covering then from the Lemma \ref{blowing_lift_composition_lem} it follows that $\overline q_\la$ is $A_0\left(\sY_\la \right)$-regular.
From the Lemma \ref{blowing_descent_lem} it follows that 
\be\label{blowing_al_en}
a_\la \bydef \bt\text{-}\sum_{	g \in \ker\left( \widehat{G}\to G\left(\widetilde \sY, \sY_\la \right)\right)  }g\overline a = A_b\left(\overline q_\la \right)  \circ \desc^c_{\overline q_\la}\left(\overline a\right) \in A_b\left(\overline q_\la \right) \left(A_\la \right) 
\ee
where $\bt\text{-}\sum$ implies the strict topology of $M\left(A_0\left(\overline \sY \right) \right)$ (cf. Definition \ref{strict_topology_defn}), $~\desc^c_{\overline q_\la}$ is the  {compactly supported $\overline q_\la$-descent} (cf. Definition \ref{blowing_descent_compactly_supported_defn}) .
and $ A_b\left(\overline q_\la \right): A_\la \hookto M\left(A_0\left( \overline \sY\right)  \right)$ is the $\overline q_\la$-lift (cf. Definition \ref{blowing_lift_hom_defn}). From the Lemma \ref{blowing_descent_pedersen_lem} it follows that 
if  $\overline A'\subset A_b\left(\overline q_\la \right) \left(A_\la \right)$ is a set of given by \eqref{blowing_al_en} elements then
$$
A_b\left(\overline q_\la \right)\left(K\left( A\left( \sY_\la\right)  \right) \right)\subset \overline A'
$$
Since $K\left( A\left( \sY_\la\right)\right) $ is dense in $ A\left( \sY_\la\right)$ the closure of $\overline A'$ equals to  $A_b\left(\overline q_\la \right)\left(  A\left( \sY_\la\right)\right)$, so the conditions of the Definition \ref{infinite_quasicovering_defn} are satisfied.

Now one should prove conditions (a) and (b) of the Definition \ref{algebraical_finite_covering_category_defn}.
\begin{enumerate}
	\item [(a)]  From the Theorem \ref{blowing_sufficient_covering_thm} it follows that $A_0\left( q^\nu_\mu\right): A_0\left( \sY_\mu\right) \hookto A_0\left( \sY_\nu\right)$ is a noncommutative finite-fold covering.
	\item[(b)] Follows from the Lemma \ref{blowing_descent_pedersen_lem}.
\end{enumerate}

\end{proof}

Under the hypotheses of the Theorem	\ref{blowing_sufficient_covering_inf_thm} set $\widehat A \bydef C^*$-$\varinjlim_{\la \in \La} A\left( \sY_\la\right) $, where $C^*$-$\varinjlim$ is the given by the Definition \ref{inductive_lim_non_defn} $C^*$-inductive limit. If $\widehat A\hookto B\left(\widehat \H \right)$ is a faithful, nondegenerate representation, and a triple  $\left(A, \overline A', \widehat{G} \right)$ is a {pre}-{covering of the algebraical finite covering category} (cf. Definition \ref{algebraical_finite_covering_category_defn}  $\mathfrak{S}_{A_0\left(q \right) }$ given by \eqref{groupoid_foli_algebraical_finite_covering_category_eqn} then 
from Lemma \ref{infinite_faithful_nondegererate_lem}  it follows that   there is a faithful, nondegenerate representation
\be\label{groupoid_folia_eqn}
\pi_{\overline A'} : \overline A' \hookto B\left(\widehat  \H \right). 
\ee
In particular one has
\be\label{groupoid_folia_c_eqn}
\pi_{A_0\left(\overline \sY \right) } : A_0\left(\overline \sY \right) \hookto B\left(\widehat  \H \right) 
\ee
(iii)
One should the conditions (a) and (b) of the  Definition \ref{good_defn} are satisfied.
\begin{enumerate}
	\item [(a)] From the equation \eqref{blowing_disconnected_eqn} and the Corollary \ref{blowing_connected_comp_cor} it follows that any connected component of $A_0\left(\overline\sY \right)$ equals to the given by 
	$$
	_{g_\iota \widetilde{\sY}} A_0\left(\overline\sY \right)_{g_\iota \widetilde{\sY}}\cong  A_0\left(g_\iota \widetilde{\sY} \right)= g_\iota  A_0\left( \widetilde{\sY} \right).
	$$
	$g_\iota \widetilde{\sY}$-algebra (cf. Definition \ref{blowing_ideals_au_ua_defn}). If both $g_{\iota_1}  A_0\left( \widetilde{\sY} \right)$ and $g_{\iota_2}  A_0\left( \widetilde{\sY} \right)$ are connected components then $g_{\iota_1}  A_0\left( \widetilde{\sY} \right)= \left(g_{\iota_1}g^{-1}_{\iota_2} \right)g_{\iota_2}  A_0\left( \widetilde{\sY} \right)$,
	\item [(b)] If
	$$
	G\left(\left. A_0\left( \widetilde{\sY} \right)~\right| A\right)\bydef 
	\left\{\left. g \in \widehat{G}\right| \forall \widetilde a^\perp \in A_0\left( \widetilde{\sY} \right)^\perp \quad g \widetilde a^\perp= \widetilde a^\perp\right\}
	$$
	is a specialization of given by the equation \eqref{infinite_covering_transformation_group_eqn} group then there is a natural isomorphism $G\left(\left. A_0\left( \widetilde{\sY} \right)~\right| A\right)\cong G\left(\left. \widetilde\sY~\right|\sY \right)$. For any $\la\in \La$ a natural homomorphism $G\left(\left. A_0\left( \widetilde{\sY} \right)~\right| A\right)\to G\left(\left. A_0\left( \sY_\la\right) ~\right|\sY \right)$ is surjective, since it is equivalent to the surjective homomorphism $G\left(\left. \widetilde\sY~\right|\sY \right)\to G\left(\left. \sY_\la~\right|\sY \right)$. 
	
\end{enumerate}
From the Definition \ref{infinite_noncommutative_covering_defn} it turns out that the triple $\left(A, A_0\left(\widetilde\sY \right), G\left(\left.\widetilde \sY\right|\sY \right) \right)$ the  {infinite noncommutative covering} of $\mathfrak{S}_{A_0\left(q \right) }$ (cf. Definition \ref{infinite_noncommutative_covering_defn})

\begin{lemma}\label{blowing_necessary_covering_inf_lem} 
	Let $C_0\left( \sY\right) \hookto M\left( A\right) $ be    Hausdorff blowing-up (cf. the Definition \ref{blowing_defn}), and
	let $q:\widetilde \sY \to \sY$ be an $A$-regular (cf. Definition \ref{blowing_a_regular_defn})  transitive covering with connected $\widetilde \sY$ and residually finite group  $G\left(\left.\widetilde \sY\right|\sY \right)$ of covering transformations.  Let   $\overline q: \overline \sY  \to  \sY$  be {disconnected covering of} $q : \widetilde \sY \to \sY$ (cf. Definition \ref{top_disconnected_defn}), and let $\widehat{G}$ be the profinite completion of  $G\left(\left.\widetilde \sY\right|\sY \right)$. If 	\be\label{blowing_algebraical_finite_covering_p_category_eqn}
	\mathfrak{S}_{A_0\left(q \right) }\bydef \left\{\left\{A_0\left( \sY_\la\right) \right\}_{\la\in \La}, \left\{A_0\left( q^\nu_\mu\right)  : A_0\left( \sY_\mu\right) \hookto A_0\left( \sY_\nu\right)\right\}_{\substack{\mu, \nu \in \La\\\mu \le \nu}}\right\}
	\ee
is the given by the Theorem \ref{blowing_sufficient_covering_inf_thm} an algebraical finite covering category (cf. Definition \ref{algebraical_finite_covering_category_defn}) and $\left(A, \overline A', \widehat{G} \right)$ is {pre}-{covering of algebraical finite covering category}	
$\mathfrak{S}_{A_0\left(q \right) }$ (cf. Definition \ref{algebraical_finite_covering_category_defn} then one has:
	\begin{enumerate}
\item[(i)] if $A_0\left(\overline\sY \right)$ is a $C^*$-subalgebra of $\overline A'$ then $A_0\left(\overline\sY \right)$ is  a hereditary $C^*$-subalgebra (cf. Definition \ref{hered_defn}) of $\overline A'$.
\item[(ii)] if $\left(\overline \sU, \overline \sV, \overline s\right)$ is a \textit{covering triple for} $\overline q$ (cf. Definition \ref{top_coveing_triple_defn}) and $$\overline C \subset K\left( A_0\left(\overline \sY \right)\right)  \cap ~ _{\overline \sU} A_0\left(\overline \sY \right) _{\overline \sU}$$ is a commutative $C^*$-subalgebra of  $A_0\left(\overline \sY \right)$ then there is the natural inclusion $\pi_{\overline C}: \overline C \hookto \overline A'$ such that
\be\label{blowing_necessary_covering_inf_eqn}
\pi_{A_0\left(\overline \sY \right) }\left( {\overline C}\right) \subset \pi_{\overline A'}\left(  \overline A'\right) 
\ee
where both $\pi_{\overline A'}$ and $\pi_{A_0\left(\overline \sY \right) }$ are given by the equations \eqref{groupoid_folia_eqn} and \eqref{groupoid_folia_c_eqn} respectively.
	\end{enumerate}
\end{lemma}
\begin{proof}
(i) From the Lemma \ref{infinite_faithful_nondegererate_lem} it follows that there  is an equivariant representation $\widehat \pi : \widehat A  \bydef C^*$-$\varinjlim_{\la \in \La} A_0\left( \sY_\la\right)  \hookto B\left(\widehat{\H}\right)$ (cf. Definition \ref{equivariant_representation_defn}) is faithful \ref{faithful_representation_defn} then there are natural faithful representations $\overline \pi: \overline A' \hookto  B\left(\widehat{\H}\right)$ and $\overline \pi_0:  A_0\left( \overline \sY\right)  \hookto  B\left(\widehat{\H}\right)$.  Under the hypotheses (iv) one has 
$$
\overline \pi_0\left(   A_0\left( \overline \sY\right)\right) \subset \overline \pi\left(  \overline A \hookto  B\left(\widehat{\H}\right)\right). 	
$$
Since Pedersen's ideal $K\left( A_0\left( \overline \sY\right) \right)$  of $C^*$-algebra $A_0\left( \overline \sY\right)$ is dense in $ A_0\left( \overline \sY\right)$ one should prove
$$
\forall \overline a \in K\left( A_0\left( \overline \sY\right) \right)\quad \forall \overline b \in  \overline A'\quad \overline b \le \overline a \quad \Rightarrow \quad \overline b \in K\left( A_0\left( \overline \sY\right) \right) 	
$$
If $\overline a \in K\left( A_0\left( \overline \sY\right) \right)$ then from the Lemma \ref{blowing_pedersen_compact_lem} it follows that the support $\supp \overline a$ of $\overline a$ (cf. Definition \ref{blowing_support_defn}) if compact.
From the Theorem \ref{top_compact_thm} it follows that there is an open set $\overline \sU$ and $\la_{\overline \sU} \in \La$ such that:
\begin{itemize}
	\item $\supp \overline a \subset  \overline \sU$,
	\item If $\overline q_{\la_{\overline \sU}} : \overline \sY \to \sY_\la$ then $\overline \sU$ is mapped homeomorphically onto $\sU_\la \bydef \overline q_{\la_{\overline \sU}}\left(  \overline \sU\right)$. 
\end{itemize}
If $\overline b \in  \overline A'$ is such that $\overline b \le \overline a$ and 
$$
\overline a_{\la_{\overline \sU}}	\bydef \sum_{	g \in \ker\left( \widehat{G}\to G_{\la_{\overline \sU}}\right) }g \overline a \in A_0\left(\sY_{\la_{\overline{   \mathcal U }}} \right) 
\overline b_{\la_{\overline \sU}}	\bydef \sum_{	g \in \ker\left( \widehat{G}\to G_{\la_{\overline \sU}}\right) }g \overline b \in A_0\left(\sY_{\la_{\overline{   \mathcal U }}} \right) 
$$
then  $ b_{\la_{\overline \sU}} \le  a_{\la_{\overline \sU}}$
From $\overline b \le \overline a$ and $\overline a =\lift^{\overline q_{\la_{\overline \sU}}}_{\overline \sU}\left( 	\overline a_{\la_{\overline \sU}}\right) $ it follows that  $\overline b =\lift^{\overline q_{\la_{\overline \sU}}}_{\overline \sU}\left( 	\overline b_{\la_{\overline \sU}}\right)\in A_0\left( \overline \sY\right)$. Taking into account $\overline a \in K\left( A_0\left( \overline \sY\right) \right)$ and $\overline b \le \overline a$ one has $\overline b \in K\left( A_0\left( \overline \sY\right) \right)$\\
	(ii) If $\desc_{\overline q} : K\left( A_0\left( \overline \sY\right) \right)\to A$ is  $\overline q$-descent (cf. Definition \ref{blowing_descent_defn}) then
\bean
C \bydef \desc_{\overline q}\left(\overline C  \right) = \bigcup \left\{\left\{\desc_{\overline q}\left(\overline c  \right)\right\}\left|\overline c \in \overline  C\right.\right\}\subset K\left(A \right) \cap ~_\sU A _\sU.
\eean
is a $C^*$-algebra  isomorphic to $\overline C$. Similarly for all $\la \in \La$ there is a $C^*$-algebra  
\bean
C_\la \bydef \desc_{\overline q_\la}\left(\overline C  \right) = \bigcup \left\{\left\{\desc_{\overline q_\la}\left(\overline c  \right)\right\}\left|\overline c \in \overline  C\right.\right\}\subset K\left(A \left(\sY_\la \right) \right) \cap ~_{\sU_\la} A \left(\sY_\la \right) _{\sU_\la}.
\eean
where $\overline q_\la: \overline \sY \to \sY_\la$ is the natural covering and $\sU_\la = \overline p_\la\left(\overline \sU \right)$. 	
If $\sX, ~\sX_\la$ and $\overline \sX$ are spectra of $C,~C_\la$ and $\overline C$ respectively then there are natural homeomorphisms $\sX_\la \cong \sX\cong \overline{   \mathcal X}$. Moreover if 
\bean
C^\oplus_\la \bydef \bigoplus_{g \in G\left(\left.\sY_\la \right|\sY \right) } g C_\la
\eean 
and $\sX^\oplus_\la$ is the spectrum of $C^\oplus_\la$ then
$$
\sX^\oplus_\la = \bigsqcup_{g \in G\left(\left.\sY_\la \right|\sY \right) } g \sX_\la
$$
where the action $G\left(\left.\sY_\la \right|\sY \right)\times \sX^\oplus_\la \to \sX^\oplus_\la$ comes from the natural action 	$G\left(\left.\sY_\la \right|\sY \right)\times C^\oplus_\la \to C^\oplus_\la$. If
$$
\overline \sX^\oplus  \bydef \bigsqcup_{g \in \widehat G } g\sX
$$
then for any $\la\in \La$ the surjective homomorphism $\widehat G \to G\left(\left.\sY_\la \right|\sY \right)$ yields the natural covering
$$
\overline p_\la : \overline \sX^\oplus \to \sX_\la
$$
If  $\overline c \in \overline C_+$  is positive element then from the condition (b) of the Definition \ref{algebraical_finite_covering_category_defn} it follows that there is $c'_{\overline A'} \in K\left( \overline A\right)_+$ such that
\be\label{blowing_ova_eqn}
c=  \b\text{-}\sum_{	g \in \widehat {G}}gc'_{\overline A'}.
\ee
where $\b\text{-}\sum$ means the convergence with respect to the strict topology of $M\left(\overline A \right)$. 
Moreover for all $\la \in \La$ there is a positive $c_\la \in C^\oplus_\la$ such that
$$
c_\la =  \b\text{-}\sum_{	g \in \ker\left(  \widehat {G}\to G\left(\left.\sY_\la \right|\sY \right)\right) }gc'_{\overline A'} . 
$$
On the other hand for any $\la \in \La$ there is 
$$
\overline c_\la \bydef \lift_{ \widetilde{p}_{\la}}\left( c_\la \right) \in C_b\left( \overline \sX^\oplus\right) 
$$
where $\lift_{ \widetilde{p}_{\la}}$ is given by the Definition \ref{top_lift_defn} $ \widetilde{p}_{\la}$-lift. The  net $\left\{\overline c_\la\right\}_{\la\in\La}$ is decreasing, it follows that there is a point-wise limit
$$
\overline c' \bydef s\text{-}\lim_{{\la \in \La}}\overline c_\la \in C\left( \overline \sX^\oplus_{\text{discr}}\right) 
$$
where  $C\left( \overline \sX^\oplus_{\text{discr}}\right)$ is   the {discontinuous extension of} $C_0\left(  \overline \sX^\oplus\right)$. From the Theorem \ref{top_discrete_p_thm} it follows that $\overline c \in C_c\left( \overline \sX^\oplus\right)$. Since the support  $\supp \overline c'$ is compact the set
$$
\widehat G_0 \bydef \left\{\left. \widehat g \in \widehat{G}\right| \widehat g \overline \sX  \cap\supp \overline c'\neq \emptyset\right\} 
$$
is finite, i.e. $\widehat G_0= \left\{\widehat g_1,..., \widehat g_n\right\}$ Since the group $\widehat{G}$ is residually finite there is $\la_0\in \La$ such that
$$
\forall g', g'' \in \widehat G_0\quad g' \neq g'' \quad \Rightarrow \quad \phi_{\la_0}\left( g'\right)\neq \phi_{\la_0}\left( g''\right). 
$$
where $\phi_{\la_0}:  \widehat {G}\to G\left(\left.\sY_\la \right|\sY \right)$ is the natural surjective homomorphism of groups. Indeed $\overline c'$ corresponds to $\overline c'_{\overline A'}$.
If $s_{\la_0}\bydef \desc_{\overline q_{\la_0}}\left(\overline s \right) $ where $ \desc_{\overline q_{\la_0}}$ is the $\overline q_{\la_0}$ descent (cf. Definition \ref{blowing_descent_defn})  then from $C_0\left( \sY\right)\subset M\left( A\left(\sY_\la \right) \right)$ and the Corollary \ref{infinite_multiplier_cor} it follows that $s_{\la_0}\in M\left(A_0\left( \sY_\la\right)\right).$ From the Corollary   \ref{infinite_multiplier_cor}  if follows that $g  s_{\la_0}\in M\left( \overline A'\right)$ for all $g \in G\left(\left.\sY_\la \right|\sY \right)$.  In result one has
$$
\overline c'_{\overline A}= \overline c'_1 + ...+ \overline c'_n
$$
where
$$
\forall j = 1,..., n \quad \overline c'_j \bydef \left( \phi_{\la_0}\left( g's_{\la_0}\right)\right) c'_{\overline A'}\in K\left( \overline A'\right).
$$
Set  $ \overline c_{\overline A'} \bydef \widehat g_1^{-1}\overline c'_1+...+ \widehat g_n^{-1}\overline c'_n \in K\left( \overline A'\right)_+$, and define a map
\bean
\pi_{\overline C}: \overline C \to \overline A',\\
\overline c \mapsto \overline c_{\overline A'}.
\eean 
We leave to the reader a proof of that $\pi_{\overline C}$ satisfies to the equation \eqref{blowing_necessary_covering_inf_eqn}. From this fact it turns out that for any $\overline{c}\in \overline C$ the operator $\pi_{\overline C}\left( \overline{c}\right)$ does not depend on choice of  $c'_{\overline A}$ satisfying to the equation \eqref{blowing_ova_eqn}.
\end{proof}
\begin{lemma}\label{blowing_universal_fg_lem} 
	Let $C\left(\sY \right) \hookto M\left(A \right)$ be a Hausdorff blowing-up (cf. Definition \ref{blowing_defn}) with compact, locally connected space $\sY$, and let 
	$q: \widetilde\sY \to \sY$  be an  $A$-regular covering (cf. Definition \ref{blowing_a_regular_defn}), with connected $\widetilde \sY$. Suppose that the $q$-lift $A_0\left( \widetilde\sY\right)$ of $A$ (cf. Definition \ref{blowing_lift_defn}) is simply connected (cf. Definition \ref{simply_connected_defn}). If $\left( A, A_0\left(\widetilde{\sY},  \right), G\left(\left.\widetilde{\sY}~\right|{\sY} \right) \right)$ is the infinite noncommutative covering \ref{infinite_noncommutative_covering_defn} of $\mathfrak{S}_{A_0\left(q \right) }$ (cf. \eqref{blowing_category_fin_eqn}) and the property $P_{\mathrm{untz}}$ of finite-fold noncommutative coverings is given by \eqref{unitization_p_eqn} then one has:
	\begin{enumerate}
		\item[(i)] the $P_{\mathrm{untz}}$-{universal  covering} 
		of $A$ (cf. Definition \ref{fundamental_group_nc_p_defn}) equals to
		$$
		\left( A, A_0\left(\widetilde{\sY},  \right), G\left(\left.\widetilde{\sY}~\right|{\sY} \right) \right), 
		$$
		\item[(ii)] if $\pi^{P_{\mathrm{untz}}}_1\left(A \right)$ is the $P_{\mathrm{untz}}$-{fundamental group} of $A$ (cf. Definition \ref{fundamental_group_nc_p_defn}) then there is a natural isomorphism $\pi^{P_{\mathrm{untz}}}_1\left(A \right)\cong  G\left(\left.\widetilde{\sY}~\right|{\sY} \right)$. 
	\end{enumerate}

\end{lemma}

\begin{proof}
From the Theorem \ref{blowing_sufficient_fin_thm} it follows that the given by the Theorem \ref{blowing_sufficient_covering_inf_thm} algebraical finite covering category (cf. Definition \ref{algebraical_finite_covering_category_defn})	
\bean
\mathfrak{S}_{A_0\left(q \right) }\bydef \left\{\left\{A_0\left( \sY_\la\right) \right\}_{\la\in \La}, \left\{A_0\left( q^\nu_\mu\right)  : A_0\left( \sY_\mu\right) \hookto A_0\left( \sY_\nu\right)\right\}_{\substack{\mu, \nu \in \La\\\mu \le \nu}}\right\}
\eean
contains {all} classes of isomorphisms of noncommutative finite-fold coverings of $A$ (cf. Definition \ref{fin_defn}) which possess the property $P_{\mathrm{untz}}$. From the Definition \ref{fundamental_group_nc_p_defn} it follows that $\pi^{P_{\mathrm{untz}}}_1\left(A \right)\cong  G\left(\left.\widetilde{\sY}~\right|{\sY} \right)$
\end{proof}
\chapter{Coverings of stable $C^*$-algebras}\label{stab_chap}

\section{Finite-fold coverings}\label{stab_fin_sec}
\paragraph*{}
Here we find a relation between  noncommutative finite-fold coverings of $C^*$-algebras and their stabilizations.
\begin{empt}\label{stab_uni_empt}
	Let $\H$ be a Hilbert space with finite or countable basis. The algebra $ \K \bydef\K\left(\H \right)$ of compact operators is isomorphic to $\mathbb{M}_n\left(\C \right)$ if  dimension of $\H$ equals to $n$ and $\K\left(\ell^2\left( \N\right)  \right)$ is $\H$ is infinite dimensional.  Denote by $\K^+$ the minimal unital $C^*$-algebra which contains $\K$. If  $\K \cong  \mathbb{M}_n\left(\C \right)$ then $\K^{+} \cong  \mathbb{M}_n\left(\C \right)$, otherwise $\K^+ = \K\left(\ell^2\left( \N\right)\right)  \bigoplus \C$ with product given by
	$$
	\left( a_1 \oplus c_1\right) \cdot  \left( a_2 \oplus c_2\right) = 
	\left(a_1a_2 + c_2a_1 + c_1a_2 \right) \oplus c_1c_2.
	$$ 
	Here we find the relation between noncommutative coverings of $C^*$-algebras and noncommutative coverings its tensor products with algebras of compact operators.
	For any $C^*$-algebra the following notation is used
	
	\break
	\begin {table}[H]
	\caption {Notations for products with compact operator} \label{stab_table} 
	\begin{center}
		\begin{tabular}{|c|c|c|c|}
			\hline
			& Notation	& Dimension $n$	 &Infinite dimension\\
			\hline
			&	&	&\\
			Hilbert space & $\H_\K$ &$\C^n$ 	&  $\ell^2\left( \N\right)$ \\
			&	&	&\\			Compact operators& $\K$ &$\mathbb{M}_n\left(\C \right)$ 	&  $\K\left(\ell^2\left( \N\right)\right)$ \\
			&	&	&\\
			Unitization &&&\\
			of compact operators& $\K^+$ &$\mathbb{M}_n\left(\C \right)$ 	&  $\K\left(\ell^2\left( \N\right)\right)\oplus \C$ \\
			&	&	&\\
			Irreducible &&&\\ representation& $\pi_\K: \K \to B\left(\H_\K \right) $ &$\H_\K=\C^n$ 	&  $\H_\K=\ell^2\left( \N\right)$ \\
			&	&	&\\
			Hilbert $C^*$-module& $X_A$ &$A^n$ 	&  $\ell^2\left(A\right) $\\
			&	&	&\\
			$A$-compact operators&$\K_A= \K\left( X_A\right) $ &$\mathbb{M}_n\left(A \right)$	&  $\K\left( \ell^2\left(A\right)\right)  $ \\ 		
			&	&	&\\
			Unitization of $\K_A$		&$\K^+_A= \K\left( X_A\right)^+ $ &$\mathbb{M}_n\left(A \right)$	&  $\K\left( \ell^2\left(A\right)\right)^+  $ \\ 		
			&	&	&\\
			\hline
		\end{tabular}
	\end{center}
	\end {table}
	If $A \otimes \K$ is the minimal or, equivalently. maximal tensor product then one has
	\be\label{stable_comp_op_eqn}
	A \otimes \K \cong \K_A.
	\ee
	Any  involutive action of $G$ on $A$ uniquely induces a involutive action  of $G$ on $A\otimes \mathcal K$.
\end{empt}
	\begin{empt}\label{stable_ap_empt}
	Note that for any $C^*$-algebra $A$ there is a natural inclusion $\K \subset M\left(A \otimes \K \right)$. If $p \in \K$ is a rank-one projector and
	\be\label{stable_ap_eqn}
	A^p \bydef \left\{ a \otimes p | a \in A\right\}\subset A \otimes  \K	
	\ee
	then $A^p$ is a hereditary $C^*$-subalgebra of $A \otimes  \K$ with the natural $*$-isomorphism
	\be\label{stable_app_eqn}
	A^p\cong A.
	\ee
\end{empt}
\begin{lemma}\label{stable_pre_lem}
	Following conditions hold:
	\begin{enumerate}
		\item [(i)] If the  quadruple $\left( A, \widetilde{A}, G, \pi\right)$ is a {noncommutative finite-fold  pre-covering} then  
		$$
		\left( A \otimes \K, \widetilde{A}  \otimes \K, G, \pi  \otimes \Id_\K \right)
		$$
		is a {noncommutative finite-fold  pre-covering}.
		\item[(ii)] If a normal subgroup $H\subset G$ is   $\left(A, \widetilde{A}, G, \pi \right)$-{proper} (cf.Definition \ref{proper_subgroup_fin_eqn}) then $H$ is $
		\left( A \otimes \K, \widetilde{A}  \otimes \K, G, \pi  \otimes \Id_\K \right)
		$-{proper}.
	\end{enumerate}
\end{lemma}
\begin{proof}
	
	(i) $G_{\K}\stackrel{\mathrm{def}}{=}\left\{ g_{\K} \in \Aut\left(\widetilde{A}  \otimes \K\right)~|~ g_{\K}a = a;~~\forall a \in A  \otimes \K\right\}$. If  $\mathfrak{e}_{1,1}$ is an elementary matrix (cf. Definition \ref{elementary_defn}) $1_{A^{\sim} } \otimes \mathfrak{e}_{1,1}$ is a multiplier of both $A \otimes \K$ and $\widetilde{A} \otimes \K$. Any *-automorphism of $\widetilde{A}\otimes \K$ can be uniquely extended to *-automorphism of $M\left( \widetilde{A}\otimes \K\right)$. If $\widetilde{a}_\K \in \widetilde{A}  \otimes  \mathfrak{e}_{1,1}$ then  $$\widetilde{a}_\K = \left( 1_{A^{\sim} } \otimes \mathfrak{e}_{1,1}\right) \widetilde{a}_\K \left( 1_{A^{\sim} } \otimes \mathfrak{e}_{1,1}\right).$$
	Hence for any $g_\K \in G_\K$ one has 
	\be\label{stab_g_act_eqn}
	\begin{split}
		g_\K\left(  \left( 1_{A^\sim} \otimes \mathfrak{e}_{1,1}\right) \widetilde{a}_\K \left( 1_{A^\sim } \otimes \mathfrak{e}_{1,1}\right)\right) = \\=\left(	g_\K  \left( 1_{A^\sim } \otimes \mathfrak{e}_{1,1}\right)\right) 	g_\K \widetilde{a}_\K \left( 	g_\K\left( 1_{A^\sim } \otimes \mathfrak{e}_{1,1}\right)\right)=\\= \left( 1_{A^\sim } \otimes \mathfrak{e}_{1,1}\right) g_\K\widetilde{a}_\K \left( 1_{A^\sim } \otimes \mathfrak{e}_{1,1}\right)
	\end{split}
	\ee
	and the above equation means
	$$
	g_\K \widetilde{A}  \otimes  \mathfrak{e}_{1,1} = \widetilde{A}  \otimes  \mathfrak{e}_{1,1}
	$$
	However there is $*$-isomorphism $\widetilde{A}  \otimes  \mathfrak{e}_{1,1}$ so $g_\K$ uniquely defines the element of the group $G\stackrel{\mathrm{def}}{=}\left\{ g \in \Aut\left(\widetilde{A}  \right)~|~ ga = a;~~\forall a \in A \right\}$. Otherwise the action of $g_\K$ on $\widetilde{A}\otimes \K$ uniquely depends on the action of $g_\K$ on $\widetilde{A}\otimes \mathfrak{e}_{1,1}$. Really if $\widetilde{b}_\K \in \widetilde{A}\otimes \mathfrak{e}_{j,k}$ then from 
	$$
	\widetilde{b}_\K = \left( 1_{A^\sim } \otimes \mathfrak{e}_{1,j}\right)\widetilde{b}_\K \left( 1_{A^\sim } \otimes \mathfrak{e}_{1,k}\right) 
	$$
	it turns out 
	\bean
	\widetilde{b}_\K = \left( 1_{A^\sim }  \otimes \mathfrak{e}_{1,j}\right)  \widetilde{c}_\K \left( 1_{A^\sim } \otimes \mathfrak{e}_{k,1}\right)\\ \text{ where }  \widetilde{c}_\K = \left(1_{A^\sim }  \otimes  \mathfrak{e}_{1,j}\right)\widetilde{b}_\K \left( 1_{A^\sim } \otimes \mathfrak{e}_{1,k}\right) \in \widetilde{A}  \otimes  \mathfrak{e}_{1,1}
	\eean
	and
	$$
	g_\K \widetilde{b}_\K = \left( 1_{A^\sim }  \otimes \mathfrak{e}_{1,j}\right)  g_\K\widetilde{c}_\K \left( 1_{A^\sim } \otimes \mathfrak{e}_{k,1}\right)
	$$
	Conversely any element of $g \in G$ uniquely defines the element $g_\K \in G_\K$. Moreover from $A = \widetilde{A}^G$ it turns out ${A}\otimes \K = \left( \widetilde{A}\otimes \K\right)^{G_\K}$.\\
	(ii) One needs check that a quadruple 	$\left(A\otimes\K, \widetilde{A}^H\otimes\K, G\left(\left.\widetilde A~\right|A \right)/H, \pi^H\otimes \Id_\K \right)$ satisfies to the conditions of the Definition \ref{fin_pre_defn}.
	\begin{enumerate}
		\item [(a)] Follows from
	$$
\left\{ \left.g \in \Aut\left(\widetilde{A}^H\ox \K \right)~\right|\forall a \in A\ox\K\quad ga = a\right\}\cong \left\{ \left.g \in \Aut\left(\widetilde{A} \right)~\right|\forall a \in A\quad ga = a\right\}.	
	$$	
		
		\item[(b)]
		Follows from
		$$
	\left( \widetilde{A}^H\ox \K\right)^{G\left(\left.\widetilde A~\right|A \right)/H} 	= A\otimes \K.
		$$
	\end{enumerate}
	
\end{proof}
\begin{lemma}\label{stable_pre_inv_lem}
Let $\left( A, \widetilde{A}, G, \pi\right)$ be a {noncommutative finite-fold  quasi-covering} (cf. Definition \ref{fin_quasi_defn}).
	\begin{enumerate}
		\item [(i)] If a  quadruple $\left( A\otimes\K, \widetilde{A}\otimes\K, G, \pi\otimes \Id_\K\right)$ is a {noncommutative finite-fold  pre-covering} (cf. Definition \ref{fin_pre_defn}) then  
		$$
		\left( A , \widetilde{A}  , G, \pi   \right)
		$$
		is a {noncommutative finite-fold  pre-covering}.
		\item[(ii)] If a normal subgroup $H\subset G$ is   $\left( A\otimes\K, \widetilde{A}\otimes\K, G, \pi\otimes \Id_\K\right)$-{proper} (cf.Definition \ref{proper_subgroup_fin_defn}) then $H$ is $
		\left( A , \widetilde{A} , G, \pi   \right)
		$-{proper}.
	\end{enumerate}
\end{lemma}
\begin{proof}(i)
If $g$ is a $*$-automorphism of $\widetilde{A}$ such that $\forall a \in \pi\left(A \right) \quad g a = a$ then there is  $*$-automorphism of $\widetilde{A}\otimes \K$ such that
$$
g_\K \left(\widetilde a \otimes x \right) \bydef  \left(g\widetilde a  \right)\otimes x
$$
such that $g_\K \left( a \otimes x\right) = \left( a \otimes x\right)$ for all $a \otimes x \in \left( \pi\otimes \Id_\K\right) \left(A \otimes \K \right)$. If $g_\K$  does not correspond to an element $g \in G$ then the  quadruple $\left( A\otimes\K, \widetilde{A}\otimes\K, G, \pi\otimes \Id_\K\right)$ is not a {noncommutative finite-fold  pre-covering}, so one has
$$
\left\{ \left.g \in \Aut\left(\widetilde{A} \right)~\right|\forall a \in \pi \left( A\right) \quad ga = a\right\}\cong G
$$
(ii) Similarly to proof of (i) one can prove that
$$
\left\{\left.g \in \Aut\left(\widetilde{A}^H \right)~\right| \forall a \in \pi \left( A\right) \quad g a = a  \right\}\cong G/ H.
$$
\end{proof}
\begin{lemma}\label{stable_fin_cov_comp_lem}
 If a  quadruple $\left( A, \widetilde{A}, G, \pi\right)$ is a  {noncommutative finite-fold  pre-covering} (cf. Definition \ref{fin_pre_defn}) then one has:
 \begin{enumerate}
 	\item [(i)] if  $\left( A, \widetilde{A}, G, \pi\right)$ is a finite-fold covering with unitization (cf. Definition \ref{fin_unitization_defn}) then  
 	$$
 	\left( A \otimes \K, \widetilde{A}  \otimes \K, G, \pi  \otimes \Id_\K \right)
 	$$
 	is a  {noncommutative finite-fold  covering} with unitization,
 	\item[(ii)] conversely if  $\left( A \otimes \K, \widetilde{A}  \otimes \K, G, \pi  \otimes \Id_\K \right)$ is a finite-fold covering with unitization (cf. Definition \ref{fin_unitization_defn}) then  
 	$$
 	\left( A, \widetilde{A}, G, \pi\right)
 	$$
 	is a  {noncommutative finite-fold  covering} with unitization
 \end{enumerate}
 \end{lemma}
\begin{proof}
	
	(i)	From the Lemma \ref{stable_pre_lem} it turns out that the quadruple
	$$
	\left( A \otimes \K, \widetilde{A}  \otimes \K, G, \pi  \otimes \Id_\K \right)
	$$	
	is a noncommutative finite-fold  pre-covering. 	From the Definition \ref{fin_unitization_defn} it follows that following conditions hold:
	
	\begin{enumerate}
		\item[(a)] 
		
		There are unitizations $A \hookto B$, $\widetilde{A} \hookto \widetilde{B}$.
		\item[(b)] There is an 
		unital  noncommutative finite-fold covering (resp. unital  noncommutative finite-fold quasi-covering)	$\left(B ,\widetilde{B}, G,\widetilde{\pi} \right)$ such that $\pi =\widetilde{\pi}|_A$ and the action $G \times\widetilde{A} \to \widetilde{A}$ is induced by $G \times\widetilde{B} \to \widetilde{B}$.
	\end{enumerate}
	We need prove that the quadruple $\left( A \otimes \K, \widetilde{A}  \otimes \K, G, \pi  \otimes 1_{\K^+} \right)$ satisfies to conditions (a), (b) the Definition \ref{fin_unitization_defn}.\\
	(a)
	Let $\mathcal K^+$ be the  unitization given by \ref{stab_uni_empt}. Algebra $A \otimes \mathcal K$ (resp. $\widetilde{A} \otimes \mathcal K$) is an essential ideal of $B \otimes \mathcal K^+$ (resp. $\widetilde{B} \otimes \mathcal K^+$), i.e. both ${A} \otimes \mathcal K \hookto B \otimes \mathcal K^+$ and $\widetilde{A} \otimes \mathcal K \hookto \widetilde{B} \otimes \mathcal K^+$ are unitizations. \\
	(b) The quadruple 
	$\left( A \otimes \K, \widetilde{A}  \otimes \K, G, \pi  \otimes \Id_\K \right)$ is a  noncommutative finite-fold  pre-covering, hence $\left( B \otimes \K^+, \widetilde{B}  \otimes \K^+, G, \widetilde{\pi}  \otimes \Id_\K \right)$ is a  noncommutative finite-fold  pre-covering. 
	Since $\left(B ,\widetilde{B}, G \right)$ is an unital  noncommutative finite-fold covering the algebra $\widetilde{B}$ is a finitely generated $B$ module, i.e. there are $\widetilde{b}_1, \dots , \widetilde{b}_n \in \widetilde{B}$ such that any $\widetilde{b} \in\widetilde{B}$ is given by
	$$
	\widetilde{b} = \sum_{j = 1}^n \widetilde{b}_jb_j; \text{	where } b_j \in B.
	$$
	From the above equation it turns out that if
	$$
	\widetilde{b}^{\mathcal K} = \sum_{k = 1}^m \widetilde{b}_k \otimes x_k \in \widetilde{B} \otimes \mathcal K^+;~ \widetilde{b}_k \in \widetilde{B}, ~ x_k \in \mathcal K^+
	$$
	and 
	$$
	\widetilde{b}_k  = \sum_{j = 1}^n \widetilde{b}_jb_{kj}
	$$
	then
	\begin{equation*}
		\begin{split}
			\widetilde{b}^{\mathcal K} = \sum_{j = 1}^n \left( \widetilde{b}_j\otimes 1_{\mathcal K^+}\right) b_j^{\mathcal K},\\
			\text{ where } b_j^{\mathcal K} = \sum_{k = 1}^m b_{kj}  \otimes x_k \in B \otimes \mathcal K.
		\end{split}
	\end{equation*}
	It turns out that the algebraic tensor product $\widetilde{B} \otimes \mathcal K^+$ is generated by $\widetilde{b}_1\otimes 1_{\mathcal K^+},..., \widetilde{b}_n\otimes 1_{\mathcal K^+}$ as $B \otimes \K$ module. The minimal (or maximal) completion the algebraic tensor product $\widetilde{B} \otimes \mathcal K^+$ is also generated by $\widetilde{b}_1\otimes 1_{\mathcal K^+},..., \widetilde{b}_n\otimes 1_{\mathcal K^+}$ so the  minimal (or maximal) tensor product is a finitely generated $B\otimes \mathcal K^+$ module. It turns out that the quadruple
	$$
	\left({B} \otimes \mathcal K^+,\widetilde{B} \otimes \mathcal K^+, G, \widetilde{   \pi} \otimes \Id_\K\right) 
	$$
	is an unital  noncommutative finite-fold covering. \\
	(ii) Let $\left\{ \widetilde{a}^\K_1, ..., \widetilde{a}^\K_n, \widetilde{b}^\K_1, ..., \widetilde{b}^\K_n\right\} \subset   M\left(  \widetilde{A}\otimes \K\right)$ be a set satisfying to the equation \eqref{cov_mult_fin_c_eqn}. If $p \in \K$ is a rank-one projector then $A^p \cong A$ and $\widetilde A^p \cong \widetilde A$ (cf. \eqref{stable_app_eqn}).
	For any $j \in \left\{1,..., n\right\}$ Let
	\bean
	L_{a_j} \widetilde A^p \to \widetilde A^p,\\
	\widetilde a \mapsto p \widetilde{a}^\K_1 p 	\widetilde a;\\
		R_{a_j} \widetilde A^p \to \widetilde A^p,\\
	\widetilde a \mapsto p \widetilde{a}^\K_1 p 	\widetilde a;\\
		L_{b_j} \widetilde A^p \to \widetilde A^p,\\
	\widetilde a \mapsto p \widetilde{b}^\K_1 p 	\widetilde a;\\
	R_{b_j} \widetilde A^p \to \widetilde A^p,\\
	\widetilde a \mapsto p \widetilde{b}^\K_1 p 	\widetilde a;\\
	\eean
	where the $p \in M\left(\widetilde A \otimes \K\right)$ is implied.
For any $j \in \left\{1,..., n\right\}$ a pair $\left( 	L_{a_j} , 	R_{a_j} \right)$ resp $\left( 	L_{b_j} , 	R_{b_j} \right)$ is a double centralizer (cf. Definition \ref{double_centralizer_defn}). From the Remark \ref{double_centralizer_rem} it follows that 
$$
\left\{ \widetilde{a}_1, ..., \widetilde{a}_n, \widetilde{b}_1, ..., \widetilde{b}_n\right\} \subset   M\left(  \widetilde{A}^p\right)
$$
From the equation \eqref{cov_mult_fin_c_eqn} one can deduce
\bean
\forall \widetilde{a} \in M\left( \widetilde{A}^p\right) \quad \widetilde{a}=	\sum_{j=1}^{n}\widetilde{b}_j a_j \\
\text{where} \quad \forall j\in \left\{1,...,n\right\}\quad a_j\bydef  \sum_{g \in G} g\left(\widetilde{a}^*_j\widetilde{a} \right)\in M\left(\pi \right) \left(  M\left( A^p\right) \right) 
\eean 

\end{proof}

\begin{theorem}\label{stable_fin_cov_thm}
	If a  quadruple $\left( A, \widetilde{A}, G, \pi\right)$ is a  {noncommutative finite-fold  pre-covering} (cf. Definition \ref{fin_pre_defn}) then one has:
	\begin{enumerate}
		\item [(i)] If the  quadruple $\left( A, \widetilde{A}, G, \pi\right)$ is a  {noncommutative finite-fold  covering}(cf. Definition \ref{fin_defn}) then  
		$$
		\left( A \otimes \K, \widetilde{A}  \otimes \K, G, \pi  \otimes 1_{\K^+} \right)
		$$
		is a  {noncommutative finite-fold  covering}.
		\item [(ii)] If the  quadruple $\left( A \otimes \K, \widetilde{A}  \otimes \K, G, \pi  \otimes 1_{\K^+} \right)$ is a  {noncommutative finite-fold  covering} then  
$$
\left( A, \widetilde{A}, G, \pi\right)
$$
is a  {noncommutative finite-fold  covering}.
		
	\end{enumerate}

\end{theorem}

\begin{proof}
	
	(i)
	From the Lemma \ref{stable_pre_lem} it turns out that $\left(A  \otimes \mathcal K, \widetilde{A} \otimes \mathcal K, G, \pi \otimes 1_{\K^+ } \right)$ is a a noncommutative finite-fold  pre-covering.
	From the Definition \ref{fin_defn} it turns out that if there is an indexed by a directed set $\La$ family $\left\{A_\la \subset A\right\}_{\la\in\La}$ of $\left(A, \widetilde{A}, G, \pi \right)$-{strictly proper} hereditary $C^*$-subalgebras (cf. the Definition \ref{strictly_proper_defn}) such that 
	\begin{enumerate}
		\item[(a)] 
		$$
		\forall \mu, \nu \in \La \quad \mu \le \nu \quad\Rightarrow\quad A_\mu \subset A_\nu,
		$$	
		\item[(b)]
		a union 
		$
		\bigcup_{\la\in\La} A_\la
		$
		is dense in $A$.
	\end{enumerate}
	We leave to the reader a proof of that  $A_\la \otimes \K$ is a $\left(A  \otimes \mathcal K, \widetilde{A} \otimes \mathcal K, G, \pi \otimes 1_{\K^+ } \right)$-strictly proper hereditary $C^*$-subalgebra (cf. the Definition \ref{strictly_proper_defn}). On the other hand one has:
	\begin{enumerate}
		\item[(a)] 
		$$
		\forall \mu, \nu \in \La \quad \mu \le \nu \quad\Rightarrow\quad A\otimes \K_\mu \subset A\otimes \K_\nu,
		$$	
		\item[(b)]
		a union 
		$
		\bigcup_{\la\in\La} A_\la\otimes\K
		$
		is dense in $A\otimes\K$.
	\end{enumerate}	
	for all $\la\in \La$.\\
	(ii) 	From the Definition \ref{fin_defn} it turns out that if there is an indexed by a directed set $\La$ family $\left\{A^\K_\la \subset A\otimes \K \right\}_{\la\in\La}$ of $\left(A\otimes \K, \widetilde{A}\otimes \K , G, \pi\otimes \Id_\K \right)$-{strictly proper} hereditary $C^*$-subalgebras (cf. the Definition \ref{strictly_proper_defn}).
If
	$$
\La' \bydef \left\{\la \in \La \left|A^\K_\la\cap A^p \neq \{0\} \right.\right\}	
	$$
	then there is a family $\left\{A^p_\la\bydef A^\K_\la\cap A^p  \right\}_{\la\in \La'}$ of hereditary subalgebras of $A^p$. From (ii) of the Lemma \ref{stable_fin_cov_comp_lem} one can deduce that $A^p_\la$ is a $\left(A^p, \widetilde{A}^p, G, \pi^p \right)$-{strictly proper} hereditary $C^*$-subalgebra (cf. the Definition \ref{strictly_proper_defn}) for any $\la \in \La'$. The union $\bigcup_{\la\in \La} A^\K_\la$ is dense in $A\otimes \K$ it follows that the union  $\bigcup_{\la\in \La'} A^p_\la$ is dense in $A^p$.
\end{proof}
\begin{exercise}\label{stable_fin_cov_exer}
	Proof that under the hypotheses of the Theorem \ref{stable_fin_cov_thm} if  a  {noncommutative finite-fold  covering} $\left( A, \widetilde{A}, G, \pi\right)$  is reduced (cf. Definition \ref{fin_red_defn}) then a  {noncommutative finite-fold  covering}
	$
	\left( A \otimes \K, \widetilde{A}  \otimes \K, G, \pi  \otimes 1_{K^+} \right)
	$ is also reduced.
\end{exercise}

\begin{empt}
Let $A$ be a $C^*$-algebra, and let $\mathfrak{FinCov}$-$A$ be e category 
\end{empt}

\begin{corollary}\label{stable_functor_cor}
	There is a  full functor $-\otimes \K: \mathfrak{FinCov}$-$A\to\mathfrak{FinCov}$-$A\otimes \K$ between categories of finite-fold coverings (cf.  Definitions \ref{fin_category_defn}  \ref{functor_defn} and \ref{funct_full_faithfull_defn}).
\end{corollary}
\begin{proof}
	Follows from the Lemma \ref{stable_pre_lem} and the Theorem \ref{stable_fin_cov_thm}.
\end{proof}
\begin{remark}\label{stab_non_iso_rem}
	From the Example \ref{nt_stable_exm}  it follows that the functor $-\otimes \K$ is not always faithful (cf. Definition \ref{funct_full_faithfull_defn}).
\end{remark}

\section{Infinite coverings}

\paragraph*{} Here we find a relation between infinite noncommutative coverings of $C^*$-algebras and their stabilizations.

\begin{lemma}\label{stable_lem}
If $A$ is a $C^*$-algebra then the algebraic  tensor product $K\left( A\right) \otimes K\left(\K \right)$ is there Pedersen's ideal $K\left( A\otimes \K\right) $ (cf. Definition \ref{pedersen_ideal_defn}) of $A\otimes \K$.
\end{lemma}
\begin{proof}
	If $p\in \K$ is a rank-one projector and $f_\eps$ is given by \eqref{f_eps_eqn} then from 
	$$
	f_\eps\left( a\otimes p  \right) =f_\eps  \left(a \right)\otimes p 
	$$
	one can deduce that
	$$
\forall b \in K\left(A \right)\quad  b \otimes p \in K\left( A\right) \otimes K\left(\K \right)
	$$
	Using this circumstance one can proof that $K\left( A\right) \otimes K\left(\K \right)\subset K\left( A\otimes \K\right)$. We leave to the reader the proof of the following statements:
	\begin{itemize}
		\item $K\left( A\right) \otimes K\left(\K \right)$ is a two sided ideal,
		\item  $K\left( A\right) \otimes K\left(\K \right)$ is dense in $A\otimes \K$. 
	\end{itemize}
	Now from the Theorem \ref{pedersen_ideal_thm} it follows that  $K\left( A\right) \otimes K\left(\K \right)$ is there Pedersen's  of $A\otimes \K$
\end{proof}

\begin{empt}\label{stable_empt}
Let $A$ be a $C^*$-algebra, and let $\widehat G$ be a profinite group (cf. Example \ref{profinite_exm}) of $*$-automorphisms of $\overline A$. Let $\widehat G \times \left( \overline A\otimes\K \right)\to \overline  A\otimes\K$ be an action which comes from the action $\widehat G \times \overline A \to \overline A$. Let  $\left\{G_\la\right\}_{\la\in \La}$ be the indexed by an ordered set  $\La$ family of all finite factor-groups of  $\widehat G$ (cf. Definition \ref{g_category_defn}).
Similarly to  \ref{infinite_quasicovering_empt}) suppose that for any element $\overline a^\K \in K\left(\overline A\otimes \K \right)$ of the Pedersen's ideal of $\overline A\otimes K$ (cf. Definition \ref{pedersen_ideal_defn}) a series 
\be\label{stab_infinite_covering_basic_eqn}
\sum_{	g \in \widehat{G}}g \overline a^\K
\ee
is convergent with respect to the strict topology of $M\left(\overline A\otimes \K\right)$ (cf. Definition \ref{strict_topology_defn}). For any $\la \in \La$ denote by $A^\K_\la$ a generated by elements
\be\label{stable_basic_cov_cl_eqn}
a^\K_\la =\bt\text{-} \sum_{	g \in \ker\left( \widehat{G}\to G_\la\right) }g \overline a^\K,
\ee
where  $\bt\text{-} \sum$ means a convergence  with respect to the strict topology of $M\left(\overline A\otimes \K\right)$. 
\end{empt}
\begin{theorem}\label{stable_thm}
Under hypotheses \ref{stable_empt} one has:
\begin{enumerate}
	\item[(i)] for any $\la\in \La$ there is a $C^*$-algebra $A_\la$ and  a natural $*$-isomorphism $A^\K_\la \cong A_\la \otimes \K$,
	\item[(ii)]  if 			
	\bean
	\mathfrak{S}^\K\bydef \left\{\left\{A^\K_\la\right\}_{\la\in \La}, \left\{\pi^\mu_\nu\otimes 1_{\K^+} : A^\K_\mu \hookto A^K_\nu\right\}_{\substack{\mu, \nu \in \La\\\mu \le \nu}}\right\}
	\eean
	is an algebraical finite covering category (cf. Definition \ref{algebraical_finite_covering_category_defn}) then  
		\bean
	\mathfrak{S}\bydef \left\{\left\{A_\la\right\}_{\la\in \La}, \left\{\pi^\mu_\nu : A_\mu \hookto A_\nu\right\}_{\substack{\mu, \nu \in \La\\\mu \le \nu}}\right\}
	\eean
	is  an algebraical finite covering category,
	\item[(iii)] if the  triple $\left( A\otimes \K, \overline A\otimes \K, \widehat G\right)$  is a {pre}-{covering of the algebraical finite covering category}   $\mathfrak{S}^\K$ (cf. Definition \ref{algebraical_finite_covering_category_defn}) then the  triple  $\left( A, \overline A, \widehat G\right)$ {pre}-{covering of the algebraical finite covering category} $\mathfrak{S}$,
	\item[(iv)] if the  triple $\left( A\otimes \K, \overline A\otimes \K, \widehat G\right)$  is the {disconnected infinite noncommutative covering} of $\mathfrak{S}^\K$ (cf. Definition \ref{disconnected_infinite_noncommutative_covering_defn}) then the  triple  $\left( A, \overline A, \widehat G\right)$ is the {disconnected infinite noncommutative covering} $\mathfrak{S}$,

\item[(v)] if a disconnected infinite noncommutative covering  $\left( A\otimes \K, \overline A\otimes \K, \widehat G\right)$  is {good} (cf. Definition \ref{good_defn}) then $\left( A, \overline A, \widehat G\right)$ is good, 
\item[(vi)] if the triple $\left(A\otimes \K, \widetilde{A}\otimes\K , G\left(\left.\widetilde{A}\otimes \K~\right| A\otimes \K\right)\right)$
is  the  {infinite noncommutative covering} of $\mathfrak S^\K$ (cf. Definition \ref{infinite_noncommutative_covering_defn}) then the triple $$\left(A, \widetilde{A}, G\left(\left.\widetilde{A}~\right| A\right)\cong G\left(\left.\widetilde{A}\otimes \K~\right| A\otimes \K\right)\right)$$
is  the  {infinite noncommutative covering} of $\mathfrak S$.
\end{enumerate}
\end{theorem}
\begin{proof}
	From the Lemma \ref{stable_lem} one can deduce that any $a^\K\in K\left(A \otimes K \right)$ can be represented by the following finite sum:
\be\label{stable_bab_eqn}
a^\K = \sum_{j=1}^n b'_j a_j b''_j; \quad \forall j \in \left\{1,...,n\right\}\quad a_j \in K\left(A^p \right)  \quad b'_j, b''_j \in K\left(\K \right).
\ee
where $A^p$ is given by the equation \eqref{stable_ap_eqn}.
(i) For any $\la\in \La$ we denote by $A_\la\subset A_\la^\K$ the $C^*$-algebra generated by elements 
$$
a^p_\la \bydef \bt\text{-} \sum_{	g \in \ker\left( \widehat{G}\to G_\la\right) }g \overline a^p\quad \in \overline A^p\quad \text{(cf. \eqref{stable_ap_eqn})}.
$$
From the equations \eqref{stable_app_eqn} and \eqref{stable_bab_eqn} it follows that
$
A^\K_\la \cong A_\la \otimes \K$.\\
(ii) and (iii) One should proof that the category $\mathfrak S$ satisfies to the conditions of the Definition \ref{algebraical_finite_covering_category_eqn}.
	\begin{enumerate}
	\item [(a)] 
Follows from (ii) of the Theorem \ref{stable_fin_cov_thm}
	\item[(b)] If $a \in  K\left( A^p \right)_+$ then from the Definition \ref{algebraical_finite_covering_category_defn}  it follows that there is $\widetilde a \in K\left( \overline A \otimes \K\right)_+$ such that $a =  \desc_{\la}\left(\widetilde a  \right)$. If $\widetilde a^p  \bydef p \widetilde a p$ then $$a =  \desc_{\la}\left(\widetilde a^p  \right).$$
\end{enumerate}
(iv) If $\phi: \overline A' \hookto  \overline A''$ is a $\mathfrak{Cov}\left(\mathfrak{S} \right)$-morphism (cf. \ref{infinite_covering_empt}) then
$$
\phi \otimes \Id_\K : \overline A' \otimes\K \hookto\overline A'' \otimes\K
$$
is a $\mathfrak{Cov}\left(\mathfrak{S}^K \right)$-morphism. So a terminal object of the category $\mathfrak{Cov}\left(\mathfrak{S}^K \right)$ corresponds to the terminal object of the category $\mathfrak{Cov}\left(\mathfrak{S} \right)$.\\
(v) One needs check (a) and (b) of the Definition \ref{good_defn}.
\begin{enumerate}
	\item[(a)] If $\widetilde A^K$ is a connected component of $\overline A\otimes \K$ then $\widetilde A^K\cap  \overline A^p\cong \widetilde A$ and $\widetilde A^K \cong \widetilde A\otimes \K$. If both $\widetilde A_1$ and $ \widetilde A_2$ are connected components of $\overline A$ then both $\widetilde A_1\otimes \K$ and $ \widetilde A_2\otimes \K$ are connected component of $\overline A\otimes \K$. From (a) of the Definition \ref{good_defn} it follows that there is $g \in \widehat g$ such that $\widetilde A_1\otimes \K = g\left(\widetilde A_2\otimes \K \right)$. From this circumstance  one has $\widetilde A_1 = g\widetilde A_2$.
	\item[(b)] From the above construction one has 
	\bean
	\overline{A} = \widetilde A \oplus \widetilde A^\perp,\\
		\overline{A}\otimes \K = \widetilde A \otimes \K \oplus \widetilde A^\perp\otimes \K;\\
		 	G\left(\left.\widetilde{A}~\right| A\right)\bydef		\left\{\left. g \in \widehat{G}\right| \forall \widetilde a^\perp \in \widetilde A^\perp \quad g \widetilde a^\perp= \widetilde a^\perp\right\}\cong\\\cong  \left\{\left. g \in \widehat{G}\right| \forall \widetilde a^\perp \in \widetilde A^\perp\otimes \K \quad g \widetilde a^\perp= \widetilde a^\perp\right\}= 	G\left(\left.\widetilde{A}\otimes \K~\right| A\otimes \K\right)
		\eean 
(cf. equation \eqref{infinite_covering_transformation_group_eqn}). For any $\la \in \La$ the natural homomorphism 	$G\left(\left.\widetilde{A}\otimes \K~\right| A\otimes \K\right)\to G_\la$ is surjective, so the homomorphism $G\left(\left.\widetilde{A}~\right| A\right)\to G_\la$ is surjective.
\end{enumerate}
(vi) An ideal  $\widetilde{A}\subset \overline A$ is a connected component of  $\overline A$ if and only if an ideal  $\widetilde{A}\otimes \K\subset \overline A\otimes \K$ is a connected component of  $\overline A\otimes \K$.
\end{proof}


\chapter{Hurewicz homomorphism. A sketch}\label{h_chap}
\section{Preliminaries}\label{h_prel}
\paragraph{}
If $\sX$ is a compact space and $K_n\left(C \left( \sX\right) \right)$ is a group of $K$-homology (cf. \cite{blackadar:ko}) then the for all $n\in \Z$ the map $\sX \mapsto K_n\left(C \left( \sX\right) \right)$ is a functor which be denote by $K_n$. The sequence of functors $K_n$ from the category of topological spaces to category of Abelian groups satisfy to the Definition \ref{top_red_homol_defn}, i.e. it is a reduced homology theory. It turns out that  there is the  {Hurewicz homomorphism} (cf. Definition \ref{top_hurewitz_defn}) 	$h_n:\pi_n\left(\sX, x_0 \right)\to K_n\left(\sX \right)$ such that
\be\label{h_hk1_eqn}
\begin{split}
	h_{K_1}	:\pi_1\left(\sX, x_0 \right)\cong \left[S^1, s_0; \sX, x_0\right]\xrightarrow{K_1}\Hom\left(K_1\left(S^1 \right),  K_1\left(\sX \right)\right) \xrightarrow{\Hom\left(\sigma^1, 1 \right) }\\
	\Hom\left(K_0\left(S^0 \right),  K_1\left(\sX \right)\right)\cong K_1\left(\sX \right).
\end{split}
\ee
From \ref{cirle_k1_empt} it follows that
$$
h_{K_1}	\left(\left[f \right]  \right) \bydef K_1\left( \left[f \right]\right) \left(\left[ 	K^1_{S^1}\right]  \right)
$$
where the following notation is used
\begin{itemize}
	\item $f: \left(S^1, s_0\right)\to \left(\sX, x_0\right)$ is a representative of $\left[f \right]\in \pi_1\left(\sX, x_0 \right)/\left[\pi_1\left(\sX, y_0 \right), \pi_1\left(\sX, y_0 \right)\right]$.
	\item $	K^1_{S^1}:C\left( u\right)\to  Q\left(\K \right)$ is a given by \eqref{s1_k_1_eqn} representative of the generator of $K_1\left(S^1 \right)$.  
\end{itemize}
It turns out that
\be\label{h_k1_f_eqn}
h_{K_1}	\left(\left[f \right]  \right)= \left[K^1_{S^1}\circ f\right] \in K_1\left(\sX \right)
\ee
where right part of \eqref{h_k1_f_eqn} means the represented by $K^1_{S^1}\circ f: \sX \to Q\left( \K\right)$  element of $K_1\left(\sX \right)$ (cf. \eqref{busby_k_eqn}). We would like to  generalize $h_{K_1}$. For any $C^*$-algebra $A$ and any (infinite) noncommutative covering $\left(A, \widetilde A, G\left( \left. \widetilde A\right| A\right) \right)$, we replace a homomorphism
\eqref{h_k1_f_eqn} by a pair $\left( h^{\widetilde A}_{\text{tors}}, h^{\widetilde A}_{\text{tors}}\right)_{\text{free}}  $ of homomorphisms
\be\label{h_main_eqn}
\begin{split}
	h^{\widetilde A}_{\text{tors}}: G\left( \left. \widetilde A\right| A\right)\to K_1\left(A \right)_{\text{tors}} ,\\
	h^{\widetilde A}_{\text{free}}: G\left( \left. \widetilde A\right| A\right)\to K_1\left( A\right)_{\text{free}} \bydef  K_1\left( A\right)/ K_1\left( A\right)_{\text{tors}}
\end{split}
\ee
where $K_1\left( A\right)$ is a group of $K$-homology (cf. Definition \ref{k_hom_defn}) and $K_1\left(A \right)_{\text{tors}}$ is its torsion part.

\section{Basic idea}

\paragraph*{}
If $A$ is a stable $C^*$-algebra (cf. Definition \ref{stable_ca_defn}) then according to \cite{blackadar:ko} any unitary element $u\in M\left( A\right)$ or $u\in M\left( A\right)/A$ yields an element of  $K$-groups of $A$ (cf. Corollary \ref{stable_k_cor}). Otherwise sometimes the element $u$ can define a noncommutative covering of $A$. This circumstance yields an interaction between noncommutative coverings and $K$-theory. This interaction is similar   Hurewicz homomorphism from homotopic groups to homology ones.
\begin{empt}\label{h_pre_empt}
	Let $A$ be a stable  $C^*$-algebra (cf. Definition \ref{stable_ca_defn}), and let $u \in M\left( A\right) $ be an unitary element.  If $A \hookto B\left(\H \right)$ is a faithful representation then $u \in B\left(\H \right)$. Let $\phi_n$ is a Borel $n^{\mathrm{th}}$ root $\phi_n$ of identity map on the set $\left\{\left. z \in \C\right| \left|z \right|=1\right\}$, i.e.
	\be\label{h_n_root}
	\left( \phi_n\right)^n = \Id_{\left\{\left. z \in \C\right| \left|z \right|=1\right\}}
	\ee
	In particular $\phi_n$ can be given by	
	\be\label{h_n_n_root}
	\phi_n\left( \varphi\right)= e^{\frac{i \varphi}{n}} 
	\ee
	where $\varphi\in \left(0, 2\pi\right]$ is the angular parameter on $\left\{\left. z \in \C\right| \left|z \right|=1\right\}$. If  $v \bydef \phi_n\left( u\right)\in   B\left(\H \right)$ then $u = v^n$. Let $\widetilde A \subset B\left(\H \right)$ be a minimal $C^*$- algebra such that $Av^n \subset \widetilde A$ for all $n \in \Z$.
	
 If  $v \bydef \phi_n\left( u\right)\in   B\left(\H \right)$ then $u = v^n$. Let $\widetilde A \subset B\left(\H \right)$ be a $C^*$- algebra which is isomorphic to the left $A$ module given by
	$$
	\widetilde A \cong \sum_{j = 0}^{n-1} Av^i
	$$
	Assume that there is the natural action $\Z_n \times \widetilde A\to \widetilde A$ such that 
	$$
	\forall \overline k \in \Z_n \quad \overline k v = e^{\frac{2\pi i k}{n}}
	$$
Suppose that the natural $A$-module homomorphism $\pi: A \hookto \widetilde A$  is a $*$-homomorphism of $C^*$-algebras.
\end{empt}
\begin{definition}\label{hurewicz_u_n_defn}
	In  situation \ref{h_pre_empt} if $\left(A, \widetilde A , \Z_n , \pi  \right)$  is a  noncommutative finite-fold   covering $\left(A, \widetilde A , \Z_n , \pi  \right)$  (cf. Definition \ref{fin_defn} then we say that $\left(A, \widetilde A , \Z_n , \pi  \right)$ is a $\left(u, n\right)$-\textit{covering}.
\end{definition}
\begin{example}\label{s_1_covering_sample}
	If  $A \bydef C\left( S^1\right)\otimes \K$ then there is the natural inclusion $C\left( S^1\right)\subset M\left(A \right)$ it is known that $K_1\left(C\left( S^1\right) \right) = \Z\left[u\right]$ where $u \in C\left( S^1\right)\in M\left( A\right)$ is an unital element. If $\phi_n$ satisfying  to \eqref{h_n_root} and 	$v = \phi\left( u\right)$ then $C\left( S^1\right)\left[v\right]\cong C\left( \widetilde \sY \right)$ where $\widetilde \sY\cong S^1$ and an inclusion  $C\left( S^1\right)\hookto C\left( S^1\right)\left[v\right]$ corresponds to an $n$-listed covering $\widetilde \sY\cong S^1 \to S^1$. In result one has an $\left(u, n\right)$-{covering} $\left(C\left( S^1\right)\otimes \K, C\left( \widetilde \sY\right)\otimes \K, \Z_n\right)$  of $C\left( S^1\right)\otimes \K$ (cf. Theorem \ref{stable_fin_cov_thm} )
\end{example}
\begin{exm}\label{s_3_covering_sample}
	It is known that $S^3$ is homeomorphic to $SU(2)$ and $K^1\left(SU(2) \right)\cong K_1(C(SU(2))\cong \Z$. The group $K_1(C(SU(2)))$ is generated by $\left[u\right]$ where   $u \in C(SU(2)) \otimes \mathbb{M}_2(\mathbb{C})\subset M\left(C(SU(2))\otimes \K \right)$ is an unitary element which commutes with $C(SU(2))$. Let $C\cong C\left( \sY\right)$ be a generated by $C\left( SU\left(2 \right) \right) \cup \{u\}$ commutative $C^*$-algebra. According to the construction \ref{h_pre_empt} there is a finite-fold covering $\widetilde \sY \to \sY$. If $\widetilde \sY$ is connected then there is a nontrivial  connected noncommutative finite-fold covering of $C(SU(2)))\otimes \K$ (cf. Definition \ref{fin_defn}) then From the Theorem \ref{stable_fin_cov_thm} it follows that a nontrivial  connected noncommutative finite-fold covering of $C(SU(2)))$. But  $\pi_1\left( SU(2)\right)$ is trivial, so from the Corollary  \ref{comm_comm_cor} it follows that $C(SU(2)))$ has no nontrivial  connected noncommutative finite-fold coverings. So there is no $\left(u, n\right)$-{coverings}. Since $\pi_1\left( SU(2)\right)$ it trivial there are no nontrivial elements of $K_1\left( SU(2) \right)$ which lie in $h_{K_1}\left( \pi_1\left( SU(2)\right)\right)$.
\end{exm}
\begin{remark}
	The Example \ref{s_1_covering_sample} demonstrates an element of $K_1$ group which gives nontrivial  connected noncommutative finite-fold coverings.  	The Example \ref{s_3_covering_sample} proves that there are elements of $K_1$ group which do not give nontrivial  connected noncommutative finite-fold coverings.
\end{remark}

\section{Building blocks}
\subsection{Torsion special case}\label{torsion_special_case}
\paragraph{} Here we construct homomorphism for torsion part of $K^1\left( A\right)$. 
\begin{empt}\label{tors_special_case_general}
	Let $A$ be a stable  $C^*$-algebra such that $K^1(A)= G \oplus \mathbb{Z}_n$, where $G$ is an Abelian group. From (\ref{uct_c}) it follows that $K_0(A) \approx G' \oplus \mathbb{Z}_n$. So there is an element $x  \in K_0(A)$ such that $K_0\left(A \right)= G' \oplus \Z_n x$. From \ref{part_index_empt} it follows that there is an unitary $w\in M\left( A\right)/A$ and a Fredholm partial isometry $v \in M\left( A\right)$ such that $w$ lifts to $v$ and
	\bean
	x = \partial \left[w\right]= \mathrm{Index}~v
	\eean
	(cf. equation \eqref{d_u_iv_eqn}).
	From $nx = 0$ it turns out that the index of $v^n$ equals to $0$. It means that $\ker ~v^n \cong \coker ~v^n$. There is unital operator $u \in M\left( A\right)$ such that $\left.u \right|_{\mathrm{coim}~v^n }= \left.v^n \right|_{\mathrm{coim}~w^n }$ and $\left.u \right|_{\ker w^n}$ isomoprhically maps $\ker w^n$ onto  $\coker w^n$. If there is an $\left(u, n\right)$-covering of $A$ then one can construct the finite fold nocommutative covering $\left(A, \widetilde A, \Z_n\right)$. Using it one can construct a homomorphism from the torsion part of $K_1\left( A\right)$ to $G\left( \left. \widetilde A\right| A\right)  /\left[G\left( \left. \widetilde A\right| A\right) , G\left( \left. \widetilde A\right| A\right) \right] \cong G\left( \left. \widetilde A\right| A\right) \cong \Z_n$.

\end{empt}
\begin{problem}\label{h_cuntz_cov}
	Let $n \in \N$, and let $O_n$  be an universal $C^*$-algebra which is generated by isometries $\left\{s_1,...,s_n\right\}$ having following relations
	\bean
	s^*_js_j = 1,\\
	s_js^*_j = p_j,\\
	p_1 + ... + p_n = 1.
	\eean 
	It is proven in \cite{cuntz:k_c_a} that $K_0\left(O_n \right)\cong \Z_{n-1}$. Suppose that $\left[p_1\right]$ is presented by   $u \in U\left(M\left(O_n\otimes\K \right)/\left( O_n\otimes\K\right)  \right)$ such that $\left[u\right] \neq 0 \in K_1\left( M\left(O_n\otimes\K \right)/\left( O_n\otimes\K\right)\right)$ and $\left( n-1\right) \left[u\right]=0$. \\ Question. Does an $\left(u, n-1\right)$-covering (cf. Definition \ref{hurewicz_u_n_defn})  $\left(O_n\otimes\K, \widetilde A, \Z_{n-1}\right)$ exist?  
\end{problem}

\subsection{Free special case}\label{free_special_case}
\paragraph{}
Let $A$ be a stable  $C^*$ - algebra such that $K^1(A)\approx G \oplus \mathbb{Z}$. From (\ref{uct_c}) it follows that
\begin{equation}\label{k1_direct_sum} 
	K_1(A)=G' \oplus \mathbb{Z}[u]
\end{equation}
where $u \in M\left(A\right)$ is an unitary representative  of $u$. Suppose that for all $n \in \N$ there is an $\left(u, n\right)$ {covering} (cf. Definition \ref{hurewicz_u_n_defn}).
Let $n_1, ..., n_j, ... \in \N$ be an increasing sequence of natural numbers such that $n_j$ divides $n_{j+1}$ for all $j\in \N$. For any $j \in \N$ let $A\hookto A_j$ be an $\left(u, n_j\right)$ covering. Suppose that there is a good  	algebraical  finite covering category 
\be\label{free_hur_cat_eqn}
\mathfrak{S}= \left\{A \hookto A_1\hookto ... \hookto A_n\hookto... \right\}
\ee
(cf. Definition \ref{good_defn}) and there is an infinite noncommutative covering  $\left(A, \widetilde A, 	\Z\right)$ of $\mathfrak{S}$. This covering defines a group $\Z \cong  G\left( \left. \widetilde A\right| A\right) $. If we map this group onto $\Z[u]$ then one has a  homomorphism 
\be\label{free_hur_hom_eqn}
G\left( \left. \widetilde A\right| A\right) \to K^1\left( A\right)
\ee

\begin{exercise}
	\begin{enumerate}
		\item If $\T^n$ is an $n$-torus then $\pi_1\left( \T^n\right) \cong \pi_1\left(C\left(\T^n \right)  \right)  \cong K^1\left( C\left(\T^n \right) \right)\cong \Z^n$. Prove that the above construction yields a natural isomorphism $\pi_1\left(C\left(\T^n \right)  \right) \xrightarrow{ \approx} K^1\left( C\left(\T^n \right) \right)$ 
		\item One has $\pi_1\left( S^1\times S^3\right)= \pi_1\left( C\left( S^1\times S^3\right) \right)= \Z$, and $K^1\left(  C\left( S^1\times S^3\right)\right)= \Z[u] \oplus \Z[v]$.
		Prove that the above construction yields a homomorphism
		\bean
		\Z \to \Z[u] \oplus \Z[v],\\
		\left( a, b \right) \mapsto a[u]\oplus 0.	
		\eean 
	\end{enumerate}
\end{exercise}

\section{Assembly}
\paragraph*{}
Here we construct Hurewiscs homomorphism for $C^*$-algebras with finitely generated $K^1$-homology group.
\subsection{Torsion special case}\label{torsion_sc_sec}
\paragraph*{}
Let $A$ be a stable $C^*$-algebra (cf. Definition \ref{stable_ca_defn}) such that the torsion subgroup  $K_0\left( A\right)_{\mathrm{tors}}\subset K_0\left( A\right)$ is finite.   Any element   $\left[p\right] \in K_0\left( A\right)$ such that $n\left[p\right]=0$ can be represented by an unitary element $u \in M\left(A \right)/A$. We way that $u$ is \textit{torsion special}  if there is an $\left(u, n\right)$ covering $\left( A, \widetilde A_u,\Z_n \right)$ (cf. Definition \ref{hurewicz_u_n_defn}). Denote by $K_0\left( A\right)_{\mathrm{tors~spec}}\subset K_0\left( A\right)_{\mathrm{tors}}$  the generated by torsion special elements subgroup. Suppose  $\left\{u_1, ,..., u_k\right\}\subset  M\left(A \right)/A$ be a set of unitary elements such any $x\in K_0\left( A\right)_{\mathrm{tors~spec}}$ can be uniquely represented by the following way
\be\label{h_zn_eqn}
x = c_1 \left[u_1\right]+...+c_k \left[u_k\right]\quad c_j \in \Z_{n_j}\quad n_j \bydef \mathrm{ord}\left(\left[u_j\right] \right) \quad j = 1,..., k.
\ee
Let
\be\label{h_znt_eqn}
\left( A, \widetilde A_{\mathrm{tors}},G_{\mathrm{tors}} = \Z_{n_1} \times...\times \Z_{n_k}\right)
\ee
 be the minimal covering which contains\\ $\left( A, \widetilde A_{u_j},\Z_{n_j}\right)$ for all $j = 1,..., k$. $\widetilde A$ is generated by unitary elements $v^l_j \widetilde A$ such that $v_j =\phi_{n_j}\left(  u_j\right) $ for all $j =1,..., n$, $~l\in \N$ where $\phi_n$ satisfies to \eqref{h_n_root}. For each $j = 1,...,n$ one has $v_j \in M\left(\widetilde A_{\mathrm{tors}} \right)$. For any $g \in G_{\mathrm{tors}}$ we define a character $\chi_g : K_0\left( A\right)_{\mathrm{tors~spec}}\to U\left(1 \right)$ such that
$$
\chi_g\left(\left[u_j\right] \right) = \frac{gv_j}{v_{j}}
$$
From \eqref{h_zn_eqn} it turns out that the  above equation uniquely defines $\chi_g$.
The map $g \mapsto  \chi_g$ is a homomorphism. On the other hand
$\mathrm{\Ext}\left(K_0\left( A\right)_{\mathrm{tors}}, \Z \right)$ is naturally isomorphic to the group of characters on $K_0\left( A\right)_{\mathrm{tors}}$ (cf. Section \ref{ext_char_sec}).
In result a one has homomorphism
$$
h^{\widetilde A_{\mathrm{tors}}}_{\mathrm{tors}} : G_{\text{tors}}\bydef G\left( \left. \widetilde A_{\text{tors}}\right| A\right)  \to \mathrm{\Ext}\left(K_0\left( A\right)_{\mathrm{tors-spec}}, \Z \right).
$$
Taking into account an inclusion $\mathrm{\Ext}\left(K_0\left( A\right)_{\mathrm{tors-spec}}, \Z \right)\subset \mathrm{\Ext}\left( K_0\left(A \right) , \Z\right)$ one 
obtains the \textit{torsion Hurewicz homomorphism}
\be\label{h_tors_eqn}
h^{\widetilde A_{\mathrm{tors}}}_{\mathrm{tors}} : G_{\text{tors}}=G\left( \left. \widetilde A_{\text{tors}}\right| A\right)  \to  \mathrm{\Ext}\left( K_0\left(A \right) , \Z\right).
\ee

\subsection{Free special case}\label{h_free_ass_sec}
\paragraph*{}
Denote by $K_1\left( A\right)_{\mathrm{free}} \bydef K_1\left( A\right)/ K_1\left( A\right)_{\mathrm{tors}}$.   
Any $\left[u\right] \in K_1\left( A\right)/ K_1\left( A\right)_{\mathrm{tors}}$  can be represented by an unitary element $u \in M\left(A\right)$. We say that $u$ is \textit{free special}  if for all $n\in \N$ there is an $\left(u, n\right)$ covering $\left( A, \widetilde A_u,\Z_n \right)$. Denote by $K_1\left( A\right)_{\mathrm{free~spec}}\subset K_1\left( A\right)_{\mathrm{free}}$  the generated by free special elements subgroup. Suppose  $\left\{u_1, ..., u_k\right\}\subset M\left(A \right) $ be a set of unitary elements such any $x\in K_1\left( A\right)_{\mathrm{free~spec}}$ can be uniquely represented by the following way
$$
x = c_1 \left[u_1\right]+...+c_k \left[u_k\right]\quad c_j \in \Z\quad j = 1,..., k
$$
so one has
\be\label{h_k_sum_eqn}
K_1\left( A\right)_{\mathrm{free~spec}}=  \Z \left[u_1\right]+...+\Z \left[u_k\right]\cong \Z^k
\ee\label{h_free_ass_eqn}
Let 
\be
\left( \widetilde A_{\mathrm{tors}}, \widetilde A_{\mathrm{free}},G_{\mathrm{free}} = \Z^k \right)
\ee
 be the minimal covering of  $\widetilde A_{\mathrm{tors}}$ which contains  $\left(u_j, l\right)$-coverings for all $j = 1,...,k$ and $l\in \N$. For any pair $\left(u_j, l\right)$ one has $\phi_l\left(u_j\right)  \in \widetilde A_{\mathrm{free}}$ where $\phi_l$ satisfies to \eqref{h_n_root}. The Abelain group generated by symbols $\left[ \phi_l\left(u_j\right)\right]$ with relations
\be\label{h_q_eqn}
l_2\left[ \phi_{l_1}\left(u_j\right)\right]= \left[ l_1\phi_{l_2}\left(u_j\right)\right]\quad l_1, l_2\in \Z
\ee 
is isomorphic to a vector space $\Q^k$. Moreover the image of the  natural inclusion $K_1\left( A\right)_{\mathrm{free~spec}}\subset \Q^{k}$ is  $\Z^k\subset \Q^{k}$. Let $\chi: \Q\to U\left( 1\right)$ be a nontrivial character  given by $x \mapsto e^{2\pi i x}$. For all $g\in G_{\mathrm{free}}$ define a character $\chi_g :\Q^{k}\to U\left( 1\right)$ such that 
\be\label{h_char_eqn}
\chi_g\left(\left[ \phi_{l}\left(u_j\right)\right]\right)= \frac{g\phi_{l}\left(u_j\right)}{\phi_{l}\left(u_j\right)} 
\ee
such that $K_1\left( A\right)_{\mathrm{free~spec}}\subset \ker \chi_g$. It is explained in \cite{weil:basic_number_theory} that there is a unique $x_g \in \Hom\left( \Q^k, \Q\right)$ such that $\chi_g\left(y \right) = \chi\left(x_g\left(y \right)  \right)$. Moreover  from  $\Z^k\subset \ker \chi_g$ it turns out that $x_g\in \Hom\left( \Z^k, \Z\right)   \cong \Hom\left(K_1\left( A\right)_{\mathrm{free~spec}} , \Z\right)\subset \Hom\left(K_1\left( A\right), \Z\right)$. 
In result one has the \textit{free Hurewicz homomorphism}
\be\label{h_free_eqn}
h^{\widetilde A_{\mathrm{free}}}_{\mathrm{free}} : G_{\text{free}}=G\left( \left. \widetilde A_{\text{free}}\right| A_{\text{tors}}\right) \to \Hom\left(K_1\left( A\right), \Z  \right) 
\ee
\begin{exercise}
	Prove that for any (infinite) noncommutative covering $\left(A, \widetilde A, G\left( \left. \widetilde A\right| A\right) \right)$ there is a pair  of homomorphisms
	\be\label{h_hom_ext_eqn}
	\begin{split}
		h^{ A}_{\mathrm{tors}} :G\left( \left. \widetilde A\right| A\right) \to \mathrm{\Ext}\left( K_0\left(A \right) , \Z\right) ,\\
		h^{ A}_{\mathrm{free}} : G\left( \left. \widetilde A\right| A\right) \to \Hom\left(K_1\left( A\right), \Z  \right). 
	\end{split}
	\ee

\end{exercise}

\begin{definition}\label{h_defn}
	For any (infinite) noncommutative covering $\left(A, \widetilde A, G\left( \left. \widetilde A\right| A\right)\right)$ a given by \eqref{h_hom_ext_eqn} pair $\left(h^{\widetilde A}_{\mathrm{tors}}, h^{\widetilde A}_{\mathrm{free}} \right)$ 
	of homomorphisms is said to be the \textit{Hurewicz homomorphism}. Since both groups   $\mathrm{\Ext}\left( K_0\left(A \right) , \Z\right)$ and $\Hom\left(K_1\left( A\right) , \Z \right)$ are Abelian the following notation will be used
	\be\label{h_gen_eqn}
	\begin{split}
		h^{\widetilde A}_{\mathrm{tors}} :  G\left( \left. \widetilde A\right| A\right)/\left[ G\left( \left. \widetilde A\right| A\right),  G\left( \left. \widetilde A\right| A\right)\right] \to   \mathrm{\Ext}\left( K_0\left(A \right) , \Z\right),\\
		h^{\widetilde A}_{\mathrm{free}} :  G\left( \left. \widetilde A\right| A\right)/ \left[ G\left( \left. \widetilde A\right| A\right),  G\left( \left. \widetilde A\right| A\right)\right] \to \Hom\left(K_1\left( A\right) , \Z \right).
	\end{split}
	\ee
\end{definition}
\begin{remark}
	If $A$ is a $C^*$-algebra and $A\in N$ (cf. Remark \ref{n_alg_rem}) then  equations \eqref{h_gen_eqn} are equivalent to
	\be\label{h_gen_n_eqn}
	\begin{split}
		h^{\widetilde A}_{\text{tors}}: G\left( \left. \widetilde A\right| A\right)\to K_1\left(A \right)_{\text{tors}} ,\\
		h^{\widetilde A}_{\text{free}}: G\left( \left. \widetilde A\right| A\right)\to K_1\left( A\right)_{\text{free}} \bydef  K_1\left( A\right)/ K_1\left( A\right)_{\text{tors}}
	\end{split}
	\ee
	(cf. equation \eqref{uct_c_eqn}).
\end{remark}
\section{Comparison with the classical Hurewicz homomorphism}
\paragraph*{} Consider the described in the Section \ref{h_prel}  situation. Let $\left( \sX, x_0\right)$ be a  pointed  $CW$-complex (cf. Remark \ref{top_obstr_rem}) and let 
$$
f_S : \left(S^1, s_0 \right)\to \left( \sX, x_0\right) 
$$
be a map which represents a generator $\left[f_S\right]$ of an Abelian group \\ $ \pi_1\left(\sY, y_0 \right)/\left[\pi_1\left(\sY, y_0 \right), \pi_1\left(\sY, y_0 \right)\right]$.
\begin{exercise}
	A given by the Proposition \ref{top_pi1_pi1_prop} action yields an equivalence relation (cf. Definition \ref{equivalence_relation_defn}) on $\pi_1\left(\sX \right)$ such that
	$$
	\forall \ga_1, \ga_2 , \a \in \pi_1\left(\sX \right)\quad \ga_1\cdot \a \sim \ga_2\cdot \a.
	$$
	Prove following statements.
	\begin{enumerate}
		\item $\forall  \a, \bt \in\pi_1\left(\sX \right)\quad \a \sim \bt \quad \Leftrightarrow \quad \a^{-1}\bt \in \left[\pi_1\left(\sX \right), \pi_1\left(\sX \right)\right]$.
		\item There is the natural isomorphism  $\pi_1\left(\sX \right)/\left[\pi_1\left(\sX \right), \pi_1\left(\sX \right)\right]\to \left[S^1, \sX \right]$ where $ \left[S^1, \sX \right]$ is a set of classes of homotopically equivalent continuous maps from $S$ to $\sX$.
	\end{enumerate} 
\end{exercise}
Using the above Exercise we replace $
f_S$ with 
$$
f_S : S^1 \to \sX.
$$
Suppose that $f_S$ homeomorphically maps $S^1$ onto $f_S\left(S^1 \right)$ i.e. there is a pair $\left(\sX, f_S\left(S^1 \right)\right)\cong \left(\sX, S^1\right)$ (cf. \cite{spanier:at,switzer:at}). 
From the Theorem \ref{hurewicz_iso_thm} it follows that 	$h_{\mathrm{sing}}\left(\left[f_S\right] \right) \in H_1\left( \sX\right)$ is a generator of $H_1\left( \sX\right)$. Moreover if the period of  $\left[f_S\right]$ equals to $n$ or infinite then the period of $h_{\mathrm{sing}}\left(\left[f_S\right]\right) $ equals to $n$ or infinite respectively. If $u_{S^1}: S^1\to U\left( 1\right)$ is a homeomorphism given by $x \mapsto e^{2\pi i x}$ then $\left[u_{S^1}\right]$ is a generator of $K^1\left(S^1 \right)$ and $Index \left(\left[u_{S^1}\right],  \left[ K^1_{S^1}\right] \right)\bydef Index \left( K^1_{S^1}\left( u_{S^1}\right)\right)  =  1$, i.e. index of the Fredholm operator  $  K^1_{S^1}\left(u_{S^1} \right)$ equals to 1 (cf. Equation \eqref{s1_k_1_eqn}).
\subsection{Torsion special case}
\paragraph*{} Suppose that   $\left[ f_S\right]\in  \pi_1\left(\sX, x_0 \right)/\left[\pi_1\left(\sX, x_0 \right), \pi_1\left(\sX, x_0 \right)\right]$  has period $n$, i.e.  $n \left[ f_S\right]=0$ and $m \left[ f_S\right]\neq 0$, for any $m < n$. Let
\be\label{h_a_eqn}
\a \bydef \left(S^1 \to U\left( 1\right) \quad  x \mapsto e^{2\pi i x}\right) 
\ee
be a map which corresponds to the natural homeomorhism.
\begin{lemma}
	If we consider the above situation then the map $\bt \bydef \a^n : S^1\to U\left(1 \right)$ can be extended up to a map $\varphi:\sX\to U\left(1 \right)$.   
\end{lemma}
\begin{proof}
	Indeed one should prove the existence of the dashed arrow in the below diagram.
	\newline
	\begin{tikzpicture}
		\matrix (m) [matrix of math nodes,row sep=3em,column sep=4em,minimum width=2em]
		{
			U\left( 1\right)\cong S^1	  & & \sX  \\ 
			& S^1  & \\};
		\path[-stealth]
		(m-1-1) edge[dashed]  (m-1-3)		(m-2-2) edge node [left]  {$\bt~~~$} (m-1-1)
		(m-2-2) edge node [right] {$~~\psi$} (m-1-3);
	\end{tikzpicture}
	\\ 
	The Theorem \ref{top_obstr_thm}	can be applied  to the above diagram because  $U\left( 1\right)\cong S^1$ is a $\left(\Z, 1\right)$-space (cf. Definition \ref{top_pi_n_defn}). From $\bt \bydef \a^n$ it follows that
	\bean
	H^1\left(\bt , \Z\right):  H^1\left(S^1 \cong U(1), \Z \right)\to  H^1\left(S^1, \Z \right),\\
	\iota \mapsto n \iota.
	\eean
	where $\iota \in  H^1\left(S^1 \cong U(1), \Z \right)$ is a characteristic element (cf. Definition \ref{top_char_defn}). From $nh_{\mathrm{sing}}\left(\left[f\right] \right)= 0 \in H_1\left( \sX\right)$ it follows that $H^1\left[\bt\right]\left( \iota \right)=  nH^1\left[f\right]\left( \iota \right) = 0\in H^1\left(S^1, \Z \right)$. So $\delta H^1\left[\bt\right]\left( \iota \right)  = 0 \in  H^2\left(\sX, S^1; \Z \right)$ so from the Theorem \ref{top_obstr_thm} $\bt$ it turns out that $\bt$ can be extended over $\sX$.
\end{proof}
Let $\varphi: \sX \to U(1)$ such that $\varphi|_{f\left(S^1 \right)}=\bt$.  Since $SU\left(2\right)$ is a topological group then there is a contravariant functor $\left[~\cdot~;\sY\right]$ from the category of topological spaces and  homotopy classes of continuous maps to the category of groups (cf. Theorem \ref{top_hgr_thm} and Remark \ref{top_hgr_rem}).
The map 
\bean
\phi_{S^1}: f\left( S^1\right) \to SU\left(2 \right),\\
x \mapsto \begin{pmatrix}
	
	\a \left(x \right)   & 0\\
	0  & \a^*\left(x \right) 
\end{pmatrix} 
\eean
is contractible it follows that $\left[\phi_{S^1}\right]= e\in \left[f\left( S^1\right), SU\left(2 \right)\right]$ where $e$ is a neutral element of group. If $\phi'_\sX : \sX \to  SU\left(2 \right)$ is such that $\left[\phi_\sX\right]= e\in \left[\sX, SU\left(2 \right)\right]$ then $\left[\phi_{S^1}\right]= \left[\psi, SU\left(2 \right)\right]\left( \left[\phi'_\sX\right]\right)$
where
$$
\left[\psi, SU\left(2 \right)\right]: \left[\sX, SU\left(2 \right)\right]\to \left[S^1, SU\left(2 \right)\right]
$$
is a homomorphism given by the Remark  Remark \ref{top_hgr_rem}. It turns out that there is $u'_1: \sX \to SU\left( 2\right)$ which is homotopic to $\phi'_\sX$ such that $\left.u'_1\right|_{f\left( S^1\right)}= \phi_{S^1}$, i.e.
\bean
u'_1|_{f\left(S^1 \right)}= \begin{pmatrix}
	
	\a   & 0\\
	0  & \a^*
\end{pmatrix}. 
\eean	
If 		\bean
p_1 \bydef \begin{pmatrix}
	
	1   & 0\\
	0  & 0
\end{pmatrix} 	
\eean
then $u'_1p_1u'^*$ is a projector in $\mathbb{M}_2\left(C\left(\sX\right) \right)$ and
$x=	\left[u'_1p_1u'^*\right] - \left[p_1\right]\in K_0\left( C\left( \sX\right) \right) \cong K^0\left(  \sX \right)$. 
\begin{lemma}
	In the space $\mathrm{Cont}\left( \sX, \SU\left(2 \right)\right)$  of continuous maps $\sX \to \SU\left(2 \right)$ with compact-open topology (cf. \ref{top_comp_open_empt}) there is a homotopy $H:[0,1]\to\mathrm{Cont}\left( \sX, \SU\left(2 \right)\right) $ such that $H\left(0\right) = u'_1$ and $\left( H\left(1\right)\right)^n =  \begin{pmatrix}
		\varphi & 0\\
		0& \varphi^*
	\end{pmatrix}$.
\end{lemma}
\begin{proof}
	If $e\in \left[\sX, SU\left(2\right)\right]$ is the neutral element of the group then from 
	$$
	\left[u'^n \right]= e = \left[\begin{pmatrix}
		\varphi & 0\\
		0& \varphi^*
	\end{pmatrix}\right]
	$$
	it turns out  that there is a homotopy  $G:[0,1]\to\mathrm{Cont}\left( \sX, \SU\left(2 \right)\right) $ such that $G\left(0\right) = u'^n_1$ and $ G\left(1\right) = \begin{pmatrix}
		\varphi & 0\\
		0& \varphi^*
	\end{pmatrix}$. 
	From the theory of Lie groups it follows that there is an open connected neighborhood $\sU\subset \SU\left(2 \right)$ of $\begin{pmatrix}
		1 & 0\\
		0& 1
	\end{pmatrix}$ such that the map  \bean
	\psi: \sU\xrightarrow{\approx}\sU^n; \\ g \mapsto g^n
	\eean is a homeomorphism from $\sU$ onto its image $\sU^n$. There is $\eps_1$ such that
	$G\left( 0\right)^{-1} G\left(t \right)\left( x\right) \in \sU^n$ for any $0\le t \le \eps_1$ and $x \in \sX$. We define   $H\left(t \right)\left( x\right)\bydef H(0)^{-1}\psi^{-1} \left( G\left( 0\right)^{-1} G\left(t \right)\right)$. Then there is $\eps_2$ such that
	$G\left( \eps_1\right)^{-1} G\left(t \right)\left( x\right) \in \sU^n$ for any $\eps_1\le t \le \eps_1+\eps_2$ and $x \in \sX$ and so on. Since $[0,1]$ using finite number of steps we find the required homotopy.
\end{proof}
If $u_1\bydef H\left( 1\right) $ then $x=	\left[u_1p_1u^*_1\right] - \left[p_1\right]$ because $u_1$ is homotopic to $u'_1$. There is the $*$-isomorphism $C\left( \sX\right)\otimes \K = \K\left( \ell^2\left( C\left( \sX\right)\right) \right)$.  From the Theorem \ref{kasparov_stab_thm} it follows that  $\ell^2\left( C\left( \sX\right)\right)\cong \C^2 \oplus H$ where $H \cong \ell^2\left( C\left( \sX\right)\right)$. 
If
$$
w \bydef \mathrm{diag}\left( u_1\begin{pmatrix}
	
	1   & 0\\
	0  & 0
\end{pmatrix}, 1_{B\left(H \right) } \right) \in M\left(\K\left( \ell^2\left( C\left( \sX\right)\right) \right) \right) 
$$
then $w$ is a Fredholm  partial isometry which represents $x \in  K_0\left(C\left( \sX\right)\right)$ (cf. equation \eqref{partial_p_eqn}). From $nx=0$ it follows that $w^n$ represents zero in $K_0\left(C\left( \sX\right)\right) $.  If 
$$
u  	\bydef \mathrm{diag}\left(\begin{pmatrix}
	
	\bt   & 0\\
	0  & 1_\C
\end{pmatrix}, 1_{B\left(H \right) } \right) \in M\left(\K\left( \ell^2\left( C\left( \sX\right)\right) \right) \right) 
$$
then $u$ is unitary and
$$
u \in w^n + \K\left( \ell^2\left( C\left( \sX\right)\right) \right)
$$
i.e. $u\approx w^n$ modulo $\K\left( \ell^2\left( C\left( \sX\right)\right) \right)$. Indeed this unitary is a specialization of explained in \ref{tors_special_case_general} unitary $u$.
If $\phi_n$ is satisfies to  \eqref{h_n_root}  and $\K\left( \ell^2\left( C\left( \sX\right)\right) \right)\hookto B\left( \H\right)$ is faithful representation then 
\be\label{h_z_eqn}
v \bydef \phi_n\left( u\right) = 	 \mathrm{diag}\left(\begin{pmatrix}
	
	\phi_n\left( \bt\right)    & 0\\
	0  & e^{\frac{2\pi i}{n}}
\end{pmatrix},~ e^{\frac{2\pi i}{n}} 1_{B\left(H \right) } \right) \in B\left(\H \right). 
\ee
Let $\widetilde A \subset B\left(\H \right)$ be a minimal $C^*$- algebra such that $\K\left( \ell^2\left( C\left( \sX\right)\right) \right)v^n \subset \widetilde A$ for all $n \in \Z$.
Using $\bt$ we construct $n$-listed universal covering of $\mathcal Z \to \sX$ such that
\be\label{h_zc_eqn}
\mathcal Z = \left\{\left.\left(x, z \right)  \in \sX\times U\left( 1\right)\right|\exists j \in \{1,...,n\} \quad z = \phi^j_n\left( x\right) \right\}
\ee
where $\phi_n$ satisfies to \eqref{h_n_root}. The covering $\mathcal Z \to \sX$ is given by
$$
\left(x, z \right)\mapsto \left(x, z \right)
$$
The action $\Z_n \times \mathcal Z$ is given by
$$
\forall \overline k \in \Z_n \quad \overline k \left(x, z \right)= \left(x, e^{\frac{2\pi i k}{n}}z \right)
$$
where $k\in \Z$ is a representative of $\overline{k}$.
\begin{exercise} Prove following statements
	\begin{enumerate}
		\item If $v_1 \bydef \phi_n\left( \bt\right)$. Prove that $C\left(\sX \right) \left[v_1\right] = C\left(\mathcal Z\right)$.
		\item Prove that $\widetilde A\cong C\left( \mathcal Z\right) \otimes \K$.

	\end{enumerate}
\end{exercise}
\subsection{Free special case}
\paragraph*{}
Suppose that both  $\left[ f\right]\in  \pi_1\left(\sX, x_0 \right)/\left[\pi_1\left(\sX, x_0 \right), \pi_1\left(\sX, x_0 \right)\right]$ have infinite period and are not divisible. 
\begin{lemma}
	If we consider the above situation then the map $\a: S^1\to U\left(1 \right)$ can be extended up to a map $\varphi:S^1\to U\left(1 \right)$.   
\end{lemma}
\begin{proof}
	Consider  the cohomology exact sequence \eqref{top_c_exact_eqn}, i.e.
	\be\label{h_h_eqn}
	...\to	 H^1\left(\sX ;\Z\right) \xrightarrow{H^1\left[f\right]=H^1\left[\a\right]} H^1\left(S^1; \Z \right) \xrightarrow{\delta} H^{2}\left(\sX, S^1 ;\Z\right)\to ...
	\ee
	The map $H_1\left[f\right]: H_1\left(S^1 \right)  \hookto  H_1\left(\sX \right)$ is injective. From $H_1\left(S^1 \right)\cong \Z$ and the Theorem \ref{top_h_to_k_thm} one can deduce that the homomorphism $H^1\left[f\right]: H^1\left(\sX \right)  \to  H^1\left(S^1 \right)$ is surjective. Since the sequence \eqref{h_h_eqn} is exact, one has $\delta = 0$, so $\delta H^1\left[\a\right]\left( \iota \right)  = 0 \in  H^2\left(\sX, S^1; \Z \right)$ so from the Theorem \ref{top_obstr_thm} it turns out that $\a$ can be extended over $\sX$.
\end{proof}
Let $\varphi:\sX \to U(1)$ be a continuous map
$\sX \to U(1)$ such that $\varphi|_{S^1 }=\a$.
The map 	$\varphi: \sX \to U(1)$ can be regarded as an unitary element $u \in C\left( \sX\right) \subset M\left(C\left(\sX\right)\otimes \K \right)$. One has $Index\left(\left[u \right] ,\left[K^1_{S^1}\circ f\right]  \right)= Index\left(K^1_{S^1}\circ f \right)= 1$. From the universal coefficient theorem it follows that $\left[u \right]\in \K^1\left(\sX \right)$ is not divisible. So if $n > 1$ and $\phi_n$ satisfies to \eqref{h_n_root}  then $v \bydef \phi_n\left( u\right)\notin  M\left(C\left(\sX\right)\otimes \K \right)$. Similarly to \eqref{h_zc_eqn} one can construct a covering $\mathcal Z$ such that $G\left(\left.\mathcal Z\right| \sX\right)= \Z_n$ and $C\left(\mathcal Z\right) \cong C\left( \sX\right)\left[v\right]$. 
\section{Stable Hurewicz homomorphism}

If $A$ is a connected $C^*$-algebra then one has
$$
K_j\left( A\right) \cong K_j\left( A \otimes \K\right) \quad j = 0,1
$$
Let $\pi^{\text{stab}}_1\left(A \right)\bydef \pi_1\left( A \otimes \K \right)$ be the {fundamental group}. Let 
		\bean
	\mathfrak{S}_{A\otimes \K} \bydef   \left(\left\{\widetilde\pi_\la: A \otimes \K\hookto  A_\la^\K\right\}_{\la \in \La}, \left\{\widetilde\pi^\mu_\nu: A_\mu^\K \hookto A_\nu^\K\right\}_{\substack{\mu, \nu \in \La\\ \nu \ge \mu}}\right)	
	\eean
	be an algebraical  finite covering category (cf. Definition \ref{algebraical_finite_covering_category_defn}) which corresponds to $\pi_1\left( A \otimes \K \right)$ and suppose that $	\mathfrak{S}_{A\otimes \K}$ is {represented by}
	\bean
	\mathfrak{S}_{A} \bydef   \left(\left\{\pi_\la: A\hookto A_\la\right\}_{\la \in \La}, \left\{\pi^\mu_\nu: A_\mu \hookto A_\nu\right\}_{\substack{\mu, \nu \in \La\\ \nu \ge \mu}}\right).	
	\eean
 If 	$\mathfrak{S}_{A\otimes \K}$ is good (cf. Definition \ref{good_defn}) then ohe has 
	$$
	\pi^{\text{stab}}_1\left(A \right)\bydef \pi_1\left( A \otimes \K \right)\cong G\left(\left.\widetilde{A}\otimes \mathcal{K}~\right|~ A\otimes \mathcal{K}\right)\cong G\left(\left.\widetilde{A}~\right|~ A\right)
	$$
	where $\widetilde{A}$ is the inverse noncommutative limit  of $\mathfrak{S}_A$ (cf. Definition \ref{infinite_noncommutative_covering_defn}). Taking into account homomorphisms \eqref{h_gen_n_eqn} following homomorphisms
	
	\be\label{h_gen_stab_n_eqn}
\begin{split}
	h^{\widetilde A}_{\text{tors}}: \pi^{\text{stab}}_1\left( A\right) \to K_1\left(A \right)_{\text{tors}} ,\\
	h^{\widetilde A}_{\text{free}}:\pi^{\text{stab}}_1\left( A\right)\to K_1\left( A\right)_{\text{free}} \bydef  K_1\left( A\right)/ K_1\left( A\right)_{\text{tors}}.
\end{split}
\ee
\begin{definition}
We say the given by \eqref{h_gen_stab_n_eqn} set of homomorphisms is a \textit{stable Hurewicz homomorphism}.
\end{definition}

\section{Resume}
\paragraph*{}
In the introduction of Soviet edition of the book "Basic Number Theory" \cite{weil:basic_number_theory} Andre Weil wrote the following:
\paragraph*{}
\textit{With this I would like to emphasize that [the current book] is not simply a textbook for the future number theorists. There was a time when Galois theory was considered as something abstract and complicated, demanded only by specialists. Moreover, I knew quite a few excellent mathematicians of my generation that overtly confessed in their complete ignorance about Galois theory and, perhaps, even took pride in that. Today, everyone perfectly understands that [Galois theory] is one of the 'stem' topics, which every serious student of mathematics has to get familiar with during the first years of their study. In my opinion, the same holds for the elementary theory of algebraic number fields, including the class field theory, and I hope that this book will ultimately promote this. I will be glad if the present translation is useful for the new generation of Soviet mathematicians in that regard}.
\paragraph*{}
This book uses this point. If $K$ is a local field (resp. $F$ is an $\mathbf{A}$-field) the class field theory \cite{weil:basic_number_theory} yields  homomorphisms
\be\label{cls_f_eqn}
\begin{split}
	K^\times \to \mathrm{Gal}\left(K^\mathrm{ab}/K \right),\\
	F^\times_{\mathbf{A}}/F^\times \to \mathrm{Gal}\left(F^\mathrm{ab}/F \right)
\end{split}
\ee
Using  the the explained in the Sections \ref{torsion_special_case} and \ref{free_special_case} constructions one can define the Hurewicz homomorphism, i.e similarly to \eqref{cls_f_eqn} one can define a homomorphism
\be\label{h_f_eqn}
\pi_1\left( A\right) /\left[\pi_1\left(A\right), \pi_1\left(A\right)\right]\to K^1\left(A\right).
\ee
Elements of $K^1\left(A\right)$ can be represented by Fredholm operators which are invertible modulo compact operators, otherwise both $K^\times$ and 	$F^\times_{\mathbf{A}}/F^\times$ correspond to invertible elements of fields. Overwise the Galois groups $\mathrm{Gal}\left(K^\mathrm{ab}/K \right)$ are analogs  $\mathrm{Gal}\left(F^\mathrm{ab}/F \right)$ of the  Abelian fundamental group $ \pi_1\left( A\right) /\left[\pi_1\left(A\right), \pi_1\left(A\right)\right]$. The obtaining of   \eqref{cls_f_eqn} and \eqref{h_f_eqn} uses explained in \cite{weil:basic_number_theory} technical methods.

\chapter{Coverings of $C^*$-algebras with Hausdorff spectrum}\label{ctr_chap}
\paragraph*{}
 The notion of the operator space with Hausdorff spectrum is a  generalization of $C^*$-algebra with Hausdorff spectrum. (cf. Section \ref{ctr_oa_sec}). Here we consider noncommutative coverings of these spaces. 
\section{Operator spaces with Hausdorff spectrum}\label{ctr_oa_sec}
\paragraph*{} Here  operator spaces  with Hausdorff spectrum are defined.
\begin{empt}\label{ctr_oa_empt}
	
If $A$ be a $C^*$-algebra, such that the spectrum $\sX$  of $A$ is Hausdorff then $A$ is a $CCR$-algebra (cf. Remark \ref{ctr_open_res_rem}). From the Corollary \ref{ctr_ccr_i_cor} it follows that $A$ is  $C^*$-algebra of  type I and its primitive spectrum coincides with its prime one (cf. Theorem \ref{ctr_hat_check_thm}).
There is the given by the Dauns Hofmann theorem action
\be\label{ctr_oa_act_eqn}
\begin{split}
 C_0\left(\sX\right) \times A \to A;\\
\left(x, a\right)\mapsto xa; \quad x \in C_0\left(\sX\right),\quad a \in A
\end{split}
\ee
(cf.  Theorem \ref{dauns_hofmann_thm}). For any $x\in \sX$ there is the irreducible representation $\rep_x: A\to B\left(\H_x\right)$.
\end{empt}

\begin{definition}\label{ctr_oa_defn}
Let $\left(X, Y\right)$ be a sub-unital operator space (cf. Definition \ref{operator_space_subunital_defn}). Let $A$ be a  $C^*$-algebra Hausdorff spectrum, and let $\left(\pi_X, \pi_Y \right):\left(X, Y\right)\hookto \left(A, A^+\right)$ be a  complete isometry (cf. Definition \ref{op_sum_space_defn}). We say that  $\left(X, Y\right)$ is an \textit{operator space with Hausdorff spectrum} if following conditions hold:
\begin{enumerate}
	\item[(a)] If $\sX$ is the spectrum of $A$ then one has $C_0\left(\sX\right)\times X \subset X$ where the action  $C_0\left(\sX\right) \times A \to A$ is given by  the Dauns Hofmann theorem (cf. equation \eqref{ctr_oa_act_eqn}).
	\item[(b)] For any  $x \in \sX$ if both $X_x$ and $Y_x$ are $C^*$-norm completions of $\rep_x\left(X\right)$ and $\rep_x\left(Y\right)$ respectively then  the $C^*$-algebra $\rep_x\left(A\right)$ is the $C^*$-envelope of  $\left( X_x, Y_x\right)$, i.e. $\rep_x\left(A\right)= C^*_e\left( X_x, Y_x \right)$ (cf. Definition \ref{operator_space_envelope_defn}).
\end{enumerate}
\end{definition}
\begin{lemma}\label{ctr_oa_lem}
	Under the hypotheses of the Definition \ref{ctr_oa_defn} the $C^*$-algebra $A$ is the $C^*$-envelope of $\left(X, Y \right)$ (cf. Definition \ref{operator_space_envelope_defn}), i.e. $A = C^*_e\left(X, Y \right)$. 
\end{lemma}
\begin{proof}
	Firstly we prove that the minimal unitization $A^+$ (cf. Definition \ref{multiplier_min_defn}) of $A$  is the $C^*$-envelope of $Y$ (cf. Definition \ref{c_env_sp_defn}).  Denote by $k: Y \hookto A^+$ the natural inclusion. If $\left( C^*_e\left(Y\right) , j\right)$ is the $C^*$-envelope of $Y$ then from the Definition  \ref{c_env_sp_defn} it follows that there is the surjective unital $*$-homomorphism  $\pi: A^+ \to  C^*_e\left(Y\right)$ such that $j = \pi\circ k$. If we represent $A^+ = \C~ 1_{A^+} \oplus A$ then clearly $\ker \pi \subset 0 \oplus A \cong A$. If $\ker\pi \neq \{0\}$ then $\ker \pi$ is a closed ideal which corresponds to an a nonempty open subset $\sU \subset \sX$ where $\sX$ is the spectrum of $A$ (cf. Theorem  \ref{jtop_thm}). If $x_0 \in \sU$ is any point then from (b) of the Definition \ref{ctr_oa_defn} there is $a \in X$ such that $\rep_{x_0} \left( a\right)\neq 0$. Since the space $\sX$ is Hausdorff there is $f \in C_0\left(\sX\right)$ such that $f\left(x_0\right)= 1$ and $\supp f \subset \sU$. Clearly $fa \neq 0$. From (a) of the Definition \ref{ctr_oa_defn} it follows that $fa \in Y$, otherwise one has $k\left( fa\right)  \in \ker \pi$. Hence $\pi \circ k$ is not injective, so there is a contradiction. It turns out that $\ker \pi =  \{0\}$, i.e. $\pi$ is injective. Otherwise $\pi$ is surjective  the Definition  \ref{c_env_sp_defn}, hence $\pi: A^+ \cong   C^*_e\left(Y\right)$ is an isomorphism.
	
	 Secondly we prove that $A = C^*_e\left(X, Y \right)$. Denote by $A' \stackrel{\text{def}}{=} C^*_e\left(X, Y \right) \subset  C^*_e\left(Y \right)$. Let us prove that $A'$ is a $C_0\left(\sX\right)$-module. Select $f \in C_0\left( \sX\right)$. 	 
	For any $a' \in A'$ and $\eps > 0$ there is a polynomial
	 \bean
	b'= \sum_{j = 1}^n x_{j,1}\cdot x_{j,2}\cdot...\cdot x_{j,k_j} \text{ where } x_{1,1}, ..., x_{n,k_n}\in X\cup X^*,\\
	\left\| b' - a'\right\|< \frac{\eps}{\left\| f\right\|}.
	\eean
	Otherwise $X$ is  $C_0\left(\sX\right)$-module it turns out that $fx_{j,1}\in X$ for all $j = 1,..., n$. So one has
		 \bean
	b''= \sum_{j = 1}^n \left( f x_{j,1}\right) \cdot x_{j,2}\cdot...\cdot x_{j,k_j} \in A',\\
	\left\| b'' - fa'\right\|< \eps,
	\eean
	hence $fa' \in A'$.	 
	 From the Definition \ref{ctr_oa_defn}  it follows that $X\cup X^* \subset A$, hence $A'\subset A$. Suppose that that  there is $x \in \sX$ $\rep_x\left(A'\right)\subsetneqq \rep_x\left(A\right)$. There are $a \in A$ and $\eps > 0$ such that $\left\| \rep_x\left(a \right) -\rep_x\left(a' \right)\right\| > \eps$ for all $a' \in A'$. Since $A_x$ is generated by $X_x$ as $C^*$-algebra one has $X_x \subsetneqq  \rep_x\left(A'\right)$, hence $X \subsetneqq A'$. There is a contradiction, it follows that $\rep_x\left(A'\right) = \rep_x\left(A\right)$. From the Corollary \ref{top_sub_eq_cor} it follows  that $A' = A$.
\end{proof}
\section{Basic constructions}\label{ctr_bas_constr}
\subsection{$C^*$-algebras induced by coverings}

  \paragraph*{}
   If $A$ is a $C^*$-algebra with Hausdorff spectrum $\sX$ then according to  the Example \ref{blowing_hausdorff_exm} there is  Hausdorff blowing-up $C_0\left(\sX\right)\hookto M\left(A\right)$ (cf. Definition \ref{blowing_defn}).

  \begin{exercise}
 Let $C_0\left( \sX\right) \to M\left(A\right)$ be described above Hausdorff blowing-up. Prove that for any open subset $\sU\subset\sX$ there are natural $*$-isomorphisms
 \be\label{ctr_au_iso_eqn}
 \left.A\right|_\sU \cong A_\sU\cong  ~_\sU A\cong~_\sU A_\sU
 \ee
 where an ideal  $\left.A\right|_\sU$ is given by \eqref{open_ideal_eqn}, both  ideals $ A_\sU$  and   $_\sU A$ anda hereditary $C^*$-subalgebra $~_\sU A_\sU$ are explained in the Definition \ref{blowing_ideals_au_ua_defn}.
  \end{exercise}
Our following constructions require an analog of the Dini's theorem.
\begin{lemma}\label{ctr_dini_lem}
	Let $A$ be a  $C^*$-\textit{algebra}  having a compact Hausdorff spectrum $\sX$. Let $\La$ be a directed set, and let $\left\{a_\la\right\}_{\la\in\La}\subset A$ be a net (cf. Definition \ref{top_net_defn}) such that
	\begin{itemize}
		\item 
		\be\label{ctr_sup_eqn}
		\forall \mu, \nu \in \La \quad \mu \ge\nu\quad \Rightarrow \quad  a_\mu\ge a_\nu.	
		\ee
		\item  There is $a\in A$ such that for all $x \in \sX$ there is the following limit
		$$
		\lim_{\la\in\La}\rep_x\left( a_\la\right) = \rep_x\left( a\right).
		$$
	\end{itemize}
	Then there is the $C^*$-norm limit $	\lim_{\la\in\La}a_\la = a$.
	
\end{lemma}
\begin{proof}
	Put $b_\la\bydef a-a_\la\in A$. We have to prove that $b_\la \to 0$ with respect to the $C^*$-norm topology.
	Let $\eps> 0$ be given. From (b) of the Lemma \ref{hausdorff_spectrum_lem} it turns out that for all $\la\in \La$ the map
	$$
	x\mapsto \left\|  \rep_x\left(b_\la\right)\right\|
	$$
	is continuous, so the set 
	$$
	\sX_\la \bydef \left\{x \in \sX \left|\left\|  \rep_x\left(b_\la\right)\right\|>\eps \right. \right\}
	$$
	is closed. Moreover $\sX_\la$ is compact because it is a subset of a compact space $\sX$. From \eqref{ctr_sup_eqn} it follows that 
	$$
	\forall \mu, \nu \in \La \quad \mu \ge\nu\quad \Rightarrow \quad  \sX_\mu \supset  \sX_\nu.	
	$$
	Fix $x\in \sX$. Since $\rep_x\left(b_\la\right) \to 0$
	we see that $x\notin \sX_\la$ if $\la$ is sufficiently large. Thus $x \notin \bigcap \sX_\la$. In other words, an intersection
	$\bigcap   \sX_\la$ is empty. It follows that   $\bigcup \left(\sX \setminus\sX_\la\right)= \sX$. Otherwise $\sX \setminus\sX_\la$ is open for all $\la\in \La$, taking into account that the space $\sX$ is compact there is $\la_0\in\La$ such that $\sX \setminus\sX_{\la_0}= \sX$ and  $\sX_{\la_0}$ is empty. 
	It follows that  $0\le \left\|  \rep_x\left(b_\la\right)\right\| < \eps$ for all $x \in \sX$  and for all $\la\ge \la_0$. This fact proves the theorem.
\end{proof}
\begin{remark}
	The proof of the Lemma \ref{ctr_dini_lem} is similar to the proof of the Theorem \ref{dini_thm}.
\end{remark}

\section{Finite-fold coverings}
\subsection{Coverings of $C^*$-algebras}
\begin{lemma}\label{ctr_sufficient_lem} 
	If   $A$ is a $C^*$-algebra, with locally compact, connected,  locally connected, Hausdorff   spectrum $\sX$, and $p: 	\widetilde \sX \to \sX$ is a finite-fold transitive covering then there is a noncommutative finite-fold covering (cf. Definition \ref{fin_defn})
	$$
	\left(A, A_0\left( \widetilde{\sX }\right) , G\left(\left.\widetilde\sX\right|\sX\right),A_0\left(p\right) \right)
	$$
	where $A_0\left( \widetilde{\sX }\right)$ is given by the Definitions  \ref{top_lift_a_f_defn} and/or \ref{blowing_finite_lift_defn} $p$-lift of $A$, and $*$-homomorphism  $A_0\left(p\right)$ is given by the equations \eqref{top_c0_ob_eqn} and/or \eqref{blowing_lift_fin_eqn}.
\end{lemma}
\begin{proof}
	This lemma follows from the Example \ref{blowing_hausdorff_exm} and the Theorem \ref{blowing_sufficient_covering_thm}.
\end{proof}
\begin{lemma}\label{ctr_comm_lem}
	Let $\left(A, \widetilde{A}, G, \pi \right)$ be a  noncommutative  finite-fold covering with unitization (cf. Definition \ref{fin_unital_defn}), Let  
	\be\label{ctr_comm_eqn}
	H  \bydef \ker\left(  G \to \left\{\left.g \in \mathrm{Homeo}\left(\widetilde\sX \right)\right| \forall \widetilde x \in \widetilde\sX \quad p\left( \widetilde x\right)= p\left( g\widetilde x\right) \right\}\right)
	\ee
	where $p:  \widetilde\sX \to \sX$ is a given by the Proposition \ref{spectrum_covering_finite_prop}  map from the spectrum $\widetilde\sX$ of $\widetilde A$ to the spectrum $\sX$ of $A$.
	If   $\sX$  is a  locally connected, locally compact Hausdorff space then one has:
	\begin{enumerate}
		\item [(i)] the space $\widetilde \sX$  is Hausdorff,
		\item[(ii)] the  map $p : \widetilde \sX\to \sX$ is a transitive finite-fold covering with $G\left(\left.\widetilde{\sX}~\right|\sX\right)\cong G/H$,
		\item[(iii)] there is a finite-fold covering with unitization given by
	\be\label{ctr_comm_e_eqn}
\left(A ,A_0\left( \widetilde\sX\right), G\left(\left.\widetilde{\sX}~\right|\sX\right),A_0\left(p \right) \right)\cong 	\left(A,\widetilde{A}^H, G/H, \pi^H    \right)
	\ee
		\item[(iv)] the normal subgroup $H \subset G$ is $\left(A, \widetilde{A}, G, \pi \right)$-proper (cf. Definition \ref{proper_subgroup_fin_defn})
	\end{enumerate}
\end{lemma}
\begin{proof}
	(i) Follows from the Corollary \ref{spectrum_ff_p_cor}.\\
	(ii) Follows from the equation \eqref{ctr_comm_eqn} and the Corollary \ref{spectrum_ff_p_cor}.\\
	(iii) Any $\widetilde{a}^H\in \widetilde{A}^H$ is given by the family 
	$$
	\left\{\rep_{  \widetilde{x }}\left( \widetilde{a}^H\right)\right\}_{\widetilde x \in \widetilde \sX } 
	$$
	and from the Corollary \ref{spectrum_ff_p_cor} it follows that
	$$
\forall \widetilde x \in \widetilde \sX \quad \rep_{  \widetilde{x }}\left( \widetilde{A}^H\right) = 	\rep_{  \widetilde{x }}\left(\pi^H\left(  A\right) \right). 
	$$
	From the Lemma \ref{hausdorff_spectrum_lem} it follows that  $\pi^H\left( A\right)$ are continuity structure for  $\widetilde \sX$ and the	$\left\{\rep_{  \widetilde{x }}\left( \pi^H\left(  A\right)\right)\right\}_{\widetilde x \in \widetilde \sX }$ (cf. Definition \ref{operator_fields_continuity_defn}). Let  $\widetilde a^H\in \widetilde{A}^H$ and $\widetilde x_0 \in \widetilde \sX$. If  $\widetilde f_{\widetilde x_0}\in C_0\left(\widetilde\sX \right)$ is  a $\widetilde x_0$-stump (cf. Definition \ref{top_stump_defn}) such that $\supp \widetilde f_{\widetilde x_0}$ is homeomorphically mapped onto $p\left( \supp \widetilde f_{\widetilde x_0}\right)$ and there is an open neighborhood $\widetilde \sU$ of $\widetilde x_0$ such that $\widetilde f_{\widetilde x_0}\left(\widetilde \sU \right)  = \{1\}$. If $\left\{g_1, ..., g_n\right\}\subset G$ is a set of representatives of $G/H$ then the element $\widetilde a' \bydef \sum_{g \in \left\{g_1, ..., g_n\right\}} g\left( f_{\widetilde x_0}\widetilde a^H\right) $ is $G$-invariant, i.e. $\widetilde a' \in \widetilde A^G = \pi\left(A \right)$. If $a \bydef \pi^{-1}\left(\widetilde a' \right)$ then
	$$
\forall \widetilde x \in \widetilde\sU \quad \rep_{  \widetilde{x }} \left( \widetilde{a}^H\right) = 	\rep_{  \widetilde{x }}\left(\pi^H\left(  a\right) \right)	
	$$
	i.e. the vector field $\left\{\rep_{  \widetilde{x }} \left( \widetilde{a}^H\right)\right\}_{\widetilde x \in \widetilde \sX}$ continuous with respect to $\pi^H\left( A\right) $ at $x_0$ (cf. Definition \ref{op_cont_fields_defn}). Since the point $\widetilde x_0$ is arbitrary the vector field $\left\{\rep_{  \widetilde{x }} \left( \widetilde{a}^H\right)\right\}_{\widetilde x \in \widetilde \sX}$  is {continuous} on $\widetilde \sX$, i.e. one has a natural inclusion 
$$
\phi: \widetilde{A}^H \hookto C\left(\widetilde\sX, \left\{\rep_{  \widetilde{x }}\left(   \pi^H\left( A\right) \right)\right\}, \pi^H\left( A\right) \right)
$$
(cf. equation \eqref{top_c_sec_eqn}). From the  Lemma \ref{hausdorff_spectrum_lem} it follows that for all $\widetilde{a}^H\in \widetilde{A}^H$ the map
\bean
\mathrm{norm}_{\widetilde{a}^H} :\widetilde \sX \to \left[0,  \left\|\widetilde{a}^H\right\|  \right],\\
\widetilde x \mapsto \left\|  \rep_{  \widetilde{x }}\left(   \widetilde{a}^H\right)\right\| 
\eean
lies in $C_0\left(\widetilde \sX \right)$ it follows that
$$
\phi\left(  \widetilde{A}^H \right) \subset  C_0\left(\widetilde\sX, \left\{\rep_{  \widetilde{x }}\left(   \pi^H\left( A\right) \right)\right\}, \pi^H\left( A\right) \right)
$$  
(cf. equation \eqref{top_cc0_eqn}).  From the Lemma \ref{top_cs_nc_lem} it follows that the space  $C_0\left(\widetilde\sX \right)\pi^H\left(A \right) $ is dense in 	$C_0\left(\widetilde\sX, \left\{\rep_{  \widetilde{x }}\left(   \pi^H\left( A\right) \right)\right\}, \pi^H\left(A \right) \right)$. From $C_0\left(\widetilde\sX \right)\pi^H\left(A \right)\subset \widetilde{A}^H$ we conclude that $\phi\left(  \widetilde{A}^H \right)$ is dense in $C_0\left(\widetilde\sX, \left\{\rep_{  \widetilde{x }}\left(   \pi^H\left( A\right) \right)\right\}, \pi^H\left( A\right) \right)$ and taking into account that $\widetilde{A}^H$ is norm closed we conclude that $\phi\left(  \widetilde{A}^H \right)= C_0\left(\widetilde\sX, \left\{\rep_{  \widetilde{x }}\left(   \pi^H\left( A\right) \right)\right\}, \pi^H\left( A\right) \right)$, i.e. there is a natural isomorphism.
$$
\varphi: \widetilde{A}^H \xrightarrow{\cong} C_0\left(\widetilde\sX, \left\{\rep_{  \widetilde{x }}\left(   \pi^H\left( A\right) \right)\right\}, \pi^H\left( A\right) \right)
$$
On the other hand from the Lemma \ref{top_compact_c0_lem} it turns out that 
$$
C_0\left(\widetilde\sX, \left\{\rep_{  \widetilde{x }}\left(   \pi^H\left( A\right) \right)\right\}, \pi^H\left( A\right) \right) = C_0\left(\lift_p\left[ C\left(\sX, \left\{\rep_x\left(A \right) \right\}_{x\in \sX},  A \right) \right] \right).
$$
Moreover $A = C_0\left(\sX, \left\{\rep_x\left(A \right) \right\}_{x\in \sX},  A \right)$ and from
$$
\forall a \in A \quad \rep_{\widetilde{x }}\left(   \pi^H\left( a\right)\right) =\rep_{p\left( \widetilde{x }\right) }\left(   a\right) 
$$
it follows that the homomorphism $\pi^H: A \to \widetilde A^H$ is equivalent to
\be\label{top_lift_a0_eqn}
\begin{split}
\lift_p|_{C_0\left(\sX, \left\{\rep_x\left(A \right) \right\}_{x \in \sX},  A \right)}:  C_0\left(\sX, \left\{\rep_x\left(A \right) \right\}_{x \in \sX},  A \right)\hookto\\
\hookto C_0\left(\widetilde\sX, \left\{\rep_{  \widetilde{x }}\left(   \pi^H\left( A\right) \right)\right\}_{\widetilde x \in \widetilde \sX}, \pi^H\left( A\right) \right).
\end{split}
\ee
(cf. Lemma \ref{top_compact_c0_lem}).
However the $*$-homomorphism \eqref{top_lift_a0_eqn} is equivalent to $A_0\left(p \right): A \hookto A_0\left( \widetilde\sX\right)$ (cf. equations \eqref{top_c0_ob_eqn}  and \eqref{top_c0p_ob_eqn}, i.e. $\pi^H: A \to \widetilde A^H$ is equivalent to $A_0\left(p \right): A \hookto A_0\left( \widetilde\sX\right)$, so one has
	$$
\left(A,\widetilde{A}^H, G/H, \pi^H A_0\left(p \right)   \right)\cong\left(A ,A_0\left( \widetilde\sX\right), G\left(\left.\widetilde{\sX}~\right|\sX\right),A_0\left(p \right) \right)
$$
(iv)	From (iii) it follows that the quadruple  $\left(A,\widetilde{A}^H, G/H, \pi^H A_0\left(p \right)   \right)$ is equivalent to $\left(A ,A_0\left( \widetilde\sX\right), G\left(\left.\widetilde{\sX}~\right|\sX\right),A_0\left(p \right) \right)$ one. However from the Lemma \ref{blowing_proper_group_lem} it follows that
	$
	\left(A,\widetilde{A}^H, G/H, \pi^H A_0\left(p \right)   \right)=\left(A ,A_0\left( \widetilde\sX\right), G\left(\left.\widetilde{\sX}~\right|\sX\right),A_0\left(p \right) \right)
	$
	is a {noncommutative finite-fold  pre-covering}, i.e. the normal subgroup $H \subset G$ is $\left(A, \widetilde{A}, G, \pi \right)$-proper.
\end{proof}

\begin{lemma}\label{ctr_ccr_sp_lem} 
	Let $\left(A, \widetilde{A}, G, \pi \right)$ is  a noncommutative finite-fold  quasi-covering with unitization.
	If $A$ is a $CCR$-algebra (cf. Definition \ref{ccr_defn}) then $\widetilde{A}$ is  a $CCR$-algebra. 
\end{lemma}
\begin{proof}
	From the Definition \ref{fin_unitization_defn} it turns out that
	\begin{enumerate}
		\item[(a)] 
		There are unital $C^*$-algebras $B$, $\widetilde{B}$  and inclusions 
		$A \subset B$,  $\widetilde{A}\subset \widetilde{B}$ such that $A$ (resp. $B$) is an essential ideal of $\widetilde{A}$ (resp. $\widetilde{B}$) and $A = B\bigcap \widetilde{A}$,
		\item[(b)] There is a 
		unital  noncommutative finite-fold covering $\left(B ,\widetilde{B}, G, \widetilde{\pi} \right)$ such that $\pi = \widetilde{\pi}_A$ and the action $G \times\widetilde{A} \to \widetilde{A}$ is induced by $G \times\widetilde{B} \to \widetilde{B}$.
	\end{enumerate}
	Any irreducible representation  $\widetilde{\rho} : \widetilde{A}\to B\left(\widetilde{\H} \right)$ can be uniquely extended up to $\widetilde{\rho} : \widetilde{B}\to B\left(\widetilde{\H} \right)$. If $\widetilde \xi \in \widetilde{\H}$ then from the Theorem \ref{irred_thm} it follows that $\widetilde \H = \widetilde \rho\left(\widetilde B \right) \widetilde \xi$. If $\H \bydef \widetilde \rho\left( B\right)\widetilde \xi$ then there is a representation $\rho : B \to B\left(\H\right)$ and taking into account the Proposition   \ref{spectrum_covering_finite_prop} we conclude that the representation is irreducible. According to the
	Definition \ref{fin_unital_defn} and the 	
	Remark \ref{fin_cov_kasp_rem} $\widetilde{B}$ is finitely generated projective module $B$-module, i.e. 
	$$
	\widetilde{B}\cong p_{\widetilde{B}} B^n
	$$
	where $p_{\widetilde{B}}: B^n \to B^n$ is a projector. So from $\widetilde \H = \widetilde \rho\left(\widetilde B \right) \widetilde \xi$ and $ \H = \widetilde \rho\left( B \right) \widetilde \xi$ it follows that 
	$$
	\widetilde{\H}\cong p_{\widetilde{\H}} \H^n
	$$
	where $p_{\widetilde{\H}}: \H^n \to \H^n$ is a projector.
	There is a diagonal action $\rho': B \to B\left( \H^n\right)$ such that
	$$
	b \in B \quad \Rightarrow\quad \rho' \left( b\right) \bydef \begin{pmatrix}
		\rho\left(b \right)  & 0 &\ldots & 0\\
		0& \rho\left(b \right)   &\ldots&0 \\
		\vdots& \vdots &\ddots&\vdots\\
		0& 0 & \ldots &\rho\left(b \right) 
	\end{pmatrix}. 
	$$
	and 
	$\widetilde \rho \left( b\right) = p_{\widetilde{\H}}\rho'\left( b\right)p_{\widetilde{\H}}$. 
	If $a \in A$ then $\rho\left(a \right)$ is compact (cf. Definition \ref{ccr_defn}). So  both operators  $\rho'\left(a \right)$ and $\widetilde \rho \left( a\right) = p_{\widetilde{\H}}\rho'\left( a\right)p_{\widetilde{\H}}$ are compact, i.e.
	\bean
	\forall a \in A \quad \widetilde \rho \left( a\right) \in \K\left( \widetilde\H\right), 
	\eean
	From \eqref{fin_comp_eqn}  it follows that there are  $\widetilde b_1, ..., \widetilde b_n \in \widetilde B$ such that 
	\be\label{ctr_i_comp_eqn}
	\widetilde A = \widetilde b_1A + ... + \widetilde b_nA,
	\ee 
	hence one has
	$$
	\widetilde \rho \left(	\widetilde A \right) = \widetilde \rho \left(\widetilde b_1\right) \widetilde \rho \left(\widetilde A\right) + ... +\widetilde \rho \left( \widetilde b_n\right) \widetilde \rho \left(A\right)	
	$$
	and taking into account that 	
	$\widetilde{\rho}\left(A \right) \subset \K\left( \widetilde{\H}\right)$ we conclude that $\widetilde{\rho}\left(\widetilde A \right) \subset \K\left( \widetilde{\H}\right)$. From the Corollary \ref{ctr_ccr_i_cor} and the Theorem \ref{ctr_hat_check_thm} it follows that $\widetilde{\rho}\left(A \right) = \K\left( \widetilde{\H}\right)$. 
\end{proof}
\begin{lemma}\label{ctr_ess_dense_lem}
	If $A$ is a $C^*$-algebra with the Hausdorff spectrum $\sX$ then the bijective map \eqref{ctr_open_id_eqn} yields the one-to-one correspondence between essential ideals and open dense subsets of $\sX$. 
\end{lemma}
\begin{proof}
	The given by \eqref{ctr_open_id_eqn} map $	\sU \leftrightarrow 	A|_{\mathcal U}$ yields the one-to-one correspondence between  ideals and open  subsets of $\sX$ (cf. Remark \ref{ctr_open_res_rem}). 
	
	Suppose that $\sU \subset \sX$ is a dense subset of $\sX$ and
	$a \in A$ be such that $a\neq 0$ and $a A|_{\mathcal U}=\left\{0\right\}$. From the Lemma \ref{hausdorff_spectrum_lem} it follows that there is an open subset $\sV \subset \sX$ such that $\rep_x\left( a\right) \neq 0$ for all $x \in \sV$. Since $\sU$ is dense, one has $\sV \cap \sU \neq \emptyset$. There is $x_0\in \sU$  such that  $\rep_{x_0}\left( a\right) \neq 0$. From the Exercise \ref{top_completely_regular_exer} it follows that the space $\sX$ is completely regular (cf. Definition \ref{top_completely_regular_defn}). There is a continuous function $f_\sU: \mathcal X \to \left[0,1 \right]$ such that $f_\sU\left(x_0 \right)= 1$ and $f_\sU\left(\mathcal X\setminus\sU \right)= \left\{0\right\} $.  Otherwise from the Lemma \ref{ctr_fact_lem} if follows that $f_\sU a^* \in A|_{\mathcal U}$.  From $\left\|\rep_{x_0}\left( a\left(  f_\sU a^*\right) \right) \right\| = \left\|\rep_{x_0}\left( a \right) \right\|^2\neq 0$ it follows that $a A|_{\mathcal U}\neq\left\{0\right\}$.  From the Lemma \ref{essential_lem} it turns out that $A|_{\mathcal U}$ is an essential ideal. 
	
	If $\sU$ is not dense then there is an open subset $\sV \subset \sX$ such that $\sU\cap\sV= \emptyset$. There is $a \in A$ such that and $x_0\in \sV$ such that $\rep_{x_0}\left(a\right)\neq 0$. There is a continuous function $f_\sV: \mathcal X \to \left[0,1 \right]$ such that $f_\sV\left(x_0 \right)= 1$ and $f_\sV\left(\mathcal X\setminus\sV \right)= \left\{0\right\}$, and $\rep_{x_0}\left( f_\sV a \right)= \rep_{x_0}\left(  a \right)\neq 0$. So $f_\sV a\neq 0$. Otherwise for all $f \in C_0\left(\sU \right)$ from $f_\sV f= 0$ it follows that $\left(f_\sV a \right)\left(fb\right)= 0$ for all $b \in A$.  On the other hand form the Lemma \ref{ctr_fact_lem} it follows that  $A|_{\mathcal U}$ is the $C^*$-norm closure of $C_0\left(\sU \right)A$, so one has $\left(f_\sV a \right)A|_{\mathcal U}=\left\{0\right\}$, i.e. the ideal  $A|_{\mathcal U}$ is not essential.
\end{proof}

\begin{empt}\label{ctr_fin_empt}
	Let $A$ be a   $C^*$-algebra  with Hausdorff spectrum $\sX$, and suppose that $A$ is $CCR$ (cf. Definition \ref{ccr_defn}). If $\left( A, \widetilde A, G, \pi\right)$ is a noncommutative finite-fold covering with unitization  (cf. Definition \ref{fin_unitization_defn}) then from the Lemma \ref{ctr_ccr_sp_lem} it follows that the $C^*$-algebra $\widetilde A$ is $CCR$. If $\widetilde \sX$ is the spectrum of $\widetilde A$ then from the Proposition \ref{spectrum_covering_finite_prop} it follows that there is a natural homomorphism
	$$
	\phi: G \to \mathrm{Homeo}\left(\widetilde \sX \right).  
	$$
	If $H \bydef \ker \phi$ then from the Lemma \ref{proper_subgroup_fin_lem} it turns out that  a quadruple 
	\be\label{ctr_proper_subgroup_finh_eqn}
	\left(\widetilde{A}^H, \widetilde{A}, H, \left.\Id_{\widetilde{A}}\right|_{\widetilde{A}^H} \right).
	\ee
	a {noncommutative finite-fold  pre-covering} (cf. Definition \ref{fin_pre_defn}).

\end{empt}

\begin{lemma}\label{ctr_fin_h_lem}
Under the hypotheses of \ref{ctr_fin_empt} the quadruple $\left( A, \widetilde A, G, \pi\right)$ is equivalent to $\left(A ,A_0\left( \widetilde\sX\right), G\left(\left.\widetilde{\sX}~\right|\sX\right),A_0\left(p \right) \right)$ one.
\end{lemma}
\begin{proof}
	Suppose that the group $H\bydef \ker\left(G \to \mathrm{Homeo}\left(\widetilde \sX \right) \right)$ is not trivial and for any $\widetilde a \in \widetilde A$ denote by
	\bean
	\widetilde a^\parallel \bydef \frac{1}{\left|H \right| } \sum_{	g \in H}g \widetilde a,\\
	\widetilde a^\perp \bydef \widetilde a-\widetilde a^\parallel,\\
	\widetilde A^\parallel \bydef \left\{\left.\widetilde a^\parallel \in \widetilde A\right| \forall g \in H \quad g \widetilde a^\parallel = \widetilde a^\parallel \right\},\\
	\widetilde A^\perp \bydef \left\{\widetilde a^\perp \in \widetilde A\left| \sum_{	g \in H}g \widetilde a^\perp = 0 \right.\right\},\\	\widetilde A = \widetilde A^\parallel \oplus \widetilde A^\perp,\\
	\forall \widetilde x \in \widetilde\sX \quad \rep_{ \widetilde{x}}\left(\widetilde{A}^\parallel \right)=\rep_{ \widetilde{x}}\left(\pi\left( A\right)  \right),\quad \rep_{ \widetilde{x}}\left(\widetilde{A}^\perp \right)\cap\rep_{ \widetilde{x}}\left(\pi\left( A\right)  \right)= \{0\}
	\eean
	where $\oplus$ means a direct sum of $A$-$A$-bimodules. 
	From the Corollary \ref{ctr_ccr_i_cor} it turns out that $\widetilde A$ is a $C^*$-algebra is of type $I$. 	From the Theorem \ref{ctr_big_thm} it follows that  $\widetilde A$ contains an essential ideal $\widetilde I$ which is a continuous trace continuous trace $C^*$-algebra. From the Lemma \ref{ctr_ess_dense_lem} it follows that the spectrum $\widetilde\sX_{\widetilde I}$ of $\widetilde I$ is a dense open subset of $\widetilde\sX$  and $\rep_{ \widetilde{x}}\left(\widetilde{I}\right)=\rep_{ \widetilde{x}}\left(\widetilde{A}\right)$ for all $\widetilde x \in\widetilde\sX_{\widetilde I}$. 
	
	Suppose  $\widetilde x_0 \in \widetilde\sX_{\widetilde I}$ is such that $\rep_{  \widetilde{x }_0}\left( \pi\left(A \right) \right)\subsetneqq\rep_{  \widetilde{x }_0}\left(\widetilde{A} \right)$. Then $H$ non trivially acts of  $\rep_{  \widetilde{x }_0}\left(\widetilde{A} \right)$. According to  the construction \ref{ctr_imprim_empt} there is a compact neighborhood $\widetilde \sV$ of $\widetilde  x_0$ such that $\widetilde\sV\subset  \widetilde\sX_{\widetilde I}$ and  $\left.\widetilde A\right|^{\widetilde\sV}\left.\widetilde e\right|^{\widetilde\sV}$
	is an $\left.\widetilde A\right|^{\widetilde\sV}-\left. \widetilde e\right|^{\widetilde\sV}\left.A\right|^{\widetilde\sV}\left.\widetilde e\right|^{\widetilde\sV}\cong C\left({\widetilde\sV}\right)$-imprimitivity bimodule. If 	 we denote it by $X_{C\left( \widetilde\sV\right)}\bydef \left.\widetilde A\right|^{\widetilde\sV}\left.\widetilde e\right|^{\widetilde\sV}$ then from $\widetilde A = \widetilde A^\parallel \oplus \widetilde A^\perp$ it follows that three is a direct sum  $	X_{C\left( \widetilde\sV\right)}= 	X^\parallel_{C\left( \widetilde\sV\right)}\oplus 	X^\perp_{C\left( \widetilde\sV\right)}$ where
	\bean
	X^\parallel_{C\left( \widetilde\sV\right)}\bydef \left.\widetilde A^\parallel\right|^{\widetilde\sV}\left.\widetilde e\right|^{\widetilde\sV},\\
	X^\perp_{C\left( \widetilde\sV\right)}\bydef \left.\widetilde A^\perp \right|^{\widetilde\sV}\left.\widetilde e\right|^{\widetilde\sV}.
	\eean
	We can suppose that $\widetilde\sV$ is compact, so from the 
Theorems \ref{comp_normal_thm} 	and \ref{urysohn_lem} it follows that there is   an open neighborhood of $\widetilde x_0$ such that
\begin{itemize}
	\item $\widetilde{\sU}\subset \widetilde{\sV}$,
	\item there is a continuous map  $\widetilde f:  \to \left[0, 1\right]$ with $\widetilde f\left(\widetilde \sU \right) = \{1\}$ and $\widetilde f\left(\widetilde \sX\setminus \widetilde\sV\right) = \{0\}$,
	\item there is a continuous map  $\varphi: \widetilde\sX \to [0,1]$ such that $\varphi\left(\widetilde x_0  \right) = 1$.
\end{itemize}
 For any $\epsilon \in \left(0, 2\pi \right)$ denote by  
	\bean
	\varphi_{\epsilon}\in C_b\left( \widetilde\sX\right) ,\\
	\widetilde x \mapsto e^{i\eps \varphi\left( \widetilde x\right) }.
	\eean	
	There is the automorphism $\chi_\varepsilon\in \End^*_{C_0\left((\widetilde\sV \right) }\left(X_{C_0\left(\widetilde\sV \right)} \right) $ of the Hilbert $C_0\left(\widetilde\sV \right)$-module $X_{C_0\left(\widetilde\sV \right)}$ given by $\widetilde\xi^\parallel + \widetilde\xi^\perp\mapsto \widetilde\xi^\parallel + \varphi_{\epsilon}\left( \widetilde\xi^\perp\right) $. This automorphism yields the *-automorphism $\psi_\epsilon \in \Aut\left(\K\left(X_{C_0\left(\widetilde\sV \right)} \right)\right)\cong \Aut\left( \left.\widetilde A\right|^{\widetilde\sV}\right) $ given by
	$$
	\sum_{j = 0}^\infty \widetilde\xi_j \left\rangle \right\langle\widetilde \eta_j \mapsto \sum_{j = 0}^\infty \chi_\varepsilon\left( \widetilde\xi_j\right)  \left\rangle \right\langle \chi_\varepsilon\left( \widetilde\eta_j\right).
	$$
There the natural inclusion $\left.\widetilde A\right|_{\widetilde\sU}\subset \left.\widetilde A\right|^{\widetilde\sV}$ such that for all $\epsilon$ the *-automorphism  $\psi_{\epsilon}$ corresponds to the *-automorphism $\theta_\epsilon \in \Aut\left(\left.\widetilde A\right|_{\widetilde\sU}\right)$. Moreover 
	if $\epsilon_1\neq \epsilon_2$ then $\th_{\epsilon_1}\neq \th_{\epsilon_2}$.  For any  $\epsilon \in \left(0, 2\pi \right)$ there is a bijective $\C$-linear map $g_\epsilon: \widetilde A\xrightarrow{\approx} \widetilde A$ given by
	$$
	\widetilde a \mapsto \left(1-\widetilde f \right) \widetilde a + \th_{\epsilon}\left(\widetilde f\widetilde a \right) .
	$$
	Since for any $\widetilde x \in \widetilde\sX$ the map 
	\be\label{ctr_inf_eqn}
	\rep_{ \widetilde{x}}\left(\widetilde A\right)  \to \rep_{ \widetilde{x}}\left(\widetilde A\right),\\
	\rep_{ \widetilde{x}}\left(\widetilde a\right)  \to \rep_{ \widetilde{x}}\left(g_{\epsilon}\left( \widetilde a\right) \right)\quad\forall  \widetilde a\in  \widetilde A
	\ee
	is *-automorphism, the  map $g_\epsilon$ is $*$-automorphism. From 
	$$
	\forall\epsilon_1,\epsilon_2 \in \left[0, 2\pi \right)\quad  \epsilon_1,\neq\epsilon_2\quad \Rightarrow\quad \varphi_{\epsilon_1}\left(\widetilde x_0 \right) \neq \varphi_{\epsilon_2}\left(\widetilde x_0 \right) 
	$$
	and taking into account \eqref{ctr_inf_eqn} one has infinitely many different automorphisms $g_\eps$, so the group $H$ is not finite, so there is a contradiction with the Definition \ref{fin_pre_defn}. From this contradiction it follows that the hypotheses about the existence of $\widetilde x_0 \in \widetilde\sX_{\widetilde I}$ such that $\rep_{  \widetilde{x }_0}\left( \pi\left(A \right) \right)\subsetneqq\rep_{  \widetilde{x }_0}\left(\widetilde{A} \right)$ is not right, i.e.
	$$
	\forall \widetilde x_0 \in \widetilde\sX_{\widetilde I}\quad \rep_{  \widetilde{x }_0}\left( \pi\left(A \right) \right)=\rep_{  \widetilde{x }_0}\left(\widetilde{A} \right).
	$$
	If $\widetilde x' \in \widetilde\sX\setminus \widetilde\sX_{\widetilde I}$ is such that $\rep_{  \widetilde{x }}\left( \pi\left(A \right) \right)\subsetneqq\rep_{  \widetilde{x }}\left(\widetilde{A} \right)$ then there is $\widetilde a \in \widetilde a$ and  $g \in H$  such that $\rep_{  \widetilde{x }'}\left(\widetilde{a} \right)\neq \rep_{  \widetilde{x }'}\left(g\widetilde{a} \right)$. So $\left\| \rep_{  \widetilde{x }'}\left(\widetilde{a} \right)-g \widetilde{a}\right\| > 0$.
	 Otherwise 
	$$
	\forall \widetilde x_0 \in \widetilde\sX_{\widetilde I}\quad	\rep_{  \widetilde{x }}\left(\widetilde{a} \right)=  \rep_{  \widetilde{x }'}\left(g\widetilde{a} \right)\quad \Rightarrow\quad \left\| \rep_{  \widetilde{x }_0}\left(\widetilde{a}-g \widetilde a \right)\right\| = 0,
	$$
	Taking into account that $\widetilde\sX_{\widetilde I}$ is an open subset of $\widetilde\sX$ we conclude that the map
	\bean
	\widetilde\sX \to \R,\\
	\widetilde x \mapsto \left\| \rep_{  \widetilde{x }}\left(\widetilde{a}-g \widetilde a \right)\right\|
	\eean
	is not lower semi-continuous (cf. Definition \ref{top_lower_semi_defn}), but this fact contradicts with the Theorem  \ref{top_lower_semi_thm}. From this contradiction it follows that 
$$
	\forall \widetilde x \in \widetilde\sX \rep_{  \widetilde{x }_0}\left( \pi\left(A \right) \right)=\rep_{  \widetilde{x }_0}\left(\widetilde{A} \right)
$$
and $H$ is a trivial group so $\left(A,\widetilde{A}, G, \pi    \right)\cong 	\left(A,\widetilde{A}^H, G/H, \pi^H    \right)
$. 
On the other hand from the equation \eqref{ctr_comm_e_eqn} it follows that 
\bean
\left(A ,A_0\left( \widetilde\sX\right), G\left(\left.\widetilde{\sX}~\right|\sX\right),A_0\left(p \right) \right)\cong 	\left(A,\widetilde{A}^H, G/H, \pi^H    \right)\cong \left(A,\widetilde{A}, G, \pi    \right).
\eean
\end{proof}
\begin{lemma}\label{ctr_necessary_lem} 
	Let $A$ be a connected   $C^*$-algebra  with Hausdorff spectrum, locally compact, locally connected spectrum $\sX$, and suppose that $A$ is $CCR$ (cf. Definition \ref{ccr_defn}). If $\left(A,\widetilde{A}, G, \pi    \right)$ is a noncommutative finite-covering (cf. Definition \ref{fin_defn}) and $p: \widetilde \sX \to \sX$ is a given by the Proposition \ref{spectrum_covering_finite_prop} map form the spectrum $\widetilde \sX$ to the spectrum of $A$, then one has:
	\begin{enumerate}
		\item[(i)] the space $\widetilde \sX$ is Hausdorff,'
		\item[(ii)] the map $p: \widetilde \sX \to \sX$ is a transitive covering,
		\item[(iii)] the quadruple $\left(A,\widetilde{A}, G, \pi    \right)$ is equivalent to the given by the Lemma \ref{ctr_sufficient_lem} noncommutative finite-fold covering 
		$
		\left(A, A_0\left( \widetilde{\sX }\right) , G\left(\left.\widetilde\sX\right|\sX\right),A_0\left(p\right) \right)
		$.
	\end{enumerate}
\end{lemma}
\begin{proof}
If  $\left\{A_\la \subset A\right\}_{\la\in\La}$ be a required by the Definition \ref{fin_defn} family of  $\left(A,\widetilde{A}, G, \pi    \right)$-{strictly proper} $C^*$-subalgebras (cf. the Definition \ref{strictly_proper_defn}) and $\sX_\la$ is a spectrum of $A_\la$ then from the Proposition  \ref{hered_spectrum_prop} it follows that there is an inclusion $\sX_\la \hookto \sX$. If $x_0 \in \sX_\la$ and  $a \in A_\la$ is such that $\left\| \rep_{x_0}\left( a\right)\right\|= 2$ then from the Theorem  \ref{lower_norm_thm} it follows that the set
$$
\left\{x \in \sX_\la |\left\| \rep_x\left( a\right)\right\| > 1  \right\}
$$ 
is open in both $\sX_\la$ and $\sX$, so the set $\sX_\la$ is open in $\sX$. If $\sX \setminus \bigcup_{\la\in \La}\sX_\la \neq \emptyset$ then there is $x_0 \in \sX \setminus \bigcup_{\la\in \La}\sX_\la$ and $a \in A$ such that $\left\| \rep_{x_0}\left( a\right)\right\|= 1$. It turns out that
$$
\forall a' \in \bigcup_{\la\in \La} A_\la \quad \left\| a -a'\right\|\ge 1,
$$ 
i.e. the intersection $ \bigcup_{\la\in \La} A_\la$ is not dense in $A$. It contradicts with the Definition \ref{fin_defn}, so $\sX \setminus \bigcup_{\la\in \La}\sX_\la = \emptyset$ and $\sX = \bigcup_{\la\in \La}\sX_\la$. If $\widetilde A_\la$ is a {hereditary} $\left(A, \widetilde{A}, G, \pi \right)$-{lift} of $A_\la$ (cf. Definition \ref{hereditary_lift_defn}) then from the Definition \ref{strictly_proper_defn} it follows that there is a noncommutative finite-fold covering with unitization (cf. Definition \ref{fin_unitization_defn}) $\left(A_\la,\widetilde{A}_\la, G, \pi|_{A_\la}    \right)$. Similarly to above proof the spectrum $\widetilde\sX_\la$ of $\widetilde A_\la$ is an open subset of $\widetilde \sX$. Applying the Lemma \ref{fin_def_lem} we conclude that a union $\cup_{\la \in \La} \widetilde{A}_\la$ is dense in $\widetilde{A}$, and similarly to the above proof we conclude that $\widetilde \sX= \bigcup_{\la\in \La}\widetilde\sX_\la$.\\
(i) From the Lemma  \ref{ctr_fin_h_lem} it follows that for all $\la \in \La$ the set $\widetilde \sX_\la$ is Hausdorff. Applying the Theorem \ref{zorn_thm} we can prove that the space  $\widetilde \sX$ is Hausdorff.\\
(ii) For all $\la \in \La$ the restriction $p_{ \widetilde \sX_\la}: \widetilde \sX_\la\to \sX_\la$ is a transitive finite-fold covering (cf. Lemma \ref{ctr_fin_h_lem}), so from $ \sX= \bigcup_{\la\in \La}\sX_\la$ and $\widetilde \sX= \bigcup_{\la\in \La}\widetilde\sX_\la$ we conclude that the map $p : \widetilde \sX \to  \sX$  is a transitive finite-fold covering.\\
(iii) For any $\la \in \La$ there is an inclusion $A_\la \subset A|_{\sX_\la}=A_{\sX_\la}=~_{\sX_\la} A= ~_{\sX_\la} A_{\sX_\la}$ (cf. Definition \ref{blowing_ideals_au_ua_defn}). From this fact one can deduce that $\widetilde A_\la \left.\subset A_0\left( \widetilde{\sX }\right)\right|_{\widetilde \sX_\la}\subset A_0\left( \widetilde{\sX }\right)$, and $\pi|_{A_\la} = \left.A_0\left(p\right)\right|_{A_\la}$. If $\widetilde a \in A_0\left( \widetilde \sX\right)$ then for each $\eps > 0$ there is a compact set  $\widetilde \sY$ such that
$$
\forall\widetilde x \in \widetilde \sX \setminus \widetilde\sY \quad  \left\| \rep_{\widetilde x}\left(\widetilde a\right)\right\| < \frac{\eps}{2}.
$$
There is a {covering sum} $\sum_{j = 1}^n \widetilde f_j$ for $\widetilde\sY$ (cf. Definition \ref{top_covering_sum_defn}) such that the set
$$
\widetilde\sU_j \bydef \left\{\left.\widetilde x \in \widetilde\sX\right| \widetilde f\left( \widetilde x\right)\neq 0 \right\}
$$
is mapped homeomorphically onto $p\left(\widetilde\sU_j  \right)$.  If $a_j \bydef \pi^{-1}\left(\sum_{	g \in G}g \left( \sqrt{\widetilde f_j} \widetilde a \right)  \right)$ then there is $\la_j \in \La$ and $a_{\la_j}\in A_{\la_j}$ such that
$$
\left\| a_{\la_j}- a_j\right\|< \frac{\eps}{2n}.
$$
From our construction it follows that
\be\label{ctr_ee_eqn}
\left\| \widetilde a - \sum_{j=1}^n\pi\left( a_{\la_j}\right) \sqrt{\widetilde f_j}  \right\|< \eps.
\ee
If $\la \in \La$ is such that $\la \ge \la_j$ for every $j = 1,...,n$ then $\sum_{j=1}^n\pi\left( a_{\la_j}\right) \sqrt{\widetilde f_j}\in \widetilde A_\la$. From \eqref{ctr_ee_eqn} it follows that the union $\bigcup_{\la\in \La} \widetilde A_\la$ is dense in $ A_0\left( \widetilde{\sX }\right)$ so $\widetilde A \cong A_0\left( \widetilde{\sX }\right)$ (cf. Lemma \ref{fin_def_lem}). From $\pi|_{A_\la} = \left.A_0\left(p\right)\right|_{A_\la}$ and taking into account that the union $\bigcup_{\la\in \La}A_\la$ is dense in $A$ we conclude that $\pi =A_0\left(p\right)$.
\end{proof}

\begin{thm}\label{ctr_fin_thm}
	Let $A$ be a connected   $C^*$-algebra  with Hausdorff, locally compact, locally connected spectrum $\sX$, and suppose that $A$ is $CCR$ (cf. Definition \ref{ccr_defn}). The  given by the Definition \ref{blowing_functor_defn} $\sX$-$A$-functor $A_0$ is an equivalence between the {category of finite-fold coverings} of $\sX$ (cf. Definition \ref{top_fin_cov_defn}) and the {category of finite-fold coverings} of $A$ (cf. Definition \ref{fin_category_defn}), i.e.
	$$
	\mathfrak{FinCov}\text{-}\sX \xrightarrow[\cong]{A_0} \mathfrak{FinCov}\text{-}A. 
	$$
\end{thm}
\begin{proof}
	Follows from the Lemmas \ref{ctr_sufficient_lem} and \ref{ctr_necessary_lem}.
\end{proof}

\subsection{Coverings of operator spaces}
\begin{empt}\label{ctr_fin_op_empt}
Let $\left(X, Y\right)$ be a  {sub-unital operator space} with Hausdorff spectrum (cf. Definition \ref{ctr_oa_defn}) and let $A\bydef C^*_e\left(X, Y \right)$ be the $C^*$-envelope of $\left(X, Y\right)$ (cf. Definition \ref{operator_space_envelope_defn}). Suppose that $A$ is a $CCR$-algebra (cf. Definition \ref{ccr_defn}), hence the spectrum $\sX$ of $A$ is the locally compact Hausdorff space. Let $p: \widetilde \sX \to \sX$ be a transitive covering, denote by
\be\label{ctr_fin_op_eqn}
\begin{split}
\widetilde X \stackrel{\text{def}}{=} C_0\left(\lift_p\left[X\right]\right),\\
\widetilde A \stackrel{\text{def}}{=}  C_0\left(\lift_p\left[A\right]\right)= A_0\left(\widetilde \sX \right),\\
\widetilde Y \stackrel{\text{def}}{=} \widetilde X + \C \cdot 1_{\widetilde A^+}.
\end{split}
\ee
(cf. equations \ref{top_c0_eqn}, \ref{top_cont_s_lift_eqn}).
If the covering $p$ is finite-fold then there are the natural injective $*$-homomorphism $\rho: A \to \widetilde A$ and the complete unital  isometry $\left(\pi_X, \pi_Y \right):  \left(X, Y \right) \hookto \left(\widetilde X, \widetilde Y \right)$ from $ \left(X, Y \right)$ to  $\left(\widetilde X, \widetilde Y \right)$ (cf. Definition \ref{op_sum_space_defn}).
\end{empt}

\begin{theorem}\label{ctr_fin_op_thm}
Let us consider the given by \ref{ctr_fin_op_empt} situation. If $\sX$ is a connected, locally connected, then there is the natural 1-1 map from the set of finite-fold transitive coverings $p: \widetilde\sX \to\sX$ 
\eqref{fin_op_defn} and noncommutative  finite-fold  coverings of the sub-unital operator space $\left(X, Y\right)$ given by
\bea\label{ctr_from_oa_eqn}
\left(p: \widetilde\sX \to\sX \right) \mapsto \left(\left(X, Y \right),\left(\widetilde X, \widetilde Y \right), G\left(\left.\widetilde\sX \right|\sX\right), \left(\pi_X, \pi_Y \right) \right).
\eea 
\end{theorem}
\begin{proof}
	Firstly we proof that $\left(\left(X, Y \right),\left(\widetilde X, \widetilde Y \right), G\left(\left.\widetilde\sX \right|\sX\right), \left(\pi_X, \pi_Y \right) \right)$ is a finite-fold  covering  for each transitive finite-fold covering $p: \widetilde\sX \to\sX$. One needs check (a) and (b) of the Definition \ref{fin_op_defn}.\\
(a)	From the Theorem \ref{ctr_fin_thm} it turns out that
	$$
	\left( C^*_e\left( X,  Y \right), \widetilde{A}, G\left(\left.\widetilde\sX \right|\sX\right), C_0\left(p \right)  \right)
	$$  a finite-fold  noncommutative covering. From $X_x \subset A_x$ for any $x \in \sX$, $\widetilde X_{\widetilde x}\cong X_{p\left( \widetilde x\right) }$ and $\widetilde A_{\widetilde x}\cong A_{p\left( \widetilde x\right) }$ it follows that $C_0\left(\widetilde X_{\widetilde x}, \widetilde Y_{\widetilde x} \right) = \rep_{ \widetilde{x}}\left(\widetilde{A}\right)$ for all $\widetilde{x}\in \widetilde{\sX}$ (cf. notation of the Definition \ref{ctr_oa_defn}).
	Applying  the Lemma \ref{ctr_oa_lem} one has $\widetilde A = C^*_e\left(\widetilde X, \widetilde Y \right)$, clearly the following condition holds $\pi_Y = \left.\rho\right|_Y$, $\pi_X = \left.\rho\right|_X$.\\
(b)
If $\widetilde X' \subset C^*_e\left(\widetilde X, \widetilde Y \right)$ is a $\C$-linear space  such that $X = \widetilde X' \cap  C^*_e\left( X,  Y \right)$ and $G\widetilde X'= \widetilde X'$ then from $\widetilde X_{ \widetilde x }\cong X_{p\left( \widetilde x\right) }$ it follows that $\widetilde X'_{ \widetilde x }\subset \widetilde X_{ \widetilde x }$ for all $ \widetilde x \in \widetilde \sX$.  From the Lemma \ref{top_sub_eq_lem} it follows that $\widetilde X' \subset\widetilde X$.

Secondly  we proof that any noncommutative finite-fold  covering $$\left(\left(X, Y \right),\left(\widetilde X, \widetilde Y \right), G, \left(\pi_X, \pi_Y \right) \right)$$  of the sub-unital operator space $\left(X, Y\right)$ yields  the transitive finite-fold covering $p: \widetilde\sX \to\sX$. From  (a) of the Definition \ref{fin_op_defn} it follows that $$\left(\left(X, Y \right),\left(\widetilde X, \widetilde Y \right), G, \left(\pi_X, \pi_Y \right) \right)$$ corresponds to noncommutative covering $\left( C^*_e\left( X,  Y \right), C^*_e\left(\widetilde X, \widetilde Y \right), G, \rho \right)$ of $C^*$-algebras. However  $C^*_e\left( X,  Y \right)$ a $CCR$-algebra (cf. Definition \ref{ccr_defn}), so and taking into account the Theorem \ref{ctr_fin_thm} on can find required covering $p: \widetilde \sX\to \sX$.
\end{proof}

\subsection{Induced representations}\label{ctr_induced_finite_sec}
\paragraph*{}
Let $A$ be a $C^*$-algebra with Hausdorff spectrum $\sX$. If $\tau : C_c\left(\sX \right) \to \C$ be a functional such that $\tau\left(a \right) > 0$ for any positive $a$ there is a functional 
\be\label{ctr_ka_t_eqn}
\begin{split}
	\tr_\tau : K\left( A\right)\to \C,\\
	a \mapsto \tau\left(\tr_a \right).  
\end{split} 
\ee
where $K\left( A\right)$ is Pedersen's ideal (cf. Definition \ref{pedersen_ideal_defn}) and  $\tr\left( a\right)  \bydef \hat a$ is given by \eqref{ctr_hat_eqn}. Note that $\tr_a\in C_c\left(\sX \right)$ for all  $a \in K\left( A\right)$. There is a $\C$-valued product on $ K\left( A\right)$ given by
\be\label{ctr_k_p_eqn}
\begin{split}
	\left(\cdot, \cdot\right): K\left( A\right)\times K\left( A\right)\to \C,\\
	\left(a,b\right)\mapsto \tr_\tau\left(a^*b\right),
\end{split}
\ee
so $K\left( A\right)$ becomes a pre-Hilbert space. If  $L^2\left(A, \tr_\tau\right)$ is the Hilbert norm completion of $K\left( A\right)$, then the action $A \times K\left(A \right)\to K\left( A\right)$ induces an action $A \times K\left(A \right)\to K\left( A\right)$ which corresponds to a faithful representation
\be\label{ctr_a_eqn}
\rho_\tau: A \hookto B\left(L^2\left(A, \tr_\tau\right) \right). 
\ee   
If  $\left(A, \widetilde{{A}}, G, \pi \right)$ is a connected finite-fold noncommutative covering  (cf. Definition \ref{fin_defn}), and $\widetilde \sX$ is the spectrum of $\widetilde A$ then from the Theorem \ref{ctr_fin_thm} it follows that $\widetilde \sX$ is Hausdorff, the given by the Proposition \ref{spectrum_covering_finite_prop} map $p:\widetilde \sX\to  \sX$ such that $\left(A, \widetilde{{A}}, G, \pi \right)$ is equivalent to  $\left(A, A_0\left(\widetilde \sX \right) , G\left(\left.\widetilde \sX \right|\sX \right) , A_0\left( p\right) \right)$.
Similarly to the equation \eqref{top_wt_eqn} one has a functional
\be\label{ctr_wt_eqn}
\begin{split}
	\widetilde	\tau : C_c\left(\widetilde \sX \right)\to \C.
\end{split}
\ee
Now one can construct a faithful representation 
\be\label{ctr_wa_eqn}
\rho_{\widetilde\tau}: \widetilde A \hookto B\left(L^2\left(\widetilde A, \tr_{\widetilde\tau}\right) \right). 
\ee
which is an analog of the representation \ref{ctr_a_eqn}.
\begin{exercise}
	Prove that $\rho_{\widetilde\tau}:\widetilde A \hookto B\left(L^2\left(\widetilde A, \tr_{\widetilde\tau}\right) \right)$  {is induced by the pair} $\left(\rho_{\tau},\left(A, \widetilde{A}, G, \pi\right)  \right)$ (cf. Definition \ref{induced_repr_fin_defn}).
\end{exercise}
\begin{exercise}\label{ctr_hilbert_exer}
	Prove following statements:
	\begin{enumerate}
		\item There are natural representations
		\be
		\begin{split}
			\pi: C^*\text{-}\varinjlim_{\la\in\La} A_0\left(\sX_\la \right)\hookto B\left(L^2\left(A_0\left( \overline\sX\right), \tr_{\overline\tau}  \right) \right) ,\\
			\overline \pi: A_0\left( \overline\sX\right)\hookto B\left(L^2\left( A_0\left( \overline\sX\right), \tr_{\overline\tau}  \right) \right).
		\end{split}
		\ee	
		\item The representation $\pi$ is equivariant (cf. Definition \ref{equivariant_representation_defn}).
		\item Suppose that  $\mathcal U\subset \mathcal X$  is a connected open subset evenly {covered} by $\overline{\mathcal U}\subset \overline{\mathcal X}$. If $\overline a \in A_0\left(\overline\sX \right)$ is  an  element such that 
		$\supp \overline a \subset \overline \sU$ then the following conditions hold:
		\begin{enumerate}
			\item [(i)] The series 	
			\be\nonumber
			\sum_{g \in G\left(\left.\overline{\sX}~\right|\sX\right) } \overline \pi\left( g\overline{a}\right)
			\ee 
			is convergent in the strong  topology of $B\left(L^2\left(A_0\left( \overline\sX\right), \tr_{\overline\tau}  \right) \right)$  (cf. Definition \ref{strong_topology_defn}).
			\item[(ii)] 
			\be\nonumber
			\begin{split}
				\sum_{g \in G\left(\left.\overline{\sX}~\right|\sX\right) } \overline \pi\left( g\overline{a}\right)  = \pi\left(  \desc_{\overline p}\left( \overline{a}\right)\right)
			\end{split}
			\ee
			where $\desc_{\overline p}$ means the descent (cf. Definition \ref{top_lift_desc_defn}).
		\end{enumerate}
		\item A pointed algebraical  finite covering category $\mathfrak{S}^{\mathrm{pt}}_{A }$ is $\pi$-good (cf. Definition \ref{good_defn}).
	\end{enumerate}
	\end{exercise}

\begin{empt}\label{str_spin_c_empt}
	Suppose that $A$ is a {homogeneous of order} $n$ $C^*$-algebra (cf. Definition \ref{ctr_homo_defn}), and $E \to \sX$ is a   vector bundle with a sesquilinear form (cf. Definition \ref{top_herm_bundle_form_defn}). Assume that for any $x \in \sX$ there is an action
	\be\label{top_repxs_eqn}
	\rep_x\left(A\right)\times E_x\to E_x
	\ee
	where $\rep_x$ is given by \eqref{rep_x_eqn}, and $E_x$ is the fiber of $E$ at the point $x$ (cf. Definition \ref{top_vb_fiber_defn}). Suppose that $E_x$ is a Hermitian  $\rep_x\left(A\right)$-module for all $x\in E_x$ (cf. Definition \ref{herm_mod_defn}). Assume that a family $$\left\{\rep_x\left(A\right)\times E_x\to E_x\right\}_{x\in\sX}$$ induces an action
	$A \times \Ga_c\left(\sX, E \right)\to \Ga_c\left(\sX, E\right)$. 
	There is a $\C$ valued scalar product
	\be\label{ctr_gac_prod}
	\begin{split}
		\left(\cdot, \cdot \right): \Ga_c\left(\sX, E \right)\times \Ga_c\left(\sX, E \right)\to \C,\\
		\left(\xi, \eta \right)\bydef \tau\left(\left\langle \xi, \eta \right\rangle_c \right) 
	\end{split}
	\ee
	where $	\left\langle \cdot, \cdot \right\rangle_c$ is given by the equation \eqref{top_ggc_eqn}. It follows that $\Ga_c\left(\sX, E \right)$ is a pre-Hilbert space, denote by $L^2\left(\sX, E, \tau \right)$ (or shortly $L^2\left(\sX, E \right)$) the Hilbert norm completion of $\Ga_c\left(\sX, E \right)$. An action $A \times \Ga_c\left(\sX, E \right)\to \Ga_c\left(\sX, S \right)$ induces an action $A \times L^2\left(\sX, S \right)\to L^2\left(\sX, E \right)$.
\end{empt}

\begin{exercise}
	Prove that the action $A \times L^2\left(\sX, S \right)\to L^2\left(\sX, S \right)$  corresponds to a faithful *-representation
	\be\label{ctr_rho_eqn}
	\phi_\tau: A \to B\left(L^2\left(\sX, S \right) \right).
	\ee
\end{exercise}
\begin{exercise}\label{ctr_bounded_exer}
	Let $\widetilde S$   be the {inverse image} of $S$ by $p$ (cf. Definition \ref{vb_inv_img_funct_defn}).
	Define a faithful representation 
	\be\label{ctr_wrho_eqn}
	\phi_{\widetilde\tau}: \widetilde A \to  B\left(L^2\left(\widetilde\sX, \widetilde S \right) \right).
	\ee
	Prove that $\phi_{\widetilde\tau}:\widetilde A \hookto B\left(L^2\left(\widetilde \sX, \widetilde S\right) \right)$  {is induced by the pair} $\left(\phi_{\tau},\left(A, \widetilde{A}, G, \pi\right)  \right)$ (cf. Definition \ref{induced_repr_fin_defn}).	
\end{exercise}

\subsection{Coverings of $*$-algebras}
\subsubsection{Coverings of bounded operator $*$-algebras}
\paragraph{}
Let $A$ be a  $CCR$-algebra such that the spectrum $\sX$ of $A$ is locally connected, locally compact, Hausdorff space. Let $R\subset C_0\left(\sX\right)$ be a $c$-soft *-subalgebra (cf. Definition \ref{top_soft_r_defn}). If $\A\subset A$ is a dense *-subalgebra which is an $R$-module then $\A$-sheaf $\mathscr S^\A$ (cf. Definition \ref{top_x_sheaf_defn}) is $c$-soft (cf. Theorem  \ref{sheaf_soft_m_thm}). If $p : \widetilde \sX \to \sX$ is a finite-fold transitive covering then there from the Lemma \ref{ctr_sufficient_lem}   it follows that there is a noncommutative finite-fold covering 
$$
\left(A, A_0\left( \widetilde{\sX }\right) , G\left(\left.\widetilde\sX\right|\sX\right),A_0\left(p\right) \right)
$$
where $A_0\left( \widetilde{\sX }\right)$ is given by the Definition  \ref{top_lift_a_f_defn} $p$-lift of $A$, and $*$-homomorphism  $A_0\left(p\right)$ is given by the equation \eqref{top_c0p_ob_eqn} and/or \eqref{blowing_lift_fin_eqn} one. If $p^{-1}\mathscr S^\A$ is an inverse image sheaf (cf. Definition \ref{sheaf_inv_im_defn}) then an intersection
\be\label{ctr_csoft_lift_eqn}
\A_0\left( \widetilde \sX \right)  \bydef  A_0\left(\widetilde \sX\right) \cap p^{-1}\mathscr S^\A \left(\widetilde \sX\right)
\ee
is a dense *-subalgebra of $A_0\left( \widetilde{\sX }\right)$.
\begin{exercise}
Similarly to the Theorem \ref{top_oa_cov_thm} prove that a quadruple
$$
\left(\A, \A_0\left( \widetilde{\sX }\right) , G\left(\left.\widetilde\sX\right|\sX\right),\left.A_0\left(p\right)\right|_{\A} \right)
$$
is a	finite-fold covering of bounded operator *-algebras  (cf. Definition \ref{fin_oa_defn}).
\end{exercise}

\subsubsection{Coverings of $O^*$-algebras}\label{ctr_fin_o_sec}
\paragraph*{}

Results of this section a closed to the  Section \ref{top_pd_fin_sec}. Exercises of this Section can be solved similarly to proofs of Lemmas and Theorems of the Section \ref{top_pd_fin_sec}.
\begin{empt}\label{ctr_fin_o_empt}
	Let $M$ be a smooth manifold, and let $A$ be  a {homogeneous of order} $n$ $C^*$-algebra (cf. Definition \ref{ctr_homo_defn}). If $M$ is a spectrum of $A$ than from the Theorem  \ref{ctr_as_field_thm} it follows that  $A$ corresponds to a vector bundle $\left( E, q, M\right)$ (cf. Definition \ref{top_vb_defn}), such that any fiber is isomorphic to $\C^{n^2}\cong \mathbb{M}_n\left(\C\right)$, and there is the natural isomorphism $A \cong \Ga_0\left(M, E \right)$ of $C_b\left(M \right)$ modules 
	.   Suppose that $E$ has a structure of smooth manifold such that the map $q: E \to M$ is smooth. Any fiber of $E$ is isomorphic to $\mathbb{M}_n\left(\C \right)$  so similarly to the Lemma \ref{op_cont_module_lem} one can proof that the space of smooth sections $\Ga^\infty\left(M, E \right)$ is *-algebra.	
	If $\Ga_0^\infty\left(M, E \right)$ is a space of  smooth vanishing at infinity sections of  $\left( E, q,  M\right)$ then 
	 $\Ga^\infty_0\left(M, E \right)$ is a *-subalgebra of $\Ga^\infty\left(M, E \right)$.
	 If $p: \widetilde M \to M$ is a transitive finite-fold  covering then  from the Proposition \ref{top_cov_mani_prop} it follows that $\widetilde M$ has a (unique) structure of smooth manifold such that $p$ is smooth. 
	On the other hand  from 
	the Theorem  \ref{ctr_fin_thm} it follows that \\ $\left(A, A_0\left(\widetilde M \right) , G\left(\left.\widetilde M \right|M \right) , A_0\left( p\right) \right)$ a connected noncommutative finite-fold covering (cf. Definition \ref{fin_defn}). If $\left(\widetilde E, \widetilde q,  \widetilde M\right)$  is the {inverse image} of $\left( E, q, M\right)$ by $p$ (cf. Definition \ref{vb_inv_img_funct_defn}) then there is a transitive finite-fold covering  $\widetilde E \to E$  so from the Proposition \ref{top_cov_mani_prop} it turns out that there is the unique smooth structure on $\widetilde E$ such that map $\widetilde E\to E$ is smooth. 
	Similarly to the above construction one has a space of smooth sections of  $\left(\widetilde E,  \widetilde q,\widetilde  M\right)$ then $\Ga^\infty\left(\widetilde M, \widetilde E \right)$ is an   *-algebra. Since smooth sections of $E$ are mapped to smooth sections of $\widetilde E$ one has $*$-homomorphism $\pi^\infty: \Ga^\infty\left( M,  E \right)\hookto\Ga^\infty\left(\widetilde M, \widetilde E \right)$. Let $\tau: C_c\left(M \right)\to \C$ be a discussed in \ref{ctr_induced_finite_sec} functional, let both \be\label{ctr_sa_eqn}
	\begin{split}
		\rho_\tau: A \hookto B\left(L^2\left(A, \tr_\tau\right) \right),\\ 
		\rho_{\widetilde\tau}: \widetilde A \hookto B\left(L^2\left(\widetilde A, \tr_{\widetilde\tau}\right) \right). 
	\end{split}
	\ee
	are given by equations \eqref{ctr_a_eqn} and \eqref{ctr_wa_eqn} representations. If 
	\bean
	\D \bydef \Ga^\infty_c\left( M,  E \right),\\
	\widetilde \D \bydef \Ga^\infty_c\left(\widetilde M, \widetilde E \right)
	\eean
	then there are dense inclusions $ \D \subset  L^2\left( A, \tr_{\tau}\right) $ and $\widetilde \D \subset  L^2\left(\widetilde A, \tr_{\widetilde\tau}\right)$ such that there are natural inclusions 
	\be\label{ctr_smooth_inc_eqn}
	\begin{split}
	\Ga^\infty\left( M,  E \right)\hookto	\L^\dagger\left( \D\right),\\
	\Ga^\infty\left(\widetilde M, \widetilde E \right) \hookto	\L^\dagger\left( \widetilde\D\right)
	\end{split}
	\ee
where both $\L^\dagger\left( \D\right)$ and $\L^\dagger\left( \widetilde\D\right)$ are given by \eqref{l_dag_eqn}).
\end{empt}{
	\begin{exercise}\label{ctr_smooth_exer}
		Using equations \eqref{ctr_smooth_inc_eqn} prove that a quadruple 
		$$
		\left(\Ga^\infty\left( M,  E \right), \Ga^\infty\left(\widetilde M, \widetilde E \right), G\left(\left.\widetilde M \right|M \right) , \pi^\infty \right)
		$$
		is a noncommutative finite-fold covering of $O^*$-algebras (cf. Definition \ref{fino*_defn}).
	\end{exercise}

	\begin{exercise}\label{ctr_smoothc_exer}
		Let both  $D\left( M,  E \right)$ and $D\left(\widetilde M, \widetilde E \right)$ be *-algebras of differential operators 
		$\Ga^\infty \left(M, E \right)  \to \Ga^\infty \left(M, E \right)$ and $\Ga^\infty \left(\widetilde M,  \widetilde E \right)  \to \Ga^\infty \left(\widetilde M, \widetilde E \right)$ respectively (cf. Definition \ref{do_man_defn}).
		Construct  a noncommutative finite-fold covering of $O^*$-algebras 	$\left(D\left(\widetilde M, \widetilde E \right) ,D\left(\widetilde M, \widetilde E \right)  , G\left(\left.\widetilde M \right|M \right) , \rho \right)$ such that $\pi^\infty= \left.\rho\right|_{\Ga^\infty\left( M,  E \right)}$.
	\end{exercise}

\begin{empt}\label{str_spinc_c_empt}
	Suppose that  $\left( E, r, M\right)$ is a  bundle with a sesquilinear form (cf. Definition \ref{top_herm_bundle_form_defn}) such that if   $E_x$ is the {fiber} of $\left( E, r, M\right)$ at the point $x$ (cf. Definition \ref{top_vb_fiber_defn}) 
there is an action
\be\label{foli_rep_eqn}
\rep_x\left(A \right)\times E_x\to E_x 
\ee
such that $E_x$ is a Hermitian  	$	\rep_x\left(A \right)$-module (cf. Definition \ref{herm_mod_defn})
Suppose that a given by \eqref{foli_rep_eqn} family of actions induces the following action
\be\label{foli_me_eqn}
\Ga\left( M, E\right) \times \Ga\left( M, E\right) \to \Ga\left( M, E\right) 
\ee
If $E$ has structure of smooth manifold then similarly to the above construction one can obtain spaces of smooth sections $\Ga_c^\infty\left(M, E \right)$ and $\Ga_c^\infty\left(\widetilde M, \widetilde E \right)$ which are dense in $L^2\left( M, E \right)$ and $L^2\left(\widetilde M, \widetilde E \right)$ respectively (cf. Equations \ref{ctr_rho_eqn} and \ref{ctr_wrho_eqn}). There is the natural inclusion $\Ga^\infty\left(M, E \right)\subset \Ga^\infty\left(\widetilde M, \widetilde E \right)$
If the action $\Ga\left( M,  E \right)\times\Ga\left( M,  E \right)\to  \Ga\left( M,  E \right)$ is such that 
\be\label{foli_sm_eqn}
\Ga^\infty\left( M,  E \right)\times\Ga^\infty\left( M,  E \right)\subset  \Ga^\infty\left( M,  E \right)
\ee
then one has
$$
\Ga^\infty\left(\widetilde M, \widetilde E \right)\times\Ga^\infty\left(\widetilde M, \widetilde E \right)\subset  \Ga^\infty\left(\widetilde M, \widetilde E \right).
$$
There are dense subspaces
\bean
\D_E \bydef  \Ga_c^\infty\left( M,  E \right)\subset  L^2\left( M,  E \right),\\
\widetilde \D_{\widetilde E} \bydef  \Ga_c^\infty\left(\widetilde M, \widetilde E \right)\subset  L^2\left(\widetilde M, \widetilde E \right).
\eean
\end{empt}

	\begin{exercise}\label{ctr_smoothdiff_exer}
 Let both  $D^*\left( M, E\right)$ and $D^*\left(\widetilde M, \widetilde E\right)$ are algebras of differential operators 
		$\Ga^\infty \left(M, E \right)  \to\Ga^\infty \left(M, E \right)$ and $\Ga^\infty \left(\widetilde M,  \widetilde E \right)  \to \Ga^\infty \left(\widetilde M, \widetilde E \right)$ respectively (cf. Definition \ref{do_man_defn}) such that
\bean
\forall X \in  D^*\left( M, E\right)\quad X^*\in D^*\left( M, E\right),\\
\forall \widetilde X\in  D^*\left( \widetilde M, \widetilde E\right)\quad  \widetilde X^*\in D^*\left(\widetilde M,\widetilde E\right),
\eean	
i.e.
\bean
\forall X \in  D^*\left( M, E\right)\quad X\in \L^\dagger\left( \Ga^\infty \left(\widetilde M, \widetilde E \right)\right) ,\\
\forall \widetilde X\in  D^*\left( \widetilde M, \widetilde E\right)\quad \widetilde X\in \L^\dagger\left( \Ga^\infty \left(\widetilde M, \widetilde E \right)\right) ,\\
\eean 	
(cf. equation \eqref{l_dag_eqn}).
	Using inclusions $D^*\left( M, E\right)\subset \L^\dagger\left( \Ga^\infty \left(\widetilde M, \widetilde E \right)\right)$ and $D^*\left(\widetilde M,\widetilde E\right)\subset \L^\dagger\left( \Ga^\infty \left(\widetilde M, \widetilde E \right)\right)$ 	prove that there is a noncommutative finite-fold covering of $O^*$-algebras 	$$\left(D^*\left( M, E\right) ,D^*\left(  \widetilde M, \widetilde E\right)  , G\left(\left.\widetilde M \right|M \right) , \phi \right).$$
	\end{exercise}

\subsection{Coverings of spectral triples}\label{ctr_d_sec}.

\paragraph*{}
Let $M$ be a Riemannian compact manifold, and 	let 
$
\left(\Coo\left( M\right), L^2\left(M, S \right), \Dslash, J_{\mathrm{comm}}   \right)
$
be a commutative spectral triple (cf. Equation \ref{comm_sp_tr_eqn}). 
Let $\left(\mathbb{M}_n\left( \C\right), \H_{\mathrm{fin}},  D_{\mathrm{fin}}, J_{\mathrm{fin}}\right)$ be an explained  in  the Section \ref{ctr_fin_sp_tr_sec} finite spectral triple. Note that $\H_{\mathrm{fin}}\cong \C^N$ is a finite-dimensional  Hermitian  $\mathbb{M}_n\left( \C\right)$-module and $D_{\mathrm{fin}}\in  B\left( \H_{\mathrm{fin}}\right)$. If $\left(\A, \H, D,J\right)$ is an explained  in the section \ref{sp_tr_prod_sec} product of spectral triples\\ $\left(\Coo\left( M\right), L^2\left(M, S \right), \Dslash, J_{\mathrm{comm}}   \right)$ and $\left(\mathbb{M}_N\left( \C\right), \H_{\mathrm{fin}},  D_{\mathrm{fin}}, J_{\mathrm{fin}}\right)$ then one has
\be\label{ctr_st_eqn}
\begin{split}
\A \bydef \Coo\left( M\right)\otimes \mathbb{M}_n\left( \C\right),\\
\H \bydef L^2\left(M, S \right)\otimes \C^N \cong L^2\left(M, S \right)^N,\\
D \bydef \Dslash\otimes \Ga_{\mathrm{fin}}+ \Id_{L^2\left(M, \sS \right)} \otimes D_{\mathrm{fin}},\\
J \bydef J_{\mathrm{comm}}\otimes J_{\mathrm{fin}}
\end{split}
\ee
where $\Ga_{\mathrm{fin}}\in B\left( \H_{\mathrm{fin}}\right)$ is the grading operator of $\left(\mathbb{M}_n\left( \C\right), \H_{\mathrm{fin}},  D_{\mathrm{fin}}\right)$ (cf. Definition \ref{df:spt-even}).
If $p: \widetilde M\to M$ is a finite-fold covering then from the Theorem \ref{ctr_fin_thm} it follows that the quadruple
$$
\mathfrak T\bydef\left(C\left(M \right)\otimes  \mathbb{M}_{n}\left( \C\right) , C_0\left(\widetilde M \right)\otimes  \mathbb{M}_{n}\left( \C\right) , G\left(\left.\widetilde M \right|M \right) , C_0\left( p\right) \otimes \Id_{ \mathbb{M}_{n}\left( \C\right) } \right).
$$
is a noncommutative finite-fold covering (cf. Definition \ref{fin_defn}). 
 \begin{exercise}
	Prove that following statements:
	\begin{itemize}
	\item 
	there is a spectral triple
	$$
	\left( \Coo\left( M\right)\otimes \mathbb{M}_n\left( \C\right), L^2\left(M, S \right)^N,  D,  J\right) ,
	$$
	\item if both  $\widetilde D \bydef \lift_p\left( D\right)$  and $\widetilde J \bydef \lift_p\left( J\right)$ are $p$-lifts of $D$ and $J$ (cf. Definition \ref{top_sheaf_lift_defn}) 
	then the spectral triple $\left( \Coo\left( \widetilde M\right)\otimes \mathbb{M}_n\left( \C\right), L^2\left(\widetilde M,  \widetilde S \right)^N,\widetilde D, \widetilde J\right) $ is a
$\mathfrak T$-lift of 	$\left( \Coo\left( M\right)\otimes \mathbb{M}_n\left( \C\right), L^2\left(M, S \right)^N, D, J\right)$ (cf. Definition \ref{spectral_triple_fin_lift_defn}). 
		\end{itemize}

\end{exercise} 
!
\subsection{Unoriented spectral triples}\label{ctr_sp_tr} 

Consider the described in the Section \ref{comm_sp_tr_sec} situation, i.e. one has:

\begin{itemize}
	\item the Riemannian manifold $\widetilde M$ with the spinor bundle $\widetilde\SS$ such that there is the spectral triple $\left(\Coo\left(\widetilde{M}  \right) , L^2\left(\widetilde{M},\widetilde{\sS} \right), \widetilde\Dslash , \widetilde J_{\mathrm{comm}}\right)$,
	\item the unoriented  Riemannian manifold $M$ with the two listed covering $p:\widetilde M\to M$ and the bundle $\SS \to M$ such that $\widetilde{\sS}$ is the $p$-inverse image of $\SS$.
	\item the unoriented  spectral triple $\left(\Coo\left({M}  \right) , L^2\left({M},{\sS} \right), \Dslash, J_{\mathrm{comm}} \right)$ given by the Exercise \ref{top_unori_exer}.
\end{itemize}

\begin{theorem}\label{ctr_unori_thm}
In the described above situation there is the given by
	\be\label{ctr_equ}
	\left(\Coo\left(M \right) \otimes\mathbb{M}_m\left(\C \right) , L^2\left({M}, {\SS}\right)^k, \slashed D \otimes \Ga_{\mathrm{fin}}+ \Id_{L^2\left( M, \sS \right)^k} \otimes D_{\mathrm{fin}}, J_{\mathrm{comm}}\otimes J_{\mathrm{fin}} \right).
	\ee
 unoriented  spectral triple (cf. Definition \ref{unoriented_defn}).
\end{theorem}
\begin{proof}
	The following table shows the specialization  of the construction  \ref{unoriented_empt}. 
	\\ \\
	\begin{tabular}{|c|c|c|}
		\hline
		&Definition \ref{unoriented_defn}& This theorem specialization\\ 
		\hline
		&	&\\
		1 & $\A$   &  $\Coo\left(M \right)\otimes\mathbb{M}_m\left(\C \right)$ \\ & & \\
		2	& $\left(A, \widetilde A, \Z_2\right)$ & $\left(C\left(M \right)\otimes\mathbb{M}_m\left(\C \right), C\left(\widetilde M\right)\otimes\mathbb{M}_m\left(\C \right), \Z_2 \right)$   \\  & & \\
		3 & $\rho: A \to B\left(\H \right)$ & $C\left(M \right)\otimes\mathbb{M}_m\left(\C \right)\to B\left( L^2\left({M}, {\SS}\right)^N\right) $  \\ & & \\
		4	& $D$ & $\slashed D \otimes \Ga_{\mathrm{fin}}+ \Id_{L^2\left( M, \sS \right)^N} \otimes D_{\mathrm{fin}}$ \\  & & \\
		5	& $J$ & $J_{\mathrm{comm}}\otimes J_{\mathrm{fin}}$ \\  & & \\
	6	& $\left(\widetilde{\A}, \widetilde{\H}, \widetilde{D} \right)$ & $	\left(\Coo\left(\widetilde M \right) \otimes\mathbb{M}_m\left(\C \right) , L^2\left(\widetilde{M}, \widetilde{\SS}\right)^N, \widetilde D, \widetilde J  \right)$ \\
		&& where $\widetilde D \bydef \widetilde\Dslash \otimes \Ga_{\mathrm{fin}}+ \Id_{L^2\left(\widetilde M, \widetilde\sS \right)^N} \otimes D_{\mathrm{fin}}$,\\ &&
			 $\widetilde J \bydef  \widetilde J_{\mathrm{comm}}\otimes J_{\mathrm{fin}}$\\ && \\
		\hline
	\end{tabular}
	\\
	\\
	\\
The details of the proof are left  to the reader.
\end{proof}

\section{Infinite coverings}\label{ctr_case_sec}
\subsection{Coverings of $C^*$-algebras. Fundamental groups}

\paragraph*{}
Here we prove that sometimes an infinite covering of $C^*$-algebra with Hausdorff spectrum $\sX$ corresponds to a transitive covering of $\sX$.
\begin{empt}\label{hausdorff_covering_inf_empt} 
Let $A$ be a $CCR$-algebra (cf. Definition \ref{ccr_defn}) with connected, locally connected (cf. Definition \ref{top_locally_connected_defn}), locally compact, Hausdorff spectrum. Let $p: \widetilde \sX \to \sX$ be a transitive covering with connected $\widetilde{\sX}$ and residually finite group (cf. Definition \ref{residually_finite_defn}) $ G\left(\left. \widetilde\sX~\right|\sX \right)$ covering group  (cf. Definition \ref{top_group_of_covering_transformations_defn}). If $\overline p : \overline \sX \to \sX$ is the {disconnected covering of} $p: \widetilde \sX \to \sX$ and 
	\bean
\mathfrak{S}_p \bydef \left\{\left\{\sX_\la\right\}_{\la \in \La}, \left\{p^\mu_\nu:\sX_\mu\to \sX_\nu\right\}_{\substack{\mu,\nu \in \La\\\mu\ge\nu}}\right\}.
\eean 
is the {finite covering category of} $p : \widetilde \sX \to \sX$ (cf. Definition \ref{top_disconnected_defn}) then from Theorem \ref{blowing_sufficient_covering_inf_thm} it follows that there is  an algebraical finite covering category 
		\be\label{hausdorff_algebraical_finite_covering_category_eqn}
\mathfrak{S}_{A_0\left(p \right) }\bydef \left\{\left\{A_0\left( \sX_\la\right) \right\}_{\la\in \La}, \left\{A_0\left( p^\nu_\mu\right)  : A_0\left( \sX_\mu\right) \hookto A_0\left( \sX_\nu\right)\right\}_{\substack{\mu, \nu \in \La\\\mu \le \nu}}\right\}
\ee		
(cf. Definition \ref{algebraical_finite_covering_category_defn}. Moreover if $\widehat G$ is the profinite completion (cf. Example \ref{profinite_exm}) the triple $\left(A, A_0\left(\overline\sX \right), \widehat{G} \right)$   is the  {pre-covering} of $\mathfrak{S}_{A_0\left(p \right) }$ (cf. Theorem \ref{blowing_sufficient_covering_inf_thm}).
\end{empt}

\begin{theorem}\label{hausdorff_covering_inf_thm} 
Under the hypotheses \ref{hausdorff_covering_inf_thm} if we 	 suppose  that $A$ is a $C^*$-algebra of type $I_0$ (cf. Definition \ref{type_I_0_defn}), then  algebraical finite covering category  $\mathfrak{S}_{A_0\left(p \right) }$ is good (cf. Definition \ref{good_defn}) and a triple
$$
\left(A, A_0\left(\widetilde\sX \right), G\left(\left. \widetilde\sX~\right|\sX \right) \right)
$$
is the {infinite noncommutative covering} of $\mathfrak{S}_{A_0\left(p \right) }$ (cf. Definition \ref{infinite_noncommutative_covering_defn})
	
\end{theorem}

\begin{proof}
	If $\widehat{A}\bydef C^*$-$\varinjlim_{\la\in \La}A_0\left( \sX_\la\right)$ is the $C^*$-inductive limit of $\left\{A_0\left( \sX_\la\right)\right\}_{\la\in \La}$ (cf. Definition \ref{inductive_lim_non_defn}) and $\widehat{A}\hookto B\left( \overline \H\right)$ is a faithful, nondegenerate, equivariant representation (cf. Definitions \ref{faithful_representation_defn}, \ref{nondegenerate_repr_defn} and \ref{equivariant_representation_defn}) then from the Lemma \ref{infinite_faithful_nondegererate_lem} it follows that there is a natural faithful, nondegenerate, equivariant representation
	$$
\overline \pi: A_0\left(\overline\sX \right)\to 	B\left( \overline \H\right)
	$$
If $\left(A, \overline A', \widehat{G} \right)$ then from the Lemma \ref{infinite_faithful_nondegererate_lem} it follows that there is a natural faithful, nondegenerate, equivariant representation
$$
\overline \pi': \overline A'\to 	B\left( \overline \H\right)
$$
such that $C^*_{\overline\rho}\left(\overline M , \overline \F \right)$ is a hereditary subalgebra of $\overline A'$. We leave to the reader a proof of following statements:
\begin{itemize}
	\item $A_0\left(\overline\sX \right)$ is a $C^*$-algebra of type $I_0$ (cf. Definition \ref{type_I_0_defn}),
	\item $A_0\left(\overline\sX \right)$ is generated by by its commutative $C^*$-subalgebras explained in (ii) of the Lemma \ref{blowing_necessary_covering_inf_lem}
\end{itemize}
From (ii) of the Lemma \ref{blowing_necessary_covering_inf_lem} it follows that 
$$
\overline \pi\left(  A_0\left(\overline\sX \right)\right) \subset \overline \pi'\left( \overline A'\right),
$$
or there is the natural inclusion
$$
A_0\left(\overline\sX \right)\subset \overline A',
$$
From (i) of the Lemma \ref{blowing_necessary_covering_inf_lem} it follows that $A_0\left(\overline\sX \right)$ is a hereditary $C^*$-subalgebra of $\overline A'$. From the Theorem \ref{hered_spectrum_prop} it turns out that the spectrum of $A_0\left(\overline\sX \right)$ is an open subset of the spectrum of $\overline A'$. Taking into account the Proposition \ref{infinite_spectrum_limit_prop}  we conclude that the spectrum $\overline\sX$ of $A_0\left(\overline\sX \right)$ is homeomorphic to the spectrum of $\overline A'$, i.e. the spectrum of  $\overline A'$ is a  Hausdorff space. From the Lemma \ref{oa_haus_alg_lem} it follows that $\overline A'$ is a full algebra of operator fields (cf. Definition \ref{full_algebra_operator_fields_defn}) on $\overline \sX$. On the other hand $A_0\left( \overline\sX\right) $ is also a full algebra of operator fields (cf. Definition \ref{full_algebra_operator_fields_defn}) on $\overline\sX$. So  one has an isomorphism
$$
A_0\left( \overline\sX\right) \xrightarrow{\cong } \overline A'
$$
i.e. any  {pre}-{covering of algebraical finite covering category} $\mathfrak{S}_{A_0\left(p \right) }$ is equivalent to $\left(A, A_0\left(\overline\sX \right), \widehat{G} \right)$, so the triple $\left(A, A_0\left(\overline\sX \right), \widehat{G} \right)$ is the {disconnected infinite noncommutative covering} of $\mathfrak{S}_{A_0\left(p \right) }$ (cf. Definition \ref{disconnected_infinite_noncommutative_covering_defn}). Now this theorem becomes a consequence of (iii) of the Theorem \ref{blowing_sufficient_covering_inf_thm}.

\end{proof}
\begin{corollary}
	If $A$ is a $CCR$-algebra of type $I_0$ with locally connected, locally compact, Hausdorff spectrum then one has:
\begin{enumerate}
		\item [(i)]  	a triple 
\be\label{ctr_uni_lim_eqn}
\left(A, A_0\left( \widetilde{\sX}_{\mathrm{res~fin}}\right),G\left( \left.{\widetilde{\sX}}_{\mathrm{res~fin}}~\right|\sX\right)\right)	
\ee	
with given by \eqref{top_x_fr_eqn}	$\widetilde{\sX}_{\mathrm{res~fin}}$, is an universal covering of $C_0\left( {\sX}\right)$ (cf. Definition \ref{fundamental_group_nc_defn});
\item[(ii)]	there is a  natural  
group isomorphism 
\be\label{ctr_fg_iso_eqn}
\pi_1\left( C_0\left(\sX \right)\right)\cong G\left( \left.{\widetilde{\sX}}_{\mathrm{res~fin}}~\right|\sX\right) 
\ee
where  $\pi_1\left( A\right)$ 
is the {fundamental group}	of  $A$  (cf. Definition \ref{fundamental_group_nc_defn}.
\end{enumerate}
\end{corollary}

\begin{proof}
If $p: \widetilde{\sX}_{\mathrm{res~fin}}\to \sX$ is the natural covering and  an algebraical finite covering category 
\bean
\mathfrak{S}_{A_0\left(p \right) }\bydef \left\{\left\{A_0\left( \sX_\la\right) \right\}_{\la\in \La}, \left\{A_0\left( p^\nu_\mu\right)  : A_0\left( \sX_\mu\right) \hookto A_0\left( \sX_\nu\right)\right\}_{\substack{\mu, \nu \in \La\\\mu \le \nu}}\right\}
\eean
is given by \eqref{hausdorff_algebraical_finite_covering_category_eqn} then from the Theorem \ref{ctr_fin_thm} it follows that $\mathfrak{S}_{A_0\left(p \right) }$  contains \textit{all} classes of isomorphisms of  noncommutative finite-fold coverings of $A$. Now this corollary becomes a consequence of the Theorem \ref{hausdorff_covering_inf_thm}.		
\end{proof}
\begin{corollary}
	Let $A$ be a $CCR$-algebra of type $I_0$ with locally connected, locally compact, Hausdorff spectrum $\sX$. If the space $\mathcal X$ is weakly semilocally 1-connected (cf. Definition \ref{top_weakly_semi1_defn}) then one has
		\bean
\pi_1\left( A\right)\cong \mathfrak{ResFin}\left( 	\pi_1^{\mathrm{w}}\left( \sX, x_0\right)\right) 
\eean
where  $\pi_1\left( A\right)$ 
is the {fundamental group}	of  $A$  (cf. Definition \ref{fundamental_group_nc_defn} and $\pi_1^{\mathrm{w}}\left( \sX, x_0\right)$ is a weak fundamental group (cf. Definition \ref{top_weak_fundamental_group_defn}) and $\mathfrak{ResFin}$ is given by the equation \eqref{top_fin_eqn}.

\end{corollary}
\begin{proof}
The proof of this corollary is similar to the \ref{comm_uni_lim_cor} one.
\end{proof}

	Let $\tau: C_0\left(\sX \right)\to \C$ be a faithful state given by  \eqref{meafunc_eqn}, i.e. there exists a Borel
measure $\mu$ on $\sX$ with values in $\left[0, +\infty\right]$ such that
\be\label{ctr_meafunc_eqn}
\tau\left(a \right) = \int_{\sX} a~ d\mu \quad \forall a\in C_0\left(\sX \right).
\ee
Let  $A$ be a continuous trace $C^*$-algebra (cf. Definition \ref{continuous_trace_c_alt_defn}) such that $\sX$ is the spectrum of $A$.
there is a given by \eqref{ctr_ka_t_eqn} functional 
\be\label{ctr_tau_eqn}
\begin{split}
	\tr_\tau : K\left( A\right)\to \C,\\
	a \mapsto \tau\left(\tr_a \right).  
\end{split} 
\ee
where $K\left( A\right)$ is Pedersen's ideal (cf. Definition \ref{pedersen_ideal_defn}) and  $\tr\left( a\right)  \bydef \hat a$ is given by \eqref{ctr_hat_eqn}. There is $\C$-valued product
\be
\begin{split}
	\left(\cdot, \cdot\right): K\left( A\right)\times  K\left( A\right)\to\C,\\
	\left(a, b\right)\mapsto \tr_\tau\left( a^*b\right). 
\end{split}
\ee
Denote by $L^2\left( A, \tr_\tau  \right)$ the Hilbert norm completion of $ K\left( A\right)$.	
If $\overline p: \overline \sX  \to  \sX$  is the {disconnected covering of} $p:\widetilde {\mathcal X}\to \sX$ 	(cf. Definition \ref{top_disconnected_defn})
then there is a given by \eqref{top_ot_eqn} functional
\be\label{ctr_ot_eqn}
\overline\tau:C_c\left( \overline{\mathcal X}\right) \to \C.
\ee
Similarly to \eqref{ctr_tau_eqn} there is a functional
\be\label{ctr_otau_eqn}
\begin{split}
	\tr_{\overline\tau} : K\left( A_0\left( \overline\sX\right)\right) \to \C,\\
	\overline	a \mapsto \tau\left(\tr_{\overline a} \right).  
\end{split} 
\ee
There is $\C$-valued product
\be\label{ctr_tx_eqn}
\begin{split}
	\left(\cdot, \cdot\right): K\left( A_0\left(  \overline\sX\right)\right) \times K\left( A_0\left(  \overline\sX\right)\right)\to\C ,\\
	\left(\overline a,\overline b\right)\mapsto \tr_{\overline\tau}\left(\overline a^*\overline b\right). 
\end{split}
\ee
Denote by $L^2\left(  A_0\left( \overline\sX\right) , \tr_{\overline\tau}  \right)$ the Hilbert norm completion of $ K\left( \overline A\right)$. 

\subsection{Coverings of operator spaces}
\paragraph*{} Here we consider the generalization of the Theorem  \ref{hausdorff_covering_inf_thm}.
\begin{empt}\label{ctr_op_main_empt}
	Let $\left(X, Y\right)$ be a  {sub-unital subspace} (cf. Definition \ref{operator_space_subunital_defn}) such that $C^*$-envelope $C^*_e\left(X, Y\right)$ of $\left(X, Y\right)$ (cf. Definition \ref{operator_space_envelope_defn})
	is a $CCR$-algebra (cf. Definition \ref{ccr_defn}) of type $I_0$ with locally connected, locally compact> Hausdorff spectrum $\sX$. From the Theorem \ref{blowing_sufficient_covering_inf_thm}  it turns out that is an {algebraical  finite covering category} (cf. Definition \ref{algebraical_finite_covering_category_defn}).
	\be\nonumber
	\begin{split}
		\mathfrak{S}_{A_0\left( p\right)  } = \\
		=\left\{A_0\left( p_\la\right) : A  \hookto A_0\left(\sX_\la \right)\right\}, \left\{A_0\left( p_\nu^\mu\right) : A_0\left(\sX_\nu \right) \hookto A_0\left(\sX_\mu \right)\right\}
	\end{split}
	\ee
(cf. the equation \eqref{hausdorff_algebraical_finite_covering_category_eqn})
	From the Theorem \ref{hausdorff_covering_inf_thm} it follows that $\mathfrak{S}_{A_0\left(p \right)  }$ is good. Moreover   the triple 
	$$
	\left(A, A_0\left(\widetilde\sX\right),G\left(\left.\widetilde\sX~\right|~ \mathcal X\right)\right)
	$$
	the  { infinite noncommutative covering} of $\mathfrak{S}_{A}$(cf. Definition \ref{infinite_noncommutative_covering_defn}). From the Theorem \ref{ctr_fin_op_thm} it follows that there is a algebraical  finite covering category of operator spaces (cf. Definition \ref{comp_op_pt_defn})
	\bean
	\mathfrak{S}_{\left(X, Y \right)}
	= \left(\left\{\left( \pi_{X_\la}, \pi_{Y\la}\right) : \left(X, Y \right) \to \left(X_\la, Y_\la \right) \right\}_{\la \in \La},\left\{\left( \pi^{\nu}_{X_\mu}, \pi^{\nu}_{Y_\mu}\right)  \right\}_{\substack{\mu, \nu \in \La\\ \nu > \mu}}\right)
	\eean
	where $X_\la = C_0\left(\lift_{p_\la}\left[X\right] \right)$,  $\pi_{X_\la} = \left.C_0\left(p_\la \right) \right|_{X}$, $\pi^{\nu}_{X_\mu}= \left.C_0\left(p_\mu \right) \right|_{X_\mu}$.
\end{empt}
\begin{theorem}\label{ctr_op_main_thm}
	Consider the situation of \ref{ctr_op_main_empt}.  If $\widetilde X$  the {inverse noncommutative limit} (cf. Definition \ref{spec_lim_defn})  of $\mathfrak{S}_{\left(X, Y \right)}$ then there is the natural complete isomorphism  $\widetilde X \cong C_0\left(\lift_{\widetilde p}
	\left[X\right]\right)$  where $\widetilde p : \widetilde \sX \to \sX$ is the natural transitive covering.  
\end{theorem}
\begin{proof}
	Let $\widehat A \stackrel{\text{def}}{=} C^*$-$\varinjlim_{\la \in \La}A_0\left( \sX_\la\right)$, and let  $\pi_a:\widehat{A} \to B\left(\H_a \right)$ be the atomic representation,
	Denote by  $\widetilde p_\la =  \widetilde \sX \to \sX_\la$ for all $\la \in \La$ the natural coverings. Let $\widetilde \sU \subset \widetilde\sX$ be an open subset which is homeomorphically mapped onto $\sU =  \widetilde p\left(\widetilde\sU\right)$. Let 
	\be\label{ctr_bt_eqn}
	\widetilde b \in C_0\left(\lift_{\widetilde p}
	\left[X\right]\right) \quad \supp\widetilde b \subset \widetilde\sU.
	\ee
	There is the net $\left\{\desc_{\widetilde p_\la}\left( \widetilde b\right)\in X_\la \right\}_{\la \in \La} \subset \widehat A$. 
	One has
	$$
	\widetilde b = \lim_{\la \in \La }\pi_a\left(\desc_{\widetilde p_\la}\left( \widetilde b\right) \right) 
	$$
	where the strong limit is implied. So $\widetilde b$  is {subordinated} to 	$\mathfrak{S}_{\left(X, Y \right)}$ (cf. Definition \ref{spec_op_defn}). According to the Definition \ref{spec_lim_defn} one has $\widetilde b\in \widetilde X$. Following proof contains two parts:
	\begin{itemize}
		\item [(i)]  $C_0\left(\lift_{\widetilde p}
		\left[X\right]\right)\subset \widetilde X$.
		\item [(ii)]  $\widetilde X\subset C_0\left(\lift_{\widetilde p}
		\left[X\right]\right)$
	\end{itemize}
	(i)	The $\C$-linear space of the given by \eqref{ctr_bt_eqn} is dense in   $C_0\left(\lift_{\widetilde p}
	\left[X\right]\right)\subset \widetilde X$ with respect to $C^*$-norm topology.\\
	(ii)	For any $x \in \sX$ denote by $S_x$ the linear space of all continuous linear functionals $s : \rep_x\left( C^*_e\left(X_x, Y_x \right)\right)  \to \C$ such that $s\left(X_x\right)=\{0\}$. Since  $C^*_e\left(X_x, Y_x \right)\cong \K$ and $\K$ is a reflexive Banach space, i.e. $\K^{**}\cong\K$, and taking into account that
	$X_x \subset \rep_x\left( C^*_e\left(X, Y \right)\right)$ is a closed subset one has
	$$
	X_x = \left\{\left.t \in \rep_x\left( C^*_e\left(X, Y \right)\right)~\right|~ s\left( t\right)= 0\quad \forall s\in S_x \right\}
	$$
	If $\H_x$ is the space of the representation $\rep_x$ then for any $s \in S$ there are $\xi, \eta \in \H_x$ such that $s\left(a \right)= \left(\xi, a \eta  \right)_{\H_x}$ for each $a \in  \rep_x\left( C^*_e\left(X, Y \right)\right)$. Denote by
	$$
	T_x = \left\{\left.\left(\xi, \eta\right)\in \H_x\times \H_x\right| a \mapsto \left(\xi, \rep_x\left( a\right)  \eta  \right)_{\H_x} \in S_x\right\} 
	$$
	If $\H_{\widetilde x}$ the space of the representation $\rep_{ \widetilde{x}}:A_0\left(\widetilde \sX\right)\to B\left(\H_{\widetilde x}\right)$ then there if the natural isomorphism $c_{\widetilde x}:\H_{\widetilde x} \cong \H_{\widetilde p\left(\widetilde x\right)}$. So for any $\widetilde x \in \widetilde \sX$ there is the the natural representation $\rep_{\widetilde x}: A \to B\left(\H_{\widetilde x}\right)$. There is the set
	$$
	T_{\widetilde x} = \left\{\left.\left(c^{-1}_{\widetilde x}\left( \xi\right) , c^{-1}_{\widetilde x}\left( \eta\right)\right)\in \H_{\widetilde x}\times \H_{\widetilde x}\right|  \left(\xi, \eta  \right)\in T_x\right\} 
	$$
	such that
	$$
	X = \left\{ a \in A \left| \left(\widetilde\xi, \rep_{\widetilde x},\left( a\right)  \widetilde\eta  \right)_{\H_{\widetilde x}}= 0\quad \forall \left(\widetilde\xi,  \widetilde\eta  \right)\in T_{\widetilde x}\right. \right\}
	$$
	From the Theorem \ref{ctr_fin_op_thm} it follows that $X_\la = C_0\left(\lift_p\left[X\right]\right)$ (cf. equation \eqref{ctr_fin_op_eqn}) so one has
	$$
	X_\la = \left\{ a_\la \in A_\la \left| \left(\widetilde\xi, \rep_{\widetilde x}\left( a_\la\right)  \widetilde\eta  \right)_{\H_{\widetilde x}}= 0\quad \forall \left(\widetilde\xi,  \widetilde\eta  \right)\in T_{\widetilde x} \right.\right\}
	$$
	If $ a_\la \in X_\la$ then $\left(\widetilde\xi, \rep_{\widetilde x}\left( a_\la\right)  \widetilde\eta  \right)_{\H_{\widetilde x}}= 0$ for all  $\left(\widetilde\xi,  \widetilde\eta  \right)\in T_{\widetilde x}$. If $\widetilde a$ is a strong limit of $\left\{a_\la\right\}$ then $\widetilde a$ is a weak limit $\left\{a_\la\right\}$ so one has
	\be\label{ctr_xi_eta_eqn}
	\left(\widetilde\xi, \rep_{\widetilde x}\left(\widetilde a\right)  \widetilde\eta  \right)_{\H_{\widetilde x}}= 0; \quad \forall \left(\widetilde\xi,  \widetilde\eta  \right)\in T_{\widetilde x} \quad \forall \widetilde x \in \widetilde\sX.
	\ee
	since the linear span of $T_{\widetilde x}$ is dense in $c^{-1}_{\widetilde x}\left(S_{p\left(\widetilde x \right) } \right)
	$
	the equation \eqref{ctr_xi_eta_eqn} is equivalent to
	\bean
	\rep_{\widetilde x}\left(\widetilde a\right) \in c^{-1}_{\widetilde x}\left(X_{p\left(x\right) } \right),
	\\
	\rep_{\widetilde x}\left(\widetilde a\right) \in  \rep_{\widetilde x}\left(C_0\left(\lift_{\widetilde p}
	\left[X\right]\right)\right)\quad \forall \widetilde x \in \widetilde\sX.
	\eean
	Hence for every $\widetilde x \in \widetilde\sX$  one has $ \rep_{\widetilde x}\left(\widetilde X\right)\subset  \rep_{\widetilde x}\left(C_0\left(\lift_{\widetilde p}
	\left[X\right]\right)\right)$. From the Lemma \ref{top_sub_eq_lem} it turns out that $\widetilde X \subset C_0\left(\lift_{\widetilde p}
	\left[X\right]\right)$.
\end{proof}
\subsection{Induced representations}\label{ctr_induced_infinite_sec}

\paragraph*{}
Let $A$ be a  $C^*$-algebra with Hausdorff spectrum $\sX$. If  $\tau : C_c\left(\sX \right) \to \C$ be a functional such that $\tau\left(a \right) > 0$ for any positive $a$ then there is a given by \eqref{ctr_a_eqn} representation
\be\nonumber
\tr_\tau: A \hookto B\left(L^2\left(A, \tr_\tau\right) \right). 
\ee
If 
\be\nonumber
\begin{split}
	\mathfrak{S}_{A_0\left( p\right)  } = \\
	=\left\{A_0\left( p_\la\right) : A  \hookto A_0\left(\sX_\la \right)\right\}, \left\{A_0\left( p_\nu^\mu\right) : A_0\left(\sX_\nu \right) \hookto A_0\left(\sX_\mu \right)\right\}
\end{split}
\ee
is a given by the Theorem  \ref{blowing_sufficient_covering_inf_thm}   algebraical  finite covering category} (cf. Definition \ref{algebraical_finite_covering_category_defn}) then from the Theorem  \ref{hausdorff_covering_inf_thm} it follows that is  the  { infinite noncommutative covering} of $\mathfrak{S}_{A}$ (cf. Definition \ref{infinite_noncommutative_covering_defn}) given by
$\left(A, A_0\left(\widetilde\sX\right),G\left(\left.\widetilde\sX~\right|~ \mathcal X\right)\right)$ where $\widetilde\sX \bydef \varprojlim  \mathfrak{S}_{\mathcal{X}}$ (cf. Definition \ref{top_disconnected_defn}).  Similarly to \eqref{ctr_tx_eqn} define a scalar product $K\left( A_0\left(  \widetilde\sX\right)\right) \times K\left( A_0\left(  \widetilde\sX\right)\right)\to\C$ and denote by $L^2\left(A_0\left(\widetilde\sX \right) , \tr_{\widetilde\tau}\right)$.
\begin{exercise}
	Prove that the natural representation $A_0\left(\widetilde\sX \right)\hookto B\left(L^2\left(A_0\left(\widetilde\sX \right) , \tr_{\widetilde\tau}\right)\right)$   is {induced} by $\left( \tr_\tau, \left( A, A_0\left(\widetilde M \right) , G\left(\left.\widetilde M~\right|M \right) \right)\right)$ (cf. Definition \ref{induced_repr_inf_defn})
\end{exercise}
\begin{exercise}
	Prove that in the described in this section situation 
	the  representation $A_0\left(\widetilde\sX \right)\hookto B\left(L^2\left(\widetilde M, \widetilde S\right)\right)$   is {induced} by $\left( \rho, \left( A, A_0\left(\widetilde M \right) , G\left(\left.\widetilde M~\right|M \right) \right)\right)$ (cf. Definition \ref{induced_repr_inf_defn}) where $\rho: A\hookto B\left(L^2\left(M, S\right)\right)$ is the  natural representation.
\end{exercise}
\subsection{Coverings of $*$-algebras}
\subsubsection{Coverings of bounded operator $*$-algebras}
\paragraph{}
	Let $A$ be a $CCR$-algebra of type $I_0$ with locally connected, locally compact, Hausdorff spectrum. Let $R\subset C_0\left(\sX\right)$ be a $c$-soft *-subalgebra (cf. Definition \ref{top_soft_r_defn}). If $\A\subset A$ is a dense *-subalgebra which is an $R$-module then $\A$-sheaf $\mathscr S^\A$ (cf. Definition \ref{top_x_sheaf_defn}) is $c$-soft (cf. Theorem  \ref{sheaf_soft_m_thm}). 
	Let 
\bean
	\mathfrak{S}_p \bydef \left\{\left\{\sX_\la\right\}_{\la \in \La}, \left\{p^\mu_\nu:\sX_\mu\to \sX_\nu\right\}_{\substack{\mu,\nu \in \La\\\mu\ge\nu}}\right\}
\eean 
be a {finite covering category of} $p : \widetilde \sX \to \sX$ (cf, Definition \ref{top_disconnected_defn}), and let
\bean
\mathfrak{S}_{A_0\left(p \right) }\bydef \left\{\pi_\la=A_0\left( p_{\la}\right) :A  \hookto A_0\left( \mathcal{X}_\la\right) \right\}_{\la \in \La}\in \mathfrak{FinAlg}. 
\eean
be a given by the Theorem \eqref{blowing_sufficient_covering_inf_thm} algebraical  finite covering category (cf. Definition \ref{algebraical_finite_covering_category_defn}).
For any $\la\in \La$ let $\A_\la\bydef  A_0\left(\sX_\la \right)\cap p^{-1}_\la \mathscr S^{\A}\left( \sX_\la\right)$,  and let $\widetilde \A\bydef  A_0\left(\widetilde \sX \right)\cap \widetilde p^{-1} \mathscr S^{\A}\left(\widetilde \sX\right)$ where a notation of the Definition \ref{top_x_sheaf_defn} is used.	
The Exercise \ref{ctr_bounded_exer} enables us to construct a algebraical  finite covering category of  bounded operator *-algebras (cf. Definition \ref{comp_oa_defn}) given by
\be\label{ctr_inf_cov_oae_eqn}
\begin{split}
	\mathfrak{S}_{\A } =\\
	\left(\left\{ \left.\A_0\left( p_\la\right)\right|_{\A} :  \A  \hookto  \A_\la\right\}, \left\{ \left.A_0\left( p^\mu_{ \nu}\right)\right|_{A_\nu} :\A_\nu \hookto \A_\mu\right\}\right).
\end{split}
\ee
The formula \eqref{ctr_inf_cov_oae_eqn} is a specialization of the given by \eqref{comp_oa_eqn} one. From the Theorem \ref{blowing_sufficient_covering_inf_thm} it follows that 
the triple $$\left(A, A_0\left( \widetilde \sX\right) , G\left(\left.\widetilde \sX\right|\sX\right)\right)$$ is an inverse noncommutative limit of $\mathfrak{S}_{A_0(p)}$ (cf. Definition \ref{infinite_noncommutative_covering_defn}).
Now we have all described in the \ref{inf_oa_empt} ingredients.

\begin{exercise}
	Similarly to the Theorem \ref{top_inf_cov_oa_thm} 
prove that the triple
$$
\left(\A,\widetilde \A, G\left(\left.\widetilde \sX \right|  \sX\right)\right) 
$$ is an {infinite noncommutative covering}  of $\mathfrak{S}_{\A_0\left( p\right) }$ in the sense of the  Definition \ref{inf_cov_oa_defn}.
\end{exercise}

\subsubsection{Coverings of $O^*$-algebras}
\paragraph*{}
Discussed below results close to  the Section \ref{top_inf_o_sec} ones. Exercises of this Section can be solved similarly to proofs of Lemmas and Theorems of the Section \ref{top_inf_o_sec}.
\begin{empt}
	Similarly to \ref{ctr_fin_o_empt} 	let $M$ be a smooth manifold, and let $A$ be  a {homogeneous of order} $n$ $C^*$-algebra (cf. Definition \ref{ctr_homo_defn}) such that  $M$ is the spectrum of $A$. There is a linear bundle $\left(E, q, M\right)  $ with fibers $\C^{n^2}$ such that there is $C_b\left( M\right)$-linear isomorphism $A \cong 
	\Ga_0\left(M, E \right)$ (cf. \ref{ctr_fin_o_empt}). If $E$ has structure of smooth manifold such that the map $E\to M$ is smooth.  	From  \eqref{ctr_fin_o_empt} it turns out that  $\Ga^\infty\left(M, E \right)$ is *-algebra.  Let $p: \widetilde M \to M$ be a transitive  covering with a residually finite covering group $G\left(\left. \widetilde M \right| M \right)$ of  (cf. Definitions \ref{residually_finite_defn} and \ref{top_group_of_covering_transformations_defn}). The Theorem  \ref{blowing_sufficient_covering_inf_thm} yields  an  algebraical  finite covering category (cf. Definition \ref{algebraical_finite_covering_category_defn}) given by 		
	\be\label{ctr_bases_point_eqn}
	\begin{split}
		\mathfrak{S}_{A_0\left( p\right)  } = \\
		=\left\{A_0\left( p_\la\right) : A  \hookto A_0\left(M_\la \right)\right\}, \left\{A_0\left( p_\nu^\mu\right) : A_0\left(M_\nu \right) \hookto A_0\left(M_\mu \right)\right\}.
	\end{split}
	\ee
	If $\overline p: \overline M\to M$ is the  disconnected covering of $p: \widetilde M \to M$ (cf. Definition \ref{top_disconnected_defn}) then there the natural covering $\overline p: \overline M \to M$. From the Proposition \ref{top_cov_mani_prop} it follows that $\overline M$ has a (unique) structure of smooth manifold such that $\overline p$ is smooth. If $\left(\overline E, \overline q,  \overline M\right)$  is the {inverse image} of $\left( E, q, M\right)$ by $p$ (cf. Definition \ref{vb_inv_img_funct_defn}) then there is a transitive  covering  $\overline E \to E$  so from the Proposition \ref{top_cov_mani_prop} it turns out that there is the unique smooth structure on $\overline E$ such that map $\overline E\to E$ is smooth. So there is a *-algebra $\Ga^\infty\left(\overline M, \overline E \right)$. 
	If $\tau: C_c\left( M\right) \to \C$ is a functional
	$$
	a \mapsto \int_M a ~d\mu. 
	$$
	such that $a > 0\quad\Rightarrow\quad \tau\left(a\right)> 0$ then similarly to \eqref{ctr_a_eqn} one can define Hilbert spaces $L^2\left(  A , \tr_{\tau}  \right)$  and $L^2\left(  A_0\left( \overline M\right) , \tr_{\overline\tau}  \right)$. If
	\bean
	\D \bydef \Ga^\infty_c\left( M,  E \right),\\
	\widehat \D \bydef \Ga^\infty_c\left(\overline M, \overline E \right),
	\eean 
	are spaces of smooth sections with compact support then  both $ \D$ and $\widehat \D$ are dense subspaces of $L^2\left( A , \tr_{\tau}  \right)$ and $L^2\left(  A_0\left( \overline M\right) , \tr_{\overline\tau}  \right)$ respectively. If both  $D\left( M,  E \right)$ and $D\left(\overline M, \overline E \right)$ be *-algebras of differential operators 
	$\Ga^\infty \left(M, E \right)  \to \Ga^\infty \left(M, E \right)$ and $\Ga^\infty \left(\overline M,  \overline E \right)  \to \Ga^\infty \left(\overline M, \overline E \right)$ respectively (cf. Definition \ref{do_man_defn}) then there are inclusions
	\be
	\begin{split}
		\Ga^\infty \left( M,   E \right)\subset D \left( M,   E \right)\subset \L^\dagger\left( D \right),\\ 
		\Ga^\infty \left(\overline M,  \overline E \right)\subset D\left(\overline M,  \overline E \right)\subset \L^\dagger\left(\widehat \D \right)
	\end{split}
	\ee
	where both $\L^\dagger\left(\D \right)$  and $\L^\dagger\left(\widehat\D \right)$ are given by \eqref{l_dag_eqn}.
	For all $\la\in \La$ one can construct *-algebras $\Ga^\infty \left( M_\la,   E_\la \right)$ and $D \left( M_\la,   E_\la \right)$, such that there are following inclusions
	\bean
	\forall \la\in\La \quad \Ga^\infty \left( M_\la,   E_\la \right)\subset D \left( M_\la,   E_\la \right),\\
	\forall \mu,\nu\in\La \quad \nu \le\mu \quad\Rightarrow\quad  \Ga^\infty \left( M_\nu,   E_\nu \right)\subset \Ga^\infty \left( M_\mu,   E_\mu \right) \text{ AND }\\
	\text{ AND } D \left( M_\nu,   E_\nu \right)\subset D \left( M_\mu,   E_\mu \right)
	\eean
	There are pointed algebraic finite covering categories of $O^*$-algebras given by
	\bean
	\mathfrak{S}^{\text{pt}}_{D\left( M,   E \right)}\bydef \left(\left\{\pi_\la: D\left( M,   E \right) \hookto D\left( M_\la,   E_\la \right) \right\}, \left\{\pi^\mu_\nu\right\}\right),\\
	\mathfrak{S}^{\text{pt}}_{\Ga^\infty\left( M,   E \right)} = \left(\left\{\left.\pi_\la\right|_{\Ga^\infty\left( M,   E \right)}: \Ga^\infty\left( M,   E \right) \hookto \Ga^\infty\left( M_\la,   E_\la \right) \right\}, \left\{\left.\pi^\mu_\nu\right|_{\Ga^\infty\left( M_\mu,   E_\mu \right)}\right\}\right).
	\eean
	(cf. Definition \ref{comp_o*_defn}).
	There is the natural inclusion
	$$
	\pi: \bigcup_{\la \in \La} D \left( M_\la,   E_\la \right)\hookto \L^\dagger\left(\widehat\D \right).
	$$

\end{empt}
%
	
	\begin{exercise}
 	Using the Exercise \ref{ctr_hilbert_exer}
	prove following statements:
	\begin{enumerate}
		\item The triple $\left(D\left( M,   E \right),D\left(\widetilde M, \widetilde  E \right), G\left(\left.\widetilde M\right|  M\right)\right) $ is the $\pi$-inverse noncommutative limit of $\mathfrak{S}^{\text{pt}}_{D\left( M,   E \right)}$ (cf. Definition \ref{inv_o*_lim_defn}).
		\item The triple $\left(\Ga^\infty\left( M,   E \right),\Ga^\infty\left(\widetilde M, \widetilde  E \right), G\left(\left.\widetilde M\right|  M\right)\right) $ is the  $\left.\pi\right|_{\cup_{\la \in \La} \Ga^\infty \left( M_\la,   E_\la \right)}$-inverse noncommutative limit of $\mathfrak{S}^{\text{pt}}_{\Ga^\infty\left( M_\la,   E_\la \right)}$.
	\end{enumerate}
	\end{exercise}

\begin{empt}\label{str_spinc_infc_empt}
	Suppose that  $\left( E, r, M\right)$ is a vector  bundle with a sesquilinear form (cf. Definition \ref{top_herm_bundle_form_defn}). If both $\widetilde  E$ and $\overline  E$ are {inverse images} of $ES$ by $\widetilde p$ and $\overline p$ (cf. Definition \ref{vb_inv_img_funct_defn}) then both $\widetilde  ES$ and $\overline  E$ are  smooth manifolds such that both maps $\widetilde E \to \widetilde M$ and $\overline E \to overline M$ are smooth (cf. Proposition \ref{top_cov_mani_prop}). Similarly to \eqref{ctr_gac_prod} there are $\C$-valued products 
	\be\nonumber
	\begin{split}
		\left(\cdot, \cdot \right): \Ga_c\left(\overline M, \overline E \right)\times \Ga_c\left(\overline M,\overline E \right)\to \C,\\
		\left(\overline \xi, \overline \eta \right)\bydef \overline \tau\left(\left\langle \overline \xi, \overline \eta \right\rangle_c \right),\\
		\left(\cdot, \cdot \right): \Ga_c\left(\widetilde M, \widetilde E \right)\times \Ga_c\left(\widetilde M,\widetilde E \right)\to \C,\\
		\left(\widetilde \xi, \widetilde \eta \right)\bydef \widetilde \tau\left(\left\langle \widetilde \xi, \widetilde \eta \right\rangle_c \right) 
	\end{split}
	\ee	
	where $\overline \tau$ is given by the equation \eqref{top_ot_eqn} and $\widetilde\tau\bydef\overline\tau|_{C_c\left(\widetilde{   M} \right)  }$.  Denote by $L^2\left( M, E \right)$ and $L^2\left(\widetilde M,\widetilde E \right)$ the Hilbert norm completions of $\Ga_c\left( M, E\right)$ and $\Ga_c\left(\widetilde M,\widetilde E \right)$. 
	
\end{empt}
\begin{exercise}
	Prove that there is an algebraic finite covering category of $O^*$-algebras given by
	\bean
	\mathfrak{S}_{D^*\left( M,   E \right)}\bydef \left(\left\{\phi_\la: D^*\left( M,   E \right) \hookto D\left( M_\la,   E_\la \right) \right\}, \left\{\phi^\mu_\nu\right\}\right)
	\eean
	(cf. Definition \ref{comp_o*_defn}), where both $D^*\left( M,   E \right)$ and $D\left( M_\la,   E_\la \right)$ are spaces of adjointable differential operators (cf. Equation \eqref{top_diff*_eqn}).
\end{exercise}
If
$$
\widehat\D_E \bydef \Ga^\infty_c\left(\overline M, \overline E\right) \subset L^2\left(\widetilde M,\widetilde E \right)
$$
then there is the natural inclusion
$$
\pi_E: \bigcup_{\la \in \La} D \left( M_\la,   E_\la \right)\hookto \L^\dagger\left(\widehat\D _E\right).
$$
where $\L^\dagger\left(\widehat\D \right)$ is given by \eqref{l_dag_eqn}.
\begin{exercise}
	Let both  $D\left( M, E \right)$ and $D\left(\widetilde M, \widetilde E \right)$ be *-algebras of differential operators 
	$\Ga^\infty \left(M, E \right)  \to \Ga^\infty \left(M, E \right)$ and $\Ga^\infty \left(\widetilde M,  \widetilde E \right)  \to \Ga^\infty \left(\overline M, \overline E \right)$. Prove that the triple $\left(D\left( M,  E \right),D\left(\widetilde M, \widetilde  E\right), G\left(\left.\widetilde M\right|  M\right)\right) $ is the $\pi_E$-inverse noncommutative limit of $\mathfrak{S}_{D\left( M,   E \right)}$ (cf. Definition \ref{inv_o*_lim_defn}).
\end{exercise}

\subsection{Coverings of spectral triples}
We leave to the reader construction of an infinite covering of spectral triples which correspond to the following mapping
the following mapping.
\\
\\
\begin{tabular}{|c|c|c|}
	\hline
	&General theory & The specialization\\ 
	\hline
	&	&\\
	Hilbert spaces & $\H$  and $\widetilde\H$ &  $L^2\left(M, \sS \right)^k$ and $L^2\left(\widetilde M, \widetilde \sS \right)^k$\\ & & \\
	Pre-$C^*$-algebra	& $\A$  & $\Coo\left(M \right)\otimes\mathbb{M}_n\left(\C \right)$   \\  & & \\
	Pedersen's ideal	& $K\left(\widetilde A\right) $  & $C_c\left(\widetilde M \right)\otimes\mathbb{M}_n\left(\C \right)$   \\  & & \\
	The space of  	& $\widetilde{W}^\infty$  & $\Coo_c\left(\widetilde M \right)\otimes\mathbb{M}_n\left(\C \right)$   \\ smooth elements & & \\& & \\
	Dirac operators & $D$ & $ \Dslash\otimes \Ga_{\mathrm{fin}}+ \Id_{L^2\left(M, \sS \right)^k} \otimes D_{\mathrm{fin}}$  \\
	& $\widetilde{D}$ & ?\\  & & \\
	&$\H^\infty\bydef \bigcap_{n =1}^\infty \Dom D^n\subset \H$ & $\Ga^\infty\left(M, \sS \right)^k= \bigcap_{n =1}^\infty \Dom D^n$ \\& & \\
	\hline
\end{tabular}
\\
\\
\\

\chapter{Groupoids, foliations  and their coverings}\label{foliations_chap}
\paragraph{}

\section{Geometric construction}
\paragraph{}
In this section we consider locally-compact groupoids $\G$ (cf. Definition \ref{groupoid_topological_defn}) such that $\G^0$ is Hausdorff.
\begin{definition}\label{groupoid_lift_defn}
Let $q: \widetilde \G^0 \to \G^0$ be transitive covering
and let $\widetilde \G^0\times_{\G^0}\G$ be pullback where the map $s : \G \to \G^0$ is implied. We say that the covering  $q$ is $\G$-\textit{regular} if there is  the groupoid $\widetilde \G^0$ with homeomorphism $\widetilde \G\cong \widetilde \G^0\times_{\G^0}\G$ such that the natural continuous map 
\be\label{groupoid_lift_eqn}
q_\G: \widetilde \G \to \G
\ee
 is a covering (cf. Definition \ref{covering_groupoid_defn}). The groupoid $\widetilde \G$ is the $q$-\textit{lift} of $\G$. 	A $p$-\textit{topology} on $\widetilde\G$ is such that both maps $q_\G$ and $s:  \widetilde\G \to \widetilde\G^0$ are continuous.
\end{definition}

\begin{remark}\label{groupoid_lift_rem}
 From the  Theorem \ref{covering_groupoid_thm} it follows that the given by the  the Definition \ref{groupoid_lift_defn} groupoid is unique.
\end{remark}
\begin{remark}\label{groupoid_liftpp_rem}
Any element $\left(\widetilde y, x \right) \in \widetilde \G^0\times_{\G^0}\G$ uniquely defines a point $r\left(\widetilde{ x} \right)\in \widetilde{\G^0}$ where  $\widetilde{ x} \in \widetilde{\G}$ comes from  $\left(\widetilde y, x \right)$.
\end{remark}
\begin{example}\label{groupoid_lift_exm}
	According to the Theorem \ref{covering_groupoid_thm}) the existence of coveting is relevant to the lifting problem
	\newline
	\begin{tikzcd}
		& \widetilde \G^0\arrow[d, ""]\\
		\F \arrow[ru,  dashrightarrow]\arrow[r]	& \G^u_u
	\end{tikzcd}
	\newline
One can find obstructions of fundamental group $\pi_1\left( \G^u_u, u \right)$. However if the space  $ \G^u_u$ is not locally path connected  then 	usage of $\pi_1\left( \G^u_u, u \right)$ does not sufficient. Application of the weak fundamental group $\pi^{\text{w}}_1\left( \G^u_u, u \right)$ (cf. Definition \ref{top_weak_fundamental_group_defn}) is more universal.

\end{example}
\begin{empt}\label{groupoid_haar_lift_empt}
	Let $\G$ be a locally compact groupoid, and let   $p: \widetilde \G^0 \to \G^0$ is a $\G$-{regular} covering. If  $\widetilde\G$ is  the $p$-{lift} of $\G$  then $\widetilde\G$ is locally compact (cf. Exercise \ref{groupoid_loc_comp_lift_exer}). Let $\left\{\la^u \left| u \in \G^0\right.\right\}$ be a left  Haar system (cf. Definition \ref{groupoid_haar_defn}) for $\G$. For any $\widetilde u \in \widetilde\G^0$ there is a homeomorphism $\widetilde \G^{\widetilde u} \cong \G^{p\left(\widetilde u \right) }$. Using it  for any $\widetilde u$ one can  naturally construct a measure $\widetilde  \la^{\widetilde u}$ on $\widetilde \G$. From the conditions (a)-(c) of the Definition \ref{groupoid_haar_defn} it follows that
	\begin{enumerate}
		\item [(a)] the support $\supp \widetilde  \la^{\widetilde u}$ of the measure $\widetilde  \la^{\widetilde u}$ is $\widetilde \G^{\widetilde u}$,
		\item [(b)]  (continuity) for any $\widetilde f \in C_c\left(\widetilde \G\right)$, $\widetilde u \mapsto \widetilde\la(\widetilde f)\left( \widetilde u\right)  = \int \widetilde f \widetilde d\la^{\widetilde u}$ is continuous, and
		\item [(c)]  (left invariance) for any $\widetilde x\in \widetilde \G$ and any $\widetilde f \in  C_c(\widetilde \G )$, $\int \widetilde  f \left(  \widetilde x \widetilde y \right)  d\la^{s(\widetilde x)}\left( \widetilde y\right)  =
		\int \widetilde f\left( \widetilde y\right) d\la^{d\left( \widetilde x\right) }\left( \widetilde y\right) $.
	\end{enumerate}
	From these circumstances it turns out that $\left\{\widetilde  \la^{\widetilde u}\left| \widetilde u \in \widetilde \G^0\right.\right\}$ is a left  Haar system for $\widetilde{\G}$.
\end{empt}
\begin{definition}\label{groupoid_haar_lift_defn}
	Under the hypotheses \ref{groupoid_haar_lift_empt} we say that the left  Haar system $\left\{\widetilde  \la^{\widetilde u}\left| \widetilde u \in \widetilde \G^0\right.\right\}$ is the $p$-\textit{lift} of $\left\{\la^u \left| u \in \G^0\right.\right\}$.
\end{definition}

\begin{empt}
	Let both $\G^n$ and $\widetilde \G^n$ be explained in \ref{groupoig_gn_empt}. If $\left(\widetilde x_0,...,  \widetilde x_{n-1}\right) \in \widetilde \G^n$ then from the Definition \ref{groupoid_lift_defn}  it follows that $\left(q_\G\left( \widetilde x_0\right) ,...,  q_\G\left( \widetilde x_{n-1}\right)\right)\in G^n$ ($q_\G$ is given by the equation \eqref{groupoid_lift_eqn}). If $\sigma$ is a continuous 2-cocycle in $Z^2\left(\G, \T\right)$ (cf. Definition \ref{groupoid_cocycle_defn}) then there is a map
	\be\label{groupoid_lift_cocyle_eqn}
	\begin{split}
		\widetilde \sigma: \widetilde \G^2 \to \T,\\
		\left(\widetilde x_0, \widetilde x_1 \right) \mapsto \left(q_\G\left( \widetilde x_0\right), q_\G\left( \widetilde x_1\right)\right). 
	\end{split}
	\ee
	Moreover from \eqref{groupoid_b_eqn} it follows that 
	$$
	\forall \left(\widetilde x_0, \widetilde x_1 \right) \in \widetilde \G^2\quad \dl^2 \widetilde \sigma \left(\widetilde x_0, \widetilde x_1 \right)=\dl^2  \sigma \left(q_\G\left( \widetilde x_0\right), q_\G\left( \widetilde x_1\right)\right)= 1,
	$$
	i.e. $\widetilde \sigma$ is  a continuous 2-cocycle (cf. Definition \ref{groupoid_cocycle_defn}).	
\end{empt}
\begin{definition}\label{groupoid_cocycle_lift_defn}
	We say that the continuous 2-cocycle $\widetilde \sigma\in Z^2\left(\G, \T\right)$ is the $p$-lift of $\sigma$.
\end{definition}
\begin{exercise}\label{groupoid_loc_comp_lift_exer}
	Prove that if a groupoid $\G$ is locally compact then any $p$-lift of $\G$ is locally compact with respect to $p$-topology (cf. Definition \ref{groupoid_lift_defn}).
\end{exercise}

\section{Algebraic construction}
\begin{empt}\label{groupoid_rho_empt}
	If $\G$ is locally compact groupoid with a continuous 2-cocycle in $Z^2\left(\G, \T\right)$ (cf. Definition \ref{groupoid_cocycle_defn}) and a left  Haar system $\left\{\la^u\left| u \in \G^0\right.\right\}$  (cf. Definition \ref{groupoid_haar_defn}) then there is a $*$-algebra $C_c\left(\G, \sigma \right)$ with given by \eqref{groupoid_*_c_eqn} operations.
	If  $\rho: 	C_c\left(\G , \sigma \right)\to B\left(\H \right)$ is a representation in the sense of the Definition \ref{groupoid_representation_defn} then there is a $C^*$-seminorm 
	\be\label{groupoid_semin_eqn}
	\begin{split}
		\left\|\cdot  \right\|_\rho :   	C_c\left(\G , \sigma \right)\to \R,\\
		a \mapsto \left\|\rho\left( a \right)  \right\|
	\end{split}
	\ee
	We suppose that $\rho$ is \textit{faithful}, i.e. $a \neq 0\quad \Rightarrow \quad \rho_c\left(a \right) \neq 0$. The completion of  $C_c\left(\G , \sigma \right)$ with respect to 	$\left\|\cdot  \right\|_\rho$ is a $C^*$-algebra denoted by $C^*_\rho\left(\G , \sigma \right)$
\end{empt}
\begin{definition}\label{groupoid_rho_defn}
	Under the hypotheses \ref{groupoid_rho_empt} we say that the $C^*$-algebra $C^*_\rho\left(\G , \sigma \right)$ is the $\rho$-\textit{completion} of  $C_c\left(\G , \sigma \right)$.
\end{definition}
\begin{remark}\label{groupoid_blowing_rem}
	From the Lemma \ref{groupoid_mult_repr_lem} it follows that there is  natural Hausdorff blowing-up (cf. Definition \ref{blowing_defn})
	\be\label{groupoid_blowing_eqn}
	C_0\left(\G^0 \right)\hookto M\left( C^*_\rho\left(\G , \sigma \right)\right)  
	\ee
\end{remark}

\begin{empt}\label{groupid_tens_empt}
	Let $\G$ be a locally groupoid, and let  $q: \widetilde \G^0 \to \G^0$ be $\G$-{regular} covering (cf. Definition \ref{groupoid_lift_defn}). Let $$\left(\	\left\{\left(\sU_\a, \sV_\a, s_\a\right)\right\}_{\a \in \mathscr A}, \left\{\left(\widetilde \sU_{\widetilde\a}, \widetilde \sV_{\widetilde\a}, \widetilde s_{\widetilde\a}\right)\right\}_{\widetilde \a \in \widetilde{\mathscr A}} \right)$$
	be a $p$-{covering} (cf. Definition \ref{top_cov_defn}).
	One can proof that following conditions hold.
	\begin{itemize}
		\item Similarly to \eqref{blowing_ajk_eqn}  any $a \in C_c\left(\G , \sigma \right)$ can be represented by following way
		\be\label{groupoid_faf_eqn}
		a = \sum^{\substack{j = m\\k=n}}_{\substack{j = 1\\k=1}}  f _{j}af_k
		\ee
		where there are open subsets $\sU'_j, \sU''_k \in \left\{{\sU}_\a\right\}_{\a\in \mathscr A}$ such that $\supp f_j \subset \sU'_j$ and $\supp f_k \subset \sU'_k$.
		\item Similarly to \eqref{blowing_tajk_eqn} any $\widetilde a \in C_c\left(\widetilde \G, \widetilde\sigma \right)$
		\be\label{groupoid_tfaf_eqn}
		\widetilde a = \sum^{\substack{j = m\\k=n}}_{\substack{j = 1\\k=1}}  \widetilde f _{j}\widetilde a\widetilde f_{k} 
		\ee
		where there are open subsets $\widetilde \sU'_j, \widetilde\sU''_k \in \left\{{\widetilde\sU}_{\widetilde\a}\right\}_{{\widetilde\a}\in \widetilde{\mathscr A}}$ such that $\supp \widetilde f_j \subset \widetilde\sU'_j$ and $\supp \widetilde f_k \subset \widetilde\sU'_k$.
	\end{itemize}
Below 	for any $a \in C_c\left(\G , \sigma \right)$ and $\widetilde f \in C_c\left( \widetilde \G^0\right)$ we will define a product $\widetilde f a \in C_c\left(\widetilde \G, \widetilde\sigma \right)$. 
 Using the decomposition  \ref{groupoid_decomp_empt} one can assume that there is an open Hausdorff subset $\sV \subset \G$ such that $a$ is represented by $f_{a}\in C_c\left( \sV\right)$. There is a subordinated  to the family $\left\{\widetilde \sU_{\widetilde \a}\right\}$ covering sum of $\supp \widetilde f$ (cf. Definition \ref{top_covering_sum_defn}), i.e.
	$$
	\widetilde f = \sum_{j=1}^m \widetilde f_j,
	$$
	where for any $j=1,..., m$ there is $\widetilde\sU_j \in \left\{\widetilde \sU_{\widetilde \a}\right\}$ such that $\supp \widetilde f_j\subset \widetilde\sU_j$. Similarly an application of the covering sum of $\supp \widetilde f$ yields the following 
	$$
	a  = \sum_{k=1}^n f_k a.
	$$
	where for any $k=1,..., n$ there is $\sU_k \in \left\{ \sU_{ \a}\right\}$.
	If $q\left( \widetilde\sU_j \right)\cap  \sU_k = \emptyset$ then we set $\widetilde f_j\left( f_k a\right) = 0$. Otherwise there is $\widetilde\sU'_k \in \left\{\widetilde \sU_{\widetilde \a}\right\}$ such that $q\left(\widetilde\sU'_k  \right) = \sU_k$ and $\widetilde\sU'_k\cap \widetilde\sU_j\neq 0$. Let $\widetilde f'_k \bydef \lift^q_{\widetilde\sU'_k}\left( f_k\right)\in C_c\left(\widetilde \G^0 \right)  $ be the  $q$-${\widetilde\sU'_k}$-lift (cf. Definition \ref{top_lift_desc_defn}). Denote by $q_\G: \widetilde \G \to \G$ the natural continuous map. If $a$ is represented by $f_a \in C_c\left(\G \right)$ and there is an open Hausdorff subspace  $\sV \subset \G$ such that $\supp f_a \subset \sV$  then the set
	$$
	\widetilde\sV \bydef q^{-1}_\G\left(\sV \right) \cap r^{-1}\left( \widetilde\sU'_k\right) \in \widetilde \G
	$$
	is mapped homeomorphically onto $q_\G\left( \widetilde\sV\right)$. If $ \widetilde f'_a \bydef \lift^{q_\G}_{\widetilde\sV}\left( f_a\right) \in C_c\left(\widetilde G\right)$ then $\widetilde f'_a$ uniquely defines $\widetilde a' \in \widetilde \G$. We set
	\be\label{groupoid_fa_eqn}
	\widetilde f_j\left( f_k a\right)\bydef\widetilde f_j  \widetilde a' 
	\ee 
\end{empt}
\begin{lemma}\label{groupoid_tens_lem}
If a covering  $q: \widetilde \G^0 \to \G^0$ be $\G$-{regular}  (cf. Definition \ref{groupoid_lift_defn}) then one has.
	\begin{enumerate}
		\item [(i)] The equation \eqref{groupoid_fa_eqn} yields a surjective map 
		\be\label{groupoid_tensor_eqn}
		\phi:	C_c\left(\widetilde \G^0 \right)\otimes_{C_0\left(\G^0 \right) }C_c\left( \G, \sigma \right)\xrightarrow{\approx} C_c\left(\widetilde \G, \widetilde\sigma \right)
		\ee
		where an algebraic tensor product is implied.
		\item[(ii)] There are natural  right and left  actions 
		\bean
		C_c\left(\widetilde \G, \widetilde\sigma \right)\times C_c\left( \G, \sigma \right)\to  C_c\left(\widetilde \G, \widetilde\sigma \right),\\
		C_c\left( \G, \sigma \right)\times  C_c\left(\widetilde \G, \widetilde\sigma \right)\to  C_c\left(\widetilde \G, \widetilde\sigma \right).
		\eean
	\end{enumerate}
\end{lemma}
\begin{proof}
	(i)  
	 From  \eqref{groupoid_tfaf_eqn}  it follows that one should prove that  any $\widetilde f' \widetilde a \widetilde f''$ with $\widetilde{\sU}', \widetilde{\sU}''\in \left\{{\widetilde\sU}_{\widetilde\a}\right\}$ such that $\supp \widetilde f' \subset \widetilde \sU'$ and $\supp \widetilde f'' \subset \widetilde \sU''$  there is $\widetilde f \otimes a \in C_c\left(\widetilde \G^0 \right)\otimes_{C_0\left(\G^0 \right) }C_c\left( \G, \sigma \right)$ which satisfies to following condition $\phi\left(\widetilde f \otimes a \right) = \widetilde f' \widetilde a \widetilde f''$. The set $\widetilde{\sU}''$ is mapped homeomorphically onto $q\left(\widetilde{\sU}''\right)$ so the set $s^{-1}\left(\widetilde{\sU}''\right)$ is mapped homeomorphically onto $q_\G\left(s^{-1}\left(\widetilde{\sU}''\right)\right)$, so if $\widetilde a \widetilde f''$ corresponds to $\widetilde h \in C_c\left(\widetilde \G \right)$ then one has a $q_\G$ descent $h \bydef \desc_{q_\G}\left( \widetilde h\right)$ (cf. Definition \ref{top_lift_desc_defn}) if $a \in  C_c\left( \G, \sigma \right)$ corresponds to $h$ then from the construction \ref{groupid_tens_empt} it follows that 
	$$
	\widetilde f' \widetilde a \widetilde f''= \widetilde f' a. 
	$$
	(ii) From (i) it follows that 
	$$
	\forall \widetilde a \in C_c\left(\widetilde \G, \widetilde\sigma \right)~~ \exists \sum_{j = 1}^n \widetilde f_j \otimes a_j \in C_c\left(\widetilde \G^0 \right)\otimes_{C_0\left(\G^0 \right) }C_c\left( \G, \sigma \right)~~ \widetilde a= \phi\left( \sum_{j = 1}^n \widetilde f_j \otimes a_j\right).
	$$
	For any $a \in  C_c\left( \G, \sigma \right)$ we define
	$$
	\widetilde a a\bydef \phi\left( \sum_{j = 1}^n \widetilde f_j \otimes a_ja\right)
	$$
	so one has a right  action  $C_c\left(\widetilde \G, \widetilde\sigma \right)\times C_c\left( \G, \sigma \right)\to  C_c\left(\widetilde \G, \widetilde\sigma \right)$. An application of $*$-operation yields the left  action $C_c\left( \G, \sigma \right)\times  C_c\left(\widetilde \G, \widetilde\sigma \right)\to  C_c\left(\widetilde \G, \widetilde\sigma \right)$.
\end{proof}

\begin{corollary}\label{groupoid_tens_cor}
	Use notation of the Definition \ref{blowing_a_regular_defn}.
	Let $\rho: {C_c\left(\G , \sigma \right)} \to B\left( \sH\right) $ be a faithful,   representation in the sense of the Definition \ref{groupoid_representation_defn}.
Let $q: \widetilde \G^0 \to \G^0$ let a $\G$-regular transitive covering (cf. Definition \ref{groupoid_lift_defn}) and $C^*_\rho\left(\G , \sigma \right)$ is the $\rho$-\textit{completion} of  $C_c\left(\G , \sigma \right)$. If $\mathscr L^2\left(\widetilde \G^0 \right)_{C^*_\rho\left(\G , \sigma \right)}$ is the $\widetilde \G^0$-$C^*_\rho\left(\G , \sigma \right)$-module (cf. Definition \ref{blowing_lift_hm_defn}) then there is the natural homomorphisms 
\be\label{groupoid_prods_eqn}
	\phi:C_c\left( \widetilde \G^0 \right)\otimes_{C_0\left(\G^0\right) } C_c\left(\G , \sigma \right)\to  \End^*_A\left( \mathscr L^2\left(\widetilde \G^0 \right)_{C^*_\rho\left(\G , \sigma \right)}\right)
\ee
of right $C_c\left(\G , \sigma \right)$-modules such that 	$\varphi\left( C_c\left( \widetilde \G^0 \right)\otimes_{C_0\left(\G^0\right) } C_c\left(\G , \sigma \right)\right)$  is a $*$-subalgebra of $\End^*_{C^*_\rho\left(\G , \sigma \right)}\left( \mathscr L^2\left(\widetilde \G^0 \right)_{C^*_\rho\left(\G , \sigma \right)}\right)$. 
\end{corollary}
\begin{proof}
From the Lemma \ref{groupoid_tens_lem} it follows that the generated by elements \ref{groupoid_prods_eqn}  $C_c\left(\G , \sigma \right)$-module is isomorphic to $*$-algebra  $C_c\left(\widetilde \G, \widetilde\sigma \right)$.
\end{proof}
\begin{theorem}\label{groupoid_lift_thm}
If   $q: \widetilde \G^0\to  \G^0$  is a $\G$-regular transitive covering (cf. Definition \ref{groupoid_lift_defn}) then  $q: \widetilde \G^0\to  \G^0$ is $C^*_\rho\left(\G , \sigma \right)$-regular (cf. Definition \ref{blowing_a_regular_defn}) 
\end{theorem}
\begin{proof}
	Since $C_c\left(\G , \sigma \right)$ is dense in $C^*_\rho\left(\G , \sigma \right)$ the map \eqref{groupoid_prods_eqn} can be uniquely extended up to a map
	$$
C_c\left( \widetilde \G^0 \right)\otimes_{C_0\left(\G^0\right) } C^*_\rho\left(\G , \sigma \right)\to  \End^*_A\left( \mathscr L^2\left(\widetilde \G^0 \right)_A\right)	
	$$
	so one has a map
	$$
\varphi : 	C_c\left( \widetilde \G^0 \right)\otimes_{C_0\left(\G^0\right) } K\left( C^*_\rho\left(\G , \sigma \right)\right) \to  \End^*_A\left( \mathscr L^2\left(\widetilde \G^0 \right)_A\right)	
	$$
	On the other hand both $C_c\left(\G , \sigma \right)$ and $K\left( C^*_\rho\left(\G , \sigma \right)\right) $ are dense in $C^*_\rho\left(\G , \sigma \right)$ so the  $C^*$-norm closure  of $\phi\left( C_c\left( \widetilde \G^0 \right)\otimes_{C_0\left(\G^0\right) } C_c\left(\G , \sigma \right)\right)$ coincides with the $C^*$-norm closure of  $\varphi\left( C_c\left( \widetilde \G^0 \right)\otimes_{C_0\left(\G^0\right) } K\left( C^*_\rho\left(\G , \sigma \right)\right)\right)$. However from the Corollary \ref{groupoid_tens_cor} it turns out that the $C^*$-norm completion of $\phi\left( C_c\left( \widetilde \G^0 \right)\otimes_{C_0\left(\G^0\right) } C_c\left(\G , \sigma \right)\right)$ is a $C^*$-subalgebra of $\End^*_A\left( \mathscr L^2\left(\widetilde \G^0 \right)_A\right)$, i.e. one has a specialization of the Definition \ref{blowing_a_regular_defn}.
\end{proof}
\begin{empt}\label{groupoid_alg_lift_empt}
Under the hypotheses of the Theorem \ref{groupoid_lift_thm} from the Lemma \ref{blowing_induced_representation_lem} it follows that there is a faithful (cf. Definition \ref{faithful_representation_defn}) 
\bean
\widetilde{\rho}\bydef\mathscr L^2\left(\widetilde \G^0 \right)_{C^*_\rho\left(\G , \sigma \right)}\text{-}\Ind^{C^*_\rho\left(\G , \sigma \right)}_{{C^*_\rho\left( \G ,  \sigma \right)}_0\left(\widetilde\G^0 \right) }\rho: C^*_\rho\left(\G , \sigma \right)_0\left(\widetilde\G^0 \right)\hookto B\left(\widetilde{\H} \right)
\eean
given by \eqref{induced_representation_eqn}. On the other hand there is an inclusion $C_c\left(\widetilde \G , \widetilde \sigma \right)\subset C^*_\rho\left(\G , \sigma \right)_0\left(\widetilde\G^0 \right)$, so there is a representation 
\be\label{groupoid_rho_c_eqn}
\widetilde{\rho} : C_c\left(\widetilde \G , \widetilde \sigma \right)\hookto B\left(\widetilde{\H} \right).
\ee

\end{empt}

\begin{definition}\label{groupoid_alg_lift_defn}
	Under the hypotheses \ref{groupoid_alg_lift_empt} we say that the completion  $C^*_{\widetilde\rho}\left(\widetilde \G , \widetilde \sigma \right)$ of $C_c\left(\widetilde \G , \widetilde \sigma \right)$ with respect to a given by
\be\label{groupoid_alg_lift_eqn}
\begin{split}
\left\| \cdot \right\|_{\widetilde\rho} : C_c\left(\widetilde \G , \widetilde \sigma \right)\to \R,\\
\widetilde a \mapsto \left\| \widetilde\rho\left(\widetilde a \right)  \right\|
\end{split}
\ee
$C^*$-norm is $q$-\textit{lift} of $C^*_\rho\left(\G , \sigma \right)$.
\end{definition}
\begin{theorem}\label{groupoid_blowing_thm} If   $q: \widetilde \G^0\to  \G^0$  is a $\G$-regular transitive covering (cf. Definition \ref{groupoid_lift_defn}) then following conditions hold:
\begin{enumerate}
	\item [(i)] the given by \eqref{groupoid_rho_c_eqn} representation $\widetilde{\rho}$ satisfies to the Definition \ref{groupoid_representation_defn},
	\item[(ii)] there is the natural $*$-isomorphism 
\be\label{groupoid_blowing_fin_eqn}
 C^*_\rho\left(\G , \sigma \right)_0\left(\widetilde\G^0 \right)\cong C^*_{\widetilde\rho}\left(\widetilde \G , \widetilde \sigma \right)
\ee
where $C^*_\rho\left(\G , \sigma \right)_0\left(\widetilde\G^0 \right)$ is explained in the Lemma \ref{blowing_lift_constr_lem}.
\end{enumerate}
\end{theorem}
\begin{proof}
(i) From the Lemma \ref{blowing_induced_representation_lem} it follows that $\widetilde{\rho}_c$ is nondegenerate. If $\widetilde \sZ \subset \widetilde \G$ is compact then from the Theorem \ref{top_compact_img_thm}  it follows that both $\widetilde \sV' \bydef s\left(\widetilde \sZ \right)$ and $\widetilde \sV'' \bydef r\left(\widetilde \sZ \right)$ are compact. Both maps $s: \widetilde \G\to \widetilde G^0$ and $r: \widetilde \G\to \widetilde G^0$  induce  right and left  actions
\bean
 C_0\left( \widetilde\sZ\right)\times C\left( \widetilde\sV'\right) \to C_0\left( \widetilde\sZ\right),\\
C\left( \widetilde\sV''\right) \times C_0\left( \widetilde\sZ\right)\to C_0\left( \widetilde\sZ\right).
\eean
If  $\left\{\widetilde f_\a\right\}_{\a \in \mathscr A}\subset  C_0\left( \widetilde\sZ\right)$ is  uniformly convergent net and both  $\sum_{j = 1}^m \widetilde f'_j$, $~\sum_{k = 1}^n \widetilde f'_j$ is  dominated  to $q$ (cf. Definition \ref{top_covering_sum_subordinated_defn}) covering sums of $\widetilde \sV'$ and  $\widetilde \sV''$  then there are families $\left\{\widetilde \sV'_j\right\}_{j\in \{1,...,m\}}$ and  $\left\{\widetilde \sV''_k\right\}_{k\in \{1,...,n\}}$ such that:
 \begin{itemize}
 	\item $\supp \widetilde f'_j \subset \widetilde \sV'_j$ and  $\supp \widetilde f''_k\subset \widetilde \sV''_k$,
 	\item $\widetilde \sV'_j$ is mapped homeomorphically onto $\sV'_j \bydef q\left( \widetilde \sV_j\right)$ then for all $j = 1,..., m$,
 	\item $\widetilde \sV''_k$ is mapped homeomorphically onto $\sV''_k \bydef q\left( \widetilde \sV''_k\right)$ then for all $k = 1,..., m$, 
 	\item the net $\left\{\widetilde f''_k \widetilde f_\a\widetilde f'_j \right\}\subset C_c\left(s^{-1}\left(\widetilde \sV'_j \right)\cap r^{-1}\left(\widetilde \sV''_k \right)\cap \widetilde\sZ  \right) $ is uniformly convergent. 
 \end{itemize}
If $j=1,...m$ and $k= 1,...,n$ then a set $\widetilde \sZ_{jk} \bydef s^{-1}\left(\widetilde \sV'_j \right)\cap r^{-1}\left(\widetilde \sV''_k \right)\cap \widetilde\sZ$ is mapped homeomorphically onto $\sZ_{jk} \bydef q_\G\left(\widetilde \sZ_{jk} \right)$. If $f^{jk}_\a \bydef \desc_{q_\G} \left(\sqrt{\widetilde f''_k} \widetilde f_\a \sqrt{\widetilde f'_j} \right)\in C_0\left(\sZ_{jk} \right) $ then the net $\left\{f^{jk}_\a\right\}$ is uniformly convergent. The net $\left\{f^{jk}_\a \right\}\subset C_0\left(\sZ_1 \right)$.  The net $\left\{f^{jk}_\a \right\}$ yields a net $\left\{a^{jk}_\a \right\} \subset C_c\left(\G , \sigma \right)$. From the Definition \ref{groupoid_representation_defn} it follows that then net $\left\{\rho\left( a^{jk}_\a\right) \right\}$ is convergent with respect to the weak topology (cf. Definition \ref{weak_topology_defn}) of $B\left( \H\right)$. One can proof that
$$
\widetilde f_\a =\sum_{\substack{j=1\\k=1}}^{\substack{j=m\\k=n}} \sqrt{\widetilde{f}''_k}a^{jk}_\a \sqrt{\widetilde{f}'_j}
$$

From the Definition \ref{induced_representation_defn} it follows that there is a dense inclusion
$$
\mathscr L^2\left(\widetilde \G^0 \right) \otimes_{C_0\left(\G^0 \right) }{C^*_\rho\left(\G , \sigma \right)}\subset \widetilde \H
$$
and since $C_c\left(\widetilde \G^0 \right)$ is dense in $\mathscr L^2\left(\widetilde \G^0 \right)$ one has a dense inclusion
$$
C_c\left(\widetilde \G^0 \right) \otimes_{C_0\left(\G^0 \right) }{C^*_\rho\left(\G , \sigma \right)}\subset \widetilde \H
$$
If $\widetilde \xi, ~\widetilde \eta \in \widetilde \H$ which correspond to $\widetilde f' \otimes \xi,~ \widetilde f'' \otimes \eta \in C_c\left(\widetilde \G^0 \right) \otimes_{C_0\left(\G^0 \right) }{C^*_\rho\left(\G , \sigma \right)}$  then
\be\label{groupoid_sum_eqn}
\left(\widetilde \eta, \widetilde a_\a \widetilde \xi\right)_{\widetilde \H} = \sum_{\substack{j=1\\k=1}}^{\substack{j=m\\k=n}}\left(\desc_q\left(\sqrt{\widetilde{f}''_k}\widetilde f'' \right)\eta,\rho\left(  a^{jk}_\a\right)  \desc_q\left(\sqrt{\widetilde{f}'_j}\widetilde f' \right)\xi\right)_\H 
\ee
For any $j = 1,..., m$ and $k=1,...,n$ the net $\left\{\rho\left( a^{jk}_\a\right) \right\}$ is convergent with respect to the weak topology of $B\left(\H\right)$ so the net
$$
\left\{\left(\desc_q\left(\sqrt{\widetilde{f}''_k}\widetilde f'' \right)\eta,\rho\left(  a^{jk}_\a\right)  \desc_q\left(\sqrt{\widetilde{f}'_j}\widetilde f' \right)\xi\right)_\H\right\}
$$
is convergent. From \eqref{groupoid_sum_eqn} it turns out that the net $\left\{\left(\widetilde \eta, \widetilde a_\a \widetilde \xi\right)_{\widetilde \H}\right\}$ is convergent. Using this fact can proof that a net $\left\{\left(\widetilde \eta, \widetilde a_\a \widetilde \xi\right)_{\widetilde \H}\right\}$ is convergent for all $\widetilde \xi,  \widetilde \eta \in C_c\left(\widetilde \G^0 \right) \otimes_{C_0\left(\G^0 \right) }{C^*_\rho\left(\G , \sigma \right)}$. Since $C_c\left(\widetilde \G^0 \right) \otimes_{C_0\left(\G^0 \right) }{C^*_\rho\left(\G , \sigma \right)}$ is dense in $\widetilde \H$  a net $\left\{\left(\widetilde \eta, \widetilde a_\a \widetilde \xi\right)_{\widetilde \H}\right\}$ is convergent for all $\widetilde \xi,  \widetilde \eta \in\widetilde \H$, i.e. a net is convergent $\left\{\widetilde{\rho} \left( \widetilde a_\a\right) \right\}$ with respect to the weak topology of $B\left(\widetilde{\H} \right)$.\\
(ii) If we consider a norm
\bean
\left\|\cdot  \right\|_{\widetilde \rho} : C_c\left(\widetilde \G , \widetilde \sigma \right)\to \R,\\
\widetilde a \mapsto \left\|\widetilde \rho_c\left(\widetilde a\right) \right\| 
\eean
then from the Definition \ref{groupoid_rho_defn} it follows that  $C^*_{\widetilde\rho}\left(\widetilde \G , \widetilde \sigma \right)$ is the completion of $C_c\left(\widetilde \G , \widetilde \sigma \right)$ with respect to the norm $\left\|\cdot  \right\|_{\widetilde \rho}$. On the other hand from the Lemma \ref{blowing_induced_representation_lem} it turns out that $C^*_\rho\left(\G , \sigma \right)_0\left(\widetilde\G^0 \right)$ is the  completions of $C_c\left(\widetilde \G , \widetilde \sigma \right)$ with respect to the norm $\left\|\cdot  \right\|_{\widetilde \rho}$.
\end{proof}
\begin{empt}\label{groupoid_inv_empt}
	According to the  construction \ref{groupoid_alg_lift_empt} a representation $\rho: C_c\left( \G ,  \sigma \right)\hookto B\left({\H}\right) $ yields a representation 
  $\widetilde{\rho} : C_c\left(\widetilde \G , \widetilde \sigma \right)\hookto B\left(\widetilde{\H} \right)$. Here we obtain inverse operation.  Suppose that the representation $\widetilde{\rho} : C_c\left(\widetilde \G , \widetilde \sigma \right)\hookto B\left(\widetilde{\H} \right)$ is known. A $q$-lift (cf. Definition \ref{blowing_lift_hom_defn}) $$C^*_\rho\left(\G , \sigma \right)_b\left( q\right) : C^*_\rho\left(\G , \sigma \right)\hookto M\left( C^*_\rho\left(\G , \sigma \right)_0\left(\widetilde\G^0 \right)=  C^*_{\widetilde \rho}\left(\widetilde \G , \widetilde \sigma \right)\right)$$ and the Definition \ref{multiplier_el_defn} yield a representation
 \be\label{groupoid_rh_eqn}
 \rho' : C_c\left( \G ,  \sigma \right)\to B\left(\widetilde \H \right) 
\ee

\end{empt}

\begin{exercise}
 Under the hypotheses \ref{groupoid_inv_empt} prove that
		$$
	\forall a \in 	C_c\left( \G ,  \sigma \right)\quad \Rightarrow \quad \left\|\rho\left(a \right)  \right\|= \left\|\rho'\left(a \right)  \right\| 
		$$
\end{exercise}
\begin{empt}\label{groupoid_red_empt}
Consider the situation \ref{groupoid_inv_empt} end suppose that $\widetilde \rho: C_c\left(\widetilde \G , \widetilde \sigma \right)\to B\left(\widetilde \H \right)$ corresponds to the  given by \ref{groupoid_red_norm_eqn} $C^*$-norm, i.e. $ C^*_{\widetilde \rho}\left(\widetilde \G , \widetilde \sigma \right)$  equals to the {reduced $C^*$ -algebra} $C^*_r\left(\widetilde \G, \widetilde \sigma \right)$ of $\G$ (cf. Definition \ref{groupoid_red_defn}). For any $\widetilde{u}\in \widetilde{G}^0$ there is a given by \eqref{groupoid_reg_eqn}representation $\widetilde L_{\widetilde{u}}: C_c\left(\widetilde \G , \widetilde \sigma \right)\to B\left(\widetilde L^2\left( \widetilde u\right) \right)$.
\end{empt}
\begin{exercise}
Under the hypothesis 
\begin{enumerate}
	\item For any $\widetilde{u}\in \widetilde{G}^0$ there is a  representation
	$$
	L_{\widetilde{u}}\bydef  C_c\left(\G, \sigma \right)\to B\left(L^2\left( {\widetilde{u}} \right) \right) 
	$$
	which complies with the inclusion $C_c\left(\G, \sigma \right) \hookto M\left( C^*_r\left(\widetilde \G, \widetilde \sigma \right)\right)$ and the representations  \bean
	M\left( C^*_r\left(\widetilde \G, \widetilde \sigma \right)\right)\to B\left( \widetilde \H\right),\\ 	L_{\widetilde{u}}:  C_c\left(\G, \sigma \right)\to B\left(L^2\left( {\widetilde{u}} \right) \right).
	\eean
	\item 
	$
	\forall \widetilde{u}\in \widetilde G^0\quad \forall  a \in  C^*_r\left( \G,  \sigma \right)\quad \left\|L_{\widetilde{u}}\left(a \right)  \right\|= \left\|L_{q\left( \widetilde{u}\right) }\left(a \right)  \right\|
	$
	where $q:   \widetilde \G^0 \to \G^0$ is the natural covering.
	\item If $ \rho' : C_c\left( \G ,  \sigma \right)\to B\left(\widetilde \H \right)$ is given by then
	$$
	\forall  a \in  C^*_r\left( \G,  \sigma \right)\quad \left\|\rho'\left(a \right)  \right\| = \left\|a \right\|_r
	$$
	where $\left\|a \right\|_r$ is given by the equation \eqref{groupoid_red_norm_eqn}.
	\item Under the hypotheses \ref{groupoid_blowing_thm} one has
	\be\label{groupoid_rho_red_eqn}
	C^*_{\widetilde \rho }\left(\widetilde \G, \widetilde \sigma \right)= C^*_r\left(\widetilde \G, \widetilde \sigma \right)\quad \Rightarrow \quad 	C^*_{ \rho }\left( \G,  \sigma \right)= C^*_r\left( \G,  \sigma \right).
	\ee
\end{enumerate}
\end{exercise}
\begin{definition}\label{groupoid_support_defn}
Let $\G$ be a Hausdorff groupoid.
If   $C^*_\rho\left(\G , \sigma \right)$ is the $\rho$-{completion} of  $C_c\left(\G , \sigma \right)_+$ (cf. Definition \ref{groupoid_rho_defn}), and  $a \in C^*_\rho\left(\G , \sigma \right)$ then we say that a set $\sU\subset \G$ is \textit{orthogonal} to $a$ if one has
$$
\forall b \in C_c\left(\G , \sigma \right)_+\quad \supp b \in \sU \quad \Rightarrow \quad ab = ba = 0.
$$
The \textit{groupoid support} of $a$ is the closure of the set
$$
\G \setminus \bigcup \left\{\left.\sU \subset \G~ \right| ~\sU  \text{ is orthogonal to }  a\right\}
$$
We denote it by $\G$-$\supp a$.
\end{definition}
\begin{remark}
The difference between both notions  supports given by the Definitions \ref{blowing_support_defn} and \ref{groupoid_support_defn} will be explained as we go along.
\end{remark}
\begin{remark}
For any $a', a'' \in C^*_\rho\left(\G , \sigma \right)$ one has
\be\label{groupoid_support_supp_eqn}
\supp \left( a' + a'' \right) = \supp a' \cup \supp a''.
\ee
\end{remark}
\begin{remark}
	For any $a \in C^*_\rho\left(\G , \sigma \right)$ one has a  disjoint union 
	$$
	G\text{-}\supp a = \bigsqcup_{\iota \in I} \mathcal W_\iota 
	$$
	where $\mathcal W_\iota$ is a quasi-component of $\supp a$ (cf. Definition \ref{top_quasi_component_defn}). From the Definition \ref{groupoid_rho_defn} one has a $C^*$-norm convergent series 
	\be\label{groupoid_support_dec_eqn}
	a = \sum_{\iota \in I} a_\iota, \quad \quad \forall \iota \in I \quad \G\text{-}\supp a_\iota = \mathcal W_\iota.
	\ee
\end{remark}

\begin{empt}\label{groupoid_decomp_empt}
Any $a \in  C_c\left( \G, \sigma \right)$ equals to a finite sum $a = \sum_{j=1}^n a_j$ such that for all $j=1,..., n$  there is an open Hausdorff subset $\sV_j \subset \G$  and $a_j$ is represented by $f_{a_j}\in C_c\left( \sV_j\right)$.
\end{empt}

\section{Noncommutative coverings of $C^*$-algebras of groupoids}
\paragraph{}
Here we apply the explained in the chapter \ref{blowing_chap} to the investigations of noncommutative coverings of $C^*$-algebras of groupoids.
\begin{empt}\label{groupoid_n_cov_empt}
	Suppose that there are following ingredients:
\begin{itemize}
	\item a connected locally compact groupoid $\G$,
	\item a left  Haar system $\left\{\left.\la^u \right| u \in \G^0\right\}$ for $\G$ (cf. Definition \ref{groupoid_haar_defn}),
	\item a continuous 2-cocycle in $\sigma \in Z^2\left(\G, \T\right)$ (cf. Definition \ref{groupoid_cocycle_defn}).
	\item $*$-algebra $C_c\left(\G, \sigma \right)$ with given by \eqref{groupoid_*_c_eqn} operations,
	\item a full   representation $\rho: 	C_c\left(\G , \sigma \right)\to B\left(\H \right)$ in the sense of the Definition \ref{groupoid_representation_defn} (cf. \ref{groupoid_rho_empt}), 
	\item a $\rho$-\textit{completion} $C^*_\rho\left(\G , \sigma \right)$ of  $C_c\left(\G , \sigma \right)$ (cf. Definition \ref{groupoid_rho_defn}),
		\item  a $\G$-regular covering $q: \widetilde \G^0 \to \G^0$ is  (cf. Definition \ref{groupoid_lift_defn}), with connected $ \widetilde \G^0$,
	\item the $q$-lift $C^*_{\widetilde\rho}\left(\widetilde \G , \widetilde \sigma \right)$ of $C^*_{\rho}\left( \G ,  \sigma \right)$ (cf. Definition \ref{groupoid_alg_lift_defn}).
\end{itemize}

\end{empt}

\begin{theorem}\label{groupoid_fin_thm}
	Under the hypothesis \ref{groupoid_n_cov_empt}  suppose that the covering $q$ is finite-fold and $C^*_{\rho}\left( \G ,  \sigma \right)_0\left( q\right) : C^*_{\rho}\left( \G ,  \sigma \right)\hookto C^*_{\rho}\left( \G ,  \sigma \right)_0\left( \widetilde{\G}^0\right)$ is the finite-fold $q$-lift (cf. Definition \ref{blowing_finite_lift_defn}). If the space $\G^0$ is locally connected  (cf. Definition \ref{top_locally_connected_defn}) the natural quadruple
	\be\label{groupoid_fin_eqn}
	\left(C^*_{\rho}\left( \G ,  \sigma \right), C^*_{\widetilde\rho}\left(\widetilde \G , \widetilde \sigma \right)\cong C^*_{\rho}\left( \G ,  \sigma \right)_0\left( \widetilde{\G}^0\right) , 	G\left(\left.\widetilde{\G}^0~~\right|~\G^0 \right), C^*_{\rho}\left( \G ,  \sigma \right)_0\left( q\right)   \right)  
	\ee
	is a noncommutative finite-fold covering (cf. Definition \ref{fin_defn}).
\end{theorem}
\begin{empt}\label{groupoid_n_cov_inf_empt}
Under the hypotheses \ref{groupoid_n_cov_empt} suppose that the covering group 	$G\left(\left.\widetilde{\G}^0~~\right|~\G^0 \right)$ (cf. Definition \ref{top_group_of_covering_transformations_defn}) is residually finite (cf. Definition \ref{residually_finite_defn}). From the equation \eqref{groupoid_blowing_fin_eqn} and the Theorems \ref{blowing_sufficient_covering_inf_thm} and \ref{groupoid_fin_thm} it follows that there is an algebraical finite covering category (cf. Definition \ref{algebraical_finite_covering_category_defn})
\be\label{groupoid_algebraical_finite_covering_category_eqn}
\begin{split}
\mathfrak{S}_{C^*_{\rho}\left( \G ,  \sigma \right)_0\left(q \right) }\bydef \\
\left\{\left\{C^*_{\rho_\la}\left( \G_\la ,  \sigma_\la \right) \right\}_{\la}, \left\{C^*_{\rho}\left( \G ,  \sigma \right)_0\left( q^\nu_\mu\right)  :C^*_{\rho_\mu}\left( \G_\mu ,  \sigma_\mu \right)\hookto C^*_{\rho_\nu}\left( \G_\nu ,  \sigma_\nu \right)\right\}_{\substack{\mu, \nu \in \La\\\mu \le \nu}}\right\}
\end{split}
\ee	
\end{empt}

\begin{theorem}\label{grupoid_intinite_covering_thm}
	Under  hypotheses \ref{groupoid_n_cov_empt} and \ref{groupoid_n_cov_inf_empt} if $\overline q: \overline \G^0  \to  \G^0$  is {disconnected covering of} $q : \widetilde \G^0 \to \G^0$ (cf. Definition \ref{top_disconnected_defn}) and algebraical finite covering category 
	$
	\mathfrak{S}_{C^*_{\rho}\left( \G ,  \sigma \right)_0\left( q\right) }$
	is given by \eqref{groupoid_algebraical_finite_covering_category_eqn} then a triple
		\be\label{groupoid_inf_eqn}
	\left(C^*_{\rho}\left( \G ,  \sigma \right), C^*_{\widetilde\rho}\left(\overline \G , \overline \sigma \right)\cong C^*_{\rho}\left( \G ,  \sigma \right)_0\left( \overline{\G}^0\right) , 	G\left(\left.\overline{\G}^0~~\right|~\G^0 \right)\right)  
	\ee 
	 is a {pre}-{covering of algebraical finite covering category}  $\mathfrak{S}_{C^*_{\rho}\left( \G ,  \sigma \right)_0\left( q\right) }$ (cf. Definition \ref{algebraical_finite_covering_category_defn}).	
	
\end{theorem}
\begin{proof}
	Follows from the Theorem \ref{blowing_sufficient_covering_inf_thm}
\end{proof}
\section{Coverings  of foliations}
\subsection{Basic constructions}
\paragraph{}
Any foliated manifold $\left(M, \F\right)$ yields a holonomy groupoid $\G\left(M, \F\right)$ (cf. Definition \ref{foliated_manifold_defn}) (cf. Definition \ref{foli_graph_defn}). 
If $\left(M, \F\right)$  is a foliated manifold then one has
\bean
\G\left( M, \mathcal{F}\right)^0 = M,\\ 
\eean
where the notation of the Definition \ref{groupoid_defn} is used. Suppose that  $\sigma\in Z^2\left(\G\left( M, \mathcal{F}\right) , \T\right) $ is a   2-cocycle on $\G\left( M, \mathcal{F}\right)$, and define
\be\label{foli_sigma_eqn}
C_c\left( M, \mathcal{F}, \sigma\right) \bydef C_c\left(\G\left( M, \mathcal{F}\right), \sigma \right).
\ee
If  $\rho: C_c\left( M, \mathcal{F}, \sigma \right)\to B\left( \H\right)$ is a faithful representation in the sense of the Definition  \ref{groupoid_representation_defn} then there is $\rho$-{completion} $C^*_\rho\left( M, \mathcal{F}, \sigma \right)$  of  $C_c\left( M, \mathcal{F},\sigma \right)$ (cf. Definition \ref{groupoid_rho_defn}).  There is Hausdorff blowing-up 
\be
C_0\left(M\right)\hookto M\left( C^*_\rho\left( M, \mathcal{F}, \sigma \right)\right) 
\ee 
(cf. Remark \eqref{groupoid_blowing_rem}). 
If $\sigma_{\text{triv}}\in Z^2\left(\G\left( M, \mathcal{F}\right) , \T\right)$ is the trivial 2-cocycle we use the following notation
\be\label{foli_triv_eqn}
C^*_\rho\left( M, \mathcal{F}\right) \bydef C^*_\rho\left( M, \mathcal{F}, \sigma_{\text{triv}}\right).
\ee
For any foliated space $\left(M, \sF\right)$ denote by $M/\sF$ the set of all leaves, and for any leaf $L \in M/\sF$ we select a point $x_L \in L$. We write 
\be\label{foli_leaf_repr_eqn}
\rho_L \stackrel{\text{def}}{=} \rho_{x_L}: C^*_r\left(M, \sF\right) \to B\left( L^2\left(\G_{x_L}\right)\right) 
\ee
where the representation $\rho_{x_L}$ is given by the equation \eqref{foli_repr_eqn}.
From the Theorem \ref{foli_irred_hol_thm} it turns out that the representation \eqref{foli_leaf_repr_eqn} is irreducible if and only if the leaf $L$ has no holonomy. If $M/\sF_{\mathrm{no~hol}} \subset M/\sF$ the subset of leaves which have no holonomy then the spectrum of $C^*_r\left(M, \sF\right)$ coincides with $M/\sF_{\mathrm{no~hol}}$.
The atomic representation $\pi_a$ of  $C^*_r\left(M, \sF\right)$ (cf. Definition \ref{atomic_repr_defn}) is given by
\be\label{foli_atomic_eqn}
\pi_a = \bigoplus_{L \in M/\sF_{\mathrm{no~hol}}} \rho_L: C^*_r\left(M, \sF\right) \hookto B\left( \bigoplus_{L \in M/\sF_{\mathrm{no~hol}}}  L^2\left(\G_{x_L}\right)\right).
\ee

\subsection{Finite-fold coverings}
\begin{exercise}\label{groupoid_foli_r_exer}
Prove that for any foliated manifold $\left(M, \F\right)$ (cf. Definition \ref{foliated_manifold_defn})   that any   regular covering $q:\left(\widetilde{M},~ \widetilde{\mathcal F} \right)\to\left({M},~ {\mathcal F} \right)$ (cf. Definition \ref{foli_reg_cov_defn}) is $\G\left(M, \F\right)$-regular (cf. Definition \ref{groupoid_lift_defn}) and vice versa.
\end{exercise}
\begin{remark}
Under the hypotheses of the Exercise \ref{groupoid_foli_r_exer} the covering $q$ is $C^*_\rho\left( M, \sF, \sigma \right)$ -regular (cf. Definition \ref{blowing_a_regular_defn}), 
\end{remark}
\begin{theorem}\label{groupoid_foli_fin_thm}
If $q:\left(\widetilde{M},~ \widetilde{\mathcal F} \right)\to\left({M},~ {\mathcal F} \right)$ is a 	regular covering of foliated manifolds (cf. Definition \ref{foli_reg_cov_defn}), and an injective $*$-homomorphism $C^*_{\rho}\left( M, \mathcal{F}, \sigma\right)_0\left( q\right) : C^*_{\rho}\left( M, \mathcal{F}, \sigma \right)\hookto C^*_{\rho}\left( M, \mathcal{F},  \sigma \right)_0\left( \widetilde{M}\right)\cong 	C^*_{\widetilde \rho}\left( \widetilde M,\widetilde \F, \widetilde \sigma\right)$ (cf. equation \eqref{groupoid_blowing_fin_eqn}) is the finite-fold $q$-lift (cf. Definition \ref{blowing_finite_lift_defn}) then the natural quadruple
	\be\label{groupoid_foli_fin_eqn}
	\left(C^*_{\rho}\left( M, \mathcal{F}, \sigma\right),  C^*_{\widetilde \rho}\left( \widetilde M,\widetilde \F, \widetilde \sigma\right) , 	G\left(\left.\widetilde{M}~~\right|~M \right), C^*_{\rho}\left( M, \mathcal{F}, \sigma\right)_0\left( q\right)   \right)  
	\ee
	is a noncommutative finite-fold covering (cf. Definition \ref{fin_defn}).
\end{theorem}
\begin{proof}
FOLLOWS FROM THE THEOREM \ref{blowing_sufficient_covering_thm}.
\end{proof}
\subsection{Infinite coverings}
\begin{empt}\label{groupoid_foli_empt}
	Let $\left(M, \F\right)$  is a foliated manifold (cf. Definition \ref{foliated_manifold_defn}),  and  $q:\left(\widetilde{M},~ \widetilde{\mathcal F} \right)\to\left({M},~ {\mathcal F} \right)$ be a   regular covering (cf. Definition \ref{foli_reg_cov_defn}), 
Suppose that the covering group 	$G\left(\left.\widetilde M~~\right|~M \right)$ (cf. Definition \ref{top_group_of_covering_transformations_defn}) is residually finite (cf. Definition \ref{residually_finite_defn}).  Let $\widehat{G}$ be the profinite completion of $G\left(\left.\widetilde M~~\right|~M \right)$ (cf. Definition \ref{profinite_exm}). if $\overline q: \overline M  \to  M$  is the {disconnected covering of} $q : \widetilde M \to M$ (cf. Definition \ref{top_disconnected_defn}From the Theorem \ref{grupoid_intinite_covering_thm} it follows that there is an algebraical finite covering category (cf. Definition \ref{algebraical_finite_covering_category_defn}).
	\be\label{groupoid_foli_algebraical_finite_covering_category_eqn}
	\begin{split}
		\mathfrak{S}_{C^*_{\rho}\left( M, \F, \sigma \right)_0\left(q \right) }\bydef 
		\left\{\left\{A_\la \right\}_{\la\in \La}, \left\{\pi^\mu_\nu   :A_\mu\hookto A_\nu\right\}\right\}_{\substack{\mu, \nu \in \La\\\mu\le \nu}},\\
 \quad A_\la \bydef C^*_{\rho_\la}\left( M_\la ,  \F_\la, \sigma_\la \right),\quad \pi^\mu_\nu \bydef C^*_{\rho}\left( M ,  \F, \sigma \right)_0\left( q^\nu_\mu\right): A_\mu \hookto A_\mu
	\end{split}
	\ee	
	which is a specialization of \ref{groupoid_algebraical_finite_covering_category_eqn} one. Moreover if $C^*_{\overline\rho}\left(\overline M, \overline\F, \overline\sigma \right)$ is the $\overline q$-lift of $C^*_{\rho}\left( M, \F, \sigma \right)$ (cf. Definition \ref{blowing_lift_defn}) then from the Theorem \ref{blowing_sufficient_covering_inf_thm} it follows that 	\be\label{groupoid_foli_inf_eqn}
	\left(C^*_{\rho}\left( M ,  \F, \sigma \right), C^*_{\overline\rho}\left(\overline M , \overline \F, \overline \sigma \right) , 	\widehat{G}  \right)  
	\ee 
	is a {pre}-{covering of algebraical finite covering category}  $\mathfrak{S}_{C^*_{\rho}\left( M ,  \F, \sigma \right)_0\left( q\right) }$ (cf. Definition \ref{algebraical_finite_covering_category_defn}.
Set  $\widehat A \bydef C^*$-$\varinjlim_{\la \in \La} A_\la$, where $C^*$-$\varinjlim$ is the given by the Definition \ref{inductive_lim_non_defn} $C^*$-inductive limit. 
If $\widehat A\hookto B\left(\widehat \H \right)$ is a faithful, nondegenerate representation, and a triple  $\left(C^*_{\rho}\left( M, \F \right), \overline A', \widehat{G} \right)$ is a {pre}-{covering of the algebraical finite covering category} (cf. Definition \ref{algebraical_finite_covering_category_defn} then similarly to the there are faithful, nondegenerate  representations
\bean
\pi_{\overline A'} : \overline A' \hookto B\left(\widehat  \H \right),\\
\pi_{C^*_{\overline\rho}\left(\overline M, \overline\F, \overline\sigma \right) } : C^*_{\overline \rho}\left(\overline M, \overline\F, \overline\sigma \right) \hookto B\left(\widehat  \H \right) 
\eean
which are generalizations of \eqref{groupoid_folia_eqn} and \eqref{groupoid_folia_c_eqn} ones.
\end{empt}

\begin{lemma}\label{grupoid_foli_intinite_covering_lem}
	Under  hypotheses \ref{groupoid_foli_empt}  there is the natural injective $*$-homomorphism 
	$$
	\pi^{C^*_{\overline\rho}\left(\overline M, \overline\F, \overline\sigma \right) }_{\overline A'}:C^*_{\overline\rho}\left(\overline M, \overline\F, \overline\sigma \right)\hookto \overline A'
	$$
	such that the image $	\pi^{C^*_{\overline\rho}\left(\overline M, \overline\F, \overline\sigma \right) }_{\overline A'}\left( C^*_{\overline\rho}\left(\overline M, \overline\F, \overline\sigma \right)\right)$ is a hereditary $C^*$-subalgebra of $\overline A'$.
\end{lemma}
\begin{proof}
Consider a $\overline q$-{covering} 
 $$\left(\	\left\{\left(\sU_\a, \sV_\a, s_\a\right)\right\}_{\a \in \mathscr A}, \left\{\left(\overline \sU_{\overline\a}, \overline \sV_{\overline\a}, \overline s_{\overline\a}\right)\right\}_{\overline \a \in \overline{\mathscr A}} \right)$$ 
 (cf. Definition \ref{top_cov_defn}) such that for any $\overline\a\in\overline{\mathscr A}$ the set $\overline \sU_{\overline\a}$ is a foliation chart (cf. Definition \ref{foli_chart_defn}), i.e. there is a homeomorphism
 $$\textbf{}
 \phi_{\overline{\a}}: \overline \sU_{\overline\a}\xrightarrow{\cong }\R^{\dim\F} \times \R^{\codim\F}.
 $$
 We also suppose that a family $\left\{\overline \sU_{\overline\a}\right\}_{\overline \a \in \overline{\mathscr A}}$ is a regular atlas (cf. Definition  \ref{foli_reg_atlas_defn}).
 From (ii) of the Lemma \ref{blowing_necessary_covering_inf_lem} it turns out that if $\overline C\subset C^*_{\overline\rho}\left(\sU_{\overline\a}, \left.\overline\F\right|_{\sU_{\overline\a}}, \overline\sigma \right)$ is a commutative $C^*$-subalgebra then
  $$
  \pi_{C^*_{\overline\rho}\left(\overline M, \overline\F, \overline\sigma \right) }\left(\overline C \right)\subset  \pi_{\overline A'}\left(  \overline A'\right),
  $$
   The homeomorphism $\phi_{\overline{\a}}$, the Proposition \ref{foli_tens_comp_prop}, and the Corollary  \ref{foli_tens_comp_cor} yield a $*$-isomorphism
  $$
  \varphi_{\overline{\a}}:\K\left( L^2\left( \R^{\dim\F} \right)\right)  \otimes C_0\left( \R^{\codim\F}\right) \xrightarrow{\cong }  C^*_{\overline\rho}\left(\sU_{\overline\a}, \left.\overline\F\right|_{\sU_{\overline\a}}, \overline\sigma \right) 
  $$
  Since  $\K\left( L^2\left( \R^{\dim\F} \right)\right)  \otimes C_0\left( \R^{\codim\F}\right)\cong \K\left( \ell^2\left( \N \right)\right)  \otimes C_0\left( \R^{\codim\F}\right)$ is a continuous trace  $C^*$-algebra it  is generated by its commutative $C^*$-subalgebras, so $C^*_{\overline\rho}\left(\sU_{\overline\a}, \left.\overline\F\right|_{\sU_{\overline\a}}, \overline\sigma \right)$ is also generated by its commutative $C^*$-subalgebras. It follows that 
  $$
    \pi_{C^*_{\overline\rho}\left(\overline M, \overline\F, \overline\sigma \right) }\left(C^*_{\overline\rho}\left(\sU_{\overline\a}, \left.\overline\F\right|_{\sU_{\overline\a}}, \overline\sigma \right) \right)\subset  \pi_{\overline A'}\left(  \overline A'\right),
  $$
  From the Corollary \ref{foli_cov_alg_cor} it follows that the $C^*$-algebra $C^*_{\overline \rho}\left(\overline M, \overline\F, \overline\sigma \right)$ is generated by a family $\left\{C^*_{\overline\rho}\left(\sU_{\overline\a}, \left.\overline\F\right|_{\sU_{\overline\a}}, \overline\sigma \right) \right\}_{\overline \a \in \overline{\mathscr A}}$, it turns out that 
   $$
 \pi_{C^*_{\overline\rho}\left(\overline M, \overline\F, \overline\sigma \right) }\left(C^*_{\overline \rho}\left(\overline M, \overline\F, \overline\sigma \right) \right)\subset  \pi_{\overline A'}\left(  \overline A'\right).
 $$
 From (i) of the Lemma \ref{blowing_necessary_covering_inf_lem} it follows that $\pi_{C^*_{\overline\rho}\left(\overline M, \overline\F, \overline\sigma \right) }\left(C^*_{\overline \rho}\left(\overline M, \overline\F, \overline\sigma \right)\right) $ is a hereditary $C^*$-subalgebra of $\overline A'$.
  \end{proof}
  \begin{lemma}\label{foli_fibration_lem}
  	Under  hypotheses \ref{groupoid_foli_empt}  suppose that  the foliated space $\left(\overline{M},~ \overline{\mathcal F} \right)$ {comes from the fibration} $\pi_{\overline B}: \overline M \to \overline B$ (cf. Definition \ref{foli_fibration_comes_defn}). If  the $C^*$-algebra  $C^*_\rho\left( M, \F\right)$ is given by \eqref{foli_triv_eqn} and  $C^*_{\overline\rho}\left( \overline M, \overline\F\right)$ is the $\overline q$-lift of $C^*_\rho\left( M, \F\right)$ (cf. Definition \ref{blowing_lift_defn}) then one has
  	\begin{enumerate}
  		\item[(i)] there is an algebraical finite covering category (cf. Definition \ref{algebraical_finite_covering_category_defn}) given by
  			\be\label{groupoid_foli_algebraical_finite_covering_category_trivial_eqn}
  		\begin{split}
  			\mathfrak{S}_{C^*_{\rho}\left( M, \F \right)_0\left(q \right) }\bydef 
  			\left\{\left\{A_\la \right\}_{\la\in \La}, \left\{\pi^\mu_\nu   :A_\mu\hookto A_\nu\right\}\right\}_{\substack{\mu, \nu \in \La\\\mu\le \nu}},\\
  			\quad A_\la \bydef C^*_{\rho_\la}\left( M_\la ,  \F_\la,  \right),\quad \pi^\mu_\nu \bydef C^*_{\rho}\left( M ,  \F, \right)_0\left( q^\nu_\mu\right): A_\mu \hookto A_\mu,
  		\end{split}
  		\ee	
 		\item[(ii)] the triple $\left(C^*_{\rho}\left( M ,  \F\right), C^*_{\overline\rho}\left(\overline M , \overline \F \right) , 	\widehat{G}  \right)$ is the   \textit{disconnected infinite noncommutative covering} of  (cf. Definition \ref{disconnected_infinite_noncommutative_covering_defn})  $\mathfrak{S}_{C^*_{\rho}\left( M, \F \right)_0\left(q \right) }$,
 		\item[(iii)] there are natural $*$-isomorphisms 
 		\bean
 		C^*_r\left( M, \F\right)\cong C^*_\rho\left( M, \F\right)\cong C^*\left( M, \F\right),\\
 		C^*_{r}\left(\overline M , \overline \F \right)\cong C^*_{\overline\rho}\left(\overline M , \overline \F \right)\cong C^*\left(\overline M , \overline \F \right)
 		\eean 
 		where $	C^*_r\left( M, \F\right), C^*_{r}\left(\overline M , \overline \F \right)$ and $C^*_r\left( M, \F\right), C^*_{r}\left(\overline M , \overline \F \right)$ are reduced and full $C^*$-algebras of foliations (cf. Definitions \ref{foli_red_defn} and \ref{foli_full_defn}).
  	\end{enumerate}
\end{lemma} 
\begin{proof}
(i) The algebraical finite covering category $\mathfrak{S}_{C^*_{\rho}\left( M, \F \right)_0\left(q \right) }$ equals to $\mathfrak{S}_{C^*_{\rho}\left( M, \F, \sigma_{\text{triv}} \right)_0\left(q \right) }$ one (cf. \eqref{foli_triv_eqn} and \eqref{groupoid_foli_algebraical_finite_covering_category_eqn}).\\
(ii) From the Theorem \ref{foli_bundle_thm} and the Remark \ref{foli_bundle_rem} it follows that
$$
	C^*_{r}\left(\overline M , \overline \F \right)\cong C^*\left(\overline M , \overline \F \right)\cong C_0\left(\overline B\right)\otimes \K\left(L^2\left(\overline F\right)\right),
$$
so if $\left\| \cdot \right\|_r$ and  $\left\| \cdot \right\|$ are $C^*$-norms of $C^*_{r}\left(\overline M , \overline \F \right)\cong C^*\left(\overline M , \overline \F \right)$ then
\bean
\forall \overline a \in 	C_c\left(\overline M , \overline \F \right)\quad \left\| \overline a \right\|_r = \left\| \overline a \right\|.
\eean
On the other hand from
\bean
\forall \overline a \in 	C_c\left(\overline M , \overline \F \right)\quad \left\| \overline a \right\|_r \le \left\| \overline a \right\|_{\overline \rho} \le  \left\| \overline a \right\|
\eean 
where  $\left\| \overline a \right\|_{\overline \rho}$ is the $C^*$-norm of $C^*_{\overline\rho}\left(\overline M , \overline \F \right)$ it turns out that
\bean
\forall \overline a \in 	C_c\left(\overline M , \overline \F \right)\quad \left\| \overline a \right\|_r = \left\| \overline a \right\|_{\overline \rho}  = \left\| \overline a \right\|,\\	
	C^*_{r}\left(\overline M , \overline \F \right)\cong C^*_{\overline\rho}\left(\overline M , \overline \F \right)\cong C^*\left(\overline M , \overline \F \right).
\eean 
From the  Theorem \ref{blowing_sufficient_covering_inf_thm} it follows that $
\left(C^*_{\rho}\left( M ,   \right), C^*_{\overline\rho}\left(\overline M , \overline \F \right) , 	\widehat{G}  \right)  
$
is a {pre}-{covering of algebraical finite covering category} $\mathfrak{S}_{C^*_{\rho}\left( M, \F \right)_0\left(q \right) }$. 
If a triple $\left(C^*_{\rho}\left( M ,  \F\right), \overline A' , 	\widehat{G}  \right)$ is a {pre}-{covering of the algebraical finite covering category}  $\mathfrak{S}_{C^*_{\rho}\left( M, \F \right)_0\left(q \right) }$ then from the Proposition  \ref{infinite_spectrum_limit_prop} it follows that there is a natural bijective set theoretic  map
$$
\varphi: \overline \sX' \xrightarrow{\cong} \overline B
$$ 
spectrum of $\overline A'$ to the spectrum of $C^*_{\overline\rho}\left(\overline M , \overline \F \right)$. From the Lemma \ref{grupoid_foli_intinite_covering_lem} it follows that there is a natural inclusion 
$
 C^*_{\overline\rho}\left(\overline M , \overline \F \right)= \overline A'
$
such that $C^*_{\overline\rho}\left(\overline M , \overline \F \right)$ is a hereditary subalgebra of $\overline A'$. From the Theorem \ref{hered_spectrum_prop} it turns out that the spectrum of $C^*_{\overline\rho}\left(\overline M , \overline \F \right)$ is an open subset of the spectrum of $C^*_{\overline\rho}\left(\overline M , \overline \F \right)$, it turns out that the map $\varphi: \overline \sX' \xrightarrow{\cong} \overline B$ is a homeomorphism, i.e. the spectrum of  $\overline A'$ is the  Hausdorff space $\overline B$. From the Lemma \ref{oa_haus_alg_lem} it follows that $\overline A'$ is a full algebra of operator fields (cf. Definition \ref{full_algebra_operator_fields_defn}) on $\overline B$. On the other hand $C_0\left(\overline B\right)\otimes \K\left(L^2\left(\overline F\right)\right)$ is also a full algebra of operator fields (cf. Definition \ref{full_algebra_operator_fields_defn}) on $\overline B$. From the Lemma \ref{top_full_oaf_lem} it follows that $C_0\left(\overline B\right)\otimes \K\left(L^2\left(\overline F\right)\right)$ is  a maximal full algebra of operator fields. Taking into account  the inclusion $C_0\left(\overline B\right)\otimes \K\left(L^2\left(\overline F\right)\right)\subset \overline A'$  one has an isomorphism
$$
C_0\left(\overline B\right)\otimes \K\left(L^2\left(\overline F\right)\right)\cong C^*_{\overline\rho}\left(\overline M , \overline \F \right)\xrightarrow{\cong } \overline A'
$$
i.e. any  {pre}-{covering of algebraical finite covering category} $\mathfrak{S}_{C^*_{\rho}\left( M, \F \right)_0\left(q \right) }$ is equivalent to $\left(C^*_{\rho}\left( M ,  \F\right), C^*_{\overline\rho}\left(\overline M , \overline \F \right) , 	\widehat{G}  \right)$.\\
(iii) It is already proven that there are isomorphisms $$ 	C^*_{r}\left(\overline M , \overline \F \right)\cong C^*_{\overline\rho}\left(\overline M , \overline \F \right)\cong C^*\left(\overline M , \overline \F \right)$$. The isomorphisms $	C^*_r\left( M, \F\right)\cong C^*_\rho\left( M, \F\right)\cong C^*\left( M, \F\right)$ follow from \eqref{groupoid_rho_red_eqn}.

\end{proof} 
\begin{theorem} \label{foli_sufficient_covering_inf_thm}
    	Under  hypotheses \ref{groupoid_foli_empt}  suppose that  the foliated space $\left(\overline{M},~ \overline{\mathcal F} \right)$ {comes from the fibration} $\pi_{\overline B}: \overline M \to \overline B$ (cf. Definition \ref{foli_fibration_comes_defn}). If  the $C^*$-algebras    $C^*_{\widetilde \rho}\left( \widetilde M, \widetilde\F\right)$ and $C^*_{\overline\rho}\left( \overline M, \overline\F\right)$ are the $q$- and  $\overline q$-lifts of $C^*_\rho\left( M, \F\right)$ (cf. Definition \ref{blowing_lift_defn}) then one has
\begin{enumerate}
	\item[(i)] the given by the Lemma  \ref{foli_fibration_lem}	{disconnected infinite noncommutative covering} $$\left(C^*_{\rho}\left( M ,  \F\right), C^*_{\overline\rho}\left(\overline M , \overline \F \right) , 	\widehat{G}  \right)$$ of    $\mathfrak{S}_{C^*_{\rho}\left( M, \F \right)_0\left(q \right) }$ is good (cf. Definition \ref{good_defn}).
	\item[(ii)] the triple   $\left(C^*_{\rho}\left( M ,  \F\right), C^*_{\widetilde\rho}\left(\widetilde M , \widetilde \F \right) , 	G\left(\left.\widetilde M~~\right|~M \right)  \right)$ the  {infinite noncommutative covering} of $\mathfrak{S}_{C^*_{\rho}\left( M, \F \right)_0\left(q \right) }$.
\end{enumerate}

\end{theorem}
\begin{proof}
	Follows from the Theorem \ref{blowing_sufficient_covering_inf_thm}.
\end{proof} 
\begin{lemma}\label{foli_simply_lem}
	Under  hypotheses \ref{groupoid_foli_empt}  suppose that  the foliated space $\left(\widetilde{M},~ \widetilde{\mathcal F} \right)$ {comes from the fibration} $\pi_{\widetilde B}: \widetilde M \to \widetilde B$ (cf. Definition \ref{foli_fibration_comes_defn}), and suppose that the manifold $M$ is compact. If $\widetilde B$ is simply connected (cf. Definition \ref{top_simply_conn_defn}) then for any noncommutative finite-fold covering with unitization $\left(C_\rho\left(M, \sF \right), \widetilde A, G, \pi \right)$ (cf. Definition \ref{fin_unitization_defn}) there is a regular finite-fold covering $ q': \left(\widetilde M',  \widetilde\F'\right)\to \left(M, \F \right)$   (cf. Definition \ref{foli_reg_cov_defn}) such that the quadruple $\left(C_\rho\left(M, \sF \right), \widetilde A, G, \pi \right)$ is equivalent to
$$
\left(C_r\left(M, \sF \right), C^*_r\left(\widetilde M',  \widetilde\sF'\right), G\left(\left.\widetilde M'\right|M \right), C_r\left(M, \sF \right)_0\left( q' \right)   \right)
$$
where $C_r\left(M, \sF \right)_0\left(q' \right)$ is the $q'$-lift (cf. Definition \ref{blowing_lift_hom_defn}).
\end{lemma}
\begin{proof}
	From (iii) of the Lemma \ref{foli_fibration_lem} it turns out that $C^*_\rho\left(M, \sF \right)\cong C_r\left(M, \sF \right)$.
From the Theorem  \ref{foli_bundle_thm} it follows that
$$
C_r\left(\widetilde M,  \widetilde\sF\right)\cong C_0\left(\widetilde B\right)\otimes \K\left(L^2\left( \widetilde\F\right)\right)
$$
and taking into account the Theorem \ref{ctr_fin_thm} we conclude that the $C^*$-algebra $C_r\left(\widetilde M,  \widetilde\sF\right)$ is simply connected (cf, Definition \ref{simply_connected_defn}). This lemma follows from the Theorem \ref{blowing_sufficient_fin_thm},
\end{proof}

\begin{lemma}\label{foli_universal_fg_lem} 
	Let $C\left(\sY \right) \hookto M\left(A \right)$ be a Hausdorff blowing-up (cf. Definition \ref{blowing_defn}) with compact, locally connected space $\sY$, and let 
	$q: \widetilde\sY \to \sY$  be an  $A$-regular covering (cf. Definition \ref{blowing_a_regular_defn}), with connected $\widetilde \sY$. Suppose that the $q$-lift $A_0\left( \widetilde\sY\right)$ of $A$ (cf. Definition \ref{blowing_lift_defn}) is simply connected (cf. Definition \ref{simply_connected_defn}). If $\left( A, A_0\left(\widetilde{\sY},  \right), G\left(\left.\widetilde{\sY}~\right|{\sY} \right) \right)$ is the infinite noncommutative covering \ref{infinite_noncommutative_covering_defn} of $\mathfrak{S}_{A_0\left(q \right) }$ (cf. \eqref{blowing_category_fin_eqn}) and the property $P_{\mathrm{untz}}$ of finite-fold noncommutative coverings is given by \eqref{unitization_p_eqn} then one has:
	\begin{enumerate}
		\item[(i)] the $P_{\mathrm{untz}}$-{universal  covering} 
		of $A$ (cf. Definition \ref{fundamental_group_nc_p_defn}) equals to
		$$
		\left( A, A_0\left(\widetilde{\sY},  \right), G\left(\left.\widetilde{\sY}~\right|{\sY} \right) \right), 
		$$
		\item[(ii)] if $\pi^{P_{\mathrm{untz}}}_1\left(A \right)$ is the $P_{\mathrm{untz}}$-{fundamental group} of $A$ (cf. Definition \ref{fundamental_group_nc_p_defn}) then there is a natural isomorphism $\pi^{P_{\mathrm{untz}}}_1\left(A \right)\cong  G\left(\left.\widetilde{\sY}~\right|{\sY} \right)$. 
	\end{enumerate}
	
\end{lemma}
\begin{proof}
	This lemma is a specialization of the \ref{blowing_universal_fg_lem} one.
\end{proof}
\section{Coverings of linear foliations of torus}

\paragraph{}
Here we consider a generalization of the Example \ref{fol_tor_exm}. Let $N \in \N$ and $n \bydef 2N$, and let $\th \in \R \setminus \Q$ be an irrational number.
Consider a vector field $\widetilde{X}_\th$ on $\R^{2N}$ given by
\[
\tilde{X}+\th=\frac{\partial}{\partial
	x_1}+....+ \frac{\partial}{\partial
	x_N}+\th \frac{\partial}{\partial
	x_{N+1}}+...+ \th \frac{\partial}{\partial
	x_{2N}}.
\]
which yields a foliation $\widetilde\F_\th$ on $\R^{2N}$.
Since $\widetilde{X}_\th$ is
invariant under all translations, it determines a vector field $X$
on the $2N$-dimensional torus ${\T}^{2N}={\R}^{2N}/{\Z}^{2N}$. The vector
field $X_\th$ determines a foliation $\mathcal{F}_\th$ on ${\T}^{2N}$.  
\begin{definition}\label{groupoid_linear_torus_fol_defn}
	We say that the foliated spaces $\left(\T^{2N}, \mathcal{F}_\th \right)$ and $\left(\R^{2N}, \widetilde\F_\th \right)$. are a \textit{linear foliation on torus} and the \textit{universal covering} of $\left(\T^{2N}, \mathcal{F}_\th \right)$. These foliations yield groupoids $\G\left( \T^{2N}, \mathcal{F}_\th\right)$ and $ \G\left(\R^{2N}, \widetilde\F_\th \right)$ respectively.  
\end{definition}

\subsection{Noncommutative coverings. Fundamental group}

\begin{lemma}
 If $\left(\T^{2N}, \mathcal{F}_\th \right)$ is a {linear foliation on torus} and
 $C_\rho\left(\T^{2N}, \mathcal{F}_\th \right)$ is a given by \eqref{foli_triv_eqn} $C^*$-algebra then one has:
 \begin{enumerate}
 	\item[(i)] there is a natural isomorphism $C^*_\rho\left(\T^{2N}, \mathcal{F}_\th \right)\cong C^*_r\left(\T^{2N}, \mathcal{F}_\th \right)$ where $C^*_r\left(\T^{2N}, \mathcal{F}_\th \right)$ is the reduced $C^*$-algebra of the foliated space $\left(\T^{2N}, \mathcal{F}_\th \right)$ (cf. Definition \ref{foli_red_defn}),
 	\item[(ii)] if $\left( C^*_r\left(\T^{2N}, \mathcal{F}_\th \right), \widetilde A, G, \pi\right)$ is a noncommutative finite-fold covering with unitization then there is a finite-fold covering $q: \widetilde \T^{2N} \to \T^{2N}$ such that the quadruple $\left( C^*_r\left(\T^{2N}, \mathcal{F}_\th \right), \widetilde A, G, \pi\right)$ is equivalent to
 	$$
 \left( C^*_r\left(\T^{2N}, \mathcal{F}_\th \right), C^*_r\left(\widetilde\T^{2N}, \widetilde\sF_{\widetilde\th } \right),G\left(\left.\widetilde\T^{2N}~\right|\T^{2N} \right)	 , C^*_r\left(\T^{2N}, \mathcal{F}_\th \right)_0\left(q \right) \right)	
 	$$
 	where the  foliated space  $\left(\widetilde\T^{2N}, \widetilde\sF_{\widetilde\th } \right)$ is the $q$-{lift} of $\left(\T^{2N}, \mathcal{F}_\th \right)$ (cf. Definition \ref{fol_cov_defn}) and $C^*_r\left(\T^{2N}, \mathcal{F}_\th \right)_0\left(q \right)$ is the $q$-\textit{lift} of $C^*_r\left(\T^{2N}, \mathcal{F}_\th \right)$ (cf. Definition \ref{blowing_lift_fin_eqn}).
 \end{enumerate}
\end{lemma}
\begin{proof}
If $\left(\R^{2N}, \widetilde\F_\th \right)$ is the universal covering (cf. Definition \ref{groupoid_linear_torus_fol_defn}) of $\left(\T^{2N}, \mathcal{F}_\th \right)$ then $\left(\R^{2N}, \widetilde\F_\th \right)$  {comes from the fibration} $\pi_{\R^{2N}}: \R^{2N} \to \R^N$ (cf. Definition \ref{foli_fibration_comes_defn}).\\
(i) The space $\R^N$ is simply connected, so from 	 (iii) of the Lemma \ref{foli_fibration_lem} it turns out that $C^*_\rho\left(\T^{2N}, \mathcal{F}_\th \right)\cong C^*_r\left(\T^{2N}, \mathcal{F}_\th \right)$.\\
(ii) Follows from the Lemma \ref{foli_simply_lem}
\end{proof}

\begin{lemma}\label{foli_t_universal_fg_lem} 
If the property $P_{\mathrm{untz}}$ of finite-fold noncommutative coverings is given by \eqref{unitization_p_eqn} then one has:
	\begin{enumerate}
		\item[(i)] the $P_{\mathrm{untz}}$-{universal  covering} 
		of $ C^*_r\left(\T^{2N}, \mathcal{F}_\th \right)$ (cf. Definition \ref{fundamental_group_nc_p_defn}) equals to
		$$
		\left( C^*_r\left(\T^{2N}, \mathcal{F}_\th \right),  C^*_r\left(\R^{2N}, \mathcal{F}_\th \right), \Z^{2N} \right), 
		$$
		\item[(ii)] if $\pi^{P_{\mathrm{untz}}}_1\left( C^*_r\left(\T^{2N}, \mathcal{F}_\th \right) \right)$ is the $P_{\mathrm{untz}}$-{fundamental group} of $ C^*_r\left(\T^{2N}, \mathcal{F}_\th \right)$ (cf. Definition \ref{fundamental_group_nc_p_defn}) then there is a natural isomorphism $\pi^{P_{\mathrm{untz}}}_1\left( C^*_r\left(\T^{2N}, \mathcal{F}_\th \right) \right)\cong  \Z^{2N}$. 
	\end{enumerate}
	
\end{lemma}
\subsection{Hurewicz homomorphism}
\paragraph{}

If  $\left(\T^{2}, \mathcal{F}_\th \right)$  is  a {linear foliation on 2-torus} and $\left(\R^{2}, \widetilde\F_\th \right)$. is the {universal covering} of $\left(\T^{2N}, \mathcal{F}_\th \right)$ (cf. Definition \ref{groupoid_linear_torus_fol_defn})
 then there is the   covering $p_2 : \R^2\to \T^2$.  If 
 $$
 \widetilde{\mathcal N}\bydef \left\{\left.\left(0, x \right)\in \R^2\right| x \in \R \right\}
 $$
 then $\widetilde{ \mathcal N}$ is a complete transversal  of $\left(\R^{2}, \widetilde\F_\th \right)$ (cf. Definition \ref{foli_complete_defn}). If
 $$
 { \mathcal N}\bydef p_2\left(  \widetilde{ \mathcal N}\right) 
 $$
then $\mathcal N$ is  a complete transversal  of $\left(\T^{2}, \widetilde\F_\th \right)$ (cf. Definition \ref{foli_complete_defn}). From the Lemma \ref{foli_stab_lem} it follows that
 $$
 C_r\left(\T^{2}, \widetilde\F_\th \right)\cong C^*_r\left(\G^{\mathcal N}_{\mathcal N} \right)\otimes \K,
 $$
 i.e. the $C^*$-algebra  $C_r\left(\T^{2}, \widetilde\F_\th \right)$ is stable (cf. Definition \ref{stable_ca_defn}).
It is proven in \cite{connes:ncg94} that 
$$
C^*_r\left(\G^{\mathcal N}_{\mathcal N} \right)\cong C\left( \T^2_\th\right) 
$$
where $C\left( \T^2_\th\right) $ is a noncommutative torus (cf. Definition \ref{nt_defn}).  
One has $\mathcal N \cong S^1$ and $C\left(S^1 \right)\cong C\left( u\right) $ is a $C^*$-algebra generated by an unitary element. This circumstance yields following inclusions
\bean
C\left(u \right)\hookto  C\left( \T^2_\th\right),\\
C\left(u \right)\hookto  M\left(  C_r\left(\T^{2}, \F_\th \right)\right) 
\eean 
so we can suppose that $u \in  C\left( \T^2_\th\right)$ and $u \in M\left(  C_r\left(\T^{2}, \F_\th \right)\right)$. From the Definition \ref{nt_defn} it follows that  $C\left( \T^2_\th\right)$ is an universal $C^*$-algebra generated by two unitary elements $u, v\in U\left( C\left( \T^2_\th\right)\right) $ such that
$$
	u v = e^{2\pi i\th}vu.
$$
It is explained in \cite{connes:ncg94} that $K_1\left(  C\left( \T^2_\th\right)\right) = \Z\left[u\right]\oplus  \Z\left[v\right]$. For any $n \in \N$ there is the natural subgroup 
\bean
A_n \bydef \left\{\left.\left(n x , y\right)\in \Z^2\right| \left(x , y\right)\in \Z^2  \right\}\subset \Z^2
\eean 
and a factor-space $\widetilde \T^2\bydef \R^2/ A_n$ with the natural covering $q_n:\left(\widetilde \T^2,  \widetilde \F_{\th_n}\right)\to \left(\T^{2}, \F_\th \right)$. If   $C\left( \T^2_{\th/n}\right)$ is generated by $u_n, v_n \in U\left( C\left( \T^2_{\th/n}\right)\right)$ then there is an injective $*$-homomorphism
\bean
\varphi_n : C\left( \T^2_{\th}\right)\hookto C\left( \T^2_{\th/n}\right);\\
u \mapsto u^n_n,  \quad v \mapsto v_n.
\eean 
such that there is the natural commutative diagram 
\\
\begin{tikzcd}
C\left( \T^2_{\th}\right)\otimes\K\arrow[rr, "\varphi_n\otimes \Id_\K" ]\arrow[d, "\parallel" ] && C\left( \T^2_{\th/n}\right)\otimes\K\arrow[d, "\parallel"]\\
C^*_r\left(\T^{2}, \F_\th \right)\arrow[rr, "A_0\left(q_n \right)  "]&& C^*_r\left(\widetilde \T^2,  \widetilde \F_{\th_n}\right)
\end{tikzcd}
\\
where $A_0\left(q_n \right)$ is the the finite-$q_n$-lift (cf. Definition \ref{blowing_finite_lift_defn}) and vertical arrows are $*$-isomorphisms. 
If $C^*_r\left(\widetilde \T^2,  \widetilde \F_{\th_n}\right)\to B\left(\H \right)$ is a faithful nondegenerate representation then from the Definition \ref{multiplier_el_defn} it follows that there are inclusions $C\left( \T^2_{\th}\right)\hookto B\left(\H \right)$ and  $C\left( \T^2_{\th/n}\right)\hookto B\left(\H \right)$. From $u \mapsto u^n_n$ it follows that there   is a Borel $n^{\mathrm{th}}$ root $\phi_n$ of identity map on the set $\left\{\left. z \in \C\right| \left|z \right|=1\right\}$ (cf. \eqref{h_n_root}) such that
$$
u_n = \phi_n \left(u \right) .
$$
 If $A \hookto B\left(\H \right)$ is a faithful representation then $u \in B\left(\H \right)$. 
 On the other hand from 
 $$
  C\left( \T^2_{\th/n}\right) = \sum_{j=0}^{n-1} C\left( \T^2_{\th}\right)u_n^j
 $$
it follows that
$$
C^*_r\left(\widetilde \T^2,  \widetilde \F_{\th_n}\right)= \sum_{j=0}^{n-1}C^*_r\left( \T^2,   \F_{\th}\right)u_n^j,
$$
i.e. 
\be\label{foli_u_ca_eqn}
\mathfrak S_u \bydef \left(C^*_r\left( \T^2,   \F_{\th}\right), C^*_r\left(\widetilde \T^2,  \widetilde \F_{\th_n}\right), \Z_n, A_0\left(q_n \right)  \right) 
\ee
 is a  $\left(u, n\right)$-{covering} (cf. Definition \ref{hurewicz_u_n_defn}).
 Similarly to \eqref{free_hur_cat_eqn} one has a category 
 \bean
 \mathfrak{S}_u= \left\{C^*_r\left( \T^2,   \F_{\th}\right)\hookto C^*_r\left( \T^2,   \F_{\th/n_1}\right)\hookto  ...\hookto C^*_r\left( \T^2,   \F_{\th/n_k}\right)\hookto ... \right\}
 \eean 
This category corresponds to a foliated $\left(S^1 \times \R, \widetilde \F_u\right)$ which comes from the fibration (cf. Definition \ref{foli_fibration_comes_defn}, so from the Lemma \ref{foli_fibration_lem} it follows that $\mathfrak S_u$ is an an algebraical finite covering category (cf. Definition \ref{algebraical_finite_covering_category_defn}). From the Theorem \ref{foli_sufficient_covering_inf_thm} it follows that $\mathfrak S_u$ is good, and the triple  $\left( C^*_r\left( \T^2,   \F_\th\right) , C^*_r\left(S^1 \times \R, \widetilde \F_u\right), \Z\right) $ is the  {infinite noncommutative covering} of $\mathfrak{S}_{u}$. Similarly to the equation \eqref{free_hur_hom_eqn} one has a homomorphism 
\bean
\Z \cong G\left( \left. C^*_r\left(S^1 \times \R, \widetilde \F_u\right)\right|  C^*_r\left( \T^2,   \F_\th\right)\right) \to K^1\left(  C^*_r\left( \T^2,   \F_\th\right)\right),\\
k \mapsto k\left[u\right]
\eean 
Also we can obtain a homomorphism
\bean
\Z \cong G\left( \left. C^*_r\left(S^1 \times \R, \widetilde \F_v\right)\right|  C^*_r\left( \T^2,   \F_\th\right)\right) \to K^1\left(  C^*_r\left( \T^2,   \F_\th\right)\right),\\
k \mapsto k\left[v\right].
\eean 
It is known that $K_0\left(C\left(\T^2_\th \right) \right)\cong K_0\left(C\left(\T^2_\th \right) \right)\cong \Z^2$, so from the Corollaries \ref{stab_k_0_cor} and \ref{stab_k_1_cor}  it follows that there is the natural isomorphism
$$
K_0\left( C_r\left(\T^{2}, \widetilde\F_\th \right)\right) \cong K_1\left( C_r\left(\T^{2}, \widetilde\F_\th \right)\right)= \Z^2,
$$
so the torsion subgroup $K_1\left( C_r\left(\T^{2}, \widetilde\F_\th \right)\right)_{\mathrm{tors}}\subset K_1\left(  C_r\left(\T^{2}, \widetilde\F_\th \right)\right)$ is trivial (cf. Section \ref{torsion_sc_sec}). So one has 
\bean
K_1\left( C_r\left(\T^{2}, \widetilde\F_\th \right)\right)_{\mathrm{free}} \bydef K_1\left( C_r\left(\T^{2}, \widetilde\F_\th \right)\right)/K_1\left( C_r\left(\T^{2}, \widetilde\F_\th \right)\right)/ K_1\left( A\right)_{\mathrm{tors}}\cong\\\cong  K_1\left( C_r\left(\T^{2}, \widetilde\F_\th \right)\right)\cong \Z\left[u\right]\oplus \Z\left[v\right]
\eean
where the notation of the Section \ref{h_free_ass_sec} is used. The given by \eqref{h_znt_eqn} noncommutative finite-fold covering is trivial, the given by \eqref{h_free_ass_eqn} equals to
$$
\left(  C_r\left(\T^{2}, \F_\th \right), C_r\left(\R^{2}, \widetilde\F_\th \right), \Z^2\right) 
$$
Similarly to \eqref{h_hom_ext_eqn} one has isomorphism
\bean
\begin{split}
	h^{ A}_{\mathrm{free}} : G\left( \left. \widetilde  C_r\left(\R^{2}, \widetilde\F_\th \right)\right|C_r\left(\T^{2}, \F_\th \right)\right) \to \Hom\left(K_1\left( C_r\left(\T^{2}, \F_\th \right)\right), \Z  \right). 
\end{split}
\eean
From the Proposition \ref{nt_khom_prop} it follows that the Hurewicz homomorphism (cf. Definition \ref{h_defn})
is an isomorphism 
$$
\pi^{P_{\mathrm{untz}}}_1\left( C^*_r\left(\T^{2N}, \mathcal{F}_\th \right) \right)\cong K^1\left(C_r\left(\T^{2}, \F_\th \right) \right). 
$$

\chapter{Coverings of  noncommutative tori}\label{nt_chap} 
\section{Basic constructions}
\begin{empt}
	For all $x = \left(x_1, ..., x_n\right)\in \R^{n}$ one has $\left\|x \right\|^2 = \left| x_1\right|^2 + ... + \left| x_{n}\right|^2$. So for any $m \in \N$ there is a polynomial map $p_m: \R^{n}\to \R$ such that
	$$
	\left\|x \right\|^{2m}= p_m\left(\left| x_1\right|, ... , \left| x_{n}\right| \right). 
	$$
If $f \in 	\SS\left(\mathbb {R} ^{n}\right)$	then from the equation \eqref{mp_sr_eqn} it follows that for all $m \in \N$ there is $C' > 0$ such that 
$$
\left| f\left(x\right)p_m\left(\left| x_1\right|,  ... , \left| x_{n}\right| \right)\right|= \left| f\left(x\right)\right|	\left\|x \right\|^{2m} < C' \quad\forall x\in \R^n.
$$
If $C'' = \sup_{\substack{x\in \R^n}} \left| f\left(x\right)\right|$ then
$
\left| f\left(x\right)\right|\left(1 + \left\|x \right\|^{2m}\right)< C' + C'',
$
i.e. 
$$
\forall m\in \N ~\quad \exists C > 0 \quad \left| f\left(x\right)\right| < \frac{C}{\left(1 + \left\|x \right\|^{2m}\right) }\quad\forall x\in \R^n.
$$
Similarly one can proof that 
\be\label{nt_ineq_eqn}
\forall f \in 	\SS\left(\mathbb {R} ^{n}\right)~ \forall m\in \N ~\quad \exists C^m_f > 0 \quad \left| f\left(x\right)\right| < \frac{C^m_f}{\left(1 + \left\|x \right\|^{m}\right) }\quad\forall x\in \R^n.
\ee
If $R > 0$ and  $y \in \R^{2N}$ is such that   $\left\|y \right\|<R$ then from $\left\|x \right\|> R$ it follows that $\left\|x + y \right\|\ge \left\|x \right\|-\left\|y \right\|$. If $\left\|x \right\|< \max\left( 1,2R\right) $ then $\left| f\left(x + y \right)\right|  < C''$ and
$$
\left| f\left(x + y \right)\right| < \frac{C^m_f}{\left(1 + \left\|x/2 \right\|^{m}\right)}.
$$
If $C^{m,R}_f \bydef C''\left(1 + \left(2R\right)\right)^m +  C^m_f 2^m$ then
\be\label{nt_ineqr_eqn}
\left\|y \right\|<R\quad \Rightarrow\quad \left| f\left(x+y\right)\right| < \frac{C^{m,R}_f}{\left(1 + \left\|x \right\|^{m}\right) }\quad\forall x\in \R^n.
\ee

\end{empt}

\section{Finite-fold coverings}
\subsection{Coverings of $C^*$-algebras}
\paragraph{}

The finite-fold coverings of  noncommutative tori are described in \cite{clarisson:phd,schwieger:nt_cov}. Here we would like to prove that these coverings comply with the general theory of noncommutative coverings described in \ref{cov_fin_bas_sec}. 
\begin{empt}
Let $\left(\T^2, \sF\right)$ be a linear foliation on torus (cf. Example \ref{fol_tor_exm}). One has $\T^2 = S^1 \times S^1$, assume that $S^1 \to S^1 \times S^1$ is the natural map to the first term of the direct product. If $N$ is the image of $S^1$ then from 
\eqref{foli_nt_eqn} it follows that
\bean
C^*_r\left( \G^N_N\right) \cong C\left(\T^2_\th\right).
\eean
On the other hand from the equation \eqref{foli_stab_eqn}  it follows that'
\bean
C\left(\T^2_\th\right)\otimes \K \cong C^*_r\left(\T^2, \sF\right).
\eean
If $\left( C\left(\T^{n}_{{\Th}} \right), \widetilde A, G, \pi\right)$ is a reduced noncommutative finite-fold covering (cf. Definitions \ref{fin_defn}, and \ref{fin_red_defn}) then from the Theorem \ref{stable_fin_cov_thm} and the Exercise \ref{stable_fin_cov_exer} it follows that  $\left( C\left(\T^{n}_{{\Th}} \right)\otimes \K=  C^*_r\left(\T^2, \sF\right), \widetilde A\otimes \K, G, \pi\right)$ is a reduced noncommutative finite-fold covering. 
\end{empt}

\begin{empt}\label{nt_ff_inclusion_empt}
	Let $C\left(\mathbb{T}^n_{\Theta}\right)$ be a noncommutative torus (cf. Definition \ref{nt_defn}).  Denote by $\left\{U_k\right\}_{k \in \Z^n}\subset C\left(\mathbb{T}^n_{\Theta}\right) $ the set of unitary elements which satisfy to \eqref{nt_unitary_product_eqn}. We would like to construct an inclusion $\pi:C\left(\mathbb{T}^n_{\Theta}\right) \hookto C\left(\mathbb{T}^n_{\widetilde{\Theta}}\right)$ where $C\left(\mathbb{T}^n_{\widetilde{\Theta}}\right)$ is another noncommutative torus. Denote by $\left\{\widetilde{U}_k\right\}_{k \in \Z^n}\subset C\left(\mathbb{T}^n_{\widetilde{\Theta}}\right)$ the set of unitary elements which satisfy to \eqref{nt_unitary_product_eqn}, i.e.
	\begin{equation}\label{nt_unitary_productc_eqn}
	\widetilde{U}_k \widetilde{U}_p = e^{-\pi ik ~\cdot~ \widetilde{\Theta} p} \widetilde{U}_{k + p}
	\end{equation}
	We suppose that there is a subgroup $\Ga \in \Z^n$ of maximal rank such that
	\be\label{nt_inc_eqn}
	\widetilde{U}_k \in \pi\left( C\left(\mathbb{T}^n_{\Theta} \right)\right) \quad \Leftrightarrow \quad k \in \Ga. 
	\ee
	Recall that, given a subgroup $\Gamma \subsetneqq \Z^n$ of maximal rank, there is an invertible $n \times n$-matrix $M$ with integer valued entries such that $\Gamma = M \mathbb{Z}^n$. It is proven in \cite{schwieger:nt_cov} that the inclusion $\pi:C\left(\mathbb{T}^n_{\Theta}\right) \hookto C\left(\mathbb{T}^n_{\widetilde{\Theta}}\right)$ exists if and only if the following condition holds
	\be\label{nt_theta_eqn}
	\Th \in M^{\text{T}} \widetilde{\Th} M + \mathbb{M}_n\left( \Z\right).
	\ee
	From $\Ga = M \mathbb{Z}^n$ it turns out that $C\left(\mathbb{T}^n_{\Theta} \right)$ is the  $C^*$-norm completion of the $\C$-linear span of of unitary elements given by
	\be\label{nt_cov_un_eqn}
	U_k = \widetilde{U}_{Mk} \in C\left(\mathbb{T}^n_{\Theta} \right); \quad k \in \Z^n.
	\ee
		For any $j = 1,..., n$ denote by
	\be\label{nt_j_eqn}
	\widetilde{u}_j \stackrel{\text{def}}{=} \widetilde{U}_{k_j },\quad \text{ where } k_j=\left(0,...,\underbrace{ 1}_{j^{\text{th}}-\text{place}},...,0 \right)\in C\left(\mathbb{T}^n_{\widetilde{\Theta}}\right).
	\ee
	There is the action 
\bean
\R^n \times  C\left(\mathbb{T}^n_{\widetilde{\Theta}}\right) \to C\left(\mathbb{T}^n_{\widetilde{\Theta}}\right);\\
s \cdot \widetilde u_j \mapsto e^{2\pi j s} \widetilde u_j.
\eean	
If $\Ga^\vee$ the dual lattice
\be\label{nt_dual_eqn}
\Ga^\vee \bydef \left\{x \in \Q^n ~|~ x \cdot \Ga \in \Z\right\}= M^{-1}\Z^n
\ee 
then the action $\R^n \times  C\left(\mathbb{T}^n_{\widetilde{\Theta}}\right) \to C\left(\mathbb{T}^n_{\widetilde{\Theta}}\right)$ can by restricted onto $\Ga^\vee\subset\R^n$ such that
$$
g a = a \quad \forall a \in  C\left(\mathbb{T}^n_{{\Theta}}\right).
$$
If $g \in \Ga^\vee \cap \Z^n$ then $g \widetilde a = a$ for all $\widetilde a \in  C\left(\mathbb{T}^n_{\widetilde{\Theta}}\right)$. It follows that there is the inclusion
\be\label{nt_group_eqn}
\Ga^\vee/\Z^n \subset \left\{ \left. g \in \Aut\left(C\left(\mathbb{T}^n_{\widetilde{\Theta}}\right)\right)~\right| ga = a;~~\forall a \in C\left(\mathbb{T}^n_{\Theta}\right) \right\}
\ee
From $\Ga = M\Z^n$ and $\Z^n = M\Ga^\vee$ one has the isomorphism $\Ga^\vee/\Z^n \cong \Z^n/\Ga$ and the inclusion.
\be\label{nt_group_ga_eqn}
 \Z^n/\Ga \subset \left\{ \left. g \in \Aut\left(C\left(\mathbb{T}^n_{\widetilde{\Theta}}\right)\right)~\right| ga = a;~~\forall a \in C\left(\mathbb{T}^n_{\Theta}\right) \right\}
\ee

\end{empt}

\begin{definition}\label{nt_ff_inclusion_m_defn}
In the described in \ref{nt_ff_inclusion_empt} situation we say that the inclusion $\pi:C\left(\mathbb{T}^n_{\Theta}\right) \hookto C\left(\mathbb{T}^n_{\widetilde{\Theta}}\right)$ \textit{corresponds} to the matrix $M$.
\end{definition}
\begin{lemma}\label{nt_cov_fin_lem}
	If $C\left(\mathbb{T}^n_{\Theta}\right) \subset C\left(\mathbb{T}^n_{\widetilde{\Theta}}\right)$ is a given by 
	\ref{nt_ff_inclusion_empt} inclusion then the group
	\be\label{nt_fin_g_eqn}
	G = \left\{ \left. g \in \Aut\left(C\left(\mathbb{T}^n_{\widetilde{\Theta}}\right)\right)~\right| ga = a;~~\forall a \in C\left(\mathbb{T}^n_{\Theta}\right) \right\}
	\ee
	is finite.
\end{lemma}
\begin{proof}
	Let $\widetilde u_j\in C\left(\mathbb{T}^n_{\widetilde{\Theta}}\right)$ be the given by \eqref{nt_j_eqn}  for any $j = 1,..., n$ the unitary element.
	From \eqref{nt_inc_eqn} it turns out that for any $j \in 1,..., n$ there is $l_j \in \N$ such that $ \widetilde{u}_j^{l_j}\in C\left(\mathbb{T}^n_{\Theta}\right)$. 
\bean
C\left(\mathbb{T}^n_{\Theta}\right)_j \bydef \left\{\left.a \in C\left(\mathbb{T}^n_{\Theta}\right) \right|a \widetilde{u}_j^{l_j}= \widetilde{u}_j^{l_j}a \right\} \subset C\left(\mathbb{T}^n_{\Theta}\right),\\
	C\left(\mathbb{T}^n_{\widetilde{\Theta}}\right)_j =  \left\{\left.\widetilde{a} \in C\left(\mathbb{T}^n_{\widetilde{\Theta}}\right)~\right|\widetilde{a}a=a\widetilde{a} \quad \forall a \in \pi\left( C\left(\mathbb{T}^n_{\Theta}\right)_j\right) \right\} 
\eean 	
	then for any $g \in G$ one has $gC\left(\mathbb{T}^n_{\widetilde{\Theta}}\right)_j\subset C\left(\mathbb{T}^n_{\widetilde{\Theta}}\right)_j$. Otherwise $C\left(\mathbb{T}^n_{\widetilde{\Theta}}\right)_j$ is a commutative $C^*$-algebra generated by $\widetilde u_j$. If
	 $\widetilde{w}= g\widetilde{u}_j$ then $\widetilde{w}^{l_j} = \widetilde{u}^{l_j}_j$, and since $C\left(\mathbb{T}^n_{\widetilde{\Theta}}\right)_j$ is commutative one has
	\be\label{nt_root_eqn}
	\widetilde{w}^{l_j}\widetilde{u}^{*l_j}_j= \left(\widetilde{w}\widetilde{u}^*_j \right)^{l_j} = 1.
	\ee
	From \eqref{nt_root_eqn} it turns out that there is $0 \le m_j < l_j$ such  $\widetilde{w} = e^{\frac{2\pi i m_j}{l_j }}\widetilde{u}^*_j$. Similarly one can prove that there are $l_j$ ($j=1, ..., n$) such that for any $g \in G$ there are $0 \le m_j < l_j$ following condition holds
	\be\label{nt_exp_eqn}
	g \widetilde{U}_{k} = \exp\left( 2\pi i \left(\frac{k_1}{l_1}+...+ \frac{k_n}{l_n}\right)\right)  \widetilde{U}_{k}; \quad \text{where} \quad k= \left(k_1, ..., k_n \right).
	\ee
	From \eqref{nt_exp_eqn} it turns out that $\left|G\right|\le l_1 \cdot ...\cdot l_n$, i.e. $G$ is finite. 
\end{proof}
\begin{empt}
	For any $x \in \R^{n}$ let us consider a continuous map 
	\bean
	f_x:[0,1] \to \Aut\left(C\left(\mathbb{T}^n_{{\Theta}}\right) \right); \\
	t \mapsto \left(  U_k \mapsto e^{2\pi i t k \cdot x}U_k\right).
	\eean
	and there is the $\pi$-lift 	$\widetilde{f}_x:[0,1] \to \Aut\left(C\left(\mathbb{T}^n_{\widetilde{\Theta}}\right) \right)$ of $f_x$ (cf. Definition \ref{upl_f_defn}) given by
	\bean
	\widetilde{f}_x\left(t \right) \widetilde{U}_{k} = e^{2\pi i t \left(M^{-1} k\right)\cdot x } \widetilde{U}_{k};\\
	\text{where } M \text{ satisfies to the equation } \eqref{nt_cov_un_eqn}.
	\eean
	From the Lemma \ref{nt_cov_fin_lem} and the Corollary \ref{lift_unique_cor} it follows that $\widetilde{f}_x$ is the unique $\pi$-lift of $f_x$. From the Lemma \ref{lift_commutes_lem} it turns out that $\widetilde{f}_x$ commutes with the given by \eqref{nt_fin_g_eqn} group $G$.  There is the continuous action 
	\be\label{nt_ract_eqn}
	\begin{split}
		\R^n \times C\left(\mathbb{T}^n_\Theta\right) \to C\left(\mathbb{T}^n_\Theta\right),\\
		x \bullet U_{k} = f_x\left(1\right) U_k= e^{2\pi i k\cdot x}U_k; \quad \forall x \in \R^n
	\end{split}
	\ee
	which can be uniquely lifted to the continuous action
	\be\label{nt_xact_eqn}
	\begin{split}
		\R^n \times C\left(\mathbb{T}^n_{\widetilde{\Theta}}\right) \to C\left(\mathbb{T}^n_{\widetilde{\Theta}}\right),\\
		x \bullet	\widetilde{U}_{k}= e^{2\pi i \left(M^{-1} k\right)\cdot x}\widetilde{U}_{k}; \quad \forall x \in \R^n
	\end{split}
	\ee
	Clearly one has
	\be\label{nt_ker_eqn}
	x \bullet \widetilde{a} = \widetilde{a};~~ \forall \widetilde{a} \in C\left(\mathbb{T}^n_{\widetilde{\Theta}}\right)\quad \Leftrightarrow\quad  x \in \Ga.
	\ee
\end{empt}

\begin{lemma}\label{nt_fin_cov_gr_lem}
	If $G$ is given by \eqref{nt_fin_g_eqn} then there is the natural group isomorphism $G \cong \Z^n/\Ga$.
\end{lemma}
\begin{proof}
	If $x \in \Z^n$ then clearly $x \bullet a = a$ for any $a \in C\left(\mathbb{T}^n_\Theta\right)$ so there is the natural group homomorphism $\phi: \Z^n \to G$. From \eqref{nt_ker_eqn} it turns out that $\ker \phi = \Ga$, i.e. there is the natural inclusion $\Z^n/\Ga \hookto G$. From the Lemma \ref{lift_commutes_lem} it turns out that the action of $G$ commutes with the action of $\Z^n/\Ga$ on $C\left(\mathbb{T}^n_{\widetilde{\Theta}}\right)$ it follows that $G$ is the direct product $G = G' \times \Z^n/\Ga$. Action of $G$ on $C\left(\mathbb{T}^n_{\widetilde{\Theta}}\right) $ induces the unitary action of $G$ on the Hilbert space $L^2\left(C\left(\mathbb{T}^n_{\widetilde{\Theta}}\right), \widetilde{\tau} \right)$. It  follows that there is the following Hilbert direct sum
	\be\label{nt_hilb_big_eqn}
	L^2\left(C\left(\mathbb{T}^n_{\widetilde{\Theta}}\right), \tau\right)  = \oplus_{\la\in \La_{G'}}\H_\la \bigoplus \oplus_{\mu\in \La_{\Z^n/\Ga}}\H_\mu
	\ee
	where $\La_{G'}$ (resp. $\La_{\Z^n/\Ga}$) is the set of irreducible representations of $G'$ (resp. $\Z^n/\Ga$), $\H_\la$ (resp. $\H_\mu$) is the Hilbert subspace which corresponds to the representation $\la$ (resp. $\mu$).
	The group	$\Z^n/\Ga$ is and commutative group, hence any irreducible representation of $\Z^n/\Ga$ has dimension 1. 
	For any  $\overline{x} \in \Ga^\vee/\Z^n$ the corresponding representation $\psi_{\overline{x}}: {\Z^n/\Ga} \times \C \to \C$  is given by
	\be\label{nt_char_r_eqn}
	\begin{split}
		\psi\left(\overline k, z \right), \quad \chi_{\overline{x}}\left( \overline{k}\right)= e^{2\pi i k x};\\
		\text{ where } z \in \C, \quad k \in \Z^n \text{ and } ~x \in \Ga^\vee \text{ are representatives of  } \\ \overline{k}\in \Z^n/\Ga \text{ and } \overline{x} \in \Ga^\vee/\Z^n.
	\end{split}
	\ee
	So there is a bijective map   $\varphi_{\Z^n/\Ga}:\Z^n/ \Ga^\vee \xrightarrow{\approx} \La_{\Z^n/\Ga}$. From \eqref{nt_xact_eqn} it turns out that
	$$
	\overline{k}' \bullet	\widetilde{U}_{k} = \chi_{\overline{M^{-1}k}}\left( \overline{k}'\right)\widetilde{U}_{k}; \text{ where } k' \text{ is representative of } \overline{k}' \in \Z^n/\Ga,
	$$
	hence one has
	\bean
	\C\widetilde{U}_k \in \H_{\varphi_{\Z^n/\Ga}\left( {M^{-1}k}\right) }\subset \bigoplus_{\mu\in \La_{\Z^n/\Ga}}\H_\mu
	\eean
	The $\C$-linear span of $\left\{\widetilde{U}_k\right\}_{k \in \Z^n}$ is dense in 	$L^2\left(C\left(\mathbb{T}^n_{\widetilde{\Theta}}\right)\right) $, so there is the following Hilbert direct sum
	\be\label{nt_hilb_little_eqn}
	L^2\left(C\left(\mathbb{T}^n_{\widetilde{\Theta}}\right), \tau\right)  =  \oplus_{\mu\in \La_{\Z^n/\Ga}}\H_\mu
	\ee
	Comparison of \eqref{nt_hilb_big_eqn} and \eqref{nt_hilb_little_eqn} gives $\oplus_{\la\in \La_{G'}}\H_\la= \{0\}$ it turns out that $G'$ is trivial, hence  $G = \Z^n/\Ga$.
\end{proof}

\begin{thm}\label{nt_fin_thm}
	In the above situation the triple $\left( C\left(\T^{n}_{{\Th}} \right), C\left(\T^{n}_{\widetilde{\Th}} \right), \Z^n / \Ga, \pi\right)$ is an unital noncommutative finite-fold  covering. 
\end{thm}
\begin{proof}
	Both $C\left(\T^{n}_{{\Th}} \right)$, $C\left(\T^{n}_{\widetilde{\Th}} \right)$ are unital $C^*$-algebras, so from the Definition \ref{fin_unital_defn} one needs check the following conditions:
	\begin{enumerate}
		\item[(i)] The quadruple $\left( C\left(\T^{n}_{{\Th}} \right), C\left(\T^{n}_{\widetilde{\Th}} \right), \Z^n / \Ga, \pi\right)$ is a {noncommutative finite-fold  pre-covering}. (cf. Definition \ref{fin_pre_defn}).
		\item[(ii)] $C\left(\T^{n}_{\widetilde{\Th}} \right)$ is a finitely generated $C\left(\T^{n}_{{\Th}} \right)$-module.
	\end{enumerate}
	(i) One needs check (a), (b)  of the Definition \ref{fin_pre_defn}.
	\begin{enumerate}
		\item[(a)] Follows from the Lemma \ref{nt_fin_cov_gr_lem}.
		\item[(b)] Suppose
		$$
		\widetilde a = \sum_{k \in \Z^n} c_k \widetilde U_k \in C\left(\T^{n}_{\widetilde{\Th}} \right)^{\Z^n / \Ga}; \quad c_k \in \C.
		$$
		If there is $k \in \Z^n \setminus \Ga$ such that $c_k \neq 0$ and  $\overline k \in \Z^n / \Ga$ is represented by $k$  then $\overline k \widetilde a \neq \widetilde a$, so one has  $\widetilde a \notin C\left(\T^{n}_{\widetilde{\Th}} \right)^{\Z^n / \Ga}$. From this contradiction we conclude that $c_k \neq 0 \Leftrightarrow k \in \Ga$. It follows that $\widetilde a\in C\left(\T^{n}_{{\Th}} \right)$, hence one has $C\left(\T^{n}_{\widetilde{\Th}} \right)^{\Z^n / \Ga}=C\left(\T^{n}_{{\Th}} \right)$.
	\end{enumerate}
	(ii) 
	Let $\widetilde a   \in C\left(\T^{n}_{\widetilde{\Th}} \right)$ be any element given by
	$$
	\widetilde a = \sum_{\widetilde x \in \Z^n} \widetilde{c}_{\widetilde x} \widetilde{U}_{\widetilde x}.
	$$
	Let $M \in \mathbb{M}_n\left(\Z\right)$ be such that $\Ga = M\Z^n$. The group $\Z^n / \Ga \cong \Z^n / M\Z^n$, let $N = \left| \Z^n / \Ga\right|$ and  $\left\{ \widetilde x_1,\dots, \widetilde x_N\right\}  \subset  \Z^n$ is a set of representatives of  $\Z^n / \Ga$. If 
	$$
	b_j = \frac{1}{\left|\Z^n / \Ga\right|}\sum_{g \in \Z^n / \Ga}g \left( \widetilde{a}\widetilde{U}^*_{\widetilde x_j}\right) \in C\left(\T^n_\Th \right)
	$$ then from the definition of $C\left(\T^n_\Th \right)$ the series
	$$
	b_j = \sum_{y \in \Z^n} d^j_y  \widetilde{U}_{My}
	$$	
	is norm convergent and $d^j_y = \widetilde{c}_{My - \widetilde x_j}$. It turns out that
	\be\label{nt_bj_eqn}
	b_j \widetilde{U}_{\widetilde x_j}= \sum_{y \in \Z^n} \widetilde{c}_{My +\widetilde x_j} \widetilde{U}_{My+ \widetilde x_j}
	\ee
	Since $\left\{ \widetilde x_1,\dots, \widetilde x_N\right\}  \subset  \Z^n$ is a set of representatives of  $\Z^n / \Ga$ from \eqref{nt_bj_eqn} it follows that
	\be\label{nt_bjn_eqn}
	\sum_{j = 1}^n b_j \widetilde{U}_{\widetilde x_j} = \sum_{\widetilde x \in \Z^n} \widetilde{c}_{\widetilde x} \widetilde{U}_{\widetilde x} = \widetilde a
	\ee
	The equation \eqref{nt_bjn_eqn} means that $C\left(\T^{n}_{\widetilde{\Th}} \right)$ is a left   $C\left(\T^{n}_{{\Th}} \right)$-module generated by the finite set $\left\{\widetilde{U}_{\widetilde x_j}\right\}$.
\end{proof}

\begin{remark}\label{nt_diag_rem}
		Let $\Th = J \th$ where $\th \in \R \setminus \Q$ and
	$$
	J = \begin{pmatrix} 0 & 1_N \\ -1_N & 0 \end{pmatrix}.
	$$
	Denote by $C\left( \T^{2N}_\th\right) \stackrel{\text{def}}{=} C\left( \T^{2N}_\Th\right)$. Let $n \in \N$ and $n > 1$ and $\widetilde \th = \th / n^2$. If $u_1, ..., u_{2N} \in C\left( \T^{2N}_\th\right)$ are generators of $C\left( \T^{2N}_\th\right)$ (cf. Definition \ref{nt_uni_defn}) and $\widetilde u_1, ...,\widetilde u_{2N} \in C\left( \T^{2N}_{\widetilde\th}\right)$ are generators of $C\left( \T^{2N}_{\widetilde\th}\right)$ then there is an injective $*$-homomorphism
	\bean
\pi:	C\left( \T^{2N}_\th\right) \hookto C\left( \T^{2N}_{\widetilde\th}\right),\\
	u_j \mapsto \widetilde u_j^n \quad j = 1,...,2N.
	\eean
	From the Theorem \ref{nt_fin_thm} it turns out that the quadruple
	$$
	\left(	C\left( \T^{2N}_\th\right) , C\left( \T^{2N}_{\widetilde\th}\right), \Z^{2N}_n, \pi\right)
	$$
	is  an unital noncommutative finite-fold  covering.
\end{remark}

\begin{corollary}
	Consider the situation of the Theorem \ref{nt_fin_thm}.  If both operator algebras $\Coo 
	\left(\T^{n}_{{\Th}} \right)$ and  $\Coo\left(\T^{n}_{\widetilde{\Th}} \right)$ are given by \eqref{nt_coo_eqn} then the quadruple $$\left( \Coo\left(\T^{n}_{{\Th}} \right), \Coo\left(\T^{n}_{\widetilde{\Th}} \right), \Z^n / \Ga, \pi|_{\Coo\left(\T^{n}_{\widetilde{\Th}} \right)}\right)$$ is a noncommutative finite-fold covering of bounded operator *-algebras (cf. Definition \ref{fin_oa_defn}).
\end{corollary}
\begin{proof}
	From the Theorem \ref{nt_fin_thm} it follows that a noncommutative finite-fold covering of bounded operator *-algebras looks like
	$$
	\left( \Coo\left(\T^{n}_{{\Th}} \right), \widetilde A, \Z^n / \Ga, \pi|_{\Coo\left(\T^{n}_{\widetilde{\Th}} \right)}\right)
	$$
	where $\widetilde A\subset C\left(\T^{n}_{\widetilde{\Th}} \right)$ is  a dense *-subalgebra which satisfies to the Definition \ref{fin_oa_defn}. From $\Coo\left(\T^{n}_{\widetilde{\Th}} \right)\cap C\left(\T^{n}_{{\Th}} \right)=\Coo\left(\T^{n}_{{\Th}} \right)$ and $\Z^n / \Ga  \Coo\left(\T^{n}_{\widetilde{\Th}} \right)= \Coo\left(\T^{n}_{\widetilde{\Th}} \right)$ it follows that $\Coo\left(\T^{n}_{\widetilde{\Th}} \right)\subseteqq  \widetilde A$. If $\Coo\left(\T^{n}_{\widetilde{\Th}} \right)\subsetneqq  \widetilde A$ then there is $\widetilde a\in \widetilde A\setminus \Coo\left(\T^{n}_{\widetilde{\Th}} \right)$. If $\left\{ \widetilde x_1,\dots, \widetilde x_N\right\}  \subset  \Z^n$ is a set of representatives of  $\Z^n / \Ga$ from \eqref{nt_bj_eqn} it follows that
	\bean
	\widetilde a= \sum_{j = 1}^n b_j \widetilde{U}_{\widetilde x_j} \quad b_1, ..., b_N\in C\left(\T^{n}_{{\Th}} \right).
	\eean
	Since $\widetilde a\notin \Coo\left(\T^{n}_{\widetilde{\Th}} \right)$ there is $j \in \left\{1,...,N\right\}$ such that $b_j \notin \Coo\left(\T^{n}_{{\Th}} \right)$. On the other hand from
	$$
	b_j = \frac{1}{\left|\Z^n / \Ga\right|}\sum_{g \in \Z^n / \Ga}g \left(\widetilde a \widetilde{U}^*_{\widetilde x_j}\right)\in \widetilde A\cap  \left(\T^{n}_{{\Th}} \right)=\Coo\left(\T^{n}_{{\Th}} \right) 
	$$
one has a contradiction, so the inclusion  $\Coo\left(\T^{n}_{\widetilde{\Th}} \right)\subsetneqq  \widetilde A$ is impossible, i.e. one has $\Coo\left(\T^{n}_{\widetilde{\Th}} \right)=  \widetilde A$.
\end{proof}

\begin{lemma}\label{nt_diag_lem}
	Let $\pi: C\left(\T^{n}_{{\Th}} \right)\hookto C\left(\T^{n}_{\widetilde{\Th}} \right)$ be an injective $*$-homomorphism which corresponds to a  noncommutative finite-fold  covering $\left( C\left(\T^{n}_{{\Th}} \right), C\left(\T^{n}_{\widetilde{\Th}} \right), \Z^n / \Ga\right)$. There is a diagonal covering $\left( C\left(\T^{n}_{{\Th}} \right), C\left(\T^{n}_{{\Th/m^2}} \right), \Z^n_m\right)$ and the injective $*$-homomorphism $\widetilde{\widetilde{\pi}}:C\left(\T^{n}_{\widetilde{\Th}} \right) \to C\left(\T^{n}_{{\Th/m^2}} \right)$ such that following condition hold:
	\begin{itemize}
		\item The homomorphism $\widetilde{\pi}:C\left(\T^{n}_{\widetilde{\Th}} \right)\hookto C\left(\T^{n}_{{\Th/m^2}} \right)$ corresponds an unial noncommutative finite-fold  covering having unique map lifting $\left( C\left(\T^{n}_{\widetilde{\Th}} \right), C\left(\T^{n}_{{\Th/k^2}} \right), \Z^n/ m\Z^n\right)$,
		\item Following digram
		\\
		\begin{tikzpicture}
		\matrix (m) [matrix of math nodes,row sep=3em,column sep=4em,minimum width=2em]
		{
			C\left(\T^{n}_{\widetilde{\Th}} \right) & &C\left(\T^{n}_{{\Th/m^2}} \right)\\ 
			& C\left(\T^{n}_{{\Th}}\right)  & \\};
		\path[-stealth]
		(m-1-1) edge node [above] {$\widetilde{\pi}$} (m-1-3)
		(m-2-2) edge node [left]  {$\pi~~~$} (m-1-1)
		(m-2-2) edge node [right] {$~~~\widetilde{\widetilde{\pi}} $} (m-1-3);
		\end{tikzpicture}
		\\
		is commutative. 
		
	\end{itemize}
\end{lemma}
\begin{proof}
	From the construction \ref{nt_ff_inclusion_empt} it turns out that:
	\bean
	\Ga = M\Z^n,\\
	\Th = M^{\mathrm{T}} \widetilde{\Th} M + P\\
	\widetilde{\Th} = \left( M^{\mathrm{T}}\right)^{-1}  {\Th} M^{-1}+ \left( M^{\mathrm{T}}\right)^{-1}  P M^{-1}
	\eean
	where $P \in \mathbb{M}_n\left( \Z\right)$. There is $m \in \N$ such that $N = mM^{-1} \in \mathbb{M}_n\left( \Z\right)$, or equivalently 
	$$
	MN = \begin{pmatrix}
	m& \ldots & 0  \\
	\vdots &\ddots  & \vdots \\
	0 &\ldots & m  \\
	\end{pmatrix}.
	$$
	If
	$$
	\widetilde{\widetilde{ \Th}}= \frac{\widetilde{\Th}-\left( M^{\mathrm{T}}\right)^{-1}  P M^{-1} }{m^2}
	$$
	then
	\bean
	\Th = mM^{\mathrm{T}} \widetilde{\widetilde{ \Th}} Mm,\\
	\widetilde{ \Th} = m 	\widetilde{\widetilde{ \Th}}m-  m\left( M^{\mathrm{T}}\right)^{-1}  P M^{-1}m =\\=m 	\widetilde{\widetilde{ \Th}}m-  N^{\mathrm{T}} P N \in m 	\widetilde{\widetilde{ \Th}}m + \mathbb{M}_n\left( \Z\right).
	\eean
	So there are noncommutative finite-fold coverings $\pi': C\left( \T^n_{{\Th}}\right) \hookto C\left( \T^n_{\widetilde{\widetilde{ \Th}}}\right)$, $\pi'':C\left( \T^n_{\widetilde{\Th}}\right) \hookto C\left( \T^n_{\widetilde{\widetilde{ \Th}}}\right)$, which correspond to triples $\left( C\left( \T^n_{{\Th}}\right), C\left( \T^n_{\widetilde{\widetilde{ \Th}}}\right), \Z^n / mM\Z^n\right)$, $\left(C\left( \T^n_{\widetilde{\Th}}\right), C\left( \T^n_{\widetilde{\widetilde{ \Th}}}\right), \Z^n_m \right)$  
	From $N = mM^{-1}$ it turns out
	\bean
	\Th= m M^{\mathrm{T}} N^{\mathrm{T}}\left( N^{\mathrm{T}}\right)^{-1}  \widetilde{\widetilde{ \Th}}  P N^{-1}NMm= m^2\left(  \left( N^{\mathrm{T}}\right)^{-1}  \widetilde{\widetilde{ \Th}}   N^{-1}\right)  m^2,\\
	\left( N^{\mathrm{T}}\right)^{-1}  \widetilde{\widetilde{ \Th}}   N^{-1} = \frac{\Th}{m^{4}},\\
	\widetilde{\widetilde{ \Th}} = N^{\mathrm{T}}\left( N^{\mathrm{T}}\right)^{-1}  \widetilde{\widetilde{ \Th}}  P N^{-1}N = N^{\mathrm{T}} \frac{\Th}{m^{4}}N
	\eean 
	From the above equations it follows that there is  a noncommutative finite-fold covering ${\widetilde{\pi}}:C\left( \T^n_{\widetilde{\widetilde{\Th}}}\right) \hookto C\left( \T^n_{\Th/m^4}\right)$ which corresponds to the triple $$\left( C\left( \T^n_{\widetilde{\widetilde{\Th}}}\right), C\left( \T^n_{\Th/m^4}\right), \Z^n/N\Z^n\right).$$ The composition $\widetilde{\widetilde{\pi}} = {\widetilde{\pi}}\circ \pi: C\left( \T^n_{\Th}\right) \to C\left( \T^n_{\Th/m^4}\right)$ corresponds to the diagonal covering $\left(C\left( \T^n_{\Th}\right), C\left( \T^n_{\Th/m^4}\right), \Z^n_{m^2} \right)$. 
\end{proof}
\begin{corollary}\label{nt_cofinal_cor}
	If $\mathfrak{S}_{C\left( \T^n_{\Th}\right)}$ is a family of all noncommutative finite-fold coverings given by the Theorem \ref{nt_fin_thm} then following conditions hold:
	\begin{enumerate}
		\item[(i)] Let $\left\{p_k\right\}_{k \in \N}\in \N$ is a sequence such that for every $s \in \N$ there are $j, m \in \N$ which satisfy to the following condition 
		\bean
		sm = \prod_{k = 1}^j p_j.
		\eean
		If $m_j = \prod_{k = 1}^j p_k$ then the ordered set 
		$$
		\mathfrak{S}_{C\left(\T^{n}_{{\Th}} \right)}^\N=\left\{ \left( C\left(\T^{n}_{{\Th}} \right), C\left(\T^{n}_{\Th/m_j^2} \right), \Z^n_{m_j}\right)  \right\}_{j \in \N}
		$$ 
		is a cofinal subset (cf. Definition \ref{cofinal_defn}) of $\mathfrak{S}_{C\left(\T^{n}_{{\Th}} \right)}$. 
		\item[(ii)] The family $\mathfrak{S}_{C\left( \T^n_{\Th}\right)}$ is a directed set with respect to the given by the Definition \ref{g_category_defn} order.
	\end{enumerate}
	
\end{corollary}

\begin{proof}
(i) Let $\left( C\left(\T^{n}_{{\Th}} \right), C\left(\T^{n}_{\widetilde{\Th}} \right), \Z^n / \Ga\right)$ be a noncommutative finite-fold  covering  having unique path  lifting. From the  Lemma \ref{nt_diag_lem} it turns out that there is a noncommutative finite-fold  covering $C\left(\T^{n}_{\widetilde{\Th}} \right) \hookto C\left(\T^{n}_{{\Th}/s^2}\right)$. Otherwise there is $j \in \N$ such $m_j = sm$, so there are following noncommutative finite-fold  coverings
	\be\label{nt_cof_eqn}
	\begin{split}
		C\left(\T^{n}_{{\Th}/s^2}\right) \hookto C\left(\T^{n}_{{\Th}/m_j^2}\right),\\
		C\left(\T^{n}_{\widetilde{\Th}}\right) \hookto C\left(\T^{n}_{{\Th}/m_j^2}\right).
	\end{split}
	\ee
	From \eqref{nt_cof_eqn} it turns out that
	$$
	\left( C\left(\T^{n}_{{\Th}} \right)\to C\left(\T^{n}_{{\Th}/m_j^2} \right)\right) \ge \left( C\left(\T^{n}_{{\Th}} \right)\to C\left(\T^{n}_{\widetilde{\Th}} \right)\right),
	$$
	i.e. $\left\{ \left( C\left(\T^{n}_{{\Th}} \right), C\left(\T^{n}_{\Th/m_j^2} \right), \Z^n_{m_j}\right)  \right\}_{j \in \N}$ is cofinal.\\
	(ii) Follows from (i).
\end{proof}

\begin{example}\label{nt_stable_exm}
	Let $\th \in \R\setminus \Q$, and 
let $C\left(\T^2_\th\right)$ be a noncommutative torus generated by unitary elements $u, v \in U\left(C\left(\T^2_\th\right) \right)$ such that $uv = e^{2\pi i \th} vu$. Similarly let both $C\left(\T^2_{\frac{\th}{mn}}\right)$ and $C\left(\T^2_{\frac{\th + 1}{mn}}\right)$  generated by unitary elements $u', v' \in U\left(C\left(\T^2_{\frac{\th}{mn}}\right) \right)$  and $u'', v'' \in U\left(C\left(\T^2_{\frac{\th +1}{mn}}\right) \right)$ noncommutative tori such that
\bean
u'v' = e^{2\pi i\frac{\th}{mn}}v'u',\\
u''v'' = e^{2\pi i\frac{\th + 1}{mn}}v''u''.
\eean
There are two non isomorphic noncommutative finite-fold coverings 
\bean 
\pi': C\left(\T^2_\th\right)\hookto C\left(\T^2_{\frac{\th}{mn}}\right),\quad 
u \mapsto u'^n, \quad v\mapsto v'^m;\\
\pi'': C\left(\T^2_\th\right)\hookto C\left(\T^2_{\frac{\th +1}{mn}}\right),\quad 
u \mapsto u''^n, \quad v\mapsto v''^m.
\eean 
However there is an $*$-isomorphism $$\pi: C\left(\T^2_{\frac{\th}{mn}}\right)\otimes \mathbb{M}_{mn}\left(\C \right)  \xrightarrow{\approx}C\left(\T^2_{\frac{\th + 1}{mn}}\right)\otimes \mathbb{M}_{mn}\left(\C \right).$$
If $U, V \in \mathbb{M}_{mn}\left(\C \right)$ are given by
\bean
U\bydef \begin{pmatrix}
	e^{2\pi i\frac{1}{mn}} & 0 & \cdots & 0\\
	0 & e^{2\pi i\frac{2}{mn}}&  \cdots & 0\\
\vdots&\vdots & \ddots & \vdots\\
0& 0& \cdots & 1
\end{pmatrix},\quad 
V\bydef \begin{pmatrix}
0 & 0 &  \cdots& 0 & 1\\
	1 & 0& \cdots& 0& 0\\
		0 & 1&  \cdots &0 & 0\\
	\vdots& \vdots & \ddots&\vdots&\vdots \\
	0& 0& \cdots & 1& 0\\
	0& 0&  \cdots & 0& 0
\end{pmatrix}
	\eean
	then one has
	\bean
	UV = e^{\frac{2\pi i}{nm}}VU
	\eean 
	and the $*$-isomorphism $\pi$ is given by
	\bean
	u' \otimes 1_{\mathbb{M}_{mn}\left(\C \right)}\mapsto u'' \otimes U,\\
	v' \otimes 1_{\mathbb{M}_{mn}\left(\C \right)}\mapsto v'' \otimes V.
	\eean
\end{example}


\subsection{Induced representations}\label{nt_i_sec}
\paragraph*{}

Let us consider an unital noncommutative finite-fold  covering $$\left( C\left(\T^{n}_{{\Th}} \right), C\left(\T^{n}_{\widetilde{\Th}} \right), \Z^n / \Ga, \pi\right)$$ given by the Theorem \ref{nt_fin_thm}. The group	$\Z^n/\Ga$ is and commutative group, hence any irreducible representation of $\Z^n/\Ga$ has dimension 1. From the proof of the Lemma \ref{nt_fin_cov_gr_lem} it follows that the moduli space of the irreducible representations $\Z^n / \Ga$ coincides with $\Ga^\vee/ \Z^n $ and any irreducible representation corresponds to a character the given by \eqref{nt_char_r_eqn} character $\chi_{\overline{x}}: \Z^n / \Ga \to \C$. It is known the orthogonality of characters (cf. \cite{hamermesh:group}), i.e.
\be\label{nt_ort_eqn}
\sum_{g \in \Z^n / \Ga} \overline{\chi}_{\overline x_1}\left( g\right) \chi_{\overline x_1}\left( g\right)= \left\{\begin{array}{c l}
	\left|\Z^n / \Ga\right|   & \overline x_1 = \overline x_2 \\
	0 & \overline x_1 \neq\overline x_2 .
\end{array}\right.
\ee

Otherwise one has
\bean
g \widetilde{U}_k = \chi_{\overline {M^{-1}k}}\left(g \right)  \widetilde{U}_k;\quad \forall g \in \Z^n / \Ga
\eean
and taking into account \eqref{nt_ort_eqn} following condition holds
\be\label{nt_hprod_eqn}
\begin{split}
	\left\langle \widetilde{U}_{k_1 }, \widetilde{U}_{k_2} \right\rangle_{C\left(\mathbb{T}^n_{{\Th}}\right)} = \sum_{g \in \Z^n / \Ga } g\left( \widetilde{U}^*_{k_1 }\widetilde{U}_{k_2 }\right)=
	\\
	=\sum_{g \in \Z^n / \Ga }\overline{\chi}_{\overline {M^{-1}k_1}}\left( g\right) \chi_{\overline {M^{-1}k_2}}\left( g\right)\widetilde{U}^*_{k_1 }\widetilde{U}_{k_2 }
	= \\
	=\sum_{g \in \Z^n / \Ga }\overline{\chi}_{\overline {M^{-1}k_1}}\left( g\right) \chi_{\overline {M^{-1}k_2}}\left( g\right)|e^{-\pi ik_1 ~\cdot~ \widetilde\Theta k_2}\widetilde{U}_{k_2-k_1 }
	= \\
	= \left\{\begin{array}{c l}
		\left|\Z^n / \Ga\right|e^{-\pi ik_1 ~\cdot~ \widetilde\Theta k_2} U_{M^{-1}\left(k_2 - k_1 \right) } \in C\left(\T^{n}_{{\Th}} \right)  & k_1 - k_2 \in \Ga \\
		0 & k_1 - k_2 \notin \Ga .
	\end{array}\right.
\end{split}
\ee
From the equation \eqref{nt_l2_eqn} it follows that
\be\label{nt_huprod_eqn}
\begin{split}
	U_{M^{-1}\left(k_2 - k_1 \right)}\xi_l = e^{-\pi iM^{-1}\left(k_2 - k_1 \right) ~\cdot~ \Theta l}\xi_{M^{-1}\left(k_2 - k_1 \right) + l}= 
	\\
	=e^{-\pi i\left(k_2 - k_1 \right) ~\cdot~ \widetilde\Theta Ml}\xi_{M^{-1}\left(k_2 - k_1 \right) + l}
\end{split}
\ee
Let $\rho: C\left(\mathbb{T}^n_{{\Th}}\right) \hookto B\left( L^2\left(C\left(\mathbb{T}^n_{\Theta}\right), \tau\right)\right)$ be the described in \ref{nt_gns_empt} faithful GNS-representation.
If $\widetilde \rho: C\left(\T^{n}_{\widetilde{\Th}} \right) \to B\left(\widetilde \H \right)$ is the induced by the pair 
$$
\left(\rho, \left( C\left(\T^{n}_{{\Th}} \right), C\left(\T^{n}_{\widetilde{\Th}} \right), \Z^n / \Ga, \pi\right)\right)
$$
representation (cf. Definition \ref{induced_repr_fin_defn}) then $\widetilde \H$ is the Hilbert norm completion of the tensor product $C\left(\T^{n}_{\widetilde{\Th}} \right) \otimes_{C\left(\T^{n}_{{\Th}} \right)} L^2\left(C\left(\mathbb{T}^n_{\Theta}\right), \tau\right)$. Denote by $\left(\cdot, \cdot \right)_{\widetilde \H}$ the scalar product on $\widetilde \H$ given by \eqref{induced_prod_equ}.
If $\xi_{l_1}, \xi_{l_2} \in L^2\left(C\left(\mathbb{T}^n_{\Theta}\right), \tau\right)$ then from \eqref{nt_hprod_eqn}, \eqref{nt_huprod_eqn} and \eqref{nt_xi_eqn} it follows that
\be\label{nt_hh_product_eqn}
\begin{split}
	\left(\widetilde{U}_{k_1 } \otimes \xi_{l_1} , \widetilde{U}_{k_2 } \otimes \xi_{l_2} \right)_{\widetilde \H} =
	\\
	=\left\{\begin{array}{c l}
		e^{-\pi ik_1 ~\cdot~ \widetilde\Theta k_2}e^{\pi i\left(k_2 - k_1 \right) ~\cdot~ \widetilde\Theta Ml_2}	\left|\Z^n / \Ga\right|  & k_1 + Ml_1 = k_2 + Ml_2  \\
		0 & k_1 + Ml_1 \neq k_2 + Ml_2 .
	\end{array}\right.
\end{split}
\ee
From the equation \eqref{nt_hh_product_eqn} it turns out that $\widetilde \H$ is a Hilbert direct sum 
$$
\widetilde\H = \bigoplus_{p \in \Z^n} \widetilde\H_p
$$
of one-dimensional Hilbert spaces such that
$$
\widetilde{U}_{k } \otimes \xi_{l}\in \widetilde\H_{k + Ml}
$$
If $k_1 + Ml_1 = k_2 + Ml_2$  and
\bean
\widetilde\xi_1 \bydef \sqrt{\left|\Z^n / \Ga\right|} e^{-\pi ik_1 ~\cdot~ \widetilde\Theta Ml_1}\widetilde{U}_{k_1 } \otimes \xi_{l_1},\\ 
\widetilde\xi_2 \bydef \sqrt{\left|\Z^n / \Ga\right|} e^{-\pi ik_2 ~\cdot~ \widetilde\Theta Ml_2}\widetilde{U}_{k_2 } \otimes \xi_{l_2}
\eean
then from \eqref{nt_hh_product_eqn} it follows that 
\bean
\begin{split}
	\left(\widetilde \xi_{1} ,\widetilde  \xi_{2} \right)_{\widetilde \H}= 	e^{-\pi ik_1 ~\cdot~ \widetilde\Theta k_2}e^{\pi i\left(k_2 - k_1 \right) ~\cdot~ \widetilde\Theta Ml_2}e^{\pi ik_1 ~\cdot~ \widetilde\Theta Ml_1} e^{-\pi ik_2 ~\cdot~ \widetilde\Theta Ml_2}=e^{\pi ix},\\
	\text{where } x = - k_1 \cdot \widetilde\Theta k_2 + \left(k_2 - k_1 \right) \cdot \widetilde\Theta Ml_2+ k_1 \cdot \widetilde\Theta Ml_1- k_2 ~\cdot~ \widetilde\Theta Ml_2=
	\\
	= - k_1 \cdot \widetilde\Theta k_2 +k_1\cdot \left(Ml_2 - Ml_1 \right)=  -k_1 \cdot \widetilde\Theta k_2  + k_1 \cdot \widetilde\Theta\left(k_2 - k_1 \right) 
\end{split}
\eean
Since $\widetilde\Theta$ is skew symmetric one has $x = 0$ it follows that $\left(\widetilde \xi_{1} ,\widetilde  \xi_{2} \right)_{\widetilde \H}=1$ and taking into account $\left\| \widetilde \xi_1\right\|_{\widetilde\H}=\left\| \widetilde \xi_2\right\|_{\widetilde\H}=1$ one has $\widetilde\xi_1 = \widetilde\xi_2$. So there is the orthonormal basis $\left\{\widetilde\xi_p\right\}_{p \in \Z^n}$ such that
\be\label{nt_basis_eqn}
\widetilde\xi_{p} \bydef \sqrt{\left|\Z^n / \Ga\right|} e^{-\pi ik ~\cdot~ \widetilde\Theta Ml}\widetilde{U}_{k} \otimes \xi_{l}\text{ where } p = Ml + k.
\ee
and for all $p\in \Z^n$ the element  $\widetilde\xi_{p}$ does not depend on the representation $p = Ml + k$. Direct check shows that
\be\label{nt_basis_act_eqn}
\widetilde U_k \widetilde\xi_{l} = 	e^{\pi ik ~\cdot~ \widetilde\Theta l}\widetilde\xi_{Mk + l}.
\ee
Using equations  \eqref{nt_basis_eqn} and \eqref{nt_l2_eqn} one has the following lemma.

\begin{lemma}
	If	$\widetilde \rho: C\left(\T^{n}_{\widetilde{\Th}} \right) \to B\left(\widetilde \H \right)$ is the induced by the pair 
	$$
	\left(\rho, \left( C\left(\T^{n}_{{\Th}} \right), C\left(\T^{n}_{\widetilde{\Th}} \right), \Z^n / \Ga, \pi\right)\right)
	$$ then $\widetilde \rho$ is equivalent to the described  \ref{nt_gns_empt} faithful GNS-representation.
\end{lemma}

\subsection{Coverings of $O^*$-algebras}\label{nt_fin_o_sec}

\paragraph*{}
There is the given by \eqref{nt_repr_eqn} representation $\pi: C\left(\T^n_\Th \right)\to B\left(L^2\left( \T^n_\Th, \tau\right)  \right)$. If  $	\Psi_\Th:C\left(\mathbb{T}^n_{\Theta}\right) \hookto L^2\left(C\left(\mathbb{T}^n_{\Theta}\right), \tau\right)$ is the given by  \eqref{nt_to_hilbert_eqn}  and
\be\label{nt_dh_eqn}
\D\left(\mathbb{T}^n_{\Theta}\right)  \bydef \Psi_\Th\left(\Coo\left(\mathbb{T}^n_{\Theta}\right) \right) 
\ee
then 	$\D\left(\mathbb{T}^n_{\Theta}\right)$ is a dense subset of $L^2\left(C\left(\mathbb{T}^n_{\Theta}\right). \tau\right)$,

\begin{lemma}
	In the above situation there is an $*$-isomorphism
	$$
	\Coo\left(\mathbb{T}^{2N}_{\Theta}\right)\cong   C\left(\mathbb{T}^{2N}_{\Theta}\right)\cap \L^\dagger\left(	\D\left(\mathbb{T}^n_{\Theta}\right)\right) .
	$$
\end{lemma}
\begin{proof}
	If $\xi \in \D\left(\mathbb{T}^n_{\Theta}\right)$ then there is $a_\xi \in \Coo\left(\mathbb{T}^n_{\Theta}\right)$ such that $\xi = \Psi_\Th\left(a_\psi \right)$. For all $a \in \Coo\left(\mathbb{T}^n_{\Theta}\right)$
	one has
	$$
	\pi\left(a \right)\xi =  \Psi_\Th\left(aa_\xi \right)
	$$ 
	and from $aa_\psi\in \Coo\left(\mathbb{T}^n_{\Theta}\right)$ it turns out that
	$$
	\pi\left(a \right)\xi \in  \Psi_\Th\left(\Coo\left(\mathbb{T}^n_{\Theta}\right) \right)= \D\left(\mathbb{T}^n_{\Theta}\right)
	$$
	It follows that
	$$
	\pi\left(\Coo\left(\mathbb{T}^n_{\Theta}\right)\right)\D\left(\mathbb{T}^n_{\Theta}\right) = \D\left(\mathbb{T}^n_{\Theta}\right),
	$$
	or, equivalently
	$$
	\pi\left(\Coo\left(\mathbb{T}^n_{\Theta}\right)\right) \subset \L^\dagger\left(	\D\left(\mathbb{T}^n_{\Theta}\right)\right).
	$$
	So from $\Coo\left(\mathbb{T}^n_{\Theta}\right)\subset C\left(\mathbb{T}^n_{\Theta}\right)$ it turns out that
	$$
	\Coo\left(\mathbb{T}^{2N}_{\Theta}\right)\subset C\left(\mathbb{T}^{2N}_{\Theta}\right)\cap \L^\dagger\left(	\D\left(\mathbb{T}^n_{\Theta}\right)\right).
	$$
	Conversely suppose that $
	C\left(\mathbb{T}^{2N}_{\Theta}\right)\cap \L^\dagger\left(	\D\left(\mathbb{T}^n_{\Theta}\right)\right)\not\subset\Coo\left(\mathbb{T}^{2N}_{\Theta}\right).
	$, i.e. there is $b\in C\left(\mathbb{T}^{2N}_{\Theta}\right)\cap \L^\dagger\left(	\D\left(\mathbb{T}^n_{\Theta}\right)\right)$ such that $b\notin \Coo\left(\mathbb{T}^{2N}_{\Theta}\right)$.
	If $\xi \bydef \Psi_\Th\left( 1_{\Coo\left(\mathbb{T}^{2N}_{\Theta}\right)}\right)\in \D\left(\mathbb{T}^n_{\Theta}\right)$ then 
	$$
	\pi\left( b\right) \xi = \Psi_\Th\left(b\right)\notin \D\left(\mathbb{T}^n_{\Theta}\right)
	$$ 
	and
	$$
	b\notin \L^\dagger\left(	\D\left(\mathbb{T}^n_{\Theta}\right)\right);
	$$
	It contradicts with $b\in C\left(\mathbb{T}^{2N}_{\Theta}\right)\cap \L^\dagger\left(	\D\left(\mathbb{T}^n_{\Theta}\right)\right)$. So one has
	$$
	C\left(\mathbb{T}^{2N}_{\Theta}\right)\cap \L^\dagger\left(	\D\left(\mathbb{T}^n_{\Theta}\right)\right)\subset\Coo\left(\mathbb{T}^{2N}_{\Theta}\right).
	$$
	
\end{proof}

Let $\left( C\left(\T^{n}_{{\Th}} \right), C\left(\T^{n}_{\widetilde{\Th}} \right), \Z^n / \Ga, \pi\right)$ be an unital noncommutative finite-fold  covering given by the Theorem  \ref{nt_fin_thm}. 
If a *-algebra  $\Coo\left(\mathbb{T}^n_{\widetilde\Theta}\right)$  is given by \eqref{nt_coo_eqn} and $\tau$ is given by \eqref{nt_state_eqn} then described in \ref{nt_gns_empt} GNS- construction yields the Hilbert space $\widetilde \H\bydef  L^2\left(C\left(\mathbb{T}^n_{\widetilde\Theta}\right), \tau\right)$, such that $\widetilde  \D\bydef \Coo\left(\mathbb{T}^n_{\widetilde\Theta}\right)$ is a dense subspace of $\widetilde \H$. If $\widetilde \delta_\mu: \Coo\left(\mathbb{T}^n_{\widetilde\Theta}\right)\to \Coo\left(\mathbb{T}^n_{\widetilde\Theta}\right)$ are given by \eqref{nt_diff_eqn} then $\widetilde\delta_\mu \in \L^\dagger\left(\widetilde\D \right)$. From  $\delta_\mu\left(\Coo\left(\mathbb{T}^n_{\Theta}\right) \right) \subset \Coo\left(\mathbb{T}^n_{\Theta}\right)$ operators can be regarded as linear maps: $\widetilde \delta_\mu: \Coo\left(\mathbb{T}^n_{\Theta}\right)\to \Coo\left(\mathbb{T}^n_{\Theta}\right)$. Denote by $D\left( \mathbb{T}^n_{\widetilde\Theta}\right)\subset \L^\dagger\left(\widetilde\D \right)$ a generated by $\Coo\left( \mathbb{T}^n_{\widetilde\Theta}\right)\sqcup\left\{\widetilde\delta_1,..., \widetilde\delta_n \right\}$ subalgebra of $\L^\dagger\left(\widetilde\D \right)$. There is the unique action $\Z^n / \Ga\times\times D\left( \mathbb{T}^n_{\widetilde\Theta}\right)\to D\left( \mathbb{T}^n_{\widetilde\Theta}\right)$ which is induced by the action $\Z^n / \Ga\times\times \Coo\left( \mathbb{T}^n_{\widetilde\Theta}\right)\to \Coo\left( \mathbb{T}^n_{\widetilde\Theta}\right)$ and such that  $\delta_\mu \in D\left( \mathbb{T}^n_{\widetilde\Theta}\right)^{\Z^n / \Ga}$. If
\be\label{nt_d_inv_eqn}
D\left( \mathbb{T}^n_{\Theta}\right)\bydef D\left( \mathbb{T}^n_{\widetilde\Theta}\right)^{\Z^n / \Ga}.
\ee
then there is an injective $*$-homomorphism $\pi: D\left( \mathbb{T}^n_{\Theta}\right)\hookto D\left( \mathbb{T}^n_{\widetilde\Theta}\right).
$
\begin{theorem}\label{nt_d_fin_thm}
The quadruple $\left( D\left( \mathbb{T}^n_{\Theta}\right), D\left( \mathbb{T}^n_{\widetilde\Theta}\right), {\Z^n / \Ga}, \pi\right)$ is a \textit{noncommutative finite-fold covering of $O^*$-algebras} (cf. Definition \ref{fino*_defn}).
\end{theorem}
\begin{proof}
Firstly we prove that  $\left( D\left( \mathbb{T}^n_{\Theta}\right), D\left( \mathbb{T}^n_{\widetilde\Theta}\right), {\Z^n / \Ga}, \pi\right)$ is a pre-covering.
If $D\left( \mathbb{T}^n_{\Theta}\right)_b$ and  $D\left( \mathbb{T}^n_{\widetilde\Theta}\right)_b$ are given by the equation \eqref{o*b_eqn} then one has $D\left( \mathbb{T}^n_{\Theta}\right)_b= \Coo\left( \mathbb{T}^n_{\Theta}\right)$ and  $D\left( \mathbb{T}^n_{\widetilde\Theta}\right)_b= \Coo\left( \mathbb{T}^n_{\widetilde \Theta}\right)$. If 
$$
G \bydef \left\{\left.g \in \Aut\left( D\left( \mathbb{T}^n_{\widetilde\Theta}\right)\right)\right| ga= a\quad \forall a \in D\left( \mathbb{T}^n_{\Theta}\right) \right\}
$$
then from the Remark \ref{o*b_rem} it follows that $G D\left( \mathbb{T}^n_{\widetilde\Theta}\right)_b= D\left( \mathbb{T}^n_{\widetilde\Theta}\right)_b$, or equivalently $G \Coo\left( \mathbb{T}^n_{\widetilde\Theta}\right)= \Coo\left( \mathbb{T}^n_{\widetilde\Theta}\right)$. So there is a surjective homomorphism of groups
\bean
\phi: G \to \left\{\left.g \in \Aut\left( \Coo\left( \mathbb{T}^n_{\widetilde\Theta}\right)\right)\right| ga= a \quad \forall a \in \Coo\left( \mathbb{T}^n_{\Theta}\right) \right\}\cong \Z^n / \Ga,\\
g \mapsto g|_{\Coo\left( \mathbb{T}^n_{\widetilde\Theta}\right)}.
\eean
From $\delta_\mu \in D\left( \mathbb{T}^n_{\widetilde\Theta}\right)^{\Z^n / \Ga}$ it turns out that $\phi$ is isomorphism. Taking into  account \eqref{nt_d_inv_eqn} one concludes that the quadruple $\left( D\left( \mathbb{T}^n_{\Theta}\right), D\left( \mathbb{T}^n_{\widetilde\Theta}\right), {\Z^n / \Ga}, \pi\right)$ is a noncommutative finite-fold pre-covering. The $C^*$-norm completions of both $D\left( \mathbb{T}^n_{\Theta}\right)_b= \Coo\left( \mathbb{T}^n_{\Theta}\right)$ and  $D\left( \mathbb{T}^n_{\widetilde\Theta}\right)_b= \Coo\left( \mathbb{T}^n_{\widetilde \Theta}\right)$ equal to $C\left( \mathbb{T}^n_{\Theta}\right)$ and  $C\left( \mathbb{T}^n_{\widetilde \Theta}\right)$ respectively. On the other hand from the Theorem \ref{nt_fin_thm}
	In the above situation the triple $\left( C\left(\T^{n}_{{\Th}} \right), C\left(\T^{n}_{\widetilde{\Th}} \right), \Z^n / \Ga, \overline \pi\right)$ is an unital noncommutative finite-fold  covering where $\overline \pi$ comes from $\pi$. Thus  $\left( D\left( \mathbb{T}^n_{\Theta}\right), D\left( \mathbb{T}^n_{\widetilde\Theta}\right), {\Z^n / \Ga}, \pi\right)$ satisfies to all conditions of the Definition \ref{fino*_defn}). 
	
\end{proof}
\begin{exercise}\label{nt_st_fin_exer}
	Prove that the quadruple $\left( \Coo\left( \mathbb{T}^n_{\Theta}\right), \Coo\left( \mathbb{T}^n_{\widetilde\Theta}\right), {\Z^n / \Ga}, \pi_{\Coo\left( \mathbb{T}^n_{\Theta}\right)}\right)$ is a \textit{noncommutative finite-fold covering of $O^*$-algebras} (cf. Definition \ref{fino*_defn}).
\end{exercise}

\begin{exercise}\label{nt_d_fin_exer}
	Prove that the quadruple $\left( \Coo\left( \mathbb{T}^n_{\Theta}\right), \Coo\left( \mathbb{T}^n_{\widetilde\Theta}\right), {\Z^n / \Ga}, \pi|_{\Coo\left( \mathbb{T}^n_{\Theta}\right)}\right)$ is a {noncommutative finite-fold covering of $O^*$-algebras} (cf. Definition \ref{fino*_defn}).
\end{exercise}

\subsection{Coverings of quasi $*$-algebras}
If $\Coo\left( \T^n\right)'$ is  the space of distribution densities (cf. Definition \ref{top_distr_dens_def}) and $ \<\cdot , \cdot >:\Coo\left( \T^n\right)'\times \Coo\left( \T^n\right)\to \C$ is the natural pairing then similarly to \ref{mp_star_ext_eqn} we define following products: 
\be\label{nt_distr_prod}
\begin{split}
\Coo\left(\mathbb{T}^{2N}_{\Theta}\right)' \times \Coo\left(\mathbb{T}^{2N}_{\Theta}\right) \to \Coo\left(\mathbb{T}^{2N}_{\Theta}\right)'\quad 
	\<T a, \xi>\bydef \<T, a  \xi>;\\
\Coo\left(\mathbb{T}^{2N}_{\Theta}\right) \times \Coo\left(\mathbb{T}^{2N}_{\Theta}\right)' \to \Coo\left(\mathbb{T}^{2N}_{\Theta}\right)'\quad 
\<aT , \xi>\bydef \<T, \xi a>.
\end{split}
\ee
The products \eqref{nt_distr_prod} yield a quasi *-algebra
\be
\left(\Coo\left(\mathbb{T}^{2N}_{\Theta}\right)', \Coo\left(\mathbb{T}^{2N}_{\Theta}\right)  \right) 
\ee
If $\left( \Coo\left( \mathbb{T}^{2N}_{\Theta}\right), \Coo\left( \mathbb{T}^{2N}_{\widetilde\Theta}\right), {\Z^{2N} / \Ga}, \pi\right)$ is given by the Exercise \ref{nt_d_fin_exer} then there is an injective *-homomoprhism
\be
\pi':\left(\Coo\left(\mathbb{T}^{2N}_{\Theta}\right)', \Coo\left(\mathbb{T}^{2N}_{\Theta}\right)  \right)\hookto \left(\Coo\left(\mathbb{T}^{2N}_{\widetilde\Theta}\right)', \Coo\left(\mathbb{T}^{2N}_{\widetilde\Theta}\right)  \right) 
\ee
of quasi *-algebras such that $\left.\pi'\right|_{\Coo\left(\mathbb{T}^{2N}_{\Theta}\right)}= \left.\pi\right|_{\Coo\left(\mathbb{T}^{2N}_{\Theta}\right)}$.

\begin{exercise}\label{nt_quasi_fin_exer}
In the described situation prove that the quadruple
\be
\left(\left(\Coo\left(\mathbb{T}^{2N}_{\Theta}\right)', \Coo\left(\mathbb{T}^{2N}_{\Theta}\right)  \right), \left(\Coo\left(\mathbb{T}^{2N}_{\widetilde\Theta}\right)', \Coo\left(\mathbb{T}^{2N}_{\widetilde\Theta}\right)  \right), \Z^{2N} / \Ga ,\pi'\right) 
\ee
is noncommutative finite-fold covering of quasi $*$-algebras (cf. Definition \ref{oq*fin_defn})
\end{exercise}

\subsection{Coverings of spectral triples}
\paragraph*{}
If  $\left( \Coo\left( \mathbb{T}^n_{\Theta}\right), \Coo\left( \mathbb{T}^n_{\widetilde\Theta}\right), {\Z^n / \Ga}, \pi|_{\Coo\left( \mathbb{T}^n_{\Theta}\right)}\right)$ is a given by Exercise \ref{nt_d_fin_exer} {noncommutative finite-fold covering of $O^*$-algebras} (cf. Definition \ref{fino*_defn}), then any $\widetilde a\in \Coo\left( \mathbb{T}^n_{\widetilde\Theta}\right)$ can be represented by $\Ga$-periodic element of $\Coo\left(\R^n \right)$. On the other hand if $u_1, ..., u_\mu$ are generators of  $C\left( \mathbb{T}^n_{\Theta}\right)$ (cf. Definition \ref{nt_uni_defn}) then  from the equation \eqref{nt_chain_eqn}  it follows that the operator $\frac{\partial}{\partial u_\mu}$ maps $\Ga$-periodic elements of $\Coo\left(\R^n \right)$ onto $\Ga$-periodic ones. So if $\delta_{\mu}\in \L^\dagger\left(\Coo\left(\mathbb{T}^n_{\Theta}\right)\right)$ is given by the equation \eqref{nt_diff_eqn}  then using the equation \eqref{nt_chaind_eqn} one can define
\be
\delta_{\mu}\in \L^\dagger\left(\Coo\left(\mathbb{T}^n_{\widetilde\Theta}\right)\right).
\ee
From the equation \eqref{mp_part_prod_eqn} the reader can deduce the following Leibniz rule.
\be\label{nt_tle_eqn}
\forall\widetilde a,\widetilde b \in C^{\infty}(\mathbb{T}^n_{\widetilde\Theta})\quad \delta_{\mu}\left(\widetilde a \widetilde b\right)  =\left( \delta_\mu \widetilde a\right)\widetilde b +\widetilde a \left( \delta_\mu\widetilde b\right).
\ee

Let $\Om^1_D$ be the {module of differential forms associated} with the given by \eqref{nt_sp_tr_eqn}  spectral triple  $\left( C^{\infty}(\mathbb{T}^n_{\Theta}),L^2\left(C\left(\mathbb{T}^n_{\Theta}\right), \tau\right)\otimes\mathbb{C}^{m},D\right)$ (cf. Definition \ref{ass_cycle_defn}). Denote a map
\be\label{nt_conn_eqn}
\begin{split}
	\nabla: \Coo\left(\mathbb{T}^n_{\widetilde\Theta}\right)\to \Coo\left(\mathbb{T}^n_{\widetilde\Theta}\right)\otimes_{\Coo\left(\mathbb{T}^n_{\Theta}\right)}\Om^1_D,\\
	\widetilde a \mapsto     = \sum_{\mu = 1}^n  \delta_\mu \widetilde  a \otimes \left(1_{\Coo\left(\mathbb{T}^n_{\Theta}\right)}\otimes \ga^\mu\right) 
\end{split}
\ee
where the equation \eqref{nt_ga_in_eqn} is used.
For any $a\in \Coo\left(\mathbb{T}^n_{\Theta}\right)$ one has
\bean
\nabla\left( \widetilde a a\right) = \sum_{\mu = 1}^n  \delta_\mu \left( \widetilde  a a\right) \otimes \left(1_{\Coo\left(\mathbb{T}^n_{\Theta}\right)}\otimes \ga^\mu\right) 
\eean
and taking into account the Leibniz rule \eqref{nt_tle_eqn} we obtain the following
\be\label{nt_conn'_eqn}
\begin{split}
	\nabla\left( \widetilde a a\right)= \sum_{\mu = 1}^n \left(   \left(\delta_\mu \widetilde  a \right) \cdot a\right) \otimes \left(1_{\Coo\left(\mathbb{T}^n_{\Theta}\right)}\otimes \ga^\mu\right)+\\+ \sum_{\mu = 1}^n  \widetilde  a \cdot \left(\delta_\mu  a\right) \otimes \left(1_{\Coo\left(\mathbb{T}^n_{\Theta}\right)}\otimes \ga^\mu\right)=\\=
	\nabla\left( \widetilde a\right)\left(a \otimes  1_{\mathbb{M}_m\left( \C\right)}\right) + \sum_{\mu = 1}^n  \widetilde  a \cdot \left(\delta_\mu  a\right) \otimes \left(1_{\Coo\left(\mathbb{T}^n_{\Theta}\right)}\otimes \ga^\mu\right).
\end{split}
\ee
On the other hand for all $b \otimes \xi \in  \Coo\left(\mathbb{T}^n_{\Theta}\right)\otimes\C^m$ one has
\bean
\left[D, a \otimes 1_{\mathbb{M}_m\left( \C\right)}\right]  \left(b \otimes \xi\right)= \sum_{\mu = 1}^n \delta_\mu \left(ab\right)\otimes \ga^\mu\xi- \sum_{\mu = 1}^n a\delta_\mu b\otimes \ga^\mu\xi=\\
= \sum_{\mu = 1}^n \left(\left(\delta_\mu a \right) b  + a\left( \delta_\mu b\right)-a\left( \delta_\mu b\right)\right)\otimes\ga^\mu\xi= \sum_{\mu = 1}^n \left( \delta_\mu a \otimes \ga^\mu \right)  \left(b \otimes \xi\right),
\eean
or, equivalently
\be\label{nt_conn''_eqn}
\left[D, a \otimes 1_{\mathbb{M}_m\left( \C\right)}\right] =  \sum_{\mu = 1}^n \delta_\mu a \otimes \ga^\mu.
\ee 
A substitution of  \eqref{nt_conn''_eqn} into \eqref{nt_conn'_eqn} yields the following equation
\bean
\nabla\left( \widetilde a a\right)=
\nabla\left( \widetilde a\right)\left(a \otimes  1_{\mathbb{M}_m\left( \C\right)}\right) + \widetilde a \left[ \left(D, a \otimes  1_{\mathbb{M}_m\left( \C\right)}\right)\right] 
\eean 
which is a specialization of \eqref{conn_triple_eqn}, i.e. $\nabla$ is a connection.
\begin{exercise} Prove following statements:
	\begin{enumerate}
		\item The given by \eqref{nt_conn_eqn} connection is $\Z^n/\Ga$-invariant (cf. equation \ref{conn_equ_conn}).
		\item If $\widetilde D \in \L^\dagger\left( \Coo\left(\mathbb{T}_{\widetilde\Th}^n \right)\otimes \C^m \right)$ is a specialization   of given by \eqref{wtd_eqn} operator then one has
		\begin{equation}\label{nt_tdirac_u_eqn}
			\widetilde	D\stackrel{\mathrm{def}}{=}\sum_{\mu =1}^{n}\delta_{\mu} \otimes\gamma^{\mu}= \sum_{\mu =1}^{n}i u_\mu \frac{\partial}{\partial u_\mu} \otimes\gamma^{\mu}\in \L^\dagger\left(\Coo\left(\mathbb{T}^n_{\widetilde\Theta}\right)\otimes\C^m\right).
		\end{equation}
		\item Let $\widetilde\H$ be a Hilbert space of the induced by a pair
		$$
		\left(\left( C\left(\mathbb{T}_{\Th}^n \right)  \hookto B \left(L^2\left(C\left(\mathbb{T}^n_{\Theta}\right), \tau\right)\right) \right) , \left( C\left(\T^{n}_{{\Th}} \right), C\left(\T^{n}_{\widetilde{\Th}} \right), \Z^n / \Ga, \pi\right)\right)
		$$
		representation $C\left(\mathbb{T}_{\widetilde\Th}^n \right)\hookto B\left(\widetilde\H\right)$ (cf. Definition \ref{induced_repr_fin_defn}). The spectral triple $\left( \Coo\left(\mathbb{T}_{\widetilde\Th}^n \right),  {\widetilde\H}, \widetilde D\right)$ is the $\left( C\left(\T^{n}_{{\Th}} \right), C\left(\T^{n}_{\widetilde{\Th}} \right), \Z^n / \Ga, \pi\right)$-lift of\\ $\left( \Coo\left(\mathbb{T}_\Th^n \right),  L^2\left(C\left(\mathbb{T}^n_{\Theta}\right), \tau\right), D\right)$ (cf. Definition \ref{spectral_triple_fin_lift_defn}).
	\end{enumerate}
	
\end{exercise}

\begin{exercise}
	Consider a given by \eqref{nt_sp_tr_eqn} spectral triple
	\begin{equation}\nonumber
		\left( \Coo\left(\mathbb{T}_\Th^n \right),  L^2\left(C\left(\mathbb{T}^n_{\Theta}\right), \tau\right), D\right).
	\end{equation} 
	(cf. Definition  \ref{df:spec-triple}). Let $\left( C\left(\T^{n}_{{\Th}} \right), C\left(\T^{n}_{\widetilde{\Th}} \right), \Z^n / \Ga, \pi\right)$ be a given by the Theorem \ref{nt_fin_thm} unital noncommutative finite-fold  covering.
	Prove that  $$\left( \Coo\left(\mathbb{T}_{\widetilde\Th}^n \right),  L^2\left(\left(\mathbb{T}_{\widetilde\Th}^n \right), \widetilde\tau\right), \widetilde D\right)$$ is the $\left( C\left(\T^{n}_{{\Th}} \right), C\left(\T^{n}_{\widetilde{\Th}} \right), \Z^n / \Ga, \pi\right)$-lift of $\left( \Coo\left(\mathbb{T}_\Th^n \right),  L^2\left(C\left(\mathbb{T}^n_{\Theta}\right), \tau\right), D\right)$ (cf. Definition \ref{spectral_triple_fin_lift_defn}).
\end{exercise}

\chapter{Isospectral deformations and their coverings}\label{isospectral_chap}
\section{Finite-fold coverings}\label{isosectral_fin_cov}
\subsection{Basic construction}
\paragraph{}
Let $M$ be Riemannian manifold which admits a Spin$^{c}$-structure (cf. Definition \ref{spin_str_defn}). Suppose that there is a    the smooth action of $\T^2\times M \to M$. 
 Let $\widetilde{x}_0 \in \widetilde{M}$ and $x_0=\pi\left(\widetilde{x}_0 \right)$. Denote by $\varphi: \R^2 \to \R^2 / \Z^2 = \T^2$ the natural covering. There are two closed paths $\om_1, \om_2: \left[0,1 \right]\to M$ given by
\begin{equation*}
\begin{split}
\om_1\left(t \right) = \varphi\left(t, 0 \right) x_0,~
\om_2\left(t \right) = \varphi\left(0, t \right) x_0.
\end{split}
\end{equation*} 
There are  lifts of these paths, i.e. maps $\widetilde{\om}_1 , \widetilde{\om}_2: \left[0,1 \right] \to\widetilde{M}$ such that
\begin{equation*}
\begin{split}
\widetilde{\om}_1\left(0 \right)= \widetilde{\om}_2\left(0 \right)=\widetilde{x}_0,~ \\
\pi\left( \widetilde{\om}_1\left(t \right)\right)  = \om_1\left(t\right),\\
\pi\left( \widetilde{\om}_2\left(t \right)\right)  = \om_2\left(t\right).
\end{split}
\end{equation*}
Since $\pi$ is a finite-fold covering there are $N_1, N_2 \in \N$ such that if 
$$
\gamma_1\left(t \right) = \varphi\left(N_1t, 0 \right) x_0,~
\gamma_2\left(t \right) = \varphi\left(0, N_2t \right) x_0.
$$
and $\widetilde{\gamma}_1$ (resp. $\widetilde{\gamma}_2$) is the lift of $\gamma_1$ (resp. $\gamma_2$) then both $\widetilde{\gamma}_1$, $\widetilde{\gamma}_2$ are closed. Let us select minimal values of $N_1, N_2$. If $\text{pr}_n: S^1 \to S^1$ is an $n$ listed covering and $\text{pr}_{N_1, N2}$ the covering given by
$$
\widetilde{\T}^2 = S^1 \times S^1 \xrightarrow{\text{pr}_{N_1}\times\text{pr}_{N_2}} \to S^1 \times S^1 = \T^2
$$
then there is the action $\widetilde{\T}^2 \times \widetilde{M} \to  \widetilde{M}$ such that

\begin{tikzpicture}
\matrix (m) [matrix of math nodes,row sep=3em,column sep=4em,minimum width=2em]
{
	\widetilde{\T}^2 \times \widetilde{M}  	&  	 & \widetilde{M} \\
	\T^2. \times M 	 & 	 & M   \\};
\path[-stealth]
(m-1-1) edge node [above] {$ \ $} (m-1-3)
(m-1-1) edge node [right] {$\mathrm{pr}_{N_1N_2} \times \pi $} (m-2-1)
(m-1-3) edge node [right] {$\pi$} (m-2-3)
(m-2-1) edge node [above] {$ \ $} (m-2-3);

\end{tikzpicture}

where  $\widetilde{\T}^2 \approx \T^2$.
Let $\widetilde{p} = \left( \widetilde{p}_1, \widetilde{p}_2\right) $ be the generator of the associated with $\widetilde{\T}^2$ two-parameters group $\widetilde{U}\left(s \right) $
so that
\begin{equation*}
\widetilde{U}\left(s \right) = \exp\left( i\left( s_1 \widetilde{p}_1 + s_2 \widetilde{p}_2\right)\right).
\end{equation*}	
The covering $\widetilde{M} \to M$ induces an involutive injective homomorphism
\begin{equation*}
\varphi:\Coo\left(M \right)\to\Coo\left( \widetilde{M} \right).
\end{equation*}

Suppose $M \to M/\T^2$ is submersion, and suppose there is  a weak fibration $\T^2 \to M \to M/\T^2$ (cf. \cite{spanier:at})
There is the exact \textit{homotopy sequence of the weak fibration}
\bean
\dots \to \pi_n\left( \T^2, e_0\right) \xrightarrow{i_{\#}} \pi_n\left( M, e_0\right)  \xrightarrow{p_{\#}} \pi_n\left( M/\T^2, b_0\right) \xrightarrow{\overline\partial} \pi_{n-1}\left( \T^2, e_0\right) \to \dots\\
\dots \to \pi_2\left( M/\T^2, b_0\right) \xrightarrow{\overline\partial} \pi_1\left( \T^2, e_0\right) \xrightarrow{i_{\#}} \pi_1\left( M, e_0\right)  \xrightarrow{p_{\#}} \pi_1\left( M/\T^2, b_0\right) \xrightarrow{\overline\partial} \pi_0\left( \T^2, e_0\right) \to \dots 
\eean
(cf. \cite{spanier:at}) where $\pi_n$ is the $n^{\mathrm{th}}$ homotopical group for any $n\in \N^0$.
If $\pi:\widetilde{M} \to M$ is a finite-fold regular covering then there is the natural surjective homomorphism $
\pi_1\left( M, e_0\right) \to G\left( \left.\widetilde{M}~\right|M\right)$. If $\pi:\widetilde{M} \to M$ induces a covering $\pi:\widetilde{M}/\T^2 \to M / \T^2$  then 
the  homomorphism $
\varphi:\pi_1\left( M, e_0\right) \to G\left( \left.\widetilde{M}~\right|M\right)$ can be included into the following commutative diagram.

\begin{tikzpicture}
\matrix (m) [matrix of math nodes,row sep=3em,column sep=4em,minimum width=2em]
{
	\pi_1\left(\T^2, e_0 \right)\cong \Z^2 &  \pi_1\left( M, e_0\right) &  \pi_1\left( M/\T^2, b_0\right) &  \{e\} \\
	G\left(\widetilde{\T}^2~|\T^2 \right)  & G\left( \left.\widetilde{M}~\right|M\right) & G\left(\widetilde{M}/\widetilde{\T}^2~|~M/\T^2\right)  & \{e\} \ \\};
\path[-stealth]
(m-1-1)  edge node [right] {$\varphi'$} (m-2-1)
(m-1-2)  edge node [right] {$\varphi$} (m-2-2)
(m-1-3)  edge node [right] {$\varphi''$} (m-2-3)
(m-1-1)  edge node [above] {$i_{\#} $} (m-1-2)
(m-1-2)  edge node [above] {$p_{\#} $} (m-1-3)
(m-2-1)  edge node [above] {$i_*$} (m-2-2)
(m-2-2)  edge node [above] {$p_*$} (m-2-3)
(m-1-3)  edge node [right] {$ $} (m-1-4)
(m-1-4)  edge node [right] {$ $} (m-2-4)
(m-2-3)  edge node [right] {$ $} (m-2-4);
\end{tikzpicture}

Denote by $G \stackrel{\mathrm{def}}{=}G\left( \left.\widetilde{M}~\right|M\right)$, $~G' \stackrel{\mathrm{def}}{=}G\left(\widetilde{\T}^2~|\T^2 \right)$, $~G'' \stackrel{\mathrm{def}}{=}G\left(\widetilde{M}/\widetilde{\T}^2~|~M/\T^2\right)$. From the above construction it turns out that $G' =G\left(\widetilde{\T}^2~|\T^2 \right) = \Z_{N_1} \times \Z_{N_2}$. Otherwise there is an inclusion of Abelian groups $G\left(\widetilde{\T}^2~|\T^2 \right) \subset \widetilde{\T}^2$. The action $\widetilde{\T}^2 \times \widetilde{M} \to  \widetilde{M}$ is free, so the action $G' \times \widetilde{M} \to $ is free, so the natural homomorphism $G' \to G$ is injective, hence there is an exact sequence of groups
\be\label{isospectral_gr_eqn}
\{e\} \to G' \to G \to G'' \to \{e\}.
\ee

Let $\th, \widetilde{\th} \in \R$ be such that
$$
\widetilde{\th}= \frac{\th +  n}{N_1N_2}, \text{ where }n \in \Z.
$$
If $\lambda= e^{2\pi i \th}$, $\widetilde{\lambda}= e^{2\pi i \widetilde{\th}}$ then
$
\lambda = \widetilde{\lambda}^{N_1N_2}.
$
There are isospectral deformations $\Coo\left(M_\th \right), \Coo\left( \widetilde{M}_{\widetilde{\th}} \right)$ and $\C$-linear isomorphisms
$l:\Coo\left(M \right) \to \Coo\left(M_\th \right)$, $\widetilde{l}:\Coo\left( \widetilde{M} \right) \to \Coo\left( \widetilde{M}_{\widetilde{\th}} \right)$.
These isomorphisms and the inclusion $\varphi$ induces the inclusion
\begin{equation*}
\begin{split}
\varphi_\th:\Coo\left(M_\th \right)\to\Coo\left( \widetilde{M}_{\widetilde{\th}} \right),
\\
\varphi_{\widetilde{\th}}\left(\Coo\left(M_\th \right) \right)_{n_1,n_2} \subset \Coo\left( \widetilde{M}_{\widetilde{\th}} \right)_{n_1N_1,~ n_2N_2}.
\end{split}
\end{equation*}
Denote by $G = G\left( \left.\widetilde{M}~\right|M\right)$ the group of covering transformations.   Since $\widetilde{l}$ is a $\C$-linear isomorphism the action of $G$ on $\Coo\left( \widetilde{M} \right)$ induces a $\C$-linear action  $G \times \Coo\left( \widetilde{M}_{ \widetilde{\th}} \right)  \to \Coo\left( \widetilde{M}_{ \widetilde{\th}} \right)$. According to the definition of the action of $\widetilde{\T}^2$ on $\widetilde{M}$ it follows that the action of $G$ commutes with the action of $\widetilde{\T}^2$.
It turns out
$$
g \Coo\left( \widetilde{M} \right)_{n_1,n_2} = \Coo\left( \widetilde{M} \right)_{n_1,n_2}
$$
for any $n_1, n_2 \in \Z$ and $g \in G$.
If $\widetilde{a} \in \Coo\left( \widetilde{M} \right)_{n_1,n_2}$, $\widetilde{b} \in \Coo\left( \widetilde{M} \right)_{n'_1,n'_2}$ then  $g\left( \widetilde{a}\widetilde{b}\right)= \left(g\widetilde{a} \right) \left(g\widetilde{b} \right)\in \Coo\left( \widetilde{M} \right)_{n_1+n'_1,n_2+n'_2} $. One has
\begin{equation*}
	\begin{split}
		\widetilde{l}\left(\widetilde{a}\right)\widetilde{l}\left(\widetilde{b}\right)= \widetilde{\la}^{n'_1n_2}\widetilde{l}\left(\widetilde{a}\widetilde{b}\right), \\
		\widetilde{\la}^{n_2\widetilde{p}_1}l\left( \widetilde{b}\right) = \widetilde{\la}^{n'_1n_2}l\left( \widetilde{b}\right) \widetilde{\la}^{n_2\widetilde{p}_1},\\
		\widetilde{l}\left(g \widetilde{a}\right)\widetilde{l}\left(g \widetilde{b}\right)= g \widetilde{a}\widetilde{\la}^{n_2\widetilde{p}_1}g \widetilde{b}\widetilde{\la}^{n'_2\widetilde{p}_1}= \widetilde{\la}^{n'_1n_2} g\left(\widetilde{a}\widetilde{b} \right) \widetilde{\la}^{\left( n_2+n_2'\right) \widetilde{p}_1}.
	\end{split}
\end{equation*}
On the other hand
\begin{equation*}
	\begin{split}
		g\left( \widetilde{l}\left(\widetilde{a}\right)\widetilde{l}\left(\widetilde{b}\right)\right) = g\left( \widetilde{\la}^{n'_1n_2}\widetilde{l}\left(\widetilde{a}\widetilde{b}\right)\right)= \widetilde{\la}^{n'_1n_2} g\left(\widetilde{a}\widetilde{b} \right) \widetilde{\la}^{\left( n_2+n_2'\right) \widetilde{p}_1}. 
	\end{split}
\end{equation*}
From above equations it turns out
$$
\widetilde{l}\left(g \widetilde{a}\right)\widetilde{l}\left(g \widetilde{b}\right) = g\left( \widetilde{l}\left(\widetilde{a}\right)\widetilde{l}\left(\widetilde{b}\right)\right),
$$
i.e. $g$ corresponds to automorphism of $\Coo\left( \widetilde{M}_{ \widetilde{\th}}\right)$. It turns out that $G$ is the group of automorphisms of $\Coo\left( \widetilde{M}_{ \widetilde{\th}}\right)$. From $\widetilde{a} \in \Coo\left( \widetilde{M}_{ \widetilde{\th}}\right)_{n_1,n_2}$ it follows that $\widetilde{a}^* \in \Coo\left( \widetilde{M}_{ \widetilde{\th}}\right)_{-n_1,-n_2}$. One has
$$
g\left(\left( \widetilde{l}\left(\widetilde{a}\right)\right)^* \right) =  g\left( \widetilde{\la}^{-n_2\widetilde{p_1}}\widetilde{a}^*\right) =
g \left(\widetilde{\la}^{n_1 n_2} \widetilde{a}^*\widetilde{\la}^{-n_2\widetilde{p_1}}\right) = \widetilde{\la}^{n_1 n_2} g\left(\widetilde{l}\left(\widetilde{a}^* \right)  \right). 
$$
On the other hand
$$
\left(g \widetilde{l}\left(\widetilde{a}\right) \right)^*= \left(\left( g \widetilde{a}\right)\widetilde{\la}^{n_2\widetilde{p_1}}  \right)^*=\widetilde{\la}^{-n_2\widetilde{p_1}}\left(ga^* \right) = \widetilde{\la}^{n_1 n_2}\left(ga^* \widetilde{\la}^{-n_2\widetilde{p_1}}\right)= \widetilde{\la}^{n_1 n_2} g\left(\widetilde{l}\left(\widetilde{a}^* \right)  \right),
$$
i.e. $g\left(\left( \widetilde{l}\left(\widetilde{a}\right)\right)^* \right) = \left(g \widetilde{l}\left(\widetilde{a}\right) \right)^*$.
It follows that $g$ corresponds to the involutive automorphism of $\Coo\left( \widetilde{M}_{ \widetilde{\th}}\right)$. Since  $\Coo\left( \widetilde{M}_{ \widetilde{\th}}\right)$ is dense in $C\left( \widetilde{M}_{ \widetilde{\th}}\right)$ there is the unique involutive action $G \times C\left( \widetilde{M}_{ \widetilde{\th}}\right) \to C\left( \widetilde{M}_{ \widetilde{\th}}\right)$.
For any $y_0 \in M/\T^2$ there is a point $x_0 \in M$ mapped onto $y_0$ and a connected submanifold $\mathcal U \subset M$ such that:
\begin{itemize}
	\item $\dim \mathcal U= \dim M-2$,
	\item $\mathcal U$ is transversal to orbits of $\T^2$-action,
	\item The  fibration $\T^2 \to \mathcal U \times \T^2 \to \mathcal U \times \T^2 / \T^2$ is the restriction of the fibration $\T^2 \to M \to M/\T^2$,
	\item The image  $\mathcal V_{y_0} \in M/\T^2$ of $\mathcal U \times \T^2$ in $M/\T^2$ is an open neighborhood of $y_0$,
	\item $\mathcal V_{y_0}$ is evenly covered by  $\widetilde{M}/ \widetilde{\T}^2 \to M / \T^2$.
\end{itemize}
It is clear that
$$
M/\T^2 = \bigcup_{y_0 \in M/\T^2} \mathcal V_{y_0}.
$$
Since $M/\T^2$ is compact there is a finite subset $I \in M/\T^2$ such that
$$
M/\T^2 = \bigcup_{y_0 \in I} \mathcal V_{y_0}.
$$
Above equation will be rewritten as
\be\label{isospectral_p_eqn}
M/\T^2 = \bigcup_{\iota \in I} \mathcal V_{\iota}
\ee
where $\iota$ is just an element of the  finite set $I$ and we denote corresponding transversal submanifold by $\mathcal{U}_\iota$. There is a smooth partition of unity subordinated to \eqref{isospectral_p_eqn}, i.e. there is a set $\left\{a_\iota \in \Coo\left(M /\T^2\right)  \right\}_{\iota \in I}$ of positive elements such that
\be\label{isospectral_part_eqn}
1_{C\left(M /\T^2\right) }= \sum_{\iota \in I} a_\iota~,\\
\ee
\be\nonumber
a_\iota\left(\left( M/\T^2\right)  \setminus \mathcal V_{\iota}\right) = \{0\}.
\ee
Denote by 
\be\label{isospectral_e_eqn}e_\iota \stackrel{\mathrm{def}}{=} \sqrt{a_\iota} \in \Coo\left( M /\T^2\right).
\ee  For any $\iota \in I$ we select an open subset $\widetilde{\mathcal V}_{\iota} \subset \widetilde{M}/\T^2$ which is homeomorphically mapped onto $\mathcal V_{\iota}$.
If $\widetilde{I}= G'' \times  \mathscr A$ then for any $\left(g'', \iota \right) \in \widetilde{I}$ we define
\be\label{isospectral_vgi_eqn}
\begin{split}
	\widetilde{\mathcal V}_{\left(g'', \iota \right)} = g''	\widetilde{\mathcal V}_{\iota}.
\end{split}
\ee
 Similarly we select a transversal submanifold
\be\nonumber
\widetilde{\mathcal U}_{\iota} \subset \widetilde{M}
\ee 
 which is homeomorphially mapped onto $\mathcal U_{\iota}$. For any $\left(g'', \iota \right) \in \widetilde{I}$ we define
 \be\label{isospectral_wgi_eqn}
 \begin{split}
 	\widetilde{\mathcal U}_{\left(g'', \iota \right)} = g	\widetilde{\mathcal U}_{\iota}.
 \end{split}
 \ee
 where $g \in G$ is an arbitrary element mapped to $g''$.
   The set $\mathcal V_{\iota}$ is evenly covered by $\pi'':\widetilde{M}/\widetilde{\T}^2\to M/\T^2$, so one has
\be\label{isostectral_empt_eqn}
g\widetilde{\mathcal V}_{\iota}\bigcap \widetilde{\mathcal V}_{\iota} = \emptyset; \text{ for any nontrivial } g \in G''.
\ee

If $\widetilde{e}_\iota \in \Coo\left(\widetilde M / \widetilde \T^2 \right) $ is given by 
$$
\widetilde{e}_\iota\left(\widetilde{x} \right) = 
\left\{\begin{array}{c l}
e_\iota\left( \pi''\left( \widetilde{x} \right) \right)  & \widetilde{x} \in \widetilde{\mathcal U}_\iota \\
0 & \widetilde{x} \notin \widetilde{\mathcal U}_\iota
\end{array}\right.
$$
then from \eqref{isospectral_part_eqn} and \eqref{isostectral_empt_eqn} it turns out
\begin{equation}
\label{isospectral_partg_eqn}
\begin{split}
1_{C\left(\widetilde M/ \widetilde \T^2 \right) } = \sum_{g \in G''} \sum_{\iota \in I} \widetilde{e}^2_{\iota},\\
\left( g\widetilde{e}_{\iota}\right)\widetilde{e}_{\iota} = 0; \text{ for any nontrivial } g \in G''. 
\end{split}
\end{equation}
If $\widetilde{I}= G'' \times I$ and $\widetilde{e}_{\left(g, \iota\right)} = g \widetilde{e}_{ \iota}$ for any $\left(g, \iota\right)  \in G'' \times I$ then from  \eqref{isospectral_partg_eqn} it turns out
\begin{equation}
\label{isospectral_parti_eqn}
\begin{split}
1_{C\left( \widetilde M / \widetilde \T^2 \right) } =\sum_{\widetilde{\iota} \in \widetilde{I}} \widetilde{e}^2_{\widetilde{\iota}},\\
\left( g\widetilde{e}_{\widetilde{\iota}}\right)\widetilde{e}_{\widetilde{\iota}} = 0; \text{ for any nontrivial } g \in G'',\\
1_{C\left( \widetilde M / \widetilde \T^2 \right) } =  \sum_{\widetilde{\iota} \in \widetilde{I}} \widetilde{e}_{\widetilde{\iota}}\left\rangle \right\langle\widetilde{e}_{\widetilde{\iota}}~.
\end{split}
\end{equation}
It is known that $C\left(\T^2 \right)$ is an universal commutative $C^*$-algebra generated by two unitary elements $u, v$, i.e. there are following relations
\be\label{isospectral_gen_eqn}
\begin{split}
uu^*=u^*u=vv^*=v^*v= 1_{C\left(\T^2 \right)},\\
uv = vu, ~u^*v = vu^*,~uv^*= v^*u, ~ u^*v^*=v^*u^*.
\end{split}
\ee
If $J = I \times \Z \times \Z$ then for any $\left({\iota}, j, k \right) \in J$ there is an element ${f}'_{\left({\iota}, j, k \right)} \in \Coo\left( {\mathcal U}_{{\iota}} \times {\T}^2\right) $ given by
\be\label{isospectral_f's_eqn}
{f}'_{\left({\iota}, j, k \right)} = {e}_{{\iota}} {u}^j{
	v}^k
\ee
where ${e}_{{\iota}}\in \Coo\left( {M}/ {\T}^2\right)$ is regarded as element of $\Coo\left( {M}\right)$.  
Let $p: {M} \to {M} / {\T}^2$.
Denote by ${f}_{\left({\iota}, j, k \right)}\in \Coo\left({M} \right)$ an element given by
\be\label{isospectral_fs_eqn}
\begin{split}
{f}_{\left({\iota}, j, k \right)}\left( {x}\right) =
\left\{\begin{array}{c l}
	{f}'_{\left({\iota}, j, k \right)}\left( {x}\right) & p\left({x} \right)  \in {\mathcal V}_{{\iota}} \\
	0 & p\left( {x}\right)  \notin {\mathcal V}_{{\iota}} 
\end{array}\right.,\\
\text{where the right part of the above equation  assumes the inclusion } \mathcal U_{{\iota}} \times {\T}^2 \hookto M.
\end{split}
\ee 
If we denote by $\widetilde{u}, \widetilde{v} \in U\left( C\left( \widetilde{\T^2} \right) \right)$ unitary generators of $C\left( \widetilde{\T^2} \right)$ then the covering $\pi':\widetilde{\T^2} \to \T^2$ corresponds to a $*$-homomorphism $C\left(\T^2 \right) \to  C\left( \widetilde{\T^2}\right)$ given by
\bean
u \mapsto \widetilde{u}^{N_1},\\
v \mapsto \widetilde{v}^{N_2}.
\eean
There is the natural action of $G\left( \widetilde{\T}^2~|~{\T^2}\right)\cong \Z_{N_1} \times \Z_{N_2}$ on $C\left(\T^2 \right)$ given by
\be\label{isospectral_uv_eqn}
\begin{split}
\left(\overline{k}_1, \overline{k}_2 \right) \widetilde{u} = e^{\frac{2\pi i k_1}{N_1}}\widetilde{u},\\
\left(\overline{k}_1, \overline{k}_2 \right) \widetilde{v} = e^{\frac{2\pi i k_2}{N_1}}\widetilde{v}
\end{split}
\ee
where $\left(\overline{k}_1, \overline{k}_2 \right) \in \Z_{N_1} \times \Z_{N_2}$. 
If we consider $C\left( \widetilde{\T}^2\right)_{C\left( {\T^2}\right)}$ as a right Hilbert module which corresponds to a finite-fold noncommutative covering then one has
\be\label{isospectral_tor_eqn}
\begin{split}
\left\langle \widetilde{u}^{j'} \widetilde{v}^{k'}, \widetilde{u}^{j''} \widetilde{v}^{k''}  \right\rangle_{C\left( \widetilde{\T^2}\right)} = N_1N_2\delta_{j'j''} \delta_{k'k''}1_{C\left( {\T^2}\right)},\\
1_{C\left( \widetilde{\T^2}\right)}= \frac{1}{N_1N_2}\sum_{\substack{j = 0\\ k = 0}}^{\substack{j = N_1\\ k = N_2}} \widetilde{u}^{j} \widetilde{v}^{k}\left\rangle \right\langle \widetilde{u}^{j} \widetilde{v}^{k}.
\end{split}
\ee
If $\widetilde{J} = \widetilde{I} \times \left\{0, \dots, N_1-1\right\} \times \left\{0, \dots, N_2-1\right\}$ then for any $\left(\widetilde{\iota}, j, k \right) \in \widetilde{J}$ there is an element $\widetilde{f}'_{\left(\widetilde{\iota}, j, k \right)} \in \Coo\left( \widetilde{\mathcal U}_{\widetilde{\iota}} \times \widetilde{\T}^2\right) $ given by
\be\label{isospectral_f'_eqn}
\widetilde{f}'_{\left(\widetilde{\iota}, j, k \right)} = \widetilde{e}_{\widetilde{\iota}} \widetilde{u}^j\widetilde{
v}^k
\ee
where $\widetilde{e}_{\widetilde{\iota}}\in \Coo\left( \widetilde{M}/ \widetilde{\T}^2\right)$ is regarded as element of $\Coo\left( \widetilde{M}\right)$.  
Let $p: \widetilde{M} \to \widetilde{M} / \widetilde{\T}^2$.
Denote by $\widetilde{f}_{\left(\widetilde{\iota}, j, k \right)}\in \Coo\left(\widetilde{M} \right)$ an element given by
\be\label{isospectral_f_eqn}
\begin{split}
\widetilde{f}_{\left(\widetilde{\iota}, j, k \right)}\left( \widetilde{x}\right) =
\left\{\begin{array}{c l}
	\widetilde{f}'_{\left(\widetilde{\iota}, j, k \right)}\left( \widetilde{x}\right) & p\left(\widetilde{x} \right)  \in \widetilde{\mathcal V}_{\widetilde{\iota}} \\
	0 & p\left( \widetilde{x}\right)  \notin \widetilde{\mathcal V}_{\widetilde{\iota}} .
\end{array}\right.\\
\text{where right the part of the above equation assumes the inclusion }  \widetilde{\mathcal U}_{ \widetilde{\iota}} \times  \widetilde{\T}^2 \hookto  \widetilde{M}.
\end{split}
\ee 
Any element $\widetilde{e}_{\widetilde{\iota}} \in C\left(\widetilde{M}/\widetilde{   \T}^2 \right)$ is regarded as element of $C\left(\widetilde{M}\right)$. From $\widetilde{e}_{\widetilde{\iota}} \in \Coo\left(\widetilde{M}\right)_{0,0}$ it turns out $\widetilde{l}\left(\widetilde{e}_{\widetilde{\iota}} \right)  = \widetilde{e}_{\widetilde{\iota}}$,  $~~\left\langle \widetilde{l}\left( \widetilde{e}_{\widetilde{\iota}'}\right) , \widetilde{l}\left( \widetilde{e}_{\widetilde{\iota}''}\right)  \right\rangle_{C\left( \widetilde{M}_{ \widetilde{\th}}\right) } = \left\langle  \widetilde{e}_{\widetilde{\iota}'} ,  \widetilde{e}_{\widetilde{\iota}''}  \right\rangle_{C\left( \widetilde{M}\right) } $.

From \eqref{isospectral_parti_eqn}-\eqref{isospectral_f_eqn} it follows that
\be\label{isospectral_undec_eqn}
1_{C\left(\widetilde{M} \right) }= \frac{1}{N_1N_2}\sum_{\widetilde{\iota} \in \widetilde{I}} \sum_{\substack{j = 0\\ k = 0}}^{\substack{j = N_1-1\\ k = N_2-1}} \widetilde{f}_{\left(\widetilde{\iota}, j, k \right)}\left\rangle \right\langle \widetilde{f}_{\left(\widetilde{\iota}, j, k \right)},
\ee
\be\nonumber
\left\langle \widetilde{f}_{\left(\widetilde{\iota}', j', k' \right)}, \widetilde{f}_{\left(\widetilde{\iota}'', j'', k'' \right)} \right\rangle_{C\left( \widetilde{M}\right) } = \frac{1}{N_1N_2} \delta_{j'j''}\delta_{k'k''}\left\langle \widetilde{e}_{\widetilde{\iota}'}, \widetilde{e}_{\widetilde{\iota}''} \right\rangle_{C\left( \widetilde{M}\right) }\in \Coo\left(M \right),
\ee
\be\label{isospectral_fprodl_eqn}
\left\langle \widetilde{l}\left( \widetilde{f}_{\left(\widetilde{\iota}', j', k' \right)}\right) , \widetilde{l}\left( \widetilde{f}_{\left(\widetilde{\iota}'', j'', k'' \right)}\right)  \right\rangle_{C\left( \widetilde{M}_{\widetilde{   \th}}\right) } =\frac{1}{N_1N_2} \delta_{j'j''}\delta_{k'k''}\left\langle \widetilde{e}_{\widetilde{\iota}'}, \widetilde{e}_{\widetilde{\iota}''} \right\rangle_{C\left( \widetilde{M}_{\widetilde{   \th}}\right) }\in  \Coo\left(M_\th \right).
\ee
From the \eqref{isospectral_undec_eqn} it turns out that $C\left(\widetilde{M} \right)$  is a right $C\left(M \right)$ module generated by  finite set of elements $\widetilde{f}_{\left(\widetilde{\iota}, j, k \right)}$ where $\left(\widetilde{\iota}, j, k \right) \in \widetilde{J}$, i.e. any $\widetilde{a} \in C\left(\widetilde{M} \right)$ can be represented as
\be\label{isospectral_wa_eqn}
\widetilde{a} = \sum_{\widetilde{\iota} \in \widetilde{I}} \sum_{\substack{j = 0\\ k = 0}}^{\substack{j = N_1-1\\ k = N_2-1}} \widetilde{f}_{\left(\widetilde{\iota}, j, k \right)} a_{\left(\widetilde{\iota}, j, k \right)}; \text{ where } a_{\left(\widetilde{\iota}, j, k \right)} \in C\left( M\right). 
\ee
Moreover if $\widetilde{a} \in \Coo\left(\widetilde{M} \right)$ then one can select $a_{\left(\widetilde{\iota}, j, k \right)} \in \Coo\left(M \right)$. However any $a_{\left(\widetilde{\iota}, j, k \right)} \in \Coo\left(M \right)$ can be uniquely
written as a doubly infinite
norm convergent sum of homogeneous elements,
\begin{equation*}
a_{\left(\widetilde{\iota}, j, k \right)} = \sum_{n_1,n_2} \, \widehat{T}_{n_1,n_2} \, ,
\end{equation*}

with $\widehat{T}_{n_1,n_2}$ of bidegree $(n_1,n_2)$ and where the sequence
of norms $||
\widehat{T}_{n_1,n_2} ||$ is of
rapid decay in $(n_1,n_2)$. One has
\bea\label{isospectral_sec_eqn}
\widetilde{l}\left(\widetilde{f}_{\left(\widetilde{\iota}, j, k \right)} a_{\left(\widetilde{\iota}, j, k \right)}\right) = \sum_{n_1, n_2}  \widetilde{f}_{\left(\widetilde{\iota}, j, k \right)} \widehat{T}_{n_1,n_2} \widetilde{\lambda}^{\left( N_2n_2+j \right) \widetilde{p}_1}= \sum_{n_1, n_2} \widetilde{f}_{\left(\widetilde{\iota}, j, k \right)} \widetilde{\lambda}^{j \widetilde{p}_1}  \widetilde{\lambda}^{kN_1n_1}\widehat{T}_{n_1,n_2} \widetilde{\lambda}^{N_2n_2 \widetilde{p}_1} 
\eea

the sequence
of norms $||
\widetilde{\lambda}^{kN_2n_2}\widehat{T}_{n_1,n_2} ||=||
\widehat{T}_{n_1,n_2} ||$ is of
rapid decay in $(n_1,n_2)$ it follows that
\begin{equation}\label{isospectral_a_eqn}
a'_{\left(\widetilde{\iota}, j, k \right)} = \sum_{n_1,n_2} \, \widetilde{\lambda}^{kN_1n_1} \widehat{T}_{n_1,n_2} \in \Coo\left(M \right) 
\end{equation}
From \eqref{isospectral_wa_eqn} - \eqref{isospectral_a_eqn} it turns out
\be\label{isospectral_sum_eqn}
\widetilde{l}\left( \widetilde{a}\right)  = \sum_{\widetilde{\iota} \in \widetilde{I}} \sum_{\substack{j = 0\\ k = 0}}^{\substack{j = N_1-1\\ k = N_2-1}} \widetilde{l}\left( \widetilde{f}_{\left(\widetilde{\iota}, j, k \right)} \right) l\left( a'_{\left(\widetilde{\iota}, j, k \right)}\right) ; \text{ where } l\left( a'_{\left(\widetilde{\iota}, j, k \right)}\right) \in \Coo\left( M_\th \right).
\ee
However $C\left(\widetilde{M}_{\widetilde{\th}} \right)$ is the norm  completion of $\Coo\left(\widetilde{M}_\th \right)$, so from \eqref{isospectral_sum_eqn} it turns out that  $C\left(\widetilde{M}_{\widetilde{\th}} \right)$ is a right Hilbert $C\left({M}_{{\th}} \right)$-module generated by a finite set
\be\label{isospectral_xi_eqn}
\Xi=\left\{\widetilde{l}\left(\widetilde{f}_{\left(\widetilde{\iota}, j, k \right)} \right) \in   \Coo\left(\widetilde{M}_{\widetilde{\th}} \right)   \right\}_{\left(\widetilde{\iota}, j, k \right) \in \widetilde{J}}
\ee  From the Corollary \ref{fin_hpro_cor} it follows that the module is projective.
So one has the following theorem.
\begin{thm}\label{isospectral_fin_thm}
	The triple $\left( C\left( M_\th\right), C\left( \widetilde{M}_{ \widetilde{\th}}\right), G\left(\widetilde{M}~|~ M \right)\right)   $ is an unital noncommutative finite-fold  covering.
\end{thm}

\begin{exercise}
	Prove the following statement
	If $\widetilde{\rho}:C\left( \widetilde{M}_{ \widetilde{\th}}\right)\to B\left(\widetilde{\H} \right) $ is induced by $\left(\rho, \left( C\left( M_\th\right), C\left( \widetilde{M}_{ \widetilde{\th}}\right), G\left(\widetilde{M}, M \right)\right)\right)$ then $\widetilde{\rho}$ can be represented by action of $C\left( \widetilde{M}_{ \widetilde{\th}}\right)$ on $ L^2\left( \widetilde{M},\widetilde{S}\right)$ by operators \eqref{l_defn}.
		The noncommutative spectral triple   $$\left( \Coo\left(\widetilde{M}_{\widetilde{\th}}\right) , L^2\left(\widetilde{M},\widetilde{S} \right), \widetilde{ \slashed D}  \right)$$ is a $\left( C\left( M_\th\right), C\left( \widetilde{M}_{ \widetilde{\th}}\right), G\left(\widetilde{M}~|~ M \right)\right)   $-lift of $\left( \Coo\left(M_\th\right) , L^2\left(M,S \right), \slashed D  \right)$. 

\end{exercise}
\subsection{Coverings of spectral triples}
\paragraph*{}
If ${f}_{\left({\iota}, j, k \right)} \in \Coo \left(M \right)$ is given by \eqref{isospectral_fs_eqn} then from then from \eqref{isospectral_f's_eqn} it turns out 
\be\nonumber
{f}'_{\left({\iota}, j, k \right)} = {e}_{{\iota}} {u}^j{
	v}^k \in C_0\left({\mathcal{U}}_{\iota} \right)= C_0\left({\mathcal{V}}_{\iota} \times \T^2 \right)
\ee
In the above formula the product ${u}^j{
	v}^k$ can be regarded as element of both $C\left( \T^2\right)$ and  $C_b\left({\mathcal{V}}_{\iota} \times \T^2 \right) =  C_b\left({\mathcal{W}}_{\iota} \right)$, where $\mathcal{W}_{\iota}\subset M$ is the homeomorphic image of ${\mathcal{V}}_{\iota} \times \T^2$. Since the Dirac operator $\slashed D$ is invariant with respect to transformations 
$
u \mapsto e^{i\varphi_u}u, ~~ v \mapsto e^{i\varphi_v}v
$
one has
\be\label{isospectral_jdb_eqn}
\left[\slashed D, u \right]= d^{\iota}_uu,~~ \left[\slashed D, v \right]= d^{\iota}_vv
\ee
where $d^\iota_u, d^\iota_v: \mathcal V_\iota \to \mathbb{M}_{\dim S} \left(\C \right)$ are continuous matrix-valued functions. I would like to avoid functions in $C_b\left({\mathcal{U}}_{\iota} \right)$, so instead \eqref{isospectral_jdb_eqn} the following evident consequence of it will be used
\be\label{isospectral_jde_eqn}
\left[\slashed D,a u^jv^k \right]= \left[\slashed D,a \right] u^jv^k  + a\left( jd^\iota_u+kd^\iota_v\right) u^jv^k; ~~ a \in C_0\left( M/ \T^2\right), ~~ \supp a \subset \mathcal{V}_{\iota}.
\ee
In contrary to \eqref{isospectral_jdb_eqn} the equation \eqref{isospectral_jde_eqn} does not operate with $C_b\left({\mathcal{U}}_{\iota} \right)$, it operates with $C_0\left({\mathcal{U}}_{\iota} \right)$. Let $\pi: \widetilde M \to M$ and let $\widetilde{\slashed D} = p^{-1} \slashed D$ be the $\pi$-inverses image of $\slashed D$ (cf. Definition  \ref{top_smooth_inv_im_defn}). Suppose $\widetilde{  \mathcal V}_{\widetilde{\iota}} \subset \widetilde{M}/\widetilde{   \T}^2$ is mapped onto ${  \mathcal V}_{{\iota}} \subset {M}/{   \T}^2$. Then we set $d^{\widetilde{\iota}}_u \stackrel{\mathrm{def}}{=} d^{{\iota}}_u$, $d^{\widetilde{\iota}}_v \stackrel{\mathrm{def}}{=} d^{{\iota}}_v$. The covering $\pi$ maps $\widetilde{\mathcal{V}}_{ \widetilde\iota} \times \widetilde{\T}^2$ onto ${\mathcal{V}}_{\iota} \times \T^2$. If $\widetilde{u}, \widetilde{v} \in C\left(\widetilde{\T}^2 \right)$ are natural generators, then the covering $\widetilde{\T}^2 \to \T^2$ is given by
\bean
C\left(\T^2 \right) \to C\left(\widetilde{\T}^2 \right),\\
u \mapsto \widetilde{u}^{N_1},~~v \mapsto \widetilde{v}^{N_2}.
\eean 
From the above equation and taking into account \eqref{isospectral_jdb_eqn} one has
\be\label{isospectral_du_eqn}
\begin{split}
\left[\widetilde{\slashed D}, \widetilde{u} \right]= \frac{d^{\widetilde{\iota}}_u}{N_1}\widetilde{u},~~ \left[\widetilde{\slashed D}, \widetilde{v} \right]= \frac{d^{\widetilde{\iota}}_v}{N_2}\widetilde{v},\\
\left[\widetilde{\slashed D}, \widetilde{a} \widetilde{u}^j\widetilde{v}^k \right]=\left[\widetilde{\slashed D},\widetilde{a} \right] \widetilde{u}^j\widetilde{v}^k  + \widetilde{a} \left(\frac{j^{\widetilde{\iota}}_u}{N_1} +\frac{k^{\widetilde{\iota}}_v}{N_2}\right) \widetilde{u}^j\widetilde{v}^k;\\ ~~ \widetilde{a} \in C_0\left( \widetilde{M}/ \widetilde{\T}^2\right), ~~ \supp \widetilde{a} \subset \widetilde{\mathcal{V}}_{\widetilde{\iota}}.
\end{split}
\ee
For any $\widetilde{a}\in \Coo\left( \widetilde{M}\right) $ such that $\supp \widetilde{a} \subset  \widetilde{\mathcal{W}}_{\widetilde{\iota}}$ following condition holds
\be\label{isospectral_da_eqn}
\left[\widetilde{\slashed D}, \widetilde{a} \right]= \mathfrak{lift}_{\widetilde{\mathcal{W}}_{\widetilde{\iota}}}\left(\left[{\slashed D}, \mathfrak{desc}\left( \widetilde{a}\right)  \right] \right) 
\ee
If $b \in \Coo\left(M/\T^2 \right)\subset \Coo\left(M \right)$ then from \eqref{isospectral_du_eqn} and \eqref{isospectral_da_eqn} it turns out
\be\nonumber
\begin{split}
\left[\widetilde{\slashed D},  \widetilde{f}_{\left(\widetilde{\iota}, j, k \right)}b u^{j'}v^{k'}\right]   = \sqrt{\widetilde{e}_{\widetilde{\iota}}}\widetilde{u}^j\widetilde{v}^k  u^{j'}v^{k'}\otimes \left[ \slashed D, \sqrt{{e}_{\widetilde{\iota}}} b  \right]+ \\ +  b\sqrt{\widetilde{e}_{\widetilde{\iota}}}\widetilde{u}^j\widetilde{v}^k  u^{j'}v^{k'}\otimes \left[ \slashed D, \sqrt{{e}_{\widetilde{\iota}}}\right]+\\+b
	\sqrt{\widetilde{e}_{\widetilde{\iota}}}\widetilde{u}^j\widetilde{v}^k  u^{j'}v^{k'}\otimes  \sqrt{{e}_{\widetilde{\iota}}}\left(j'd_u+\frac{ jd_u}{N_1} +k'd_v + \frac{kd_v}{N_2}\right)\\
\end{split}
\ee 
and taking into account $\left[ \widetilde{\slashed D}, \widetilde{l}\right]=0$ one has
\be\label{isospectral_predu_eqn}
\begin{split}
	\left[\widetilde{\slashed D}, \widetilde{l}\left(  \widetilde{f}_{\left(\widetilde{\iota}, j, k \right)}b u^{j'}v^{k'}\right) \right]   = \sqrt{\widetilde{e}_{\widetilde{\iota}}}\widetilde{l}\left(\widetilde{u}^j\widetilde{v}^k  u^{j'}v^{k'}\right)\otimes \left[ \slashed D, \sqrt{{e}_{\widetilde{\iota}}} b  \right]+ \\ +  b\sqrt{\widetilde{e}_{\widetilde{\iota}}}\widetilde{l}\left(\widetilde{u}^j\widetilde{v}^k  u^{j'}v^{k'}\right)\otimes \left[ \slashed D, \sqrt{{e}_{\widetilde{\iota}}}\right]+\\+b
	\sqrt{\widetilde{e}_{\widetilde{\iota}}}\widetilde{l}\left(\widetilde{u}^j\widetilde{v}^k  u^{j'}v^{k'}\right)\otimes  \sqrt{{e}_{\widetilde{\iota}}}\left(j'd_u+\frac{ jd_u}{N_1} +k'd_v + \frac{kd_v}{N_2}\right)\\
\end{split}
\ee
Taking into account that  
$
\widetilde{l}\left(\widetilde{u}^j\widetilde{v}^k \right)\widetilde{l}\left( u^{j'}v^{k'}\right)= \widetilde{\la}^{j'N_2k} \widetilde{l}\left(\widetilde{u}^j\widetilde{v}^k  u^{j'}v^{k'}\right)
$ the equation \eqref{isospectral_predu_eqn} is equivalent to

\be\label{isospectral_pred_eqn}
\begin{split}
	\left[\widetilde{\slashed D},  \widetilde{l}\left( \widetilde{f}_{\left(\widetilde{\iota}, j, k \right)}\right) b \widetilde{l}\left( u^{j'}v^{k'}\right) \right]   = \sqrt{\widetilde{e}_{\widetilde{\iota}}}\widetilde{l}\left(\widetilde{u}^j\widetilde{v}^k \right)\widetilde{l}\left( u^{j'}v^{k'}\right)\otimes \left[ \slashed D, \sqrt{{e}_{\widetilde{\iota}}} b  \right]+ \\ +  b\sqrt{\widetilde{e}_{\widetilde{\iota}}}\widetilde{l}\left(\widetilde{u}^j\widetilde{v}^k \right) \widetilde{l}\left( u^{j'}v^{k'}\right)\otimes \left[ \slashed D, \sqrt{{e}_{\widetilde{\iota}}}\right]+\\+b
	\sqrt{\widetilde{e}_{\widetilde{\iota}}}\widetilde{l}\left(\widetilde{u}^j\widetilde{v}^k \right) \widetilde{l}\left( u^{j'}v^{k'}\right)\otimes  \sqrt{{e}_{\widetilde{\iota}}}\left(j'd_u+\frac{ jd_u}{N_1} +k'd_v + \frac{kd_v}{N_2}\right)\\
\end{split}
\ee
\begin{exercise}
		Prove that a noncommutative spectral triple   $$\left( \Coo\left(\widetilde{M}_{\widetilde{\th}}\right) , L^2\left(\widetilde{M},\widetilde{S} \right), \widetilde{ \slashed D}  \right)$$ is a $\left( C\left( M_\th\right), C\left( \widetilde{M}_{ \widetilde{\th}}\right), G\left(\widetilde{M}~|~ M \right)\right)   $-lift of $\left( \Coo\left(M_\th\right) , L^2\left(M,S \right), \slashed D  \right)$  (cf. Definition \ref{spectral_triple_fin_lift_defn}). 
\end{exercise}

\subsection{Unoriented  twisted spectral triples}

\paragraph*{}
Suppose that $M$ is unoreintable manifold which satisfies to \eqref{isos_t_act_eqn}, i.e.
\begin{equation*}
\mathbb{T}^2 \subset \mathrm{Isom}(M) \, ,
\end{equation*}
Suppose that the natural 2-fold covering $p:\widetilde{M}\to M$ is such that $\widetilde{M}$ is a Riemannian manifold which  admits a Spin$^c$ structure (cf. Definition \ref{spin_str_defn}), so there is an oriented spectral triple
$\left( \Coo\left(\widetilde{M}\right) , L^2\left(\widetilde{M},\widetilde{\sS} \right), \widetilde{ \slashed D}  \right)$. Assume that there is a bundle $\sS \to M$ on $M$ with sesquilinear form (cf. Definition \ref{top_herm_bundle_form_defn}) such that $\widetilde{\sS}$ is the $p$-inverse image of $\SS$ (cf. Definition \ref{vb_inv_img_funct_defn}).
From \ref{comm_sp_tr} it turns out that there is an unoriented  spectral triple given by \eqref{top_unori_eqn}, i.e.

\be\nonumber
\left(\Coo\left(M \right), L^2\left(\widetilde{M}, \widetilde{\SS}\right)^{\Z_2}, \slashed D  \right).
\ee
Otherwise from \eqref{isos_twisetd_eqn} it follows that there is an oriented twisted spectral triple
\be\nonumber
\left( l\left( \Coo\left(\widetilde{M}\right)\right)  , L^2\left(\widetilde{M},\widetilde{\sS} \right), \widetilde{ \slashed D}  \right).
\ee.

Action of $G\left( \left.\widetilde{M}~\right|M\right) \cong \Z_2$ on $\widetilde{M}$ induces an action of $\Z_2$ on both $\Coo\left(\widetilde{M}\right) $ and $l\Coo\left(\widetilde{M}\right) $ such that
\be\nonumber
\begin{split}
	\Coo\left(M \right)= \Coo\left(\widetilde{M}\right)^{\Z_2},\\
	l\Coo\left(M \right)= l\Coo\left(\widetilde{M}\right)^{\Z_2},\\
\end{split}
\ee
From the above construction we have an unoriented twisted spectral triple
$$
\left(l\Coo\left(M \right), L^2\left(\widetilde{M}, \widetilde{\SS}\right)^{\Z_2}, \slashed D  \right).
$$

which satisfies to the Definition \ref{unoriented_defn}  and \ref{unoriented_empt}.

\chapter{Coverings of quantum groups}\label{qdr_chap}

\section{The double covering of the quantum group $SO_q(3)$}
\paragraph*{} 

The covering of the quantum $SO\left(3 \right)$ is described in  \cite{dijkhuizen:so_doublecov,podles:so_su}. Here we would like to prove that the covering complies with the general theory of noncommutative coverings described in \ref{cov_fin_bas_sec}. Moreover we prove that the covering gives an unoriented unoriented spectral triple (cf. Definition \ref{unoriented_defn}  and \ref{unoriented_empt}).

\subsection{Basic constructions}
\paragraph*{}
Denote by 
\be\label{su_q_2_G_eqn}
	G = \left\{\left. g \in \Aut\left(C\left( SU_q\left(2 \right)\right)  \right)~\right|~ ga = a;~~\forall a \in C\left(SO_q\left(3 \right)\right) \right\}
\ee
Denote by
\be\label{su_q_2_gr_eqn}
C\left(SU_q\left(2 \right)\right)_N \bydef \left\{ \left.\widetilde{a} \in C\left(SU_q\left(2 \right)\right)~\right|~ \widetilde{a}\bt\bt^* = q^{2N}\bt\bt^*\widetilde{a} \right\}; \quad \forall N \in \Z
\ee
and let us prove that the equation \eqref{su_q_2_gr_eqn} yields the $\Z$-grading of $C\left(SU_q\left(2 \right)\right)$. Really following conditions hold:
\begin{itemize}
	\item 
	\be\label{su_q_2_N_eqn}
	\begin{split} 
		\widetilde{a} \in C\left(SU_q\left(2 \right)\right)_N \Leftrightarrow	 \widetilde{a} =
		\left\{
		\begin{array}{c l}
			\al^N	\sum_{\substack{ j = 0\\ k = 0}}^\infty  b_{jk}\bt^j\bt^{*k} & N \ge 0 \\
			& \\
			\left( \al^{*}\right)^{-N}  	\sum_{\substack{ j = 0\\ k = 0}}^\infty  b_{jk}\bt^j\bt^{*k} & N < 0
		\end{array}\right.
	\end{split}	
	\ee
	\item From the Theorem \ref{su_q_2_bas_thm} it follows that the linear span of given by the equation \eqref{su_q_2_fin_eqn} elements is dense in $C\left( SU_q(2)\right)$, hence $C\left( SU_q(2)\right)$ is the $C^*$-norm completion of the following direct sum
	\be\label{su_q_2_nsum_eqn}
	\bigoplus_{N \in \Z}  C\left(SU_q\left(2 \right)\right)_N.
	\ee
\end{itemize}
From  $\bt\bt^*\in C\left(SO_q\left(3 \right)\right)$ it follows that the grading is $G$-invariant, i.e.
\be\label{su_q_2_ginv_eqn}
\begin{split}
	g C\left(SU_q\left(2 \right)\right)_N = C\left(SU_q\left(2 \right)\right)_N; \quad \forall g \in G, \quad \forall N \in \N.
\end{split}	
\ee	
$C\left(SU_q\left(2 \right)\right)_0$ is a commutative $C^*$-algebra generated generated by $\bt$ and $\bt^*$ i.e.
$$
\widetilde{a} \in C\left(SU_q\left(2 \right)\right)_0 \Rightarrow \widetilde{a} = \sum_{\substack{j = 0\\ k = 0}}^\infty c_{jk}\bt^j\bt^{*k}; \quad c_{jk}\in \C.
$$
It follows that $g\bt \in  C\left(SU_q\left(2 \right)\right)_0$, i.e.  
$
g\bt = \sum_{\substack{j = 0\\ k = 0}}^\infty c_{jk}\bt^j\bt^{*k},
$
and taking into account $\bt^2\in C\left( SO_q(3)\right)\Rightarrow\left(g\bt \right)^2 = \bt^2$ one concludes
\be\label{su_q_2_bpm_eqn}
g\bt = \pm \bt.
\ee
From $\al \in C\left(SU_q\left(2 \right)\right)_1$ and the equations \eqref{su_q_2_N_eqn}, \eqref{su_q_2_ginv_eqn}  it turns out that for any $g \in G$ one has
\be\label{su_q_2_gal_eqn}
g \al\in  C\left(SU_q\left(2 \right)\right)_1 \Rightarrow g\al=\al\sum_{\substack{ j = 0\\ k = 0}}^\infty b_{jk}\bt^j\bt^{*k};\quad b_{jk}\in \C
\ee
Otherwise from $\al^2 \in SO_q\left(3 \right)$ it turns out 
\be\label{su_q_2_al_eqn}
\begin{split}
\al^2 = g\al^2 = (g\al)^2 = \left(\al\sum_{\substack{ j = 0\\ k = 0}}^\infty b_{jk}\bt^j\bt^{*k} \right)^2 =
\\
= \al^2\sum_{\substack{j = 0\\ k = 0}}^\infty\sum_{\substack{l = 0\\m = 0}}^\infty q^{k+m}b_{jk}b_{lm} \bt^{j+l}\left( \bt^*\right)^{k+m}= \al^2\sum_{\substack{ j = 0\\ k = 0}}^\infty c_{rs}\bt^r\bt^{*s}; \\
\text{where} \quad c_{rs} = \sum_{l = 0}^r \sum_{m = 0}^s q^{k+m}b_{r-l, s-m}b_{lm}
\end{split}
\ee
It turns out that $c_{0,0}= 1$ and $c_{rs}=0$ for each $\left(r,s \right) \neq \left( 0,0\right)$. Otherwise $c_{0,0}=b_{0,0}^2$ it turns out that $b_{0,0}= \epsilon = \pm1$. Suppose that there are $j, k \in \N^0$ such that $b_{jk}\neq 0$ and $j + k > 0$. If $j$ and $k$ are such that 
\bean
b_{jk}\neq 0,\\
j+k =\min_{\substack{b_{mn}\neq 0\\ m+n>0}} m+n
\eean
then $c_{jk} = \epsilon b_{jk}\left(1+q^{j+k} \right) \neq 0$. There is a contradiction with $c_{rs}=0$ for each $r+s > 0$. It follows that if $j+k > 0$ then $b_{jk}= 0$, hence one has
$$
g \al = \eps\al = \pm\al.
$$
In result we have 
\bean
g \al = \pm \al,
g \bt = \pm \bt
\eean
If $g\al = \al$ and $g \bt = -\bt$ then $g\left(\al\bt \right) = - \al\bt$, it is impossible because $\al\bt \in SO_q\left(3 \right)$. It turns out that $G=\Z_2$ and if $g \in G$ is not trivial then $g\al = -\al$ and $g\bt = -\bt$.
So we proved the following lemma
 \begin{lem}\label{su_q_2_group_lem}
 If $G = \left\{ g \in \Aut\left(C\left( SU_q\left(2 \right)\right)  \right)~|~ ga = a;~~\forall a \in C\left(SO_q\left(3 \right)\right) \right\}$ then $G \approx \Z_2$. Moreover if $g \in G$ is the nontrivial element then
 \bean
g\left( \al^k\bt^n\bt^{*m}\right) = \left( -1\right)^{k+m+n}  \al^k\bt^n\bt^{*m}, \\ g\left( \al^{*k}\bt^n\bt^{*m}\right) = \left( -1\right)^{k+m+n}  \al^{*k}\bt^n\bt^{*m}.
 \eean 

 \end{lem}
\subsection{Covering of $C^*$-algebra}

\begin{lem}\label{su_q_2_fin_lem}
	$C\left( SU_q\left(2\right)\right)$ is a finitely generated projective $C\left( 	SO_q\left(3 \right) \right)$ module.
\end{lem}
\begin{proof}
	Let $A_f$ be given by the Theorem \ref{su_q_2_bas_thm}. 
	If $A_f^{\Z_2}= A_f\bigcap C\left( 	SO_q\left(3 \right) \right) $ then from \eqref{su_q_2_z2_ncomm_eqn} and the Theorem \ref{su_q_2_bas_thm} it turns out that given by \eqref{su_q_2_fin_eqn} elements
	\be\nonumber
	\al^k\bt^n\bt^{*m}~~ \text{ and }~~ \al^{*k'}\bt^n\bt^{*m}
	\ee 
	with even $k+m+n$ or $k'+m+n$ is the basis of $A_f^{\Z_2}$.	If 
	$$
	\widetilde{a} = \al^k\bt^n\bt^{*m} \notin A_f^{\Z_2}
	$$
	then $k+m+n$ is odd. If $m > 0$ then
	$$
	\widetilde{a} = \al^k\bt^n\bt^{*m-1} \bt^* = a \bt^* \text{ where } a \in A_f^{\Z_2}.
	$$
	If $m = 0$ and $n > 0$ then
	$$
	\widetilde{a} = \al^k\bt^{n-1} \bt = a \bt \text{ where } a \in A_f^{\Z_2}
	$$
	If $m = 0$ and $n= 0$ then $k > 0$ and
	$$
	\widetilde{a} = \al^{k-1}\al = a \al \text{ where } a \in A_f^{\Z_2}.
	$$
	From 
	$$
	\widetilde{a} = \al^{*k'}\bt^n\bt^{*m} \notin A_f^{\Z_2}
	$$
	it follows that  $k'+m+n$ is odd. Similarly to the above proof one has
	$$
	\widetilde{a} = a \al\text{ or } \widetilde{a} = a \al^* \text{ or } \widetilde{a} = a \bt \text{ or } \widetilde{a} = a \bt^*   \text{ where } a \in A_f^{\Z_2}.
	$$
	From the above equations it turns out that $A_f$ is a left  $A_f^{\Z_2}$-module generated by $\al, \al^*, \bt, \bt^*$. Algebra $A_f^{\Z_2}$ (resp. $A_f$) is dense in $ C\left( 	SO_q\left(3 \right) \right) $ (resp.  $C\left( 	SU_q\left(2 \right) \right)$ ) it follows that  $C\left( 	SU_q\left(2 \right) \right)$ is a left  $C\left( 	SO_q\left(3 \right) \right)$-module generated by $\al, \al^*, \bt, \bt^*$. From the Corollary \ref{fin_hpro_cor} it turns out that $C\left( SU_q\left(2\right)\right)$ is a finitely generated projective $C\left( 	SO_q\left(3 \right) \right)$ module.
	
\end{proof}
\begin{corollary}
	The	triple $\left(C\left( 	SO_q\left(3 \right) \right), C\left( 	SU_q\left(2 \right) \right), \Z_2 \right)$ is an unital noncommutative finite-fold  covering.
\end{corollary}
\begin{proof}
	Follows from $C\left( 	SO_q\left(3 \right) \right)=C\left( SU_q\left(2 \right) \right)^{\Z_2}$ and Lemmas \ref{su_q_2_group_lem}, \ref{su_q_2_fin_lem}.
\end{proof}

\begin{corollary}
	The	triple $\left(C\left( 	SO_q\left(3 \right) \right), C\left( 	SU_q\left(2 \right) \right), \Z_2\times \Z_2, \pi\right)$ is an unital noncommutative finite-fold  covering.
\end{corollary}

\begin{problem}
	It is known that there is an unitary element $u \in C\left(  SU_q\left(2 \right)\right)$ (cf. \cite{chakraborty_pal:quantum_su_2}) such that $\left[u\right]\in \K_1\left(C\left( \left( SU_q\left(2 \right)\right) \right)\right) $ is not trivial and has infinite period.
\\ Question. Does an $\left(u, n\right)$-covering   $\left(C\left(  SU_q\left(2 \right)\right)\otimes\K, \widetilde A, \Z_{n}\right)$ (cf. Definition \ref{hurewicz_u_n_defn}) exist?  
\end{problem}

\subsection{$SO_q\left( 3\right)$ as an unoriented  spectral triple }
Let  $h: C\left( SU_q\left(2\right)\right)$ be a given by the equation \eqref{su_q_2_haar_eqn} state, and $L^2\left( C\left( SU_q\left(2\right)\right), h\right)$ is the Hilbert space of the corresponding GNS representation  (cf. Section \ref{gns_constr_sec}). The state $h$ is  $\Z_2$-invariant, hence there is an action of $\Z_2$ on $L^2\left( C\left( SU_q\left(2\right)\right), h\right)$ which is naturally induced by the action  $\Z_2$ on $\Coo\left( SU_q\left(2\right)\right) $. From the above construction it follows that the unital orientable spectral triple
\be\nonumber
\left(\Coo\left( SU_q\left(2\right)\right), L^2\left( C\left( SU_q\left(2\right)\right), h\right), \widetilde{D} \right).
\ee
\begin{exercise}
Using an unital noncommutative finite-fold covering  $$\left(C\left( 	SO_q\left(3 \right) \right), C\left( 	SU_q\left(2 \right) \right), \Z_2 \right)$$
(cf. Definition \ref{fin_unital_defn})
rove that there is is an unoriented  spectral triple (cf. Definition \ref{unoriented_defn}) given by
$$
\left(\Coo\left( SO_q\left(3\right)\right), L^2\left( C\left( SU_q\left(2\right)\right), h\right)^{\Z_2},D \right) 
$$ where $\Coo\left( SO_q\left(3\right)\right)\bydef C\left( SO_q\left(3\right)\right)\bigcap \Coo\left( SU_q\left(2\right)\right)$ and $D\bydef \widetilde{ D}|_{L^2\left( C\left( SU_q\left(2\right)\right), h\right)^{\Z_2}}$.

\end{exercise}
\section{Covering of quantum group of euclidean transformations}

\paragraph*{}
Let us fix a real number $\mu$ such that $\left|\mu \right| > 1$, and  let $A \bydef A_\mu$  be the $*$-algebra generated by two elements $u$ and $n$ which satisfy to the equations \eqref{em_eqn}. Let  $k\in \N$ is such that $k > 1$ and $\widetilde\mu \bydef \mu^{1/k}$ then similarly to \eqref{em_eqn} there is a *-algebra $\widetilde A\bydef A_{\widetilde \mu}$  generated by two elements $\widetilde u$ and $\widetilde n$ such that
\bean
\begin{split}
	\widetilde v^*\widetilde v = \widetilde v\widetilde v^* = 1_{A_{\widetilde\mu}},\\
\widetilde	n^*\widetilde n = \widetilde  n\widetilde n^*,\\
	\widetilde v^*\widetilde n \widetilde v=\widetilde \mu\widetilde n.
\end{split}
\eean
There is a *- homomorphism $\pi:A\hookto \widetilde A$ such that 
\bean
\pi\left(n \right) \bydef \widetilde n,\\
\pi\left(v \right)\bydef \widetilde v^k. 
\eean 
There is the action $\Z_k \times A_{\widetilde\mu}\to A_{\widetilde\mu}$ such that for all $\overline l \in \Z_k$ one has
\bean
\overline l \widetilde n = \overline n;\\
\overline l \widetilde v = e^{\frac{2\pi i l}{k}}\widetilde v.
\eean 
where $l \in \Z$ is a representative of $\overline l$.
\begin{exercise}

Prove that a quadruple $\left(A, \widetilde A, \Z_k, \pi \right)$ is a {noncommutative finite-fold  pre-covering of *-algebras} (cf. Definition \ref{fin_pre*_defn}).
\end{exercise}
A reduced crossed product $\widetilde B\bydef C_0\left(\C_{\left( \widetilde\mu\right)}\right) \rtimes_r \Z$ (cf. Section \ref{discr_cr_prod_sec}) is a completion of $C_c\left( \Z, C_c\left(\C_{\left( \widetilde\mu\right)} \right) \right)$ with respect to the reduced crossed norm (cf. Equation \eqref{discr_red_n_eqn}). There is a *-automorphism 
\be\label{qg_a_eqn}
\begin{split}
\a \in \Aut\left( C_c\left( \Z, C_c\left(\C_{\left( \widetilde\mu\right)} \right)\right) \right) ,\\
\sum_{t \in \Z} a_t u_t\mapsto \sum_{t \in \Z}e^{\frac{2\pi  it}{k}} a_t u_t\quad \forall t \in G \quad a_t  \in A.
\end{split}
\ee
Since $\a^k$ is the identical isomorphism it yields an action $\Z_k \times C_c\left( \Z, C_c\left(\C_{\left( \widetilde\mu\right)} \right)\right)\to C_c\left( \Z, C_c\left(\C_{\left( \widetilde\mu\right)} \right)\right)$, which can be uniquely extended to the action
$$
\Z_k \times \widetilde B \to \widetilde B.
$$
The  given by \eqref{dm_act_eqn} action 
$$ A_{\widetilde \mu}\times  C_c\left( \Z, C_c\left(\C_{\left(  \widetilde\mu\right) }\right)  \right)\to  C_c\left( \Z, C_c\left(\C_{\left(  \widetilde\mu\right) }\right)  \right)
 $$
 can be extended up to 
 $$
   A_{\widetilde \mu}\times  C_c\left( \Z, C_c\left(\C_{\left(  \widetilde\mu\right) }\right)  \right)\widetilde B\to  C_c\left( \Z, C_c\left(\C_{\left(  \widetilde\mu\right) }\right)  \right) \widetilde B.
 $$
If $\widetilde{\mathfrak{B}}$ is a completion of a dense right $\widetilde B$-ideal $C_c\left(\C_{\left(  \widetilde\mu\right) }\right)   \widetilde B$ with respect to the graph norms (cf. Definition \ref{graph_norm_eqn}) then there is a representation
$$
 \widetilde A = A_{\widetilde \mu}\hookto \End^*_{\widetilde B}\left( \widetilde{\mathfrak{B}}\right).
$$
Denote by 
\bean
B \bydef \widetilde B^{\Z_k} \bydef \left\{\left.  \widetilde b \in  \widetilde B \right| \forall g \in \Z_k\quad g  \widetilde b = \widetilde b\right\},\\
\mathfrak{B} \bydef \widetilde{\mathfrak{B}}\cap B\cong  \widetilde{\mathfrak{B}}^{\Z_k}.
\eean
\begin{lemma}
In the above situation the triple $\left( B, \widetilde B, \Z_k\right)$ is a noncommutative finite-fold covering with unitization (cf. Definition \ref{fin_unitization_defn}).
\end{lemma}
 \begin{proof}
 If $C_0\left(\C_{\left( \widetilde\mu\right)}\right)^+$ is the minimal unitization (cf. Definition \ref{multiplier_min_defn}) of $C_0\left(\C_{\left( \widetilde\mu\right)}\right)$ and   $\widetilde C\bydef C_0\left(\C_{\left( \widetilde\mu\right)}\right)^+ \rtimes_r \Z$ is the reduced crossed product then there is the natural inclusion
 $$
 \widetilde B= C_0\left(\C_{\left( \widetilde\mu\right)}\right) \rtimes_r \Z\subset C_0\left(\C_{\left( \widetilde\mu\right)}\right)^+ \rtimes_r \Z=  \widetilde C
 $$
 which is a unitization (cf. Definition \ref{unitization_defn}). Moreover if
 \be\label{qg_c_eqn}
 C \bydef \widetilde C^{\Z_k} \bydef \left\{\left.  \widetilde c \in  \widetilde B \right| \forall g \in \Z_k\quad g  \widetilde c = \widetilde c\right\}
 \ee
 then the inclusion $B\subset C$ is also unitization, i.e. the condition (a) of the Definition \ref{fin_unitization_defn} holds. Let us prove (b), i.e. $\left(C, \widetilde{C}, G, \pi_C \right)$ is  a 
 unital  noncommutative finite-fold quasi-covering. From \eqref{qg_c_eqn}  it follows that $\left(C, \widetilde{C}, G, \pi_C \right)$ is  noncommutative finite-fold quasi-covering \ref{fin_quasi_defn}
  For any $t\in \Z$ denote by $\widetilde c_t \bydef 1_{C_0\left(\C_{\left( \widetilde\mu\right)}\right)^+}~u_t$ (cf equation \eqref{discr_cr_prod_cc_eqn}). We leave to the reader the proof of the following equation
 $$
 \forall \widetilde c \in \widetilde C \quad \exists~ c_0,..., c_{k-1} \in C\quad \widetilde c = \widetilde c_0c_0 +...+\widetilde c_{k-1}c_{k-1},
 $$
 or equivalently
\be\label{qg_ct_eqn}
 \widetilde C=  \widetilde c_0C \oplus ...\oplus\widetilde c_{k-1}C.
\ee
 It means that  $\widetilde C$ is an unital  noncommutative finite-fold quasi-covering	$\left(B ,\widetilde{B}, G, \widetilde{\pi} \right)$ (cf. Definition \ref{fin_unital_defn}).  Now let us prove that
 		\be\nonumber
 G\bydef  \left\{ \left.g \in \Aut\left(\widetilde{C} \right)~\right|~ gc = c;~~\forall c \in C\right\}\cong \Z_k
 \ee
 where the action $\Z_k \times \widetilde{C}\to \widetilde{C}$ comes from the given by \eqref{qg_c_eqn} periodic *-automorphism $\a$.
 From \eqref{qg_c_eqn} 
 and $\widetilde c_t = \widetilde c_1^t$ it follows that  $g \in G$ depends on $g\widetilde c_1$ only.
 From  $\widetilde c^k_1= \widetilde c_k \in C$ it follows that 
\be\label{qg_k_eqn}
\left( g\widetilde c_1\right)^k= g \widetilde c_k= \widetilde c_k
\ee
From the equations \eqref{discr_cr_prod_cc_eqn} and \eqref{discr_cr_prod_op_eqn} it follows that $\widetilde C \bydef C_0\left(\C_{\left( \widetilde\mu\right)}\right)^+ \rtimes_r \Z$ is $\Z$ graded. Taking into account \eqref{qg_k_eqn} one can deduce that the action of $G$ preserves this grading, it follows that $\left( g\widetilde c_1\right)^k = a u_1$ with $a \in C_0\left(\C_{\left( \widetilde\mu\right)}\right)^+$. Using \eqref{qg_k_eqn}  one has
$$
 \left( g\left( a u_1\right) \right)^k=  \left(a \left( g u_1\right) \right)^k = \widetilde c_k = 1_{C_0\left(\C_{\left( \widetilde\mu\right)}\right)^+} u_1^k.
$$
From the above equation it follows that 
$
\exists j \in \{0, ..., k-1\} \quad a =  e^{\frac{2\pi ij}{k}}1_{C_0\left(\C_{\left( \widetilde\mu\right)}\right)^+},
$, or equivalently $g = \a^j$
so $G$ uniquely depends on $j$ modulo $k$, so  $G\cong \Z_k$. It follows that  $\left(C, \widetilde{C}, G, \pi_C \right)$ is  
noncommutative finite-fold pre-covering (cf. Definition  \ref{fin_pre_defn}). Since both $B$ and $\widetilde B$ are essential ideals of   $C$ and $\widetilde C$ the quadruple $\left(B, \widetilde{B}, G, \left.\pi_C \right|_B\right)$ is also noncommutative finite-fold pre-covering.
  \end{proof}
\begin{exercise}
	Prove that in the above situation the quadruple  $\left(A, \widetilde{A}, \Z_k, \left.\pi_A\right|_{tA} \right)$ is an {associated with $\left(B, \widetilde{B}, \Z_k, \left.\pi_B\right|_B \right)$  noncommutative finite-fold  covering of *-algebras} (cf. Definition \ref{fin_chull_defn}).
	
\end{exercise}

\begin{appendices}
	\chapter{Categories, functors  and ordered sets}
	
	\section{Categories and limits}
	\begin{definition}\label{category_defn}\cite{goldblatt:topoi}
	A \textit{category} $\mathscr C$ comprises 
	\begin{itemize}
		\item[(a)] 	 a collection of things called $\mathscr C$-\textit{objects}; 
	\item[(b)] a collection of things called $\mathscr C$-\textit{arrows}  or $\mathscr C$-\textit{morphisms}; 
	\item[(c)] operations assigning to each  $\mathscr C$-arrow $f$  $\mathscr C$-object $\mathrm{dom}f$ (the 
		"domain" of $f$) and a $\mathscr C$-object $\mathrm{cod}f$ (the "codomain" of $f$). If $a =\mathrm{dom} f$  
		and $b =\mathrm{cod} f$ we display this as 
		\bean
		f: a \to b\quad \text{or}\quad  a \xrightarrow{f} b;
		\eean
		\item[(d)] an operation assigning to each pair $\left(g, f\right)$ of $\mathscr C$-arrows with $\mathrm{dom}g = \mathrm{cod} f$.
		 A $\mathscr{C}$-arrow  $g\circ f$, the \textit{composite of f and g}, having $\mathrm{dom}g\circ f= \mathrm{dom}f$ and  $\mathrm{cod }g\circ f= \mathrm{cod}g$, i.e. $g\circ f:  \mathrm{dom}f \to \mathrm{cod}g$  such that the following condition obtains:
		 
		 \textit{	Associative Law}: Given the configuration 
		$$
		a \xrightarrow{f}	b \xrightarrow{g}	c \xrightarrow{h}d
		$$ 
		of  $\mathscr{C}$-objects and $\mathscr{C}$-arrows  then $h \circ \left(g \circ f\right)= \left(h \circ g\right)\circ f$;
				\item[(e)] an assignment to each $\mathscr{C}$-object $b$ of a $\mathscr{C}$-arrow $\mathbb{1}_b : b \to b$, called the \textit{identity arrow} on $b$, such that 
				
				\textit{Identity Law}: For any $\mathscr{C}$-arrows $f: a \to b$ and $g : b \to c$ one has
				\bean
				\mathbb{1}_b \circ f = f \quad \text{and}\quad g\circ \mathbb{1}_b= g.
				\eean
		\end{itemize}

	\end{definition}

\begin{definition}\label{subcategory_defn}\cite{goldblatt:topoi}
If $\mathscr{C}$ is a category, and $a$ and $b$ are $\mathscr{C}$-objects, 
we introduce the symbol $\mathscr{C}\left( a, b\right) $ to denote the collection of all $\mathscr{C}$-arrows 
with $\mathrm{dom} = a$ and $\mathrm{cod} = b$, i.e. 
\be\label{category_not_eqn}
\mathscr{C}\left( a, b\right) \bydef \left\{f| f \text{ is } \mathscr{C} \text{arrow and } a \xrightarrow{f}b\right\}.
\ee
$\mathscr{C}$ is said to be a \textit{subcategory} of category $\mathscr{D}$, denoted $\mathscr{C}\subseteqq\mathscr{D}$, if
\begin{enumerate}
	\item[(i)] every $\mathscr{C}$-object is a $\mathscr{D}$-object, and 
	\item[(ii)] if $a$ and $b$ are any two $\mathscr{C}$-objects, then $\mathscr{C}\left( a, b\right)\subseteqq\mathscr{D}\left( a, b\right)$, i.e. all the 	$\mathscr{C}$-arrows $a \to b$ are present in $\mathscr{D}$. \\$\mathscr{C}$ is a \textit{full subcategory} of $\mathscr{D}$ if  $\mathscr{C}\subseteqq\mathscr{D}$, and 
	
	\item[(iii)] for any $\mathscr{C}$-objects $a$ and $b$, $\mathscr{D}\left( a, b\right)=\mathscr{C}\left( a, b\right)$, i.e. $\mathscr{D}$ has no arrows 
$a \to b$ other than the ones already in $\mathscr{C}$. 
\end{enumerate} 

\end{definition}
	\begin{definition}\label{dual_cat_defn}\cite{goldblatt:topoi}
	If  $\mathscr C$ is a category then its \textit{dual} or \textit{opposite} category $\mathscr C^{\mathrm{op}}$ is constructed as follows: 
	\begin{itemize}
		\item $\mathscr C$ and $\mathscr C^{\mathrm{op}}$ have the same objects.
		\item For each arrow $f: a\to b$ o we 
		introduce an arrow  $f^{\mathrm{op}}: b\to a$ in $\mathscr C^{\mathrm{op}}$, these being all and only the arrows.
		\item  The composite $f^{\mathrm{op}}\circ g^{\mathrm{op}}$ is defined precisely when $g\circ f$ is defined in $\mathscr C$
		and has $f\circ g\bydef \left(g\circ f \right)^{\mathrm{op}}$.
	\end{itemize}
	
\end{definition}
\begin{definition}\label{init_ob_defn}\cite{goldblatt:topoi}
	An object $0$ is \textit{initial} in a category $\mathscr C$ if for every $\mathscr C$-object $a$ there is one and only one arrow from $0$ to $a$ in $\mathscr C$.
\end{definition}
\begin{definition}\label{terminal_object_defn}\cite{goldblatt:topoi}
	An object $\mathbb{1}$ is \textit{terminal} or \textit{final} in a category $\mathscr C$ if for every $\mathscr C$-object $a$ there is one and only one arrow from $a$ to $\mathbb{1}$ in $\mathscr C$.
\end{definition}

\subsection{Limits}\label{limit_sec}
\paragraph{} Here I follow to \cite{goldblatt:topoi}. The notion of \textit{commutative diagram}, is a 
	very important aid to understanding used in category theory.
By a 
diagram we simply mean a display of some objects, together with some 
arrows (here representing functions) linking the objects. The "triangle" of 
arrows $f$, $g$, $h$ as shown is another diagram. 
		\newline
\begin{tikzpicture}
	\matrix (m) [matrix of math nodes,row sep=3em,column sep=4em,minimum width=2em]
	{
	A & B  \\ 
		& C\\};
	\path[-stealth]
		(m-1-1) edge node [above] {$f$} (m-1-2)
	(m-1-1) edge node [left]  {$h~~$} (m-2-2)
	(m-1-2) edge node [right] {$~~g$} (m-2-2);
\end{tikzpicture}
\\ 	
It will be said to \textit{commute} if $h = g\circ f$. The point is that the diagram offers 
two paths from $A$ to $C$, either by composing to follow $f$ and then $g$, or by 
following h directly. Commutativity means that the two paths amount to 
the same thing. A more complex diagram, like the previous one, is said to 
be commutative when all possible triangles that are parts of the diagram 
are themselves commutative. This means that any two paths of arrows in 
the diagram that start at the same object and end at the same object 
compose to give the same overall arrow. By a \textit{diagram} $D$ in a 
category $\mathscr C$ we simply mean a collection of $\mathscr C$-objects $d_j, d_k,...$ together 
with some $\mathscr C$-arrows $g: d_j \to d_k$ between certain of the objects in the 
diagram. (Possibly more than one arrow between a given pair of objects, 
possibly none). 
A \textit{cone} for diagram $D$ consists of a $\mathscr C$-object $c$ together with a $\mathscr C$-arrow 
$c\to d_j$ for each object $f_j$ in $D$, such that
		\newline
\begin{tikzpicture}
	\matrix (m) [matrix of math nodes,row sep=3em,column sep=4em,minimum width=2em]
	{
	d_j  & & d_k \\ 
		& c  & \\};
	\path[-stealth]
	(m-1-1) edge node [above] {$g$} (m-1-3)
	(m-2-2) edge node [left]  {$f_j~~$} (m-1-1)
	(m-2-2) edge node [right] {$~~f_k$} (m-1-3);
\end{tikzpicture}
\\ 	
 commutes, whenever g is an arrow in the diagram $D$. We use the 
symbolism $\left\{f_j: c\to d_j\right\}$ to denote a cone for $D$.

A \textit{limit} for a diagram $D$ is a $D$-cone $\left\{f_j: c\to d_j\right\}$ with the property that for 
	any other $D$-cone $\left\{f'_j: c'\to d_j\right\}$ there is exactly one arrow $f:c' \to c$ such 
			\newline
	\begin{tikzpicture}
		\matrix (m) [matrix of math nodes,row sep=3em,column sep=4em,minimum width=2em]
		{
			  & d_j &  \\ 
		c'	&  & c \\};
		\path[-stealth]
		(m-2-1) edge node [above] {$g$} (m-2-3)
		(m-2-1) edge node [left]  {$f'_j~~$} (m-1-2)
		(m-2-3) edge node [right] {$~~f_j$} (m-1-2);
	\end{tikzpicture}
	\\ 	
		commutes for every object $d_j$ in $D$. 
		This limiting cone, when it exists, is said to have the \textit{universal property} 
		with respect to $D$ -cones.
		 A limit for 
		diagram $D$ is unique up to isomorphism.

\begin{definition}\label{pull_back_defn}\cite{goldblatt:topoi}
	A \textit{pullback} of a pair $a \xrightarrow{f}c \xleftarrow{g} b$ of $\mathscr C$-arrows with a common codomain 
	is a limit in $\mathscr C$ for the diagram 
			\newline
\begin{tikzpicture}
	\matrix (m) [matrix of math nodes,row sep=3em,column sep=4em,minimum width=2em]
	{
		& b  \\ 
		a	&  c \\};
	\path[-stealth]
	(m-1-2) edge node [right]  {$g$} (m-2-2)
	(m-2-1) edge node [above] {$f$} (m-2-2);
\end{tikzpicture}
\\ 	
	
\end{definition}
\section{Functors}
\begin{definition}\label{functor_defn}\cite{goldblatt:topoi}
A \textit{functor} $F$ from category $\mathscr{C}$ to category $\mathscr{D}$ is a function that assigns 
\begin{enumerate}
	\item [(i)]
 to each $\mathscr{C}$-object $a$, a $\mathscr{D}$-object $F(a)$; 
\item[(ii)] to each $\mathscr{C}$-arrow $f:a \to b$ a $\mathscr{D}$-arrow $F(f): F(a) \to F(b)$, 
such that 
\begin{enumerate}
	\item[(a)]  $F\left(\mathbb 1_a\right) = \mathbb 1_{F\left(a\right)}$ for all  $\mathscr{C}$-objects $a$, i.e. the identity arrow on $a$ is assigned 
	the identity on $F\left(a\right)$,
	\item[(b)]  $F\left(g\circ f\right)=F\left(g\right)\circ F\left( f\right) $, whenever $g \circ f$ is defined. 
	This last condition states that the $F$-image of a composite of two arrows 
	is the composite of their $F$-images.
\end{enumerate}

\end{enumerate}
 We write $F:\mathscr{C}\to \mathscr{D}$ or  $\mathscr{C}\xrightarrow{F} \mathscr{D}$ to indicate that $F$ is a 
functor from $
\mathscr{C}$ to $\mathscr{D}$. Briefly then a functor is a transformation that 
"preserves" dom's, cod's, identities and composites. 
\end{definition}
If $a$ and $b$ are  $\mathscr{C}$-objects then a functor $\mathscr{C}\xrightarrow{F} \mathscr{D}$ yields a map
\be\label{f_ab_funct_eqn}
F_{a,b}:\mathscr{C}\left(a, b \right)  \to \mathscr{D}\left( F\left(a\right), F\left(b\right)\right)  
\ee
(cf. Notation \eqref{category_not_eqn}).
\begin{definition}\label{funct_full_faithfull_defn}\cite{bass}
A functor $\mathscr{C}\xrightarrow{F} \mathscr{D}$ is said to be \textit{faithful} (resp. \textit{full}) if the given by \eqref{f_ab_funct_eqn} map is injective (resp. surjective).
\end{definition}
\begin{defn}\label{functor_contravariant_defn}
A \textit{contravariant} functor is one that reverses direction by mapping domains 
to codomains and vice versa. 
Thus $\mathscr{C}\xrightarrow{F} \mathscr{D}$ is a contravariant functor if it assigns to $f: a\to b$ an 
arrow $F(f):F(b)\to F(a)$, so that $F\left(\mathbb{1}_a\right)= \mathbb{1}_{F(a)}$ as before, but now 
$$
F\left( g\circ f\right) = F\left( f\right)\circ  F\left( g\right). 
$$
\end{defn}

\begin{definition}\label{exact_functor_defn}
A left/right exact functor is a functor that preserves finite limits/finite colimits.
\end{definition}
\section{Natural transformations}\label{natural_transformation_sec}
\paragraph{}
Here I follow to \cite{goldblatt:topoi}.
Given two categories $\mathscr C$ and $\mathscr D$ we are going to construct a category, 
denoted $\mathrm{Funct}\left(\mathscr C, \mathscr D\right)$, or $\mathscr D^{\mathscr C}$, whose objects are the functors from $\mathscr C$ to $\mathscr D$. 
We need a definition of arrow from one functor to another. Let 
$F: \mathscr C\to \mathscr D$ and $G: \mathscr C\to \mathscr D$ be two functors. For any $\mathscr C$-object $a$ we define a $\mathscr D$-arrow $\tau_a : F\left(a\right)\to  G\left(a\right)$. We require that each  $\mathscr C$-arrow $f: a \to b$  gives rise to a diagram 
	\begin{tikzpicture}
	\matrix (m) [matrix of math nodes,row sep=3em,column sep=4em,minimum width=2em]
	{
	a	& F\left(a \right)  &   G\left(a \right)\\ 
		b	& F\left(b \right)  &   G\left(b \right) \\};
	\path[-stealth]
	(m-1-1) edge node [left]  {$f$} (m-2-1)
	(m-1-2) edge node [above] {$\tau_a$} (m-1-3)
	(m-2-2) edge node [above]  {$\tau_b$} (m-2-3)
	(m-1-2) edge node [left]  {$F\left( f \right)$} (m-2-2)
	(m-1-3) edge node [left]  {$G\left( f \right)$} (m-2-3);
\end{tikzpicture}
\\ 	
that commutes.
In summary then, a \textit{natural transformation} from functor $F: \mathscr C\to \mathscr D$ and $G: \mathscr C\to \mathscr D$  to functor $F: \mathscr C\to \mathscr D$ and $G: \mathscr C\to \mathscr D$  is an assignment $\tau$ that provides, for each $\mathscr C$-object $\mathscr D$-arrow $\tau_a :F(a) \to G(a)$, such that for any $\mathscr C$-arrow $f:a\to b$, the above diagram commutes in $\mathscr D$, i.e. $\tau_b \circ F(f)= G(f)\circ \tau_a$. We use the symbolism $\tau: F\to G$, or $F \xrightarrow{\tau}G$, to denote that $\tau$ is a natural transformation from $F$ to $G$. The arrows $\tau_a$ are called the \textit{components} of $a$. Now if each component $\tau_a$ of $a$ is an iso arrow in $\mathscr D$ then  case we call $\tau$ a \textit{natural isomorphism}. Each $\tau_a: F(a)\to G(a)$ then has an inverse $\tau_a^{-1}: G(a) \to F(a)$, and these $\tau^{-1}_a$'s form the components of a natural isomorphism $\tau^{-1}: G \to F$. We denote natural isomorphism by $\tau: F \cong G$. 
\begin{example}
The identity natural transformation $\mathbb{1}_F :F\to F$ assigns to each object $a$, the identity arrow $\mathbb{1}_{F(a)}:F(a)\to F(a)$. This is clearly a natural isomorphism. 
\end{example}

	\begin{definition}\label{category_equivalence_definition}\cite{goldblatt:topoi}
	A functor $F: \mathscr C\to  \mathscr D$ is called an \textit{equivalence of categories} if there 
	is a functor $G: \mathscr D\to  \mathscr C$ such that there are natural isomorphisms $\tau : 1_{\mathscr C} \cong G \circ F$, and $\sigma : 1_{\mathscr D} \cong F \circ G$, from the identity functor on ${\mathscr C}$ to $ G \circ F$, and 
	from the identity functor on ${\mathscr D}$ to $ F \circ G$.
	
	Categories $\mathscr C$ and $\mathscr D$ are \textit{equivalent},  $\mathscr C$ and $\mathscr D$ when there exists an equivalence $F: \mathscr C\to  \mathscr D$ .
\end{definition}

	\section{Skeleton and pre-order category}

	\begin{definition}\label{skeletal_defn}\cite{goldblatt:topoi}
		A \textit{skeletal category} is one in which "isomorphic" does actually mean the 
		same as "is", i.e. in which whenever $a\cong b$, then $a=b$. 
	\end{definition}
	\begin{definition}\label{skeleton_defn}\cite{goldblatt:topoi}
		A \textit{skeleton}
		of a category $\mathscr C$ is a full subcategory $\mathscr C_0$ of $\mathscr C$ that is skeletal, and such that each $\mathscr C$-object is isomorphic to one (and only one) $\mathscr C$-object.
	\end{definition}
	\begin{remark}
For any category there is a skeleton.
	\end{remark}
\begin{definition}\label{preordercat_defn}\cite{goldblatt:topoi}
		A category is said to be a \textit{pre-order category} if there is at most one morphism between different objects.
	\end{definition}
	\begin{definition}\label{ord_defn}\cite{goldblatt:topoi}
		A binary relation $R \subset \La\times \La$ on the set $\La$ (writing $pRq$ in place of $\left(p,q \right)\in R$ ) is said to be \textit{pre-ordering} if it is 
		\begin{enumerate}
			\item [(a)] \textit{reflexive}, i.e. for each $p$ we have $pRp$ and
			\item[(b)] \textit{transitive}, whenever $pRq$ and $qRs$, we have $pRs$.
		\end{enumerate}
	\end{definition}
	\begin{remark}\label{pre_order_rem}\cite{goldblatt:topoi}
		For any set pre-ordered $\La$ with pre-ordering $R$ there is a pre-order category such that
		\begin{enumerate}
			\item [(a)] Objects of the category are elements of $\La$.
			\item[(b)] Morphisms are pairs $\left(p, q \right) \in R$.
			\item[(c)] The composition of morphisms is given by $\left(q,s\right)\circ \left(p, q\right) =\left(p,s\right)$.
		\end{enumerate}
		
	\end{remark}
	\section{Equivalence relations}
	
	\begin{definition}\label{equivalence_relation_defn}\cite{spanier:at}
	An \textit{equivalence relation} in a set $A$ is a relation $\sim$ between elements of $A$ 
	which is \textit{reflexive} (that is, $a\sim a$ for all $a \in A$), \textit{symmetric} (that is, $a\sim a'$ 
	implies $a'\sim a$ for $a, a'\in A$), and \textit{transitive} (that is, $a\sim a'$ and $a'\sim a''$
	imply $a \sim a''$ for $a, a', a'' \in A$). The \textit{equivalence class} of $a \in A$ with 
	respect to $\sim$ is the subset $\left\{\left. a'\in A\right| a'\sim a \right\}$. The set of all equivalence classes 
	of elements of $A$ with respect to ~ is denoted by $A/\sim$ and is called a 	\textit{quotient set} of $A$. There is a \textit{projection map} $A \to A/\sim$ which sends $a\in A$ to
	its equivalence class. 
	\end{definition}

	\section{Directed sets}
	
		\begin{definition}\label{directed_set_defn}\cite{engelking:general_topology}
		Let $\La$ be a set and $\le$ is relation on $\La$. We say that $\le$ \textit{directs} $\La$ or $\La$ is \textit{directed} by $\le$, if $\le$ has following properties:
		\begin{enumerate}
			\item [(a)] If $\la \le \mu$ and $\mu \le \nu$, then $\la \le \nu$,
			\item[(b)] For every $\la \in \La$, $\la \le \la$,
			\item[(c)] For any $\mu,\nu \in \La$ there exists a $\la \in \La$ such that $\mu \le \la$ and $\nu \le \la$.
		\end{enumerate}
	\end{definition}
	\begin{definition}\label{cofinal_defn}\cite{engelking:general_topology}
		A subset $\Xi \subset \La$ is said to be \textit{cofinal} in $\La$ if for every $\la \in \La$ there is $\chi \in \Xi$ such that $\la \le \chi$. 
	\end{definition}
\section{Inverse systems}
\begin{definition}\label{inverse_limit_defn}\cite{spanier:at}
An \textit{inverse system of sets} $\left\{A_\a, f^\bt_\a\right\}$ consist of a collection of sets $\left\{A_\a\right\}$ indexed by by a directed set $\La = \left\{\a\right\}$ and a collection of functions $f^\bt_\a: A_\bt \to A_\a$ such that
\begin{enumerate}
	\item [(a)] $f^\a_\a = 1_{A_\a}: A_\a\cong A_\a$ for $\a\in\La$
	\item[(b)] $f^\ga_\a = f^\bt_\a \circ f^\ga_\bt : A_\ga \to A_\a$ for $\a \le \bt \le \ga$ in $\La$
\end{enumerate}
The \textit{inverse limit} $\varprojlim A_\a$ is the subset of $\prod A_\a$ consisting of all points $\left\{a_\a\right\}$ such that if $\a\le\bt$, then $a_\a = f^\bt_\a\left( a_\bt\right) $. For each $\a$ there is a map $p_\a: \varprojlim A_\a\to A_\a$, and if $\a \le \bt$, then $p_\a = f^\bt_\a \circ p_\bt$ . 
\end{definition}
\begin{lemma}\label{inv_lim_lem}\cite{spanier:at}
 Given an inverse system of sets  $\left\{A_\a, f^\bt_\a\right\}$ and given a set $B$ and for every $\a\in \La$ function $g_\a: B \to A_\a$ such that $g_\a = f^\bt_\a \circ g_\bt$ if $\a\le\bt$, there is a unique function $g : B \to \varprojlim A_\a$  such that $g_\a = p_\a\circ g$ for all $\a\in \La$. 
\end{lemma}

\section{Ordinals}\label{ord_sec}
\paragraph*{}
Here the book \cite{harzheim:os} is cited.
A simply ordered set is called \textit{well-ordered}, if each non-empty 
subset of it has a first (or least) element.
We recall in short, without going into the details, the construction of 
the class of ordinal numbers by the method of von Neumann. Consider 
the sequence 
$	\emptyset ,\{\emptyset\},\{\emptyset, \{\emptyset\}\},\{\emptyset, \{\emptyset\}, \{\emptyset, \{\emptyset\}\}\},....$, where from every element x of 
the sequence we create its immediate successor by forming the set $x \cup \{x\}$. 
This step can be repeated transfinitely. The formal definition is: 
A set $M$ whose elements are again sets is called an \textit{ordinal} if there 
holds: $M$, with the relation $\subseteqq$, is well-ordered, and each $x\in M$ is the 
set of all elements of M which are $\subset x$. 
This entails that the class On of all ordinals is well-ordered, and if
$\a$ and $\bt$ are ordinals, we have 
$$
\a < \bt \Leftrightarrow a \text{ is a proper subset of }\bt \Leftrightarrow\a\in \bt.
$$.
The first infinite ordinal is denoted by $\om$ or $\om_0$, and the finite  
ordinals by $0,1,2,...$. 
If $\la$ is an ordinal $\neq \emptyset$, for which the set of ordinals which are $< \la$ has 
no greatest element, then $\la$ is called a \textit{limit ordinal} or \textit{limit number}. If 
this set has a greatest element $\kappa$, $\la$ is said to be a \textit{successor ordinal} or 
\textit{successor number}, and then we also denote $\kappa$ by $\la -1$. 
In On we have operations + and $\bullet$  which can be defined by transfinite 
recursion: 
If $a$ is an ordinal we put $a + 0 = a$ and define $a + 1$ to be the least 
ordinal which is $> a$. If for an ordinal $\bt$ we have already defined $\a +  \bt$
we put $\a + (\bt + 1) = (a + \bt) + 1$.
\begin{thm}\label{zorn_thm}\cite{spanier:at} (Zorn's lemma). A partially ordered set in which every  simply ordered set has an upper bound contains maximal elements.
\end{thm}

	\subsection{Projective systems}\label{proj_sys_sec}
\paragraph*{}
Here I follow to \cite{phillips:inv_lim_app,phillips:inv_lim,spanier:at}.
An \textit{inverse system} (also called a 
\textit{projective system}) of sets consists of a directed set $D$, a set $X_d$ for each $d\in D$
and functions $\pi_{d,e} : X_d \to X_e$ for all $d,e\in D$ with $d\ge e$, satisfying the consistency 
conditions $\pi_{d,e}= \Id{X_d}$ and $\pi_{e,f}\circ \pi_{d,e} = \pi_{d,f}$. The inverse limit (or 
projective limit) of the system $\left\{X_d\right\}$ is a set $X$ together with functions $\ka_d: X \to X_d$
such that $\pi_{d,e}\circ \ka_e = \ka_d$ for $d \ge e$, and satisfying the following universal property: 
given any set $Y$ and functions $\varphi_d : Y \to X_d$ such that $\pi_{d,e}\circ \varphi_{d,e} = \varphi_d$, there is a 
unique function $\phi : Y \to X$ such that $\phi_s = \ka_d\circ \pi$ for all $d$. The inverse limit $X$ 
is denoted by $\varprojlim X_d$, and it is clearly unique up to unique isomorphism. 
The definition of an inverse limit generalizes in an obvious way to any 
category, and it is easily seen that the same construction produces an inverse limit in 
the categories of abelian groups and homomorphisms, topological spaces and continuous 
maps, and topological *-algebras over $\C$ and continuous $*$-homomorphisms. (The product in (*) is given the product topology and the pointwise algebraic operations).

\chapter{Topology}

	\section{General topology}
	\begin{definition}\label{top_base_defn}\cite{munkres:topology}
		If $\sX$ is a set, a \textit{basis} for a topology on $\sX$ is a collection $\mathscr B$ of subsets of $\sX$
		(called \textit{basis elements}) such that:
		\begin{enumerate}
			\item [(a)]  for each $x\in \sX$, there is at least one basis element $B$ containing $x$,
			\item [(b)]  if $x$ belongs to the intersection of two basis elements $B_1$ and $B_2$, then there is a
			basis element $B_3$ containing $x$ such that $B_3 \subset B_1\cap B_2$.
		\end{enumerate}
		
		If  $\mathscr B$ satisfies these two conditions, then we define the \textit{topology}  $\mathscr T$ \textit{generated by}  $\mathscr B$ as
		follows: A subset $\sU$ of $\sX$ is said to be \textit{open} in $\sX$ (that is, to be an element of $\mathscr T$) if for
		each $x \in \sU$, there is a basis element $B \in \mathscr B$ such that $x \in B$ and $B \subset \sU$.
	\end{definition}
	\begin{remark}\label{top_base_rem}It is proven that given by the Definition \ref{top_base_defn} collection $\mathscr T$ is a topology  (cf. \cite{munkres:topology}).
	\end{remark}
	\begin{proposition}\label{top_final_prop}\cite{bourbaki_sp:gt}
Let $\sX$ be a set, let $\left\{\sY_\iota\right\}_{\iota\in I}$ be a family of topological spaces, and for each $\iota\in I$ let $f_\iota$ be a mapping of $\sY_\iota$ into $\sX$. Let $\mathfrak{D}$ be a set of subsets $\sU$ of $\sX$ such that $f_\iota^{-1} \left(\sU \right)$ is open in $\sY_\iota$ for each $\iota \in I$; then  $\mathfrak{D}$ is a set of open sets in a topology $\mathfrak{T}$. In particular $\mathfrak{T}$ is the finest topology for which the mappings $f_\iota$ are continuous. In other words, if $g$ is mapping on $\sX$ into a topological space $\sZ$, then $g$ is continuous ($\sX$ carrying the topology $\mathfrak{T}$) if and only if each of the mappings $g\circ f_\iota$ is continuous.
	\end{proposition}
	\begin{definition}\label{top_final_defn}
	Under the hypotheses of Proposition \ref{top_final_prop} we say that the topology $\mathfrak{T}$ is \textit{final} (\textit{with respect to the family of maps} $\left\{f_\iota:\sY_\iota \to \sX\right\}_{\iota\in I}$).
	\end{definition}
		\begin{remark}
		The Definition \ref{top_final_defn} is a specialization of final object (cf. Definition \ref{terminal_object_defn}).
	\end{remark}
	\begin{corollary}\label{top_final_cor}
	Under the hypothesis of the Proposition \ref{top_final_prop} a subset $\sV$ of $\sX$ is closed in the topology  $\mathfrak{T}$ if and only if $f^{-1}\left(\sV \right)$ is closed in $\sY_\iota$ for each $\iota \in I$. 
	\end{corollary}
	
	\begin{definition}\label{top_pointed_defn}
		If $\sX$ is a set (resp. topological space) and $x_0 \in \sX$ is a point then we say that the pair $\left(\sX, x_0 \right)$ is a \textit{pointed  set} (resp. \textit{pointed space}) (cf. \cite{spanier:at,switzer:at}). If both $\left(\sX, x_0 \right)$ and $\left(\sY, y_0 \right)$ are  pointed  spaces and $\varphi: \sX \to \sY$ is such that $\varphi\left(x_0\right)= y_0$ then we say that $\varphi$ is a \textit{pointed map}. We write
		\be\label{top_pointed_eqn}
		\varphi: \left(\sX, x_0 \right)\to\left(\sY, y_0 \right).
		\ee
		We say that $x_0$ is the \textit{base}-\textit{point}.
	\end{definition}
\begin{definition}\label{top_hausdorff_defn}\cite{munkres:topology}
	A topological space $\sX$ is called a Hausdorff space if for each pair $x_1, x_2$
	of distinct points of $\sX$, there exist neighborhoods $\sU_1$, and $\sU_2$ of $x_1$ and , respectively,
	that are disjoint.
\end{definition}	
\begin{definition}\label{top_connected_defn}\cite{munkres:topology}
Let $\sX$ be a topological space. A \textit{separation} of $\sX$ is a pair $\sU$, $\sV$ of disjoint
nonempty open subsets of $\sX$ whose union is $\sX$. The space  $\sX$ is said to be \textit{connected}
if there does not exist a separation of  $\sX$.
\end{definition}
\begin{remark}\label{top_conn_rem}\cite{munkres:topology}
A space $\sX$ is connected if and only if the only subsets of $\sX$ that are both
open and closed in $\sX$ are the empty set and $\sX$ itself.
\end{remark}
\begin{theorem}\label{top_connected_closure_thm}\cite{munkres:topology}
Let $A$ be a connected subspace of $\sX$, and let $\overline A$ be the closure of $A$, i.e. intersection of containing $A$ closed sets. If $A \subset B \subset \overline A$, then $B$ is also
connected.
\end{theorem}
\begin{definition}\label{top_path_connected_defn}\cite{munkres:topology}
Given points $x$ and $y$ of the space $\sX$, a \textit{path} in $\sX$ from $x$ to $y$ is a continuous map $f: \left[a, b\right]\to \sX$ of some closed interval in the real line into $\sX$, such
that $f(a) = x$ and $f(b) = y$. A space $\sX$ is said to be \textit{path connected} if every pair of
points of $\sX$ can be joined by a path in $\sX$.
\end{definition}

\begin{definition}\label{top_connected_component_defn}\cite{munkres:topology}
	Given $\sX$, define an equivalence relation on $\sX$ by setting $x\sim y$ if there
	is a connected subspace of $\sX$ containing both $x$ and $y$. The equivalence classes are
	called the \textit{components} (or the \textit{connected components}) of $\sX$.
\end{definition}

\begin{theorem}\label{top_conn_comp_thm}\cite{munkres:topology}
	The components of $\sX$ are connected disjoint subspaces of $\sX$ whose
	union is  $\sX$, such that each nonempty connected subspace of  $\sX$ intersects only one of
	them.
\end{theorem}
	\begin{definition}\label{top_locally_connected_defn}\cite{munkres:topology}.
	A space $\mathcal X$ is said to be \textit{locally connected at} $x$ if for every neighborhood $\mathcal U$ of $x$ there is a connected neighborhood $\mathcal V$ of $x$  contained in $\mathcal U$. If $\mathcal X$ is locally connected at each of its points, it is said simply to be \textit{locally connected}.
\end{definition}
\begin{theorem}\label{top_loc_conn_thm}\cite{munkres:topology}
	A space $\mathcal X$ is locally connected if and only if for every open set $\mathcal U$, each component of $\mathcal U$ is open in $\mathcal X$.
\end{theorem}
\begin{proposition}\label{top_quasi_component_prop}\cite{bredon:topology_geometry}
The statement "$d(p) = d(q)$" for every discrete valued map $d$
	on $\sX$" is an equivalence relation. 
\end{proposition}
\begin{definition}\label{top_quasi_component_defn}\cite{bredon:topology_geometry}
The equivalence classes of the relation in Proposition \ref{top_quasi_component_prop}
are called the \textit{quasi}-\textit{components} of $\sX$.
\end{definition}
\begin{proposition}\label{top_quasi_component_closed_prop}\cite{bredon:topology_geometry}
Quasi-components of a space $\sX$ are closed. Each connected
set is contained in a quasi-component. (In particular, each connected component is contained
in a quasi-component.) Quasi-components are either equal or disjoint,
and fill out $\sX$.
\end{proposition}

\begin{exercise}\label{top_loc_conn_exer}\cite{bredon:topology_geometry}
Prove that If $\sX$ is locally connected, show that its components are open and equal its
quasi-components.
\end{exercise}

\begin{definition}\cite{spanier:at}
Given a set $\sX$ and an indexed collection of topological spaces $\left\{\sX_j\right\}_{j \in J}$ functions $f_j : \sX \to \sX_j$, the \textit{topology induced} on $\sX$ by the functions $\left\{f_j\right\}$ is the 
smallest or coarsest topology such that each $f_j$ is continuous. 
\end{definition}

\begin{definition}\label{top_inverse_limit_defn}\cite{spanier:at}
If $\left\{\sX_\a\right\}_{\a\in\La}$ is an inverse system of topological 
spaces (that is, $\sX_\a$ is a topological space for $\a\in \La$ and $f^\bt_\a: \sX_\bt\to \sX_\a$ is  
continuous for $\a\le\bt$) their \textit{inverse limit} $\varprojlim \sX_\a$ (cf. Definition \ref{inverse_limit_defn}) is given the topology induced by the functions $p_\a : \varprojlim \sX_\a\to\sX_\a$  for $\a\in\La$. 
\end{definition}
\begin{definition}\label{top_cover_defn}\cite{munkres:topology}
A collection $\A$ of subsets of a space $\sX$ is said to \textit{cover} $\sX$, or to be a
\textit{covering} of $\sX$, if the union of the elements of $\A$ is equal to $\sX$. 
It is called an \textit{open
covering} of $\sX$ if its elements are open subsets of $\sX$.
\end{definition}

\begin{definition}\label{top_compact_defn}\cite{munkres:topology}
A space $\sX$ is said to be \textit{compact} if every open covering $\A$ of $\sX$ contains
a finite subcollection that also covers $\sX$.
\end{definition}	
	\begin{theorem}\label{top_compact_img_thm}\cite{munkres:topology}
	The image of a compact space under a continuous map is compact.
\end{theorem}	

\begin{definition}\label{top_locally_compact_defn}
	A space $\sX$ is said to be \textit{locally compact} at $x$ if there is some compact
	subspace $\sV$  of $\sX$ that contains a neighborhood of $x$. If $\sX$ is locally compact at each of
	its points, $\sX$ is said simply to be \textit{locally compact}.
\end{definition}
\begin{empt}\cite{engelking:general_topology}
		Every set of cardinal numbers being well-ordered by $<$, the set of all cardinal numbers
		of the form $\left|B\right|$, where $B$ is a base for a topological space $\sX$, has a smallest element; this cardinal number is called the \textit{weight of the topological space} $\sX$ and is denoted by $w\left(\sX\right)$.
		
	\end{empt}

	\begin{theorem}\label{top_weght_thm}\cite{engelking:general_topology}
		If the weight of both $\sX$ and $\sX$ is not larger than $\mathfrak{m} \ge \aleph_0$ and $\sY$ is a locally
		compact space, then the weight of the space $\sY^\sX$ with the compact-open topology (cf. \ref{top_comp_open_empt})  is not larger
		than $\mathfrak{m}$.
	\end{theorem}	
	
\begin{definition}\label{top_net_defn}\cite{engelking:general_topology}
		A \textit{net in topological space} $\mathcal{ X}$ is an arbitrary function from a non-empty directed set $\La$ (cf. Definition \ref{directed_set_defn}) to the space $\mathcal{ X}$. Nets will be denoted by $S = \left\{x_\la \in \mathcal{ X}\right\}_{\la \in \La}$. 
	\end{definition}
	\begin{definition}\label{top_net_lim_defn}\cite{engelking:general_topology}
		A point $ x \in \mathcal{ X}$ is called a \textit{limit of a net} $S = \left\{x_\la \in \mathcal{ X}\right\}_{\la \in \La}$ if for every neighborhood $\mathcal U$ of $x$  there exists $\la_0 \in \La$ such that $x_\la \in \mathcal U$ for every $\la \ge \la_0$; we say that $\left\{x_\la \right\}$ \textit{converges} to $x$. The limit will be denoted by $x = \lim S$, or $x = \lim_{\la \in \La} x_\la$, or $\lim x_\la$.
	\end{definition}
	\begin{remark}
		A net can converge to many points, however if $\mathcal{ X}$ is Hausdorff then the net can have the unique limit.
	\end{remark}
	\begin{definition}\label{haudorff_defn}\cite{munkres:topology}
		A topological space $\sX$ is called a \textit{Hausdorff space} if for each pair $x_1$, $x_2$ of distinct points of $\sX$, there exist neighborhoods $\sU_1$, and $\sU_2$ of $x_1$  and $x_2$, respectively, that are disjoint.
	\end{definition}
	\begin{definition}\label{top_normal_defn}\cite{munkres:topology}
		Suppose that one-point sets are closed in $\sX$. Then $\sX$ is said to be  
		\textit{regular} if for each pair consisting of a point $x$ and a closed set $B$ disjoint from $x$, there
		exist disjoint open sets containing $x$ and $B$ respectively. The space $\sX$ is said to be
		\textit{normal} if for each pair $A$, $B$ of disjoint closed sets of $\sX$ there exist disjoint open sets containing $A$ and $B$ respectively.
	\end{definition}
	\begin{definition}\label{top_completely_regular_defn}\cite{munkres:topology}
		A topological space $\mathcal X$ is \textit{completely regular} if one-point sets are closed in  $\mathcal X$ and for each point $x_0$ and each closed $\mathcal \sY \subset \mathcal X$ not containing $x_0$, there is a continuous function $f: \mathcal X \to \left[0,1 \right]$ such that $f\left(x_0 \right)= 1$ and $f\left(\mathcal Y \right)= \left\{0\right\} $.  
	\end{definition}
	\begin{theorem}\label{top_sc_regular_normal_thm}\cite{munkres:topology}
		Every regular space with a countable basis is normal.
	\end{theorem}
	\begin{exercise}\label{top_completely_regular_exer}\cite{munkres:topology}
		Show that every locally compact, Hausdorff space is completely regular.
	\end{exercise}
	
	\begin{thm}\label{comp_normal_thm}\cite{munkres:topology}
		Every compact Hausdorff space is normal.
	\end{thm}
	
	
	\begin{thm}\label{urysohn_lem}\cite{munkres:topology}\textbf{ Urysohn lemma.}
		Let $\mathcal X$ be a normal space, let $\mathcal A$, $\mathcal B$ be disconnected closed subsets of $\mathcal X$. Let $\left[a,b\right]$ be a closed interval in the real line. Then there exist a continuous map $f: \mathcal X \to \left[a, b\right]$ such that $f(\mathcal A)=\{a\}$ and $f(\mathcal B)=\{b\}$.
	\end{thm}
	\begin{theorem}\label{tietze_ext_thm}\cite{munkres:topology} \textbf{Tietze extension theorem.} Let  $\mathcal X$ be a normal space; let  $\mathcal A$ be a closed subspace of  $\mathcal X$.
		\begin{enumerate}
			\item [(a)] Any continuous map of  $\mathcal A$ into the closed interval $[a,b]$ of $\R$ may be extended to a continuous map of all  $\mathcal X$ into $[a,b]$
			\item [(b)] Any continuous map of  $\mathcal A$ into $\R$ may be extended to a continuous map of all  $\mathcal X$ into $\R$.
		\end{enumerate}
		
	\end{theorem}
\begin{definition}\label{top_onen_map_defn}\cite{munkres:topology}
A map $f: \sX \to \sY$ is said to be an \textit{open map} if for every open set $\sU$ of $\sX$, the
set $f\left(\sU \right)$  is open in $\sY$.
\end{definition}	
	\begin{definition}\label{top_bicont_defn}\cite{kurat:topI}
A continuous map 	 	$f$ is called \textit{bicontinuous} if it is onto and if
	 $$
 \left(A \text{ is open}\right) \Leftrightarrow  \left(f^{-1}\left( A\right)  \text{ is open}\right)
	 $$
	 Equivalently 
	 $$
	 \left(A \text{ is closed}\right) \Leftrightarrow  \left(f^{-1}\left( A\right)  \text{ is closed}\right).
	 $$	
	 \end{definition}
	 \begin{theorem}\label{top_bicont_oc_thm}\cite{kurat:topI}
	Let $f: \sX \to \sY$ be open or closed onto. Then $f$ is bicontinuous.
\end{theorem}
	 \begin{theorem}\label{top_bicont_thm}\cite{kurat:topI}
	 If $f$ is bicontinuous and one-to-one, then $f$ is a homeomorphism.
	 \end{theorem}
	 	 \begin{theorem}\label{top_bicont_inv_thm}\cite{kurat:topI}
	 	If $f$ is bicontinuous and $g: \sY \to \sZ$. If $h \bydef g \circ f$ is continuous, then so is $g$. If $h$ is bicontinuous, then so is $g$.
	 	\newline
	 	\begin{tikzcd}
	 		\sX \arrow[rd, "h"] \arrow[r, "f"]
	 		& \sY\arrow[d, "g"]\\
	 		& \sZ
	 	\end{tikzcd}
	 \end{theorem}
	 
	 \begin{definition}\label{top_separable_defn}\cite{munkres:topology}
	 A space having a countable  dense subset is said to be \textit{separable}.
	 \end{definition}	
	\begin{defn}\label{top_support_defn}\cite{munkres:topology}
		If $\phi: \mathcal X \to \mathbb{C}$  is continuous then the \textit{support} of $\phi$ is defined to be the closure of the set $\phi^{-1}\left(\mathbb{\C}\setminus \{0\}\right)$ Thus if $x$ lies outside the support, there is some neighborhood of $x$ on which $\phi$ vanishes. Denote by $\supp \phi$ the support of $\phi$.
	\end{defn}

	There are two equivalent definitions of $C_0\left(\mathcal{X}\right)$ and both of them are used in this book.
	\begin{defn}\label{c_c_closure_defn}\cite{murphy}
		An algebra $C_0\left(\mathcal{X}\right)$ is the $C^*$-norm closure of the algebra $C_c\left(\mathcal{X}\right)$ of compactly supported continuous complex-valued functions.
	\end{defn}
	\begin{defn}
		\label{c_c_compact_defn}\cite{blackadar:ko,murphy}
		A $C^*$-algebra $C_0\left(\mathcal{X}\right)$ is given by the following equation
	\bean
			C_0\left(\mathcal{X}\right) \ \left\{\varphi \in C_b\left(\mathcal{X}\right) \ | \ \forall \varepsilon > 0 \ \ \exists K \subset \mathcal{X} \ ( K \text{ is compact}) \ \& \ \forall x \in \mathcal X \setminus K \ \left|\varphi\left(x\right)\right| < \varepsilon  \right\},
		\eean
		i.e.
		\begin{equation*}
			\left\|\varphi|_{\mathcal X \setminus K}\right\| < \varepsilon.
		\end{equation*}
	\end{defn}
	\begin{definition}\label{top_ind_lim_defn}\cite{candel:foliII,MW08}

	Let $\sX$ be a locally compact Hausdorff space, and let $C_c\left(\sX\right)$ denote the	space of continuous, compactly supported functions on $\sX$. The natural	topology on $C_c\left(\sX\right)$ is the \textit{inductive limit topology}, defined as follows. A net	$\left\{f_\a\right\}$ in $C_c\left(\sX\right)$ converges to $f$ if there is a compact set $K\subset\sX$ containing	the supports of all $f_\a$ and $f$, and such that $f_\a$ converges uniformly to $f$ on$K$.
	\end{definition}
	\begin{defn}\label{top_loc_fin_defn}\cite{munkres:topology}
		An indexed family of sets $\left\{A_\al\right\}$ of topological space $\sX$ is said to be \textit{locally finite} if each point $x$ in $\sX$  has a neighborhood that intersects for only finite many values of $\a$.
	\end{defn}
	\begin{defn}\label{top_part_of_unity_defn}\cite{munkres:topology}
		Let $\left\{\mathcal U_\alpha\in \mathcal X\right\}_{\alpha \in \mathscr A}$ be an indexed open covering of $\mathcal{X}$. An indexed family of functions 
		\begin{equation*}
			\phi_\alpha : \mathcal X \to \left[0,1\right]
		\end{equation*}
		is said to be a {\it partition of unity }, dominated by $\left\{\mathcal{U}_\alpha \right\}$, if:
		\begin{enumerate}
			\item $\phi_\alpha\left(\mathcal X \setminus \mathcal U_\alpha\right)= \{0\}$
			\item The family $\supp\phi_\alpha$ is locally finite.
			\item $\sum_{\alpha \in \mathscr A}\phi_\alpha\left(x\right)=1$ for any $x \in \mathcal X$.
		\end{enumerate}
	\end{defn}
 \begin{definition}\label{top_paracompact_defn}\cite{munkres:topology}
	A space $\sX$ is \textit{paracompact} if every open covering $\sX = \cup~\sU_{   \a}$ has a locally finite open refinement $\sX = \cup~\sV_{   \bt}$.
\end{definition}	
\begin{thm}\label{top_part_u_thm}\cite{munkres:topology}
	Let $\mathcal X$ be a paracompact Hausdorff space; let $\left\{\mathcal U_\alpha\in \mathcal X\right\}_{\alpha \in \mathscr A}$ be an indexed open covering of $\mathcal{X}$. Then there exists a partition of unity, dominated by $\left\{\mathcal{U}_\alpha \right\}$.  
\end{thm}

	\begin{prop}\label{top_smooth_part_unity_prop} \cite{brickell_clark:diff_m}
		A differential manifold $M$ admits a (smooth) partition of unity if and only if it is paracompact. 
	\end{prop}
	\begin{defn}\label{lindel_defn}\cite{engelking:general_topology}
	 We say that a topological space $\sX$ is a \textit{Lindel\"{o}f space}, or has  a \textit{Lindel\"{o}f property}, if $\sX$ is regular and every open cover of $\sX$ contains a countable  subcover. 
	 \end{defn}
	\begin{theorem}\label{lindel_thm}\cite{engelking:general_topology}
	Every regular second-countable  space in Lindel\"{o}f.
	\end{theorem}
	\begin{theorem}\label{lindel_para_thm}\cite{munkres:topology}
		Every regular Lindel\"{o}f space is paracompact.
	\end{theorem}
	\begin{definition}\label{top_comp_defn}\cite{munkres:topology}
		If $\sY$ is a compact Hausdorff space and $\sX$ is a proper subspace of $\sY$ whose
		closure equals $\sY$, then $\sY$ is said to be a \textit{compactification} of $\sX$. If $\sY\setminus\sX$ equals a single
		point, then $\sY$ is called the \textit{one-point compactification} of $\sX$.
	\end{definition}
	\begin{remark}
		It is shown  (cf. \cite{munkres:topology}) that $\sX$ has a one-point compactification $\sY$ if and only if $\sX$ is
		a locally compact Hausdorff space that is not itself compact. We speak of $\sY$ as "the"
		one-point compactification because $\sY$ is uniquely determined up to a homeomorphism.
	\end{remark}
	\begin{theorem}\label{sc_comp_thm}\cite{munkres:topology}
		Let $\sX$ be a completely regular space. There exists a  
		compactification $\sY$ of $\sX$ having the property that every bounded continuous map $f : \sX\to\R$
		extends uniquely to a continuous map of $\sY$ into $\R$.
	\end{theorem}
	
	\begin{definition}\label{sc_comp_defn}\cite{munkres:topology}
		For each completely regular space $\sX$, let us choose, once and for all,
		a compactification of $\sX$ satisfying the extension condition of Theorem \ref{sc_comp_thm}. We will
		denote this compactification of  by $\bt\sX$ and call it the \textit{Stone-\v{C}ech compactification}  
		of $\sX$. It is characterized by the fact that any continuous map $\sX\to\sY$ of $\sX$ into a
		compact Hausdorff space $\sY$ extends uniquely to a continuous map $\bt\sX\to\sY$.
	\end{definition}
\begin{theorem}\label{top_sc_to_compact_thm}\cite{munkres:topology}
Let $\sX$ be a completely regular space; let $\bt\sX$ be the  {Stone-\v{C}ech compactification} of $\sX$. Given any continuous map $f: \sX \to \sY$ into a compact Hausdorff space $\sY$, the map $f$ extends uniquely to a
continuous map $\bt f: \bt \sX\to\sY$.
\end{theorem}	
	
	\begin{definition}\label{f_topology_defn}\cite{rudin:fa}
		Suppose next that $\sX$ is a set and $\mathscr F$ is a nonempty family of mappings $f: \sX \to \sY_f$, where each $\sY_f$ is a topological space. (In many important cases, $\sY_f$ is the same for all $f \in \mathscr F$.) Let $\tau$ be the collection of all unions of finite intersections of sets $f^{-1}\left(\sV  \right)$ with $f\in \mathscr F$ and $\sV$ open in $\sY_f$. Then $\tau$ is a topology on $\sX$, and it is in fact the weakest topology on $\sX$ that makes every 
		$f\in \mathscr F$ continuous: If $\tau'$ is any other topology with that property, then $\tau\subset\tau'$. This $\tau$ is called the \textit{weak topology on} $\sX$ \textit{induced by} $\mathscr F$, or, more 
		shortly, the $\mathscr F$-\textit{topology of} $\sX$.
	\end{definition}
	\begin{definition}\label{top_locally_contractible_defn}\cite{spanier:at,switzer:at}
		A topological space $\sX$ is said to be 
		\textit{contractible} if the identity map of $\sX$ is homotopic to some constant map of $\sX$ 
		to itself. $\sX$ is \textit{locally contractible} if
		every neighborhood $\sU$ of a point $x$ contains a neighborhood $\sV$ of $x$ deformable to $x$ in $\sV$). 
		
	\end{definition} 
\begin{definition}\label{top_simply_conn_defn}\cite{spanier:at}
Let $E^{k+1}, S^k\subset\R^{k+1}$ be a disk ans a sphere, i.e.
\bean
E^{k+1} \bydef \left\{\left.\left(x_1, ..., x_{k+1}\right) \in \R^{k+1}~\right| x_1^2+...+x^2_{k+1}\le 1\right\},\\
S^k \bydef \left\{\left.\left(x_1, ..., x_{k+1}\right) \in \R^{k+1}~\right| x_1^2+...+x^2_{k+1}= 1\right\}
\eean
with the natural inclusion $S^k\subset E^{k+1}$. A space $\sX$ is said to be $n$-\textit{connected} if every continuous map $f: S^k \to \sX$ for $k\le n$ has a continuous extension over $E^{k+1}$. A 1-connected space is also said to be \textit{simply connected}.
\end{definition}
	
	\begin{theorem}\label{it_lim_thm}\cite{kelley:gt}
	\textit{Theorem on iterated limits}. Let $D$ be a directed set, let $E_m$ be a directed set for each $m$ in $D$, let $F$ be the product $D\times \prod_{m \in D}  E_m$, and for $(m,f)$ in $F$ let $R(m,f) = (m, f(m))$. If $S(m,n)$ is 
	a member of a topological space for each m in $D$ and each $n$ in $E_n$, 
	then $S\circ R$ converges to $\lim_m \lim_nS(m,n)$ whenever this iterated limit 
	exists.
\end{theorem}
\begin{definition}\label{top_lower_semi_defn}\cite{bourbaki_sp:gt}
A real-valued function $f$, defined on a topological space $\sX$, is said to be \textit{lower semi-continuous} (resp. \textit{upper semi-continuous} at point $a \in \sX$ if for each $h < f\left(a\right)$ (resp. each $k > f\left(a\right))$ there is an open neighborhood $\sV$ of $a$ such that $h < f\left(x\right)$ (resp. each $k > f\left(x\right))$ for each $x\in\sV$.\\
A real-valued function $f$, is said to be \textit{lower semi-continuous} (resp. \textit{upper semi-continuous} on $\sX$ if it is  {lower semi-continuous} (resp. {upper semi-continuous} at every point of $\sX$.
\end{definition}
\begin{theorem}\label{top_lower_semi_thm}\cite{bourbaki_sp:gt}
	Let $f$ be a lower semi-continuous function on a non-empty quasi-compact space $\sX$. Then there is at least one point $a \in \sX$ such that $f\left(a \right) = \inf_{x \in \sX}f\left( x\right) $ (in other words, $f$ attains its greatest lower bound in $\sX$).
\end{theorem}
\begin{remark}\label{top_lower_semi_rem}
There is inverted result about upper semi-continuous functions and supremums.
\end{remark}

\begin{empt}\label{top_cauchy_empt}\cite{engelking:general_topology}
If $\left( \sX, \rho\right) $ is a metric space then wa say that $\left\{x_j\right\}_{j \in \N}\subset \sX$ is a \textit{Cauchy sequence} if for every $\eps > 0$ there exist $k > 0$ such that $\rho\left(x_j, x_k \right) < \eps$ whenever $j> k$. If $\left( \sX, \rho\right)$ is complete then any Cauchy sequence is convergent.
\end{empt}
\begin{remark}\label{top_cauchy_rem}
Above statements can be generalized to nets (cf. Definition \ref{top_net_lim_defn}) the word \textit{Cauchy net} is also used in this book.
\end{remark}


\begin{definition}\label{metric_space_defn}\cite{rudin:pa}
	A set $X$, whose elements we shall call \textit{points}, is said to be a 
	\textit{metric space} if with any two points $p$ and $p$ of $X$ there is associated a real 
	number $(p, q)$, called the \textit{distance} from $p$ to $q$, such that
	\begin{enumerate}
		\item [(a)]
		\bean
		p \neq q \quad \Leftrightarrow \quad d\left(p, q\right)> 0;\\
		d(p,p)=0;
		\eean
		\item[(b)] $d(p,q)=d(q,p);$
		\item[(c)] $\forall r \in X \quad d(p, q) < d(p, r) + d(r, q)$.
	\end{enumerate} 
	Any function with these three properties is called a \textit{distance function}, or a \textit{metric}. 
\end{definition}

\begin{remark}\label{triangle_ineq_rem}
	The condition (b) of the Definition \ref{metric_space_defn} is said to be the \textit{triangle inequality}. If $Y$ is a normed vector space then there is a distance function on $Y$ such that
	$$
	d(x, y)\bydef \left\|x - y \right\| \quad \forall x, y \in Y.
	$$
	The triangle identity is given by
\be\label{triangle_ineq_vec_eqn}
	\forall x, y, z \in Y \quad  \left\|x - y \right\|\le  \left\|x - z \right\|+ \left\|y - z \right\|.
\ee
\end{remark}
\begin{definition}\label{top_order_top_defn}\cite{counter_topology}
	Let $\sX$ be a set which is linearly ordered by the transitive relation "<". We define the \textit{order (or interval) topology} on $\sX$ by taking as a basis the open intervals
	\be\label{top_order_top_eqn}
	\left(x,y\right)\bydef \left\{x \in \sX| y < x < z\right\}
	\ee 
\end{definition}
\begin{empt}\label{top_lex_square_empt}\cite{counter_topology,munkres:topology}
	Let $I^2_o = \left[0,1\right]^2$ be a unit square. We order $I^2_o$ lexicographically
	$$
	\forall (x, y),(u, w) \in \left[0,1\right]^2 = I^2_o \quad (x, y)< (u, w)  = \begin{cases}
		y < w & x = u\\
		x < u & \text{otherwise}
	\end{cases}
	$$
	and place the order topology on $I^2_o$. It is proven in  \cite{counter_topology,munkres:topology} $I^2_o$ is connected but it is not path connected.
	The space $I^2_o$ is said to be the \textit{the ordered square} (cf. \cite{munkres:topology}).
\end{empt}
	\section{Vector bundles}\label{top_vb_sub_sub}
\paragraph*{}
Let $k$ be the field of real or complex numbers, and let $\sX$ be a topological space.
\begin{definition}\label{top_vb_fiber_defn}\cite{karoubi:k}
	A \textit{quasi-vector bundle with base} $\mathcal X$ is given by:
	\begin{enumerate}
		\item [(a)] A finite dimensional $k$-vector space $E_x$ for every point $x$ of $\mathcal X$.
		\item[(b)] A topology on the disjoint union $E = \bigsqcup E_x$ which induces the natural topology on each $E_x$, such that the obvious projection $\pi: E \to  \mathcal X$ is continuous.
	\end{enumerate}
	The quasi-vector bundle with base will be denoted by $\xi = \left( E, \pi, \mathcal X\right)$. The space $E$ is the \textit{total space} of $\xi$ and $E_x$ is the \textit{fiber} of $\xi$ at the point $x$.
\end{definition}
\begin{empt}\label{trivial_vb_empt}\cite{karoubi:k}	Let $V$ be a finite dimensional vector space over $k$,  $E_x = V$ and the total space may be identified with $\sX \times V$ with the product topology then the quasi-vector bundle  $\left(\sX \times V, \pi, \mathcal X\right)$ is called a \textit{trivial vector bundle}.
\end{empt}
\begin{empt}\cite{karoubi:k}	Let $\xi = \left( E, \pi, \mathcal X\right)$ be a quasi-vector bundle, and let let $\sX' \subset \sX$ be a subspace of $\sX$. The triple $\xi' = \left(\pi^{-1}\left( \sX'\right) , \pi|_{\pi^{-1}\left( \sX'\right)}, \mathcal X'\right)$ is called the \textit{restriction} of $\xi$ to $\sX'$. The fibers of $\xi'$ are just fibers of $\xi$ over the subspace $\xi$. One has
	\be
	\sX'' \subset \sX' \subset \sX \quad \Rightarrow\quad \left.\left(\xi|_{\sX'} \right) \right|_{\sX''}= \xi|_{\sX''}.
	\ee
\end{empt}
\begin{definition}\label{top_vb_defn}\cite{karoubi:k}
	Let $\xi = \left( E, \pi, \mathcal X\right)$ be quasi-vector bundle. Then $\xi$ is said to be \textit{locally trivial} or a \textit{vector bundle} if for every point $x$ in $\sX$, there exists a neighborhood of $x$ such that $\xi|_{\sU}$ is isomorphic to a trivial bundle.
\end{definition}
\begin{definition}\label{top_vb_cs_defn}\cite{karoubi:k}	
	 Let $\xi \bydef \left( E, \pi, \mathcal X\right)$ be a vector bundle. Then a \textit{section} of $\xi$ is a map $s: \sX \to E$ such that $\pi \times s = \Id_{\sX}$. A section $s$ is called \textit{continuous} if $s$ is a continuous. The space $\Ga\left({\mathcal X}, {E}\right)$, of continuous sections can be regarded as both left  and right $C_b\left(\mathcal{X} \right)$-module.
\end{definition}
\begin{notation}\label{top_sec_notn}
	We denote by $\Ga\left( M, E\right)$ the $\C$-space of continuous sections. Both  $\Ga_0\left( M, E\right)$ and $\Ga_c\left( M, E\right)$ are subspaces of  sections tending to zero at infinity and having compact support.
\end{notation}
\begin{remark}
	In \cite{karoubi:k} the vector bundles over fields $\R$ and $\C$ are considered. Here we consider complex vector bundles only.
\end{remark}
\begin{definition}\label{vb_inv_img_funct_defn}\cite{karoubi:k}
	Let $f: \mathcal X' \to \mathcal X$ be a continuous map. For every point $x'$ of $\mathcal X'$, let $E'_{x'}= E_{f\left(x' \right) }$. Then the set $E' = \bigsqcup_{x' \in \mathcal X'}E'_{x'}$ may be identified with the \textit{fiber product} $\mathcal X' \times_{\mathcal X} E$ formed by the pairs $\left(x',e \right)$ such that $f\left(x' \right) = \pi\left( e\right)$. If  $\pi': E' \to \mathcal X'$ is defined by $\pi'\left(x',e \right) = x'$, it is clear that the triple $\xi = \left( E', \pi', \mathcal X'\right)$ defines a quasi-vector bundle over $ \mathcal X'$, when we provide $E'$ with the topology induced by $ \mathcal X' \times E$. We write $\xi' = f^*\left(\xi \right)$ or $f^*\left(E \right)$: this is the \textit{inverse image} of $\xi$ by $f$.  
\end{definition}
\begin{theorem}\label{serre_swan_thm}\cite{karoubi:k}.
	Theorem (Serre, Swan). Let $A = C_k(\sX)$ be the ring of continuous 
	functions on a compact space $\sX$ with values in $k$. Then the section functor $\Ga$ (cf. Notation \ref{top_sec_notn}) induces an 
	equivalence of categories of vector bundles over $\sX$ and finitely generated projecive $A$-modules.
\end{theorem}

\begin{definition}\label{top_herm_bundle_form_defn}\cite{karoubi:k}
Let $E$ be a vector bundle over $\C$. A \textit{sesquilinear form} on $E$ is a 
continuous map $\varphi: E \times_\sX E \to \C$ which has the following property. The map 
$\varphi_x : E_x \times E_x$ induced on each fiber is "sesquilinear" with respect to the $\C$-vector space structure of $E_x$. In other words, ($\varphi_x$ is $\R$-bilinear and $\varphi_x\left(\la e, e' \right)  = \varphi_x\left( e, \overline \la e' \right) = \la \varphi_x\left( e, e' \right)$ f
for $\la\in \C$, $e \in E_x$, and $e' \in E_x$. 
\end{definition}


\begin{remark}\label{top_herm_bundle_rem}
	Let $\mathcal X$ be a locally compact topological space and $E$ the complex vector bundle on $\mathcal X$ with a {sesquilinear form} $\varphi: E \times_\sX E \to \C$ then there are the following pairings
	\bea\label{top_gg_eqn}
	\left\langle \cdot, \cdot \right\rangle:  \Ga_0\left( \sX, E\right)\times \Ga_0\left( \sX, E\right) \to C_0\left( \sX\right) \quad \left(\xi, \eta \right) \mapsto \left(x \mapsto \varphi_x\left(\xi_x, \eta_x \right)  \right) ,\\
	\label{top_ggc_eqn}
	\left\langle \cdot, \cdot \right\rangle_c:  \Ga_c\left( \sX, E\right)\times \Ga_c\left( \sX, E\right) \to C_c\left( \sX\right) \quad \left(\xi, \eta \right) \mapsto \left(x \mapsto \varphi_x\left(\xi_x, \eta_x \right)  \right) 
	\eea
	(cf. Notation \ref{top_sec_notn}).
\end{remark}

\section{Algebraic topology}
\subsection{$H$-groups}
\paragraph{}
Here I follow to \cite{spanier:at}. If both $\sX$ and $\sY$ are topological spaces then we denote by $\left[\sX; \sY\right]$ the set of homotopy classes of continuous maps from $\sX$ to $\sY$. One method of obtaining a group structure on  $\left[\sX; \sY\right]$ is to start with a 
group structure on $\sY$. Thus, let $\sY$ be a topological group with identity element as base point. There is a law of composition in the set of all base-point-preserving continuous maps from $\sX$ to $\sY$ defined by point-wise multiplication 
of functions. That is, if $g_1, g_2: \sX \to\sY$, then $g_1g_2: \sX \to\sY$ is defined by $g_1g_2\left(x\right)\bydef g_1\left(x\right)g_2\left(x\right)$, where the right-hand side is the group product in $\sY$. With 
this law of composition, the set of base-point-preserving continuous maps 
from $\sX$ to $\sY$ is a group (which is Abelian if $\sY$ is Abelian). The law of  
composition carries over to give an operation on homotopy classes such that $\left[g_1\right]\left[g_2\right]\bydef \left[g_1g_2\right]$, and we have the following theorem. 
\begin{theorem}\label{top_hgr_thm}\cite{spanier:at}
If $\sY$ is a topological group, there is a contravariant functor from 
the homotopy category of  pointed  topological spaces to the category of groups 
and homomorphisms. 
\end{theorem}
\begin{remark}\label{top_hgr_rem}
For any topological group $\sY$ we denote by $\left[~\cdot~;\sY\right]$ the given by the Theorem \ref{top_hgr_thm} functor.
\end{remark}

	\subsection{Loop space and suspension in topology}


\begin{empt}\label{top_comp_open_empt} {\it Function spaces.}\cite{switzer:at} If $\sX$ and $\sY$ are topological spaces, we let $\sY^\sX$ denote the set of all continuous functions functions $f: \sX \to \sY$. We give this set a topology, called the {\it compact-open topology}, by taking as a sub-base for the topology all sets of the form $N_{K,U} =\{f:f(K)\subset \sU\}$, $K\subset \sU$ compact, $\sU \subset \sX$ open. 
\end{empt}
\begin{empt}\label{base_points}\cite{switzer:at} {\it Base points.} \cite{switzer:at} Algebraic topology often have to consider not just topological space $\sX$ but rather a space $\sX$ together with a distinguished point $x_0\in \sX$ called the {\it base point}. The pair $(\sX, x_0)$ is called a {\it  pointed space} (one also speaks of  pointed sets). When we are concerned with  pointed spaces $(\sX, x_0)$, $(\sY, y_0)$, etc. we always require that all functions $f: \sX \to \sY$ shall preserve base point, i.e. $f(x_0)=y_0$, and that all homotopies $F:\sX\times I \to \sY$ be relative to base point, i.e. $F(x_0, t) = y_0$, $\forall t \in I$ unless an explicit disclaimer to be contrary is made. We shall use the notation $[\sX, x_0, \sY, y_0]$ to denote the homotopy classes of base point preserving functions, where homotopies are $\mathrm{rel}x_0$, of course. $[\sX, x_0, \sY, y_0]$ is a  pointed space with base point $f_0$ the constant function: $f_0(x)=y_0$, $\forall x\in \sX$. If $(\sX, x_0)$, $(\sY, y_0)$ are  pointed spaces then we have the space $(\sY, y_0)^{(\sX, x_0)}$ of base point preserving functions. We use the notation $\sX \vee \sY$ for the subspace $\sZ \times \{x_0\} \cup \{z_0\} \sX$. It can be thought of as the result of taking the disjoint union $\sZ \cup \sX$ and identifying $z_0$ with $x_0$, $\sZ \vee \sX$ is again a  pointed space with base point $(z_0, x_0)$.  Given maps $f: (\sX, x_0)\to (\sX', x'_0)$ and $g: (\sY, y_0)\to (\sY', y'_0)$ the map $f \times g : \sX\times \sY \to \sX'\times \sY'$ maps $\sX \vee \sY$ into $\sX' \vee \sY'$ and so induces a map $f \wedge g: \sX\wedge \sY \to \sX'\wedge \sY'$.
\end{empt}
\begin{empt}\label{smash_pointed_sp}
	If both $(\sX, x_0)$, $\left( \sY, y_0\right) $ are  pointed spaces then we define the {\it smash product} $(\sX \wedge \sY, *)$ to be the quotient 
	\begin{equation}\nonumber
		\sX \wedge \sY = \sX \times \sY / \sX \vee \sY
	\end{equation}
	with the point $* = p(\sX\vee \sY)$ as base point. Here $p: \sX \times \sY \to X \vee Y$ is the projection. For any $(x,y)\in \sX \times \sY$ we denote $p(x,y)\in \sX \wedge \sY$ by $[x,y]$.
\end{empt}
\begin{thm}\cite{switzer:at}
	If $(\sX, x_0), (\sY, y_0), (\sZ, z_0)$ are  pointed spaces, $\sX, \sZ$ Hausdorff and $\sZ$ locally compact, then there is a natural equivalence
	\begin{equation}
		A: [\sZ \wedge \sX, *; \sY, y_0] \to [\sX, x_0, (\sY, y_0)^{(\sZ, z_0)}, f_0]
	\end{equation}
	defined by $A[f]=[\hat f]$, where if $f: \sZ \wedge \sX \to \sY$ is a map then $\hat f: \sZ \wedge \sX \to \sY$ is a map given by $(\hat f(x)(z)=f[z,x]$.
\end{thm}
\begin{empt}\label{topological_loop_spaces} \cite{switzer:at} If $(\sY, y_0)$ is a  pointed space, we define the {\it loop space} $(\Omega \sY, \omega_0)$ of $\sY$ to be function space 
	\begin{equation}\nonumber
		(\Omega \sY, \omega_0) = (\sY, y_0)^{(S^1, s_0)}
	\end{equation}
	with constant loop $\omega_0$ $(\omega_0(s))= y_0, \ \forall s \in S^1$ as base point.
\end{empt}
\begin{empt}\label{topological_suspension}\cite{switzer:at}
	If $(\sX, x_0)$ we define {\it suspension} $(\Sigma \sX, *)$ to be the smash product $(S^1 \wedge \sX, *)$ of $\sX$ with the 1-sphere.
\end{empt}
\begin{cor}\label{loop_susp_adj_cor}\cite{switzer:at}		If both $(\sX, x_0), (\sY, y_0)$ are  pointed spaces and $\sX$ is Hausdorff, then there is a natural equivalence
	\begin{equation}\nonumber
		A:[\Sigma \sX,*;\sY,y_0]\cong [\sX, x_0, \Omega \sY, y_0].
	\end{equation}
\end{cor}
\begin{prop}\label{loop_susp_adj_prop}\cite{switzer:at}
	The adjoint correspondence
	\begin{equation}\nonumber
		A:[\Sigma \sX,*;\sY, y_0]\cong [\sX,x_0, \Omega \sY, \omega_0]
	\end{equation}
	is an isomorphism of groups.
\end{prop}
\begin{remark}\cite{switzer:at}
Both sets $[\Sigma \sX,*;\sY, y_0]$ and $[\sX,x_0, \Omega \sY, \omega_0]$ posses the natural group structure.
\end{remark}

\begin{definition}\label{top_hn_bas_defn}\cite{switzer:at}
	For $n\ge1$  we define  $n^{\text{th}}$ \textit{homotopy group} $\pi_n\left(\sX, x_0 \right)$ of $\left(\sX, x_0 \right)$ as 
	$$
	\pi_n\left( \sX, x_0\right) \bydef \pi_0\left( \Om^n\sX, \om_0\right) \cong \pi_1\left( \Om^{n -1}\sX, \om'_0\right).
	$$
	where $\pi_0\left( \Om^n\sX\right) $ is a set of connected components of $\Om^n\sX$,
\end{definition}
\begin{remark}\label{top_homotopy_group_rem}\cite{switzer:at}
The group $\pi_1\left(\sX, x_0 \right)$ is a set of homotopy equivalence classes of closed paths $f: \left[0, 1\right] \to\sX$ such that $f\left(0\right)= f\left(1\right)= x_0$. For all $n \in \N$ there is the natural bijective map $\pi_n\left(\sX, x_0 \right)\cong \left[\sX, x_0; S^n, s_0\right]$. We say that $\pi_1\left(\sX, x_0 \right)$ is the \textit{fundamental group} of the  pointed space $\left(\sX, x_0 \right)$. 
\end{remark}

	\begin{defn}\label{top_semi1_defn}\cite{spanier:at}
	A space $\sX$ is said to by \textit{semilocally} 1-\textit{connected} if every point $x_0$ has a neighborhood $N$ such that $\pi_1\left( N, x_0\right) \to \pi_1\left(\sX, x_0 \right)$ is trivial. 
\end{defn} 

\begin{remark}\cite{switzer:at}
For all $n > 1$ the group $	\pi_{n} (\sX, x_0)$ is Abelian.
\end{remark}
\begin{remark}\label{top_ho_fund_rem}\cite{switzer:at}
	From the Corollary \ref{loop_susp_adj_cor} and the Proposition \ref{loop_susp_adj_prop} it follows that there is the natural isomorphism of Abelian groups
	\be\label{top_ho_fund_eqn}
	\pi_{n+m} (\sX, x_0) = \pi_n(\Omega^m(\sX, \om_m)). 
	\ee
\end{remark}

\subsection{Coverings and fundamental groups}
	\paragraph*{}
	Results of this section are copied from \cite{spanier:at}. The \textit{covering projection} word is replaced with \textit{covering} outside this section.
	
	\begin{defn}\label{top_covering_defn}\cite{spanier:at}
		Let $\widetilde{\pi}: \widetilde{\mathcal{X}} \to \mathcal{X}$ be a continuous map. An open subset $\mathcal{U} \subset \mathcal{X}$ is said to be {\it evenly covered } by $\widetilde{\pi}$ if $\widetilde{\pi}^{-1}(\mathcal U)$ is the disjoint union of open subsets of $\widetilde{\mathcal{X}}$ each of which is mapped homeomorphically onto $\mathcal{U}$ by $\widetilde{\pi}$. A continuous map $\widetilde{\pi}: \widetilde{\mathcal{X}} \to \mathcal{X}$ is called a {\it covering projection} if each point $x \in \mathcal{X}$ has an open neighborhood evenly covered by $\widetilde{\pi}$. $\widetilde{\mathcal{X}}$ is called the {
			\it covering space} and $\mathcal{X}$ the {\it base space} of the covering.
	\end{defn}
	\begin{defn}\cite{spanier:at}
		A topological space $\mathcal X$ is said to be \textit{locally path-connected} if the path components of open sets are open.
	\end{defn}
	\subsubsection{Group of covering transformations}
	\begin{defn}\label{top_group_of_covering_transformations_defn}\cite{spanier:at}
		Let $p: \mathcal{\widetilde{X}} \to \mathcal{X}$ be a covering.  A self-equivalence is a homeomorphism $f:\mathcal{\widetilde{X}}\to\mathcal{\widetilde{X}}$ such that $p \circ f = p$. This group of such homeomorphisms is said to be the {\it group of covering transformations} of $p$ or the {\it covering group}. Denote by $G\left( \left.\mathcal{\widetilde{X}}~\right|~\mathcal{X}\right)$ this group.
	\end{defn}
	\begin{definition}\label{top_properly_disc_group_defn}\cite{spanier:at} 
A group $G$ of homeomorphisms of a topological space $\sY$  is said to be  
\textit{discontinuous} if the orbits of $G$ in $\sY$ are discrete subsets of $\sY$. $G$ is \textit{properly discontinuous} if for $y \in \sY$ there is an open neighborhood $\sU$ of $y$ in $\sY$ such 
that if $g, g' \in  G$ and $g\sU$ meets $g'\sU$, then $g = g'$. $G$ \textit{acts without fixed points} if the only element of G having fixed points is the identity element.	
\end{definition}
\begin{theorem}\label{top_group_of_covering_transformations_thm}\cite{spanier:at}
	Let $G$ a properly discontinuous group of homeomorphisms 
	of a space $\sY$. Then the projection of $\sY$ to the orbit space $\sY/G$ is a covering projection. If $\sY$ is connected, this covering projection is regular and $G$  is its	group of covering transformations. 
\end{theorem}
		\subsubsection{Unique path lifting}
\begin{theorem}\label{top_fundamental_group_mor_thm}\cite{spanier:at}
	Let $p: \widetilde{\sX}\to \sX$ a fibration with unique path lifting. Let $\sX$ be path connected and let $x_0\in \sX$. Then  there is  a monomorphism $\psi$ of $G\left( \left.\widetilde{\sX} \right|\sX\right) $ to 
	the quotient group $N\left(\pi_1\left( p\right) \left( \pi_1\left(\widetilde\sX,\widetilde x_0 \right)\right)\right)   /\pi_1\left(\widetilde\sX,\widetilde x_0 \right)$ where $N\left(\pi_1\left( p\right) \left( \pi_1\left(\widetilde\sX,\widetilde x_0 \right)\right)\right)$ is a maximal among subgroups $G\subset \pi_1\left(\sX, x_0 \right)$ such that $\pi_1\left(\widetilde\sX,\widetilde x_0 \right)$ is a normal subgroup of $G$.	If $\sX$ is also locally path connected, 
	$\psi$ is an isomorphism. 
\end{theorem}

	
	\begin{theorem}\label{top_path_lift_thm}\cite{spanier:at}
		A fibration has unique path lifting if and only is every fiber has no nonconstant paths. 
	\end{theorem}

	\begin{theorem}\label{top_locally_conn_cov_thm}\cite{spanier:at}
		If $\mathcal{ X}$ is locally connected , a continuous map  $p: \widetilde{\mathcal X}\to {\mathcal X}$ is a covering projection if and only if for each component $\mathcal Y$ of $\mathcal X$ the map
		$$
		p|_{p^{-1}\left(\mathcal Y \right) }:p^{-1}\left(\mathcal Y \right) \to \mathcal Y
		$$
		is a covering projection.
	\end{theorem}	
	%
	\begin{corollary}\label{top_cov_cat_cor}\cite{spanier:at}
		Consider a commutative triangle
		\newline
		\begin{tikzpicture}
			\matrix (m) [matrix of math nodes,row sep=3em,column sep=4em,minimum width=2em]
			{
				\widetilde{\mathcal X}_1 & &\widetilde{\mathcal X}_2\\ 
				& {\mathcal X}\\};
			\path[-stealth]
			(m-1-1) edge node [above] {$p$} (m-1-3)
			(m-1-1) edge node [left]  {$p_1~~$} (m-2-2)
			(m-1-3) edge node [right] {$~~p_2$} (m-2-2);
		\end{tikzpicture}
		\\
		where $\sX$ is locally connected (cf. Definition \ref{top_locally_connected_defn}) and $p_1$, $p_2$ are covering projections. If $p$ is a surjection then $p$ is a covering projection.
	\end{corollary}
	\begin{theorem}\label{top_locally_conn_cov_com_thm}\cite{spanier:at}
		If   $p: \widetilde{\mathcal X}\to {\mathcal X}$ is a covering projection onto locally connected base space, then for any component  $\widetilde{\mathcal Y}$ of $\widetilde{\mathcal X}$ the map
		$$
		p|_{\widetilde{\mathcal Y}}:\widetilde{\mathcal Y} \to p\left(\widetilde{\mathcal Y} \right) 
		$$
		is a covering projection onto some component of $\widetilde{\mathcal X}$.
	\end{theorem}
	\begin{theorem}\label{top_conjugate_thm}\cite{spanier:at}
		Let $p_1: 	\widetilde{\sX}_1 \to \sX$, $p_2: 	\widetilde{\sX}_2 \to \sX$ be objects  in the category of connected covering spaces of a connected locally path-connected space $\sX$. The following are equivalent
		\begin{enumerate}
			\item [(a)] There is a coveting projection $f:\widetilde{\sX}_1 \to \widetilde{\sX}_2$ such that $p_2 \circ f = p_1$.
			\item[(b)] For all $\widetilde{x}_1\in \widetilde{\sX}_1$ and  $\widetilde{x}_2\in \widetilde{\sX}_2$ such that $p_1\left(\widetilde{x}_1 \right) = p_2\left(\widetilde{x}_2 \right)$, $\pi_1\left(p_1\right)\left( \pi_1 \left( \widetilde{\sX}_1, \widetilde{x}_1\right)  \right)$ is conjugate to a subgroup of  $\pi_1\left(p_2\right)\left( \pi_1 \left( \widetilde{\sX}_2, \widetilde{x}_2\right)  \right)$.
			\item[(c)] There exist $\widetilde{x}_1\in \widetilde{\sX}_1$ and  $\widetilde{x}_2\in \widetilde{\sX}_2$ such that $\pi_1\left(p_1\right)\left( \pi_1 \left( \widetilde{\sX}_1, \widetilde{x}_1\right)  \right)$ is conjugate to a subgroup of  $\pi_1\left(p_2\right)\left( \pi_1 \left( \widetilde{\sX}_2, \widetilde{x}_2\right)  \right)$.
		\end{enumerate} 
	\end{theorem}
\begin{definition}\cite{spanier:at}
A continuous map $f: \sX \to \sY$ is called a \textit{local homeomorphism} if each 
point $y \in \sY$ has an open neighborhood mapped homeomorphically by $f$ onto 
an open subset of $\sX$
\end{definition}
\begin{lemma}\label{top_cov_lem}\cite{spanier:at}
	A covering projection is a local homeomorphism.
\end{lemma}

	\begin{lemma}\label{top_om_lem}\cite{spanier:at}
		A local homeomorphism is an open map (cf. Definition \ref{top_onen_map_defn}).
	\end{lemma}

	\begin{remark}\label{top_cov_bicont_rem}
 From the Lemmas \ref{top_om_lem} and \ref{top_om_lem} it follows that any covering is a bicontinuous map (cf. Definition \ref{top_bicont_defn}).
	\end{remark}

	\begin{defn}\label{top_path_lifting_defn}\cite{spanier:at}
		A continuous map $p:\mathcal{E} \to \mathcal{B}$ is said to have the {\it unique path lifting} if, given paths $\omega$ and $\omega'$ in $E$ such that $p \circ \omega = p \circ \omega'$ and $\omega(0)=\omega'(0)$, then $\omega=\omega'$.
	\end{defn}
	\begin{thm}\cite{spanier:at}\label{spanier_thm_un}
		Let $p: \widetilde{\mathcal{X}} \to \mathcal{X}$ be a covering projection and let $f, g: \mathcal{Y} \to \widetilde{\mathcal{X}}$ be liftings of the same map (that is, $p \circ f = p \circ g$). If $\mathcal{Y}$ is connected and $f$ agrees with g for some point of $\mathcal{Y}$ then $f=g$.
	\end{thm}
	\begin{rem}
		From theorem \ref{spanier_thm_un} it follows that a covering projection has unique path lifting.
	\end{rem}
	
	\subsubsection{Regular and universal coverings}
	\begin{defn}\label{top_regular_defn}\cite{spanier:at}
		A fibration $p: \mathcal{\widetilde{X}} \to \mathcal{X}$ with unique path lifting is said to be  {\it regular} if, given any closed path $\omega$ in $\mathcal{X}$, either every lifting of $\omega$ is closed or none is closed.
	\end{defn}
	\begin{thm}\label{locally_path_thm}\cite{spanier:at}
		Let $p: \widetilde{\mathcal X} \to \mathcal X$ be a fibration with unique path lifting and assume that a nonempty $\widetilde{\mathcal X}$ is a locally path-connected space. Then $p$ is regular if and only if for some $\widetilde{x}_0 \in  \widetilde{\mathcal X}$, $\pi_1\left(p\right)\pi_1\left(\widetilde{\mathcal X}, \widetilde{x}_0\right)$ is a normal subgroup of $\pi_1\left(\mathcal X, p\left(\widetilde{x}_0\right)\right)$.
	\end{thm}
\begin{definition}\label{top_properly_disc_defn}\cite{spanier:at}
A group $G$ of homeomorphisms of a topological space $\sY$ is said to be  
\textit{discontinuous} if the orbits of $G$ in $\sY$ are discrete subsets of $\sY$. $G$ is \textit{properly 
discontinuous} if for $y \in \sY$ there is an open neighborhood $\sU$ of $y$ in $\sY$ such 
that if $g, g' \in  G$ and $g\sU$ meets $g'\sU$, then $g = g'$. $G$ acts \textit{ without fixed points} 
if the only element of $G$ having fixed points is the identity element. The 
following are clear. 
A properly discontinuous group of homeomorphisms is discontinuous 
and acts without fixed points. 
\end{definition}
\begin{remark}\label{top_properly_disc_rem}\cite{spanier:at}
A finite group of homeomorphisms acting without fixed points on a 
Hausdorff space is properly discontinuous.
\end{remark}
	\begin{theorem}\label{top_cov_fact_thm}\cite{spanier:at}
		Let $G$ be a properly discontinuous group of homeomorphisms of space $\mathcal Y$. Then the projection of $\mathcal Y$ to the orbit space $\mathcal Y/G$ is a covering projection. If $\mathcal Y$ is connected, this covering is regular and $G$ is its group of covering transformations, i.e. $G = G\left(\mathcal Y~|~\mathcal Y/G \right)$. 
	\end{theorem}
	\begin{lem}\label{top_cov_from_pi1_cor}\cite{spanier:at}
		Let $p: \widetilde{\mathcal X} \to \mathcal X$ be a fibration with a unique path lifting. If $ \widetilde{\mathcal X}$ is connected and locally path-connected and $\widetilde{x}_0 \in \widetilde{\mathcal X}$ then $p$ is regular if and only if $G\left(\widetilde{\mathcal X}~|~{\mathcal X} \right)$ transitively acts on each fiber of $p$, in which case 
		$$
		\psi: G\left(\widetilde{\mathcal X}~|~{\mathcal X} \right) \approx \pi_1\left(\mathcal X, p\left( \widetilde{x}_0\right)  \right) / \pi_1\left( p\right)\pi_1\left(\widetilde{\mathcal X}, \widetilde{x}_0 \right)  
		$$
		where $\pi_1$ denotes fundamental groups (cf. Remark \ref{top_homotopy_group_rem}).
	\end{lem}
	\begin{remark}\label{top_sim_con_reg_rem}\cite{spanier:at}
		If $\widetilde{   \sX}$ is simply connected, any fibration $p: \widetilde{   \sX} \to \sX$ is regular, and we also have the next result.
	\end{remark}
	\begin{cor}\label{top_cov_pi1_cor}\cite{spanier:at}
		Let $p: \widetilde{\mathcal X} \to \mathcal X$ be a fibration with a unique path lifting where $ \widetilde{\mathcal X}$ is simply connected locally path-connected and nonempty.  Then the group of self-equivalences of $p$ is isomorphic to the fundamental group of ${\mathcal X}$ (cf. Remark \ref{top_homotopy_group_rem}), i.e. $\pi_1\left( {\mathcal X}\right)\approx G\left(\left.\widetilde{\sX}~\right|\sX\right)$.  
	\end{cor}
	\begin{definition}\label{top_universal_covering_defn}\cite{spanier:at}
		A \textit{universal covering space} of a connected space $\sX$ is an object $p: \widetilde{\sX}\to \sX$ of the category of connected covering spaces of $\sX$ such that for any object $p': \widetilde{\sX}'\to \sX$ of this category there is a morphism 
		\newline
		\begin{tikzpicture}
			\matrix (m) [matrix of math nodes,row sep=3em,column sep=4em,minimum width=2em]
			{
				\widetilde{\sX}  &  & \widetilde{\sX}'\\
				& \sX & \\};
			\path[-stealth]
			(m-1-1) edge node [above] {$f$} (m-1-3)
			(m-1-1) edge node [left]  {$p~~$} (m-2-2)
			(m-1-3) edge node [right] {$~~p'$} (m-2-2);
		\end{tikzpicture}
		\\
		in the category.
	\end{definition}
	
	\begin{lem}\label{top_simply_con_cov_lem}\cite{spanier:at}
		A connected locally path-connected space $\mathcal X$ has a simply connected covering space if and only if  $\mathcal X$ is semilocally 1-connected (cf. Definition \ref{top_semi1_defn}).
	\end{lem}
	
	\begin{lem}\label{top_uni_exist_spa_lem}\cite{spanier:at}
		A simply connected covering space of a connected locally path-connected space $\sX$ is an universal covering space of $\sX$.
	\end{lem}

\begin{definition}\label{top_sus_defn}\cite{hatcher:at}
	For a space $\sX$, the (unreduced) suspension $\Sigma\sX$ is the quotient of $\sX\times [0,1]$ obtained by collapsing $\sX\times\{0\}$ to one point and  $\sX\times\{1\}$ to another.
point.
\end{definition}

\begin{prop}\cite{switzer:at} \label{top_pi1_pi1_prop}
	(cf. Proposition 3.30).
	If $\sX$ is a topological space then there is an action of $\pi_1\left(\sX\right)\times \pi_1\left(\sX\right)\to \pi_1\left(\sX\right)$ given conjugation; i.e. such that for $\ga, \a \in \pi_1\left(\sX\right)$ one has $\ga \cdot \a = \ga \a\ga^{-1}$.
\end{prop}
!!!\ref{pedersen_ideal_defn}

\subsection{Homology theory and Hurewicz homomorphism}

Let $\mathscr {PT}$ be a category of  pointed spaces and let $\mathscr {PT}'$ be a homotopy category, i.e. objects are as in $\mathscr {PT}$, morphisms are homotopy classes of maps. On $\mathscr {PT}'$ we have a suspension functor $S\left(\sX, x_0 \right)\bydef\left( S\sX, *\right)$ and $S\left[f \right]\bydef\left[Sf \right]   \bydef \left[1_{S^1} \wedge f\right]$ 
\begin{definition}\label{top_red_homol_defn}
A \textit{reduced homology theory} $k_*$ on $\mathscr {PT}'$ is a collection of functors $k_n: \mathscr {PT}'\to \mathscr A$ to the category of Abelian groups and natural equivalences $\sigma_n: k_n \to k_{n+1}\circ S$, satisfying
	\begin{enumerate}
			\item [(a)] \textit{Exactness}: for every  pointed pair $\left(\sX, \A, x_0\right)$ with inclusions $\iota:\left( \sX, x_0\right) \to \left( \sX_0\right)$ and $j : \left( \sX\cup C\A, *\right)$ the sequence 
		\bean
	k_n\left(\A, x_0\right)\xrightarrow{k_n\left[\iota\right]} k_n\left( \sX, x_0\right)\xrightarrow{k_n\left[j\right]} \left( \sX\cup C\A, *\right)	
	\eean
		is exact.
	\end{enumerate}

\end{definition}

\begin{proposition}
There are natural transformations $\hat \partial_n $ and for each  pointed pair $\left(\sX, \A, x_0\right)$ a long exact sequence
		\bean
	k_n\left(\A, x_0\right)\xrightarrow{k_n\left[\iota\right]} k_n\left( \sX, x_0\right)\xrightarrow{k_n\left[j\right]} \left( \sX\cup C\A, *\right)	
\xrightarrow{\hat \partial_n\left(\sX, \A, x_0\right)}\\	\xrightarrow{\hat \partial_n\left(\sX, \A, x_0\right)}	k_{n-1}\left(\A, x_0\right)\xrightarrow{k_{n-1}\left[\iota\right]}...
\eean
\end{proposition}  
\begin{definition}\label{top_hurewitz_defn}\cite{switzer:at}
	Let $k_*$ be a reduced homology theory. We define the \textit{Hurewicz homomorphism} 	$h:\pi_n\left(Y, y_0 \right)\to k_n\left(Y \right)$  to be the composite
	\be
	\begin{split}
		\pi_n\left(\sY, y_0 \right)\cong \left[S^n, s_0; \sY, y_0\right]\xrightarrow{k_n}\Hom\left(k_n\left(S^n \right),  k_n\left(\sY \right)\right) \xrightarrow{\Hom\left(\sigma^n, 1 \right) }\\
		\Hom\left(k_n\left(S^0 \right),  k_n\left(\sY \right)\right)\cong k_n\left(\sY \right).
	\end{split}
	\ee
\end{definition}

Any homology theory defines a cohomology theory. In particular it satisfies to the following exactness axiom
\be\label{top_c_exact_eqn}
...\xrightarrow{\delta} k^n\left(\sX, \A \right) \xrightarrow{k^n\left[j\right]} k^n\left(\sX \right) \xrightarrow{k^n\left[i\right]} k^n\left(\A \right) \xrightarrow{\delta} k^{n+1}\left(\sX, \A \right)\to ...
\ee
\begin{theorem}\cite{rotman:ag}
 If $\sX$ is path connected, then the Hurewicz
map $\varphi: \pi_1\left(\sX, x_0 \right) \to H_1\left( \sX\right)$  is a surjection with kernel $\pi_1\left(\sX, x_0 \right)'$, the commutator
subgroup of $\pi_1\left(\sX, x_0 \right)$. Hence
$$
\pi_1\left(\sX, x_0 \right)/\pi_1\left(\sX, x_0 \right)'\cong H_1\left( \sX\right).
$$
\end{theorem}
We shall use the following equivalent formulation of the above theorem.
\begin{theorem}\label{hurewicz_iso_thm}
	If $H_*$ is the theory of singular homology then the following Hurewicz homomorphism
	$$
	h_{\mathrm{sing}}:\pi_1\left( \sX\right) / \left[\pi_1\left( \sX\right), \pi_1\left( \sX\right)\right]\to H_1\left(\sX \right) 
	$$
	is an isomorphism.
\end{theorem}
\begin{theorem}\label{top_h_to_k_thm}\cite{switzer:at}
For every topological pair $\left(\sX, \A\right)$ and Abelian group $G$ there is a natural exact sequence
\bean
0 \to \mathrm{\Ext}\left(H_{n-1}\left(\sX, \A\right),G \right)\to H^{n}\left(\sX, \A;G \right)\xrightarrow{\mu^*}\\\mathrm{Hom}\left(H_{n}\left(\sX, \A \right), G\right) \to 0
\eean
and unnatural splitting
$$
H^{n}\left(\sX, \A;G \right)\cong \mathrm{\Ext}\left(H_{n-1}\left(\sX, \A \right),G \right)\oplus\mathrm{Hom}\left(H_{n}\left(\sX, \A\right) ,G \right).
$$
\end{theorem}

\subsection{Obstruction theory}
\begin{definition}\label{top_char_defn}\cite{spanier:at}
	Let $\pi$ an Abelian group and $\sY$ is a path-connected  pointed space. An element $v \in H^n\left(\sY, y, \pi \right)$ is said to be $n$-\textit{characteristic} for $\sY$ if the composite
	$$
	\pi_n\left(\sY, y_0 \right)\xrightarrow{\varphi} H^n\left(\sY, y, \pi \right)\xrightarrow{h(v)}\pi
	$$ 
	is isomorphism (where $\varphi$ is the Hurewicz homomorphism and $h\left( v\right)$ is the natural homomorhism which follows from the universal coefficient theorem  \cite{spanier:at}).
\end{definition}
\begin{definition}\label{top_pi_n_defn}\cite{spanier:at}
Let $\pi$ be a group and let $n$ be an integer $\ge 1$. A space of type $\left(\pi, n\right)$ is a path-connected  pointed space $\sY$ such that $\pi_q\left(\sY \right)= 0$  for $q\neq n$ and 
$\pi_n\left(\sY \right)= 0$ is isomorphic to $\pi$.
\end{definition}
\begin{theorem}\label{top_obstr_thm}\cite{spanier:at}
	Let $\sY$ be a space of type $\left(\pi, n\right)$, with $n\ge 1$ and $\pi$ Abelian, and let $\iota$ be $n$-characteristic for $\sY$. If $(\sX,\A)$ is a relative $CW$-complex, a map $f: \A\to \sY$ 
	can be extended over $\sX$ if and only if $\delta H^n\left[f\right]\left( \iota \right)  = 0 \in  H^{n + 1}\left(\sX, \sA; \pi \right)$. 
\end{theorem}
\begin{remark}\label{top_obstr_rem}
A notion of  a (relative) $CW$-complex is explained in \cite{spanier:at,switzer:at}.
\end{remark}
\begin{remark}
	The  Theorem \ref{top_obstr_thm} implies the cohomology exact sequence \eqref{top_c_exact_eqn}, i.e.
	\bean
	...\xrightarrow{\delta} H^n\left(\sX, \A ;\pi\right) \xrightarrow{H^n\left[j\right]} H^n\left(\sX ;\pi\right) \xrightarrow{H^n\left[i\right]} H^n\left(\A; \pi \right) \xrightarrow{\delta} H^{n+1}\left(\sX, \A ;\pi\right)\to ...
	\eean
	and a homomorphism $H^n\left[f\right]: H^n\left(\sY; \pi \right) \to H^n\left(\A; \pi \right)$.
\end{remark}
\chapter{Differential geometry}
	\section{Smooth manifolds and multinormed algebras}\label{smooth_op_system_sec}
	
\begin{definition}\label{diff_mani_defn}\cite{do_carmo:rg}
 A \textit{differentiable manifold} of dimension $n$ is a set 
$M$ and a family of injective mappings $\mathbf{x}_\a: \sU_\a \subset \R^n \to M$ of open 
sets $\sU_\a$ of $\R^n$ into $M$ such that: 
\begin{enumerate}
	\item  $\cup \mathbf{x}_\a\left(\sU_{{\a}} \right) = M$.
		\item for any pair $\a, \bt$ , with $\mathbf{x}_\a\left(\sU_{{\a}} \right) \cap \mathbf{x}_\bt\left(\sU_{{\bt}} \right)= \mathcal{W}\neq \emptyset$, the sets 
	$\mathbf{x}^{-1}_\a\left( \mathcal{W}\right)$  and $\mathbf{x}^{-1}_\bt\left( \mathcal{W}\right)$ are open sets in $\R^n$ and the mappings 
		$\mathbf{x}^{-1}_\bt \circ \mathbf{x}_\a$ are differentiable
				\item  The family $\left\{\left(\mathbf{x}_\a, \sU_\a \right) \right\}$ is maximal relative to the conditions 
				1) and 2).
		
\end{enumerate}
\end{definition}	
\begin{definition}\label{diff_mani_tan_defn}\cite{do_carmo:rg}
Let $M$ be a differentiable manifold. A differentiable function $\a: \left(-\eps, \eps\right)\to M$ is called a (differentiable) \textit{curve} in 
$M$. Suppose that $\a\left(0\right)= p \in M$, and let $\D$ be the set of functions 
on $M$ that are differentiable at $p$. The \textit{tangent vector to the curve} 
$\a$ at $t=0$ is a function $\a'\left(0 \right)$  given by
$$
\a'\left(0 \right)f = \left. \frac{d\left(f \circ \a\right)}{dt}\right|_{t=0}.
$$ 

A \textit{tangent vector at} $p$ is the tangent vector at $t = 0$ of some curve 
 $\a: \left(-\eps, \eps\right)\to M$ with $\a\left(0\right)=p$. The set of all tangent vectors to $M$ 
at $p$ will be indicated by $T_pM$. 

\end{definition}
	\begin{definition}\label{ori_man_defn}\cite{do_carmo:rg}
Let $M$ be a differentiable manifold. We say that 
$M$ is \textit{orientable} if $M$ admits a differentiable structure $\left\{\left(\mathbf{x}_\a, \sU_\a \right) \right\}$ 
such that  for every pair $\a, \bt$ , with $\mathbf{x}_\a\left(\sU_{{\a}} \right) \cap \mathbf{x}_\bt\left(\sU_{{\bt}} \right)= \mathcal{W}\neq \emptyset$, the 
	differential of the change of coordinates $\mathbf{x}^{-1}_\bt \circ \mathbf{x}_\a$ has positive 
	determinant. 
	In the opposite case, we say that $M$ is \textit{unorientable}. 
	\end{definition}

\begin{definition}\label{comm_smooth_usual_defn} \cite{cinfty_manifolds}
	If the  topology of $M$ has a countable  basis,
	and let us consider a countable  family $\left\{K_j\right\}_{j\in \N}$ of compact subsets such that
	$$
	M = \bigcap_{j\in \N} \mathring{K}_j 
	$$
	and each compact set $K_j$ is contained in some coordinate open set $\left(\sU_j, u_1, ..., u_n\right)$. The \textit{usual topology} of $\Coo\left(M \right)$  is defined by the following
	submultiplicative seminorms
	\be\label{comm_smooth_usual_eqn}
	p_{j,k}\left( f\right) \bydef \max 2^k\left|\frac{\partial^{\left|\a\right|}}{\partial u^\a}\right|,
	\ee
	where the maximum is considered over all points $x \in K_j$ and all orders of derivation
	$\a = \left(\a_1, ..., \a_n\right)$ such that $\left|\a\right| = \a_1 + . . . + \a_n \le k$.
	
\end{definition}
It is a basic fact that $\Coo\left(M \right)$ is complete, so that it is a Fr\'echet algebra.
\begin{rem}\label{comm_smooth_usual_rem}
	Under the hypotheses of the Definition \ref{comm_smooth_usual_defn} $\Coo\left(M\right)$  is a local operator system (cf. Remark \ref{la_los_rem}).
\end{rem}

	\begin{prop}\label{top_cov_mani_prop}\cite{kobayashi_nomizu:diff_geom,lee:smooth_manifolds})
	If $M$ is a manifold, and $p:\widetilde M\to M$ is a covering then $\widetilde M$ has a (unique) structure of manifold such that $p$ is differentiable.
\end{prop}
\section{Riemannian manifolds}
\begin{definition}\label{riemann_mani_defn}\cite{do_carmo:rg}
	A \textit{Riemannian metric} (or \textit{Riemannian structure}) 
	on a differentiable manifold $M$ is a correspondence which associates 
	to each point $p$ of $M$ an inner product $\left\langle \cdot, \cdot\right\rangle_p$ (that is, a symmetric, 
	bilinear, positive-definite form) on the tangent space $T_pM$, which 
	varies differentiably in the following sense: If $\mathbf{x}: \sU \cong \R^n \to M$ 
	is a system of coordinates around $p$, with $\mathbf{x}\left(x_1,..., x_n\right) = q \in  \mathbf{x}\left(\sU\right)$
	and $\frac{\partial q}{\partial x_j}= d\mathbf{x}\left( 0,..., 1,...,0\right) $, 
	then
	\be\label{riemann_mani_eqn}
	\left\langle \frac{\partial}{\partial x_j} q, \frac{\partial}{\partial x_k} q\right\rangle_p= g_{jk}\left(x_1,..., x_n \right) 
	\ee
	is a differentiable function on $\sU$. 	It is clear this definition does not depend on the choice of 
	coordinate system. 
	A differentiable manifold with a given Riemannian metric will be called a \textit{Riemannian manifold}. 
	
\end{definition}
\begin{remark}\label{riemann_mani_rem}
	The Riemannian metric is completely defined by the \textit{metric tensor} $\left[g_{jk}\right]$ (cf. equation \eqref{riemann_mani_eqn})
\end{remark}

\begin{prop}\label{comm_cov_mani_prop}(Proposition 5.9 \cite{kobayashi_nomizu:diff_geom})
	\begin{enumerate}
		\item[(i)] Given a connected manifold $M$ there is a unique (unique up to isomorphism) universal covering manifold, which will be denoted by $\widetilde{M}$.
		\item[(ii)] The universal covering manifold $\widetilde{M}$ is a principal fibre bundle over $M$ with group $\pi_1(M)$ and projection $p: \widetilde{M} \to M$, where $\pi_1(M)$ is the first homotopy group of $M$.
		\item[(iii)] The isomorphism classes of covering spaces over $M$ are in 1:1 correspondence with the conjugate classes of subgroups of $\pi_1(M)$. The correspondence is given as follows. To each subgroup $H$ of $\pi_1(M)$, we associate $E=\widetilde{M}/H$. Then the covering manifold $E$ corresponding to $H$ is a fibre bundle over $M$ with fibre $\pi_1(M)/H$ associated with the principal bundle  $\widetilde{M}(M, \pi_1(M))$. If $H$ is a normal subgroup of $\pi_1(M)$, $E=\widetilde{M}/H$ is a principal fibre bundle with group $\pi_1(M)/H$ and is called a regular covering manifold of $M$.
	\end{enumerate}
\end{prop}

\begin{empt}\label{covering_metric_empt}
	If $\widetilde{M}$ is a covering space of a Riemannian manifold $M$ then it is possible to give $\widetilde{M}$ a Riemannian structure such that $\pi: \widetilde{M} \to M$ is a local isometry (this metric is called the {\it covering metric}) cf. \cite{do_carmo:rg} for details.
\end{empt}

\section{Smooth vector bundles}

\begin{definition}\label{top_sm_bundle_defn}\cite{lee:smooth_manifolds}
	Let $M$ be a smooth manifold.
A \textit{smooth vector bundle of rank $k$ over $M$} is a smooth manifold $E$ together
with a smooth surjective map $\pi: E\to M$ satisfying:
\begin{enumerate}
	\item [(i)] for each $p\in M$ the set $E_p \bydef \pi^{-1}\left( p\right)$  (called the \textit{fiber} of $E$ over
$p$) is endowed with the structure of a real vector space,
\item[(ii)] for each $p\in M$ there exists a neighborhood $\sU$ of $p$ in $M$ and a diffeomorphism $\Phi:\pi^{-1}\left( p\right)\to \sU \times \R^k$  such that the following diagram
		\newline
\begin{tikzpicture}
	\matrix (m) [matrix of math nodes,row sep=3em,column sep=4em,minimum width=2em]
	{
		\pi^{-1}\left(\sU \right)  & & \sU \times \R^k \\ 
		& \sU & \\};
	\path[-stealth]
	(m-1-1) edge node [above] {$\Phi$} (m-1-3)
	(m-1-1) edge node [left]  {$\pi~~$} (m-2-2)
	(m-1-3) edge node [right] {$~~\pi_1$} (m-2-2);
\end{tikzpicture}
\\ 	
and the restriction of $\Phi$ to $E_p$ is a linear isomorphism from $E_p$ to $\{p\}\times \R^k \xrightarrow{\approx } \R^k$. (Here $\pi_1$ is the projection on the first factor).
\end{enumerate}
The manifold $E$ is called the \textit{total space} of the bundle, $M$ is called its \textit{base}, and $\pi$ is its \textit{projection}. Each map $\Phi$ in the above definition is called
a \textit{local trivialization} of $E$ over $\sU$. 

\end{definition}

\begin{definition}\label{top_sm_sec_defn}\cite{lee:smooth_manifolds}
Let $E$ be a smooth vector bundle over a smooth manifold $M$, with projection $\pi: E\to M$. A section of $E$ is a \textit{section} of the map $\pi$, i.e., a continuous map $\sigma: M \hookto E$ such that $\pi\circ\sigma = \Id_M$. If $\sU \subset M$ is an open subset, a
section of the restricted bundle $E|_{\sU}$ is called a \textit{local section} of $E$. A \textit{smooth
section} is a (local or global) section that is smooth as a map between manifolds.
\end{definition}
\begin{remark}
The definitions \ref{top_sm_bundle_defn} and \ref{top_sm_sec_defn} are specializations of \ref{top_vb_defn} and \ref{top_vb_cs_defn} ones.
\end{remark}

\section{Differential operators on manifolds}
\begin{definition}\label{do_man_defn}\cite{kahn:glo_an}
	Given an $n$-dimensional smooth manifold $M$ and two
	smooth, complex vector bundles over $M$, say
	$\left(E_1, M, p_1\right)$ and $\left(E_1, M, p_1\right)$,
	let $ \Ga^\infty\left(M, E_j \right)$  be the complex vector space of sections of the bundle
	$\left(E_j, M, p_j\right)$, i.e., smooth maps $s:M \to E_j$ such that $\pi_j\circ s=\Id_{M}$.
	A \textit{linear differential operator} means a linear map of vector spaces
	$$
	P: \Ga^\infty\left(M, E_1 \right)\to \Ga^\infty\left(M, E_2 \right)
	$$
	such that $\supp P(s) \subset  \supp s$. We denote by 
	\be\label{do_man_eqn}
	D\left(M, E_1, E_2 \right)
	\ee
the space of all differential operators $\Ga^\infty\left(M, E_1 \right)\to \Ga^\infty\left(M, E_2 \right)$.
\end{definition}
\begin{remark}\label{do_man_rem}
If $E_1 = E_2= E$ then we denote
\be\label{do_man_alg_eqn}
	D\left(M, E \right)\bydef 	D\left(M, E, E \right).
\ee
The space $D\left(M, E \right)$ is an algebra over $\C$ with the given by the composition of operators product.
\end{remark}
\begin{definition}\label{do_man_order_defn}\cite{kahn:glo_an}
	Let $M$ be a smooth manifold, $x_0\in M$. Let $E_1$, $E_2$ be two
	smooth vector bundles over $M$, and let
	$$
	P:  \Ga^\infty\left(M, E_1 \right)\to  \Ga^\infty\left(M, E_2 \right)
	$$
	be a linear differential operator, in the sense of Definition \ref{do_man_defn}.
	Then we say that $P$ \textit{has order} $m$ at the point $x_0$ if $m$ is the largest nonnegative
	integer such that there is some $s \in  \Coo\left(M, E_1 \right)$ and some smooth function
	$f$ defined in an open neighborhood of $x_0$ and vanishing at $x_0$ such that $P(f^ms)(x_0)\neq0$.
	The order of $P$ is the maximum of the orders of $P$ at all the points of $M$. 
\end{definition}

\section{Distributions on manifolds}
\paragraph*{}
The notion of a distribution on smooth manifold $M$ is explained in 6.3.3 \cite{horm:I} and the space of distributions  is denoted by $\mathscr{D}'\left( M\right)$.
However $\mathscr{D}'\left( M\right)$ is not equal to the space of linear forms  $C_c\left( M\right) \cap \Coo\left(M \right)\to \C$.
\begin{definition}\label{top_distr_dens_def}
	If $M$ is a smooth manifold then a $\C$-valued linear form  $\Coo_c \left(M  \right)\bydef C_c\left( M\right) \cap \Coo\left(M \right)\to\C$ is said to be a \textit{distribution density}. The space of distribution densities will be denoted by $\Coo_c \left(M  \right)'$. Denote by 
	\be\label{top_distr_dens_eqn}
	\begin{split}
	\left\langle \cdot, \cdot \right\rangle : \Coo_c \left(M  \right)'\times \Coo_c \left(M  \right) \to \C,\\
\left\langle \phi, a \right\rangle \bydef \phi\left(a \right) 
	\end{split}
	\ee
	the natural pairing.
\end{definition}
\begin{remark}
The space $\mathscr{D}'\left( M\right)$ is relevant to $\Coo_c \left(M  \right)'$ one (cf. \cite{horm:I})
\end{remark}
\begin{remark}
	If there is a positive measure $\mu$ on $M$ which is locally presented by a $\dim M$-dimensional external differential form then there is the natural inclusion
	\be
	\begin{split}
		C\left(M \right) \subset \Coo_c \left(M  \right)',\\
		f \mapsto \left( \phi \mapsto  \int_M f \phi d \mu \right) .
	\end{split}
	\ee
\end{remark}
\begin{remark}
	From $\Coo \left(M  \right) \Coo_c \left(M  \right)\subset \Coo_c \left(M  \right)$ it follows that there is the action $\Coo \left(M  \right)\times \Coo_c \left(M  \right)' \to \Coo_c \left(M  \right)'$ given by
	\bean
	\langle ab, c\rangle \bydef \langle b, ac\rangle\quad a \in \Coo \left(M  \right), \quad b \in \Coo_c \left(M  \right)', \quad c \in \Coo_c \left(M  \right).
	\eean
	where 	$\langle \cdot , \cdot \rangle: \Coo_c \left(M  \right)'\times \Coo_c \left(M  \right)\to\C$ is the natural pairing.
\end{remark}
\begin{remark}\label{top_distr_dens_q_rem}
	From the pairing $\Coo \left(M  \right)\times \Coo_c \left(M  \right)' \to \Coo_c \left(M  \right)'$ it turns out that a pair $\left(\Coo_c \left(M  \right)', \Coo \left(M  \right)\right)$ is a quasi *-algebra (cf. Definition \ref{quasi_defn}).
\end{remark}


	\section{Flat connections}\label{geom_flat_subsec}
\paragraph*{}
Here I follow to \cite{kobayashi_nomizu:diff_geom}. Let $M$ be a manifold and $G$ a Lie group. A (\textit{differentiable}) \textit{principal bundle over M with group} $G$ consists of a manifolfd $P$ and an action of $G$ on $P$ satisfying the following conditions:
\begin{enumerate}
	\item [(a)] $G$ acts freely on $P$ on the right: $\left(u, a \right) \in P \times G \mapsto ua = R_au \in P$;
	\item[(b)] $M$ is the quotient space of $P$ by the equivalence relation induced by $G$, i.e. $M = P/G$, and the canonical projection $\pi: P \to M$ is differentiable;
	\item[(c)] $P$ is locally trivial, that is, every point $x$ of $M$ has an open neighborhood $U$ such that $\pi^{-1}\left( U\right)$ is isomophic to  $U\times G$ in the sense that there is a diffeomorphism $\psi:  \pi^{-1}\left( U\right) \to U \times G$ such that $\psi\left( u\right) = \left( \pi\left( u\right), \varphi\left(u \right) \right) $ where $\varphi$ is a mapping of $\pi^{-1}\left(U \right)$ into $G$ satisfying  $\psi\left(ua \right)= \left( \psi\left( u\right)\right) a$  for all $u \in \pi^{-1}\left(U \right)$ and $a \in G$. 
\end{enumerate}
A principal fibre bundle will be denoted by $P\left( M, G, \pi\right), ~ P\left(M, G \right)$ or simply $P$.
\paragraph*{} Let $P\left(M, G \right)$ be a principal fibre bundle over a manifold with group $G$. For each $u \in P$ let $T_u\left(P \right)$ be a tangent space of $P$ at $u$ and $G_u$ the subspace of $T_u\left( P\right)$ consisting of vectors tangent to the fibre through $u$. A \textit{connection} $\Ga$   in $P$ is an assignment of a subspace $Q_u$ of $T_u\left(P \right)$ to each $u \in P$ such that
\begin{enumerate}
	\item [(a)] $T_u\left(P \right) = G_u \oplus Q_u$ (direct sum);
	\item[(b)] $Q_{ua}= \left(R_a \right)_*Q_u$ for every $u \in P$ and $a \in G$, where $R_a$ is a transformation of $P$ induced by $a \in G, ~ R_au=ua$.
\end{enumerate}
\paragraph*{}
Let $P = M \times G$ be a trivial principal bundle. For each $a \in G$, the set $M \times \{a\}$ is a submanifold of $P$. The \text{canonical flat connection} in $P$ is defined by taking the tangent space to $M \times \{a\}$ at $u = \left(x, a \right)$ as the horizontal tangent subspace at $u$. A connection in any principal bundle is called \textit{flat} if every point has a neighborhood such that the induced connection in $P|_U = \pi^{-1}\left(U \right)$ is isomorphic with the canonical flat connection.

\begin{cor}\label{dg_flat_cor}(Corollary II 9.2 \cite{kobayashi_nomizu:diff_geom})
	Let $\Ga$ be a connection in $P\left(M, G \right)$ such that the curvature vanishes identically. If $M$ is paracompact and simply connected, then $P$ is isomorphic to the trivial bundle and $\Ga$ is isomorphic to the canonical flat connection in $M \times G$.
\end{cor}
\begin{remark}
	The \textit{curvature} notion is explained in \cite{kobayashi_nomizu:diff_geom}.
\end{remark}
\paragraph*{}

If $\widetilde{\pi}: \widetilde{M} \to M$ is a covering then the $\widetilde{\pi}$-\textit{lift} of $P$ is a principal $\widetilde{P}\left(\widetilde{M}, G \right)$  bundle, given by
\be\nonumber
\widetilde{P} = \left\{\left(u, \widetilde{x}\right) \in P \times \widetilde{M}~|~ \pi\left(u \right) = \widetilde{\pi}\left( \widetilde{x}\right) \right\}.
\ee
If $\Ga$ is a  connection on $P\left( M, G\right)$ and $\widetilde{M} \to M$ is a covering then is a canonical connection $\widetilde{\Ga}$ on $\widetilde{P}\left(\widetilde{M}, G \right)$ which is the \textit{lift} of $\Ga$, that is, for any $\widetilde{u} \in \widetilde{P}$ the horizontal space $\widetilde{Q}_{\widetilde{u}}$ is isomorphically  mapped onto the horizontal space $Q_{\widetilde{\pi}\left(\widetilde{u} \right) }$ associated with the connection $\Ga$.
If $\Ga$ is flat then from the Proposition (II 9.3 \cite{kobayashi_nomizu:diff_geom}) it turns out that there is a covering $\widetilde{M} \to M$ such that $\widetilde{P}\left(\widetilde{M}, G \right)$ (which is the lift of $P\left(M,G\right)$) is a trivial bundle, so the lift $\widetilde{\Ga}$ of $\Ga$ is a canonical flat connection (cf. Corollary \ref{dg_flat_cor}). From the the Proposition (II 9.3 \cite{kobayashi_nomizu:diff_geom}) it follows that for any flat connection $\Ga$  on $P\left(M, G \right)$ there is a group homomorphism $\varphi: G\left( \widetilde{M} ~|~ M\right) \to G$ such that
\begin{enumerate}
	\item [(a)] There is an action $G\left( \widetilde{M} ~|~ M\right) \times \widetilde{P} \to \widetilde{P} \approx \widetilde{M} \times G$ given by
	\be\label{nontriv_buble_eqn}
	g \left(\widetilde{x}, a \right) = \left( g\widetilde{x}, \varphi\left( g\right) a\right); \quad \forall \widetilde{x} \in \widetilde{M}, ~ a \in G, 
	\ee
	\item[(b)] There is the canonical diffeomorphism  $P = \widetilde{P}/G\left( \widetilde{M} ~|~ M\right)$,
	\item[(c)] The lift $\tilde{\Ga}$ of $\Ga$ is a canonical flat connection.
\end{enumerate}

\begin{defn}
	In the above situation we say that the flat connection $\Ga$ is \textit{induced} by the covering $\widetilde{M}\to M$ and the homomorphism $G\left( \left.\widetilde{M}~\right|M\right) \to G$, or we say that $\Ga$ \textit{comes from} $G\left( \left.\widetilde{M}~\right|M\right) \to G$.
\end{defn}
\begin{remark}
	The  Proposition (II 9.3 \cite{kobayashi_nomizu:diff_geom}) assumes that $\widetilde{M} \to M$ is the universal covering however it is not always necessary requirement.
\end{remark}
\begin{remark}
	If $\pi_1\left(M, x_0 \right)$ is the fundamental group \cite{spanier:at} then there is the canonical surjective homomorphism $\pi_1\left(M, x_0 \right) \to G\left( \left.\widetilde{M}~\right|M\right)$. So there exist the composition $\pi_1\left(M, x_0 \right) \to G\left( \left.\widetilde{M}~\right|M\right) \to G$. It follows that any flat connection comes from the homomorphisms $\pi_1\left(M, x_0 \right) \to G$. 
\end{remark} 
\paragraph*{}
Suppose that there is the right action of $G$ on $P$ and suppose that $F$ is a manifold with the left  action of $G$. There is an action of $G$ on $P \times F$ given by $a\left( u, \xi\right) = \left(u a, a^{-1}\xi \right)$ for any $a \in G$ and $\left( u, \xi\right) \in P\times F$. The quotient space $P \times_G F = \left(P \times F \right)/G$ has the natural structure of a manifold and if $E =  P \times_G F$ then\\ $E\left(M, F, G, P \right)$ is said to be the \textit{fibre bundle over the base $M$, with (standard) fibre $F$, and (structure) group G which is associated with the principal bundle P} (cf. \cite{kobayashi_nomizu:diff_geom}). If $P = M \times G$ is the trivial bundle then $E$ is also trivial, that is, $E = M \times F$. If $F = \C^n$ is a vector space and the action of $G$ on $\C^n$ is a linear representation of the group then $E$ is the linear bundle. Denote by $T\left(M \right)$ (resp. $T^*\left(M \right)$) the tangent  (resp. contangent) bundle, and denote by $\Ga\left( E\right)$, $\Ga\left(T\left(M \right)\right)$, $\Ga\left(T^*\left(M \right)\right)$ the spaces of sections of $E$, $T\left(M \right)$, $T^*\left(M \right)$ respectively. Any connection $\Ga$ on $P$ gives a covariant derivative  on $E$, that is,
for any section  $X \in \Ga\left( T\left(M \right)\right) $ and any section $\xi \in \Ga\left( E\right)$ there is the derivative given by
\be\nonumber
\nabla_X\left( \xi\right) \in \Ga\left( E\right).
\ee
If $E = M \times \C^n$, $\Ga$ is the canonical flat connection and $\xi$ is a trivial section, that is, $\xi = M \times \{x\}$ then 
\be\label{comm_triv_eqn}
\nabla_X \xi = 0,~~ \forall X \in T\left(M \right).
\ee
For any connection there is the unique map
\be\label{comm_alg_conn}
\nabla : \Ga\left( E\right) \to \Ga\left( E \otimes T^*\left( M\right) \right)
\ee
such that
\be\nonumber
\nabla_X \xi = \left(\nabla \xi, X \right)
\ee
where the pairing $\left(\cdot, \cdot\right) : \Ga\left( E \otimes T^*\left( M\right) \right) \times \Ga\left( T\left( M\right)\right)  \to \Ga\left( E \right) $ is induced by the pairing $\Ga\left( T^*\left(M \right)\right)   \times \Ga\left( T\left( M\right)\right)  \to \Coo\left(M \right) $.

\chapter{Algebra}
\section{Groups}\label{grops_sec}
\paragraph{}
We suppose that the reader is familiar  with the notion of group. Here I follow to \cite{lang}.
Let $G$ be a group and $H$ a subgroup. A \textit{left  coset} of $H$ in $G$ is a subset of 
$G$ of type $aH$, for some element $a \in G$. 
$G$ is the disjoint union of the left  cosets of $H$. A similar 
remark applies to right cosets (i.e. subsets of $G$ of type $Ha$). The number of left  
cosets of $H$ in $G$ is denoted by $(G : H)$, and is called the \textit{(left) index} of $H$ in $G$. The kernel of a group-homomorphism is a 
subgroup. We now wish to characterize such subgroups. 
Let $f: G' \to G$ be a group-homomorphism, and let $H$ be its kernel. If $x$ is an 
element of $G$, then
$xH = Hx$ .

	\begin{definition}\label{normal_subgroup_defn}\cite{kurosh:lga}
	A subgroup $H$ is called a \textit{normal subgroup} (or \textit{invariant subgroup}) of the group $G$, if the left-sided partition of the group $G$
	with respect to the subgroup $H$ coincides with the right-sided
	partition, i.e. if for every $g\in G$ the equation
	\be\label{normal_subgroup_eqn}
	gH = Hg
	\ee
	or equivalently
	\be\label{normal_subgroup_ghg_eqn}
	gHg^{-1} = H
	\ee
	holds (understood in the sense that the two subsets coincide in $G$).
	Thus we can speak simply of \textit{the partition of the group} $G$ \textit{with
		respect to the normal subgroup} $H$.
\end{definition}
The left index of any normal subgroup coincides with the right one.
\begin{definition}\label{group_coset_defn}\cite{lang}
Let $A$ be a group and $H$ a subgroup. A \textit{left  coset} of $H$ in $G$ is a subset of $G$ of type $aH$ for some $a\in G$. An element of $aH$ is called a \textit{coset representative} of $aH$. Similarly we define \textit{right cosets} $Ha$.
\end{definition}
\begin{definition}\label{group_inv_lim_defn}\cite{spanier:at}
If $\left\{G_\a\right\}$ is an inverse system of
groups  (that is, $G_\a$ is a group for each $\a$ and $f^\bt_\a: G_\bt\to  G_\a$ is a homomorphism 
for $\a\le \bt$), their \textit{inverse limit} $\varprojlim G_\a$ (which is a set (cf. Definition \ref{inverse_limit_defn})) is a subgroup of $\prod G_\a$. 
\end{definition}

\begin{example}\label{profinite_exm}\cite{lang}
Let $G$ be a group. Let $\F$ be the family of normal subgroups of 
finite index. If $H$, $K$ are normal of finite index, then so is $H \cap K$ so ($\F$ is a 
directed family. We may then form the inverse limit $\varprojlim G/H$ with $H \in \F$
A group which is an inverse limit of finite groups is called \textit{profinite}. It is the \textit{profinite completion} of $G$.

\end{example}
	\section{Finitely generated modules}
	\begin{theorem}\label{fin_gen_thm}\cite{kasch:mr}
		The module $M_R$ over a (unital) ring $R$ is finitely generated if and only if for every set $\left\{A_\la\right\}_{\la\in\La}$ of submodules $A_\la\hookto  M$ with 
		$$
		M=\sum_{\la \in \La} A_\la
		$$
		there is a finite subset $\La_0\subset \La$  such that 
		$$
		M = \sum_{\la \in \La_0} A_\la.
		$$
	\end{theorem}

	\section{Algebraic Morita equivalence}
	
	\begin{defn}\label{morita_ctx_defn}
		A \textit{Morita context} $\left( A,B,P,Q,\varphi,\psi\right)$ or, in some authors (e.g. Bass \cite{bass}) the \textit{pre-equivalence data} is a generalization of Morita equivalence between categories of modules. In the case of right modules, for two associative $\mathbf{k}$-algebras (or, in the case of $\mathbf{k} = \mathbb{Z}$, rings) $A$ and $B$, it consists of bimodules $_AP_B$, $_BQ_A$ and bimodule homomorphisms $\varphi: P\otimes_B Q\to A$, $\psi: Q\otimes_A P\to B$ satisfying mixed associativity conditions, i.e. for any $p,p' \in P$ and $q,q' \in Q$ one  has:
		\begin{equation}\label{morita_ctx_eqn}
			\begin{split}
				\varphi\left(p\otimes q \right) p' = p \psi\left(q\otimes p' \right),\\
				\psi\left(q \otimes p \right) q' = q\varphi\left(p\otimes q' \right).  
			\end{split}
		\end{equation}
		A Morita context is a \textit{Morita equivalence} if both $\varphi$ and $\psi$ are isomorphisms of bimodules. 
	\end{defn}
	
	\begin{rem}\label{morita_rem}
		The Morita context $\left( A,B,P,Q,\varphi,\psi\right)$  is a Morita equivalence if and only if $A$-module $P$ is a finitely generated projective generator (cf. \cite{bass} II 4.4).
	\end{rem}

	\section{Finite Galois coverings}\label{fin_gal_cov_sec}
	\paragraph*{} Here I follow to \cite{auslander:galois}. Let $A \hookto \widetilde{A}$ be an injective homomorphism of unital algebras, such that
	\begin{itemize}
		\item $\widetilde{A}$ is a projective finitely generated $A$-module,
		\item There is an action $G \times \widetilde{A} \to \widetilde{A}$ of a finite group $G$ such that $$A = \widetilde{A}^G=\left\{\widetilde{a}\in \widetilde{A}~|~g\widetilde{a}=\widetilde{a}; ~\forall g \in G\right\}.$$
	\end{itemize}
	Let us consider the category $\mathscr{M}^G_{\widetilde{A}}$ of $G-\widetilde{A}$ modules, i.e.  any object $M \in \mathscr{M}^G_{\widetilde{A}}$ is a $\widetilde{A}$-module with equivariant action of $G$, i.e. for any $m \in M$ a following condition holds
	$$
	g\left(\widetilde{a}m \right)=  \left(g\widetilde{a} \right) \left(gm \right) \text{ for any } \widetilde{a} \in \widetilde{A}, ~ g \in G.
	$$
	Any morphism $\varphi: M \to N$ in the category $\mathscr{M}^G_{\widetilde{A}}$ is $G$- equivariant, i.e.
	$$
	\varphi\left( g m\right)= g \varphi\left( m\right)   \text{ for any } m \in M, ~ g \in G.
	$$
	Let $\widetilde{A}\left[ G\right]$ be an algebra such that $\widetilde{A}\left[ G\right] \approx \widetilde{A}\times G$ as an Abelian group and a multiplication law is given by
	$$
	\left( a, g\right)\left( b, h\right) =\left(a\left(gb \right), gh  \right).
	$$
	The category $\mathscr{M}^G_{\widetilde{A}}$ is equivalent to the category $\mathscr{M}_{\widetilde{A}\left[ G\right]}$ of $\widetilde{A}\left[ G\right]$ modules. Otherwise in \cite{auslander:galois} it is proven that if $\widetilde{A}$ is a finitely generated, projective $A$-module then there is an  equivalence between a category $\mathscr{M}_{A}$ of $A$-modules and the category $\mathscr{M}_{\widetilde{A}\left[ G\right]}$ (cf. Definition \ref{category_equivalence_definition}). It turns out that the category $\mathscr{M}^G_{\widetilde{A}}$ is equivalent to the category $\mathscr{M}_{A}$. The equivalence is given by mutually inversed functors $\left( -\right) \otimes \widetilde{A}: \mathscr{M}_{A} \to \mathscr{M}^G_{\widetilde{A}}$ and $\left( -\right) ^G: \mathscr{M}^G_{\widetilde{A}}\to \mathscr{M}_{A}$.

	\section{Profinite and residually finite groups}\label{profinite_section}
	\begin{empt}\label{prof_empt}\cite{dix:profinite}
		Let $G$ be a group and $\left\{G_\al \right\}$ be the set of all normal subgroups of finite index in $G$. Then the set $\left\{G/G_\al, \phi_{\al\bt}\right\}$ of finite quotients $G/G_\al$, of $G$ together with the canonical projections $\phi_{\al\bt}:G/G_\al \to G/G_\bt$  whenever $G_\al \subset G_\bt$ is an inverse system. The inverse limit $\varprojlim G/G_\a$ of this system is called the \textit{profinite completion} of $G$ and is denoted by $\widehat{G}$. The group $\widehat{G}$ can also be described as the closure of the image of $G$ under the
		the diagonal mapping $\Delta : G \to \prod \left( G/G_\a\right)$  where $G/G_\a$, is given the discrete topology and $\prod \left( G/G_\a\right)$ has the product topology. In this description the elements of $\widehat{G}$ are the
		elements $\left(g_\a \right) \in \prod \left( G/G_\a\right)$ which satisfy $\phi_{\al\bt}\left( g_\al\right)= \bt$  whenever $G_\a \subset G_\bt$.
	\end{empt} 
	
	\begin{definition}\label{residually_finite_defn}\cite{bogopolsky:group_theory}
		A group $G$ is said to be \textit{residually finite} if for each nontrivial element from $G$ there exists a finite group $K$ and a homomorphism $\varphi: G\to K$ such that $\varphi\left( g\right)$ is not trivial.
	\end{definition}
	
	


	\chapter{Functional analysis}
\section{Weak topologies}
	
	\begin{definition}\label{w_topology_defn}\cite{rudin:fa}
		Suppose $X$ is a topological vector space (with topology $\tau$) whose dual $X^*$ separates points on $X$.   The $X^*$-topology (cf. Definition \ref{f_topology_defn}) of $X$ is called the 
	\textit{weak topology} of $X$.
	\end{definition}
	
	\begin{definition}\label{w*_topology_defn}\cite{rudin:fa}
		Let $X$ again be a  	topological vector space whose dual is $X^*$. For the definitions that follow, it is irrelevant whether $X^*$ separates points on X or not. The important  observation to make is that every $x\in X$ induces a linear functional $f_x$ on $X^*$, 
		defined by 
		$$
		f_x\La= \La x
		$$
		The $X$-topology of $X^*$ (cf. Definition \ref{f_topology_defn}) is called the \textit{weak}*-\textit{topology of} $X^*$.
		
	\end{definition}
	\begin{remark}\label{weak_convergence_rem}\cite{rudin:fa}
		In the above situation the net $\left\{x_\a\right\}_{\a\in \mathscr A}\in X$ is convergent with respect to weak topology if the net $\left\{\La x_\a\right\}\in \C$ is convergent for all $\La \in X^*$. Similarly the net $\left\{\La_\a\right\}_{\a\in \mathscr A}\in X^*$ is convergent with respect to weak*-topology if the net $\left\{\La_\a x\right\}\in \C$ is convergent for all $x\in X$.
	\end{remark}
	\section{Unbounded operators}
	
	\begin{definition}\label{unb_defn}(cf. \cite{kreyszig:fa})
		Let $\H$ be a Hilbert space, and let $\D\subset \H$ is a dense subspace. A $\C$-linear map $T: \D \to \H$ is an \textit{unbounded operator} if there is no $\widetilde T \in B\left( \H\right)$ such that $T = \left.\widetilde T\right|_\D= T$. Denote by
		\be\label{unb_eqn}
		\D\left(T \right) \bydef \D = \Dom T.
		\ee 
	\end{definition}
	
	\begin{definition}\label{unb_ext_defn}(cf. \cite{kreyszig:fa})
		Let $T, S$ be unbounded operator. We say that $T$ is an \textit{extension of} $S$ if one has
		\be\label{unb_ext_eqn}
		\begin{split}
			\D\left(S\right)\subset	\D\left(T \right),\\
			\left.T\right|_{\D\left(S\right) }= S.
		\end{split}
		\ee 
	\end{definition}
	
	\begin{definition}\cite{kreyszig:fa}
		Let $T: \D\left( T\right)\to \H$  be
		a (possibly unbounded) densely defined linear operator in a complex
		Hilbert space $\H$. Then the \textit{Hilbert-adjoint operator} (or simply adjoint) $T^*: \D\left(T^*\right) \to \H$ of $T$ is defined as follows. The domain $\D\left(T^*\right)$ of $T^*$ consists of all $y\in \H$
		such that there is a $y\in \H$ satisfying
		\be
		\left( Tx, y\right)_\H = \left(x, y^*\right)_\H 
		\ee
		for all $x\in\D\left(T^*\right)$. For each such $y\in \D\left(T^*\right)$ the Hilbert-adjoint operator
		$T^*$ is then defined in terms of that $y$ by
		\be
		y^* \bydef T^*y.
		\ee
		
	\end{definition}
	\begin{theorem}\cite{kreyszig:fa}
		Let $S: \D\left( S\right)\to \H$ and
		$T: \D\left( T\right)\to \H$ be linear operators which are densely defined in a
		complex Hilbert space $\H$. Then:
		\begin{enumerate}
			\item[(a)] If $S \subset T$ then $T^*\subset S^*$.
			\item[(b)]  If $\D\left( T\right)$ is dense in $\H$ then $T\subset T^{**}$.
		\end{enumerate}
	\end{theorem}
	\begin{definition}\cite{kreyszig:fa}
		Let $T: \D\left( T\right)\to\H$
		be a linear operator which is densely defined in a complex Hilbert
		space $\H$. Then $T$ is called a \textit{symmetric linear operator} if for all $x, y \in \D\left( T\right)$
		\be
		(Tx, y)_\H=(x, Ty)_\H.
		\ee
	\end{definition}
	\begin{lemma}\cite{kreyszig:fa}
		A densely defined linear operator
		$T$ in a complex Hilbert space H is symmetric if and only if
		\be
		T \subset T^*.
		\ee
	\end{lemma}
	\begin{definition}\cite{kreyszig:fa}
		Let $T: \D\left( T\right)\to\H$ be a linear operator which is densely defined in a complex Hilbert
		space $\H$. Then $T$ is called a \textit{self-adjoint linear operator} if
		\be
		T = T^*
		\ee
		
		Every self-adjoint linear operator is symmetric.
	\end{definition}
	\begin{definition}\cite{kreyszig:fa}
		Let $T: \D\left( T\right)\to\H$ be a
		linear operator, where $\D\left( T\right)\subset\H$) and $\H$ is a complex Hilbert space.
		Then $T$ is called a \textit{closed linear operator} if its graph
		$$
		\G\left(T\right)\bydef\left\{\left(x, y \right) \in \H\times \H\left|\right.x \in \D\left( T\right), \quad y = Tx\right\}
		$$
		is closed in $\H\times \H$, where the norm on $\H\times \H$ is defined by
		$$
		\left\| \left(x, y \right) \right\| \bydef \sqrt{\left\|x \right\|^2+ \left\|y \right\|^2}
		$$
		and results from the inner product defined by
		$$
		\left(\left(x_1, y_1 \right), \left(x_2, y_2 \right)\right)_{\H\times \H} = \left(x_1, x_2 \right)_\H + \left(y_1, y_2 \right)_\H
		$$
	\end{definition}
	\begin{theorem}\cite{kreyszig:fa}
		Let $T: \D\left( T\right)\to\H$ be a
		linear operator, where $\D\left( T\right)\subset \H$ and $\H$ is a complex Hilbert space.
		Then:
		\begin{enumerate}
			\item [(a)] $T$ is closed if and only if
			\be
			x_\a  \to x \quad \{x_\a\}\subset \D\left( T\right) \quad \text{and} \quad 	Tx_\a  \to y
			\ee
			together imply that $x \in \D\left( T\right)$ and $Tx = y$,
			\item[(b)]  if $T$ is closed and $\D\left( T\right)$ is closed, then $T$ is bounded,
			\item[(c)] If $T$ be bounded then $T$ is closed if and only if $\D\left( T\right)$ is
			closed.
		\end{enumerate}
	\end{theorem}
	\begin{definition}\cite{kreyszig:fa}
		If a linear operator $T$
		has an extension $T_1$ which is a closed linear operator, then $T$ is said to
		be \textit{closable}, and $T_1$ is called a closed linear extension of $T$.
		A closed linear extension $\widetilde{T}$ of a closable linear operator T is said
		to be \textit{minimal} if every closed linear extension $T_1$ of $T$ is a closed linear
		extension of  $\widetilde{T}$. This minimal extension $\widetilde{T}$ of $T$ -if it exists-is called
		the \textit{closure} of $T$.
	\end{definition}
	\begin{theorem}\cite{kreyszig:fa}
		Let $T: \D\left( T\right)\to\H$  be a linear operator,
		where $\H$ is a complex Hilbert space and $\D\left( T\right)$ is dense in $\H$. Then if $T$ is
		symmetric, its closure $\widetilde{T}$ exists and is unique.
	\end{theorem}

	\section{Fourier transformation}
	\paragraph*{}
	There is a norm on $\mathbb{Z}^n$ given by
	\begin{equation}\label{mp_znorm_eqn}
		\left\|\left(k_1, ..., k_n\right)\right\|= \sqrt{k_1^2 + ... + k^2_n}.
	\end{equation}
	The space of complex-valued Schwartz  functions on $\Z^n$ is given by 
	\begin{equation}\label{schwartz_z_eqn}
		\begin{split}
	\sS\left(\mathbb{Z}^n\right)\bydef\\\bydef  \left\{\left.a = \left\{a_k\right\}_{k \in \mathbb{Z}^n} \in \mathbb{C}^{\mathbb{Z}^n}~\right|~ \mathrm{sup}_{k \in \mathbb{Z}^n}\left(1 + \|k\|\right)^s \left|a_k\right| < \infty, ~ \forall s \in \mathbb{N} \right\}.
		\end{split}
		\end{equation}
	Let $\mathbb{T}^n$ be an ordinary $n$-torus. We will often use real coordinates for $\mathbb{T}^n$, that is, view $\mathbb{T}^n$ as $\mathbb{R}^n / \mathbb{Z}^n$. Let $\Coo\left(\mathbb{T}^n\right)$ be an algebra of infinitely differentiable complex-valued functions on $\mathbb{T}^n$. 
	There is the bijective Fourier transformations  $\mathcal{F}_\T:\Coo\left(\mathbb{T}^n\right)\xrightarrow{\approx}\sS\left(\mathbb{Z}^n\right)$;  $f \mapsto \widehat{f}$ given by
	\begin{equation}\label{nt_fourier_eqn}
		\widehat{f}\left(p\right)= \mathcal F_\T (f) (p)= \int_{\mathbb{T}^n}e^{- 2\pi i x \cdot p}f\left(x\right)dx
	\end{equation}
	where $dx$ is induced by the Lebesgue measure on $\mathbb{R}^n$ and   $\cdot$ is the  scalar
	product on the Euclidean space $\R^n$.
	The Fourier transformation carries multiplication to convolution, i.e.
	\begin{equation*}
		\widehat{fg}\left(p\right) = \sum_{r +s = p}\widehat{f}\left(r\right)\widehat{g}\left(s\right).
	\end{equation*}
	The inverse Fourier transformation $\mathcal{F}^{-1}_\T:\sS\left(\mathbb{Z}^n\right)\xrightarrow{\approx} \Coo\left(\mathbb{T}^n\right)$;  $\quad \widehat{f}\mapsto f$ is given by
	$$
	f\left(x \right) =\mathcal{F}^{-1}_\T \widehat f\left( x\right)  = \sum_{p \in \Z^n} \widehat f\left( p\right)   e^{ 2\pi i x \cdot p}.
	$$
	There  is the $\C$-valued scalar product  on $\Coo\left( \T^n\right)$ given by
	$$
	\left(f, g \right) = \int_{\T^n}fg dx =\sum_{p \in \Z^n}\widehat{f}\left( -p\right) \widehat{g}\left(p \right).  
	$$ 
	Denote by $\SS\left( \R^{n}\right) $ be the space of
	complex Schwartz (smooth, rapidly decreasing) functions on $\R^{n}$. 
	\be\label{mp_sr_eqn}
	\begin{split}
		\SS\left(\mathbb {R} ^{n}\right)\bydef\\
	\bydef\left\{f\in C^{\infty }(\mathbb {R} ^{n})\left|\|f\|_{\alpha  ,\beta )}<\infty \quad \forall \alpha 
		=\left( \al_1,...,\al_n\right) ,\beta =\left( \bt_1,...,\bt_n\right)\in \mathbb {Z} _{+}^{n}\right.\right\},\\
		\|f\|_{{\alpha ,\beta }}=\sup_{{x\in {\mathbb  {R}}^{n}}}\left|x^{\alpha }D^{\beta }f(x)\right|
	\end{split}
	\ee
	where 
	\bean
	x^\al = x_1^{\al_1}\cdot...\cdot x_n^{\al_n},\\
	D^{\beta} = \frac{\partial}{\partial x_1^{\bt_1}}~...~\frac{\partial}{\partial x_n^{\bt_n}}.
	\eean
	The topology on $\SS\left(\mathbb {R} ^{n}\right)$ is given by seminorms $\|\cdot\|_{{\alpha ,\beta }}$.
	Let $\mathcal F$ and $\mathcal F^{-1}$ be the ordinary and inverse Fourier transformations given by
	\begin{equation}\label{intro_fourier}
		\begin{split}
			\left(\mathcal{F}f\right)(u) = \int_{\mathbb{R}^{2N}} f(t)e^{-2\pi it\cdot u}dt,~\left(\mathcal F^{-1}f\right)(u)=\int_{\mathbb{R}^{2N}} f(t)e^{2\pi it\cdot u}dt 
		\end{split}
	\end{equation}
	which satisfy  following conditions
	$$
	\mathcal{F}\circ\mathcal{F}^{-1}|_{\SS\left( \R^{n}\right)} = \mathcal{F}^{-1}\circ\mathcal{F}|_{\SS\left( \R^{n}\right)} = \Id_{\SS\left( \R^{n}\right)}.
	$$
	There is the $\C$-valued scalar product  on $\SS\left( \R^n\right)$ given by
	\begin{equation}\label{fourier_scalar_product_eqn}
		\left(f, g \right) = \int_{\R^n}fg dx =\int_{\R^n}\mathcal{F}f\mathcal{F}g dx. 
	\end{equation}
	
	which if $\mathcal{F}$-invariant, i.e.
	\be\label{mp_inv_eqn}
	\left(f, g \right)_{L^2\left( \R^n\right) } = \left(\mathcal{F}f, \mathcal{F}g \right)_{L^2\left( \R^n\right) }.
	\ee
		There is the action of $\Z^n$ on $\R^n$ such that
	$$
	g x = x + g; ~ x \in \R^n,~ g \in \Z^n
	$$
	and $\T^n \approx \R^n / \Z^n$. Any $f \in \Coo\left( \T^n\right)$ can be regarded as $\Z^n$- invariant and smooth function on $\R^n$. On the other hand if $f \in \SS\left( \R^n\right)$ then the series 
	$$
	h = \sum_{g \in \Z^n} g f 
	$$ 
	is point-wise convergent and $h$ is a smooth $\Z^n$ - invariant function. So we can assume that $h \in \Coo\left(\T^n \right)$. This construction provides a map
	\begin{equation}\label{mp_sooth_sum_eqn}
		\begin{split}
			\SS\left(\R^n\right) \to  \Coo\left(\T^n\right), \\
			f \mapsto h = \sum_{g \in \Z^n} g f.
		\end{split}
	\end{equation}
	If $\widehat f = \mathcal F f$,  $\widehat h = \mathcal F_{\T} h$ then for any $p \in \Z^n$ a following condition holds
	\begin{equation}\label{fourier_from_r_to_z_eqn}
		\widehat h\left(p\right) = \widehat f\left( p\right). 
	\end{equation}  
	
	\section{Measure theory}

	\begin{theorem}\label{meafunc_thm}\cite{bogachev_measure_v2}
		Let $\sX$ be a locally compact space and let $L$ be a linear
		function on $C_c\left(\sX \right)$  such that $L\left(f\right)\ge 0$ if $f\ge 0$. Then, there exists a Borel
		measure $\mu$ on $\sX$ with values in $\left[0, +\infty\right]$ such that
		\be\label{meafunc_eqn}
		L\left(f \right) = \int_{\sX} f~ d\mu \quad \forall f \in C_0\left(\sX \right).
		\ee
		In addition, one can choose $\mu$ in such a way that it will be Radon on all sets of finite measure 
		 and there is only one measure with this property).
	\end{theorem}
\begin{remark}
In \cite{bogachev_measure_v2} the space of $C_c\left( \sX\right)$ is denoted by  $C_0\left( \sX\right)$.
\end{remark}
	
	\begin{defn}\label{measure_image_defn} Definition \ref{top_lift_measure_defn}\cite{bogachev_measure_v2}
		Let $\mu$ be a Borel measure on a topological space $\mathcal{X}$ and let $f$ be a $\mu$- measurable mapping from $\mathcal{X}$ to a topological space $\mathcal{Y}$. Then on $\mathcal{Y}$ we obtain
		the Borel measure $\nu = \mu \circ f^{-1}: B \mapsto \mu \circ f^{-1}\left( B\right)$. The measure $\nu$ is called
		the \textit{image} of $\mu$ under the mapping $f$, and $\mu$ is called a \textit{preimage} of $\nu$. The  same terms are used in the case of general measurable mappings of measurable  spaces.
	\end{defn}

	\chapter{Operator algebras}
	
	\section{$C^*$-algebras and von Neumann algebras}
	\paragraph{}In this section I follow to \cite{apt_mult,blackadar:ko,murphy,pedersen:ca_aut,rae:ctr_morita}.
	\begin{definition}\cite{murphy}
		A \textit{Banach *-algebra} is a *-algebra $A$ together with a complete submultiplicative norm such that $\left\| a^*\right\|=\left\| a\right\| ~\forall a\in A$. If, in addition, $A$ has a unit
		such that $\left\| 1\right\|=1$, we call $A$ a \textit{unital Banach *-algebra}.
		A \textit{$C^*$-algebra} is a Banach *-algebra such that
		\be
		\left\| a^*a\right\|=\left\| a\right\|^2\quad \forall a\in A.
		\ee
	\end{definition}
	\begin{example}\cite{murphy}
		If $\H$ is a Hilbert space, then the algebra of bounded operators $B\left(\H\right)$ is a $C^*$-algebra. 
	\end{example}
	\begin{definition}\label{ideal_left_right_defn}\cite{murphy}
A \textit{left} (respectively, \textit{right}) \textit{ideal} in an algebra $A$ is a vector subspace $I$
of $A$ such that
$$
\forall a \in A\quad \forall b \in I \quad ab \in I\quad \text{(respectively } ba \in I \text{)}.
$$
	\end{definition}
\begin{definition}\label{opposite_algebra_defn}
If $A$ is an associative algebra then an \textit{	opposite algebra} $A^\opp$, consisting of elements
$
\left\{\left. a^0~ \right| a\in A \right\}
$
with product $a^0b^0\bydef \left(ba\right)^0$. If $A$ is a *-algebra then one can define $a^0 \bydef a^*$.
\end{definition}
	\begin{definition}\label{essential_defn}\cite{rae:ctr_morita}
		A two sided ideal $I$ in a $C^*$-algebra $A$
		is essential if $I$ 
		has nonzero intersection 
		with every other nonzero ideal $A$. 
	\end{definition}
	Alternatively the essential ideal can be given by the following lemma.
	\begin{lemma}\label{essential_lem}\cite{rae:ctr_morita}
		An ideal $I$ 
		is essential if and only if $aI= \left\{0\right\}$
		implies $a 
		=0$.
	\end{lemma}
	\begin{example}\label{comm_ess_exm}\cite{rae:ctr_morita}
		Let $A\bydef C_0\left( \sX\right)$ and let $\sU$ be an open subset of $\sX$. Then 
		$$
		I\bydef \left(\left.f\in A\right| f\left(x\right)= 0 \quad \forall x \in \sX\setminus\sU\right)
		$$
		is an essential ideal in $A$ 
		if and only if $\sU$ 
		is dense 
		in $\sX$. 
	\end{example}

\begin{definition}\label{faithful_representation_defn}\cite{murphy}
A representation $\rho : A\to B\left( \H\right)$ is called \textit{faithful} if the $*$-homomorphism $\rho$ is injective.
\end{definition}
	
	\begin{definition}\label{nondegenerate_repr_defn}\cite{matro:hcm}
		A representation $\rho : A\to B\left( \H\right)$ is called \textit{nondegenerate} if for any $\xi \in \H$  there exists an element $a \in A$ such that $\rho\left(a \right)\xi \neq 0$. 
	\end{definition}
	\begin{lemma}\label{nondegenerate_repr_lem}\cite{matro:hcm}
		A representation $A\to B\left( \H\right)$ is {nondegenerate} if $\rho\left(A\right)\H$ is dense in $\sH$. 
	\end{lemma}
	\begin{defn}\label{unitization_defn}\cite{rae:ctr_morita}
		A \textit{unitization} of a $C^*$-algebra $A$ is  a $C^*$-algebra $B$ with identity and an injective $*$-homomorphism $\iota: A \hookto B$ such that $\iota\left(A\right)$ is an essential ideal of $B$. 
	\end{defn}
	
	\begin{example}\label{unitization_exm}
		Suppose $A$ is  $C^*$-algebra which has no identity. Then $A^+ = A \oplus \C$
		is a *-algebra with 
	\be\label{a_plus_eqn}
		\left( a \oplus \la\right)\left( b \oplus \mu\right) = \left(ab + \la b + \mu a\right)  \oplus \la \mu, \quad \left(a \oplus \la \right)^* = a^*\oplus \overline{\la}. 
\ee
		It is proven in \cite{rae:ctr_morita} that there is the natural unique $C^*$-norm $\left\| \cdot\right\|_{A^+}$  on $A^+$ such that 
		$$
		\left\|a \oplus 0 \right\|_{A^+}=\left\|a \oplus 0 \right\|_{A}
		$$
		where $\left\|\cdot \right\|_{A}$ is the $C^*$-norm on $A$. Thus $A^+$ is an unital $C^*$-algebra, and the natural map $A\hookto A^+$ is a unitization.
	\end{example}
	\begin{definition}\label{multiplier_min_defn}
		Let $A$ be a $C^*$-algebra. The described in the Example \ref{unitization_exm} unitization   $\iota: A \hookto B$  is called \textit{minimal}.
	\end{definition}
	\begin{definition}\label{multiplier_max_defn}\cite{rae:ctr_morita}
		A unitization   $\iota: A \hookto B$  is called \textit{maximal} if for every embedding $j: A\hookto C$ of $A$ as an essential ideal of a $C^*$-algebra there is a $*$-homomorphism $\phi: C\to B$ such that $\phi \circ j = \iota$. 
	\end{definition}
	
	\begin{rem}\label{multiplier_rem}
		It is proven in \cite{rae:ctr_morita} that for any $C^*$-algebra $A$ there unique maximal unitization.
	\end{rem}
	
	\begin{definition}\label{multiplier_defn}\cite{murphy,pedersen:ca_aut}
		We say that the maximal unitization of $A$ is the \textit{multiplier algebra} of $A$ and denote it by $M\left( A\right)$. 
	\end{definition}
	\begin{defn}\label{strict_topology_defn}\cite{pedersen:ca_aut}
		Let $A$ be a $C^*$-algebra.  The {\it strict topology} on the multiplier algebra $M(A)$ is the topology generated by seminorms 
		\be\label{strict_topology_norm_eqn}
		\vertiii{x}_a\bydef \|ax\| + \|xa\|,\quad a\in A.
		\ee
 If $\La$ is a directed set and $\left\{a_\la\in M\left( A\right) \right\}_{\la\in \La}$ is a net the we denote by $\bt\text{-}\lim_{\la\in\La }a_\la$ the limit of $\left\{a_\la \right\}$ with respect to the strict topology.
		If $x \in M(A)$  and a sequence of partial sums $\sum_{i=1}^{n}a_i$ ($n = 1,2, ...$), ($a_i \in A$) tends to $x$ in the strict topology then we shall write
		\begin{equation}\label{strict_topology_eqn}
			x = \bt\text{-}\sum_{i=1}^{\infty}a_i.
		\end{equation}
	\end{defn}
We also use an alternative equivalent definition of the multiplier algebra.
\begin{definition}\label{multiplier_el_defn}\cite{matro:hcm}
	Let $\rho: A\hookto B\left( \H\right)$ be a faithful {nondegenerate} (cf. Definitions \ref{faithful_representation_defn}, \ref{nondegenerate_repr_defn}) representation, so we assume $A \subset B\left( \H\right)$. An operator $x \in B\left(\H\right)$ is called (two-sided) \textit{multiplier} if 
	\be\label{multiplier_el_eqn}
	xa \in A, \quad ax\in A
\ee
	for each $a\in A$. Denote by $M\left(A\right)$ the set of all multipliers. It is easy to see that $M\left(A\right)$ is an involutive unital algebra.
\end{definition}
\begin{lemma}\label{multiplier_al_lem}\cite{matro:hcm}
	The set $M\left(A\right)$ is an unital $C^*$-algebra 
	$$
	A \subset M\left( A\right) \subset A^{**}
	$$	
	and $A$ is an essential ideal of $M\left( A\right)$. If $A$ has no unit then $A^+\subset M\left(A \right)$. 
\end{lemma}
\begin{definition}\label{double_centralizer_defn}\cite{matro:hcm}
A pair $\left(L, R\right)$ of maps
\be\label{double_centralizer_eqn}
L: A \to A, \quad R: A \to A\quad R\left(a \right) b = a L\left( b\right) \quad \forall a, b \in A
\ee
is called a \textit{double centralizer}. Let us denote $\mathbf{DC}\left(A \right) $ the space of all double  centralizers of $A$.
\end{definition}

\begin{proposition}\label{dc_prop}\cite{matro:hcm}
	If $\left(L, R\right)\in \mathbf{DC}\left(A \right)$ then
	\begin{enumerate}
		\item [(i)] $L\left( ab\right) = L\left( a\right)b$ and $R\left( ab\right) = aR\left( b\right)$;
		\item [(ii)] $L$ and $R$ are linear;
		\item [(iii)] $L$ and $R$ are bounded and $\left\|L \right\|= \left\|R \right\|$;
	\end{enumerate}
The set $\mathbf{DC}\left(A \right) $ with operations
\bean
\left(L_1, R_1\right) + \left(L_2, R_2\right)\bydef \left(L_1 + L_2, R_1+R_2\right),\quad z\left(L,R\right)\bydef \left(zL, zR\right) \quad \forall z \in \C,\\
\left(L_1, R_1\right)  \left(L_2, R_2\right)\bydef \left(L_1  L_2, R_2R_2\right),\\
\left(L, R\right)^*\bydef \left(R^*, L^*\right),\\
L^*\left( a\right)\bydef \left( L\left( a^*\right)\right)^*, \quad R^*\left( a\right)\bydef \left( L\left( a^*\right)\right)^*\quad \forall a \in A.
\eean 
is a normed involutive algebra with respect to the norm
$$
\left\|\left(L, R\right) \right\|\bydef \left\|L \right\|= \left\|R \right\|.
$$
\end{proposition}
\begin{remark}\label{double_centralizer_rem}
	In \cite{matro:hcm} it is proven  there is a natural $*$-isomorphism $\mathbf{DC}\left(A \right)\cong M\left(A\right)$. 
\end{remark}

	\begin{defn}\label{approximate_unit_defn} \cite{pedersen:ca_aut}
		Let $A$ be a $C^*$-algebra. A net $\left\{u_\la \right\}_{\la \in \La}$ in $A_+$ with $\left\|u_\la \right\| \le 1$ for all $\la \in \La$ is called an \textit{approximate unit} for $A$ if $\la < \mu$ implies $u_\la < u_\mu$ and if $\lim \left\|x- xu_\la \right\| = 0$ for each $x$ in $A$. Then, of course, $\lim \left\|x- u_\la x \right\| = 0$ as well.
	\end{defn}
\begin{theorem}\label{left_ideal_thm}\cite{murphy}
 If $L$ is a closed left ideal in a $C^*$-algebra $A$, then there
	is an increasing net $\left\{u_\la\right\}_{\la\in\La}$ of positive elements in the closed unit ball of
$L$ such that $a = \lim_{\la\in \La}au_\la $ for all $a\in L$.
\end{theorem}
\begin{proof}
Set $B \bydef L\cap L^*$. Since $B$ is a $C^*$-algebra, it admits an approximate
unit $\left\{u_\la\right\}_{\la\in \La}\subset B$. If $a \in L$, then $a^*a\in B$, so $0 =
\lim_{\la\in \La} a^*a\left(1_{A^\sim }- u_\la \right)$ . Hence,\\ $\lim_{\la\in \La}\left\|  a - au_\la \right\|^2= \lim_{\la\in \La}\left\| \left(1-u_\la \right)a^*a \left(1-u_\la \right) \right\|
\le \lim_{\la\in \La}\left\|a^*a \left(1-u_\la \right) \right\|=0$, and therefore $\lim_{\la\in \La}\left\|a - au_\la \right\|=0 $.
In the preceding proof we worked in the unitization $A^\sim$ of $A$. 
frequently do this tacitly.
\end{proof}

	\begin{thm}\label{approximate_unit_thm} \cite{pedersen:ca_aut}
		Each $C^*$-algebra contains an \textit{approximate unit}.
	\end{thm}
	\begin{proposition}\label{mult_str_pos_prop}\cite{apt_mult}
		If $B$ is a $C^*$-subalgebra of $A$ containing an
		approximate unit for $A$, then  $M\left(B \right) \subset M\left( A\right)$  (regarding $B''$ as a subalgebra	of $A''$).
	\end{proposition}
	\begin{definition}\label{hered_defn}\cite{pedersen:ca_aut}
		A cone $M$ in the positive part of $C^*$-algebra $A$ is said to be \textit{hereditary} if $0 \le x \le y$, $y \in M$ implies $x \in M$ for each $x \in A$. A *-subalgebra $B$ of $A$ is \textit{hereditary} if $B_+$ is hereditary in $A_+$.
	\end{definition}
	\begin{lemma}\label{hered_bab_lem}\cite{murphy}
		Let $B$ be a $C^*$-subalgebra of $C^*$-algebra $A$. Then $B$ is hereditary in $A$ if and only if $bab' \in B$ for all $b, b' \in B$ and $a \in A$.
	\end{lemma}
	\begin{lemma}\label{hered_ideal_lem}\cite{murphy}
		Let $A$ be a $C^*$-algebra.
		\begin{enumerate}
			\item[(i)] If $L$ is a closed left ideal in $A$ then $L\cap L^*$ is a hereditary $C^*$-subalgebra of $A$. The map $L \mapsto L\cap L^*$ is the bijection from the set of closed left deals of $A$ onto the the set of hereditary $C^*$-subalgebras of $A$.
			\item[(ii)] If $L_1, L_2$ are closed left ideals, then $L_1 \subseteq L_2$ is and only if $L_1\cap L_1^* \subset L_2\cap L_2^*$.
			\item[(iii)] If $B$ is a hereditary $C^*$-subalgebra of $A$, then the set 
			$$
			L\left(B \right) = \left\{\left.a \in A~\right| a^*a \in B \right\}
			$$
			is the unique closed left ideal of $A$ corresponding to $B$.
		\end{enumerate}
	\end{lemma}
\begin{remark}\cite{murphy}
Obviously, $0$ and
$A$ are hereditary $C^*$-subalgebras of $A$, and any intersection of hereditary
$C^*$-subalgebras is one also. 

\end{remark}
\begin{definition}\label{hered_gen_defn}\cite{murphy}
The hereditary $C^*$-subalgebra \textit{generated} by a
subset $S$ of $A$ is the smallest hereditary $C^*$-subalgebra of $A$ containing $S$.
\end{definition}

	\begin{defn}\label{center_defn}\cite{murphy}
		If $A$ is a $C^*$-algebra, its \textit{center} $C$ is the set of elements of $A$ commuting with every $a\in A$.
	\end{defn}
	\begin{defn}
		\label{strong_topology_defn}\cite{pedersen:ca_aut} Let $\H$ be a Hilbert space. The {\it strong} topology on $B\left(\H\right)$ is the locally convex vector space topology associated with the family of seminorms of the form $x \mapsto \|x\xi\|$, $x \in B(\H)$, $\xi \in \H$.
	\end{defn}
	\begin{defn}\label{weak_topology_defn}\cite{pedersen:ca_aut} Let $\H$ be a Hilbert space. The {\it weak} topology on $B\left(\H\right)$ is the locally convex vector space topology associated with the family of seminorms of the form $x \mapsto \left|\left(x\xi, \eta\right)\right|$, $x \in B(\H)$, $\xi, \eta \in \H$.
	\end{defn}
\begin{defn}\label{commutant_defn}
For each subset $M$ of $B(\H)$ let $M'$ denote the \textit{commutant} of $M$, i.e. 
\bean
M'\bydef \left\{b \in B\left(\sH\right)|\forall a \in M \quad ab = ba \right\}
\eean
The $C^*$-algebra 
\bean
M''\bydef (M')'
\eean
is said to be a \textit{bicommutant} of $M$.
\end{defn}
	
	\begin{thm}\label{von_Neumann_thm}\cite{pedersen:ca_aut}
		Let $M$ be a $C^*$-subalgebra of $B(\H)$, containing the identity operator. The following conditions are equivalent:
		\begin{itemize}
			\item $M=M''$ where $M''$ is the bicommutant of $M$;
			\item $M$ is weakly closed;
			\item $M$ is strongly closed.
		\end{itemize}
	\end{thm}

	\begin{defn}\cite{pedersen:ca_aut}
		Any $C^*$-algebra $M$ is said to be a {\it von Neumann algebra} or a {\it $W^*$- algebra} if $M$ satisfies to the conditions of the Theorem \ref{von_Neumann_thm}.
	\end{defn}
	\begin{definition}\label{factor_defn}\cite{pedersen:ca_aut}
		We say that a  von Neumann algebra $M$ is 	a \textit{factor} if the center of consists only of scalar multiplies of $1_M$.
	\end{definition}

	\begin{lemma}\label{increasing_convergent_w_lem}\cite{pedersen:ca_aut} Let $\Lambda$ be an increasing in the partial ordering.  Let $\left\{x_\lambda \right\}_{\la \in \La}$ be an increasing of self-adjoint operators in $B\left(\H\right)$, i.e. $\la \le \mu$ implies $x_\la \le x_\mu$. If $\left\|x_\la\right\| \le \ga$ for some $\ga \in \mathbb{R}$ and all $\la$ then $\left\{x_\lambda \right\}$ is strongly convergent to a self-adjoint element $x \in B\left(\H\right)$ with $\left\|x_\la\right\| \le \ga$.
	\end{lemma}
\begin{defn}\cite{pedersen:ca_aut}
	For any $x\in B(\H)$ element $|x| \stackrel{\text{def}}{=} (xx^*)^{1/2}$ is said to be the {\it absolute value of} $x$.
\end{defn}
For each $x\in B(\H)$ we define the { range projection} of $x$ (denoted by $[x]$) as projection on the closure of $x\H$. If $M$ is a von Neumann algebra and $x \in M$ then $[x]\in M$.

	\begin{prop}\label{polar_decomposition_prop}\cite{pedersen:ca_aut}
		For each element $x$ in   a von Neumann algebra $M$ there is a unique partial isometry $u\in M$ and positive $\left|x\right| \in M_+$ with $uu^*=[|x|]$ and  $x=|x|u$.
	\end{prop}
	\begin{defn}\label{polar_decomposition_defn}
		The formula $x=|x|u$ in the Proposition \ref{polar_decomposition_prop} is said to be the \textit{polar decomposition}.
	\end{defn}

	\begin{definition}\label{stritly_pos_defn}\cite{pedersen:ca_aut}
		We say that an element $h$ is a $C^*$-algebra is \textit{strictly positive} if $\phi\left( h\right)>0$ for any nonzero positive linear functional $\phi$ on $A$.
	\end{definition}
	\begin{thm}\label{stritly_pos_thm}\cite{pedersen:ca_aut}
		Let $A$ be a $C^*$-algebra. The following conditions are equivalent:
		\begin{enumerate}
			\item [(i)] there is a strictly positive element $h$ in $A_+$,
			\item[(ii)] there is an element $h$ in $A_+$ such that $\left[ h\right] = 1$ in $A''$,
			\item[(iii)] there is a countable approximate unit for $A$.
		\end{enumerate}
	\end{thm}
	\begin{proposition}\label{stritly_pos_prop}\cite{blackadar:ko}
		Let $A$ be a $C^*$-algebra, and $h \in A_+$. Then $h$ is strictly
		positive if and only if $hA$ is dense in $A$.
	\end{proposition}
	\begin{definition}\label{simple_ca_defn}\cite{pedersen:ca_aut}
		A $C^*$-algebra $A$ is said to be \textit{simple} if $0$ and $A$ are its only closed ideals. In this book we consider only simple $C^*$-algebras of compact operators $\K\left(\H \right)$ where $\H = \C^n$ ($n \in \N$) or $\H = \ell^2\left( \N\right)$. 
	\end{definition}
\begin{lemma}\label{haus_simple_lem}\cite{rae:ctr_morita}
Suppose 
$A$
is 
a 
$C^*$-algebra 
with 
Hausdorff 
spectrum 
$\sX$, 
and 
$P_x$
is 
the 
primitive 
ideal 
corresponding 
to $x\in\sX$.
Then 
the 
quotient $A\left(x\right)\bydef A/ P_x$
is 
a 
simple 
$C^*$-algebra 
with, 
up 
to 
equivalence, 
a 
unique 
irreducible 
representation. 
\end{lemma}

	\begin{thm}\label{pedersen_ideal_thm}  \cite{pedersen:ca_aut} 
		For each $C^*$-algebra $A$ there is a dense hereditary ideal $K(A)$,
		which is minimal among dense ideals.
		
	\end{thm}
	\begin{proof}
		Let $K(]0, \infty [)$ denote the set of continuous functions on $]0, \infty [$ with 
		compact support and define 
		\be\label{pedersen_k0_eqn}
		K\left( A \right)_0 \bydef \left\{f\left(x\right) \left|x \in A_+, \quad f \in K(]0, \infty [) \right.\right\}.
		\ee
		Let 
		\be\label{pedersen_k_plus_eqn}
		K\left( A \right)_+ \bydef \left\{x \in A_+ \left|x \le \sum_{j = 1}^nx_j, \quad x_j \in  	K\left( A \right)_0\right.\right\}, 	
		\ee
	so that $	K\left( A \right)_+$ is the smallest hereditary cone  containing $K\left( A \right)_0$. If $K(A)$ 
		is  the algebraic  $\C$-linear span of $K(A)_+$ then $K(A)$,
		which is minimal among dense ideals. The full  proof is  described in \cite{pedersen:ca_aut}.
			\end{proof}
			
	\begin{defn}\label{pedersen_ideal_defn}\cite{blackadar:ko}
		The ideal $K\left( A\right) $ from the theorem \ref{pedersen_ideal_thm} is said to be the {\it Pedersen's ideal of $A$}. 
	\end{defn}
	\begin{remark}\label{pedersen_ideal_rem}
	There is an explained in \cite{pedersen:mea_c} alternative definition  Pedersen's ideal. We think of the elements of $A$ as operators on its universal Hilbert space and denote by $A''$ the weak closure of $A$ (cf. \ref{weak_topology_defn}). A projection $p \in A''$ is called \textit{majorized} (relative to $A$) if there exists $b \in A_+$ such that $p \le b$. We let $\left[a\right]$ denote the range projection of any operator $a$ on a Hilbert space and define
	\be\label{pedersen_ideal_eqn}
	\begin{split}
	K\left(A \right)_0  \bydef \left\{\left. a \in A_+ \right|  \left[a\right] \le b\right\},\\
		K\left(A \right)_+  \bydef \left\{\left. a \in A_+ \right|\exists a_j \in K\left( A\right)_0 , ~ j = 1,...,n, ~ a \le \sum a_j \right\}.
	\end{split}
	\ee
	$K\left(A \right)_0$ is a set of operators in $A_+$ with majorized supports, and $K\left(A \right)_+$ is the smallest ideal containing $K\left(A \right)_0$. The definition of $K\left(A \right)_0$ may be rephrased without reference to $A''$ as follows: $a \in K\left(A \right)_0$ if there exists $b\in A_+$ such that for all functions $\varphi$ continuous functions on the spectrum of $a$, $~0\le\varphi\le 1$ we have $\varphi\left( a\right) \le b$.
Pedersen's ideal $K\left(A\right)$ is the $\C$-linear span of $K\left(A \right)_0$.
	\end{remark}

\begin{remark}\cite{pedersen:ca_aut} 
One has
\bea\label{peder_k_eqn}
K\left( \K\right) = \left\{\left. a \in \K\right| a  \text{ is a finite rank operator}\right\},\\
\label{peder_c0_eqn}
K\left(C_0\left(\sX \right)  \right) = C_c\left(\sX \right).
\eea

	\end{remark}	

	\begin{theorem}\label{noncom_tietze_thm}\cite{apt_mult}
		Let $\pi$ be a surjective morphism between separable $C^*$-algebras $A$ and $B$. Then $\pi$ extends to a surjective morphism of $M\left(A \right)$ onto $M\left( B\right)$.  
	\end{theorem}
	\begin{remark}\label{noncom_tietze_rem}
		The Theorem \ref{noncom_tietze_thm} can be regarded as noncommutative Tietze's extension theorem \ref{tietze_ext_thm}, see \cite{apt_mult} for details.
	\end{remark}
	\begin{empt}\cite{pedersen:ca_aut}
		Let $\H$ be a Hilbert space and $M$ a subset of $B\left(\H \right)_{\text{sa}}$ (where $B\left(\H \right)_{\text{sa}}\subset B\left(\H \right)$ is the $\R$-space of self-adjoint operators). The monotone  sequential closure of $M$ is defined as the smallest class $\mathscr B\left(M\right)$ in $B\left(\H \right)_{\text{sa}}$ that  contains $M$ and contains the strong limit of each monotone (increasing or  decreasing) sequence of elements from $\mathscr B\left(M\right)$.
	\end{empt}
	\begin{lemma}\cite{pedersen:ca_aut}
		Each countable  subset of  $\mathscr B\left(M\right)$ lies in the monotone sequential  closure of a separable subset of $M$. 
	\end{lemma}
	\begin{theorem}\cite{pedersen:ca_aut}
		Let $A$ be a $C^*$-subalgebra of $B\left(\H \right)$. Then $\mathscr B\left(A_{\mathrm{sa}}\right)$ is the self-adjoint part of a C*-algebra. 
	\end{theorem}

	
	\begin{empt}\label{von_nemann_op_empt}\cite{pedersen:ca_aut}
		For each $x\in B(\H)$ we define the {\it range projection} of $x$ (denoted by $[x]$) as projection onto a closure of $x\H$. If $x\ge 0$ then the sequence $\left(\left((1 /n) +x\right)^{-1}x\right)$ is monotone increasing to $[x]$.  If $p$ and $q$ are projections then $p \vee q = [p + q]$ and thus $p \wedge q = 1 - \left[2 - \left(p+q\right)\right]$. Similarly we have $p \setminus q = p - p\wedge q$. Since $[x]H$ is the orthogonal complement of the null space of $x^*$ we have $[x]=[xx^*]$. If $M$ is a von Neumann algebra in $B(\H)$ then $[x]\in M$ for any $x\in M$. \end{empt}

	\begin{definition}\label{central_supp_defn}\cite{blackadar:oa}
	Every projection $p$ in a von Neumann algebra $M$ has a central
	\textit{carrier} (or \textit{central support projection}) $z_p$, the smallest projection in the center
	$\mathscr Z(M)$ containing $p$ as a subprojection ($z_p$ exists since $\mathscr Z(M)$ is itself a von
	Neumann algebra).
	$$
	1- z_p = \bigvee\left\{q\in M | qMp = 0\right\}
	$$
	Projections $p$ and $q$ have nonzero equivalent subprojections if and only if
	$z_pz_q = 0$. In particular, if $q \preceq p$, then $z_q \le z_p$.
	\end{definition}

	\section{States and representations}\label{st_rep_sec}
\begin{definition}\label{representation_defn}\cite{pedersen:ca_aut}
A representation of a $C^*$-algebra $A$ is a *- homomorphism  consisting $\pi : A \to B\left( \H\right)$ where $B\left(\H \right)$ is a $C^*$-algebra of bounded operators  on a Hilbert space $\H$. 

\end{definition}
\begin{definition}\label{subrepr_defn}\cite{blackadar:oa}
	A \textit{subrepresentation} of a representation $\pi: A \to B\left(\widetilde \H\right)$ on $\H$ is the restriction of $\pi$ to a
	closed invariant subspace of $\H \subset \widetilde \H$.
\end{definition}	
	\begin{definition}\label{representation_equivalence_defn}\cite{pedersen:ca_aut}
We say that two representations $\pi_1: A \to B\left( \H_1\right)$  and $\pi_2: A \to B\left( \H_2\right)$ of a $C^*$-algebra $A$ are \textit{spatially equivalent} (or \textit{unitarily equivalent}) if there is an isometry $u$ of $\H_1$ 	onto $\H_2$ such that $u\pi_1\left(x \right)u^*  = \pi_2\left(x \right)$  for all $x$ in $A$. We say that the representations are \textit{equivalent} (or \textit{quasi-equivalent}) if there is an isomorphism $\rho$ of	$\pi_1\left( A\right)''$ onto 	$\pi_2\left( A\right)''$ such that $\rho\left(  \pi_1\left( x\right)\right) =\pi_2\left( x\right)$ for all $x$ in $A$. 
	\end{definition}
	\begin{definition}\label{state_defn}\cite{pedersen:ca_aut}
		A state of a $C^*$ -algebra $A$ is a positive functional of norm one. The set of 
		states of $A$ is denoted by $SA$ (or just $S$ if no confusion can arise).
	\end{definition}
	\begin{definition}\label{ps_defn}\cite{murphy}
		We say a state $\tau$ on a $C^*$-algebra $A$ is \textit{pure} if it has the property that
		whenever $p$ is a positive linear functional on $A$ such that $p \le \tau$, necessarily
		there is a number $t\in \left[0,1\right]$ such that $p=t\tau$.
		The set of pure states on $A$ is denoted by $PS(A)$.
	\end{definition}
	
	\subsection{GNS construction}\label{gns_constr_sec}
	\paragraph*{}
	
	Any state $\tau$ of  $C^*$-algebra $A$  induces a GNS representation  \cite{murphy}. There is a $\mathbb{C}$-valued product on $A$ given by
	\begin{equation}\label{tau_prod_eqn}
		\left(a, b\right)\bydef\tau\left(a^*b\right).
	\end{equation}
	This product induces a product on $A/\mathcal{I}_\tau$ where 
	\be\label{tau_ideal_eqn}
	\mathcal{I}_\tau\bydef\left\{\left.a \in A \ \right| \ \tau(a^*a)=0\right\}
	\ee 
	So $A/\mathcal{I}_\tau$ is a pre-Hilbert space. Denote by $L^2\left(A, \tau\right)$ the Hilbert  completion of $A/\mathcal{I}_\tau$.  The Hilbert space  $L^2\left(A, \tau\right)$ is a space of a  GNS representation  $A\to B\left(L^2\left(A, \tau\right) \right)$ (or equivalently $A\times L^2\left(A, \tau\right)\to L^2\left(A, \tau\right) $) which comes from the Hilbert norm completion of the natural action $A \times A/\mathcal{I}_\tau \to A/\mathcal{I}_\tau$. The natural map  $A \to A/\mathcal{I}_\tau$ induces the homomorphism of left  $A$-modules
	\be\label{from_a_to_l2_eqn}
	\begin{split}
		f_\tau : A \to L^2\left(A, \tau\right),\\
		a \mapsto a + \mathcal{I}_\tau
	\end{split}
	\ee
	such that $f_\tau\left(A \right)$ is a dense subspace of $L^2\left(A, \tau\right)$.
	\begin{theorem}\label{state_repr_thm}\cite{pedersen:ca_aut} 
		For each positive functional $\tau$ on a $C^*$-algebra $A$ there is a cyclic representation $\pi_{\tau}: A \to B\left( L^2\left(A, \tau\right)\right) $ with a cyclic vector $\xi_{\tau}\in L^2\left(A, \tau\right)$ such that $\left( \pi_{\tau}\left( x\right) \xi_{\tau}, \xi_{\tau} \right)= \tau\left(a \right)$ for all $x \in A$.   
	\end{theorem}
	\begin{defn}\label{gns_defn}\cite{murphy,pedersen:ca_aut}
		The given by the Theorem \ref{state_repr_thm} representation is said to be a \textit{GNS representation}. We say that the representation 
		$\pi_{\tau}$ given by the Theorem \ref{state_repr_thm} is the cyclic representation \textit{associated with} $\tau$.
	\end{defn}
	\begin{proposition}\label{state_repr_prop}\cite{pedersen:ca_aut}
		Let $\phi$ be a positive functional on a $C^*$-algebra A and let 
	$\pi_\phi$ its associated representation. For each positive functional $\psi \le \phi$, 
		there is a unique element $a \in \pi\left(A \right)'$  with $0\le a \le 1$ such that 
	$$
\pi_\psi\left( x\right) = \left(\pi_\phi\left(x \right)  a \xi_\phi, \xi_\phi \right).	
	$$
	for all $x \in A$.
	\end{proposition}
	
	\begin{empt}\label{l2_mu}
	If $C^*$-algebra $A$ is commutative then	from the Theorem \ref{meafunc_thm} it follows that the state $\tau: A = C_0\left(\sX\right)\to \C$ can be represented by the following integral
		\begin{equation}\label{hilb_integral}
			\tau\left(a\right)= \int_{\mathcal X}a \ d\mu
		\end{equation}
		where $\mu$ is a positive Borel measure.   
		In analogy with the Riemann integration, one can define the integral of a 	bounded continuous function $a$ on $\mathcal{X}$. There is a $\mathbb{C}$-valued product on $C_0\left(\mathcal X\right)$ given by
		\begin{equation}\label{mu_prod_eqn}
			\left(a, b\right)=\tau\left(a^*b\right)= \int_{\mathcal X}a^*b \ d\mu,
		\end{equation}
		hence $C_0\left(\mathcal X\right)$ is a pre-Hilbert space. Denote by $L^2\left(C_0\left(\mathcal X\right), \tau\right)$ or $L^2\left(\mathcal X, \mu\right)$ the Hilbert space completion of $C_0\left(\mathcal X\right)$. 
	\end{empt}
	\begin{lemma}\label{domi_mult_repr_lem}\cite{pedersen:ca_aut}
		Let $\psi$ and $\phi$ be positive functionals on a $C^*$-algebra $A$ and $\psi$ 	is dominated by a multiple of $\phi$. If both representations $\pi_\phi : A \to B\left( \H_\phi\right)$ and  $\pi_\psi : A \to B\left( \H_\psi\right)$  correspond $\phi$ and $\psi$ respectively   then the representation $\pi_\psi : A \to B\left( \H_\psi\right)$  is spatially 	equivalent to a subrepresentation of $\pi_\phi : A \to B\left( \H_\phi\right)$ (cf. Definition \ref{subrepr_defn}).  
	\end{lemma}
	\begin{definition}\label{orth_func_defn}(cf. \cite{pedersen:ca_aut})
	Let $A$ be  $C^*$-algebra. Two positive functionals $\phi: A \to \C$ and  $\psi: A \to \C$ are \textit{orthogonal} if
	$$
	 \left\|\phi - \psi \right\| = \left\|\phi \right\| +  \left\|\psi \right\|. 
	$$
	We write
	\bean
	\phi \perp \psi.
	\eean
	If $\pi_\phi : A \to B\left( \H_\phi\right)$ and $\pi_\psi : A \to B\left( \H_\psi\right)$ are representations which correspond to $\phi$ and $\psi$ respectively then we write
		\be\label{orth_func_eqn}
	\pi_\phi \perp \pi_\psi.
	\ee
	 
	\end{definition}

\begin{definition}\label{orth_repr_defn}\cite{blackadar:oa}
	Let $\pi$ and $\rho$ be two (nondegenerate) representations
	of a $C^*$-algebra $A$, on Hilbert spaces $\H_\pi$ and $\H_\rho$.\\
	(i) An \textit{intertwiner} of $\pi$ and $\rho$ is an operator $T\in B\left(\H_\pi, \H_\rho \right)$  with
	$$
	T\pi(x) = \rho(x)T
	$$
	for all $x \in A$. Denote by $Int\left(\pi,\rho\right)$ the set of intertwiners of $\pi$ and $\rho$.\\
	(ii) The representations $\pi$ and $\rho$ are \textit{disjoint}, written $\pi\perp\rho$ if no nonzero
	subrepresentation of $\pi$ is equivalent to a subrepresentation of $\rho$.\\
	(iii) The representation $\pi$ is \textit{subordinate} to $\rho$, written $\pi\preceq\rho$, if no nonzero
	subrepresentation of $\pi$ is disjoint from $\rho$ (i.e. every nonzero subrepresentation of $\pi$ contains a nonzero subrepresentation equivalent to a subrepresentation of $\rho$).\\
	(iv) The representations  $\pi$ and $\rho$ are \textit{quasi}-\textit{equivalent}, written  $\pi\sim\rho$, if $\pi\preceq\rho$
	and $\rho\preceq\pi$.
\end{definition}

\begin{proposition}\label{orth_repr_prop}\cite{blackadar:oa}
	Let $\pi$ be a representation of $A$ on $\H$, and $p,q$ projections in $\pi\left(A \right)'$ 
	, with central supports $z_p,z_q$. Let $\rho, \sigma$ be the subrepresentations of $\pi$ on $p\H$ and $q\H$ respectively. Then\\
	(i) $Int\left(\rho, \sigma\right)= q\pi\left(A\right)'p$,\\
	(ii) $\rho\perp\sigma$ if and only if $z_p \perp z_q$.\\
	(iii) $\rho \preceq\sigma$ if and only if $z_p\le z_q$.\\
	(iv) $\rho\sim\sigma$ if and only if $z_p = z_q$.
\end{proposition}
\begin{corollary}\label{domin_rep_cor}\cite{pedersen:ca_aut}
If $\varphi$ and $\psi$ are positive functionals on a $C^*$-algebra $A$ and $\psi$  
is dominated by a multiple of $\varphi$ then the representation $\pi_\psi: A \to B\left( \H_\psi\right)$  is spatially 
equivalent to a subrepresentation of $\pi_\varphi: A \to B\left( \H_\varphi\right)$. 
\end{corollary}
	\subsection{Irreducible representations}

	\begin{theorem}\label{irred_thm}\cite{pedersen:ca_aut}
		Let $\pi: A \to B\left(\H \right)$ be a nonzero representation of $C^*$-algebra $A$. The following conditions are equivalent:
		\begin{enumerate}
			\item [(i)] there are no non-trivial $A$-subspaces for $\pi$,
			\item[(ii)] the commutant of $\pi\left(A \right)$ is the scalar multipliers of 1,
			\item[(iii)] $\pi\left(A \right)$ is strongly dense in   $B\left(\H \right)$,
			\item[(iv)] for any two vectors $\xi, \eta \in \H$ with $\xi \neq 0$ there is $a \in A$ such that $\pi\left(a \right)\xi = \eta$,
			\item[(v)] each nonzero vector in $\H$ is cyclic for  $\pi\left(A \right)$,
			\item[(vi)]  $A \to B\left(\H \right)$ is spatially equivalent to a cyclic representation associated with a pure state of $A$.
		\end{enumerate} 
	\end{theorem}
	\begin{definition}\label{irred_defn}\cite{pedersen:ca_aut}
		Let $A \to B\left(\H \right)$ be a nonzero representation of $C^*$-algebra $A$. The representation is said to be \textit{irreducible} if it satisfies to the Theorem \ref{irred_thm}.
	\end{definition}
	
	\begin{remark}\label{irr_func_rem}
		From the condition (i) of the Theorem \ref{irred_thm} it turns out that irreducibility is the categorical property, i.e. if there is the equivalence between category of representations of $C^*$-algebras then any irreducible representation is mapped to irreducible one.  The equivalence between category of representations corresponds to the strong Morita equivalence (cf. \ref{strong_morita_sec}).
	\end{remark}
	
	\begin{definition}\label{equivalent_representation_defn}\cite{pedersen:ca_aut}
		Let $A$ be a $C^*$-algebra.
		We say that two representations $\pi_1: A \to B\left( \H_1\right)$ and $\pi_2: A \to B\left( \H_2\right)$  are \textit{spatially equivalent} (or \textit{unitary equivalent}) if there is an isometry $u$ of $\H_1$ onto $\H_2$ such that $u\pi_1\left(a \right)u^*=\pi_2\left(a \right)$ for all $a \in A$. By the \textit{spectrum} or  of $A$ we understand the set $\hat A$ of spatially equivalence classes of irreducible representations. For any $x \in \hat A$ we denote by
		\be\label{rep_x_eqn}
		\begin{split}
			\rep_x : A \to B\left( \H\right) \quad \text{OR} \quad 	\rep^A_x : A \to B\left( \H\right)
		\end{split}
		\ee
		a representation which corresponds to $x$.  Sometimes we use alternative notation of the spectrum $A\hat~\stackrel{\text{def}}{=}\hat A$.
	\end{definition}
	
	\begin{lem}\label{disj_repr_lem}\cite{pedersen:ca_aut}
Two irreducible representations $\pi_1: A \to B\left(\H_1 \right)$   and  $\pi_2: A \to B\left(\H_2 \right)$   of a 
$C^*$-algebra $A$ are either disjoint (cf. Definition \ref{orth_repr_defn}) or spatially equivalent (cf. Definition \ref{equivalent_representation_defn}). 
	\end{lem}
	
	\begin{defn}\label{universal_rep_defn}\cite{pedersen:ca_aut} 
		\cite{pedersen:ca_aut}
		Let $A$ be a $C^*$-algebra, and let $S$ be the state space of $A$. For any $s \in S$ there is an associated representation $\pi_s: A \hookto B\left( \H_s\right)$. The representation $\prod_{s \in S} \pi_s: A \hookto \prod_{s \in S} B\left(\H_s \right)$ is said to be the \textit{universal representation}. The universal representation can be regarded as $A \hookto B\left( \widetilde \H\right)$ where $\widetilde \H$ is a norm completion of the algebraic direct sum $\bigoplus_{s \in S}\H_s$. 
	\end{defn} 
	\begin{defn}\label{env_alg_defn}\cite{pedersen:ca_aut}
		Let   $A$ be a $C^*$-algebra, and let $A \to B\left(\H \right)$ be the universal representation $A \to B\left(\H \right)$. The strong closure of $\pi\left( A\right)$ is said to be   the  {\it enveloping von Neumann algebra} or  the {\it enveloping $W^*$-algebra}  of $A$. The enveloping  von Neumann algebra will be denoted by $A''$.
	\end{defn}
	\begin{prop}\label{env_alg_sec_dual_prop}\cite{pedersen:ca_aut}
		The enveloping von Neumann algebra $A''$ of a $C^*$-algebra $A$ is isomorphic, as a Banach space, to the second dual of $A$, i.e. $A'' \approx A^{**}$.
	\end{prop}
	
	\begin{thm}\label{env_alg_thm}\cite{pedersen:ca_aut}
		For each nondegenerate representation $\pi: A \to B\left(\H \right)$ of a $C^*$-algebra $A$ there is a unique  normal morphism of  $A''$ onto $\pi\left( A\right)''$ which extends $\pi$.   
	\end{thm}
	
	\begin{definition}\label{spectrum_prime_primtive_defn}\cite{pedersen:ca_aut}
		An ideal $I$ in a $C^*$-algebra $A$ is \textit{prime}  if $xAy \subset I$  implies $x\in I$  or $y\in I$
		for all $x$, $y$ in $A$.  Equivalently, $I$ is prime if $I_1I_2\subset I$ implies  $I_1\subset I$ or 
	 $I_2\subset I$ for any two (left, right, or two-sided) ideals $I_1$ and $I_2$ of $A$. 
		We say that $I$ is a \textit{primitive} ideal if $I= \ker\pi$ for some irreducible 
		representation $\pi: A \to B\left(\H \right)$. The set of prime ideals will be denoted by $\check{A}$ or $\mathrm{Prim}\left(A \right)$  and the set of primitive ideals will be denoted by $\hat A$. We say that $\check{A}$ is a \textit{prime spectrum} of $A$. The set  $\hat A$ is said to by a \textit{primitive  spectrum} or simply a \textit{spectrum} of $A$. For any $x \in \hat A$ denote  by $\rep_x: A \to B\left(\H_x\right)$ a corresponding irreducible representation. This definition complies with \ref{equivalent_representation_defn} one.
	\end{definition}
\begin{proposition}\label{primitive_prime_prop}\cite{pedersen:ca_aut}
Each primitive ideal of a $C^*$-algebra is prime. 
\end{proposition}
\begin{theorem}\label{prime_primitive_separable_thm}\cite{rae:ctr_morita}
Every prime ideal in a separable (cf. Definition \ref{top_separable_defn}) $C^*$-algebra is primitive.
\end{theorem}
	\begin{empt}\cite{pedersen:ca_aut}
		For each set $F$ in $\check{A}$ define a closed ideal 
		\be\label{ker_f_eqn}
		\ker\left( F\right) = \bigcap_{t \in F} t.
		\ee
		For each subset $I$ of $A$ define a set
		\be\label{hull_eqn}
		\mathrm{hull}\left( I\right) = \left\{t \in \check{A}~|~I \subset t \right\}.
		\ee
		Using a natural map
		\be\label{check_hat_eqn}
		\hat{A}\to \check{A}
		\ee
		one can define the $\mathrm{hull}\left( I\right)\subset \hat A$ as the preimage of $\mathrm{hull}\left( I\right)$ in $\check{A}$.
	\end{empt}
	\begin{theorem}\label{jtop_thm}\cite{pedersen:ca_aut}
		The class $\left\{\mathrm{hull}\left( I\right)~|~I \subset A\right\}$ form the closed sets for a topology on $\check{A}$. There is a bijective order preserving isomorphism between open sets in this topology and the two-sided closed ideals in $A$.
	\end{theorem}
	\begin{definition}\label{jtop_defn}\cite{pedersen:ca_aut}
		The topology on  $\check{A}$ defined in the Theorem \ref{jtop_thm} is called be the \textit{Jacobson topology}. We define the \textit{Jacobson topology} on $\hat A$ as the topology for which the 
		natural map $
		\hat{A}\hookto \check{A}
		$ is open and continuous.
	\end{definition}
	\begin{remark}
		In the following text we consider the only Jacobson topology on both $\hat{A}$ and $\check{A}$.
	\end{remark}
\begin{definition}\label{primitive_prime_spectrum_defn}\cite{pedersen:ca_aut}
The spaces $\hat A$ and $\check A$ supplied with Jacobson topology are said to be the \textit{primitive spectrum} and the \textit{prime spectrum} of $A$ respectively. Sometimes we use the \textit{spectrum} world instead of the \textit{primitive spectrum} one.
\end{definition}
	\begin{proposition}\label{oa_sc_prim_thm}\cite{pedersen:ca_aut}
		If $A$ is a separable $C^*$-algebra then every closed prime ideal is primitive. 
	\end{proposition}
\begin{lemma}\label{hered_repr_lem}\cite{pedersen:ca_aut}
Let $B$ be a hereditary $C^*$-subalgebra of $A$. For each irreducible 
representation $\pi: A \to B\left( \H\right)$  such that $B \not\subset\ker\pi$ the map $\pi|_B: B \to B\left( \pi\left(B \right)\H\right)$  is an irreducible 
representation of $B$. 

\end{lemma}
\begin{proof}
Let $\left\{u_\la\right\}$ be an approximate unit for $B$ and let $p$ be the projection on the 
closure of $\pi\left(B \right)\H$ . Then $\left\{\pi\left( u_\la\right) \right\}$ is strongly convergent to $p$. For any pair of vectors $\xi, \eta\in p\H$ with $\xi \neq 0$ there is by  an element $x \in A$ with $\pi\left(x \right)\xi   = \eta$ (cf. Theorem \ref{irred_thm}). But $u_\la x u_\la \in B$ and 
$$
\left\|\pi\left(u_\la x u_\la \right) \xi - \eta \right\|\to \left\|p\pi\left( x  \right)p \xi - \eta \right\| = 0.
$$
		Consequently $\pi\left( B\right)$  acts topologically irreducibly on $p\H$. But then it also acts 
		algebraically irreducibly, so there must be a $y$ in $B$ for which $\pi\left( y\right)\xi  = \eta$. In 
		particular,  $\pi\left( B\right)\H$  is closed and $\pi|_B: B \to \pi\left(B \right)\H$ is irreducible. 
\end{proof}
	
	\begin{proposition}\label{hered_spectrum_prop}\cite{pedersen:ca_aut}
	If $B$ is a hereditary $C^*$-subalgebra of $A$ then the map $t \mapsto t\cap B$ is a homeomorphism between $\check A\setminus \mathrm{hull}\left( B\right)$ and $\check B$, where 
		$$
		\mathrm{hull}\left( B\right) = \left\{\left. x \in \hat A~\right|~ \rep_x\left(B \right)= \{0\} \right\} .
		$$ 
		Moreover we have a commutative diagram:
		\\
		\begin{tikzcd}
\hat A \setminus \mathrm{hull}\left(B \right)\arrow[d]\arrow[r, "\approx"] & \hat B\arrow[d]\\
\check A \setminus \mathrm{hull}\left(B \right)\arrow[r, "\approx"] & \check B
\end{tikzcd}
		\\ 	
	\end{proposition}
	\begin{cor}\label{hered_closed_cor}\cite{murphy}
		Every closed ideal  of $C^*$-algebra is a hereditary $C^*$-subalgebra.
	\end{cor}
	\begin{remark}\label{ideal_open_inc_rem}
	From the Theorem \ref{jtop_thm}, the Proposition \ref{hered_spectrum_prop} and the Corollary \ref{hered_closed_cor} any closed ideal $I$ corresponds of $C^*$-algebra $A$ corresponds to  the open subset $\mathcal U \subset \hat A$ such that there is the natural homeomorphism $\hat I \cong  \mathcal U$. Moreover if $I$ is an essential ideal then $\mathcal U$ is a dense subset of $\hat A$.
	\end{remark}
\begin{theorem}\label{ideal_spectrum_thm}\cite{pedersen:ca_aut}
Let $I$ be a closed ideal of $A$.
\begin{enumerate}
	\item [(i)] The maps $\left(\pi: A \to B\left(\H \right)\right)   \mapsto \left(\pi|_I: I \to B\left(\H \right)\right)$ and \\$\left(\pi: A \to B\left(\H \right)\right)   \mapsto \left(\pi: A/I \to B\left(\H \right)\right)$ from 
	$\mathrm{Irr} A \setminus \mathrm{hull}\left( I\right)$  to $\mathrm{Irr} I$ and from $\mathrm{hull}\left( I\right)$ to $\mathrm{Irr} A/I$, respectively, induce 
	isomorphisms 	$\hat A \setminus \mathrm{hull}\left( I\right)$  onto $\hat I$ and from $\mathrm{hull}\left( I\right)$ onto $\left( A/I\right)\hat~ $ 
	\item[(ii)] The maps $t \mapsto t\cap I$ and $t \mapsto t/I$ are homeomorphisms from $\check{A}\setminus \mathrm{hull}\left(I \right)$  onto $\check{I}$ and from $\mathrm{hull}\left( I\right)$ onto $\left( A/I\right)\check~$.
	\item[(iii)] The resulting diagrams, below, are commutative:
			\newline
\begin{tikzpicture}
	\matrix (m) [matrix of math nodes,row sep=3em,column sep=4em,minimum width=2em]
	{
		\hat I   & \hat A\setminus \mathrm{hull}\left( I\right),  &\mathrm{hull}\left( I\right) &	\left( A/I\right)\hat~ \\ 
		\check I  &  \check A\setminus \mathrm{hull}\left( I\right), & \mathrm{hull}\left( I\right) & 	\left( A/I\right)\check~\\};
	\path[-stealth]
	(m-1-2) edge node [above] {$\approx $} (m-1-1)
	(m-1-3) edge node [above] {$\approx$} (m-1-4)
	(m-2-2) edge node [above] {$\approx$} (m-2-1)
	(m-2-3) edge node [above] {$\approx$} (m-2-4)
	(m-1-1) edge node [right]  {} (m-2-1)
	(m-1-2) edge node [right]  {} (m-2-2)
	(m-1-3) edge node [right] {} (m-2-3)
	(m-1-4) edge node [right] {} (m-2-4);
\end{tikzpicture}
\\ 	
	
\end{enumerate}

\end{theorem}
	\begin{proposition}\label{lift_prop}\cite{murphy}
		Let $B$ be a $C^*$-algebra. For each positive functional $\phi$ on $B$ there is a norm preserving extension of $\phi$ to a positive functional on $A$. If $B$ is hereditary this extension is unique.
	\end{proposition}
	\begin{proposition}\label{state_prop}\cite{pedersen:ca_aut}
		If $B$ is a $C^*$-subalgebra of $A$ then each pure state of $B$ can be extended to a pure state of $A$.
	\end{proposition}
	\begin{proposition}\label{sur_prop}\cite{pedersen:ca_aut}
		If $B$ is a $C^*$-subalgebra of $A$ then for each irreducible representation $\rho: B \to B\left( \H_1\right)$ of $B$ there is an irreducible representation $A \to B\left(\H \right)$ of $A$ with a closed subspace $\H_1 \subset \H$ such that $\rho_B  :B \to B\left(\H_1 \right)$ is spatially equivalent to $\rho: B \to B\left( \H_1\right)$.  
	\end{proposition}	
	\begin{thm}(Dauns Hofmann)\label{dauns_hofmann_thm}\cite{pedersen:ca_aut}
		For each $C^*$-algebra $A$ there is the natural isomorphism from the center of $M\left( A\right)$ onto the class of bounded continuous  functions on $\check{A}$. 
	\end{thm}
	\begin{theorem}\label{lower_norm_thm}\cite{rae:ctr_morita}
		Suppose that $A$ is a $C^*$-algebra and that $a \in A$.
		\begin{enumerate}
			\item [(a)] The function $x \mapsto \left\| \rep_x\left(a \right) \right\|$ is lower semi-continuous on $\hat A$; that is $$\left\{x \in \hat A~|~  \left\| \rep_x\left(a \right) \right\|\le k \right\}$$ is closed for all $k \ge 0$.
			\item[(b)] For each $k > 0$, $\left\{x \in \hat A~|~  \left\| \rep_x\left(a \right) \right\|\ge k \right\}$ is compact.
		\end{enumerate}
	\end{theorem}
\begin{definition}\label{ctr_homo_defn}\cite{fell:operator_fields} 
	A $C^*$-algebra is \textit{homogeneous of order} $n$ if every irreducible *-representation is of the same finite dimension $n$.
\end{definition}
\begin{remark}\label{ctr_homo_rem}\cite{fell:operator_fields} 
If $\sX$ is a spectrum of a {homogeneous of order} $n$ algebra $A$ then $A$ corresponds  some fibre bundle  with base space $\sX$,
\end{remark}
	\begin{proposition}\label{less_n_pi_prop}\cite{pedersen:ca_aut}
		The subset $_n\check{A}$ of $\check{A}$ corresponding to irreducible representations of $A$ with finite dimension less or equal to $n$ is closed. The set $_n\check{A}\setminus_{n-1}\check{A}$ of $n$-dimensional representations is a Hausdorff space in its relative topology.
	\end{proposition}

	\begin{defn}\label{spectral_proj_defn}\cite{reed_simon:mp_1}
		Let $A$ be a bounded self-adjoint operator and $\Om$ a Borel set of $\mathbb{R}$. If $\chi_\Om$ is the characteristic function of $\Om$ (cf. Definition \ref{top_char_f_defn}) then $P_\Om=\chi_\Om\left(A\right)$ is called the \textit{spectral projection} of $A$.
	\end{defn}
	
	\begin{empt} \cite{rae:ctr_morita}
		Let $\hat A$ be the space of primitive ideals of $A$.
		The topology on $\hat A$ always determines the ideal structure of $A$: the open
		sets $\sU$ in $\mathrm{Prim}~ A$ are in one-to-one correspondence with the ideals
		\be\label{open_ideal_eqn}
		\left.A\right|_\sU \stackrel{\mathrm{def}}{=} \bigcap \left\{\left.P\in \hat A~\right| P \notin \sU\right\}
		\ee
		and there are natural homeomorphisms $P \mapsto P\cap \sU$ of $\sU$ onto the set of primitive ideals of  $ \left.A\right|_\sU$, and
		$P \mapsto P/A_\sU$ of $\hat A \setminus \sU$ onto  the set of primitive ideals of  $A/A_\sU$. 
		When $\hat A$
		is a (locally compact) Hausdorff space $T$, we can localize at a point $t$ ‚ that is,
		examine behavior in a neighborhood of $t$ ‚ either by looking at the ideal $A_\sU$
		corresponding to an open neighborhood $\sU$ of $t$, or by passing to the quotient
		\be\label{closed_ideal_eqn}
		\left.A\right|^F \stackrel{\text{def}}{=} A/\left( 
		\left.A\right|_{\hat A \setminus F} \right) 
		\ee
		corresponding to a compact neighborhood $F$ of $t$. 
	\end{empt}
	\begin{remark}
		Clearly there are natural injective  $\left.A\right|_\sU \to A$ and surjective $A \hookto \left.A\right|^F$ $*$-homomorphisms,
		
	\end{remark}
	
	\begin{remark}
		The notations  \ref{open_ideal_eqn} and \ref{closed_ideal_eqn} slightly differs from  \cite{rae:ctr_morita}. Here we write $	\left.A\right|_\sU $ and $\left.A\right|^F$ instead of $\left.A\right._\sU $ and $\left.A\right.^F$ respectively.
	\end{remark}
	
	
	\begin{definition}\label{atomic_repr_defn}\cite{pedersen:ca_aut}
		Let $A$ be a $C^*$-algebra with the spectrum $\hat A$. We choose for any $t \in \hat A$ a pure state $\phi_t$ and  associated representation $\pi_t: A \to B\left(\H_t\right)$.
		The representation 
		\be
		\pi_a = \bigoplus_{t \in \hat A} \pi_t \quad \text{on the closure } \H_a \text{ of an algebraic direct sum}\quad  \bigoplus_{t \in \hat A} \H_t
		\ee
		is called the (reduced) \textit{atomic representation} of $A$. Any two atomic representations are unitary equivalent and any atomic representation of $A$ is faithful and nondegenerate  (cf.  Definitions \ref{faithful_representation_defn}, \ref{nondegenerate_repr_defn} and \cite{pedersen:ca_aut}).
	\end{definition}

\begin{lemma}\label{atomic_uni_lem}\cite{pedersen:ca_aut}
If  $\pi_a: A \hookto B\left( \H_a\right)$ is the atomic representation then one has
$$
\pi_a\left( A\right)'' = \prod_{x\in \hat A} B\left( \H_x\right). 
$$
\end{lemma}
\begin{theorem}\label{vn_ext_thm}\cite{pedersen:ca_aut}
For each nondegenerate representation $A \to B\left(\H \right)$  of a $C^*$-algebra 
$A$ there is a unique normal morphism $\pi$ of $A''$ onto $\pi\left(A\right)''$ which extends $\pi$. 
\end{theorem}
	
		\section{Hilbert modules and compact operators}\label{hilbert_modules_chap}
	\begin{definition}\label{banach_non_defn}\cite{rae:ctr_morita}
		A left  $A$-module $X$ \textit{Banach} $A$-\textit{module} if $X$ is a Banach space and $\left\| a \cdot x\right\| \le \left\|a \right\| \left\| x\right\|$ for all $a\in A$ and $x\in\sX$.
		A  \textit{Banach} $A$-{module} is \textit{nondegenerate}  
		is nondegenerate if $\text{span}\left\{\left.a \cdot x\right|a\in A\quad x\in \sX\right\}$
		is dense in $X$. We then have \\$a_\la\cdot x \to x$
		whenever $x\in\sX$ and $\left\{a_\la\right\}$
		is a bounded approximate identity for $A$. 
	\end{definition}
	
	\begin{proposition}\label{banach_non_prop}\cite{rae:ctr_morita}
		Suppose that $X$ 
		is a nondegenerate Banach $A$-module. Then 
		every element of $X$ 
		is of the form $a\cdot x$
		for some $a\in A$ and $x\in X$.
	\end{proposition}
	
	\begin{definition}\label{hilbert_module_defn}Paschke \cite{Paschke:73}, Rieffel \cite{Rieffel:74a}
		Let~$B$ be a $C^*$-algebra.  A \emph{pre-Hilbert $B$-module} is a right
		$B$-module~$X$ (with a compatible $\C$-vector space structure),
		equipped with a conjugate-bilinear map (linear in the second variable)
		$\left\langle{\blank},{\blank}\right\rangle_B\colon X\times X\to B$ satisfying
		\begin{enumerate}
			\item[(a)] $\left\langle{x},{x}\right\rangle_B\ge0$ for all $x\in X$;
			\item[(b)] $\left\langle{x},{x}\right\rangle_B=0$ only if $x=0$;
			\item[(c)] $\left\langle{x},{y}\right\rangle_B=\left\langle{y},{x}\right\rangle_B^\ast$ for all $x,y\in X$;
			\item[(d)] $\left\langle{x},{y\cdot a}\right\rangle_B=\left\langle{x},{y}\right\rangle_B\cdot a$ for all $x,y\in X$, $a\in B$.
		\end{enumerate}
		The map $\left\langle{\blank},{\blank}\right\rangle_B$ is called a \emph{$B$-valued inner product
			on~$X$}.
		Following equation
		\be\label{hilbert_module_norm_eqn}
		\|x\|=\|\left\langle{x},{x}\right\rangle_B\|^{\nicefrac{1}{2}}
		\ee defines a norm on~$X$.
		If~$X$ is complete with respect to this norm, it is called a $C^*$-\emph{Hilbert
			$B$-module} or simply a \emph{Hilbert
			$B$-module}.  
	\end{definition}
	\begin{definition}\label{full_hilb_defn}\cite{rae:ctr_morita}
		A Hilbert  $A$-module $X_A$  is a \textit{full} Hilbert $A$-module if the ideal 
		$$
		I\bydef \text{span}\left\{\left.\left\langle\xi, \eta \right\rangle_A\right| \xi, \eta\in X_A \right\} 
		$$
		is dense in $A$. 
	\end{definition}
	\begin{remark}\label{full_hilb_rem}
		Suppose that $A$ is a $C^*$-algebra and that $p$
		is a projection in $A$ 
		(or $M(A)$ ). Following facts are proven in the Example 2.12 of \cite{rae:ctr_morita}.
		Then $Ap \bydef \left\{\left.ap\right| a\in A\right\}$
		is a Hilbert $pAp$ module 
		with inner product 
		\be\label{full_hilb_eqn}
		\left\langle ap, bp \right\rangle_{pAp} \bydef pa^*bp
		\ee
		This Hilbert module is full. Similarly, $pA$ 
		is a Hilbert $A$-module which is full over the ideal $\overline{ApA}\bydef \overline{\text{span}}\left\{\left.a p b\right| a,b \in A\right\}$
		generated by $p$, and $Ap$ 
		is itself a full left  Hilbert $\overline{ApA}$-module. 
	\end{remark}

	\begin{remark}\label{polarization_equality_rem}
		For any $C^*$-pre-Hilbert $X$ module, or more  precisely, for any sesquilinear form $\left\langle \cdot , \cdot \right\rangle$ the \textit{polarization equality}
		\be\label{polarization_equality_eqn}
		\begin{split}
			4 \left\langle \xi , \eta \right\rangle=\sum_{k = 0}^3i^k\left\langle \xi + i^k\eta, \xi + i^k\eta \right\rangle
		\end{split}
		\ee 
		is obviously satisfied for all $\xi, \eta \in X$.  If $\H$ is a Hilbert space then there is the following analog of the identity \eqref{polarization_equality_eqn}
		\be\label{polarization_hilb_equality_eqn}
		\\
		\left( \xi , \eta \right)_{\H} = \frac{\sum_{k = 0}^3i^k\left\| \xi + i^k\eta \right\|}{4}, \quad\forall \xi, \eta \in X
		\ee
	\end{remark}
	\begin{definition}\label{adjointable_operator_defn}\cite{matro:hcm}
		Let $X$, $Y$ be Hilbert modules over $C^*$-algebra $A$. A bounded $\C$-linear $A$-homomorphism from $X$ to $Y$ is called an \textit{operator} from $X$ to $Y$. Let $\Hom_A\left(X, Y \right)$ denote the set of all {operators} from $X$ to $Y$. If $Y = X$ then $\End_A\left(X\right)\bydef \Hom_A\left(X, X \right)$ is a Banach algebra. We say that $L\in \Hom_A\left(X, Y \right)$ 
		\textit{adjointable} if there is $L^*\in\Hom_A\left(Y, X\right)$ such that 
		$$
		\left\langle \eta, L \xi\right\rangle_A=  \left\langle L^*\eta,  \xi\right\rangle_A \quad  \forall \xi\in X,~ \eta \in Y.
		$$
		Denote by $\Hom^*_A\left(X, Y \right)\subset \Hom_A\left(X, Y \right)$ of all  {adjointable} operators. The set $\End^*_A\left(X \right)\bydef \Hom^*_A\left(X, X \right)$ is a  $C^*$-algebra.
	\end{definition}
	\begin{defn}\label{compact_a_operator_defn}\cite{pedersen:ca_aut}
		If $X$ is a $C^*$ Hilbert $A$-rigged module then
		denote by $\theta_{\xi, \zeta} \in \End^*_A\left(X \right)$   such that
		\begin{equation}\label{rank_one_eqn}
	\forall  \xi, \eta, \zeta \in X\quad		\theta_{\xi, \zeta} (\eta) = \zeta \langle\xi, \eta \rangle_X.
		\end{equation}
		The $C^*$-norm closure of  a generated by such endomorphisms ideal is said to be the {\it algebra of compact operators} which we denote by $\mathcal{K}(X)$. The $\mathcal{K}(X)$ is an ideal of  $\End^*_A(X)$.
	\end{defn}
\begin{remark}
 Also we shall use a following notation 
\be\label{rank_one_notation_eqn}
\begin{split}
	\xi\rangle \langle \zeta \stackrel{\text{def}}{=} \theta_{\xi, \zeta}: X_A \to X_A,\\
	\eta \mapsto \xi \langle\zeta, \eta \rangle_X.
\end{split}
\ee
\end{remark}
	\begin{theorem}\label{comp_mult_thm}\cite{blackadar:ko}
		Let $X_A$ is a Hilbert $A$-module.
		The $C^*$-algebra of adjointable maps  $\End^*_A\left( X_A\right)\bydef \Hom^*_A\left( X_A, X_A\right)$ is naturally isomorphic to the algebra $M\left(\K\left( X_A\right) \right)$ of multiplies of compact operators $\K\left( X_A\right)$. 
	\end{theorem}
	
	\begin{remark}
		The text of the Theorem \ref{comp_mult_thm} differs from the text of the Theorem 13.4.1 \cite{blackadar:ko}. It is made for a compatibility with this book.
	\end{remark}
	\begin{defn}\label{standard_h_m_defn}(cf. \cite{matro:hcm}) 
		The direct sum of countable  number of copies of a Hilbert module $X$ is denoted by $\ell^2\left(X \right)$. The Hilbert module  
		$\ell^2\left( A\right)$ is said to be the \textit{standard Hilbert $A$-module} over $A$. If $A$ is unital then $\ell^2\left( A\right)$ possesses the standard basis $\left\{\xi_j\right\}_{j \in \N}$.
		\begin{equation}\label{st_hilb_eqn}
			\begin{split}
				\ell^2\left( A\right) = \left\{\left\{a_n\right\}_{n \in \N}\in A^{\N}~|~\sum_{n =1}^\infty a^*_na_n < \infty \right\},\\
				\left\langle\left\{a_n\right\}, \left\{b_n\right\}\right\rangle_{\ell^2\left( A\right)}=\sum_{n =1}^\infty a^*_nb_n.
			\end{split}
		\end{equation}
	\end{defn}
	\begin{thm}\label{kasparov_stab_thm}\textbf{Kasparov Stabilization or Absorption Theorem.}\cite{blackadar:ko}
		If $X_A$ is a countably generated Hilbert $A$-module, then $X_A \oplus \ell^2\left( A\right) \cong \ell^2\left( A\right)$.
	\end{thm}
	\begin{theorem}\label{fin_hpro_thm}\cite{wegge_olsen}
		Let $A$ be a unital C*-algebra. The following conditions 
		on a right $A$-module $\E$ are equivalent: 
		\begin{enumerate}
			\item [(a)] $\E$  is projective and algebraically finitely generated.
			\item [(b)] $\E$ s isomorphic as a module to a direct summand in $A^n$ for some $n\in N$. 
			\item [(c)] $\E$  is isomorphic as a module to a closed, algebraically finitely generated submodule of $\ell^2\left( A\right)$. 
			\item [(d)] $\E\cong Q\ell^2\left( A\right) $ for some compact projection $Q\in \K\left( \ell^2\left( A\right)  \right)$. 
			\item [(e)] $\E\cong Q\ell^2\left( A\right) $ for some compact projection $Q\in \K\left( \ell^2\left( A\right)  \right)$ with $Q < P_n$.
			
		\end{enumerate}
		In particular, every algebraically finitely generated projective $A$-module can 
		be endowed with a structure as a Hilbert $A$-module.
	\end{theorem}
	\begin{corollary}\label{fin_hpro_cor}\cite{wegge_olsen}
		Any algebraically finitely generated Hilbert module over 
		a unital $C^*$-algebra is projective. 
	\end{corollary}
	
	\begin{remark}\label{fin_gen_un_equ_rem}\cite{wegge_olsen}
		If $A$ is an unital $C^*$-algebra and $X_A$ is a Hilbert $A$-module then the following conditions are equivalent
		\begin{itemize}
			\item [(a)] $X_A$ is algebraically finitely generated.
			\item[(b)]  The $C^*$-algebra $\K\left(X_A \right)$ of compact operators is unital. (cf. \cite{wegge_olsen} Remark 15.4.3).
		\end{itemize}
	\end{remark}
%
%
		

	\begin{empt}\label{hm_dual_empt}\cite{matro:hcm}
		For a Hilbert $C^*$-module $X$ over $C^*$-algebra $A$ let us denote by $X'$ the set of  all bounded $A$-linear maps  from  $X$ to $A$. 
		The formula
		\be
		\left( f \cdot a\right)\left( \xi\right) \bydef a^*f\left( x\right); \quad a \in A, \quad f \in X' , \quad \xi \in X
		\ee
		introduced the structure if right $A$-module on $X'$.
		The elements of $X'$ are called \textit{functionals} on a Hilbert module $X$. Note that there is an obvious isometric inclusion
		\be\label{to_dual_eqn}
		\begin{split}
			X \subset X',\\
			x \mapsto\left\langle x, \cdot \right\rangle.
		\end{split}
		\ee
		The space $X'$ is called the \textit{dual Banach module} of $X$.
	\end{empt}	
	\begin{definition}\label{hm_selfdual_defn}\cite{matro:hcm}
		A Hilbert module $X$ is called  \textit{self-dual} if $X\cong X'$.
	\end{definition}
	\begin{definition}\label{hm_MI_defn}\cite{matro:hcm}
		A $C^*$-algebra is said to be \textit{module}-\textit{infinite} (MI) if each countably generated Hilbert $A$-module is projective and finitely generated if and only if it is self-dual.
	\end{definition}
	
	\begin{definition}\label{hm_di_empt}\cite{matro:hcm}
		A commutative $C^*$-algebra $A = C\left(\sY \right)$ is said to be DI (\textit{divisible infinite}) if for any infinite sequence $\left\{u_j\right\}_{j\in \N}$ of norm $1\ge\left\|u_j \right|\ge C > 0$ in $A$ there exists a subsequence $j\left(k\right)$ and elements $0 < b_k \in A$ of norm 1 such that
		\begin{enumerate}
			\item [(i)] The supports $b_k$ in $\sY$ are pairwise disjoint.
			\item[(ii)] For each $k$ there exists points $y_k, ~y'_k$ such that $b_k\left(y_k\right)= 1$, $y'_k \notin \supp b_k$, $~\left\|u_{j\left(y_k\right)} \right|\ge\delta$, $~\left\|u_{j\left(y'_k\right)} \right|\ge\delta$ and the sequences $\left\{y_k\right\}, ~ \left\{y'_k\right\}$ have a common accumulation point. In particular $\sum_k b_k^s$ is a discontinuous function for any integer $s\ge 1$.
		\end{enumerate} 
	\end{definition}
	\begin{theorem}\label{hm_di_mi_thm}\cite{matro:hcm}
		If a commutative unital $C^*$-algebra $A$ is $\mathrm{DI}$, then it is  $\mathrm{MI}$.
	\end{theorem}
	
	\begin{theorem}\label{hm_sep_di_thm}\cite{matro:hcm}
		A commutative separable  unital $C^*$-algebra $A$ is $\mathrm{DI}$ if and only if its Gelfand spectrum has no isolated points.
	\end{theorem}	
	
	
	\section{Hermitian  modules and functors}
	\paragraph*{}
	In this section we consider an analogue of the $A \otimes_B - : \ _B\mathcal{M}\to _A\mathcal{M}$ functor or an algebraic generalization of continuous maps. Following text is in fact a citation of \cite{rieffel_morita}.
	
	\begin{defn}\label{herm_mod_defn}
		\cite{rieffel_morita} Let $B$ be a $C^*$-algebra. By a (left) {\it Hermitian  $B$-module} we will mean the Hilbert space $\H$ of a nondegenerate *-representation $A \rightarrow B(\H)$. Denote by $\mathbf{Herm}(B)$ the category of Hermitian  $B$-modules.
	\end{defn}
	\begin{empt}
		\cite{rieffel_morita} Let $A$, $B$ be $C^*$-algebras. In this section we will study some general methods for construction of functors from  $\mathbf{Herm}(B)$ to  $\mathbf{Herm}(A)$.
	\end{empt}

	\begin{defn}\label{corr_defn}\cite{rieffel_morita}
		Let $A$ and $B$ be $C^*$-algebras. By a {\it Hermitian  $B$-rigged $A$-module} we mean a $C^*$-Hilbert $B$ module, which is a left  $A$-module by means of $*$-homomorphism of $A$ into $\mathcal{L}_B(X)$.
	\end{defn}
	
	\begin{empt}\label{herm_functor_empt}
		Let $X$ be a Hermitian $B$-rigged $A$-module. If $V\in \mathbf{Herm}(B)$ then we can form the algebraic tensor product $X \otimes_{B_{\mathrm{alg}}} V$, and equip it with an ordinary pre-inner-product which is defined on elementary tensors by
		\begin{equation}\label{inital_hilb_prod_eqn}
			\langle x \otimes v, x' \otimes v' \rangle = \langle \langle x',x \rangle_B v, v' \rangle_V.
		\end{equation}
		Completing the quotient $X \otimes_{B_{\mathrm{alg}}} V$ by subspace of vectors of length zero, we obtain an ordinary Hilbert space, on which $A$ acts by 
		\be\label{comp_hilb_act_eqn}
		a(x \otimes v)=ax\otimes v
		\ee to give a  *-representation of $A$. We will denote the corresponding Hermitian  module by $X \otimes_{B} V$. The above construction defines a functor $X \otimes_{B} -: \mathbf{Herm}(B)\to \mathbf{Herm}(A)$ if for $V,W \in \mathbf{Herm}(B)$ and $f\in \mathrm{Hom}_B(V,W)$ we define $f\otimes X \in \mathrm{Hom}_A(V\otimes X, W\otimes X)$ on elementary tensors by $(f \otimes X)(x \otimes v)=x \otimes f(v)$.	We can define action of $B$ on $V\otimes X$ which is defined on elementary tensors by
		\begin{equation}\nonumber\label{comp_hilb_pre_act_eqn}
			b(x \otimes v)= (x \otimes bv) = x b \otimes v.
		\end{equation}
		The complete proof of above facts is contained in the Proposition 2.66 \cite{rae:ctr_morita}
	\end{empt}
	\begin{defn}\label{induced_representation_defn}\cite{rae:ctr_morita}
		Let us consider the situation \ref{herm_functor_empt}. 	The functor $X \otimes_{B} -: \mathbf{Herm}(B)\to \mathbf{Herm}(A)$ is said to be the \textit{Rieffel correspondence}. If $\pi: B \to B\left( \H_B\right)$ a representation then the representation $\rho :  A \to B\left( \H_A\right)$  given by the functor $X \otimes_{B} -: \mathbf{Herm}(B)\to \mathbf{Herm}(A)$ is said to be the \textit{induced representation}. We use following notation
		\be\label{induced_representation_eqn}
		X\text{-}\Ind^A_B\pi\bydef \rho\quad \text{ or } \quad \Ind^A_B\pi\bydef \rho\quad \text{ or } \quad	 X\text{-}\Ind\pi\bydef\rho
		\ee
		(cf. Proposition 2.66 of \cite{rae:ctr_morita}).
	\end{defn}
	
	\begin{remark}\label{left_right_rem}
		If $A$ is a $C^*$-algebra and $A \to B\left(\H \right)$ is a representation then $\left(a \xi, \eta \right)_{\H}  = \left( \xi, a^* \eta \right)_{\H}$ for each $a \in A$ and $\xi, \eta$. Taking into account (c) of the Definition \ref{hilbert_module_defn}  the equation \eqref{inital_hilb_prod_eqn} can be rewritten by following way
		
		\begin{equation}\label{hilb_prod_eqn}
			\langle x \otimes v, x' \otimes v' \rangle = \langle v,  \langle x,x '\rangle_B v' \rangle_V.
		\end{equation}
		The equation \eqref{hilb_prod_eqn} is more convenient for our purposes than \eqref{inital_hilb_prod_eqn} one. 
		
	\end{remark}

	\section{$C^*$-correspondences}\label{correspondeces_chap}
Here I follow to \cite{uni_groupoid_ca}.
	A $C^*$-\textit{correspondence} from a $C^*$-algebra~$A$
	to another $C^*$-algebra~$D$
	consists of a (right) Hilbert $D$-module~$\F$
	with a nondegenerate $*$-homomorphism~\(\varphi\) from~$A$
	to~$\mathbb{B}(\F)$,
	the $C^*$-algebra of adjointable operators on~$\F$.
	We view a $C^*$-correspondence from~$A$
	to~$D$
	as an arrow \(A\to D\)
	and usually write $A\xrightarrow{\F} D$.
	We also view~\(\varphi\)
	as a \emph{representation} of~$A$
	on~$\F$.
	Two $C^*$-correspondences
	$\F_1$)
	and~$\F_2$)
	from~$A$
	to~$D$
	are \emph{isomorphic} if there is a unitary bimodule map
	$U: \F_1 \xrightarrow{\sim} \F_2$.
	
	We write~\(\otimes\)
	for suitably completed tensor products of $C^*$-correspondences,
	and~\(\odot\)
	for the tensor product of vector spaces without any completion.
	In particular, the composite of two $C^*$-correspondences
	$A\xrightarrow{\E} B$ and $B\xrightarrow{\F} D$ is their
	(balanced) tensor product~$\E \ox_B \F$
	(see~\cite{Lance:Hilbert_modules}
	). This is
	the completion of the algebraic (balanced) tensor product
	~$\E \odot_B \F$) with respect to the $D$-valued inner product
	$$
\left\langle \xi_1\ox \eta_1, \xi_2\ox \eta_1 \right\rangle\bydef  \left\langle \eta_1, \varphi\left(\xi_1, \xi_2 \right) \right\rangle 
$$
	where $\varphi\colon B\to \mathbb{B}\left(\F \right)$  is the underlying
	homomorphism that gives the left  $B$-module structure of~$\F$;
	the notation $\E \ox_\varphi \F$ is used instead of
$\E \ox_B \F$ to highlight~$\varphi$.

	\section{Strong Morita equivalence for $C^*$-algebras}\label{strong_morita_sec}

	\paragraph*{}
	The notion of the strong Morita
	equivalence was introduced by Rieffel.
	\begin{definition}\label{strong_morita_defn}[Rieffel \cite{Rieffel:74a,rieffel_morita,Rieffel:76}]
		Let $A$ and~$B$ be $C^*$-algebras.  By an \emph{$A$-$B$-equivalence
			bimodule} (or \emph{$A$-$B$-imprimitivity
			bimodule}) we mean an $\left(B,A\right)$-bimodule which is equipped with $A$- and
		$B$-valued inner products with respect to which~$X$ is a right Hilbert
		$A$-module and a left  Hilbert $B$-module such that
		\begin{enumerate}
			\item[(a)] $\left\langle{x},{y}\right\rangle_B z = x\left\langle{y},{z}\right\rangle_A$ for all $x,y,z\in X$;
			\item[(b)] $\left\langle{X},{X}\right\rangle_A$ spans a dense subset of~$A$ and $\left\langle{X},{X}\right\rangle_B$ spans a dense
			subset of~$B$, i.e. $X$ is a full Hilbert $A$-module and a full Hilbert $B$-module.
		\end{enumerate}
		We call $A$ and~$B$ \emph{strongly Morita equivalent} if there is an
		$A$-$B$-equivalence bimodule.
	\end{definition}
	
	\begin{example}\label{imp_p_exm}\cite{rae:ctr_morita}
		Let $p$ be a projection in $M 
		(A)$. We saw in the Remark \ref{full_hilb_rem} that $Ap$
		is a full right Hilbert $pAp$-module and a full left  Hilbert $\overline{ApA}$-module. (Recall that 
		$\overline{ApA}$
		denotes the ideal generated by $p$, which is the closed span of the set $ApA$.) 
		Thus $Ap$ is an $\overline{ApA}$-$pAp$-imprimitivity 
		bimodule. 
	\end{example}
	
	\begin{rem}
			If~$_AX_B$ be an $A$-$B$ imprimitivity bimodule then there is a \textit{dual} $B$-$A$ bimodule $_B\widetilde X_A$ which is a $B$–$A$-imprimitivity bimodule (cf. \cite{rae:ctr_morita}).
	\end{rem}
	\begin{theorem}\label{morita_herm_thm} \cite{rae:ctr_morita,Rieffel:74a}
		Let~$X$ be an $A$-$B$ equivalence bimodule.  Then given by \ref{induced_representation_defn} functor $X\otimes_B\blank$
		induces an equivalence (cf. Definition \ref{category_equivalence_definition}) between the category of Hermitian  $B$ modules and
		the category of Hermitian  $A$-modules, the inverse being given by
		$\tilde{X} \otimes_A\blank$.
		
	\end{theorem}
	\begin{theorem}\label{rieffel_equiv_thm}\cite{rae:ctr_morita}
		Suppose that $X$ is an $A$–$B$-imprimitivity bimodule, and $\pi$, $\rho$
		are nondegenerate representations of $B$ 
		and $A$, respectively. Then $\widetilde X-\Ind\left(  X-\Ind \pi \right)$ is 
		naturally unitarily equivalent to $\pi$, and $X-\Ind\left(\widetilde  X-\Ind \rho \right)$ is 
		naturally unitarily equivalent to $\rho$.
		
	\end{theorem}
	\begin{corollary}\cite{rae:ctr_morita}
		If $X$ is an $A$–$B$-imprimitivity bimodule, then the inverse of the Rieffel correspondence $X-\Ind$ is  $\widetilde X-\Ind$.
	\end{corollary}
	\begin{corollary}\label{morita_irred}\cite{rae:ctr_morita}
		Suppose that $X$ is an $A$–$B$-imprimitivity bimodule, and that $\pi$
		is 
		a nondegenerate representation of $B$. Then $\widetilde  X-\Ind \pi$
		is irreducible if and only if $\pi$ 
		is irreducible. 
	\end{corollary}
	\begin{corollary}\label{rieffel_homeo_cor}\cite{rae:ctr_morita}
		If $X$ is $A$-$B$-imprimitivity bimodule then the Rieffel correspondence restricts to a homeomorphism $h_X: \mathrm{Prim}~B \xrightarrow{\approx}\mathrm{Prim}~A$.
	\end{corollary}
	\begin{definition}\label{rieffel_homeo_defn}\cite{rae:ctr_morita}
		Under the hypotheses of the Corollary \ref{rieffel_homeo_cor} the homeomorphism $h_X: \mathrm{Prim}~B \xrightarrow{\approx}\mathrm{Prim}~A$
		is said to be the \textit{Rieffel homeomorphism}.
	\end{definition}
\begin{definition}\label{stable_ca_defn}\cite{wegge_olsen}
 A $C^*$-algebra A is said to be \textit{stable} when $A\otimes \K \cong  A$. 
The \textit{stabilization} of $A$ is $A^s\bydef A\otimes \K$. 
$C^*$-algebras $A$ and $B$ are said to be \textit{stably equivalent} when $A\otimes \K\cong B\otimes \K$.
\end{definition}
	\begin{theorem}\label{stable_morita_thm}\cite{brown:stable}
		Let $B$ and $E$ be $C^*$-algebras. If $B$ and $E$ are 
		stably isomorphic, then they are strongly Morita equivalent. Conversely, 
		if $B$ and $E$ are strongly Morita equivalent and if they both 
		possess strictly positive elements, then they are stably isomorphic.	
	\end{theorem}
	
		\begin{proposition}\label{top_allows_morita_prop}\cite{connes:ncg94}
		Let a group $\Ga$  act freely and properly on a topological space  $\widetilde \sX$ and let $\sX\bydef \widetilde \sX/\Ga$. Then
		the $C^*$-algebra $C_0\left(\sX \right)$  is strongly Morita equivalent to the  product $C^*$-algebra
		$C_0\left(\widetilde \sX\right)\rtimes \Ga$.
	\end{proposition}
	\begin{remark}\label{top_allows_morita_rem}\cite{connes:ncg94}
The equivalence
	$C^*$-bimodule $E$ is easy to describe; it is given by the bundle of Hilbert spaces $\left\{\H_x\right\}_{x\in \sX}$
	over $\sX$  whose fiber at $x\in\sX$ is the $\ell^2$-space of the orbit $x\in\sX$. This yields
	the required $\left(C_0\left(\widetilde \sX \right)\rtimes \Ga, C_0\left(\sX \right) \right)$  $C^*$-bimodule.
	\end{remark}

	

\section{$C^*$-algebras of type I}
\subsection{Basic facts}
\paragraph*{}   Let $A$ be a $C^*$-algebra. For each positive $x\in A_+$ and irreducible representation $\pi: A \to B\left( \H\right)$   the (canonical) trace of $\pi(x)$ depends only on the equivalence class of $\pi$, so that we may define a function 
\be\label{ctr_hat_eqn}
\begin{split}
	\hat x : \hat A \to [0,\infty];\\
	\hat x\left(t \right) =\tr\left( \pi\left(x \right) \right) 
\end{split}
\ee
where $\hat A$ is the spectrum of $A$ (the space of equivalence classes of irreducible representations) and $\tr$ is the trace of the operator. From the Proposition 4.4.9 \cite{pedersen:ca_aut} it follows that $\hat x$ is lower semi-continuous function  in the Jacobson topology (cf. Definition \ref{jtop_defn}).

\begin{lemma}\label{ctr_para_sum}\cite{rae:ctr_morita}
	Let $A$ be a $C^*$-algebra whose spectrum $\sX$ is paracompact (and hence Hausdorff). Suppose that $\left\{\sU_\a\right\}_{\a\in\mathscr A}$ is a locally finite cover of $\sX = \cup_{\a\in\mathscr A} \sU_\a$ by relatively 
	compact open sets, $\left\{\phi_\a\right\}$
	is a partition of unity subordinate to $\left\{\sU_\a\right\}$, and $\left\{a_\a\right\}$ is a set of elements in $A$ parameterized by $\mathscr A$. If there is a function $f\in C_0\left(\sX\right)$ such 
	that $\left\| \rep_{ x}\left(a_\a \right) \right\| \le f\left( x\right)$ for all $x$ and $\a$, then there is a unique element $a$ of $A$ such that 
	$\rep_{ x}\left(a \right)= \sum_{ {\a}\in {\mathscr A}}\phi_\a \rep_{ x}\left(a_\a \right)$
	for every $x\in\sX$
	; we write $a =\sum_{ {\a}\in {\mathscr A}}\phi_\a a_\a$.
\end{lemma}
\begin{defn}\label{continuous_trace_c_a_defn}\cite{pedersen:ca_aut} We say that element $x\in A_+$ has {\it continuous trace } if $\hat x \in C_b(\hat A)$. We say that $C^*$-algebra has {\it continuous trace } if a set of elements with continuous trace  is dense in $A_+$.
\end{defn}
\begin{remark}
	If a $C^*$-algebra $A$ has {continuous trace } then for simplicity we say that $A$ is \textit{continuous-trace} $C^*$-\textit{algebra},
\end{remark}
There are alternative equivalent definitions {continuous-trace} $C^*$-{algebras}. One of them is presented below.
\begin{definition}\label{continuous_trace_c_alt_defn}\cite{rae:ctr_morita}
	A \textit{continuous-trace} $C^*$-\textit{algebra} is a $C^*$-algebra $A$ with Hausdorff
	spectrum $\sX$ such that, for each $x_0\in\sX$ there are a neighborhood $\sU$ of $t$ and $a\in A$ such that $\rep_{ x}\left( a\right) $ is a rank-one projection for all $x \in \sU$.
\end{definition}
\begin{defn}\label{abelian_element_defn}\cite{pedersen:ca_aut}
	A positive element in $C^*$ - algebra $A$ is {\it Abelian} if subalgebra $xAx \subset A$ is commutative.
\end{defn}
\begin{defn}\label{type_I_0_defn}\cite{pedersen:ca_aut}
	We say that a $C^*$-algebra $A$ is of type $I$ if each non-zero quotient of $A$ contains a non-zero
	Abelian element. If $A$ is even generated (as $C^*$-algebra) by its Abelian elements we say
	that it is \textit{of type} $I_0$.
\end{defn}
\begin{lemma}\label{typeI_lem}\cite{pedersen:ca_aut}
	If $A$ is a $C^*$-algebra of operators acting irreducibly on a Hilbert space $\H$ such that $A\cap \K\left( \H\right) \neq \{0\}$ then $\K\left( \H\right)\subset A$ and each faithful representation is unitary equivalent with the identity representation.
\end{lemma}

\begin{prop}\label{abelian_element_proposition}\cite{pedersen:ca_aut}
	A positive element $x$ in $C^*$-algebra $A$ is Abelian if $\mathrm{dim}~\pi(x) \le 1$ for every irreducible representation $\pi:A \to B(\H)$.
\end{prop}
\begin{proposition}\label{ctr_pedersen_prop}\cite{pedersen:ca_aut}
If $A$ is a $C^*$-algebra with continuous trace  there is for 
each $x$ in $K(A)_+$ an $n < \infty$ such that $\dim \pi\left(x \right)< n$  for every irreducible 
representation $\pi: A \to B\left(\H \right)$ . Moreover, the map $x \mapsto \left( \pi \mapsto \tr \circ \pi\left( x\right) \right)$  is a faithful, positive linear 
surjection of $K\left(A \right)$ onto  $K\left(\hat A \right)$. 
\end{proposition}

\begin{defn}\label{ccr_defn}\cite{pedersen:ca_aut}
	A $C^*$-algebra is called \textit{liminaly} (or $CCR$) if $\rho\left( A\right)= \K$ for each irreducible representation $\rho: A \to B\left( \H\right)$.  
\end{defn}

\begin{theorem}\label{ctr_hat_check_thm}\cite{pedersen:ca_aut}
	Let $A$ be a $C^*$-algebra of type I. Then 
	\begin{enumerate}
		\item [(i)] $\K \subset \pi\left(A \right)$ for each irreducible representation $\pi: A \to B\left(\H \right)$  of $A$.
		\item[(ii)] $\hat A = \check{A}$ (cf. Definition \ref{spectrum_prime_primtive_defn}).
	\end{enumerate}
\end{theorem}
\begin{prop}\label{ctr_hered_prop}\cite{pedersen:ca_aut}
	Each hereditary $C^*$-subalgebra and each quotient of a 
	$C^*$-algebra which is of type $I_0$ or has continuous trace  is of type $I_0$  or has 
	continuous trace , respectively. 
\end{prop}

\begin{corollary}\label{ctr_ccr_i_cor}\cite{pedersen:ca_aut}
	Any $CCR$ $C^*$-algebra is of type I.
\end{corollary}
\begin{definition}\label{oa_haus_defn}\cite{rae:ctr_morita}
	A topological space is called \textit{almost Hausdorff} if every closed subset
	has a dense open subset which is Hausdorff in the relative topology.
\end{definition}

\begin{theorem}\label{oa_haus_prim_thm}\cite{rae:ctr_morita}
	Suppose that $A$ 
	is a $C^*$-algebra such that $\check{A}$
	is either second-countable  or almost Hausdorff (cf. Definition \ref{oa_haus_defn}). Then every prime ideal is primitive (cf. Definition \ref{spectrum_prime_primtive_defn}). 
\end{theorem}

\begin{proposition}\label{foli_ot_prop}\cite{rae:ctr_morita}
	Suppose that $A$ and $B$ are $C^*$-algebras and $X$ is an $A$-$B$-imprimitivity bimodule (cf. Definition \ref{strong_morita_defn}). 
	\begin{itemize}
		\item [(a)] If $h_X: \mathrm{Prim}~B\xrightarrow{\approx}\mathrm{Prim}~A
		$
		is the Rieffel homeomorphism and $f\in C_b(\mathrm{Prim}~A)$ then $f\cdot x = x \cdot \left( f\circ h_X\right)$  for all $x\in \sX$.
		\item[(b)]  If $A$ and $B$ have Hausdorff spectrum $\sX$, then $X$ is an imprimitivity bimodule over $\sX$ if and only if $f\cdot x=x\cdot f$
		for all $x\in X$ and $f\in C_0\left(\sX\right)$.
	\end{itemize}
\end{proposition}

\begin{proposition}\label{ctr_morita_prop}\cite{rae:ctr_morita} 
	A $C^*$-algebra $A$ with Hausdorff spectrum $\sX$ has continuous
	trace if and only if $A$ is locally Morita equivalent to $C_0\left(\sX \right)$, in the sense that each
	$x \in \sX$ has a compact neighborhood $\sV$ such that $\left.A\right|^\sV$ (cf. equation \eqref{closed_ideal_eqn}) is Morita equivalent to $C\left( \sV\right)$ 
	over $\sV$.
\end{proposition}
\begin{empt}\label{ctr_imprim_empt}
	For our research it is important the proof of the Proposition \ref{ctr_morita_prop}. Here is  the fragment of the proof described in \cite{rae:ctr_morita}.
	Suppose that $A$ 
	has continuous trace , and fix $x_0\in \sX$. Choose a compact neighbourhood $\sV$ of $x_0$ and $p\in A$ such that $\rep_x\left(p \right)$  is a rank-one projection for all $x\in \sV$ Then $\left.p\right|^\sV$  is a projection in $\left.A\right|^\sV$
	, and $\left.A\right|^\sV\left.p\right|^\sV$
	is an $\left.A\right|^\sV-\left.p\right|^\sV\left.A\right|^\sV\left.p\right|^\sV$
	-imprimitivity bimodule (cf. Example \ref{imp_p_exm}). Note that the map $f\mapsto 
	f\left.p\right|^\sV$ is an isomorphism of $C_0\left( \sV\right)$ into  $\left.p\right|^\sV\left.A\right|^\sV\left.p\right|^\sV$. On the other hand, if $a\in A$ and $x\in \sV$,  then $\rep_x\left(p\right)$ is a rank-one 
	projection in the algebra $\rep_x\left(A\right)$ of compact operators, and $\rep_x\left(pap\right) = \rep_x\left(p\right)\rep_x\left(a\right)(pap)\rep_x\left(p\right)$ 
	must be a scalar multiple $f_a\left(x\right)\rep_x\left(p\right)$ of $x \mapsto \rep_x\left(p\right)$. We claim that $f_a$ 
	is continuous, so that $f \mapsto f\left.p\right|^\sV$ is an isomorphism of $C\left(\sV\right)$ onto $\left.p\right|^\sV\left.A\right|^\sV\left.p\right|^\sV$.
	Well, for $x,y\in\sV$ we have 
	$$
	\left|f_a\left(x\right)-f_a\left(y\right)\right|= \left\| \rep_x\left(f_a\left(x\right)-f_a\left(y\right)\right) p\right\| = \left\|\rep_x\left( pap- f_a\left(y\right)p\right) \right\|.
	$$
	Since $\hat A$ is Hausdorff, for fixed $y$ the right-hand side is a continuous function of $x$ by Lemma \ref{lower_norm_thm}; since it vanishes at $y$, the left-hand side goes to 0 as $x\to y$ In other words, $f_a$ 
	is continuous, as claimed, and $\left.A\right|^\sV\left.p\right|^\sV$ 
	is an $\left.A\right|^\sV$-$C(\sV)$-imprimitivity $\sV$-bimodule. Because the actions of $C\left(\sV\right)$ on the left  and right of $\left.A\right|^\sV\left.p\right|^\sV$ coincide, it 
	is actually an $\left.A\right|^\sV$-$C(\sV)$-imprimitivity bimodule by Proposition \ref{foli_ot_prop}. 
\end{empt}
\begin{proposition}\label{ctr_loc_morita_prop}\cite{rae:ctr_morita}	  
	 Suppose that $A$ is a continuous-trace
	    $C^*$-algebra with paracompact spectrum $\sX$. Then there are compact subsets $\left\{\sV_\a \subset \sX\right\}_{\a \in \mathscr A}$ whose interiors form a
	cover $\left\{\sU_\a\right\}$ of $\sX$, such that for each $\a \in \mathscr A$, there is an $\left.A\right|^{\sV_{   \a}}-C\left( \sV_{   \a}\right)$-imprimitivity bimodule $X_\a$.
\end{proposition}

\begin{proposition}\label{ctr_ped_prop}\cite{pedersen:ca_aut}
	If $A$ is a $C^*$-algebra with continuous trace  there is for each $x \in K\left(A \right)_+$ and $n < \infty$ such that $\dim \pi\left(A\right)  \le n$ for every irreducible representation $\pi: A \to B\left(\H \right)$. Moreover, the map $x \mapsto \hat x$ (cf. \eqref{ctr_hat_eqn}) is a faithful, positive linear surjection of $K\left(A\right)$ onto $K\left( \hat A\right)$.  
\end{proposition}




\begin{lem}\label{sep_sec_cou_lem}\cite{pedersen:ca_aut}
	If $A$ is a separable algebra then $\check{A}$ is second-countable .
\end{lem}
\begin{prop}\label{continuous_trace_c_a_proposition}\cite{pedersen:ca_aut}
	Let $A$ be a $C^*$ - algebra with continuous trace . Then
	\begin{enumerate}
		\item[(i)] $A$ is of type $I_0$;
		\item[(ii)] $\hat A$ is a locally compact Hausdorff space;
		\item[(iii)] For each $t \in \hat A$ there is an Abelian element $x \in A$ such that $\hat x \in K(\hat A)$ and $\hat x(t) = 1$.
	\end{enumerate}
	The last condition is sufficient for $A$ to have continuous trace .
\end{prop}

\begin{theorem}\label{ctr_big_thm}\cite{pedersen:ca_aut}
	Let $A$ be a $C^*$-algebra of type I. Then $A$ contains an essential ideal which has continuous trace . Moreover, $A$ has an essential composition series $\left\{ I_\al~|~0\le \al \le \bt \right\}$ such that $I_{\al+1} / I_\al$ has continuous trace  for each $\al < \bt$.
\end{theorem}

\begin{theorem}\label{ctr_fin_bundle_thm}\cite{fell:operator_fields}
	Every {homogeneous $C^*$-algebra of order} $n$ (cf. Definition \ref{ctr_homo_defn}) is isomorphic with some $C_0\left(B \right)$, where $B$ is a fibre bundle with base space $\hat A$, fibre space $\mathbb{M}_n\left(\C\right)$ , and group $PU\left(n \right)$. 
\end{theorem}
\begin{lemma}\label{ctr_sep_haus_lem}\cite{chun-yen:separability}
Let $A$ be a separable C*-algebra. Then its spectrum $\sX$
is metrizable if and only if $\sX$ is Hausdorff. In this case, A is a $CCR$.
\end{lemma}
\begin{thm}\label{ctr_hausdoff_thm}\cite{kaplansky:certain}
	Let $A$ be a $C^*$-algebra in which for every primitive ideal $P$,
	$P$ is finite-dimensional and of order independent of $P$. Then the structure space of $A$ is Hausdorff.
\end{thm}
\begin{remark}\label{ctr_open_res_rem}
	It is  shown n \cite{rae:ctr_morita} that the Jacobson topology on of the primitive spectrum (cf. Definition \ref{jtop_defn}) always determines the ideal structure of $A$: the open
	sets $\sU$ in $\check{A}$ are in one-to-one correspondence with the ideals
	\be\label{ctr_open_id_eqn}
	\sU \leftrightarrow 	A|_{\mathcal U}.
	\ee
	On the other hand any separable $C^*$-algebra $A$ with Hausdorff spectrum is $CCR$-algebra (cf. Lemma \ref{ctr_sep_haus_lem}). From the Theorem \ref{ctr_hat_check_thm} it turns out that the prime spectrum of $A$ coincides with its primitive one (cf. Definition \ref{spectrum_prime_primtive_defn}).  So  the equation \ref{ctr_open_id_eqn} establishes the one-to-one correspondence with the ideals of $A$ and open sets of the spectrum of $A$.
\end{remark}

\begin{lemma}\label{ctr_fact_lem}\cite{rae:ctr_morita}
	Suppose that $A$ is a $C^*$-algebra with Hausdorff spectrum $\mathcal{X}$ and that $\mathcal U$ is an open subset of $\hat A$. If $\left.A\right|_\sU$ is given by \eqref{open_ideal_eqn} then $\left.A\right|_\sU$  is the closure of the space
	$$
	\left\{f\cdot a~|~a\in A \text{ and } f \in C_0\left(\mathcal{X}\right) \text{ vanishes off }~ \mathcal U \right\},
	$$
	or, equivalently the closure of the space
	$$
	C_0\left(\mathcal{U}\right)A.
	$$
\end{lemma}
\begin{lemma}\label{hausdorff_spectrum_lem}\cite{rae:ctr_morita}
	Suppose $A$ is a $C^*$-algebra with Hausdorff spectrum $\mathcal{X}$.
	\begin{itemize}
		\item [(a)] If $a, b \in A$ and $\mathfrak{rep}_x\left(a \right)=  \mathfrak{rep}_x\left(b \right)$ for every $x \in  \mathcal{X}$, then $a = b$.
		\item[(b)] For each $a \in A$ the function $x \mapsto \left\|\mathfrak{rep}_x\left(a \right) \right\|$ is continuous on  $\mathcal{X}$, vanishes at infinity and has sup-norm equal to $\left\| a\right\|$. 
	\end{itemize}
\end{lemma}

\begin{proposition}\label{ctr_bundle_prop}\cite{rae:ctr_morita}
Let $\H$ be a separable infinite-dimensional Hilbert space. If $A$
is a stable separable continuous-trace $C^*$-algebra with spectrum $\sX$, there is a locally
trivial bundle $\left( X, \pi,\sX\right) $ with fibre $\K(\H)$ and structure group $\Aut \left( \K(\H)\right)$ such that $A$ is $C_0\left( \sX\right)$ -isomorphic to the space of sections $\Ga_0\left(X \right)$ (cf. equation \eqref{top_ub_eqn}).
\end{proposition}

\subsection{Fields of operators}\label{f_op_sec}
\paragraph*{}
Let $\sX$ be a locally compact Hausdorff space called the base space; and for each $x$ in
$\sX$, let $A_x$ be a (complex) Banach space.
\begin{definition}\label{operator_fields_defn}\cite{fell:operator_fields}
 A \textit{vector field (with values in the $\left\{A_x\right\}$)} is a function
$a$ on $\sX$ given by $x \mapsto a_x$ such that $a_x\in A_x$ for each $x$ in $\sX$. Obviously the vector fields form a complex linear
space. If each $A_x$ is a $*$-algebra, then the vector fields form a $*$-algebra under the point-wise
operations; in that case the vector fields will usually be referred to as \textit{operator fields}.
Here, either each $A_x$ will be a Hilbert space or each $A_x$ will be a $C^*$-algebra.
\end{definition}

\begin{definition}\label{operator_fields_continuity_defn}\cite{fell:operator_fields}
	A \textit{continuity structure for} $\sX$ \textit{and the} $\left\{A_x\right\}$ is a linear space $\sF$ of vector fields on $\sX$, with values in the $\left\{A_x\right\}$ , satisfying:
	\begin{enumerate}
		\item[(a)] If $a \in \sF$, the real-valued function $x \mapsto \left\| a_x\right\|$  is continuous on $\sX$.
		\item[(b)] For each $x$ in $\sX$,  $\left\{\left. a_x~\right|~ a \in \sF\right\}$ is dense in $A_x$.
		\item[(c)] If each $A_x$ is a $C$*-algebra, we require also that
		$\sF$ is closed under pointwise multiplication and involution.
	\end{enumerate}
	
\end{definition} 
\begin{definition}\label{op_cont_fields_defn}\cite{fell:operator_fields}
	A vector field $a$ is \textit{continuous (with respect to $\sF$)} at $x_0$, if for each $\eps>0$
	there is an element $b$ of $\sF$ and a neighborhood $\sU$ of $x_0$ such that $\left\|a_x- b_x\right\|< \eps$ for all
	$x$ in $\sU$. We say that $a$ is \textit{continuous} on $\sX$ if it is continuous at all points of $\sX$.
\end{definition} 
\begin{lemma}\label{op_cont_uni_lem}\cite{fell:operator_fields}
	If a sequence of vector fields $\left\{a_n\right\}$ continuous (with respect to $\sF$) at $x_0$ converges uniformly on $\sX$ to a vector field $a$, then $a$ is continuous at  $x_0$ (with respect to $\sF$). 
\end{lemma}
\begin{lemma}\label{op_cont_con_lem}\cite{fell:operator_fields}
	If a vector field $a$ is continuous with respect to $\sF$ at $x_0$, then $x \mapsto  \left\| a_x\right\|$ is
	continuous at $x_0$.
\end{lemma}
\begin{lemma}\label{op_cont_module_lem}\cite{fell:operator_fields}
	The vector fields continuous (with respect to $\sF$) at $x_0$ form a linear space,
	closed under multiplication by complex-valued functions on $\sX$ which are continuous at $x_0$.
	If each $A_x$ is a $C^*$-algebra, the vector fields continuous at $x_0$ are also closed under pointwise multiplication and involution.
\end{lemma}
\begin{lemma}\cite{fell:operator_fields}
A vector field $x$ is continuous (with respect to $\sF$) at $t_0$ if and only if, for each $y \in \sF$, the function $t\mapsto \left\|x\left(t\right) - y\left(t\right)\right\|$ is continuous at $t_0$. 
\end{lemma}
\begin{lemma}\cite{fell:operator_fields}
If a sequence of vector fields $\left\{x_n\right\}$ continuous (with respect to $\sF$) at $t_0$ converges uniformly on $T$ to a vector field $x$, then $x$ is continuous at $t_0$ (with respect to $\sF$). 
\end{lemma}
\begin{lemma}\cite{fell:operator_fields}
For every $t$ in $T$, and every $\a$ in $A_t$, there is a vector field $x$, continuous on $T$ with respect to $\sF$, such that $x\left(t\right)=\a$.
\end{lemma}
\begin{defn}\cite{fell:operator_fields}
If $\sF'$ is another continuity structure for $T$ and the $\left\{A_t\right\}$, then we shall say that $\sF$ and $\sF'$ are \textit{strictly equivalent} $\left(\sF \sim \sF'\right)$ if, for all $t$ in $T$, a vector field is continuous at $t$ with respect to $\sF$ if and only if it is so with respect to $\sF'$. 
\end{defn}
\begin{lemma}\cite{fell:operator_fields}
If $\sF'$ is another continuity structure for $T$ and the $\{A_t\}$, and if there exists a family $\sG$ of vector fields such that 
\begin{enumerate}
	\item [(i)] each $x$ in $\sG$ is continuous on $T$ with respect to both $\sF$ and $\sF'$, and
	\item[(ii)]   for each $t$ in $T$, $\left\{\left. x(t)\right| x \in \sG\right\}$ is dense in $A_t$,
\end{enumerate}
then $\sF \sim \sF'$.
\end{lemma}
\begin{defn}\label{full_algebra_operator_fields_defn}\cite{fell:operator_fields}
	A \textit{full algebra of operator fields} is a family $A$ of operator fields on $T$ 	satisfying: 
		\begin{enumerate}
			\item[(a)] $A$ is a *-algebra, i.e., it is closed under all the point-wise algebraic operations;
			\item[(b)] 	for each $x$ in A, the function $t\mapsto \left\|x\left(t\right) \right\|$ is continuous on $T$ and vanishes at infinity; 
			\item[(c)] for each $t$ , $\left\{x(t)| x\in A\right\}$ is dense in $A_t$;
			\item[(d)] A is complete in the norm $\left\| x\right\|\bydef \sup_{t\in T} \left\| x(t)\right\|$.
		\end{enumerate}
\end{defn}

\begin{empt}\label{top_cs_not_empt}\cite{fell:operator_fields}
	A full algebra of operator fields is evidently a continuity structure. If $\sF$ is any continuity 
	structure, let us define $C_0\left(\sF\right)$ to be the family of all vector fields $x$ which are continuous on $T$ with respect to $\sF$, and for which $t\mapsto \left\|x\left(t\right) \right\|$ vanishes at infinity. 
\end{empt}
\begin{lemma}\label{top_full_oaf_lem}\cite{fell:operator_fields}
For any full algebra $A$ of  operator fields on $T$, the following three conditions are equivalent: 
\begin{enumerate}
	\item [(i)] $A$ is a maximal full algebra of operator fields;
	\item [(ii)] $A = C_0\left(\sF\right)$ for some continuity structure $\sF$;
	\item [(iii)] $A = C_0(A)$.
\end{enumerate}

\end{lemma}
It is proven in \cite{structure_of_standard} that any $C^*$-algebra with Hausdorff spectrum is a full algebra of operator fields (cf. Definition \ref{full_algebra_operator_fields_defn}). Below we consider a simple proof of this fact.
\begin{lemma}\label{oa_haus_alg_lem}\cite{fell:operator_fields}
Any $C^*$-algebra $A$ with Hausdorff spectrum $\sX$ is a full algebra of operator fields (cf. Definition \ref{full_algebra_operator_fields_defn}.
\end{lemma}
\begin{proof}
There is a family $\left\{\rep_x\left( A\right) \right\}_{x\in \sX}$ of Banach spaces. From the Lemma \ref{hausdorff_spectrum_lem} it follows that $A$ is a continuity structure for $\sX$ and $\left\{\rep_x\left( A\right) \right\}$ (cf. Definition \ref{op_cont_fields_defn}), such that
\be\label{top_a_c0a_eqn}
A \subset C_0\left(A \right). 
\ee
If $\left\{a_x \right\}_{x\in \sX} \in  C_0\left(A \right)$ and $\eps > 0$ then there is a compact set $\sY$ such that
$$
\forall x \in \sX \setminus\sY \quad \left\|a_x  \right\|< \eps.
$$
From the Definition \ref{op_cont_fields_defn} it follows that  there is a family $\left\{\sU_x\right\}_{x \in \sX}$ of open subsets of $\sX$ such that for all $x \in \sX$ one has $x \in \sU_x$ and there is $b_x \in A$ which satisfies to the following condition
$$
\forall y \in \sU_x \quad \left\|a_y - \rep_y\left( b_x\right)  \right\| < \eps.
$$
If $\sum_{j=1}^n f_{x_j}$ is {covering sum} for $\sY$ {dominated} by the family $\left\{\sU_x\right\}$ (cf. Definition \ref{top_covering_sum_defn}) and 
$$
b \bydef \sum_{j=1}^n f_{x_j}b_{x_j}\in A
$$
then 
$$
\forall x \in \sX \quad \left\|a_x - \rep_x\left( b\right)  \right\| < \eps,
$$
or, equivalently $\left\|\left\{a_x \right\} - b  \right\| < \eps$. Since $A$ is $C^*$-norm closed one has
$$
\left\{a_x \right\} \in A.
$$
It turns out that
$
C_0\left(A \right) \subset A$ and taking into account \eqref{top_a_c0a_eqn} we conclude that $A = C_0(A)$. From the Lemma \ref{top_full_oaf_lem} it follows that $A$ is a full algebra of operator fields.
\end{proof}

\subsection{$C^*$-algebras as cross sections and their multipliers}\label{cross_sections_sec}
\paragraph*{}	
Here I follow to \cite{apt_mult}. The following definition is a specialization of the Definition \ref{op_cont_fields_defn}.

For any $C^*$-algebra $A$ denote by $M\left( A\right)_\bt$  the algebra $M\left( A\right)$ of multipliers with the strict topology (cf. Definition \ref{strict_topology_defn}).
\begin{definition}\label{ctr_crooss_alg_defn}\cite{apt_mult}
	A bounded  section $a$ of the fibered space $\left\{\mathcal X, M\left( A_x \right) \right\}$ is said to be \textit{strictly continuous (with respect to)} $\mathscr F$ if for each $\eps > 0$ and for each $c \in \mathscr F$ there is an element $b \in \mathscr{F}$ and an neighborhood $\mathcal U$ of $x_0$ such that $\left\|c _x \left( a _x - b_x\right)   \right\|+\left\|\left( a _x - b_x\right)c _x    \right\|< \eps$ for every $x$ in $\mathcal U$.
\end{definition}

We denote by $C_b\left(\mathcal X, M\left( A_x\right)_\bt,  \mathscr F\right)$ the set of all bounded strictly continuous cross sections in $\left\{\mathcal X, M\left( A_x \right) \right\}$
\begin{theorem}\label{cross_mult_thm}\cite{apt_mult}
	There is the natural $*$-isomorphism $$M\left( C_0\left(\mathcal X, A_x, \mathscr F\right)\right) \cong C_b\left(\mathcal X, M\left( A_x\right)_\bt, \mathscr F\right).$$
\end{theorem}
 
\begin{empt}\label{ctr_crooss_alg_empt}\cite{dixmier_ca}
	Let $\sX$ be a locally compact space,  let $\left\{A_x\right\}_{x \in \sX}$ be a family of $C^*$-algebras, and let $\sF$ be {continuity structure for} $\sX$ {and the} $\left\{A_x\right\}$. Let 	$A= C_0\left(\mathcal X, A\left( x\right), \mathscr F\right)$ be the set of all continuous sections of $\left\{\mathcal X, A\left( x\right) \right\}$ vanishing at infinity. Then $A$ is an involutive algebra. For $a \in A$, put
	$$
	\lVert a \rVert = \sup_{x \in \sX} \lVert a\left(x\right) \rVert.
	$$
	It is immediate that, for this norm, $A$ is a C*-algebra which we will call  \textit{	the $C^*$-algebra defined by $\left(\mathcal X, A\left( x\right), \mathscr F\right)$}.
\end{empt}
For any $C^*$-algebra $A$ denote by $M\left( A\right)_\bt$ is the algebra $M\left( A\right)$ with the strict topology.

We denote by $C_b\left(\mathcal X, M\left( A\left( x\right)\right)_\bt,  \mathscr F\right)$ the set of all bounded strictly continuous cross sections in $\left\{\mathcal X, M\left( A\left( x\right) \right) \right\}$

\begin{theorem}\label{ctr_as_field_thm}\cite{dixmier_ca}
	Let $A$ be a liminal $C^*$-algebra whose spectrum $\sX$ is Hausdorff. Let  $\left(\sX, \left\{A_x\right\}_{x \in \sX}, \sF\right)$ be the continuous field of non-zero elementary $C^*$-algebras over $\sX$ defined by $A$. Let $A'$ be the  $C^*$-algebra
	defined by $\left(\sX, \left\{A_x\right\}_{x \in \sX}, \sF\right)$. For every $a \in A$, let $a'  = \left\{a\left(x\right)\in A_x\right\}_{x \in \sX}$ be the element of  $C_0\left( \sX, \left\{A_x\right\}_{x \in \sX}, \sF\right)$  defined by $a$.  Then $a' \in A'$ and $a \mapsto a'$ is an isomorphism of $A$ onto $A'$.
\end{theorem}
\begin{empt}\label{cross_mult_empt}\cite{apt_mult}
	Theorem \ref{cross_mult_thm} allows us (in principle) to determine the multipliers of any liminal algebra $C^*$-algebra $A$ with Hausdorff spectrum. Because in this case $A$ can be represented as an algebra of continuous cross sections 
\be\label{cross_mult_eqn}
	A = C_0\left(\hat A, A\left(t \right), A  \right) 
\ee
	where $A\left( t\right) = \K\left(\H_t \right)$. Since $M\left(\K\left(\H_t \right) \right) $ is equal to $B\left( \H\right)$ equipped with the strong* topology, denoted by $B\left( \H\right)_{s^*}$, we can state the following.  
\end{empt}
\begin{remark}\cite{apt_mult}
	Even in case where $A$ has only one and two dimensional representations it is not necessary true that $M\left( A\right)$ has Hausdorff spectrum. 
\end{remark}
\begin{theorem}\label{ctr_cont_thm} \cite{fell:operator_fields}
	Let $\sX$ be a locally compact Hausdorff space, and let $\left\{A_x\right\}_{x \in \sX}$ be a family of $C^*$-algebras. If $\sF$ is a {continuity structure for} $\sX$ {and the} $\left\{A_x\right\}$. If $A_x \cong \K\left(\ell^2\left( \N\right) \right)$ for all $x \in \sX$ then $C_0\left( \sX,  \left\{A_x\right\}, \sF\right)$ is a $C^*$-algebra with continuous trace .
\end{theorem}
\begin{remark}\label{ctr_cont_rem}\cite{fell:operator_fields}
	There are $C^*$-algebras with continuous trace  which do not match to the Theorem \ref{ctr_cont_thm}, i.e. these algebras cannot be represented as  $C_0\left( \sX,  \left\{A_x\right\}, \sF\right)$.
\end{remark}

%
%

\subsection{$K$-theory of $C^*$-algebras with Hausdorff spectrum}
\begin{theorem}\cite{sudo:k_ctr}\label{ctr_k_thm}
	Let 	$A \bydef C_0\left(\mathcal X, \left\{A_x\right\}_{x \in \sX}, \mathscr F\right)$ be a $C^*$-algebra of continuous fields on locally compact, paracompact Hausdorff space $\sX$ with fibers $A_x$, elementary $C^*$-algebras acting on $\C^n$ or $\ell^2\left(\C\right)$. Then for $*=0,1$,
	\be\label{ctr_k_eqn}
	\varphi_{K_*\left(C_0\left(\mathcal X\right)\right)}:K_*\left(C_0\left(\mathcal X, \left\{A_x\right\}_{x \in \sX}, \mathscr F\right)\right)\cong K_*\left(C_0\left(\mathcal X\right)\right).
	\ee
\end{theorem} 

\section{Inductive limits of $C^*$-algebras}
	\begin{definition}\label{principal_defn}\cite{takeda:inductive}
		An  injective $*$-homomorphism  $\phi: A \hookto B$ of unital $C^*$-algebras is said to be \textit{unital} if and only if $\phi\left(1_A\right)= 1_B$. 
	\end{definition}
	\begin{remark}
		In the cited text \cite{takeda:inductive} the word "principal" used instead "unital". In this book all entries of the word "principal $*$-homomorphism" are replaced with "unital $*$-homomorphism".
	\end{remark}
	\begin{definition}\label{inductive_lim_defn}\cite{takeda:inductive}
		Let $\La$ be an increasingly directed set and $A_\la$ be a $C^*$-algebra
		having an identity $1_\la$ associated with $\la \in \La$ If there exists a $C^*$-algebra $A$ with the identity $1$ and a unital isomorphism $f_\la$ of $A_\la$ into $A$ for every $\la \in \La$ such that
		$$
		f_\mu\left(A_\mu\right)\subset f_\nu\left(A_\nu\right) \quad \text{ if } \quad \nu < \mu; \quad \mu, \nu \in \La
		$$
		and that the join of $f_\la\left(A_\la\right)$ ($\la \in \La$) is uniformly dense in $A$, $A$ is called the $C^*$-\textit{inductive limit} of $A_\la$, and is denoted by $C^*\text{-}\varinjlim_\La A_\la$ or  $C^*\text{-}\varinjlim A_\la$
	\end{definition}
	\begin{theorem}\label{inductive_lim_thm}\cite{takeda:inductive}
		Let $\left\{A_\la\right\}_{\la \in \La}$ be a family of $C^*$-algebras where $\La$  denotes an increasingly directed set. If, for every $\mu, \nu$ with $\mu < \nu$, there exists a unital injective $*$-homomorphism $f_{\mu\nu}: A_\mu \to A_\nu$ satisfying
		$$
		f_{\mu\nu} = f_{\mu\la}\circ f_{\la\nu}\quad~\mathrm{where} \quad \mu < \la < \nu.
		$$
		then there exists the $C^*$-inductive limit of $A_\la$.
	\end{theorem}
\begin{remark}\label{inductive_lim_rem}
A proof of the Theorem is presented in \cite{takeda:inductive}. Roughly speaking the union $\bigcup A_\la$ has the natural structure of $*$-algebra with $C^*$-norm. The $C^*$-algebra $A$ is the $C^*$-norm completion of $\bigcup A_\la$. Moreover one has
\be\label{inductive_lim_eqn}
\forall a \in A \quad \exists \left(a_\la\right)_{\la\in \La} \in \prod A_\la\quad  a = \lim_{\la \in \La}a_\la
\ee
where the $C^*$-norm limit is implied.
\end{remark}

	\begin{empt}\label{inductive_lim_empt}
		Let $\Om$ and $\Om_{\la}$ be state spaces of $A$ and $A_\la$ respectively. When $A$ is a $C^*$-
		inductive limit of $\left\{A_\la\right\}$, every state $\tau$ of $A$ defines a state $\tau_{\la}$ on $A_\la$.
		Then, for every $\mu, \nu \in \La$ such as $\mu < \nu$, we put $f^*_{\mu\nu}$ the conjugate mapping of 	the principal isomorphism $f_{\mu\nu}$ of $A_\mu$ into $A_\nu$ which maps $\Om_{\nu}$ onto  $\Om_{\nu}$
		has the following properties
		\be\label{ind_7_eqn}
		\tau_\mu = f^*_{\mu\nu}\left(\tau_\nu \right) 
		\ee
		\be\label{ind_8_eqn}
		f^*_{\mu\nu}= f^*_{\mu\la}\circ f^*_{\la\nu}\quad \text{ if }\quad \mu < \la < \nu.
		\ee
		Conversely, a system of states $\left\{\tau_\la \in \Om_{\la} \right\}_{\la \in \La}$ which satisfies the condition \eqref{ind_7_eqn}  defines a state on $A$ since every positive bounded linear functional on the algebraic inductive limit $A^0$ of  $\left\{A_\la\right\}$, is uniquely extended over $A$.
	\end{empt}
	\begin{prop}\label{inductive_lim_prop}\cite{takeda:inductive}
		Let a $C^*$-algebra $A$ be a $C^*$-inductive limit of $\left\{A_\la\right\}_{\la \in \La}$, and $\La'$ be a cofinal subset in $\La$, then $A$ is the $C^*$-inductive limit of $\left\{A_{\la'})\right\}_{\la' \in \La'}$.
	\end{prop}
	\begin{theorem}\label{inductive_lim_state_thm}\cite{takeda:inductive}
		If a $C^*$-algebra $A$ is a $C^*$-inductive limit of $A_\la$ ($\la \in \La$), the
		state space $\Om$ of A is homeomorphic to the projective limit of the state spaces $\Om_\la$ of $A_\la$.
	\end{theorem}
	\begin{corollary}\label{inductive_lim_state_cor}\cite{takeda:inductive}
		If a commutative $C^*$-algebra $A$ is a $C^*$-inductive limit of the
		commutative  $C^*$-algebras $A_\la$ ($\la \in \La$), the spectrum $\mathcal X$ of $A$ is the projective limit of spectra $\mathcal X_\la$ of $A_\la$ ($\la \in \La$).
	\end{corollary}

\section{Operator spaces and algebras}\label{oa_sec}
	\begin{definition}\label{op_space_defn}\cite{blecher_merdy}
		A concrete \textit{operator space} is a (usually closed) linear subspace $X$ of $B(\H_1,\H_2)$,
		for Hilbert spaces $\H_1,\H_2$ (indeed the case $\H_1=\H_2$ usually suffices, via the canonical
		inclusion $B(\H_1,\H_2)\subset B(\H_1\oplus\H_2)$.
	\end{definition}
	\begin{definition}\label{c_op_alg_defn}\cite{blecher_merdy}
		A \textit{concrete operator algebra} is a closed subalgebra of $B(\H)$, for some Hilbert space $\H$.
	\end{definition}
\begin{remark}\label{c_op_alg_n_rem}
Sometimes we consider operator algebras which are dense subalgebras of operator algebras.
\end{remark}
	\begin{remark}\label{c_op_alg_rem}
		There are  abstract definitions of  \textit{operator algebras} and and an \textit{operator spaces}. It is proven that the classes of operator algebras and operator spaces are equivalent to the class of concrete operator algebras and spaces respectively (cf. \cite{blecher_merdy}.)
	\end{remark}
	\begin{definition}\label{op_u_space_defn}\cite{blecher_merdy}
		(Unital operator spaces). We say that an operator space $X$ is \textit{unital} if it has a distinguished element usually written as $e$ or 1, called the \textit{identity} of $X$, 
		such that there exists a complete isometry $u: X \hookto A$ into a unital $C^*$-algebra with $u(e) = 1_A$. A \textit{unital-subspace} of such $X$ is a subspace containing $e$. 
	\end{definition}
	
	\begin{empt}\label{complete_maps_empt}\cite{blecher_merdy}
		(Completely bounded maps). Suppose that $X$ and $Y$ are operator spaces 
		and that $u: X \to Y$ is a linear map. For a positive integer $n$, we write $u_n$ for 
		the associated map $[x_{jk} ] \mapsto \left[u\left(x_{jk}\right)\right]$ from $\mathbb{M}_n(X)$ to $\mathbb{M}_n(Y )$. This is often called the ($n^{\text{th}}$) amplification of $u$, and may also be thought of as the map $\Id_{\mathbb{M}_n\left( \C\right)} \otimes u$ on 
		$\mathbb{M}_n\left( \C\right)\otimes X$. Similarly one may define $u_{m,n} : \mathbb{M}_{m,n}(X) \to \mathbb{M}_{m,n}(Y)$. If each matrix 
		space $\mathbb{M}_n(X)$ and $\mathbb{M}_n(Y)$ has a given norm $\left\|\cdot \right\|_n$ , and if un is an isometry for 
		all $n\in \N$, then we say that $u$ is 
		\textit{completely isometric}, or is a \textit{complete isometry}. 
		Similarly, $u$ is \textit{completely contractive} (resp. is a \textit{complete quotient map}) if each 
		$u_n$ is a contraction (resp. takes the open ball of $\mathbb{M}_n(X)$ onto the open ball of $\mathbb{M}_n(Y )$). A map $u$ is completely bounded if 
		$$
		\left\| u\right\|_{\mathrm{cb}}  \stackrel{\mathrm{def}}{=}\sup\left\{\left.\left\| \left[u\left(x_{jk}\right)\right]\right\|~\right|~\left\| \left[ x_{jk}\right] \right\|< 1 \quad\forall n \in \N\right\} 
		$$
		Compositions of completely bounded maps are completely bounded, and one 
		has the expected relation $\left\| u\circ v\right\|_{\mathrm{cb}}\le \left\| u\right\|_{\mathrm{cb}}\left\| v\right\|_{\mathrm{cb}}$. If $u: X \to Y$ is a completely 
		bounded linear bijection, and if its inverse is completely bounded too, then we 
		say that $u$ is a complete isomorphism. In this case, we say that $X$ and $Y$ are 
		completely isomorphic and we write $X \approx Y$. 
		The matrix norms satisfy to the following \textit{Ruan's axioms}
	\end{empt}
	\begin{definition}
		The following conditions:
		\begin{itemize}
			\item[(R1)]  $\left\|\a x\bt \right\|_n \le \left\|\a \right\|\left\|x \right\|_n\left\|\bt \right\|$, for all  $n\in \N$ and all  $\a, \bt \in \mathbb{M}_n\left(\C\right)$ and $x \in \mathbb{M}_n\left(X\right)$
			(where multiplication of an element of $\mathbb{M}_n\left(X\right)$ by an element of $\mathbb{M}_n\left(\C\right)$ is
			defined in the obvious way).
			\item[(R2)]  For all $x\in\mathbb{M}_n\left(X\right)$ and $y\in\mathbb{M}_n\left(Y\right)$, we have
			$$
			\left\|\begin{pmatrix} x & 0\\
				0& y
			\end{pmatrix}\right\|_{n + m} = \max\left(\left\|x \right\|_n, \left\|y \right\|_m\right).
			$$
			
		\end{itemize}
		are said to be \textit{Ruan's axioms}.
	\end{definition}
	\begin{remark}\label{c_star_iso_rem}\cite{blecher_merdy}
		Any $*$-homomorphism $u: A \to B$ of $C^*$-algebras is injective if and only if $u$ is completely isometric (cf. \cite{blecher_merdy}).
	\end{remark}
	\begin{empt}\label{blecher_merdy_empt}\cite{blecher_merdy}
		(Extensions of operator spaces). An \textit{extension} of an operator space $X$ is 
		an operator space $Y$ , together with a linear completely isometric map $j : X\hookto Y$. 
		Often we suppress mention of $j$, and identify $X$ with a subspace of $Y$. We say that $Y$ is a \textit{rigid} extension of $X$ if $\Id_Y$ is the only linear completely contractive 
		map $Y\to Y$ which restricts to the identity map on $j\left( X\right)$. We say $Y$ is an \textit{essential} extension of $X$ if whenever $u: Y\to Z$ is a completely contractive map 
		into another operator space $Z$ such that $u\circ j$ is a complete isometry, then $u$ is a complete isometry. We say that $(Y, j)$ is an \textit{injective envelope} of $X$ if $Y$ is 
		injective, and if there is no injective subspace of $Y$ containing $j(X)$. 
	\end{empt}
	\begin{lem}\label{os_eqv_lem}
		Let $(Y, j)$ be an extension of an operator space $X$ such that $Y$ is 
		injective. The following are equivalent: 
		\begin{enumerate}
			\item[(i)] $Y$ is an injective envelope of $X$, 
			\item[(ii)] $Y$ is a rigid extension of $X$, 
			\item[(iii)] $Y$ is an essential extension of $X$.
		\end{enumerate}
	\end{lem}
	\begin{definition}\label{c_ext_defn}\cite{blecher_merdy}
		Thus we define a $C^*$-\textit{extension} of a unital operator space 
		$X$  to be a pair $(A, j)$ consisting of a unital $C^*$-algebra $A$, and a unital 
		complete isometry $j : X\to  A$, such that $j(X)$ generates $A$ as a $C^*$-algebra. 
	\end{definition}
	
	\begin{corollary}\cite{blecher_merdy} Following conditions hold:
		$ $
		\begin{itemize}
			\item[(i)] If $X$ is a unital operator space (resp. unital operator algebra, approximately 
			unital operator algebra), then there is an injective envelope $(I(X), j)$ for $X$, 
			such that $I(X)$ is a unital $C^*$-algebra and $j$ is a unital map (resp. $j$ is a 
			unital homomorphism, $j$ is a homomorphism). 
			\item[(ii)] If $A$ is an approximately unital operator algebra, and if $(Y, j)$ is an injective 
			envelope for $A^+$, then $(Y, j|_A)$ is an injective envelope for $A$. 
			\item[(iii)] If $A$ is an approximately unital operator algebra which is injective, then $A$ 
			is a unital $C^*$-algebra. 
		\end{itemize}
		
	\end{corollary}

	\begin{theorem}\label{c_env_sp_thm}\cite{blecher_merdy}
		(Arveson‚ Hamana) If $X$ is a unital operator space, then there 
		exists a $C^*$-extension $(B, j)$ of $X$ with the following universal property: Given 
		any  $C^*$-extension $(A, k)$ of $X$, there exists a (necessarily unique and surjective) 
		$*$-homomorphism $\pi: A \to B$, such that $\pi \circ k = j$. 
	\end{theorem}
	\begin{definition}\label{c_env_sp_defn}\cite{blecher_merdy}
		If $X$ is a unital operator space, then the given by the Theorem \ref{c_env_sp_thm} pair  $(B, j)$ is said to be the $C^*$-\textit{envelope} of $X$. Denote by
		\be\label{c_env_sp_eqn}
		C^*_e\left( X\right)\stackrel{\text{def}}{=} B.
		\ee 
	\end{definition}

	\begin{definition}\label{c_env_alg_defn}\cite{blecher_merdy}
		The $C^*$-\textit{envelope} of a nonunital operator algebra $A$ is  a pair $(B, j)$, where $B$ is the $C^*$-subalgebra generated by the copy $j(A)$ of $A$ inside a $C^*$-envelope $(C^*_e (A^+), j)$ of the unitization $A^+$ of $A$. Denote by 
		\be\label{c_env_alg_eqn}
		C^*_e\left( A\right)\stackrel{\text{def}}{=} B.
		\ee
	\end{definition}

	\subsection{Real operator spaces}
	
	\begin{definition}\label{real_os_defn}\cite{ruan:real_os, ruan:real_comp}
		A \textit{real operator space} on a real Hilbert space $\H$ is a norm closed subspace $V$ of $B\left(\H\right)$ 
		together with the canonical matrix norm inherited from  $B\left(\H\right)$ . Then every real $C^*$-algebra,
		which is defined to be a norm closed $*$-subalgebra of some $B\left(\H\right)$, is a real operator space with
		a canonical matrix norm (actually, a real $C^*$-algebra matrix norm) inherited from $B\left(\H\right)$ . 
	\end{definition}
	
	\begin{remark}
		It is easy to verify
		that there the family of matrix norms $\left\{\left\|\cdot \right\| \right\}_{n \in \N} $ which satisfies
		\begin{enumerate}
			\item [(M1)]  $\left\|x \oplus y \right\|_{n + m}= \max \left(\left\|x \right\|_{n}, \left\| y \right\|_m\right)$.
			\item[(M2)] $\left\|\a x\bt\right\|_{n}\le\left\|\a \right\|\left\|x \right\|_{n}\left\|\bt \right\|$
		\end{enumerate}
		for all $x\in \mathbb{M}_n\left( V\right)$,   $y\in \mathbb{M}_m\left( V\right)$, $\a, \bt\in \mathbb{M}_n\left( \R\right)$,
		
		As in the complex case, we can define an \textit{abstract real operator space} to be
		a real space $V$ together with a Banach space norm $\left\|\cdot \right\|_n$ on each matrix space
		$\left\{\left\|\cdot \right\| \right\}_{n \in \N} $ such that the conditions (M1) and (M2) are satisfied.
	\end{remark}
	
	\begin{example}\label{op_real_exm}\cite{ruan:real_os}
		If $\Om$ is a compact Hausdorff space, then the space of all real-valued continuous functions
		$C\left(\Om, \R \right)$  is a commutative real $C^*$-algebra, and we may obtain a canonical (real $C^*$-algebra)
	\end{example}
		
	\begin{empt}\label{complexification_empt}\cite{ruan:real_comp}
		Let $V$ be a real vector space. The complexification $\C V$ of $V$ is defined as
		the direct sum
		$$
		\C V \bydef V \oplus iV = \left\{\left. x +iy\right|x, y \in V\right\}
		$$
		
		There is a natural complex linear structure on $\C V$ given by
		$$(x_1 + iy_1) + (x_2 + iy_2) = (x_1 + x_2) + i(y_1 + y_2)$$
		and
		$\left(a + i\bt \right)\left(x + iy\right)\bydef \left(\a x - \bt y\right)+ i\left(\bt x + \a y\right)$
		$(\a + i\bt)(x + iy) = (\a x -\bt y) + i(\bt x + 
\a y)$
		We can also define a natural conjugation on $\C V$ by letting
		$\overline {x + iy} \bydef x - iy$,
		Then up to the identification $x = x+i0$, we may identify $V$ with the real part
		of $\C V$ since an element $z = x + iy\in \C V$ is contained in $V$ if and only if $\overline z = z$.
	\end{empt}

	\begin{thm}\label{ros_contr_iso_thm}\cite{ruan:real_comp}
		Let $V$ and W be real operator spaces on real Hilbert spaces
		$H$ and $K$ and let $T : V \to W$ be a complete contraction (respectively, a
		complete isometry from $V$ into $W$). If $\C V$ and $\C W$ are equipped with the
		canonical complex operator space matrix norms from $B(H_c)$ and $B(K_c)$, then
		$\C T : \C V \to \C W$ is a complete contraction with $\left\| \C T\right\|_{\mathrm{cb}}  = \left\| T\right\|_{\mathrm{cb}}$ (respectively, a
		complete isometry from $\C V$ into $\C W$).
	\end{thm}
	\begin{empt}
		
		The \textit{canonical complex operator space structure}
		on $\C V$ is determined by the identification
		\be\label{ros_can_eqn}
		\C V = \left\{\left. \begin{pmatrix}
			x & -y \\
			y& x
		\end{pmatrix}\right| x, y \in V\right\}
		\ee
		where the latter space is a real subspace of $\mathbb{M}_2(V )$.
	\end{empt}
	\begin{definition}\label{ros_reason_defn}\cite{ruan:real_comp}
		The matrix norm on its complexification $\C V = V\oplus
		iV$ satisfies the \textit{reasonable condition} if one has
		\be\label{ros_reason_eqn}
		\left\|x + iy\right\|_n= \left\| x - iy\right\|_n
		\ee
		for all $x + iy \in  \mathbb{M}_n\left( \C V\right) = \mathbb{M}_n\left( V\right) + i\mathbb{M}_n\left( V\right)$ and $n\in \N$.
	\end{definition}
	\begin{theorem}\label{rs_e_thm}\cite{ruan:real_comp}
		Let $V$ be a real operator space. If a complex operator space
		matrix norm $\left\{\left\|\cdot \right\| \right\}_{n \in \N} $ on $\C V$ is a reasonable complex extension of the original
		matrix norm on $V$ , then $\left\{\left\|\cdot \right\| \right\}_{n \in \N} $ must be equal to the canonical matrix norm
		$\left\{\left\|\cdot \right\| \right\}_{n \in \N} $ on $\C V$ given by \eqref{ros_can_eqn}.
	\end{theorem}
	\section{$K$-theory and $K$-homology}
	\paragraph*{} The notions of $K$-theory, $K$-homology and $KK$-theory are explained in \cite{blackadar:ko,wegge_olsen}.
	\begin{definition}\label{stable_mult_defn}\cite{blackadar:ko}
		The \textit{stable multiplier algebra} $M^s\left( A\right)$  is the multiplier algebra
		of $A\otimes\K$. The \textit{stable outer multiplier algebra} $Q^s\left( A\right)$ is the quotient
		$M\left(A\otimes\K\right) /A\otimes\K$.
	\end{definition}
	\begin{corollary}\label{stable_k_cor}\cite{blackadar:ko}
		If $A$ is any $C^*$-algebra, then $K_n\left(A \right) \cong K_{n-1}\left(Q^s\left(A \right)  \right)$  for
		$n = 0; 1$.
	\end{corollary}
	
	\begin{corollary}\label{stab_k_0_cor}\cite{wegge_olsen}
The morphism $A \to A\otimes \K$ sending $a \mapsto a\otimes \mathfrak{e}_{1,1}$, where $\mathfrak{e}_{1,1}$ is a rank 1 
projection in $\K$, induces an isomorphism $K_0\left(A\right)\cong K_0\left(A\otimes\K\right)$. 
In particular, if $A$ and $B$ are stably isomorphic
then $K_0\left(A\right)\cong K_0\left(B\right)$. 
	\end{corollary}
\begin{cor}\label{stab_k_1_cor}\cite{wegge_olsen}
$K_1\left(A\right)\cong K_1\left(A\otimes\K\right)$) for all $C^*$-algebras $A$. 
		In particular, $K_1(A)\cong K_1(B)$ when $A$ and $B$ are stably isomorphic. 
\end{cor}

\begin{corollary}\label{kk_stable_cor}\cite{blackadar:ko}
For any $C^*$-algebras $A$ and $B$, and for any $m$ and $n$, there are natural
	isomorphisms
	\bean
	KK\left(A, B \right) \cong 	KK\left(A\otimes \mathbb{M}_n\left(\C \right) , B \otimes \mathbb{M}_m\left(\C \right)\right)	\cong KK\left(A, B\otimes\K \right)\cong \\\cong	KK\left(A, B\otimes\K \right)\cong 	KK\left(A\otimes\K, B \otimes\K\right).
	\eean
	
\end{corollary}

	\subsection{Fredholm index}
	
	\begin{definition}\label{fred_a_defn}\cite{matro:hcm}
		Let $A$ be a $C^*$-algebra. A bounded $A$ linear operator $F: \ell^2\left(A \right)\to \ell^2\left(A \right)$ is called $A$-\textit{Fredholm} if
		\begin{enumerate}
			\item [(a)] It is adjointable.
			\item[(b)] There exists a decomposition of the domain $\ell^2\left(A \right) = \M_1 \oplus \mathcal N_1$ and the range $\ell^2\left(A \right) = \M_2 \oplus \mathcal N_2$ (where $\M_1,~\M_2,~\mathcal N_1, ~\mathcal N_2$ are closed $A$-modules and $\mathcal N_1, ~\mathcal N_2$ have a finite number of generators), such that $F$ has the matrix form $F = \begin{pmatrix}
				F_1& 0\\
				0 & F_2
			\end{pmatrix}$ with respect to these decompositions and $F_1: \M_1 \cong \M_2$ is an isomorphism. 
		\end{enumerate} 
	\end{definition}
	\begin{rem}\cite{matro:hcm}
		If the conditions of the Definition \ref{fred_a_defn}  hold, then both $\mathcal N_1$ and $\mathcal N_2$ are projective $A$-modules and $\left[ \mathcal N_1\right] - \left[\mathcal  N_2\right] \in K\left(A \right) $ is well defined.
	\end{rem}
	\begin{definition}\label{fred_index_defn}\cite{matro:hcm}
		If the conditions of the Definition \ref{fred_a_defn} we define the \textit{index} of $F$ by
		$$
		\mathrm{Index}~F \bydef \left[ \mathcal N_1\right]-  \left[ \mathcal  N_2\right] \in K\left(A \right).
		$$
	\end{definition}
	
	\subsection{Long Exact Sequence of K-Theory}
	\paragraph*{}
If $A$ is a $C^*$-algebra and $J \subset A$ is a closed two-sided ideal then there there a $*$-homomorphisms $\iota: J \hookto A$, $~\pi: A \to A / J$ and a long exact sequence
\be\label{long_exact_seq_eqn}
K_1\left(J \right) \xrightarrow{\iota_*}K_1\left(A \right) \xrightarrow{\pi} K_1\left(A/J \right) \xrightarrow{\partial}K_0\left(J \right) \xrightarrow{\iota_*}K_0\left(A \right) \xrightarrow{\pi} K_0\left(A/J \right)
\ee
(cf. \cite{blackadar:ko})
\begin{definition}\cite{blackadar:ko}
	 $~$Let $u \in  GL_n(A/J)$, and let $w \in GL_{2n}(A)$ be a lift of $\mathrm{diag}(u,u^{-1})$. Define $\partial([u]) = \left[wp_nw^{-1}\right]- \left[p_n\right]\in K_0(J)$. The map $\partial$ is called the \textit{index map}.
\end{definition}
\begin{empt}\label{v_lift_empt}\cite{blackadar:ko}
Suppose A is a unital $C^*$-algebra and $u$ is a unitary in $\mathbb{M}_n(A/J)$. If $u$ lifts to a partial isometry $v \in \mathbb{M}_n(A)$, then $\mathrm{diag}\left(u, u^{-1} \right)$ lifts to the unitary
\be\label{partial_p_eqn}
w = \begin{pmatrix}
	v & 1 - vv^*\\
	1- v^*v& v^*
\end{pmatrix},
\ee
so 
\be\label{partial_u_eqn}
\begin{split}
\partial([u]) = \left[wp_nw^{-1}\right]- \left[p_n\right]=\left[\mathrm{diag}\left( v^*v, 1-vv^*\right) \right]- \left[p_n\right] = \\=\left[ 1-vv^*\right]- \left[1-v^*v\right].
\end{split}
\ee
\end{empt}

\begin{empt}\label{part_index_empt}
It is proven in \cite{blackadar:ko} that both groups $K_0\left( M\left( A\otimes\K\right)\right)$ and $K_1\left( M\left( A\otimes\K\right)\right)$ are trivial. From the equation \eqref{long_exact_seq_eqn} it follows that 
$$
\partial : K_1\left(  M\left( A\otimes\K\right)/\left( A\otimes\K\right)\right)\cong  K_0\left( A\otimes\K\right) 
$$
is an isomorphism. If $u \in M\left( A\otimes\K\right)/\left( A\otimes\K\right)$ is an unitary which lifts to a partial isometry $v \in M\left( A\otimes\K\right)$ then similarly  \ref{v_lift_empt} $~v$ is a partial isometry. Moreover one has
\be\label{d_u_iv_eqn}
\partial \left[u\right]= \mathrm{Index}~v
\ee
(cf. Definition \ref{fred_index_defn}).
\end{empt}


	\subsection{Universal coefficient theorem}
	\paragraph*{}	The universal coefficient theorem \cite{blackadar:ko} establishes (in particular) a relationship between $K$ - theory and $K$- homology. 
		\begin{theorem}
			\label{uct_thm}\cite{blackadar:ko}
		(Universal 	Coefficient 	Theorem 	(UCT)).	
		Let 
		$A$ 
		and 
		$B$ 
		be 
		separable 
		$C^*$-algebras, with $A\in N$. Then 
		there 
		is 
		a 
		short 
		exact 
		sequence 
		\begin{equation}\label{uct_eqn}
			\mathrm{\Ext}^1_\Z(K_*(A), K_*(B)) \xrightarrow{\delta} KK^*(A, B) \xrightarrow{\ga} \mathrm{Hom}(K_*(A), K_*(B))
		\end{equation}
		The 
		map 
		$\ga$ 
		has 
		degree 
		0 and 
		$\delta$
		has 
		degree 
		1. The 
		sequence 
		is 
		natural 
		in 
		each 
		variable, and 
		splits 
		unnaturally. So 
		if 
		$K\left( A\right)$  is 
		free 
		or 
		$K\left( A\right)$ is 
		divisible, then 
		$\ga$
		is 
		an 
		isomorphism. 
	\end{theorem}
\begin{remark}\label{n_alg_rem}
	The class $N$ of $C^*$-algebras is defined in \cite{blackadar:ko} (Definition 22.3.4).
\end{remark}
	If $\tau \in KK^1(A, B)$ is represented by extension
	\begin{equation}\nonumber
		0 \to  B \to D \to A \to  0
	\end{equation}
	then $\gamma$ is given as connecting maps $\partial$ in the associated six-term exact sequence of $K$-theory
	\newline
	\newline
	\begin{tikzpicture}
		\matrix (m) [matrix of math nodes,row sep=3em,column sep=4em,minimum width=2em]
		{
			K_0(B) &    K_0(D)&  {K_0(A)}\\
			K_1(A)	& K_1(D)   & K_1(B) \\};
		\path[-stealth]
		(m-1-1) edge node [above] {} (m-1-2)
		(m-1-2) edge node [above] {} (m-1-3)
		(m-2-3) edge node [above] {} (m-2-2)
		(m-2-1) edge node [right] {$\partial$} (m-1-1)
		(m-1-3) edge node [right] {$\partial$} (m-2-3)
		(m-2-2) edge node [above] {}  (m-2-1);
		
	\end{tikzpicture}
	\newline
	
	If $\gamma(\tau)=0$ for an extension $\tau$ then the six-term $K$-theory exact sequence degenerates into two short exact sequences 
	\begin{equation}\nonumber
		0 \to K_j(A) \to K_j(D) \to K_j(B) \to 0 \quad (j=0,1)
	\end{equation}
	and thus determines an element $\kappa(\tau)\in \mathrm{\Ext}^1(K_*(A), K_*(B))$. Homomorphism $\delta$ is inverse to $\kappa$.
	\begin{definition}\label{k_hom_defn}\cite{blackadar:ko}
	If $A$ is a $C^*$-algebra then $K$-\textit{homology} functor $K_*$ is given by $K_*\left( A\right)\bydef KK^*\left( A, \C\right)$. 
	\end{definition}
\begin{definition}\label{busby_defn}\cite{blackadar:ko}
If 
	\be\label{ext_eqn}
0 \to 	B \to D \to A \to 0
\ee
is an extension of $C^*$-algebra, then natural $*$-homomorphism 
	\be\label{busby_eqn}
A \to Q\left(B \right) \bydef M\left(B \right)/B 
\ee
is said to be  the \textit{Busby invariant}.
\end{definition}
\begin{remark}
	The {Busby invariant} uniquely defines the extension up to isomorphism.
\end{remark}
\begin{empt}\label{busby_k_1_empt}\cite{blackadar:ko}  
	It is proven that any $x \in K_1\left(A \right)$ corresponds to a class of isomorphism of extension
	\be\label{ext_k_eqn}
0 \to 	\K \to D \to  A \to 0.
	\ee
	So $x$ can be represented by the Busby invariant
	\be\label{busby_k_eqn}
A \to Q\left(\K \right) \bydef M\left(\K\right)/\K. 
\ee
\end{empt}
\begin{empt}\label{cirle_k1_empt}\cite{blackadar:ko}
Let $C\left( u\right) \cong C\left( S^1\right)$ be a $C^*$-algebra generated by a single unitary element $u$.
If $w \in B\left( \H\cong \ell\left(\N\right)\right)$ is  an unilateral shift then $w$ represents an unitary element $v \in Q\left(\K \right)$. 
It is known \cite{blackadar:ko} that $K^1\left( C\left( u\right)\right)\cong \Z$.  and a generator of $K^1\left( C\left( u\right)\right)$ is represented by the homomorphism
\be\label{s1_k_1_eqn}
\begin{split}
K^1_{S^1}:C\left( u\right)\to  Q\left(\K \right),\\
u \mapsto v.
\end{split}
\ee
\end{empt}
	
 For any $C^*$-algebra $A$ there is a natural homomorphism
	\begin{equation}\label{free_spec_eqn}
		\gamma: K^1(A) \to \mathrm{Hom}(K_1(A), K_0(\mathbb{C})) \approx \mathrm{Hom}(K_1(A), \mathbb{Z})
	\end{equation}
	which is the adjoint of following pairing
	\begin{equation}\nonumber
		KK(\mathbb{C}, A) \otimes KK(A, \mathbb{C}) \to KK(\mathbb{C}, \mathbb{C})\cong\Z.
	\end{equation}
On the other hand  the equation \eqref{free_spec_eqn} comes from \eqref{uct_eqn}.
	If $\tau \in KK^1(A, \mathbb{C})$ is represented by extension
	\begin{equation}\nonumber
		0 \to  \mathbb{C} \to D \to A \to  0
	\end{equation}
	then $\gamma$ is given as connecting maps $\partial$ in the associated six-term exact sequence of $K$ theory
		\newline
	\newline
	\begin{tikzpicture}
		\matrix (m) [matrix of math nodes,row sep=3em,column sep=4em,minimum width=2em]
		{
		 K_0(\mathbb{C}) &    K_0(D)&  {K_0(A)}\\
	K_1(A)	& K_1(D)   & K_1(\mathbb{C}) \\};
		\path[-stealth]
		(m-1-1) edge node [above] {} (m-1-2)
		(m-1-2) edge node [above] {} (m-1-3)
		(m-2-3) edge node [above] {} (m-2-2)
		(m-2-1) edge node [right] {$\partial$} (m-1-1)
		(m-1-3) edge node [right] {$\partial$} (m-2-3)
		(m-2-2) edge node [above] {}  (m-2-1);
		
	\end{tikzpicture}
	\newline

	If $\gamma(\tau)=0$ for an extension $\tau$ then the six-term $K$-theory exact sequence degenerates into two short exact sequences 
	\begin{equation}\nonumber
		0 \to K_j(A) \to K_j(D) \to K_j(\mathbb{C}) \to 0 \quad (j=0,1)
	\end{equation}
	and thus determines an element $\kappa(\tau)\in \mathrm{\Ext}^1(K_*(A), K_*(\mathbb{C})$. 
	In result we have a sequence of abelian group homomorphisms
	\begin{equation}\label{uct_c_eqn}
		\mathrm{\Ext}^1(K_0(A), K_0(\mathbb{C})) \to KK^1(A, \mathbb{C}) \to \mathrm{Hom}(K_1(A), K_0(\mathbb{C}))
	\end{equation}
	such that composition of the homomorphisms is trivial. Above sequence can be rewritten by following way
	\begin{equation}\label{uct_c}
		\mathrm{\Ext}^1(K_0(A), \mathbb{Z}) \to K^1(A) \to \mathrm{Hom}(K_1(A), \mathbb{Z})).
	\end{equation}
	If $G$ is an abelian group that 
	\begin{equation}\nonumber
		\mathrm{\Ext}^1(G, \mathbb{Z}) = \mathrm{\Ext}^1(G_{tors}, \mathbb{Z}),
	\end{equation}
	\begin{equation}\nonumber
		\mathrm{Hom}(G, \mathbb{Z}) =  \mathrm{Hom}(G / G_{tors}, \mathbb{Z})).
	\end{equation}
	From (\ref{uct_c}) it follows that $K^1(A)$ depends on $K_0(A)_{tors}$ and $K_1(A)/K_1(A)_{tors}$. We say that dependence \eqref{uct_c} on $K_0(A)_{tors}$ is a {\it torsion special case} and dependence (\ref{free_spec_eqn}) of $K^1(A)$ on $K_1(A)/K_1(A)_{tors}$ is a {\it free special case}.
	


\subsection{Miscellany}	

\begin{proposition}\label{nt_khom_prop}\cite{had:ntk,sudo:ntk}
For $0\le \th \le 1$, the $K$-homology groups of the $C\left(\T^2_\th \right)$ are $K^j\left(C\left(\T^2_\th \right) \right)\bydef  KK^j\left( C\left(\T^2_\th \right)\right) \cong \Z^2$.

\end{proposition}
\section{Inverse limits of $C^*$-algebras}

\subsection{Basic constructions}
\paragraph*{}
Inverse limits of $C^*$-algebras are specializations of the explained  in \ref{proj_sys_sec} theory.
	\begin{prop}\label{pro_c_prop}\cite{phillips:inv_lim_app}
		Let $A$ be a topological $*$ -algebra over 
		$\C$. Then $A$ is isomorphic, as a topological  $*$ -algebra, to an inverse limit of $C^*$- algebras in the above sense if and only if $A$ is Hausdorff, complete, and its topology 
		is determined by the set of all continuous $C^*$-seminorms on $A$.
	\end{prop}	
	
	\begin{definition}\label{pro_c_defn}\cite{phillips:inv_lim_app}
		A topological $*$-algebra satisfying the 
		conditions of the previous proposition is called a \textit{pro}-$C^*$ -\textit{algebra}. If it is a countable  
		inverse limit of $C^*$-algebras (equivalently, if its topology is determined by countably 
		many continuous $C^*$-seminorms), then it is called a $\sigma$-$C^*$-algebra. 
		
	\end{definition}
	
	\begin{notn}\label{ap_notn}
		If $A$ is a pro-$C^*$-algebra, then $S(A)$ is the set of all	continuous $C^*$-seminorms on $A$, ordered as in the proof of Proposition \ref{pro_c_prop}. For 
		$p\in S\left(A\right)$, we let $\ker(p) \bydef \left\{\left.a\in A\right|p(a)=0\right\}$ and we let $A_p$ be the completion of 
		$A/\ker(p)$ in the norm determined by $p$. It is proven in \cite{phillips:inv_lim} that $A/\ker(p)$ is in fact already complete, so $A/\ker(p)\cong A_p$ is a $C^*$-algebra.  
	\end{notn} 
	Note that $A_p$ is a $C^*$-algebra, and that the proof of Proposition \ref{pro_c_prop} 
	gives a canonical representation of $A$ as an inverse limit, namely $A = \varprojlim_{p\in S(A)} A_p$. 
	\begin{definition}\label{pro_bound_defn}\cite{phillips:inv_lim}
		Let $A$ be a pro-$C^*$-algebra. Then the set of \textit{bounded elements} is the set
		$$
		b\left(A\right)\bydef \left\{ a \in A\left| \left\| a\right\|_\infty = \sup\left\{p\left( a\right)\left| p\in S\left(a\right) \right.\right\}< \infty\right. \right\}.
		$$
		
	\end{definition}

	\begin{proposition}\label{pro_bound_prop}\cite{phillips:inv_lim_app}
		Let $A$ be a pro-$C^*$-algebra. Then:
		
		\begin{enumerate}
			\item $b\left(A\right)$ is a $C^*$-algebra in the norm $\left\| \cdot\right\|_\infty$.
			\item If $a \in A$ is normal and $f\in C\left(Sp\left(a\right)\right)$ is bounded then $f\left(a\right)\in b\left(A\right)$.
			\item If $a \in A$ is normal then $a\in b\left(A\right)$ if and only if $sp\left(a\right)$ is bounded.
			\item $b\left(A\right)$ is dense in $A$.
			\item For $a\in b\left(A\right)$, we have $sp_{b\left(A\right)}\left(a \right)=\overline{sp_A\left(a \right)}$.
			\item If $q \in S\left(A\right)$ then the map from $b\left(A\right)$ to $A_q$ is surjective.
		\end{enumerate}
		
	\end{proposition}

	\begin{prop}\label{pro_dir_lim_prop}\cite{phillips:inv_lim}
		Direct limits exist in the category of pro-$C^*$-algebras.
	\end{prop}
	\begin{proof}
		If $\left\{A_{\a}\right\}_{\a \in \mathscr A}$ is a direct system of pro-$C^*$-algebras with homomorphisms $\varphi_{\a\bt}: A_\a \to A_\bt$ for $\a\le\bt$ then the direct limit is constructed as follows. Let
		\be\label{seminorms_inv_lim}
		D\bydef\left\{\left. p \in \prod_{\a \in \mathscr A} S\left(A_\a \right)  \right|\forall \a, \bt\in \mathscr A \quad \bt \ge \a \quad p_\bt \circ \varphi_{\a\bt} \le p_\a \right\}
		\ee
		be ordered by $p\le q$ if $p_\a\le q_\a$ for all $\a$. If $p$ denote by the set $B_p \bydef \varinjlim \left( A_\a\right)_p$ and $B = \varprojlim B_p$ then $\varinjlim A_\a$ is the closure of the images of $A_\a$ in $B$.  The complete proof is omitted, because direct limits badly behaved that they do not seem to be of much use.
	\end{proof}
	\begin{remark}\label{pro_dir_lim_rem}
		However direct limits of pro-$C^*$-algebras can be effectively used in the theory of noncommutative coverings (cf. Definition \ref{inv_pro_lim_defn}).
	\end{remark}
	\begin{empt}\label{pro_fin_comp_empt}
		There are explained in \cite{phillips:inv_lim} Hilbert modules over pro-$C^*$-algebras. It is shown that if $E_A$ is a finitely generated module over pro-$C^*$-algebra  $A$  then the  pro-$C^*$-algebra  $\mathcal{L}_A\left( X_A\right)$ adjointable is naturally isomorphic to the algebra $\K\left( X_A\right)$ of  compact operators.
	\end{empt}
	\subsection{Generators and relations}\label{generators_and_relations}
	\paragraph{}
Here I follow to \cite{phillips:inv_lim_app}.	We begin with the construction of the $C^*$-algebra on a properly chosen
	set of generators and relations. Our treatment is based that of \cite{blackadar:shape_theory}, but includes
	changes to make it more suitable for our purposes. If $G$ is any set, we denote by
	$F(G)$ the free associative complex *-algebra (without identity) on the set $G$. Thus,
	$F(G)$ consists of all polynomials in the noncommuting variables $C \coprod G^*$ (disjoint union), with complex coefficients and no constant term. By definition, any function
	$p$ from $G$ to a $C^*$-algebra $A$ extends to a unique $*$-homomorphism, which we also
	call $p$, from $F(G)$ to $A$.
	
	A set $R$ of relations on $G$ is a collection of statements about the
	elements of $G$ which make sense for elements of a $C^*$ -algebra. Possible relations
	include statements of the form "$\|x\| \in S$"', where $x \in F(G)$ and $S \in \mathbb{R}$, "$x$ is
	positive," or the statement that some equation in the variables $G \coprod G^*$ and some unknowns has a solution, or that some function from a topological space into $G$ is
	continuous. Note that Blackadar considers only relations of the form $\|x\| \le \eta$ for $\eta > 0$ and $x$ in the unitization $F(G)^+$ of $F(G)$. (We do not allow relations involving
	elements of $F(G)^+=F(G)$ because they do not make sense in a nonunital $C^*$-
	algebra. However, it is perfectly possible for $R$ to include the relations $eg = ge Ñ g$
	for some fixed $e \in G$ and all $g \in G \coprod G^*$).
	\begin{definition}\label{gr_ca_defn}\cite{phillips:inv_lim_app}
		Let $(G,R)$ be a set of \textit{generators and relations}.
		(That is, $G$ is a set and $R$ is a set of relations on $G$). Then a \textit{representation} of
		$(G,R)$ in a $C^*$-algebra $A$ is a function $\rho : G \to A$ such that the elements $\rho(g)$ for $g \in G$ satisfy the relations $R$ in $A$. A representation on a Hilbert space $\H$ is a
		representation in $B(\H)$.
	\end{definition} 
	\begin{definition}\label{adm_gr_defn}
		(Compare \cite{blackadar:shape_theory}). A set $(G,R)$ of generators and
		relations is \textit{admissible} if the following conditions hold:
		\begin{enumerate}
			\item The function from $G$ to the zero $C^*$-algebra is a representation of $(G,R)$.
			\item If $\rho$ is a representation of $(G,R)$ in a $C^*$-algebra $A$, and if $B$ is a $C^*$-subalgebra of $A$ which contains $\rho(G)$, then $\rho$ is a representation of $(G,R)$ in $B$.
			\item If $\rho$ is a representation of $(G, R)$ in a $C^*$-algebra $A$, and if $\psi: A\to B$ is a surjective homomorphism, then $\psi \circ \rho$ is a representation of $(G, R)$ in $B$.
			\item For every $g \in G$ there is a constant $M(g)$ such that $\|\rho (g) \| < M(g)$ for all representations $p$ of $(G, R)$.
			\item If $\{\rho_{\alpha}\}$ is a family of representations of $(G, R)$ on Hilbert spaces $\H_{\alpha}$ then $g \mapsto \rho(g)=\oplus_{\alpha}\rho_{\alpha}(g)$   is a representation of $(G, R)$ on $\H= \oplus \H_{\alpha}$. (That is, the elements $\rho(g)$, which are in $B(\H)$ by (4), in fact satisfy the relations $R$.)
		\end{enumerate}
	\end{definition}
	Note that, in the presence of (3), condition (1) is equivalent to "there
	exists a representation of $(G, R)$." Also note that, for relations of the sort considered
	by Blackadar, (2) and (3) are automatic and (5) follows from (4).
	\begin{definition}\label{universal_ca_defn}\cite{phillips:inv_lim_app}
	The \textit{universal $C^*$-algebra on the generators $G$ and relations $R$} is a $C^*$- algebra $C^*(G,R)$ with a representation $\rho$ of $(G,R)$ in $C^*(G,R)$, such that, given any representation $\sigma$ of $(G,R)$ in a $C^*$-algebra $B$, then there is a unique homomorphism
$\psi : C^*(G,R) \to B$ such that $\sigma=\psi \circ \rho$.	
\end{definition}
 If $(G,R)$ is admissible, then $C^*(G,R)$ exists
	and, following Blackadar, can be obtained as the Hausdorff completion of $F(G)$ in the $C^*$ -seminorm $\|x\| = \mathrm{sup}|\rho(x)\|$ where $\rho$ is a representation of $(G,R)$.
	Note that condition (4) guarantees that $\|x\| < \infty$ for $x \in G$, and that condition
	(5) guarantees that the obvious map from $G$ to $C^*(G, R)$ is in fact a representation.
	Admissibility, or something close to it, is also necessary for the existence
	of $C^*(G,R)$. Condition (1) is needed, since otherwise there may be no
	representations at all; conditions (2) and (3) are needed to ensure that the notion of a universal $C^*$-algebra is sensible, and without conditions (4) and (5) it will not be
	possible to construct a universal $C^*$ -algebra. However, if the relations $R$ also make
	sense in a pro-$C^*$ -algebra, then a pro-$C^*$-algebra with the required properties will
	exist under much weaker conditions, as in the following definition. A representation
	of $(G, R)$ in a pro-$C^*$-algebra has the obvious meaning.
	\begin{defn}\label{weak_adm_defn}\cite{phillips:inv_lim_app}
		
		A set $(G, R)$ of generators and relations is called
		\textit{weakly admissible} if the following conditions are satisfied:
		\begin{enumerate}
			
			\item The function from $G$ to the zero $C^*$-algebra is a representation of $(G,R)$.
			\item If $\rho$ is a representation of $(G,R)$ in a $C^*$-algebra $A$, and if $B$ is a $C^*$-subalgebra of $A$ which contains $\rho(G)$, then $\rho$ is a representation of $(G,R)$ in $B$.
			\item If $\rho$ is a representation of $(G, R)$ in a $C^*$-algebra $A$, and if $\varphi: A\to B$ is a surjective homomorphism, then $\varphi \circ \rho$ is a representation of $(G, R)$ in $B$.
			\item For every $g \in G$ there is a constant $M(g)$ such that $\|\rho (g) \| < M(g)$ for all representations $p$ of $(G, R)$.
			\item If $\{\rho_{\alpha}\}$ is a family of representations of $(G, R)$ on Hilbert spaces $\H_{\alpha}$ then $g \mapsto \rho(g)=\oplus_{\alpha}\rho{\alpha}(g)$   is a representation of $(G, R)$ on $\H= \oplus \H_{\alpha}$. (That is, the elements $\rho(g)$, which are in $B(\H)$ by (4), in fact satisfy the relations $R$.)
			
			\item If $A$ is a pro-$C^*$-algebra, and $\rho : G \to A$ is a function such that, for every $p \in S(A)$, the composition of $\rho$ with $A \to A_p$ is a representation of $(G, R)$ in $A_p$, then $\rho$ is a representation of $(G,R)$.
			\item If $\rho_1, ..., \rho_n$ are representations of $(G,R)$ in $C^*$-algebras $A_1,...,A_n$
			then $g \mapsto (\rho_1(g), ..., \rho_n(g))$ is a representation of $(G,R)$ in $A_1 \oplus ... \oplus A_n$.
		\end{enumerate}
	\end{defn}
	\begin{prop}\label{uni_pro_c_rep_prop} \cite{phillips:inv_lim_app}
		Let $(G,R)$ be a weakly admissible set of generators
		and relations. Then there exists a pro-$C^*$ -algebra C*(G,R), equipped with
		a representation $\rho : G \to C^*(G,R)$ of $(G,R)$, such that, for any representation $\sigma$ of $(G,R)$ in a pro-$C^*$-algebra $B$, there is a unique homomorphism $\varphi : C^*(G,R) \to B$
		satisfying $\sigma = \varphi \circ \rho$. If $(G,R)$ is admissible, then $C^*(G,R)$ is a $C^*$-algebra.
	\end{prop}
	\begin{remark}\label{uni_pro_c_rep_rem}
		The construction of $C^*(G,R)$ is described in the proof of the Proposition \ref{uni_pro_c_rep_prop}.  Let $D$ be the set of all $C^*$-seminorms on $F(G)$ of the form 
		$p(x) = \left\|\sigma\left(x \right)  \right\|$  for some representation $\sigma$ of $G$ in a $C^*$-algebra. Then take $C^*(G,R)$ 
		to be the Hausdorff completion of $F(G)$ in the family of $C^*$-seminorms $D$.
	\end{remark}

	\begin{example}\label{wa_exm}\cite{phillips:inv_lim_app}
		Any combination (including the empty set) of the 
		following kinds of relations is weakly admissible: 
		\begin{enumerate}
			\item 	Any algebraic relation among  the elements of $F(G)$, or the elements 
			of the $n\times n$ matrix algebra $\mathbb{M}_n\left(F\left( G\right)  \right)$.
			\item Any norm inequality of the form $\left\| x\right\|  < \eta \quad  (\eta \ge 0)$ or $\left\|x\right\| < \eta \quad  (\eta > 0)$, where $x \in F\left( G\right) $. or $x\in \mathbb{M}_n\left(F\left( G\right)  \right)$, and where norm relations are interpreted as 
			applying to all continuous C*-seminorms. (Thus, $\left\|x\right\| < 1$ means $p(x) < 1$ for all $p$. 
			Note that this is a weaker condition than $\left\|a\right\|_\infty< 1$).
			\item Any operator inequality $x > 0$ or $x > y$ for $x, y\in \mathbb{M}_n\left(F\left( G\right)  \right)$.
			\item  The assertion that a given function from an appropriate space to 
			$F(G)$ is continuous, Lipschitz, differentiable, continuously differentiate, or r times 
			(continuously) differentiate for $0 < r < \infty$. 
		\end{enumerate}
	\end{example}

	\subsection{Commutative pro-$C^*$-algebras}
	\begin{definition}\label{comm_pro_dist_defn}\cite{phillips:inv_lim_app}
		Let $\sX$ be a topological space. Then a family $F$ 
		of compact subsets is said to be \textit{distinguished} if it contains all one point sets, is closed 
		under finite unions and passage to compact subsets, and determines the topology in 
		the sense that a subset $C$ of $\sX$ is closed if and only if $C\cap K$ is closed for all $K \in F$. 
		If $\left(\sX_1, F_1\right)$ and  $\left(\sX_2, F_2\right)$ are spaces with distinguished families of compact subsets, 
		then a morphism from  $\left(\sX_1, F_1\right)$ to $\left(\sX_2, F_2\right)$ is a continuous function $f:\sX_1\to\sX_2$ 
		such that $f\left(K\right)\in F_2$ for every $K\in F_1$. 
	\end{definition}
	\begin{empt}\label{comm_pro_dist_empt}\cite{phillips:inv_lim_app}
		If $\sX$ is a topological space, and $F$ is a set of compact subsets of $\sX$, then we 
		write $C_K\left(\sX \right)$  for the topological *-algebra of all continuous functions from $\sX$ to $\C$, 
		with the topology of uniform convergence on the members of $F$. If $F$ is omitted, 
		it is understood to be the set of all compact subsets of $F$. Of course, in general 
		$Cont\left( \sX\right)$  can fail to be a pro-$C^*$-algebra by not being complete. We need one more definition. 
	\end{empt}

	\begin{definition}\label{comm_pro_ch_defn}\cite{phillips:inv_lim_app}
		We call a topological space $\sX$ \textit{completely Hausdorff}
		for any two distinct points $x,y\in \sX$ there is a continuous function $f: \sX \to \left[0,1\right]$ such that $f(x) = 0$ and $f(y) = 1$. 
		This condition lies between Hausdorff and completely regular. 
	\end{definition}

	\begin{theorem}\label{comm_pro_thm}\cite{phillips:inv_lim_app}
		The assignment $\left(\sX, F \right)\mapsto Cont\left( \sX, F\right)$   is a 
		contravariant category equivalence from the category of completely Hausdorff spaces 
		with distinguished families of compact subsets to the category of commutative unital 
		pro-$C^*$-algebras and unital homomorphisms. 
	\end{theorem}
	
	\begin{lemma}\label{comm_pro_lem}\cite{phillips:inv_lim}
		Let $\sX$ be a completely Hausdorff quasitopological space. Then
		$$
		Cont\left(\sX \right) \cong \varprojlim_{K \in F} C\left( K\right). 
		$$
	\end{lemma}

	\begin{remark}\label{comm_pro_rem}
 Consider  the situation of the Lemma \ref{comm_pro_lem}. If $\sX$  is a topological Hausdorff space then $	Cont\left(\sX \right)$ is a completion of $	C_0\left(\sX \right)$ with respect to a family $\left\{p_K\right\}_{K \in F}$ of seminorms given by
 \be\label{comm_pro_eqn}
p_K\left(a\right)\bydef\max_{ x \in  K}\left| a\left( x \right)  \right|. 
 \ee
\end{remark}
	
	\subsection{Noncommutative suspension and loop space}\label{noncom_loop_space_section}
\paragraph*{}	
	We would like to generalize notion of a loop space.
		\begin{definition}\label{abelianization_defn}\cite{phillips:inv_lim_app}
Let $A$ be a pro-$C^*$-algebra. Then the \textit{abelianization} of $A$ is the Hausdorff completion of $A$ in the topology determined by all 
continuous $C^*$-seminorms which vanish on the closed ideal $[A, A]$ generated by all 
commutators $xy - yx$ for $x,y \in  A$. 
The abelianization can also be described as the completion of $A/[A, A]$ 
with respect to the quotient topology, or as the inverse limit $\varprojlim A_p/[A_p, A_p]$. It is 
the largest abelian quotient of $A$. 	
\end{definition}
\begin{remark}\label{abelianization_rem}
	If $B$ is a commutative pro-$C^*$-algebra then there is the natural bijection $\Hom\left(A, B\right)\cong \Hom\left(A^{\text{ab}}, B\right)$ where $A^{\text{ab}}$ is the abelianization of $A$.
\end{remark}

	\begin{definition}\label{pointed_alg_defn}\cite{phillips:inv_lim_app}
A \textit{pointed} (pro-)$C^*$-algebra is a pair $\left(A,\a)\right)$, where $A$ is a unital (pro-)$C^*$-algebra and $\a : A \to \C$ is a unital homomorphism.
	\end{definition}
	
	\begin{defn}\label{pointed_hom_defn}\cite{phillips:inv_lim_app}
		If $(A, \alpha)$ is a  pointed pro-$C^*$-algebra, than its
		\textit{suspension} is the pro-$C^*$-algebra $\Sigma A = \{f : S^1 \to A \ \mathrm{continuous} \ : \alpha(f(\zeta)) \cdot 1 = f(1), \ \forall \zeta \in S^1\}$, together with the homomorphism $ev_1 : \Sigma A \to \C$ of evaluation at 1. Here $S^1$ is
		identified with $\{\zeta \in \mathbb{C}:|\zeta| = 1\}$. (Note that $ev_1$ makes sense as a homomorphism
		to $\mathbb{C}$ since $f(1) \in \mathbb{C} \cdot 1$ for $f \in \Sigma A$). We also denote by $\left(\Sigma A, ev_1 \right)$ the corresponding  pointed pro-$C^*$-algebra.   
	\end{defn}
	Note that $\Sigma(A^+) = (S A)^+$ , where $SA$ is the conventional suspension $C_0(\mathbb{R} \otimes A$. A left  adjoint for $\Sigma$ in the  pointed category immediately gives a left
	adjoint for $\Sigma$ in the category of pro-$C^*$-algebras and arbitrary homomoprhisms,
	simply by taking the kernel of the homomorphism to $\mathbb{C}$ which comes with the
	 pointed pro-$C^*$-algebra. 
	\begin{defn}\label{pointed_noncommutative_loop_defn}\cite{phillips:inv_lim_app}
		Let $(A, \alpha)$ be a  pointed pro-$C^*$-algebra. We
		construct a  pointed pro-$C^*$-algebra $\Omega A$ in terms of generators and relations as
		follows. Let the generating set $G$ consist of the symbols $z(a, \zeta)$ for a $a \in A$ and $\zeta \in S^1$, and let the relations $R$ be as follows:
		\begin{enumerate}
			\item The map $(a, \zeta) \mapsto z(a, \zeta)$ continuous.
			\item For each fixed $\zeta \in S^1$, the elements $z(a, \zeta)$ satisfy all the algebraic
			relations satisfied by the corresponding elements of $A$.
			\item $z(1,\zeta) = z(1,1)$ $\forall S^1$ .
			\item $z(a, 1) = z(\alpha(a)1, 1)$ $\forall \in  A$.
		\end{enumerate}
		Then set $\Omega A\bydef C^*(G, R)$. The required homomorphism from  $\om_1: \Omega A\to\C$ is given by $z(a, \zeta) \mapsto \alpha(a)$  for $a \in A$ and $\zeta \in S^1$. We say that the pair $\left(\Omega A, \om_1\right)$ is the \textit{ pointed loop algebra} of $(A, \alpha)$.
	\end{defn}
	\begin{thm}\label{nc_s_loop_iso}\cite{phillips:inv_lim_app}
		The map $\Phi : \Hom_+(\Omega A, B) \to \Hom_+(A, \Sigma B)$ defined by $\Phi(\varphi)(a)(\zeta) = \varphi(z((a, \zeta))))$, is a natural bijection, and also defines a natural bijection $[\Omega A,B]^+ \to [A, \Sigma B]^+$ . In particular, $\Omega$ is a left  adjoint for the functor
		$\Sigma$.
	\end{thm}
	\begin{remark}\label{nc_loop_rem}\cite{phillips:inv_lim_app}
We should note that if $A \bydef C(\sX)$ for some compactly generated
completely Hausdorff space $\sX$, and if $\a = ev_{x_0}$, evaluation at $x_0$, for some $x_0\in\sX$,
then the abelianization of $\Om A$ (Definition \ref{abelianization_defn}) is $C(\Om(\sX,x_0))$, where $\left( \Om\sX, \om_0\right) $ is
the usual loop space of $\sX$ relative to the basepoint (cf. \ref{topological_loop_spaces}), and is equipped
with the compactly generated topology (cf. \ref{top_comp_open_empt}). Note,
however, that $\Om C\left(\sX\right)$ is essentially never commutative, because if $\varphi :\Om A\to \Sigma B$
 is a homomorphism, then there is generally no reason for the ranges of the
homomorphisms $ev_\xi : A \to \Sigma B$ to commute with each other.
\end{remark}
	\section{Local operator spaces}
	\begin{definition}\label{loc_op_sp_defn}\cite{dosi:multi,effros:loc_conv}
		Let $X$ be a linear space and let
		$p^{n} :\mathbb{M}_n(X)\to\left[0, \infty\right]$, $n \in \N$, be gauges (respectively, seminorms) over all matrix spaces. The
		family $p \bydef \left\{p^{n}\right\}_{n\in \N}$ is said to be a \textit{ matrix gauge} (respectively, \textit{matrix seminorm})  on $X$ if
		$p$ possesses the following properties:
		\begin{itemize}
			\item [M1.] $p^{(m+n)}(v \oplus w)= \max\left( p^{(m)}\left(v \right), p^{(n)}\left(w \right)\right) $.
			\item [M2.] $ p^{(n)}\left(avb \right)\le \left\|a \right\|   p^{(m)}\left(v \right)\left\|b \right\| $.
		\end{itemize}
		for all $v = [v_{ij}] \in\mathbb{M}_m(X),~w = [v_{ij}] \in\mathbb{M}_n(E), ~a = [a_{ij}] \in\mathbb{M}_{n,m}(\C),~b = [b_{ij}] \in\mathbb{M}_{m,n}(\C),~ n,m \in\N$. A linear space $X$ with a (separated) family of matrix seminorms $\left\{p_\a\right\}_{\a\in \mathscr A}$ is called an
		\textit{abstract local operator space}. 
	\end{definition}
\begin{definition}\label{multinormed_a_defn}\cite{dosi:multi}
A seminorm $p$ on a unital associative algebra $A$ is called a \textit{multiplicative seminorm} if $p(1_A) = 1$ and $p(ab) \le p(a)p(b)$ for all $a, b \in A$. A multiplicative seminorm on an
associative *-algebra A is said to be a \textit{$C^*$-seminorm} if $p\left(a^*\right)= p(a)$ and $p\left(a^*a\right)= p(a)^2$ for
all $p\in A$. A complete polynormed  algebra with a defining family of multiplicative seminorms (respectively, $C^*$-seminorms) is called an \textit{Arens–Michael algebra} (respectively, \textit{multinormed
$C^*$-algebra}).
\end{definition}
	\begin{remark}\label{la_los_rem}
		Any local operator $C^*$-algebra is an example of a local operator system (cf. \cite{dosi:multi}).
	\end{remark}
		\begin{remark}\label{la_pro_rem}
		Any pro-$C^*$-algebra is an example of a local operator system.
	\end{remark}
		
	\begin{empt}\label{complete_loc_maps_empt}\cite{dosi:multi} 
		Similarly to \ref{complete_maps_empt}  we define a \textit{complete isometry}. 
		\textit{completely contractive} and \textit{complete quotient map}) of abstract local operator spaces.
	\end{empt}
	\begin{definition}\cite{dosi:multi} \label{conc_loc_op_space}
		A \textit{(concrete) local operator space} is a  $\C$-linear subspace $X$ of a multinormed $C^*$-algebra $A$ such that $1_A \in X$. 
	\end{definition}
	\begin{remark}\label{loc_op_abs_conc_rem}\cite{dosi:multi}
		Any local operator space naturally possesses the  structure of an abstract  local operator space. This structure comes from the structure of an abstract  local operator space of multinormed $C^*$-algebra $A$. Conversely  any abstract  local operator space is a concrete  local operator space, i.e. $\C$-subspace of a multinormed $C^*$-algebra.
	\end{remark}

	\begin{remark}\label{loc_op_sp_rem}
		Any abstract operator local space $V$ is an inverse limit of operator spaces (cf. \cite{dosi:multi} Remark 4.1).
	\end{remark}

	\section{Locally convex quasi *-algebras and their representations}
	
	\paragraph*{}
	Here I follow to \cite{quasi_star}.
	
	\begin{definition}
		A \textit{partial *-algebra} is a complex vector space $\mathfrak{A}$, endowed with	an involution $a\mapsto a^*$(that is, a bijection, such that $a^{**}=a$, for all $a\in\mathfrak{A}$) and
		a partial multiplication defined by a set $\Ga\subset \mathfrak{A} \times \mathfrak{A}$ (a binary relation), with the	following properties
		\begin{enumerate}
			\item [(a)]	$\left(a,b\right)\in\Ga$ implies $\left(b^*,a^*\right)\in\Ga$;
			\item[(b)] $\left(a,b_1\right),\left(a,b_2\right)\in\Ga$ implies $\left(a,\la b_1 +  \mu bb_2\right)\in\Ga\forall \la,\mu\in\C$;
			\item[(c)] for any $\left(a,b\right)\in\Ga$ , a product $a \cdot b \in \mathfrak{A}$ is defined, which is distributive with 		respect to the addition and satisfies the relation $\left( a^* \cdot b^*\right) = b^*a^*$
		\end{enumerate}
		We shall assume that the partial *-algebra $\mathfrak{A}$ contains a unit $e$, if
		$$
		e= e^*, \left( e, a\right) \in \Ga, \forall a\in \mathfrak{A}  \quad \text{and} \quad e \cdot a = a \cdot e = a,\quad  \forall a \in  \mathfrak{A}.
		$$
		If $\mathfrak{A}$ has no unit, it may always be embedded into a larger partial *-algebra with unit
		(the so-called unitization of $\mathfrak{A}$), in the standard fashion (cf. \cite{antoine:part_s}).
	\end{definition}

	\begin{definition}\label{qousi_star_defn}\label{quasi_defn}
		A \textit{quasi *-algebra} $\left(\mathfrak{A}, \mathfrak{A}_0\right)$ is a pair consisting of a vector space $\mathfrak{A}$
		and a *-algebra $\mathfrak{A}_0$ contained in $\mathfrak{A}$ as a subspace, such that
		\begin{itemize}
			\item [(a)] $\mathfrak{A}$ carries an involution $a\mapsto a^*$
			extending the involution of $\mathfrak{A}_0$;
			\item [(b)] $\mathfrak{A}$ is a bimodule over $\mathfrak{A}$ and the module multiplications extend the multiplication
			of $\mathfrak{A}_0$; In particular, the following associative laws hold:
			\be\label{guasi_ass_eqn}
			(xa)y = x(ay);\quad  a(xy) = (ax)y, \quad a \in \mathfrak{A},\quad x,y \in \mathfrak{A}_0; 
			\ee
			\item [(c)]
			$\left(ax\right)^*=x^*a^*$, for every $a \in \mathfrak{A}$ and $x\in \mathfrak{A}_0$.
		\end{itemize}
		We say that a quasi *-algebra $\left(\mathfrak{A}, \mathfrak{A}_0\right)$  is \textit{unital}, if there is an element $e \in \mathfrak{A}_0$,
		such that $ae = a = ea$, for all $a \in \mathfrak{A}$; $e$ is unique and called unit of $\left(\mathfrak{A}, \mathfrak{A}_0\right)$.
		We say that $\left(\mathfrak{A}, \mathfrak{A}_0\right)$ has a \textit{quasi-unit} if there exists an element $q \in \mathfrak{A}$, such that
		$qx = xq = x$, for every $x \in \mathfrak{A}_0$. It is clear that the unit $e$, if any, is a quasi-unit but
		the converse is false, in general.
	\end{definition}
	
	\begin{remark}
		If $\mathfrak{A}$ has no unit, it may always be embedded into a larger partial *-algebra with unit
		(the so-called unitization of $\mathfrak{A}$), in the standard fashion \cite{antoine:part_s}.
	\end{remark}
	
	\begin{definition}\label{qousi_star_tau_defn}\cite{quasi_star}
		Let  $\left(\mathfrak{A}, \mathfrak{A}_0\right)$ be a quasi *-algebra and $\tau$  a locally convex topology on A. We say that $\left(\mathfrak{A}\left[\tau\right], \mathfrak{A}_0\right)$ is a locally
		convex quasi *-algebra if
		\begin{enumerate}
			\item [(a)]  the map $a\in \mathfrak{A} \mapsto a^*\in\mathfrak{A}$ is continuous;
			\item[(b)] for every $x\in \mathfrak{A}_0$, the maps $a \mapsto ax$, $a\mapsto xa$ are continuous in $\mathfrak{A}\left[\tau\right]$;\\
			(c) $\mathfrak{A}_0$ is $\tau$-dense in $\mathfrak{A}$.
		\end{enumerate}
		
	\end{definition}
	\begin{definition}\label{quasi_hom_defn}
		Let $\left(\mathfrak{A}, \mathfrak{A}_0\right)$ be a quasi *-algebra and $\mathfrak{B}$ a partial *-algebra (cf. \cite{quasi_star}). A
		linear map 
		from $\mathfrak{A}$ into $\mathfrak{B}$ is called a homomorphism of $\left(\mathfrak{A}, \mathfrak{A}_0\right)$ into $\mathfrak{A}$ if,
		for $a \in\mathfrak{A}$ and $x\in\mathfrak{A}_0$ 
		$\Phi\left(a \right) \Phi\left( x\right)$ and 
		$\Phi\left( x\right) \Phi\left(a \right)$ are well defined in $\mathfrak{B}$ and $\Phi\left( a\right) \Phi\left( x\right)= 
		\Phi(ax)$, 
		$\Phi\left(x \right) \Phi\left(a \right)= 
		\Phi(xa)$, respectively. The homomorphism is called a *-\textit{homomorphism} if 
		$\Phi\left(a^*\right)
		) = \Phi\left(a\right)^*$, for every $a\in\mathfrak{A}$.
	\end{definition}
	
	\subsection{Partial $O^*$-algebras}
	\paragraph*{}
	Here I follow to \cite{antoine:part_o}. Define partial *-algebras of closable operators in Hilbert spaces. Let $\H$ be a Hilbert space with inner product $\left(\cdot, \cdot\right)_\H$ and $\D$ a dense subspace of $\H$. We denote by $\L\left(\D, \H \right)$ the set of all closable linear 
	operators $X$ in $\H$ such that $\D\left(X \right)  =\D$ and put
	\be\label{l_dag_eqn}
	\begin{split}
		\L\left( \D\right) \bydef\left\{\left.X\in \L\left(\D, \H \right)\right|X\D\subset \D \right\},\\
		\L^\dagger\left( \D, \H\right) \bydef\left\{\left.X\in \L\left(\D, \H \right)\right|\D\left(X^* \right)\supset \D \right\},\\
		\L^\dagger\left( \D\right) \bydef\left\{\left.X\in \L^\dagger\left(\D, \H \right)\cap \L\left( \D\right) \right|X^*\D\subset \D \right\}.
	\end{split}
	\ee
	Then $\L^\dagger\left( \D, \H\right)$ is a vector space 
	with the usual operations: $X + Y$, $\la X$, and $\L\left(\D \right)$ is a subspace of $\L^\dagger\left( \D, \H\right)$ and 
	an algebra with the usual 
	multiplication $XY$
	For $\L^\dagger\left(\D, \H \right)$ 
	and $\L^\dagger\left(\D\right)$, we 
	have 
	the following
	\begin{proposition}\cite{antoine:part_o}
		$\L^\dagger\left(\D, \H \right)$ is a partial *-algebra with respect to the following operations: the sum $X + Y$, the scalar multiplication $\la X$, the involution $X\mapsto X^\dagger\bydef X^*|_\D$  and the (weak) partial multiplication $X\Hsquare Y\bydef X^{\dagger *} Y$, defined whenever $X$ is a left  multiplier of $Y$, ($X \in L^{\mathrm{w}}\left(Y \right)$ or $Y \in R^{\mathrm{w}}\left(X \right)$), that is, iff $Y\D\subset\D\left(X^{\dagger *} \right)$ and $X^\dagger\D\subset \D\left( Y^*\right)$. $\L^\dagger\left(\D\right)$ is a *-algebra 
		with the usual multiplication $XY$ (which here coincides with the weak multiplication $\Hsquare$) and the involution $X \mapsto X^\dagger$.
	\end{proposition} 
	\begin{remark}\label{o*b_rem}
		If 
		\be
		\begin{split}\label{o*b_eqn}
			\L^\dagger\left( \D, \H\right)_b\bydef\left\{\left.X\in \L^\dagger\left( \D, \H\right)\right|\overline{X}\in B\left(\H \right) \right\},\\
			\L^\dagger\left( \D\right)_b\bydef \L^\dagger\left( \D, \H\right)_b\cap \L^\dagger\left(\D \right) 
		\end{split}
		\ee
		then the  inclusion $\L^\dagger\left( \D, \H\right)_b\subset B\left(\H \right)$ induces a norm  on $\L^\dagger\left( \D, \H\right)_b$ which is a $C^*$-norm on $\L^\dagger\left( \D\right)_b$.
	\end{remark}
	
	\begin{definition}\label{o*alg_defn}\cite{antoine:part_o}
		A subset (subspace) of $\L^\dagger\left( \D, \H\right)$ is 
		called an $O$-\textit{family} ($O$-\textit{vector space}) 
		on $\D$, and a subalgebra of $\L\left( \D\right)$
		is called an $O$-\textit{algebra} on $\D$. A $\dagger$-
		invariant subset (subspace) of subset (subspace) of $\L^\dagger\left( \D, \H\right)$ is 
		called an $O^*$-\textit{family} ($O^*$-\textit{vector space}) 
		on $\D$. A partial *-subalgebra of $\L^\dagger\left( \D, \H\right)$ is called a \textit{partial} $O^*$-\textit{algebra} on $\D$, and a 
		*-subalgebra of  $\L^\dagger\left( \D\right)$ is 
		called an $O^*$-\textit{algebra} on $\D$. 
	\end{definition}
	\begin{remark}\label{*_bound_rem}
		If $A$ is a  $O^*$-\textit{algebra} on $\D$, i.e, 
		$A \subset \L^\dagger\left( \D\right)$ then the subalgebra
		\be\label{*_bound_eqn}
		A_b \bydef A\cap \L^\dagger\left( \D\right)_b
		\ee
		and $C^*$-norm on 	it depend on multiplication and *-operation on $A$. The *-algebra	$A_b$
		and $C^*$-norm on 	it do not depend in the inclusion $A \subset \L^\dagger\left( \D\right)$
	\end{remark}

	\subsection{Quasi *-algebras and partial *-algebras of operators}
	\paragraph{}
	Let $\D$ be a dense subspace of a Hilbert space $\H$.
	We denote by $\L^\dagger\left( \D, \H\right)$  the set of all (closable) linear operators $X$ in $\H$, such that $D\left( X\right) = \D$, $D\left( X^*\right)\supseteqq \D$ where $D\left( X\right)$ denotes the domain of $X$. The set $\L^\dagger\left( \D, \H\right)$ is a partial *-algebra with respect to the following operations:
	the usual sum $X_1 + X_2$, the scalar multiplication $\la X$, the involution $X \mapsto X^\dagger\bydef X^*\D$ and the (weak) partial multiplication $X1\square X_2\bydef X_1^{\dagger *}
	X_2$ (where $X_1^{\dagger *}\bydef  \left( X_1^{\dagger}\right)^*$). The latter is defined whenever $X_2$ is a \textit{weak right multiplier} of  $X_2$ (for this, we shall write $X_2\in R^w\left( X_1\right) $ or  $X_1\in L^w\left( X_1\right) $, that is, if and only if, $X_2\D\subset D\left(X_1^{\dagger *}\right)$ and $X_1^\dagger\D\in D\left(X_2^* \right)$. The operator $1_\D$, restriction to $\D$ of the identity
	operator $1_\H$ on $\H$, is the unit of the partial *-algebra $\L^\dagger\left( \D, \H\right)$. By $\L^\dagger\left( \D, \H\right)_b$ we shall denote the bounded part of $\L^\dagger\left( \D, \H\right)$; i.e.,
	$$
	\L^\dagger\left( \D, \H\right)_b\bydef\left\{\left.X\in \L^\dagger\left( \D, \H\right)\right|\overline{X}\in B\left(\H \right) \right\}
	$$
	where $\overline{X}$ is the closure of $X$, i.e., a minimal closed extension of $X$.
	Let us denote by $\L^\dagger\left(\D \right)$  the space of all linear operators $X:\D\to\D$, having an adjoint $X^\dagger : \D\to\D$, by which we simply mean that 
	$$
	\left(X\xi, \eta\right)= \left(\xi, X^\dagger\eta\right)\quad\forall \xi, \eta\in \D.
	$$
	If $\L^\dagger\left( \D, \H\right)$ is endowed with the strong *-topology $\mathrm{t}_{s^*}$ , defined by the set of seminorms
	$$
	p_\xi(X)\bydef \left\| X\xi\right\|  + \left\| X^\dagger\xi\right\| \quad \xi\in \D, \quad  X\in \L^\dagger\left( \D, \H\right),
	$$
	then $\left(\L^\dagger\left( \D, \H\right)\left[\mathrm{t}_{s^*}\right],\L\dagger(D)_b)\right)$ is a locally convex quasi *-algebra or, more precisely, a \textit{locally convex quasi $C^*$-normed algebra}.
	If $\L^\dagger\left( \D, \H\right)$ is endowed with the weak topology $\mathrm{t}_w$ defined by the set of seminorms
	$$
	p_{\xi,\eta}(X)\bydef \left| \left(X\xi, \eta \right) \right|   \quad \xi, \eta\in \D, \quad  X\in \L^\dagger\left( \D, \H\right),
	$$
	then, again, $\left(\L^\dagger\left( \D, \H\right)\left[\mathrm{t}_{w}\right],\L\dagger(D)_b)\right)$. is a locally convex quasi *-algebra.
	\begin{definition}
		Let $\left(\mathfrak{A}, \mathfrak{A}_0\right)$ be a quasi *-algebra and $\D_\pi$ a dense domain in a
		certain Hilbert space $\H_\pi$ . A linear map $\pi$ from $\mathfrak{A}$ into $\L^\dagger\left(\H_\pi, \D_\pi \right)$ is called a
		*-representation of $\left(\mathfrak{A}, \mathfrak{A}_0\right)$, if the following properties are fulfilled:
		\begin{enumerate}
			\item [(a)] $\pi\left( a^*\right) = \pi\left( a\right)^\dagger\quad \forall a \in pi\left( a^*\right)$.
			\item[(b)]  for $a\in \mathfrak{A}$ and $x\in \mathfrak{A}$ $\pi\left( a\right)\square \pi\left( x\right)$  is well defined and $\pi\left( a\right)\square \pi\left( x\right)=\pi\left(ax \right)$ 
		\end{enumerate}
		
		In other words,  $\pi$ is a $*$-homomorphism of $\left(\mathfrak{A}, \mathfrak{A}_0\right)$ into the partial *-algebra
		$\L^\dagger\left(\D_\pi, \H_\pi \right)$
		
		If $\left(\mathfrak{A}, \mathfrak{A}_0\right)$ has a unit $e \in \mathfrak{A}_0$, we assume $\pi\left(e \right)  = \Id_{\D_\pi}$, where $\Id_{\D_\pi}$ is the identity operator on the space $\D_\pi$. If $\pi_o\bydef\pi|_{ \mathfrak{A}_0}$  is a *-representation of the *-algebra $\mathfrak{A}_0$ into $\L^\dagger\left(\D_\pi \right) $ we say that $\pi$ is a qu*-representation.
	\end{definition}
	
	\begin{empt}If
		$$
		\mathfrak{A}^+_0\bydef\left\{\sum_{k=1}^n x^*_kx_k\quad x_k \in  \mathfrak{A}_0\quad n\in \N \right\}
		$$
		
		then, $\mathfrak{A}^+_0$ is a wedge in $\mathfrak{A}_0$ and we call the elements of $\mathfrak{A}^+_0$ \textit{positive elements} of $\mathfrak{A}_0$. We call \textit{positive elements} of $\mathfrak{A}\left[\tau\right]$ the elements of closure (with respect to $\tau$) of  $\mathfrak{A}^+_0$ in 
		and we denote them by $\mathfrak{A}^+$.
		
	\end{empt}
	\begin{proposition}
		If $a \ge 0$, then $\pi\left(a \right)\ge 0$ , for every $\left(\tau,\mathrm{t}_{w} \right)$-continuous *-representation of $\left(\mathfrak{A}\left[\tau\right], \mathfrak{A}_0\right)$.
	\end{proposition}
	\begin{theorem}
		Assume that $\mathfrak{A}^+\cap \left(-\mathfrak{A}^+\right)= \{0\}$. Let $a \in \mathfrak{A}^+$, $a\neq 0$. Then, there
		exists a continuous linear functional $\om$ on $\mathfrak{A}\left[\tau\right]$ with the properties:
		\begin{enumerate}
			\item[(i)] $\om\left(b \right)\ge 0\quad \forall b \in \mathfrak{A}^+$;
			\item[(ii)] $\om\left(a \right)> 0$.
		\end{enumerate}
		
	\end{theorem}
	
	If $\left(\mathfrak{A}\left[\tau\right], \mathfrak{A}_0\right)$ is an arbitrary locally convex quasi *-algebra then there is a natural order related to the topology $\tau$. This order can be  used to define bounded elements. In what follows, we shall assume that $\left(\mathfrak{A}\left[\tau\right], \mathfrak{A}_0\right)$
	has a unit $a$. Let $a\in \mathfrak{A}$; put  $\mathfrak{R}\left( a\right)\bydef \frac{1}{2}\left(a + a^* \right)$ $\mathfrak{J}\left( a\right)\bydef \frac{1}{2i}\left(a - a^* \right)$
	Then, $\mathfrak{R}\left( a\right),\mathfrak{J}\left( a\right) \in\mathfrak{A}_h$ and $a = \mathfrak{R}\left( a\right)+i\mathfrak{J}\left( a\right)$
	\begin{definition}
		An element $a\in $ is called \textit{order bounded} if there exists $\ga \ge 0$, such that 
		$$\pm \mathfrak{R}\left( a\right) \le \ga e,\quad \pm \mathfrak{J}\left( a\right) \le \ga e$$.
		We denote by $\mathfrak{A}^{\mathrm{or}}_{\mathrm{b}}$
		b the set of all order bounded elements of $\mathfrak{A}\left[\tau\right]$.
	\end{definition}

	We recall that an unbounded $C^*$-seminorm $p$ on a partial *-algebra $\mathfrak{A}$ is a seminorm defined on a partial $*$-subalgebra $D(p)$ of $\mathfrak{A}$, the domain of $p$, with the properties:
	\begin{itemize}
		\item $p(ab) ≤ p(a)p(b)$, whenever $ab$ is well-defined;
		\item $p\left(a^*a \right) = p\left(a \right)^2$  whenever $a^*a$ is well-defined
	\end{itemize}
	
	\begin{proposition}
		$\left\|\cdot \right\|^{\mathrm{or}}_{\mathrm{b}}$ is an unbounded $C^*$-norm on $\mathfrak{A}$ with domain $\mathfrak{A}_b$.
	\end{proposition}
\section{Representations by unbounded operators on Hilbert modules}
\label{sec:rep_Hilbert_module}
\paragraph*{} 
Let $A $ be a unital $*$-algebra,  $D $ a $C^*$-algebra,
and $\mathcal{E} $ a Hilbert  $D $-module.  

\begin{definition}
	\label{def:rep_Hilbert_module}\cite{meyer:unb_repr}
	A \emph{representation} of $A $ on $\mathcal{E} $ is a
	pair $(\mathfrak{E},\pi) $, where  $\mathfrak{E}\subseteq\mathcal{E} $ is a dense
	$D $-submodule and  $\pi\colon A\to\End_D(\mathfrak{E}) $ is a unital
	algebra homomorphism to the algebra of  $D$-module endomorphisms
	of $\mathfrak{E} $, such that
	\[
	\braket{\pi(a)\xi}{\eta}_D = \braket{\xi}{\pi(a^*)\eta}_D
	\qquad
	\text{for all }a\in A,\ \xi,\eta\in\mathfrak{E}.
	\]
	
	We call $\mathfrak{E} $ the \emph{domain} of the representation.  We
	may drop $\pi $ from our notation by saying that $\mathfrak{E} $ is an
	$A,D $-bimodule with the right module structure inherited
	from $\mathcal{E} $, or we may drop $\mathfrak{E} $ because it is the common
	domain of the partial linear maps $\pi(a) $ on $\mathcal{E} $ for all
	$a\in A $.
	We equip $\mathfrak{E} $ with the \emph{graph topology}, which is
	generated by the \emph{graph norms}
	\be	\label{graph_norm_eqn}
	\left\| {\xi}\right\| _a \bydef \left\langle \xi, \xi\right\rangle  + \left\langle  \pi(a)\xi,\pi(a)\xi\right\rangle^{\nicefrac12}	= \left\langle  \xi,\pi(1+a^*a)\xi\right\rangle ^{\nicefrac12}
	\ee
	for  $a\in A $.  The representation is \emph{closed} if $\mathfrak{E} $
	is complete in this topology.  A \emph{core} for $(\mathfrak{E},\pi) $ is
	an  $A,D $-subbimodule of $\mathfrak{E} $ that is dense in $\mathfrak{E} $ in
	the graph topology.
\end{definition}

Definition~\ref{def:rep_Hilbert_module}  for  $D=\C $ is the usual
definition of a representation of a $*$-algebra on a Hilbert
space by unbounded operators.   For  $\mathcal{E}=D $ with the
canonical Hilbert  $D $-module structure, we get representations
of $A $ by \emph{densely defined unbounded multipliers}.  The
domain of such a representation is a dense right ideal
$\mathfrak{E}\subseteq D $.  

Given two norms  $p,q $, we write  $p \preceq q $ if there is a
scalar  $c>0 $ with  $p \le c q $.

\begin{lemma}
	\label{lem:graph_norms_directed}\cite{meyer:unb_repr}
	The set of graph norms partially ordered by $\preceq $ is
	directed: for all  $a_1,\dotsc,a_n\in A $ there are  $b\in A $ and
	$c\in\R_{>0} $ so that  $\mathrm{norm}{\xi}_{a_i} \le c\mathrm{norm}{\xi}_b $ for
	any representation  $(\mathfrak{E},\pi) $, any  $\xi\in\mathfrak{E} $, and
	$i=1,\dotsc,n $.
\end{lemma}

\begin{definition}
	\label{def:regular}
	A densely defined operator $t $ on a Hilbert module $\mathcal{E} $ is
	\emph{semiregular} if its adjoint is also densely defined; it is
	\emph{regular} if it is closed, semiregular and  $1+t^*t $ has dense
	range.  An \emph{affiliated multiplier} of a
	$C^*$-algebra $D $ is a regular operator on $D $ viewed as
	a Hilbert  $D $-module.
\end{definition}

The usual norm on $\mathcal{E} $ is the graph norm for  $0\in A $.  Hence
the inclusion map  $\mathfrak{E}\hookto\mathcal{E} $ is continuous for the graph
topology on $\mathfrak{E} $ and extends continuously to the
completion $\overline{\mathfrak{E}} $ of $\mathfrak{E} $ in the graph topology.

\begin{proposition}
	\label{pro:closure_rep}
	The canonical map  $\overline{\mathfrak{E}}\to\mathcal{E} $ is injective, and its
	image is
	\begin{equation}
		\label{eq:domain_closure}
		\overline{\mathfrak{E}} = \bigcap_{a\in A} \Dom \overline{\pi(a)}.
	\end{equation}
	Thus $(\mathfrak{E},\pi) $ is closed if and only if  $\mathfrak{E} =
	\bigcap_{a\in A} \Dom \overline{\pi(a)} $.  Each  $\pi(a) $ extends
	uniquely to a continuous operator $\overline{\pi(a)} $
	on $\overline{\mathfrak{E}} $.  This defines a closed
	representation $(\overline{\mathfrak{E}},\overline{\pi}) $ of $A $, called the
	\emph{closure} of $(\mathfrak{E},\pi) $.
\end{proposition}

We shall need a generalisation of~\eqref{eq:domain_closure} that
replaces $A $ by a sufficiently large subset.

\begin{definition}
	\label{def:strong_generating_set}
	A subset  $S\subseteq A $
	is called a \emph{strong generating} set if it generates $A $
	as an algebra and the graph norms for  $a\in S $
	generate the graph topology in any representation.  That is, for any
	representation on a Hilbert module, any vector $\xi $
	in its domain and any  $a\in A $,
	there are  $c\ge1 $
	in $\R $
	and  $b_1,\dotsc,b_n\in S $
	with  $\mathrm{norm}{\xi}_a \le c \sum_{i=1}^n \mathrm{norm}{\xi}_{b_i} $.
\end{definition}

An estimate  $\mathrm{norm}{\xi}_a \le c \sum_{i=1}^n \mathrm{norm}{\xi}_{b_i} $
is usually shown by finding  $d_1,\dotsc,d_m\in A $
with
$a^* a + \sum_{j=1}^m d_j^* d_j = c\cdot \sum_{i=1}^n b_i^* b_i $,
compare the proof of Lemma~\ref{lem:graph_norms_directed}.

\begin{proposition}
	\label{pro:equality_if_closure_equal}
	Let  $S\subseteq A $ be a strong generating set.  Two closed
	representations  $(\mathfrak{E}_1,\pi_1) $ and $(\mathfrak{E}_2,\pi_2) $
	of $A $ on the same Hilbert module $\mathcal{E} $ are equal if and
	only if  $\overline{\pi_1(a)} = \overline{\pi_2(a)} $ for all  $a\in S $.
\end{proposition}

Here we consider 

\begin{corollary}
	\label{cor:bounded_rep}
	Let $S $ be a strong generating set of $A $ and
	let $(\mathfrak{E},\pi) $ be a closed
	representation of $A $ with  $\Dom \overline{\pi(a)} = \mathcal{E} $ for each
	$a\in S $.  Then  $\mathfrak{E}=\mathcal{E} $ and $\pi $ is a
	$*$-homomorphism to the $C^*$-algebra  $\mathbb{B}(\mathcal{E}) $ of
	adjointable operators on $\mathcal{E} $.
\end{corollary}

An \emph{isometry}  $I\colon \mathcal{E}_1\hookto\mathcal{E}_2 $
between two Hilbert  $D $-modules  $\mathcal{E}_1 $
and $\mathcal{E}_2 $
is a right  $D $-module
map with  $\braket{I\xi_1}{I\xi_2} = \braket{\xi_1}{\xi_2} $
for all  $\xi_1,\xi_2\in\mathcal{E}_1 $.

\begin{definition}
	\label{def:isometric_intertwiner}
	Let  $(\mathfrak{E}_1,\pi_1) $
	and  $(\mathfrak{E}_2,\pi_2) $
	be representations on Hilbert  $D $-modules
	$\mathcal{E}_1 $
	and $\mathcal{E}_2 $,
	respectively.  An \emph{isometric intertwiner} between them is an
	isometry  $I\colon \mathcal{E}_1\hookto\mathcal{E}_2 $
	with  $I(\mathfrak{E}_1)\subseteq \mathfrak{E}_2 $
	and  $I\circ\pi_1(a) (\xi) = \pi_2(a) \circ I (\xi) $
	for all  $a\in A $,
	$\xi\in\mathfrak{E}_1 $;
	equivalently,  $I\circ\pi_1(a)\subseteq \pi_2(a)\circ I $
	for all  $a\in A $,
	that is, the graph of $\pi_2(a)\circ I $
	contains the graph of $I\circ \pi_1(a) $.
	We neither ask $I $
	to be adjointable nor  $I(\mathfrak{E}_1)=\mathfrak{E}_2 $.
	Let  $\mathrm{Rep}(A,D) $
	be the category with closed representations of $A $
	on Hilbert  $D $-modules
	as objects, isometric intertwiners as arrows, and the usual
	composition.  The unit arrow on $(\mathfrak{E},\pi) $
	is the identity operator on $\mathcal{E} $.
\end{definition}

\begin{lemma}
	\label{lem:closure_functorial}
	Let  $(\mathfrak{E}_1,\pi_1) $ and  $(\mathfrak{E}_2,\pi_2) $ be
	representations on Hilbert  $D $-modules  $\mathcal{E}_1 $
	and $\mathcal{E}_2 $, respectively, and let  $I\colon
	\mathcal{E}_1\hookto\mathcal{E}_2 $ be an isometric intertwiner.  Then $I $ is
	also an intertwiner between the closures of  $(\mathfrak{E}_1,\pi_1) $
	and  $(\mathfrak{E}_2,\pi_2) $.
\end{lemma}
\begin{proposition}
	\label{pro:intertwiner_strong_generators}
	Let  $(\mathfrak{E}_1,\pi_1) $ and  $(\mathfrak{E}_2,\pi_2) $ be closed
	representations of $A $ on Hilbert  $D $-modules  $\mathcal{E}_1 $
	and $\mathcal{E}_2 $, respectively.  Let  $S\subseteq A $ be a strong
	generating set.  An isometry  $I\colon \mathcal{E}_1\hookto \mathcal{E}_2 $ is
	an intertwiner from  $(\mathfrak{E}_1,\pi_1) $ to  $(\mathfrak{E}_2,\pi_2) $
	if and only if  $I\circ \overline{\pi_1(a)} \subseteq \overline{\pi_2(a)}\circ
	I $ for all  $a\in S $.
\end{proposition}

Now we relate the categories  $\mathrm{Rep}(A,D) $ for different
$C^*$-algebras $D $.

\begin{definition}
	\label{def:Cstar-correspondence}
	Let  $D_1 $ and $D_2 $ be two $C^*$-algebras.  A
	\emph{$C^*$-correspondence} from $D_1 $ to $D_2 $ is a
	Hilbert  $D_2 $-module with a representation of $D_1 $ by
	adjointable operators (representations of $C^*$-algebras are
	tacitly assumed nondegenerate).  An \emph{isometric intertwiner}
	between two correspondences from $D_1 $ to $D_2 $ is an
	isometric map on the underlying Hilbert  $D_2 $-modules that
	intertwines the left   $D_1 $-actions.  Let  $\mathrm{Rep}(D_1,D_2) $
	denote the category of correspondences from $D_1 $ to $D_2 $
	with isometric intertwiners as arrows and the usual composition.
\end{definition}

Let $\mathcal{E} $ be a Hilbert  $D_1 $-module and $\F $ a
correspondence from $D_1 $ to $D_2 $.  The interior tensor product
$\mathcal{E}\otimes_{D_1} \F $ is the (Hausdorff) completion of the
algebraic tensor product  $\mathcal{E}\odot \F $ to a Hilbert
$D_2 $-module, using the inner product
\begin{equation}
	\label{eq:interior_tensor}
	\braket{\xi_1\otimes\eta_1}{\xi_2\otimes\eta_2}
	= \braket{\eta_1}{\braket{\xi_1}{\xi_2}_{D_1}\cdot\eta_2}_{D_2}.
\end{equation}
We may use the balanced tensor product
$\mathcal{E}\odot_{D_1} \F $ instead of  $\mathcal{E}\odot \F $
because the inner product~\eqref{eq:interior_tensor} descends to this
quotient.  If we want to emphasise the left  action  $\varphi\colon
D_1\to\mathbb{B}(\F) $ in the
$C^*$-correspondence $\F $, we write  $\mathcal{E}\otimes_\varphi
\F $ for  $\mathcal{E}\otimes_{D_1} \F $.

In addition, let $(\mathfrak{E},\pi) $ be a closed representation of $A $
on $\mathcal{E} $.  We are going to build a closed representation
$(\mathfrak{E}\otimes_{D_1}\F, \pi\otimes_{D_1}1) $ of $A $ on
$\mathcal{E}\otimes_{D_1}\F $.  First let  $X\subseteq
\mathcal{E}\otimes_{D_1}\F $ be the image of  $\mathfrak{E}\odot \F $
or  $\mathfrak{E}\odot_{D_1} \F $ under the canonical map to
$\mathcal{E}\otimes_{D_1}\F $.

\begin{lemma}
	\label{lem:tensor_rep_with_corr}
	For  $a\in A $, there is a unique linear operator
	$\pi(a)\otimes1\colon X\to X $ with  $(\pi(a)\otimes1)
	(\xi\otimes\eta) = \pi(a)(\xi)\otimes \eta $ for all
	$\xi\in\mathfrak{E} $,  $\eta\in\F $.  The map  $a\mapsto
	\pi(a)\otimes 1 $ is a representation of $A $ with domain $X $.
\end{lemma}

\begin{definition}
	Let  $(\mathfrak{E}\otimes_{D_1}\F, \pi\otimes_{D_1}1) $ be the
	closure of the representation on $\mathcal{E}\otimes_{D_1}\F $
	defined in Lemma~\ref{lem:tensor_rep_with_corr}.
\end{definition}

\begin{lemma}
	\label{lem:rep_tensor_corr_functor}
	Let  $I\colon \mathcal{E}_1\hookto\mathcal{E}_2 $ be an isometric intertwiner
	between two representations  $(\mathfrak{E}_1,\pi_1) $ and
	$(\mathfrak{E}_2,\pi_2) $, and let  $J\colon \F_1\hookto \F_2 $
	be an isometric intertwiner of $C^*$-correspondences.  Then
	$I\otimes_{D_1} J\colon \mathcal{E}_1\otimes_{D_1}\F_1 \hookto
	\mathcal{E}_2\otimes_{D_1}\F_2 $ is an isometric intertwiner between
	$(\mathfrak{E}_1\otimes_{D_1}\F_1, \pi_1\otimes1) $ and
	$(\mathfrak{E}_2\otimes_{D_1}\F_2, \pi_2\otimes1) $.
\end{lemma}

The lemma gives a bifunctor
\begin{equation}
	\label{eq:interior_tensor_bifunctor}
	\otimes_{D_1}\colon \mathrm{Rep}(A,D_1)\times \mathrm{Rep}(D_1,D_2) \to
	\mathrm{Rep}(A,D_2).
\end{equation}
The corresponding bifunctor
\[
\otimes_{D_1}\colon \mathrm{Rep}(B,D_1)\times \mathrm{Rep}(D_1,D_2) \to \mathrm{Rep}(B,D_2)
\]
for a $C^*$-algebra $B $ is the usual composition of
$C^*$-correspondences.  This composition is associative up to
canonical unitaries
\begin{equation}
	\label{eq:tensor_associative}
	\mathcal{E} \otimes_{D_1} (\F \otimes_{D_2} \mathcal{E}) \xrightarrow{\cong}
	(\mathcal{E} \otimes_{D_1} \F) \otimes_{D_2} \mathcal{E},\qquad
	\xi \otimes (\eta\otimes \zeta) \mapsto
	(\xi \otimes \eta)\otimes \zeta,
\end{equation}
for all triples of composable $C^*$-correspondences.

\begin{lemma}
	\label{lem:tensor_associative}
	If $\mathcal{E} $ carries a representation $(\mathfrak{E},\pi) $ of a
	$*$-algebra $A $, then the unitary
	in~\eqref{eq:tensor_associative} is an intertwiner  $(\mathfrak{E},\pi)
	\otimes_{D_1} (\F\otimes_{D_2} \mathcal{E}) \xrightarrow{\cong}
	\bigl((\mathfrak{E},\pi) \otimes_{D_1} \F\bigr) \otimes_{D_2}
	\mathcal{E}$.
\end{lemma}

\begin{definition}
	\label{def:star-intertwiner}
	Let  $(\mathfrak{E}_1,\pi_1) $ and $(\mathfrak{E}_2,\pi_2) $ be two
	representations of $A $ on Hilbert  $D $-modules  $\mathcal{E}_1 $
	and $\mathcal{E}_2 $.  An adjointable operator  $x\colon
	\mathcal{E}_1\to\mathcal{E}_2 $ is an \emph{intertwiner} if
	$x(\mathfrak{E}_1)\subseteq \mathfrak{E}_2 $ and  $x\pi_1(a)\xi =
	\pi_2(a)x\xi $ for all  $a\in A $,  $\xi\in\mathfrak{E}_1 $.  It is
	a \emph{$*$-intertwiner} if both  $x $ and $x^* $ are
	intertwiners.
	\label{def:adjointable_intertwiner}
\end{definition}

Any adjointable intertwiner between two representations of a
$C^*$-algebra $B $ is a $*$-intertwiner.  In contrast, for
a general $*$-algebra, even the adjoint of a unitary
intertwiner $u $ fails to be an intertwiner if
$u(\mathfrak{E}_1)\subsetneq \mathfrak{E}_2 $.

\begin{example}
	\label{exa:Friedrichs_extension}
	Let $t $ be a positive symmetric operator on a Hilbert
	space $ \H $.  Assume that  $\bigcap_{n\in\N} \Dom t^n $ is
	dense in $ \H $, so that $t $ generates a
	representation $\pi $ of the polynomial algebra $\C[x] $
	on $ \H $.  The Friedrichs extension of $t $ is a positive
	self-adjoint operator $t' $ on $ \H $.  It generates
	another representation $\pi' $ of $\C[x] $ on $ \H $.  The
	identity map on $ \H $ is a unitary intertwiner
	$\pi\hookto\pi' $.  It is not a $*$-intertwiner unless
	$t=t' $.
\end{example}

The following proposition characterises when an adjointable isometry
$I\colon \mathcal{E}_1\hookto \mathcal{E} $ between two representations
on Hilbert  $D $-modules is a $*$-intertwiner.  Since
$\mathcal{E}\cong \mathcal{E}_1 \oplus \mathcal{E}_1^\bot $ if $I $ is adjointable, we
may as well assume that $I $ is the inclusion of a direct summand.

\begin{proposition}
	\label{pro:isometry_star-intertwiner}
	Let  $\mathcal{E}_1 $ and $\mathcal{E}_2 $ be Hilbert modules over a
	$C^*$-algebra $D $ and let  $(\mathfrak{E}_1,\pi_1) $ and
	$(\mathfrak{E},\pi) $ be representations of $A $ on  $\mathcal{E}_1 $
	and $\mathcal{E}_1\oplus\mathcal{E}_2 $, respectively.  The following are
	equivalent:
	\begin{enumerate}
		\item \label{enum:isometry_star-intertwiner1} the canonical
		inclusion  $I\colon \mathcal{E}_1\hookto\mathcal{E}_1\oplus\mathcal{E}_2 $ is a
		$*$-intertwiner from $\pi_1 $ to $\pi $;
		\item \label{enum:isometry_star-intertwiner2} the canonical
		inclusion  $I\colon \mathcal{E}_1\hookto\mathcal{E}_1\oplus\mathcal{E}_2 $ is an
		intertwiner from $\pi_1 $ to $\pi $ and  $\mathfrak{E} =
		\mathfrak{E}_1 + (\mathfrak{E}\cap \mathcal{E}_2) $;
		\item \label{enum:isometry_star-intertwiner3} there is a
		representation  $(\mathfrak{E}_2,\pi_2) $ on $\mathcal{E}_2 $ such that  $\pi
		= \pi_1 \oplus \pi_2 $.
	\end{enumerate}
\end{proposition}

\section{Discrete crossed product of $C^*$-algebras}\label{discr_cr_prod_sec}	
Here I follow to \cite{discrete_crossed}.
Let $G$ be a discrete group and $A$ a $C^*$-algebra. An action of $G$ on $A$ is a group homomorphism from $G$ into $\Aut(A)$ — the group of *-automorphisms on $A$. The action is denoted by a dot in the following way
$$
t \mapsto (a \mapsto t. a),\quad t \in G,\quad a \in A.
$$
A $C^*$-dynamical system $(A,G)$ consist of a discrete group $G$, a $C^*$-algebra $A$ and an action of $G$ on $A$. For a $C^*$-dynamical system $(A,G)$ with $G$ discrete let $C_c(G,A)$ be the linear span of finitely supported functions on G with values in $A$. A typical element a in $C_c(G,A)$ is written as a sum
\be\label{discr_cr_prod_cc_eqn}
\sum _{t \in G} a_t u_t\quad \forall t \in G \quad a_t  \in A,
\ee
where only finitely many $a_t$’s are non-zero. One equips $C_c(G,A)$ with a *-operation and a twisted product, in such
a way that the action becomes ”inner”. By this we mean that $t. a = u_tau^*_t$ for every $t \in G$, $a \in A$. More precisely we have the following product and *-operation
\be\label{discr_cr_prod_op_eqn}
\begin{split}
	ab \bydef \sum_{s,t\in G}a_t(t\cdot b_s)u_{ts}, \quad
	a^* =\sum_{t\in G} (t^{-1}.a^*_t )u_{t^{-1}} ,\\ a =\sum _{t \in G} a_t u_t, \quad b =\sum _{s \in G} a_s u_s
\end{split}
\ee
Fix a $C^*$-dynamical system $(A,G)$ with $G$ discrete. A covariant representation $\left(\pi, u, \H \right)$  of $(A,G)$ consists of an unitary representation $u: G \to B(\H)$ of $G$ and a representation $\pi  : A \to B(\H)$ of $A$ such that $u(t)π(a)u(t)^* = \pi(t. a)$
for every $t \in  G$ and $a \in A$. For a covariant representation $\left(\pi, u, \H \right)$ let $\pi\times u$ denote the associated representation of $C_c(G,A)$ on $\H$. The \textit{full $C^*$-algebra norm} on $C_c(G,A)$ is given by
$$
\left\|\cdot \right\|  \bydef \sup\left\| (\pi \times u)(\cdot )\right\| 
$$
where the supremum is taken over all covariant representations  $\left(\pi, u, \H \right)$ of $(A,G)$. The \textit{full crossed product}, denoted $A \rtimes  G$, is the completion of
$C_c(G,A)$ with respect to the full $C^*$-algebra norm.
The \textit{reduced $C^*$-algebra norm} on $C_c(G,A)$ is given by
\be\label{discr_red_n_eqn}
\left\|\cdot \right\|_\la  \bydef \sup\left\| \widetilde\pi \times \la(\cdot )\right\|
\ee
where the regular representation $\widetilde\pi : C_c(G,A) \to B(\widetilde \H )$ is the representation associated to the covariant representation given by 
\be\label{discr_red_repr_eqn}
\widetilde\pi\left(\delta_{t, \xi} \right) \bydef \delta_{t, \pi\left(t^{-1}.a\right) \xi}\quad \la(s)\delta_{st, \xi} \bydef \delta_{t, \xi}, \quad a \in A, \quad t, s \in G, \quad \xi \in \H
\ee
where $\pi: A \hookto B\left(\H\right)$ is any faithful representation, $\widetilde \H\bydef \ell^2\left(G, \H\right) $ and $\delta_{t, \xi}$ is the map $s \mapsto \delta_{t, s}\xi \in \H$ ($ \delta_{t, s}$ is the Kronecker delta). In this way $\widetilde\pi$ becomes faithful. The reduced $C^*$-algebra norm does not depend on the choice of the faithful representation $\pi$ (cf.  \cite{discrete_crossed}). The \textit{reduced crossed product}, denoted by $A\rtimes_rG$ is the completion of $C_c(G,A)$ with respect to the reduced $C^*$-algebra norm. 


\chapter{Spectral triples}
	
	\paragraph{}
	This section contains citations of  \cite{hajac:toknotes}. 
	\section{Definition of spectral triples}\label{sp_tr_defn_sec}
	\begin{defn}
		\label{df:spec-triple}\cite{hajac:toknotes}
		A (unital) {\it {spectral triple}} $(\A, \H, D)$ consists of:
		\begin{itemize}
			\item
			an unital pre-$C^*$-algebra $\A$ with an involution $a \mapsto a^*$, equipped
			with a faithful representation on:
			\item
			a \emph{Hilbert space} $\H$; and also
			\item
			a \emph{self-adjoint operator} $D$ on $\mathcal{H}$, with dense domain
			$\Dom D \subset \H$, such that $a(\Dom D) \subseteq \Dom D$ for all 
			$a \in \mathcal{A}$.
		\end{itemize}
	\end{defn}
	\begin{defn}
		\label{df:spt-real_defn}\cite{hajac:toknotes}
		A {\it real  spectral triple} is a  spectral triple $(\sA, \sH, D)$,
		together with an antiunitary operator $J\:\sH\to\sH$ such that
		$J(\Dom D)\subset \Dom D$, and 
	\be\label{df:spt-real_eqn}
	[a, Jb^*J^\dagger] = 0
	\ee for all
		$a, b \in \sA$.
	\end{defn}
	
	\begin{defn}
		\label{df:spt-even}\cite{hajac:toknotes}
		A  spectral triple $(\sA, \sH, D)$ is {\it even} if there is
		a selfadjoint unitary \textit{grading operator} $\Ga$ on $\sH$ such that $a\Ga = \Ga a$
		for all $a \in \sA$, $\Ga(\Dom D) = \Dom D$, and $D\Ga = -\Ga D$. If
		no such $\bZ_2$-grading operator $\Ga$ is given, we say that
		the spectral triple is {\it odd}.
	\end{defn}
	There is a set of axioms for  spectral triples described in \cite{hajac:toknotes,varilly:noncom}. In this article the following axioms are used only.
	\begin{axiom}\label{regularity_axiom}\cite{varilly:noncom} (Regularity). 
		For any $a \in \A$, $[D,a]$ is a bounded operator on~$\H$, and both
		$a$ and $[D,a]$ belong to the domain of smoothness
		$\bigcap_{k=1}^\infty \Dom(\delta^k)$ of the derivation $\delta$
		on~$B(\H)$ given by $\delta(T) \stackrel{\mathrm{def}}{=} [\left|D\right|,T]$.
	\end{axiom}
	\begin{axiom}\label{finiteness_axiom}
		(Finiteness)
		The subspace of smooth vectors
		$\sH^\infty \bydef \bigcap_{k\in\bN} \Dom D^k$ is a \emph{finitely
			generated projective} left  $\sA$-module.
		
		This is equivalent to saying that, for some $N \in \bN$, there is a 
		projector $p = p^2 = p^*$ in~$M_N(\sA)$ such that 
		$\sH^\infty \isom \sA^N p$ as left  $\sA$-modules. 
	\end{axiom}
	\begin{axiom}[Real structure]
		The antiunitary operator $J : \sH \to \sH$ satisfying
		$J^2 = \pm 1$, $JDJ^\dagger = \pm D$, and $J\Ga = \pm \Ga J$ in the even
		case, where the signs depend only on $n \bmod 8$ (and thus are given
		by the table of signs for the standard commutative examples).
		\[
		\begin{array}[t]{|c|cccc|}
			\hline
			n \bmod 8               & 0 & 2 & 4 & 6 \rule[-5pt]{0pt}{17pt} \\
			\hline
			J^2 = \pm 1             & + & - & - & + \rule[-5pt]{0pt}{17pt} \\
			JD = \pm D J & + & + & + & + \rule[-5pt]{0pt}{17pt} \\
			J\Ga = \pm\Ga J         & + & - & + & - \rule[-5pt]{0pt}{17pt} \\
			\hline
		\end{array}
		\qquad\qquad
		\begin{array}[t]{|c|cccc|}
			\hline
			n \bmod 8               & 1 & 3 & 5 & 7 \rule[-5pt]{0pt}{17pt} \\
			\hline
			J^2 = \pm 1             & + & - & - & + \rule[-5pt]{0pt}{17pt} \\
			JD = \pm D J & - & + & - & + \rule[-5pt]{0pt}{17pt} \\
			\hline
		\end{array}
		\]
		Moreover, $b \mapsto J b^* J^\dagger$ is an antirepresentation of $\sA$
		on~$\sH$ (that is, a representation of the opposite algebra
		$\sA^\opp$), which commutes with the given representation of~$\sA$:
		\be\label{st_comn_eqn}
		[a, J b^* J^\dagger] = 0,  \word{for all} a,b \in \sA,
		\ee
		(cf. Definition \ref{df:spt-real_defn}).
	\end{axiom}

	\begin{axiom}\label{fist_order_st_ax}(First order).
		For each $a,b \in \sA$, the following relation holds:
		\begin{equation}\label{fist_order_eqn}
			[[D, a], J b^* J^\dagger] = 0,  \word{for all} a,b \in \sA.
		\end{equation}
		This generalizes, to the noncommutative context, the condition that
		$D$ be a first-order differential operator.
		Since 
		\[
		[[D, a], Jb^*J^\dagger]
		= [[D, Jb^*J^\dagger], a] + [D, \underbrace{[a, Jb^*J^\dagger]}_{=0}],
		\]
		this is equivalent to the condition that 
		\be\label{first_order_dual_eqn}
		[a, [D, Jb^*J^\dagger]] = 0.
		\ee
	\end{axiom}

\begin{axiom}[Orientation]\label{orientation_st_ax}
	There is a Hochschild $n$-cycle
	$$
	\cc = \tsum_j (a_j^0 \ox b_j^0) \ox a_j^1 \oxyox a_j^n 
	\in Z_n(\sA, \sA \ox \sA^\opp),
	$$
	such that
	\begin{equation}\label{eq:vol-cond}
		\pi_D(\cc)
		\equiv \tsum_j a_j^0 (J b_j^{0*} J^{-1}) \,[D,a_j^1] \dots [D,a_j^n]
		= \begin{cases} \Ga, &\text{if $n$ is even}, \\
			1, &\text{if $n$ is odd}. \end{cases}
		\end{equation}
\end{axiom}

	\begin{lem}
		\label{lm:proj-approx}\cite{hajac:toknotes}
		Let $\sA$ be an unital Fr\'echet pre-$C^*$-algebra, whose
		$C^*$-completion is~$A$. If $\tilde{q} = \tilde{q}^2 = \tilde{q}^*$ is
		a projection in $A$, then for any $\eps > 0$, we can find a projection
		$q = q^2 = q^* \in \sA$ such that $\|q - \tilde{q}\| < \eps$.
	\end{lem}
	\begin{theorem}\label{smooth_k_iso_thm}\cite{varilly_bondia}
		If $\sA$ is  a Fr\'echet pre-$C^*$-algebra with $C^*$-completion $A$, then the inclusion $j: \A \to A$ induces an isomorphism $K_0\left( j\right) : K_0\left(\A\right)\to K_0\left(A\right)$.
	\end{theorem}
	\begin{remark}\label{smooth_k_iso_rem}
		The $K_0$-symbol in the above Theorem is the $K_0$-functor of $K$-theory (cf. \cite{bass,blackadar:ko}). The $K_0$-functor is related to projective finitely generated modules. Otherwise any projective finitely generated module corresponds to the idempotent of the matrix algebra. The idea of the proof of the Theorem \ref{smooth_k_iso_thm} contains following ingredients:
		\begin{itemize}
			\item For any idempotent $\tilde e \in \mathbb{M}_n\left(A\right)$ one constructs and idempotent $e\in \mathbb{M}_n\left(\A\right)$ such that there is the isomorphisms $\tilde e A^n\cong e A^n$ of $A$-modules.
			\item The inverse to $K_0j$ homomorphism is roughly speaking given by
			$$
			\left[e\A^n\right]\mapsto \left[ \tilde e A^n\right] 
			$$
			where $\left[\cdot \right]$ means the $K$-theory class  of projective finitely generated module. From the isomorphism $\tilde e A^n\cong e A^n$ the above equation can be replaced with
			\be\label{smooth_k_iso_eqn}
			\left[e\A^n\right]\mapsto \left[  e A^n\right]
			\ee
		\end{itemize}
	\end{remark}
	\section{Representations of  smooth algebras}\label{s_repr}
	
	\paragraph*{}
	Let $(\A, \H, D)$ be a spectral triple. Similarly to \cite{bram:atricle} we  define a representation of $\pi^1:\A \to B(\H^2)$ given by
	\begin{equation}\label{s_diff1_repr_equ}
		\pi^1(a) =  \begin{pmatrix} a & 0\\
			[D,a] & a\end{pmatrix}.
	\end{equation}
	We can inductively construct  representations $\pi^s: \A \to B\left(\H^{2^s}\right)$ for any $s \in \mathbb{N}$. If $\pi^s$ is already constructed then  $\pi^{s+1}: \A \to B\left(\H^{2^{s+1}}\right)$ is given by
	\begin{equation}\label{s_diff_repr_equ}
		\pi^{s+1}(a) =  \begin{pmatrix}  \pi^{s}(a) & 0 \\ \left[D,\pi^s(a)\right] &  \pi^s(a)\end{pmatrix}
	\end{equation}
	where we assume diagonal action of $D$ on $\H^{2^s}$, i.e.
	\begin{equation*}
		D \begin{pmatrix} x_1\\ ... \\ x_{2^s}
		\end{pmatrix}= \begin{pmatrix} D x_1\\ ... \\ D x_{2^s}
		\end{pmatrix}; \ x_1,..., x_{2^s}\in \H.
	\end{equation*}
	For any $s \in \N^0$ there is a seminorm $\left\|\cdot \right\|_s$  on $\A$ given by
	\begin{equation}\label{s_semi_eqn}
		\left\|a \right\|_s = \left\| \pi^{s}(a) \right\|.
	\end{equation}
	The definition of spectral triple requires that $\A$ is a Fr\'echet algebra with respect to seminorms $\left\|\cdot \right\|_s$.

	\section{Noncommutative differential forms}\label{ass_cycle_sec}

	\paragraph*{} 
	\begin{defn}\cite{connes:ncg94}\label{cycle_defn}
		\begin{enumerate}
			\item [(a)]  A \textit{cycle} of dimension $n$ is a triple $\left(\Om, d, \int \right)$ where $\Om = \bigoplus_{j=0}^n\Om^j$  is a graded algebra over $\C$, $d$ is a graded derivation of degree 1 such that $d^2=0$, and $\int :\Om^n \to \C$ is a closed graded trace on $\Om$,
			\item[(b)] Let $\A$  be an algebra over $\C$. Then a \textit{cycle over} $\A$ is given by a cycle $\left(\Om, d, \int \right)$	and a homomorphism $\A \to \Om^0$.
		\end{enumerate}
	\end{defn}
	\begin{empt}\label{connes_cycle_empt}
		We  let $\Om^*\A$ be the reduced universal differential graded algebra over $\A$ (cf. \cite{connes:ncg94} Chapter III.1).
		It is by definition equal to $\A$ in degree 0 and is generated by symbols $da$ ($a \in \A$) of
		degree 1 with the following presentation:
		\begin{enumerate}
			\item[$(\a)$] $d(ab) = (da)b + adb \quad \forall a, b \in \A$,
			\item[$(\bt)$] $d1 = 0$.
		\end{enumerate}
		
		One can check that $\Om^1\A$ is isomorphic as an $\A$-bimodule to the kernel $\ker(m)$ of the
		multiplication mapping $m : \A\otimes\A \to A$, the isomorphism being given by the mapping
		$$
		\sum a_j \otimes b_j \in \ker(m) \mapsto \sum a_j d b_j\in \Om^1\A
		$$
		The involution * of $\A$ extends uniquely to an involution on * with the rule
		$$
		\left(da \right)^* \stackrel{\text{def}}{=} - da^*.
		$$
		The differential $d$ on $A$ is defined unambiguously by
		$$
		d(a_0da_1 ... da_n) =  da_0da_1 ... da_n \quad \forall a_j \in \A,
		$$
		and it satisfies the relations
		\bean
		d^2\om = 0 \quad \forall \om \in \Om^*\A,\\
		d(\om_1\om_2) = (d\om_1)\om_2 + \left( -1\right)^{\partial \om_1} \om_1 d\om_2.
		\eean
		
	\end{empt}
	\begin{proposition}\label{st_cycle_connes_prop}\cite{connes:ncg94}
		Let  $\left( \A, \H, D\right)$  be a spectral triple.
		\begin{enumerate}
			\item 	The following equality defines a *-representation $\pi$ of the reduced universal
			algebra $\Om^*\A$ on $\H$:
			\be\label{diff_repr_eqn}
			\pi(a_0da_1 ... da_n) = a_0[D, a_1] ... [D, a_n] \quad \forall a_j \in \A.
			\ee
			\item Let $J_0 = \ker \pi$  be the graded two-sided ideal of $\Om^*\A$ given by 
			\be\label{junk_grad_eqn}
			J_0^k = \left\{\left.\om \in \Om^k~\right| \pi\left(\om \right)=0 \right\}
			\ee
			then 
			\be\label{junk_cycle_eqn}
			J = J_0 + dJ_0
			\ee
			is a graded differential two-sided ideal of
			$\Om^*\A$.
		\end{enumerate}
	\end{proposition}
	\begin{remark}
		Using 2) of Proposition \ref{st_cycle_connes_prop}, we can now introduce the graded differential algebra
		\be\label{connes_cycle_eqn}
		\Om_D \bydef \Om^*\A/J.
		\ee
	\end{remark}
	
	Thus any spectral triple $\left( \A, \H, D\right)$  naturally defines a cycle $\rho : \A \to \Om_D$ (cf. Definition \ref{cycle_defn}). 
	In particular for any spectral triple there is an $\A$-bimodule $\Om^1_D\subset B\left(\H \right) $ of differential forms which is the $\C$-linear span of operators given by
	\begin{equation}\label{dirac_d_module}
		a\left[D, b \right];\quad a,b \in \A.
	\end{equation}
	There is the differential map
	\begin{equation}\label{diff_map_eqn}
		\begin{split}
			d: \A \to \Om^1_D, \\
			a \mapsto \left[D, a \right].
		\end{split}
	\end{equation}
	
	\begin{definition}\label{ass_cycle_defn}
		We say that that both the cycle $\rho : \A \to \Om_D$ and the differential \eqref{diff_map_eqn} are \textit{associated} with the triple  $\left( \A, \H, D\right)$. We say that  $\A$-bimodule $\Om^1_D$ is the \textit{module of differential forms associated} with the spectral triple  $\left( \A, \H, D\right)$.
	\end{definition}
	
	\subsection{Noncommutative connections and curvatures}
	
	\begin{defn}\label{connection_defn}\cite{connes:ncg94}
		Let $\A\xrightarrow{\rho} \Om$ be a cycle over $\A$, and $\E$ a finite projective module over $\A$.
		Then a \textit{connection} $\nabla$ on $\E$ is a linear map  $\nabla: \E \to \E \otimes_{\A} \Om^1$ such that
		\be\label{conn_prop_eqn}
		\nabla\left(\xi x \right) =  \nabla\left(\xi \right) x =  \xi \otimes d\rho\left(x \right) ; ~ \forall \xi \in \E, ~ \forall x \in \A.
		\ee
		Here $\E$ is a right module over $\A$ and $\Om^1$ is considered as a bimodule over $\A$.
	\end{defn}
	\begin{remark}
		In case of associated cycles (cf. Definition \ref{connection_defn}) the connection equation \eqref{conn_prop_eqn} has the following form
		\be\label{conn_triple_eqn}
		\nabla \left(\xi a \right) = \nabla \left(\xi \right) a + \xi \left[ D,a\right] .
		\ee
	\end{remark}

	\begin{rem}
		The map $\nabla: \E \to \E \otimes_{\A} \Om^1$ is an algebraic analog of the map $\nabla : \Ga\left( E\right) \to \Ga\left( E \otimes T^*\left( M\right) \right)$ given by \eqref{comm_alg_conn}. 
	\end{rem}
	
	\begin{prop}\label{conn_prop}\cite{connes:ncg94}
		Following conditions hold:	
		\begin{enumerate} 
			\item[(a)] 	Let $e \in \End_{\A}\left( \E\right)$ be an idempotent and $\nabla$ is a connection on $\E$; then 
			\be\label{idem_conn}
			\xi \mapsto \left(e \otimes 1 \right) \nabla \xi
			\ee
			is a connection on $e\E$,
			\item[(b)] Any finitely generated projective module $\E$ admits a connection,
			\item[(c)]  The space of connections is an affine space over the vector space $\Hom_{\sA}\left(\E, \E \otimes_{\A} \Om^1 \right)$, 
			\item[(d)] Any connection $\nabla$ extends uniquely up to a linear map  of  $\widetilde{\mathcal E}= \mathcal E \otimes_{\A} \Om$ into itself
			such that
			\be\label{eop_into_eqn}
			\nabla\left(\xi \otimes \om \right) \bydef \nabla\left(\xi \right) \om + \xi \otimes d\om; \quad\forall \xi \in \mathcal E, ~ \om \in \Om. 
			\ee
		\end{enumerate}
	\end{prop}
	\begin{remark}\label{conn_rem}
		The Proposition \ref{conn_prop} implicitly assumes that $\mathcal E\cong  \mathcal E \otimes_{\A} \A\cong \mathcal E \otimes_{\A} \Om^0$ so  there is the natural inclusion
		\be\label{conn_eqn}
		\mathcal E \subset \mathcal E \otimes_{\A} \Om
		\ee
		of right $\A$-modules.
		
	\end{remark}

	\subsection{Connection and curvature}
	\begin{definition}\label{curvature_definition}
		A \textit{curvature} of a connection $\nabla$ is a (right $\A$-linear) map
		\be
		F_\nabla : \mathcal E \to \mathcal E \otimes_{\A} \Om^2
		\ee
		defined as the restriction of $\nabla \circ \nabla$ to $\mathcal E$, that is, $F_{\nabla}\bydef \left.\nabla \circ \nabla \right|_{\mathcal E}$ (cf. \eqref{eop_into_eqn},\eqref{conn_eqn}).
		A connection is said to be \textit{flat} if
		its curvature is identically equal to $0$ (cf. \cite{brzezinsky:flat_co}).
	\end{definition}

	\begin{remark}
		Above algebraic notions of curvature and flat connection are generalizations of corresponding geometrical notions explained in \cite{kobayashi_nomizu:diff_geom} and the Section \ref{geom_flat_subsec}.
	\end{remark}
	\begin{defn}\label{triv_conn_defn}
		
		For any projective $\A$ module $\mathcal E$ there is the unique \textit{trivial connection}
		\bean
		\nabla:  \mathcal E \otimes_{\A} \Om \to \mathcal E \otimes_{\A} \Om, \\
		\nabla = \Id_{\mathcal E} \otimes d.
		\eean
	\end{defn}
	\begin{remark}\label{triv_conn_rem}
		From $d^2 = d \circ d = 0$ it follows that $\left(\Id_{\mathcal E} \otimes d \right) \circ   \left(\Id_{\mathcal E} \otimes d \right)$ = 0, i.e. any trivial connection is flat.
	\end{remark}

	\section{Commutative spectral triples}\label{comm_sp_tr_sec}
	
	\paragraph*{} This section contains citation of \cite{hajac:toknotes,varilly:noncom}.

	\subsection{Spin$^{c}$ manifolds}\label{spin_mani_sec}
	\paragraph*{}
	Let  $M$ be a compact $n$-dimensional orientable Riemannian manifold with a metric~$g$ on  its tangent bundle~$TM$.  For any section  $X \in \Ga\left(M, TM \right)$ of the tangent bundle there is the \textit{derivative}  (cf. \cite{kobayashi_nomizu:diff_geom})
	$$
	X: \Coo\left(M \right) \to C\left( M\right) 
	$$
	We say that $X$ is \textit{smooth} if $X \left(  \Coo\left(M \right)\right) \subset\Coo\left(M \right)$.
	Denote by $\Ga^\infty\left(M, TM \right)$ the space of smooth vector fields. 
	A  \textit{tensor bundle}  corresponds to a tensor products
	$$
	\Ga\left(M, TM \right)\otimes_{C\left( M\right) } ...\otimes_{C\left( M\right) } \Ga\left(M, TM \right)
	$$
	\textit{Smooth sections} of tensor bundle correspond to elements of
	\bean
	\Ga^\infty\left(M, TM \right)\otimes_{C^\infty\left( M\right) } ...\otimes_{C^\infty\left( M\right) } \Ga^\infty\left(M, TM \right)\subset \\ \subset\Ga\left(M, TM \right)\otimes_{C\left( M\right) } ...\otimes_{C\left( M\right) } \Ga\left(M, TM \right).
	\eean
	The metric  tensor $g\bydef \left[g_{jk}\right]$ (cf. Remark \ref{riemann_mani_rem}) is a section of the tensor  bundle, or equivalently $g\in\Ga\left(M, TM \right)\otimes_{C\left( M\right) }\Ga\left(M, TM \right)$,  we suppose that the section is smooth, i.e $g\in\Ga^\infty\left(M, TM \right)\otimes_{\Coo\left( M\right) }\Ga^\infty\left(M, TM \right)$. On the other hand $g$ yields the isomorphism
	$$
	\Ga\left(M, TM \right)\cong \Ga\left(M, T^*M \right)
	$$
	between sections of the tangent and the cotangent bundle. The section of cotangent bundle is said to be \textit{smooth} if it is an image of the smooth section of the tangent bundle. Below we consider tensor bundles which correspond to tensor products of exemplars of $\Ga\left(M, T^*M \right)$ and/or $\Ga\left(M, T^*M \right)$.
	The metric tensor  yields the given by
	\be\label{top_vol_eqn}
	v = \sqrt{\left[ g_{jk}\right] }~dx_1...dx_n
	\ee
	\textit{volume element} (cf. \cite{do_carmo:rg}. The volume element can be regarded as an element of $\Ga\left(M, T^*M \right)\otimes_{C\left( M\right) } ...\otimes_{C\left( M\right) } \Ga\left(M, T^*M \right)$, i.e. $v$ is a section of the tensor field. On the other hand element $v$ defines the unique \textit{Riemannian measure} $\mu$ on $M$.
	If $E \to M$ is a complex bundle then from the Serre-Swan Theorem \ref{serre_swan_thm} it turns out that there is an idempotent $e \in \mathbb{M}_n\left(C\left(M \right)  \right)$ such that there is the  $C(M)$-module isomorphism
	$$
	\Ga\left(M, E \right) \cong e C\left(M \right)^n.
	$$
	From the Theorem \ref{smooth_k_iso_thm} one can suppose that $e \in \mathbb{M}_n\left(\Coo\left(M \right)  \right)$ (cf. Remark \ref{smooth_k_iso_rem}).
	\begin{definition}\label{top_smooth_m_defn}
		The isomorphic to $e \Coo\left(M \right)^n$ group  with the induced by
		\be\label{top_smooth_m_eqn}
		\Ga^\infty\left(M, E \right) \cong e \Coo\left(M \right)^n \subset e C\left(M \right)^n\cong \Ga\left(M, E \right)\quad e \in \mathbb{M}_n\left(\Coo\left(M \right)  \right)
		\ee
		inclusion $\Ga^\infty\left(M, E \right) \subset \Ga\left(M, E \right)$ is said to be the subgroup of \textit{smooth sections}.
	\end{definition} 
	It is clear that 
	\be
	\Coo\left(M \right)\Ga^\infty\left(M, E \right)\subset \Ga^\infty\left(M, E \right)
	\ee
	and the map
	$$
	\Ga^\infty\left(M, E \right)\mapsto \Ga\left(M, E \right)
	$$
	yields the given by the Theorem \ref{smooth_k_iso_thm} isomorphism $K_0\left(\Coo\left(M \right) \right) \cong K_0\left( C\left(M \right)\right) $. 
	It can be proved that   the Definition \ref{top_smooth_m_defn} complies with the above definition of smooth tensor fields.
	We build a Clifford algebra bundle
	$\C\ell(M) \to M$ whose fibres are full matrix algebras (over~$\C$) as
	follows. If $n$ is even, $n = 2m$, then
	$\C\ell_x(M) \bydef \Cl(T_xM,g_x) \ox_\R \C \simeq M_{2^m}(\C)$ is the
	complexified Clifford algebra over the tangent space $T_xM$. If $n$ is
	odd, $n = 2m + 1$, the analogous fibre splits as
	$\mathbb{M}_{2^m}(\C) \oplus \mathbb{M}_{2^m}(\C)$, so we take only the \textit{even}
	part of the Clifford algebra:
	$\C\ell_x(M) \bydef \Cl^\mathrm{even}(T_xM) \ox_\R \C \simeq \mathbb{M}_{2^m}(\C)$. 
	What we gain is that in all cases, the bundle
	\be\label{st_cliffod_eqn}
	 \C\ell(M) \to M
	 \ee
	is a locally trivial field of (finite-dimensional) elementary
	$C^*$-algebras, so $B = \Ga\left( M, \C\ell(M)\right)$ is a $C^*$-algebra 
	$S_x$ such that 
\be\label{st_spinor_eqn}
	\C\ell_x(M) \simeq \End(S_x);
	\ee the
	class~$\delta(\C\ell(M))$ is precisely the obstruction to patching them
	together (there is no obstruction to the existence of the algebra
	bundle $\C\ell(M)$). It was shown by Plymen~\cite{plymen:mor_s} that
	$\delta(\C\ell(M)) = W_3(TM)\in  H^3(M,\Z)$, the integral class that is the
	obstruction to the existence of a \textit{Spin$^c$ structure} in the
	conventional sense of a lifting of the structure group of~$TM$ {}from
	$SO(n)$ to $\Spin^c(n)$: see \cite{lawson_m} for more
	information on~$W_3(TM)$.
	
	Thus $M$ admits Spin$^c$ structures if and only if
	$\delta(\C\ell(M)) = 0$. But in the Dixmier--Douady theory,
	$\delta(\C\ell(M))$ is the obstruction to constructing (within the
	$C^*$-category) a $B$-$A$-bimodule $\SS$ that implements a Morita
	equivalence between $A = C(M)$ and $B = C(M,\C\ell(M))$. Let us
	paraphrase Plymen's redefinition of a Spin$^c$ structure, in the
	spirit of noncommutative geometry:
	\begin{definition}\label{spin_str_defn}
		Let $M$ be a Riemannian manifold, $A = C(M)$ and\\
		$B = C(M,\C\ell(M))$. We say that the tangent bundle $TM$
		\textit{admits a }Spin$^c${ structure} if and only if it is orientable
		and $\delta(\C\ell(M)) = 0$. In that case, a \textit{Spin$^{\bf c}$
			structure} on~$TM$ is a pair $(\eps,S)$ where $\eps$ is an
		orientation on~$TM$ and $S$ is a $B$-$A$-equivalence bimodule.
	\end{definition}
	
	
	What is this equivalence bimodule~$\SS$? By the Serre-Swan theorem  \ref{serre_swan_thm}, it
	is of the form $\Ga(M, S)$ for some complex vector bundle $S \to M$ that also
	carries an irreducible left  action of the Clifford algebra bundle
	$\C\ell(M)$. This is the \textit{spinor bundle} whose existence displays
	the Spin$^c$ structure in the conventional picture. We call
	$\Ga^\infty(M,S)$ the \textit{spinor module}; it is an irreducible
	Clifford module in the terminology of~\cite{atiyah_b_s}, and has
	rank~$2^m$ over $C(M)$ if $n = 2m$ or $2m + 1$.
	\begin{remark}
		If $\B \bydef \Ga^\infty \left(M, \C\ell(M) \right)$ then $\B$ is an be an unital Fr\'echet pre-$C^*$-algebra, such that $\B$ is dense in $B$. Similarly if $\A\bydef \Coo\left(M \right)$ then $\A$ is dense in  $A\bydef C\left(M \right)$. The Morita equivalence between $B$ and $A$ is given by the projective $B$-$A$ bimodule $\Ga\left(M,S \right)$. Similarly  the Morita equivalence between $\B$ and $\A$ is given by the projective $\B$-$\A$ bimodule $\Ga^\infty\left(M,S \right)$.
	\end{remark}
\begin{empt}
	Let $\sS^\sharp\bydef \Hom_{C\left(M\right)}\left(\Ga\left(M,S \right), C\left(M\right) \right)$
	\end{empt}

\begin{prop}\label{pr:charge-conj_prop}\cite{hajac:toknotes}
If $\sS^\sharp\bydef \Hom_{C\left(M\right)}\left(\Ga\left(M,S \right), C\left(M\right) \right)$ then 	there is a $B$-$A$-bimodule isomorphism $\sS^\sharp \isom \sS$ if and only
	if there is an \textit{antilinear} endomorphism $C$ of $\sS$ such that
	\begin{enumerate}
		\item[(a)]
		$C(\psi\,a) = C(\psi) \,\bar a$\quad
		for $\psi \in \sS$, $a \in A$;
		\item[(b)]
		$C(b\,\psi) = \chi(\bar b)\, C(\psi)$\quad
		for $\psi \in \sS$, $b \in B$;
		\item[(c)]
		$C$ is \emph{antiunitary} in the sense that
		$\pairing{C\phi}{C\psi} = \pairing{\psi}{\phi} \in A$,\quad
		for $\phi, \psi \in \sS$;
		\item[(d)]
		$C^2 = \pm 1$ on $\sS$ whenever $M$ is connected.
	\end{enumerate}
\end{prop}
\begin{remark}\label{pr:charge-conj_rem}
	Operator $C$ from the Proposition \ref{pr:charge-conj_prop} is a bijective isomorphism $\Ga\left(M,S \right)\cong \Ga\left(M,S \right)$ of Abelian groups. On the other hand since $\Ga\left(M,S \right)$ is dense in $L^2\left(M, S, \mu \right)$ it can be extended up to antiunitary  operator $J:L^2\left(M, S, \mu \right)\cong L^2\left(M, S, \mu \right)$. According to the construction of commutative real spectral triples $J$ is an operator which defines the structure of a real spectral triple (cf. Definition \ref{df:spt-real_defn}).
\end{remark}

\begin{remark}\label{top_spin_prod_rem}\cite{varilly:noncom}
	Any spinor bundle $S$ has a sesquilinear form $S \times_\sX S \to \C$ (cf. Definition \ref{top_herm_bundle_form_defn}). In \cite{varilly:noncom} a bundle with  sesquilinear form is said to be a \textit{Hermitian bundle}.
\end{remark}


	\subsection{The Dirac operator}\label{comm_dirac_sec}
	\paragraph{}
	As soon as a spinor module makes its appearance, one can introduce the
	\textit{Dirac operator}. Let $\mu$ be the Riemannian measure given by the volume element (cf. equation \eqref{top_vol_eqn}). If  $\H \bydef L^2(M,S, \mu)$ is the space $\H \bydef L^2(M,S, \mu)$
	of square-integrable spinors then $\Ga^\infty\left(M, \SS \right)\subset \H$.
	This is a selfadjoint first-order
	differential operator $\Dslash$ defined on the space $\H$
	of square-integrable spinors, whose domain includes the space of smooth spinors
	$\SS = \Ga^\infty\left(M, S \right)$.  The Riemannian metric $g = [g_{ij}]$
	defines isomorphisms $T_xM \simeq T^*_xM$ and induces a metric
	$g^{-1} = [g^{ij}]$ on the cotangent bundle $T^*M$. Via this
	isomorphism, we can redefine the Clifford algebra as the bundle with
	fibres $\C\ell_x(M) \bydef \Cl(T^*_xM, g_x^{-1}) \ox_\R \C$ (replacing $\Cl$
	by $\Cl^\mathrm{even}$ when $\dim M$ is odd). Let $\Ga(M, T^*M)$ be
	the $C\left( M\right) $-module of \hbox{$1$-forms} on~$M$. 
	The spinor module $\SS$ is
	then a $B$-$A$-bimodule on which the algebra $B = \Ga\left( M,\C\ell\left( M\right) \right) $
	acts irreducibly. If 
\be\label{spin_gamma_eqn}
\begin{split}
\ga: B \cong \End_A\left( \SS\right) \\
\text{or, equivalently } \ga: \Ga(M,\C\ell(M)) \cong \End_{C\left(M \right) }\left( \Ga\left(M,S \right) \right)
\end{split}
\ee	
is	 the natural isomorphism then $\ga$ obeys the anticommutation rule
	$$
	\{\ga(\a), \ga(\b)\} = -2 g^{-1}(\a,\b) = -2 g^{ij} \a_i \b_j\in C\left( M\right) 
	\sepword{for}  \a,\b \in \Ga(M, T^*M).
	$$
	The action~$\ga$ of $\Ga(M, \C\ell(M))$ on the Hilbert-space completion
	$\H$ of~$\SS$ is called the \textit{spin representation}.
	
	The metric $g^{-1}$ on $T^*M$ gives rise to a canonical
	\textit{Levi-Civita connection}\\
	$\nabla^g \: \Ga^\infty(M, T^*M) \to \Ga^\infty(M, T^*M) \otimes_\A \Ga^\infty(M, T^*M)$ that, as well as
	obeying the Leibniz rule
	$$
	\nabla^g(\om a) = \nabla^g(\om)\,a + \om \ox da,
	$$
	preserves the metric and is torsion-free. The \textit{spin connection}
	is then a linear operator
	\be\label{spin_conn_defn_eqn}
\begin{split}
	\nabla^S \: \Ga^\infty(M,S) \to \Ga^\infty(M, S) \otimes_{\Coo\left(M \right) } \Ga^\infty(M, TM)
\end{split}
\ee
	satisfying two Leibniz rules, one for the right action of~$\A$ and the
	other, involving the Levi-Civita connection, for the left  action of
	the Clifford algebra:
	\be\label{spin_conn_eqn}
	\begin{split}
		\nabla^S(\psi a) \bydef \nabla^S(\psi)a + \psi \otimes  da,\\
		\nabla^S(\ga(\om)\psi)= \ga(\nabla^g\om)\psi + \ga(\om)\nabla^S\psi,
	\end{split}
	\ee
	$$
	$$
	for $a \in \A$, $\om \in \Ga^\infty(M, T^*M)$, $\psi \in \Ga^\infty(M, S)$.
	
	Once the spin connection is found, we define the Dirac operator as the
	composition $\ga \circ \nabla^S$; more precisely, the local expression
	\be\label{comm_dirac_eqn}
	\Dslash \bydef \sum_{j=1}^n\ga(dx^j) \, \nabla^S_{\del/\del x^j}
	\ee
	is independent of the local coordinates and defines $\Dslash$ on the domain
	$\SS \subset \H$.
	One can check that this operator is symmetric; it extends to an 	unbounded selfadjoint operator on~$\H$, also called $\Dslash$.
	In result one has the commutative spectral triple
	\be\label{comm_sp_tr_eqn}
	\left(\Coo\left( M\right), L^2\left(M, S \right), \Dslash , J  \right).
	\ee
	It is shown in \cite{varilly:noncom,krajewski:finite} that the first order condition (cf. Axiom \ref{fist_order_st_ax}) is equivalent to the following equation
	\begin{equation}\label{comm_matr_x_eqn}
		\begin{split}
			\left[ {\slashed D}, {a}\right]b= b\left[ {\slashed D}, {a}\right] \quad\forall a, b \in \Coo\left( M\right).
		\end{split}
	\end{equation}


	\section{Finite spectral triples}\label{ctr_fin_sp_tr_sec}
	
	\paragraph*{}
	Here I follow to \cite{krajewski:finite}.
	Finite spectral triples are particular cases of spectral triples of dimension 0. The latter are rigorously defined within the axioms of noncommutative geometry and yield a general theory of discrete spaces. Among all discrete spaces, we focus on finite ones, thus the algebra is finite dimensional.  Futhermore, we will also assume that the Hilbert space is finite dimensional, an infinite dimensional one corresponding to a theory with an infinite number of elementary fermions. Accordingly, a finite spectral triple $\left(\A, \H, D , J\right) $ is defined as a spectral triple of dimension 0 such that both $\A$ and $\H$ are finite dimensional. Using such a triple, it is possible to a construct Yang-Mills theory with spontaneous symmetry breaking whose gauge group is the group of unitary elements of  $\A$, $\H$ is the fermionic Hilbert space and $D$ is the mass matrix.
	\par
	It is  known that finite dimensional real involutive algebras which admit a faithful representation on a finite dimensional Hilbert space are just direct sums of matrix algeras over the fields of real numbers, complex numbers and quaternions. Therefore, we write the algebra as a direct sum
	\be\label{fin_oa_eqn}
	\A=\bigoplus_{j=1}^{N}\mathbb{M}_{n_{j}}\left( \mathbb{K}\right),
	\ee
	where $\mathbb{M}_{n}\left( \mathbb{K}\right)$ denotes the algebra of square matrices of order $n$ with entries in the field $\mathbb{K}=\R$, $\C$ or $\mathbb{H}$.  This remark simplifies considerably finite noncommutative geometry and it becomes possible to give a detailed account of all finite spectral triples.
	We suppose that there is a faithful representation 
	$$
	\A\hookto \mathbb{M}_N\left(\C\right).
	$$
	The Dirac operator $D$ is represented by the mass matrix $M\in \mathbb{M}_{N}\left( \C\right)$, i.e.
	\be\label{ctr_fin_dirac_eqn}
	D \xi = M \xi;~~\xi \in \C^N.
	\ee
	The charge conjugation operator $J$ is described in \cite{krajewski:finite}.
	\begin{remark}
		Any  finite spectral triple $\left(\A, \H, D, J \right)$ does not satisfies to all described in \ref{sp_tr_defn_sec} conditions because the closure of  $\A$ in $\mathbb{M}_N\left(\C\right)$ is not a $C^*$-algebra in general. However this circumstance does not contradict with constructions of this book.
	\end{remark}

	\section{Product of spectral triples}\label{sp_tr_prod_sec}

	\paragraph*{}
	Here I follow to \cite{dabrowski:product}. Let both $\left(\A_1, \H_1, D_1, J_1\right)$ and  $\left(\A_2, \H_2, D_2, J_2\right)$ be spectral triples, we would like to find their direct product $\left(\A, \H, D, J\right) =\left(\A_1, \H_1, D_1, J_1\right)\times\left(\A_2, \H_2, D_2, J_2\right)$. It is known that there are even and odd spectral triples, so
	there are following cases:
	
	\begin{enumerate}
		\item[(i)]\textit{Even-even case}. In this case $\A \stackrel{\text{def}}{=} \A_1 \otimes \A_2$ and  $\H = \H_1 \otimes \H_2$. Moreover if $\rho_1: \A_1 \to B\left(  \H_1\right)$ and $\rho_2: \A_2 \to B\left(  \H_2\right)$ are representations of spectral triples then the representation of the product is given by $\rho \stackrel{\text{def}}{=} \rho_1 \otimes \rho_1: \A_1 \otimes \A_2 \to B\left( \H_1 \otimes \H_2 \right)$. The Dirac operator is given by
		\begin{equation}
			\begin{split}
				D\stackrel{\text{def}}{=} D_1\otimes \id_{\mathcal{H}_2} + \chi_1\otimes D_2\;,\\
				D'\stackrel{\text{def}}{=}D_1\otimes \chi_2 + \id_{\mathcal{H}_1}\otimes D_2\, ,
			\end{split}
		\end{equation}
		where $\chi_1$ (resp. $\chi_2$) is the chirality operator of $\left(\A_1, \H_1, D_1\right)$ (resp.  $\left(\A_2, \H_2, D_2\right)$) (cf.  \cite{dabrowski:product})
		Operators $D$ and $D'$ are unitary equivalent. A  reality structure is given by $J \bydef J_1\otimes J_2$.
		\item[(ii)]\textit{Even-odd case}. The algebra, the *-representation, Dirac operator and reality structure are the same as in the even-even case.
		\item[(iii)] \textit{Odd-odd case}. This case is not considered it this work, so it is not described here.
	\end{enumerate}

\section{Morita equivalence and Hermitian connections}

We start with any geometry $(\A,\H,D,\Ga,J)$ and a finite projective
right $\A$-module~$\E$. Using the representation $\pi\: \A \to B(\H)$
and the antirepresentation $\pi^0 \: b \mapsto J\pi(b^*)J^\dagger$, we can
regard the space $\H$ as an $\A$-bimodule. This allows us to introduce
the vector space
\be\label{morita_prod_eqn}
\widetilde\H \bydef \E \otimes_A \H \otimes_A \overline\E.
\ee
If $\E = p\A^m$, then $\overline\E = \overline\A^m\,p$ and
$\widetilde\H = \pi(p) \pi^0(p) [\H \ox \C^{m^2}]$, so that $\widetilde\H$
becomes a Hilbert space under the scalar product
\be\label{morita_sc_prod_eqn}
\<r \ox \eta \ox \bar q, s \ox \xi \ox \bar t\,>
\bydef \<\eta, \pi(r,s)\, \pi^0(t,q) \,\xi>.
\ee
If $\H = \H^{\mathrm{op}} \H^-$ is $\Z_2$-graded, there is an obvious
$\Z_2$-grading of $\widetilde\H$.

The antilinear correspondence $s \mapsto \bar s$ between $\E$
and~$\overline\E$ also gives an obvious way to extend $J$ to~$\widetilde\H$:
\be\label{morita_J_eqn}
\widetilde J(s \ox \xi \ox \bar t\,) \bydef t \ox J\xi \ox \bar s.
\ee
Let $\B \bydef \End_\A \E$, and recall that $\E$ is a left  $\B$-module.
Then
\be\label{st_left_eqn}
\rho(b) : s \ox \xi \ox \bar t \longmapsto b\,s \ox \xi \ox \bar t
\ee
yields a representation $\rho$ of $\B$ on~$\widetilde\H$, satisfying
\be\label{st_right_eqn}
\rho^0(b) \bydef \widetilde J \rho(b^*) \widetilde J^\7
: s \ox \xi \ox \bar t \longmapsto s \ox \xi \ox \bar t\,b,
\ee
where $\bar t\,b \bydef \overline{b^*t}$, of course. 

\begin{remark}\label{sp_tr_commute_rem}\cite{varilly:noncom}
The action $\rho$,
$\rho^0$ of $\B$ on~$\widetilde\H$ obviously commute.
\end{remark}

\subsubsection{Where connections come from}
The nontrivial part of the construction of the new geometries
$(\B,\widetilde\H,\widetilde D,\widetilde\Ga,\widetilde J)$ is the determination of an
appropriate operator $\widetilde D$ on~$\widetilde\H$. Guided by the
differential properties of Dirac operators, the most suitable
procedure is to postulate a \textit{Leibniz rule}:
\be\label{morita_d_eqn}
\widetilde D(s \ox \xi \ox \bar t\,) \bydef (\nabla s)\xi \ox \bar t
+ s \ox D\xi \ox \bar t + s \ox \xi(\overline{\nabla t}),
\ee
where $\nabla s$, $\nabla t$ belong to some space whose elements can
be represented on~$\H$ by suitable extensions of~$\pi$ and $\pi^0$.

Consistency of~\eqref{morita_d_eqn} with the actions of $\A$ on~$\E$ and~$\H$
demands that $\nabla$ itself comply with a Leibniz rule. Indeed, since
$$
sa \ox \xi \ox \bar t = s \ox a\xi \ox \bar t
\sepword{for all} a \in \A,
$$
we get from \eqref{morita_d_eqn}
$$
\nabla(sa)\xi \ox \bar t  + s \ox a\,D\xi \ox \bar t
= (\nabla s)a\xi \ox \bar t  + s \ox D\,a\xi \ox \bar t,
$$
so we infer that
\be\label{morita_nabla_eqn}
\nabla(sa) = (\nabla s)a + [D,a],
\ee
or more pedantically, $\nabla(sa) = (\nabla s)\pi(a) + [D,\pi(a)]$ as
operators on~$\H$.

To satisfy these requirements, we introduce the space of bounded
operators
\be\label{om_d_eqn}
\Om_D^1 \bydef \text{ a linear span of }\set{a\,[D,b] : a,b \in \A} \subseteq B(\H),
\ee
which is evidently an $\A$-bimodule, the right action of~$\A$ being
given by $a\,[D,b]\.c \bydef a\,[D,bc] - ab\,[D,c]$. The notation is
chosen to remind us of differential $1$-forms; indeed, for the
commutative geometry $(\Coo(M), L^2(M,S), \Dslash, \chi, J)$, we get
$$
\Om_{\Dslash}^1 = \set{\ga(\a) : \a \in \A^1(M)},
$$
i.e., conventional $1$-forms on~$M$, represented on spinor space as
(Clifford) multiplication operators.

\begin{definition}
We can now form the right $\A$-module $\E \otimes_A \Om_D^1$. A
\textbf{connection} on~$\E$ is a linear mapping
$$
\nabla \: \E \to \E \otimes_A \Om_D^1
$$
that satisfies the Leibniz rule \eqref{morita_nabla_eqn}.

It is worth mentioning that only \textit{projective} modules admit
connections (cf. \cite{connes:ncg94}). In the present case, if we define linear
maps
$$
0 \xrightarrow \E \otimes_A \Om_D^1 \xrightarrow{j}  \E \ox_\C \A \xrightarrow{m} \E \xrightarrow 0
$$
by $j(s\,[D,a]) \bydef sa \ox 1 - s \ox a$ and $m(s \ox a) \bydef sa$, we get
a short exact sequence of right $\A$-modules (think of $\E \ox_\C \A$
as a free $\A$-module generated by a vector-space basis of~$\E$). Any
linear map $\nabla \: \E \to \E \otimes_A \Om_D^1$ gives a linear section
of~$m$ by $f(s) \bydef s \ox 1 - j(\nabla s)$. Then
$f(sa) - f(s)a = j(s\,[D,a] - \nabla(sa) + (\nabla s)a)$, so $f$ is an
$\A$-module map precisely when $\nabla$ satisfies the Leibniz
rule~\eqref{morita_nabla_eqn}. If that happens, $f$ splits the exact sequence and
embeds $\E$ as a direct summand of the free $\A$-module
$\E \ox_\C \A$, so $\E$ is projective.
\end{definition}

\subsubsection{Hermitian connections}
The operator $\widetilde D$ must be selfadjoint on~$\widetilde\H$. If
$\xi,\eta \in \Dom(D)$, we get
\bean
\left\langle r \otimes \eta \otimes  \overline q\left| \widetilde D\left(  s \otimes \xi \otimes  \overline t \right) \right. \right\rangle =  \left\langle \eta\left|\pi_D\left(\left(r, \nabla s \right)  \right)\pi^0\left(\left( t, q\right)\right) \xi    \right.\right\rangle+\\+ \left\langle \eta\left|\pi\left(\left(r,  s \right)  \right)\pi^0\left(\left( t, q\right)\right) D \xi    \right.\right\rangle
+ \left\langle \eta\left|\pi\left(\left(r,  s \right)  \right)\pi^0_D\left(\left(\nabla  t, q\right)\right)  \xi    \right.\right\rangle,\\
\\
\left\langle\left.\widetilde{D} \left( r \otimes \eta \otimes  \overline q\right) \right|   s \otimes \xi \otimes  \overline t  \right\rangle =  \left\langle \eta\left|\pi_D\left(\left(\nabla r,  s \right)  \right)\pi^0\left(\left( t, q\right)\right) \xi    \right.\right\rangle+\\+ \left\langle \eta\left|D\pi\left(\left(r,  s \right)  \right)\pi^0\left(\left( t, q\right)\right) \xi    \right.\right\rangle
+ \left\langle \eta\left|\pi\left(\left(r,  s \right)  \right)\pi^0_D\left(\left(  t, \nabla q\right)\right) \xi    \right.\right\rangle.
\eean
This reduces to the condition that
\be\label{herm_conn_eqn}
(r,\nabla s) - (\nabla r,s) = [D, (r,s)]
\sepword{for all} r,s \in \E.
\ee
where the order one condition ensures commutation of $\pi_D(\Om_D^1)$
with $\pi^0(\A)$.
\begin{definition}\label{herm_conn_defn}\cite{varilly:noncom}.
We call the connection $\nabla$ \textit{Hermitian} (with respect
to~$D$) if \eqref{herm_conn_eqn} holds.
\end{definition}

To sum up: two \textit{geometries} $(\A,\H,D,\Ga,J)$ and
$(\B,\widetilde\H,\widetilde D,\widetilde\Ga,\widetilde J)$ are \textit{Morita-equivalent}
if there exist a finite projective right $\A$-module $\E$ and an
$\Om_D^1$-valued Hermitian connection $\nabla$ on~$\E$, such that:
$\B = \End_\A \E$, $\widetilde\H$ and $\widetilde\Ga$ are given by \eqref{morita_prod_eqn},
$\widetilde J$ by \eqref{morita_J_eqn}, and $\widetilde D$ by \eqref{morita_d_eqn}.

\chapter{Noncommutative Parallel Transport}\label{sec:fuzzy}
\section{Parallel transports}
\paragraph*{}
Definitions of the Appendix \ref{comm_prot} cannot be directly used in the noncommutative case because the noncommutative geometry does contain closed paths. However paths can be replaced with module parallel transports.
\begin{defn}\label{defi:transport}\cite{schenkel:nc_parallel}
	Let $A$ be an associative and unital algebra and $\mathcal{E}$ a right $A$-module.
	\begin{itemize}
		\item[1.)]
		A {\it one-parameter group of automorphisms} of $A$ is a map
		$\varphi: \mathbb{R} \times A \to A\,,~(\tau,a) \mapsto \varphi(\tau,a) =\varphi_\tau(a)$, such that
		\begin{itemize}
			\item[(i)] $\varphi_\tau(a\,b) = \varphi_\tau(a) \,\varphi_\tau(b)$, for all $\tau\in\mathbb{R}$ and $a,b\in A$
			\item[(ii)] $\varphi_0 =\mathrm{id}_A$
			\item[(iii)] $\varphi_{\tau+\sigma} = \varphi_\tau\circ\varphi_\sigma$, for all $\tau,\sigma\in\mathbb{R}$
		\end{itemize}
		\item[2.)] Let $\varphi:\mathbb{R} \times A \to A $ be a one-parameter group of automorphisms of $A$. A {\it module parallel transport} on
		$\mathcal{E}$ along $\varphi$ is a map $\Phi: \mathbb{R} \times \mathcal{E} \to \mathcal{E}\,,~(\tau,s)\mapsto \Phi(\tau,s)=\Phi_\tau(s)$, such that
		\begin{itemize}
			\item[(i)] $\Phi_\tau(s\,a) = \Phi_\tau(s)\,\varphi_\tau(a)$, for all $\tau\in\mathbb{R}$, $s\in\mathcal{E}$ and $a\in A$
			\item[(ii)] $\Phi_0=\mathrm{id}_\mathcal{E}$
			\item[(iii)] $\Phi_{\tau+\sigma}= \Phi_\tau\circ \Phi_\sigma$, for all $\tau,\sigma\in\mathbb{R}$.
		\end{itemize}
	\end{itemize}
	\paragraph{} If $A$ and $\mathcal{E}$ are equipped with a smooth structure, the maps $\varphi$ and $\Phi$ are required to be smooth. Denote by $\mathrm{Trans}_{\mathcal{E}}$ the set of module parallel transports.
\end{defn}

\section{Parallel transform in noncommutative differential geometry}\label{par_trans_conn_sect}
\paragraph{} In this section I follow to \cite{connes:c_alg_dg}.
\begin{empt}\label{c_din_empt}\cite{connes:c_alg_dg}
	Let $(A,G,\alpha)$ be a $C^*$-\textit{dynamical system}, where $G$ is a Lie group.
	We shall say that $x \in A$ is of $C^{\infty}$ class if and only if
	the map $g \mapsto \alpha_g (x)$ from $G$ to the normed space $A$ is in
	$C^{\infty}$. The involutive algebra 
\be\label{par_inf_eqn}	
	A^{\infty} \bydef \{ x \in A ,
	\, x \ \hbox{of class} \ C^{\infty} \}
\ee	
	 is norm dense in $A$.
	
	Let $\mathcal{E}^{\infty}$ be a finite projective module on $A^{\infty}$,
	(we shall write it as a right module); $\mathcal{E} = \mathcal{E}^{\infty}
	\otimes_{A^{\infty}}$ $A$ is then  a finite projective module on
	$A$.
\end{empt}

\begin{lem}\cite{connes:c_alg_dg}
	For every finite projective module $\mathcal{E}$ on $A$, there exists a
	finite projective module $\mathcal{E}^{\infty}$ on $A^{\infty}$, unique up to
	isomorphism, such that $\mathcal{E}$ is isomorphic to $\mathcal{E}^{\infty}
	\otimes_{A^{\infty}} \, A$.
\end{lem}
\begin{empt}
	In the sequel we let $\E^{\infty}$ be a finite projective module 
	on $A^{\infty}$. An Hermitian  structure on 
	$\E^{\infty}$ is given by a positive Hermitian  form  $\langle 
	\xi , \eta \rangle \in A^{\infty}$, $ \, \xi , \eta \in \E^{\infty}$ 
	such that
	$$
	\langle \xi \cdot x , \eta \cdot y \rangle = y^* \langle \xi , \eta 
	\rangle \, x \ , \qquad  \, \xi , \eta \in \E^{\infty} \, , \quad \, 
	x,y \in A^{\infty} \, .
	$$
	For $n \in \N$, $\E^{\infty} \otimes \C^n$ is a finite projective module
	on $\mathbb{M}_n (A^{\infty}) = A^{\infty} \otimes \mathbb{M}_n(\C)$, this allows, 
	replacing $A$ by $\mathbb{M}_n  (A) = A \otimes \mathbb{M}_n (\C)$ (and the 
	$G$-action by $\a \otimes {\rm id}$) to assume the existence of a  
	selfadjoint idempotent $e \in A^{\infty}$ and of an isomorphism $F$ with the module 
	$e \, A^{\infty}$ on $\E^{\infty}$. We then endow $\E^{\infty}$ with
	the following Hermitian  structure~:
	$$
	\langle \xi , \eta \rangle = F^{-1} (\eta)^* \, F^{-1} (\xi) \in 
	A^{\infty} \, .
	$$
	Let $\delta$ be the representation of Lie $G$ in the Lie-algebra of derivations of
	$A^{\infty}$ given by
	$$
	\delta_X (x) = \lim_{t \rightarrow 0} \ \frac{1}{t} \, (\alpha_{g_t} (x) - x) \, ,
	\qquad \hbox{where} \quad \dot{g}_0 = X \, , \ x \in A^{\infty} \, .
	$$
\end{empt}
\medskip

\begin{defn}\label{connes_dg_connection_defn}\cite{connes:c_alg_dg} $\mathcal{E}^{\infty}$ be a finite projective module on $A^{\infty}$,
	a $G$-{\it connection} (on $\mathcal{E}^{\infty}$) is a linear map
	$\nabla : \mathcal{E}^{\infty} \to \mathcal{E}^{\infty} \otimes(\mathrm{Lie} ~ G)^*$ such that,
	for all $X \in \mathrm{Lie} ~ G$ and $\xi \in \mathcal{E}^{\infty}$, $x \in
	A^{\infty}$ one has
	$$
	\nabla_X (\xi \cdot x) = \nabla_X (\xi) \cdot x + \xi \cdot \delta_X (x) \, .
	$$
\end{defn}
\begin{remark}
We use the $G$-{\it connection} word instead \textit{connection} (cf. \cite{connes:c_alg_dg}) to distinguish Definitions \ref{connection_defn} and \ref{connes_dg_connection_defn}.
\end{remark}


\medskip

We shall say that $\nb$ is compatible with the Hermitian  structure 
iff~:
$$
\langle \nb_X \, \xi , \xi' \rangle + \langle \xi , \nb_X \, \xi' 
\rangle = \delta_X \langle \xi , \xi' \rangle \, , \qquad ~ \, \xi , 
\xi' \in \mathcal{E}^{\infty} \, , \quad ~ \, X \in \hbox{Lie} \ G \, .
$$
Every finite projective module, $\mathcal{E}^{\infty}$, on $A^{\infty}$ 
admits a connection~; on the module $e \, A^{\infty}$ the following 
formula defines the  {\it Grassmannian connection }
$$
\nb_X^0 (\xi) = e \, \delta_X (\xi) \in e \, A^{\infty} \, , \qquad ~ \, \xi 
\in e \, A^{\infty} \, , \quad ~ \, X \in \hbox{Lie} \ G \, .
$$
This connection is compatible with the Hermitian  structure
$$
\langle \xi , \eta \rangle = \eta^* \, \xi \in A^{\infty} \, , \qquad 
~ \, \xi , \eta \in e \, A^{\infty} \, .
$$

To the representation $\d$ of Lie $G$ in the Lie-algebra of derivations of
$A^{\infty}$ corresponds the complex $\Om = A^{\infty} \otimes 
\Lambda \ (\hbox{Lie} \ G)^*$ of left-invariant differential forms on $G$ 
with coefficients in $A^{\infty}$. We endow $\Om$ with the algebra 
structure given by the tensor product of $A^{\infty}$ by the exterior algebra of 
$(\hbox{Lie} \ G)_{\C}^*$
(we use the notation $\om_1 \wedge \om_2$ for the product 
of $\om_1$ with $\om_2$, one no longer has, of course, 
the equality
$\om_2 \wedge \om_1 = (-1)^{\partial \om_1 \, \partial \om_2} \, \om_1 
\wedge \om_2$).
The exterior differential
$d$ 
is such that~:
\begin{itemize}
	\item[1$^{\rm o}$] for $a \in A^{\infty}$ and $X \in \hbox{Lie} \ G$ one 
	has $\langle X , da \rangle = \d_X (a)$~;
	\item[2$^{\rm o}$] $d (\om_1 \wedge \om_2) = d \om_1 \wedge \om_2 + 
	(-1)^p \, \om_1 \wedge d \om_2$, $~ \, \om_1 \in \Om^p$, $~ \, \om_2 
	\in \Om$~;
	\item[3$^{\rm o}$] $d^2 \, \om = 0$, $~ \, \om \in \Om$.
\end{itemize}
As $A^{\infty} \subset \Om$, $\Om$ is a bimodule on $A^{\infty}$.

Every connection on $e \, A^{\infty}$ is of the form $\nb_X (\xi) = 
\nb_X^0 (\xi) + \th_X \, \xi$, $~ \, \xi \in e \, A^{\infty}$, $X \in 
\hbox{Lie} \ G$, where the form $\th \in e \, \Om^1 \, e$ is uniquely 
determined by $\nabla$, one has $\th_X^* = -\th_X$, $~ \, X \in 
\hbox{Lie} 
\ G$ iff $\nb$ is compatible with the Hermitian  structure of $e \, 
A^{\infty}$.

\medskip

\begin{defn}\cite{connes:c_alg_dg}
	Let $\nabla$ be a connection on the finite projective module
	$\mathcal{E}^{\infty}$ (on $A^{\infty}$), the \textit{curvature} of $\nabla$
	is the element
	$\mathcal{T}$ of $\mathrm{End}_{A^{\infty}} (\mathcal{E}^{\infty}) \otimes\mathbb{L}^2 (\mathrm{Lie}~
	G)^*$ given by
	$$
	\mathcal{T} (X,Y) = \nabla_X \nabla_Y - \nabla_Y \nabla_X - \nabla_{[X,Y]} \in \mathrm{End}_{A^{\infty}} (\mathcal{E}^{\infty}) \, , \qquad \forall \, X,Y \in \mathrm{Lie}~G
	\, .
	$$
\end{defn}
\begin{defn}\cite{connes:c_alg_dg} A connection with zero curvature is said to be {\it flat}.
	Denote by $\mathrm{Con}_A(\mathcal{E})_0$ a space of flat connections.
\end{defn}

\begin{empt}\label{connection_construction}\cite{connes:c_alg_dg}
	If $e \in A^{\infty}$ is an idempotent then every connection on $e \, A^{\infty}$ is of the form $\nabla_X (\xi) =
	\nabla_X^0 (\xi) + \theta_X \, \xi$, $\forall \, \xi \in e \, A^{\infty}$, $X \in
	\hbox{Lie} ~ G$, where the form $\theta \in e \, \Om^1 \, e$ is uniquely
	determined by $\nabla$, one has $\theta_X^* = -\theta_X$, $\forall \, X \in
	\hbox{Lie}
	~ G$ iff $\nabla$ is compatible with the Hermitian  structure of $e \,
	A^{\infty}$. We identify $\mathrm{End} (e \, A^{\infty})$ with $e \, A^{\infty} \, e \subset
	A^{\infty}$, the {curvature} $\mathcal{T}_0$ of the grassmannian connection is
	the 2-form $e(de \wedge de) \in \Om^2$, the curvature of $\nabla = \nabla^0 +
	\theta \wedge$ equals to
	\begin{equation}\label{curvature_formula}
		\mathcal{T}_0 + e (d\theta + \theta \wedge \theta) \, e \in \Om^2.
	\end{equation}
	
\end{empt}
\begin{empt}\label{connes_transp_procedure}
	There is a natural correspondence between elements of Lie algebra and one-parameter transformation groups \cite{kobayashi_nomizu:diff_geom}.
	Let $\nabla$ be a connection (on $\mathcal{E}^{\infty}$). If $X\in \mathrm{Lie}~G$  defines a one-parameter group of automorphisms $\varphi$ of $A^{\infty}$  then $\nabla_X$ defines a  module parallel transport $\Phi$ on
	$\mathcal{E}^{\infty}$ along $\varphi$. So we have a connection transport procedure
	\begin{equation*}
		\mathbf{Transport} : \mathrm{Con}_{A^{\infty}}(\mathcal{E}^{\infty}) \times \mathrm{Paths}_A^{\infty}  \to \mathrm{Trans}_{\mathcal{E}^{\infty}}.
	\end{equation*}
	
\end{empt}
\begin{defn}
	Let $A$ be a $\mathbb{C}$-algebra, and let $\mathcal{E}$ be a finite projective $A$ module. Suppose that $\mathrm{Paths}_A$ is a set of one-parameter group of $A$ automorphisms.   A {\it connection transport procedure} is a natural map
	\begin{equation*}
		\mathbf{Transport} : \mathrm{Con}_A(\mathcal{E}) \times \mathrm{Paths}_A  \to \mathrm{Trans}_{\mathcal{E}}
	\end{equation*}
	such that
	$\mathbf{Transport}(\nabla, \varphi)$ is a transport along $\varphi$  for any $\nabla \in \mathrm{Con}_A(\mathcal{E})$ and $\varphi \in \mathrm{Paths}_A$.
\end{defn}
\begin{rem}
	A connection transport procedure should have a good  math and/or physical sense. Such procedures are known in following cases:
	\begin{enumerate}
		\item Commutative differential geometry \cite{kobayashi_nomizu:diff_geom,schenkel:nc_parallel}.
		\item $A = \mathbb{M}_n(\mathbb{C})$ \cite{schenkel:nc_parallel}.
		\item Noncommutative torus \cite{alekseev_bytsko:wilson_nc_tori,Ambjorn:2000cs,Driver,Sengupta}.
		\item Noncommutative differential geometry \cite{connes:c_alg_dg}.
	\end{enumerate}
\end{rem}

	\chapter{Presheaves and sheaves}
\section{Basic constructions}\label{sheaves_bas_sec}
\begin{definition}\label{presheaf_defn}\cite{hartshorne:ag}
	Let $\sX$ be a topological space. A \textit{presheaf} $\mathscr F$ of Abelian groups on  $\sX$ consists of the data
	\begin{itemize}
		\item[(a)] for every open subset $\sU \subseteq \sX$, an Abelian group $\mathscr F\left(\sU\right)$, and 
		\item[(b)] for every inclusion $\sV \subseteq \sU$ of open subsets of $\sX$, a morphism of Abelian groups $\rho_{\sU \sV}:\mathscr F\left(\sU\right) \to \mathscr F\left(\sV\right)$,\\
		subject to conditions
		\begin{itemize}
			\item [(0)] $\mathscr F\left(\sV\right)= 0$, where $\emptyset$ is the empty set,
			\item[(1)] $\rho_{\sU \sU}$ is the identity map, and
			\item[(2)] if $\mathcal W \subseteq \sV \subseteq \sU$ are three open sets, then $\rho_{\sU \mathcal W} = \rho_{\sV \mathcal W }\circ \rho_{\sU \sV}$.
		\end{itemize}
	\end{itemize}
\end{definition}

\begin{definition}\label{sheaf_defn}\cite{hartshorne:ag}
	A \textit{presheaf} $\mathscr F$ on  
	a topological space $\sX$ is a \textit{sheaf}  if it satisfies the following supplementary conditions:
	\begin{itemize}
		\item[(3)] If $\sU$ is an open set, if $\left\{\sV_{\a}\right\}$ is an open covering of $\sU$, and if $s \in \mathscr F\left(\sU\right)$ is an element such that $\left.s\right|_{\sV_{\a}}= 0$ for all $\a$, then $s = 0$;
		\item[(4)] If $\sU$ is an open set, if $\left\{\sV_{\a}\right\}$ is an open covering of $\sU$ (i.e. $\sU = \cup\sV_\a$), and we have elements $s_\a$ for each $\a$, with property that for each $\al, \bt, \left.s_\a\right|_{\sV_{\a}\cap \sV_{\bt}}= \left.s_\bt\right|_{\sV_{\a}\cap \sV_{\bt}}$, then there is an element $s \in \mathscr F\left(\sU\right)$ such that $\left.s\right|_{\sV_\a} = s_\a$ for each $\a$.
	\end{itemize}
	(Note condition (3) implies that $s$ is unique.)
\end{definition}
\begin{definition}\label{sheaf_stalk_defn}\cite{hartshorne:ag}
	If $\mathscr F$ is a {presheaf} on $\sX$, and if $x$ is a point of $\sX$ we define the \textit{stalk} or the \textit{germ} $\mathscr F_x$ 
	\textit{of} $\mathscr F$ \textit{at} $x$ to be the direct limit of groups $\mathscr F\left(\sU\right)$ for all open sets $\sU$ containing $x$, via restriction maps $\rho$.
\end{definition}

\begin{definition}\cite{hartshorne:ag}
If $\mathscr F$  and $\mathscr G$ are presheaves on $\sX$, a \textit{morphism} $\varphi:\mathscr F\to\mathscr G$  consists 
of a morphism of Abelian groups $\varphi_\sU:\mathscr F\left( \sU\right) \to\mathscr F\left( \sU\right)$ for each open set 
	$\sU$, such that whenever $\sV\subset\sU$ is an inclusion, the diagram 
	\newline
\begin{tikzpicture}
	\matrix (m) [matrix of math nodes,row sep=3em,column sep=4em,minimum width=2em]
	{
\mathscr F\left( \sU\right)   & \mathscr G\left( \sU\right)\\
\mathscr F\left( \sV\right)    & \mathscr G\left( \sV\right)\\
	};
	\path[-stealth]
	(m-1-1) edge node [above] {$\rho_{\sU\sV}$} (m-1-2)
	(m-1-1) edge node [right] {$\varphi_\sU$} (m-2-1)
	(m-1-2) edge node [right] {$\varphi_\sV$} (m-2-2)
	(m-2-1) edge node [above] {$\rho'_{\sU\sV}$} (m-2-2);
\end{tikzpicture}
\\
is commutative, where $\rho_{\sU\sV}$ and $\rho'_{\sU\sV}$ are the restriction maps in $\mathscr F$  and $\mathscr G$. If $\mathscr F$  and $\mathscr G$ are sheaves on $\sX$, we use the same definition for a morphism 
of sheaves. An isomorphism is a morphism  which has a two-sided inverse. 
\end{definition}

\begin{prdf}\label{sheaf_prdf}\cite{hartshorne:ag}
	Given a presheaf $\mathscr F$, there is a sheaf  $\mathscr F^+$ and a morphism $\th: \mathscr F \to \mathscr F^+$, with the property that for any sheaf  $\mathscr G$, and any morphism $\varphi: \mathscr F \to \mathscr G$, there is a unique morphism $\psi:\mathscr F^+\to \mathscr G$ such that $\varphi = \psi \circ \th$. Furthermore the pair $\left(\mathscr F^+, \th\right)$ is unique up to unique isomorphism. $\mathscr F^+$ is called the $\mathrm{sheaf~associated}$ to the presheaf $\mathscr F$. 
\end{prdf}

\begin{empt}\label{sheaf_empt}
	Following text is the citation of the proof of \ref{sheaf_prdf} (cf. \cite{hartshorne:ag}). For any open set $\sU$, let $\mathscr F^+\left(\sU\right)$ be set of functions $s$ from $\sU$ to the union $\bigcup_{x \in \sU}  \mathscr F_x$ of stalks of $\mathscr F$ over points of $\sU$, such that
	\begin{itemize}
		\item [(1)] for each $x \in \sU$, $s\left(x\right)\in \mathscr F_x$, and
		\item[(2)] for each $x \in \sU$, there is a neighborhood $\sV$ of $x$ contained in $\sU$ and an element $t \in \mathscr F\left(\sV\right)$, such that for all $y \in \sV$ the stalk (germ) $t_y$ of $t$ at $y$ is equal to $s\left(y\right)$.
	\end{itemize}
\end{empt}

\begin{exercise}\label{sheaf_etale_exer}\cite{hartshorne:ag}
	\textit{
		\'Espace Etal\'e of a Presheaf}. 
	Given a presheaf $\mathscr F$ on $\sX$, we define a topological space $\mathrm{Sp\acute{e}}\left(\mathscr F \right)$ , called the \textit{
		\'espace etal\'e} of a presheaf of $\mathscr F$ as 
	follows. As a set, $\mathrm{Sp\acute{e}}\left(\mathscr F \right)= \bigcup_{x\in\sX} \mathscr F_x$. We define a projection map $p: \mathscr \mathrm{Sp\acute{e}}\left(\mathscr F \right)$
	by sending $s_x\in \mathscr F_x$ to $x$. For each open set $\sU\subset\sX$ and each section $s\in \mathscr F\left(\sU \right)$  we 
	obtain a map: $\overline{s}: \sU \to  \mathrm{Sp\acute{e}}\left(\mathscr F \right)$ by sending $x \mapsto s_x$, its germ at $x$. This map has the property that 
	$p\circ \overline{s}= \Id_\sX$, in other words, it is a "section" of $p$ over $\sU$. We now 
	make $\mathrm{Sp\acute{e}}\left(\mathscr F \right)$ into a topological space by giving it the strongest topology such that 
	all the maps $\overline{s}: \sU \to  \mathrm{Sp\acute{e}}\left(\mathscr F \right)$ for all $\sU$ and all $s\in \mathscr F\left(\sU \right)$ , are continuous. Now 	show that the sheaf $\mathscr F^+$ associated to $\mathscr F$ can be described as follows: for any 
	open set $\sU\subset \mathscr F$, $\mathscr F\left(\sU \right)$ is the set of continuous sections of $\mathrm{Sp\acute{e}}\left(\mathscr F \right)$ over $\sU$. In 
	particular, the original presheaf $\mathscr F$ was a sheaf if and only if for each $\sU\subset \sX$,  $\mathscr F\left(\sU \right)$ is 
	equal to the set of all continuous sections of $\mathrm{Sp\acute{e}}\left(\mathscr F \right)$ over $\sU$. 
	
\end{exercise}
\begin{exercise}
Let $\mathscr F$  and $\mathscr G$ by presheaves. Prove that any morphism corresponds to a continuous map $\phi:\mathrm{Sp\acute{e}}\left(\mathscr F \right)\to \mathrm{Sp\acute{e}}\left(\mathscr G \right)$ such that
\begin{itemize}
	\item $p_{\mathscr G}= p_{\mathscr F} \circ \phi$ where both $p_{\mathscr G}: \mathrm{Sp\acute{e}}\left(\mathscr G\mathscr F \right)\to \sX$ are natural projections.
\item For all $x \in \sX$ a restriction $\phi|_{\mathscr F_x}: \mathscr F_x \to \mathscr G_x$ is a homomorphism of Abelian groups, rings, modules, etc.
\end{itemize}
\end{exercise}

\begin{remark}\label{sheaf_rem}
	If $\sV \subset \sX$ is any subset  then  one can define a group $\mathscr F\left(\sV\right)$  which is a group of continuous maps  $\iota :\sV\to \mathrm{Sp\acute{e}}\left(\mathscr F \right)$ such that $p\circ\iota = \Id_\sV$. Denote by  $\mathscr F\left(\sV \right)$ the Abelian group of such maps.
\end{remark}

\begin{definition}\label{sheaf_inv_im_defn}\cite{hartshorne:ag}
	Let $f: \sX\to \sY$ be a continuous map of topological spaces. For any sheaf  $\mathscr F$ on $\sX$, we define the \textit{direct image} sheaf  $f_*\mathscr F$ on $\sY$ by $\left(f_*\mathscr F\right)\left(\sV\right)= \mathscr F\left(f^{-1}\left(\sV\right)\right)$ for any open set $\sV \subseteq \sY$. For any sheaf  $\mathscr G$ on $\sY$, we define the \textit{inverse image} sheaf  $f^{-1}\mathscr G$ on $\sX$ be the sheaf  associated to the presheaf  $\sU \mapsto \lim_{\sV \supseteq f\left(\sU\right)} \mathscr G\left(\sV\right)$, where $\sU$ is any open set in $\sX$, and the limit is taken over all open sets $\sV$ of $\sV$ containing $f\left(\sU\right)$.
\end{definition}

\begin{definition}\label{sub_sheaf_defn}\cite{hartshorne:ag}
	A \textit{subsheaf}  of a sheaf $\mathscr F$ is a sheaf $\mathscr F'$ such that for every open set 
	$\sU \subset \sX$, $\quad \mathscr F'\left(\sU \right)$  is a subgroup of $\mathscr F\left(\sU \right)$, and the restriction maps of the 
	sheaf $\mathscr F'$ are induced by those of $\mathscr F$. It follows that for any point $P$, the 
	stalk $\mathscr F'_P$ is a subgroup of $\mathscr F_P$. 
		
	\end{definition}

\begin{definition}\label{sheaf_hom_defn}\cite{hartshorne:ag}
	Let $\mathscr F$, $\mathscr G$ be sheaves of Abelian  groups  on $\sX$. For any open set $\sU \subseteq \sX$ the set of morphisms
	 $\Hom\left(\left.\mathscr F\right|_{\sU}, \left.\mathscr G\right|_{\sU}\right)$ has the natural structure of Abelian group. 
	It is a sheaf  (cf. \cite{hartshorne:ag}). It is called the \textit{sheaf of local morphisms} of $\mathscr F\to \mathscr G$, "sheaf hom" for short, and is denoted by $\mathscr Hom \left(\mathscr F, \mathscr G\right)$.
\end{definition}

\begin{definition}\label{sheaf_ringed_space_defn}\cite{hartshorne:ag}
	A \textit{ringed space} is a pair $\left(\sX, \mathcal O_\sX \right)$  consisting of a topological space 
	$\sX$ and a sheaf of rings $\mathcal O_\sX$  on $\sX$.
			
		\end{definition}
		\begin{example} \label{sheaf_ringed_space_exm}
			One has a following ringed spaces (cf.Definition \ref{sheaf_ringed_space_exm})
			\begin{enumerate}
				\item [(i)] If $\sX$ is topological space and a sheaf   $\mathscr C_\sX$ associated with a presheaf
				$$
				\sU \mapsto C_0\left( \sU\right) 
				$$
				then $\left(\sX, \mathscr C_\sX \right)$ is a ringed space.
				\item [(ii)]If $M$ is smooth manifold and a sheaf   $\mathscr C^\infty_M$ associated with a presheaf
				$$
				\sU \mapsto \Coo\left( \sU\right) 
				$$
				then $\left(M, \mathscr C^\infty_\sX \right)$ is a ringed space.
			\end{enumerate}
		\end{example}
		
		\begin{definition}\label{sheaf_of_modules_den}\cite{hartshorne:ag}
			Let $\left(\sX, \mathcal O_\sX \right)$  be a ringed space. A \textit{sheaf  of $ \mathcal O_\sX $-modules} 
			(or simply an $\mathcal O_\sX$ -\textit{module}) is a sheaf  $\mathscr F$ on $\sX$, such that for each open set 
			$\sU\subset\sX$, the group 
			$\mathscr F\left(\sU \right)$ is an $\mathcal O_\sX$-module, and for each inclusion of 
			open sets $\sV \subset\sU$, the restriction homomorphism $\mathscr F\left(\sU \right) \to \mathscr F\left(\sV \right)$  is compatible with the module structures via the ring homomorphism $\mathcal O_\sX\left(\sU \right) \to \mathcal O_\sX \left(\sV \right)$.
		\end{definition}
\section{Flasque (flabby) and soft sheaves}
\begin{definition}\label{sheaf_flasque}\cite{hartshorne:ag}
	A sheaf  $\mathscr F$ on a topological space is \textit{flasque} or \textit{flabby} if for every inclusion $\sV \subseteq \sU$ of open sets, the restriction map $\mathscr F\left(\sU \right)\to\mathscr F\left(\sV \right)$ is surjective. 
\end{definition}
Recall that for $s \in \mathscr{F}\left(\sX \right)$, 
\be\label{sheaf_supp_eqn}
\supp s= \left\{\left.x \in \sX \right|s\left(x \right)\neq 0  \right\}
\ee
 denotes the \textit{support} of the section $s$.
\begin{definition}\label{phi_supp_defn}\cite{godement:sheaf}.
	Let $\sX$ be a topological space. A \textit{family of supports} on $\sX$ is a family $\Phi$ of closed
	subsets of $\sX$  such that:
	\begin{enumerate}
		\item a closed subset of a member of $\Phi$ is a member of $\Phi$;
		\item $\Phi$ is closed under finite unions.
	\end{enumerate}
	
	$\Phi$  is said to be a  \textit{paracompactifying} family of supports if in addition:
	\begin{enumerate}
		\item  each element of $\Phi$ is paracompact;
		\item each element of $\Phi$ has a (closed) neighborhood which is in $\Phi$.
	\end{enumerate}
\end{definition}\label{sheaf_soft_defn}
The family of all compact subsets of $\sX$ is denoted by $c$.
\begin{definition}\label{soft_sheaf_defn}\cite{bredon:sheaf}
	A sheaf  $\mathscr{F}$ on $\sX$ is called $\Phi$-\textit{soft} if the restriction map	$\mathscr{F}\left(\sX \right) \to  \mathscr{F}\left(K \right)$  is surjective for all $K \in \Phi$. If $\Phi$ is a set of all closed sets  then $\Phi$ is simply
	called \textit{soft}. If $\Phi= c$ is a set of all compact  sets  then $\mathscr{F}$ is called $c$-\textit{soft}.
\end{definition}
\begin{theorem}\label{sheaf_soft_m_thm}\cite{bredon:sheaf}
 Let $\left(\sX, \mathcal O_\sX \right)$ be a ringed space (cf. Definition \ref{sheaf_ringed_space_defn}), and let $\mathscr F$ be an $\mathcal O_\sX$-module (cf. Definition \ref{sheaf_of_modules_den}). Suppose that the space $\sX$ is paracompact and $\Phi$ is a {paracompactifying} family of supports(cf. Definition \ref{phi_supp_defn}). If the sheaf $\mathcal O_\sX$ is $\Phi$-soft (cf. Definition \ref{soft_sheaf_defn}) then  $\mathscr F$ is also $\Phi$-soft.
\end{theorem}
Now is $\mathscr P$ is a presheaf  on $\sX$ such and $s \in \mathscr P\left(\sX \right)$ we put
$\left|s\right| = \left|\th\left( s\right) \right|$, where $\th: \mathscr{P}\left(\sX \right)\to \mathscr{P}^+\left(\sX \right) $ is the canonical map, $\mathscr{P}^+$ being the sheaf  generated by $\mathscr{P}$.
Note that for $s \in \mathscr{P}$ one has $x \notin \left|s\right|\Leftrightarrow\left.s \right|_{\sU}$ for some neighborhood $\sU$ of $x$. If $\mathscr{F}$ is a sheaf, we put
$$
\Ga_\Phi = \left\{\left.s\in \mathscr{F}\left( \sX\right)\right| \left|s\right| \in \Phi \right\}.
$$Maybe   the general case sheaf cohomology do not preserve homotopy equivalence. I read a book written by Bredon, page 80.

\begin{theorem}\label{sheaf_homototy_thm} Any two properly homotopic maps (with respect to $\Phi$ and $\Psi$)
	of a space $X$  into a space $Y$ induce identical homomorphisms
	$$
	H^*_\Psi(Y;G)\to  H^*_\Phi(X;G),
	$$
	where $G$ is any constant coefficient group.
	
\end{theorem}
Note the special cases:
\begin{itemize}
	\item [(a)]  $\Phi = cld = \Psi$. In this case, "properly homotopic" is the same as
	"homotopic."
	\item [(b)] $X, Y$ locally compact Hausdorff, $\Phi = c = \Psi$.
\end{itemize}
There is a theorem in Bredon's book (page 80).
\begin{theorem} Any two properly homotopic maps (with respect to $\Phi$ and $\Psi$)
	of a space $X$  into a space $Y$ induce identical homomorphisms
	$
	H^*_\Psi(Y;G)\to  H^*_\Phi(X;G),
	$
	where $G$ is any constant coefficient group.
	
\end{theorem}

For a presheaf $\mathscr{P}$  on $\sX$ we put  $\mathscr{P}_\Phi\left(\sX \right) = \left\{\left.s\in \mathscr{P}\left( \sX\right)\right| \left|s\right| \in \Phi \right\}$.
\begin{lemma}\label{schv_soft_n_lem}\cite{godement:sheaf,torsten:sheaves}
	Let $\sX$ be a paracompact space, $\mathscr{F}$ a sheaf on $\sX$. Then $\mathscr{F}$ is soft if and only if every point $x \in \sX$ has a closed neighborhood $\overline{   \mathcal U }$ such that $\left.\mathscr{F}\right|_{\overline{   \mathcal U }}$ is soft.
\end{lemma}

\begin{theorem}\label{sheaf_neigh_thm}\cite{godement:sheaf}
	Let $\mathscr{F}$ be a sheaf  of set over the space $\sX$, $\mathcal S$ is a subset of $\sX$ and $s$ is a section of the sheaf  $\mathscr{F}$ over $\mathcal S$. If $\mathcal S$ allows a fundamental system of paracompact neighborhoods then $s$ can be extended to a neighborhood of $\mathcal S$ in $\sX$. 
\end{theorem}
\begin{theorem}
Soit $\mathscr F$ un faisceau surun espace $\sX$ paracompact. Supposons que 
tout point de $\sX$ possède un  voisinage $\sU$ vérifiant la condition suivante : toute section de $\mathscr F$
au-dessus d'un sous-ensemble fermé de $\sX$ contenu dans $\sU$, se prolonge à $\sU$. Alors $\mathscr F$ est mou.
\end{theorem}
Following text is an English translation of the above theorem
\begin{theorem}
	Let  $\mathscr{ F}$ be a sheaf on  a paracompact space $\sX$. Suppose that for any $x \in\sX$ there is a neighborhood  $\sU$ which satisfies to the following condition; every section of $\mathscr{ F}$ over closed subset  $\sY\subset \sX$ such that $\sY\subset \sU$ can be extended up to $\sU$. Then the sheaf $\mathscr{ F}$ is soft.
\end{theorem}
\begin{thm}\cite{godement:sheaf}
Soit $\mathscr{ A}$ un faisceau d'anneaux avec unité sur un espace paracompact $\sX$. Pour que $\mathscr{ A}$ soit mou, il faut et il suffit que tout point de $\sX$ possède un  voisinage $\sU$ tel que, 	étant donnés des fermé disjoints $S,T \subset \sU$, il existe une section de $\mathscr{ A}$ au-dessus de $\sU$, égale à 1 sur $S$ et à 0 sur $T$.
\end{thm}
Following text is an English translation of the above theorem.
\begin{theorem}\label{sheaf_a_soft_thm}
Let $\mathscr{ A}$ be a sheaf of unital algebras on a paracompact apace $\sX$. The sheaf $\mathscr{ A}$ is soft if and only if any $x\in\sX$ have a neighborhood $\sU$ such that for any closed sets $S, T\subset \sU$ one has
\bean
S\cap T = \emptyset \quad \Rightarrow \quad \exists s \in \mathscr{ A}\left( \sU\right) \quad s|_S = 1\quad \mathrm{AND} \quad \quad s|_T = 0.
\eean
\end{theorem}

\begin{proposition}\label{flabby_soft_prop}\cite{bredon:sheaf}
	If $\sU$ is a compact relatively Hausdorff subspace of $\sX$,
	then $\left.\mathscr{F}\right|_{\sU} $ is soft for any flabby sheaf  $\mathscr{F}$ on $\sX$.
\end{proposition}

\section{Sheaves over non Hausdorff spaces}\label{sheaves_nh_sec}
\begin{empt}\label{nh_csoft_empt}\cite{cra_moe:nhaus}
	Here we assume that the space $\sX$ has an open cover by subsets $\sU \subset \sX$ which are each paracompact, Hausdorff, locally compact and of cohomological dimension bounded a number $d$ (depending  on $\sX$ but not $\sU$).	
\end{empt}

\begin{definition}\label{nh_csoft_defn}\cite{cra_moe:nhaus}
	Let $\sX$ be a space satisfying the assumptions of \ref{nh_csoft_empt}. An Abelian sheaf $\mathscr{F}$ is said to be $c$-\textit{soft} if for any open $\sU \subset \sX$ its restriction $\left.\mathscr{F}\right|_{\sU}$ is a $c$-soft on $\sU$ in the usual sense. By the same proprety of for Hausdorff space, it follows that $c$-softness is a local property, i.e. a sheaf   $\mathscr{F}$ is $c$-soft if and only if there is an open cover $\sX= \cup \sU_\a$  such that $\left.\mathscr{F}\right|_{\sU_\a}$ is a $c$-soft sheaf  on $\sU$.
\end{definition}
\begin{definition}\label{nh_csoft_gc_defn}\cite{cra_moe:nhaus}
	Let $\mathscr{F}$ be a $c$-soft sheaf  and let $\mathscr{F}'$  be its Godement resolution (i.e. $\mathscr{F}'= \Ga\left(\sU_{\text{distr}}, \mathscr{F} \right)$ is the set of all (not necessary continuous) sections for any open $\sU \subset \sX$). For any Hausdorff open set $\mathcal W$, let $\Ga_c\left(\mathcal W, \mathscr{F} \right)$ be the usual set of compactly supported sections. If $\mathcal W \subset \sU$, there is an evident homomorphism "extension by 0" 
	\be\label{sheaf_inc_eqn}
	\Ga_c\left(\mathcal W, \mathscr{F} \right) \hookto \Ga_c\left(\mathcal U, \mathscr{F} \right)\subset \Ga_c\left(\mathcal U, \mathscr{F}' \right).
	\ee
 For any (not necessary Hausdorff) open set $\sU \subset \sX$ we define $\Ga_c \left(\sU, \mathscr{F} \right)$ to be the image of the map:
	$$
	\bigoplus_{\mathcal W} \Ga_c\left(\mathcal W, \mathscr{F}' \right) \hookto \Ga\left(\mathcal U, \mathscr{F}' \right),
	$$ 
	where $\mathcal W$ ranges over all open subsets $\mathcal W \subset \mathcal U$. 
\end{definition}
\begin{proposition}\label{nh_csoft_gc_lem}
	Let  $\mathscr{F}$ be a $c$-soft sheaf. For any open cover $\sU = \cup \mathcal W_\a$ where each $\mathcal W_\a$ is Hausdorff the sequence $\oplus_\a \Ga_c\left(\mathcal W_\a,  \mathscr{F}\right) \to \Ga_c\left(\mathcal U,  \mathscr{F}\right)\to 0$ is exact.  
\end{proposition}

\chapter{Groupoids foliations, pseudogroups and operator algebras}\label{foliations_sec}
		\section{Groupoids}
	\paragraph*{}
	A groupoid is a small category with inverses, or more explicitly:
	\begin{definition}\label{groupoid_defn}\cite{connes:ncg94}
		A \textit{groupoid} consists of a set $\G$, a distinguished subset $\G^0\subset\G$, two maps
		$r, s : \G\to \G^0$ and a law of composition
		$$
		\circ: \G^2\bydef\left\{\left.\left(\ga_1,\ga_2 \right) \in \G\times\G~\right| s\left(\ga_1\right)= r\left(\ga_2\right)\right\}\to \G
		$$
		such that
		\begin{enumerate}
			\item $s\left(\ga_1\circ\ga_2\right)=s\left(\ga_2\right), \quad r\left(\ga_1\circ\ga_2\right)=r\left(\ga_1\right)\quad \forall\left(\ga_1, \ga_2 \right) \in \G$
			\item $s\left(x\right)=r\left(x\right)=x \quad\forall x\in\G^0$
			\item $\ga\circ s\left(\ga\right)= r\left(\ga\right)\circ\ga = \ga\quad \forall\ga\in\G$
			\item $\left( \ga_1\circ\ga_2\right) \circ\ga_3=\ga_1\circ\left( \ga_2\circ\ga_3\right) $
			\item Each $\ga \in\G$ has a two-sided inverse $\ga^{-1}$, with $\ga\circ\ga^{-1}=r\left(\ga\right)$, $\ga^{-1}\circ\ga=r\left(\ga\right)$.
		\end{enumerate}
		The maps $r$, $s$ are called the \textit{range} and \textit{source} maps.
	\end{definition}
	\begin{definition}\label{groupoid_sets_defn}
If $A$ and $B$ are subsets of $\G$, one may form the following subsets o f $\G$ :
\bean
A^{-1} \bydef \left\{x \in \G \left| x^{-1}\in A\right.\right\},\\
AB \bydef \left\{z \in \G | x \in A, ~y \in B \quad  z = x y \right\}.
\eean
A groupoid $\G$ is said to be \textit{principal} if the map $( r , s )$ from $\G$ into $\G^0\times \G^0$ is one-to-one, it is said to be \textit{transitive}  the map ( r , d ) is onto.
For $u, v, \in \G^0$ , $G^u \bydef r^{-1}(u),\quad G_v \bydef s^{-1}(v), \quad G^u_v \bydef G^u\cap G_v $ and
$G(u) = G^u_u$ which is a group, is called  the \textit{isotropy group} at $u$.
The relation  $u \sim v$ if and only if $G^u\cap G_v \neq \emptyset$ is an equivalence relation on the unit space $\G^0$. its equivalence classes are called \textit{orbits} and the \textit{orbit} of  $u$ is denoted $[ u ]$ . 
$G^0/G$
denotes the \textit{orbit space}. A groupoid is transitive if and only if it has a single orbit.	
\end{definition} 
	\begin{example}
The set $\G^2$ of composable elements may be given the following groupoid structure:
$( x , y )$ and $( y ' , z )$ are composable if and only if  $y' = xy, ~ ( x , y ) ( x y , z ) = ( x , y z )$, and $( x , y )^{-1} =( x y , y ^{-1} )$.
Then $r^2 ( x , y ) = ( x , r ( y ) ) = ( x , d ( x ) )$ and $d^2 ( x , y ) = ( x y , d ( x y ) )$. The map
$x\mapsto ( x , d ( x ) )$ identified the unit space of $G^2$ with $G$. The groupoid $G^2$ is principal.
One may notice that it  comes from the action of $\G$ G on itself. It is transitive if and only if $\G$  is a group.
	\end{example}
		
	\begin{definition}\label{groupoid_hom_defn}\cite{renault:gropoid_ca}
	Let $\G$ and $\mathcal H$ be groupoids a map $\phi: \G \to \H$ is \textit{homomorphism} if one has:
	\begin{itemize}
		\item
	\bean
	\left(x, y \right) \in \G^2 \quad \Rightarrow \quad \left(\phi\left( x\right), \phi\left(y\right)  \right) \in \H^2,\\
	\phi\left(\G^0 \right) \subset \H^0, 
	\eean
	\item the  map
	\bean
	\phi^2 : \G^2 \to \H^2,\\
	\left(x, y\right)\mapsto \left(\phi(x), \phi(y) \right) 
	\eean 
is a homeomorphism. 
\end{itemize}
Two homomorphism are \textit{similar} (write $\phi\sim \psi$) if there exists a function $\th: \G^0 \to \H$ such that $\left(\th\circ r\right)(x)\phi\left(x\right)= \psi\left(x\right)\left( \th\circ s\right)(x)$. Groupoids $\G$ and $\H$ are called \textit{similar} (write $\G\sim \H$) if there exists homomorphisms $\phi: \G \to \H$ and $\psi : \H\to \G$  such that $\phi \circ \psi$ and $\psi \circ \phi$ are similar to identity isomorphisms. 
	\end{definition}
	\begin{definition}\label{groupoid_reduction_defn}\cite{renault:gropoid_ca}
	Let $\G$ be a groupoid, and let $E$ be a subset of $\G^0$.
	A subgroupoid 
	$$
	\G^E_E \bydef \left\{x \in G | r(x), s(x)\in E\right\}
	$$
	with unit space $E$ is said to be  the \textit{reduction} 
	of $\G$ by $E$.
\end{definition}
\begin{defn}\label{groupoid_isotropy_defn}\cite{renault:gropoid_ca}
For $u, v\in \G^0$, $~\G^u\bydef r^{-1}\left( u\right)$,  $~\G_v\bydef s^{-1}\left( v\right)$  $~\G^u_v\bydef \G^u\cap \G_v$ and
$\G(u) = G^u_u$ which is a group, is called  the \textit{isotropy  group} at $u$.
\end{defn}
\begin{prop}\label{groupoid_reduction_prop}\cite{renault:gropoid_ca}
	Let $\G$ be a groupoid, $E$ a subset of $\G^0$ which meets each orbit in $\G$; then  	$\G^E_E$ and $\G$ are similar (cf. Definition \ref{groupoid_hom_defn}).
\end{prop}
\begin{definition}\cite{renault:gropoid_ca}
Let $\G$ be a groupoid, $A$ a group and $c: \G\to A$ a homomorphism, the
\textit{skew-product} $\G(c)$ is the groupoid $\G\times A$ where : $( x , a )$ and $(y,b)$ are composable if and only if  $x$ and $y$ are composable and $b = a c ( x )$, $( x , a ) ( y , a c ( x ) )\bydef ( x y , a )$, and $( x , a )^{-1} =
\left(  x^{-1}  , a c ( x )\right)$; $\quad ( x , a )  \bydef( r ( x ) , a )$ , $\quad s ( x , a )\bydef ( d ( x ) , a c ( x ) )$ . Its unit space is $\G^0\times A$.
\end{definition}

\begin{definition}\cite{renault:gropoid_ca}
 Let $\G$ be a groupoid, let $A$ be a group and let $\a: A \to \Aut\left(\G\right)$ be a
homomorphism. We write $x*a \bydef \left[\a\left(a^{-1}\right)\right]$ for $a\in A$ and $x\in \G$. The \textit{semi-direct product} $G\rtimes_\a A$ is the groupoid $G\times A$ where $( x , a )$ and $( z , b )$ are composable if and only if $z = y*a$ with  $x$ and $y$ composable,$ ( x , a ) ( y * a,b) = ( x y , a b )$ , and $(x,a)^{ -1}\bydef \left( x ^{-1} * a, a^{ - l} \right)$ .
Then, $r ( x , a ) = ( r ( x ) , e )$ and $s ( x , a ) = (d(x) • a , e )$ . The unit  space may be identified 
with $\G^0$.
\end{definition}

\begin{proposition}
With above notation,
\begin{enumerate}
	\item[(i)] $\G(c) \rtimes_\a A$  is similar to $\G$ and
	\item[(ii)] $\left( \G\rtimes_a A\right) (c)$ is similar to $\G$.
\end{enumerate}

.\end{proposition}
\begin{definition}
 An \textit{inverse semi}-\textit{group} is a set $\mathscr G$  endowed with an associative binary operation , noted multiplication, and an inverse map
\bean
 \mathscr G\to  \mathscr G,\\
 s \mapsto s^{-1}
\eean
such that the following relations  are satisfied  $ss^{-1}s = s$ and $s^{-1}ss^{-1}=s^{-1}$ .
\end{definition}
\begin{definition}
Let $\G$ be a groupoid. A subset $s$ of $\G$ will be called a $\G$-\textit{set} if 
the restriction of $r$ and $s$ to it are one-to-one . Equivalently, $s$ is a $\G$-set if and only is $s^{-1}s$
and $ss^{-1}$  are contained in $\G^0$.

\end{definition}

\begin{definition}
Suppose that $\mathscr C$ is some category. A map $p$ from a set $A$ onto a set $A^0$ such that
each fiber $p^{-1} ( u )$ is an object of  $\mathscr C$ will be called a  $\mathscr C$-\textit{bundle} map and $A$ will be called  $\mathscr C$-bundle.  Let $A$ be
a  $\mathscr C$ - bundle with the bundle map $p: A \to A_0$. Write $A_u \bydef p^{-1} ( u )$.
$$
 \text{Iso}(A) = \left\{\text{isomorphisms }\phi_{u,v}| A_u\to A_v\quad u,v \in A^0\right\} 
 $$	
	has a natural structure of groupoid. 
	\end{definition}
\begin{definition}
	Let $\G$ be a groupoid. A $\G$-\textit{bundle} $(A,L)$ is a $\mathscr C$ - bundle $A$ together
	with a homomorphism $L : G \to \text{Iso} (A)$ such that $L^0: \G^0\to A^0$ is a bijection . (We will often identify $\G^0$ and $A^0$). When  $\mathscr C$  is the category of Abelian groups, one speaks of a
	$\G$-\textit{module bundle}.
\end{definition}
\begin{empt}\label{groupoig_gn_empt}
Given a $\G$-module bundle $( A , L )$, one can form the following cochain complex. Let
us first define $\G^n$ for any $n\in\N$. The sets $\H^0$ , $\G^1\bydef \G$ and $\G^2$ have already been defined. For $n> 2$, $\quad \G^n$ is the set of $n$-tuples $(x_0 . . . . . x_{n-1}) \in \G\times...\times \G$ such that for $
j = 1 , . . . , n - 1$ , $\quad x_j$ is composable with its left  neighbor. A $n$-cochain is a function from $G^n$ to $A$ which satisfies the conditions
\begin{enumerate}
	\item[(i)] $p\circ f(x_0 . . . . . x_{n-1}) = d(x_0)$ and
	\item[(ii)] if $n > 0$ and for some $j = 0, . . . , n - 1$ , $\quad x_0 \in \G^0$, then $f ( x_0, ..., x_j , ..., x_{n-1})\in A^0$.
\end{enumerate}
The set $C^n\left(\G, A\right)$ of $n$-cochains is an Abelian group under point-wise addition. The
sequence 
\bean
0 \to C^0\left( \G, A\right)\to C^0\left( \G, A\right)\to C^1\left( \G, A\right)\to...\to  C^n\left( \G, A\right) \xrightarrow{\delta^n} C^{n + 1}\left( \G, A\right)\to ...
\eean
 where 
 \be\label{groupoid_b_eqn}
 \begin{split}
 \delta^0f(x)\bydef L(x)\quad  f\circ s(x) - f \circ r ( x ),\\
\delta^n(f(x_0 ,..., x_n) = L ( x_0 ) f ( x_1,..., x_n) +\\+ \sum_{j=1}^n (-1)^j
f ( x_0 ,..., x_{j-1}x_j . . . . . x_{n-1})+(-1)^{n-1} f ( x_0 , . . . , x_{n-1}) \quad  n > 0,
\end{split}
\ee

 is a cochain complex.
 \end{empt}	
\begin{definition}\label{groupoid_cocycle_defn}
The group of $n$-cocycles of this complex will be denoted by $Z^n(\G,A)$,
the group of $n$-coboundaries will be denoted by $B^n(\G,A)$ and the $n$-th cohomology group
$Z^n(\G,A)/B^n(\G,A)$ will be denoted by $H^n(\G,A)$.
\end{definition}

	\begin{definition}\label{groupoid_topological_defn}\cite{renault:gropoid_ca}
	A \textit{topological groupoid} consists of a groupoid $\G$ and a topology compatible with the groupoid structure:
	\begin{enumerate}
		\item [(a)] $\G \to \G \quad x \mapsto x^{-1}$ is continuous,
		\item [(b)] $\G^2\to \G\quad \left(x,y\right)\mapsto xy$ is continuous where $\G^2$ has the induced topology from $\G \times \G$.
	\end{enumerate}
		\end{definition}
	
	\begin{remark}\cite{renault:gropoid_ca}
		One has:
	\begin{itemize}
		\item the map $x \mapsto x^{-1}$ is a homeomorphism,
		\item if $\G$ is Hausdorff then $\G^0$ is closed in $\G$,
		\item if $\G^0$ is Hausdorff then  $\G^2$ is closed in $\G \times \G$, $\G^0$ is both a subspace of $\G$ and a quotient of $\G$ (by the map $r$), the induced and the quotient topology coincide.
	\end{itemize}
	\end{remark}
	\begin{definition}\label{groupoid_haar_defn}\cite{renault:gropoid_ca}
 Let $\G$ be a locally compact groupoid. A \textit{left  Haar system} for $\G$
consists of measures $\left\{\left.\la^u \right| u \in \G^0\right\}$ on $\G$ such that
\begin{enumerate}
	\item [(a)] the support $\supp\la^u$ of the measure $\la^u$ is $\G^u$,
	\item [(b)]  (continuity) for any $f \in C_c\left(\G\right)$, $u \mapsto \la(f)(u) = \int f d\la^u$ is continuous, and
	\item [(c)]  (left invariance) for any $x\in \G$ and any $f \in  C_c(\G )$, $\int  f ( x y ) d\la^{s(x)}(y) =
\int f(y)d\la^{r(x)}(y)$.

\end{enumerate}
\end{definition}
\section{Covering groupoids}\label{covering_groupoid_sec}
\paragraph{}
The theory of covering groupoids is explained in \cite{zhi:cov_group}. Given $u \in \G^0$, a \textit{sieve} on $\G$ is a set $S\subset G^u$ (cf. Definition \ref{groupoid_sets_defn}) G such that
\bean
f \in S \text{ and the composite } fh \text{ is defined implies} fh \in S.
\eean
Since $x \in \G$ is invertible, then every nonempty sieve on $\G$ coincides with $G^u$.
\begin{definition}\label{covering_groupoid_defn}\cite{zhi:cov_group}
 Let $p : \widetilde{\G}\to G$ be a morphism of groupoids. An ordered pair
$\left( \widetilde{\G}, p\right)$ is a \text{covering groupoid} if for each object $\widetilde u \in \widetilde{\G}$ the restriction of $p$
$$
\widetilde{\G}^{\widetilde{u}}\to {\G}^u
$$
is bijection. The morphism $p$ is called the \textit{covering}. 
\end{definition}

\begin{theorem}\label{covering_groupoid_thm}\cite{zhi:cov_group}
Let $\left(\widetilde \G, p \right) $ be a covering groupoid of $\G$
with $p\left(\widetilde{u} \right) = u$  where $\widetilde u \in \widetilde \G$ and $u \in \G$, and $f:  \F \to \G$ a groupoid
morphism with $f\left( v\right)= u$  such that$\F$ is connected. Then $f$ lifts
to a morphism $\widetilde f:  \F \to \widetilde\G$ with $\widetilde f\left(v \right)  =\widetilde u$ if and only if $f^*\F^v_v \subset p^*\widetilde\G^{\widetilde u}_{\widetilde u}$
and if this lifting exists, then it is unique.
\end{theorem}

There are following motivations this definition:
\begin{itemize}
	\item analog of unique lifting theorem,
	\item  analog of fundamental group,
	\item analog of Galois theory.
\end{itemize}

\section{Groupoid $C^*$-algebras}
\paragraph{}
	Let $\G$ be a locally compact Hausdorff groupoid with left  Haar system $\left\{\la^u\right\}$ and let $\sigma$ be a continuous 2-cocycle in $Z^2\left(\G, \T\right)$. For $f ,g \in C_c(\G, \sigma )$, let us define
\be\label{groupoid_*_c_eqn}
\begin{split}
	f * g \left(x\right)\bydef 
	\int f ( x y ) g \left( y^{-1}\right)\sigma\left(xy, y^{-1} \right) d\la^{d(x)}(y),\\
		f^* ( x ) \bydef \overline{f ( x^{ -1})}~\overline{\sigma\left(x, x^{-1} \right)}, 	
\end{split}
\ee	
In particular is $\sigma$ is trivial then one has a *-algebra
\be\label{groupoid_*_eqn}
\begin{split}
	f * g \left(x\right)\bydef 
	\int f ( x y ) g \left( y^{-1}\right) d\la^{d(x)}(y),\\
	f^* ( x ) \bydef\overline{f ( x^{ -1})} 	
\end{split}
\ee	
which is a specialization of \eqref{groupoid_*_c_eqn}
	
\begin{empt} 
Let $\G$ be a locally compact groupoid with left  Haar system $\left\{\la^u\right\}$  For $f$ and $g\in C_c\left(\G\right)$, let  us define
\be\label{groupoid_*__defn}
\begin{split}
f * g \left(x\right)\bydef 
 \int  f ( x y ) g \left( y^{-1}\right) d\la^{d(x)}(y),\\
f^* ( x ) \bydef\overline{f ( x^{ -1})} 	
\end{split}
\ee
It is proven in \cite{renault:gropoid_ca} that the equations yield a *-algebra. 
\end{empt}
\begin{definition}\label{groupoid_representation_defn}\cite{renault:gropoid_ca}
A \textit{representation} of $C_c\left(\G, \sigma\right)$ on a Hilbert space $\H$ is a $*$-homomorphism $L : C_c\left(\G, \sigma\right) \to B\left(\H \right)$ which is continuous when $C_c\left(\G, \sigma\right)$ has the inductive limit 
topology (cf. Definition \ref{top_ind_lim_defn}) and $B\left(\H \right)$  the weak operator topology (cf. Definition \ref{weak_topology_defn}), and is such that the linear  span of
$$
L\left(f \right) \xi , \quad f \in C_c\left(\G, \sigma\right), \quad \xi \in \H
$$
is dense in $\H$.
\end{definition}
\begin{lemma}\label{groupoid_mult_repr_lem}\cite{renault:gropoid_ca}
If $L$ is a representation of $C_c\left(\G, \sigma\right)$, there exists a unique representation
$M$ of $C_c\left(\G^0\right)$ such that for every $h \in C_c\left(\G^0\right)$ and every $f\in C_c\left(\G, \sigma\right)$, $L\left(h f \right)= M\left(h \right)L\left(f \right)$  and
$L\left(fh \right)= L\left(f \right)M\left(h \right)$. 
\end{lemma}

	\begin{definition}\label{foli_groupoid_red_defn}\cite{renault:gropoid_ca}
	If $\Pi$ is the set of irreducible representations of $\C\left[\G \right]$ then the completion of $\C\left[\G \right]$ with respect to $C^*$-norm
	\be\label{groupoid_red_norm}
	\left\| a\right\|_r = \sup_{\pi\in \Pi}\left\|\pi\left( a\right)\right\| 
	\ee
	is said to be the \textit{reduced algebra} of $\G$. It will be denoted by $C^*_r\left(\G\right)$.
\end{definition}
\begin{empt}\label{groupoid_reg_empt}\cite{renault:gropoid_ca}
Let $\sigma$ be a 2-cocycle and $\mu$ a quasi-invariant measure. Consider the measurable
field of Hilbert space $\left\{L^2\left(\G, \la^u \right) \right\}_{u \in \G^0}$ with square integrable sections
$\int^\oplus L^2\left(\G, \la^u \right) d\mu\left( u\right)= L^2\left(u \right)$. For $x\in \G$, define $L_u\left( x\right)$  mapping $L^2\left( \G, \la^{d\left( x\right) }\right)$  to
$L^2\left( \G, \la^{r\left( x\right) }\right)$ by $L_u\left( x\right)\xi\left( y\right) \bydef \sigma\left( x, x^{-1}\right) \xi\left( x^{-1}y\right)$. This yields a given by 
\be\label{groupoid_reg_eqn}
L_u :  C_c\left(\G, \sigma\right)\to B\left(L^2\left(u \right)  \right) 
\ee
$\sigma$-representation of $\G$,
\end{empt}
\begin{definition}\label{groupoid_reg_defn}\cite{renault:gropoid_ca} 
The above $\sigma$-representation of $\G$ will be called the $\sigma$-\textit{regular representation
of $\G$ on $\mu$}. Its integrated form is the \textit{regular representation on $\mu$ of
$C_c\left(\G, \sigma \right)$}. 
\end{definition}
\begin{proposition}\label{groupoid_reg_prop}\cite{renault:gropoid_ca}
	$C_c\left(\G, \sigma \right)$ has a faithful family of bounded representations, consisting
of regular representations.
\end{proposition}
It results from Proposition \ref{groupoid_reg_prop} that the function defined by 
\be\label{groupoid_red_norm_eqn}
\begin{split}
\left\|\cdot  \right\|_r : C_c\left(\G, \sigma \right)\to \R,\\
f \mapsto \sup_{u \in \G^0} \left\|L_u\left( f\right)   \right\|,
\end{split}
\ee
where $L_u$ ranges over all representations induced from the unit space, is
a $C^*$ -norm on $ C_c\left(\G, \sigma \right)$ dominated by the $C^*$ -norm $\left\|f  \right\|$.
\begin{definition}\label{groupoid_red_defn}	\cite{renault:gropoid_ca}
 The \textit{reduced $C^*$ -algebra} $C^*_r\left(\G, \sigma \right)$ of $\G$ is the completion of
$C^*_r\left(\G, \sigma \right)$ for the reduced norm $\left\|\cdot  \right\|_r$.
\end{definition}
\subsection{$C^*$-algebras of non-Hausdorff \'etale groupoids}
\paragraph{} Here I follow to \cite{neshv:non_haudorff}. 
\begin{definition}
Assume $\G$ is a locally compact, not necessarily Hausdorff, \'etale groupoid.
By this, we mean that $\G$ is a groupoid endowed with a locally compact topology
such that
\begin{itemize}
	\item  the groupoid operations are continuous;
	\item the unit space $\G^0$ is a locally compact Hausdorff space in the relative
	topology;
	\item the range map $r : \G \to \G^0$ and the source map $r : \G \to \G^0$ are local
	homeomorphisms.
	\end{itemize}
\end{definition}

	For an open Hausdorff subset $\sV \subset \G$, consider the usual space $C_c\left( \sV\right)$  of continuous
	compactly supported functions on $\sV$ . Every such function can be extended
	by zero to $\G$; in general, this extension is not a continuous function on $\G$.
	This way, we can view  $C_c\left( \sV\right)$ as a subspace of the space of functions $\text{Func}\left(\G \right)$ 
	on $\G$. For arbitrary open subsets $\sU\subset G$ we denote by $C_c\left( \sU\right)$  the
	linear span of the subspaces $C_c\left( \sV\right) \subset \text{Func}\left(\G \right)$ for all open Hausdorff subsets
	$\sV \subset \sU$ Instead of all possible $\sV$ , it suffices to take a collection of open bisections
	covering $\sU$.
	
\section{Strong Morita equivalence of groupoid $C^*$-algebras}

\begin{empt}\cite{renault:gropoid_equiv}
	The definition of a $\G$-\textit{space} $\sX$ is a straightforward generalization of that for a group action. Here we require a continuous open map from the locally compact space $\sX$ onto $\G^0$, which we call $\rho$ or $\sigma$, according to the side on which $\G$ acts. For example a \textit{left} $G$-\textit{space} is given by a continuous map $\G*\sX\to\sX$ where $\G*\sX$ denotes the set of composable pairs $\left(\ga, x\right)$ with $s(\ga)=\rho(x)$. 
	
	We say that the action is \textit{free} id $\ga \cdot x = x$ only when $\ga$ is a unit. We say that the action is \textit{proper} if the map
	\bean
	\G*\sX \to \sX \times \sX,\\
	(\ga, x)\mapsto (\ga\cdot x, x)
	\eean 
	is proper.	 The space is a \textit{principal} $\G$-space if the action is both free and proper.
\end{empt}
\begin{definition}\label{groupoid_equiv_defn}\cite{renault:gropoid_equiv}
	Let $\G$ and $\H$ be locally compact groupoids. We say that a locally compact space $\sZ$ is a $(\G,\H)$-\textit{equivalence} if
	\begin{enumerate}
		\item [(a)] $\sZ$ is a left  principal $\G$-space,
		\item [(b)] $\sZ$ is a right principal $\H$-space,
		\item [(c)] the $\G$ and $\H$ actions commute.
		\item [(d)] the map $\rho$ induces a bijection of $\sZ/\H$ onto $\G^0$, and
		\item [(e)]  the map $\sigma$ induces a bijection of $\G\backslash\sZ$ onto $\G^0$.
	\end{enumerate}
\end{definition}
\begin{example}\cite{renault:gropoid_equiv}
Let $\G$ be locally compact Hausdorff groupoid and let $N$ be a closed subset of $\G^0$ that meets each orbit in $\G^0$. Then as easy to see, $\G_N$ is a principal left  $\G$-space and a principal right $\G^N_N$-space. The maps $\sigma\bydef  s|_{\G_N}$ and $\rho\bydef  r|_{\G_N}$ satisfy to (d) and (e) of Definition \ref{groupoid_equiv_defn}, so if they are open then $\G_N$ is a $\left(\G,\G^N_N \right)$-equivalence.
\end{example}
\begin{theorem}\label{groupoid_morita_defn}\cite{renault:gropoid_equiv}
Suppose that $(\G, \la)$ and $\left(\H, \bt\right)$ are second countable , locally compact groupoids with Haar systems $\la$ and $\bt$. Then for any $\left(\G, \H\right)$-equivalence $\sZ$, $C_c\left(\sZ\right)$ can naturally be completed into $C^*(\G, \la)$-$C^*\left(\H, \bt\right)$ imprimitivity bimodule. In particular $C^*(\G, \la)$ and $C^*\left(\H, \bt\right)$ are strongly Morita equivalent.
\end{theorem}	
	
	\begin{lemma}\cite{renault:gropoid_equiv}
		Let $\Om$ be a principal left  $\G$-space.
		\begin{enumerate}
			\item [(i)] If $F \in C_c\left(\Om\times \G\right)$ then
			$$
		\varphi\left(\om, u\right)	\bydef \int_{\G} F\left(\om, \ga \right)d\la^u\left(\ga\right)
			$$
			defines an element of $C_c\left(\Om\times \G^0\right)$.
			\item[(ii)] If $f \in C_c\left(\Om\right)$ then
			$$
			\la(f)\left(\left[\om\right]\right)\bydef \int_{\G} f\left(\la^{-1}\cdot \om\right)d\la^{\rho(\om)}(\ga)
			$$
			defines a surjection of $C_c(\Om)$ onto $C_c\left(\G\backslash\Om\right)$.
		\end{enumerate}
	\end{lemma}
\begin{empt}
There are two pre-$C^*$-algebras $A \bydef C_c(\G, \la)$ and $B \bydef C_c\left(\H, \bt\right)$. We define the left  $A$-action and the right $B$-action as follows:
\bean
f \cdot \varphi(z)\bydef\int_{\G}f(\ga)\varphi\left(\ga^{-1}\cdot z\right)d\la^{\rho(z)}(\ga),\\
\varphi \cdot g(z)\bydef\int_{\G}\varphi(z\cdot \eta)g \left(\eta^{-1}\right)d\bt^{\sigma(z)}(\eta)
\eean
where $\varphi \in C_c\left(\sZ\right)$, $\quad f\in A$, and $g \in B$. Also we define $B$ and $A$-valued inner products
\bean
\left\langle\varphi, \psi \right\rangle_B(\eta)\bydef \int_{\G} \overline{\varphi\left(\ga^{-1}\cdot z\right)}\psi\left(\ga^{-1}\cdot z\cdot \eta\right)\d^\la{\rho(z)}(\ga)\in B,\quad \sigma(z)=r(\eta);\\
\left\langle\varphi, \psi \right\rangle_A(\ga)\bydef \int_{\G} \varphi\left(\ga^{-1}\cdot z\cdot \eta\right)\overline{\psi\left(\ga^{-1}\cdot \eta\right)}\d^\bt{\sigma(z)}(\eta)\in B,\quad \rho(z)=r(\ga);
\eean 
\end{empt}
	\begin{definition}\label{groupoid_discete_defn}\cite{renault:gropoid_ca}
A locally compact groupoid is $r$ - \textit{discrete} if its unit space is an
open subset.
\end{definition}
\begin{lemma}\label{groupoid_discete_lem}\cite{renault:gropoid_ca}
 Let $\G$ be an $r$-discrete groupoid.
 \begin{enumerate}
 	\item [(i)] For any $u\in \G^0$, both $\G^u$ and $G_u$ are discrete spaces.
 \item[(ii)] If a Haar system exists, it is essentially the counting measures system.
 \item[(iii)] If a Haar system exists, $s$ and $r$ are local homeomorphisms.	
 \end{enumerate}

\end{lemma}
	\begin{empt}
		Let $\G$ be a groupoid. 
		Consider an involutive algebra $\C\left[\G \right]$ over $\C$ generated by $\G$ which satisfies to the following relations
		\be\label{groupoid_a_eqn}
		\begin{split}
			\left(a\cdot b\right)\left(\ga\right)= \sum_{\ga_1\circ\ga_2=\ga}a\left( \ga_1\right)a\left( \ga_1\right),\\
			a^*\left(\ga\right)= \overline{a\left(\ga^{-1}\right)}.
		\end{split}
		\ee
	\end{empt}-
	
\section{Foliations and pseudogroups}
	\begin{definition}\cite{connes:ncg94}		
		Let $M$ be a smooth manifold and $TM$ its tangent bundle, so that
		for each $x \in M$, $T_x M$ is the tangent space of $M$ at $x$. A
		smooth subbundle $\mathcal{F}$ of $TM$ is called {\it integrable} if and only if one of
		the following equivalent conditions is satisfied:
		
		\smallskip
		
		\begin{enumerate}
			
			\item[(a)] Every $x \in M$ is contained in a submanifold $W$ of $M$ such that
			$$
			T_y (W) = \mathcal{F}_y \qquad \forall \, y \in W \, ,
			$$
			
			\smallskip
			
			\item[(b)] Every $x \in M$ is in the domain $U \subset M$ of a
			submersion $p : U \to {\mathbb R}^q$ ($q = {\rm codim} \, \mathcal{F}$) with
			$$
			\mathcal{F}_y = {\rm Ker} (p_*)_y \qquad \forall \, y \in U \, ,
			$$
			
			\smallskip
			\item[(c)] $C^{\infty} \left( \mathcal{F}\right)  = \{ X \in C^{\infty} \left(TM\right) \, , \ X_x \in
			\mathcal{F}_x \quad \forall \, x \in M \}$ is a Lie algebra,
			
			\smallskip
			
			\item[(d)] The ideal $J\left( \mathcal{F}\right) $ of smooth exterior differential forms which
			vanish on $\mathcal{F}$ is stable by exterior differentiation.
		\end{enumerate}
		
	\end{definition}

	\begin{empt}\label{foli_leaf_empt}\cite{connes:ncg94}
		A foliation of $M$ is given by an integrable subbundle $\mathcal{F}$ of $TM$.
		The \textit{leaves} of the foliation $\left(M , \mathcal F\right)$ are the maximal connected
		submanifolds $L$ of $M$ with $T_x (L) = \mathcal{F}_x $, $\forall \, x \in L$,
		and the partition of $M$ in leaves $$M = \cup
		L_{\alpha}\,,\quad\alpha \in X$$ is characterized geometrically by
		its ``local triviality'': every point $x \in M$ has a neighborhood
		$\mathcal U$ and a system of local coordinates
		$(x^j)_{j = 1 , \ldots , \dim V}$ called
		{\it foliation charts}, so
		that the partition of $\mathcal U$ in connected components of
		leaves corresponds to the partition of 
		\begin{equation*}
			{\mathbb
				R}^{\dim M} = {\mathbb R}^{\dim \mathcal F} \times {\mathbb R}^{\text{codim}
				\, \mathcal F}
		\end{equation*}
		in the parallel affine subspaces 
		$
		{\mathbb R}^{\dim \mathcal F}
		\times {\rm pt}$.
		The corresponding foliation will be denoted by
		\begin{equation}\label{fol_chart_eqn}
			\left(\R^n, \mathcal{F}_p \right) 
		\end{equation}
		where $p = \dim   \mathcal{F}_p$.
		To each foliation $\left(M, \mathcal{F}\right)$ is canonically associated a $C^*$- algebra
		$C^*_r (M, ~\mathcal{F})$ which encodes the topology of the space of leaves.  To
		take this into account one first constructs a manifold $\mathcal G$, $\dim
		\, \mathcal G = \dim \,M + \dim \,\mathcal F$. 
	\end{empt}
	\begin{definition}\label{foli_trans_defn}\cite{candel:foliI}
		Let $N \subset M$ be a smooth submanifold. We say that $\sF$ is \textit{transverse} to $N$ (and write $\sF\pitchfork N$) if, for each leaf $L$ of $\sF$ and each point $x \in L\cap N$, $T_x\left(L \right)$ ans $T_x\left(N \right)$ together span $T_x\left( M\right)$. At the other extreme At the other extreme, we say that $\sF$ is tangent to $N$ if, for each leaf $L$ of $\sF$, either $L \cap N = \emptyset$ or $L \subset N$.
	\end{definition}
	The symbol $\mathbb{F}^p$ denotes either the full Euclidean space $\R^p$ or Euclidean half space $\mathbb{H}^p = \left\{\left.\left(x_1,..., x_n \right) \in \mathbb{R}^p\right| x_1 \le 0 \right\}$.
	\begin{definition}\label{foli_rect_defn}\cite{candel:foliI}
		A rectangular neighborhood in $\mathbb{F}^n$ is an open subset of the form $B = J_1\times...\times J_n$, where each $J_j$ is a (possibly unbounded) relatively open interval in the $j^{\text{th}}$ coordinate axis. If $J_1$ is of the form $\left( a,0\right]$, we say that $B$ has boundary $\partial B\left\{\left(0, x_2,..., x_n \right)\right\}\subset B$.	
	\end{definition}
	\begin{definition}\label{foli_chart_defn}\cite{candel:foliI}
		Let $M$ be an $n$-manifold. A \textit{foliated chart} on $M$ of codimension $q$ is a pair $\left(\sU, \varphi)\right)$, where $\sU\subset M$ is open and $\varphi : \sU \xrightarrow{\approx} B_\tau\times B_\pitchfork$ is a diffeomorphism, $B_\pitchfork$ being a rectangular neighborhood in $\mathbb{F}^q$ and $B_\tau$ a rectangular neighborhood in $\mathbb{F}^{n-q}$. The set $P_y = \varphi^{-1}\left(B_\tau \times \left\{y\right\} \right)$ , where $y \in B_\pitchfork$, is called a \textit{plaque} of this foliated chart. For each $x \in B_\tau$, the set  $S_x=\varphi^{-1}\left(\left\{x\right\} \times B_\pitchfork \right)$  is called a \textit{transversal} of the foliated chart. The set $\partial_{\tau}\sU = \varphi^{-1}\left(B_\tau \times \left(\partial B_\pitchfork \right)  \right)$ is called the \textit{tangential boundary} of $\sU$ and $\partial_{\pitchfork}\sU = \varphi^{-1}\left(\partial \left( B_\tau\right)  \times \partial B_\pitchfork \right)$ is called the \textit{transverse boundary} of $\sU$.
	\end{definition}
	
	\begin{definition}\label{foliated_manifold_defn}\cite{candel:foliI}
		Let $M$ be an $n$-manifold, possibly with boundary and corners, and let $\sF= \left\{L_\la\right\}_{\la \in \La}$ be a decomposition  of $M$ into connected, topologically immersed submanifolds of dimension $k=n-q$. Suppose that $M$ admits an atlas $\left\{\sU_\a \right\}_{\a \in \mathfrak A}$ of foliated charts of codimension $q$ such that, for each $\a \in \mathscr A$ and each $\la \in \La$, $L_\la \cap \sU_\a$ is a union of plaques. Then $\sF$ is said to be a \textit{foliation} of $M$ of codimension $q$ (and dimension $k$) and $\left\{\sU_\a \right\}_{\a \in \mathscr A}$ l is called a \textit{foliated atlas} associated to $\sF$. Each $L_x$ is called a leaf of the foliation and the pair $\left(M, \sF \right)$  is called a \textit{foliated manifold}. If the foliated atlas is of class $C^r$ ($0 \le r \le \infty$ or $r=\om$), then the foliation $\sF$ and the foliated manifold $\left(M, \sF \right)$. is said to be \textit{of class} $C^r$.
	\end{definition}
	\begin{definition}\label{fol_res_defn}\cite{candel:foliI}
		If $\left(M,\mathcal F \right)$ is a foliation and $\mathcal{U} \subset M$ be an open subset $\mathcal F|_{\mathcal{U}}$ is the restriction of $\mathcal F$ on ${\mathcal{U}}$ then we say that  $\left(\mathcal{U},\mathcal F|_{\mathcal{U}}\right) $ is the \textit{restriction}   of $\left(M,\mathcal F \right)$ \textit{to}  $\mathcal{U}$. (cf. \cite{connes:ncg94})
	\end{definition}
	\begin{definition}\label{foli_atlas_defn}\cite{candel:foliI}
		A \textit{foliated atlas} of codimension $q$ and class $C^r$ on the $n$-manifold $M$ is a $C^r$-atlas $\mathfrak{A}\bydef\left\{\sU_\a \right\}_{\a \in \mathscr A}$  of foliated charts of codimension $q$ which are \textit{coherently foliated} in the sense that, whenever $P$ and $Q$ are plaques in distinct charts of $\mathfrak{A}$, then $P\cap Q$ is open both in $P$ and $Q$. 
	\end{definition}
	\begin{definition}\label{foli_coh_atlas_defn}\cite{candel:foliI}
		Two foliated atlases  and $\mathfrak{A}$ on $\mathfrak{A}'$ of the same codimension and smoothness class $C^r$ are \textit{coherent} ($\mathfrak{A}\approx\mathfrak{A}'$) if $\mathfrak{A}\cup\mathfrak{A}'$ is a foliated $C^*$-atlas.
	\end{definition}
	\begin{lemma}\label{foli_coh_atlas_eq_lem}\cite{candel:foliI}
		Coherence of foliated atlases is an equivalence relation.
	\end{lemma}
	\begin{lemma}\label{foli_coh_atlas_ass_lem}
		Let  $\mathfrak{A}$ and $\mathfrak{A}'$ be foliated atlases on $M$ and suppose that $\mathfrak{A}$ is associated to a foliation $\sF$. Then $\mathfrak{A}$ and $\mathfrak{A}'$ are coherent if and only if $\mathfrak{A}'$ is also associated to $\sF$.
	\end{lemma}
	\begin{definition}\label{foli_reg_atlas_defn}\cite{candel:foliI}
		A foliated atlas $\mathfrak{A}\bydef\left\{\sU_\a \right\}_{\a \in \mathscr A}$ of class $C^r$ is said to be \textit{regular} if
		\begin{enumerate}
			\item [(a)] For each $\al \in \mathscr A$, the closure $\overline{\sU}_\al$ of $\sU_\al$ is a compact subset of a foliated chart  $\left\{\sV_\a \right\}$ and $\varphi_\a = \psi|_{\sU_\a }$.
			\item[(b)] The cover $\left\{\sU_\a \right\}$ is locally finite.
			\item[(c)] if $\sU_\a$ and $\sU_\bt$ are elements of $\mathfrak{A}$, then the interior of each closed plaque $P \in \overline \sU_\a$ meets at most one plaque in $\overline \sU_\bt$.
		\end{enumerate}
	\end{definition}
	\begin{lemma}\label{foli_reg_atlas_ref_lem}\cite{candel:foliI}
		Every foliated atlas has a coherent refinement that is regular.
	\end{lemma}
	\begin{thm}\label{foli_thm}\cite{candel:foliI}
		The correspondence between foliations on $M$ and their associated foliated atlases induces a one-to-one correspondence between the set of foliations on $M$.
	\end{thm}
	
	We now have an alternative definition of the term "foliation". 
	\begin{defn}\label{foli_alt_defn}\cite{candel:foliI}
		A \textit{foliation} $\sF$ of codimension $q$ and class $C^r$ on $M$ is a coherence class of foliated atlases of codimension $q$ and class $C^r$ on $M$.
	\end{defn}
	By Zorn's lemma \ref{zorn_thm}, it is obvious that every coherence class of foliated atlases contains a unique maximal foliated atlas. 
	\begin{defn}\label{foli_max_defn}\cite{candel:foliI}
		A \textit{foliation of codimension} $q$ and class $C^r$ on $M$ is a maximal foliated $C^r$-atlas of codimension $q$ on $M$.
	\end{defn}
	\begin{empt}\label{foli_topology_empt}\cite{candel:foliII}
	We construct a topology for $\G = \G\left( M, \F\right)$  as follows.
	Let $\ga \in \G$ be an element with $s\left(\ga\right)= x$ and $r\left(\ga\right)= x$ in the leaf
	$L$. Represent $\ga$ by an immersion $\a \left[0, 1\right]\to L$. Embed $\left[0, 1\right]$ in $\R^{\dim L}$
	as a closed straight line segment $I$ from the point $p_0$ to the point $p_1$ in
	$\R^{\dim L}$. Then, for some $\eps > 0$, $\a$ can be extended to an immersion $f$ of the
	open set $B \bydef  \left\{\left.w \in \R^{\dim L} |\right| \left\|w - I \right\| , \eps \right\}$ so that $f|_{B_\eps\left(p_0 \right) }$, respectively
$f|_{B_\eps\left(p_1 \right) }$, is a coordinate chart around $x$, respectively $y$, in $L$. 
There is a space $Z$
	homeomorphic to a transversal through $x$ in $M$ so that $f$ can be extended
	to an immersion $F$ of the trivial foliated space $B \times Z$ into $M$, and so that
	$F|_{B_\eps\left(p_0 \right)}\times Z$, respectively 	$F|_{B_\eps\left(p_1 \right)}\times Z$, is a foliation chart around $x$, respectively
	$y$, in $M$. Let $\a_{x,w}$ denote the straight path from $0$ to $w$ in $B_\eps\left(p_0 \right)$),
	and let $\a_{y,w}$ be the straight path from $p_1$ to $w$ in $B_\eps\left(p_1 \right)$. The collection of
	elements $\ga \in G$ that have a representative of the form
	$$
	t \in I \mapsto  F\left(\a_{x,u} \#\a\#  \a_{x,u}\left( t\right), z  \right) 	\in M
	$$
	for some $u \in B_\eps\left(p_0 \right)$, $v \in B_\eps\left(p_1 \right)$ and $z \in Z$, is a neighborhood $U\left( \ga\right)$  of $\ga$ in
	$\G$. ($\a\#\b$denotes the path $\a$ followed by the path $\b$.) Note that
	this makes sense because an element of G can be represented by at most one
	such path. This neighborhood $\sU$ is diffeomorphic to $B_\eps\left(p_0 \right)\times B_\eps\left(p_1 \right)\times Z$
	(as a foliated space with leaves $B_\eps\left(p_0 \right)\times B_\eps\left(p_1 \right)\times {z}$).
	\end{empt}

	\begin{example}\label{foli_fibration_exm}\cite{candel:foliI}
Let $M$ be an $n$-manifold, $B$ a manifold of dimension $n-k$, $\partial M = \emptyset = \partial B$, and let  $\pi: M \to B$ be a smooth \textit{fiber bundle}. That is there is a $k$-manifold $F$ and, for each $x \in  B$, an open neighborhood $\sU$ of $x$ in $B$ and a commutative diagram
\\
\begin{tikzcd}
\pi^{-1}\left(\sU \right) \arrow[r, "\varphi"] \arrow[d, "\pi"] & \sU\times F  \arrow[d, "p"]\\
\sU\arrow[r, "\Id"]& \sU
\end{tikzcd}
\\ 
with $\varphi$ a diffeomorphism  and $p$ the canonical projection 	onto the first factor. Each subspace $\pi^{-1}\left(x \right)$, $x\in B$, is clearly an embedded $k$-manifold diffeomorphic to $F$. The manifold is the \text{base} of the fiber bundle and $F$ is the \textit{fiber}. The \textit{total space} of the bundle is $M$.
\end{example}

	\begin{definition}\label{foli_fibration_comes_defn}\cite{candel:foliI} 
Under the hypotheses of the Example \ref{foli_fibration_exm} we say that the foliated space $\left(M, \F \right)$  \textit{comes from the fibration} $\pi: M \to B$.
\end{definition}
	
	\begin{empt}\label{foli_graph_empt}\cite{candel:foliII}
		Let  $\Pi\left( M,\sF\right)$ be the space of paths on leaves, that is, maps $\a : [0,1] \to M$ that are continuous with respect to the leaf topology on $M$. For such a path  let $s\left(\a \right) = \a\left( 0\right)$  be its source or initial point and let  $r\left(\a \right) = \a\left( 1\right)$ be its range or terminal point. The space $\Pi\left( M,\sF\right)$ has a partially defined multiplication: the product $\a\cdot \bt$ of two elements $\a$ and $\bt$ is defined if the terminal point of $\bt$ is the initial point of $\a$, and the result is the path $\bt$ followed by the path $\a$. (Note that this is the opposite to the usual composition of paths  $\al\#\bt = \bt \cdot \a$ used in defining the fundamental group of a space.)
	\end{empt}
	\begin{definition}\label{foli_path_space_defn}\cite{candel:foliI}
		Under the hypotheses of \ref{foli_graph_empt} we say that the topological space $\Pi\left( M,\sF\right)$ is the \textit{space of path on leaves}.
	\end{definition}
	\begin{defn}\label{foli_graph_defn}\cite{candel:foliII}
		The \textit{graph}, or \textit{holonomy groupoid}, of the foliated space $\left( M,\sF\right)$  is the quotient space of $\Pi\left( M,\sF\right)$ by the equivalence relation that identifies two paths $\a$ and $\bt$ if they have the same initial and terminal points, and the loop $\a \cdot \bt$ has trivial germinal holonomy.
		The graph of $\left( M,\sF\right)$ will be denoted by $\G\left(M, \sF\right)$, or simply by  $\G\left( M\right)$ or by $\G$ when all other variables are understood.
	\end{defn}
	\begin{remark}\label{foli_graph_rem}
		The {holonomy groupoid} is a locally compact topological  groupoid (cf. Definitions \ref{groupoid_defn} and \ref{groupoid_topological_defn}).
	\end{remark}
	\begin{remark}
		There is the natural surjective continuous map
		\be\label{foli_cov_map_eqn}
		\Phi : \Pi\left( M,\sF\right)\to \G\left( M,\sF\right)
		\ee
		from the space of path on leaves to the foliation graph.
	\end{remark}
	\begin{proposition}\label{foli_chart_prop}\cite{candel:foliI}
		Let $\mathfrak{A}= \left\{\sU_\iota\right\}$ be a regular foliated atlas of $M$. For each finite sequence of indices $\left\{\a_0 ,...,\a_k\right\}$, the product
		\be	\label{foli_chart_eqn}
		\sV_{   \a} =	 \G\left(\sU_{\iota_0} \right)~...~\G\left(\sU_{\iota_k} \right) \in \G\left(M, \sF\right) \quad \a = \left({\iota_0},...,{\iota_k}\right)
		\ee
		is either empty or a foliated chart for the graph $\G$. The collection of all such finite products is a covering of $\G$ by foliated charts. 
	\end{proposition}
	\begin{theorem}\label{foli_graph_thm}\cite{candel:foliII}
		The graph $\G$ of $\left( M,\sF\right)$ is a groupoid with unit space $\G_0 = M$, and this algebraic structure is compatible with a foliated structure on $\G$ and $M$. Furthermore, the following properties hold.
		\begin{enumerate}
			\item [(i)] The range and source maps $r, s : \G \to M$ are topological submersions. 
			\item[(ii)] The inclusion of the unit space $M \to\G$ is a smooth map. 
			\item[(iii)] The product map $\G\times_M \G \to G$, given by $\left( \ga_1 , \ga_2\right) \mapsto\ga_1 \cdot \ga_2$, is smooth.
			\item[(iv)]  There is an involution $j: \G \to \G$, given by $j\left( \ga\right) = \ga^{-1}$, which is a diffeomorphism of $\G$, sends each leaf to itself, and exchanges the foliations given by the range. 
		\end{enumerate}
		
	\end{theorem}
	
	Above definitions refines the equivalence relation coming from
	the partition of $M$ in leaves $M = \cup L_{\alpha}$. 
	An element $\gamma$ of $\mathcal G$ is given by two points $x = s(\gamma)$,
	$y = r(\gamma)$ of $M$ together with an equivalence class of smooth
	paths: $\gamma (t)\in M$, $t \in [0,1]$; $\gamma (0) = x$, $\gamma
	(1) = y$, tangent to the bundle $\mathcal{F}$ ( i.e. with $\dot\gamma (t)
	\in \mathcal{F}_{\gamma (t)}$, $\forall \, t \in {\mathbb R}$) up to the
	following equivalence: $\gamma_1$ and $\gamma_2$ are equivalent if and only if
	the {\it my} of the path $\gamma_2 \circ \gamma_1^{-1}$ at the
	point $x$ is the {\it identity}. The graph $\mathcal G$ has an obvious
	composition law. For $\gamma , \gamma' \in G$, the composition
	$\gamma \circ \gamma'$ makes sense if $s(\gamma) = r(\gamma')$. If
	the leaf $L$ which contains both $x$ and $y$ has no my, then
	the class in $\mathcal G$ of the path $\gamma (t)$ only depends on the pair
	$(y,x)$. In general, if one fixes $x = s(\gamma)$, the map from $\mathcal G_x
	= \{ \gamma , s(\gamma) = x \}$ to the leaf $L$ through $x$, given
	by $\gamma \in \mathcal G_x \mapsto y = r(\gamma)$, is the my covering
	of $L$.
	Both maps $r$ and $s$ from the manifold $\mathcal G$ to $M$ are smooth
	submersions and the map $(r,s)$ to $M \times M$ is an immersion
	whose image in $M \times M$ is the (often singular) subset
	\begin{equation*}\label{subset}
		\{ (y,x)\in M \times M: \, \text{ $y$ and $x$ are on the same leaf}\}.
	\end{equation*}
	For
	$x\in M$ one lets $\Omega_x^{1/2}$ be the one dimensional complex
	vector space of maps from the exterior power $\wedge^k \,  \mathcal{F}_x$, $k =
	\dim F$, to ${\mathbb C}$ such that
	$$
	\rho \, (\lambda \, v) = \vert \lambda \vert^{1/2} \, \rho \, (v)
	\qquad \forall \, v \in \wedge^k \,  \mathcal{F}_x \, , \quad \forall \,
	\lambda \in {\mathbb R} \, .
	$$
	Then, for $\gamma \in\mathcal G$, one can identify $\Omega_{\gamma}^{1/2}$ with the one
	dimensional complex vector space $\Omega_y^{1/2} \otimes
	\Omega_x^{1/2}$, where $\gamma : x \to y$. In other words
	\be\label{foli_om_g_eqn}
	\Omega_{\mathcal G}^{1/2}=\, r^*(\Omega_M^{1/2})\otimes s^*(\Omega_M^{1/2})\,.
	\ee
	
	
	\begin{empt}\label{foli_sc_haus_empt}\cite{candel:foliII}
	The  groupoid of a foliated space all leaves of which are simply connected is Hausdorff.
	\end{empt}
	
	\begin{exercise}\label{foli_haus_exer}\cite{candel:foliII}
	Prove or decide the following.
\begin{enumerate}
	\item The graph of a foliated space all leaves of which are simply connected is Hausdorff.
	\item The graph of a foliated space all leaves of which have trivial holonomy is Hausdorff.
\end{enumerate}	
\end{exercise}

	\begin{definition}\label{foli_pseudo_defn}\cite{candel:foliI}.
		Let $N$ be a $q$-manifold. A $C^r$ pseudogroup $\Ga$ on $N$ is a collection of $C^r$ diffeomorphisms $h : D(h) \xrightarrow{\approx} R(h)$ between open subsets of $N$ satisfying the following axioms. 
		\begin{enumerate}
			\item  If $g, h \in \Ga$  and $R(h) \subset G(g)$, then $g \circ h \in \Ga$
			\item	 If $h \in \Ga$, then $h^{-1} \in \Ga$. 
			\item $\Id_N \in \Ga$. 
			\item If $h \in \Ga$ and $W \subset  D(h)$ is an open subset, then $\left.h\right|_W \in \Ga$. 
			\item If $h:D(h) \xrightarrow{\approx} R(h)$ is a $C^r$ diffeomorphism between open subsets of $N$ and if, for each $w \in  D(h)$, there is a neighborhood $W$ of $w\in  D(h)$ such that $\left.h\right|_W \in \Ga$, then $h \in \Ga$. 
		\end{enumerate}
		If $\Ga' \subset \Ga$ is also a pseudogroup, it is called a subpseudogroup of $\Ga$.
	\end{definition}

	\begin{remark}\cite{candel:foliI}
		Any pseudogroup is a groupoid (cf. Definition \ref{groupoid_defn}).
	\end{remark}
\begin{remark}\cite{candel:foliI}\label{foli_pseudo_rem}
In the case of general foliations, the total my group of a foliated bun dle must be replaced by a local analogue called the \textit{my pseudogroup}.
\end{remark}
	\begin{remark}\label{foli_groupoid_n_red_defn}\cite{candel:foliI}
		If $\left(M, \sF\right)$ is a foliated manifold and $N$ is a tranversal then
		\be\label{foli_gnn_eqn}
		\G^N_N \bydef \left\{\left. \ga \in \G\left(M, \sF\right)\right| s\left(\ga\right), r\left(\ga\right)\in N\right\}
		\ee
		is	a pseudogroup.

	\end{remark}
	\begin{definition}\cite{candel:foliI}
	In the above situation we say that $\G^N_N$ is a \textit{reduced groupoid}.
	\end{definition}

	\section{Operator algebras of foliations}\label{foli_alg_subsec}
	\paragraph*{}
	Here I follow to \cite{candel:foliII,connes:ncg94}.  Since the bundle $\Om^{1/2}$ is trivial (because $\G\left(M, \sF\right)$ admits partitions of unity), a choice of an everywhere positive density $\nu$ allows us to identify $\Ga_c\left(\G\left(M, \sF\right),\Om^{1/2}  \right)$  with $\Coo_c\left( \G\left(M, \sF\right)\right)$. 
	The definition of foliated space makes sense even when the underlying topological space fails to satisfy the Hausdorff separation axiom. Non-Hausdorff spaces appear naturally in the theory of foliations. In a graph of a foliated space is not necessary Hausdorff. It will be necessary to use functions with compact support on such spaces. However, a non-Hausdorff space may not have sufficiently many such functions, the basic reason being that compact subsets of a Hausdorff space are not necessarily closed. The non-Hausdorff spaces that will appear here have a particularly simple local structure, and even when it is possible to construct appropriate functions using this local structure, the standard operation of “extension by 0” of local objects to the full space does not pro duce continuous functions. M. Crainic and I. Moerdijk \cite{cra_moe:nhaus} proposed a very natural way of dealing with this problem, and this preliminary section describes it. (That paper develops an extended sheaf  theory for non-Hausdorff manifolds (cf. Appendix \ref{sheaves_nh_sec})).
	Here $\sX$ will denote a separable topological space having the structure of a foliated space, but it is not required that the topology be Hausdorff. It is only required that $\sX$ can be covered by countably many open sets homeomorphic to a product $\D \times \mathcal Z$, where $D$ is an open ball in Euclidean $n$-space and $\mathcal Z$ is a separable locally compact Hausdorff space. Let $\mathcal C\infty$ denote the structure sheaf  of the foliated space $\sX$, that is, the sheaf  of smooth functions on $\sX$. Let $\A$ be a sheaf of $\mathcal C\infty$-modules over $\sX$, for instance, the sheaf  of differential forms or other sheaves of smooth sections of foliated vector bundles. For such a sheaf  $\A$ over $\sX$, let $\A$ denote its Godement resolution: $A'\left(\sU\right)$ is the set of all sections (continuous or not) of $\A$ over $\sU\subset\sX$. For a Hausdorff open subset $\mathcal W$ of $\sX$, let $\Ga_c\left( \mathcal W, \A\right) $ denote the set of continuous compactly supported sections of $\A$ over $\mathcal W$. If $\mathcal W\subset\sU$, then “extension by 0” induces a well defined homomorphism $\Ga_c\left( \mathcal W, \A\right)\to \A'\left(\sU \right)$. For an open subset $\sU$ of $\sX$, let $\Ga_c\left( \mathcal U, \A\right)$ denote the image of the homomorphism $\oplus \Ga_c\left( \mathcal W, \A\right)\to \A'\left(\sU\right)$ (cf. Definition \ref{nh_csoft_gc_defn}). From the equation \eqref{sheaf_inc_eqn} it follows that there is the inclusion
	\be\label{foli_incc_eqn}
\Ga_c\left(\mathcal W, \A \right) \hookto \Ga_c\left(\mathcal U, \A \right).
\ee
Let $\G\bydef \G\left(M, \sF \right)$ be a foliation chart. 	The bundle $\Omega_M^{1/2}$ is trivial on $M$, and we
	could choose once and for all  a trivialization $\nu$ turning
	elements of $\Ga_c \left(\mathcal G , \Omega_{\mathcal G}^{1/2}\right)$ into functions.
	Let us
	however stress that the usage of half densities makes all the
	construction completely canonical.
	For $f,g \in \Ga_c \left(\mathcal G , \Omega_{\mathcal G}^{1/2}\right)$, the convolution
	product $f * g$ is defined by the equality
	\be\label{foli_prod_eqn}
	(f * g) (\gamma) = \int_{\gamma_1 \circ \gamma_2 = \gamma}
	f(\gamma_1) \, g(\gamma_2) \, .
	\ee
	This makes sense because, for fixed $\gamma : x \to y$ and fixing $v_x
	\in \wedge^k \,  \mathcal{F}_x$ and $v_y \in \wedge^k \,  \mathcal{F}_y$, the product
	$f(\gamma_1) \, g(\gamma_1^{-1} \gamma)$ defines a $1$-density on
	$G^y = \{ \gamma_1 \in G , \, r (\gamma_1) = y \}$, which is smooth
	with compact support (it vanishes if $\gamma_1 \notin\supp f$),
	and hence can be integrated over $G^y$ to give a scalar, namely $(f * g)
	(\gamma)$ evaluated on $v_x , v_y$.
	The $*$ operation is defined by $f^* (\gamma) =
	\overline{f(\gamma^{-1})}$,  i.e. if $\gamma : x \to y$ and
	$v_x \in \wedge^k \, \mathcal{F}_x$, $v_y \in \wedge^k \, \mathcal{F}_y$ then $f^*
	(\gamma)$ evaluated on $v_x , v_y$ is equal to
	$\overline{f(\gamma^{-1})}$ evaluated on $v_y , v_x$. We thus get a
	$*$-algebra $\Ga_c \left(\mathcal G , \Omega_{\mathcal G}^{1/2}\right)$. 
	where $\xi$ is a square integrable half density on $\mathcal G_x$. 
	For each leaf $L$ of
	$\left(M, \mathcal{F}\right)$ one has a natural representation of this $*$-algebra on the
	$L^2$ space of the my covering $\tilde L$ of $L$. Fixing a
	base point $x \in L$, one identifies $\tilde L$ with $\mathcal G_x = \{
	\gamma , s(\gamma) = x \}$ and defines
	\begin{equation}\label{foli_repr_eqn}
		(\rho_x (f) \, \xi) \, (\gamma) = \int_{\gamma_1 \circ \gamma_2 =
			\gamma} f(\gamma_1) \, \xi (\gamma_2) \qquad \forall \, \xi \in L^2
		(\mathcal G_x),\
	\end{equation}

	Given
	$\gamma : x \to y$ one has a natural isometry of $L^2 (\mathcal G_x)$ on $L^2
	(G_y)$ which transforms the representation $\rho_x$ in $\rho_y$.
	\begin{lemma}\cite{candel:foliII} 
		If $f_1 \in \Ga_c\left(\sU_{   \a_1},\Om^{1/2} \right)$ and $f_2 \in \Ga_c\left(\sU_{   \a_2},\Om^{1/2} \right)$ then their convolution is a well-defined element $f_1*f_2 \in \Ga_c\left(\sU_{   \a_1}\cdot\sU_{   \a_2},\Om^{1/2} \right)$
	\end{lemma}
	
	\begin{proposition}\label{foli_repr_prop}\cite{candel:foliII} 
		If $\sV \subset \G$ is a foliated chart for the graph of $\left(M, \sF\right)$ and $f \in \Ga_c\left(\sV, \Om^{1/2}\right)$ , then $\rho_x\left( f\right)$, given by \eqref{foli_repr_eqn}, is a bounded integral operator on $L^2\left(\G_x \right)$.
	\end{proposition}
	\begin{empt}\cite{candel:foliII} 
		The space of compactly supported half-densities on $\G$ is taken as given by the exact sequence 
		\be\label{foli_ga_p_eqn}
		\bigoplus_{\a_0\a_1}\Ga_c\left(\sU_{   \a_0\a_1}, \Om^{1/2} \right) \to \bigoplus_{\a_0}\Ga_c\left(\sU_{   \a_0},\Om^{1/2} \right) \xrightarrow{\Ga_\oplus}  \Ga_c\left(\G,\Om^{1/2}\right) 
		\ee
		associated to a regular cover for $\left((M, \sF)\right)$ as above. The first step for defining a convolution is to do it at the level of $\bigoplus_{\a_0}\Ga_c\left(\sU_{   \a_0}\Om^{1/2} \right)$, as the following lemma indicates. 
	\end{empt}

	\begin{defn}\label{foli_red_defn}\cite{candel:foliII} 
		The \textit{reduced} $C^*$-algebra of the foliated space $\left(M,\sF\right)$ is the completion of $\Ga_c\left( \G,\Om^{1/2}\right)$ with respect to the pseudonorm \be\label{foli_pseudo_norm_eqn}
		\left\|f \right\| =\sup_{x \in M}\left\|  \rho_x\left(f\right)\right\|
		\ee where $\rho_x$ is given by  \eqref{foli_repr_eqn}.
		This $C^*$-algebra is denoted by $C^*_r\left(M,\sF\right)$.
	\end{defn} 
	An obvious consequence of the construction of 
	$C^*_r\left(M,\sF\right)$ is the following. 
	
	\begin{cor}\label{foli_cov_alg_cor}\cite{candel:foliII} 
		Let $M$ be a foliated space and let $\mathfrak{A}$ be a regular cover by foliated charts. Then the algebra generated by the convolution algebras $\Ga_c\left( \G\left(\sU \right), \Om^{1/2}\right)$, $~\sU\in\mathfrak{A}$ is dense in  $C^*_r\left(M,\sF\right)$.
	\end{cor}
	\begin{definition}\label{foli_full_defn}\cite{candel:foliII}
The \textit{full} $C^*$-algebra of the foliated space $(M,\F)$, denoted
by $C^*(M,\F)$, is the completion of $\Ga_c\left(\G,\Om^{1/2}\right)$ with respect to the
pseudonorm
\be\label{groupoid_full_norm_eqn}
\left\|f \right\|_f =\sup_{\pi}\left\|  \rho_x\left(f\right)\right\|
\ee
where $\pi$ runs through all the involutive representations of $\Ga_c\left(\G,\Om^{1/2}\right)$ on a
separable Hilbert space whose restrictions to the graph $\G\left( \sU\right)$  of each foliated
chart $\sU$ for $(M,\F)$ are weakly continuous for the inductive limit topology
on $\Ga_c\left(\G,\Om^{1/2}\right)$.	
\end{definition}
	\begin{empt}\label{foli_res_inc_empt}\cite{candel:foliII} 
		Let $\left(M,\sF\right)$ be an arbitrary foliated space and let $\sU\subset  M$ be an open subset.  Then $\left(\sU,\sF|_\sU\right)$ is a foliated space and the inclusion  $\sU\hookto  M$ induces a homomorphism of groupoids $\G\left(\sU \right)\hookto \G$ , hence a mapping
		\be\label{foli_inc_gc_eqn}
		j_\sU : \Ga_c\left(  \G\left(\sU \right) ,\Om^{1/2}\right)  \hookto \Ga_c\left( \G\left( M\right) ,\Om^{1/2}\right)
		\ee
		that is an injective homomorphism of involutive algebras. 
	\end{empt}
	\begin{prop}\label{foli_res_inc_prop} \cite{candel:foliII} 
		Let $\sU$ be an open subset of the foliated space $M$. Then the inclusion $\sU \hookto M$ induces an isometry of $C^*_r\left(\sU,\sF|_\sU\right)$  into $C^*_r\left(M,\sF\right)$.
	\end{prop}
	\begin{lemma}\label{foli_leaf_lem}\cite{candel:foliII} 
		Each element $\ga \in \G$ induces a unitary operator $\rho_\ga : L^2\left(s\left( \ga\right)\right)   \xrightarrow{\approx}  L^2\left(r\left( \ga\right)\right) $   that conjugates the operators $\rho_{s\left(\ga \right)} \left(f \right) $ and $\rho_{r\left(\ga \right)} \left(f \right)$. In particular, the norm of $\rho_x\left(f \right)$  is independent of the point in the leaf through $x$.
	\end{lemma}
	
	\begin{lem}\label{foli_point_lem}\cite{candel:foliII} 
		If $f \in \Ga^\infty_c\left(\G, \Om^{1/2}\right)$ does not evaluate to zero at each $\ga \in \G$, then there exists a point $x$ in $M$ such that $\rho_x\left( f\right)  \neq  0$.
	\end{lem}
	\begin{definition}\label{foli_fibration_defn}\cite{candel:foliI}
		A foliated space $\left(M, \sF\right)$ is a \textit{fibration} if for any $x$ there is an open transversal $N$ such that $x\in N$ and for every  leaf $L$ of $\left(M, \sF\right)$ the intersection $L\cap  N$ contains no more then one point.
	\end{definition}
	\begin{prop}\label{foli_one_leaf_prop}
		The reduced $C^*$-algebra of a foliated space $M$ consisting of exactly one leaf is the algebra $\K\left( L^2\left(M \right) \right)$  of compact operators on $L^2\left(M \right)$.
	\end{prop}

	\begin{prop}\label{foli_tens_comp_prop}\cite{candel:foliII}
		The reduced  $C^*$-algebra $C^*_r\left(\mathcal N\times \mathcal Z \right)$ of the trivial foliated space $\mathcal N \times \mathcal Z$ is the tensor product $\K\otimes C_0\left(\mathcal Z\right)$, where $\K$ is the algebra of compact operators on $L^2\left(\mathcal N \right)$  and $ C_0\left(\mathcal Z\right)$ is the space of continuous functions on $\mathcal Z$ that vanish at infinity.
	\end{prop}
	\begin{corollary}\label{foli_tens_comp_cor}\cite{candel:foliII}
The full $C^*$-algebra of a trivial foliated space $\mathcal N \times \sZ$ is
the tensor product $\K\left( L^2\left( \mathcal N\right) \right) \otimes C_0\left( \sZ\right)$. 	
\end{corollary}
	\begin{thm}\label{foli_bundle_thm}\cite{candel:foliII}
		Assume that $\left(M, \sF\right)$ is given by the fibers of a fibration $p : M \to B$ with fiber $F$. Then $C^*\left(M, \sF\right)$ is isomorphic to $C_0\left(B\right)\otimes \K\left(L^2\left(\F\right)\right)$.
	\end{thm}
\begin{proof}	
	It is clear  that the $C^*$-algebra $C^*\left(M, \F \right)$  is as follows. There is a locally trivial
	bundle $E$ over $B$ with fiber $\K = \K\left(L^2\left(\F \right)  \right)$  so that  $C^*_r\left(M, \sF\right)$ is isomorphic
	to the algebra of continuous sections $s : B \to  E$ vanishing at infinity on $B$
	with the norm $\left\| s\right\| = \sup_{x\in B}\left\|s\left( x\right) \right\|$
 The structure of this bundle $E$ will
	now be examined.

	The bundle $E$ is associated to a locally trivial Hilbert bundle $H\to B$
	over $B$ with fiber $L^2\left(\sF \right)$  and transition functions given by the cocycle that
	defines $M \to B$. Indeed, if $\left\{B_\iota\right\}$ is a cover of $B$ that trivializes $M$, then
	there is a cocycle
	$$
	g_{\iota j} : B_\iota \cap  B_j \to \Diff\left(\sF \right) 
	$$
	that is continuous for the topology on $\Diff\left(\sF \right)$ of convergence on compact
	subsets. By the continuity of the representation $\Diff\left(\sF \right)\to U\left(L^2\left(\F \right)  \right)$,  $\left\{	g_{\iota j} \right\}$
	induces a cocycle for the Hilbert bundle $H$.

	Since the group of unitary operators on an infinite dimensional Hilbert
	space is contractible, the Hilbert bundle $H$ is, in
	fact, trivial, which means that, perhaps after passing to a refinement of the
	covering, there are continuous maps
	$$
	h_\iota : B_\iota  \to  U\left(L^2\left(\F \right)  \right)
	$$
	such that $g_{\iota j}= h_\iota^* h_j$.
 This in turn implies that the operator bundle $E\to B$ 
	is, in fact, trivial; that is, it is isomorphic to $C_0\left(B\right)\otimes \K\left(L^2\left(\F\right)\right)$, and thus the $C^*$-
	algebra $C^*_r\left(M, \sF\right)$ is isomorphic to the tensor product $C_0\left(B\right)\otimes \K\left(L^2\left(\F\right)\right)$.
\end{proof}	
		\begin{remark}\label{foli_bundle_rem}
		From the proof of the Theorem \ref{foli_bundle_thm} it follows that the  $C^*_r\left(M, \sF\right)$ is isomorphic to $C_0\left(B\right)\otimes \K\left(L^2\left(\F\right)\right)$.
	\end{remark}
	

	
	\begin{theorem}\label{foli_irred_hol_thm}\cite{candel:foliII}
		Let $(M,\sF)$ be a foliated space and let $x \in M$. Then the representation $\rho_x$ is irreducible if and only if the leaf through $x$ has no holonomy.
	\end{theorem}
	\begin{empt}\label{foli_pseudo_empt}\cite{connes:ncg94}
		If $\G^N_N$ is a  given by \eqref{foli_gnn_eqn} pseudogroup then $\G^N_N$ is a manifold.  One can define a structure of $*$-algebra on $\Coo_c\left( \G^N_N\right)$ such that for any $a, b \in \Coo_c\left( \G^N_N\right)$ one has:
		\be\label{foli_pseudo_eqn}
		\begin{split}
			\left(a\cdot b\right)\left(\ga\right)= \sum_{\ga_1\circ\ga_2=\ga}a\left( \ga_1\right)a\left( \ga_1\right),\\
			a^*\left(\ga\right)= \overline{a\left(\ga^{-1}\right)} 
		\end{split}
		\ee
		Denote by $C^*_r\left(	\G^N_N \right)$ the reduced algebra of $\G^N_N$ (cf. Definition \ref{foli_groupoid_red_defn}). The equation \eqref{foli_pseudo_eqn} is a specialization of the \eqref{groupoid_a_eqn} one. 
	\end{empt} 	
	\begin{definition}\label{foli_complete_defn}\cite{hilsum_scandalis:stab}
	A transversal $\mathcal N$ of a foliated space $\left(M, \F \right)$   is said to be \textit{complete} if $\mathcal N$ meets every leaf. 
\end{definition}
	\begin{remark}
	In \cite{hilsum_scandalis:stab} the French word "fidèle" is used instead the "complete" one. However the word "complete" is used in \cite{candel:foliII}.
\end{remark}
	
	\begin{theorem}\label{foli_mor_thm}\cite{candel:foliII,hilsum_scandalis:stab}
		Let $(M,\F)$ be a foliated space and let $(Z, \mathscr H)$ be the holonomy pseudogroup associated to a complete transversal $Z$. Then the foliation $C^*$-algebra is Morita equivalent to the pseudogroup $C^*$-algebra.
	\end{theorem}
	\begin{remark}\label{foli_pseudo_alg_rem}\cite{hilsum_scandalis:stab}
		A pseudogroup $C^*$-algebra is $C^*_r\left(	\G^N_N \right)$.
	\end{remark}

	The following  lemma is a more strong version of the Theorem \ref{foli_mor_thm}.
	\begin{lemma}\label{foli_stab_lem}\cite{hilsum_scandalis:stab}
		Si $N$ est une transversale fidèle de $M$
		on a: 
	\be\label{foli_stab_eqn}
	C^*_r\left(M, \sF \right)\cong C^*_r\left(\G^N_N \right)\otimes \K
	\ee	
	\end{lemma}
	\begin{remark}\label{foli_stab_rem}
		The English translation of the Theorem \ref{foli_stab_lem} is "if $N$ is a complete transversal then 
	$C^*_r\left(M, \sF \right)\cong C^*_r\left(\G^N_N \right)\otimes \K$".
	\end{remark}
	\begin{empt}\label{foli_stab_empt}
		Besides the Lemma \ref{foli_stab_lem} we need some details of its proof. If $N$ is a complete transversal and
		$$
		E^M_N \bydef \Ga_c\left(r^{-1}\left(M \right), s^*\left(\Om^{1/2} \right)   \right) 
		$$
		then there is a $\Ga_c\left(\G, \Om \right)$-valued product  on $E^M_N$ given by
		$$
		\left\langle \xi, \eta \right\rangle\left( \ga\right) = \sum_{\substack{ s\left(\ga\right)=s\left(\ga'\right)\\r\left(\ga'\right)\in N}}\overline\xi\left(\ga'\circ\ga^{-1}\right)\eta\left(\ga'\right), \quad 
		$$
		If $\E^M_N$ is a completion of $E^M_N$ with respect to the norm $\left\| \xi\right|  \bydef \sqrt{\left\langle \xi, \xi \right\rangle}$ then $\E^M_N$ is $C^*_r\left(M,\sF\right)$-$C^*_r\left(\G^N_N\right)$-imprimitivity bimodule. (cf. Definition \ref{strong_morita_defn}). It follows that
		\be\label{foli_k1_eqn}
		C^*_r\left(M,\sF\right)\cong \K\left(\E^M_N\right).
		\ee
		On the other hand in \cite{hilsum_scandalis:stab} it is proven that there is an isomorphism 
		\be\label{foli_k2_eqn}
		\E^M_N\cong\ell^2\left(C^*_r\left(\G^N_N\right) \right)\cong C^*_r\left(\G^N_N\right)\otimes \ell^2\left( \N\right).  
		\ee
		The combination of equations \eqref{foli_k1_eqn}, \eqref{foli_k2_eqn} yields the Lemma  \ref{foli_stab_lem}.
		Below we briefly remind  the idea described in \cite{hilsum_scandalis:stab} proof of the equation  \ref{foli_k2_eqn}. There is a tubular neighborhood $V_N$ of $N$ which is homeomorphic to $N \times [0,1]^q$. This neighborhood yields the isomorphism
		\be\label{foli_e_sum_eqn}
		\E^M_N\cong C^*_r\left(\G^N_N\right)\otimes L^2\left( [0,1]^q\right) \oplus  \E^W_N
		\ee
		where $\E^W_N$ is a countably generated Hilbert $C^*_r\left(\G^N_N\right)$-module.
		From $ L^2\left( [0,1]^q\right)\cong\ell^2\left(\N\right)$ and the Kasparov stabilization theorem \ref{kasparov_stab_thm}, it follows that\\ $\E^M_N\cong\ell^2\left(C^*_r\left(\G^N_N\right) \right)$. If we  select any
		\be\label{foli_vec_eqn}
		\xi \in L^2\left( [0,1]^q\right)
		\ee
		such that $\left\|\xi\right\|= 1$ then  is an orthogonal basis
		\be\label{foli_basis_eqn}
		\left\{\eta_j\right\}_{j\in \N}\subset \ell^2\left(\N \right), \quad\eta_1\bydef\xi
		\ee
		In this basis the element  $a\otimes \xi\left\rangle \right\langle \xi\in C^*_r\left(\G^N_N\right)\otimes \K$ is represented by the following matrix
		\be\label{foli_mat_eqn}
		a\otimes \xi\left\rangle \right\langle \xi  =  	\begin{pmatrix}
			a& 0 &\ldots \\
			0& 0 &\ldots \\
			\vdots& \vdots &\ddots\\
		\end{pmatrix}.
		\ee
		
	\end{empt}

	\begin{theorem}\cite{candel:foliII, hilsum_scandalis:stab}
	Let $(M, \sF)$ be a foliated space and let $\left(Z, \mathcal H\right)$ be the holonomy
	pseudogroup associated to a complete transversal $Z$. Then the foliation
	$C^*$-algebra is Morita equivalent to the pseudogroup $C^*$-algebra.
	\end{theorem}	
	\begin{example}\label{fol_tor_exm}\emph{Linear foliation on torus}. Here I follow to \cite{candel:foliI,connes:ncg94}.
		Consider a vector field $\tilde{X}$ on $\R^2$ given by
		\[
		\tilde{X}=\alpha\frac{\partial}{\partial
			x}+\beta\frac{\partial}{\partial y}
		\]
		with constant $\alpha$ and $\beta $. Since $\tilde{X}$ is
		invariant under all translations, it determines a vector field $X$
		on the two-dimensional torus ${\T}^2={\R}^2/{\Z}^2$. The vector
		field $X$ determines a foliation $\mathcal{F}$ on ${\T}^2$. The leaves of
		$\mathcal{F}$ are the images of the parallel lines
		$\tilde{L}=\{(x_0+t\alpha, y_0+t\beta): t\in\R\}$ with the slope
		$\theta=\beta/\alpha $ under the projection $\R^2\to \T^2$.
		In the case when $\theta$ is rational, all leaves of $\mathcal{F}$ are
		closed and are circles, and the foliation $\mathcal{F}$ is determined by
		the fibers of a fibration $\T^2\to S^1$. In the case when $\theta$
		is irrational, all leaves of $\mathcal{F}$ are everywhere dense in $\T^2$. We say that  $\left(\T^2, \mathcal{F}_\th \right)$ is  the \textit{Kronecker foliation} $dy = \th dx$ of the 2-torus $ \T^2 \bydef \R^2/\Z^2$ 
	with natural coordinates $\left((x, y \right)\in \R^2$. Here $\th\in (0,1)$ is an irrational number.
	The graph $\G$ of this foliation is the manifold $\G \cong \T\times \R$ with range and source maps  $\G \to \T^2$ given by
	\bean
	r((x, y), t) = (x + t, y + \th t).\\
	s((x, y), t) = (x, y)
	\eean
	and with composition given by $((x, y), t)((x_0, y_0), t_0) = ((x_0, y_0), t + t_0)$ for any pair of	composable elements.
	Every closed geodesic of the at torus $\T^2$ yields a compact transversal. More precisely,
	for each pair $(p, q)$ of relatively prime integers we let $N_{p,q}$ be the submanifold of $\T^2$
	given by
	$$
	N_{p,q} = \left\{(ps, qs) | s\in \R/\Z\right\}
$$
	The graph $\G$ reduced by $N = N_{p,q}$, i.e. $\G^N_N \bydef \left\{\left.\ga \in \G \right| r(\ga)\in N, \quad s(\ga)\in N\right\}$
	the manifold $\G^N_N=\T\times \Z$ with range and source maps given by:
	\bean
	r(x, n) = x + n\th' ;\\
	 s(x, n) = x
	\eean
	where $\th'$ is determined uniquely by any pair $(p_0, q_0)$ of integers such that $pq_0 - p_0q= 1$, $\quad \th' = \frac{p'\th - q'}{p\th - q}$
If $N = N_{0,1}$ then 
\bean
r(x, n) = x + n\th ;\\
s(x, n) = x
\eean
Let $U$ be an element of $C\left(\T\right)$ which comes from
\bean
U:\R \to \C,\\
t \mapsto e^{2\pi i t}.
\eean
If both  $u, v\in C^\infty_c\left(\G^N_N \right)$ are such that
\be\label{foli_uv_eqn}
\begin{split}
u\left(x, n \right) = \begin{cases}
	U\left(x \right) & n = 0\\
	0 & n\neq 0
\end{cases};\\
v\left(x, n \right) = \begin{cases}
	1 & n = 1\\
	0 & n\neq 1
\end{cases},\\
\end{split}
\ee
then from  the equations \eqref{foli_pseudo_eqn}
one has
\bean
u^*=u^{-1},\\
v^*=v^{-1},\\
	vu = e^{2\pi i\th } uv.
\eean
Indeed above equations describe a noncommutative 2-torus $C\left(\T^2_\th\right)$ (cf Definition \ref{nt_defn}). Using this fact one deduce that 
\be\label{foli_nt_eqn}
C^*_r\left( \G^N_N\right) \cong C\left(\T^2_\th\right).
\ee
	\end{example}
	
		\begin{exercise}\cite{candel:foliII}
		Let $\mathcal W \bydef \G\left(\sU_1 \right)....\  \G\left(\sU_1 \right)$. Then the collection of linear
		combinations of elements of the form $f = f_1*...* f_n,~ f_j \in \G\left(\sU_1 \right)$, is dense
		in the space of compactly supported continuous half densities on $\mathcal W$, in the
		inductive limit topology (cf. Definition \ref{top_ind_lim_defn}).	
		\end{exercise}
	
	\subsection{Restriction of foliation}

	\begin{lemma}\label{fol_res_lem}\cite{connes:ncg94}
		If $\sU \subset M$ is an open set and  $\left(\mathcal{U},\mathcal F|_{\mathcal{U}}\right) $ is the {restriction}   of $\left(M,\mathcal F \right)$ {to}  $\mathcal{U}$ then the graph 
		$\G\left(\mathcal{U},\mathcal F|_{\mathcal{U}}\right)$ is an open set in the graph $\G\left(M,\mathcal F \right)$, and the inclusion
		$$
		\Coo\left(\mathcal{U},\mathcal F|_{\mathcal{U}}\right)\hookto\Coo\left(M,\mathcal F \right)
		$$
		extends to an isometric $*$-homomorphism of $C^*$-algebras
		$$
		C^*_r\left(\mathcal{U},\mathcal F|_{\mathcal{U}}\right)\hookto C^*_r\left(M,\mathcal F \right).
		$$
	\end{lemma}
	
	\begin{remark}\label{fol_res_rem}\cite{connes:ncg94}
		This lemma, which is still valid in the non-Hausdorff case \cite{connes:foli_survey}, allows one to reflect
		algebraically the local triviality of the foliation.
		Thus one can cover the manifold $M$ by
		open sets $\sU_\la$ such that $\sF$ restricted to $\sU_\la$ has a Hausdorff space of leaves, $\sV_{\la} = \sU_\la/\sF$. and hence such that the C*-algebras $C^*_r\left(\mathcal{U}_\la,\mathcal F|_{\mathcal{U}_\la}\right)$  are strongly Morita equivalent to the commutative $C^*$-algebras $C_0\left( B_\la\right)$. These subalgebras $C^*$-algebras $C^*_r\left(\mathcal{U}_\la,\mathcal F|_{\mathcal{U}_\la}\right)$ generate $ C^*_r\left(M,\mathcal F \right)$. 
		but of course they fit together in a very complicated way which is related to the global
		properties of the foliation.
	\end{remark}
	
	\subsection{Lifts of foliations}
	\paragraph*{}
	Let $M$ be a smooth manifold and let
	is an  $\mathcal F\subset TM$ be an integrable subbundle.
	If $p:\widetilde M \to M$ is a covering and $\widetilde{ \mathcal F} \subset T\widetilde M$ is the lift of $\mathcal F$ given by a following diagram
	\newline
	\newline
	\begin{tikzpicture}
		\matrix (m) [matrix of math nodes,row sep=3em,column sep=4em,minimum width=2em]
		{
			\widetilde{\mathcal F} &    T\widetilde{M} \\
			\mathcal F	& TM    \\};
		\path[-stealth]
		(m-1-1) edge node [above] {$\hookto$} (m-1-2)
		(m-1-1) edge node [right] {} (m-2-1)
		(m-1-2) edge node [right] {} (m-2-2)
		(m-2-1) edge node [above] {$\hookto$}  (m-2-2);
		
	\end{tikzpicture}
	\newline
	then $\widetilde{\mathcal F}$ is integrable.
	
	\begin{definition}\label{fol_cov_defn}
		In the above situation we say that a foliation $\left(\widetilde{M},~ \widetilde{\mathcal F} \right)$ is the \textit{induced by} $p$ \textit{covering} of $\left(M,\mathcal F \right)$ or the $p$-\textit{lift} of $\left(M,\mathcal F \right)$. 
	\end{definition}
	\begin{remark}
		The $p$-lift of a foliation is described in  \cite{ouchi:cov_fol, xiaolu:foli_cov}.
	\end{remark}
	\begin{empt}
		If $\gamma: \left[0,1\right]\to M$ is a path which corresponds to an element of the holonomy groupoid then we denote by $\left[\gamma\right]$ its equivalence class, i.e. element of groupoid.
		There is the space of half densities $\Omega_{\widetilde{M}}^{1/2}$ on $\widetilde{M}$ which is a lift  the space of half densities $\Omega_{M}^{1/2}$ on $M$. If $L$ is a leaf of $\left(M,\mathcal F \right)$, $L'=\pi^{-1}\left( L\right)$   then a space $\widetilde{L}$ of holonomy covering of  $L$ coincides with the space of the holonomy covering of $L'$. It turns out that $L^2\left( \widetilde{\mathcal G}_{\widetilde{x}}\right)\approx L^2\left(\mathcal G_{ \pi\left(\widetilde{x}\right)} \right)$ for any $\widetilde{x} \in \widetilde{M}$.
		If $\mathcal G$ (resp. $\widetilde{\mathcal G}$) is a holonomy groupoid of $\left(M,\mathcal F \right)$ (resp. $\left(\widetilde{M},~ \widetilde{\mathcal F} \right)$) then there is the surjective map $p_{\mathcal G}:\widetilde{\mathcal G} \to \mathcal G$ given by
		\begin{equation*}
			\begin{split}
				\left[\widetilde{\gamma}\right] \mapsto \left[p \circ\widetilde{\gamma}\right].
			\end{split}
		\end{equation*}
		If the covering is finite-fold then the map $p_{\mathcal G}:\widetilde{\mathcal G} \to \mathcal G$ induces 
		a natural involutive homomorphism
		\begin{equation*}
			\begin{split}
				\Coo_c\left(\mathcal G,  \Omega_{M}^{1/2}\right) \hookto  \Coo_c\left(\widetilde{\mathcal G},  \Omega_{\widetilde{M}}^{1/2}\right).
			\end{split}
		\end{equation*}	
		Completions of $\Coo_c\left(\mathcal G,  \Omega_{M}^{1/2}\right)$ and $\Coo_c\left(\widetilde{\mathcal G},  \Omega_{\widetilde{M}}^{1/2}\right)$ with respect to given by \eqref{foli_pseudo_norm_eqn} norms gives an injective *- homomorphism 
		\be\label{foli_inc_eqn}
		\pi: C^*_r\left( M,\mathcal F\right) \hookto C^*_r\left(\widetilde{M},~ \widetilde{\mathcal F} \right)
		\ee
		of $C^*$-algebras. The action of the group $G\left( \left.\widetilde{M}~\right|M\right)$ of covering transformations on $\widetilde{M}$ naturally induces an action of $G\left( \left.\widetilde{M}~\right|M\right)$ on $\left(\widetilde{M},~ \widetilde{\mathcal F} \right)$.  It follows that there is the natural action   $C^*_r\left(\widetilde{M},~ \widetilde{\mathcal F} \right)$ such that 
		\be\label{foli_exp_g_eqn}
		C^*_r\left( M, \mathcal{F}\right)   = C^*_r\left(\widetilde{M},~ \widetilde{\mathcal F} \right)^{G\left( \left.\widetilde{M}~\right|M\right)}
		\ee
	\end{empt}
	
	Let $\G\left({M},~ {\mathcal F} \right)$ and $\G\left(\widetilde{M},~ \widetilde{\mathcal F} \right)$ be the holonomy groupoids of $\left({M},~ {\mathcal F} \right)$ and $\left(\widetilde{M},~ \widetilde{\mathcal F} \right)$ respectively. The natural surjective map $\G\left(\widetilde{M},~ \widetilde{\mathcal F} \right)\to\G\left({M},~ {\mathcal F} \right)$ induces the injective $*$-homomorphism $C_r^*\left(\widetilde{M},~ \widetilde{\mathcal F} \right)\hookto C^*_r\left(\widetilde{M},~ \widetilde{\mathcal F} \right)$.
	Assume both $\G\left({M},~ {\mathcal F} \right)$ and $\G\left(\widetilde{M},~ \widetilde{\mathcal F} \right)$ are Hausdorff. Let $G\left( \left.\widetilde{M}~\right|M\right)$  be the covering 	group of  $p:\widetilde M \to M$ . The $G\left( \left.\widetilde{M}~\right|M\right)$-action on $\widetilde{M}$ can be naturally extended to $\G\left(\widetilde{M},~ \widetilde{\mathcal F} \right)$ by sending $\widetilde{\ga}$ to $g\widetilde{\ga}$ for a representative path $\widetilde{\ga}$ in $ \widetilde{M}$ and $g \in G\left( \left.\widetilde{M}~\right|M\right)$. 
	\begin{lemma}\cite{xiaolu:foli_cov}
		Let $p:\widetilde M \to M$ be a regular covering manifold with covering group $G\left( \left.\widetilde{M}~\right|M\right)$
		$N\subset M$ a connected submanifold, and $\widetilde N$ a connected component of $p^{-1}\left( N\right) $
		Then the restriction $p_{\widetilde N}$ of $p$ to $N$ is also regular, with the covering group $G\left( \left.\widetilde{N}~\right|N\right)$ being a  	subgroup  $G\left( \left.\widetilde{M}~\right|M\right)$.
	\end{lemma}
	In particular, if $L_{\widetilde x}$ is the leave in $\left(\widetilde{M},~ \widetilde{\mathcal F} \right)$ containing $\widetilde x \in p^{-1}\left( x\right)$ , where $x \in M$,
	then $L_{\widetilde x}$ is a regular cover of $L_x$. We denote the covering group by $G\left( \left.\widetilde{M}~\right|M\right)_x$. On the  other hand, for each $x \in M$, there is a holonomy group $\G^x_x$, and we have the holonomy group bundle $\left\{\G^x_x\right\}$ over $M$. If $x_1$, $x_2$ are on the same leaf, then any path $\ga$  connecting $x_1$ and $x_2$ induces an isomorphism $\ga^*: \G^{x_1}_{x_1}\xrightarrow{\cong}\G^{x_2}_{x_2}$  by mapping $\left[\ga_1\right]$ to $\left[\ga\ga_1\ga^{-1}\right]$.
	As a local homeomorphism, the covering map $p$ induces an embedding $\overline{p} : \G^{\widetilde x}_{\widetilde x}\to \G^x_x$ for each $\widetilde x \in p^{-1}\left( x\right)$.
	\begin{lemma}\cite{xiaolu:foli_cov}
		The group $\overline{p}^*_x\left(\G^{\widetilde x}_{\widetilde x}\right)$ is a normal subgroup of $G^x_x$. Equivalently,
		$\overline{p}^*_x\left( \G^{\widetilde x_1}_{\widetilde x_1}\right) =\overline{p}^*_x\left( \G^{\widetilde x_2}_{\widetilde x_2}\right)$ for $\widetilde x_1, \widetilde x_2 \in p^{-1}\left(x\right)$ if $L_{\widetilde x_1}=L_{\widetilde x_2}$,
	\end{lemma}
	Thus we may form the quotient holonomy group bundle $\left\{\G^{\widetilde x}_{\widetilde x}\right\}$ over $M$. There is an obvious group homomorphism $\phi_x: G\left( \left.\widetilde{M}~\right|M\right)_x \to \G^{ x}_{ x}/\G^{\widetilde x}_{\widetilde x}$ defined as follows. An element $g\in G\left( \left.\widetilde{M}~\right|M\right)_x$ corresponds to a point $x_g \in p^{-1}\left( x\right)\cap L_{\widetilde x}$ if we fix $\widetilde x$ corresponding
	to the unit $e$. A path $\ga_g$ starting at  $\widetilde x$ and ending at $g \widetilde x$ gives a loop $\pi\left(\widetilde \ga_g \right)$  in $M$ representing an element $\phi_x\left(g \right)$  in $G_x$, whose class in $\G^{ x}_{ x}/\G^{\widetilde x}_{\widetilde x}$ is uniquely defined
	by $g$. Given any $\left[\ga\right]$ in $\G^{ x}_{ x}$ there is a preimage $\widetilde \ga$ in $L_{\widetilde x}$ starting at $\widetilde x$. The point $r\left( \widetilde \ga\right) \in p^{-1}\left( x\right)\cap  L_{\widetilde x}$ corresponding to some $g\in G\left( \left.\widetilde{M}~\right|M\right)_x$. So $\phi_x$ is onto.
	\begin{definition}\cite{xiaolu:foli_cov}\label{foli_reg_cov_defn}
		The covering map $p:\left(\widetilde{M},~ \widetilde{\mathcal F} \right)\to\left({M},~ {\mathcal F} \right)$ of foliations is said to be \textit{regular} if the map $\phi$ is an isomorphism from the leaf covering group bundle to the quotient holonomy group bundle.
	\end{definition}
	
	\begin{remark}\label{foli_action_rem}
		Every regular covering map $p:\left(\widetilde{M},~ \widetilde{\mathcal F} \right)\to\left({M},~ {\mathcal F} \right)$ of foliations induces nontrivial action of $G\left( \left.\widetilde{M}~\right|M\right)$ on $C^*_r\left(\widetilde{M},~ \widetilde{\mathcal F}  \right)$ such that $$C^*_r\left(\widetilde{M},~ \widetilde{\mathcal F}  \right)^{ G\left( \left.\widetilde{M}~\right|M\right)}\cong C^*_r\left(\widetilde{M},~ \widetilde{\mathcal F}  \right).$$
	\end{remark}
	\begin{remark}\label{foli_trans_rem}
		If a map $p:\left(\widetilde{M},~ \widetilde{\mathcal F} \right)\to\left({M},~ {\mathcal F} \right)$ is regular covering of foliations  and a leaf $\widetilde L \in \widetilde M$ has no holonomy then  $g\widetilde L \neq \widetilde L$ for all nontrivial $g \in G\left( \left.\widetilde{M}~\right|M\right)$, i.e. $G\left( \left.\widetilde{M}~\right|M\right)$ transitively acts on leaves having no holonomy.
	\end{remark}

	\chapter{Noncommutative torus and Moyal plane}\label{nt_descr_sec}
	\section{Noncommutative torus $\mathbb{T}^n_{\Theta}$}\label{nt_descr_subsec}
	\subsection{Definition of noncommutative torus $\mathbb{T}^n_{\Theta}$}
	\begin{definition}\label{nt_qirr_defn}\cite{wagner:pb}
		Denote by "$\cdot$" the scalar product on $\R^n$.
		The matrix $\Theta$ is called \emph{quite irrational} if, for all $\lambda \in \Z^n$, the condition $\exp(2\pi i \, \lambda\cdot \Theta  \mu) = 1$ for all $
		\la,\mu \in \Z^n$ implies $\lambda = 0$.
	\end{definition}
	\begin{definition}\label{nt_defn}\cite{wagner:pb}
		Let $\Theta$ be an invertible, real  skew-symmetric quite irrational $n \times n$ matrix. A \textit{noncommutative torus} $C\left(\mathbb{T}^n_{\Theta}\right)$ is the universal $C^*$-algebra generated by the set $\left\{U_k\right\}_{k \in \Z^n}$ of unitary elements which satisfy to the following relations.
		\begin{equation}\label{nt_unitary_product_eqn}
			U_k U_p = e^{-\pi ik ~\cdot~ \Theta p} U_{k + p}; ~~~   \end{equation}
	\end{definition}
	Following condition holds
	\begin{equation}\label{nt_unitary_product_comm_eqn}
		U_k U_p = e^{-2\pi ik ~\cdot~ \Theta p}U_p U_k.
	\end{equation}
	An alternative description of $\C\left(\mathbb{T}^n_{\Theta}\right)$ is such that if
	\begin{equation}\label{nt_th_eqn}
		\Th = \begin{pmatrix}
			0& \th_{12} &\ldots & \th_{1n}\\
			\th_{21}& 0 &\ldots & \th_{2n}\\
			\vdots& \vdots &\ddots & \vdots\\
			\th_{n1}& \th_{n2} &\ldots & 0
		\end{pmatrix}=\begin{pmatrix}
		0& \th_{12} &\ldots & \th_{1n}\\
		-\th_{12}& 0 &\ldots & \th_{2n}\\
		\vdots& \vdots &\ddots & \vdots\\
		-\th_{1n}& -\th_{2n} &\ldots & 0
		\end{pmatrix}
	\end{equation}
	then $C\left(\mathbb{T}^n_{\Theta}\right)$ is the universal $C^*$-algebra generated by unitary elements   $u_1,..., u_n \in U\left( C\left(\mathbb{T}^n_{\Theta}\right)\right) $ such that following condition holds
	\begin{equation}\label{nt_com_eqn}
		\begin{split}
			u_j u_k = e^{-2\pi i \theta_{jk} }u_k u_j.
		\end{split}
	\end{equation}
	Unitary  operators $u_1,..., u_n$ correspond to the standard basis of $\mathbb{Z}^n$ and they are given by
	\be\label{nt_gen_eqn}
	u_j = U_{k_j },\quad \text{ where } k_j=\left(0,...,\underbrace{ 1}_{j^{\text{th}}-\text{place}},...,0 \right)
	\ee
	\begin{defn}\label{nt_uni_defn}
		The unitary elements 
		$u_1,..., u_n \in U\left(C\left(\mathbb{T}^n_{\theta}\right)\right)$ which satisfy the relations \eqref{nt_com_eqn}, \eqref{nt_gen_eqn}
		are said to be \textit{generators} of $C\left(\mathbb{T}^n_{\Theta}\right)$. The set $\left\{U_l\right\}_{l \in \Z^n}$ is said to be the \textit{basis} of $C\left(\mathbb{T}^n_{\Theta}\right)$.
	\end{defn}

	If $a \in C\left(\mathbb{T}^n_{\Th}\right)$ is presented by a series
	\be\label{nt_series_eqn}
	a = \sum_{l \in \mathbb{Z}^{n}}c_l U_l;~~ c_l \in \mathbb{C}
	\ee
	and the series $\sum_{l \in \mathbb{Z}^{n}}\left| c_l\right| $ is convergent then from the triangle inequality it follows that the series is $C^*$-norm convergent and the following condition holds.
	\begin{equation}\label{nt_norm_estimation_eqn}
		\left\|a \right\| \le \sum_{l \in \mathbb{Z}^{n}}\left| c_l\right|.
	\end{equation}
	In particular if $\mathrm{sup}_{l \in \mathbb{Z}^n}\left(1 + \|l\|\right)^s \left|c_l\right| < \infty, ~ \forall s \in \mathbb{N}$ then $\sum_{l \in \mathbb{Z}^{n}}\left| c_l\right|< \infty $ and taking into account the equation \eqref{schwartz_z_eqn} there is the natural inclusion $\phi_\infty: \sS\left(\mathbb{Z}^n\right) \subset C\left(\mathbb{T}^n_{\Theta}\right)$ of vector spaces. If
	\be\label{nt_coo_eqn}
	\Coo\left(\mathbb{T}^n_{\Theta}\right)\bydef \phi_\infty \left( \sS\left(\mathbb{Z}^n\right)\right) \subset C\left(\mathbb{T}^n_{\Theta}\right)
	\ee
	then $\Coo\left(\mathbb{T}^n_{\Theta}\right)$ is a pre-$C^*$-algebra, and the Fourier transformation \ref{nt_fourier_eqn} yields the $\C$-isomorphism
	\be\label{nt_coo_iso_eqn}
	\Coo\left(\mathbb{T}^n\right)\cong \Coo\left(\mathbb{T}^n_{\Theta}\right).
	\ee
	\begin{empt}\label{nt_gns_empt}
		There is a  state 
		\be\label{nt_state_eqn}
		\begin{split}
			\tau: C\left(\mathbb{T}^n_{\Theta}\right) \to \C;\\
			\sum_{k \in \Z^n} a_k U_k \mapsto a_{\left(0,...,0\right) }; \quad \text{ where } a_k \in \C,
		\end{split}
		\ee
		which induces the faithful GNS representation. 
		The $C^*$-norm completion  $C\left(\mathbb{T}^n_{\Theta}\right)$ of $\Coo\left(\mathbb{T}^n_{\Theta}\right)$ is a $C^*$-algebra and there is a faithful representation
		\begin{equation}\label{nt_repr_eqn}
			C\left(\mathbb{T}^n_{\Theta}\right) \to B\left( L^2\left(C\left(\mathbb{T}^n_{\Theta}\right), \tau\right)\right) .
		\end{equation}
		(cf. Definition \ref{gns_defn}).  Similarly to the equation \eqref{from_a_to_l2_eqn} there is a $\C$-linear map 
		\begin{equation}\label{nt_to_hilbert_eqn}
			\Psi_\Th:C\left(\mathbb{T}^n_{\Theta}\right) \hookto L^2\left(C\left(\mathbb{T}^n_{\Theta}\right), \tau\right).
		\end{equation}
		If 
		\be\label{nt_xik_eqn}
		\xi_k \bydef \Psi_\Th\left(U_k \right)
		\ee 
		then from \eqref{nt_unitary_product_eqn}, \eqref{nt_state_eqn} it turns out
		\begin{equation}\label{nt_h_product_eqn}
			\tau\left(U^*_k  U_l \right) = \left(\xi_k, \xi_l \right)  = \delta_{kl},   
		\end{equation} 
		i.e. the subset $\left\{\xi_k\right\}_{k \in \mathbb{Z}^n}\subset L^2\left(C\left(\mathbb{T}^n_{\Theta}\right), \tau\right)$ is an orthogonal basis of  $L^2\left(C\left(\mathbb{T}^n_{\Theta}\right), \tau\right)$.
		Hence the Hilbert space  $L^2\left(C\left(\mathbb{T}^n_{\Theta}\right), \tau\right)$ is naturally isomorphic to the Hilbert space $\ell^2\left(\mathbb{Z}^n\right)$ given by
		\begin{equation*}
			\ell^2\left(\mathbb{Z}^n\right) = \left\{\xi = \left\{\xi_k \in \mathbb{C}\right\}_{k\in \mathbb{Z}^n} \in \mathbb{C}^{\mathbb{Z}^n}~|~ \sum_{k\in \mathbb{Z}^n} \left|\xi_k\right|^2 < \infty\right\}
		\end{equation*}
		and the $\C$-valued scalar product on $\ell^2\left(\mathbb{Z}^n\right)$ is given by
		\begin{equation}\label{nt_xi_eqn}
			\left(\xi,\eta\right)_{ \ell^2\left(\mathbb{Z}^n\right)}= \sum_{k\in \mathbb{Z}^n}    \overline{\xi}_k\eta_k.
		\end{equation}
		From \eqref{nt_unitary_product_eqn} and \eqref{nt_xik_eqn} it follows that the representation $C\left(\mathbb{T}^n_{\Theta}\right)\hookto B\left( L^2\left(C\left(\mathbb{T}^n_{\Theta}\right), \tau\right)\right)$ corresponds to the following action
		\be\label{nt_l2_eqn}
		\begin{split}
			C\left(\mathbb{T}^n_{\Theta}\right)\times  L^2\left(C\left(\mathbb{T}^n_{\Theta}\right), \tau\right)\to  L^2\left(C\left(\mathbb{T}^n_{\Theta}\right), \tau\right);\\
			U_k\xi_l = e^{-\pi ik ~\cdot~ \Theta l}\xi_{k + l} .
		\end{split}
		\ee
	\end{empt}

	\subsection{Geometry of noncommutative tori}\label{nt_geom_sec}
	\paragraph*{}
	In the below text we imply that  $\Theta$ is {quite irrational}. 
	The restriction of given by \eqref{nt_state_eqn}  state on $\Coo\left(\mathbb{T}^n_{\Theta}\right)$ satisfies to the following equation
	\begin{equation}\label{nt__smotth_state_eqn}
		\tau\left(f\right)= \widehat{f}\left(0\right)
	\end{equation}
	where $\widehat{f}$ means the Fourier transformation.
	From  $\Coo\left(\mathbb{T}^n_{\Theta} \right) \approx \SS\left( \Z^n\right)$ it follows  that there is a $\C$-linear isomorphism 
	\begin{equation}\label{nt_varphi_inf_eqn}
		\varphi_\infty: \Coo\left(\mathbb{T}^n_{\Theta} \right) \xrightarrow{\approx}  \Coo\left(\mathbb{T}^n \right).
	\end{equation} 
	such that following condition holds
	\begin{equation}\label{nt_state_integ_eqn}
		\tau\left(f \right)=  \frac{1}{\left( 2\pi\right)^n }\int_{\mathbb{T}^n} \varphi_\infty\left( f\right) ~dx.
	\end{equation}
	From \eqref{nt_state_integ_eqn} it follows that for any $a, b \in \Coo\left(\mathbb{T}^n_{\Theta}\right)$  the scalar product on $L^2\left(C\left(\mathbb{T}^n_{\Theta}\right), \tau\right)$ is given by
	\begin{equation}\label{nt_int_sc_pr_eqn}
		\left(a, b \right)= \int_{\T^n} a_{\text{comm}}^*b_{\text{comm}}dx 
	\end{equation}
	where $a_{\text{comm}}\in \Coo\left(\T^n \right)$ (resp. $b_{\text{comm}})$ is a commutative function which corresponds to $a$ (resp. $b$).
		\begin{empt}\label{nt_k_1_empt}
		If		$u_1,..., u_n \in U\left(C\left(\mathbb{T}^n_{\theta}\right)\right)$ are {generators} of $C\left(\mathbb{T}^n_{\Theta}\right)$ (cf. Definition \ref{nt_uni_defn}) then for all $\left( k_1, ..., k_n\right) \in\Z^n $ one has
		$$
		\left( k_1, ..., k_n\right)\neq 0 \quad \Rightarrow\quad	\tau \left( u_1^{k_1}\cdot ...\cdot u_n^{k_n}\right)= 0. 
		$$
		where $\tau$ is given by \eqref{nt_state_integ_eqn}. On the other hand from
		$
		\tau\left(1_{C\left(\mathbb{T}^n_{\theta}\right)}\right)= 1
		$
	 it follows that any nontrivial product $ u_1^{k_1}\cdot ...\cdot u_n^{k_n}$ is not homotopic to $1_{C\left(\mathbb{T}^n_{\theta}\right)}$ in
			 $U\left(C\left(\mathbb{T}^n_{\theta}\right)\right)$. Similarly to \cite{varilly:noncom} one can deduce the following
			\be\label{nt_k_1_eqn}
		K_1\left(C\left(\mathbb{T}^n_{\theta}\right)\right)= \Z\left[u_1\right]\oplus ... \oplus \Z\left[u_n\right]
			\ee
			where $K_1$ means the $K_1$-functor (cf. \cite{blackadar:ko}).
	\end{empt}
	\begin{defn}\label{nt_symplectic_defn}
		If  $\Theta$ is non-degenerate, that is to say,
		$\sigma(s,t) \stackrel{\mathrm{def}}{=} s\.\Theta t$ to be \textit{symplectic}. This implies even
		dimension, $n = 2N$. One then selects
		\begin{equation}\label{nt_simpectic_theta_eqn}
			\Theta = \theta J
			\stackrel{\mathrm{def}}{=} \th \begin{pmatrix} 0 & 1_N \\ -1_N & 0 \end{pmatrix}
		\end{equation}
		where  $\th > 0$ is defined by $\th^{2N} \stackrel{\mathrm{def}}{=} \det\Theta$.
		Denote by
		\be\label{nt_symplectic_eqn}
		\begin{split}
	\Coo\left(\mathbb{T}^{2N}_\th\right)\stackrel{\mathrm{def}}{=}\Coo\left(\mathbb{T}^{2N}_\Th\right),\\
	C\left(\mathbb{T}^{2N}_\th\right)\stackrel{\mathrm{def}}{=}C\left(\mathbb{T}^{2N}_\Th\right).
		\end{split}
\ee
	\end{defn}
	\paragraph*{}
	The space $\Coo\left(\T^n_{\Theta}\right)$ is a dense subspace of  the Hilbert space $L^2\left(C\left(\mathbb{T}^n_{\Theta}\right), \tau\right)$. Let  $\L^\dagger\left(\Coo\left(\mathbb{T}^n_{\Theta}\right)\right)$ be a given by  the equation \eqref{l_dag_eqn} $O^*$-algebra. Denote by $\delta_{\mu}\in \L^\dagger\left(\Coo\left(\mathbb{T}^n_{\Theta}\right)\right)$ ($\mu = 1,\dots, n$) the analogues of the partial derivatives  on $C^{\infty}(\mathbb{T}^n)$ which are derivations on the algebra $C^{\infty}(\mathbb{T}^n_{\Theta})$ given by
	\be\label{nt_diff_eqn}
k = \left(k_1,..., k_n\right)\quad \Rightarrow\quad	\delta_{\mu}(U_k)=ik_{\mu} U_k.
	\ee
	These derivations have the following property
	$$
	\delta_{\mu}(a^*)=-(\delta_{\mu}a)^*,
	$$
	and also satisfy  the integration by parts formula
	$$\tau(a\delta_{\mu}b)=-\tau((\delta_{\mu}a)b),\quad a,b\in C^{\infty}(\mathbb{T}^n_{\Theta}).$$
From 	the equation \eqref{nt_unitary_product_eqn}	one has
	\bean
		\delta_\mu\left(U_k U_l \right)=  e^{-\pi ik ~\cdot~ \Theta l}	\delta_\mu\left(U_{k+l} \right)= i\left(k_\mu + l_\mu \right) e^{-\pi ik ~\cdot~ \Theta l}U_{k+l} =\\=
		i\left(k_\mu + l_\mu \right)U_k U_l= \left( ik_\mu U_k\right) U_l + U_k \left( il_\mu U_l\right)= \left( \delta_\mu U_k\right) U_l + U_k \left( \delta_\mu U_l\right),
	\eean
	and from the above equation one has a Leibniz rule
	\be\label{nt_le_eqn}
\forall a, b \in C^{\infty}(\mathbb{T}^n_{\Theta})\quad \delta_{\mu}\left( a\cdot b\right)  =\left( \delta_\mu a\right)\cdot b + a\cdot \left( \delta_\mu b\right)
	\ee
	where $\cdot$ means the product of the involutive algebra $\Coo\left(\mathbb{T}^n_{\Theta}\right)$.
	The spectral triple describing the noncommutative geometry of noncommutative $n$-torus consists of the algebra $C^{\infty}(\mathbb{T}^n_{\Theta})$, the Hilbert space $\mathcal{H}=L^2\left(C\left(\mathbb{T}^n_{\Theta}\right), \tau\right)\otimes\mathbb{C}^{m}$, where $m=2^{[n/2]}$ with  the representation  $\pi\otimes 1_{\mathbb{M}_n\left(\C\right) }:C^{\infty}(\mathbb{T}^n_{\Theta})\to B\left( \H\right) $  where $\pi: C^{\infty}(\mathbb{T}^n_{\Theta})\to B\left( L^2\left(C\left(\mathbb{T}^n_{\Theta}\right), \tau\right)\right)$ is given by \eqref{nt_repr_eqn}.
	The Dirac operator is given by
	\begin{equation}\label{nt_dirac_eqn}
		D\stackrel{\mathrm{def}}{=}\sum_{\mu =1}^{n}\delta_{\mu} \otimes\gamma^{\mu}\in \L^\dagger\left(\Coo\left(\mathbb{T}^n_{\Theta}\right)\otimes\C^m\right),
	\end{equation}
	where $\delta_{\mu}$ is given by \eqref{nt_diff_eqn}, seen as an unbounded self-adjoint operator on $L^2\left(C\left(\mathbb{T}^n_{\Theta}\right), \tau\right) $ and $\gamma^{\mu}$s are Clifford (Gamma) matrices in $\mathbb{M}_m(\mathbb{C})$ satisfying the relation
	\be\label{nt_gam_eqn}\gamma^\mu\gamma^\nu+\gamma^\nu\gamma^\mu=2\delta_{\mu\nu} 1_{\mathbb{M}_m(\mathbb{C})}.\ee
	There is a spectral triple 
	\begin{equation}\label{nt_sp_tr_eqn}
		\left( C^{\infty}(\mathbb{T}^n_{\Theta}),L^2\left(C\left(\mathbb{T}^n_{\Theta}\right), \tau\right)\otimes\mathbb{C}^{m},D\right).
	\end{equation}
There is an alternative description of $D$. Any $b \in \Coo\left(\T^n\right)$ can be regarded as $f_b \in  \Coo\left(\R^n\right)$ which is $\Z^n$-periodic, i.e. it is invariant with respect to $\Z^n$-shifts. If $x_1,..., x_n$ are coordinates of $\R^n$ then for all $\mu= 1,..., n$ 
the partial derivation $\frac{\partial f_b}{\partial x_\mu}$ is smooth and  $\Z^n$-periodic. We write $\frac{\partial b}{\partial x_\mu}$ instead of $\frac{\partial f_b}{\partial x_\mu}$ so for all $\mu= 1,..., n$ one can assume $\frac{\partial b}{\partial x_\mu}\in \Coo\left(\T^n\right)$, i.e. there is an operator $\frac{\partial}{\partial x_\mu}\in \L^\dagger\left(  \Coo\left(\T^n\right)\right) $. Using the given by \eqref{nt_varphi_inf_eqn} $\C$-linear isomorphism  $\varphi_\infty:\Coo\left(\mathbb{T}^n_{\Theta} \right) \cong \Coo\left(\mathbb{T}^n \right)$   one can set 
	\be\label{nt_ddx_eqn}
	\frac{\partial}{\partial x_\mu}\in \L^\dagger\left( \Coo\left(\mathbb{T}^n_{\Theta} \right) \right)
	\ee
Elements of the basis $U_k$ and generators $u_\mu$ (cf. Definition  \ref{nt_uni_defn})  $U_k$ correspond to  $\Z^n$-periodic maps  $\R^n \to \C$ such that 
\bean
\\
U_k: x \mapsto e^{2\pi i k \cdot x},\\
u_\mu : x \mapsto e^{2\pi i x_\mu}\quad x = \left(x_1, ..., x_n\right)
\eean
Using chain rule one has
\be\label{nt_chain_eqn}
\frac{\partial}{\partial x_\mu} = \frac{1}{2\pi i}u^*\frac{\partial}{\partial u_\mu}, 
\ee	
so $\frac{\partial}{\partial u_\mu}$ can be regarded as a $\C$=linear map 
\be\label{nt_chainu_eqn}
\frac{\partial}{\partial u_\mu}: \Coo\left(\R^n \right)\to \Coo\left(\R^n \right). 
\ee
Otherwise one has
\be\label{nt_chaind_eqn}
\delta_\mu =  i u_\mu \frac{\partial}{\partial u_\mu} = 2\pi \frac{\partial}{\partial x_\mu} ,
\ee	
so the equation \eqref{nt_dirac_eqn}  can be rewritten by the following way
	\begin{equation}\label{nt_dirac_u_eqn}
	D\stackrel{\mathrm{def}}{=}\sum_{\mu =1}^{n}\delta_{\mu} \otimes\gamma^{\mu}= i \sum_{\mu =1}^{n} u_\mu \frac{\partial}{\partial u_\mu} \otimes\gamma^{\mu}\in \L^\dagger\left(\Coo\left(\mathbb{T}^n_{\Theta}\right)\otimes\C^m\right).
\end{equation}
Form the equation  \eqref{nt_diff_eqn} it follows that 
Let $b \in \Coo\left(\mathbb{T}^n_{\Theta} \right)$, $~\xi\in \C^m$
\be\label{nt_ga_eqn}
\begin{split}
\left( \left(u^*_\mu \otimes 1_{\mathbb{M}_m\left(\C\right)} \right)  \left[ D, u_\mu \otimes 1_{\mathbb{M}_m\left(\C\right)} \right]\right)\left( b\otimes \xi\right) =\\
\left(u^*_\mu \otimes 1_{\mathbb{M}_m\left(\C\right)} \right)\left(  D\left( u_\mu b\otimes \xi\right) - \left(u_\mu \otimes 1_{\mathbb{M}_m\left(\C\right)} \right) D\left( b\otimes \xi\right) \right).
\end{split}
 \ee
From the definition \eqref{nt_dirac_u_eqn} of the Dirac operator and tha chain rule \eqref{nt_le_eqn} it follows that 
\bean
D\left( u_\mu b\otimes \xi\right)= \left( \sum_{\mu =1}^{n}\delta_{\mu} \otimes\gamma^{\mu}\right)\left( u_\mu b\otimes \xi\right)= i \left( u_\mu \cdot b\right) \otimes\xi +  \sum_{\mu =1}^{n}u_\mu \delta_\mu b \otimes \ga^\mu\xi,\\
\left(u_\mu \otimes 1_{\mathbb{M}_m\left(\C\right)} \right) D\left( b\otimes \xi\right)= \sum_{\mu =1}^{n}u_\mu \delta_\mu b \otimes \ga^\mu\xi
\eean
The substitution of the above equation into \eqref{nt_ga_eqn} yields the following
\bean
\left( \left(u^*_\mu \otimes 1_{\mathbb{M}_m\left(\C\right)} \right)  \left[ D, u_\mu \otimes 1_{\mathbb{M}_m\left(\C\right)} \right]\right)\left( b\otimes \xi\right) =\\
\left(iu^*_\mu  \otimes 1_{\mathbb{M}_m\left(\C\right)}\right) \left( u_\mu b\otimes \ga \xi\right) = i\left(b \otimes \ga^\mu \xi \right).
\eean
or, equivalently
\be\label{nt_ga'_eqn}
\left(u^*_\mu \otimes 1_{\mathbb{M}_m\left(\C\right)} \right)  \left[ D, u_\mu \otimes 1_{\mathbb{M}_m\left(\C\right)} \right]= i 1_{\Coo\left(\mathbb{T}^n_{\Theta} \right)}\otimes \ga^\mu.
\ee
If $\Om^1_D$ is a {module of differential forms associated} with the given by \eqref{nt_sp_tr_eqn}  spectral triple (cf. Definition \ref{ass_cycle_defn}) then $\left(u^*_\mu \otimes 1_{\mathbb{M}_m\left(\C\right)} \right)  \left[ D, u_\mu \otimes 1_{\mathbb{M}_m\left(\C\right)} \right]\in \Om^1_D$ and taking into account \eqref{nt_ga'_eqn} one has
\be\label{nt_ga_in_eqn}
1_{\Coo\left(\mathbb{T}^n_{\Theta} \right)}\otimes \ga^\mu\in \Om^1_D.
\ee
\section{Moyal plane}
	
	\begin{defn}
		Denote the \textit{Moyal plane} product $\star_\th$ on $\SS\left(\R^{2N} \right)$ given by
		\begin{equation}\label{mp_prod_eqn}
			\left(f \star_\th h \right)\left(u \right)= \int_{y \in \R^{2N} } f\left(u - \frac{1}{2}\Th y\right) g\left(u + v \right)e^{2\pi i y \cdot v }  dydv
		\end{equation}
		
		where $\Th$ is given by \eqref{nt_simpectic_theta_eqn}.
	\end{defn}
	There is the tracial property \cite{moyal_spectral} of the Moyal product
	\begin{equation}\label{nt_tracial_prop}
		\int_{\R^{2N}} \left( f\star_\th g\right) \left(x \right)dx =  \int_{\R^{2N}}  f\left(x \right) g\left(x \right)dx.
	\end{equation}
	The Fourier transformation of the star product satisfies to the following condition.
	\begin{equation}\label{mp_fourier_eqn}
		\mathcal{F}\left(f \star_\th g\right) \left(x \right) =    	\int_{\R^{2N}}\mathcal{F}{f}\left(x-y \right) \mathcal{F}{g}\left(y\right)e^{ \pi i  y \cdot \Th x }~dy.
	\end{equation}
	\begin{prop}\label{mp_factor_prop}\cite{moyal_spectral}
		The algebra $\SS\left(\R^{2N}, \star_\th \right)$  has the (nonunique) factorization property:
		for all $h \in  \SS\left(\R^{2N} \right)$ there exist $f,g \in  \SS\left(\R^{2N} \right)$ that $h = f \star_\th g$.
	\end{prop}
	\begin{defn}\label{mp_mult_defn}\cite{varilly_bondia:phobos}
		\label{df:Moyal-alg}	Denote by $\SS'\left( \R^{n}\right) $ the vector space dual to $\SS\left( \R^{n}\right) $, i.e. the space of continuous functionals on $\SS\left( \R^{n}\right)$.
		The Moyal product can be defined, by duality, on larger sets than
		$\SS\left(\R^{2N}\right)$. For $T \in \SS'\left(\R^{2N}\right)$, write the evaluation on $g \in \SS\left(\R^{2N}\right)$ as
		$\<T, g> \in \C$; then, for $f \in \SS$ we may define $T \star_{\theta} f$ and
		$f \star_{\theta} T$ as elements of~$\SS'\left(\R^{2N}\right)$ by
		\begin{equation}\label{mp_star_ext_eqn}
			\begin{split}
				\<T \star_{\theta} f, g> \stackrel{\mathrm{def}}{=} \<T, f \star_{\theta} g>\\
				\<f \star_{\theta} T, g> \stackrel{\mathrm{def}}{=} \<T, g \star_{\theta} f>	\end{split}
		\end{equation}  using the continuity of the
		star product on~$\SS\left(\R^{2N}\right)$. Also, the involution is extended to  by
		$\<T^*,g> \stackrel{\mathrm{def}}{=} \overline{\<T,g^*>}$.	
		Consider the left  and right multiplier algebras:
		\begin{equation}\label{mp_multi_eqn}
			\begin{split}
				\M_L^\th
				&\stackrel{\mathrm{def}}{=} \set{T \in \SS'(\R^{2N}) : T \star_{\theta} h \in \SS(\R^{2N})
					\text{ for all } h \in \SS(\R^{2N})},
				\\
				\M_R^\th
				&\stackrel{\mathrm{def}}{=} \set{T \in \SS'(\R^{2N}) : h \star_{\theta} T \in \SS(\R^{2N})
					\text{ for all } h \in \SS(\R^{2N})},\\
				\M^\th &\stackrel{\mathrm{def}}{=} \M_L^\th \cap \M_R^\th.\\
			\end{split}
		\end{equation} 
	\end{defn}
	In \cite{varilly_bondia:phobos} it is proven that
	\bea\label{mp_mult_distr}
	\M_R^\th \star_{\theta} \SS'\left(\R^{2N}\right) = \SS'\left(\R^{2N}\right) \text{ and }
	\SS'\left(\R^{2N}\right) \star_{\theta} \M_L^\th = \SS'\left(\R^{2N}\right),\\
	\label{mp_parial_mult_eqn}
	\partial_j S \times f = 	\partial_j \left( S \times f \right) - f \times \partial_j S \quad \forall S \in \M_L^\th \quad \forall f \in \SS(\R^{2})
	\eea

	It is known \cite{moyal_spectral} that the domain of  the Moyal plane product can be extended up to $L^2\left(\R^{2N} \right)$. 
	
	\begin{lem}\label{nt_l_2_est_lem}\cite{moyal_spectral}
		If $f,g \in L^2 \left(\R^{2N} \right)$, then $f\star_\th g \in L^2 \left(\R^{2N} \right)$ and $\left\|f\right\|_{\mathrm{op}} < \left(2\pi\th \right)^{-\frac{N}{2}} \left\|f\right\|_2$.
		where	$\left\|\cdot\right\|_{2}$ be the $L^2$-norm given by
		\begin{equation}\label{nt_l2_norm_eqn}
			\left\|f\right\|_{2} \stackrel{\mathrm{def}}{=} \left|\int_{\R^{2N}} \left|f\right|^2 dx \right|^{\frac{1}{2}}.
		\end{equation}
		and the operator norm
		\be\label{mp_op_eqn}
		\|T\|_{\mathrm{op}} \stackrel{\mathrm{def}}{=}\sup\set{\|T \star g\|_2/\|g\|_2 : 0 \neq g \in L^2\left( \R^{2N})\right) }
		\ee 
	\end{lem}
	
	\begin{defn}\label{mp_star_alg_defn}
		Denote by $\SS\left(\R^{2N}_\th \right)$  (resp. $L^2\left(\R^{2N}_\th \right)$ ) the operator algebra  which is $\C$-linearly isomorphic to $\SS\left(\R^{2N} \right)$  (resp. $L^2\left(\R^{2N} \right)$ ) and product coincides with  $\star_\th$. Both  $\SS\left(\R^{2N}_\th \right)$ and $L^2\left(\R^{2N}_\th \right)$ act on the Hilbert space $L^2\left(\R^{2N} \right)$. Denote by
		\begin{equation}\label{mp_psi_th_eqn}
			\Psi_\th:  \SS\left(\R^{2N} \right)\xrightarrow{\approx}\SS\left(\R^{2N}_\th \right)
		\end{equation}
		the natural $\C$-linear isomorphism.  
	\end{defn}
	\begin{defn}\label{nt_r_2_N_repr_defn}\cite{moyal_spectral} Let $\SS'\left(\R^{2N} \right)$ be a vector space dual to $\SS\left(\R^{2N} \right)$. Denote by $C_b\left(\R^{2N}_\th\right)\stackrel{\mathrm{def}}{=} \set{T \in \SS'\left(\R^{2N}\right) : T \star_\th g \in L^2\left(\R^{2N}\right) \text{ for all } g \in L^2(\R^{2N})}$, provided with the given by \eqref{mp_op_eqn} operator norm $	\|\cdot\|_{\mathrm{op}}$.
		Denote by $C_0\left(\R^{2N}_\th \right)$ the operator norm completion of $\SS\left(\R^{2N}_\th \right).$  
	\end{defn}
	\begin{rem}
		Obviously $\SS\left(\R^{2N}_\th\right)  \hookto C_b\left(\R^{2N}_\th\right)$. But $\SS\left(\R^{2N}_\th\right)$ is not dense in $C_b\left(\R^{2N}_\th\right)$, i.e. $C_0\left(\R^{2N}_\th\right) \subsetneq C_b\left(\R^{2N}_\th\right)$  (cf. \cite{moyal_spectral}).
	\end{rem}
	
	\begin{rem}
		$L^2\left(\R^{2N}_\th\right)$ is the $\|\cdot\|_2$ norm completion of $\SS\left(\R^{2N}_\th\right)$ hence 
		from the Lemma \ref{nt_l_2_est_lem} it follows that 
		\begin{equation}\label{mp_2_op_eqn}
			L^2\left(\R^{2N}_\th\right) \subset C_0\left(\R^{2N}_\th\right).
		\end{equation} 
	\end{rem}
	\begin{rem}
		Notation of the Definition \ref{nt_r_2_N_repr_defn} differs from \cite{moyal_spectral}. Here symbols $A_\th, \A_\th, A^0_\th$ are replaced with $C_b\left(\R^{2N}_\th\right), \SS\left(\R^{2N}_\th\right), C_0\left(\R^{2N}_\th\right)$ respectively.
	\end{rem}
	\begin{rem}
		The $\C$-linear space $C_0\left(\R^{2N}_\th \right)$ is not isomorphic to $C_0\left(\R^{2N}\right)$. 
	\end{rem}

	There are elements $\left\{f_{nm}\in\SS\left(\R^{2} \right)\right\}_{m,n \in \N^0}$, described in \cite{varilly_bondia:phobos}, which satisfy to the following Lemma. 
	\begin{lem}\label{mp_osc_lem} {\rm\cite{varilly_bondia:phobos}}
		\label{lm:osc-basis}
		Let $m,n,k,l \in \N$. Then $f_{mn} \star_\th f_{kl} = \delta_{nk}f_{ml}$
		and $f_{mn}^* = f_{nm}$. Thus $f_{nn}$ is an orthogonal projection and
		$f_{mn}$ is nilpotent for $m \neq n$. Moreover,
		$\left\langle f_{mn}, f_{kl}\right\rangle  = 
		\delta_{mk}\,\delta_{nl}$. The family
		$\set{f_{mn} : m,n\in \N^0} \subset \sS\left( \mathbb{R}^{2}\right) \subset L^2(\R^{2})$ is an
		orthogonal basis.
	\end{lem}
	\begin{remark}
		One has 
		\be\label{nt_dmn_eqn}
		\int_{\R^2} f_{mn} = \delta_{mn}.
		\ee
	\end{remark}

	\begin{prop}\label{mp_fmn_prop}\cite{moyal_spectral,varilly_bondia:phobos}
		Let $N = 1$. Then $\SS\left(\R^{2N}_\th\right)=\SS\left(\R^{2}_\th\right) $ has a Fr\'echet algebra isomorphism with
		the matrix algebra of rapidly decreasing double sequences
		$c = (c_{mn})$ of complex numbers such that, for each $k \in \N$,
		\begin{equation}\label{mp_matr_norm}
			r_k(c) \stackrel{\mathrm{def}}{=} \biggl( \sum_{m,n=0}^\infty
			\th^{2k}  \left( m+\half\right)^k \left( n+\half\right)^k |c_{mn}|^2 \biggr)^{1/2}
		\end{equation}
		
		is finite, topologized by all the seminorms $(r_k)$; via the
		decomposition $f = \sum_{m,n=0}^\infty c_{mn} f_{mn}$ of~$\SS(\R^2)$ in
		the $\{f_{mn}\}$ basis.
		The twisted product $f \star_\th g$ is
		the matrix product $ab$, where
		\begin{equation}\label{mp_mult_eqn}
			\left( ab\right)_{mn} \stackrel{\mathrm{def}}{=} \sum_{k= 0}^{\infty} a_{\mu\nu}b_{kn}.
		\end{equation}
		
		For $N > 1$, $\Coo\left(\R^{2N}_\th\right)$ is isomorphic to the (projective) tensor product
		of $N$ matrix algebras of this kind, i.e.
		\begin{equation}\label{mp_tensor_prod}
			\SS\left(\R^{2N}_\th\right) \cong \underbrace{\SS\left(\R^{2}_\th\right)\otimes\dots\otimes\SS\left(\R^{2}_\th\right)}_{N-\mathrm{times}}
		\end{equation}
		with the projective topology induced by seminorms $r_k$ given by \eqref{mp_matr_norm}.	
	\end{prop}
	\begin{rem}
		If $A$ is  $C^*$-norm completion of the matrix algebra with the norm \eqref{mp_matr_norm} then $A \approx \mathcal K$, i.e.
		\begin{equation}\label{mp_2_eqn}
			C_0\left(\R^{2}_\th\right) \approx \mathcal K.
		\end{equation}
		Form \eqref{mp_tensor_prod} and \eqref{mp_2_eqn} it follows that
		\begin{equation}\label{mp_2N_eqn}
			C_0\left(\R^{2N}_\th\right) \cong \underbrace{C_0\left(\R^{2}_\th\right)\otimes\dots\otimes C_0\left(\R^{2}_\th\right)}_{N-\mathrm{times}} \approx \underbrace{\mathcal K\otimes\dots\otimes\mathcal K}_{N-\mathrm{times}} \approx \mathcal K
		\end{equation}
		where $\otimes$ means minimal or maximal tensor product ($\mathcal{K}$ is nuclear hence both products coincide).
		
	\end{rem}
	
	\begin{definition}\cite{varilly_bondia:phobos}
		For $s, t \in\R$ we denote by $\G_{s,t}$ the Hilbert space obtained by completing $\sS\left( \R^2 \right)$  with
		respect to the norm
		$$
		\left\|\sum_{m,n = 0}^\infty c_{mn} f_{mn}\right\|= \sqrt{\sum_{m,n = 0}^\infty \left(2m + 1\right)^s \left(2n + 1\right)^t \left| c_{mn}\right|^2 }
		$$
		
	\end{definition}
	
	Denote by 
	\be
	\mathcal{I}_{s,t}\bydef \G_{s,0}\times \G_{0,t} = \left\{f\times g~ \left|~ f\in \G_{s,0}, \quad g\in \G_{0,t}  \right|\right\}
	\ee

	\begin{empt}\label{nt_m_planes_empt}\cite{moyal_spectral}
		By plane waves we understand all functions of the form
		$$
		x \mapsto \exp(ik\cdot x) 
		$$
		for $k\in \R^{2N}$.  One obtains for the Moyal
		product of plane waves:
		\begin{equation}\label{mp_wave_prod_eqn}
			\begin{split}
				\exp\left(ik\cdot\right) \star_{\Theta}\exp\left(ik\cdot\right)=\exp\left(ik\cdot\right) \star_{\theta}\exp\left(ik\cdot\right)= \exp\left(i\left( k+l\right) \cdot\right) e^{-\pi i k \cdot \Th l}.
			\end{split}
		\end{equation}
		It is proven in \cite{moyal_spectral} that plane waves lie in $C_b\left(\R^{2N}_\th \right)$. 	
	\end{empt}

	\begin{empt}\label{mp_scaling_empt}

		Let us consider the unitary dilation operators $E_a$ given
		by
		$$
		E_af(x) \stackrel{\mathrm{def}}{=} a^{N/2} f(a^{1/2}x),
		$$
		It is proven in \cite{moyal_spectral} that
		\begin{equation}\label{eq:starscale}
			f {\star_{\theta}} g =
			(\th/2)^{-N/2} E_{2/\th}(E_{\th/2}f \star_2 E_{\th/2}g).
		\end{equation}
		We can simplify our construction by setting $\th = 2$. Thanks to
		the scaling relation~\eqref{eq:starscale} any qualitative result can is true if it is true in case of 
		$\th = 2$. We use the following notation
		\begin{equation}\label{mp_times_eqn}
			f {\times} g\stackrel{\mathrm{def}}{=}f {\star_{2}} g
		\end{equation}

		Introduce the symplectic Fourier transform $F$ by
		\begin{equation}
			Ff(x) \stackrel{\text{def}}{=} (2\pi)^{-N} \int f(t) e^{ix \cdot Jt} \,d^{2N}t; \quad 	\widetilde Ff(x) \stackrel{\text{def}}{=} (2\pi)^{-N} \int f(t) e^{ix \cdot Jt} \,d^{2N}t; 
			\label{eq:Fsympl}
		\end{equation}
		The \textit{twisted convolution} \cite{varilly_bondia:phobos} $f\diamond g$ is defined by
		\begin{equation}\label{nt_diamond_eqn}
			f\diamond g \left(u \right) \stackrel{\text{def}}{=} \int f\left(u - t\right)  g\left(t \right)e^{-iu \cdot J t}  dt.
		\end{equation}
		Following conditions hold \cite{varilly_bondia:phobos}:
		\bea\label{nt_prod_duality_eqn}
		\mathcal{F}\left(f \times g \right) = \mathcal{F} f \diamond  \sF g ;\quad \mathcal{F}\left(f \diamond g \right) = \mathcal{F} f \times \mathcal{F} g;\\
		\label{nt_prodf_duality}
		f \times g = Ff \diamond g = f \diamond  \widetilde Fg;  \quad f \diamond g = Ff \times g = f \times \widetilde F g;\\
		\label{nt_prod_ass_eqn}
		\left( f \times g \right) \times h= f \times \left(g \times h\right);  \quad \left( f \diamond g \right) \diamond h= f \diamond \left(g \diamond h\right);\\
		\label{nt_prod_*_eqn}
		\left( f \times g \right)^* = g^*\times h^*;  \quad \left( f \diamond g \right)^* = g^*\diamond h^*.
		\eea
		
	\end{empt} 
	\	
	
	\begin{defn}
		\label{df:Gst}\cite{moyal_spectral}
		We may as well introduce more Hilbert spaces $\G_{st}$ (for
		$s,t \in \R$) of those 
		$$
		f \in \SS'(\R^2) = \sum_{m,n = 0}^\infty c_{mn} f_{mn}
		$$ for which the following sum
		is finite:
		$$
		\|f\|_{st}^2 \stackrel{\mathrm{def}}{=} \sum_{m,n=0}^\infty
		(m+\half)^s (n+\half)^t |c_{mn}|^2.
		$$
		for~$\G_{st}$.
	\end{defn}
	
	\begin{rem}\label{mp_l2_rem}
		It is proven in \cite{varilly_bondia:phobos} $f, g \in L^2\left( \R^2\right)$, then $f\times g \in L^2\left( \R^2\right)$ and $\left\|f\times g  \right\|\le\left\|f \right\|\left\| g  \right\|$.
		Moreover, $f \times g$ lies in $C_0\left( \R^2\right)$ : the continuity follows by adapting the analogous argument for (ordinary) convolution. 
	\end{rem}
	\begin{rem}\label{mp_ss_rem}
		It is shown in \cite{varilly_bondia:phobos} that 
		\begin{equation}\label{mp_ssin_eqn}
			\SS\left( \R^2\right) = \bigcap_{s, t \in \R} \G_{st}.
		\end{equation}
	\end{rem}
	
	\begin{empt}\label{mp_coord_constr}
		This part contains a useful equations proven in \cite{varilly_bondia:phobos}.
		There are coordinate functions $p,q$ on $\R^2$ such that for any $f \in \SS\left( \R^2\right)$ following conditions hold
		\begin{equation}\label{mp_pq_mult_eqn}
			\begin{split}
				q \times f = \left(q + i \frac{\partial}{\partial p} \right)  f; ~~p \times f = \left(p - i \frac{\partial}{\partial q} \right)  f;\\
				f \times q = \left(q - i \frac{\partial}{\partial p} \right)  f; ~~f \times p = \left(p + i \frac{\partial}{\partial q} \right)  f.
			\end{split}
		\end{equation}
		From $q \times f, f \times q,p  \times f, f \times p \in \SS\left(\R^{2N} \right)$ it follows that $p, q \in \M^2$ (cf. \eqref{mp_multi_eqn}). From 
		\eqref{mp_mult_distr} it follows that
		\begin{equation}\label{mp_pq_mult}
			\begin{split}
				q \times \SS'\left(\R^{2N} \right) \subset \SS'\left(\R^{2N} \right); ~~p \times \SS'\left(\R^{2N} \right) \subset \SS'\left(\R^{2N} \right);\\
				\SS'\left(\R^{2N} \right) \times q \subset \SS'\left(\R^{2N} \right); ~~\SS'\left(\R^{2N} \right) \times p \subset \SS'\left(\R^{2N} \right).
			\end{split}
		\end{equation}
		If $f \in \SS'\left( \R^2\right)$ then from \eqref{mp_pq_mult_eqn} it follows that
		\begin{equation}\label{mp_partial_eqn}
			\frac{\partial}{\partial p} f = -iq \times f + i f \times q, ~~  \frac{\partial}{\partial q} f = ip \times f - i f \times p
		\end{equation}

		If 
		\begin{equation}\label{mp_ham_eqn}
			\begin{split}
				a \stackrel{\mathrm{def}}{=} \frac{q + ip}{\sqrt{2}}, ~~\overline{a} \stackrel{\mathrm{def}}{=} \frac{q - ip}{\sqrt{2}},\\
				\frac{\partial}{\partial a}\stackrel{\mathrm{def}}{=} \frac{\partial_q + i\partial_p}{\sqrt{2}}, ~~ \frac{\partial}{\partial \overline{a}}\stackrel{\mathrm{def}}{=} \frac{\partial_q - i\partial_p}{\sqrt{2}},\\
				H \stackrel{\mathrm{def}}{=}a \overline{a}= \frac{1}{2}\left(p^2 + q^2 \right) , \\
				\overline{a}\times a = H - 1, ~~ a \times \overline{a} = H + 1
			\end{split}
		\end{equation}
		then
		\begin{equation}\label{mp_ap_eqn}
			\begin{split}
				a \times f = af + \frac{\partial f}{\partial \overline{a}}, ~~ f \times a = af - \frac{\partial f}{\partial \overline{a}}, \\
				\overline{a} \times f = \overline{a}f - \frac{\partial f}{\partial a}, ~~ f \times a = \overline{a}f + \frac{\partial f}{\partial a}, \\
			\end{split}
		\end{equation} \begin{equation}\label{mp_ham_act_eqn}
			\begin{split}
				H\times f_{mn} = (2m + 1) f_{mn};~~f_{mn} \times H = 2(n + 1)f_{mn}
			\end{split}
		\end{equation} 
		\begin{equation}\label{mp_hamd_act_eqn}
			\begin{split}
				a \times f_{mn} = \sqrt{2m}f_{m-1,n};~~ f_{mn}\times a = \sqrt{2n + 2} f_{m, n+1};\\
				\overline{a} \times f_{mn} = \sqrt{2m+2}f_{m+1,n};~~ f_{m+1,n}\times \overline{a} = \sqrt{2n} f_{m, n-1}.
			\end{split}
		\end{equation}

		It is proven in \cite{varilly_bondia:phobos} that
		\begin{equation}\label{mp_part_prod_eqn}
			\partial_j\left(f \times g \right)  = \partial_j f \times g + f \times \partial_jg;
		\end{equation}
		where $\partial_j = \frac{\partial}{ \partial x_{j}}$ is the partial derivation in $\SS\left(\R^{2N} \right)$. 
	\end{empt}
	
	From \eqref{mp_ham_eqn} it follows that 
	$$
	q =  \sqrt{2}\left( a   +  \overline{a}\right), \quad p =  \sqrt{2} i\left( a   -  \overline{a}\right),
	$$
	so taking into account \eqref{mp_partial_eqn} one has
	\be\label{mp_dpq_e_eqn}
	\begin{split}
		\frac{\partial}{\partial p} f = \sqrt{2}\left( -i\left( a   +  \overline{a}\right) \times f +if \times \left( a   +  \overline{a}\right)\right) \\
		\frac{\partial}{\partial q} f = \sqrt{2}\left( -\left( a   -  \overline{a}\right) \times f -f \times \left( a   +  \overline{a}\right)\right).
	\end{split}
	\ee

\chapter{Miscellany}
	
\section{Some quantum groups}	
	\subsection{Quantum $SU\left(2 \right)$ and $SO\left(3 \right)$}
	\paragraph*{}
	There is a quantum generalization of $SU\left(2 \right)$ and  we will introduce a quantum analog of $SO\left(3 \right)$.
	Let $q$ be a real number such that $0<q<1$. 
	A quantum group $C\left( \SU_q(2)\right) $ is the universal $C^*$-algebra algebra generated by two elements $\al$ and $\beta$ satisfying the following relations:
	\begin{equation}\label{su_q_2_rel_eqn}
		\begin{split}
			\al^*\al + \beta^*\beta = 1, ~~ \al\al^* + q^2\beta\beta^* =1,
			\\
			\al\bt - q \bt\al = 0, ~~\al\bt^*-q\bt^*\al = 0,
			\\
			\bt^*\bt = \bt\bt^*.
		\end{split}
	\end{equation}
	From  $C\left( SU_1\left(2 \right)\right) \approx C\left(SU\left(2 \right)  \right)$ it follows that  $C\left( \SU_q(2)\right) $ can be regarded as a noncommutative deformation of the space $SU(2)$. 
	The dense pre-$C^*$-algebra $\Coo\left( SU_q(2)\right) \subset C\left( SU_q(2)\right)$ is defined in  \cite{chakraborty_pal:quantum_su_2}. 
	Let $Q, S \in B\left( \ell_2\left(\N^0 \right)\right) $ be given by 
	\begin{equation*}
		\begin{split}
			Qe_k= q^ke_k, \\
			Se_k = \left\{
			\begin{array}{c l}
				e_{k-1} & k > 0 \\
				0 & k = 0
			\end{array}\right..
		\end{split}
	\end{equation*}
	and let $R \in B\left( \ell_2\left(\Z \right)\right) $ be given by $e_k \mapsto e_{k+1}$.
	There is a faithful representation \cite{woronowicz:su2}  $C\left(\SU_q\left( 2\right) \right) \to B\left(\ell_2\left(\N^0 \right) \otimes \ell_2\left(\Z \right) \right)  $ given by
	
	\begin{equation}\label{su_q_2_repr_eqn}
		\begin{split}
			\al \mapsto S\sqrt{1 - Q^2} \otimes 1, \\
			\bt \mapsto Q \otimes R.
		\end{split}
	\end{equation}
	There is a faithful state $h:C\left( SU_q(2)\right) \to \C$ given by
	\be\label{su_q_2_haar_eqn}
	h\left(a \right)  = \sum_{n = 0}^\infty q^{2n}\left(e_n \otimes e_0, a  e_n \otimes e_0\right) 
	\ee
	where $a \in C\left( SU_q(2)\right)$ and $e_0 \otimes e_n \in \ell_2\left(\N^0 \right) \otimes \ell_2\left(\Z \right)$ (cf. \cite{woronowicz:su2}).
	\begin{defn}
		The state $h$ is said to be the \textit{Haar measure}.
	\end{defn}
	
	Denote by $L^2\left( C\left( SU_q\left(2\right)\right), h\right) $ the GNS space associated with the state $h$. 
	The representation theory of $SU_q(2)$ is strikingly similar to its
	classical counterpart.  In particular, for each $l\in\{0,\frac{1}{2},
	1,\ldots\}$, there is a unique irreducible unitary representation
	$t^{(l)}$ of dimension $2n+1$.  Denote by $t^{(l)}_{jk}$ the
	$jk$\raisebox{.4ex}{th} entry of $t^{(l)}$. These are all elements of
	$\A_f$ and they form an orthogonal basis for $L^2\left( C\left( SU_q\left(2\right)\right), h\right)$. Denote by
	$e^{(l)}_{jk}$ the normalized $t^{(l)}_{jk}$'s, so that
	$\{e^{(l)}_{jk}: n=0,\frac{1}{2},1,\ldots, i,j=-n,-n+1,\ldots, n\}$ is
	an orthonormal basis. The definition of equivariant operators (with respect to action of quantum groups) is described in \cite{chakraborty_pal:inv_hom}. It is proven in \cite{chakraborty_pal:quantum_su_2} that any unbounded equivariant operator $\widetilde{D}$ satisfies to the following condition
	\be \label{su_q_2_dirac_eqn}
	\widetilde{D}: e^{(l)}_{jk}\mapsto d(l,j)e^{(l)}_{jk},
	\ee
	Moreover if
	\be \label{su_q_2_genericd_eqn}
	d(l,j)=\begin{cases}2l+1 &  l \neq j,\cr
		-(2l+1) & l=j,\end{cases}
	\ee
	
	then there is a 3-summable spectral triple 
	\be \label{su_q_2_spt_eqn}
	\left(\Coo\left( SU_q\left(2\right)\right), L^2\left( C\left( SU_q\left(2\right)\right), h\right), \widetilde{D} \right) 
	\ee
	described in \cite{chakraborty_pal:quantum_su_2}. 
	
	According to \cite{kl-sch} (equations (4.40)-(4.44) ) following condition holds
	\begin{equation}\label{su_q_2_tij_eqn}
		\begin{split}
			t^{\left(l \right) }_{jk}= M^l_{jk} \al^{-j-k}\bt^{k-j}p_{l + k}\left( \bt\bt^*; q^{-2\left(k-j \right) }q^{2\left(j + k \right)}~|~q^2 \right); ~~ j + k \le 0 ~\&~ k \ge j,
			\\ 
			t^{\left(l \right) }_{jk}= M^l_{jk} \al^{-j-k}\bt^{*\left( k-j\right) }p_{l + k}\left( \bt\bt^*; q^{-2\left(k-j \right) }q^{2\left(j + k \right)}~|~q^2 \right); ~~ j + k \le 0 ~\&~ k \le j,
			\\ 
			t^{\left(l \right) }_{jk}= M^l_{j,k}p_{l - k}\left( \bt\bt^*; q^{-2\left(k-j \right) }q^{2\left(j + k \right)}~|~q^2 \right)\bt^{k-j}\al^{*j+k}; ~~ j + k \ge 0 ~\&~ k \ge j,\\ 
			t^{\left(l \right) }_{jk}= M^l_{k,j}p_{l - j}\left( \bt\bt^*; q^{-2\left(k-j \right) }q^{2\left(j + k \right)}~|~q^2 \right)\bt^{*j-k}\al^{*j + k}; ~~ j + k \ge 0 ~\&~ j \ge k\\ 
		\end{split}
	\end{equation}
	where $M^l_{jk} \in \R$ for any $l, j, k$ and $p_{l - k}\left( x; q^{-2\left(k-j \right) }q^{2\left(j + k \right)}~|~q^2 \right)$ is little Jacobi polynomial (cf. \cite{kl-sch}).  
	Denote by $g \in \Z_2$ the unique nontrivial element.
	There is a surjective group homomorphism
	$$
	\Phi : SU\left(2 \right) \to SO\left(3 \right) ,~ \ker \Phi = \Z_2 = \{\pm 1\}
	$$
	and the natural action of $\Z_2$ on $SU(2)$ such that
	\be\label{su_q_2_z2_comm_eqn}
	\begin{split}
		SO\left(3 \right) \cong SU\left(2\right)/\Z_2,\\
		g\begin{pmatrix} \al & -\overline{\bt}\\
			\bt & \overline{\al}\end{pmatrix}= \begin{pmatrix} -\al & \overline{\bt}\\
			-\bt & -\overline{\al}\end{pmatrix};\quad \forall \begin{pmatrix} \al & -\overline{\bt}\\
			\bt & \overline{\al}\end{pmatrix}\in SU\left(2\right).
	\end{split}
	\ee	
	This action induces an action of $\Z_2$  on a $C^*$-algebra $C\left(SU\left(2\right) \right)$ given by
	\be\nonumber
	\begin{split}
		g\al = -\al, ~~ g\bt = - \bt
	\end{split}
	\ee
	where $\al, \bt$ are regarded as functions $\SU\left(2 \right) \to \C$.
	Indeed 	$ SU\left(2 \right)$ is an oriented manifold, $SO\left(3 \right)$ is an unoriented  one, and $ SU\left(2 \right) \to SO\left(3 \right)$ is a two-fold covering
	There is a quantum generalization of $SU\left(2 \right)$ and  we will introduce a quantum analog of $SO\left(3 \right)$.
	Let $q$ be a real number such that $0<q<1$. 
	There is a noncommutative analog of the action \eqref{su_q_2_z2_comm_eqn} described in \cite{dijkhuizen:so_doublecov,lance:so,podles:so_su}
	\begin{equation}\label{su_q_2_z2_ncomm_eqn}
		\begin{split}
			\Z_2 \times C\left( SU_q\left(2\right)\right)\to C\left( SU_q\left(2\right)\right),\\
			g\al = -\al,~g \bt = -\bt.
		\end{split}
	\end{equation} 
	In \cite{dijkhuizen:so_doublecov,lance:so,podles:so_su} the quantum group 	$SO_q\left(3 \right)$ is defined, moreover in 
	it is proven in \cite{dijkhuizen:so_doublecov,podles:so_su} following 
	$$
	C\left( 	SO_q\left(3 \right) \right) \stackrel{\mathrm{def}}{=} C\left( SU_q\left(2\right)\right) ^{\Z_2}\cong \left\{\widetilde{a}\in C\left( SU_q\left(2\right)\right) ,~ g\widetilde{a} = \widetilde{a} \right\}.
	$$
	Above equation can be used as the definition of $SO_q\left(3 \right)$
	\begin{definition}\label{su_q_2_so3_defn}
		Denote by
		\begin{equation}\label{su_q_2_so3_eqn}
			C\left( 	SO_q\left(3 \right) \right) \stackrel{\mathrm{def}}{=} C\left( SU_q\left(2\right)\right) ^{\Z_2}\cong \left\{\widetilde{a}\in C\left( SU_q\left(2\right)\right) ,~ g\widetilde{a} = \widetilde{a} \right\}.
		\end{equation}
		The $C^*$-algebra $C\left( 	SO_q\left(3 \right)\right) $ is said to be the \textit{quantum 	$ SO\left(3 \right)$}.
	\end{definition}
	\begin{remark}
		Our definition of $SO_q\left(3 \right)$  differs from the Definition given in other sources. For example the quantum group defined in \cite{rtf:qlie} is a
		quantization of $O\left( 3\right)$  and not of $SO(3)$ (cf. \cite{dijkhuizen:so_doublecov}). 
	\end{remark}
	
	\begin{thm}\label{su_q_2_bas_thm}\cite{woronowicz:su2}
		Let  $q \neq 0$, and let $\A_f$ the dense involutive subalgebra of $C\left( SU_q(2)\right)$ generated  by $\alpha, \beta$. The set of elements of the form
		\be\label{su_q_2_fin_eqn}
		\al^k\bt^n\bt^{*m}~~ \text{ and }~~ \al^{*k'}\bt^n\bt^{*m}
		\ee 
		where $k, m, n = 0, \dots;~k'=1, 2,\dots$ forms a basis in $\A_f$: any element
		of $\A_f$ can be written in the unique way as a finite linear combination of
		elements of \eqref{su_q_2_fin_eqn}.
	\end{thm}
	
\subsection{The group of motions of the Euclidean quantum plane}

\paragraph{}The group of motions of the Euclidean quantum plane is represented by a *-algebra of unfounded operators. Here I follow to \cite{woronowicz:unb_affil}. Let us fix a real number $\mu$ such that $\left|\mu \right| > 1$; $\mu=1$ corresponds to the classical (i.e. non-quantum) case. In \cite{woronowicz:unb_affil} the quantum group $G_\mu$ is introduced. Let $A_\mu$  be the $*$-algebra generated by two elements $u$ and $n$ such that
\be\label{em_eqn}
\begin{split}
	v^*v = vv^* = 1_{A_\mu},\\
	n^*n = nn^*,\\
	v^*nv= \mu n
\end{split}
\ee
(cf. \cite{woronowicz:unb_affil} equation (3.3)). Elements of $A_\mu$ are called polynomials on $G_\mu$.  Let $u$ and $n$ be (unbounded) operators acting on a Hubert space $\H$. We say that the pair $\left(v, n\right)$ is a
representation of commutation relations \eqref{em_eqn} if $v$ is unitary, $n$ is normal and
$v^*nv= \mu n$. We recall that the normality of n means that $\Dom\left(n^*\right)= \Dom\left(n\right)$ and $\left\|n^*\psi \right\|= \left\|n\psi \right\|$ for any $\psi \in \Dom\left(n\right)$. One can
easily verify the following facts: As usual any representation is a direct integral of
irreducible ones. Any representation is either infinite or one-dimensional. All one-dimensional
representations are of the form $\left(c 1_{\mathscr A}, 0\right)$, where $c\in\C$, $\left|c \right|=1$. If $\left( v, n\right)$ is an irreducible infinite-dimensional representation then $\mathrm{Sp}~ n = \left\{0\right\}\cup \Ga$,
where $\Ga$ is of the form
\be\label{em_sp_eqn}
\Ga \bydef \left\{\left.t_0\mu^k \right|k\in \Z\right\}
\ee
and $t_0\in \C\setminus\{0\}$. Moreover, one can find an orthonormal basis $\left\{\left|t\right\rangle| t\in\Ga\right\}$ such
that
\bea\label{em_vb_eqn}
v\left|t\right\rangle = \mu\left|t\right\rangle,\\
\label{em_nb_eqn} 
n\left|t\right\rangle = t\left|t\right\rangle 
\eea
for any $t \in \Ga$. It means that $\left( v, n\right)$ is uniquely determined (up to a unitary
equivalence) by $\mathrm{Sp}~ n$ In \cite{woronowicz:unb_affil} the following condition is introduced
\be\label{em_sc_eqn}
\mathrm{Sp}~n\subset \C_{\left( \mu\right) } \bydef \{0\} \cup \left\{\left.\xi\in \C\right|\exists k \in \Z\quad \left|\xi \right|= \mu^k\right\}
\ee
We call \eqref{em_sc_eqn} the \textit{spectral condition}.
Let $C_c\left(\C_{\left( \mu\right) } \right)$  be the algebra of all continuous functions with compact support on
$\C_{\left( \mu\right) }$. There is a natural action of 
\bean
\begin{split}
	\Z\times C_c\left(\C_{\left( \mu\right) } \right)\to \C_{\left( \mu\right) },\\
	\forall k \in \Z\quad \forall f \in C_c\left(\C_{\left( \mu\right) } \right)\quad \forall z \in \C_{\left( \mu\right) }\quad \left( k. f\right) \left(z \right) \bydef f\left(\mu^k z \right) 
\end{split}
\eean
If   
$C_c\left( \Z, C_c\left(\C_{\left( \mu\right) } \right) \right)$
is given by \eqref{discr_cr_prod_cc_eqn} then $C_c\left( \Z, C_c\left(\C_{\left( \mu\right) } \right) \right)$ is a *-algebra with given by \eqref{discr_cr_prod_op_eqn} operations. If $Cont\left( \C_{\left( \mu\right) }\right)$ is a *-algebra of all  (possibly unbounded) continuous maps on $\C_{\left( \mu\right) }$ then there is an action $Cont\left( \C_{\left( \mu\right) }\right)\times C_c\left( \C_{\left( \mu\right) }\right)\to C_c\left( \C_{\left( \mu\right) }\right)$. Denote by 
\bean
N \in Cont\left( \C_{\left( \mu\right) }\right),\\
N\left(z \right)= z.
\eean
From the above equation if follows that $\C_{\left( \mu\right) }\subset C_c\left( \Z, C_c\left(\C_{\left( \mu\right) } \right) \right)$, using this circumstance one can define
\be\label{em_Na_eqn}
\forall a \in C_c\left( \Z, C_c\left(\C_{\left( \mu\right) }\right)  \right)\quad N a\in  C_c\left( \Z, C_c\left(\C_{\left( \mu\right) } \right) \right) 
\ee

\begin{exercise}\label{dm_exer}
	Prove that there is an action $A_\mu \times  C_c\left( \Z, C_c\left(\C_{\left( \mu\right) }\right)  \right)\to  C_c\left( \Z, C_c\left(\C_{\left( \mu\right) }\right)  \right)$ such that for any $a \in  C_c\left( \Z, C_c\left(\C_{\left( \mu\right) }\right)  \right)$ one has
	\be\label{dm_act_eqn}
	\begin{split}
		na \bydef Na \quad \text{(cf. \eqref{em_Na_eqn})},\\
		v^k a \bydef k. a.
	\end{split}
	\ee
\end{exercise}

\section{Isospectral deformations}
	\paragraph*{}A very general construction of isospectral
	deformations
	of noncommutative geometries is described in \cite{connes_landi:isospectral}. The construction
	implies in particular that any
	compact Riemannian manifold $M$ which  admits a Spin$^c$ structure (cf. Definition \ref{spin_str_defn}), whose isometry group has rank
	$\geq 2$ admits a
	natural one-parameter isospectral deformation to noncommutative geometries
	$M_\theta$.
	We let $(\Coo\left(M \right)  , \H = L^2\left(M,S \right)  , \slashed D)$ be the canonical spectral triple associated with a
	compact spin-manifold $M$. We recall that $\mathcal{A} = C^\infty(M)$ is
	the algebra of smooth
	functions on $M$, $S$ is the spinor bundle and $\slashed D$
	is the Dirac operator.
	Let us assume that the group $\mathrm{Isom}(M)$ of isometries of $M$ has rank
	$r\geq2$.
	Then, we have an inclusion
	\begin{equation}\label{isos_t_act_eqn}
		\mathbb{T}^2 \subset \mathrm{Isom}(M) \, ,
	\end{equation}
	with $\mathbb{T}^2 = \mathbb{R}^2 / 2 \pi \mathbb{Z}^2$ the usual torus, and we let $U(s) , s \in
	\mathbb{T}^2$, be
	the corresponding unitary operators in $\H = L^2(M,S)$ so that by construction
	\begin{equation*}
		U(s) \, \slashed D = \slashed D \, U(s).
	\end{equation*}
	Also,
	\begin{equation}\label{isospectral_sym_eqn}
		U(s) \, a \, U(s)^{-1} = \alpha_s(a) \, , \, \, \, \forall \, a \in \mathcal{A} \, ,
	\end{equation}
	where $\alpha_s \in \mathrm{Aut}(\mathcal{A})$ is the action by isometries on the
	algebra of functions on
	$M$.
	
	\noindent
	We let $p = (p_1, p_2)$ be the generator of the two-parameters group $U(s)$
	so that
	\begin{equation*}
		U(s) = \exp(i(s_1 p_1 + s_2 p_2)) \, .
	\end{equation*}
	The operators $p_1$ and $p_2$ commute with $D$.
	Both $p_1$ and $p_2$
	have integral spectrum,
	\begin{equation*}
		\mathrm{Spec}(p_j) \subset \mathbb{Z} \, , \, \, j = 1, 2 \, .
	\end{equation*}
	
	\noindent
	One defines a bigrading of the algebra of bounded operators in $\H$ with the
	operator $T$ declared to be of bidegree
	$(n_1,n_2)$ when,
	\begin{equation*}
		\alpha_s(T) = \exp(i(s_1 n_1 + s_2 n_2)) \, T \, , \, \, \, \forall \, s \in
		\mathbb{T}^2 \, ,
	\end{equation*}
	where $\alpha_s(T) = U(s) \, T \, U(s)^{-1}$ as in \eqref{isospectral_sym_eqn}.
	\paragraph{}
	Any operator $T$ of class $C^\infty$ relative to $\alpha_s$ (i. e. such that
	the map $s \rightarrow \alpha_s(T) $ is of class $C^\infty$ for the
	norm topology) can be uniquely
	written as a doubly infinite
	norm convergent sum of homogeneous elements,
	\begin{equation*}
		T = \sum_{n_1,n_2} \, \widehat{T}_{n_1,n_2} \, ,
	\end{equation*}
	with $\widehat{T}_{n_1,n_2}$ of bidegree $(n_1,n_2)$ and where the sequence
	of norms $||
	\widehat{T}_{n_1,n_2} ||$ is of
	rapid decay in $(n_1,n_2)$.
	Let $\lambda = \exp(2 \pi i \theta)$. For any operator $T$ in $\H$ of
	class $C^\infty$ we define
	its left  twist $l(T)$ by
	\begin{equation}\label{l_defn}
		l(T) = \sum_{n_1,n_2} \, \widehat{T}_{n_1,n_2} \, \lambda^{n_2 p_1} \, ,
	\end{equation}
	and its right twist $r(T)$ by
	\begin{equation*}
		r(T) = \sum_{n_1,n_2} \, \widehat{T}_{n_1,n_2} \, \lambda^{n_1 p_2} \, ,
	\end{equation*}
	Since $|\lambda | = 1$ and $p_1$, $p_2$ are self-adjoint, both series
	converge in norm. Denote by $\Coo\left(M \right)_{n_1, n_2} \subset \Coo\left(M \right) $ the $\C$-linear subspace of elements of bidegree $\left( n_1, n_2\right) $. \\
	One has,
	\begin{lem}\label{conn_landi_iso_lem}\cite{connes_landi:isospectral}
		\begin{itemize}
			\item[{\rm a)}] Let $x$ be a homogeneous operator of bidegree $(n_1,n_2)$
			and $y$ be
			a homogeneous operator of  bidegree $(n'_1,n'_2)$. Then,
			\begin{equation}
				l(x) \, r(y) \, - \,  r(y) \, l(x) = (x \, y \, - y \, x) \,
				\lambda^{n'_1 n_2} \lambda^{n_2 p_1 + n'_1 p_2}
			\end{equation}
			In particular, $[l(x), r(y)] = 0$ if $[x, y] = 0$.
			\item[{\rm b)}] Let $x$ and $y$ be homogeneous operators as before and
			define
			\begin{equation*}
				x * y = \lambda^{n'_1 n_2} \, x y \, ; \label{star}
			\end{equation*}
			then $l(x) l(y) = l(x * y)$.
		\end{itemize}
	\end{lem}

	\noindent
	The product $*$ defined in (\ref{star}) extends by linearity
	to an associative product on the linear space of smooth operators and could
	be called a $*$-product.
	One could also define a deformed `right product'. If $x$ is homogeneous of
	bidegree
	$(n_1,n_2)$ and $y$ is homogeneous of bidegree $(n'_1,n'_2)$ the product is
	defined by
	\begin{equation*}
		x *_{r} y = \lambda^{n_1 n'_2} \, x y \, .
	\end{equation*}
	Then, along the lines of the previous lemma one shows that $r(x) r(y) = r(x
	*_{r} y)$.
	
	We can now define a new spectral triple where both $\H$ and the operator
	$\slashed D$ are unchanged while the
	algebra $\Coo\left(M \right)$  is modified to $l(\Coo\left(M \right))$ . By
	Lemma~{\ref{conn_landi_iso_lem}}~b) one checks that  $l\left( \Coo\left(M \right)\right) $ is still an algebra. Since $\slashed D$ is of bidegree $(0,0)$ one has,
	\begin{equation}\label{isospectral_ld_eqn}
		[\slashed D, \, l(a) ] = l([\slashed D, \, a]) 
	\end{equation}
	which is enough to check that $[\slashed D, x]$ is bounded for any $x \in l(\mathcal{A})$. There is a spectral triple $\left(l\left( \Coo\left(M \right)\right) , \H, \slashed D\right)$.
	\paragraph{} Denote by $\Coo\left( M_\th\right)$ (resp.  $C\left(M_\th \right)$) the algebra $l\Coo\left( M\right)$ (resp. the operator norm completion of $l\left(\Coo\left( M\right)  \right)  $). Denote by $\rho: C\left(M\right) \to L^2\left( M, S\right) $ (resp. $\pi_\th: C\left(M_\th\right) \to B\left( L^2\left( M, S\right)\right) $ ) natural representations.
	There is an oriented twisted spectral triple 
	\be\label{isos_twisetd_eqn}
	\left(l\left( \Coo\left(M \right)\right) , \H, \slashed D\right).
	\ee
	described in \cite{connes_landi:isospectral}.

\section{Miscellany}
	
	\paragraph*{} Our constructions require an analog of Dini's theorem. Below there are  citations of  Dini's theorem it and its proof.
	\begin{theorem}\label{dini_thm}\cite{rudin:pa}
		Suppose $K$ is compact, and 
		\begin{enumerate}
			\item[(a)]  $\left\{f_n\right\}_{n \in \N}$ is a sequence of continuous functions on $K$, 
			\item[(b)]  $\left\{f_n\right\}$ converges point-wise to a continuous function $f$ on $K$, 
			\item[(c)] $f_n(x)\ge f_{n+l}(x)$ for all $x\in K$, $n = 1, 2, 3,...$. 
		\end{enumerate}
		Then $f_n \to f$ uniformly on $K$. 
	\end{theorem}
	\begin{proof}
		Put $g_n\bydef f_n-f$. Then $g_n$ is continuous, $g_n\to 0$ pointwise, and $g_n\ge g_{n+1}$. We have to prove that $g_n \to 0$ uniformly on $K$. 
		Let $\eps> 0$ be given. Let $K_n$ be the set of all $x\in K$ with $g_n\left( x\right) > \eps$.  
		Since $g_n$ is continuous, $K_n$ is closed 
		hence compact.
		Since $g_n > g_{n+1}$, we have $K_n\supset K_{n+`1}$. Fix $x\in K$. Since $g_n\left(x\right) \to x$
		we see that $x\notin K_n$ if $n$ is sufficiently large. Thus $x \notin \bigcap K_n$. In other words, 
		the intersection $\bigcap K_n$ is empty. Hence $K_N$ is empty for some $N$. 
		It follows that  $0\le g\left( x\right) < \eps$ for all $x \in K$  and for all $n\ge N$. This proves the theorem.
	\end{proof}
	
	\
	
	
	
	\end{appendices}

\part{Applications to the Theoretical Physics}\label{phys_part}

\chapter{Hausdorff non-Hausdorff duality}\label{haus_non_haus_dul_sec}
\paragraph*{}

The most impressive directions of present day theoretical physics are superstring theory and noncommutative geometry. Both directions  are aimed to describe fundamental physical laws, i.e. these directions have the single purpose. However these directions seem very different. The superstring theory operates with the classical geometry of Hausdorff spaces (cf. Definition \ref{haudorff_defn}). Noncommutative geometry replaces spaces with $C^*$-algebras, which simultaneously describe both "outer space" (stars and galaxies) and "inner space" (fundamental particles and forces). The spectra of $C^*$-algebras should or  should not be Hausdorff. In \cite{connes_lott:particle} the Standard model is described by an operator algebra which generates a $C^*$-algebra with Hausdorff spectrum (cf. Section \ref{oa_sec}). In \cite{connes_rieffel:nc_ym} a physical model is described by a $C^*$-algebra with non-Hausdorff spectrum.
\paragraph*{}
The present day theoretical physics operates with models having  unbounded dimension of the internal physical space. The relationship between large $N$ matrix models and noncommutative geometry in
string theory is described in \cite{wittenD,bfss,bss,cds,ikkt,sw1}.
Here examples   of a duality between Hausdorff and non-Hausdorff geometry are presented. The Hausdorff and non-Hausdorff
geometry can be  manifestations of the single physical theory. All examples come from  having  unbounded dimension
 models,

\section{Prototype. Matrix model and noncommutative torus}\label{phys_mm_sec}
\paragraph{} 
The paper \cite{landi:nm_nt} can be regarded as a manifestation of  Hausdorff and non-Hausdorff duality. The infinite sequence of $C^*$-algebras $\Coo\left(\T_{M, N} \right)$ yield $C\left(\T^2_\th \right)$ where $\th \in \R \setminus \Q$.
The algebra $\Coo\left(\T^2_\th\right)$ of smooth functions on the `noncommutative
two-torus'
${\bf T}_\theta^2$ is the unital $*$-algebra generated by two unitary
elements $U_1,
U_2$ with the relation
\be U_1U_2=e^{2\pi i\theta}\,U_2U_1 \ .
\label{Tthetarel}
\ee
A generic element $a\in \Coo\left(\T^2_\th\right)$ is written as a convergent series
of the form
\be
a = \sum_{(m,n)\in \\Z^2} a_{mn} ~(U_1)^m (U_2)^n
\ee
where $a_{mn}$ is a complex-valued Schwarz function on $\Z^2$, i.e. a
sequence
of complex numbers $\{a_{mn} \in  \C~| ~ (m,n) \in \Z^2 \}$
which decreases rapidly at `infinity'. When the deformation parameter
$\theta=M/N$ is a rational number, with $M$ and $N$
positive integers which we take to be relatively prime, the algebra
$\Coo\left(\T^2_{M/N}\right)$ is intimately
related to the algebra $C^\infty({\bf T}^2)$ of smooth functions on the
ordinary torus
${\bf T}^2$. Precisely, $\Coo\left(\T^2_{M/N}\right)$ is Morita equivalent to
$C^\infty({\T}^2)$, i.e. $\Coo\left(\T^2_{M/N}\right)$ is a twisted matrix bundle over
$C^\infty({\bf T}^2)$ of topological charge $M$ whose fibers are $N\times N$
complex matrix algebras.
Physically, this  implies that noncommutative $U(1)$ Yang-Mills theory with
rational deformation parameter
$\theta=M/N$ is dual to a conventional $U(N)$ Yang-Mills theory with $M$
units of 't~Hooft flux.

The algebra $\Coo\left(\T^2_{M/N}\right)$ has a
`huge' center ${\cal C}(\Coo\left(\T^2_{M/N}\right))$ which is generated by the elements
$(U_1)^N$ and $(U_2)^N$. One identifies ${\cal C}(\Coo\left(\T^2_{M/N}\right))$ with the algebra
$C^\infty({\T}^2)$, while the appearence of finite dimensional matrix
algebras can be seen as follows. With $\omega=e^{2\pi iM/N}$, one introduces
the $N\times N$ clock and shift matrices
\be\label{Uarational}
\widetilde{U}_1=\left({\begin{array}{lllll}
		1& & & & \\ &\omega& & & \\
		& &\omega^2& & \\& &
		&\ddots& \\ & & & &
		\omega^{N-1}
\end{array}}\right)~~~~~~,~~~~~~
\widetilde{U}_2=\left({\begin{array}{lllll}
		0&1& & &0\\ &0&1& & \\
		& &\ddots&\ddots& \\
		& & &\ddots&1\\ 1& & & &0\end{array}}\right) \ .
\ee
These matrices are traceless (since $\sum_{k=0}^{N-1} \omega^k = 0$), they obey
the relation \ref{Tthetarel}, and they satisfy
\be
\left(\widetilde{U}_1\right)^N =\left(\widetilde{U}_2\right)^N = \id_N~.
\ee
Since $M$ and $N$ are relatively prime, the matrices \ref{Uarational} generate
the finite dimensional algebra $\mathbb{M}_N(\C)$ of $N\times N$ complex
matrices. In \cite{landi:nm_nt} the following representation
\be
\theta=\lim_{n\to\infty}\theta_n=\lim_{n\to\infty}p_n/q_n;\quad \th \in \R\setminus \Q,\quad p_n \in \Z, \quad q_n \in \N 
\label{thetalim}
\ee
is used.
Then, for any integer $n$, the semi-simple algebra
\be
A_n=\mathbb{M}_{q_n}(\C)\oplus\mathbb{M}_{q_{n-1}}(\C)
\label{Andef}
\ee
is considered.
At each finite level labelled by the integer $n$, let $A_{\theta_n}$ be
the algebra of
the noncommutative two-torus with rational deformation parameter
$\theta_n= p_n / q_n$ given in \ref{thetalim}, and generators
$U_a^{(n)}~, ~a=1,2$ obeying the relation
\be\label{hnh_comm_matr_eqn}
U_1^{(n)}U_2^{(n)}=e^{2\pi ip_n/q_n}\,U_2^{(n)}U_1^{(n)} 
\ee
Having the sequence of $C^*$-algebras $A_n$ and $A_{\th_n}$ one obtains $C\left(\T^2_n \right)$. The complete description of this procedure is presented in \cite{landi:nm_nt}. 
\section{Application of noncommutative coverings}
\paragraph{}
Here we prove that described in \cite{landi:nm_nt}  results come from noncommutative coverings.
\subsection{Finite-fold coverings}\label{phys_fin_sec}
\paragraph{}  The described in \cite{landi:nm_nt} construction of based on the Morita equivalence between $C^*$-algebras $\Coo\left(\T^2_{M/N}\right)$ and
$C^\infty({\T}^2)$. We shall use the Morita equivalence given by the Theorem \ref{finite_morita_main_theorem}. If $\sX \cong S^1$ and $p: \sX\to S^1$ is an $N$-listed covering then there is an unital noncommutative finite-fold covering $\left(C\left(S^1\right), C\left(\sX\right), \Z_N, C_0\left(p\right)\right)$. If $u\in C\left(S^1\right)$ is an untary such that $C\left(S^1\right)\cong C\left(u\right)$ then is an untary $v\in C\left(\sX\right)$ such that $C\left(\sX\right)\cong C\left(v\right)$ and $C_0\left(p\right)$ is given by $u \mapsto v^N$. From the Theorem 
\ref{finite_morita_main_theorem} it follows that the algebra $C\left(S^1\right) $ is Morita equivalent to the algebra $C\left(\sX\right) \rtimes \Z_N$. The action $\Z_N \times C\left(\sX\right)\to C\left(\sX\right)$ is given by
$$
\forall m \in \Z\quad \overline k\cdot v^m = e^{\frac{2\pi i km}{N}} v^m\quad k \in \Z \text{ is a repesentative of } \overline k.
$$ If $\a\in \Aut\left(C\left(\sX\right) \right) $ corresponds to represented by $M$ element $\overline M \in \Z_N$ then one has 
\be\label{phys_Uarational}
v \a  = e^{\frac{2\pi i M}{N}}\a v
\ee
On the other hand if $\widetilde{U}_1, ~\widetilde{U}_2$ are given by \eqref{Uarational} then one has
$$
\widetilde{U}_1\widetilde{U}_2  = e^{\frac{2\pi i M}{N}}\widetilde{U}_2\widetilde{U}_1.
$$
The above equation is a counterpart of \eqref{phys_Uarational}. If one replaces $\widetilde{U}_1, ~\widetilde{U}_2$ with $v, ~ \a\in C\left(\sX\right) \rtimes \Z_N$ then one can give a based on the theory of coverings version of explained in \cite{landi:nm_nt} theory.
\subsection{Infinite coverings}
\paragraph{}

Here we would like reformulate the given by \ref{phys_mm_sec} formalism in terms of noncommutative coverings. Consider a sequences $\left\{U_1^{(n)}\right\}_{n\in \N}$ $~\left\{U_2^{(n)}\right\}_{n \in \N}$ of matrices which satisfy to the equation
\eqref{hnh_comm_matr_eqn}, i.e. 
\bean
U_1^{(n)}U_2^{(n)}=e^{2\pi ip_n/q_n}\,U_2^{(n)}U_1^{(n)}
\eean
Suppose that 

$$
\theta=\lim_{n\to\infty}\theta_n=\lim_{n\to\infty}p_n/q_n;\quad \th \in \R\setminus \Q,\quad p_n \in \Z, \quad q_n \in \N.
$$

The application of finite-coverings for obtaining $C\left( \T^2_\th\right)$ requires an infinite procedure  for obtaining the limit \eqref{thetalim}, i.e.
$$
\theta=\lim_{n\to\infty}\theta_n = \lim_{n\to\infty}p_n/ q_n \quad p_n,q_n \in \N.
$$
However an infinite covering can perform this operation at once. The application of infinitely many finite-fold coverings of $C\left(S^1 \right)$ yields the infinite covering $\left(C\left(S^1 \right), C_0\left( \R\right), \Z  \right)$.
If $\a \in C_0\left( \R\right)$ corresponds to the 1 shift of $\R$ and 
\bean
v \in M\left( C_0\left( \R\right)\right) = C_b\left( \R\right),
x \mapsto e^{\pi i \th x} 
\eean
then one has
$$
v \a  = e^{2\pi i \th}\a v.
$$
The above equation is a counterpart of \eqref{Tthetarel}, so the noncommutative torus is obtained.

\section{Unbounded  internal space dimension and Hausdorff non-Hausdorff duality}
 \paragraph*{}
 In the Section \ref{phys_mm_sec} it is proven that the algebra of  physical system with internal space of unbounded dimension yields inclusion of the noncommutative torus.  Here we prove that the infinite covering of physics with internal space of unbounded dimension yields the noncommutative torus itself.  Consider an one dimensional physics on $\R$ with gauge groups $U\left( N\right)$ for all $N\in \N$. For each $N \in \N$  the group $U\left( N\right)$ acts on $\C^N$. If $\C_\infty \bydef \cup_{N \in \N} \C^N$ then $\C_\infty$ is a dense subspace of $\ell_2\left(\C\right)$. If we use an isomorphism $\ell_2\left(\C\right)\to L^2\left(\R \right)$ then the generated by gauge elements group is isomorphic to
 $$
 C_0\left( \R\right) \otimes \K\left( \ell_2\left(\C\right)\right) \cong C_0\left( \R\right) \otimes \K\left(L^2\left(\R \right)\right),
 $$
 On the other hand from the Theorem  \ref{foli_bundle_thm} we have an isomorphism
 $$
 C_0\left( \R\right) \otimes \K\left(L^2\left(\R \right)\right)\cong C^*_r\left( \R^2, \widetilde{\sF}\right) 
 $$
 Otherwise
 one has
 an infinite noncommutative covering
 $$
 \left(C^*_r\left(\T^2, \sF_\th\right), C^*_r\left(\R^2, \widetilde \sF\right), \Z^2 \right)\cong \left(C\left( \T^2_\th\right)\otimes \K, C^*_r\left(\R^2, \widetilde \sF\right), \Z^2 \right)
 $$
 $C^*$-algebra $C\left( \T^2_\th\right)\otimes \K$ is Morita equivalent to $C\left( \T^2_\th\right)$. The spectrum  $C_0\left( \R\right) \otimes \K\left(L^2\left(\R \right)\right)\cong C^*_r\left( \R^2, \widetilde{\sF}\right)$  is the Hausdorff space $\R$ and the spectrum of  $C\left( \T^2_\th\right)\otimes \K$ is not Hausdorff, i.e. infinite noncommutative covering yields an example of the Hausdorff non-Hausdorff duality.
 
 \chapter{Wilson lines and noncommutative covering spaces}\label{wl_chap}
 \section{Preliminaries}
 \paragraph{} Commutative geometry has a lot of local structures, for example local sections of bundles. There are  bundles such that they have no global sections. Let $p: P \to M$ is a bundle such that $p$ has no global sections. Let $\pi : \widetilde{M} \to M$ be the universal covering projection, and $\widetilde{p}: \widetilde{M} \to M$ be the pullback \cite{spanier:at} of $p$ by $\pi$. Then $\widetilde{p}$ can have global sections. Any locally pure gauge field can be locally represented by \eqref{pure_gauge}, but cannot be represented by \eqref{pure_gauge} globally in general case. However pullback of this field can be globally repesented  by \eqref{pure_gauge}. Noncommutative geometry has no local sections. However in this book the noncommutative generalisation of covering projections are described. Local pure gauge fields can be regarded as global pure gauge fields on universal noncommutative covering projections. Locally gauge fields satisfy \eqref{pure_gauge} which can be rewritten in physical notation \cite{alekseev_bytsko:wilson_nc_tori}
 \begin{equation} \label{Ag}
 	\partial_i g = i A_i * g \,.
 \end{equation}
 This equation has a noncommutative analog.  Let $A^{\infty}$ be a smooth algebra and $\mathcal{E}^{\infty}$ be a finitely generated projective smooth . Spaces $A^{\infty}$,  $\widetilde{A}^{\infty}$, $\mathcal{E}^{\infty}$ are operator spaces. Let denote $\widetilde{\mathcal{E}}^{\infty}=\mathcal{E}^{\infty} \otimes _{A^{\infty}}\widetilde{A}^{\infty}$ where $\otimes$ means the Haagerup tensor product. The $\widetilde{\mathcal{E}}^{\infty}$ module can be regarded as pullback of $\mathcal{E}^{\infty}$. Any vector $X \in \mathrm{Lie}G$ can be lifted to the vector $\widetilde{X}\in \mathrm{Lie}\widetilde{G}$ where $\mathrm{Lie}\widetilde{G}$ is a Lie algebra of infinitesimal transformations of the $\widetilde{A}^{\infty}$. Any connection $\nabla : \mathcal{E}^{\infty} \to \mathcal{E}^{\infty} \otimes \Omega^1A^{\infty}$ can be lifted to $G$ equivariant connection $\widetilde{\nabla}: \widetilde{\mathcal{E}}^{\infty} \to \widetilde{\mathcal{E}}^{\infty} \otimes \Omega^1\widetilde{A}^{\infty}$. Otherwise any connection $\widetilde{\nabla}: \widetilde{\mathcal{E}}^{\infty} \to \widetilde{\mathcal{E}}^{\infty} \otimes \Omega^1\widetilde{A}^{\infty}$ naturally induces a map $\widetilde{\nabla}':\mathrm{End}_{\widetilde{A}^{\infty}}\left(\widetilde{\mathcal{E}}^{\infty}\right) \to \mathrm{End}_{\widetilde{A}^{\infty}}\left(\widetilde{\mathcal{E}}^{\infty}\right) \otimes \Omega^1\widetilde{A}^{\infty}$.
 Noncummutative analog of \eqref{Ag} is given by
 \begin{equation}\label{nc_Ag}
 	\widetilde{X}U = \widetilde{\nabla}_{\widetilde{X}}U; \ \forall X \in \mathrm{Lie}G, \ U\in \mathrm{Aut} \left(\widetilde{\mathcal{E}}^{\infty}\right).
 \end{equation}
 $U$ is a noncommutative analog of a global gauge transformation on the universal covering space.
 Wilson line can be regarded as a group homomorphisms $G \to \mathrm{Aut}\left(\widetilde{\mathcal{E}}^{\infty}\right)$ given by
 \begin{equation}\label{wilson_rel}
 	g \mapsto (gU) \cdot U^{-1}
 \end{equation}
 where $g \in G(\widetilde{A}^{\infty}|A^{\infty})$ and $U$ satisfies \eqref{nc_Ag}.
 
 \section{Alternative Description of  Wilson Lines}
 \paragraph{} As well as in \ref{alter_prot} (See \cite{green_schwarz_witten:superstring}) we can define an alternative approach to Wilson lines. Suppose that there is a spectral triple $(\mathcal{A}, \H, D)$ \cite{varilly:noncom}. Using the Section \ref{flat_sec} one can find a projective $\A$ bimodule $\E$ and a given by \eqref{nc_flat_conn}
 \be\nonumber
 \nabla : \mathcal E \to \mathcal E \otimes_{\A} \Om^1_D
 \ee
 (cf. Lemma \ref{flat_lem}). This construction is a noncommutative version of a geometrical canonical connection explained in the Appendix \ref{geom_flat_subsec}. The mapping between algebra and geometry is presented in the Table \ref{geo_alg_table}.
 
  	\chapter{Commutative Wilson Lines}\label{comm_prot}
 \section{Basic construction}
 \begin{empt}{\it Standard description of Wilson lines} \cite{green_schwarz_witten:superstring}. Here the physical notation is used.
 	Let $M$ be a compact Riemann manifold, and let $A$ be a locally pure gauge field, i.e. $A$ can be locally represented by following way
 	\begin{equation}\label{pure_gauge}
 		A_i = \partial_i U \cdot U^{-1}
 	\end{equation}
 	where $U$ is the gauge transformation. Whether \eqref{pure_gauge} is also true globally is another story. Condition \eqref{pure_gauge} is equivalent to that the field stength equals to zero \cite{green_schwarz_witten:superstring}.
 	If the manifold is not simply connected, $\pi_1(M) \neq \{e\}$ then there is a more general possibility that appears in electrodynamics as the Bohm-Aharonov effect. Let $\gamma$ be a noncontractible loop in $M$, beginning and ending at the same point $x$. Then "Wilson line"
 	\begin{equation}\label{gauge_transformation}
 		U_{\gamma} = P \ \mathrm{exp} \oint\limits_{\gamma} A \cdot dx
 	\end{equation}
 	is gauge invariant, and if $U_{\gamma}\neq 1$ it cannot be set to one by gauge transformation, $U_{\gamma}$ depends only on $[\gamma] \in \pi_1(M)$. If $G$ is the gauge transformation group then is a group homomorphism
 	\begin{equation}\label{from_fund_to_gauge}
 		\varphi: \pi_1(M)\to G,
 	\end{equation}
 	\begin{equation*}
 		[\gamma] \mapsto U_{\gamma}.
 	\end{equation*}
 \end{empt}
 \begin{empt}\label{alter_prot} {\it Alternative description of Wilson lines}.
 	In this construction I follow to \cite{green_schwarz_witten:superstring} 16.4.1. Let $\widetilde{M}$ be a simply connected manifold, and let $F$ be a discrete symmetry group which acts freely. Let $M = \widetilde{M}/F$. An ordinary field $\psi(x)$ is equivalent to a field on $\widetilde{M}$ that obeys
 	\begin{equation}\label{invariant_simple}
 		\psi(fx)=\psi(x); \ \forall f \in F.
 	\end{equation}
 	Then we can generalize \eqref{invariant_simple} as follows. There is a natural isomorphism $\pi_1(M) \approx F$, and from \eqref{from_fund_to_gauge} it follows that there is a natural homomorphism
 	\begin{equation}\label{from_trans_to_field}
 		F\to G,
 	\end{equation}
 	\begin{equation*}
 		f \mapsto U_f.
 	\end{equation*}
 	Now we require that $\psi$ obey not \eqref{invariant_simple} but
 	\begin{equation}\label{invariant_twist}
 		\psi(fx)=U_f\psi(x); \ \forall f \in F.
 	\end{equation}
 	This operation enables us replace guage field \eqref{gauge_transformation} by "twist" in boundary conditions \eqref{invariant_twist}. The gauge field $A$ which obeys \eqref{gauge_transformation} (resp. field $\psi$) is replaced with the trivial gauge field (resp. field $\psi'$ which obey \eqref{gauge_transformation}). Field $\psi'$ is given by
 	\begin{equation*}
 		\psi'(y) = \left(P \ \mathrm{exp} \int\limits_{\omega} A' \cdot dx\right) \psi(\pi(y))
 	\end{equation*}
 	where
 	\begin{itemize}
 		\item $\pi: \widetilde{M} \to M$ is a covering projection.
 		\item $\omega : [0,1] \to \widetilde{M}$ is such that $\omega(0) = y_0$ is a fixed point, $\omega(1)=y$.
 		\item $A'$ is a lift of $A$ by $\pi$.
 	\end{itemize}
 	This construction is similar to passive/active approach to physical transformations. Change of a gauge field is similar to a passive transformation of a coordinat system, a change of field is similar to an active transformation of point's position. Both transformations describe the same physical phenomenon.
 \end{empt}
 
 \begin{empt}{\it Replacement of closed paths by covering projections}.
 	Construction from \cite{green_schwarz_witten:superstring} can be generalized such that the universal covering is replaced by a covering  $\pi:\widetilde{M} \to M$ such that the image of the composition
 	\begin{equation*}
 		\pi_1(\widetilde{M})\xrightarrow{\pi_1(\pi)} \pi_1(M) \xrightarrow{\varphi} G
 	\end{equation*}
 	coincides with image of $\varphi$. In this case homomorphism \eqref{from_fund_to_gauge} can be replaced with
 	\begin{equation}\label{from_trans_to_gauge}
 		G(\widetilde{M}|M) \to G
 	\end{equation}
 	where $G(\widetilde{M}|M)$ is a group of covering transformations \cite{spanier:at}.
 \end{empt}

 \chapter*{Acknowledgment}

\paragraph*{}
I am very grateful to Prof. Joseph C Varilly  and Arup Kumar Pal
 for advising me on the properties of Moyal planes and equivariant spectral triples respectively. Author would like to acknowledge members of the Moscow State University Seminars
 "Noncommutative geometry and topology", "Algebras in analysis" leaded by professors A. S. Mishchenko and  A. Ya. Helemskii for a discussion
 of this work. I am grateful to my grandfather Petr who forested my mind and my grandson Petr who  watched my work during the quarantine.

 \end{document}